\documentclass[12pt, a4paper]{article}

\usepackage[latin9]{inputenc}
\usepackage{hyperref}
\usepackage{amsmath}
\usepackage{amssymb}
\usepackage{amsthm}
\usepackage{enumitem}
\usepackage{geometry}

\usepackage[arrow]{xy}
\usepackage{tikz-cd}
\usepackage{tikz-3dplot}

\tdplotsetmaincoords{50}{50}

\theoremstyle{plain}
\newtheorem{prop}{Proposition}[subsection]
\newtheorem{conj}[prop]{Conjecture}
\newtheorem{lem}[prop]{Lemma}
\newtheorem{thm}[prop]{Theorem}
\newtheorem{cor}[prop]{Corollary}

\theoremstyle{definition}
\newtheorem{defn}[prop]{Definition}
\newtheorem{ex}[prop]{Example}
\newtheorem{rem}[prop]{Remark}

\newtheorem{claim}[prop]{Claim}

\theoremstyle{plain}
\newtheorem{lem0}{Lemma}[section]
\theoremstyle{definition}
\newtheorem{defn0}[lem0]{Definition}
\newtheorem{ex0}[lem0]{Example}
\newtheorem{rem0}[lem0]{Remark}

\theoremstyle{plain}

\theoremstyle{definition}

\newcommand{\Z}{\mathbb{Z}}
\newcommand{\Zp}{\mathbb{Z}_p}
\newcommand{\bQ}{\mathbb{Q}}
\newcommand{\Qp}{\mathbb{Q}_p}
\newcommand{\bR}{\mathbb{R}}
\newcommand{\bH}{\mathbb{H}}
\newcommand{\bP}{\mathbb{P}}

\newcommand{\GL}{\mathrm{GL}}
\newcommand{\U}{\mathrm{U}}

\newcommand{\fb}{\mathfrak{b}}

\newcommand{\fg}{\mathfrak{g}}

\newcommand{\ft}{\mathfrak{t}}
\newcommand{\fl}{\mathfrak{l}}
\newcommand{\fp}{\mathfrak{p}}
\newcommand{\fu}{\mathfrak{u}}
\newcommand{\fn}{\mathfrak{n}}
\newcommand{\fz}{\mathfrak{z}}

\newcommand{\cB}{\mathcal{B}}
\newcommand{\cC}{\mathcal{C}}
\newcommand{\cD}{\mathcal{D}}
\newcommand{\cF}{\mathcal{F}}
\newcommand{\cH}{\mathcal{H}}
\newcommand{\cI}{\mathcal{I}}
\newcommand{\cJ}{\mathcal{J}}
\newcommand{\cL}{\mathcal{L}}
\newcommand{\cM}{\mathcal{M}}

\newcommand{\cO}{\mathcal{O}}
\newcommand{\cP}{\mathcal{P}}
\newcommand{\cT}{\mathcal{T}}

\newcommand{\cX}{\mathcal{X}}

\newcommand{\al}{\alpha}
\newcommand{\val}{\mathrm{val}}
\newcommand{\Id}{\mathrm{id}}
\newcommand{\Hom}{\mathrm{Hom}}
\newcommand{\Ext}{\mathrm{Ext}}
\newcommand{\soc}{\mathrm{soc}}
\newcommand{\Ker}{\mathrm{ker}}
\newcommand{\Dim}{\mathrm{dim}}
\newcommand{\dom}{\mathrm{dom}}

\newcommand{\un}[1]{\underline{#1}}
\newcommand{\tld}[1]{\widetilde{#1}}
\newcommand{\smat}[1]{\left( \begin{smallmatrix} #1 \end{smallmatrix} \right)}
\newcommand{\defeq}{\stackrel{\textrm{\tiny{def}}}{=}}

\title{Splitting and making explicit the de Rham complex of the Drinfeld space}

\author{Christophe Breuil\footnote{CNRS, B\^atiment 307, Facult\'e d'Orsay, Universit\'e Paris-Saclay, 91405 Orsay Cedex, France}\\
\and
Zicheng Qian\footnote{Morningside Center of Mathematics, Academy of Mathematics and Systems Science, Chinese Academy of Sciences, Beijing 100190, China; University of the Chinese Academy of Sciences, Beijing 100049, China}}

\date{ }

\begin{document}

\maketitle

\begin{abstract}
Let $p$ be a prime number, $K$ a finite extension of $\Qp$ and $n$ an integer $\geq 2$. We completely and explicitly describe the global sections $\Omega^\bullet$ of the de Rham complex of the Drinfeld space over $K$ in dimension $n-1$ as a complex of (duals of) locally $K$-analytic representations of $\GL_n(K)$. Using this description, we construct an explicit section in the derived category of (duals of) finite length admissible locally $K$-analytic representations of $\GL_n(K)$ to the canonical morphism of complexes $\Omega^\bullet \twoheadrightarrow H^{n-1}(\Omega^\bullet)[-(n-1)]$.
\end{abstract}

\tableofcontents

\newpage

\section{Introduction}\label{sec: intro}

Let $p$ be a prime number and $K$ a finite field extension of the field of $p$-adic numbers $\Qp$. This monograph has to do with certain $p$-adic representations of ${\rm GL}_n(K)$. One of its main aims is to completely describe the complex of differential forms of the Drinfeld space of dimension $n-1$ as a complex of representations of ${\rm GL}_n(K)$, where $n\geq 2$. The Drinfeld spaces are very important $p$-adic spaces over $K$ introduced by Drinfeld in the seventies, and there is one Drinfeld space (of dimension $n-1$) for each integer $n\geq 2$. The group ${\rm GL}_n(K)$ acts on the Drinfeld space of dimension $n-1$, hence on its complex of differential forms, yielding representations of ${\rm GL}_n(K)$ that mathematicians have started studying in the eighties. Understanding these representations was one of the main motivations for the development of the theory of locally analytic representations of ${\rm GL}_n(K)$ in the 2000's, which can be seen as a $p$-adic analogue of Harish-Chandra's $(\mathfrak{gl}_n,K)$-modules where $K$ here is not the above field $K$ but a maximal compact subgroup of ${\rm GL}_n(\bR)$. One difference is that, whereas ${\rm GL}_n(\bR)$ does not act on $(\mathfrak{gl}_n,K)$-modules (though it acts on analytic globalizations, see e.g.~\cite{GKS11}, \cite{S24}), the whole group ${\rm GL}_n(K)$ acts on locally analytic representations. We review below the history of that subject (which is intimately related to the Drinfeld space), and then we explain the main results of this book.

\subsection{Historical background}\label{history}

We fix $n\geq 2$ an integer and we let ${\bP^{n-1}_{\rm{rig}}}_{\!\!\!\!/K}$ be the rigid analytic projective space of dimension $n-1$ over $K$. Fifty years ago, Drinfeld introduced in \cite{Dr74} what is now called the Drinfeld space $\bH_{/K}$ over $K$ of dimension $n-1$, which is the rigid analytic (admissible) open subspace of $\bP^{n-1}_{\rm{rig}}$ defined as
\begin{equation}\label{defDrinfeld}
\bH\defeq \bP^{n-1}_{\rm{rig}}\setminus\bigcup_{\cH}\cH
\end{equation}
where $\cH$ runs through the $K$-rational hyperplanes inside $\bP^{n-1}_{\rm{rig}}$. It has become such a familiar space in $p$-adic arithmetic geometry that it is hard to imagine a time when it was not defined. In \cite{Dr76} Drinfeld gave a second definition of $\bH$ as a moduli space of certain $p$-divisible groups, and used this definition to define a tower of \'etale coverings of $\bH$. These coverings have now been vastly generalized into the Rapoport-Zink spaces (\cite{RZ96}), and more recently the local Shimura varieties (\cite{RV14}, \cite{SW20}). Their cohomology plays a fundamental role in the Langlands program.\bigskip

In this work, we only use the definition of $\bH$ given in (\ref{defDrinfeld}). The group $G\defeq {\rm GL}_n(K)$ naturally acts on $\bP^{n-1}_{\rm{rig}}$ and preserves the open subspace $\bH$. By functoriality, it follows that any cohomology group of $\bH$ is naturally endowed with a left action of $G$. More than thirty years ago, Schneider and Stuhler in \cite{SS91} computed this action on any abstract cohomology theory satisfying certain axioms (see \cite[\S2]{SS91} for more details, we won't need these axioms). The first example of such a cohomology theory is the $\ell$-adic \'etale cohomology
\[H^\bullet_{\textrm{\'et}}\big(\bH\times_K \widehat{\overline K}, {\mathbb Q}_\ell\big)\]
where $\widehat{\overline K}$ is the $p$-adic completion of an algebraic closure $\overline K$ of $K$ and $\ell$ is a prime number distinct from $p$\footnote{The case $\ell=p$ was treated more recently in \cite{CDN20}.}. The second example is the de Rham cohomology
\[H^\bullet_{\rm dR}(\bH).\]
We recall their result. For $P\subseteq G$ the $K$-points of a lower standard parabolic subgroup of ${\rm GL}_n$, we write $\mathrm{Ind}_P^G (1)$ for the smooth parabolic induction of the trivial representation of $P$. For $k\in \{0,\dots, n-1\}$ we let $P_{[1,n-k-1]}$ be the $K$-points of the lower standard parabolic subgroup of $G$ of Levi $\mathrm{GL}_{n-k}\times \mathrm{GL}_1\times \cdots \times \mathrm{GL}_1$. Following our notation in the text (see (\ref{gensteinberg}) and (\ref{extranot})), we define:
\begin{equation}\label{steinbergintro}
V_{[1,n-k-1],\Delta}^{\infty}\defeq (\mathrm{Ind}_{P_{[1,n-k-1]}}^G 1)^\infty\ \ /\!\sum_{P_{[1,n-k-1]}\subsetneq P}(\mathrm{Ind}_P^G 1)^\infty
\end{equation}
(where $(\mathrm{Ind}_{P}^G (-))^\infty$ is the usual smooth parabolic induction and $\Delta$ the set of simple roots of ${\rm GL}_n$) which is called a smooth generalized Steinberg representation of $G$. The representations $V_{[1,n-k-1],\Delta}^{\infty}$ are absolutely irreducible, and note that $V_{[1,n-1],\Delta}^{\infty}=1$ (the trivial representation of $G$) while $V_{\emptyset,\Delta}^{\infty}$ is the smooth Steinberg representation of $G$, that we also denote by $\mathrm{St}_n^\infty$.

\begin{thm}[\cite{SS91}]\label{SSintro}
We have $H^k_{\mathrm{\acute et}}\big(\bH\times_K \widehat{\overline K}, {\mathbb Q}_\ell\big)=H^k_{\rm dR}(\bH)=0$ for $k\geq n$, and for $k\in \{0,\dots, n-1\}$ we have $G$-equivariant isomorphisms
\[H^k_{\mathrm{\acute et}}\big(\bH\times_K \widehat{\overline K}, {\mathbb Q}_\ell\big)\cong (V_{[1,n-k-1],\Delta}^{\infty})^\vee\ \ \textrm{and}\ \ \ H^k_{\rm dR}(\bH)\cong (V_{[1,n-k-1],\Delta}^{\infty})^\vee\]
where the first $V_{[1,n-k-1],\Delta}^{\infty}$ is seen with ${\mathbb Q}_\ell$-coefficients, the second with $K$-coefficients, and where $(-)^\vee$ is the corresponding algebraic dual.
\end{thm}

More precisely the above theorem follows from \cite[\S3 Thm.~1]{SS91} (with $A={\mathbb Q}_\ell$ or $A=K$) together with \cite[\S4 Lemma~1]{SS91}. In particular we see that the $\ell$-adic and de Rham cohomology groups of $\bH$ are duals of \emph{smooth} representations of $G$, which is not obvious \emph{a priori}.\bigskip

By \cite[\S1 Prop.~4]{SS91} the rigid analytic space $\bH$ is \emph{quasi-Stein}, by which we mean that $\bH$ admits an admissible covering given by the union of an ascending sequence $U_1\subseteq U_2\subseteq \cdots \subseteq U_n\subseteq \cdots $ of affinoid open subspaces such the restriction maps of Banach spaces $\Gamma(U_n,\mathcal{O}_{U_n})\rightarrow \Gamma(U_{n-1},\mathcal{O}_{U_{n-1}})$ have dense image. It follows that $H^k(\bH, \mathcal F)=0$ for $k\geq 1$ and $\mathcal F$ any coherent sheaf on $\bH$ (\cite[Satz 2.4]{Kie67}). In particular the de Rham cohomology of $\bH$ (which for an arbitrary rigid space is defined as the hypercohomology of its de Rham complex) is just here the cohomology of the complex $\Omega^\bullet$ of its \emph{global sections}, i.e.
\[H^k_{\rm dR}(\bH) = H^k(\Omega^\bullet)\]
where
\[\Omega^\bullet\defeq [\Omega^0\longrightarrow \Omega^1\longrightarrow \cdots \longrightarrow \Omega^{n-1}]\]
with
\begin{equation}\label{Omegai}
\Omega^k\defeq \Gamma(\bH, \Omega^k_{\bH/K}).
\end{equation}
By functoriality each $\Omega^k$ is also a representation of $G$ over $K$ and the complex $\Omega^\bullet$ is $G$-equivariant. Moreover, as we have seen in Theorem \ref{SSintro}, $H^k(\Omega^\bullet)$ is the (algebraic) dual of a smooth representation of $G$. However, the individual representations $\Omega^i$ are far from being (duals of) smooth representations of $G$.\bigskip

Investigating which type of $G$-action they carry was one of the main motivations for Schneider and Teitelbaum's theory of (admissible) locally analytic representations of $G$, and completely determining the representations $\Omega^k$ of $G$ and the differential maps in $\Omega^\bullet$ is one of the main goals of this monograph. Indeed, motivated by the beginnings of the $p$-adic Langlands program, the first author has been fascinated during years (in the early 2000's) by the idea of fully understanding the internal structure of $\Omega^\bullet$ as a complex of representations of $G$, with the hope that it was maybe hiding ``secrets''. In the end, this internal structure is not that complicated and does not hide that many secrets, but completely unravelling $\Omega^\bullet$ remains an interesting and challenging (though sometimes technical) task.\bigskip

Writing $\bH = \cup_n U_n$ for affinoids $U_n\subseteq U_{n+1}$ as above, one has $\Omega^k\buildrel\sim\over\rightarrow \varprojlim_{n}\Gamma(U_n, \Omega^k_{\bH/K})$. Since all $\Gamma(U_n, \Omega^k_{\bH/K})$ are ($p$-adic) Banach spaces, the projective limit topology gives a natural ($p$-adic) Fr\'echet topology on each $\Omega^k$ (recall that a $p$-adic Fr\'echet space is a projective limit of countably many $p$-adic Banach spaces), and it is not very hard to check that the map $G \times \Omega^k\rightarrow \Omega^k$, $(g, v)\mapsto gv$ is continuous. In the eighties, Morita pioneered the study of the continuous $G$-representations $\Omega^0$ and $\Omega^1$ when $n=2$ (\cite{Mor82}, \cite{Mor84}, \cite{Mor85}). In particular in \cite[\S 5]{Mor84} he proved the following theorem:

\begin{thm}[\cite{Mor84}]\label{morita}
Let $B\subseteq G={\rm GL}_2(K)$ be the $K$-points of the lower Borel of ${\rm GL}_2$, then one has $G$-equivariant topological isomorphisms
\[\Omega^0/1\cong \big((\mathrm{Ind}_B^G \ t^{-1}\!\boxtimes t)^{\rm{an}}\big)^\vee\ \ \textrm{and}\ \ \ \Omega^1\cong \big((\mathrm{Ind}_B^G \ 1)^{\rm{an}}/1\big)^\vee.\]
Here $((\mathrm{Ind}_B^G \chi_1\boxtimes \chi_2)^{\rm{an}})^\vee$ is the topological dual of the locally $K$-analytic principal series $(\mathrm{Ind}_B^G \chi_1\boxtimes \chi_2)^{\rm{an}}$ where $\chi_i:K^\times \rightarrow K$ are locally $K$-analytic characters and $t$ is the identity character $K^\times \rightarrow K$, $t\mapsto t$.
\end{thm}

Recall that $(\mathrm{Ind}_B^G \chi_1\boxtimes \chi_2)^{\rm{an}}$ is the $K$-vector space of locally $K$-analytic fonctions $f:G\longrightarrow K$ such that
\[f(gb)=b^{-1}\cdot f(g)=(\chi_1(t_1)\chi_2(t_2))^{-1}f(g),\ \ \ g\in G,\ b=\smat{\!t_1 & 0\\ * & t_2\!}\in B\]
with the (left) action of $G$ given by
\begin{equation}\label{actionG}
(g(f))(g')\defeq f(g^{-1}g'),\ \ g,g'\in G.
\end{equation}
It is naturally endowed with a topology which makes it a locally convex $K$-vector space of compact type for which the action of $G$ is continuous (see \cite[\S 6]{ST05} and \cite[Lemma 2.1]{ST02a}). The representation $(\mathrm{Ind}_B^G 1)^{\rm{an}}/1$ is called the locally $K$-analytic Steinberg representation of $G$. Note that Morita more generally proves Theorem \ref{morita} for (what is now called) Schneider's holomorphic discrete series (\cite[\S 3]{S92}). In \cite[Thm.~1]{Mor85} and \cite[Thm.~2]{Mor84} (see also \cite[Thm.~1]{Mor84} and \cite[Thm.~2]{Mor84}), Morita moreover proves that $(\mathrm{Ind}_B^G t^{-1}\boxtimes t)^{\rm{an}}$ is (topologically) irreducible as a $G$-representation and that $(\mathrm{Ind}_B^G 1)^{\rm{an}}/1$ is an extension of $(\mathrm{Ind}_B^G \ t^{-1}\boxtimes t)^{\rm{an}}$ by the (irreducible) smooth Steinberg representation of $G$.\footnote{Morita's results are proven for ${\rm SL}_2(K)$ instead of ${\rm GL}_2(K)$. Moreover the proof of the irreducibility is actually flawed, see \cite[p.443-444]{ST02a}.}\bigskip

For many years, there was no progress on and (maybe) no interest in these $p$-adic questions (except \cite{S92}) until Schneider and Teitelbaum, using F\'eaux de Lacroix's thesis \cite{Lac99}, decided to start the theory from scratch in \cite{ST02a}, \cite{ST02b}, \cite{ST01}, \cite{ST03} and \cite{ST05}. In particular in \cite{ST03} they defined an important abelian category of admissible locally $K$-analytic representations of $G$ (more generally of a locally $K$-analytic group) on locally convex $K$-vector spaces of compact type, or equivalently taking continuous duals an abelian category of (so called) coadmissible $D(G)$-modules. Here $D(G)$ is the $K$-algebra of locally $K$-analytic distributions on $G$, i.e.~the continuous dual of 
\[C^{\rm{an}}(G)\defeq \{f: G \longrightarrow K,\ \ f\textrm{ locally $K$-analytic}\}\]
endowed with its natural locally convex topology (\cite[\S 2]{ST02a}). Shortly after, Emerton gave his own account of the theory in \cite{Em17} (which was published years later). All this work ultimately lead, a decade or so later, to the beautiful theory of Orlik-Strauch representations in \cite{OS15}, which can be seen as one of the key achievements of the theory.\bigskip

Orlik-Strauch representations $\cF_{P}^G(M,\pi^{\infty})$ are now widely used and are, essentially, the only admissible locally $K$-analytic representations of $G$ which are so far well understood. We briefly review their main properties in the text (see Theorem \ref{prop: OS property}). Suffice it here to say that, if $P$ is the $K$-points of a (lower standard) parabolic subgroup of ${\rm GL}_n$, $M$ a $U(\fg)$-module in Bernstein-Gelfand-Gelfand category's ${\mathcal O}^{\fp}$ (\cite[\S 9.3]{Hum08}) such that all its weights are integral (here $\fg$, $\fp$ are the respective $K$-Lie algebras of $G$, $P$ and $U(-)$ the enveloping algebra) and $\pi^\infty$ a smooth admissible representation of the Levi factor $L_P$ of $P$, then 
\begin{equation}\label{OSintro}
\cF_{P}^G(M,\pi^{\infty}) \defeq \Big(\big(\mathrm{Ind}_P^G \ W^\vee\!\otimes \pi^\infty\big)^{\rm{an}}\Big)^{\Ker =0}
\end{equation}
where $W$ is any finite dimensional algebraic representation of $P$ such that one has a surjection of $U(\fg)$-modules $U(\fg) \otimes_{\fp}W \twoheadrightarrow M$ and where $(-)^{\Ker =0}$ is the (closed) subrepresentation of the locally $K$-analytic parabolic induction $(\mathrm{Ind}_P^G \ W^\vee \!\otimes \pi^\infty)^{\rm{an}}$ of vectors killed by
\[\Ker \defeq \Ker\big(U(\fg) \otimes_{\fp}W \twoheadrightarrow M\big).\]
Here $(\mathrm{Ind}_P^G \ W^\vee \!\otimes \pi^\infty)^{\rm{an}}$ is defined similarly to $(\mathrm{Ind}_B^G \chi_1\boxtimes \chi_2)^{\rm{an}}$ above with left action (\ref{actionG}) of $G$. We refer the reader to \cite{OS15} for more details on the admissible representations $\cF_{P}^G(M,\pi^{\infty})$. In particular, when $\pi^\infty$ is of finite length, they are moreover (topologically) of finite length.\bigskip

Going back to the $G$-representations $\Omega^i$ of (\ref{Omegai}), the first result after Morita's Theorem \ref{morita} came in \cite{ST02b} where the authors could describe explicitly the graded pieces of a filtration on $\Omega^{n-1}$ (for any $n\geq 2$). More precisely in \cite[Thm.~8.6]{ST02b} they proved:

\begin{thm}[\cite{ST02b}]\label{boundary}
The $D(G)$-module $\Omega^{n-1}$ admits a filtration by closed $D(G)$-submodules
\[\Omega^{n-1}=\mathrm{Fil}^0(\Omega^{n-1})\supsetneq \mathrm{Fil}^1(\Omega^{n-1})\supsetneq \cdots \supsetneq\mathrm{Fil}^{n}(\Omega^{n-1})=0\]
such that for $0\leq j\leq n-1$:
\begin{equation}\label{filtrationintro}
\mathrm{gr}^{j}(\Omega^{n-1})\defeq \mathrm{Fil}^j(\Omega^{n-1})/\mathrm{Fil}^{j+1}(\Omega^{n-1})\cong \cF_{P_{\widehat{j}}}^{G}\big(L(s_j s_{j-1}\cdots s_1\cdot 0),1_{{\rm GL}_j}\boxtimes \mathrm{St}^\infty_{n-j}\big)^{\!\vee}\footnote{The reference \cite{OS15} was not available at the time of \cite{ST02b} but they directly used the description on the right hand side of (\ref{OSintro}).}
\end{equation}
where $P_{\widehat{j}}$ is the $K$-points of the lower parabolic subgroup of ${\rm GL}_n$ of Levi factor ${\rm GL}_j\times {\rm GL}_{n-j}$, $s_1,\dots,s_{n-1}$ are the simple reflections of ${\rm GL}_n$ and $L(s_j s_{j-1}\cdots s_1\cdot 0)$ is the unique irreducible $U(\fg)$-module in ${\mathcal O}^{\fp_{\widehat{j}}}$ of highest weight $s_j s_{j-1}\cdots s_1\cdot 0$ (with $P_{\widehat j}= G$, $s_j s_{j-1}\cdots s_1= 1$ when $j=0$, and the dot action $\cdot$ being relative to the lower Borel subgroup of ${\rm GL}_n$).
\end{thm}

Since all graded pieces in Theorem \ref{boundary} are finite length coadmissible $D(G)$-modules, and since coadmissibility and finite length are preserved under extensions (for coadmissibility see for instance the proof of \cite[Lemme 2.1.1]{Bre19}), we see in particular that $\Omega^{n-1}$ is also a finite length coadmissible $D(G)$-module. To prove some of the results of this monograph (see \S\ref{introresults} below), we use a weak variant of Theorem \ref{boundary} as a \emph{key} ingredient (see Theorem \ref{thm: dR complex filtration} in the text together with the comment before Theorem \ref{thm: main dR}).\bigskip

A few years later, Pohlkamp in \cite{Po04} proved a result analogous to Theorem \ref{boundary} but where $\Omega^{n-1}$ is replaced by the global sections $\Omega^0$ of the structural sheaf of $\bH$. Finally, a few more years later Orlik considerably generalized both statements to $\Omega^k$ for all $k$ in \cite{Or08} (see also \cite{Or15} or Theorem \ref{thm: dR complex filtration} in the text). In particular all $D(G)$-modules $\Omega^k$ are coadmissible and (topologically) of finite length. Note that, when $0<k<n-1$, $\Omega^k$ is more envolved because it has (essentially) twice as many irreducible constituents as $\Omega^{n-1}$ or $\Omega^0$.\bigskip

Although the above theorems describe the graded pieces of a filtration on $\Omega^k$, and in particular the irreducible constituents of $\Omega^k$, this does not give its full internal structure. For instance we do not know the extensions between the graded pieces (some extensions as subquotients could be split). Or it could be that there are several analogous finite length coadmissible $D(G)$-modules with the same graded pieces but with different extensions as subquotients. As an example, let us go back to the case $n=2$ and Theorem \ref{morita}, and denote ${\rm SP}\defeq (\mathrm{Ind}_B^G \ t^{-1}\!\boxtimes t)^{\rm{an}}$. Although $\Omega^1$ is fully determined there, this is not the case of $\Omega^0$: we only know it is an extension of ${\rm SP}$ by the trivial representation $1$. To obtain the full structure of $\Omega^0$, we need to know that it is the unique such non-split extension\footnote{Curiously the authors could not find this classical result explicitly and clearly stated in the literature (it follows for instance from the baby case $n=2$, $k=0$, $\mu_0=(0,0)$ of Theorem \ref{thm: main dR} in the text).}. In particular for $n=2$ we can make the complex $\Omega^\bullet = [\Omega^0\longrightarrow \Omega^1]$ transparent by rewriting it
\begin{equation}\label{explicitcomplex}
\Omega^\bullet = \left[(\!\begin{xy} (0,0)*+{1}="a"; (13,0)*+{{\rm SP}^\vee}="b"; {\ar@{-}"a";"b"}\end{xy}\!) \longrightarrow (\!\begin{xy} (0,0)*+{{\rm SP}^\vee}="a"; (18.5,0)*+{(\mathrm{St}_2^\infty)^\vee}="b"; {\ar@{-}"a";"b"}\end{xy}\!) \right]
\end{equation}
where as usual a line means a non-split extension between two irreducible constituents (with the socle on the left and the cosocle on the right) and where we use that $(\mathrm{Ind}_B^G \ 1)^{\rm{an}}/1$ is the (unique) non-split extension of ${\rm SP}$ by the smooth Steinberg $\mathrm{St}_2^\infty$. An analogous complete description of $\Omega^\bullet$ for $G=\GL_3(\Qp)$ was given by Schraen in his thesis (\cite[\S 6.4]{Schr11}).\bigskip

The history of $\Omega^\bullet$ (so far) did not quite stop there. By another result of Schraen (\cite[Thm.~6.1]{Schr11}\footnote{Although the result is stated there only for $K=\Qp$, its proof works for arbitrary $K$ using Lemma \ref{erratum}.}) crucially based on results of Orlik (\cite[Thm.~1]{Or05}) and on a theorem of Dat (\cite[Cor.~A.1.3]{Dat06}) which itself is an elaboration of Deligne's splitting result (\cite{De68}), the complex $\Omega^\bullet$ splits in the bounded derived category of all (abstract) $D(G)$-modules, i.e.~there exists an isomorphism in this derived category:
\begin{equation}\label{qiintro}
\Omega^\bullet \cong \bigoplus_{k=0}^{n-1} (V_{[1,n-k-1],\Delta}^{\infty})^\vee[-k].
\end{equation}
However, trying to unravel the abstract proofs of \cite[Cor.~A.1.3]{Dat06} and \cite[Thm.~6.1]{Schr11} to produce an explicit such isomorphism, say in the bounded derived category of finite length coadmissible $D(G)$-modules (instead of all $D(G)$-modules), seems seriously challenging. Yet, when $n=2$, one can easily produce an explicit such isomorphism as follows. Choosing a $p$-adic logarithm $\log : K^\times \rightarrow K$, one can ``glue'' $\Omega^0$ and $\Omega^1$ in (\ref{explicitcomplex}) into one length $3$ $D(G)$-module $\!\begin{xy} (0,0)*+{1}="a"; (11,0)*+{{\rm SP}^\vee}="b"; (27.2,0)*+{(\mathrm{St}_2^\infty)^\vee}="c"; {\ar@{-}"a";"b"}; {\ar@{-}"b";"c"}\end{xy}\!$ (see for instance \cite[\S 3.2]{Bre19}). One then has an explicit section in the derived category $(\mathrm{St}_2^\infty)^\vee[-1] \dashrightarrow \Omega^\bullet$ to the canonical morphism of complexes $\Omega^\bullet \rightarrow H^1(\Omega^\bullet)[-1]\cong (\mathrm{St}_2^\infty)^\vee[-1]$ provided by (see \cite[\S 5.1]{Schr10}):
\begin{equation}\label{splittingcomplex}
(\mathrm{St}_2^\infty)^\vee[-1]\ \longleftarrow \ \left[(\!\begin{xy} (0,0)*+{1}="a"; (11,0)*+{{\rm SP}^\vee}="b"; {\ar@{-}"a";"b"}\end{xy}\!) \rightarrow (\!\begin{xy} (0,0)*+{1}="a"; (11,0)*+{{\rm SP}^\vee}="b"; (27.3,0)*+{(\mathrm{St}_2^\infty)^\vee}="c"; {\ar@{-}"a";"b"}; {\ar@{-}"b";"c"}\end{xy}\!)\right] \ \longrightarrow \ \Omega^\bullet
\end{equation}
(where the morphisms of complexes are easily guessed and are quasi-isomorphisms). Adding up the (trivial) morphism of complexes $1[0]\rightarrow \Omega^\bullet$, we deduce an isomorphism as in (\ref{qiintro}) $1[0]\oplus (\mathrm{St}_2^\infty)^\vee[-1] \buildrel\sim\over \dashrightarrow \Omega^\bullet$. In fact, there exists a slightly better variant which, in the case $K=\Qp$, is more directly related to the $p$-adic local Langlands correspondence for $\GL_2(\Qp)$. Let $\varepsilon:K^\times \rightarrow K$ be the $p$-adic cyclotomic character (which factors as $K^\times \rightarrow \Qp^\times \twoheadrightarrow \Zp^\times\subset K$ where the first map is the norm and the second the projection sending $p$ to $1$) and let ${\rm SP}'\defeq (\mathrm{Ind}_B^G \varepsilon^{-1}\!\boxtimes \varepsilon)^{\rm{an}}$. Then one can add in a unique way the constituent ${\rm SP}'^\vee$ to each length $3$ $D(G)$-module as above and obtain a length $4$ $D(G)$-module
\[\begin{xy} (0,0)*+{1}="a"; (10,0)*+{{\rm SP}^\vee}="b"; (25,0)*+{(\mathrm{St}_2^\infty)^\vee}="c"; (-10.5,0)*+{{\rm SP}'^\vee}="d"; {\ar@{-}"a";"b"}; {\ar@{-}"b";"c"}; {\ar@{-}"d";"a"}\end{xy}\!.\]
The isomorphism classes of such $D(G)$-modules are in non-canonical bijection with ${\mathbb A}^{1}(K)$ (which corresponds to the choice of $\log(p)\in K$ in the $p$-adic logarithm). When $K=\Qp$, such length $4$ representations precisely correspond to $2$-dimensional semi-stable non-crystalline representations of ${\rm Gal}(\overline\Qp/\Qp)$ (over $\Qp$) of Hodge-Tate weights $(0,1)$ (see \cite{CDP14} or \cite[\S 3.1]{Bre19}, the value of $\log(p)\in K$ is related to Fontaine's $\cL$-invariant on the Galois side). One then has a section (in the derived category) similar to (\ref{splittingcomplex})
\begin{equation}\label{splittingcomplexbis}
(\mathrm{St}_2^\infty)^\vee[-1] \longleftarrow \left[(\!\begin{xy} (0,0)*+{1}="a"; (10,0)*+{{\rm SP}^\vee}="b"; (-10.5,0)*+{{\rm SP}'^\vee}="d"; {\ar@{-}"a";"b"}; {\ar@{-}"d";"a"}\end{xy}\!) \rightarrow (\!\begin{xy} (0,0)*+{1}="a"; (10,0)*+{{\rm SP}^\vee}="b"; (25,0)*+{(\mathrm{St}_2^\infty)^\vee}="c"; (-10.5,0)*+{{\rm SP}'^\vee}="d"; {\ar@{-}"a";"b"}; {\ar@{-}"b";"c"}; {\ar@{-}"d";"a"}\end{xy}\!)\right]
 \longrightarrow \Omega^\bullet.
\end{equation}
Following the notation of this text, the constituent ${\rm SP}^\vee$ is denoted $X_1$ and the constituent ${\rm SP}'^\vee$ is denoted $Y_1$.\bigskip

With an explicit description of $\Omega^\bullet$ at hand for $n\geq 2$ (as we will soon have), it becomes tempting to look for a generalization of (\ref{splittingcomplexbis}) to $n\geq 3$. We provide such explicit sections in this work, which also apply to Schneider's holomorphic discrete series.

\subsection{The main results}\label{introresults}

The first aim of this monograph is to finish the work started by Schneider-Teitelbaum and continued by Orlik (and Schraen) by giving a transparent description of $\Omega^\bullet$ for $n\geq 2$ analogous to (\ref{explicitcomplex}) when $n=2$.\bigskip

We first need a bit of notation. For $\un{j}=(j_0,j_1,j_2)$ such that
\begin{equation}\label{Jintro}
1\leq j_0, j_1\leq n-1,~1\leq j_2\leq n\ \ {\rm and}\ \ 0\leq j_2-j_1\leq n-1
\end{equation}
we set using the notation in (\ref{OSintro}) (see (\ref{cj}))
\begin{equation}\label{cjintro}
C_{\un{j}}=C_{(j_0,j_1,j_2)}\defeq \cF_{P_{\widehat{j}_1}}^{G}(L(w_{j_1,j_0}\cdot 0),\pi_{j_1,j_2}^{\infty})
\end{equation}
where $P_{\widehat{j}_1}$ is as in Theorem \ref{boundary}, $w_{j_1,j_0}\defeq \left\{\begin{array}{ccc} s_{j_1}s_{j_1-1}\cdots s_{j_0}&{\rm if}&j_1\geq j_0\\
s_{j_1}s_{j_1+1}\cdots s_{j_0}&{\rm if}&j_1\leq j_0
\end{array}\right.$, $L(w_{j_1,j_0}\cdot 0)$ is as in Theorem \ref{boundary} and where $\pi_{j_1,j_2}^{\infty}$ is an explicit irreducible smooth representation of ${\rm GL}_{j_1}(K)\times {\rm GL}_{n-j_1}(K)$ defined in (\ref{pij_1,_2}). The $C_{\un{j}}$ are irreducible admissible locally analytic representations of $G$ over $K$. For $k\in \{1,\dots,n-1\}$ one first proves that there exists a unique finite length coadmissible $D(G)$-module $X_k$ of the following form (see (\ref{definitionofxyz}) with Theorem \ref{thm: unique rep} and recall $(-)^\vee$ means the continuous dual):
\begin{equation}\label{defineXk}
X_k = \!\begin{xy} (0,0)*+{C_{(n-k,n-1,n-1)}^\vee}="a"; (32,0)*+{C_{(n-k,n-2,n-2)}^\vee}="b"; (54.5,0)*+{\cdots}="c"; (73,0)*+{C_{(n-k,2,2)}^\vee}="d"; (96.5,0)*+{C_{(n-k,1,1)}^\vee}="e"; {\ar@{-}"a";"b"}; {\ar@{-}"b";"c"}; {\ar@{-}"c";"d"}; {\ar@{-}"d";"e"}\end{xy}
\end{equation}
(we say that $X_k$ is uniserial).

\begin{thm}[Theorem \ref{thm: main dR}, \ref{it: unicity rep 1} of Corollary \ref{unicitystuff}]\label{jolicomplexe}
For $k\in \{0,\dots,n-1\}$ $\Omega^k$ is the unique coadmissible $D(G)$-module of the form 
\[\begin{xy}(0,0)*+{X_{k}}="a"; (25.2,0)*+{(V_{[1,n-k-1],\Delta}^\infty)^\vee}="b"; (52,0)*+{X_{k+1}}="c";
{\ar@{-}"a";"b"}; {\ar@{-}"b";"c"}\end{xy}\!\]
with $X_0=X_n\defeq 0$. It is indecomposable multiplicity free with an irreducible socle and cosocle. Moreover the $k$-th differential map is the unique (up to non-zero scalar) non-zero map of $D(G)$-modules $\Omega^k \rightarrow \Omega^{k+1}$.
\end{thm}

The first statement of Theorem \ref{jolicomplexe} follows from Theorem \ref{thm: main dR} and \ref{it: unicity rep 1} of Corollary \ref{unicitystuff}, while the second statement follows from (\ref{definitionofd0}) and (\ref{definitionofdk}). Note that the differential map in Theorem \ref{jolicomplexe} obviously factors as
\[\begin{xy}(0,0)*+{X_{k}}="a"; (22,0)*+{(V_{[1,n-k-1],\Delta}^\infty)^\vee}="b"; (46,0)*+{X_{k+1}}="c";
{\ar@{-}"a";"b"}; {\ar@{-}"b";"c"}\end{xy}\!\twoheadrightarrow X_{k+1}\hookrightarrow \!\begin{xy}(0,0)*+{X_{k+1}}="a"; (23.7,0)*+{(V_{[1,n-k-2],\Delta}^\infty)^\vee}="b"; (47.4,0)*+{X_{k+2}}="c";{\ar@{-}"a";"b"}; {\ar@{-}"b";"c"}\end{xy}\!.\]
In fact one exactly knows the form of $\Omega^k\rightarrow \Omega^{k+1}$, which looks like (for $2\leq k\leq n-4$):
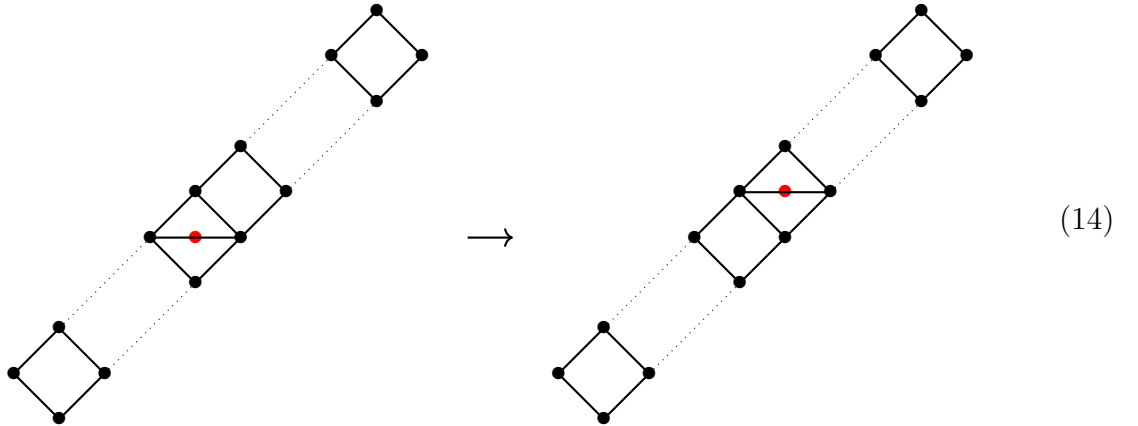
\begin{equation}\label{dessin}
\begin{gathered}
\begin{tikzpicture}[scale=0.6]

\draw (1,1) node{$\bullet$};
\draw (2,2) node{$\bullet$};
\draw (4,4) node{$\bullet$};
\draw (5,5) node{$\bullet$};
\draw (6,6) node{$\bullet$};
\draw (8,8) node{$\bullet$};
\draw (9,9) node{$\bullet$};

\draw[red] (4,5) node{$\bullet$};

\draw (0,2) node{$\bullet$};
\draw (1,3) node{$\bullet$};
\draw (3,5) node{$\bullet$};
\draw (4,6) node{$\bullet$};
\draw (5,7) node{$\bullet$};
\draw (7,9) node{$\bullet$};
\draw (8,10) node{$\bullet$};

\draw[thick](1,1) -- (2,2);
\draw[dotted](2,2) -- (4,4);
\draw[thick](4,4) -- (5,5);
\draw[thick](5,5) -- (6,6);
\draw[dotted](6,6) -- (8,8);
\draw[thick](8,8) -- (9,9);

\draw[thick](0,2) -- (1,3);
\draw[dotted](1,3) -- (3,5);
\draw[thick](3,5) -- (4,6);
\draw[thick](4,6) -- (5,7);
\draw[dotted](5,7) -- (7,9);
\draw[thick](7,9) -- (8,10);

\draw[thick](0,2) -- (1,1);
\draw[thick](1,3) -- (2,2);
\draw[thick](3,5) -- (4,4);
\draw[thick](4,6) -- (5,5);
\draw[thick](5,7) -- (6,6);
\draw[thick](7,9) -- (8,8);
\draw[thick](8,10) -- (9,9);

\draw[thick](3,5) -- (4,5);
\draw[thick](4,5) -- (5,5);

\draw[thick, ->] (10,5) -- (11,5);

\draw (13,1) node{$\bullet$};
\draw (14,2) node{$\bullet$};
\draw (16,4) node{$\bullet$};
\draw (17,5) node{$\bullet$};
\draw (18,6) node{$\bullet$};
\draw (20,8) node{$\bullet$};
\draw (21,9) node{$\bullet$};

\draw[red] (17,6) node{$\bullet$};

\draw (12,2) node{$\bullet$};
\draw (13,3) node{$\bullet$};
\draw (15,5) node{$\bullet$};
\draw (16,6) node{$\bullet$};
\draw (17,7) node{$\bullet$};
\draw (19,9) node{$\bullet$};
\draw (20,10) node{$\bullet$};

\draw[thick](13,1) -- (14,2);
\draw[dotted](14,2) -- (16,4);
\draw[thick](16,4) -- (17,5);
\draw[thick](17,5) -- (18,6);
\draw[dotted](18,6) -- (20,8);
\draw[thick](20,8) -- (21,9);

\draw[thick](12,2) -- (13,3);
\draw[dotted](13,3) -- (15,5);
\draw[thick](15,5) -- (16,6);
\draw[thick](16,6) -- (17,7);
\draw[dotted](17,7) -- (19,9);
\draw[thick](19,9) -- (20,10);

\draw[thick](12,2) -- (13,1);
\draw[thick](13,3) -- (14,2);
\draw[thick](15,5) -- (16,4);
\draw[thick](16,6) -- (17,5);
\draw[thick](17,7) -- (18,6);
\draw[thick](20,8) -- (19,9);
\draw[thick](21,9) -- (20,10);

\draw[thick](16,6) -- (17,6);
\draw[thick](17,6) -- (18,6);

\end{tikzpicture}
\end{gathered}
\end{equation}
where each bullet is an irreducible constituent, where the $4$ long diagonals are (from left to right) $X_k$, $X_{k+1}$, $X_{k+1}$, $X_{k+2}$, where the left (resp.~right) red bullet is $(V_{[1,n-k-1],\Delta}^\infty)^\vee$ (resp.~$(V_{[1,n-k-2],\Delta}^\infty)^\vee$) and where the socle (resp.~cosocle) of $\Omega^k$ or $\Omega^{k+1}$ is the leftmost (resp.~rightmost) bullet. All this follows from Theorem \ref{thm: main dR} and the description of $D_k$ in (\ref{definitionofdk}). See also the full complex $\Omega^\bullet$ for $n=4$ and $n=5$ in Figure $1$ and Figure $7$ of Appendix \ref{sec: appendix2} where we see that the dual of the smooth constituent (i.e.~the red bullet) goes up and up when moving from left to right.\bigskip

Theorem \ref{jolicomplexe} is more generally proven for holomorphic discrete series (\cite{S92}), where for instance $(V_{[1,n-k-1],\Delta}^\infty)^\vee$ is replaced by $(V_{[1,n-k-1],\Delta}^\infty)^\vee\otimes_K L(\mu_0)$ for an arbitrary dominant weight $\mu_0$ (with respect to the lower Borel).\bigskip

Let $\mathrm{St}_n^{\rm{an}}\defeq (\mathrm{Ind}_B^G 1)^{\rm{an}}/\!\sum_{B\subsetneq P}(\mathrm{Ind}_P^G 1)^{\rm{an}}$ be the locally $K$-analytic Steinberg representation of $G$ (where $B$ is the $K$-points of the lower Borel of $\GL_n$). Then Theorem \ref{jolicomplexe} is used in \cite{Qi25} to show that the $K$-vector space of homomorphisms $(\mathrm{St}_n^{\rm{an}})^\vee[1-n] \dashrightarrow \Omega^\bullet$ in the derived category of $D(G)$-modules has a natural structure of an admissible filtered $(\varphi,N)$-module in the sense of Fontaine (\cite{Fon94}) corresponding to a certain explicit ``universal'' semi-stable non-crystalline $p$-adic representation of ${\rm Gal}(\overline K/K)$.\bigskip

The second aim of this monograph is to generalize the section (\ref{splittingcomplexbis}) to $n\geq 2$.\bigskip

But there is a problem to fix. As soon as $n>2$, it is \emph{impossible} in general to ``glue'' consecutive $\Omega^k$ as was done to obtain (\ref{splittingcomplex}). For instance, already for $n=3$, one cannot ``glue'' $\Omega^0$ and $\Omega^1$ as it turns out that the $D(G)$-module (using the notation of (\ref{dessin}) for $n=3$)
\begin{equation*}
\begin{tikzpicture}[scale=0.65]

\draw (-5.5,2) node{\!\begin{xy} (-1,0)*+{X_1}="a"; (15.15,0)*+{(V_{\{1\},\Delta}^\infty)^\vee}="b"; (31.3,0)*+{X_2}="c"; (-10.4,0)*+{1}="d"; {\ar@{-}"a";"b"}; {\ar@{-}"b";"c"}; {\ar@{-}"d";"a"}\end{xy}\!};
\draw (-1.6,2) node{$\simeq$};

\draw (1,1) node{$\bullet$};
\draw (2,2) node{$\bullet$};

\draw[red] (1,2) node{$\bullet$};
\draw[red] (-1,2) node{$1$};

\draw (0,2) node{$\bullet$};
\draw (1,3) node{$\bullet$};

\draw[thick](1,1) -- (2,2);
\draw[thick](0,2) -- (1,3);
\draw[thick](0,2) -- (1,1);
\draw[thick](1,3) -- (2,2);

\draw[thick](-1,2) -- (0,2);
\draw[thick](0,2) -- (1,2);
\draw[thick](1,2) -- (2,2);

\end{tikzpicture}
\end{equation*}
(where the red bullet is $(V_{\{1\},\Delta}^\infty)^\vee$) does not exist. When $K=\Qp$ and $n=3$, the first author proved in 2019 that there exists a coadmissible length $2$ $D(G)$-module $Y_2^\flat$ such that there exist a unique coadmissible $D(G)$-module of the form
\begin{equation}\label{bemolintro}
\widetilde \Omega^{1}{}^\flat\defeq \begin{xy}
(0,0)*+{1}="b"; (9,-5)*+{X_{1}}="c"; (25,0)*+{(V_{\{1\},\Delta}^\infty)^\vee}="d"; (9,5)*+{Y_2^\flat}="e"; (42,0)*+{X_{2},}="f";
{\ar@{-}"b";"c"}; {\ar@{-}"c";"d"}; {\ar@{-}"b";"e"}; {\ar@{-}"e";"d"}; {\ar@{-}"d";"f"};
\end{xy}
\end{equation}
and non-unique coadmissible $D(G)$-modules of the form
\begin{equation*}
\widetilde\Omega^2{}^\flat\defeq \begin{xy}
(0,0)*+{Y_2{}^\flat}="a"; (16,0)*+{(V_{\{1\},\Delta}^\infty)^\vee}="b"; (32.5,0)*+{X_{2}}="c"; (47.4,0)*+{(\mathrm{St}_3^\infty)^\vee}="d";
{\ar@{-}"a";"b"}; {\ar@{-}"b";"c"}; {\ar@{-}"c";"d"}
\end{xy}
\end{equation*}
where the non-unicity is similar to the non-unicity of $\!\begin{xy} (0,0)*+{1}="a"; (10,0)*+{X_1}="b"; (25.5,0)*+{(\mathrm{St}_2^\infty)^\vee}="c"; {\ar@{-}"a";"b"}; {\ar@{-}"b";"c"}\end{xy}\!$ in (\ref{splittingcomplex}). Hence, using $Y_2{}^\flat$ and remembering that $\Omega^0\simeq \!\begin{xy} (0,0)*+{1}="a"; (9,0)*+{X_1}="b";{\ar@{-}"a";"b"}\end{xy}\!$, one still has a surjection of complexes $[\Omega^0\rightarrow \widetilde\Omega^{1}{}^\flat \rightarrow \widetilde \Omega^{2}{}^\flat]\twoheadrightarrow (\mathrm{St}_3^\infty)^\vee[-2]$ which is a quasi-isomorphism and a natural morphism of complexes $[\Omega^0\rightarrow \widetilde\Omega^{1}{}^\flat\rightarrow \widetilde\Omega^{2}{}^\flat]\rightarrow \Omega^\bullet$ which give an explicit section $(\mathrm{St}_3^\infty)^\vee[-2] \dashrightarrow \Omega^\bullet$ to the morphism $\Omega^\bullet \rightarrow H^2(\Omega^\bullet)[-2]\cong (\mathrm{St}_3^\infty)^\vee[-2]$ in the derived category of $D(G)$-modules:
\begin{equation*}
(\mathrm{St}_3^\infty)^\vee[-2]\ \longleftarrow \ [\Omega^0\rightarrow \widetilde \Omega^{1}{}^\flat \rightarrow \widetilde \Omega^{2}{}^\flat] \longrightarrow \Omega^\bullet.
\end{equation*}
Just like (\ref{splittingcomplexbis}) is better than (\ref{splittingcomplex}) when $n=2$, we can also define $\widetilde\Omega^0\defeq \!\begin{xy} (-8,0)*+{Y_1}="c"; (0,0)*+{1}="a"; (9,0)*+{X_1}="b";{\ar@{-}"a";"b"}; {\ar@{-}"c";"a"}\end{xy}\!$ where $Y_1$ is a certain finite length coadmissible $D(G)$-module $Y_1$, and modify accordingly $Y_2{}^\flat$ and $\widetilde\Omega^{1}{}^\flat$, $\widetilde\Omega^{2}{}^\flat$ into slightly larger (finite length coadmissible) $D(G)$-modules $Y_2$ and
\[\widetilde\Omega^1\defeq \begin{xy}
(-9.1,0)*+{Y_{1}}="a"; (0,0)*+{1}="b"; (9,-5)*+{X_{1}}="c"; (25,0)*+{(V_{\{1\},\Delta}^\infty)^\vee}="d"; (9,4)*+{Y_{2}}="e"; (42,0)*+{X_{2}}="f";{\ar@{-}"a";"b"}; {\ar@{-}"b";"c"}; {\ar@{-}"c";"d"}; {\ar@{-}"b";"e"}; {\ar@{-}"e";"d"}; {\ar@{-}"d";"f"}\end{xy}, \ \ \ \widetilde\Omega^2\defeq \begin{xy}
(0,0)*+{Y_2}="a"; (16,0)*+{(V_{\{1\},\Delta}^\infty)^\vee}="b"; (32.5,0)*+{X_{2}}="c"; (47.4,0)*+{(\mathrm{St}_3^\infty)^\vee}="d";{\ar@{-}"a";"b"}; {\ar@{-}"b";"c"}; {\ar@{-}"c";"d"}\end{xy}.\]
This has two advantages. The first is that, as for $n=2$, the continuous dual of $\widetilde\Omega^2$ is then similar when $K=\Qp$ to the representations of $\GL_3(\Qp)$ appearing in the completed cohomology for $3$-dimensional semi-stable non-crystalline representations of ${\rm Gal}(\overline\Qp/\Qp)$ with Hodge-Tate weights $(0,1,2)$ (for instance the dual of $\widetilde\Omega^2$ is exactly like the ``upper branch'' of the representation in \cite[(1.1)]{BD20}). The second (which was discovered in the present work) is that the $D(G)$-modules $\widetilde\Omega^0$, $\widetilde\Omega^1$, $\widetilde\Omega^2$ behave remarkably well with respect to wall-crossing functors (as will be explained in \S\ref{intermediateintro} below) contrary to the $D(G)$-modules $\Omega^0$, $\widetilde \Omega^{1}{}^\flat$, $\widetilde \Omega^{2}{}^\flat$.\bigskip

The coadmissible $D(G)$-modules $Y_k$ and $\widetilde\Omega^k$ turn out to nicely generalize to any $n\geq 2$.

\begin{thm}\label{concatenation}
There exist indecomposable multiplicity free finite length coadmissible \!$D(\!G)$-modules $Y_1, \dots, Y_{n-1}$ with irreducible socle and cosocle satisfying the following properties.
\begin{enumerate}[label=(\roman*)]
\item\label{generick}
For $k\in \{0,\dots,n-2\}$ there exists a unique coadmissible $D(G)$-module $\widetilde\Omega^k$ of the form
\begin{equation*}
\begin{xy}
(-3,0)*+{Y_{k}}="a"; (20,0)*+{(V_{[1,n-k],\Delta}^\infty)^\vee}="b"; (40,-9)*+{X_{k}}="c"; (60,0)*+{(V_{[1,n-k-1],\Delta}^\infty)^\vee}="d"; (40,9)*+{Y_{k+1}}="e"; (87,0)*+{X_{k+1}}="f";
{\ar@{-}"a";"b"}; {\ar@{-}"b";"c"}; {\ar@{-}"c";"d"}; {\ar@{-}"b";"e"}; {\ar@{-}"e";"d"}; {\ar@{-}"d";"f"};
\end{xy}
\end{equation*}
where $Y_0=X_0=(V_{[1,n],\Delta}^\infty)^\vee\defeq 0$. Moreover $\widetilde\Omega^k$ is indecomposable multiplicity free with irreducible socle and cosocle.
\item\label{genericn-1}
The set of isomorphism classes of coadmissible $D(G)$-modules $\widetilde\Omega^{n-1}$ of the form
\begin{equation*}
\begin{xy}
(0,0)*+{Y_{n-1}}="a"; (20,0)*+{(V_{\{1\},\Delta}^\infty)^\vee}="b"; (41,0)*+{X_{n-1}}="c"; (61,0)*+{(\mathrm{St}_n^\infty)^\vee}="d";
{\ar@{-}"a";"b"}; {\ar@{-}"b";"c"}; {\ar@{-}"c";"d"};
\end{xy}
\end{equation*}
is in non-canonical bijection with ${\mathbb A}^{n-1}(K)$\footnote{One has $\Dim_K\mathrm{Ext}_{D(G)}^1\big((\mathrm{St}_n^{\rm{alg}})^\vee,\!\begin{xy}(0,0)*+{Y_{n-1}}="a"; (16,0)*+{(V_{\{1\},\Delta}^\infty)^\vee\!}="b"; (32.4,0)*+{X_{n-1}}="c";{\ar@{-}"a";"b"}; {\ar@{-}"b";"c"};\end{xy}\!\big)=n$, see \ref{it: main Ext 1} of Theorem \ref{thm: main Ext}.}. Moreover any such $\widetilde\Omega^{n-1}$ is indecomposable multiplicity free with irreducible socle and cosocle. 
\end{enumerate}
\end{thm}

Contrary to the $D(G)$-modules $X_k$ in (\ref{defineXk}), the $D(G)$-modules $Y_k$ are not uniserial, they look like ``triangles'', and contrary to $X_k$ they contain duals of smooth irreducible constituents (when $k>1$), see for instance Figure $2$ and Figure $8$ in Appendix \ref{sec: appendix2} when $n=4$ and $n=5$. To represent the $D(G)$-modules $\widetilde\Omega^k$ one needs $3$-dimensional drawings (at least when $k\in \{1,\dots,n-2\}$), see for instance Figure $3$ to $6$ and Figure $9$ to $13$ in Appendix \ref{sec: appendix2} when $n=4$ and $n=5$ (where $\widetilde\Omega^k$ is denoted $\widetilde D_k$).\bigskip

The proof of Theorem \ref{concatenation} is dissiminated throughout the paper. The first statement about the $Y_k$ follows from (\ref{definitionofxyz}). The first statement of \ref{generick} of Theorem \ref{concatenation} follows from \ref{it: main rep 4} of Theorem \ref{thm: main rep} with \ref{it: unicity rep 0},\ref{it: unicity rep 2} of Corollary \ref{unicitystuff}, while the second statement follows from (\ref{definitionofdk}). The first statement of \ref{genericn-1} of Theorem \ref{concatenation} follows from \ref{it: main Ext 1} of Theorem \ref{thm: main Ext} with \ref{it: main rep 3} of Theorem \ref{thm: main rep} (for $k=n-1$) and the discussion around (\ref{equ: branch filtration 1}), while the second statement follows from \ref{it: Ext1 with alg 2} of Lemma \ref{lem: Ext1 factor 2}, the definition of $X_{n-1}$ in (\ref{definitionofxyz}) (or (\ref{defineXk})) and the definition and unicity of $Z_{n-1}$ in (\ref{definitionofxyz}) and \ref{it: unicity rep 0} of Corollary \ref{unicitystuff}. Moreover, like Theorem \ref{jolicomplexe}, Theorem \ref{concatenation} is also proved in the setting of holomorphic discrete series.\bigskip

As for $n=2$ and $n=3$ when $K=\Qp$, one can expect that, for all $n\geq 2$ and all $K$, the duals of some of the $D(G)$-modules $\widetilde\Omega^{n-1}$ in \ref{genericn-1} of Theorem \ref{concatenation} (i.e.~for some values of the parameter in ${\mathbb A}^{n-1}(K)$) occur in the completed cohomology for $n$-dimensional semi-stable non-crystalline representations of ${\rm Gal}(\overline K/K)$ with Hodge-Tate weights $(0,1,\dots,n-1)$ in all directions. The values of the parameter should then be determined by the Hodge filtration on Fontaine's filtered $(\varphi,N)$-module associated to the ${\rm Gal}(\overline K/K)$-representation, see for instance \cite[Thm.~1.1]{BD20} or \cite[Thm.~1.9]{BD23} when $n=3$ and $K=\Qp$. More generally one can wonder if the duals of \emph{all} the $D(G)$-modules $\widetilde\Omega^{k}$ for $k\in \{0,\dots,n-2\}$ in \ref{generick} of Theorem \ref{concatenation} do not also occur as \emph{subquotients} in the completed cohomology for such Galois representations (for instance this is obvious for $n=2$ and one can check it for $\widetilde\Omega^{1}$, and thus $\widetilde\Omega^{0}$, when $n=3$ and $K=\Qp$).\bigskip

Using Theorem \ref{concatenation} (and Theorem \ref{jolicomplexe}), we immediately obtain a section to $\Omega^\bullet \rightarrow H^{n-1}(\Omega^\bullet)[-(n-1)]\cong (\mathrm{St}_n^\infty)^\vee[-(n-1)]$ as follows. There are unique (up to scalar) non-zero $D(G)$-equivariant morphisms $\tld{\Omega}^{k}\rightarrow \tld{\Omega}^{k+1}$ for $k\in \{0,\dots,n-2\}$. We can thus define the complex of finite length coadmissible $D(G)$-modules (fixing an arbitrary choice for $\widetilde\Omega^{n-1}$)
\[\widetilde\Omega^\bullet\defeq [\widetilde \Omega^0\rightarrow \widetilde \Omega^1\rightarrow \cdots \rightarrow \widetilde \Omega^{n-2} \rightarrow \widetilde \Omega^{n-1}]\]
which is exact in degree $<n-1$ and has cohomology $(\mathrm{St}_n^\infty)^\vee$ in degree $n-1$. In particular there is a quasi-isomorphism $\widetilde\Omega^\bullet\rightarrow (\mathrm{St}_n^\infty)^\vee[-(n-1)]$. There are also unique (up to non-zero scalar) non-zero $D(G)$-equivariant morphisms $\tld{\Omega}^{k}\rightarrow {\Omega}^{k}$ for $k\in \{0,\dots,n-1\}$ (which are surjections) which give a morphism of complexes $\widetilde\Omega^\bullet\rightarrow \Omega^\bullet$. 

\begin{cor}[Theorem \ref{thm: main split}, Corollary \ref{cor: split dR}]\label{qiintro+}
There is an explicit section
\[(\mathrm{St}_n^\infty)^\vee[-(n-1)] \dashrightarrow \Omega^\bullet\]
in the derived category of finite length coadmissible $D(G)$-modules to the canonical morphism of complexes $\Omega^\bullet \rightarrow (\mathrm{St}_n^\infty)^\vee[-(n-1)]$ given by 
\[(\mathrm{St}_n^\infty)^\vee[-(n-1)]\longleftarrow \widetilde \Omega^\bullet \longrightarrow \Omega^\bullet.\]
\end{cor}

We see that this explicit section exists in the derived category of finite length coadmissible $D(G)$-modules with all irreducible constituents being (duals of) Orlik-Strauch representations (\ref{OSintro}). Corollary \ref{qiintro+} and (\ref{qiintro}) raise the question of the existence of analogous explicit sections to the morphisms of complexes $\tau_{\leq \ell}\Omega^\bullet \twoheadrightarrow H^{\ell}(\Omega^\bullet)[-\ell]\cong (V_{[1,n-\ell-1],\Delta}^{\infty})^\vee[-\ell]$ for $\ell\in \{0,\dots,n-2\}$. This is obvious for $\ell=0$ and not too difficult for $\ell=1$ (Proposition \ref{prop: little split}). In particular we obtain a full explicit splitting of $\Omega^\bullet$ when $n=3$ (Corollary \ref{cor: split3}). For some time we thought we could use the $D(G)$-modules $\tld{\Omega}^{k}$ to also obtain sections for all $\ell\in \{2,\dots,n-2\}$ (in the derived category of finite length coadmissible $D(G)$-modules with Orlik-Strauch irreducible constituents), but this does not seem to work. Nevertheless, we do expect such sections to exist:

\begin{conj}[Conjecture \ref{conj: split}]
For $\ell\in \{2,\dots,n-2\}$ the morphism of complexes
\[\tau_{\leq \ell}\Omega^\bullet \twoheadrightarrow H^{\ell}(\Omega^\bullet)[-\ell]\]
admits a section in the derived category of finite length coadmissible $D(G)$-modules with Orlik-Strauch irreducible constituents.
\end{conj}

We could prove the above conjecture for $n=4$ and $\ell=2$, but the complex we construct for that (analogous to $\widetilde\Omega^\bullet$ when $\ell=3$ or to the complex in Proposition \ref{prop: little split} when $\ell=1$) looks much more complicated. In particular it is not clear to us if the above complexes $\widetilde \Omega^\bullet$ admit nice generalizations to $\tau_{\leq \ell}\Omega^\bullet$ and $(V_{[1,n-\ell-1],\Delta}^{\infty})^\vee[-\ell]$ for $\ell \in \{1,\dots,n-2\}$.

\subsection{Some intermediate results and ideas of proofs}\label{intermediateintro}

If $V_0$, $V_1$ are admissible locally $K$-analytic representations of $G$ over $K$, it is not difficult to check that $\mathrm{Ext}_G^i(V_0,V_1)\cong \mathrm{Ext}_{D(G)}^i(V_1^\vee,V_0^\vee)$ for $i=0,1$ where the first $\mathrm{Ext}^i$ is computed in the abelian category of admissible locally $K$-analytic representations of $G$ (\`a la Yoneda for $i=1$) and the second $\mathrm{Ext}^i$ is computed in the category of all (abstract) $D(G)$-modules (see \cite[\S 6]{ST03} when $i=0$, \cite[Lemma 2.1.1]{Bre19} when $i=1$). Hence, in order to prove Theorem \ref{jolicomplexe} and Theorem \ref{concatenation}, it is enough to control the dimensions of the $K$-vector spaces $\mathrm{Ext}_{D(G)}^i(V_1^\vee,V_0^\vee)$ for $i=0,1$ and certain $V_0$, $V_1$. By d\'evissage, it is enough to control $\mathrm{Ext}_{D(G)}^i(V_1^\vee,V_0^\vee)$ for $i=0,1,2$ where $V_0$, $V_1$ are Orlik-Strauch representations (\ref{OSintro}) such that $\pi^\infty$ is of finite length. However, this is still a non-trivial task. For instance, when $i\ne 0$ it is not even clear that such $K$-vector spaces are finite dimensional.\bigskip

To explain our method, we need some more notation. For $i=0,1$ write $V_i=\cF_{P_{I_i}}^G(M_i,\pi_i^{\infty})$ where $I_i\subseteq \Delta$, $P_{I_i}$ is the $K$-points of the associated lower parabolic subgroup of $\GL_n$, $L_{I_i}$ its Levi factor, $M_i$ an object of ${\mathcal O}^{\fp_{I_i}}$ with integral weights and $\pi_i^\infty$ a smooth finite length representation of $L_{I_i}$. By d\'evissage, we can always reduce ourselves to the case where $M_1$ is a generalized Verma module, i.e.~$M_1\cong U(\fg)\otimes_{U(\fp_{I_1})}L^{I_1}(\mu)$ where $\mu$ is an (integral) dominant weight with respect to the lower Borel of $L_{I_{1}}$ and $L^{I_1}(\mu)$ the irreducible finite dimensional algebraic representation of $L_{I_1}$ of highest weights $\mu$. In this case $V_1$ is the locally $K$-analytic parabolic induction $(\mathrm{Ind}_{P_{I_1}}^G L^{I_1}(\mu)^\vee\!\otimes \pi_1^\infty)^{\rm{an}}$.\bigskip

By a result of Schneider-Teitelbaum when $K=\Qp$ completed by Schmidt when $K$ is arbitrary (proof of \cite[Lemma 6.3(ii)]{ST05} replacing \cite[Lemma~6.2]{ST05} by \cite[Prop.~2.6]{Schm09}), we have in this case isomorphisms for $i\geq 0$
\begin{equation}\label{schneiderschmidt}
\mathrm{Ext}_{D(G)}^i(V_1^\vee,V_0^\vee)\cong \mathrm{Ext}_{D(P_{I_1})}^i(L^{I_1}(\mu)\otimes_E (\pi_1^{\infty})^\vee, V_0^\vee)
\end{equation}
where $D(P_{I_1})$ is the $K$-algebra of locally $K$-analytic distributions on $P_{I_1}$ and $\mathrm{Ext}^i_{D(P_{I_1})}$ is computed in the category of $D(P_{I_1})$-modules. So we are reduced to computing the (hopefully finite) dimension of $\mathrm{Ext}_{D(P_{I_1})}^i(L^{I_1}(\mu)\otimes_E (\pi_1^{\infty})^\vee, V_0^\vee)$. The method is then to use a spectral sequence, but there are two possible choices.\bigskip

The first choice is to write $P_{I_1}=N_{I_1}L_{I_1}$, where $N_{I_1}$ is the unipotent radical of $P_{I_1}$, and use the spectral sequence:
\[\mathrm{Ext}_{D(L_{I_1})}^i\big(L^{I_1}(\mu)\otimes_E (\pi_1^{\infty})^\vee, \mathrm{Ext}_{D(N_{I_1})}^i(1,V_0^\vee)\big)\implies \mathrm{Ext}_{D(P_{I_1})}^{i+j}\big(L^{I_1}(\mu)\otimes_E (\pi_1^{\infty})^\vee, V_0^\vee\big).\]
This is for instance the method used in \cite[\S 4]{Schr11} or in \cite[\S 5]{Bre19} when $K=\Qp$ and $n=3$. Though one can maybe proceed this way, it is not the spectral sequence that we use in this work. One reason is that, already when $K=\Qp$ and $n=3$, using this spectral sequence proved quite laborious in \emph{loc.~cit.}\bigskip

The second choice, which is ours, is to use the spectral sequence (\cite[\S 3]{ST05}, see (\ref{equ: ST seq})):
\begin{equation}\label{equ: ST seqintro}
\mathrm{Ext}_{D^{\infty}(P_{I_1})}^i\big((\pi_1^{\infty})^\vee, \mathrm{Ext}_{U(\fp_{I_1})}^j(L^{I_1}(\mu),V_0^\vee)\big)\implies \mathrm{Ext}_{D(P_{I_1})}^{i+j}\big(L^{I_1}(\mu)\otimes_E (\pi_1^{\infty})^\vee, V_0^\vee\big)
\end{equation}
where $D^{\infty}(P_{I_1})$ is the (algebraic) dual of the locally constant $K$-valued functions on $P_{I_1}$, $\mathrm{Ext}_{D^{\infty}(P_{I_1})}^i$ is computed in the category of $D^{\infty}(P_{I_1})$-modules and $\mathrm{Ext}_{U(\fp_{I_1})}^j$ in the category of $U(\fp_{I_1})$-modules. This spectral sequence seems more suited to our dimension calculations because (as we will see in (\ref{spectralbisintro}) below) it turns out we can ``separate'' the smooth part and the Lie part on the left hand side of (\ref{equ: ST seqintro}), so that we are essentially reduced to computing dimensions of extensions groups either in the world of smooth representations or in the world of modules over Lie algebras, which is much easier (and where there is no more topology).\bigskip

We thus need to compute the left hand side of (\ref{equ: ST seqintro}). However, we do not compute it directly. Rather we define a filtration on $V_0^\vee$ and first compute
\[\mathrm{Ext}_{D^{\infty}(P_{I_1})}^i\big((\pi_1^{\infty})^\vee, \mathrm{Ext}_{U(\fp_{I_1})}^j(L^{I_1}(\mu),\mathrm{graded\ pieces})\big).\]
Denote by $W^{I_0,I_1}\subseteq W(G)$ the subset of minimal length representatives of the double coset $W(L_{I_0})\backslash W(G)/W(L_{I_1})$ where $W(-)$ are the respective (finite) Weyl groups. The Bruhat decomposition $G=\bigsqcup_{w\in W^{I_0,I_1}} P_{I_1}w^{-1}P_{I_0}$ induces a filtration indexed by $W^{I_0,I_1}$ on $V_0^\vee$ by closed $D(P_{I_1})$-submodules (see (\ref{equ: Bruhat filtration OS})), and we denote by $\mathrm{gr}_w(V_0^\vee)$ its graded pieces. We first prove the following key description of $\mathrm{Ext}_{U(\fp_{I_1})}^j(L^{I_1}(\mu), \mathrm{gr}_w(V_0^\vee))$.

\begin{thm}[Theorem \ref{prop: p coh graded} with Lemma \ref{lem: dual of compact induction}]\label{prop: p coh gradedintro}
For $w\in W^{I_0,I_1}$ and $j\geq 0$ there is a $D^{\infty}(P_{I_1})$-equivariant isomorphism of Fr\'echet spa\-ces
\begin{equation}\label{equ: p coh gradedintro}
\mathrm{Ext}_{U(\fp_{I_1})}^j(L^{I_1}(\mu), \mathrm{gr}_w(V_0^\vee))\cong \mathrm{Ext}_{U(\fp_{I_1})}^j(L^{I_1}(\mu), M_0^w)\otimes_K\big((\mathrm{ind}_{P_{I_1}\cap w^{-1}P_{I_0}w}^{P_{I_1}}\pi_0^{\infty,w})^\infty\big)^\vee
\end{equation}
where the $K$-vector space $\mathrm{Ext}_{U(\fp_{I_1})}^j(L^{I_1}(\mu), M_0^w)$ is finite dimensional and has trivial action of $D^{\infty}(P_{I_1})$. Here $(-)^w$ means that $\mathfrak x\in U(\mathfrak g)$ (resp.~$x\in P_{I_1}\cap w^{-1}P_{I_0}w$) acts by $w{\mathfrak x}w^{-1}$ (resp.~$wxw^{-1}$) and $(\mathrm{ind}_{P_{I_1}\cap w^{-1}P_{I_0}w}^{P_{I_1}}\pi_0^{\infty,w})^\infty$ is the smooth induction with compact support.
\end{thm}

The first main ingredient for the proof of Theorem \ref{prop: p coh gradedintro} is an explicit description of the continous dual $\mathrm{gr}_w(V_0^\vee)$ (Proposition \ref{lem: graded OS tensor}), where the \emph{canonical Fr\'echet completion} $\cM_0$ of $M_0$ defined in \cite{Schm13} (see also Proposition \ref{prop: category O completion}) shows up, or more precisely its twist $\cM_0^w$. The second main ingredient is the following statement.

\begin{thm}[Lemma \ref{lem: p coh M isom} and Lemma \ref{lem: p coh M separated}]\label{separateintro}
For $w\in W^{I_0,I_1}$ and $j\geq 0$ the natural morphism
\[\mathrm{Ext}_{U(\fp_{I_1})}^j(L^{I_1}(\mu), M_0^w)\longrightarrow \mathrm{Ext}_{U(\fp_{I_1})}^j(L^{I_1}(\mu), \cM_0^w)\]
is an isomorphism of finite dimensional $K$-vector spaces. Moreover $\mathrm{Ext}_{U(\fp_{I_1})}^j(L^{I_1}(\mu), \cM_0^w)$ is separated for its natural topology (coming from the topology on $\cM_0^w$ and from the Chevalley-Eilenberg complex).
\end{thm}

We remark that Theorem \ref{separateintro} is also proved in \cite{BCGP25}.\bigskip

Using these two key ingredients, the proof of Theorem \ref{prop: p coh gradedintro} then consists in carefully analysing the Chevalley-Eilenberg complex of $\mathrm{gr}_w(V_0^\vee)$. It is given in \S\ref{subsec: loc an spectral seq}, and is quite long and tedious because we give all the (topological) technical details. Note that, had the topological vector space $\mathrm{Ext}_{U(\fp_{I_1})}^j(L^{I_1}(\mu), \cM_0^w)$ been non-separated, the solid techniques of \cite{RR22} would probably have been necessary.\bigskip

Plugging in (\ref{equ: p coh gradedintro}) inside the left hand side of (\ref{equ: ST seqintro}) with $\mathrm{gr}_w(V_0^\vee)$ instead of $V_0^\vee$, we deduce canonical isomorphisms for $w\in W^{I_0,I_1}$ and $i,j\geq 0$ (see Corollary \ref{cor: Ext P graded})
\begin{multline}\label{equ: extp coh gradedintro}
\mathrm{Ext}_{D^{\infty}(P_{I_1})}^i\big((\pi_1^{\infty})^\vee, \mathrm{Ext}_{U(\fp_{I_1})}^j(L^{I_1}(\mu), \mathrm{gr}_w(V_0^\vee))\big)\\
\cong \mathrm{Ext}_{U(\fp_{I_1})}^j(L^{I_1}(\mu), M_0^w)\otimes_K \mathrm{Ext}_{L_{I_1}}^i\Big(\big((\mathrm{ind}_{P_{I_1}\cap w^{-1}P_{I_0}w}^{P_{I_1}}\pi_0^{\infty,w})^\infty\big)_{N_{I_1}}, \pi_1^{\infty}\Big)^{\infty}
\end{multline}
where $\mathrm{Ext}_{L_{I_1}}^i(-)^\infty$ means extensions in the category of smooth representations of $L_{I_1}$. Moreover the $K$-vector space on the right hand side of (\ref{equ: extp coh gradedintro}) is finite dimensional for all $w,i,j$ (using Theorem \ref{separateintro} for $\mathrm{Ext}_{U(\fp_{I_1})}^j(L^{I_1}(\mu), M_0^w)$ and the finite length of $\pi_0^{\infty}$, $\pi_1^{\infty}$ for the other factor). Hence we obtain from the analogue of (\ref{equ: ST seqintro}) with $\mathrm{gr}_w(V_0^\vee)$ instead of $V_0^\vee$ that the $K$-vector space $\mathrm{Ext}_{D(P_{I_1})}^{i}(L^{I_1}(\mu)\otimes_E (\pi_1^{\infty})^\vee, \mathrm{gr}_w(V_0^\vee))$ is also finite dimensional for all $i\geq 0$. By d\'evissage on the filtration on $V_0^\vee$ and by (\ref{schneiderschmidt}), we deduce the following nice by-product result, which is new.

\begin{thm}[Theorem \ref{thm: finitedim}]\label{finiteintro}
Let $\Pi_0$, $\Pi_1$ be finite length admissible locally $K$-analytic representations of $G$ over a finite extension of $K$ with all (topological) irreducible constituents being Orlik-Strauch representations (\ref{OSintro}). Then the $K$-vector space $\mathrm{Ext}_{D(G)}^i(\Pi_0^\vee,\Pi_1^\vee)$ is finite dimensional for $i\geq 0$.
\end{thm}

We mention here another aside result, purely on the Lie algebra side, that we need in the proofs and that we couldn't find in the literature (it is also proved in \cite{BCGP25}):

\begin{prop}[Proposition \ref{lem: n coh O}]
Let $I$ a subset of $\Delta$. Denote by $\fb$ the $K$-Lie algebra of the lower Borel of $G$, $\fn_I$, $\fl_I$ the $K$-Lie algebras of $N_I$, $L_I$, and $\fb_I$ the $K$-Lie algebra of the lower Borel of $L_I$. Then for any $M$ in $\cO^{\fb}$ with integral weights and any $j\geq 0$ the Lie algebra cohomology group $H^j(\fn_{I}, M)$ is an object of the category $\cO^{\fb_I}$ relative to $\fl_I$ (with integral weights). 
\end{prop}

Going back to (\ref{equ: extp coh gradedintro}), in practice we only need to apply it with $\pi_0^{\infty}$ such that the smooth $L_{I_1}$-representa\-tions $((\mathrm{ind}_{P_{I_1}\cap w^{-1}P_{I_0}w}^{P_{I_1}}\pi_0^{\infty,w})^\infty)_{N_{I_1}}$ lie in distinct Bernstein blocks when $w$ varies in $W^{I_0,I_1}$. Assuming that $\pi_1^{\infty}$ lies in a unique Bernstein block (which will be our case), we thus deduce that there is \emph{at most one} $w\in W^{I_0,I_1}$ only depending on $\pi_0^{\infty}$ and $\pi_1^{\infty}$ such that 
\[\mathrm{Ext}_{L_{I_1}}^i\!\Big(\!\big((\mathrm{ind}_{P_{I_1}\cap w^{-1}P_{I_0}w}^{P_{I_1}}\pi_0^{\infty,w})^\infty\big)_{N_{I_1}}, \pi_1^{\infty}\Big)^{\infty}\ne 0\]
and hence such that (\ref{equ: extp coh gradedintro}) is possibly non-zero. By an obvious d\'evissage on the filtration on $V_0^\vee$, we thus see that, in order to compute the left hand side of (\ref{equ: ST seqintro}), it is sufficient (in our case) to compute (\ref{equ: extp coh gradedintro}) because there is at most one $w$ for which (\ref{equ: extp coh gradedintro}) is non-zero! In particular by (\ref{schneiderschmidt}) the spectral sequence (\ref{equ: ST seqintro}) becomes
\begin{multline}\label{spectralbisintro}
\mathrm{Ext}_{U(\fp_{I_1})}^j(L^{I_1}(\mu), M_0^w)\otimes_K \mathrm{Ext}_{L_{I_1}}^i\Big(\big((\mathrm{ind}_{P_{I_1}\cap w^{-1}P_{I_0}w}^{P_{I_1}}\pi_0^{\infty,w})^\infty\big)_{N_{I_1}}, \pi_1^{\infty}\Big)^{\infty}\\
\implies \mathrm{Ext}_{D(G)}^{i+j}(V_1^\vee,V_0^\vee)
\end{multline}
for this unique $w$ (with possibly all vector spaces being $0$).\bigskip

Now, the big advantage of (\ref{spectralbisintro}) is that the term on its left hand side is computable because $\mathrm{Ext}_{U(\fp_{I_1})}^j(L^{I_1}(\mu), M_0^w)$ \ is \ purely \ in \ the \ category \ of \ $U(\fp_{I_1})$-modules \ while $\mathrm{Ext}_{L_{I_1}}^i(((\mathrm{ind}_{P_{I_1}\cap w^{-1}P_{I_0}w}^{P_{I_1}}\pi_0^{\infty,w})^\infty)_{N_{I_1}}, \pi_1^{\infty})^{\infty}$ is purely in the category of smooth representations of $L_{I_1}$ (in other words we have ``separated'' the Lie part and the smooth part). We give in \S\ref{sec: smooth rep} all the material needed to compute the dimensions of the $K$-vector spaces $\mathrm{Ext}_{L_{I_1}}^\bullet(((\mathrm{ind}_{P_{I_1}\cap w^{-1}P_{I_0}w}^{P_{I_1}}\pi_0^{\infty,w})^\infty)_{N_{I_1}}, \pi_1^{\infty})^{\infty}$, at least for the $\pi_i^{\infty}$ we need, and we give in \S\ref{sec: n coh} all the material needed to compute the dimensions of $\mathrm{Ext}_{U(\fp_{I_1})}^\bullet(L^{I_1}(\mu), M_0^w)$. Very often, either all terms on the left hand side of (\ref{spectralbisintro}) are $0$, and we obtain the vanishing of $\mathrm{Ext}_{D(G)}^i(V_1^\vee,V_0^\vee)$, or only one of them is non-zero, and we obtain a useful description of $\mathrm{Ext}_{D(G)}^i(V_1^\vee,V_0^\vee)$. Combining this with some d\'evissage, we can compute the dimensions of lots of $\mathrm{Ext}_{D(G)}^i(V_1^\vee,V_0^\vee)$ (for $M_1$ in ${\mathcal O}^{\fp_{I_1}}$ not necessarily a generalized Verma module), see the various results in \S\S\ref{subsec: Ext OS}, \ref{subsec: square}. As a sample, let us mention the following statement:

\begin{cor}[Proposition \ref{prop: Ext1 OS}]
Let $\un{j}=(j_0,j_1,j_2)$ and $\un{j}'=(j'_0,j'_1,j'_2)$ as in (\ref{Jintro}) and assume $(j_0,j_1)\ne (j'_0,j'_1)$, there is an isomorphism of (finite dimensional) $K$-vector spaces
\begin{multline}\label{sampleintro}
\mathrm{Ext}_{D(G)}^1(C_{\un{j}}^\vee, C_{\un{j}'}^\vee)\\
\cong \mathrm{Ext}_{U(\fg)}^1(L(w_{j_1,j_0}\cdot 0),L(w_{j_1',j_0'}\cdot 0))\otimes_K \Hom_{L_{\widehat{j}_1}}\Big(\big((\mathrm{ind}_{P_{\widehat j_1}\cap P_{\widehat j'_1}}^{P_{\widehat j_1}}\pi_{j_1',j_2'}^{\infty})^\infty\big)_{N_{\widehat j_1}},\pi_{j_1,j_2}^{\infty}\Big).
\end{multline}
\end{cor}

One then deduces $\Dim_K\mathrm{Ext}_{D(G)}^1(C_{\un{j}}^\vee, C_{\un{j}'}^\vee)$ by computing the dimension of the right hand side of (\ref{sampleintro}), see for instance Lemma \ref{lem: Ext1 factor 1}.\bigskip

Pushing further this kind of arguments, we prove the existence and unicity of the $D(G)$-modules $X_k$ in (\ref{defineXk}), of the $D(G)$-modules $\!\begin{xy}(0,0)*+{X_{k}}="a"; (22,0)*+{(V_{[1,n-k-1],\Delta}^\infty)^\vee}="b"; (46.3,0)*+{X_{k+1}}="c";
{\ar@{-}"a";"b"}; {\ar@{-}"b";"c"}\end{xy}\!$ in Theorem \ref{jolicomplexe}, and of all the $D(G)$-modules in Theorem \ref{concatenation}. Note that, apart from the $X_k$, most of these $D(G)$-modules are not uniserial: they involve subquotients which look like ``squares'' or ``cubes''. We call them Ext-squares or Ext-cubes and we construct them in \S\ref{subsec: square} (and \S\ref{U(g)square} for the Lie counterpart). For instance the Ext-squares involved in the $D(G)$-modules (\ref{dessin}) are computed in \ref{it: easy square 2} of Proposition \ref{prop: easy square}, the Ext-squares involved in the $D(G)$-modules $Y_k$ of Theorem \ref{concatenation} are computed in Proposition \ref{prop: hard square}, while the $D(G)$-modules $\widetilde\Omega^k$ of Theorem \ref{concatenation} also involve Ext-cubes which are computed in Proposition \ref{prop: factor cube}. These computations can be long and technical, but since we give all details they should not be hard to follow by an interested reader.\bigskip

(It is worth mentioning here that the smooth generalized Steinberg $V_{[1,n-k-1],\Delta}^{\infty}$ of (\ref{steinbergintro}) and the smooth representations $\pi_{j_1,j_2}^{\infty}$ in (\ref{cjintro}) are instances of more general smooth representations that we call \emph{$G$-basic} (see Definition \ref{def: basic rep}) and that we completely study in \S\ref{sec: smooth rep} (see \S\ref{subsec: sm example} for the special case of the $\pi_{j_1,j_2}^{\infty}$). These $G$-basic representations are quite convenient because one can easily compute explicitly everything that is needed for this work: Jacquet functors, constituents of smooth parabolic inductions, smooth extension groups, etc.)\bigskip

But this is still not the end of the proof of Theorem \ref{jolicomplexe}.\bigskip

For \ $k\in \{0,\dots,n-1\}$ \ denote \ by \ $D_k$ \ the \ unique \ coadmissible \ $D(G)$-module $\!\begin{xy}(0,0)*+{X_{k}}="a"; (22,0)*+{(V_{[1,n-k-1],\Delta}^\infty)^\vee}="b"; (46.3,0)*+{X_{k+1}}="c";{\ar@{-}"a";"b"}; {\ar@{-}"b";"c"}\end{xy}\!$ of Theorem \ref{jolicomplexe}. The $D(G)$-modules $D_k$ and all the $D(G)$-modules of Theorem \ref{concatenation} are constructed (in \S\ref{subsec: final}) \emph{independently} of the Drinfeld space. This is enough for Theorem \ref{concatenation} where there is no mention of the Drinfeld space. But in order to complete the proof of Theorem \ref{jolicomplexe}, we need to show $\Omega^k\cong D_k$. This is non-trivial and uses two more key ingredients that we explain now.\bigskip

The first ingredient is a unicity statement which strengthens Schneider-Teitelbaum's Theorem \ref{boundary}.

\begin{thm}\label{boundaryunicity}
\hspace{2em}
\begin{enumerate}[label=(\roman*)]
\item\label{unicityD}
The coadmissible $D(G)$-module $D_{n-1}=\!\begin{xy}(0,0)*+{X_{n-1}}="a"; (18,0)*+{(\mathrm{St}_n^\infty)^\vee}="b"; 
{\ar@{-}"a";"b"}\end{xy}\!$ of Theorem \ref{jolicomplexe} (for $k=n-1$) is the unique coadmissible $D(G)$-module $D$ which admits a filtration by closed $D(G)$-submodules
\[D=\mathrm{Fil}^0(D)\supsetneq \mathrm{Fil}^1(D)\supsetneq \cdots \supsetneq\mathrm{Fil}^{n}(D)=0\]
satisfying (\ref{filtrationintro}) and such that $H^0(N_{\widehat {n-1}}, D)\cong s_{n-1} s_{n-2}\cdots s_1\cdot 0$ (recall $N_{\widehat {n-1}}$ is the unipotent radical of $P_{\widehat{n-1}}$).
\item\label{Omegan-1}
We have $H^0(N_{\widehat {n-1}}, \Omega^{n-1})\cong s_{n-1} s_{n-2}\cdots s_1\cdot 0$.
\end{enumerate}
\end{thm}

Part \ref{unicityD} is proven in Theorem \ref{prop: top deg}. Its proof is a bit long but it is remarkable that the condition $H^0(N_{\widehat {n-1}}, D)\cong s_{n-1} s_{n-2}\cdots s_1\cdot 0$ is enough to ``rigidify everything'' and ensure unicity. In particular all extensions between the graded pieces of the filtration $\mathrm{Fil}^\bullet(D)$ are then non-split. The easier part \ref{Omegan-1} is proven in Lemma \ref{lem: dR complex N coh} for any $k\in \{0,\dots,n-1\}$ (not just $k=n-1$). Note that both parts were already proven in \cite{Schr11} when $G=\GL_3(\Qp)$ (see \cite[Prop.~6.3]{Schr11} and \cite[\S 6.4]{Schr11}), and the strategy in \emph{loc.~cit.}~actually inspired Theorem \ref{boundaryunicity}.\bigskip

From Theorem \ref{boundary} and Theorem \ref{boundaryunicity} we immediately obtain Theorem \ref{jolicomplexe} for $k=n-1$. We now need to go from $k=n-1$ to smaller $k$, and this is where we use the second ingredient: translation functors.\bigskip

Translation functors were first defined in the setting of representations of real Lie algebras and of $(\fg, K)\footnote{Here $K$ is of course a maximal compact subgroup, not the field $K$!}$-modules in the late seventies. Together with wall-crossing functors, they quickly became a major tool in the study of BGG's category ${\mathcal O}^{\fp}$, see for instance \cite[\S 7]{Hum08}. Strangely, it is only recently that these functors have been introduced in the framework of locally analytic representations (\cite{JLS24}), and they are only beginning to be used (\cite{Di24}). Let us quickly recall their definition.\bigskip

Let $Z(\fg)$ be the center of $U(\fg)$, a $D(G)$-module $D$ is said to be $Z(\fg)$-finite, if for any $v\in M$, the $K$-vector subspace $\langle Z(\fg)v \rangle$ of $D$ generated by elements in $Z(\fg)v$ is finite dimensional. Let $\lambda$, $\mu$ be integral weights, $\xi_\lambda, \xi_{\mu}: Z(\fg)\rightarrow K$ the infinitesimal characters of the respective Verma modules $U(\fg)\otimes_{U(\fb)}\lambda$, $U(\fg)\otimes_{U(\fb)}\mu$, and $\overline{\lambda-\mu}$ the unique highest weight (relative to the lower Borel of $\GL_n$) in the $W(G)$-orbit of $\lambda-\mu$ for the standard action of $W(G)$ (not the dot action). We denote by $L(\overline{\lambda-\mu})$ the unique finite dimensional algebraic representation of $G$ over $K$ with highest weight $\overline{\lambda-\mu}$. Then for any $Z(\fg)$-finite $D(G)$-module $D$ we define (following \cite{JLS24})
\begin{equation}\label{translationintro}
\cT_\lambda^\mu (D) \defeq \mathrm{pr}_{\mu}\big(L(\overline{\lambda-\mu})\otimes_K \mathrm{pr}_{\lambda}(D)\big)
\end{equation}
where $\mathrm{pr}_{\lambda}$, $\mathrm{pr}_{\mu}$ is the projection onto the generalized eigenspace for the infinitesimal character $\xi_\lambda$, $\xi_\mu$ (which is well-defined thanks to the $Z(\fg)$-finiteness). The $D(G)$-module $\cT_\lambda^\mu(D)$ is still $Z(\fg)$-finite and the functor $D\rightarrow \cT_\lambda^\mu(D)$ is exact on $Z(\fg)$-finite $D(G)$-modules. It is called a translation functor. Moreover with the notation of (\ref{OSintro}) one has $\cT_\lambda^\mu(\cF_{P}^G(M,\pi^{\infty}))\cong \cF_{P}^G(T_\lambda^\mu M,\pi^{\infty})$ where $T_\lambda^\mu$ is the usual translation functor on the category ${\mathcal O}^{\fp}$ (\cite[\S 7.1]{Hum08}, in fact $T_\lambda^\mu$ is defined exactly as in (\ref{translationintro}) replacing $D$ by $M$).\bigskip

Translation functors are already interesting. For instance one can use them to revisit the construction of \cite{S92} and show that Schneider's holomorphic discrete series are translations of $\Omega^\bullet$, see Lemma \ref{lem: dR to discrete}. However, there exist even more interesting functors.\bigskip

Let $\rho$ be half the sum of the roots of the lower Borel of $\GL_n$ and fix a simple reflection $s_k$ of $\GL_n$ for $1\leq k\leq n-1$. Let $\mu$ be an integral weight such that $\langle\mu+\rho,\alpha^\vee\rangle\geq 0$ for $\al\in \Delta$ and the stabilizer of $\mu$ in the Weyl group $W(G)$ for the dot action is $\{1,s_k\}$. Then we define the \emph{wall-crossing} functor on $Z(\fg)$-finite $D(G)$-modules as
\[\Theta_\mu\footnote{Experts on wall-crossing functors are more used to the notation $\Theta_{s_k}$. However we do not know if the functor $\Theta_\mu$ defined on all $Z(\fg)$-finite $D(G)$-modules only depends on $s_k$ and not on the choice of $\mu$ (as it does when defined on the category ${\mathcal O}^{\fb}$), see Remark \ref{independence}. This doesn't affect this work.}\defeq \cT_{w_0\cdot 0}^\mu \circ \cT_\lambda^{w_0\cdot 0}\]
where $w_0$ is the longest element of $W(G)$ and $w_0\cdot 0=w_0(\rho)-\rho$. The terminology comes from that fact that the stabilizer of $w_0\cdot 0$ in $W(G)$ (for the dot action) is trivial, which is not the case for $\mu$, see \cite[\S 7.15]{Hum08}. By the argument in \cite[\S 7.2]{Hum08} one has canonical and functorial adjunction morphisms of ($Z(\fg)$-finite) $D(G)$-modules $D\rightarrow \Theta_\mu(D)$ and $\Theta_\mu(D) \rightarrow D$ which are non-zero when $D$ and $\Theta_\mu(D)$ are non-zero (see (\ref{adjunctionlocadm})).\bigskip

Although their definition may look puzzling at first, wall-crossing functors have remarkable properties. For instance they behave very well on the $D(G)$-modules $D_k$ and $\widetilde\Omega^k$:

\begin{thm}\label{wallcrossingintro}
Let $1\leq k\leq n-1$ and let $\mu$ be an integral weight such that $\langle\mu+\rho,\alpha^\vee\rangle\geq 0$ for $\al\in \Delta$ and the stabilizer of $\mu$ in the Weyl group $W(G)$ for the dot action is $\{1,s_k\}$.
\begin{enumerate}[label=(\roman*)]
\item \label{it: D wall crossingintro0} 
We have non-split short exact sequences of coadmissible $D(G)$-modules $0\rightarrow D_k\rightarrow \Theta_{\mu}(D_k)\rightarrow D_{k-1}\rightarrow 0$ where $D_k\rightarrow \Theta_{\mu}(D_k)$ is the canonical adjunction map.
\item \label{it: D wall crossingintro1} 
We have short exact sequences of coadmissible $D(G)$-modules $0\rightarrow \Omega^k\rightarrow \Theta_{\mu}(\Omega^k)\rightarrow \Omega^{k-1}\rightarrow 0$.
\item \label{it: D wall crossingintro2} 
If $k\ne n-1$ we have non-split short exact sequences of coadmissible $D(G)$-modules $0\rightarrow \tld{\Omega}^k\rightarrow \Theta_{\mu}(\tld{\Omega}^k)\rightarrow \tld{\Omega}^{k-1}\rightarrow 0$ where $\tld{\Omega}^k\rightarrow \Theta_{\mu}(\tld{\Omega}^k)$ is the canonical adjunction map.
\end{enumerate}
\end{thm}

Part \ref{it: D wall crossingintro0} of Theorem \ref{wallcrossingintro} is proven in \ref{it: D wall crossing 2} of Theorem \ref{thm: D wall crossing}, part \ref{it: D wall crossingintro1} of Theorem \ref{wallcrossingintro} is proven in Lemma \ref{lem: discrete wall crossing} and part \ref{it: D wall crossingintro2} of Theorem \ref{wallcrossingintro} is proven in \ref{it: D wall crossing 3} of Theorem \ref{thm: D wall crossing}. When $k=n-1$, \ref{it: D wall crossingintro2} of Theorem \ref{wallcrossingintro} is not true anymore, but it is almost true, see \ref{it: D wall crossing 4} of Theorem \ref{thm: D wall crossing}. Note that, going back to $n=3$ and $\Omega^0$, $\widetilde \Omega^{1}{}^\flat$ in (\ref{bemolintro}), we do \emph{not} have a short exact sequence $0\rightarrow \widetilde \Omega^{1}{}^\flat \rightarrow \Theta_\mu(\widetilde \Omega^{1}{}^\flat)\rightarrow \Omega^0\rightarrow 0$ (some constituents are missing in $\Omega^0$).\bigskip

Now, when $k=n-1$, we know that $\Omega^{n-1}\cong D_{n-1}$. Using \ref{it: D wall crossingintro1} of Theorem \ref{wallcrossingintro} and \ref{it: D wall crossingintro0} of Theorem \ref{wallcrossingintro} for $k=n-1$, we can deduce $\Omega^{n-2}\cong D_{n-2}$. Applying again \emph{loc.~cit.}~for $k=n-2$ we obtain $\Omega^{n-3}\cong D_{n-3}$ and so on, see Theorem \ref{thm: main dR}. We finally obtain $\Omega^k\cong D_k$ for all $k$, which finishes the proof of Theorem \ref{jolicomplexe}. Note that \emph{a posteriori} the exact sequence in \ref{it: D wall crossingintro1} of Theorem \ref{wallcrossingintro} is also non-split and the injection is also the adjunction map.\bigskip

Other nice properties are satisfied. For instance, for $1\leq k\leq n-1$, the composition $\Theta_\mu(\Omega^k)\twoheadrightarrow \Omega^{k-1}\rightarrow \Omega^k$, where the surjection is in \ref{it: D wall crossingintro1} of Theorem \ref{wallcrossingintro} and the second map is the differential map on $\Omega^\bullet$, is nothing other than the other adjunction morphism $\Theta_\mu(\Omega^k)\rightarrow \Omega^k$, and likewise with $\widetilde\Omega^k$, see Remark \ref{secondadjunction}. Moreover the wall-crossing functors $\Theta_\mu$ are important tools for constructing Ext-squares and Ext-cubes of $D(G)$-modules, see for instance (among many other statements) Proposition \ref{prop: Jantzen middle}, Lemma \ref{lem: g square as wall crossing}, Lemma \ref{lem: square as wall crossing}, Lemma \ref{lem: cube step 2}, Proposition \ref{prop: general wall crossing}, etc.\bigskip

We finally briefly give the contents of each section. In \S\ref{sec: smooth rep} we introduce (smooth) $G$-basic representations, among which are the $V_{[1,n-k-1],\Delta}^{\infty}$ and the $\pi_{j_1,j_2}^{\infty}$, and we prove all the necessary material on smooth representations for the later sections. In \S\ref{sec: n coh} we prove all results on $U(\fg)$-modules and $\mathrm{Ext}$ groups of $U(\fg)$-modules needed in \S\ref{sec: spectral seq} and \S\ref{sec: extension}. In \S\ref{sec: spectral seq} we prove Theorem \ref{prop: p coh gradedintro} and deduce many consequences (e.g.~Theorem \ref{finiteintro}). Finally, in \S\ref{sec: extension} we construct all the relevant finite length coadmissible $D(G)$-modules and we prove the main results of \S\ref{introresults}. Appendix \ref{sec: appendix} is devoted to technical combinatorial lemmas on certain elements of $W(G)$ while Appendix \ref{sec: appendix2} gives pictures for most of the previous $D(G)$-modules when $G=\GL_4(K)$ and $G=\GL_5(K)$, which is fairly representative of the general case.

\subsection{Some general notation}\label{generalnotation}

We end up this introduction with general notation which will be used throughout this work. More specialized notation will be gradually introduced within the text.\bigskip

We fix $E$ an arbitrary finite extension of $K$ which can be $K$. If $x\in K$ we let $|x|_K\defeq q^{-e\val(x)}$ where $q$ is the cardinality of the residue field of $K$, $e$ the ramification index of $K$ and val is normalized by $\val(p)=1$.\bigskip

We fix an integer $n\geq 2$ and, \emph{unless otherwise stated}, $G$ is either the algebraic group $\GL_n/K$ or its $K$-points $\GL_n(K)$ (the context being clear). We denote by $T$ the diagonal matrices in $G$, $B$ the lower triangular matrices of $G$ and $U$ the lower unipotent matrices of $B$ (so $B=UT$). We let $\fg=\fg\fl_n$, $\ft$ the Lie algebra of $T$, $\fb$ the one of $B$, $\fu\subseteq \fb$ its radical, $\fb^+$ the upper Borel, $\fu^+\subseteq \fb^+$ its radical and $Z(\fg)$ the center of the enveloping algebra $U(\fg)$. By extension of scalars from $K$ to $E$, all Lie algebras are considered as $E$-vector spaces.\bigskip

We let $\Lambda\defeq X(T)\simeq \Z^n$ be the integral weights, $\Lambda^{\rm dom}$ the set of (integral) dominant weights \emph{with respect to $\fb$}, $\Phi^+$ the set of positive roots in $\fu^+$ and $\Delta\subseteq \Phi^+$ the subset of positive simple roots. We fix a bijection $\{1,\cdots,n-1\}\cong \Delta$ sending $j$ to $e_j-e_{j+1}$ and we use both notation $\al$ or $j\in \{1, \dots, n-1\}$ for an element of $\Delta$. We set $\widehat{j}\defeq \Delta\setminus\{j\}$ for $j\in\Delta$.\bigskip

For $I\subseteq \Delta$ a subset, we let $P_I$ (resp.~$P_I^+$) be the ($K$-points of) the lower (resp.~upper) standard parabolic subgroup of $\GL_n$ associated to $I$, $N_I$ the unipotent radical of $P_I$ and $L_I\cong P_I/N_I$ its Levi factor (so $P_I=N_IL_I$). We denote by $B_I=B\cap L_I$ the lower Borel of $L_I$ and $U_I \subseteq B_I$ its unipotent radical. We let $\fp_I$, $\fp_I^+$, $\fn_I\subseteq \fp_I$, $\fl_I\cong \fp_I/\fn_I$, $\fb_I=\fb\cap\fl_I$ and $\fu_I\subseteq \fb_I$ be the respective Lie algebras of $P_I$, $P_I^+$, $N_I$, $L_I$, $B_I$ and $U_I$. We let $Z(\fl_{I})$ be the center of the enveloping algebra $U(\fl_{I})$, $\fn_I^+$ the radical of $\fp_I^+$, $\fb_I^+=\fb^+\cap \fl_I$ the upper Borel of $\fl_I$ and $\fu_I^+\subseteq \fb_I^+$ the radical of $\fb_I^+$. We denote by $\Lambda_I^{\rm dom}$ the set of (integral) dominant weights with respect to $\fb_I$ and $\Phi_I^+\subseteq \Phi^+$ the roots of $\fu_I^+$.\bigskip

We let $W(G)$ be Weyl group of $G$, $W(L_I)$ the Weyl group of the Levi $L_I$ corresponding to $I$, $\ell(w)\in \Z_{\geq 0}$ the length of $w\in W(G)$ and $W^{I_0,I_1}(L_I)$ for $I_0,I_1\subseteq I$ the set of minimal length representatives of $W(L_{I_0})\backslash W(L_I)/W(L_{I_1})$ (see \cite[Lemma~5.4]{DM91}). When $L_I=G$ we write $W^{I_0,I_1}$. We let $\rho_I$ be half the sum of the roots in $\fb_I$ (so $\langle\rho_I, \alpha^\vee\rangle =-1$ for $\alpha \in I$ and $\rho_I$ is dominant with respect to $\fb_I$) and we define the dot action $w\cdot \mu \defeq w(\mu + \rho_I) - \rho_I$ for $w\in W(L_I)$ and $\mu\in X(T)$. Note that, since $\rho - \rho_I$ is invariant under $W(L_I)$, the choice of $\rho$ or $\rho_I$ doesn't change the above dot action. When $I=\emptyset$, we forget the index $I$ in the notation. We denote by $w_I$ the longest element in $W(L_I)$ and $w_0\defeq w_\Delta$. We endow $W(L_I)$ with the Bruhat order $<$ and for $w\in W(L_I)$, we set
\begin{eqnarray}
\label{equ: left set} D_L(w)&\defeq &\{\al\in I\mid s_\alpha w<w\}=\{\al\in I\mid - w^{-1}(\alpha)\in \Phi_I^+\}\\
\label{equ: right set} D_R(w)&\defeq &\{\al\in I\mid w s_\alpha<w\}=\{\al\in I\mid - w(\alpha)\in \Phi_I^+\}
\end{eqnarray}
(so $D_L(1)=D_R(1)=\emptyset$). As $x<w$ if and only if $x^{-1}<w^{-1}$, we have $D_L(w)=D_R(w^{-1})$.\bigskip

If $\cC$ is an abelian category, we recall that a finite length object in $\cC$ is an object $M$ of $\cC$ such that there exists finitely many subobjects $M_1\subseteq M_2\subseteq \cdots \subseteq M_m=M$ in $\cC$ such that $M_i/M_{i-1}$ is a simple object of $\cC$. We then write $\mathrm{JH}_{\cC}(M)$ for the (finite) set of isomorphism classes of the simple subquotients of $M$. If $M$ is multiplicity free, we equip $\mathrm{JH}_{\cC}(M)$ with the following partial order: given $M_1$, $M_2$ in $\mathrm{JH}_{\cC}(M)$ we write $M_1 \leq M_2$ if $M_1$ is a subquotient of the unique subobject of $M$ with cosocle isomorphic to $M_2$. If $\cC$ is a full subcategory of the category of left $A$-modules (for $A$ any unital associative ring) or of representations of $G$ (for $G$ any group), we write $\mathrm{JH}_A(M)$, $\mathrm{JH}_G(M)$ instead of $\mathrm{JH}_{\cC}(M)$.\bigskip

If $A$ is a unital associative ring, we denote by $\mathrm{Mod}_{A}$ the (abelian) category of all abstract left $A$-modules. If $M$ is a finite length object in $\mathrm{Mod}_{A}$, we write
\begin{equation*}
\cdots\subseteq \mathrm{Rad}^{k+1}(M)\subseteq \mathrm{Rad}^k(M)\subseteq \cdots\subseteq \mathrm{Rad}^1(M)\subseteq \mathrm{Rad}^0(M)=M
\end{equation*}
for the radical filtration of $M$, where $\mathrm{Rad}^{k+1}(M)$ is the minimal submodule of $\mathrm{Rad}^k(M)$ such that $\mathrm{Rad}_k(M)\defeq \mathrm{Rad}^k(M)/\mathrm{Rad}^{k+1}(M)$ is semi-simple. We also write $\mathrm{rad}(M)=\mathrm{Rad}^1(M)$ (the radical of $M$). The Loewy length of $M$, written $\ell\ell(M)$ is by definition the minimal integer $k\geq 0$ such that $\mathrm{Rad}^k(M)=0$. Similarly we define the socle filtration $\cdots\supseteq \mathrm{Soc}^k(M)\supseteq \cdots\supseteq \mathrm{Soc}^0(M)=0$ of $M$ where $\mathrm{Soc}^{k+1}(M)$ is the maximal submodule of $M$ containing $\mathrm{Soc}^k(M)$ such that $\mathrm{Soc}^{k+1}(M)/\mathrm{Soc}^k(M)$ is semi-simple. Note that $\ell\ell(M)$ is also the minimal integer $k\geq 0$ such that $\mathrm{Soc}^{k+1}(M)/\mathrm{Soc}^k(M)=M$. The module $M$ is called \emph{rigid} if $\mathrm{Rad}^k(M)=\mathrm{Soc}^{\ell\ell(M)-k}(M)$ for $0\leq k\leq \ell\ell(M)$.\bigskip

If $\pi$ is a representation of a subgroup $H'$ of some group $H$ and $h\in H$, we denote by $\pi^h$ the representation of $h^{-1}H'h$ with same underlying space as $\pi$ where $h'\in h^{-1}H'h$ acts by $hh'h^{-1}$. If $H$ is a locally compact group, $\mathrm{Ext}_{H}^\bullet(-, -)^{\infty}$ means extensions in the category of smooth representations of $H$ over $E$-vector spaces. If $\pi$ is a smooth representation of $L_I$ (for some $I\subseteq \Delta$), we denote by $\pi^\sim$ its smooth contregredient.\bigskip

This monograph is the culmination of a long-term work. At various stages of its gestation, we benefited from discussions with several people. In particular we would like to thank Laurent Clozel, Jean-Fran\c cois Dat, Yiwen Ding, Gabriel Dospinescu, Guy Henniart, Valentin Hernandez, Florian Herzig, Joaqu\'\i n Rodrigues Jacinto, Stefano Morra, Vincent Pilloni and Benjamin Schraen. We also thank Benchao Su for pointing out an inaccuracy. Z.~Q.~is partially supported by National Key R\&D Program of China 2025YFA 1018000.

\newpage

\section{Preliminaries on smooth representations}\label{sec: smooth rep}

We prove all results on smooth representations of $G$ needed in \S\ref{sec: spectral seq} and especially in \S\ref{sec: extension}. Most results are not really new, but we provide complete proofs. In particular, using Bernstein-Zelevinsky's theory, we define and study the convenient notion of $G$-basic smooth representation of $L_I$ for some $I\subseteq \Delta$ and we give several results on certain $G$-basic representations which are crucially used afterwards.

\subsection{\texorpdfstring{$G$}{G}-basic representations and Bernstein-Zelevinsky's theory}\label{subsec: BZ}

We define $G$-basic representations (Definition \ref{def: basic rep}) and use Bernstein-Zele\-vinsky's geometric lemma and segment theory (\cite{BZ77}, \cite{Z80}) to prove several useful (and presumably well-known) results on them.\bigskip

Let $G$ be a locally compact topological group with left Haar measure $\mu_{G}$, recall its \emph{modulus character} $\delta_G: G \rightarrow \bR_{>0}$ is the unique character satisfying $\mu_{G}(A\cdot x)=\delta_{G}(x)\mu_{G}(A)$ for any Borel subset $A\subseteq G$. When $G$ is moreover $p$-adic analytic, $\delta_G$ is $\bQ^\times$-valued and we will see $\delta_G$ as an $E^\times$-valued character.\bigskip

Let $G$ be a locally profinite group. We write $\mathrm{Rep}^{\infty}(G)$ for the abelian category of all smooth representations of $G$ over $E$ and $\mathrm{Rep}^{\infty}_{\rm{adm}}(G)$ for the full abelian subcategory of admissible ones. Let $H\subseteq G$ be a closed subgroup and $\pi^{\infty}\in \mathrm{Rep}^{\infty}(H)$. We define $(\mathrm{Ind}_{H}^{G}\pi^{\infty})^\infty$ to be the $E$-vector space of uniformly locally constant functions $f: G\rightarrow \pi^{\infty}$ such that $f(xh)=h^{-1}\cdot f(x)$ for $x\in G$ and $h\in H$, which is naturally a (left) smooth $G$-representation via $(g(f))(x)\defeq f(g^{-1}x)$ ($g,x\in G$, $f\in (\mathrm{Ind}_{H}^{G}\pi^{\infty})^\infty$). We also consider the subspace $(\mathrm{ind}_{H}^{G}\pi^{\infty})^\infty\subseteq (\mathrm{Ind}_{H}^{G}\pi^{\infty})^\infty$ consisting of those $f$ for which there exists a compact open subset $C_f$ of $G$ such that $f(x)=0$ for $x\notin C_fH$. They give the so-called (unnormalized) induction and compact induction functors
\[(\mathrm{Ind}_H^G)^\infty,\ (\mathrm{ind}_H^G)^\infty: \mathrm{Rep}^{\infty}(H) \rightarrow \mathrm{Rep}^{\infty}(G)\]
which are both exact. Note that they do not send $\mathrm{Rep}^{\infty}_{\rm{adm}}(H)$ to $\mathrm{Rep}^{\infty}_{\rm{adm}}(G)$ in general.

\begin{rem}
We add $(-)^\infty$ as exponent to avoid possible confusion with locally analytic inductions. Moreover our convention in the definition of $(\mathrm{Ind}_H^G)^\infty$ and $(\mathrm{ind}_H^G)^\infty$ (which is the one used in \cite{ST03}, \cite{ST05}, \cite{OS10}, \cite{OS15}, etc.) is different from the more commonly one used for instance in \cite{Bu90} or \cite[\S III.2.2]{Re10}, where $f\in (\mathrm{Ind}_{H}^{G}\pi^{\infty})^\infty$ satisfies $f(hx)=h\cdot f(x)$ and $g\in G$ acts on $(\mathrm{Ind}_{H}^{G}\pi^{\infty})^\infty$ by $(g(f))(x)\defeq f(xg)$. But the isomorphism $f\mapsto [g\mapsto f(g^{-1})]$ gives an isomorphism between our $(\mathrm{Ind}_{H}^{G}\pi^{\infty})^\infty$ and theirs, and we can freely use the results of \cite{Re10}.
\end{rem}

We start with a general lemma (which will be used in Lemma \ref{lem: coinvariantinside} below).

\begin{lem}\label{lem: induction coinvariant}
Let $P$ be a locally profinite group and $N\subseteq H\subseteq P$ be closed subgroups such that $N$ is \emph{normal} in $P$.
Assume that there exists a continuous section $s:P/H\hookrightarrow P$ of the canonical surjection $P\twoheadrightarrow P/H$ which induces a homeomorphism of locally profinite topological spaces $P/H \times H\buildrel\sim\over\rightarrow P$.
Then for each $\pi^{\infty}$ in $\mathrm{Rep}^{\infty}(H)$, we have a canonical isomorphism in $\mathrm{Rep}^{\infty}(P/N)$
\begin{equation*}
\big((\mathrm{ind}_H^P\pi^{\infty})^\infty\big)_{N}\buildrel\sim\over\longrightarrow \big(\mathrm{ind}_{H/N}^{P/N}(\pi^{\infty})_{N}\big)^\infty
\end{equation*}
where $(-)_N$ means the usual $N$-coinvariants.
\end{lem}
\begin{proof}
As $N$ acts trivially on $\pi^{\infty}_{N}$, we have a natural isomorphism $(\mathrm{ind}_{H/N}^{P/N}(\pi^{\infty})_{N})^\infty\cong (\mathrm{ind}_{H}^{P}(\pi^{\infty})_{N})^\infty$. For a $E$-vector space $M$ equipped with a smooth $N$-action, let $V(M)\subseteq M$ be the subspace spanned by vectors $n\cdot v-v$ for $v\in M$ and $n\in N$. Then $(\pi^{\infty})_{N}=\pi^{\infty}/V(\pi^{\infty})$ by definition, and thus $(\mathrm{ind}_H^P(\pi^{\infty})_{N})^\infty=(\mathrm{ind}_H^P\pi^{\infty})^\infty/(\mathrm{ind}_H^PV(\pi^{\infty}))^\infty$ by the exactness of $(\mathrm{ind}_H^P)^\infty$. As $((\mathrm{ind}_H^P\pi^{\infty})^\infty)_{N}=(\mathrm{ind}_H^P\pi^{\infty})^\infty/V((\mathrm{ind}_H^P\pi^{\infty})^\infty)$ by definition, it suffices to show that we have
\begin{equation}\label{equ: induction coinvariant subspace}
V((\mathrm{ind}_H^P\pi^{\infty})^\infty)=(\mathrm{ind}_H^PV(\pi^{\infty}))^\infty
\end{equation}
as subspaces of $(\mathrm{ind}_H^P\pi^{\infty})^\infty$. Let $f\in (\mathrm{ind}_H^P\pi^{\infty})^\infty$, then for $x\in P$ and $n\in N$, we have $x^{-1}nx\in N\subseteq H$ and
\begin{equation}\label{equ: conjugation N action}
(n\cdot f)(x)=f(n^{-1}x)=f(x(x^{-1}n^{-1}x))=x^{-1}nx\cdot f(x).
\end{equation}
In particular, $(n\cdot f-f)(x)=x^{-1}nx\cdot f(x)-f(x)\in V(\pi^{\infty})$, and thus we have an inclusion $V((\mathrm{ind}_H^P\pi^{\infty})^\infty)\subseteq (\mathrm{ind}_H^PV(\pi^{\infty}))^\infty$. Let us prove that it is a surjection. Write $X\defeq s(P/H)\subseteq P$ and define $C^{\infty}_{c}(X,\pi^{\infty})$ as the $E$-vector space of locally constant function $h: X\rightarrow \pi^{\infty}$ with compact support, and similarly with $C^{\infty}_{c}(X,V(\pi^{\infty}))$. By definition of $(\mathrm{ind}_H^P\pi^{\infty})^\infty$, the map $f\mapsto f|_{X}$ induces an isomorphism of $E$-vector spaces
\begin{equation}\label{equ: pass to X induction}
(\mathrm{ind}_H^P\pi^{\infty})^\infty\buildrel\sim\over\longrightarrow C^{\infty}_{c}(X,\pi^{\infty}),
\end{equation}
which induces an isomorphism of $E$-vector spaces
\begin{equation}\label{equ: pass to X V}
(\mathrm{ind}_H^PV(\pi^{\infty}))^\infty\buildrel\sim\over\longrightarrow C^{\infty}_{c}(X,V(\pi^{\infty})).
\end{equation}
Let $n\in N$, $v\in\pi^{\infty}$ and $x\in X$. As the $N$-action on $\pi^{\infty}$ is smooth and the $P$-action on $N$ by conjugation is continuous, there exists a compact open subset $C_x\subseteq X$ such that $x\in C_x$ and $(y^{-1}ny)\cdot v=(x^{-1}nx)\cdot v$ for each $y\in C_x$. Let $h_{x,v}\in C^{\infty}_{c}(X,\pi^{\infty})$ be the function defined by $h_{x,v}(y)=v$ for $y\in C_x$ and $h_{x,v}(y)=0$ for $y\in X\setminus C_x$, and let $f_{x,v}\in (\mathrm{ind}_H^P\pi^{\infty})^\infty$ correspond to $h_{x,v}$ under (\ref{equ: pass to X induction}). Similarly to (\ref{equ: conjugation N action}), we have $(n\cdot f_{x,v})(y)=(y^{-1}ny)\cdot f_{x,v}(y)=(y^{-1}ny)\cdot h_{x,v}(y)=(y^{-1}ny)\cdot v=(x^{-1}nx)\cdot v$ for $y\in C_x$ and $(n\cdot f_{x,v})(y)=(y^{-1}ny)\cdot f_{x,v}(y)=0$ for each $y\in X\setminus C_x$. Hence, we have $(n\cdot f_{x,v}-f_{x,v})(y)=(x^{-1}nx)\cdot v-v$ for $y\in C_x$ and $(n\cdot f_{x,v}-f_{x,v})(y)=0$ for $y\in X\setminus C_x$. As $n\in N$, $v\in\pi^{\infty}$, $x\in X$ are arbitrary and as $C_x$ can be an arbitrarily small neighborhood of $x$ in $X$, elements of the form $(n\cdot f_{x,v}-f_{x,v})|_X$ span $C^{\infty}_{c}(X,V(\pi^{\infty}))$ over $E$, and thus elements of the form $n\cdot f_{x,v}-f_{x,v}$ span $(\mathrm{ind}_H^PV(\pi^{\infty}))^\infty$ over $E$ by (\ref{equ: pass to X V}). Since we clearly have $n\cdot f_{x,v}-f_{x,v}\in V((\mathrm{ind}_H^P\pi^{\infty})^\infty)$, we see that (\ref{equ: induction coinvariant subspace}) holds.
\end{proof}

Let $G$ be the $K$-points of a $p$-adic reductive algebraic group over $K$. For (the $K$-points of) a parabolic subgroup $P=L_PN_P\subseteq G$, the unnormalized parabolic induction functor $(\mathrm{Ind}_{P}^{G})^\infty: \mathrm{Rep}^{\infty}(L_P)\rightarrow \mathrm{Rep}^{\infty}(G)$ restricts to a functor $(\mathrm{Ind}_{P}^{G})^\infty: \mathrm{Rep}^{\infty}_{\rm{adm}}(L_P)\rightarrow \mathrm{Rep}^{\infty}_{\rm{adm}}(G)$ (see \cite[\S III.2.3]{Re10}). The functor $(\mathrm{Ind}_{P}^{G})^\infty$ admits a left adjoint functor $J_{N_P}: \mathrm{Rep}^{\infty}(G)\rightarrow \mathrm{Rep}^{\infty}(L_P)$ which is exact and restricts to a functor $J_{N_P}: \mathrm{Rep}^{\infty}_{\rm{adm}}(G)\rightarrow \mathrm{Rep}^{\infty}_{\rm{adm}}(L_P)$ (cf.\cite[\S\S VI.1.1, VI.6.1]{Re10} but note that our $J_{N_P}$ is the unnormalized Jacquet functor). Arguing as in the beginning of \cite[\S VI.9.6]{Re10}, the functor $(\mathrm{Ind}_{P}^{G})^\infty$ also admits a right adjoint $J_{N_P}': \mathrm{Rep}^{\infty}(G)\rightarrow \mathrm{Rep}^{\infty}(L_P)$ which is also exact and restricts to $J_{N_P}': \mathrm{Rep}^{\infty}_{\rm{adm}}(G)\rightarrow \mathrm{Rep}^{\infty}_{\rm{adm}}(L_P)$ (cf.~\cite[\S VI.9.6]{Re10} taking care again that we are unnormalized).\bigskip

We also need the \emph{normalized parabolic induction functor} $i_{P}^{G}: \mathrm{Rep}^{\infty}(L_P)\rightarrow \mathrm{Rep}^{\infty}(G)$ defined by $i_{P}^{G}(-)\defeq (\mathrm{Ind}_{P}^{G}((-)\otimes_E \delta_{P}^{-1/2}))^\infty$, and the \emph{normalized parabolic restriction functor} $r_{P}^{G}: \mathrm{Rep}^{\infty}(G)\rightarrow \mathrm{Rep}^{\infty}(L_P)$ defined by $r_{P}^{G}(-)\defeq J_{N_P}(-)\otimes_E \delta_{P}^{1/2}$ (cf.~\cite[\S VI.1.2]{Re10} and recall that $\delta_{P}: P\rightarrow E^\times$ factors through $L_P$). Note that, here, we might need to extend scalars to $E'=E(\sqrt{q})$ so that $\delta_{P}^{1/2}$ is $E'$-valued. But these normalized functors will only play a minor intermediate role in that paper, and everything is ultimately $K$-rational (see for instance Remark \ref{rem: rational dot} and Remark \ref{rem: rational geometric lemma} below).\bigskip

For parabolic subgroups $Q=L_QN_Q\subseteq P=L_PN_P$, $Q\cap L_P$ is a parabolic subgroup of $L_p$ with reductive quotient $L_Q$, and we have the formula
\begin{equation}\label{equ: modulus product}
\delta_{Q}=\delta_{Q\cap L_P}\cdot (\delta_{P}|_{L_Q}): L_Q \rightarrow E^\times.
\end{equation}

From now on $G=\mathrm{GL}_n(K)$ as in \S\ref{generalnotation}. For $I\subseteq I_0\subseteq \Delta$, we use the shortened notation
\begin{equation}\label{extranot}
i_{I,I_0}^{\infty}\defeq \big(\mathrm{Ind}_{P_I \cap L_{I_0}}^{L_{I_0}}\big)^\infty,\ J_{I_0,I}\defeq J_{N_I\cap L_{I_0}},\ J_{I_0,I}'\defeq J_{N_I\cap L_{I_0}}'.
\end{equation}
Then for each $I'\subseteq I$, we clearly have
\begin{equation*}
i_{I',I_0}^{\infty}\cong i_{I,I_0}^{\infty}(i_{I',I}^{\infty}),\ J_{I_0,I'}\cong J_{I,I'}(J_{I_0,I}),\ J_{I_0,I'}'\cong J_{I,I'}'(J_{I_0,I}').
\end{equation*}
For $\pi^{\infty}$ in $\mathrm{Rep}^{\infty}_{\rm{adm}}(L_I)$ and $\pi_0^{\infty}$ in $\mathrm{Rep}^{\infty}_{\rm{adm}}(L_{I_0})$, we have canonical isomorphisms for $k\geq 0$:
\begin{equation}\label{equ: first adjunction}
\mathrm{Ext}_{L_{I_0}}^k(\pi_0^{\infty}, i_{I,I_0}^{\infty}(\pi^{\infty}))^{\infty}\cong \mathrm{Ext}_{L_{I}}^k(J_{I_0,I}(\pi_0^{\infty}), \pi^{\infty})^{\infty}
\end{equation}
and
\begin{equation}\label{equ: second adjunction}
\mathrm{Ext}_{L_{I_0}}^k(i_{I,I_0}^{\infty}(\pi^{\infty}), \pi_0^{\infty})^{\infty}\cong \mathrm{Ext}_{L_{I}}^k(\pi^{\infty}, J_{I_0,I}'(\pi_0^{\infty}))^{\infty}.
\end{equation}
Recall that, for $I,I'\subseteq I_0$, $W^{I',I}(L_{I_0})\subseteq W(L_{I_0})$ is the set of minimal length representatives of the double coset $W(L_{I'})\backslash W(L_{I_0})/W(L_{I})$. Let $w_{I_0,I}$ be the \emph{longest} element inside $W^{\emptyset,I}(L_{I_0})$, we have $w_{I_0,I}(I)\subseteq I_0$~and
\begin{equation}\label{interI_0}
w_{I_0,I}^{-1}(\Phi_{I_0}^+)\cap \Phi_{I_0}^+ = \Phi_{I}^+.
\end{equation}
Note that $P_I^+\cap L_{I_0}$ is the parabolic of $L_{I_0}$ opposite to $P_I\cap L_{I_0}$. It thus follows from \cite[(VI.9.6.1)]{Re10}, $w_{I_0,I}^{-1}L_{w_{I_0,I}(I)}w_{I_0,I}=L_I$ and $\delta_{P_I^+\cap L_{I_0}}=\delta_{P_I\cap L_{I_0}}^{-1}$ that
\begin{eqnarray*}
\Hom_{L_{I_0}}(i_{I,I_0}^{\infty}(\pi^{\infty}), \pi_0^{\infty})& = & \Hom_{L_{I_0}}(i_{P_I\cap L_{I_0}}^{L_{I_0}}(\pi^{\infty}\otimes_E \delta_{P_I\cap L_{I_0}}^{1/2}), \pi_0^{\infty})\\
&\cong &\Hom_{L_I}(\pi^{\infty}\otimes_E \delta_{P_I\cap L_{I_0}}^{1/2}, r_{P_I^+\cap L_{I_0}}^{L_{I_0}}(\pi_0^{\infty}))\\
&\cong &\Hom_{L_I}\Big(\pi^{\infty}\otimes_E \delta_{P_I\cap L_{I_0}}^{1/2}, J_{N_I^+\cap L_{I_0}}(\pi_0^{\infty})\otimes_E\delta_{P_I^+\cap L_{I_0}}^{1/2}\Big)\\
&\cong &\Hom_{L_I}\Big(\pi^{\infty}, (J_{I_0,w_{I_0,I}(I)}(\pi_0^{\infty}))^{w_{I_0,I}}\otimes_E\delta_{P_I^+\cap L_{I_0}}^{1/2}\otimes_E \delta_{P_I\cap L_{I_0}}^{-1/2}\Big)\\
&\cong &\Hom_{L_I}\Big(\pi^{\infty}, (J_{I_0,w_{I_0,I}(I)}(\pi_0^{\infty}))^{w_{I_0,I}}\otimes_E \delta_{P_I\cap L_{I_0}}^{-1}\Big).
\end{eqnarray*}
By (\ref{equ: second adjunction}) for $k=0$ and since this holds for arbitrary $\pi^{\infty}$ in $\mathrm{Rep}^{\infty}_{\rm{adm}}(L_I)$ we deduce
\begin{equation}\label{equ: explicit twist Jacquet}
J_{I_0,I}'(\pi_0^{\infty})\cong (J_{I_0,w_{I_0,I}(I)}(\pi_0^{\infty}))^{w_{I_0,I}}\otimes_E \delta_{P_I\cap L_{I_0}}^{-1}.
\end{equation}
We write $\widehat{T}^{\infty}$ for the set of smooth $E$-valued characters of $T$, which is naturally an abelian group under multiplication. We let $W(G)$ act on the left on $\widehat{T}^{\infty}$ via the following dot action
\begin{equation}\label{equ: smooth dot action}
w\cdot \chi\defeq w(\chi\otimes_E\delta_B^{1/2})\otimes_E \delta_B^{-1/2}=(\chi\otimes_E\delta_B^{1/2})(w^{-1} \cdot w)\otimes_E \delta_B^{-1/2}
\end{equation}
and recall that $\delta_{B}$ sends $(t_1,\dots,t_n)\in T$ to $\prod_{i=1}^n|t_i|_K^{(n+1)-2i}$, which can be rewritten as $|\sum_{\alpha \in \Phi^+} \alpha |_K=|-2 \rho |_K$. If $w\in W(L_I)$ for some $I\subseteq \Delta$ one can check that $w\cdot \chi=w(\chi\otimes_E\delta_{B_I}^{1/2})\otimes_E \delta_{B_I}^{-1/2}$ (recall $B_I=B \cap L_I$).

\begin{rem}\label{rem: rational dot}
Note that the character $\delta_{B}$ is $K^\times$-valued but $\delta_{B}^{1/2}$ is not in general. Nevertheless, $w(\delta_{B}^{1/2})\otimes_E\delta_{B}^{-1/2}$ is $K^\times$-valued for any $w\in W(G)$. Similarly, $w(\delta_{B_I}^{1/2})\otimes_E\delta_{B_I}^{-1/2}$ is $K^\times$-valued for $I\subseteq \Delta$ and $w\in W(L_I)$. Consequently, $\chi\in\widehat{T}^{\infty}$ is $K^\times$-valued if and only if $w\cdot\chi=w(\chi)\otimes_E(w(\delta_{B}^{1/2})\otimes_E\delta_{B}^{-1/2})$ is $K^\times$-valued for each $w\in W(G)$.
\end{rem}

\begin{defn}\label{def: basic rep}
Let $I\subseteq \Delta$, $\pi^{\infty}$ a finite length representation in $\mathrm{Rep}^{\infty}_{\rm{adm}}(L_I)$ and recall that $J_{I,\emptyset}(\pi^{\infty})$ has finite length (cf.~\cite[\S VI.6.4]{Re10}). We write $\cJ(\pi^{\infty})\subseteq \widehat{T}^{\infty}$ for the subset of $\chi$ such that
\begin{equation*}
\Hom_{L_I}(\pi^{\infty}, i_{\emptyset,I}^{\infty}(\chi))\cong\Hom_{T}(J_{I,\emptyset}(\pi^{\infty}),\chi)\neq 0.
\end{equation*}
\begin{enumerate}[label=(\roman*)]
\item \label{it: basic 0} A character $\chi\in \widehat{T}^{\infty}$ is called \emph{$G$-regular} if $w\cdot \chi\neq \chi$ for $1\neq w\in W(G)$.
\item \label{it: basic 1} If $\pi^{\infty}$ is irreducible, it is called \emph{$G$-regular} if $\cJ(\pi^{\infty})$ contains a $G$-regular element. In general, $\pi^{\infty}$ is called \emph{$G$-basic} if there exists $I_1\subseteq I$ and an irreducible $G$-regular $\pi_1^{\infty}$ in $\mathrm{Rep}^{\infty}_{\rm{adm}}(L_{I_1})$ such that $\pi^{\infty}\cong i_{I_1,I}^{\infty}(\pi_1^{\infty})$.
\item \label{it: basic 3} For $i=0,1$ let $I_i\subseteq I\subseteq\Delta$ and $\pi_i^{\infty}$ in $\mathrm{Rep}^{\infty}_{\rm{adm}}(L_{I_i})$, we define the \emph{$L_I$-distance from $\pi_0^{\infty}$ to $\pi_1^{\infty}$} as
\[d_I(\pi_0^{\infty},\pi_1^{\infty})\defeq \inf\{k\mid\mathrm{Ext}_{L_I}^k(i_{I_0,I}^{\infty}(\pi_0^{\infty}),i_{I_1,I}^{\infty}(\pi_1^{\infty}))^{\infty}\neq 0\}\leq \infty.\]
We write $d(\pi_0^{\infty},\pi_1^{\infty})\defeq d_\Delta(\pi_0^{\infty},\pi_1^{\infty})$.
\end{enumerate}
\end{defn}

Note that the trivial character $1_T\in \widehat{T}^{\infty}$ is $G$-regular. By (\ref{equ: first adjunction}) for $k=0$ applied to $\pi_1^\infty$ as in \ref{it: basic 1} of Definition \ref{def: basic rep} and the left-exactness of induction we see that $\pi^{\infty}\cong i_{I_1,I}^{\infty}(\pi_1^{\infty})$ embeds into $i_{\emptyset,I}^{\infty}(\chi)\cong i_{I_1,I}^{\infty}(i_{\emptyset,I_1}^{\infty}(\chi))$ for $\chi\in \cJ(\pi_1^{\infty})$, and hence $\chi\in \cJ(\pi^{\infty})$. In particular we have $\cJ(\pi_1^{\infty})\subseteq \cJ(\pi^{\infty})$.
Note also that any $G$-basic representation of $L_I$ is (admissible) of finite length (as so is $i_{I_1,I}^{\infty}(\pi_1^{\infty})$, see \cite[\S VI.6.2]{Re10}) and that $i_{I,I'}^{\infty}(\pi^{\infty})\in\mathrm{Rep}^{\infty}_{\rm{adm}}(L_{I'})$ is again $G$-basic for any $I\subseteq I'\subseteq \Delta$.

\begin{rem}
In Definition~\ref{def: basic rep}, as $J_{I,\emptyset}(\pi^{\infty})$ is a finite length smooth representation of $T$, we have $J_{I,\emptyset}(\pi^{\infty})\cong \bigoplus_{\chi\in\cJ(\pi^{\infty})}J_{I,\emptyset}(\pi^{\infty})_{\chi}$ with $J_{I,\emptyset}(\pi^{\infty})_{\chi}$ having only $\chi$ as Jordan-H\"older factor. In particular, the set $\cJ(\pi^{\infty})$ equals $\mathrm{JH}_{T}(J_{I,\emptyset}(\pi^{\infty}))$ as well as the set of $\chi$ such that
\[\Hom_{T}(\chi,J_{I,\emptyset}(\pi^{\infty}))\neq 0.\]
\end{rem}

\begin{rem}\label{rem: twist Jacquet support}
For $I\subseteq I_0\subseteq \Delta$, by an easy calculation using (\ref{interI_0}) one can check
\[\delta_{P_I\cap L_{I_0}}^{-1}|_{T}=\delta_{w_{I_0,I}^{-1}B_{I_0}w_{I_0,I}}^{1/2}\otimes_E\delta_{B_{I_0}}^{-1/2}=\delta_{B_{I_0}}^{1/2}(w_{I_0,I}\cdot w_{I_0,I}^{-1})\otimes_E\delta_{B_{I_0}}^{-1/2}.\]
Using this and (\ref{equ: explicit twist Jacquet}) we deduce isomorphisms of $T$-representations for $\pi_0^{\infty}$ in $\mathrm{Rep}^{\infty}_{\rm{adm}}(L_{I_0})$
\begin{eqnarray*}
J_{I,\emptyset}(J_{I_0,I}'(\pi_0^{\infty})) &\cong &(J_{w_{I_0,I}(I),\emptyset}(J_{I_0,w_{I_0,I}(I)}(\pi_0^{\infty})))^{w_{I_0,I}}\otimes_E \delta_{P_I\cap L_{I_0}}^{-1}|_{T}\\
&\cong &J_{I_0,\emptyset}(\pi_0^{\infty})^{w_{I_0,I}}\otimes_E \delta_{P_I\cap L_{I_0}}^{-1}|_{T}\\
&\cong &(J_{I_0,\emptyset}(\pi_0^{\infty})\otimes_E \delta_{B_{I_0}}^{1/2})^{w_{I_0,I}}\otimes_E \delta_{B_{I_0}}^{-1/2}.
\end{eqnarray*}
From $w_{I_0,I}^{-1}\cdot \chi=(\chi\otimes_E\delta_{B_{I_0}}^{1/2})(w_{I_0,I}\cdot w_{I_0,I}^{-1})\otimes_E \delta_{B_{I_0}}^{-1/2}$ we finally obtain
\begin{equation}\label{equ: twist Jacquet support}
\cJ(J_{I_0,I}'(\pi_0^{\infty}))=w_{I_0,I}^{-1}\cdot\cJ(\pi_0^{\infty})
\end{equation}
where the set $\cJ(-)$ is as in Definition~\ref{def: basic rep}.
\end{rem}

Let $I\subseteq \Delta$. We refer to \cite[\S VI.5.2]{Re10} for the definition of equivalence classes of cuspidal data (for the group $L_I$). We consider here equivalence classes of cuspidal data $(T, \chi\otimes_E \delta_{B_I}^{1/2})$ for $\chi\in\widehat{T}^{\infty}$. Recall that it consists of the cuspidal data $(gTg^{-1},(\chi\otimes_E \delta_{B_I}^{1/2})(g^{-1} \cdot g)$ for $g\in L_I$. In particular for $\chi'\in W(L_I)\cdot\chi$ all cuspidal data $(T,\chi'\otimes_E \delta_{B_I}^{1/2})$ are equivalent. To a left $W(L_I)$-coset $\Sigma$ (under the dot action (\ref{equ: smooth dot action})) we can thus associate an element in the set $\Omega(L_I)$ of all equivalence classes of cuspidal data.\bigskip

For each finite subset $\Sigma\subseteq \widehat{T}^{\infty}$ which is stable under the left dot action of $W(L_I)$, we write $\cB^I_{\Sigma}$ for the category of finite length representations $\pi^{\infty}$ in $\mathrm{Rep}^{\infty}_{\rm{adm}}(L_I)$ with each $\sigma^{\infty}\in\mathrm{JH}_{L_I}(\pi^{\infty})$ satisfying $\emptyset\neq \cJ(\sigma^{\infty})\subseteq \Sigma$. By \cite[\S VI.7.2]{Re10} and Remark \ref{rem: rational block} below, we know that $\cB^I_{\Sigma}$ is the direct sum of the $\cB^I_{\Sigma'}$ for $\Sigma'$ running through the $W(L_I)$-cosets contained in $\Sigma$. In particular, for each object $\pi^{\infty}$ in $\cB^I_{\Sigma}$, we have a canonical decomposition
\begin{equation}\label{equ: sm block decomposition}
\pi^{\infty}\cong \bigoplus_{\Sigma'}\pi^{\infty}_{\cB^{I}_{\Sigma'}}
\end{equation}
with $\pi^{\infty}_{\cB^{I}_{\Sigma'}}\in \cB^I_{\Sigma'}$ and $\Sigma'$ running through $W(L_I)$-cosets contained in $\Sigma$. The exactness of $J_{I,\emptyset}$ implies that $\cJ(\pi^{\infty}_{\cB^{I}_{\Sigma'}})=\cJ(\pi^{\infty})\cap \Sigma'$.
We say that $\Sigma$ is \emph{$G$-regular} if each $\chi\in \Sigma$ is $G$-regular.

\begin{rem}\label{rem: rational block}
The coefficient field in \cite[\S VI.7.2]{Re10} being algebraically closed, the decomposition (\ref{equ: sm block decomposition}) deserves some justification. Let $\overline{E}$ an algebraic closure of $E$ and $\mathrm{Rep}^{\infty}_{\rm{adm}}(L_I,\overline{E})$ the abelian category of smooth admissible representations of $L_I$ over $\overline{E}$-vector spaces. Let $I\subseteq \Delta$ and $\Sigma\subseteq \widehat{T}^{\infty}$ a finite subset which is stable under the left dot action of $W(L_I)$ (and that we also see inside the smooth $\overline{E}$-valued characters of $T$). We write $\cB^I_{\Sigma,\overline{E}}$ for the category of finite length representations $\pi^{\infty}_{\overline{E}}$ in $\mathrm{Rep}^{\infty}_{\rm{adm}}(L_I,\overline{E})$ with each $\sigma^{\infty}_{\overline{E}}\in\mathrm{JH}_{L_I}(\pi^{\infty}_{\overline{E}})$ satisfying $\emptyset\neq \cJ(\sigma^{\infty}_{\overline{E}})\subseteq \Sigma$. Now let $\pi^{\infty}$ in $\cB^I_{\Sigma}\subseteq \mathrm{Rep}^{\infty}_{\rm{adm}}(L_I)$ and $\pi^{\infty}_{\overline{E}}\defeq \pi^{\infty}\otimes_E\overline{E}$. By \cite[\S VI.7.2]{Re10} we have a decomposition
\begin{equation}\label{equ: sm block decomposition bar}
\pi^{\infty}_{\overline{E}}\cong \bigoplus_{\Sigma'}(\pi^{\infty}_{\overline{E}})_{\cB^{I}_{\Sigma',\overline{E}}}
\end{equation}
where $(\pi^{\infty}_{\overline{E}})_{\cB^{I}_{\Sigma',\overline{E}}}\in \cB^I_{\Sigma',\overline{E}}$ and $\Sigma'$ runs through the $W(L_I)$-cosets contained in $\Sigma$. There exists an obvious $E$-linear action of $\mathrm{Gal}(\overline{E}/E)$ on $\pi^{\infty}_{\overline{E}}$ such that $(\pi^{\infty}_{\overline{E}})^{\mathrm{Gal}(\overline{E}/E)}\cong \pi^{\infty}$, which induces by (\ref{equ: sm block decomposition bar}) an $E$-linear action of $\mathrm{Gal}(\overline{E}/E)$ on $\oplus_{\Sigma'}(\pi^{\infty}_{\overline{E}})_{\cB^{I}_{\Sigma',\overline{E}}}$, and also by scalar extension from $E$ to $\overline{E}$ an $\overline{E}$-linear action of $\mathrm{Gal}(\overline{E}/E)$ on
\[(\pi^{\infty}_{\overline{E}})\otimes_E\overline{E}\cong \bigoplus_{\Sigma'}\big((\pi^{\infty}_{\overline{E}})_{\cB^{I}_{\Sigma',\overline{E}}}\big)\otimes_E\overline{E}.\]
But the smooth reducible representations $\big((\pi^{\infty}_{\overline{E}})_{\cB^{I}_{\Sigma',\overline{E}}}\big)\otimes_E\overline{E}$ do not share any irreducible constituent over $\overline{E}$ for \emph{distinct} $\Sigma'$ since the $\Sigma'$ consist of $E$-valued characters, hence the action of $\mathrm{Gal}(\overline{E}/E)$ on $\pi^{\infty}_{\overline{E}}$ must stabilize each $(\pi^{\infty}_{\overline{E}})_{\cB^{I}_{\Sigma',\overline{E}}}$. Let $\pi^{\infty}_{\cB^{I}_{\Sigma'}}\defeq ((\pi^{\infty}_{\overline{E}})_{\cB^{I}_{\Sigma',\overline{E}}})^{\mathrm{Gal}(\overline{E}/E)}$, we deduce from (\ref{equ: sm block decomposition bar})
\[\pi^{\infty}\cong \bigoplus_{\Sigma'}\pi^{\infty}_{\cB^{I}_{\Sigma'}}\]
and $\pi^{\infty}_{\cB^{I}_{\Sigma'}}\otimes_E\overline{E}\cong (\pi^{\infty}_{\overline{E}})_{\cB^{I}_{\Sigma',\overline{E}}}$ in $\mathrm{Rep}^{\infty}_{\rm{adm}}(L_I,\overline{E})$ (extending scalars to $\overline{E}$ again). This is the decomposition (\ref{equ: sm block decomposition}).
\end{rem}

Let $I_0,I_1\subseteq I\subseteq \Delta$ and recall $W^{I_0,I_1}(L_I)\subseteq W(L_I)$ is the subset of minimal length representatives of $W(L_{I_0})\backslash W(L_I)/W(L_{I_1})$. Recall also the Bruhat decomposition (\cite[Lemma 5.5]{DM91})
\[\bigsqcup_{w\in W^{I_0,I_1}(L_I)}(P_{I_0}\cap L_I)w(P_{I_1}\cap L_I)=L_I=\bigsqcup_{w\in W^{I_0,I_1}(L_I)}(P_{I_1}\cap L_I)w^{-1}(P_{I_0}\cap L_I).\]
For $w\in W^{I_0,I_1}(L_I)$ we have
\begin{equation}\label{caca}
\begin{array}{ccccccc}
B_{I_0}&=&B\cap L_{I_0}&=&wBw^{-1}\cap L_{I_0}&\subseteq & wP_{I_1}w^{-1}\cap L_{I_0}\\
B_{I_1}&=&B\cap L_{I_1}&=&w^{-1}Bw\cap L_{I_1}&\subseteq & w^{-1}P_{I_0}w\cap L_{I_1}.
\end{array}
\end{equation}
Hence $wP_{I_1}w^{-1}\cap L_{I_0}$ is a parabolic subgroup of $L_{I_0}$ with Levi quotient $wL_{I_1}w^{-1}\cap L_{I_0}=L_{w(I_1)\cap I_0}$, and $w^{-1}P_{I_0}w\cap L_{I_1}$ is a parabolic subgroup of $L_{I_1}$ with Levi quotient $w^{-1}L_{I_0}w\cap L_{I_1}=L_{w^{-1}(I_0)\cap I_1}$. By (\ref{caca}) and the equality $L_{w(I_1)\cap I_0}=wL_{w^{-1}(I_0)\cap I_1}w^{-1}$ we have
\begin{multline}\label{equ: Borel conj}
B_{w(I_1)\cap I_0}=B_{I_0}\cap L_{w(I_1)\cap I_0}=wBw^{-1}\cap L_{I_0}\cap L_{w(I_1)\cap I_0}=w(B\cap L_{w^{-1}(I_0)\cap I_1})w^{-1}\\=wB_{w^{-1}(I_0)\cap I_1}w^{-1}.
\end{multline}
As $w(I_1)\cap I_0\subseteq I_0$, we have $P_{w(I_1)\cap I_0}\subseteq P_{I_0}$ and thus we have the decomposition
\begin{equation*}
P_{w(I_1)\cap I_0}\cap L_I=(N_{I_0}\cap L_I)\cdot(P_{w(I_1)\cap I_0}\cap L_{I_0})=(N_{I_0}\cap L_I)\cdot(wP_{I_1}w^{-1}\cap L_{I_0}),
\end{equation*}
which together with (\ref{equ: modulus product}) (taking $P=P_{I_0}\cap L_I$ and $Q=P_{w(I_1)\cap I_0}\cap L_I\subseteq P$) implies
\begin{equation}\label{equ: outer modulus}
\delta_{wP_{I_1}w^{-1}\cap L_{I_0}} \cdot (\delta_{P_{I_0}\cap L_I}|_{wL_{I_1}w^{-1}\cap L_{I_0}}) =\delta_{P_{w(I_1)\cap I_0}\cap L_I}.
\end{equation}
Likewise, as $w^{-1}(I_0)\cap I_1\subseteq I_1$, we have $P_{w^{-1}(I_0)\cap I_1}\subseteq P_{I_1}$ and thus the decomposition
\begin{equation*}
P_{w^{-1}(I_0)\cap I_1}\cap L_I=(N_{I_1}\cap L_I)\cdot(P_{w^{-1}(I_0)\cap I_1}\cap L_{I_1})=(N_{I_1}\cap L_I)\cdot(w^{-1}P_{I_0}w\cap L_{I_1}),
\end{equation*}
which together with (\ref{equ: modulus product}) (taking $P=P_{I_1}\cap L_I$ and $Q=P_{w^{-1}(I_0)\cap I_1}\cap L_I \subseteq P$) implies
\begin{equation}\label{equ: inner modulus}
\delta_{w^{-1}P_{I_0}w\cap L_{I_1}}\cdot (\delta_{P_{I_1}\cap L_I}|_{L_{I_1}\cap w^{-1}L_{I_0}w})=\delta_{P_{w^{-1}(I_0)\cap I_1}\cap L_I}.
\end{equation}

We now consider the two functors (see (\ref{extranot}) for the notation)
\begin{equation}\label{equ: Bruhat Jacquet}
J_{I_0,I_1,w}(-)\defeq (J_{I_0,w(I_1)\cap I_0}(-))^w\otimes_E \delta_{I_0,I_1,w}: \mathrm{Rep}^{\infty}_{\rm{adm}}(L_{I_0})\longrightarrow \mathrm{Rep}^{\infty}_{\rm{adm}}(L_{w^{-1}(I_0)\cap I_1})
\end{equation}
where
\begin{equation}\label{equ: Jacquet twist}
\delta_{I_0,I_1,w}\defeq \delta_{w^{-1}P_{w(I_1)\cap I_0}w\cap L_I}^{1/2} \cdot \delta_{P_{w^{-1}(I_0)\cap I_1}\cap L_I}^{-1/2},
\end{equation}
(note that $\delta_{I_0,I_1,1}=1$) and
\begin{equation}\label{equ: Bruhat induction}
i_{I_0,I_1,w}^{\infty}(-)\defeq i_{w^{-1}(I_0)\cap I_1,I_1}^{\infty}(-): \mathrm{Rep}^{\infty}_{\rm{adm}}(L_{w^{-1}(I_0)\cap I_1})\longrightarrow \mathrm{Rep}^{\infty}_{\rm{adm}}(L_{I_1}).
\end{equation}
It follows from (\ref{equ: Borel conj}) (and (\ref{equ: modulus product})) that
\begin{eqnarray}\label{equ: twist Borel}
\nonumber \delta_{I_0,I_1,w}|_{T}&=&\delta_{w^{-1}P_{w(I_1)\cap I_0}w\cap L_I}^{1/2}|_{T}\cdot \delta_{w^{-1}B_{w(I_1)\cap I_0}w}^{1/2}\cdot \delta_{B_{w^{-1}(I_0)\cap I_1}}^{-1/2}\cdot \delta_{P_{w^{-1}(I_0)\cap I_1}\cap L_I}^{-1/2}|_{T}\\
&=&\delta_{w^{-1}B_I w}^{1/2} \otimes_E \delta_{B_I}^{-1/2}.
\end{eqnarray}
In particular $\delta_{I_0,I_1,w}$ is $K^\times$-valued (as in Remark \ref{rem: rational dot}). Moreover if $I_0=\emptyset$ (so that $L_{I_0}=L_{w^{-1}(I_0)\cap I_1}=T$), we have for $\chi\in\widehat{T}^{\infty}$ and $w\in W^{\emptyset,I_1}(L_I)$
\begin{equation}\label{equ: char case}
J_{\emptyset,I_1,w}(\chi)=(w^{-1}(\chi)\otimes_E\delta_{w^{-1}B_Iw}^{1/2})\otimes_E \delta_{B_I}^{-1/2}=w^{-1}\cdot\chi.
\end{equation}

\begin{lem}\label{lem: coinvariantinside}
Let $I_0,I_1\subseteq I\subseteq \Delta$, $w\in W^{I_0,I_1}(L_I)$ and $\pi_0^{\infty}$ in $\mathrm{Rep}^{\infty}_{\rm{adm}}(L_{I_0})$. We have an isomorphism in $\mathrm{Rep}^{\infty}_{\rm{adm}}(L_{I_1}\cap w^{-1}L_{I_0}w)$
\begin{equation}\label{equ: N coinvariant}
J_{I_0,I_1,w}(\pi_0^{\infty})\cong \big((\mathrm{ind}_{N_{I_1}\cap w^{-1}P_{I_0}w}^{N_{I_1}}\pi_0^{\infty,w})^{\infty}\big)_{N_{I_1}}
\end{equation}
and an isomorphism in $\mathrm{Rep}^{\infty}_{\rm{adm}}(L_{I_1})$
\begin{equation}\label{equ: N coinvariant induction}
i_{I_0,I_1,w}^{\infty}(J_{I_0,I_1,w}(\pi_0^{\infty}))\cong \big((\mathrm{ind}_{P_{I_1}\cap w^{-1}P_{I_0}w}^{P_{I_1}}\pi_0^{\infty,w})^\infty\big)_{N_{I_1}}
\end{equation}
(where we view $\pi_0^{\infty,w}\in \mathrm{Rep}^{\infty}_{\rm{adm}}(w^{-1}L_{I_0}w)\subseteq \mathrm{Rep}^{\infty}_{\rm{adm}}(w^{-1}P_{I_0}w)$ as a smooth representation of $P_{I_1}\cap w^{-1}P_{I_0}w$ by restriction).
\end{lem}
\begin{proof}
The isomorphism (\ref{equ: N coinvariant}) easily follows from \cite[(VI.5.1.3)]{Re10} (applied with $J=N=N_{I_1}$ and $H=N_{I_1}\cap w^{-1}P_{I_0}w$) and from $(J_{I_0,w(I_1)\cap I_0}(\pi_0^\infty))^w \cong J_{N_{I_1 \cap w^{-1}(I_0)} \cap L_{w^{-1}(I_0)}}(\pi_0^{\infty,w})$, noting that the action of $N_{I_1}\cap w^{-1}P_{I_0}w$ on $\pi_0^{\infty,w}$ factors through $N_{I_1} \cap w^{-1}L_{I_0}w=N_{I_1 \cap w^{-1}(I_0)} \cap L_{w^{-1}(I_0)}$. Note that $N_{I_1}$ is a normal subgroup of $P_{I_1}$ and thus $N_{I_1}(P_{I_1}\cap w^{-1}P_{I_0}w)$ is a subgroup of $P_{I_1}$ satisfying
\begin{equation}\label{equ: change gp 1}
N_{I_1}(P_{I_1}\cap w^{-1}P_{I_0}w)/(P_{I_1}\cap w^{-1}P_{I_0}w)\cong N_{I_1}/(N_{I_1}\cap w^{-1}P_{I_0}w)
\end{equation}
and
\begin{equation}\label{equ: change gp 2}
P_{I_1}/N_{I_1}(P_{I_1}\cap w^{-1}P_{I_0}w) \cong L_{I_1}/(L_{I_1}\cap w^{-1}P_{I_0}w).
\end{equation}
The isomorphism (\ref{equ: N coinvariant induction}) then follows from
\begin{eqnarray*}
\big((\mathrm{ind}_{P_{I_1}\cap w^{-1}P_{I_0}w}^{P_{I_1}}\pi_0^{\infty,w})^{\infty}\big)_{N_{I_1}}&\cong &\big((\mathrm{ind}_{N_{I_1}(P_{I_1}\cap w^{-1}P_{I_0}w)}^{P_{I_1}}(\mathrm{ind}_{P_{I_1}\cap w^{-1}P_{I_0}w}^{N_{I_1}(P_{I_1}\cap w^{-1}P_{I_0}w)}\pi_0^{\infty,w})^{\infty})^{\infty}\big)_{N_{I_1}}\\
&\cong &\big(\mathrm{ind}_{N_{I_1}(P_{I_1}\cap w^{-1}P_{I_0}w)}^{P_{I_1}}((\mathrm{ind}_{P_{I_1}\cap w^{-1}P_{I_0}w}^{N_{I_1}(P_{I_1}\cap w^{-1}P_{I_0}w)}\pi_0^{\infty,w})^{\infty})_{N_{I_1}}\big)^{\infty}\\
&\cong &\big(\mathrm{Ind}_{L_{I_1}\cap w^{-1}P_{I_0}w}^{L_{I_1}}((\mathrm{ind}_{N_{I_1}\cap w^{-1}P_{I_0}w}^{N_{I_1}}\pi_0^{\infty,w})^{\infty})_{N_{I_1}}\big)^{\infty}\\
&\cong &i_{I_0,I_1,w}^{\infty}(J_{I_0,I_1,w}(\pi_0^{\infty}))
\end{eqnarray*}
where the second isomorphism follows from Lemma~\ref{lem: induction coinvariant} applied with $P=P_{I_1}$, $H=N_{I_1}(P_{I_1}\cap w^{-1}P_{I_0}w)$, $N=N_{I_1}$ and $\pi^{\infty}=(\mathrm{ind}_{P_{I_1}\cap w^{-1}P_{I_0}w}^{N_{I_1}(P_{I_1}\cap w^{-1}P_{I_0}w)}\pi_0^{\infty,w})^{\infty}$ (using that $P_{I_1}\twoheadrightarrow P_{I_1}/N_{I_1}(P_{I_1}\cap w^{-1}P_{I_0}w)$ admits a continuous section), the third isomorphism follows from (\ref{equ: change gp 1}) and (\ref{equ: change gp 2}), and the last from (\ref{equ: N coinvariant}).
\end{proof}

Now we recall the classical Bernstein-Zelevinsky geometric lemma.

\begin{lem}\label{lem: general smooth geometric lemma}
Let \ $I_0,I_1\subseteq I\subseteq \Delta$ \ and \ $\pi_0^{\infty}$ \ in \ $\mathrm{Rep}^{\infty}_{\rm{adm}}(L_{I_0})$ \ of \ finite \ length. \ Then \ $J_{I,I_1}(i_{I_0,I}^{\infty}(\pi_0^{\infty}))\in \mathrm{Rep}^{\infty}_{\rm{adm}}(L_{I_1})$ admits a canonical decreasing filtration with graded pieces $i_{I_0,I_1,w}^{\infty}(J_{I_0,I_1,w}(\pi_0^{\infty}))$ for $w\in W^{I_0,I_1}(L_I)$.
\end{lem}
\begin{proof}
This is \cite[\S 2.12]{BZ77} (see also \cite[\S VI.5.1]{Re10}), however since we use different normalizations, we need to make a translation. It follows from \cite[Thm.~VI.5.1]{Re10} that
\[J_{I,I_1}(i_{I_0,I}^{\infty}(\pi_0^{\infty}))\cong r_{P_{I_1}\cap L_I}^{L_I}(i_{P_{I_0}\cap L_I}^{L_I}(\pi_0^{\infty}\otimes_E\delta_{P_{I_0}\cap L_I}^{1/2}))\otimes_E \delta_{P_{I_1}\cap L_I}^{-1/2}\]
admits a decreasing filtration with graded pieces ($w\in W^{I_0,I_1}(L_I)$):
\begin{multline*}
\Big(i_{w^{-1}P_{I_0}w\cap L_{I_1}}^{L_{I_1}}\circ w\circ r_{wP_{I_1}w^{-1}\cap L_{I_0}}^{L_{I_0}}(\pi_0^{\infty}\otimes_E\delta_{P_{I_0}\cap L_{I}}^{1/2})\Big)\!\otimes_E \delta_{P_{I_1}\cap L_{I}}^{-1/2}\\
\begin{array}{rrl}
&\cong &\!\!i_{I_0,I_1,w}^{\infty}\Big(\!\big(J_{wP_{I_1}w^{-1}\cap L_{I_0}}(\pi_0^{\infty}\otimes_E\delta_{P_{I_0}\cap L_I}^{1/2})\otimes_E \delta_{wP_{I_1}w^{-1}\cap L_{I_0}}^{1/2}\big)^w\!\!\otimes_E \delta_{w^{-1}P_{I_0}w\cap L_{I_1}}^{-1/2}\Big)\!\otimes_E \delta_{P_{I_1}\cap L_{I}}^{-1/2}\\
&\cong &\!\!i_{I_0,I_1,w}^{\infty}\Big(\!\big(J_{wP_{I_1}w^{-1}\cap L_{I_0}}(\pi_0^{\infty})\otimes_E \delta_{P_{w(I_1)\cap I_0}\cap L_I}^{1/2}\big)^w \!\otimes_E \delta_{P_{w^{-1}(I_0)\cap I_1}\cap L_{I}}^{-1/2}\Big)\\
&\cong &\!\!i_{I_0,I_1,w}^{\infty}\Big(\!J_{wP_{I_1}w^{-1}\cap L_{I_0}}(\pi_0^{\infty})^w \otimes_E \big(\delta_{w^{-1}P_{w(I_1)\cap I_0}w\cap L_{I}}^{1/2} \!\cdot \delta_{P_{w^{-1}(I_0)\cap I_1}\cap L_{I}}^{-1/2}\big)\Big)\\
&\cong &\!\!i_{I_0,I_1,w}^{\infty}(J_{I_0,I_1,w}(\pi_0^{\infty}))
\end{array}
\end{multline*}
where the second isomorphism follows from (\ref{equ: outer modulus}) and (\ref{equ: inner modulus}) and the last from (\ref{equ: Jacquet twist}) in the definition of $J_{I_0,I_1,w}$.
\end{proof}

\begin{rem}\label{rem: rational geometric lemma}
We emphasize that the statement of Lemma \ref{lem: general smooth geometric lemma} holds for $E=K$, though we cannot take $E=K$ in its proof because of characters such as $\delta_{P_{I_0}\cap L_I}^{1/2}$ and $\delta_{P_{I_1}\cap L_I}^{-1/2}$ which are not $K^\times$-valued in general. The usual formulation of the geometric lemma (\cite[\S 2.12]{BZ77}, see also \cite[\S VI.5.1]{Re10}) uses \emph{normalized} induction and restriction functors. But if one reformulates the geometric lemma using \emph{unnormalized} induction and restriction functors (as we do in Lemma~\ref{lem: general smooth geometric lemma}), one sees that its proof works for any $E$, including $E=K$.
\end{rem}

\begin{lem}\label{lem: support of induction}
Let $I_0\subseteq I\subseteq \Delta$ and $\pi_0^{\infty}$ in $\mathrm{Rep}^{\infty}_{\rm{adm}}(L_{I_0})$ of finite length. Then we have
\begin{equation}\label{equ: support of induction}
\cJ(i_{I_0,I}^{\infty}(\pi_0^{\infty}))=\bigcup_{w\in W^{I_0,\emptyset}(L_I)}w^{-1}\cdot\cJ(\pi_0^{\infty})\subseteq W(L_I)\cdot\cJ(\pi_0^{\infty}).
\end{equation}
\end{lem}
\begin{proof}
By Lemma~\ref{lem: general smooth geometric lemma}, we know that $J_{I,\emptyset}(i_{I_0,I}^{\infty}(\pi_0^{\infty}))$ admits a canonical filtration with graded pieces $J_{I_0,\emptyset,w}(\pi_0^{\infty})$ for $w\in W^{I_0,\emptyset}$, which implies that
\begin{equation}\label{equ: support of induction 1}
\cJ(i_{I_0,I}^{\infty}(\pi_0^{\infty}))=\bigcup_{w\in W^{I_0,\emptyset}}\cJ(J_{I_0,\emptyset,w}(\pi_0^{\infty})).
\end{equation}
It follows from (\ref{equ: Bruhat Jacquet}) and (\ref{equ: twist Borel}) that
\[\cJ(J_{I_0,\emptyset,w}(\pi_0^{\infty}))=w^{-1}\cdot\cJ(J_{I_0,\emptyset}(\pi_0^{\infty}))=w^{-1}\cdot \cJ(\pi_0^{\infty}),\]
which together with (\ref{equ: support of induction 1}) gives (\ref{equ: support of induction}).
\end{proof}

\begin{rem}\label{rem: support of PS}
By taking $I_0=\emptyset$ and $\pi_0^{\infty}$ to be some $\chi\in\widehat{T}^{\infty}$ in Lemma~\ref{lem: support of induction}, we deduce that $\cJ(i_{\emptyset,I}^{\infty}(\chi))=W(L_I)\cdot\chi$. Suppose $\pi^{\infty}\in\mathrm{Rep}^{\infty}_{\rm{adm}}(L_I)$ is irreducible with $\chi\in \cJ(\pi^{\infty})$, then $\pi^{\infty}$ embeds into $i_{\emptyset,I}^{\infty}(\chi)$ (using (\ref{equ: first adjunction})) and $W(L_I)\cdot\cJ(\pi^{\infty})=W(L_I)\cdot\chi$.
\end{rem}

\begin{lem}
Let $I_0,I_1\subseteq I\subseteq \Delta$, $w\in W^{I_0,I_1}(L_I)$ and $\pi_0^{\infty}$ in $\mathrm{Rep}^{\infty}_{\rm{adm}}(L_{I_0})$ of finite length. Then we have
\begin{equation}\label{jjj}
\cJ(J_{I_0,I_1,w}(\pi_0^{\infty}))=w^{-1}\cdot\cJ(J_{I_0, w(I_1)\cap I_0}(\pi_0^{\infty}))=w^{-1}\cdot\cJ(\pi_0^{\infty})\subseteq w^{-1}W(L_{I_0})\cdot\cJ(\pi_0^{\infty})
\end{equation}
and
\begin{equation}\label{equ: support bound}
\cJ(i_{I_0,I_1,w}^{\infty}(J_{I_0,I_1,w}(\pi_0^{\infty})))\subseteq W(L_{I_1})\cdot \cJ(J_{I_0,I_1,w}(\pi_0^{\infty}))\subseteq W(L_{I_1})w^{-1}W(L_{I_0}) \cdot\cJ(\pi_0^{\infty}).
\end{equation}
\end{lem}
\begin{proof}
The first claim (\ref{jjj}) follows directly from (\ref{equ: Bruhat Jacquet}) and (\ref{equ: twist Borel}). The second claim (\ref{equ: support bound}) follows from (\ref{jjj}) together with (\ref{equ: support of induction}) (replacing $I_0$, $I$ and $\pi_0^{\infty}$ there by $w^{-1}(I_0)\cap I_1$, $I_1$ and $J_{I_0,I_1,w}(\pi_0^{\infty})$).
\end{proof}

\begin{lem}\label{lem: PS geometric lemma}
For a $G$-regular $\chi\in\widehat{T}^{\infty}$ (see \ref{it: basic 0} of Definition \ref{def: basic rep}) and $I\subseteq \Delta$, we have a canonical isomorphism in $\mathrm{Rep}^{\infty}_{\rm{adm}}(T)$
\begin{equation}\label{equ: PS geometric lemma}
J_{I,\emptyset}(i_{\emptyset,I}^{\infty}(\chi))\cong\bigoplus_{w\in W(L_I)}w^{-1}\cdot\chi.
\end{equation}
\end{lem}
\begin{proof}
By (\ref{equ: char case}) we have $i_{\emptyset,\emptyset,w}^{\infty}(J_{\emptyset,\emptyset,w}(\chi))=J_{\emptyset,\emptyset,w}(\chi)\cong w^{-1}\cdot\chi$. Then the statement follows from Lemma~\ref{lem: general smooth geometric lemma} applied with $I_0=I_1=\emptyset$ and $\pi_0^{\infty}=\chi$, noting that, as $\chi$ is $G$-regular, the $w\cdot\chi$ are distinct for $w\in W(L_I)$, so the canonical filtration in Lemma~\ref{lem: general smooth geometric lemma} must split.
\end{proof}

\begin{lem}\label{lem: Jacquet of PS}
Let $I\subseteq \Delta$.
\begin{enumerate}[label=(\roman*)]
\item \label{it: PS 1} Let $\pi^{\infty}$ in $\mathrm{Rep}^{\infty}_{\rm{adm}}(L_I)$ be $G$-basic, then the $T$-representation $J_{I,\emptyset}(\pi^{\infty})$ is semi-simple and multiplicity free and $W(L_I)\cdot\cJ(\pi^{\infty})$ is a single regular left $W(L_I)$-coset.
\item \label{it: PS 3} Let $\pi^{\infty}$ in $\mathrm{Rep}^{\infty}_{\rm{adm}}(L_I)$ be irreducible and $G$-regular, then the $L_I$-representation $\pi^{\infty}$ is uniquely determined by the set $\cJ(\pi^{\infty})$. Moreover non-isomorphic irreducible $G$-regular $\pi_0^{\infty}, \pi_1^{\infty}\in \mathrm{Rep}^{\infty}_{\rm{adm}}(L_I)$ satisfy $\cJ(\pi_0^{\infty})\cap \cJ(\pi_1^{\infty})=\emptyset$.
\item \label{it: PS 2} Let $\chi\in\widehat{T}^{\infty}$ be $G$-regular, then $i_{\emptyset,I}^{\infty}(\chi)$ is multiplicity free and any irreducible constituent $\pi^{\infty}$ of $i_{\emptyset,I}^{\infty}(\chi)$ is such that $\cJ(\pi^{\infty})\neq 0$. Moreover for $w\in W(L_I)$ the semi-simplification of $i_{\emptyset,I}^{\infty}(w\cdot\chi)$ doesn't depend on $w$.
\end{enumerate}
\end{lem}
\begin{proof}
We prove \ref{it: PS 1}. For a $G$-regular $\chi$, Lemma~\ref{lem: PS geometric lemma} implies that $J_{I,\emptyset}(i_{\emptyset,I}^{\infty}(\chi))$ is multiplicity free and
\begin{equation}\label{equ: disjoint union support}
W(L_I)\cdot\chi=\cJ(i_{\emptyset,I}^{\infty}(\chi))=\bigsqcup_{\pi^{\infty}}\cJ(\pi^{\infty})
\end{equation}
where $\pi^{\infty}$ runs through the Jordan-H\"older factors of $i_{\emptyset,I}^{\infty}(\chi)$. This implies
\ref{it: PS 1} as each $G$-basic $\pi^{\infty}$ embeds into $i_{\emptyset,I}^{\infty}(\chi)$ for some $G$-regular $\chi$ by \ref{it: basic 1} of Definition~\ref{def: basic rep}. The first half of \ref{it: PS 2} follows from \cite[Prop.~2.1(b),(c)]{Z80}. We prove \ref{it: PS 3}. The first half of \ref{it: PS 2} together with (\ref{equ: disjoint union support}) easily imply the first half of \ref{it: PS 3}. For the second half of \ref{it: PS 3}, any $\pi_0^{\infty}, \pi_1^{\infty}\in \mathrm{Rep}^{\infty}_{\rm{adm}}(L_I)$ satisfying $\chi\in\cJ(\pi_0^{\infty})\cap \cJ(\pi_1^{\infty})$ both inject into $i_{\emptyset,I}^{\infty}(\chi)$, again contradicting (\ref{equ: disjoint union support}). Finally, the second half of \ref{it: PS 2} follows from \ref{it: PS 3} and the fact that $\cJ(i_{\emptyset,I}^{\infty}(w\cdot\chi))=W(L_I)\cdot \chi$ (by (\ref{equ: disjoint union support})) is independent of the choice of $w\in W(L_I)$.
\end{proof}

In particular it follows from Lemma \ref{lem: Jacquet of PS} that if $\pi^{\infty}$ is a $G$-basic representation of $L_I$ then $\pi^{\infty}$ is in $\cB^I_{\Sigma}$ for $\Sigma=W(L_I)\cdot\cJ(\pi^{\infty})$. We give below several other useful remarks on $G$-basic representations.

\begin{rem}\label{rem: basic PS intertwine}
\hspace{2em}
\begin{enumerate}[label=(\roman*)]
\item \label{it: basic as quotient}
Let $I\subseteq \Delta$ and $\pi^{\infty}\in\mathrm{Rep}^{\infty}_{\rm{adm}}(L_I)$ be $G$-basic, then by \ref{it: basic 1} of Definition~\ref{def: basic rep} there exist $I_1\subseteq I$ and an irreducible $G$-regular $\pi_1^{\infty}\in\mathrm{Rep}^{\infty}_{\rm{adm}}(L_{I_1})$ such that $\pi^{\infty}\cong i_{I_1,I}^{\infty}(\pi_1^{\infty})$. As $J_{I_1,\emptyset}(\pi_1^{\infty})$ is semi-simple and multiplicity free by \ref{it: PS 1} of Lemma~\ref{lem: Jacquet of PS}, we deduce from (\ref{equ: explicit twist Jacquet}) that $J_{I_1,\emptyset}'(\pi_1^{\infty})$ is also semi-simple and multiplicity free with
\begin{equation}\label{equ: JH Jacquet twist}
\mathrm{JH}_{T}(J_{I_1,\emptyset}'(\pi_1^{\infty}))=w_{I_1}\cdot \mathrm{JH}_{T}(J_{I_1,\emptyset}(\pi_1^{\infty}))
\end{equation}
where we recall from \S\ref{generalnotation} that $w_{I_1}$ is the longest element of $W(L_{I_1})$ (so that $w_{I_1}(\delta_{B_{I_1}})=w_{I_1}^{-1}(\delta_{B_{I_1}})=\delta_{B_{I_1}}^{-1}$) and $\mathrm{JH}_{T}(-)$ the set of irreducible constituents. Hence, by combining (\ref{equ: JH Jacquet twist}) with (\ref{equ: first adjunction}) and (\ref{equ: second adjunction}) for $k=0$, we see that $\pi_1^{\infty}$ embeds into $i_{\emptyset,I_1}^{\infty}(\chi)$ if and only if $\pi_1^{\infty}$ is a quotient of $i_{\emptyset,I_1}^{\infty}(w_{I_1}\cdot\chi)$. In particular, we deduce that $\pi^{\infty}\cong i_{I_1,I}^{\infty}(\pi_1^{\infty})$ is a quotient of $i_{\emptyset,I}^{\infty}(w_{I_1}\cdot\chi)$ for any $\chi\in\cJ(\pi_1^{\infty})$.
\item \label{it: PS socle cosocle}
Let $I\subseteq \Delta$, $\pi^{\infty}\in\mathrm{Rep}^{\infty}_{\rm{adm}}(L_I)$ be irreducible and $G$-regular, and $\chi\in\widehat{T}^{\infty}$ be $G$-regular. By (\ref{equ: first adjunction}) for $k=0$ (resp.~by (\ref{equ: second adjunction}) for $k=0$ with (\ref{equ: JH Jacquet twist})) $\pi^{\infty}$ is in the socle (resp.~cosocle) of $i_{\emptyset,I}^{\infty}(\chi)$ if and only if $\chi\in\cJ(\pi^{\infty})$ (resp.~$\chi\in w_I\cdot\cJ(\pi^{\infty})$). The second half of \ref{it: PS 3} of Lemma~\ref{lem: Jacquet of PS} then implies that $i_{\emptyset,I}^{\infty}(\chi)$ has simple socle and cosocle.
\item \label{it: PS intertwine}
For $G$-regular $\chi,\chi'\in\widehat{T}^{\infty}$ such that $W(L_I)\cdot\chi=W(L_I)\cdot\chi'$, we have canonical isomorphisms by (\ref{equ: first adjunction}) (for $k=0$) and Lemma~\ref{lem: PS geometric lemma}:
\[\Hom_{L_I}(i_{\emptyset,I}^{\infty}(\chi'),i_{\emptyset,I}^{\infty}(\chi))\cong \Hom_{T}(J_{I,\emptyset}(i_{\emptyset,I}^{\infty}(\chi')),\chi)\cong \Hom_{T}\bigg(\bigoplus_{w\in W(L_I)}w^{-1}\cdot\chi',\chi\bigg)\]
of one dimensional spaces. In particular, there exists a unique (up to scalar) non-zero map $i_{\emptyset,I}^{\infty}(\chi')\rightarrow i_{\emptyset,I}^{\infty}(\chi)$ for such $\chi,\chi'$.
\item \label{it: basic as image}
Let $I\subseteq \Delta$ and $\pi^{\infty}\in\mathrm{Rep}^{\infty}_{\rm{adm}}(L_I)$ be $G$-basic. By \ref{it: basic 1} of Definition~\ref{def: basic rep} and \ref{it: basic as quotient} above there exist $G$-regular $\chi,\chi'\in\widehat{T}^{\infty}$ such that $\pi^{\infty}$ is a subrepresentation of $i_{\emptyset,I}^{\infty}(\chi)$ and a quotient of $i_{\emptyset,I}^{\infty}(\chi')$. As both $i_{\emptyset,I}^{\infty}(\chi)$ and $i_{\emptyset,I}^{\infty}(\chi')$ are multiplicity free by \ref{it: PS 2} of Lemma~\ref{lem: Jacquet of PS}, we deduce from \ref{it: PS intertwine} above that $\pi^{\infty}$ is the image of the unique (up to scalar) non-zero map $i_{\emptyset,I}^{\infty}(\chi')\rightarrow i_{\emptyset,I}^{\infty}(\chi)$, is multiplicity free and (using \ref{it: PS socle cosocle} above) has simple socle and cosocle.
\end{enumerate}
\end{rem}

\begin{lem}\label{lem: G basic twist}
Let $I\subseteq \Delta$, $w\in W(G)$ such that $w^{-1}(I)\subseteq \Delta$, $\delta: L_{w^{-1}(I)}\rightarrow E^\times$ a smooth character and $\pi^{\infty}\in\mathrm{Rep}^{\infty}_{\rm{adm}}(L_I)$ a $G$-basic representation. Assume that $\cJ((\pi^{\infty})^w\otimes_E\delta)\subseteq W(G)\cdot\cJ(\pi^{\infty})$. Then $(\pi^{\infty})^w\otimes_E\delta\in\mathrm{Rep}^{\infty}_{\rm{adm}}(L_{w^{-1}(I)})$ is again $G$-basic.
\end{lem}
\begin{proof}
By \ref{it: basic 1} of Definition~\ref{def: basic rep} there exist $I_1\subseteq I$ and an irreducible $G$-regular $\pi_1^{\infty}\in\mathrm{Rep}^{\infty}_{\rm{adm}}(L_{I_1})$ such that $\pi^{\infty}\cong i_{I_1,I}^{\infty}(\pi_1^{\infty})$, which implies that
\begin{equation}\label{equ: induction twist}
(\pi^{\infty})^w\otimes_E\delta\cong i_{w^-1(I_1),w^-1(I)}^\infty((\pi_1^{\infty})^w\otimes_E\delta|_{L_{w^{-1}(I_1)}}).
\end{equation}
As $W(L_I)\cdot\cJ(\pi^{\infty})$ is a $G$-regular left $W(L_I)$-coset by \ref{it: PS 1} of Lemma~\ref{lem: Jacquet of PS}, the sets
\[\cJ((\pi_1^{\infty})^w\otimes_E\delta|_{L_{w^{-1}(I_1)}})\subseteq \cJ((\pi^{\infty})^w\otimes_E\delta)\subseteq W(G)\cdot \cJ(\pi^{\infty})\]
consist of $G$-regular elements. This forces $(\pi_1^{\infty})^w\otimes_E\delta|_{L_{I_1}}$ to be $G$-regular (irreducible) and thus (\ref{equ: induction twist}) implies that $(\pi^{\infty})^w\otimes_E\delta$ is $G$-basic.
\end{proof}

\begin{lem}\label{lem: smooth geometric lemma}
Let $I_0,I_1\subseteq I\subseteq \Delta$, $\Sigma_0\subseteq \widehat{T}^{\infty}$ be a left $W(L_{I_0})$-stable finite subset which is $G$-regular, and $\pi_0^{\infty}\in\cB^{I_0}_{\Sigma_0}$.
\begin{enumerate}[label=(\roman*)]
\item \label{it: sml1}
We have a canonical decomposition in $\mathrm{Rep}^{\infty}_{\rm{adm}}(L_{I_1})$
\begin{equation}\label{equ: smooth geometric lemma}
J_{I,I_1}(i_{I_0,I}^{\infty}(\pi_0^{\infty}))\cong \bigoplus_{w\in W^{I_0,I_1}(L_I)}i_{I_0,I_1,w}^{\infty}(J_{I_0,I_1,w}(\pi_0^{\infty})).
\end{equation}
\item \label{it: sml2}
If $\Sigma_0$ is a $G$-regular left $W(L_{I_0})$-coset, we have a canonical decomposition for $w\in W^{I_0,I_1}(L_I)$ induced by (\ref{equ: sm block decomposition})
\begin{equation}\label{equ: refine geometric lemma}
i_{I_0,I_1,w}^{\infty}(J_{I_0,I_1,w}(\pi_0^{\infty}))\cong \bigoplus_{\Sigma}i_{I_0,I_1,w}^{\infty}(J_{I_0,I_1,w}(\pi_0^{\infty}))_{\cB^{I_1}_{\Sigma}}
\end{equation}
where $\Sigma$ runs through the ($G$-regular) left $W(L_{I_1})$-cosets contained in
\[W(L_{I_1})w^{-1}W(L_{I_0})\cdot\cJ(\pi_0^{\infty}).\]
In particular $J_{I,I_1}(i_{I_0,I}^{\infty}(\pi_0^{\infty}))_{\cB^{I_1}_{\Sigma}}\cong i_{I_0,I_1,w}^{\infty}(J_{I_0,I_1,w}(\pi_0^{\infty}))_{\cB^{I_1}_{\Sigma}}$.
\end{enumerate}
\end{lem}
\begin{proof}
Using (\ref{equ: sm block decomposition}), it suffices to treat $(\pi_0^{\infty})_{\cB^{I_0}_{\Sigma_0'}}$ separately for each left $W(L_{I_0})$-coset $\Sigma_0'$ contained in $\Sigma_0$. Hence, we assume from now that $\Sigma_0$ is a single $G$-regular left $W(L_{I_0})$-coset, namely $\Sigma_0=W(L_{I_0})\cdot\cJ(\pi_0^{\infty})=W(L_{I_0})\cdot\chi$ for an arbitrary $\chi\in\cJ(\pi_0^{\infty})$.
As $\chi\in \cJ(\pi_0^{\infty})$ is $G$-regular, we know that
\begin{multline*}
W(L_I)\cdot\cJ(\pi_0^{\infty})=W(L_I)\cdot\chi=\bigsqcup_{w\in W^{I_0,I_1}(L_I)}W(L_{I_1})w^{-1}W(L_{I_0})\cdot\chi\\
=\bigsqcup_{w\in W^{I_0,I_1}(L_I)}W(L_{I_1})w^{-1}W(L_{I_0})\cdot\cJ(\pi_0^{\infty}),
\end{multline*}
which together with (\ref{equ: support bound}) and (\ref{equ: sm block decomposition}) forces the canonical filtration on $J_{I,I_1}(i_{I_0,I}^{\infty}(\pi_0^{\infty}))$ described in Lemma~\ref{lem: general smooth geometric lemma} to be split, and thus \ref{it: sml1} follows. \ref{it: sml2} follows directly from \ref{it: sml1}, (\ref{equ: support bound}) and (\ref{equ: sm block decomposition}).
\end{proof}

\begin{rem}\label{rem: general PS Jacquet}
If we take $I_0=\emptyset$ and $\pi_0^{\infty}=\chi$ in Lemma~\ref{lem: smooth geometric lemma}, then by (\ref{equ: char case}) we deduce a canonical isomorphism in $\mathrm{Rep}^{\infty}_{\rm{adm}}(L_{I_1})$
\begin{equation*}
J_{I,I_1}(i_{\emptyset,I}^{\infty}(\chi))\cong \bigoplus_{w\in W^{\emptyset,I_1}(L_I)}i_{\emptyset,I_1}^{\infty}(w^{-1}\cdot\chi).
\end{equation*}
\end{rem}

Let $I\subseteq \Delta$ and $\Sigma\subseteq \widehat{T}^{\infty}$ be a $G$-regular left $W(L_I)$-coset. We now attach an undirected graph $\Gamma_{\Sigma}$ to $\Sigma$ following \cite[\S 2.2]{Z80}. We choose an arbitrary $\chi\in\Sigma$ and write
\begin{equation}\label{equ: cuspidal data}
\chi\otimes_E \delta_{B}^{1/2}=\rho_1\boxtimes\cdots\boxtimes\rho_n
\end{equation}
where the $\rho_k$, $1\leq k\leq n$, are smooth distinct characters of $K^\times$ (as $\chi$ is $G$-regular). We define the set of vertices of $\Gamma_{\Sigma}$ as $V(\Gamma_{\Sigma})\defeq \{\rho_1,\dots,\rho_n\}$ and note that the set $V(\Gamma_{\Sigma})$ is independent of the choice of $\chi\in\Sigma$. Two vertices $\rho$, $\rho'\in V(\Gamma_{\Sigma})$ are connected by one edge if and only if $\rho'\rho^{-1}\in \{|\cdot|_K,|\cdot|_K^{-1}\}$ and $\rho=\rho_k$, $\rho'=\rho_{k'}$ for $k,k'$ in the same Levi block of $L_I$. So each connected component of $\Gamma_{\Sigma}$ (there are at least as many as the number of blocks in $L_I$) has its vertices of the form $\{\rho,|\cdot|_K\rho,\dots,|\cdot|_K^{\ell-1}\rho\}$ for some $\ell\geq 1$, which is called a segment $[\rho, |\cdot|_K^{\ell-1}\rho]$ of length $\ell$. An \emph{orientation} $\vec{\Gamma_{\Sigma}}$ (on $\Gamma_{\Sigma}$) is a directed graph whose underlying undirected graph is $\Gamma_{\Sigma}$. Each $\chi\in\Sigma$ as above determines an orientation denoted $\vec{\Gamma_{\Sigma}}(\chi)$ by requiring that an edge connecting $\rho$ and $\rho'$ has direction $\rho\rightarrow \rho'$ if and only if $\rho=\rho_k$ and $\rho'=\rho_{k'}$ for some $k<k'$.\bigskip

Recall from the last statement in \ref{it: PS 2} of Lemma~\ref{lem: Jacquet of PS} that the set of irreducible constituents of $i_{\emptyset,I}^{\infty}(\chi)$ is independent of $\chi\in\Sigma$, and we denote it by $\mathrm{JH}_{\Sigma}$.

\begin{thm}[\cite{Z80}, Thm.~2.2]\label{thm: JH of PS}
There exists a unique bijection $\vec{\Gamma_{\Sigma}}\mapsto \omega(\vec{\Gamma_{\Sigma}})$ between the set of all orientations on $\Gamma_{\Sigma}$ and the set $\mathrm{JH}_{\Sigma}$ such that
\begin{equation}\label{equ: JH of PS}
\cJ(\omega(\vec{\Gamma_{\Sigma}}))=\{\chi\mid \vec{\Gamma_{\Sigma}}(\chi)=\vec{\Gamma_{\Sigma}}\}.
\end{equation}
\end{thm}

\begin{rem}
Note that (\ref{equ: JH of PS}) differs from the statement in \cite[Thm.~2.2]{Z80} by a twist~$\delta_{B}^{1/2}$, and that the unicity of $\omega$ in Theorem \ref{thm: JH of PS} follows from \ref{it: PS 3} of Lemma~\ref{lem: Jacquet of PS}.
\end{rem}

\begin{rem}\label{rem: rational JH}
Let $I\subseteq \Delta$ and $\Sigma\subseteq \widehat{T}^{\infty}$ be a $G$-regular left $W(L_I)$-coset.
We say that two elements of $\Sigma$ are \emph{equivalent} if they correspond to the same orientation on $\Gamma_{\Sigma}$.
As $\delta_{B}^{1/2}$ might not be $K^\times$-valued, $\Gamma_{\Sigma}$ might not be defined as in (\ref{equ: cuspidal data}) in general. However, the equivalence relation on $\Sigma$ discussed above is always well-defined regardless of whether $\delta_{B}^{1/2}$ is $K^\times$-valued or not. In particular, Theorem~\ref{thm: JH of PS} still makes sense when $E=K$.
We also observe that the equivalence relation on $\Sigma$ remains unchanged if we view $\Sigma$ as a set of $(E')^\times$-valued characters for some finite extension $E'$ of $E$. This together with Theorem~\ref{thm: JH of PS} implies that each constituents in $\mathrm{JH}_{\Sigma}$ is absolutely irreducible.
\end{rem}

Given $I\subseteq \Delta$ and a multiplicity free finite length representation $\pi^{\infty}$ in $\mathrm{Rep}^{\infty}_{\rm{adm}}(L_I)$, recall that we have defined in \S\ref{generalnotation} a partial order on the set $\mathrm{JH}_{L_I}(\pi^{\infty})$ of constituents of $\pi^{\infty}$. We slightly reformulate \cite[Thm.~2.8]{Z80} as follows.

\begin{thm}[\cite{Z80}, Thm.~2.8]\label{thm: structure of PS}
Let $\chi\in\Sigma$. The partial order on the set of orientations on $\Gamma_{\Sigma}$ given via the bijection in Theorem \ref{thm: JH of PS} by the partial order on $\mathrm{JH}_{L_I}(i_{\emptyset,I}^{\infty}(\chi))$ is the following: two orientations $\vec{\Gamma_{\Sigma}}$, $\vec{\Gamma_{\Sigma}}'$ satisfy $\vec{\Gamma_{\Sigma}}\leq \vec{\Gamma_{\Sigma}}'$ if and only if each edge of $\Gamma_{\Sigma}$ which has the same direction in $\vec{\Gamma_{\Sigma}}'$ and $\vec{\Gamma_{\Sigma}}(\chi)$ also has the same direction in $\vec{\Gamma_{\Sigma}}$.
\end{thm}

Concretely, $\vec{\Gamma_{\Sigma}}(\chi)$ is the orientation on $\Gamma_{\Sigma}$ corresponding to the socle $\omega(\vec{\Gamma_{\Sigma}}(\chi))$ of $i_{\emptyset,I}^{\infty}(\chi)$ (which is irreducible by \ref{it: basic as image} of Remark~\ref{rem: basic PS intertwine}), and one gets the orientations corresponding to the Jordan-H\"older factors in higher layers by successively reversing (more and more) arrows in $\vec{\Gamma_{\Sigma}}(\chi)$, until all arrows of $\vec{\Gamma_{\Sigma}}(\chi)$ are reversed which gives the orientation corresponding to the (irreducible by \emph{loc.\ cit.}) cosocle of $i_{\emptyset,I}^{\infty}(\chi)$.\bigskip

For $I_1\subseteq I\subseteq \Delta$ we define
\begin{equation}\label{gensteinberg}
V_{I_1,I}^{\infty}\defeq i_{I_1,I}^{\infty}(1_{L_{I_1}})/\sum_{I_1\subsetneq I_1'}i_{I_1',I}^{\infty}(1_{L_{I_1'}})
\end{equation}
which is an (absolutely irreducible) smooth generalized Steinberg representation of $L_I$ (the smooth Steinberg being $V_{\emptyset,I}^{\infty}$). Note that it is $G$-regular (for instance by \ref{it: PS 2} of Lemma \ref{lem: Jacquet of PS} and (\ref{equ: PS geometric lemma})). The following corollary is classical, we provide a short proof for the reader's convenience.

\begin{cor}\label{cor: trivial PS}
Let $I\subseteq \Delta$ and $w\in W(L_I)$.
\begin{enumerate}[label=(\roman*)]
\item \label{it: trivial PS 1} For $I_1\subseteq I$ we have
\begin{equation}\label{equ: Jacquet Steinberg}
\cJ(V_{I_1,I}^{\infty})=\{x\cdot 1_{T}\mid x\in W(L_I), I_1=I\setminus D_R(x)\}.
\end{equation}
\item \label{it: trivial PS 2} The representation $i_{\emptyset,I}^{\infty}(w\cdot 1_T)$ is multiplicity free with socle $V_{I\setminus D_R(w), I}^{\infty}$, cosocle $V_{I\cap D_R(w),I}^{\infty}$ and constituents
\begin{equation}\label{equ: JH trivial PS}
\mathrm{JH}_{L_I}(i_{\emptyset,I}^{\infty}(w\cdot 1_T))=\{V_{I_1,I}^{\infty}\mid I_1\subseteq I\}.
\end{equation}
\item \label{it: trivial PS 3} The partial order $\leq_w$ on $\{V_{I_1,I}^{\infty}\mid I_1\subseteq I\}$ induced from the one on $\mathrm{JH}_{L_I}(i_{\emptyset,I}^{\infty}(w\cdot 1_T))$ via (\ref{equ: JH trivial PS}) is: $V_{I_1,I}^{\infty}\leq_w V_{I_2,I}^{\infty}$ (with $I_1,I_2\subseteq I$) if and only if $I_2\cap (I\setminus D_R(w))\subseteq I_1\subseteq I_2\cup(I\setminus D_R(w))$.
\end{enumerate}
\end{cor}
\begin{proof}
We only prove the case $I=\Delta$, the general case follows by treating each Levi block of $L_I$ separately. The graph $\Gamma \defeq\Gamma_{W(G)\cdot 1_{T}}$ attached to $1_T$ (or equivalently to any $w\cdot 1_T$ with $w\in W(G)$) has vertices $\{1, |\cdot|_K,\dots, |\cdot|_K^{n-1}\}$. The set of orientations $\vec{\Gamma}$ is in natural bijection to the set of subsets of $\Delta$ by sending $\vec{\Gamma}$ to the subset $I_{\vec{\Gamma}}$ of $j\in\Delta$ such that $\vec{\Gamma}$ contains an arrow $|\cdot|_K^{n-j}\rightarrow |\cdot|_K^{n-j-1}$. For $w\in W(G)$, the attached orientation $\vec{\Gamma}(w)\defeq\vec{\Gamma}(w\cdot 1_T)$ is characterized as follows: for $j\in\Delta$, there exists an arrow $|\cdot|_K^{n-j}\rightarrow |\cdot|_K^{n-j-1}$ if $j\notin D_R(w)$, and an arrow $|\cdot|_K^{n-j-1}\rightarrow |\cdot|_K^{n-j}$ if $j\in D_R(w)$, i.e.~$I_{\vec{\Gamma}(w)}=\Delta\setminus D_R(w)$. This together with (\ref{equ: JH of PS}) and the discussion after Theorem~\ref{thm: structure of PS} imply
\begin{equation}\label{equ: Jacquet soc PS}
\begin{array}{rcl}
\cJ(\mathrm{soc}_{G}(i_{\emptyset,\Delta}^{\infty}(w\cdot 1_{T})))&=&\{x\cdot 1_{T}\mid x\in W(G), D_R(x)=D_R(w)\}\\
\cJ(\mathrm{cosoc}_{G}(i_{\emptyset,\Delta}^{\infty}(w\cdot 1_{T})))&=&\{x\cdot 1_{T}\mid x\in W(G), D_R(x)=\Delta\setminus D_R(w)\}.
\end{array}
\end{equation}
We now prove \ref{it: trivial PS 1}. By (\ref{equ: support of induction}) with $I_0$, $I$, $\pi_0^{\infty}$ there being $I_1$, $\Delta$, $1_{L_{I_1}}$, we have for $I_1\subseteq \Delta$
\[\cJ(i_{I_1,\Delta}^{\infty}(1_{L_{I_1}}))=\{x^{-1}\cdot 1_{T}\mid x\in W^{I_1,\emptyset}\}=\{x\cdot 1_{T}\mid I_1\subseteq \Delta\setminus D_R(x)\}.\]
Together with $V_{I_1,\Delta}^{\infty}=i_{I_1,\Delta}^{\infty}(1_{L_{I_1}})/\sum_{I_1\subsetneq I_1'}i_{I_1',\Delta}^{\infty}(1_{L_{I_1'}})$ this implies
\begin{multline}\label{equ: Jacquet Steinberg definition}
\cJ(V_{I_1,\Delta}^{\infty})=\{x\cdot 1_{T}\mid I_1\subseteq \Delta\setminus D_R(x)\}\setminus \{x\cdot 1_{T}\mid I_1\subsetneq \Delta\setminus D_R(x)\}\\
=\{x\cdot 1_{T}\mid I_1=\Delta\setminus D_R(x)\}
\end{multline}
which is \ref{it: trivial PS 1}. Then \ref{it: trivial PS 2} follows from (\ref{equ: Jacquet Steinberg definition}), (\ref{equ: Jacquet soc PS}) and
\[\{x\cdot 1_{T}\mid x\in W(G)\}=\bigsqcup_{I_1\subseteq \Delta} \{x\cdot 1_{T}\mid \Delta\setminus D_R(x)=I_1\}.\]

We prove \ref{it: trivial PS 3}. We write $\vec{\Gamma}_i$ for the orientation on $\Gamma$ with $I_{\vec{\Gamma}_i}=I_i$, $i=1,2$. It follows from Theorem~\ref{thm: structure of PS} that $V_{I_1,\Delta}^{\infty}\leq_w V_{I_2,\Delta}^{\infty}$ if and only if each edge of $\Gamma$ which has the same direction in $\vec{\Gamma}_2$ and $\vec{\Gamma}(w)$ also has this direction in $\vec{\Gamma}_1$. Using the above definition of $I_{\vec{\Gamma}_i}$ and $I_{\vec{\Gamma}(w)}=\Delta\setminus D_R(w)$, this can be easily translated to $I_2\cap (\Delta\setminus D_R(w))\subseteq I_1$ and $(\Delta\setminus I_2)\cap D_R(w)\subseteq \Delta\setminus I_1$, or equivalently $I_2\cap (I\setminus D_R(w))\subseteq I_1\subseteq I_2\cup(I\setminus D_R(w))$ which is \ref{it: trivial PS 3}.
\end{proof}

\begin{cor}\label{cor: twist of St}
Let $I\subseteq \Delta$, $\pi^{\infty}$ in $\mathrm{Rep}^{\infty}_{\rm{adm}}(L_I)$ irreducible $G$-regular and $\Gamma_{\Sigma}$ the undirected graph attached to $\Sigma\defeq W(L_I)\cdot\cJ(\pi^{\infty})$ above Theorem~\ref{thm: JH of PS} (recall that $\Sigma$ is a single regular left $W(L_I)$-coset by \ref{it: PS 1} of Lemma~\ref{lem: Jacquet of PS}). If $\Gamma_{\Sigma}$ has one connected component for each Levi block of $L_I$, then there exists $I_1\subseteq I$ and a smooth character $\delta: L_I\rightarrow E^\times$ such that $\pi^{\infty}\cong V_{I_1,I}^{\infty}\otimes_E \delta$.
\end{cor}
\begin{proof}
It is harmless to treat each Levi block of $L_I$ separately, so we may assume $I=\Delta$ and $\Gamma_{\Sigma}$ connected. By definition of $\Gamma_{\Sigma}$, it is connected if and only if its set of vertices has the form $\{\rho,\rho\otimes_E|\cdot|_K,\dots,\rho\otimes_E|\cdot|_K^{n-1}\}$ for some smooth character $\rho: K^\times\rightarrow E^\times$. Let $\delta\defeq \rho\circ\mathrm{det}$ where $\mathrm{det}: G\rightarrow E^\times$ is the determinant character. Then the graph $\Gamma_{\Sigma'}$ attached to $\Sigma'\defeq W(G)\cdot\cJ(\pi^{\infty}\otimes_E \delta^{-1})$ has vertices $\{1,|\cdot|_K,\dots,|\cdot|_K^{n-1}\}$, which forces $\pi^{\infty}\otimes_E \delta^{-1}$ to be a Jordan-H\"older factor of $i_{\emptyset,\Delta}^{\infty}(1_T)$, and thus to be $V_{I_1,\Delta}^{\infty}$ for some $I_1\subseteq \Delta$ by (\ref{equ: JH trivial PS}).
\end{proof}

\begin{cor}\label{cor: basic subquotient}
Let $I\subseteq \Delta$ and $\chi\in\widehat{T}^{\infty}$ $G$-regular. Let $\pi^{\infty}$ be a subquotient of $i_{\emptyset,I}^{\infty}(\chi)$ with simple socle and cosocle. Then $\pi^{\infty}$ is $G$-basic and there exist $I_0\subseteq I$ and an irreducible $G$-regular $\sigma^{\infty}\in\mathrm{Rep}^{\infty}_{\rm{adm}}(L_{I_0})$ such that $\pi^{\infty}\cong i_{I_0,I}^{\infty}(\sigma^{\infty})$.
\end{cor}
\begin{proof}
Write $\sigma_1^{\infty}$ (resp.~$\sigma_2^{\infty}$) for the socle (resp.~cosocle) of $\pi^{\infty}$. For $\chi'\in\cJ(\sigma_1^{\infty})$, the image of the unique (up to scalar) non-zero map $i_{\emptyset,I}^{\infty}(\chi)\rightarrow i_{\emptyset,I}^{\infty}(\chi')$ between multiplicity free objects (see \ref{it: PS 2} of Lemma~\ref{lem: Jacquet of PS} and \ref{it: PS intertwine} of Remark~\ref{rem: basic PS intertwine}) is the unique quotient of $i_{\emptyset,I}^{\infty}(\chi)$ with socle $\sigma_1^{\infty}$, and thus contains $\pi^{\infty}$ (as $\pi^{\infty}$ has simple socle $\sigma_1^{\infty}$). Upon replacing $\chi$ with $\chi'$, we assume from on that $\pi^{\infty}$ is a subrepresentation of $i_{\emptyset,I}^{\infty}(\chi)$, hence $\sigma_1^{\infty}$ is the socle of $i_{\emptyset,I}^{\infty}(\chi)$.

Let $\Gamma_{\Sigma}$ be the undirected graph attached to $\Sigma\defeq W(L_I)\cdot\chi$ (see before Theorem~\ref{thm: JH of PS}), $\vec{\Gamma_{\Sigma,i}}\defeq \omega^{-1}(\sigma_i^{\infty})$, $i=1,2$ (see (\ref{equ: JH of PS})), and recall that $\vec{\Gamma_{\Sigma,1}}=\vec{\Gamma_{\Sigma}}(\chi)$. Write $\chi\otimes_E \delta_{B}^{1/2}=\rho_1\boxtimes\cdots\boxtimes\rho_n$, modifying $\chi$ within $W(L_I)\cdot\chi$ without changing $\vec{\Gamma_{\Sigma}}(\chi)$, we can assume that each pair of vertices $\rho_k$, $\rho_{k'}$ that are connected by an edge of $\Gamma_{\Sigma}$ are such that $|k-k'|=1$. For such a $\chi$, the set of edges of $\Gamma_{\Sigma}$ having the same orientation in $\vec{\Gamma_{\Sigma,1}}$ and $\vec{\Gamma_{\Sigma,2}}$ naturally determines a subset $I_0\subseteq I$. Let $\Gamma_{\Sigma_0}$ be the undirected graph attached to $\Sigma_0\defeq W(L_{I_0})\cdot\chi$, and note that $\Gamma_{\Sigma_0}$ is obtained from $\Gamma_{\Sigma}$ by exactly deleting the edges corresponding to elements \emph{not} in $I_0$. Each orientation on $\Gamma_{\Sigma}$ induces one on $\Gamma_{\Sigma_0}$, and by definition of $I_0$ $\vec{\Gamma_{\Sigma,1}}$ and $\vec{\Gamma_{\Sigma,2}}$ induce the same orientation $\vec{\Gamma_{\Sigma_0}}(\chi)$ on $\Gamma_{\Sigma_0}$, hence a well-defined Jordan-H\"older factor $\sigma^{\infty}$ of $i_{\emptyset,I_0}^{\infty}(\chi)$ by Theorem~\ref{thm: JH of PS}, which is actually its socle. Note that $\chi\in \cJ(\sigma^{\infty})$ from (\ref{equ: first adjunction}) (for $k=0$), hence $\sigma^{\infty}$ is $G$-regular.

We claim that $i_{I_0,I}^{\infty}(\sigma^{\infty})\cong \pi^{\infty}$. As $i_{I_0,I}^{\infty}(\sigma^{\infty})$ is a subrepresentation of $i_{\emptyset,I}^{\infty}(\chi)$, it suffices to show that $i_{I_0,I}^{\infty}(\sigma^{\infty})$ has cosocle $\sigma_2^{\infty}$. Let $w_{I,I_0}$ be the longest element inside $W^{\emptyset,I_0}(L_{I})$, by (\ref{equ: twist Jacquet support}) we have
\begin{equation}\label{equ: cosocle twist support}
\cJ(J_{I,I_0}'(\sigma_2^{\infty}))=w_{I,I_0}^{-1}\cdot\cJ(\sigma_2^{\infty}).
\end{equation}
The orientation on $\Gamma_{\Sigma}$ associated to $w_{I,I_0}\cdot\chi$ is easily checked to be obtained from the one associated to $\chi$ by reversing the orientation of the edges of $\Gamma_{\Sigma}$ corresponding to $I\setminus I_0$, which is $\vec{\Gamma_{\Sigma,2}}= \omega^{-1}(\sigma_2^{\infty})$ by definition of $I_0$. In particular we have $w_{I,I_0}\cdot\chi\in \cJ(\sigma_2^{\infty})$, which together with (\ref{equ: cosocle twist support}) implies $\chi\in \cJ(J_{I,I_0}'(\sigma_2^{\infty}))$. By (\ref{equ: second adjunction}) for $k=0$ we have for $\tau^{\infty}\defeq \mathrm{cosoc}_{L_I}(i_{I_0,I}^{\infty}(\sigma^{\infty}))$
\[0\neq \Hom_{L_I}(i_{I_0,I}^{\infty}(\sigma^{\infty}),\tau^{\infty})\cong \Hom_{L_{I_0}}(\sigma^{\infty},J_{I,I_0}'(\tau^{\infty})).\]
Thus the irreducible $\sigma^\infty$ injects into $J_{I,I_0}'(\tau^{\infty})$ and $\cJ(\sigma^{\infty})\subseteq \cJ(J_{I,I_0}'(\tau^{\infty}))$, which implies $\chi\in \cJ(J_{I,I_0}'(\tau^{\infty}))$ since $\chi\in \cJ(\sigma^{\infty})$. But it follows from (\ref{equ: explicit twist Jacquet}) and Lemma~\ref{lem: PS geometric lemma}, and the exactness of $J_{I,I_0}'$ and $J_{w_{I,I_0}(I_0),\emptyset}$, that $\sigma_2^{\infty}$ is the only Jordan-H\"older factor of $i_{\emptyset,I}^{\infty}(\chi)$ satisfying $\chi\in \cJ(J_{I,I_0}'(\sigma_2^{\infty}))$. Hence we must have $\sigma_2^{\infty}=\tau^{\infty}$.
\end{proof}

\begin{cor}\label{cor: basic reducible}
Let $I\subseteq \Delta$ and $\pi^{\infty}$ in $\mathrm{Rep}^{\infty}_{\rm{adm}}(L_I)$ $G$-basic and \emph{reducible}. Then there exists $G$-basic $\pi_i^{\infty}$ in $\mathrm{Rep}^{\infty}_{\rm{adm}}(L_I)$ for $i=0,1$ such that $\pi^{\infty}$ fits into $0 \rightarrow \pi_1^{\infty}\rightarrow \pi^{\infty} \rightarrow \pi_0^{\infty} \rightarrow 0$ which is non-split.
\end{cor}
\begin{proof}
We continue to use the notation from the proof of Corollary~\ref{cor: basic subquotient} and write $\sigma_1^{\infty}$ (resp.~$\sigma_2^{\infty}$) for the socle (resp.~cosocle) of $\pi^{\infty}$, with $\vec{\Gamma_{\Sigma,i}}= \omega^{-1}(\sigma_i^{\infty})$ for $i=1,2$. Since $\pi^{\infty}$ is reducible, we have $\sigma_1^{\infty}\neq\sigma_2^{\infty}$ and thus $\vec{\Gamma_{\Sigma,1}}\neq \vec{\Gamma_{\Sigma,2}}$. We fix an arbitrary edge $e$ of $\Gamma_{\Sigma}$ on which $\vec{\Gamma_{\Sigma,1}}$ and $\vec{\Gamma_{\Sigma,2}}$ have opposite direction. Let $\tau_i^{\infty}$ be the Jordan-H\"older factor of $\pi^{\infty}$ whose attached orientation on $\Gamma_{\Sigma}$ differs from $\vec{\Gamma_{\Sigma,i}}$ by changing only the direction of the fixed edge $e$. We define $\pi_1^{\infty}$ (resp.~$\pi_0^{\infty}$) as the unique sub (resp.~quotient) of $\pi^{\infty}$ with socle $\sigma_1^{\infty}$ and cosocle $\tau_2^{\infty}$ (resp.~with socle $\tau_1^{\infty}$ and cosocle $\sigma_2^{\infty}$). By applying Theorem~\ref{thm: structure of PS} to $i_{\emptyset,I}^{\infty}(\chi)$ for some $\chi\in\cJ(\sigma_1^{\infty})$ (with $\pi^{\infty}$ being the unique sub of $i_{\emptyset,I}^{\infty}(\chi)$ with cosocle $\sigma_2^{\infty}$ as in the proof of Corollary~\ref{cor: basic subquotient}), we see that $\mathrm{JH}_{L_I}(\pi_1^{\infty})$ (resp.~$\mathrm{JH}_{L_I}(\pi_0^{\infty})$) consists of exactly those Jordan-H\"older factors whose attached orientation on $\Gamma_{\Sigma}$ have the same direction as $\vec{\Gamma_{\Sigma,1}}$ (resp.~$\vec{\Gamma_{\Sigma,2}}$) on the fixed edge $e$, and in particular, $\pi^{\infty}$ fits into a non-split short exact sequence $0 \rightarrow \pi_1^{\infty}\rightarrow \pi^{\infty} \rightarrow \pi_0^{\infty} \rightarrow 0$.
\end{proof}

\begin{lem}\label{lem: basic as induction}
Let $I\subseteq \Delta$ and $\pi^{\infty}$ in $\mathrm{Rep}^{\infty}_{\rm{adm}}(L_I)$ $G$-basic. Then there exist $I_1\subseteq I_0\subseteq I$ and a smooth character $\delta: L_{I_0}\rightarrow E^\times$ such that $V_{I_1,I_0}^{\infty}\otimes_E \delta$ is irreducible $G$-regular and $\pi^{\infty}\cong i_{I_0,I}^{\infty}(V_{I_1,I_0}^{\infty}\otimes_E \delta)$.
\end{lem}
\begin{proof}
By the definition of $G$-basic (\ref{it: basic 1} of Definition \ref{def: basic rep}), we can assume $\pi^{\infty}$ irreducible $G$-regular. Recall that $\Sigma\defeq W(L_I)\cdot\cJ(\pi^{\infty})$ is a single left $W(L_I)$-coset (\ref{it: PS 1} of Lemma~\ref{lem: Jacquet of PS}). As in the proof of Corollary \ref{cor: basic subquotient}, we fix an arbitrary $G$-regular $\chi\in \Sigma$ such that two vertices $\rho_k$, $\rho_{k'}$ of the graph $\Gamma_{\Sigma}$ that are connected by an edge are adjacent, i.e.~such that $|k-k'|=1$. Then there exists $I_0\subseteq I$ such that there is a bijection between the Levi blocks of $L_{I_0}$ and the connected components of $\Gamma_{\Sigma}$, hence also a bijection between the set of orientations on $\Gamma_{\Sigma_0}$ ($\defeq$ the graph of $\Sigma_0\defeq W(L_{I_0})\cdot\chi$) and on $\Gamma_{\Sigma}$. Using Theorem~\ref{thm: JH of PS} applied to both $\Gamma_{\Sigma}$ and $\Gamma_{\Sigma_0}$, we deduce a bijection between the set of constituents of $i_{\emptyset,I}^{\infty}(\chi)$ and of $i_{\emptyset,I_0}^{\infty}(\chi)$. In particular they have same length, and the exactness of $i_{I_0,I}^{\infty}(-)$ then implies that $\sigma_0^{\infty}\mapsto i_{I_0,I}^{\infty}(\sigma_0^{\infty})$ induces a bijection between both sets of constituents. Using \ref{it: PS 2} of Lemma \ref{lem: Jacquet of PS}, we see that $\pi^{\infty}$ is a constituent of $i_{\emptyset,I}^{\infty}(\chi)$, hence $\pi^{\infty}\cong i_{I_0,I}^{\infty}(\pi_0^{\infty})$ for a (unique) constituent $\pi_0^{\infty}$ of $i_{\emptyset,I_0}^{\infty}(\chi)$. Applying Corollary~\ref{cor: twist of St} to $\pi_0^{\infty}$, we deduce $\pi_0^{\infty}\cong V_{I_1,I_0}^{\infty}\otimes_E \delta$ for a smooth character $\delta: L_{I_0}\rightarrow E^\times$ and $I_1\subseteq I_0$, and thus $\pi^{\infty}\cong i_{I_0,I}^{\infty}(V_{I_1,I_0}^{\infty}\otimes_E \delta)$.
\end{proof}

\begin{lem}\label{lem: Jacquet basic}
Let $I_0\subseteq I\subseteq \Delta$ and $\pi^{\infty}$ in $\mathrm{Rep}^{\infty}_{\rm{adm}}(L_I)$ $G$-basic. Then we have a canonical decomposition in $\mathrm{Rep}^{\infty}_{\rm{adm}}(L_{I_0})$ induced by (\ref{equ: sm block decomposition})
\begin{equation}\label{equ: general Jacquet block}
J_{I,I_0}(\pi^{\infty})\cong \bigoplus_{\Sigma}J_{I,I_0}(\pi^{\infty})_{\cB^{I_0}_{\Sigma}}
\end{equation}
where $\Sigma$ runs through left $W(L_{I_0})$-cosets in $W(L_I)\cdot \cJ(\pi^{\infty})$ and $J_{I,I_0}(\pi^{\infty})_{\cB^{I_0}_{\Sigma}}$ is $G$-basic if non-zero. Moreover, if $\pi^{\infty}$ is simple, then $J_{I,I_0}(\pi^{\infty})_{\cB^{I_0}_{\Sigma}}$ is simple if non-zero.
\end{lem}
\begin{proof}
Note first that we can indeed apply (\ref{equ: sm block decomposition}) to $J_{I,I_0}(\pi^{\infty})$ since $\cJ(J_{I,I_0}(\pi^{\infty}))=\cJ(\pi^{\infty})$ consists of $G$-regular weights by the last assertion in \ref{it: PS 1} of Lemma \ref{lem: Jacquet of PS}. By \ref{it: basic as image} of Remark~\ref{rem: basic PS intertwine} there exist $G$-regular $\chi,\chi'\in\widehat{T}^{\infty}$ such that $\pi^{\infty}$ is the image of the unique (up to scalar) non-zero map $i_{\emptyset,I}^{\infty}(\chi')\rightarrow i_{\emptyset,I}^{\infty}(\chi)$. By the exactness of $J_{I,I_0}(-)$ we deduce that $J_{I,I_0}(\pi^{\infty})$ is the image of the induced map $J_{I,I_0}(i_{\emptyset,I}^{\infty}(\chi'))\rightarrow J_{I,I_0}(i_{\emptyset,I}^{\infty}(\chi))$, which by Remark \ref{rem: general PS Jacquet} is the same as the map
\begin{equation*}
\bigoplus_{w'\in W^{\emptyset,I_0}(L_I)}i_{\emptyset,I_0}^{\infty}((w')^{-1}\cdot\chi') \longrightarrow \bigoplus_{w\in W^{\emptyset,I_0}(L_I)}i_{\emptyset,I_0}^{\infty}(w^{-1}\cdot\chi).
\end{equation*}
Consequently, for any left $W(L_{I_0})$-coset $\Sigma$ in $W(L_I)\cdot \cJ(\pi^{\infty})$, $J_{I,I_0}(\pi^{\infty})_{\cB^{I_0}_{\Sigma}}$ is necessarily the image of a (possibly zero) map
\begin{equation}\label{equ: component map PS}
i_{\emptyset,I_0}^{\infty}((w')^{-1}\cdot\chi')\longrightarrow i_{\emptyset,I_0}^{\infty}(w^{-1}\cdot\chi)
\end{equation}
for some $w,w'\in W^{\emptyset,I_0}(L_I)$ such that $\Sigma=W(L_{I_0})(w')^{-1}\cdot\chi'=W(L_{I_0})w^{-1}\cdot\chi$. But the image of (\ref{equ: component map PS}) is either zero or has simple socle and cosocle (since both representations in (\ref{equ: component map PS}) have simple socle and cosocle by the last assertion in \ref{it: PS socle cosocle} of Remark~\ref{rem: basic PS intertwine}), hence is either zero or $G$-basic by Corollary~\ref{cor: basic subquotient}. It follows that $J_{I,I_0}(\pi^{\infty})_{\cB^{I_0}_{\Sigma}}$ is either zero or $G$-basic.

Now we assume that $\pi^{\infty}$ is simple and prove that $J_{I,I_0}(\pi^{\infty})_{\cB^{I_0}_{\Sigma}}$ is either zero or simple. We can assume $\sigma_0^{\infty}\defeq J_{I,I_0}(\pi^{\infty})_{\cB^{I_0}_{\Sigma}}\neq 0$, and we write $\sigma^{\infty}$ for its cosocle. The natural surjections $J_{I,I_0}(\pi^{\infty})\twoheadrightarrow \sigma_0^{\infty}\twoheadrightarrow \sigma^{\infty}$ induce by (\ref{equ: first adjunction}) (for $k=0$) maps $\pi^{\infty}\rightarrow i_{I,I_0}^{\infty}(\sigma_0^{\infty}) \twoheadrightarrow i_{I,I_0}^{\infty}(\sigma^{\infty})$ the composition of which is non-zero. Hence $\pi^{\infty}$ (which is simple) appears in the socle of both $i_{I,I_0}^{\infty}(\sigma_0^{\infty})$ and $i_{I,I_0}^{\infty}(\sigma^{\infty})$. As $\sigma_0^{\infty}$ is $G$-basic, so is $i_{I,I_0}^{\infty}(\sigma_0^{\infty})$. Moreover $i_{I,I_0}^{\infty}(\sigma^{\infty})$ is clearly $G$-basic. Hence by the end of \ref{it: basic as image} of Remark~\ref{rem: basic PS intertwine} both $i_{I,I_0}^{\infty}(\sigma_0^{\infty})$ and $i_{I,I_0}^{\infty}(\sigma^{\infty})$ have simple socle and cosocle, and thus have same socle $\pi^{\infty}$. As $i_{I,I_0}^{\infty}(\sigma^{\infty})$ is a quotient of $i_{I,I_0}^{\infty}(\sigma_0^{\infty})$ this forces $i_{I,I_0}^{\infty}(\sigma^{\infty})=i_{I,I_0}^{\infty}(\sigma_0^{\infty})$ and thus $\sigma^{\infty}=\sigma_0^{\infty}$. In particular $\sigma_0^{\infty}=J_{I,I_0}(\pi^{\infty})_{\cB^{I_0}_{\Sigma}}$ is simple.
\end{proof}

\begin{rem}\label{rem: twist Jacquet basic}
Let $I_0\subseteq I\subseteq \Delta$ and $\pi^{\infty}$ in $\mathrm{Rep}^{\infty}_{\rm{adm}}(L_I)$ $G$-basic. As $J_{I,w_{I,I_0}(I_0)}(\pi^{\infty})$ is a direct sum of $G$-basic objects in distinct Bernstein blocks by Lemma~\ref{lem: Jacquet basic}, we deduce from (\ref{equ: explicit twist Jacquet}), (\ref{equ: twist Jacquet support}) and Lemma~\ref{lem: G basic twist} that an analogous statement as Lemma \ref{lem: Jacquet basic} holds for $J_{I,I_0}'(\pi^{\infty})$. Similarly, for $I_0,I_1\subseteq \Delta$, $w\in W^{I_0,I_1}$ and $\pi_0^{\infty}$ in $\mathrm{Rep}^{\infty}_{\rm{adm}}(L_{I_0})$ $G$-basic, as $J_{I_0,w(I_1)\cap I_0}(\pi_0^{\infty})$ is a direct sum of $G$-basic objects in distinct Bernstein blocks by Lemma~\ref{lem: Jacquet basic}, we deduce from (\ref{equ: Bruhat Jacquet}), (\ref{equ: Jacquet twist}), (\ref{equ: twist Borel}) and Lemma~\ref{lem: G basic twist} that the statement of Lemma \ref{lem: Jacquet basic} holds for $J_{I_0,I_1,w}(\pi_0^{\infty})$. All these statements remain true if we replace everywhere $G$-basic by irreducible $G$-regular.
\end{rem}

\begin{lem}\label{lem: structure of induction}
Let $I_1\subseteq I\subseteq \Delta$, $w\in W(L_I)$ and $\pi^{\infty}\defeq \mathrm{cosoc}_{L_{I_1}}(i_{\emptyset,I_1}^{\infty}(w\cdot 1_{T}))$. Then $i_{I_1,I}^{\infty}(\pi^{\infty})$ is isomorphic to the unique quotient of $i_{\emptyset,I}^{\infty}(w\cdot 1_{T})$ with socle $V_{I_2,I}^{\infty}$ (and cosocle $V_{D_R(w),I}^{\infty}$) where $I_2\defeq I\setminus D_R(w_{I_1}w)$.
\end{lem}
\begin{proof}
Note first that $i_{I_1,I}^{\infty}(\pi^{\infty})$ is clearly a quotient of $i_{\emptyset,I}^{\infty}(w\cdot 1_{T})$. As $i_{I_1,I}^{\infty}(\pi^{\infty})$ is $G$-basic (recall $\pi^{\infty}$ is irreducible by \ref{it: PS socle cosocle} of Remark~\ref{rem: basic PS intertwine}), it has simple socle by \ref{it: basic as image} of Remark~\ref{rem: basic PS intertwine}. By (\ref{equ: JH trivial PS}) this socle has the form $V_{I_2,I}^{\infty}$ for some $I_2\subseteq I$. It follows from \ref{it: basic as quotient} of Remark~\ref{rem: basic PS intertwine} that $\pi^{\infty}\cong \mathrm{soc}_{L_{I_1}}(i_{\emptyset,I_1}^{\infty}(w_{I_1}w\cdot 1_{T}))$, and hence $w_{I_1}w\cdot 1_{T}\in \cJ(\pi^{\infty})$ ((\ref{equ: first adjunction}) with $k=0$). Since by \emph{loc.~cit.}
\[\Hom_{L_{I_1}}(J_{I,I_1}(V_{I_2,I}^{\infty}),\pi^{\infty})\cong \Hom_{L_I}(V_{I_2,I}^{\infty}, i_{I_1,I}^{\infty}(\pi^{\infty})) \ne 0\]
we see that $\cJ(\pi^{\infty})\subseteq \cJ(J_{I,I_1}(V_{I_2,I}^{\infty}))=\cJ(V_{I_2,I}^{\infty})$ as $\pi^{\infty}$ is irreducible and thus $w_{I_1}w\cdot 1_{T}\in \cJ(V_{I_2,I}^{\infty})$. In particular (\ref{equ: Jacquet Steinberg}) implies $I_2=I\setminus D_R(w_{I_1}w)$.
\end{proof}

\begin{lem}\label{lem: basic subquotient of Ext}
Let $I\subseteq \Delta$ and $\pi_0^{\infty},\pi_1^{\infty}$ in $\mathrm{Rep}^{\infty}_{\rm{adm}}(L_I)$ both $G$-basic with $\mathrm{JH}_{L_I}(\pi_0^{\infty})\cap \mathrm{JH}_{L_I}(\pi_1^{\infty})=\emptyset$. Let $\pi^{\infty}$ in $\mathrm{Rep}^{\infty}_{\rm{adm}}(L_I)$ which fits into a non-split short exact sequence
\begin{equation}\label{equ: basic subquotient of Ext}
0\rightarrow \pi_1^{\infty} \rightarrow \pi^{\infty} \rightarrow \pi_0^{\infty} \rightarrow 0.
\end{equation}
Then $\pi^{\infty}$ admits a unique subquotient $\sigma^{\infty}$ which is $G$-basic with simple socle $\mathrm{soc}_{L_I}(\pi_1^{\infty})$ and simple cosocle $\mathrm{cosoc}_{L_I}(\pi_0^{\infty})$.
\end{lem}
\begin{proof}
As $\pi_i^{\infty}$ is $G$-basic, by \ref{it: basic as image} of Remark~\ref{rem: basic PS intertwine} it is multiplicity free with simple socle and cosocle. In particular, $\pi^{\infty}$ is multiplicity free and we write $\pi_2^{\infty}$ for its unique subrepresentation with cosocle $\mathrm{cosoc}_{L_I}(\pi_0^{\infty})$. In particular the composition of $\pi_2^{\infty}\rightarrow \pi^{\infty} \rightarrow \pi_0^{\infty}$ is a surjection. If $\pi_2^{\infty}\cap \pi_1^{\infty}=0$ inside $\pi^{\infty}$, then (\ref{equ: basic subquotient of Ext}) induces an isomorphism $\pi_2^{\infty}\buildrel\sim\over\rightarrow \pi_0^{\infty}$ and thus (\ref{equ: basic subquotient of Ext}) splits, a contradiction to our assumption. If $\pi_2^{\infty}\cap \pi_1^{\infty}\neq 0$, as $\pi_1^{\infty}$ has simple socle, we must have $\mathrm{soc}_{L_I}(\pi_1^{\infty})\subseteq \pi_2^{\infty}$, so we can define $\sigma^{\infty}$ as the unique quotient of $\pi_2^{\infty}$ with socle $\mathrm{soc}_{L_I}(\pi_1^{\infty})$. Since $\pi_0^{\infty}$ and $\pi_1^{\infty}$ are $G$-basic with $\mathrm{JH}_{L_I}(\pi_0^{\infty})\cap \mathrm{JH}_{L_I}(\pi_1^{\infty})=\emptyset$, by \ref{it: PS 1} and \ref{it: PS 2} of Lemma~\ref{lem: Jacquet of PS} we know that $J_{I,\emptyset}(\pi_0^{\infty})$ and $J_{I,\emptyset}(\pi_1^{\infty})$ are multiplicity free, semi-simple and share no common constituent. Since we have a short exact sequence $0\rightarrow J_{I,\emptyset}(\pi_1^{\infty}) \rightarrow J_{I,\emptyset}(\pi^{\infty}) \rightarrow J_{I,\emptyset}(\pi_0^{\infty}) \rightarrow 0$ (as $J_{I,\emptyset}$ is exact), it follows that $J_{I,\emptyset}(\pi^{\infty})$ is also semi-simple and multiplicity free, and thus so is $J_{I,\emptyset}(\sigma^{\infty})$. Now we choose an arbitrary $\chi\in\cJ(\mathrm{soc}_{L_I}(\pi_1^{\infty}))\subseteq \cJ(\sigma^{\infty})$, which by (\ref{equ: first adjunction}) (with $k=0$) gives a non-zero map $\sigma^{\infty}\rightarrow i_{\emptyset,I}^{\infty}(\chi)$. Since $i_{\emptyset,I}^{\infty}(\chi)$ is multiplicity free with socle $\mathrm{soc}_{L_I}(\pi_1^{\infty})\cong \mathrm{soc}_{L_I}(\sigma^{\infty})$ by \ref{it: PS 2} of Lemma~\ref{lem: Jacquet of PS} and \ref{it: PS socle cosocle} of Remark~\ref{rem: basic PS intertwine}, the map $\sigma^{\infty}\rightarrow i_{\emptyset,I}^{\infty}(\chi)$ is an injection. Since $\sigma^{\infty}$ has simple socle and cosocle by definition, Corollary~\ref{cor: basic subquotient} finally implies that $\sigma^{\infty}$ is $G$-basic.
\end{proof}

\subsection{Results on smooth Ext groups}

We prove several useful results on smooth $\mathrm{Ext}$ groups of $G$-basic representations.\bigskip

We start with some preliminaries. For $I_0,I_1\subseteq \Delta$, we define
\begin{equation}\label{di0i1}
d(I_0,I_1)\defeq \#(I_0\setminus I_1)+\#(I_1\setminus I_0)
\end{equation}
and $[I_0,I_1]\defeq \{I\subseteq \Delta\mid d(I_0,I_1)=d(I_0,I)+d(I,I_1)\}$. One easily checks $d(I_0,I_1)=d(I_0,I)+d(I,I_1)$ if and only if
\[(I_0\setminus I_1)\sqcup (I_1\setminus I_0)=(I_0\setminus I)\sqcup (I\setminus I_1)\sqcup (I_1\setminus I)\sqcup (I\setminus I_0),\]
if and only if $(I_0\setminus I)\cap (I_1\setminus I)=\emptyset$ and $(I\setminus I_0)\cap (I\setminus I_1)=\emptyset$, if and only if $I_0\cap I_1\subseteq I\subseteq I_0\cup I_1$. In other words, we have
\begin{equation}\label{equ: same interval}
[I_0,I_1]=\{I\subseteq \Delta\mid I_0\cap I_1\subseteq I\subseteq I_0\cup I_1\}.
\end{equation}
Given another $I_1'\subseteq \Delta$, it is clear that we have
\begin{equation}\label{equ: intersection of interval}
[I_0,I_1]\cap [I_0,I_1']=[I_0,I_1'']
\end{equation}
with $I_1''\defeq (I_0\setminus (I_1\cup I_1'))\sqcup ((I_1\cap I_1')\setminus I_0)$.

\begin{lem}\label{lem: Hom cube}
Let $I\subseteq \Delta$.
\begin{enumerate}[label=(\roman*)]
\item \label{it: Hom cube 1} For $I_0,I_1\subseteq I$, there exists a unique $G$-basic $Q_I(I_0,I_1)\in\mathrm{Rep}^{\infty}_{\rm{adm}}(L_I)$ with socle $V_{I_0,I}^{\infty}$ and cosocle $V_{I_1,I}^{\infty}$, and it has set of Jordan-H\"older factors $\{V_{I',I}^{\infty}\mid I'\in [I_0,I_1]\}$.
\item \label{it: Hom cube 2} For $I_0,I_1,I_0',I_1'\subseteq I$, there exists a non-zero map $Q_I(I_0',I_1')\rightarrow Q_I(I_0,I_1)$ if and only if $I_1'\in[I_0,I_1]$ and $I_0\in[I_0',I_1']$, in which case the map is unique (up to scalar) with image isomorphic to $Q_I(I_0,I_1')$.
\item \label{it: Hom cube 3} For $I_0,I_1\subseteq I$ and $I_0',I_1'\in [I_0,I_1]$, $Q_I(I_0',I_1')$ is a subquotient of $Q_I(I_0,I_1)$ if and only if $I_0'\in[I_0,I_1']$ if and only if $I_1'\in[I_0',I_1]$.
\end{enumerate}
\end{lem}
\begin{proof}
We prove \ref{it: Hom cube 1}. By \ref{it: PS socle cosocle} and \ref{it: basic as image} of Remark~\ref{rem: basic PS intertwine} any $G$-basic representation of $L_I$ with socle $V_{I_0,I}^{\infty}$ and cosocle $V_{I_1,I}^{\infty}$ must be the image of the unique (up to scalar) non-zero map $i_{\emptyset,I}^{\infty}(\chi')\rightarrow i_{\emptyset,I}^{\infty}(\chi)$ for any $\chi\in \cJ(V_{I_0,I}^{\infty})$ and $\chi'\in w_I\cdot\cJ(V_{I_1,I}^{\infty})$. This implies the unicity of such a representation. For its existence, note that by \ref{it: PS intertwine} of Remark~\ref{rem: basic PS intertwine} and (\ref{equ: Jacquet Steinberg}) there is a non-zero map $i_{\emptyset,I}^{\infty}(\chi')\rightarrow i_{\emptyset,I}^{\infty}(\chi)$ for any $\chi\in \cJ(V_{I_0,I}^{\infty})$ and $\chi'\in w_I\cdot\cJ(V_{I_1,I}^{\infty})$. Its image, which has socle $V_{I_0,I}^{\infty}$ and cosocle $V_{I_1,I}^{\infty}$, is $G$-basic by Corollary~\ref{cor: basic subquotient}. For the last statement of \ref{it: Hom cube 1}, choose $w\in W(L_I)$ such that $I_0=I\setminus D_R(w)$. One can check that $Q_I(I_0,I\setminus I_0)\cong i_{\emptyset,I}^{\infty}(w\cdot 1_T)$ and (using \ref{it: trivial PS 2} of Corollary \ref{cor: trivial PS}) that $Q_I(I_0,I_1)$ is the unique subrepresentation of $i_{\emptyset,I}^{\infty}(w\cdot 1_T)$ with cosocle $V_{I_1,I}^{\infty}$. By \emph{loc.~cit.}~the constituents of $Q_I(I_0,I_1)$ are the $V_{I',I}^{\infty}$ such that $V_{I',I}^{\infty}\leq_w V_{I_1,I}^{\infty}$ (where $\leq_w$ is the partial order defined by $i_{\emptyset,I}^{\infty}(w\cdot 1_T)$). The last statement in \ref{it: Hom cube 1} follows then from \ref{it: trivial PS 3} of Corollary \ref{cor: trivial PS} and (\ref{equ: same interval}).

We prove \ref{it: Hom cube 2}. By \ref{it: PS socle cosocle} and \ref{it: basic as image} of Remark~\ref{rem: basic PS intertwine} again, we can choose $\chi\in \cJ(V_{I_0,I}^{\infty})$ and $\chi'\in w_I\cdot\cJ(V_{I_1',I}^{\infty})$ such that $Q_I(I_0,I_1)\hookrightarrow i_{\emptyset,I}^{\infty}(\chi)$ and $i_{\emptyset,I}^{\infty}(\chi')\twoheadrightarrow Q_I(I_0',I_1')$. Assume that there is a non-zero map $Q_I(I_0',I_1')\rightarrow Q_I(I_0,I_1)$, then the composition
\begin{equation*}
i_{\emptyset,I}^{\infty}(\chi')\twoheadrightarrow Q_I(I_0',I_1')\rightarrow Q_I(I_0,I_1)\hookrightarrow i_{\emptyset,I}^{\infty}(\chi)
\end{equation*}
is also non-zero and has same image. By the existence part in \ref{it: Hom cube 1} there is a unique (up to scalar) non-zero map $i_{\emptyset,I}^{\infty}(\chi')\rightarrow i_{\emptyset,I}^{\infty}(\chi)$ with image $Q_I(I_0,I_1')$. Hence $Q_I(I_0,I_1')$ must be a quotient of $Q_I(I_0',I_1')$ and a subrepresentation of $Q_I(I_0,I_1)$, forcing $I_1'\in[I_0,I_1]$ and $I_0\in[I_0',I_1']$ by the last statement in \ref{it: Hom cube 1}. Conversely, if $I_1'\in[I_0,I_1]$, then $V_{I_1',I}^{\infty}$ shows up in $Q_I(I_0,I_1)$ by \ref{it: Hom cube 1}, hence $Q_I(I_0,I_1')$ is a subrepresentation of $Q_I(I_0,I_1)$ (by unicity of $Q_I(I_0,I_1')$). Similarly, $I_0\in[I_0',I_1']$ implies that $Q_I(I_0,I_1')$ is a quotient of $Q_I(I_0',I_1')$. So if $I_1'\in[I_0,I_1]$ and $I_0\in[I_0',I_1']$ there is a non-zero map $Q_I(I_0',I_1')\rightarrow Q_I(I_0,I_1)$ (with image $Q_I(I_0,I_1')$).

Finally, \ref{it: Hom cube 3} follows from \ref{it: Hom cube 2} and the observation that $Q_I(I_0',I_1')$ for $I_0',I_1'\in [I_0,I_1]$ is a subquotient of $Q_I(I_0,I_1)$ if and only if it is a subrepresentation of $Q_I(I_0',I_1)$ if and only it is a quotient of $Q_I(I_0,I_1')$.
\end{proof}

\begin{rem}\label{rem: dual cube}
Let $I\subseteq \Delta$ and $w\in W(L_I)$. It follows from \cite[Thm.~III.2.7]{Re10} (with $W$ there being $(w\cdot 1_{T})^{\sim}\cong \delta_{B_I}^{1/2}\cdot w(\delta_{B_I})^{-1/2}$ and recalling that $(-)^\sim$ is the smooth contragredient):
\begin{multline*}
i_{\emptyset,I}^{\infty}(w\cdot 1_{T})\cong i_{\emptyset,I}^{\infty}((w\cdot 1_{T})^{\sim}\cdot \delta_{B_I}^{-1})^{\sim}\cong i_{\emptyset,I}^{\infty}(w(\delta_{B_I})^{-1/2}\cdot \delta_{B_I}^{-1/2})\\
\cong i_{\emptyset,I}^{\infty}(w(\delta_{B_I}^{-1/2})\cdot \delta_{B_I}^{-1/2})\cong i_{\emptyset,I}^{\infty}(ww_I(\delta_{B_I}^{1/2})\cdot \delta_{B_I}^{-1/2})=i_{\emptyset,I}^{\infty}(ww_I\cdot 1_{T}).
\end{multline*}
Together with (the proof of) \ref{it: Hom cube 1} of Lemma~\ref{lem: Hom cube} this implies for $I_0,I_1\subseteq I$
\begin{equation}\label{equ: dual cube}
Q_I(I_0,I_1)^{\sim}\cong Q_I(I_1,I_0).
\end{equation}
By \cite[Prop.~5]{Vig97} and (\ref{equ: dual cube}) we have for $k\geq 0$
\begin{multline}\label{equ: dual sm Ext}
\mathrm{Ext}_{L_I}^k(Q_I(I_0,I_1),Q_I(I_0',I_1'))^\infty\cong \mathrm{Ext}_{L_I}^k(Q_I(I_0',I_1')^{\sim}, Q_I(I_0,I_1)^{\sim})^\infty\\
\cong \mathrm{Ext}_{L_I}^k(Q_I(I_1',I_0'), Q_I(I_1,I_0))^\infty
\end{multline}
and in particular for $I_0,I_1,I_0',I_1'\subseteq I$
\begin{equation*}
d_I(Q_I(I_0,I_1),Q_I(I_0',I_1'))=d_I(Q_I(I_1',I_0'), Q_I(I_1,I_0)).
\end{equation*}
\end{rem}

We now recall the following classical result (see \ref{it: basic 3} of Definition \ref{def: basic rep} for $d(-,-)$).

\begin{lem}\label{lem: Ext sm St}
Let $I_0,I_1\subseteq I\subseteq \Delta$. Then we have $d(V_{I_0,I}^{\infty},V_{I_1,I}^{\infty})=d(I_0,I_1)$ and
\[\Dim_E\mathrm{Ext}_{L_I}^{d(I_0,I_1)}(V_{I_0,I}^{\infty},V_{I_1,I}^{\infty})^{\infty}=1.\]
\end{lem}
\begin{proof}
This follows directly from \cite[Cor.~2]{Or05} or \cite[Thm.~1.3]{Dat06}.
\end{proof}

\begin{lem}\label{lem: dim one}
Let $I_i\subseteq I\subseteq \Delta$ and $\pi_{i}^{\infty}$ in $\mathrm{Rep}^{\infty}_{\rm{adm}}(L_{I_i})$ $G$-basic for $i=0,1$. Assume that $d_I(\pi_0^{\infty},\pi_1^{\infty})<\infty$. Then $d_{I'}(\pi_0^{\infty},\pi_1^{\infty})=d_I(\pi_0^{\infty},\pi_1^{\infty})$ for $I'\supseteq I$ and
\begin{equation}\label{equ: dim one}
\mathrm{Ext}_{L_I}^{d_I(\pi_0^{\infty},\pi_1^{\infty})}(i_{I_0,I}^{\infty}(\pi_0^{\infty}),i_{I_1,I}^{\infty}(\pi_1^{\infty}))^{\infty}
\end{equation}
is one dimensional.
\end{lem}
\begin{proof}
Let \ $\Sigma\defeq W(L_I)\cdot\cJ(\pi_1^{\infty})=W(L_I)\cdot\cJ(i_{I_1,I}^{\infty}(\pi_1^{\infty}))$ \ (see \ Lemma \ \ref{lem: support of induction}), \ then $d_I(\pi_0^{\infty},\pi_1^{\infty})\!<\infty$ implies in particular
\begin{equation*}
i_{I_0,I}^{\infty}(\pi_0^{\infty})\in \cB^{I}_{\Sigma}.
\end{equation*}
Consequently, for $I'\supseteq I$, we have by the last assertion in \ref{it: sml2} of Lemma \ref{lem: smooth geometric lemma} applied with $\pi_0^{\infty}$ there being $i_{I_0,I}^{\infty}(\pi_0^{\infty})$ and with $w=1$
\[J_{I',I}(i_{I,I'}^{\infty}(i_{I_0,I}^{\infty}(\pi_0^{\infty})))_{\cB^{I}_{\Sigma}}\cong i_{I,I,1}^{\infty}(J_{I,I,1}(i_{I_0,I}^{\infty}(\pi_0^{\infty}))) \cong i_{I_0,I}^{\infty}(\pi_0^{\infty}).\]
Together with (\ref{equ: first adjunction}) this implies canonical isomorphisms for $k\geq 0$
\begin{multline*}
\mathrm{Ext}_{L_{I'}}^{k}(i_{I_0,I'}^{\infty}(\pi_0^{\infty}),i_{I_1,I'}^{\infty}(\pi_1^{\infty}))^{\infty}\cong \mathrm{Ext}_{L_{I}}^{k}(J_{I',I}(i_{I_0,I'}^{\infty}(\pi_0^{\infty})), i_{I_1,I}^{\infty}(\pi_1^{\infty}))^\infty\\
\cong \mathrm{Ext}_{L_{I}}^{k}(i_{I_0,I}^{\infty}(\pi_0^{\infty}),i_{I_1,I}^{\infty}(\pi_1^{\infty}))^{\infty},
\end{multline*}
and thus $d_{I'}(\pi_0^{\infty},\pi_1^{\infty})=d_I(\pi_0^{\infty},\pi_1^{\infty})$ for $I'\supseteq I$.

Now we prove the second assertion by induction on $I$. We assume inductively that
\[\Dim_E\mathrm{Ext}_{L_{I'}}^{d_{I'}(\sigma_0^{\infty},\sigma_1^{\infty})}(i_{I_0',I'}^{\infty}(\sigma_0^{\infty}),i_{I_1',I'}^{\infty}(\sigma_1^{\infty}))^{\infty}=1\]
for any $I_0',I_1'\subseteq I'\subsetneq I$ and $G$-basic $\sigma_{i}^{\infty}$ in $\mathrm{Rep}^{\infty}_{\rm{adm}}(L_{I_i'})$, $i=0,1$, such that $d_{I'}(\sigma_0^{\infty},\sigma_1^{\infty})<\infty$. Note that the induction hypothesis trivially holds when $I=I_0=I_1=\emptyset$ since we only have smooth characters then.\bigskip

\textbf{Case $1$}: If there exist $I_2\subsetneq I$ and an irreducible $G$-regular $\pi_2^{\infty}$ in $\mathrm{Rep}^{\infty}_{\rm{adm}}(L_{I_2})$ such that $i_{I_1,I}^{\infty}(\pi_1^{\infty})\cong i_{I_2,I}^{\infty}(\pi_2^{\infty})$, then by (\ref{equ: first adjunction}) we have isomorphisms for $k\geq 0$
\begin{multline}\label{isoajoutes}
\mathrm{Ext}_{L_I}^{k}(i_{I_0,I}^{\infty}(\pi_0^{\infty}),i_{I_1,I}^{\infty}(\pi_1^{\infty}))^{\infty}\cong
\mathrm{Ext}_{L_I}^{k}(i_{I_0,I}^{\infty}(\pi_0^{\infty}),i_{I_2,I}^{\infty}(\pi_2^{\infty}))^{\infty}\\
\cong \mathrm{Ext}_{L_{I_2}}^{k}(J_{I,I_2}(i_{I_0,I}^{\infty}(\pi_0^{\infty})),\pi_2^{\infty})^{\infty}
\cong \mathrm{Ext}_{L_{I_2}}^{k}(\pi_3^{\infty},\pi_2^{\infty})^{\infty}
\end{multline}
where $\Sigma\defeq W(L_{I_2})\cdot\cJ(\pi_2^{\infty})$ and $\pi_3^{\infty}\defeq J_{I,I_2}(i_{I_0,I}^{\infty}(\pi_0^{\infty}))_{\cB^{I_2}_{\Sigma}}$. By Lemma~\ref{lem: Jacquet basic} $\pi_3^{\infty}$ is either zero or $G$-basic. The assumption $d_I(\pi_0^{\infty},\pi_1^{\infty})<\infty$ forces $\pi_3^{\infty}\neq 0$, and thus $\pi_3^{\infty}$ is $G$-basic with $d_{I_2}(\pi_3^{\infty},\pi_2^{\infty})=d_I(\pi_0^{\infty},\pi_1^{\infty})$ by (\ref{isoajoutes}). Hence, (\ref{equ: dim one}) is isomorphic to $\mathrm{Ext}_{L_{I_2}}^{d_{I_2}(\pi_3^{\infty},\pi_2^{\infty})}(\pi_3^{\infty},\pi_2^{\infty})^{\infty}$ which is one dimensional by our induction hypothesis as $I_2\subsetneq I$.\bigskip

\textbf{Case $2$}: If there exist $I_2\subsetneq I$ and an irreducible $G$-regular $\pi_2^{\infty}$ in $\mathrm{Rep}^{\infty}_{\rm{adm}}(L_{I_2})$ such that $i_{I_0,I}^{\infty}(\pi_0^{\infty})\cong i_{I_2,I}^{\infty}(\pi_2^{\infty})$, then by (\ref{equ: second adjunction}) we have isomorphisms for $k\geq 0$
\begin{multline*}
\mathrm{Ext}_{L_I}^{k}(i_{I_0,I}^{\infty}(\pi_0^{\infty}),i_{I_1,I}^{\infty}(\pi_1^{\infty}))^{\infty}\cong
\mathrm{Ext}_{L_I}^{k}(i_{I_2,I}^{\infty}(\pi_2^{\infty}),i_{I_1,I}^{\infty}(\pi_1^{\infty}))^{\infty}\\
\cong \mathrm{Ext}_{L_{I_2}}^{k}(\pi_2^{\infty},J_{I,I_2}'(i_{I_1,I}^{\infty}(\pi_1^{\infty})))^{\infty}\cong
\mathrm{Ext}_{L_{I_2}}^{k}(\pi_2^{\infty},\pi_3^{\infty})^{\infty}
\end{multline*}
where $\Sigma\defeq W(L_{I_2})\cdot\cJ(\pi_2^{\infty})$ and $\pi_3^{\infty}\defeq J_{I,I_2}'(i_{I_1,I}^{\infty}(\pi_1^{\infty}))_{\cB^{I_2}_{\Sigma}}$. We deduce from Remark~\ref{rem: twist Jacquet basic} that $\pi_3^{\infty}$ is either zero or $G$-basic. We conclude by induction as in Case $1$ since $I_2\subsetneq I$.\bigskip

\textbf{Case $3$}: If we are in none of the above two cases, then we must have $I_0=I_1=I$ and, using Lemma~\ref{lem: basic as induction}, subsets $I_2,I_3\subseteq I$ and smooth characters $\delta_2,\delta_3: L_I\rightarrow E^\times$ such that $\pi_0^{\infty}\cong V_{I_2,I}^{\infty}\otimes_E\delta_2$ and $\pi_1^{\infty}\cong V_{I_3,I}^{\infty}\otimes_E\delta_3$. Our assumption $d_I(\pi_0^{\infty},\pi_1^{\infty})<\infty$ forces $\pi_0^{\infty}$ and $\pi_1^{\infty}$ to lie in the same Bernstein block and thus
\[W(L_I)\cdot\delta_2|_{T}=W(L_I)\cdot\cJ(V_{I_2,I}^{\infty}\otimes_E\delta_2)=W(L_I)\cdot\cJ(V_{I_3,I}^{\infty}\otimes_E\delta_3)=W(L_I)\cdot\delta_3|_{T},\]
which implies $\delta_2=\delta_3$ as $\delta_2|_{T}$ is the unique element in $W(L_I)\cdot\delta_2|_{T}$ which extends to a smooth character $L_I\rightarrow E^\times$. Consequently, we obtain isomorphisms $\mathrm{Ext}_{L_I}^{k}(\pi_0^{\infty},\pi_1^{\infty})^{\infty}\cong \mathrm{Ext}_{L_I}^{k}(V_{I_2,I}^{\infty},V_{I_3,I}^{\infty})$ for $k\geq 0$ which by Lemma~\ref{lem: Ext sm St} forces (\ref{equ: dim one}) to be one dimensional.
\end{proof}

Recall that $\widehat{j}=\Delta\setminus\{j\}$ for $j\in\Delta$.

\begin{lem}\label{lem: length two trivial block}
Let $I_1,I_2\subseteq \Delta$ with $d(I_1,I_2)=1$, then there exists a unique non-split extension $0\rightarrow V_{I_1,\Delta}^{\infty}\rightarrow\pi^{\infty}\rightarrow V_{I_2,\Delta}^{\infty}\rightarrow 0$ in $\mathrm{Rep}^{\infty}_{\rm{adm}}(G)$. More precisely $\pi^{\infty}\cong i_{\widehat{j},\Delta}^{\infty}(\tau^{\infty})$ (hence $\pi^{\infty}$ is $G$-basic) for $j\in \Delta$ and $\tau^{\infty}\in \mathrm{Rep}^{\infty}_{\rm{adm}}(L_{\widehat{j}})$ irreducible $G$-regular which are as follows:
\begin{enumerate}[label=(\roman*)]
\item \label{it: length two 1} if $I_2=I_1\setminus\{j_1\}$, then $j=j_1$ and $\tau^{\infty}=V_{I_2,\widehat{j}}^{\infty}$;
\item \label{it: length two 2} if $I_2=I_1\sqcup\{j_2\}$, then $j=n-j_2$ and $\tau^{\infty}=V_{w^{\widehat{j}}(I_1),\widehat{j}}^{\infty}\otimes_E\delta$ where $w^{\widehat{j}}$ is the longest element of $W^{\widehat{j},\emptyset}$ and $\delta: L_{\widehat{j}}\rightarrow E^\times$ is the unique character such that $\delta|_{T}\cong w^{\widehat{j}}\cdot 1_{T}$ (see~(\ref{equ: smooth dot action})).
\end{enumerate}
\end{lem}
\begin{proof}
As $d(I_0,I_1)=1$, by Lemma~\ref{lem: Ext sm St} there exists a unique (up to isomorphism) length two $\pi^{\infty}$ in $\mathrm{Rep}^{\infty}_{\rm{adm}}(G)$ with socle $V_{I_1,I}^{\infty}$ and cosocle $V_{I_2,I}^{\infty}$.\bigskip

\textbf{Case $1$}: If $I_2=I_1\setminus \{j_1\}$, set $j\defeq j_1$ and $\tau^{\infty}\defeq V_{I_2,\widehat{j}}^{\infty}$. It follows from \ref{it: trivial PS 2} of Corollary \ref{cor: trivial PS} that $\tau^{\infty}\cong \mathrm{cosoc}_{L_{\widehat{j}}}(i_{\emptyset,\widehat{j}}(w_1\cdot 1_{T})^{\infty})$ for any $w_1\in W(L_{\widehat{j}})$ such that $D_R(w_1)=I_1\setminus\{j\}=I_2$. Let $w_{\widehat{j}}\in W(L_{\widehat{j}})$ be the longest element, we have $D_R(w_{\widehat{j}}w_1)=\widehat{j}\setminus D_R(w_1)=\Delta\setminus I_1$, hence $\Delta\setminus D_R(w_{\widehat{j}}w_1)=I_1$. Lemma~\ref{lem: structure of induction} then implies that $i_{\widehat{j},\Delta}^{\infty}(\tau^{\infty})$ is the unique quotient of $i_{\emptyset,\Delta}^{\infty}(w_1\cdot 1_{T})$ with socle $V_{I_1,\Delta}^{\infty}$ and cosocle $V_{I_2,\Delta}^{\infty}$, which has length $2$ by \ref{it: trivial PS 3} of Corollary~\ref{cor: trivial PS} and is thus isomorphic to $\pi^{\infty}$.\bigskip

\textbf{Case $2$}: If $I_2=I_1\sqcup \{j_2\}$, set $j\defeq n-j_2$ and $w^{\widehat{j}}\defeq w_{\widehat{j}}w_0=w_0w_{\widehat{j_2}}$. Then $w^{\widehat{j}}$ is the longest element in $W^{\emptyset, \widehat{j_2}}$ and we have $w^{\widehat{j}}(\Delta\setminus \{j_2\})= \Delta\setminus \{j\}$, or equivalently $w^{\widehat{j}}L_{\widehat {j_2}}{w^{\widehat{j}}}^{-1}=L_{\widehat {j}}$ (note then that $w^{\widehat{j}}$ is also the longest element in $W^{\widehat{j},\emptyset}$). In particular $w^{\widehat{j}}(I_1)\subseteq \Delta\setminus \{j\}$ since $j_2\notin I_1$ and there exists $w_2\in W(L_{\widehat{j}})$ such that $D_R(w_2)=w^{\widehat{j}}(I_1)$. We choose any such $w_2$ and define $\tau^{\infty}\defeq \mathrm{cosoc}_{L_{\widehat{j}}}(i_{\emptyset,\widehat{j}}(w_2w^{\widehat{j}}\cdot 1_{T})^{\infty})\cong V_{w^{\widehat{j}}(I_1),\widehat{j}}^{\infty}\otimes_E\delta$ (using \ref{it: trivial PS 2} of Corollary \ref{cor: trivial PS}). By Lemma~\ref{lem: structure of induction} $i_{\widehat{j},\Delta}^{\infty}(\tau^{\infty})$ is the unique quotient of $i_{\emptyset,\Delta}^{\infty}(w_2w^{\widehat{j}}\cdot 1_{T})$ with socle $V_{\Delta\setminus D_R(w_{\widehat{j}}w_2w^{\widehat{j}}),\Delta}^{\infty}$ and cosocle $V_{D_R(w_2w^{\widehat{j}}),\Delta}^{\infty}$. Let us first compute $\Delta\setminus D_R(w_{\widehat{j}}w_2w^{\widehat{j}})$. Let $w_2'\defeq (w^{\widehat{j}})^{-1}w_2w^{\widehat{j}}\in W(L_{\widehat{j_2}})$, we have $D_R(w_2')=(w^{\widehat{j}})^{-1}(D_R(w_2))=(w^{\widehat{j}})^{-1}w^{\widehat{j}}(I_1)=I_1$ and hence $D_R(w_0w_2')=\Delta\setminus I_1$ by (\ref{equ: right set}). This implies $\Delta\setminus D_R(w_{\widehat{j}}w_2w^{\widehat{j}})=\Delta\setminus D_R(w_0w_2')=I_1$. Let us now compute $D_R(w_2w^{\widehat{j}})$. Let $j'\in \Delta\setminus\{j_2\}$, then $w^{\widehat{j}}(j')\in \Delta\setminus \{j\}$ and we have $\ell(w_2w^{\widehat{j}}s_{j'})=\ell(w_2s_{w^{\widehat{j}}(j')}w^{\widehat{j}})=\ell(w_2s_{w^{\widehat{j}}(j')})+\ell(w^{\widehat{j}})$ as $w^{\widehat{j}}\in W^{\widehat{j},\emptyset}$ and $w_2s_{w^{\widehat{j}}(j')}\in W(L_{\widehat j})$. By (\ref{equ: right set}) we see that $j'\in D_R(w_2w^{\widehat{j}})$ if and only if $w^{\widehat{j}}(j')\in D_R(w_2)$. If $j'=j_2$, then as $w^{\widehat{j}}>w^{\widehat{j}}s_{j_2}$ (using that $w^{\widehat{j}}s_{j_2}\in W^{\widehat{j},\emptyset}$ and that $w^{\widehat{j}}\in W^{\widehat{j},\emptyset}$ has the longest possible length) we have $\ell(w_2w^{\widehat{j}}s_{j'})=\ell(w_2)+\ell(w^{\widehat{j}}s_{j'})=\ell(w_2w^{\widehat{j}})-1$, hence $j'\in D_R(w_2w^{\widehat{j}})$. Thus $D_R(w_2w^{\widehat{j}})=(w^{\widehat{j}})^{-1}(D_R(w_2))\sqcup\{j_2\}=I_1\sqcup\{j_2\}=I_2$. We deduce that $i_{\widehat{j},\Delta}^{\infty}(\tau^{\infty})$ is the (unique) quotient of $i_{\emptyset,\Delta}^{\infty}(w_2w^{\widehat{j}}\cdot 1_{T})$ with socle $V_{I_1,\Delta}^{\infty}$ and cosocle $V_{I_2,\Delta}^{\infty}$. It has length $2$ by \ref{it: trivial PS 3} of Corollary~\ref{cor: trivial PS} and is thus isomorphic to $\pi^{\infty}$.
\end{proof}

\begin{lem}
Let $I\subseteq I'\subseteq \Delta$ and $I_0,I_1\subseteq I$. Then we have the following isomorphism in $\mathrm{Rep}^{\infty}_{\rm{adm}}(L_{I'})$
\begin{equation*}
i_{I,I'}^{\infty}(Q_I(I_0,I_1))\cong Q_{I'}(I_0\sqcup(I'\setminus I),I_1).
\end{equation*}
\end{lem}
\begin{proof}
By choosing a sequence of subsets $I'=J_0\supsetneq J_1\supsetneq \cdots \supsetneq J_t=I$ with $t=\#I'\setminus I$ (and thus $\#J_{t'-1}\setminus J_{t'}=1$ for $1\leq t'\leq t$), we easily reduce to the case when $I'\setminus I=\{j\}$ for some $j\in\Delta$.
As $Q_I(I_0,I_1)$ is $G$-basic, so is $i_{I,I'}^{\infty}(Q_I(I_0,I_1))$. By \ref{it: basic as image} of Remark~\ref{rem: basic PS intertwine} we know that $i_{I,I'}^{\infty}(Q_I(I_0,I_1))$ is multiplicity free with simple socle and cosocle. As $i_{I,I'}^{\infty}(-)$ is exact and $Q_I(I_0,I_1)$ has socle $V_{I_0,I}^{\infty}$ and cosocle $V_{I_1,I}^{\infty}$, by \ref{it: length two 1} of Lemma~\ref{lem: length two trivial block} we know that $V_{I_0\sqcup\{j\},I'}^{\infty}$ (resp.~$V_{I_1,I'}^{\infty}$) occurs in the socle (resp.~cosocle) of $i_{I,I'}^{\infty}(Q_I(I_0,I_1))$, which forces the latter to have socle $V_{I_0\sqcup\{j\},I'}^{\infty}$ (resp.~cosocle $V_{I_1,I'}^{\infty}$). By \ref{it: Hom cube 1} of Lemma~\ref{lem: Hom cube} we conclude that $i_{I,I'}^{\infty}(Q_I(I_0,I_1))\cong Q_{I'}(I_0\sqcup(I'\setminus I),I_1')$.
\end{proof}

\begin{lem}\label{lem: smooth devissage St}
Let $I_0,I_1\subseteq I\subseteq \Delta$ with $I_0\neq I_1$, and $I_2\subseteq I$ such that $d(I_0,I_2)=d(I_0,I_1)-1$ and $d(I_2,I_1)=1$. Let $\pi^{\infty}$ in $\mathrm{Rep}^{\infty}_{\rm{adm}}(L_I)$ the unique non-split extension $0\rightarrow V_{I_1,I}^{\infty}\rightarrow \pi^{\infty}\rightarrow V_{I_2,I}^{\infty}\rightarrow 0$ (Lemma \ref{lem: Ext sm St}). We have $d_I(V_{I_0,I}^{\infty},V_{I_2,I}^{\infty})=d_I(V_{I_0,I}^{\infty},V_{I_1,I}^{\infty})-1$ and $d_I(V_{I_0,I}^{\infty},\pi^{\infty})=\infty$.
\end{lem}
\begin{proof}
Note first that $d_I(V_{I_0,I}^{\infty},\pi^{\infty})=\infty$ and the short exact sequence $0\rightarrow V_{I_1,I}^{\infty}\rightarrow \pi^{\infty}\rightarrow V_{I_2,I}^{\infty}\rightarrow 0$ imply $\mathrm{Ext}_{L_I}^{k-1}(V_{I_0,I}^{\infty}, V_{I_2,I}^{\infty})^\infty\cong \mathrm{Ext}_{L_I}^{k}(V_{I_0,I}^{\infty}, V_{I_1,I}^{\infty})^\infty$ for $k\geq 1$, and in particular $d_I(V_{I_0,I}^{\infty},V_{I_2,I}^{\infty})=d_I(V_{I_0,I}^{\infty},V_{I_1,I}^{\infty})-1$. Hence it is enough to prove $d_I(V_{I_0,I}^{\infty},\pi^{\infty})=\infty$, i.e.~$\mathrm{Ext}_{L_I}^{k}(V_{I_0,I}^{\infty}, \pi^{\infty})^\infty=0$ for $k\geq 0$.\bigskip

\textbf{Case $1$}: If $I_2=I_1\setminus\{j_1\}$ for some $j_1\in I_1\setminus I_0$, then we have $\pi^{\infty}\cong i_{I\setminus\{j_1\},I}^{\infty}(V_{I_2,I\setminus\{j_1\}}^{\infty})$ by \ref{it: length two 1} of Lemma~\ref{lem: length two trivial block} (which extends verbatim with $\Delta$ replaced by $I$ and $G$ by $L_I$ when $I_0,I_1\subseteq I$). Let $I^-\defeq I\setminus\{j_1\}$, $\Sigma\defeq W(L_{I^-})\cdot 1_{T}$ and recall that $V_{I_2,I^-}^{\infty}\in \cB^{I^-}_{\Sigma}$ by (\ref{equ: Jacquet Steinberg}). By (\ref{equ: first adjunction}) and (\ref{equ: general Jacquet block}) we deduce for $k\geq 0$
\[\mathrm{Ext}_{L_I}^{k}(V_{I_0,I}^{\infty},\pi^{\infty})^\infty\cong \mathrm{Ext}_{L_{I^-}}^{k}(J_{I,I^-}(V_{I_0,I}^{\infty})_{\cB^{I^-}_{\Sigma}}, V_{I_2,I^-}^{\infty})^\infty.\]
Hence, to prove $d_I(V_{I_0,I}^{\infty},\pi^{\infty})=\infty$, it suffices to show that $J_{I,I^-}(V_{I_0,I}^{\infty})_{\cB^{I^-}_{\Sigma}}=0$, or equivalently $\cJ(V_{I_0,I}^{\infty})\cap \Sigma=\emptyset$ (see the discussion below (\ref{equ: sm block decomposition})). Let $w\cdot 1_{T}\in \cJ(V_{I_0,I}^{\infty})\cap \Sigma$ (using (\ref{equ: Jacquet Steinberg})), then $w\in W(L_{I^-})$ by definition of $\Sigma$ and $I_0=I\setminus D_R(w)$ by (\ref{equ: Jacquet Steinberg}). However, $w\in W(L_{I^-})$ implies $w(j_1)\in \Phi^+$ and thus $j_1\notin D_R(w)$ by (\ref{equ: right set}). Hence $j_1\in I\setminus D_R(w)=I_0$ which contradicts $j_1\in I_1\setminus I_0$. Consequently we have $\cJ(V_{I_0,I}^{\infty})\cap \Sigma=\emptyset$.\bigskip

\textbf{Case $2$}: If $I_2=I_1\sqcup\{j_2\}$ for some $j_2\in I_0\setminus I_1$, then we deduce from \ref{it: length two 2} of Lemma~\ref{lem: length two trivial block} (applied with $L_I$ instead of $G$) that $\pi^{\infty}\cong i_{I\setminus\{j\},I}^{\infty}(V_{w^{J^-}(I_1),I\setminus\{j\}}^{\infty}\otimes_E\delta)$ where $J\subseteq I$ is the subset corresponding to the Levi block of $L_I$ containing $j_2$, $j$ is the unique element of $J$ such that, for $J^-\defeq J\setminus\{j\}$ and $w^{J^-}$ the longest element of $W^{J^-,\emptyset}(L_J)$, we have $j_2\in D_R(w^{J^-})$, and $\delta: L_I\rightarrow E^\times$ is the unique smooth character such that $\delta|_{T}\cong w^{J^-}\!\!\cdot 1_{T}$. Let $I^-\defeq I\setminus\{j\}$, $\Sigma\defeq W(L_{I^-})w^{J^-}\!\!\cdot 1_{T}$, and note that $W^{J^-,\emptyset}(L_J)\subseteq W^{I^-,\emptyset}(L_I)$. As in Case $1$ we have $V_{w^{J^-}(I_1),I^-}^{\infty}\otimes_E\delta\in \cB^{I^-}_{\Sigma}$ and for $k\geq 0$:
\[\mathrm{Ext}_{L_I}^{k}(V_{I_0,I}^{\infty},\pi^{\infty})^\infty\cong \mathrm{Ext}_{L_{I^-}}^{k}(J_{I,I^-}(V_{I_0,I}^{\infty})_{\cB^{I^-}_{\Sigma}}, V_{w^{J^-}(I_1),I^-}^{\infty}\otimes_E\delta)^\infty.\]
Again it is enough to prove $\cJ(V_{I_0,I}^{\infty})\cap \Sigma=\emptyset$. Let $w\cdot 1_{T}\in \cJ(V_{I_0,I}^{\infty})\cap \Sigma$, then $w\in W(L_{I^-})w^{J^-}$ and $I_0=I\setminus D_R(w)$ as in Case $1$. However, since $w\in W(L_{I^-})w^{J^-}$ and $w^{J^-}\in W^{I^-,\emptyset}(L_I)$ we have $D_R(w^{J^-})\subseteq D_R(w)$. As $j_2\in D_R(w^{J^-})$, we have $j_2\in D_R(w)$ and thus $j_2\notin I\setminus D_R(w)=I_0$, which contradicts $j_2\in I_0\setminus I_1$. Hence, we must have $\cJ(V_{I_0,I}^{\infty})\cap \Sigma=\emptyset$.
\end{proof}

For $I_0,I_1,I_2\subseteq \Delta$, we define $d(I_2,[I_0,I_1])\defeq \min\{d(I_2,I_1')\mid I_1'\in[I_0,I_1]\}$.

\begin{lem}\label{lem: minimal distance}
For $I_0,I_1,I_2\subseteq \Delta$, there exists a unique $I_3\in [I_0,I_1]$ such that
\begin{equation}\label{equ: minimal distance}
d(I_2,I_3)=d(I_2,[I_0,I_1]).
\end{equation}
\end{lem}
\begin{proof}
Let us first prove that
\begin{equation}\label{equ: explicit minimal distance}
I_3\defeq (I_0\cap I_1)\cup (I_2\cap (I_0\cup I_1))
\end{equation}
satisfies (\ref{equ: minimal distance}). In fact, we prove that it satisfies the stronger statement: for all $I_1'\in [I_0,I_1]$
\begin{equation}\label{equ: orthogonal distance}
d(I_2,I_3)=d(I_2,I_1')-d(I_3,I_1').
\end{equation}
Recall from (\ref{equ: same interval}) that $I_1'\in [I_0,I_1]$ if and only if
\begin{equation}\label{equ: set bound}
I_0\cap I_1\subseteq I_1'\subseteq I_0\cup I_1.
\end{equation}
It follows from (\ref{equ: set bound}) and (\ref{equ: explicit minimal distance}) that, for any $I'_1\in [I_0,I_1]$, $I_3\setminus I_2=(I_0\cap I_1)\setminus I_2\subseteq I_1'\setminus I_2$ and $I'_1\setminus I_3\subseteq I'_1\setminus I_2$, which implies $I'_1\setminus I_2=I'_1\setminus I_3 \sqcup I_3\setminus I_2$. Likewise we have $I_2\setminus I_3=I_2\setminus (I_0\cup I_1)\subseteq I_2\setminus I_1'$ and $I_3\setminus I'_1\subseteq I_2\setminus I_1'$, hence $I_2\setminus I'_1=I_2\setminus I_3 \sqcup I_3\setminus I'_1$. This implies $d(I_2,I_1')=d(I_2,I_3)+d(I_3,I_1')$, i.e.~(\ref{equ: orthogonal distance}), and therefore (\ref{equ: minimal distance}).

Now we prove unicity. Let $I'_3\in [I_0,I_1]$ satisfying (\ref{equ: minimal distance}), which implies $d(I_2,I'_3)=d(I_2,I_3)$. By the previous paragraph we have $I_3\setminus I_2\subseteq I'_3\setminus I_2$ and $I_2\setminus I_3\subseteq I_2\setminus I'_3$, hence we must have $I_3\setminus I_2=I'_3\setminus I_2$ and $I_2\setminus I_3= I_2\setminus I'_3$. Assume there is $x\in I'_3\setminus I_3$, then necessarily $x\notin I_2$ by (\ref{equ: explicit minimal distance}), hence $x\in I'_3\setminus I_2$ but $x\notin I_3\setminus I_2$, which is a contradiction. Assume there is $x\in I_3\setminus I'_3$, then necessarily $x\notin I_0\cap I_1$, hence $x\in I_2$. Thus $x\in I_2\setminus I'_3$ but $x\notin I_2\setminus I_3$, which is again a contradiction. We therefore must have $I_3=I'_3$.
\end{proof}

\begin{lem}\label{lem: Ext of cube}
Let $I\subseteq \Delta$ and $I_0,I_1,I_2\subseteq I$. We have $d_I(Q_I(I_0,I_1),V_{I_2,I}^{\infty})<\infty$ if and only~if
\begin{equation}\label{equ: minimize distance}
d(I_2,I_1)=d(I_2,[I_0,I_1]),
\end{equation}
in which case $d_I(Q_I(I_0,I_1),V_{I_2,I}^{\infty})=d(I_2,I_1)$.
\end{lem}
\begin{proof}
Let $I_3$ as in Lemma~\ref{lem: minimal distance}, so (\ref{equ: minimize distance}) is equivalent to $I_3=I_1$. If $I_3=I_1$, then by (\ref{equ: orthogonal distance}) we have $d(I_2,I_1')=d(I_2,I_1)+d(I_1,I_1')>d(I_2,I_1)$ for $I_1'\in[I_0,I_1]\setminus\{I_1\}$. By \ref{it: Hom cube 1} of Lemma \ref{lem: Hom cube}, Lemma~\ref{lem: Ext sm St} and an obvious d\'evissage we have $\mathrm{Ext}_{L_I}^k(\mathrm{ker}(Q_I(I_0,I_1)\rightarrow V_{I_1,I}^{\infty}), V_{I_2,I}^{\infty})^{\infty}=0$ for $k\leq d(I_2,I_1)$. Hence, the surjection $Q_I(I_0,I_1)\twoheadrightarrow V_{I_1,I}^{\infty}$ induces an isomorphism for $k\leq d(I_2,I_1)$
\[\mathrm{Ext}_{L_I}^k(V_{I_1,I}^{\infty}, V_{I_2,I}^{\infty})^{\infty}\buildrel\sim\over\longrightarrow \mathrm{Ext}_{L_I}^k(Q_I(I_0,I_1), V_{I_2,I}^{\infty})^{\infty}.\]
As $d_I(V_{I_1,I}^{\infty},V_{I_2,I}^{\infty})=d(I_1,I_2)=d(I_2,I_1)$ (Lemma~\ref{lem: Ext sm St}) we deduce $d_I(Q_I(I_0,I_1),V_{I_2,I}^{\infty})=d(I_2,I_1)$.

We assume from now $I_3\neq I_1$ and prove that $d_I(Q_I(I_0,I_1),V_{I_2,I}^{\infty})=\infty$. As $I_3\neq I_1$, one easily sees that there exists $I_4\in [I_3,I_1]$ such that $d(I_3,I_4)=1$ and $d(I_3,I_1)=d(I_4,I_1)+1$ ($I_4$ is obtained by either adding to $I_3$ an element of $I_1\setminus I_3$ or withdrawing from $I_3$ an element of $I_3\setminus I_1$ (which is hence in $I_3\cap I_0$)). By applying (\ref{equ: orthogonal distance}) (with $I_4$ replacing $I_1'$ there) we deduce $d(I_4,I_2)=d(I_3,I_2)+1$, and thus we have either $\emptyset\neq I_3\setminus I_4\subseteq I_2$, or $I_4\setminus I_3\ne \emptyset$ with $(I_4\setminus I_3)\cap I_2=\emptyset$. Now consider any $I_3',I_4'\in[I_0,I_1]$ such that $I_3\setminus I_4=I_3'\setminus I_4'$ and $I_4\setminus I_3=I_4'\setminus I_3'$. This implies $d(I_3',I_4')=1$ and one easily checks that one still has $d(I_3',I_1)=d(I_4',I_1)+1$, and that we have either $\emptyset\neq I'_4\setminus I'_3\subseteq I_1$, or $I'_3\setminus I'_4\ne \emptyset$ with $(I'_3\setminus I'_4)\cap I_1=\emptyset$. In all cases $I'_4\in [I'_3,I_1]$ and by \ref{it: Hom cube 3} of Lemma~\ref{lem: Hom cube} we deduce that $Q_I(I_3',I_4')$ is a subquotient of $Q_I(I_0,I_1)$. Moreover we also have either $\emptyset\neq I_3'\setminus I_4'\subseteq I_2$, or $I'_4\setminus I'_3\ne \emptyset$ with $(I'_4\setminus I'_3)\cap I_2=\emptyset$, and in both cases $d(I_2,I_4')=d(I_2,I_3')+1$. We then deduce by (\ref{equ: dual sm Ext}) for $k\geq 0$:
\begin{equation}\label{equ: dual length two}
\mathrm{Ext}_{L_I}^k(Q_I(I_3',I_4'), V_{I_2,I}^{\infty})^{\infty}\!\cong \!\mathrm{Ext}_{L_I}^k((V_{I_2,I}^{\infty})^{\sim}, Q_I(I_3',I_4')^{\sim})^{\infty}\!\cong \!\mathrm{Ext}_{L_I}^k(V_{I_2,I}^{\infty}, Q_I(I_4',I_3'))^{\infty}\!=\!0
\end{equation}
where the last equality follows from the last assertion in Lemma~\ref{lem: smooth devissage St} together with \ref{it: Hom cube 1} of Lemma \ref{lem: Hom cube}. Now, for a given $j\in I_1 \setminus I_0$ (resp.~$j\in I_0 \setminus I_1$), one can write $[I_0,I_1]$ as a disjoint union of $[I_0',I_1']$ for $I_0',I_1'\in[I_0,I_1]$ such that $I_1'=I_0' \sqcup \{j\}$ (resp.~$I_0'=I_1' \sqcup \{j\}$). It follows that $[I_0,I_1]$ is the disjoint union of $[I_3',I_4']$ for some $I_3',I_4'\in [I_0,I_1]$ as discussed above, and hence that $Q_I(I_0,I_1)$ admits a filtration whose graded pieces are the $Q_I(I_3',I_4')$ (using \ref{it: Hom cube 3} of Lemma~\ref{lem: Hom cube}). By (\ref{equ: dual length two}) and an obvious d\'evissage we deduce $\mathrm{Ext}_{L_I}^k(Q_I(I_0,I_1), V_{I_2,I}^{\infty})^{\infty}=0$ for $k\geq 0$, i.e.~$d_I(Q_I(I_0,I_1),V_{I_2,I}^{\infty})=\infty$.
\end{proof}

\begin{lem}
Let $I\subseteq \Delta$ and $I_0,I_1,I_0',I_1'\subseteq I$ with $[I_0,I_1]\cap [I_0',I_1']=\emptyset$. Then we have
\begin{equation}\label{equ: sm cube Ext1}
\mathrm{Ext}_{L_I}^1(Q_I(I_0',I_1'),Q_I(I_0,I_1))^{\infty}\neq 0
\end{equation}
if and only if $[I_0,I_1']\subseteq [I_0,I_1]\sqcup [I_0',I_1']$, in which case (\ref{equ: sm cube Ext1}) is one dimensional.
\end{lem}
\begin{proof}
It follows from (\ref{equ: intersection of interval}) that we have $[I_0,I_1']\cap [I_0,I_1]=[I_0,I_1'']$ with $I_1''\defeq (I_0\cap (I_1\cup I_1'))\sqcup ((I_1\cap I_1')\setminus I_0)$, and similarly $[I_0,I_1']\cap [I_0',I_1']=[I_0'',I_1']$ with $I_0''\defeq (I_1'\cap (I_0\cup I_0'))\sqcup ((I_0\cap I_0')\setminus I_1')$.

Assume first that $[I_0,I_1']\subseteq [I_0,I_1]\sqcup [I_0',I_1']$, which together with $[I_0,I_1]\cap [I_0',I_1']=\emptyset$ implies that
\begin{equation}\label{equ: sm cube Ext1 factor}
[I_0,I_1']=([I_0,I_1']\cap [I_0,I_1])\sqcup ([I_0,I_1']\cap [I_0',I_1'])=[I_0,I_1'']\sqcup [I_0'',I_1'].
\end{equation}
By \ref{it: Hom cube 1} and \ref{it: Hom cube 2} of Lemma~\ref{lem: Hom cube} and since $I_0''\in [I_0,I_1']$ (resp.~$I_1''\in [I_0,I_1']$), there exists a (unique up to scalar) surjection $Q_I(I_0,I_1')\twoheadrightarrow Q_I(I_0'',I_1')$ (resp.~injection $Q_I(I_0,I_1'')\hookrightarrow Q_I(I_0,I_1')$). By \ref{it: Hom cube 1} of Lemma~\ref{lem: Hom cube} and (\ref{equ: sm cube Ext1 factor}), we deduce a short exact sequence
\[0\rightarrow Q_I(I_0,I_1'') \rightarrow Q_I(I_0,I_1') \rightarrow Q_I(I_0'',I_1') \rightarrow 0\]
which has to be non-split as $Q_I(I_0,I_1')$ has simple socle and cosocle. In particular, we have
\begin{equation}\label{equ: sm cube Ext1 1}
\mathrm{Ext}_{G}^1(Q_I(I_0'',I_1'),Q_I(I_0,I_1''))\neq 0.
\end{equation}
Our assumption $[I_0,I_1]\cap [I_0',I_1']=\emptyset$ implies that $Q_I(I_0',I_1')$ and $Q_I(I_0,I_1)$ have no common Jordan-H\"older factor. Thus the surjection $Q_I(I_0',I_1')\twoheadrightarrow Q_I(I_0'',I_1')$ and the injection $Q_I(I_0,I_1'')\hookrightarrow Q_I(I_0,I_1)$) from \ref{it: Hom cube 2} of Lemma~\ref{lem: Hom cube} induce an injection
\begin{equation*}
\mathrm{Ext}_{G}^1(Q_I(I_0'',I_1'),Q_I(I_0,I_1''))\hookrightarrow \mathrm{Ext}_{G}^1(Q_I(I_0',I_1'),Q_I(I_0,I_1)),
\end{equation*}
which together with (\ref{equ: sm cube Ext1 1}) gives (\ref{equ: sm cube Ext1}). It is clear from $[I_0,I_1]\cap [I_0',I_1']=\emptyset$ and (\ref{equ: sm cube Ext1}) that $d(Q_I(I_0',I_1'),Q_I(I_0,I_1))=1$, which by Lemma~\ref{lem: dim one} implies that (\ref{equ: sm cube Ext1}) is one dimensional.

Now we assume that (\ref{equ: sm cube Ext1}) holds and let $V$ be a representationt which fits into a non-split extension
\begin{equation*}
0\rightarrow Q_I(I_0,I_1) \rightarrow V \rightarrow Q_I(I_0',I_1') \rightarrow 0.
\end{equation*}
By Lemma~\ref{lem: basic subquotient of Ext} we know that $V$ admits unique subquotient $V'$ which is $G$-basic with socle $\mathrm{soc}_{L_I}(Q_I(I_0,I_1))\cong V_{I_0,I}^{\infty}$ and cosocle $\mathrm{cosoc}_{L_I}(Q_I(I_0',I_1'))\cong V_{I_1',I}^{\infty}$. By \ref{it: Hom cube 1} of Lemma~\ref{lem: Hom cube} we must have $V'\cong Q_I(I_0,I_1')$ and in particular
\[\mathrm{JH}_{L_I}(Q_I(I_0,I_1'))\subseteq \mathrm{JH}_{L_I}(V)=\mathrm{JH}_{L_I}(Q_I(I_0,I_1))\sqcup \mathrm{JH}_{L_I}(Q_I(I_0',I_1'))\]
which (again by \emph{loc.~cit.}) gives $[I_0,I_1']\subseteq [I_0,I_1]\sqcup [I_0',I_1']$.
\end{proof}

\begin{lem}\label{lem: Hom vers induction}
For $i=0,1$ let $I_i\subseteq \Delta$ and $\pi_i^{\infty}$ in $\mathrm{Rep}^{\infty}_{\rm{adm}}(L_{I_i})$ irreducible $G$-regular with $\Sigma_i\defeq W(L_{I_i})\cdot\cJ(\pi_i^{\infty})$. Let $w\in W^{I_0,I_1}$.
\begin{enumerate}[label=(\roman*)]
\item \label{it: Hom vers induction1} We have
\begin{equation}\label{equ: Hom induction cell}
\Hom_{L_{I_1}}(i_{I_0,I_1,w}^{\infty}(J_{I_0,I_1,w}(\pi_0^{\infty})),\pi_1^{\infty})\neq 0
\end{equation}
if and only if we have an isomorphism of irreducible $G$-regular representations
\begin{equation}\label{equ: Hom vers induction Jacquet}
J_{I_0,I_1,w}(\pi_0^{\infty})_{\cB^{w^{-1}(I_0)\cap I_1}_{w^{-1}\cdot\Sigma_0\cap\Sigma_1}}\cong J_{I_1,w^{-1}(I_0)\cap I_1}'(\pi_1^{\infty})_{\cB^{w^{-1}(I_0)\cap I_1}_{w^{-1}\cdot\Sigma_0\cap\Sigma_1}},
\end{equation}
in which case (\ref{equ: Hom induction cell}) is one dimensional.
If we denote by $\pi^{\infty}$ the representation in (\ref{equ: Hom vers induction Jacquet}), we have $\pi_1^{\infty}\cong \mathrm{cosoc}_{L_{I_1}}(i_{w^{-1}(I_0)\cap I_1,I_1}^{\infty}(\pi^{\infty}))$ and
\begin{multline*}
\pi_0^{\infty}\cong \mathrm{soc}_{L_{I_0}}(i_{w^{-1}(I_0)\cap I_1,w^{-1}(I_0)}^{\infty}(\pi^{\infty}\otimes_E \delta_{I_0,I_1,w}^{-1})^{w^{-1}})\\
\cong \mathrm{soc}_{L_{I_0}}(i_{I_0\cap w(I_1),I_0}^{\infty}((\pi^{\infty}\otimes_E \delta_{I_0,I_1,w}^{-1})^{w^{-1}})).\end{multline*}
\item \label{it: Hom vers induction2} Assume (\ref{equ: Hom induction cell}) and let $\sigma_0^{\infty}$ be an irreducible constituent of $i_{I_0\cap w(I_1),I_0}^{\infty}((\pi^{\infty}\otimes_E \delta_{I_0,I_1,w}^{-1})^{w^{-1}})$ and $\sigma_1^{\infty}$ an irreducible constituent of $i_{w^{-1}(I_0)\cap I_1,I_1}^{\infty}(\pi^{\infty})$. Then we have
\[J_{I_0,I_1,w}(\sigma_0^{\infty})_{\cB^{w^{-1}(I_0)\cap I_1}_{w^{-1}\cdot\Sigma_0\cap\Sigma_1}}=0\text{ if }\sigma_0^{\infty}\neq \pi_0^{\infty},\ \ J_{I_1,w^{-1}(I_0)\cap I_1}'(\sigma_1^{\infty})_{\cB^{w^{-1}(I_0)\cap I_1}_{w^{-1}\cdot\Sigma_0\cap\Sigma_1}}=0\text{ if }\sigma_1^{\infty}\neq \pi_1^{\infty},\]
and we also have
\begin{equation}\label{equ: Hom vers induction JH}
\mathrm{Ext}_{L_{I_1}}^k(i_{I_0,I_1,w}^{\infty}(J_{I_0,I_1,w}(\sigma_0^{\infty})),\sigma_1^{\infty})\neq 0\ {for\ some\ }k\geq 0
\end{equation}
if and only if $\sigma_0^{\infty}=\pi_0^{\infty}$ and $\sigma_1^{\infty}=\pi_1^{\infty}$.
\end{enumerate}
\end{lem}
\begin{proof}
Recall that by (\ref{equ: second adjunction}) and (\ref{equ: Bruhat induction}) we have canonical isomorphisms for $k\geq 0$
\begin{equation}\label{equ: Hom induction adjunction}
\mathrm{Ext}_{L_{I_1}}^k(J_{I_0,I_1,w}(\pi_0^{\infty}),J_{I_1,w^{-1}(I_0)\cap I_1}'(\pi_1^{\infty}))\cong \mathrm{Ext}_{L_{I_1}}^k(i_{I_0,I_1,w}^{\infty}(J_{I_0,I_1,w}(\pi_0^{\infty})),\pi_1^{\infty})
\end{equation}
and that any $G$-regular left $W(L_{w^{-1}(I_0)\cap I_1})$-coset $\Sigma\subseteq \widehat{T}^{\infty}$ gives a block $\cB^{w^{-1}(I_0)\cap I_1}_{\Sigma}$ (see the discussion around (\ref{equ: sm block decomposition})). We write $I\defeq w^{-1}(I_0)\cap I_1$ and $\delta\defeq \delta_{I_0,I_1,w}$.

We prove \ref{it: Hom vers induction1}. Assume that (\ref{equ: Hom induction cell}) holds, thus the vector spaces in (\ref{equ: Hom induction adjunction}) are non-zero for $k=0$. Then by Lemma~\ref{lem: Jacquet basic} and Remark \ref{rem: twist Jacquet basic} there exists a left $W(L_{I})$-coset $\Sigma$ such that
\begin{equation}\label{sigmanonzero}
J_{I_0,I_1,w}(\pi_0^{\infty})_{\cB^{I}_{\Sigma}}\neq 0\neq J_{I_1,I}'(\pi_1^{\infty})_{\cB^{I}_{\Sigma}}.
\end{equation}
Recall that $\Sigma_i=W(L_{I_i})\cdot \cJ(\pi_i^{\infty})$ are single left regular $W(L_{I_i})$-cosets for $i=0,1$ by (the last statement in) \ref{it: PS 1} of Lemma \ref{lem: Jacquet of PS}. An easy exercise using that all characters here are $G$-regular shows that $w^{-1}\cdot\Sigma_0\cap\Sigma_1$ is again a single (regular) left $W(L_{w^{-1}(I_0)\cap I_1})$-coset. As we have $\cJ(J_{I_0,I_1,w}(\pi_0^{\infty}))\subseteq w^{-1}\cdot\Sigma_0$ and $\cJ(J_{I_1,I}'(\pi_1^{\infty}))\subseteq \Sigma_1$, we necessarily have $\Sigma=w^{-1}\cdot\Sigma_0\cap\Sigma_1$. In other words, $\Sigma$ as in (\ref{sigmanonzero}) is uniquely determined by $\pi_0^{\infty}$, $\pi_1^{\infty}$ and $w$, and thus we have a canonical isomorphism
\begin{equation}\label{equ: Hom vers induction block}
0\neq\Hom_{L_{I_1}}(J_{I_0,I_1,w}(\pi_0^{\infty}),J_{I_1,I}'(\pi_1^{\infty}))
\cong \Hom_{L_{I_1}}(J_{I_0,I_1,w}(\pi_0^{\infty})_{\cB^{I}_{\Sigma}},J_{I_1,I}'(\pi_1^{\infty})_{\cB^{I}_{\Sigma}}).
\end{equation}
As $\pi_0^{\infty}$ and $\pi_1^{\infty}$ are irreducible $G$-regular, both $J_{I_0,I_1,w}(\pi_0^{\infty})_{\cB^{I}_{\Sigma}}$ and $J_{I_1,I}'(\pi_1^{\infty})_{\cB^{I}_{\Sigma}}$ are irreducible $G$-regular by (the last statement in) Lemma~\ref{lem: Jacquet basic} and Remark \ref{rem: twist Jacquet basic}. Hence they must be isomorphic by (\ref{equ: Hom vers induction block}) and we denote them by $\pi^{\infty}$. Note in particular that (\ref{equ: Hom vers induction block}) is one dimensional and thus (\ref{equ: Hom induction cell}) is one dimensional.
We have by (\ref{equ: second adjunction})
\[\Hom_{L_{I_1}}(i_{I,I_1}^{\infty}(\pi^{\infty}),\pi_1^{\infty})\cong \Hom_{L_{I}}(\pi^{\infty}, J_{I_1,I}'(\pi_1^{\infty}))\neq 0\]
and by (\ref{equ: first adjunction}) (together with (\ref{equ: Bruhat Jacquet}))
\[\Hom_{L_{I_0}}((\pi_0^{\infty})^{w}, i_{I,w^{-1}(I_0)}^{\infty}(\pi^{\infty}\otimes_E\delta^{-1}))\cong \Hom_{L_I}(J_{I_0,I_1,w}(\pi_0^{\infty}),\pi^{\infty})\neq 0.\]
As $\pi_0^{\infty}$ and $\pi_1^{\infty}$ are irreducible $G$-regular, we deduce $\pi_0^{\infty}\cong \mathrm{soc}_{L_{I_0}}(i_{w(I),I_0}^{\infty}(\pi^{\infty}\otimes_E\delta^{-1})^{w^{-1}})$ and $\pi_1^{\infty}\cong \mathrm{cosoc}_{L_{I_1}}(i_{I,I_1}^{\infty}(\pi^{\infty}))$.

We prove \ref{it: Hom vers induction2}. We borrow $\Sigma$ and $\pi^\infty$ from \ref{it: Hom vers induction1} and we fix $\sigma_0^{\infty}\in\mathrm{JH}_{L_{I_0}}(i_{w(I),I_0}^{\infty}((\pi^{\infty}\otimes_E \delta^{-1})^{w^{-1}}))$ and $\sigma_1^{\infty}\in\mathrm{JH}_{L_{I_1}}(i_{I,I_1}^{\infty}(\pi^{\infty}))$.
It follows from Lemma~\ref{lem: smooth geometric lemma} that the natural injection
\begin{multline*}
J_{I_0,I_1,w}(\pi_0^{\infty})_{\cB^{I}_{\Sigma}}\cong \pi^{\infty}\cong i_{I,I,1}^{\infty}(J_{I,I,1}(\pi^{\infty}\otimes_E\delta^{-1}))\otimes_E\delta\\
\hookrightarrow (J_{w^{-1}(I_0),I}(i_{I,w^{-1}(I_0)}^{\infty}(\pi^{\infty}\otimes_E\delta^{-1}))\otimes_E\delta)_{\cB^{I}_{\Sigma}}
\end{multline*}
is an isomorphism, and hence that $J_{I_0,I_1,w}(i_{I,w^{-1}(I_0)}^{\infty}(\pi^{\infty}\otimes_E\delta^{-1})^{w^{-1}}/\pi_0^{\infty})_{\cB^{I}_{\Sigma}}=0$, or equivalently $J_{I_0,I_1,w}(\sigma_0^{\infty})_{\cB^{I}_{\Sigma}}=0$, for $\sigma_0^{\infty}\neq \pi_0^{\infty}$. Similarly, for $w_1$ the longest element of $W^{\emptyset, I}(L_{I_1})$, it follows from (\ref{equ: explicit twist Jacquet}) and Lemma~\ref{lem: smooth geometric lemma} that
\begin{eqnarray*}
J_{I_1,I}'(i_{I,I_1}^{\infty}(\pi^{\infty}))_{\cB^{I}_{\Sigma}}&\cong &
(J_{I_1,w_1(I)}(i_{I,I_1}^{\infty}(\pi^{\infty}))^{w_1}\otimes_E \delta_{P_{I}\cap L_{I_1}}^{-1})_{\cB^{I}_{\Sigma}}\\
&\cong &(i_{I,w_1(I),w_1}^{\infty}(J_{I,w_1(I),w_1}(\pi^{\infty}))^{w_1}\otimes_E \delta_{P_{I}\cap L_{I_1}}^{-1})_{\cB^{I}_{\Sigma}}\\
&\cong &((\pi^{\infty})^{w_1^{-1}}\otimes_E \delta_{P_{w_1(I)}\cap L_{I_1}})^{w_1}\otimes_E \delta_{P_{I}\cap L_{I_1}}^{-1}\\
&\cong &\pi^{\infty} \cong J_{I_1,I}'(\pi_1^{\infty})_{\cB^{I}_{\Sigma}}
\end{eqnarray*}
which implies $J_{I_1,I}'(\sigma_1^{\infty})_{\cB^{I}_{\Sigma}}=0$ for $\sigma_1^{\infty}\neq \pi_1^{\infty}$. Furthermore, arguing as in \ref{it: Hom vers induction1} with $\sigma_i^{\infty}$ instead of $\pi_i^{\infty}$ (and $\mathrm{Ext}_{L_I}^k$ instead of $\Hom_{L_I}$), we have seen that (\ref{equ: Hom vers induction JH}) forces
\[J_{I_0,I_1,w}(\sigma_0^{\infty})_{\cB^{I}_{\Sigma'}}\neq 0\neq J_{I_1,I}'(\sigma_1^{\infty})_{\cB^{I}_{\Sigma'}}\]
where $\Sigma'\defeq w^{-1}\cdot W(L_{I_0})\cdot\cJ(\sigma_0^{\infty}) \cap W(L_{I_1})\cdot\cJ(\sigma_1^{\infty})$. But it easily follows from \ref{it: PS 1} of Lemma~\ref{lem: Jacquet of PS} that
\[W(L_{I_1})\cdot\cJ(\sigma_1^{\infty})=W(L_{I_1})\cdot \cJ(i_{I,I_1}^{\infty}(\pi^{\infty}))=W(L_{I_1})\cdot\cJ(\pi_1^{\infty})=\Sigma_1\]
and similarly $W(L_{I_0})\cdot\cJ(\sigma_0^{\infty})=\Sigma_0$, which forces $\Sigma'=\Sigma$. By what we have proven before on $J_{I_0,I_1,w}(\sigma_0^{\infty})_{\cB^{I}_{\Sigma}}$ and $J_{I_1,I}'(\sigma_1^{\infty})_{\cB^{I}_{\Sigma}}$, we deduce that (\ref{equ: Hom vers induction JH}) forces $\sigma_0^{\infty}=\pi_0^{\infty}$ and $\sigma_1^{\infty}=\pi_1^{\infty}$.
\end{proof}

\subsection{Examples of \texorpdfstring{$G$}{G}-basic representations}\label{subsec: sm example}

We study specific $G$-basic representations: the representations $\pi_{j_1,j_2}^{\infty}$ below. We prove several technical lemmas on these representations which will be mainly used in \S\ref{subsec: square} below.\bigskip

Given $j,j'\in\Z$, we write $[j,j']\defeq \{j''\in\Z\mid j\leq j''\leq j'\}$ (hence $[j,j']=\emptyset$ if $j'<j$). Recall that $\widehat{j}= \Delta\setminus\{j\}$ for $j\in\Delta$. We also use the convenient notation $\widehat{n}\defeq \Delta$. We define
\[\mathbf{J}^{\infty}\defeq \{(j_1,j_2)\mid 1\leq j_1\leq j_2\leq n, \ j_1\leq n-1\}\]
equipped with the partial order $(j_1,j_2)\leq (j_1',j_2')$ if and only if $j_2\leq j_2'$ and $j_2-j_1\leq j_2'-j_1'$. If $j\geq 1$, we write $\mathrm{St}_j^{\infty}$ (resp.~$1_j$) for the smooth Steinberg (resp.~the trivial representation) of $\GL_j(K)$. For $(j_1,j_2)\in \mathbf{J}^{\infty}$, we set
\[\sigma_{j_1,j_2}^{\infty}\defeq |\mathrm{det}_{j_1}|_K^{j_1-j_2}\boxtimes_E\big(|\mathrm{det}_{j_2-j_1}|_K^{j_1}\otimes_E\mathrm{St}_{j_2-j_1}^{\infty}\big)\boxtimes_E\mathrm{St}_{n-j_2}^{\infty}\in \mathrm{Rep}^{\infty}_{\rm{adm}}(L_{\widehat{j}_1\cap\widehat{j}_2}),\]
which is irreducible $G$-regular, and we define the $G$-basic representation
\begin{equation}\label{pij_1,_2}
\left\{
\begin{array}{lll}
\pi_{j_1,j_2}^{\infty}&\defeq &i_{\widehat{j}_1\cap\widehat{j}_2,\widehat{j}_1}^{\infty}(\sigma_{j_1,j_2}^{\infty})\ \ \mathrm{if}\ 1\leq j_1<j_2\leq n-1\\
\pi_{j_1,j_2}^{\infty}&\defeq &\sigma_{j_1,j_2}^{\infty}\ \ \mathrm{if}\ j_1=j_2\ \mathrm{or}\ j_2=n.
\end{array}
\right.
\end{equation}
In particular, we have $\pi_{j_1,j_1}^{\infty}=\sigma_{j_1,j_1}^{\infty}=1_{j_1}\boxtimes_E \mathrm{St}_{n-j_1}^{\infty}=V_{[1,j_1-1],\widehat{j}_1}^{\infty}$ when $j_2=j_1$, and $\pi_{j_1,n}^{\infty}=\sigma_{j_1,n}^{\infty}=|\mathrm{det}_{j_1}|_K^{j_1-n}\boxtimes_E\big(|\mathrm{det}_{n-j_1}|_K^{j_1}\otimes_E\mathrm{St}_{n-j_1}^{\infty}\big)$ when $j_2=n$.
Recall that
\begin{equation}\label{equ: main coset}
\Sigma_{j_1,j_2}\defeq W(L_{\widehat{j}_1})\cdot \cJ(\sigma_{j_1,j_2}^{\infty})=W(L_{\widehat{j}_1})\cdot \cJ(\pi_{j_1,j_2}^{\infty})
\end{equation}
and $\Sigma_{j_1,j_2}'\defeq W(L_{\widehat{j}_1\cap\widehat{j}_2})\cdot \cJ(\sigma_{j_1,j_2}^{\infty})$ are single (regular) left cosets by \ref{it: PS 1} of Lemma \ref{lem: Jacquet of PS}. Let $\Gamma_{\Sigma_{j_1,j_2}}$ and $\Gamma_{\Sigma_{j_1,j_2}'}$ be their respective associated undirected graphs (see above Theorem \ref{thm: JH of PS}). They have the same vertices but $\Gamma_{\Sigma_{j_1,j_2}}$ possibly has more edges. However, when $j_1<j_2<n$ or equivalently when $L_{\widehat{j}_1\cap\widehat{j}_2}$ has $3$ Levi blocks, one can explicitly check that both $\Gamma_{\Sigma_{j_1,j_2}'}$ and $\Gamma_{\Sigma_{j_1,j_2}}$ have $3$ connected components, which implies $\Gamma_{\Sigma_{j_1,j_2}}=\Gamma_{\Sigma_{j_1,j_2}'}$. Thus any orientation on $\Gamma_{\Sigma_{j_1,j_2}'}$ defines a unique orientation on $\Gamma_{\Sigma_{j_1,j_2}}$, which means by Theorem~\ref{thm: JH of PS} that, for any $\chi\in \cJ(\sigma_{j_1,j_2}^{\infty})$, the two principal series $i_{\emptyset,\widehat{j}_1\cap\widehat{j}_2}^{\infty}(\chi)$ and $i_{\emptyset,\widehat{j}_1}^{\infty}(\chi)$ have the same number of constituents. This implies in particular that $\pi_{j_1,j_2}^{\infty}= i_{\widehat{j}_1\cap\widehat{j}_2,\widehat{j}_1}^{\infty}(\sigma_{j_1,j_2}^{\infty})$ remains irreducible ($G$-regular).\bigskip

For $(j_1,j_2)\in\mathbf{J}^{\infty}$, let $x_{j_1,j_2}$ be the longest element in
\begin{equation}\label{equ: special Weyl element}
\{x\in W(G)\mid D_L(x)=\{j_1\},\ \mathrm{Supp}(x)\subseteq [1,j_2-1]\}
\end{equation}
with $x_{j_1,j_2}\defeq 1$ if $j_1=j_2$. Recall that the condition $D_L(x)=\{j_1\}$ is equivalent to $x\in W^{\widehat{j}_1,\emptyset}$ (use for instance (\ref{equ: left set})). The following lemma gives the structure of $i_{\widehat{j}_1\cap \widehat{j}_2,\Delta}^{\infty}(\sigma_{j_1,j_2}^{\infty})$.

\begin{lem}\label{lem: explicit smooth induction}
For $(j_1,j_2)\in \mathbf{J}^{\infty}$, we have $\Sigma_{j_1,j_2}=W(L_{\widehat{j}_1})x_{j_1,j_2}\cdot 1_{T}$ and the explicit structure of $i_{\widehat{j}_1,\Delta}^{\infty}(\pi_{j_1,j_2}^{\infty})\cong i_{\widehat{j}_1\cap \widehat{j}_2,\Delta}^{\infty}(\sigma_{j_1,j_2}^{\infty})$ is given as follows.
\begin{enumerate}[label=(\roman*)]
\item \label{it: explicit induction 1} If $1\leq j_1<j_2\leq n-1$, then $i_{\widehat{j}_1\cap \widehat{j}_2,\Delta}^{\infty}(\sigma_{j_1,j_2}^{\infty})$ has Loewy length $3$, with socle $V_{[j_2-j_1+1,j_2],\Delta}^{\infty}$, cosocle $V_{[j_2-j_1,j_2-1],\Delta}^{\infty}$, and middle layer $V_{[j_2-j_1+1,j_2-1],\Delta}^{\infty}\oplus V_{[j_2-j_1,j_2],\Delta}^{\infty}$.
\item \label{it: explicit induction 2} If $1\leq j_1=j_2\leq n-1$, then $i_{\widehat{j}_1\cap \widehat{j}_2,\Delta}^{\infty}(\sigma_{j_1,j_2}^{\infty})$ has Loewy length $2$, with socle $V_{[1,j_1],\Delta}^{\infty}$ and cosocle $V_{[1,j_1-1],\Delta}^{\infty}$.
\item \label{it: explicit induction 3} If $1\leq j_1<j_2=n$, then $i_{\widehat{j}_1\cap \widehat{j}_2,\Delta}^{\infty}(\sigma_{j_1,j_2}^{\infty})$ has Loewy length $2$, with socle $V_{[n-j_1+1,n-1],\Delta}^{\infty}$ and cosocle $V_{[n-j_1,n-1],\Delta}^{\infty}$.
\end{enumerate}
\end{lem}
\begin{proof}
Let $z'_{j_1,j_2}\in W(L_{\widehat{j}_1\cap\widehat{j}_2})\subseteq W(G)$ such that $z'_{j_1,j_2}\cdot 1_{T}\in \cJ(\sigma_{j_1,j_2}^{\infty})$, $w_{\widehat{j}_1\cap\widehat{j}_2}\in W(L_{\widehat{j}_1\cap\widehat{j}_2})$ the longest element and $z_{j_1,j_2}\defeq w_{\widehat{j}_1\cap\widehat{j}_2}z'_{j_1,j_2}$, then by \ref{it: PS socle cosocle} of Remark~\ref{rem: basic PS intertwine} we have
\[\sigma_{j_1,j_2}^{\infty}\cong \mathrm{soc}_{L_{\widehat{j}_1\cap\widehat{j}_2}}(i_{\emptyset,\widehat{j}_1\cap \widehat{j}_2}^{\infty}(z'_{j_1,j_2}\cdot 1_{T}))\cong \mathrm{cosoc}_{L_{\widehat{j}_1\cap\widehat{j}_2}}(i_{\emptyset,\widehat{j}_1\cap \widehat{j}_2}^{\infty}(z_{j_1,j_2}\cdot 1_{T})).\]
By Lemma~\ref{lem: structure of induction} this implies that $i_{\widehat{j}_1\cap \widehat{j}_2,\Delta}^{\infty}(\sigma_{j_1,j_2}^{\infty})$ is the unique quotient of $i_{\emptyset,\Delta}^{\infty}(z_{j_1,j_2}\cdot 1_{T})$ with socle $V_{\Delta\setminus D_R(z'_{j_1,j_2}),\Delta}^{\infty}$ (and cosocle $V_{D_R(z_{j_1,j_2}),\Delta}^{\infty}$). Using \ref{it: trivial PS 3} of Corollary~\ref{cor: trivial PS}, we only need to find an explicit $z'_{j_1,j_2}$ as above and compute $D_R(z_{j_1,j_2})$ and $\Delta\setminus D_R(z'_{j_1,j_2})$.

Let $y_{j_1,j_2}$ be the longest element in $\{x\in W(G)\mid \mathrm{Supp}(x)\subseteq [1,j_1-1]\}$. Note that $y_{j_1,j_2}\in W(L_{[1,j_1-1]})\cong W(\GL_{j_1})\subseteq W(L_{\widehat{j}_1\cap \widehat{j}_2})$ and that $D_R(y_{j_1,j_2})=[1,j_1-1]$. We set
\begin{equation*}
\delta_{j_1,j_2}\defeq |\mathrm{det}_{j_1}|_K^{j_1-j_2}\boxtimes_E |\mathrm{det}_{j_2-j_1}|_K^{j_1}\boxtimes_E 1_{\mathrm{GL}_{n-j_2}}: L_{\widehat{j}_1\cap \widehat{j}_2}\rightarrow E^\times.
\end{equation*}
We observe that $x_{j_1,j_2}\cdot 1_{T}=\delta_{j_1,j_2}|_{T}$ and that $y_{j_1,j_2}x_{j_1,j_2}\cdot 1_{T}=(y_{j_1,j_2}\cdot 1_{T})\otimes_E(\delta_{j_1,j_2}|_{T})$ (use $y_{j_1,j_2}\in W(L_{\widehat{j}_1\cap \widehat{j}_2})$). This implies
\begin{equation}\label{equ: twist induction}
i_{\emptyset,\widehat{j}_1\cap\widehat{j}_2}^{\infty}(y_{j_1,j_2}x_{j_1,j_2}\cdot 1_{T})\cong i_{\emptyset,\widehat{j}_1\cap\widehat{j}_2}^{\infty}((y_{j_1,j_2}\cdot 1_{T})\otimes_E\delta_{j_1,j_2}|_{T})\cong
i_{\emptyset,\widehat{j}_1\cap\widehat{j}_2}^{\infty}(y_{j_1,j_2}\cdot 1_{T})\otimes_E\delta_{j_1,j_2}.
\end{equation}
It follows from \ref{it: trivial PS 2} of Lemma~\ref{cor: trivial PS} that
\[\mathrm{cosoc}_{L_{\widehat{j}_1\cap\widehat{j}_2}}(i_{\emptyset,\widehat{j}_1\cap\widehat{j}_2}^{\infty}(y_{j_1,j_2}\cdot 1_{T}))\cong 1_{\mathrm{GL}_{j_1}}\boxtimes_E \mathrm{St}_{j_2-j_1}^{\infty}\boxtimes_E \mathrm{St}_{n-j_2}^{\infty},\]
which together with (\ref{equ: twist induction}) implies
\begin{equation}\label{equ: explicit induction cosocle}
\mathrm{cosoc}_{L_{\widehat{j}_1\cap\widehat{j}_2}}(i_{\emptyset,\widehat{j}_1\cap\widehat{j}_2}^{\infty}(y_{j_1,j_2}x_{j_1,j_2}\cdot 1_{T}))\cong \sigma_{j_1,j_2}^{\infty}.
\end{equation}
Hence we may take $z_{j_1,j_2}=y_{j_1,j_2}x_{j_1,j_2}$. Note that (\ref{equ: explicit induction cosocle}) and \ref{it: PS 1} of Lemma~\ref{lem: Jacquet of PS} imply
\begin{equation*}
\Sigma_{j_1,j_2}=W(L_{\widehat{j}_1})\cdot\cJ(\sigma_{j_1,j_2}^{\infty})=W(L_{\widehat{j}_1})y_{j_1,j_2}x_{j_1,j_2}\cdot 1_{T}=W(L_{\widehat{j}_1})x_{j_1,j_2}\cdot 1_{T}.
\end{equation*}
Let $y_{j_1,j_2}'\defeq w_{\widehat{j}_1\cap \widehat{j}_2}y_{j_1,j_2}$, it is the longest element of $W(\mathrm{GL}_{j_2-j_1}\times\mathrm{GL}_{n-j_2})$ (with $\mathrm{GL}_{j_2-j_1}\times\mathrm{GL}_{n-j_1}$ seen as Levi blocks of $L_{\widehat{j}_1\cap \widehat{j}_2}$), and we have $D_R(y_{j_1,j_2}')=[j_1+1,j_2-1]\sqcup [j_2+1,n]$. We have the following two cases.\bigskip

\textbf{Case $1$}: If $j_1=j_2$, then $z_{j_1,j_2}=y_{j_1,j_2}$ and we have $D_R(z_{j_1,j_2})=D_R(y_{j_1,j_2})=[1,j_1-1]$ and $\Delta\setminus D_R(z'_{j_1,j_2})=\Delta\setminus D_R(y'_{j_1,j_2})=\Delta\setminus[j_1+1,n]=[1,j_1]$. This gives \ref{it: explicit induction 2}.\bigskip

\textbf{Case $2$}: We assume $j_2>j_1$ and write $x$, $y$, $y'$, $z=yx$, $z'=y'x$ for $x_{j_1,j_2}$, $y_{j_1,j_2}$, $y_{j_1,j_2}'$, $z_{j_1,j_2}$, $z'_{j_1,j_2}$. As $x\in W^{\widehat{j}_1,\emptyset}$, we have $\ell(ux)=\ell(u)+\ell(x)$ for each $u\in W(L_{\widehat{j}_1})$. In particular, as $y,y'\in W(L_{\widehat{j}_1})$, we have $\ell(z)=\ell(y)+\ell(x)$ and $\ell(z')=\ell(y')+\ell(x)$.
We also have $D_R(z)\subseteq \mathrm{Supp}(yx)\subseteq [1,j_2-1]$ and $D_R(z')\subseteq \mathrm{Supp}(y'x)\subseteq \widehat{j_2}$. If $1\leq j\leq j_2-j_1-1$, we have $j+j_1>j_1$ and thus $ys_{j+j_1},y's_{j+j_1}\in W(L_{\widehat{j}_1})$. This together with $ys_{j+j_1}>y$ and $y's_{j+j_1}<y$ gives $yxs_j=ys_{j+j_1}x>yx$ and $y'xs_j=y's_{j+j_1}x<yx$, and thus $j\in D_R(z')\setminus D_R(z)$. If $j_2-j_1+1\leq j\leq j_2-1$, we have $j-(j_2-j_1)<j_1$ and thus $ys_{j-(j_2-j_1)},y's_{j-(j_2-j_1)}\in W(L_{\widehat{j}_1})$. This together with $ys_{j-(j_2-j_1)}<y$ and $y's_{j-(j_2-j_1)}>y'$ gives $yxs_j=ys_{j-(j_2-j_1)}x<yx$ and $y'xs_j=y's_{j-(j_2-j_1)}x>y'x$, and thus $j\in D_R(z)\setminus D_R(z')$. If $j_2+1\leq j\leq n-1$, we have $y'xs_j=y's_jx<y'x$ and $j\notin\mathrm{Supp}(z)$, and thus $j\in D_R(z')\setminus D_R(z)$.
To finish, it suffices to prove that $j_2-j_1\in D_R(x)$, which implies $j_2-j_1\in D_R(z)\cap D_R(z')$. Let $I\defeq [1,j_2-1]$, $I'\defeq [1,j_2-1]\setminus\{j_1\}$ and recall that $w_I$, $w_{I'}$ is the longest element of respectively $W(L_I)$, $W(L_{I'})$. By definition of $x$ we have $x=w_{I'}w_I$ (with $D_R(x)\subseteq \mathrm{Supp}(x)\subseteq I$) and thus $x s_{j'}=w_{I'}s_{j_2-j'}w_I$ for $j'\in I$. Note that $u<u'$ if and only if $uw_I>u'w_I$ for each $u,u'\in W(L_I)$. It follows that $j'\in D_R(x)\subseteq I$ if and only if $x s_{j'}<x$ (see (\ref{equ: right set})) if and only if $w_{I'}s_{n-j'}>w_{I'}$ if and only if $j'\in I\setminus w_I(D_R(w_{I'}))=I\setminus w_I(I')=\{j_2-j_1\}$.\bigskip

To sum up, we have proven $D_R(z)=[j_2-j_1,j_2-1]$, $D_R(z')=[1,j_2-j_1]\sqcup[j_2+1,n-1]$ and $\Delta\setminus D_R(z')=[j_2-j_1+1,j_2]\cap \Delta$. This gives \ref{it: explicit induction 1} and \ref{it: explicit induction 3}.
\end{proof}

For $(j_1,j_2)\in\mathbf{J}^{\infty}$, we define subsets $I_{j_1,j_2}^+,I_{j_1,j_2}^-\subseteq \Delta$ by
\begin{equation}\label{i+-chiant}
\mathrm{soc}_{G}(i_{\widehat{j}_1,\Delta}^{\infty}(\pi_{j_1,j_2}^{\infty}))\cong V_{I_{j_1,j_2}^+,\Delta}^{\infty},\ \ \mathrm{cosoc}_{G}(i_{\widehat{j}_1,\Delta}^{\infty}(\pi_{j_1,j_2}^{\infty}))\cong V_{I_{j_1,j_2}^-,\Delta}^{\infty}
\end{equation}
and we note that $I_{j_1,j_2}^+, I_{j_1,j_2}^-$ are explicitly given in Lemma~\ref{lem: explicit smooth induction}. We have $i_{\widehat{j}_1,\Delta}^{\infty}(\pi_{j_1,j_2}^{\infty})\cong Q_{\Delta}(I_{j_1,j_2}^+,I_{j_1,j_2}^-)$ by \ref{it: Hom cube 1} of Lemma~\ref{lem: Hom cube}. The following lemma will be used in \S\ref{subsec: square} to study extensions between (certain) non-locally algebraic Orlik-Strauch representations.

\begin{lem}\label{lem: connect Hom}
Let $(j_1,j_2),(j_1',j_2')\in\mathbf{J}^{\infty}$ with $(j_1,j_2)<(j_1',j_2')$. Then we have
\begin{enumerate}[label=(\roman*)]
\item \label{it: connect 1} $d(\pi_{j_1',j_2'}^{\infty},\pi_{j_1,j_2}^{\infty})=0$ if and only if $(j_1',j_2')\in\{(j_1+1,j_2+1), (j_1-1,j_2), (j_1,j_2+1)\}$;
\item \label{it: connect 2} $d(\pi_{j_1',j_2'}^{\infty},\pi_{j_1,j_2}^{\infty})=1$ if and only if $(j_1',j_2')\in\{(j_1+2,j_2+2), (j_1-2,j_2)\}$;
\item \label{it: connect 3} $\cJ(i_{\widehat{j}_1',\Delta}^{\infty}(\pi_{j_1',j_2'}^{\infty}))\cap \cJ(\pi_{j_1,j_2}^{\infty})=\emptyset$ if $(j_1',j_2')\notin\{(j_1+1,j_2+1), (j_1-1,j_2), (j_1,j_2+1)\}$.
\end{enumerate}
\end{lem}
\begin{proof}
We prove \ref{it: connect 1}. By \ref{it: Hom cube 1} of Lemma~\ref{lem: Hom cube}, it suffices to check the conditions $I_{j_1',j_2'}^-\!\in [I_{j_1,j_2}^+,I_{j_1,j_2}^-]$ and $I_{j_1,j_2}^+\in [I_{j_1',j_2'}^+,I_{j_1',j_2'}^-]$. This is a straightforward check using Lemma~\ref{lem: explicit smooth induction}:
\begin{itemize}
\item If $(j_1',j_2')\notin\{(j_1+1,j_2+1), (j_1-1,j_2), (j_1,j_2+1)\}$, then one can check that $[I_{j_1,j_2}^+,I_{j_1,j_2}^-]\cap [I_{j_1',j_2'}^+,I_{j_1',j_2'}^-]\ne \emptyset$ only when $j_1=j'_1=1$ and $j'_2>j_2+1$, in which case $[I_{j_1,j_2}^+,I_{j_1,j_2}^-]\cap [I_{j_1',j_2'}^+,I_{j_1',j_2'}^-]=\{\emptyset\}$ (i.e.~the empty set is the only subset of $\Delta$ in this intersection). Since $I_{j_1,j_2}^+, I_{j_1',j_2'}^-\ne \emptyset$, the conditions $I_{j_1',j_2'}^-\!\in [I_{j_1,j_2}^+,I_{j_1,j_2}^-]$ and $I_{j_1,j_2}^+\in [I_{j_1',j_2'}^+,I_{j_1',j_2'}^-]$ are never satisfied.
\item If $(j_1',j_2')\in\{(j_1+1,j_2+1), (j_1-1,j_2), (j_1,j_2+1)\}$, then one can check that $[I_{j_1,j_2}^+,I_{j_1,j_2}^-]\cap [I_{j_1',j_2'}^+,I_{j_1',j_2'}^-]=[I_{j_1,j_2}^+,I_{j_1',j_2'}^-]$, more precisely $[I_{j_1,j_2}^+,I_{j_1,j_2}^-]\cap [I_{j_1',j_2'}^+,I_{j_1',j_2'}^-]=\{[j_2-j_1+1,j_2],[j_2-j_1,j_2]\}$ if $(j_1',j_2')=(j_1+1,j_2+1)$, $[I_{j_1,j_2}^+,I_{j_1,j_2}^-]\cap [I_{j_1',j_2'}^+,I_{j_1',j_2'}^-]=\{[j_2-j_1+1,j_2],[j_2-j_1+1,j_2-1]\}$ if $(j_1',j_2')=(j_1-1,j_2)$, and $I_{j_1,j_2}^+=I_{j_1',j_2'}^-=[j_2-j_1+1,j_2]$ if $(j_1',j_2')=(j_1,j_2+1)$.
\end{itemize}

We prove \ref{it: connect 2}. By \ref{it: connect 1} it suffices to find all $(j_1,j_2),(j_1',j_2')\in\mathbf{J}^{\infty}$ satisfying $(j_1,j_2)< (j_1',j_2')$, $(j_1',j_2')\notin\{(j_1+1,j_2+1),(j_1-1,j_2),(j_1,j_2+1)\}$ and
\begin{equation}\label{equ: explicit sm Ext1}
\mathrm{Ext}_{G}^1(i_{\widehat{j}_1',\Delta}^{\infty}(\pi_{j_1',j_2'}^{\infty}),i_{\widehat{j}_1,\Delta}^{\infty}(\pi_{j_1,j_2}^{\infty}))^{\infty}
\end{equation}
is non-zero. We have two cases.\bigskip

\textbf{Case $1$}: If $[I_{j_1,j_2}^+,I_{j_1,j_2}^-]\cap [I_{j_1',j_2'}^+,I_{j_1',j_2'}^-]\neq \emptyset$, as $(j_1',j_2')\notin\{(j_1+1,j_2+1),(j_1-1,j_2),(j_1,j_2+1)\}$, we must have $j_1=j_1'=1$, $j_2'>j_2+1$ and $[I_{j_1,j_2}^+,I_{j_1,j_2}^-]\cap [I_{j_1',j_2'}^+,I_{j_1',j_2'}^-]=\{\emptyset\}$ (see the proof of \ref{it: connect 1}). Then $i_{\widehat{j}_1,\Delta}^{\infty}(\pi_{j_1,j_2}^{\infty})$ contains a length two subrepresentation $\sigma^{\infty}$ with socle $V_{\{j_2\},\Delta}^{\infty}$ and socle $\mathrm{St}_n^{\infty}=V_{\emptyset,\Delta}^{\infty}$. It follows from Lemma~\ref{lem: smooth devissage St} that $\mathrm{Ext}_{G}^1(V_{I,\Delta}^{\infty},\sigma^{\infty})^\infty=0$ for each $I\subseteq \Delta$ satisfying $d(I,\{j_2\})=d(I,\emptyset)+1$, or equivalently $j_2\notin I$. As any $I\in [I_{j_1',j_2'}^+,I_{j_1',j_2'}^-]$ satisfies $j_2\notin I$, we deduce by d\'evissage
\begin{equation}\label{equ: explicit sm Ext1 vanishing 1}
\mathrm{Ext}_{G}^1(i_{\widehat{j}_1',\Delta}^{\infty}(\pi_{j_1',j_2'}^{\infty}),\sigma^{\infty})^{\infty}=0.
\end{equation}
If $i_{\widehat{j}_1,\Delta}^{\infty}(\pi_{j_1,j_2}^{\infty})/\sigma^{\infty}\neq 0$, then it has length two with socle $V_{[j_2-1,j_2],\Delta}^{\infty}$ and cosocle $V_{\{j_2-1\},\Delta}^{\infty}$, and it follows from Lemma~\ref{lem: smooth devissage St} that
\begin{equation}\label{equ: explicit sm Ext1 vanishing 2}
\mathrm{Ext}_{G}^1(V_{\emptyset,\Delta}^{\infty},i_{\widehat{j}_1,\Delta}^{\infty}(\pi_{j_1,j_2}^{\infty})/\sigma^{\infty})^{\infty}=0.
\end{equation}
As any $\emptyset\neq I\in [I_{j_1',j_2'}^+,I_{j_1',j_2'}^-]$ satisfies $d(I,[j_2-1,j_2])\geq 3$ and $d(I,\{j_2-1\})\geq 2$, we deduce from the first statement of Lemma~\ref{lem: Ext sm St} that $\mathrm{Ext}_{G}^1(V_{I,\Delta}^{\infty}, i_{\widehat{j}_1,\Delta}^{\infty}(\pi_{j_1,j_2}^{\infty})/\sigma^{\infty})^{\infty}=0$ for such $I$, which together with (\ref{equ: explicit sm Ext1 vanishing 2}) implies $\mathrm{Ext}_{G}^1(i_{\widehat{j}_1',\Delta}^{\infty}(\pi_{j_1',j_2'}^{\infty}),i_{\widehat{j}_1,\Delta}^{\infty}(\pi_{j_1,j_2}^{\infty})/\sigma^{\infty})^{\infty}=0$. With (\ref{equ: explicit sm Ext1 vanishing 1}) we deduce that (\ref{equ: explicit sm Ext1}) is zero.\bigskip

\textbf{Case $2$}: We assume $[I_{j_1,j_2}^+,I_{j_1,j_2}^-]\cap [I_{j_1',j_2'}^+,I_{j_1',j_2'}^-]=\emptyset$, then we must have $\max\{j_1,j_1'\}\geq 2$. If (\ref{equ: explicit sm Ext1}) is non-zero, then $i_{\widehat{j}_1,\Delta}^{\infty}(\pi_{j_1,j_2}^{\infty})$ (resp.~$i_{\widehat{j}_1',\Delta}^{\infty}(\pi_{j_1',j_2'}^{\infty})$) must contain a Jordan-H\"older factor $V_{I_0,\Delta}^{\infty}$ (resp.~$V_{I_1,\Delta}^{\infty}$) such that $I_0\neq I_1$ and
\begin{equation}\label{equ: explicit sm Ext1 JH}
\mathrm{Ext}_{G}^1(V_{I_0,\Delta}^{\infty},V_{I_1,\Delta}^{\infty})^{\infty}\neq 0,
\end{equation}
which by Lemma~\ref{lem: Ext sm St} implies $d(I_0,I_1)=1$. We have the possibilities.
\begin{itemize}
\item If $(j_1',j_2')\in \{(j_1+2,j_2+2), (j_1-2,j_2)\}$, then one can check that $[I_{j_1,j_2}^+,I_{j_1,j_2}^-]\sqcup [I_{j_1',j_2'}^+,I_{j_1',j_2'}^-]=[I_{j_1,j_2}^+,I_{j_1',j_2'}^-]$, and thus $Q_\Delta(I_{j_1,j_2}^+,I_{j_1',j_2'}^-)$ gives a non-zero element of (\ref{equ: explicit sm Ext1}).
\item If $(j_1',j_2')\in \{(j_1+1,j_2+2), (j_1-1,j_2+1)\}$, then the only possible choice of $I_0,I_1$ for (\ref{equ: explicit sm Ext1 JH}) to hold is $I_0=I_{j_1',j_2'}^+$ and $I_1=I_{j_1,j_2}^-$, and there always exists $I_2\in [I_{j_1,j_2}^+,I_{j_1,j_2}^-]\setminus\{I_1\}$ such that $d(I_1,I_2)=1$ and $d(I_0,I_2)=2$. Write $\pi^{\infty}$ for the unique length two quotient of $i_{\widehat{j}_1,\Delta}^{\infty}(\pi_{j_1,j_2}^{\infty})$ with socle $V_{I_2,\Delta}^{\infty}$ and cosocle $V_{I_1,\Delta}^{\infty}$, then it follows from Lemma~\ref{lem: smooth devissage St} that $\mathrm{Ext}_{G}^1(V_{I_0,\Delta}^{\infty},\pi^{\infty})=0$. By a d\'evissage on both $i_{\widehat{j}_1',\Delta}^{\infty}(\pi_{j_1',j_2'}^{\infty})$ and $i_{\widehat{j}_1,\Delta}^{\infty}(\pi_{j_1,j_2}^{\infty})$ (using a filtration on $i_{\widehat{j}_1',\Delta}^{\infty}(\pi_{j_1',j_2'}^{\infty})$ with simple graded pieces and a filtration on $i_{\widehat{j}_1,\Delta}^{\infty}(\pi_{j_1,j_2}^{\infty})$ with $\pi^{\infty}$ being the unique non-simple graded piece), we deduce that (\ref{equ: explicit sm Ext1}) is zero.
\item If $(j_1',j_2')\notin\{(j_1+2,j_2+2), (j_1-2,j_2), (j_1+1,j_2+2), (j_1-1,j_2+1)\}$, then a pair $I_0,I_1$ as in (\ref{equ: explicit sm Ext1 JH}) does not exist since we have $d(I_0,I_1)\geq 2$ for any $I_0\in [I_{j_1,j_2}^+,I_{j_1,j_2}^-]$ and $I_1\in [I_{j_1',j_2'}^+,I_{j_1',j_2'}^-]$. Hence (\ref{equ: explicit sm Ext1}) is zero by d\'evissage.
\end{itemize}

Finally we prove \ref{it: connect 3}. By (\ref{equ: first adjunction}) we have
\[\Hom_{L_{\widehat{j}_1}}(J_{\Delta,\widehat{j}_1}(V_{I_{j_1,j_2}^+,\Delta}^{\infty}),\pi_{j_1,j_2}^{\infty})\cong \Hom_{G}(V_{I_{j_1,j_2}^+,\Delta}^{\infty}, i_{\widehat{j}_1,\Delta}^{\infty}(\pi_{j_1,j_2}^{\infty}))\neq 0\]
and thus $\cJ(\pi_{j_1,j_2}^{\infty})\subseteq \cJ(J_{\Delta,\widehat{j}_1}(V_{I_{j_1,j_2}^+,\Delta}^{\infty}))=\cJ(V_{I_{j_1,j_2}^+,\Delta}^{\infty})$ since $\pi_{j_1,j_2}^{\infty}$ is irreducible. Hence, to prove \ref{it: connect 3}, it suffices to show that $\cJ(i_{\widehat{j}_1',\Delta}^{\infty}(\pi_{j_1',j_2'}^{\infty}))\cap \cJ(V_{I_{j_1,j_2}^+,\Delta}^{\infty})=\emptyset$, or equivalently by \ref{it: PS 3} of Lemma~\ref{lem: Jacquet of PS} that $V_{I_{j_1,j_2}^+,\Delta}^{\infty}$ is not a constituent of $i_{\widehat{j}_1',\Delta}^{\infty}(\pi_{j_1',j_2'}^{\infty})\cong Q_{\Delta}(I_{j_1',j_2'}^+,I_{j_1',j_2'}^-)$, i.e.~$I_{j_1,j_2}^+\notin [I_{j_1',j_2'}^+,I_{j_1',j_2'}^-]$ by \ref{it: Hom cube 1} of Lemma~\ref{lem: Hom cube}. By the proof of \ref{it: connect 1}, we already see that $(j_1,j_2)< (j_1',j_2')$, $(j_1',j_2')\notin\{(j_1+1,j_2+1),(j_1-1,j_2),(j_1,j_2+1)\}$ and $[I_{j_1,j_2}^+,I_{j_1,j_2}^-]\cap [I_{j_1',j_2'}^+,I_{j_1',j_2'}^-]\neq \emptyset$ happens only when $j_1=j_1'=1$ and $j_2'>j_2+1$, in which case $[I_{j_1,j_2}^+,I_{j_1,j_2}^-]\cap [I_{j_1',j_2'}^+,I_{j_1',j_2'}^-]=\{\emptyset\}$. But we never have $I_{j_1,j_2}^+=\emptyset$, which proves the statement.
\end{proof}

The following lemma will be used in \S\ref{subsec: square} to study extensions between (certain) non-locally algebraic Orlik-Strauch representations and (certain) locally algebraic representations.

\begin{lem}\label{lem: distance from St}
Let $(j_1,j_2)\in\mathbf{J}^{\infty}$ and $I\subseteq \Delta$.
\begin{enumerate}[label=(\roman*)]
\item \label{it: distance St 1} We have $d(\pi_{j_1,j_2}^{\infty}, V_{I,\Delta}^{\infty})=0$ if and only if $I=I_{j_1,j_2}^-$, and $d(V_{I,\Delta}^{\infty}, \pi_{j_1,j_2}^{\infty})=0$ if and only if $I=I_{j_1,j_2}^+$.
\item \label{it: distance St 2} If $I\neq I_{j_1,j_2}^-$, then we have $d(\pi_{j_1,j_2}^{\infty}, V_{I,\Delta}^{\infty})=1$ if and only if $I\notin[I_{j_1,j_2}^+,I_{j_1,j_2}^-]$ and $d(I_{j_1,j_2}^-,I)=1$. Similarly, if $I\neq I_{j_1,j_2}^+$, then $d(V_{I,\Delta}^{\infty}, \pi_{j_1,j_2}^{\infty})=1$ if and only if $I\notin[I_{j_1,j_2}^+,I_{j_1,j_2}^-]$ and $d(I,I_{j_1,j_2}^+)=1$.
\end{enumerate}
\end{lem}
\begin{proof}
\ref{it: distance St 1} follows from $i_{\widehat{j}_1,\Delta}^{\infty}(\pi_{j_1,j_2}^{\infty})\cong Q_{\Delta}(I_{j_1,j_2}^+, I_{j_1,j_2}^-)$. \ref{it: distance St 2} follows from this and Lemma~\ref{lem: Ext of cube}, using moreover the isomorphisms
\[\mathrm{Ext}_{G}^1(i_{\widehat{j}_1}^{\infty}(\pi_{j_1,j_2}^{\infty}),V_{I,\Delta}^{\infty})^{\infty}\cong
\mathrm{Ext}_{G}^1(V_{I,\Delta}^{\infty},i_{\widehat{j}_1}^{\infty}(\pi_{j_1,j_2}^{\infty})^{\sim})^{\infty}
\cong \mathrm{Ext}_{G}^1(V_{I,\Delta}^{\infty},Q_{\Delta}(I_{j_1,j_2}^-,I_{j_1,j_2}^+))^{\infty}\]
for the second statement (see Remark \ref{rem: dual cube}).
\end{proof}

Recall that $\Sigma_{j_1,j_2}$ is defined in (\ref{equ: main coset}). The following lemma will be used in Lemma \ref{lem: sm cell} below (and in \S\ref{subsec: square}).

\begin{lem}\label{lem: coset intersection}
Let $(j_1,j_2),(j_1',j_2')\in\mathbf{J}^{\infty}$ with $(j_1,j_2)<(j_1',j_2')$.
\begin{enumerate}[label=(\roman*)]
\item \label{it: coset 1} We have $\Sigma_{j_1,j_2}\cap \Sigma_{j_1',j_2'}\neq \emptyset$ if and only if either $j_2'=j_2$ or $j_2'-j_1'=j_2-j_1$.
\item \label{it: coset 2} If $j_1'=j_1$ then $\Sigma_{j_1,j_2}\cap (w^{-1}\cdot\Sigma_{j_1,j_2'})\neq \emptyset$ for some $w\in W^{\widehat{j}_1,\widehat{j}_1}$ if and only if $w$ is the representative of $w_{[1,j_2'-1]}w_{[1,j_2-1]}^{-1}$, i.e.~$w_{[1,j_2'-1]}w_{[1,j_2-1]}^{-1}\in W(L_{\widehat{j}_1})w W(L_{\widehat{j}_1})$. Moreover $\Sigma_{j_1,j_2}\cap (s_{j_1}\cdot\Sigma_{j_1',j_2'})\neq \emptyset$ (i.e.~$w=s_{j_1}$) if and only if either $j_1=1$ or $j_2'=j_2+1$.
\end{enumerate}
\end{lem}
\begin{proof}
We prove \ref{it: coset 1}. Recall that $\Sigma_{j_1,j_2}=W(L_{\widehat{j}_1})x_{j_1,j_2}\cdot 1_{T}$, $\Sigma_{j'_1,j'_2}=W(L_{\widehat{j}'_1})x_{j'_1,j'_2}\cdot 1_{T}$ (see (\ref{equ: special Weyl element}) and Lemma \ref{lem: explicit smooth induction}). We need to find all pairs of $(j_1,j_2)<(j_1',j_2')$ such that
$\Sigma_{j_1,j_2}\cap \Sigma_{j_1',j_2'}\neq \emptyset$. As in the proof of \ref{it: Hom vers induction1} of Lemma \ref{lem: Hom vers induction}, one checks that $\Sigma_{j_1,j_2}\cap \Sigma_{j_1',j_2'}$ is a single left $W(L_{\widehat{j}_1\cap \widehat{j}_1'})$-coset if non-empty, and thus contains a unique element $x'\cdot 1_T$ with $x'\in W(L_{\widehat{j}_1})x_{j_1,j_2}\cap W(L_{\widehat{j}'_1})x_{j'_1,j'_2}$ of miminal length. We write $x'=x_{j_1,j_2}'x_{j_1,j_2}$ and $x'=x_{j_1',j_2'}'x_{j_1',j_2'}$ with $x_{j_1,j_2}'\in W(L_{\widehat{j}_1})$, $x_{j'_1,j'_2}'\in W(L_{\widehat{j}'_1})$ and both expressions reduced. Then it is clear that $x_{j_1,j_2}'\in W(L_{\widehat{j}_1})$ is minimal in $W(L_{\widehat{j}_1\cap \widehat{j}_1'})x_{j_1,j_2}'$, and that $x_{j_1',j_2'}'\in W(L_{\widehat{j}_1'})$ is minimal in $W(L_{\widehat{j}_1\cap \widehat{j}_1'})x_{j_1',j_2'}'$. We have the following cases.
\begin{itemize}
\item If $j_1'=j_1$, then $\Sigma_{j_1,j_2}$ and $\Sigma_{j_1,j_2'}$ are both $W(L_{\widehat{j}_1})$-cosets, and thus $\Sigma_{j_1,j_2}\cap \Sigma_{j_1',j_2'}\neq \emptyset$ if and only if $\Sigma_{j_1,j_2}=\Sigma_{j_1',j_2'}$ if and only if $x_{j_1,j_2}=x_{j_1',j_2'}$, which contradicts $(j_1,j_2)<(j_1',j_2')$. Thus $\Sigma_{j_1,j_2}\cap \Sigma_{j_1',j_2'}=\emptyset$.
\item If $j_1'<j_1$ and $j_2'=j_2$, then $W(L_{\widehat{j}_1})x_{j_1,j_2}$ and $W(L_{\widehat{j}_1'})x_{j_1',j_2'}=W(L_{\widehat{j}_1'})x_{j_1',j_2}$ both contain the (unique) maximal length element in the set $\{x\in W(G)\mid \mathrm{Supp}(x)\subseteq [1,j_2-1]\}$, and thus $\Sigma_{j_1,j_2}\cap \Sigma_{j_1',j_2'}\neq \emptyset$.
\item If $j_1'<j_1$ and $j_2'>j_2$, assume $\Sigma_{j_1,j_2}\cap \Sigma_{j_1',j_2'}\ne \emptyset$. The minimal length element $x'\in W(L_{\widehat{j}_1})x_{j_1,j_2}\cap W(L_{\widehat{j}_1'})x_{j_1',j_2'}$ satisfies $x'=x_{j_1',j_2'}'x_{j_1',j_2'}\geq x_{j_1',j_2'}$. Thus $j_2'-1\in[1,j_2'-1]=\mathrm{Supp}(x_{j_1',j_2'})\subseteq \mathrm{Supp}(x')$, which together with $x'=x_{j_1,j_2}'x_{j_1,j_2}$ and $j_2'-1\notin[1,j_2-1]=\mathrm{Supp}(x_{j_1,j_2})$ forces $j_2'-1\in\mathrm{Supp}(x_{j_1,j_2}')$. As $j_1'<j_1$ and $x_{j_1,j_2}'\in W(L_{\widehat{j}_1})$ is minimal in $W(L_{\widehat{j}_1\cap \widehat{j}_1'})x_{j_1,j_2}'$, we deduce $\mathrm{Supp}(x_{j_1,j_2}')\subseteq [1,j_1-1]$ (use that $L_{\widehat{j}_1}$ and $L_{\widehat{j}_1\cap \widehat{j}_1'}$ have the lower right block $\GL_{n-j_1}$ in common when $j'_1<j_1$). But $[1,j_1-1]$ doesn't contain $j'_2-1$ as $j_1<j'_2$, a contradiction. Thus $\Sigma_{j_1,j_2}\cap \Sigma_{j_1',j_2'}=\emptyset$.
\item If $j_1'>j_1$ and $j_2'-j_1'=j_2-j_1$, then one checks that $x_{j_1',j_2'}\in W(L_{\widehat{j}_1})x_{j_1,j_2}$ and thus $\Sigma_{j_1,j_2}\cap \Sigma_{j_1',j_2'}\neq \emptyset$.
\item If $j_1'>j_1$ and $j_2'-j_1'>j_2-j_1$, assume $\Sigma_{j_1,j_2}\cap \Sigma_{j_1',j_2'}\ne \emptyset$. The minimal length element $x'\in W(L_{\widehat{j}_1})x_{j_1,j_2}\cap W(L_{\widehat{j}_1'})x_{j_1',j_2'}$ satisfies $x'=x_{j_1',j_2'}'x_{j_1',j_2'}\geq x_{j_1',j_2'}$. On the other hand, we have $x'\leq x_{j_1,j_2}''x_{j_1,j_2}$ where $x_{j_1,j_2}''$ is the element of maximal length in $W^{\widehat{j}_1\cap \widehat{j}_1',\emptyset}(L_{\widehat{j}_1})$, thus $x_{j_1',j_2'}\leq x_{j_1,j_2}''x_{j_1,j_2}$. But it follows from \ref{it: saturated cell 1} of Lemma~\ref{lem: saturated cell} (together with $j_1'>j_1$ and $j_2'-j_1'>j_2-j_1$) that $x_{j_1',j_2'}\not\leq x_{j_1,j_2}''x_{j_1,j_2}$, a contradiction. Thus $\Sigma_{j_1,j_2}\cap \Sigma_{j_1',j_2'}=\emptyset$.
\end{itemize}
We have shown that $\Sigma_{j_1,j_2}\cap \Sigma_{j_1',j_2'}\neq \emptyset$ for $(j_1,j_2)<(j_1',j_2')$ if and only if either $j_2=j_2'$ or $j_2'-j_1'=j_2-j_1$, which is \ref{it: coset 1}.

We prove \ref{it: coset 2} and assume from now on $j_1=j_1'$. Note first that $\Sigma_{j_1,j_2}\cap (w^{-1}\cdot\Sigma_{j_1',j_2'})\neq \emptyset$ is equivalent to $\Sigma_{j_1,j_2}\subseteq W(L_{\widehat{j}_1})w^{-1}\Sigma_{j_1',j_2'}$, and hence can hold for at most one $w\in W^{\widehat{j}_1,\widehat{j}_1}$ (by the double coset decomposition, see for instance the end of the proof of Lemma \ref{lem: smooth geometric lemma}). From the definitions we have $w_{[1,j_2-1]}=w'x_{j_1,j_2}$ where $w'$ is the longest element in $W(L_{[1,j_2-1]\setminus\{j_1\}})$, in particular $w_{[1,j_2-1]}\in W(L_{\widehat{j}_1})x_{j_1,j_2}$. Likewise we have $w_{[1,j_2'-1]}\in W(L_{\widehat{j}_1'})x_{j_1,j_2'}$. Let $w\in W^{\widehat{j}_1,\widehat{j}_1}$ such that $w_{[1,j_2'-1]}w_{[1,j_2-1]}^{-1}\in W(L_{\widehat{j}_1})w W(L_{\widehat{j}_1})$, i.e.~$w_{[1,j_2'-1]}w_{[1,j_2-1]}^{-1}=y'wy$ for some $y,y'\in W(L_{\widehat{j}_1})$. Then $yw_{[1,j_2-1]}=w^{-1}(y')^{-1}w_{[1,j_2'-1]}$~with
\[\begin{array}{rcl}
yw_{[1,j_2-1]}&\in &W(L_{\widehat{j}_1})w_{[1,j_2-1]}=W(L_{\widehat{j}_1})x_{j_1,j_2}\\
w^{-1}(y')^{-1}w_{[1,j_2'-1]}&\in & w^{-1}W(L_{\widehat{j}_1})w_{[1,j_2'-1]}=w^{-1}W(L_{\widehat{j}_1})x_{j_1,j_2'},
\end{array}\]
and thus $\Sigma_{j_1,j_2}\cap (w^{-1}\cdot \Sigma_{j_1,j_2'})\neq \emptyset$. Finally, it follows from \ref{it: saturated cell 2} of Lemma~\ref{lem: saturated cell} that $w_{[1,j_2'-1]}w_{[1,j_2-1]}^{-1}\in W(L_{\widehat{j}_1})s_{j_1}W(L_{\widehat{j}_1})$ if and only if either $j_1=1$ or $j_2'=j_2+1$.
\end{proof}

Let $(j_1,j_2),(j_1',j_2')\in\mathbf{J}^{\infty}$ with $(j_1,j_2)<(j_1',j_2')$.
From (\ref{equ: first adjunction}) and (\ref{equ: smooth geometric lemma}) we have canonical isomorphisms for $k\geq 0$
\begin{equation}\label{equ: explicit sm cell}
\mathrm{Ext}_{G}^k(i_{\widehat{j}_1',\Delta}^{\infty}(\pi_{j_1',j_2'}^{\infty}),i_{\widehat{j}_1,\Delta}^{\infty}(\pi_{j_1,j_2}^{\infty}))^{\infty}\\
\cong \bigoplus_{w\in W^{\widehat{j}_1',\widehat{j}_1}}\mathrm{Ext}_{L_{\widehat{j}_1}}^k(i_{\widehat{j}_1',\widehat{j}_1,w}^{\infty}(J_{\widehat{j}_1',\widehat{j}_1,w}(\pi_{j_1',j_2'}^{\infty})),\pi_{j_1,j_2}^{\infty}).
\end{equation}
Using (\ref{equ: refine geometric lemma}) we see that (\ref{equ: explicit sm cell}) is non-zero for some $k\geq 0$ if only if there exists $w\in W^{\widehat{j}_1',\widehat{j}_1}$ (necessarily unique) such that $\Sigma_{j_1,j_2}\subseteq W(L_{\widehat{j}_1})w^{-1}\Sigma_{j_1',j_2'}$ (equivalently $\Sigma_{j_1,j_2}\cap (w^{-1}\cdot\Sigma_{j_1',j_2'})\neq \emptyset$) and such that \ref{equ: explicit sm cell} induces an isomorphism
\begin{equation}\label{equ: explicit single sm cell}
0\neq \mathrm{Ext}_{G}^k(i_{\widehat{j}_1',\Delta}^{\infty}(\pi_{j_1',j_2'}^{\infty}),i_{\widehat{j}_1,\Delta}^{\infty}(\pi_{j_1,j_2}^{\infty}))^{\infty}\\
\cong \mathrm{Ext}_{L_{\widehat{j}_1}}^k(i_{\widehat{j}_1',\widehat{j}_1,w}^{\infty}(J_{\widehat{j}_1',\widehat{j}_1,w}(\pi_{j_1',j_2'}^{\infty}))_{\cB^{\widehat{j}_1}_{\Sigma_{j_1,j_2}}},\pi_{j_1,j_2}^{\infty}).
\end{equation}

\begin{lem}\label{lem: sm cell}
Let $(j_1,j_2),(j_1',j_2')\in\mathbf{J}^{\infty}$ with $(j_1,j_2)<(j_1',j_2')$.
\begin{enumerate}[label=(\roman*)]
\item \label{it: sm cell 1}
We have
\begin{equation}\label{equ: sm cell 1}
\Hom_{L_{\widehat{j}_1}}\big(i_{\widehat{j}_1\cap\widehat{j}_1',\widehat{j}_1}^{\infty}(J_{\widehat{j}_1',\widehat{j}_1\cap\widehat{j}_1'}(\pi_{j_1',j_2'}^{\infty})),\pi_{j_1,j_2}^{\infty}\big)\neq 0
\end{equation}
if and only if $(j_1',j_2')\in \{(j_1+1,j_2+1), (j_1-1,j_2)\}$.
\item \label{it: sm cell 2}
When $j_1'=j_1$, we have
\begin{equation*}
\Hom_{L_{\widehat{j}_1}}\big(i_{\widehat{j}_1,\widehat{j}_1,s_{j_1}}^{\infty}(J_{\widehat{j}_1,\widehat{j}_1,s_{j_1}}(\pi_{j_1,j_2'}^{\infty})),\pi_{j_1,j_2}^{\infty}\big)\neq 0
\end{equation*}
if and only if $j_2'=j_2+1$.
\end{enumerate}
\end{lem}
\begin{proof}
We prove \ref{it: sm cell 1}. By definition (see (\ref{equ: Bruhat Jacquet}) and (\ref{equ: Bruhat induction})) we have
\[i_{\widehat{j}_1\cap\widehat{j}_1',\widehat{j}_1}^{\infty}(J_{\widehat{j}_1',\widehat{j}_1\cap\widehat{j}_1'}(\pi_{j_1',j_2'}^{\infty}))=i_{\widehat{j}_1',\widehat{j}_1,1}^{\infty}(J_{\widehat{j}_1',\widehat{j}_1,1}(\pi_{j_1',j_2'}^{\infty})),\]
hence by (\ref{equ: explicit single sm cell}) (and the line above it) we see that (\ref{equ: sm cell 1}) holds if and only if $\Sigma_{j_1,j_2}\cap \Sigma_{j_1',j_2'}\neq \emptyset$ and $d(\pi_{j_1',j_2'}^{\infty},\pi_{j_1,j_2}^{\infty})=0$. We thus deduce \ref{it: sm cell 1} from \ref{it: coset 1} of Lemma~\ref{lem: coset intersection} and \ref{it: connect 1} of Lemma~\ref{lem: connect Hom}.

We prove \ref{it: sm cell 2}. Likewise, (\ref{equ: sm cell 1}) holds if and only if $\Sigma_{j_1,j_2}\cap s_{j_1}\Sigma_{j_1,j_2'}\neq \emptyset$ and $d(\pi_{j_1,j_2'}^{\infty},\pi_{j_1,j_2}^{\infty})=0$. Thus \ref{it: sm cell 2} follows from the last statement in \ref{it: coset 2} of Lemma~\ref{lem: coset intersection} and \ref{it: connect 1} of Lemma~\ref{lem: connect Hom}.
\end{proof}

Let $(j_1,j_2)\in\mathbf{J}^{\infty}$ with $1<j_1<n-1$ and $j_2<n$ and thus $(j_1+1,j_2+1),(j_1-1,j_2),(j_1,j_2+1)\in\mathbf{J}^{\infty}$. We write
\begin{equation}\label{equ: I+-}
I_{+}\defeq \Delta\setminus\{j_1,j_1+1\},\ \ I_{-}\defeq \Delta\setminus\{j_1,j_1-1\},\ \ I_{\pm}\defeq \Delta\setminus\{j_1-1,j_1,j_1+1\}
\end{equation}
and we set
\begin{equation}\label{equ: sigmainter}
\left\{\begin{array}{lclrcl}
\Sigma_{+,0}&\defeq &\Sigma_{j_1+1,j_2+1}\cap \Sigma_{j_1,j_2+1}&\Sigma_{-,0}&\defeq &\Sigma_{j_1-1,j_2}\cap \Sigma_{j_1,j_2+1}\\
\Sigma_{+,1}&\defeq &\Sigma_{j_1+1,j_2+1}\cap \Sigma_{j_1,j_2}&\Sigma_{-,1}&\defeq &\Sigma_{j_1-1,j_2}\cap \Sigma_{j_1,j_2}\\
\Sigma_{\pm}&\defeq &\Sigma_{j_1,j_1}\cap s_{j_1}\cdot\Sigma_{j_1,j_2+1}.&&&
\end{array}\right.
\end{equation}
It follows from Lemma~\ref{lem: coset intersection} that $\Sigma_{\ast,0}$ and $\Sigma_{\ast,1}$ are single left $W(L_{I_{\ast}})$-cosets for $\ast\in\{+,-\}$, and $\Sigma_{\pm}$ is a single left $W(L_{I_{\pm}})$-coset. Note that $\pi_{j_1,j_2}^{\infty}$, $\pi_{j_1+1,j_2+1}^{\infty}$, $\pi_{j_1-1,j_2}^{\infty}$ and $\pi_{j_1,j_2+1}^{\infty}$ are all irreducible $G$-regular. In all cases (i) to (v) below, (both parts of) Lemma \ref{lem: smooth geometric lemma} together with (\ref{equ: first adjunction}) (for $k=0$) and \ref{it: connect 1} of Lemma~\ref{lem: connect Hom} imply that (\ref{equ: Hom induction cell}) holds, which gives us by the first statement in \ref{it: Hom vers induction1} of Lemma~\ref{lem: Hom vers induction} isomorphisms of irreducible $G$-regular representations:
\begin{enumerate}[label=(\roman*)]
\item $\pi_{+,0}^{\infty}\defeq J_{\widehat{j}_1,I_{+}}(\pi_{j_1,j_2+1}^{\infty})_{\cB^{I_{+}}_{\Sigma_{+,0}}}\cong J_{\Delta\setminus\{j_1+1\},I_{+}}'(\pi_{j_1+1,j_2+1}^{\infty})_{\cB^{I_{+}}_{\Sigma_{+,0}}}$;
\item $\pi_{-,0}^{\infty}\defeq J_{\widehat{j}_1,I_{-}}(\pi_{j_1,j_2+1}^{\infty})_{\cB^{I_{-}}_{\Sigma_{-,0}}}\cong J_{\Delta\setminus\{j_1-1\},I_{-}}'(\pi_{j_1-1,j_2}^{\infty})_{\cB^{I_{-}}_{\Sigma_{-,0}}}$;
\item $\pi_{+,1}^{\infty}\defeq J_{\Delta\setminus\{j_1+1\},I_{+}}(\pi_{j_1+1,j_2+1}^{\infty})_{\cB^{I_{+}}_{\Sigma_{+,1}}}\cong J_{\widehat{j}_1,I_{+}}'(\pi_{j_1,j_2}^{\infty})_{\cB^{I_{+}}_{\Sigma_{+,1}}}$;
\item $\pi_{-,1}^{\infty}\defeq J_{\Delta\setminus\{j_1-1\},I_{-}}(\pi_{j_1-1,j_2}^{\infty})_{\cB^{I_{-}}_{\Sigma_{-,1}}}\cong J_{\widehat{j}_1,I_{-}}'(\pi_{j_1,j_2}^{\infty})_{\cB^{I_{-}}_{\Sigma_{-,1}}}$;
\item $\pi_{\pm}^{\infty}\defeq J_{\widehat{j}_1,\widehat{j}_1,s_{j_1}}(\pi_{j_1,j_2+1}^{\infty})_{\cB^{I_{\pm}}_{\Sigma_{\pm}}}\cong J_{\widehat{j}_1,I_{\pm}}'(\pi_{j_1,j_2}^{\infty})_{\cB^{I_{\pm}}_{\Sigma_{\pm}}}$.
\end{enumerate}
Using the irreducibility of $\pi_{j_1,j_2}^{\infty}$, $\pi_{j_1,j_2+1}^{\infty}$, $\pi_{j_1+1,j_2+1}^{\infty}$ and $\pi_{j_1-1,j_2}^{\infty}$ as well as (\ref{equ: first adjunction}) and (\ref{equ: second adjunction}), we observe that $\pi_{j_1,j_2}^{\infty}\cong\mathrm{cosoc}_{L_{\widehat{j}_1}}(i_{I_{\ast},\widehat{j}_1}^{\infty}(\pi_{\ast,1}^{\infty}))$ and $\pi_{j_1,j_2+1}^{\infty}\cong\mathrm{soc}_{L_{\widehat{j}_1}}(i_{I_{\ast},\widehat{j}_1}^{\infty}(\pi_{\ast,0}^{\infty}))$
for $\ast\in\{+,-\}$. Similarly, we have $\mathrm{soc}_{L_{\widehat{j}_1}}(i_{I_{+},\Delta\setminus\{j_1+1\}}^{\infty}(\pi_{+,1}^{\infty}))\cong\pi_{j_1+1,j_2+1}^{\infty}\cong \mathrm{cosoc}_{L_{\widehat{j}_1}}(i_{I_{+},\Delta\setminus\{j_1+1\}}^{\infty}(\pi_{+,0}^{\infty}))$ and
$\mathrm{soc}_{L_{\widehat{j}_1}}(i_{I_{-},\Delta\setminus\{j_1-1\}}^{\infty}(\pi_{-,1}^{\infty}))\cong\pi_{j_1-1,j_2}^{\infty}\cong \mathrm{cosoc}_{L_{\widehat{j}_1}}(i_{I_{-},\Delta\setminus\{j_1-1\}}^{\infty}(\pi_{-,0}^{\infty}))$.
As $\pi_{\pm}^{\infty}$ is irreducible $G$-regular, for each $I\supseteq I_{\pm}$ the representation $i_{I_{\pm},I}^{\infty}(\pi_{\pm}^{\infty})$ is $G$-basic and thus multiplicity free with simple socle and cosocle by the last statement of \ref{it: basic as image} of Remark~\ref{rem: basic PS intertwine}.\bigskip

The last lemma below will be used in Lemma \ref{lem: square Ext1 vanishing} and in Lemma \ref{lem: conj Ext1 vanishing OS}, themselves used in the important Proposition \ref{prop: hard square}.

\begin{lem}\label{lem: coset square}
Let $(j_1,j_2)\in\mathbf{J}^{\infty}$ with $1<j_1<n-1$ and $j_2<n$.
\begin{enumerate}[label=(\roman*)]
\item \label{it: coset square 1} We have $s_{j_1}\cdot\Sigma_{j_1+1,j_2+1}=\Sigma_{j_1+1,j_2+1}$ and $s_{j_1}\cdot\Sigma_{j_1-1,j_2}=\Sigma_{j_1-1,j_2}$, which induce the following equalities for $\ast\in\{+,-\}$
\begin{equation}\label{equ: coset square intersection}
\Sigma_{\pm}=\Sigma_{j_1,j_2}\cap s_{j_1}\cdot\Sigma_{\ast,0}=\Sigma_{\ast,1}\cap s_{j_1}\cdot\Sigma_{j_1,j_2+1}.
\end{equation}
\item \label{it: coset square 2} We have
\[\pi_{\pm}^{\infty}\cong J_{I_{\ast},I_{\pm}}'(\pi_{\ast,1}^{\infty})_{\cB^{I_{\pm}}_{\Sigma_{\pm}}}\cong J_{I_{\ast},\widehat{j}_1,s_{j_1}}(\pi_{\ast,0}^{\infty})_{\cB^{I_{\pm}}_{\Sigma_{\pm}}}\]
and $\pi_{\ast,1}^{\infty}\cong \mathrm{cosoc}_{L_{I_{\ast}}}(i_{I_{\pm},I_{\ast}}^{\infty}(\pi_{\pm}^{\infty}))$ for $\ast\in\{+,-\}$. We also have
\[\pi_{j_1+1,j_2+1}^{\infty}\in\mathrm{JH}_{L_{\Delta\setminus\{j_1+1\}}}(i_{I_{\pm},\Delta\setminus\{j_1+1\}}^{\infty}(\pi_{\pm}^{\infty}))\text{ and }\pi_{j_1-1,j_2}^{\infty}\in\mathrm{JH}_{L_{\Delta\setminus\{j_1-1\}}}(i_{I_{\pm},\Delta\setminus\{j_1-1\}}^{\infty}(\pi_{\pm}^{\infty})).\]
\item \label{it: coset square 3} We have $J_{\Delta\setminus\{j_1+1\},I_{+}}'(\tau^{\infty})_{\cB^{I_{+}}_{\Sigma_{+,0}}}=0$ for each $\tau^{\infty}\in\mathrm{JH}_{L_{\Delta\setminus\{j_1+1\}}}(i_{I_{\pm},\Delta\setminus\{j_1+1\}}^{\infty}(\pi_{\pm}^{\infty}))$ satisfying $\tau^{\infty}< \pi_{j_1+1,j_2+1}^{\infty}$, and similarly $J_{\Delta\setminus\{j_1-1\},I_{-}}'(\tau^{\infty})_{\cB^{I_{-}}_{\Sigma_{-,0}}}=0$ for each \ \ $\tau^{\infty}\in\mathrm{JH}_{L_{\Delta\setminus\{j_1-1\}}}(i_{I_{\pm},\Delta\setminus\{j_1-1\}}^{\infty}(\pi_{\pm}^{\infty}))$ satisfying $\tau^{\infty}< \pi_{j_1-1,j_2}^{\infty}$.
\item \label{it: coset square 4} We have $J_{\widehat{j}_1,\widehat{j}_1,s_{j_1}}(\tau^{\infty})_{\cB^{I_{\pm}}_{\Sigma_{\pm}}}=0$ for each $\tau^{\infty}\in\mathrm{JH}_{L_{\widehat{j}_1}}(i_{I_{\ast},\widehat{j}_1}^{\infty}(\pi_{\ast,0}^{\infty})/\pi_{j_1,j_2+1}^{\infty})$ and $\ast\in \{+,-\}$.
\end{enumerate}
\end{lem}
\begin{proof}
We prove \ref{it: coset square 1}.
It follows from $s_{j_1}\in W(L_{\Delta\setminus\{j_1+1\}})\cap W(L_{\Delta\setminus\{j_1-1\}})$ that we have $s_{j_1}\cdot\Sigma_{j_1+1,j_2+1}=\Sigma_{j_1+1,j_2+1}$ and $s_{j_1}\cdot\Sigma_{j_1-1,j_2}=\Sigma_{j_1-1,j_2}$, which together with (\ref{equ: sigmainter}) give $s_{j_1}\cdot\Sigma_{+,0}=\Sigma_{j_1+1,j_2+1}\cap s_{j_1}\cdot\Sigma_{j_1,j_2+1}$ and $s_{j_1}\cdot\Sigma_{-,0}=\Sigma_{j_1-1,j_2}\cap s_{j_1}\cdot\Sigma_{j_1,j_2+1}$. Hence, in order to prove (\ref{equ: coset square intersection}), it suffices to prove that
\begin{equation}\label{equ: coset square inclusion}
\Sigma_{\pm}\subseteq \Sigma_{j_1+1,j_2+1}\cap \Sigma_{j_1-1,j_2}.
\end{equation}
Following the proof of $\Sigma_{j_1,j_2}\cap s_{j_1}\Sigma_{j_1,j_2+1}\neq \emptyset$ in \ref{it: coset 2} of Lemma~\ref{lem: coset intersection}, we write $w_{[1,j_2]}w_{[1,j_2-1]}^{-1}=y's_{j_1}y$ for some $y,y'\in W(L_{\widehat{j}_1})$. We have
\[yw_{[1,j_2-1]}=s_{j_1}(y')^{-1}w_{[1,j_2]}\in W(L_{\widehat{j}_1})x_{j_1,j_2}\cap s_{j_1}W(L_{\widehat{j}_1})x_{j_1,j_2+1}.\]
More precisely, as $w_{[1,j_2]}w_{[1,j_2-1]}^{-1}=s_1\cdots s_{j_2}$, we have $y=s_{j_1+1}\cdots s_{j_2}\in W(L_{\Delta\setminus\{j_1-1\}})$ (resp. $s_{j_1}(y')^{-1}=s_{j_1}s_{j_1-1}\cdots s_1\in W(L_{\Delta\setminus\{j_1+1\}})$),\ and \ since \ $w_{[1,j_2-1]}\in W(L_{\Delta\setminus\{j_1-1\}})x_{j_1-1,j_2}$ \ (resp. $w_{[1,j_2]}\in W(L_{\Delta\setminus\{j_1+1\}})x_{j_1+1,j_2+1}$)\ this \ implies \ $yw_{[1,j_2-1]}\in W(L_{\Delta\setminus\{j_1-1\}})x_{j_1-1,j_2}$ (resp.~$s_{j_1}(y')^{-1}w_{[1,j_2]}\in W(L_{\Delta\setminus\{j_1+1\}})x_{j_1+1,j_2+1}$). In other words, we have shown
\[yw_{[1,j_2-1]}\cdot 1_{T}=s_{j_1}(y')^{-1}w_{[1,j_2]}\cdot 1_{T}\in \Sigma_{j_1,j_2}\cap s_{j_1}\Sigma_{j_1,j_2+1}\cap \Sigma_{j_1-1,j_2}\cap \Sigma_{j_1+1,j_2+1}.\]
In \ particular, \ $\Sigma_{\pm}$ \ (resp. $\Sigma_{j_1+1,j_2+1}$, \ resp. $\Sigma_{j_1-1,j_2}$) \ is \ the \ unique \ $W(L_{I_{\pm}})$-coset (resp.~$W(L_{\Delta\setminus\{j_1+1\}})$-coset, resp.~$W(L_{\Delta\setminus\{j_1-1\}})$-coset) containing $yw_{[1,j_2-1]}\cdot 1_{T}$, which forces (\ref{equ: coset square inclusion}) and thus (\ref{equ: coset square intersection}).
This finishes the proof of \ref{it: coset square 1}.

We prove \ref{it: coset square 2}. On one hand, we have $J_{\widehat{j}_1,I_{\pm}}'(-)\cong J_{I_{\ast},I_{\pm}}'(J_{\widehat{j}_1,I_{\ast}}'(-))$ which together with $\Sigma_{\pm}\subseteq \Sigma_{\ast,1}$ (and the definition of $\pi_{\pm}^{\infty}$ and $\pi_{\ast,1}^{\infty}$) gives the isomorphisms
\begin{multline*}
\pi_{\pm}^{\infty}\cong J_{\widehat{j}_1,I_{\pm}}'(\pi_{j_1,j_2}^{\infty})_{\cB^{I_{\pm}}_{\Sigma_{\pm}}} \cong J_{I_{\ast},I_{\pm}}'(J_{\widehat{j}_1,I_{\ast}}'(\pi_{j_1,j_2}^{\infty}))_{\cB^{I_{\pm}}_{\Sigma_{\pm}}}\\
\cong J_{I_{\ast},I_{\pm}}'(J_{\widehat{j}_1,I_{\ast}}'(\pi_{j_1,j_2}^{\infty})_{\cB^{I_{\ast}}_{\Sigma_{\ast,1}}})_{\cB^{I_{\pm}}_{\Sigma_{\pm}}} \cong J_{I_{\ast},I_{+,-}}'(\pi_{\ast,1}^{\infty})_{\cB^{I_{\pm}}_{\Sigma_{\pm}}},
\end{multline*}
which together with (\ref{equ: second adjunction}) (and the irreducibility of $\pi_{\pm}^{\infty}$) give $\pi_{\ast,1}^{\infty}\cong \mathrm{cosoc}_{L_{I_{\ast}}}(i_{I_{\pm},I_{\ast}}^{\infty}(\pi_{\pm}^{\infty}))$
for $\ast\in \{+,-\}$. On the other hand, we have the equalities
\begin{equation}\label{equ: intersection of index}
\widehat{j}_1\cap s_{j_1}(\widehat{j}_1)=I_{\pm}=I_{\ast}\cap s_{j_1}(\widehat{j}_1),
\end{equation}
which together with (\ref{equ: Bruhat Jacquet}), (\ref{equ: Jacquet twist}) and $\Sigma_{\pm}\subseteq s_{j_1}\cdot\Sigma_{\ast,0}$ (and the definitions of $\pi_{\pm}^{\infty}$, $\pi_{\ast,0}^{\infty}$) give isomorphisms
\begin{multline*}
\pi_{\pm}^{\infty}\cong J_{\widehat{j}_1,\widehat{j}_1,s_{j_1}}(\pi_{j_1,j_2+1}^{\infty})_{\cB^{I_{\pm}}_{\Sigma_{\pm}}} \cong J_{I_{\ast},\widehat{j}_1,s_{j_1}}(J_{\widehat{j}_1,I_{\ast}}(\pi_{j_1,j_2+1}^{\infty}))_{\cB^{I_{\pm}}_{\Sigma_{\pm}}}\\
\cong J_{I_{\ast},\widehat{j}_1,s_{j_1}}(J_{\widehat{j}_1,I_{\ast}}(\pi_{j_1,j_2+1}^{\infty})_{\cB^{I_{\ast}}_{\Sigma_{\ast,0}}})_{\cB^{I_{\pm}}_{\Sigma_{\pm}}} \cong J_{I_{\ast},\widehat{j}_1,s_{j_1}}(\pi_{\ast,0}^{\infty})_{\cB^{I_{\pm}}_{\Sigma_{\pm}}}.
\end{multline*}
Finally, recall from the discussion right before this lemma that
\[\pi_{j_1+1,j_2+1}^{\infty}\cong\mathrm{soc}_{L_{\widehat{j}_1}}(i_{I_{+},\Delta\setminus\{j_1+1\}}^{\infty}(\pi_{+,1}^{\infty}))\text{ \ and\ }\pi_{j_1-1,j_2}^{\infty}\cong\mathrm{soc}_{L_{\widehat{j}_1}}(i_{I_{-},\Delta\setminus\{j_1-1\}}^{\infty}(\pi_{-,1}^{\infty})),\]
which together with $\pi_{\ast,1}^{\infty}\cong \mathrm{cosoc}_{L_{I_{\ast}}}(i_{I_{\pm},I_{\ast}}^{\infty}(\pi_{\pm}^{\infty}))$
for $\ast\in \{+,-\}$ implies that $\pi_{j_1+1,j_2+1}^{\infty}\in\mathrm{JH}_{L_{\Delta\setminus\{j_1+1\}}}(i_{I_{\pm},\Delta\setminus\{j_1+1\}}^{\infty}(\pi_{\pm}^{\infty}))$
and
$\pi_{j_1-1,j_2}^{\infty}\in\mathrm{JH}_{L_{\Delta\setminus\{j_1-1\}}}(i_{I_{\pm},\Delta\setminus\{j_1-1\}}^{\infty}(\pi_{\pm}^{\infty}))$.
This finishes the proof of \ref{it: coset square 2}.

We prove the first half of \ref{it: coset square 3} and leave the second half, which is similar, to the reader. We write $\delta\defeq \delta_{\widehat{j}_1,\widehat{j}_1,s_{j_1}}: L_{I_{\pm}}\rightarrow E^\times$ (see (\ref{equ: Jacquet twist})) and note that $\delta=\delta_{I_{\ast},\widehat{j}_1,s_{j_1}}$ for $\ast\in\{+,-\}$ (using (\ref{equ: intersection of index})).
Recall from the discussion right before this lemma that $i_{I_{\pm},\Delta\setminus\{j_1+1\}}^{\infty}(\pi_{\pm}^{\infty})$ is $G$-basic and multiplicity free with simple socle and cosocle, and from \ref{it: coset square 2} that $\pi_{j_1+1,j_2+1}^{\infty}\in\mathrm{JH}_{L_{\Delta\setminus\{j_1+1\}}}(i_{I_{\pm},\Delta\setminus\{j_1+1\}}^{\infty}(\pi_{\pm}^{\infty}))$.
\ We \ write \ $\pi_{j_1+1,j_2+1,-}^{\infty}$ \ for \ the \ unique \ subrepresentation \ of $i_{I_{\pm},\Delta\setminus\{j_1+1\}}^{\infty}(\pi_{\pm}^{\infty})$ with cosocle $\pi_{j_1+1,j_2+1}^{\infty}$.
Note that $\pi_{j_1+1,j_2+1,-}^{\infty}$ is $G$-basic by Corollary~\ref{cor: basic subquotient}, and that we have an injection $q_1: \pi_{j_1+1,j_2+1,-}^{\infty}\hookrightarrow i_{I_{\pm},\Delta\setminus\{j_1+1\}}^{\infty}(\pi_{\pm}^{\infty})$ and a surjection $q_2: \pi_{j_1+1,j_2+1,-}^{\infty}\twoheadrightarrow \pi_{j_1+1,j_2+1}^{\infty}$. We write $J_{+}(-)\defeq J_{\Delta\setminus\{j_1+1\},I_{+}}'(-)_{\cB^{I_{+}}_{\Sigma_{+,0}}}$ and recall that $J_{+}(-)$ is an exact functor. 
By the exactness of $J_{+}(-)$, the injection $q_1$ induces an injection $q_1': J_{+}(\pi_{j_1+1,j_2+1,-}^{\infty})\hookrightarrow J_{+}(i_{I_{\pm},\Delta\setminus\{j_1+1\}}^{\infty}(\pi_{\pm}^{\infty}))$, and the surjection $q_2$ induces a surjection $q_2': J_{+}(\pi_{j_1+1,j_2+1,-}^{\infty})\twoheadrightarrow J_{+}(\pi_{j_1+1,j_2+1}^{\infty})$.
For any $G$-basic representation $\sigma^{\infty}$ in $\mathrm{Rep}^{\infty}_{\rm{adm}}(L_{\Delta\setminus\{j_1+1\}})$, it follows from Lemma~\ref{lem: Jacquet basic} and Remark~\ref{rem: twist Jacquet basic} that $J_{+}(\sigma^{\infty})$ is either zero or $G$-basic. As $J_{+}(\pi_{j_1+1,j_2+1}^{\infty})=\pi_{+,0}^{\infty}$ by definition of $\pi_{+,0}^{\infty}$, we have $J_{+}(\pi_{j_1+1,j_2+1,-}^{\infty})\ne 0$, hence $J_{+}(\pi_{j_1+1,j_2+1,-}^{\infty})$ and $J_{+}(i_{I_{\pm},\Delta\setminus\{j_1+1\}}^{\infty}(\pi_{\pm}^{\infty}))$ are $G$-basic, and thus multiplicity free with simple socle and cosocle (last statement in \ref{it: basic as image} of Remark \ref{rem: basic PS intertwine}). In particular, the surjection $q_2'$ implies that $J_{+}(\pi_{j_1+1,j_2+1,-}^{\infty})$ has cosocle $\pi_{+,0}^{\infty}$.
It follows from (\ref{equ: second adjunction}) and then (\ref{equ: first adjunction}) that
\begin{multline}\label{equ: + - adjunction}
\Hom_{L_{I_{+}}}(\pi_{+,0}^{\infty}, J_{+}(i_{I_{\pm},\Delta\setminus\{j_1+1\}}^{\infty}(\pi_{\pm}^{\infty})))\\
\cong\Hom_{L_{\Delta\setminus\{j_1+1\}}}(i_{I_{+},\Delta\setminus\{j_1+1\}}^{\infty}(\pi_{+,0}^{\infty}),i_{I_{\pm},\Delta\setminus\{j_1+1\}}^{\infty}(\pi_{\pm}^{\infty}))\\
\cong\Hom_{L_{I_{\pm}}}(J_{\Delta\setminus\{j_1+1\},I_{\pm}}(i_{I_{+},\Delta\setminus\{j_1+1\}}^{\infty}(\pi_{+,0}^{\infty}))_{\cB^{I_{\pm}}_{\Sigma_{\pm}}},\pi_{\pm}^{\infty}).
\end{multline}
Applying \ref{it: sml2} of Lemma~\ref{lem: smooth geometric lemma} with $I_0=I_{+}$ and $I_1=I_{\pm}$, we see that 
\[J_{\Delta\setminus\{j_1+1\},I_{\pm}}(i_{I_{+},\Delta\setminus\{j_1+1\}}^{\infty}(\pi_{+,0}^{\infty}))_{\cB^{I_{\pm}}_{\Sigma_{\pm}}}\cong J_{I_{+},I_{\pm},s_{j_1}}(\pi_{+,0}^{\infty})_{\cB^{I_{\pm}}_{\Sigma_{\pm}}}.\]
Since $J_{I_{+},I_{\pm},s_{j_1}}(-)=J_{I_{+},\widehat{j}_1,s_{j_1}}(-)$ from (\ref{equ: Bruhat Jacquet}) and $\pi_{\pm}^{\infty}\cong J_{I_{+},\widehat{j}_1,s_{j_1}}(\pi_{+,0}^{\infty})_{\cB^{I_{\pm}}_{\Sigma_{\pm}}}$ from \ref{it: coset square 2}, this implies by (\ref{equ: + - adjunction})
\begin{equation}\label{equ: + - Hom}
\Hom_{L_{I_{+}}}(\pi_{+,0}^{\infty}, J_{+}(i_{I_{\pm},\Delta\setminus\{j_1+1\}}^{\infty}(\pi_{\pm}^{\infty})))\neq 0.
\end{equation}
As $J_{+}(i_{I_{\pm},\Delta\setminus\{j_1+1\}}^{\infty}(\pi_{\pm}^{\infty}))$ has simple socle by the previous discussion, by (\ref{equ: + - Hom}) it must have socle $\pi_{+,0}^{\infty}$, hence the same holds for its subrepresentation $J_{+}(\pi_{j_1+1,j_2+1,-}^{\infty})$. But recall that $J_{+}(\pi_{j_1+1,j_2+1,-}^{\infty})$ is multiplicity free with cosocle $\pi_{+,0}^{\infty}$ by the previous discussion, hence we have $J_{+}(\pi_{j_1+1,j_2+1,-}^{\infty})\cong \pi_{+,0}^{\infty}$ and in particular $J_{+}(\mathrm{ker}(q_2))=0$ by the exactness of $J_{+}(-)$. As $\tau^{\infty}\in\mathrm{JH}_{L_{\Delta\setminus\{j_1+1\}}}(\mathrm{ker}(q_2))$ if and only if $\tau^{\infty}\in\mathrm{JH}_{L_{\Delta\setminus\{j_1+1\}}}(i_{I_{\pm},\Delta\setminus\{j_1+1\}}^{\infty}(\pi_{\pm}^{\infty}))$ and $\tau^{\infty}<\pi_{j_1+1,j_2+1}^{\infty}$ (from the definition of the partial order on $\mathrm{JH}_{L_{\Delta\setminus\{j_1+1\}}}(i_{I_{\pm},\Delta\setminus\{j_1+1\}}^{\infty}(\pi_{\pm}^{\infty}))$ in \S\ref{generalnotation}), we obtain the first half of \ref{it: coset square 3} using $J_{+}(\mathrm{ker}(q_2))=0$.

We prove \ref{it: coset square 4}. As $\pi_{\pm}^{\infty}\cong J_{I_{\ast},\widehat{j}_1,s_{j_1}}(\pi_{\ast,0}^{\infty})_{\cB^{I_{\pm}}_{\Sigma_{\pm}}}$ by \ref{it: coset square 2}, we have $\pi_{\ast,0}^{\infty}\cong \mathrm{soc}_{L_{I_{\ast}}}(i_{I_{\pm},s_{j_1}(I_{\ast})}^{\infty}(\pi_{\pm}^{\infty}\otimes_E\delta^{-1})^{s_{j_1}})$ for $\ast\in \{+,-\}$ by the last statement in \ref{it: Hom vers induction1} of Lemma~\ref{lem: Hom vers induction} applied with $I_0=I_\ast$, $I_1=\widehat{j}_1$ and $w=s_{j_1}$, and thus $i_{I_{\ast},\widehat{j}_1}^{\infty}(\pi_{\ast,0}^{\infty})$ is a subrepresentation of
\[i_{I_{\ast},\widehat{j}_1}^{\infty}(i_{I_{\pm},s_{j_1}(I_{\ast})}^{\infty}(\pi_{\pm}^{\infty}\otimes_E\delta^{-1})^{s_{j_1}})\cong i_{I_{\pm},\widehat{j}_1}^{\infty}((\pi_{\pm}^{\infty}\otimes_E\delta^{-1})^{s_{j_1}}).\]
By the first statement of \ref{it: Hom vers induction2} of Lemma~\ref{lem: Hom vers induction} (which uses the last statement of \ref{it: Hom vers induction1} of Lemma~\ref{lem: Hom vers induction}) applied with $I_0=I_1=\widehat{j}_1$, $w=s_{j_1}$, $\pi^{\infty}=\pi_{\pm}^{\infty}$ and $\sigma_0^{\infty}=\tau^{\infty}$, we have $J_{\widehat{j}_1,\widehat{j}_1,s_{j_1}}(\tau^{\infty})_{\cB^{I_{\pm}}_{\Sigma_{\pm}}}=0$ for $\tau^{\infty}\in\mathrm{JH}_{L_{\widehat{j}_1}}(i_{I_{\pm},\widehat{j}_1}^{\infty}((\pi_{\pm}^{\infty}\otimes_E\delta^{-1})^{s_{j_1}})/\pi_{j_1,j_2+1}^{\infty})$, and in particular for $\tau^{\infty}\in\mathrm{JH}_{L_{\widehat{j}_1}}(i_{I_{\ast},\widehat{j}_1}^{\infty}(\pi_{\ast,0}^{\infty})/\pi_{j_1,j_2+1}^{\infty})$ as $i_{I_{\ast},\widehat{j}_1}^{\infty}(\pi_{\ast,0}^{\infty})\subseteq i_{I_{\pm},\widehat{j}_1}^{\infty}((\pi_{\pm}^{\infty}\otimes_E\delta^{-1})^{s_{j_1}})$.
\end{proof}

\newpage

\section{Results on Lie algebra cohomology groups}\label{sec: n coh}

We prove all results on $U(\fg)$-modules needed in \S\ref{sec: spectral seq} and especially in \S\ref{sec: extension}. In particular we prove many statements on unipotent cohomology groups and on various $\mathrm{Ext}$ groups of $U(\fg)$-modules, and we construct important explicit finite length $U(\fg)$-modules. We use the notation in \S\ref{generalnotation} and fix throughout a weight $\mu_0\in \Lambda^{\dom}$. For standard facts on the Bernstein-Gelfand-Gelfand category and on Kazhdan-Lusztig theory, we use \cite{Hum08} (the reader can find in \emph{loc.\ cit.}\ the original references where the results we use are actually proven).

\subsection{Categories of \texorpdfstring{$U(\fg)$}{U(g)}-modules}\label{subsec: category}

We introduce various abelian categories of $U(\fg)$-modules and prove several basic results on unipotent cohomology groups and $\mathrm{Ext}$ groups of $U(\fg)$-modules in these categories, and on the relations between the two.\bigskip

For $M$ in $\mathrm{Mod}_{U(\ft)}$ and $\mu\in \Lambda=X(T)$, we define $M_{\mu}\subseteq M$ as the maximal $U(\ft)$-submodule of $M$ on which $t-\mu(t)$ acts {\it nilpotently} for each $t\in\ft$. So $M_{\mu}$ is the generalized weight space attached to the weight $\mu$. Hence, we always have a $U(\ft)$-equivariant embedding
\begin{equation}\label{equ: sum of generalized wt space}
\bigoplus_{\mu\in \Lambda}M_{\mu}\hookrightarrow M.
\end{equation}
We define $\cC_{\rm{alg}}$ as the full subcategory of $\mathrm{Mod}_{U(\fg)}$ of those $M$ such that the embedding (\ref{equ: sum of generalized wt space}) is an isomorphism, and $\cC_{\rm{alg}}^{\rm{fin}}\subset \cC_{\rm{alg}}$ as the full subcategory of $M$ satisfying moreover $\Dim_E M_{\mu}<\infty$ for each $\mu\in \Lambda$. In particular each object of $\cC_{\rm{alg}}^{\rm{fin}}$ has countable dimension as $E$-vector space. We define $\tld{\cO}_{\rm{alg}}^{\fb}\subset \cC_{\rm{alg}}$ as the full subcategory of $M$ which are \emph{locally $\fb$-finite}, i.e.\ $M$ is the union of its finite dimensional $U(\fb)$-submodules. Note that the full subcategory of the category $\cO$ of \cite{Hum08} of objects with integral (equivalently algebraic) weights is the full subcategory $\cO_{\rm{alg}}^{\fb}\subset \tld{\cO}_{\rm{alg}}^{\fb}\cap \cC_{\rm{alg}}^{\rm{fin}}$ consisting of those $M$ which are moreover $U(\ft)$-semi-simple and finitely generated as $U(\fg)$-modules (\cite[\S 1.1]{Hum08}). We also write $\cO_{\rm{alg}}^{\fb,\infty}\subseteq \tld{\cO}_{\rm{alg}}^{\fb}$ for the full abelian subcategory consisting of finite length objects.\bigskip

For each $\mu\in \Lambda$, we have a Verma module $M(\mu)\defeq U(\fg)\otimes_{U(\fb)}\mu \in \cO_{\rm{alg}}^{\fb}$, which has an irreducible cosocle denoted by $L(\mu)$ (\cite[Thm.~1.2(f)]{Hum08}). Moreover recall that each simple object of $\cO_{\rm{alg}}^{\fb}$ has the form $L(\mu)$ for some $\mu\in \Lambda$ (\cite[\S 1.3]{Hum08}) and that each object of $\cO_{\rm{alg}}^{\fb}$ has finite length (\cite[\S 1.11]{Hum08}), and thus $\cO_{\rm{alg}}^{\fb}$ is the full subcategory of $\cO_{\rm{alg}}^{\fb,\infty}$ consisting of those $M$ which are moreover $U(\ft)$-semi-simple. As each $M$ in $\tld{\cO}_{\rm{alg}}^{\fb}$ contains at least one $\fb$-stable $E$-line with weight $\mu\in \Lambda$ (since $M$ is locally $\fb$-finite), there exists a non-zero map $M(\mu)\rightarrow M$. It follows that $\tld{\cO}_{\rm{alg}}^{\fb}$, $\cO_{\rm{alg}}^{\fb,\infty}$ and $\cO_{\rm{alg}}^{\fb}$ all share the same simple objects, namely the $L(\mu)$ for $\mu\in \Lambda$. (It is thus clear that $\cO_{\rm{alg}}^{\fb,\infty}\subseteq \cC_{\rm{alg}}^{\rm{fin}}$.) In the sequel, we write $N(\mu)$ for the kernel of the surjection $M(\mu)\twoheadrightarrow L(\mu)$.\bigskip

We say that a full subcategory $\cC\subseteq \mathrm{Mod}_{U(\fg)}$ is \emph{stable under extensions} if for each short exact sequence $0\rightarrow M_1\rightarrow M \rightarrow M_2\rightarrow 0$ in $\mathrm{Mod}_{U(\fg)}$ with $M_1, M_2$ in $\cC$, we necessarily have $M$ in $\cC$. We observe that both $\cC_{\rm{alg}}$ and $\cC_{\rm{alg}}^{\rm{fin}}$ are stable under extensions. It is not difficult to check that $\tld{\cO}_{\rm{alg}}^{\fb}$ and $\cO_{\rm{alg}}^{\fb,\infty}$ are also stable under extensions (by pull-back one can replace $M_2$ by a finite dimensional $U(\fb)$-submodule $M'_2\subseteq M_2$ , and then by induction on $\Dim_E M'_2$ one can reduce to $M'_2=E$ by taking a $\fb$-stable $E$-line in $M'_2$ and twisting, and then use that $\Ext_{U(\fb)}^1(E, M_1)=H^1(\fb, M_1)=\varinjlim H^1(\fb, M'_1)$ where the limit runs along the finite dimensional $U(\fb)$-submodules $M'_1$ of $M_1$). Recall however that the extension in $\mathrm{Mod}_{U(\fg)}$ of two objects of $\cO_{\rm{alg}}^{\fb}$ is not an object of $\cO_{\rm{alg}}^{\fb}$ in general.\bigskip

For each object $M$ of $\cC_{\rm{alg}}^{\rm{fin}}$, we set
\begin{equation}\label{tauduality}
M^\tau\defeq \bigoplus_{\mu\in \Lambda}M_\mu^\ast
\end{equation}
with $M_\mu^\ast\defeq \Hom_E(M_\mu, E)$ and we make $\fg$ act on $M^\tau$ by (see \cite[\S 3.2]{Hum08})
\[(x\cdot f)(v)\defeq f(\tau(x)\cdot v)\]
where $x\in \fg$ and $\tau: \fg\rightarrow \fg$ is Chevalley's anti-involution introduced at the end of \cite[\S 0.5]{Hum08} (for $\GL_n$ is it induced by the transpose map). Then $M\mapsto M^\tau$ defines a (contravariant) endo-functor of $\cC_{\rm{alg}}^{\rm{fin}}$ which is an exact involution, and thus a self-equivalence. As $\cC_{\rm{alg}}^{\rm{fin}}$ is stable under extensions, the functor $\tau$ induces an isomorphism for $M_1, M_2$ in $\cC_{\rm{alg}}^{\rm{fin}}$
\begin{equation}\label{equ: dual Ext1}
\mathrm{Ext}_{U(\fg)}^1(M_1,M_2)\cong \mathrm{Ext}_{U(\fg)}^1(M_2^\tau, M_1^\tau).
\end{equation}
Moreover, $\tau$ restricts to an exact involution of $\cO_{\rm{alg}}^{\fb}$ which satisfies
\begin{equation}\label{equ: simple self dual}
L(\mu)^\tau\cong L(\mu)
\end{equation}
for each $\mu\in \Lambda$ (\cite[Thm.\ 3.2]{Hum08}).\bigskip

For $I\subseteq \Delta$, we consider the full subcategory $\tld{\cO}_{\rm{alg}}^{\fp_{I}}\subseteq \tld{\cO}_{\rm{alg}}^{\fb}$ of those $M$ which are locally $\fp_{I}$-finite, i.e.\ equal to the union of their finite dimensional $U(\fp_{I})$-submodules. As $\fp_{I}=\fl_{I}\oplus \fn_{I}$ and the category of finite dimensional $U(\fl_{I})$-modules is semi-simple, $M\in \cO_{\rm{alg}}^{\fb}$ is locally $\fp_{I}$-finite if and only if the (underlying) $U(\fl_I)$-module $M$ is a direct sum of (simple) finite dimensional $U(\fl_{I})$-modules. We also define $\cO_{\rm{alg}}^{\fp_{I},\infty}\defeq \tld{\cO}_{\rm{alg}}^{\fp_{I}}\cap \cO_{\rm{alg}}^{\fb,\infty}$ and $\cO_{\rm{alg}}^{\fp_{I}}\defeq \tld{\cO}_{\rm{alg}}^{\fp_{I}}\cap \cO_{\rm{alg}}^{\fb}$. Replacing $\fg$ with $\fl_{I}$, we can define analogous full subcategories $\cC_{\fl_{I},\rm{alg}}$, $\cC_{\fl_{I},\rm{alg}}^{\rm{fin}}$, $\tld{\cO}_{\fl_{I},\rm{alg}}^{\fb_{I}}$, $\cO_{\fl_{I},\rm{alg}}^{\fb_{I},\infty}$, $\cO_{\fl_{I},\rm{alg}}^{\fb_{I}}$ of $\mathrm{Mod}_{U(\fl_{I})}$ with fully faithful embeddings $\cO_{\fl_{I},\rm{alg}}^{\fb_{I}}\hookrightarrow \cC_{\fl_{I},\rm{alg}}^{\rm{fin}} \hookrightarrow \cC_{\fl_{I},\rm{alg}}$ and $\cO_{\fl_{I},\rm{alg}}^{\fb_{I}}\hookrightarrow\cO_{\fl_{I},\rm{alg}}^{\fb_{I},\infty}\hookrightarrow \tld{\cO}_{\fl_{I},\rm{alg}}^{\fb_{I}}$. We also define $\tld{\cO}_{\fl_{I},\rm{alg}}^{\fl_{I}\cap \fp_{I'}}\hookrightarrow \tld{\cO}_{\fl_{I},\rm{alg}}^{\fb_{I}}$, $\cO_{\fl_{I},\rm{alg}}^{\fl_{I}\cap \fp_{I'},\infty}$ and $\cO_{\fl_{I},\rm{alg}}^{\fl_{I}\cap \fp_{I'}}$ for each $I'\subseteq I$ (note that $\fl_{I}\cap \fp_{I'}$ is a lower parabolic in $\fl_{I}$). We write $L^I(\mu)\in \cO_{\fl_{I},\rm{alg}}^{\fb_{I}}$ for the unique simple quotient of $U(\fl_{I})\otimes_{U(\fb_{I})}\mu$, and set
\begin{equation}\label{belongtoOb}
M^I(\mu)\defeq U(\fg)\otimes_{U(\fp_{I})}L^I(\mu).
\end{equation}
Be careful that we allow $M^I(\mu)$ and $L^I(\mu)$ to be defined for any $\mu\in \Lambda$, in particular $L^I(\mu)$ can be infinite dimensional. As $M^I(\mu)$ is a quotient of
\[M(\mu)=U(\fg)\otimes_{U(\fb)}\mu\cong U(\fg)\otimes_{U(\fp_{I})}(U(\fl_{I})\otimes_{U(\fb_{I})}\mu),\]
we see that $M^I(\mu)$ is in $\cO_{\rm{alg}}^{\fb}$ and has $L(\mu)$ as unique simple quotient. We write $N^I(\mu)$ for the kernel of the surjection $M^I(\mu)\twoheadrightarrow L(\mu)$. Moreover it follows from \cite[Prop.\ 9.3(e)]{Hum08} and \cite[Thm.\ 9.4]{Hum08} that $M^I(\mu)$ is in $\cO_{\rm{alg}}^{\fp_{I}}$ if and only if $L(\mu)$ is in $\cO_{\rm{alg}}^{\fp_{I}}$ if and only if $\mu\in \Lambda_I^{\dom}$ if and only if $L^I(\mu)$ is finite dimensional.\bigskip

For $w\in W(G)$ and $I\subseteq \Delta$ we set
\[\begin{array}{rclrclrcl}
M(w)&\defeq &M(w\cdot \mu_0),\ \ &L(w)&\defeq &L(w\cdot \mu_0),\ \ &N(w)&\defeq &N(w\cdot \mu_0)\\
M^I(w)&\defeq &M^I(w\cdot \mu_0),\ \ &L^I(w)&\defeq &L^I(w\cdot \mu_0),\ \ &N^I(w)&\defeq &N^I(w\cdot \mu_0).
\end{array}\]
Note that all Jordan-H\"older factors of $N(w)$ have the form $L(w')$ for some $w'>w$ (\cite[\S\S 5.1,5.2]{Hum08} and \cite[\S 8.3(a)]{Hum08}). For $w\in W(G)$, recall the sets $D_L(w),D_R(w)\subseteq \Delta$ from (\ref{equ: left set}) and (\ref{equ: right set}).

\begin{lem}\label{lem: dominance and left set}
Let $w\in W(G)$ be an element and $I\subseteq \Delta$ be a subset. Then $L(w)$ is in $\cO_{\rm{alg}}^{\fp_{I}}$ if and only if $I\cap D_L(w)=\emptyset$.
\end{lem}
\begin{proof}
We have $L(w)\in \cO_{\rm{alg}}^{\fp_{I}}$ if and only if $w\cdot \mu_0 \in \Lambda_I^{\dom}$ if and only if $\langle w\cdot \mu_0, \alpha^\vee\rangle \leq 0$ for $\alpha \in I$ if and only if $\langle w(\mu_0+\rho), \alpha^\vee\rangle \leq \langle \rho, \alpha^\vee\rangle = -1$ for $\alpha \in I$ if and only if $\langle \mu_0+\rho, w^{-1}(\alpha)^\vee\rangle < 0$ for $\alpha \in I$ if and only if $w^{-1}(\alpha)\in \Phi^+$ for $\alpha \in I$ (as $\mu_0\in \Lambda^{\dom}$ and hence $\mu_0+\rho\in \Lambda^{\dom}$) if and only if $s_{\alpha} w > w$ for $\alpha \in I$ (\cite[\S 0.3(4)]{Hum08}) if and only if $I\cap D_L(w)=\emptyset$.
\end{proof}

For $M_1, M_2$ in $\mathrm{Mod}_{U(\fg)}$ we write $\mathrm{Ext}_{U(\fg)}^k(M_1,M_2)$ for the extension groups computed in the category $\mathrm{Mod}_{U(\fg)}$. When $M_1, M_2$ are in $\cO_{\rm{alg}}^{\fb}$, we write $\mathrm{Ext}_{\cO_{\rm{alg}}^{\fb}}^k(M_1,M_2)$ for the extension groups computed in $\cO_{\rm{alg}}^{\fb}$ (which still has enough projective and injective objects by \cite[Thm.\ 3.8]{Hum08}). Given two objects $M_1, M_2$ in $\cO_{\rm{alg}}^{\fb}$, the fully faithful embedding $\cO_{\rm{alg}}^{\fb}\hookrightarrow \mathrm{Mod}_{U(\fg)}$ induces an injection
\[\mathrm{Ext}_{\cO_{\rm{alg}}^{\fb}}^1(M_1,M_2)\hookrightarrow \mathrm{Ext}_{U(\fg)}^1(M_1,M_2),\]
but the comparison between $\mathrm{Ext}_{\cO_{\rm{alg}}^{\fb}}^k(M_1,M_2)$ and $\mathrm{Ext}_{U(\fg)}^k(M_1,M_2)$ for $k\geq 2$ is more complicated in general. Since the dual functor $\tau: \cC_{\rm{alg}}^{\rm{fin}}\rightarrow \cC_{\rm{alg}}^{\rm{fin}}$ in (\ref{tauduality}) restricts to an exact involution of $\cO_{\rm{alg}}^{\fb}$, we have a canonical isomorphism for $M_1, M_2$ in $\cO_{\rm{alg}}^{\fb}$ and $k\geq 0$
\begin{equation}\label{equ: dual category O}
\mathrm{Ext}_{\cO_{\rm{alg}}^{\fb}}^k(M_1,M_2)\cong \mathrm{Ext}_{\cO_{\rm{alg}}^{\fb}}^k(M_2^\tau,M_1^\tau).
\end{equation}

For $M$ in $\mathrm{Mod}_{U(\fg)}$ and $I\subseteq \Delta$, we consider the Chevalley-Eilenberg complex (see for instance \cite[\S 3]{ST05})
\begin{equation}\label{equ: CE complex}
M\rightarrow \Hom_E(\fn_{I}, M)\rightarrow\cdots\rightarrow \Hom_E(\wedge^k\fn_{I},M)\rightarrow\cdots
\end{equation}
with $M$ in degree zero, and we define $H^k(\fn_{I}, M)$ as the cohomology group of this complex in degree $k\geq 0$. As the complex (\ref{equ: CE complex}) is $U(\fl_{I})$-equivariant, $H^k(\fn_{I},M)$ is naturally a $U(\fl_{I})$-module, and thus in particular a $U(\ft)$-module. We will use the following lemma.

\begin{lem}\label{lem: n coh wt}
Let $\mu,\mu'\in \Lambda$, $I\subseteq \Delta$ and $k\geq 0$. Assume that
\[H^k(\fn_{I},L(\mu))_{\mu'}\neq 0.\]
Then there exists $k$ distinct roots $\al_1,\dots,\al_k\in \Phi^+\setminus \Phi_I^+$ (no roots if $k=0$) such that
\begin{equation}\label{equ: n coh wt gap}
\mu'-\mu-\sum_{\ell=1}^k\al_{\ell}\in \Z_{\geq 0}\Phi^+.
\end{equation}
\end{lem}
\begin{proof}
The statement is obvious for $k=0$, hence we can assume $k\geq 1$. Let $\mu''\in \Lambda$ be an (integral) weight. We observe that $(\fn_{I})_{\mu''}\neq 0$ if and only if $\mu''=-\al$ for some $\al\in\Phi^+\setminus \Phi_I^+$, and thus $(\wedge^k\fn_{I})_{\mu''}\neq 0$ if and only if $\mu''=-\sum_{\ell=1}^k\al_{\ell}$ for $k$ distinct roots $\al_1,\dots,\al_k$ in $\Phi^+\setminus \Phi_I^+$. Consequently, $\Hom_E(\wedge^k\fn_{I}, L(\mu))_{\mu'}\neq 0$ if and only if there exists $\mu''\in \Lambda$ such that $(\wedge^k\fn_{I})_{\mu''}\neq 0$ and $L(\mu)_{\mu'+\mu''}\neq 0$. Note that $L(\mu)_{\mu'+\mu''}\neq 0$ implies $\mu'+\mu''-\mu\in\Z_{\geq 0}\Phi^+$. Hence, $\Hom_E(\wedge^k\fn_{I}, L(\mu))_{\mu'}\neq 0$ implies the existence of $k$ distinct roots $\al_1,\dots,\al_k\in \Phi^+\setminus \Phi_I^+$ such that (\ref{equ: n coh wt gap}) holds. As (\ref{equ: CE complex}) is $U(\ft)$-equivariant and $H^k(\fn_{I},L(\mu))$ is a subquotient of $\Hom_E(\wedge^k\fn_{I}, L(\mu))$ as $U(\ft)$-module, we obtain the statement.
\end{proof}

It is not difficult to check that, if $M$ is in $\cC_{\rm{alg}}$ (resp.~$\cC_{\rm{alg}}^{\rm{fin}}$, resp.~$\tld{\cO}_{\rm{alg}}^{\fb}$), then the $U(\fl_I)$-module $H^k(\fn_{I},M)$ is also in $\cC_{\fl_{I},\rm{alg}}$ (resp.~$\cC_{\fl_{I},\rm{alg}}^{\rm{fin}}$, resp.~$\tld{\cO}_{\fl_{I},\rm{alg}}^{\fb_{I}}$, for the latter recall that $\fn_I\subseteq \fb$). If $M$ is in $\cO_{\rm{alg}}^{\fb}$, we will prove in Proposition~\ref{lem: n coh O} below that $H^k(\fn_{I},M)\in \cO^{\fb_{I}}_{\fl_{I},\rm{alg}}$.\bigskip

For $I,I'\subseteq \Delta$, we have in $\fg$
\[\fu_{I}\cap \fn_{I'}=\fu_{I}\cap \fn_{I\cap I'}=\fl_{I}\cap \fn_{I\cap I'}=\fl_{I}\cap \fn_{I'},\]
which is the nilpotent radical of $\fl_{I}\cap \fp_{I'}=\fl_{I}\cap \fp_{I\cap I'}$, and
$H^k(\fl_{I}\cap \fn_{I'}, M_I)$ is naturally a $U(\fl_{I\cap I'})$-module for $M_I$ in $\mathrm{Mod}_{U(\fl_{I})}$. Since $\fn_{I}$ is an ideal in $\fn_{I\cap I'}$ with quotient naturally identified with $\fl_{I}\cap \fn_{I'}$, we have a $U(\fl_{I\cap I'})$-equivariant spectral sequence for $M$ in $\mathrm{Mod}_{U(\fg)}$ (see for instance \cite[\S 7.5]{Wei94})
\begin{equation}\label{equ: double n coh}
H^{\ell_1}(\fl_{I}\cap \fn_{I'}, H^{\ell_2}(\fn_{I}, M))\implies H^{\ell_1+\ell_2}(\fn_{I\cap I'}, M).
\end{equation}
When $I'=\emptyset$ and $I=\{j\}$ for some $j\in\Delta$, $\fl_{I}\cap \fn_{I'}=\fl_{\{j\}}\cap \fn =\fu_{\{j\}}$ is a $1$-dimensional Lie algebra, and (\ref{equ: double n coh}) induces the short exact sequence for $k\geq 1$ (see e.g.\ \cite[Exercise~5.2.1]{Wei94})
\begin{equation}\label{equ: special double n coh}
0\rightarrow H^1(\fu_{\{j\}}, H^{k-1}(\fn_{\{j\}}, M))\rightarrow H^k(\fu, M)\rightarrow H^0(\fu_{\{j\}}, H^k(\fn_{\{j\}}, M))\rightarrow 0.
\end{equation}

For $I\subseteq \Delta$, $M_I$ an $U(\fl_{I})$-module and $\xi: Z(\fl_{I})\rightarrow E$ a character (i.e.~an $E$-algebras homomorphism), we write $M_{I,\xi}$ for the maximal $U(\fl_{I})$-submodule of $M_I$ on which $z-\xi(z)$ acts nilpotently for each $z\in Z(\fl_{I})$. For a simple object $L^I(\mu)$ of $\tld{\cO}_{\fl_{I},\rm{alg}}^{\fb_{I}}$, we denote by $\xi_{\mu}: Z(\fl_{I})\rightarrow E$ the unique character such that $L^I(\mu)_{\xi_{\mu}}\neq 0$ (\cite[\S 1.7]{Hum08}). By Harish-Chandra's theorem (\cite[Thm.\ 1.10]{Hum08}) $\xi_{\mu}=\xi_{\mu'}$ if and only if there exists $w\in W(L_I)$ such that $\mu'=w\cdot \mu$. We consider the endo-functor
\[\mathrm{pr}_{\xi}: \cC_{\fl_{I},\rm{alg}}^{\rm{fin}}\rightarrow \cC_{\fl_{I},\rm{alg}}^{\rm{fin}},~ M_I\mapsto M_{I,\xi}.\]

\begin{lem}\label{lem: central component}
For $\xi: Z(\fl_{I})\rightarrow E$ the endo-functor $\mathrm{pr}_{\xi}$ has the following properties.
\begin{enumerate}[label=(\roman*)]
\item \label{it: block 2} The functor $\mathrm{pr}_{\xi}$ is exact and an idempotent, and we have a natural transformation
\[\mathrm{id}\cong \bigoplus_{\xi}\mathrm{pr}_{\xi}.\]
\item \label{it: block 3} Let $M_{I,1}, M_{I,2}$ in $\cC_{\fl_{I},\rm{alg}}^{\rm{fin}}$ such that $\mathrm{Ext}_{U(\fl_{I})}^k(M_{I,1},M_{I,2})\neq 0$ for some $k\geq 0$. Then there exists $\xi$ such that $M_{I,1,\xi}\neq 0$ and $M_{I,2,\xi}\neq 0$. The same statement holds for $M_{I,1}, M_{I,2}$ in $\cO_{\fl_{I},\rm{alg}}^{\fb_{I}}$ when $\mathrm{Ext}_{\cO_{\fl_{I},\rm{alg}}^{\fb_{I}}}^k(M_{I,1},M_{I,2})\neq 0$.
\item For each indecomposable object $M_I$ in $\cC_{\fl_{I},\rm{alg}}^{\rm{fin}}$, there exists a unique $\xi: Z(\fl_{I})\rightarrow E$ such that $M_I=M_{I,\xi}$.
\end{enumerate}
\end{lem}
\begin{proof}
It suffices to prove $M_I\cong \bigoplus_{\xi} M_{I,\xi}$ for $M_I$ in $\cC_{\fl_{I},\rm{alg}}^{\rm{fin}}$ and the rest is abstract non-sense. For $\mu\in \Lambda$, the action of $Z(\fl_I)$ on $M_I$ stabilizes $M_{I,\mu}$. As $M_{I,\mu}$ is finite dimensional by assumption, we deduce $M_{I,\mu}\cong \bigoplus_{\xi} M_{I,\xi,\mu}$. As $M_I\cong \bigoplus_{\mu\in \Lambda}M_{I,\mu}$ by assumption, we obtain $M_I\cong \bigoplus_{\xi} M_{I,\xi}$.
\end{proof}

Recall that there is an $E$-linear projection $\mathrm{pr}: U(\fg)\twoheadrightarrow U(\ft)$ obtained by sending to $0$ all monomials of $U(\fg)$ (in a standard Poincar\'e-Birkhoff-Witt basis associated to the decomposition $\fg = \fu^+ \oplus \ft \oplus \fu$) containing factors which are not in $\ft$. Using the decompositions $\fg = \fn_I^+ \oplus \fl_I \oplus \fn_I$ and $\fl_I = \fu_I^+ \oplus \ft \oplus \fu_I$ for $I\subseteq \Delta$, we see that pr uniquely factors as $U(\fg)\twoheadrightarrow U(\fl_I)\twoheadrightarrow U(\ft)$. It then follows from Harish-Chandra's theory (see \cite[\S\S 1.7, 1.9, 1.10]{Hum08}) and from the fact that $\rho-\rho_I$ is invariant under $W(L_I)$ that the above surjections restrict to injective morphisms of commutative $E$-algebras $Z(\fg) \hookrightarrow Z(\fl_I)\hookrightarrow Z(\ft)=U(\ft)$. We denote by $\psi_I:Z(\fg) \hookrightarrow Z(\fl_I)$ the first injection.\bigskip

We recall the following version of a classical theorem by Casselman-Osborne \cite[Thm.~2.6]{CO75}.

\begin{lem}\label{lem: HC n coh}
Let $M$ in $\cC_{\rm{alg}}^{\rm{fin}}$, $I\subseteq \Delta$ and $\xi: Z(\fl_I)\rightarrow E$. Then we have for $k\geq 0$
\begin{equation}\label{equ: center n coh}
H^k(\fn_I, M)_\xi\subseteq H^k(\fn_I, M_{\xi\circ\psi_I}).
\end{equation}
\end{lem}
\begin{proof}
Recall first that $H^k(\fn_I, M_{\xi\circ\psi_I})$ is an $U(\fl_I)$-submodule of $H^k(\fn_I, M)$ by \ref{it: block 2} of Lemma \ref{lem: central component} (applied to $I=\emptyset$ and $M$). Making explicit the action of $U(\fp_I)$ on the Cheval\-ley-Eilenberg complex (\ref{equ: CE complex}) (using $[\fp_I, \fn_I]\subseteq \fn_I$) and using that $\fn_I$ acts trivially on $H^k(\fn_I, \!M_{\xi\circ\psi_I})$, it is not difficult to check that the characters $\xi':Z(\fl_I)\rightarrow E$ such that $H^k(\fn_I, M_{\xi\circ\psi_I})_{\xi'}\ne 0$ satisfy $\xi'\circ \psi_I = \xi \circ \psi_I$. As any character of $Z(\fg)$ can be written $\xi''\circ \psi_I$ for some character $\xi''$ of $Z(\fl_I)$, decomposing $M$ using \ref{it: block 2} of Lemma \ref{lem: central component} and applying $H^k(\fn_I, -)_\xi$ implies the inclusion (\ref{equ: center n coh}).
\end{proof}

The following important proposition is probably known, but we couldn't find it in the published literature (it is also proved in \cite{BCGP25}).

\begin{prop}\label{lem: n coh O}
Let $M$ in $\cO^{\fb}_{\rm{alg}}$ and $I\subseteq \Delta$. Then we have for $k\geq 0$.
\begin{equation*}
H^k(\fn_{I}, M)\in \cO^{\fb_I}_{\fl_I,\rm{alg}}.
\end{equation*}
\end{prop}
\begin{proof}
As $M$ is in $\cO^{\fb}_{\rm{alg}}$, $M$ is locally $\fb_{I}$-finite and $U(\ft)$-semi-simple with finite dimensional weight spaces. Since $\wedge^{k}\fn_{I}$ is a finite dimensional $U(\fl_{I})$-module which is semi-simple as $U(\ft)$-module, $\Hom_{E}(\wedge^{k}\fn_{I},M)$ is still locally $\fb_{I}$-finite and $\ft$-semi-simple with finite dimensional weight spaces. Hence so is $H^k(\fn_{I}, M)$ for $k\geq 0$ by $U(\fl_{I})$-equivariance of the Chevalley-Eilenberg complex (\ref{equ: CE complex}). It follows that any finitely generated $U(\fl_{I})$-submodule of $H^k(\fn_{I}, M)$ is necessarily in $\cO^{\fb_{I}}_{\fl_{I},\rm{alg}}$, and thus has finite length by \cite[\S 1.11]{Hum08}.

It remains to show that $H^k(\fn_{I}, M)$ itself has finite length. Since $M$ is in $\cO^{\fb}_{\rm{alg}}$, $M$ has finite length and hence by \ref{it: block 2} of Lemma~\ref{lem: central component} there exists a \emph{finite} set $\Omega$ of characters $\xi: Z(\fg)\rightarrow E$ such that $M\cong \bigoplus_{\xi\in\Omega}M_{\xi}$ wich $M_{\xi}\neq 0$. Let $\xi_I: Z(\fl_I)\rightarrow E$ be a character such that $H^k(\fn_{I}, M)_{\xi_I}\neq 0$. By Lemma~\ref{lem: HC n coh} we know that $M_{\xi_I\circ \psi_I}\neq 0$ and thus $\xi_I\circ \psi_I \in\Omega$. By \ref{it: block 2} of Lemma~\ref{lem: central component} applied to the object $H^k(\fn_{I}, M)$ of $\cC_{\fl_{I},\rm{alg}}^{\rm{fin}}$ we deduce an isomorphism of $U(\fl_{I})$-modules where $\Omega_I\defeq \{\xi_I: Z(\fl_I)\rightarrow E\mid \xi_I\circ \psi_I\in \Omega\}$:
\begin{equation*}
H^k(\fn_{I}, M)\cong \bigoplus_{\xi_I\in\Omega_I}H^k(\fn_{I}, M)_{\xi_I}.
\end{equation*}
Write $U(\ft)$ as a polynomial algebra $E[t_1,\dots,t_n]$, then there are constants $\rho_1,\dots,\rho_n\in E$ such that, for any $\mu: U(\ft)\rightarrow E$, the polynomial $\prod_{i=1}^n(X-\rho_i-\mu(t_i))\in E[X]$ only depends on $\mu\vert_{Z(\fg)}$ by \cite[Thm.~1.10(a)]{Hum08} (the $\rho_i$ being related to the shift by $\rho$ in \emph{loc.\ cit.}). Since it has finitely many roots, we deduce that there is only a finite number of characters $\mu : U(\ft)\rightarrow E$ with a given $\mu\vert_{Z(\fg)}$. Since any character of $Z(\fl_I)$ is the restriction of a character of $U(\ft)$, it follows that the set $\{\xi_I\mid \xi_I\circ \psi_I=\xi\}$ is \emph{a fortiori} finite for each $\xi\in\Omega$, hence $\Omega_I$ is again finite (as $\Omega$ is). Now, assume on the contrary that $H^k(\fn_{I}, M)$ has infinite length. Then there exists $\xi_I\in \Omega_I$ such that $H^k(\fn_{I}, M)_{\xi_I}$ has infinite length. Using \cite[Thm.~1.10(b)]{Hum08}, we deduce that there exists at least one $\mu\in\Lambda$ such that $L^I(\mu)$ appears infinitely many times as a subquotient of (finitely generated $U(\fl_{I})$-submodules of) $H^k(\fn_{I}, M)_{\xi_I}$. But $H^k(\fn_{I}, M)_{\xi_I}$ is $\ft$-semi-simple, so the infinite multiplicity of $L^I(\mu)$ in $H^k(\fn_{I}, M)_{\xi_I}$ forces the $\mu$-weight space of $H^k(\fn_{I}, M)_{\xi_I}$ to have infinite dimension, a contradiction. Consequently $H^k(\fn_{I}, M)$ has finite length and thus lies in $\cO^{\fb_{I}}_{\fl_{I},\rm{alg}}$.
\end{proof}

For $M$ in $\mathrm{Mod}_{U(\fg)}$ and $M_I$ in $\mathrm{Mod}_{U(\fl_{I})}$, recall the Hochschild-Serre spectral sequence (\cite[\S 7.5]{Wei94})
\begin{equation}\label{equ: g spectral seq}
\mathrm{Ext}_{U(\fl_{I})}^{\ell_1}(M_I, H^{\ell_2}(\fn_{I},M))\implies \mathrm{Ext}_{U(\fp_{I})}^{\ell_1+\ell_2}(M_I, M)\cong \mathrm{Ext}_{U(\fg)}^{\ell_1+\ell_2}(U(\fg)\otimes_{U(\fp_{I})}M_I, M)
\end{equation}
where the last isomorphism is Shapiro's lemma for Lie algebra cohomology. In particular, we have a canonical isomorphism
\begin{equation}\label{equ: g spectral seq 0}
\Hom_{U(\fl_{I})}(M_I, H^0(\fn_{I},M))\cong \Hom_{U(\fg)}(U(\fg)\otimes_{U(\fp_{I})}M_I, M),
\end{equation}
and an exact sequence
\begin{multline}\label{equ: 5terms}
0\rightarrow\mathrm{Ext}_{U(\fl_{I})}^{1}(M_I, H^{0}(\fn_{I},M))\rightarrow \mathrm{Ext}_{U(\fg)}^{1}(U(\fg)\otimes_{U(\fp_{I})}M_I, M)\\
\rightarrow\Hom_{U(\fl_{I})}(M_I, H^{1}(\fn_{I},M))\rightarrow \mathrm{Ext}_{U(\fl_{I})}^{2}(M_I, H^{0}(\fn_{I},M)).
\end{multline}

\begin{lem}
Let $M_I$ in $\cO_{\fl_{I},\rm{alg}}^{\fb_{I}}$ and $M$ in $\cO_{\rm{alg}}^{\fb}$. Then we have a spectral sequence
\begin{equation}\label{equ: category O HS}
\mathrm{Ext}_{\cO^{\fb_{I}}_{\fl_{I},\rm{alg}}}^{\ell_1}(M_I,H^{\ell_2}(\fn_{I},M))\implies \mathrm{Ext}_{\cO_{\rm{alg}}^{\fb}}^{\ell_1+\ell_2}(U(\fg)\otimes_{U(\fp_{I})}M_I,M).
\end{equation}
\end{lem}
\begin{proof}
Note first that $H^{\ell_2}(\fn_{I},M)$ is indeed in $\cO^{\fb_{I}}_{\fl_{I},\rm{alg}}$ by Lemma \ref{lem: n coh O}. For $I\subseteq \Delta$ we use in this proof the extension groups $\mathrm{Ext}_{\fl_I,\ft}^k(-,-)$ computed in the abelian category of $(\fl_I,\ft)$-modules as defined in \cite[\S I.2]{BW00}, i.e.\ the full subcategory of $\mathrm{Mod}_{U(\fl_I)}$ of $U(\fl_I)$-modules which are semi-simple as $U(\ft)$-modules (i.e.~such that (\ref{equ: sum of generalized wt space}) is an isomorphism with all generalized weight spaces $M_\mu$ being weight spaces). We have analogous groups replacing $(\fl_I, \ft)$ by $(\fp_I, \ft)$. We first need to recall a little background.

Let $H^k(\fl_I,\ft,-)\defeq \mathrm{Ext}_{\fl_I,\ft}^k(E,-)$ and $H^k(\fp_I,\ft,-)\defeq \mathrm{Ext}_{\fp_I,\ft}^k(E,-)$. For $M, N$ in $\mathrm{Mod}_{U(\fl_I)}$, we endow $\Hom_E(N,M)$ with the unique structure of left $U(\fl_I)$-module such that $(\mathfrak{x}\cdot f)(n)\defeq \mathfrak{x}(f(n)) - f(\mathfrak{x}(n))$ for $\mathfrak{x} \in \fl_I$, $f\in \Hom_E(N,M)$ and $n\in N$. If $M$ is in $\mathrm{Mod}_{U(\fp_I)}$, by the same formula $\Hom_E(N,M)$ is naturally a left $U(\fp_I)$-module (with $U(\fp_I)$ acting on $N$ via $U(\fp_I)\twoheadrightarrow U(\fl_I)$). If $M$ and $N$ are moreover semi-simple as $U(\ft)$-modules (in the above sense), then by standard homological arguments we have canonical isomorphisms $H^{\ell}(\fl_{I},\ft, \Hom_E(N, M)) \cong \mathrm{Ext}_{\fl_{I},\ft}^{\ell}(N,M)$ (or $H^{\ell}(\fp_{I},\ft, \Hom_E(N, M)) \cong \mathrm{Ext}_{\fp_{I},\ft}^{\ell}(N,M)$ if $M$ is in $\mathrm{Mod}_{U(\fp_I)}$) for $\ell\geq 0$.

Let $M_I$ in $\mathrm{Mod}_{U(\fl_I)}$, $M$ in $\mathrm{Mod}_{U(\fg)}$ and assume that $M$ and $M_I$ are semi-simple as $U(\ft)$-modules. As $\fn_{I}$ is an ideal in $\fp_{I}$ (with $\fl_{I}\cong \fp_{I}/\fn_{I}$) and $\fn_{I}\cap \ft=0$, it then follows from \cite[Thm.~I.6.5]{BW00} and \cite[Rem.~I.6.7]{BW00} that there is a Hochschild-Serre type spectral sequence
\begin{equation}\label{equ: relative HS sequence}
H^{\ell_1}(\fl_{I},\ft, H^{\ell_2}(\fn_{I}, \Hom_E(M_I,M)))\implies H^{\ell_1+\ell_2}(\fp_{I},\ft, \Hom_E(M_I,M)).
\end{equation}
As $M_I$ is semi-simple as $U(\ft)$-module, so is $U(\fg)\otimes_{U(\fp_{I})}M_I$, and we have canonical isomorphisms (the second being the usual Shapiro's lemma)
\begin{eqnarray}\label{equ: relative inflation}
\nonumber H^{\ell_1+\ell_2}(\fp_{I},\ft, \Hom_E(M_I,M))&\cong &\mathrm{Ext}_{\fp_I,\ft}^{\ell_1+\ell_2}(M_I,M)\\
&\cong &\mathrm{Ext}_{\fg,\ft}^{\ell_1+\ell_2}(U(\fg)\otimes_{U(\fp_{I})}M_I,M).
\end{eqnarray}
As $\fn_{I}$ acts trivially on $M_I$, we have a $U(\fl_{I})$-equivariant isomorphism
\[H^{\ell_2}(\fn_{I}, \Hom_E(M_I,M))\cong \Hom_E(M_I, H^{\ell_2}(\fn_{I},M)),\]
which together with (\ref{equ: relative inflation}) implies that (\ref{equ: relative HS sequence}) can be rewritten
\begin{equation}\label{equ: HS g t}
\mathrm{Ext}_{\fl_{I},\ft}^{\ell_1}(M_I,H^{\ell_2}(\fn_{I},M))\\
\implies \mathrm{Ext}_{\fg,\ft}^{\ell_1+\ell_2}(U(\fg)\otimes_{U(\fp_{I})}M_I,M).
\end{equation}

Now take $M_I$ in $\cO_{\fl_{I},\rm{alg}}^{\fb_{I}}$, $M$ in $\cO_{\rm{alg}}^{\fb}$. By Delorme's theorem (\cite[Thm.\ 6.15]{Hum08}) we have in that case for $k\geq 0$
\[\mathrm{Ext}_{\cO_{\rm{alg}}^{\fb}}^k(U(\fg)\otimes_{U(\fp_{I})}M_I,M)\cong \mathrm{Ext}_{\fg,\ft}^k(U(\fg)\otimes_{U(\fp_{I})}M_I,M),\]
and for $\ell_1,\ell_2\geq 0$
\[\mathrm{Ext}_{\cO^{\fb_{I}}_{\fl_{I},\rm{alg}}}^{\ell_1}(M_I,H^{\ell_2}(\fn_{I},M))\cong \mathrm{Ext}_{\fl_{I},\ft}^{\ell_1}(M_I,H^{\ell_2}(\fn_{I},M)).\]
Then (\ref{equ: HS g t}) gives the spectral sequence (\ref{equ: category O HS}).
\end{proof}

\begin{lem}\label{lem: HC}
Let $w\in W(G)$, $I\subseteq \Delta$ and $M$ in ${\cO}_{\fl_I,\rm{alg}}^{\fb_I}$. Assume there is $k\geq 0$ such that
\begin{equation}\label{equ: same block as simple}
\mathrm{Ext}_{U(\fl_I)}^k(L^I(w),M)\neq 0
\end{equation}
or
\begin{equation}\label{equ: same O block as simple}
\mathrm{Ext}_{{\cO}_{\fl_I,\rm{alg}}^{\fb_I}}^k(L^I(w),M)\neq 0
\end{equation}
or
\begin{equation}\label{equ: Hom u coh}
\Hom_{U(\ft)}(w\cdot\mu_0, H^k(\fu_I,M))\neq 0.
\end{equation}
Then there exists $x\in W(L_I)w$ such that
\[\Hom_{U(\fl_I)}(L^I(x),M)\neq 0.\]
\end{lem}
\begin{proof}
Let $\xi: Z(\fl_I)\rightarrow E$ be the unique homomorphism such that $L^I(w)_\xi\neq 0$.

We first prove that the three hypothesis all imply $M_\xi\ne 0$. If either (\ref{equ: same block as simple}) or (\ref{equ: same O block as simple}) holds, this follows from \ref{it: block 3} of Lemma~\ref{lem: central component}. Replacing $(\fg, \fl_I, \ft)$ by $(\fl_I, \ft, \ft)$, we have a spectral sequence analogous to (\ref{equ: category O HS}). But since $\Ext_{{\cO}_{\ft,\rm{alg}}^{\ft}}^k\!\!=0$ if $k>0$, it gives isomorphisms for $k\geq 0$
\[\Hom_{U(\ft)}(w\cdot\mu_0, H^k(\fu_I,M))\cong \mathrm{Ext}_{{\cO}_{\fl_I,\rm{alg}}^{\fb_I}}^k(U(\fl_I)\otimes_{U(\fb_I)}w\cdot\mu_0,M).\]
If (\ref{equ: Hom u coh}) holds, we thus obtain $\mathrm{Ext}_{{\cO}_{\fl_I,\rm{alg}}^{\fb_I}}^k(L^I(w'),M)\neq 0$ for some Jordan-H\"older factor $L^I(w')$ of $U(\fl_I)\otimes_{U(\fb_I)}w\cdot\mu_0$ (with $w'\in W(L_I)w$ and $L^I(w')_\xi\neq 0$ since all constituents of $U(\fl_I)\otimes_{U(\fb_I)}w\cdot\mu_0$ have the same infinitesimal character). As above with (\ref{equ: same O block as simple}), this implies $M_\xi\neq 0$.

We now prove the statement. Take $\mu\in \Lambda$ such that $\Hom_{U(\fl_I)}(L^I(\mu),M_\xi)\neq 0$, thus in particular $L^I(\mu)_\xi\neq 0$. From Harish-Chandra's theorem (\cite[Thm.~1.10]{Hum08}) we deduce $\mu=x\cdot\mu_0$ for some $x\in W(L_I)$. This finishes the proof.
\end{proof}

For $\mu, \mu'\in \Lambda$, (\ref{equ: g spectral seq 0}) applied with $\fl_I=\ft$ gives an isomorphism
\[\Hom_{U(\ft)}(\mu', H^0(\fu, L(\mu))) \cong \Hom_{U(\fg)}(M(\mu'), L(\mu)).\]
As $H^0(\fu, L(\mu))$ is a semi-simple $U(\ft)$-module and $\Hom_{U(\fg)}(M(\mu'), L(\mu))\ne 0$ if $\mu'\ne \mu$ and has dimension $1$ if $\mu'=\mu$, we deduce a $U(\ft)$-equivariant isomorphism $H^0(\fu, L(\mu))\cong \mu$. Similarly, we have a $U(\ft)$-equivariant isomorphism $H^0(\fu_{I}, L^I(\mu))\cong \mu$ for each $I\subseteq \Delta$.

\begin{lem}\label{lem: n coh dominance}
Let $I\subseteq \Delta$.
\begin{enumerate}[label=(\roman*)]
\item \label{it: dominance 1} For $I'\subseteq \Delta$, $M$ in $\cO_{\rm{alg}}^{\fp_{I'}}$ and $k\geq 0$, the $U(\fl_{I})$-module $H^k(\fn_{I}, M)$ is locally $\fl_{I}\cap \fp_{I'}$-finite.
\item \label{it: dominance 2} For $\mu\in \Lambda$ we have $H^0(\fn_{I},L(\mu))\cong L^I(\mu)$.
\item \label{it: dominance 3} For $\mu\in \Lambda$, the unique (by \ref{it: dominance 2}) $\xi: Z(\fl_{I})\rightarrow E$ such that $H^0(\fn_{I},L(\mu))_\xi\neq 0$ is such that $H^k(\fn_{I},L(\mu))_\xi=0$ for $k\geq 1$.
\end{enumerate}
\end{lem}
\begin{proof}
As $M$ is in $\cO_{\rm{alg}}^{\fp_{I'}}$, $M$ is locally $\fp_{I'}$-finite and {\it a fortiori} locally $\fl_{I}\cap \fp_{I'}$-finite. Hence so is $\Hom_E(\wedge^k\fn_{I},M)$ for $k\geq 0$ (as $\wedge^k\fn_{I}$ is finite dimensional), and \ref{it: dominance 1} follows by the $U(\fl_{I})$-equivariance of (\ref{equ: CE complex}).

We prove \ref{it: dominance 2}. Note first that for $\mu'\in \Lambda$
\begin{multline*}
\Hom_{U(\fl_{I})}(U(\fl_{I})\otimes_{U(\fb_{I})}\mu', H^0(\fn_{I}, L(\mu)))\cong \Hom_{U(\ft)}(\mu', H^0(\fu_{I}, H^0(\fn_{I}, L(\mu))))\\
\cong \Hom_{U(\ft)}(\mu', H^0(\fu, L(\mu))) \cong \Hom_{U(\ft)}(\mu',\mu)
\end{multline*}
is non-zero if and only if $\mu'=\mu$. This implies that $H^0(\fn_{I}, L(\mu))$ has simple socle $L^I(\mu)$. It is thus enough to prove that $H^0(\fn_{I}, L(\mu))$ is a highest weight $U(\fl_I)$-module of weight $\mu$. The following argument is due to Florian Herzig (note that we know that $H^0(\fn_{I}, L(\mu))$ is in $\cO^{\fb_I}_{\fl_I,\rm{alg}}$ by Proposition \ref{lem: n coh O}, but we need the above more precise statement). Recall that $-\Phi^+$ (resp.\ $-\Phi_I^+$) are the roots of $\fu$ (resp.\ $\fu_I$) and $L(\mu)=\sum_{\lambda \in \Z_{\geq 0}\Phi^+}L(\mu)_{\mu+\lambda}$. Consider the following $U(\ft)$-submodules of $L(\mu)$
\[L(\mu)'\defeq \!\!\!\sum_{\lambda \in \Z_{\geq 0}\Phi_I^+}L(\mu)_{\mu+\lambda} \mathrm{\ \ and\ \ }
L(\mu)''\defeq \!\!\!\sum_{\lambda \notin \Z_{\geq 0}\Phi_I^+}L(\mu)_{\mu+\lambda}.\]
We have $L(\mu) = L(\mu)' \oplus L(\mu)''$, $L(\mu)_\mu = L(\mu)'_\mu$ and $L(\mu)'=U(\fu_I^+)\cdot L(\mu)_\mu$. Since the action of $\fl_I$ modifies a weight by a character in $\Z \Phi_I^+$ and since $\Z_{\geq 0}\Phi^+\cap \Z\Phi_I^+=\Z_{\geq 0}\Phi_I^+$, we see that $L(\mu')$ and $L(\mu'')$ are $U(\fl_I)$-submodules of $L(\mu)$. Since the action of $\fn_I$ modifies a weight by a character in $-\Z_{\geq 0}(\Phi^+\setminus \Phi_I^+)$, we see that $\fn_I$ necessarily acts by $0$ on $L(\mu)_{\mu+\lambda}$ for $\lambda \in \Z_{\geq 0}\Phi_I^+$, i.e.\ $L(\mu)'\subseteq H^0(\fn_I, L(\mu))$. Assume that $L(\mu)'\subsetneq H^0(\fn_I, L(\mu))$, or equivalently $L(\mu)''\cap H^0(\fn_I, L(\mu))\ne 0$. The action of $\fb$ on $L(\mu)''\cap H^0(\fn_I, L(\mu))$ factors through $\fb_I$ and since $L(\mu)$ is locally $\fb$-finite, then $L(\mu)''\cap H^0(\fn_I, L(\mu))$ is locally $\fb_I$-finite. In particular it contains a non-zero maximal vector for the action of $\fu_I$, or equivalently $\fu$, i.e.\ $H^0(\fu, L(\mu)'')\ne 0$. However $H^0(\fu, L(\mu)'')\subseteq H^0(\fu, L(\mu))=L(\mu)_\mu = L(\mu)'_\mu$ which contradicts $L(\mu)' \cap L(\mu)''=0$. It follows that $L(\mu)'=H^0(\fn_I, L(\mu))$, hence $H^0(\fn_I, L(\mu))=U(\fu_I^+)\cdot L(\mu)_\mu$ is a highest weight $U(\fl_I)$-module of weight $\mu$, which finishes the proof of \ref{it: dominance 2}.

We prove \ref{it: dominance 3}. Recall first that $H^k(\fn_{I}, L(\mu))$ is in $\cO^{\fb_I}_{\fl_I,\rm{alg}}$ by Proposition~\ref{lem: n coh O}. Let $\xi: Z(\fl_{I})\rightarrow E$ be the unique homomorphism such that $L^I(\mu)_\xi\neq 0$ and assume that there exists $k\geq 1$ and a Jordan-H\"older factor $L^I(\mu')$ of $H^k(\fn_{I},L(\mu))$ such that $L^I(\mu')_\xi\neq 0$. Then we have $\mu'=w\cdot \mu$ for some $w\in W(L_I)$ by Harish-Chandra's theorem. As the weight space $H^k(\fn_{I},L(\mu))_{\mu'}$ is non-zero, by Lemma~\ref{lem: n coh wt} there exists distinct $\al_1,\dots,\al_k\in \Phi^+\setminus \Phi_I^+$ such that $\mu'-\mu-\sum_{\ell=1}^k\al_{\ell}\in \Z_{\geq 0}\Phi^+$. However, $\mu'=w\cdot \mu$ (with $w\in W(L_I)$) implies $\mu'-\mu\in \Z\Phi_I^+$, which is a contradiction as
\[\bigg(\Big(\sum_{\ell=1}^k\al_{\ell}\Big)+\Z_{\geq 0}\Phi^+\bigg) \cap \Z\Phi_I^+=\emptyset.\]
It follows that $H^k(\fn_{I},L(\mu))_\xi=0$ for $k\geq 1$.
\end{proof}

\begin{lem}\label{lem: Ext 1 degenerate}
Let $\mu\in \Lambda$, $I\subseteq \Delta$ and $M_I$ in $\cC_{\fl_{I},\rm{alg}}^{\rm{fin}}$ such that $M_I=M_{I,\xi}$ for some $\xi: Z(\fl_{I})\rightarrow E$.
\begin{enumerate}[label=(\roman*)]
\item \label{it: degenerate 1} The map
\begin{equation*}
d_2^{k,1}:\mathrm{Ext}_{U(\fl_{I})}^k(M_I, H^1(\fn_{I},L(\mu)))\rightarrow \mathrm{Ext}_{U(\fl_{I})}^{k+2}(M_I, H^0(\fn_{I},L(\mu)))
\end{equation*}
in (\ref{equ: g spectral seq}) is zero for $k\geq 0$. In fact either the source or the target of $d_2^{k,1}$ is $0$.
\item \label{it: degenerate 2} If $H^0(\fn_{I},L(\mu)))_\xi\neq 0$ then $\mathrm{Ext}_{U(\fl_{I})}^k(M_I, H^\ell(\fn_{I},L(\mu)))=0$ for $k\geq 0$, $\ell\geq 1$ and we have isomorphisms for $k\geq 0$
\begin{equation*}
\mathrm{Ext}_{U(\fl_{I})}^k(M_I, H^0(\fn_{I},L(\mu)))\cong \mathrm{Ext}_{U(\fg)}^k(U(\fg)\otimes_{U(\fp_{I})}M_I, L(\mu)).
\end{equation*}
\item \label{it: degenerate 3} If $H^0(\fn_{I},L(\mu))_\xi=0$ then $\mathrm{Ext}_{U(\fl_{I})}^k(M_I, H^0(\fn_{I},L(\mu)))=0$ for $k\geq 0$ and we have an isomorphism
\begin{equation*}
\Hom_{U(\fl_{I})}(M_I, H^1(\fn_{I},L(\mu)))\cong \mathrm{Ext}_{U(\fg)}^1(U(\fg)\otimes_{U(\fp_{I})}M_I, L(\mu))
\end{equation*}
and an injection
\begin{equation}\label{equ: Ext2 H1 embedding}
\mathrm{Ext}_{U(\fl_{I})}^1(M_I, H^1(\fn_{I},L(\mu)))\hookrightarrow \mathrm{Ext}_{U(\fg)}^2(U(\fg)\otimes_{U(\fp_{I})}M_I, L(\mu)).
\end{equation}
\item \label{it: degenerate 4} Analogous statements as \ref{it: degenerate 1}, \ref{it: degenerate 2}, \ref{it: degenerate 3} hold when $M_I$ is in $\cO_{\fl_{I},\rm{alg}}^{\fb_{I}}$ and replacing $\mathrm{Ext}_{U(\fl_{I})}^k$ and $\mathrm{Ext}_{U(\fg)}^k$ by respectively $\mathrm{Ext}_{\cO^{\fb_{I}}_{\fl_{I},\rm{alg}}}^k$ and $\mathrm{Ext}_{\cO_{\rm{alg}}^{\fb}}^k$.
\end{enumerate}
\end{lem}
\begin{proof}
Combining \ref{it: dominance 2} and \ref{it: dominance 3} of Lemma~\ref{lem: n coh dominance} with \ref{it: block 3} of Lemma~\ref{lem: central component} and using (\ref{equ: g spectral seq}), (\ref{equ: 5terms}) and (\ref{equ: category O HS}), we obtain \ref{it: degenerate 1}, \ref{it: degenerate 2} and the first two statements in \ref{it: degenerate 3} in both cases $M_I\in\cC_{\fl_{I},\rm{alg}}^{\rm{fin}}$ and $M_I\in\cO_{\fl_{I},\rm{alg}}^{\fb_{I}}$. For (\ref{equ: Ext2 H1 embedding}), note that the bottom line of the $E_2$-terms in the spectral sequence (\ref{equ: g spectral seq}) is identically $0$. This implies that $E_2^{1,1}=E_\infty^{1,1}=\mathrm{Ext}_{U(\fl_{I})}^1(M_I, H^1(\fn_{I},L(\mu)))$ is the first non-zero graded piece of the abutment filtration on $E_\infty^2=\mathrm{Ext}_{U(\fg)}^2(U(\fg)\otimes_{U(\fp_{I})}M_I, L(\mu))$, whence (\ref{equ: Ext2 H1 embedding}). When $M_I$ is in $\cO_{\fl_{I},\rm{alg}}^{\fb_{I}}$ the proof is analogous using (\ref{equ: category O HS}) instead.
\end{proof}

\subsection{Results on \texorpdfstring{$\mathrm{Ext}^1$}{Ext1} groups in the category \texorpdfstring{$\cO_{\rm{alg}}^{\fb}$}{Obalg}}\label{subsec: Ext O}

We prove results on $\mathrm{Ext}$ groups (mainly $\mathrm{Ext}^1$ groups) in the category $\cO_{\rm{alg}}^{\fb}$, in the category of all $U(\fg)$-modules, and on the comparison between the two.\bigskip

We denote by $w_0\in W(G)$ the element of maximal length. For $x,w\in W(G)$ such that $x\leq w$ we let $P_{x,w}(q)\in \Z[q]$ be the associated Kazhdan-Lusztig polynomial (\cite{KL79}) and let $P_{x,w}(q)\defeq 0$ if $x\not\leq w$. When $x\leq w$ recall that $\mathrm{deg}P_{x,w}(q)\leq \frac{1}{2}(\ell(w)-\ell(x)-1)$ if $x\ne w$ and $P_{x,w}=1$ if $\ell(w)\leq \ell(x)+2$. Moreover the Kazhdan-Lusztig conjectures (independently proved in the 1980s by Beilinson-Bernstein and Brylinski-Kashiwara, and reproved by several people ever since) imply the following formula in the Grothendieck group of $\cO_{\rm{alg}}^{\fb}$ (with obvious notation)
\[[L(w)]=\sum_{w\leq x}(-1)^{\ell(x)-\ell(w)}P_{xw_0,ww_0}(1)[M(x)].\]
Following \cite[Def.\ 1.2]{KL79}, we write $x\prec w$ if $x<w$ and $\mathrm{deg}P_{x,w}(q)= \frac{1}{2}(\ell(w)-\ell(x)-1)$ (so $x\prec w$ implies that $\ell(w)-\ell(x)$ is odd) and we define $\mu(x,w)$ as the leading coefficient of $P_{x,w}$ if $x\prec w$, and $\mu(x,w)\defeq 0$ otherwise. For instance if $x<w$ and $\ell(w)=\ell(x)+1$, we have $x\prec w$ (and $\mu(x,w)=1$).

\begin{lem}\label{lem: Ext with Verma}
Let $x,w\in W(G)$ and $\mu\in \Lambda$.
\begin{enumerate}[label=(\roman*)]
\item \label{it: Ext Verma 1} We have an isomorphism for $k\geq 0$
\begin{equation}\label{equ: Ext with Verma}
\mathrm{Ext}_{\cO_{\rm{alg}}^{\fb}}^k(M(\mu),L(w))\cong \Hom_{U(\ft)}(\mu, H^k(\fu,L(w))).
\end{equation}
Moreover, if (\ref{equ: Ext with Verma}) is non-zero then there exists $x'\geq w$ in $W(G)$ such that $\mu=x'\cdot\mu_0$.
\item \label{it: Ext Verma 2} The dimension of $\mathrm{Ext}_{\cO_{\rm{alg}}^{\fb}}^k(M(x),L(w))$ (for $k\geq 0$) is equal to the coefficient of the monomial of degree $\frac{1}{2}(\ell(x)-\ell(w)-k)$ in $P_{xw_0,ww_0}$.
\item \label{it: Ext Verma 3} If $k\geq 0$ is the minimal integer such that (\ref{equ: Ext with Verma}) is non-zero and if $x'$ is as in \ref{it: Ext Verma 1}, then we have for $k'\leq k$
\[\mathrm{Ext}_{U(\fg)}^{k'}(M(x'),L(w))\cong \mathrm{Ext}_{\cO_{\rm{alg}}^{\fb}}^{k'}(M(x'),L(w)).\]
\end{enumerate}
\end{lem}
\begin{proof}
We prove \ref{it: Ext Verma 1}. The isomorphism (\ref{equ: Ext with Verma}) is simply \cite[Thm.\ 6.15(b)]{Hum08}. If (\ref{equ: Ext with Verma}) is non-zero, then by \ref{it: block 3} of Lemma~\ref{lem: central component} $Z(\fg)$ acts on $M(\mu)$ and $L(w)$ by the same character, which together with Harish-Chandra's theorem gives some $x'\in W(G)$ such that $\mu=x'\cdot\mu_0$. Moreover the non-vanishing of (\ref{equ: Ext with Verma}) forces $x'\geq w$ by \cite[Thm.\ 6.11(a)]{Hum08}. \ref{it: Ext Verma 2} follows from (the second statement in) \ref{it: Ext Verma 1} and \cite[Thm.\ 8.11(b),(c)]{Hum08}. We prove \ref{it: Ext Verma 3}. By assumption we have
\[H^{k'}(\fu,L(w))_{x'\cdot\mu_0}\cong \Hom_{U(\ft)}(x'\cdot\mu_0, H^{k'}(\fu,L(w)))=0\mathrm{\ for\ }k'<k.\]
Hence, from \ref{it: block 3} of Lemma~\ref{lem: central component} applied with $I=\emptyset$ (recall from Proposition \ref{lem: n coh O} that $H^{k'}(\fu,L(w))$ is in $\cO_{\ft,\rm{alg}}^{\ft}\subset \cC_{\ft,\rm{alg}}^{\rm{fin}}$) we obtain $\mathrm{Ext}_{U(\ft)}^{k''}(x'\cdot\mu_0,H^{k'}(\fu,L(w)))=0$ for $k''\geq 0$ and $k'<k$. By (\ref{equ: g spectral seq}) (applied with $I=\emptyset$) we deduce for $k'\leq k$
\begin{equation*}
\mathrm{Ext}_{U(\fg)}^{k'}(M(x'),L(w))\cong \Hom_{U(\ft)}(x'\cdot\mu_0, H^{k'}(\fu,L(w)))
\end{equation*}
(which is zero if $k'<k$).
\end{proof}

We collect the following standard results on $\mathrm{Ext}_{\cO_{\rm{alg}}^{\fb}}^1$.

\begin{lem}\label{lem: Ext 1 category O}
Let $x,w\in W(G)$.
\begin{enumerate}[label=(\roman*)]
\item \label{it: Ext O 1} If $\mathrm{Ext}_{\cO_{\rm{alg}}^{\fb}}^1(L(x),L(w))\neq 0$, then either $x<w$ or $w<x$.
\item \label{it: Ext O 2} If $x<w$, then $\mathrm{Ext}_{\cO_{\rm{alg}}^{\fb}}^1(L(x),L(w))\neq 0$ if and only if $x\prec w$, and is $\mu(x,w)$-dimensional in that case.
\item \label{it: Ext O 3} If $x<w$ and $\ell(w)=\ell(x)+1$, then $\Dim_E\mathrm{Ext}_{\cO_{\rm{alg}}^{\fb}}^1(L(x),L(w))=1$.
\end{enumerate}
\end{lem}
\begin{proof}
We start with \ref{it: Ext O 1}. It follows from (\ref{equ: dual category O}) and $L(x)^{\tau}\cong L(x)$, $L(w)^\tau\cong L(w)$ that
\begin{equation}\label{equ: exchange x w}
\mathrm{Ext}_{\cO_{\rm{alg}}^{\fb}}^1(L(x),L(w))\cong \mathrm{Ext}_{\cO_{\rm{alg}}^{\fb}}^1(L(w),L(x)).
\end{equation}
So upon exchanging $x$ and $w$, it suffices to treat the case $x\not\geq w$. Note that $x\not\geq w$ is equivalent to $x\cdot \mu_0- w\cdot \mu_0\notin \Z_{\geq 0}\Phi^+$. Hence, by \cite[Prop.\ 3.1(a)]{Hum08} (and our conventions)
\begin{equation}\label{equ: Ext 1 O vanishing}
\mathrm{Ext}_{\cO_{\rm{alg}}^{\fb}}^1(M(x),L(w))=0.
\end{equation}
As $L(x)$ is the cosocle of $M(x)$, we obviously have
\begin{equation}\label{equ: O head}
\Hom_{\cO_{\rm{alg}}^{\fb}}(L(x),L(w))\cong \Hom_{\cO_{\rm{alg}}^{\fb}}(M(x),L(w)).
\end{equation}
The short exact sequence $0\rightarrow N(x)\rightarrow M(x)\rightarrow L(x)\rightarrow 0$ together with (\ref{equ: Ext 1 O vanishing}) and (\ref{equ: O head}) induce an isomorphism
\begin{equation}\label{equ: O Hom radical}
\Hom_{\cO_{\rm{alg}}^{\fb}}(N(x),L(w))\cong \mathrm{Ext}_{\cO_{\rm{alg}}^{\fb}}^1(L(x),L(w))\neq 0
\end{equation}
which in particular implies $w>x$. We prove \ref{it: Ext O 2}. If $x<w$ we have by (\ref{equ: exchange x w}), \cite[Thm.~8.15(c)]{Hum08} and \cite[p.9]{Brenti}
\begin{equation}\label{equ: Ext 1 O dim}
\Dim_E\mathrm{Ext}_{\cO_{\rm{alg}}^{\fb}}^1(L(x),L(w))=\mu(ww_0,xw_0)=\mu(x,w),
\end{equation}
which implies \ref{it: Ext O 2}. We prove \ref{it: Ext O 3}. If $x<w$ and $\ell(w)=\ell(x)+1$, then $P_{x,w}=1$ and $\mathrm{deg}P_{x,w}=\frac{1}{2}(\ell(w)-\ell(x)-1)=0$, so we have $x\prec w$ and $\mu(x,w)=1$, and thus $\Dim_E\mathrm{Ext}_{\cO_{\rm{alg}}^{\fb}}^1(L(x),L(w))=1$ by (\ref{equ: Ext 1 O dim}).
\end{proof}

\begin{lem}\label{lem: Ext1 O g}
Let $x\in W(G)$ and $M$ in $\cO_{\rm{alg}}^{\fb}$ with all irreducible constituents isomorphic to $L(x')$ for some $x'\not\leq x$. Then we have a canonical isomorphism
\begin{equation}\label{equ: Ext1 O g}
\mathrm{Ext}_{\cO_{\rm{alg}}^{\fb}}^1(L(x),M)\buildrel\sim\over\longrightarrow \mathrm{Ext}_{U(\fg)}^1(L(x),M).
\end{equation}
\end{lem}
\begin{proof}
As (\ref{equ: Ext1 O g}) is clearly injective, it is enough to show that it is also surjective. Given a $U(\fg)$-module $M'$ representing a non-split extension of $M$ by $L(x)$, it is enough to show that $M'$ lies in $\cO_{\rm{alg}}^{\fb}$ (note that it lies in $\cC_{\rm{alg}}^{\rm{fin}}$). For $x'\not\leq x$ and $\mu\in \Z_{\geq 0}\Phi^+$, we have $L(x')_{x\cdot\mu_0-\mu}=0$. Indeed, $x'\not\leq x$ if and only if $x\cdot \mu_0- x'\cdot \mu_0\notin \Z_{\geq 0}\Phi^+$ and $L(x')_{\mu'}= 0$ if and only if $\mu'-x'\cdot\mu_0\notin\Z_{\geq 0}\Phi^+$. By the assumption and an obvious d\'evissage, we deduce $M_{x\cdot\mu_0-\mu}=0$ and thus $M'_{x\cdot\mu_0-\mu}\buildrel\sim\over\rightarrow L(x)_{x\cdot\mu_0-\mu}$ for $\mu\in \Z_{\geq 0}\Phi^+$. In particular, $\dim_E M'_{x\cdot\mu_0}=1$ and $M'_{x\cdot\mu_0-\al}=0$ for $\al\in\Phi^+$, forcing $M'[\fu_\alpha]_{x\cdot\mu_0}\buildrel\sim\over\rightarrow M'_{x\cdot\mu_0}$ for $\al\in\Phi^+$ ($\fu_\alpha$ being the one dimensional Lie subalgebra of $\fu$ corresponding to $\al\in\Delta$) and thus $M'[\fu]_{x\cdot\mu_0}\buildrel\sim\over\rightarrow M'_{x\cdot\mu_0}$. We deduce isomorphisms of one dimensional vector spaces
\[\Hom_{U(\ft)}(x\cdot\mu_0, M'[\fu])\cong \Hom_{U(\fb)}(x\cdot\mu_0, M')\buildrel\sim\over\longrightarrow \Hom_{U(\fb)}(x\cdot\mu_0, L(x)),\]
giving a non-zero map $M(x)\rightarrow M'$ whose image $M''$ is necessarily distinct from $M$. In particular $M'$ is the amalgamate sum of $M$ and $M''$ over $M\cap M''\subseteq M$. But $M''$ is in $\cO_{\rm{alg}}^{\fb}$ as it is a quotient of $M(x)$ and $M$ is in $\cO_{\rm{alg}}^{\fb}$ by assumption. It follows that the quotient $M'$ of $M\oplus M''$ is in also in $\cO_{\rm{alg}}^{\fb}$.
\end{proof}

We will extensively use the following consequences of Lemma \ref{lem: Ext 1 category O} and Lemma \ref{lem: Ext1 O g}:

\begin{lem}\label{rabiotext}
Let $x,w\in W(G)$.
\begin{enumerate}[label=(\roman*)]
\item \label{it: rabiotext 1} We have $\mathrm{Ext}_{\cO_{\rm{alg}}^{\fb}}^1(L(x),L(w))\neq 0$ if and only if $x\prec w$ or $w\prec x$, in which case it has dimension $\mu(x,w)$ or $\mu(w,x)$ respectively and $\vert\ell(x) - \ell(w)\vert$ is odd.
\item \label{it: rabiotext 2} Assume $x\neq w$, then $\mathrm{Ext}_{\cO^{\fb}_{\rm{alg}}}^1(L(x),L(w))\cong \mathrm{Ext}_{U(\fg)}^1(L(x),L(w))$. In particular for $x\neq w$ we have $\mathrm{Ext}_{U(\fg)}^1(L(x),L(w))\neq 0$ if and only if $x\prec w$ or $w\prec x$, in which case it has dimension $\mu(x,w)$ or $\mu(w,x)$ respectively.
\end{enumerate}
\end{lem}
\begin{proof}
\ref{it: rabiotext 1} follows from parts \ref{it: Ext O 1} and \ref{it: Ext O 2} of Lemma~\ref{lem: Ext 1 category O} together with (\ref{equ: exchange x w}). \ref{it: rabiotext 2} follows from Lemma~\ref{lem: Ext1 O g} together with (\ref{equ: dual Ext1}) and (\ref{equ: dual category O}), and from \ref{it: rabiotext 1}.
\end{proof}

\begin{lem}\label{lem: dominance control}
If $x\prec w$ and $\ell(w)>\ell(x)+1$, then we have $D_L(w)\subseteq D_L(x)$ and $D_R(w)\subseteq D_R(x)$.
\end{lem}
\begin{proof}
The inclusion $D_R(w)\subseteq D_R(x)$ is \cite[Prop.\ 5.1.9]{BB05}, and the inclusion $D_L(w)\subseteq D_L(x)$ is its symmetric version which follows from $P_{x,w}=P_{x^{-1},w^{-1}}$ (\cite[p.9]{Brenti}).
\end{proof}

Recall that $W^{I,\emptyset}$ is the set of minimal length representatives of $W(L_{I})\backslash W(G)$.

\begin{lem}\label{lem: Levi reduction}
Let $I\subseteq \Delta$, $y\in W^{I,\emptyset}$, $x,w\in W(L_I)$ with $x<w$ and $P^I_{x,w}$ the Kazhdan-Lusztig polynomial for the Coxeter group $W(L_I)$. Then we have $P^I_{x,w}=P_{x,w}=P_{xy,wy}$.
\end{lem}
\begin{proof}
Recall that the multiplication on the right by $y$ gives a bijection preserving the Bruhat order between the interval $[x,w]$ in $W(G)$ and the interval $[xy,wy]$ in $W(G)$, and that $\ell(x'y)=\ell(x')+\ell(y)$ for $x'\in W(L_I)$ (one reference for such facts is \cite[\S 2]{BB05}). Using this, the equality $P_{x,w}=P_{xy,wy}$ easily follows by induction on $\ell(y)$ from \cite[Thm.\ 2(iii)]{Brenti} and \cite[Thm.\ 4(iv)]{Brenti} (arguing as in the last paragraph of \cite[p.8]{Brenti}). By the same combinatorics, the equality $P^I_{x,w}=P_{x,w}$ comes from the fact that the interval $[x,w]$ in $W(G)$ is the same as the interval $[x,w]$ in $W(L_I)$.
\end{proof}


\begin{lem}\label{lem: H0 and H1}
Let $x,w\in W(G)$ and $I\subseteq \Delta$.
\begin{enumerate}[label=(\roman*)]
\item \label{it: coh 1} If $w\in W(L_I)x$, then we have canonical isomorphisms for $k\geq 0$
\begin{equation*}
\mathrm{Ext}_{\cO^{\fb_{I}}_{\fl_{I},\rm{alg}}}^k(L^I(w), H^0(\fn_{I},L(x)))\cong \mathrm{Ext}_{\cO_{\rm{alg}}^{\fb}}^k(M^I(w),L(x)).
\end{equation*}
\item \label{it: coh 3} The $E$-vector space $\Hom_{U(\fl_{I})}(L^I(w),H^1(\fn_{I},L(x)))$ is non-zero if and only if $x\prec w$ and $w\notin W(L_I)x$, and has dimension $\mu(x,w)$ if non-zero.
\item \label{it: coh 4} For $k\geq 0$ the canonical map
\begin{equation}\label{equ: n coh sequence}
H^k(\fu_{I}, H^1(\fn_{I}, L(x)))\rightarrow H^{k+2}(\fu_{I},H^0(\fn_{I},L(x)))
\end{equation}
induced from the spectral sequence (\ref{equ: double n coh}) (with $I'=\emptyset$) is zero.
\item \label{it: coh 5} The inclusion $\mathrm{soc}_{U(\fl_{I})}(H^1(\fn_{I},L(x)))\subseteq H^1(\fn_{I},L(x))$ induces a $U(\ft)$-equivariant isomorphism
\begin{equation}\label{equ: coh of socle}
H^0\big(\fu_{I}, \mathrm{soc}_{U(\fl_{I})}(H^1(\fn_{I},L(x)))\big)\buildrel\sim\over\longrightarrow H^0(\fu_{I}, H^1(\fn_{I},L(x)))
\end{equation}
with both $U(\ft)$-modules in (\ref{equ: coh of socle}) isomorphic to $\bigoplus_w(w\cdot\mu_0)^{\oplus\mu(x,w)}$ where $w$ runs through those $w\in W(G)$ such that $w\notin W(L_I)x$ and $x\prec w$.
\end{enumerate}
\end{lem}
\begin{proof}
\ref{it: coh 1} follows from \ref{it: dominance 2} of Lemma~\ref{lem: n coh dominance} and \ref{it: degenerate 4} of Lemma \ref{lem: Ext 1 degenerate} (together with Harish-Chandra's theorem).

We prove \ref{it: coh 3}. Let $\mu\in \Lambda$, by Lemma \ref{lem: n coh wt} if $H^1(\fn_{I},L(x))_\mu\!\neq 0$ then $\mu-x\cdot\mu_0-\al\in \Z_{\geq 0}\Phi^+$ for some $\al\in \Phi^+\setminus \Phi_I^+$. In particular
\begin{equation}\label{w>x}
\Hom_{U(\fl_{I})}(L^I(w),H^1(\fn_{I},L(x)))\neq 0\implies w\cdot \mu_0-x\cdot\mu_0\in \al + \Z_{\geq 0}\Phi^+ \implies w>x.
\end{equation}
On the other hand, each Jordan-H\"older factor $L(w')$ of $N(w)$ satisfies $w'>w$. So for $w>x$ the surjections $M(w)\twoheadrightarrow M^I(w)\twoheadrightarrow L(w)$ induce embeddings
\begin{equation}\label{equ: Ext mu}
\mathrm{Ext}_{\cO_{\rm{alg}}^{\fb}}^1(L(w),L(x))\hookrightarrow \mathrm{Ext}_{\cO_{\rm{alg}}^{\fb}}^1(M^I(w),L(x))\hookrightarrow \mathrm{Ext}_{\cO_{\rm{alg}}^{\fb}}^1(M(w),L(x)).
\end{equation}
But since the first and third vector spaces both have dimension $\mu(x,w)$ by \cite[p.9]{Brenti}, \cite[Thm.\ 8.15(c)]{Hum08} and \cite[Thm.\ 8.11(b)]{Hum08}, we deduce that the embeddings in (\ref{equ: Ext mu}) are all isomorphisms. In particular, when $x<w$ each $E$-vector space in (\ref{equ: Ext mu}) has dimension $\mu(x,w)$ and is non-zero if and only if $x\prec w$.

Now, from \ref{it: degenerate 2}, \ref{it: degenerate 3} and \ref{it: degenerate 4} of Lemma~\ref{lem: Ext 1 degenerate} we have $\Hom_{U(\fl_{I})}(L^I(w),H^1(\fn_{I},L(x)))\neq 0$ if and only if $\mathrm{Ext}_{\cO_{\rm{alg}}^{\fb}}^1(M^I(w),L(x))\neq 0$ and $H^0(\fn_I, L(x))_\xi=0$ where $\xi$ is the infinitesimal character of $L^I(w)$. By Harish-Chandra's theorem and \ref{it: dominance 2} of Lemma~\ref{lem: n coh dominance}, $H^0(\fn_I, L(x))_\xi=0$ if and only if $w\notin W(L_I)x$. By the previous paragraph $x\prec w$ if and only if $x<w$ and $\mathrm{Ext}_{\cO_{\rm{alg}}^{\fb}}^1(M^I(w),L(x))\neq 0$. Using (\ref{w>x}), we get \ref{it: coh 3}.

We prove \ref{it: coh 4}. Since $H^0(\fn_{I},L(x))\cong L^I(x)$ by \ref{it: dominance 2} of Lemma~\ref{lem: n coh dominance}, for each $\mu\in \Lambda$ satisfying $H^{k+2}(\fu_{I},H^0(\fn_{I},L(x)))_\mu\neq 0$, we deduce in particular from Lemma~\ref{lem: n coh wt} (applied with $\fl_I$ instead of $\fg$) that
\begin{equation}\label{equ: n coh wt 1}
\mu-x\cdot\mu_0\in\Z_{\geq 0}\Phi_I^+.
\end{equation}
Let $\mu'\in \Lambda$ such that $H^k(\fu_{I}, H^1(\fn_{I},L(x)))_{\mu'}\neq 0$.
By a d\'evissage on $H^1(\fn_{I},L(x))$
(using Proposition~\ref{lem: n coh O})
and by applying Lemma~\ref{lem: n coh wt} twice, first to $H^1(\fn_{I},L(x))$ and then to $H^k(\fu_{I},L^I(\mu''))$ for each Jordan-H\"older factor $L^I(\mu'')$ of $H^1(\fn_{I},L(x))$, we deduce the existence of $\al\in\Phi^+\setminus \Phi_I^+$ such that
\begin{equation}\label{equ: n coh wt 2}
\mu'-x\cdot\mu_0-\al\in \Z_{\geq 0}\Phi^+.
\end{equation}
If the map (\ref{equ: n coh sequence}) is non-zero, then there must exists $\mu=\mu'$ satisfying both (\ref{equ: n coh wt 1}) and (\ref{equ: n coh wt 2}), which is impossible as
\[\Z_{\geq 0}\Phi_I^+\cap (\al+\Z_{\geq 0}\Phi^+)=\emptyset.\]

We prove \ref{it: coh 5}. Note first that \ref{it: coh 4} together with (\ref{equ: double n coh}) imply a short exact sequence
\begin{equation}\label{equ: n coh two terms}
0\rightarrow H^1(\fu_{I},H^0(\fn_{I},L(x)))\rightarrow H^1(\fu, L(x))\rightarrow H^0(\fu_{I}, H^1(\fn_{I},L(x)))\rightarrow 0.
\end{equation}
Secondly, from \ref{it: Ext Verma 1} and \ref{it: Ext Verma 2} of Lemma~\ref{lem: Ext with Verma} (and $\mu(ww_0,xw_0)=\mu(x,w)$) we have
\begin{equation}\label{equ: full H1}
H^1(\fu,L(x))\cong \bigoplus_{x\prec w}(w\cdot\mu_0)^{\oplus \mu(x,w)}.
\end{equation}
Thirdly, write $x=x_Ix^I$ where $x_I\in W(L_I)$ and $x^I\in W^{I, \emptyset}$ and note that $x^I\cdot \mu_0\in \Lambda_I^{\dom}$, so that, when dealing with the reductive group $L_I$, we can replace $\mu_0$ by $x^I\cdot \mu_0$. Then a proof analogous to the proof of \ref{it: Ext Verma 1} and \ref{it: Ext Verma 2} of Lemma~\ref{lem: Ext with Verma} replacing $G$ by $L_I$ and using Lemma \ref{lem: Levi reduction} gives
\begin{equation}\label{equ: full HI}
H^1(\fu_{I},L^I(x))\cong \bigoplus_{x\prec w, w\in W(L_I)x}(w\cdot\mu_0)^{\oplus \mu(x,w)}.
\end{equation}
As $H^0(\fn_{I},L(x))\cong L^I(x)$ by \ref{it: dominance 2} of Lemma~\ref{lem: n coh dominance}, we deduce from (\ref{equ: n coh two terms}), (\ref{equ: full H1}) and (\ref{equ: full HI})
\begin{equation*}
H^0(\fu_{I}, H^1(\fn_{I},L(x)))\cong \bigoplus_{x\prec w, w\notin W(L_I)x}(w\cdot\mu_0)^{\oplus \mu(x,w)}.
\end{equation*}
But by \ref{it: coh 3} we also have
\[\mathrm{soc}_{U(\fl_{I})}(H^1(\fn_{I},L(x)))\cong \bigoplus_{x\prec w, w\notin W(L_I)x}L^I(w)^{\oplus \mu(x,w)}\]
which clearly implies (\ref{equ: coh of socle}).
\end{proof}

\begin{rem}\label{rem: Levi Ext}
Let $x,w\in W(G)$ with $x\neq w$ and $I\defeq \Delta\setminus D_L(w)$. The short exact sequence $0\rightarrow N^I(w)\rightarrow M^I(w)\rightarrow L(w)\rightarrow 0$ together with $x\neq w$ give an exact sequence
\begin{equation}\label{equ: Levi Ext devissage}
0\rightarrow \Hom_{U(\fg)}(N^I(w),L(x))\rightarrow \mathrm{Ext}_{U(\fg)}^1(L(w),L(x))\rightarrow \mathrm{Ext}_{U(\fg)}^1(M^I(w),L(x)).
\end{equation}
By \cite[Thm.~9.4(c)]{Hum08} and Lemma \ref{lem: dominance and left set} we have $\Hom_{U(\fg)}(N^I(w),L(x))\neq 0$ if and only if $\Hom_{U(\fg)}(N(w),L(x))\neq 0$ and $D_L(x)\subseteq D_L(w)$, and by (\ref{equ: O Hom radical}) and \ref{it: Ext O 2} of Lemma \ref{lem: Ext 1 category O} that $\Hom_{U(\fg)}(N(w),L(x))\neq 0$ if and only if $w\prec x$. Since $w\prec x$ and $D_L(x)\subseteq D_L(w)$ if and only if $w\prec x$ and $x\notin W(L_I)w$ (using Lemma~\ref{lem: dominance control} to deal with the implication $\Leftarrow$ when $\ell(x)>\ell(w)+1$), we finally deduce (with \ref{it: Ext O 2} of Lemma \ref{lem: Ext 1 category O}) that $\Hom_{U(\fg)}(N^I(w),L(x))\neq 0$ if and only if $w\prec x$ and $x\notin W(L_I)w$, in which case it has dimension $\mu(w,x)$.

Moreover it follows from \ref{it: degenerate 2}, \ref{it: degenerate 3} of Lemma \ref{lem: Ext 1 degenerate} (with Harish-Chandra's theorem) and \ref{it: dominance 2} of Lemma~\ref{lem: n coh dominance} that one has isomorphisms
\begin{equation}\label{equ: Levi Ext H0}
\mathrm{Ext}_{U(\fg)}^1(M^I(w),L(x))\cong \mathrm{Ext}_{U(\fl_I)}^1(L^I(w),H^0(\fn_I,L(x)))\cong \mathrm{Ext}_{U(\fl_I)}^1(L^I(w),L^I(x))
\end{equation}
if $w\in W(L_I)x$, and
\begin{equation}\label{equ: Levi Ext H1}
\mathrm{Ext}_{U(\fg)}^1(M^I(w),L(x))\cong \Hom_{U(\fl_I)}(L^I(w),H^1(\fn_I,L(x)))
\end{equation}
if $w\notin W(L_I)x$.
When $w\in W(L_I)x$, it then follows from Lemma~\ref{lem: Levi reduction} (upon writing $x=x'y$, $w=w'y$ for $w'<x'$ in $W(L_I)$ and $y\in W^{I,\emptyset}$) and \ref{it: rabiotext 2} of Lemma \ref{rabiotext} (applied with $\fl_I$ instead of $\fg$) that (\ref{equ: Levi Ext H0}) has dimension $\mu(x,w)$ (resp.~$\mu(w,x)$) when $x<w$ (resp.~when $w<x$), which is also $\Dim_E\mathrm{Ext}_{U(\fg)}^1(L(w),L(x))$ by \ref{it: rabiotext 2} of Lemma \ref{rabiotext}.
When $w\notin W(L_I)x$, it then follows from \ref{it: coh 3} of Lemma~\ref{lem: H0 and H1} that (\ref{equ: Levi Ext H1}) is non-zero if and only if $x\prec w$ (and $w\notin W(L_I)x$) and has dimension $\mu(x,w)=\Dim_E\mathrm{Ext}_{U(\fg)}^1(L(w),L(x))$ in that case.

By (\ref{equ: Levi Ext devissage}) and the above discussion on $\Hom_{U(\fg)}(N^I(w),L(x))$ and $\mathrm{Ext}_{U(\fg)}^1(M^I(w),L(x))$, we see that $\mathrm{Ext}_{U(\fg)}^1(L(w),L(x))\neq 0$ if and only if exactly one of the following holds:
\begin{itemize}
\item $w\prec x$, $x\notin W(L_I)(w)$, $\mathrm{Ext}_{U(\fg)}^1(M^I(w),L(x))=0$ and (\ref{equ: Levi Ext devissage}) induces an isomorphism
\[\Hom_{U(\fg)}(N^I(w),L(x))\buildrel\sim\over\longrightarrow \mathrm{Ext}_{U(\fg)}^1(L(w),L(x));\]
\item $\Hom_{U(\fg)}(N^I(w),L(x))=0$ and (\ref{equ: Levi Ext devissage}) induces an isomorphism
\[\mathrm{Ext}_{U(\fg)}^1(L(w),L(x))\buildrel\sim\over\longrightarrow \mathrm{Ext}_{U(\fg)}^1(M^I(w),L(x))\]
(as it is then an embedding between two vector spaces of the same dimension).
\end{itemize}
\end{rem}

\begin{lem}\label{lem: Ext with dominant}
Let $w\in W(G)$. We have canonical isomorphisms for $\ell\leq \ell(w)$ induced by $M(w)\twoheadrightarrow L(w)$
\begin{equation}\label{equ: Ext with dominant}
\mathrm{Ext}_{\cO_{\rm{alg}}^{\fb}}^\ell(L(w),L(1))\buildrel\sim\over\rightarrow \mathrm{Ext}_{\cO_{\rm{alg}}^{\fb}}^\ell(M(w),L(1))\cong \Hom_{U(\ft)}(w\cdot\mu_0,H^\ell(\fu,L(1)))
\end{equation}
where all $E$-vector spaces in (\ref{equ: Ext with dominant}) are $0$ if $\ell<\ell(w)$ and $1$-dimensional if $\ell=\ell(w)$.
\end{lem}
\begin{proof}
Recall that by definition $L(1)=L(\mu_0)$.
By Bott's formula (see \cite[\S 6.6]{Hum08}) there is a $U(\ft)$-equivariant isomorphism for $\ell\geq 0$
\begin{equation}\label{equ: Bott formula}
H^\ell(\fu, L(1))\cong \bigoplus_{\ell(w)=\ell}w\cdot \mu_0.
\end{equation}
We deduce from (\ref{equ: category O HS}) applied with $I=\emptyset$ (using $\Ext_{\cO^{\ft}_{\ft,\rm{alg}}}^{\ell_1}=0$ if $\ell_1>0$) that for $\ell\geq 0$
\begin{equation*}
\mathrm{Ext}_{\cO_{\rm{alg}}^{\fb}}^\ell(M(w),L(1))\cong \Hom_{U(\ft)}\big(w\cdot\mu_0, H^\ell(\fu, L(1))\big),
\end{equation*}
which is $1$-dimensional if $\ell=\ell(w)$ and $0$ otherwise by (\ref{equ: Bott formula}). This gives the second isomorphism in (\ref{equ: Ext with dominant}). We prove the first isomorphism in (\ref{equ: Ext with dominant}) when $\ell\leq \ell(w)$ by decreasing induction on $\ell(w)$. If $w=w_0$, then $L(w_0)=M(w_0)$ and there is nothing more to prove. If $w<w_0$, we have by induction $\mathrm{Ext}_{\cO_{\rm{alg}}^{\fb}}^\ell(L(w'),L(1))=0$ for $\ell\leq \ell(w)$ and $w<w'$, which implies $\mathrm{Ext}_{\cO_{\rm{alg}}^{\fb}}^\ell(N(w),L(1))=0$ for $\ell\leq \ell(w)$ by an obvious d\'evissage on $N(w)$. The short exact sequence $0\rightarrow N(w)\rightarrow M(w)\rightarrow L(w)\rightarrow 0$ then induces an isomorphism for $\ell\leq \ell(w)$
\[\mathrm{Ext}_{\cO_{\rm{alg}}^{\fb}}^\ell(L(w),L(1))\buildrel\sim\over\longrightarrow \mathrm{Ext}_{\cO_{\rm{alg}}^{\fb}}^\ell(M(w),L(1)).\qedhere\]
\end{proof}

\begin{lem}\label{lem: Ext1 with socle}
Let $M$ in $\cO^{\fb}_{\rm{alg}}$. Assume that the inclusion $\mathrm{soc}_{U(\fg)}(M)\subseteq M$ induces an isomorphism
\begin{equation}\label{equ: socle u inv}
H^0(\fu, \mathrm{soc}_{U(\fg)}(M))\buildrel\sim\over\longrightarrow H^0(\fu, M)
\end{equation}
and that
\begin{equation}\label{equ: exclude dominant}
\Hom_{U(\ft)}(\mu_0, H^0(\fu, M))=0,
\end{equation}
then it also induces an isomorphism
\begin{equation}\label{equ: Ext with socle}
\mathrm{Ext}_{U(\fg)}^1(L(1),\mathrm{soc}_{U(\fg)}(M))\buildrel\sim\over\longrightarrow \mathrm{Ext}_{U(\fg)}^1(L(1),M).
\end{equation}
\end{lem}
\begin{proof}
We have the decomposition $M=\bigoplus_\xi M_\xi$ from \ref{it: block 2} of Lemma~\ref{lem: central component}. Since both vector spaces in (\ref{equ: Ext with socle}) are $0$ when $M$ is replaced by $M_\xi$ with $\xi\ne \xi_{\mu_0}$ by \ref{it: block 3} of Lemma~\ref{lem: central component}, we can assume $M=M_{\xi_{\mu_0}}$, i.e.\ that all Jordan-H\"older factors of $M$ are of the form $L(x)$ for some $x\in W(G)$.

For $\al\in\Phi^+$, we have $L(x)_{\mu_0-\al}=0$ for each $x\in W(G)$ and thus $M_{\mu_0-\al}=0$. Consequently, any $v\in M_{\mu_0}$ must be killed by $\fu$, i.e.\ $M_{\mu_0}\subseteq H^0(\fu, M)$. Then (\ref{equ: exclude dominant}) forces $M_{\mu_0}=0$, i.e.\ $L(1)$ is not a constituent of $M$. Likewise, for $\al\in\Phi^+$ and $x,x'\in W(G)$ with $x'\neq 1$ and $\ell(x)=1$, we have $L(x')_{x\cdot\mu_0-\al}=0$. Consequently $M_{x\cdot\mu_0-\al}=0$ for $\al\in\Phi^+$ and $x\in W(G)$ with $\ell(x)=1$ (using that $L(1)$ does not appear in $M$) and thus $M_{x\cdot\mu_0}\subseteq H^0(\fu,M)$. From (\ref{equ: socle u inv}) we then obtain $M_{x\cdot\mu_0}\subseteq H^0(\fu,\mathrm{soc}_{U(\fg)}(M))\subseteq \mathrm{soc}_{U(\fg)}(M)$. Since any constituent $L(x)$ of $M$ contributes to $M_{x\cdot\mu_0}$, we deduce that all constituents $L(x)$ of $M$ with $\ell(x)=1$ can only appear in $\mathrm{soc}_{U(\fg)}(M)$. Since, by \ref{it: rabiotext 2} of Lemma \ref{rabiotext}, for $x\ne 1$ we have $\mathrm{Ext}_{U(\fg)}^1(L(1),L(x))\neq 0$ if and only if $\ell(x)=1$, an easy d\'evissage implies (\ref{equ: Ext with socle}).
\end{proof}

\subsection{\texorpdfstring{$W(G)$}{W(G)}-conjugates of objects of \texorpdfstring{$\cO_{\rm{alg}}^{\fb}$}{Obalg}}

We study unipotent cohomology groups and $\mathrm{Ext}^1$, $\mathrm{Ext}^2$ groups of $U(\fg)$-modules which are conjugates of $U(\fg)$-modules in the category $\cO_{\rm{alg}}^{\fb}$ by elements of $W(G)$.\bigskip

Recall we defined full subcategories $\cC_{\rm{alg}}$, $\cC_{\rm{alg}}^{\rm{fin}}$, $\tld{\cO}_{\rm{alg}}^{\fb}$ and $\cO_{\rm{alg}}^{\fb}$ of $\mathrm{Mod}_{U(\fg)}$ in \S\ref{subsec: category}. For $M$ in $\mathrm{Mod}_{U(\fg)}$ and $g\in G$, we define $M^g$ in $\mathrm{Mod}_{U(\fg)}$ as the same underlying $U(\fg)$-module as $M$ but where $x\in U(\fg)$ acts by ${\rm ad}(g)(x)=gxg^{-1}$. For $g_1, g_2\in G$ we have $(M^{g_1})^{g_2}\cong M^{g_1g_2}$. For $M_1,M_2$ in $\mathrm{Mod}_{U(\fg)}$ and $g\in G$, we (clearly) have isomorphisms for $k\geq 0$
\begin{equation}\label{equ: twist isom}
\mathrm{Ext}_{U(\fg)}^k(M_1,M_2)\cong \mathrm{Ext}_{U(\fg)}^k(M_1^g,M_2^g).
\end{equation}
Recall that an algebraic action of $\ft$, resp.~of $\fb$, on a finite dimensional $E$-vector space lifts uniquely to an (algebraic) action of $T$, resp.~of $B$, see for instance \cite[Lemma~3.2]{OS15}. It follows that, for $M$ in $\cC_{\rm{alg}}^{\rm{fin}}$ and $t\in T$, the resulting action of $t$ on $M$ induces an isomorphism $M \buildrel\sim\over\rightarrow M^t$ in $\cC_{\rm{alg}}^{\rm{fin}}$. Likewise, for $M$ in $\tld{\cO}_{\rm{alg}}^{\fb}$ and $b\in B$, we have $M \buildrel\sim\over\rightarrow M^b$ in $\tld{\cO}_{\rm{alg}}^{\fb}$. In particular, for $M$ in $\cC_{\rm{alg}}^{\rm{fin}}$, the isomorphism class of $M^{g}$ is independent of the choice of $g\in N_G(T)$ lifting $w\in N_G(T)/T$ and we denote it by $M^w$. For $I\subseteq \Delta$ and $M\in \cO_{\rm{alg}}^{\fp_I}$, by \cite[Lemma~3.2]{OS15} the action of $U(\fp_I)$ on $M$ lifts uniquely to an action of $P_I$ on $M$, and thus as above $M \buildrel\sim\over\rightarrow M^g$ as $U(\fg)$-modules when $g\in P_I$. In particular, $M^w\cong M$ for each $w\in W(L_I)$.

\begin{lem}\label{lem: H0 Weyl conjugate}
Let $I,I'\subseteq \Delta$, $M$ in $\cO_{\rm{alg}}^{\fp_{I'}}$ and $w\in W(G)$. We write $w=w_1w_2w_3$ for (unique) $w_1\in W(L_{I'})$, $w_3\in W(L_I)$ and $w_2\in W^{I',I}$.
\begin{enumerate}[label=(\roman*)]
\item \label{it: H0 conjugate 1} For $k \geq 0$ the $U(\fl_I)$-module $H^k(\fn_I, M^{w_2})$ lies in $\cO^{\fb_{I}}_{\fl_{I},\rm{alg}}$ and we have $H^k(\fn_{I},M^w)\cong (H^k(\fn_I, M^{w_2}))^{w_3}$ in $\mathrm{Mod}_{U(\fl_I)}$.
\item \label{it: H0 conjugate 2} For \ $\mu\in \Lambda$ \ such \ that \ $I'$ \ is \ maximal \ for \ the \ condition \ $L(\mu)\in \cO_{\rm{alg}}^{\fp_{I'}}$, \ we \ have $H^0(\fn_{I},L(\mu)^w)\neq 0$ if and only if $w_2=1$.
\item \label{it: H0 conjugate 3} For $\mu\in \Lambda$ such that $I'$ is maximal for the condition $L(\mu)\in \cO_{\rm{alg}}^{\fp_{I'}}$, we have $L(\mu)^w\in \cO_{\rm{alg}}^{\fb}$ if and only if $L(\mu)^w\cong L(\mu)$ if and only if $w\in W(L_{I'})$.
\end{enumerate}
\end{lem}
\begin{proof}
We prove the second statement in \ref{it: H0 conjugate 1}. Since any lift of $w_1\in W(L_{I'})$ in $N_G(T)$ lies in $L_{I'}\subseteq P_{I'}$, we have $M^{w_1}\cong M$ as $U(\fg)$-modules by the sentence just before Lemma \ref{lem: H0 Weyl conjugate}. Since $L_I$ normalizes $N_I$, we have $w_3^{-1}\fn_Iw_3\cong \fn_I$ and thus (using (\ref{equ: CE complex}) for the last isomorphism)
\begin{equation}\label{twisttwist}
H^k(\fn_I, M_1^{w_3})\cong H^k(w_3^{-1}\fn_Iw_3,M_1^{w_3})\cong H^k(\fn_I, M_1)^{w_3}
\end{equation}
for any $U(\fg)$-module $M_1$ and $k\geq 0$. We apply this to $M_1=M^{w_1w_2}\cong M^{w_2}$, which gives the second statement. We prove the first statement in \ref{it: H0 conjugate 1}. The minimal length assumption on $w_2$ implies in particular $\fb_{I}=\fb\cap \fl_{I}\subseteq w_2^{-1}\fb w_2$, and thus $M^{w_2}$ is locally $\fb_{I}$-finite, which together with the $U(\fl_{I})$-equivariance of the Chevalley-Eilenberg complex (\ref{equ: CE complex}) (and the fact $\Dim_E\fn_I<+\infty$) implies that the $U(\fl_I)$-module $H^k(\fn_{I},M^{w_2})$ is locally $\fb_{I}$-finite for $k\geq 0$. A similar argument shows that $H^k(\fn_{I},M^{w_2})$ is also $\ft$-semi-simple. Then the argument to show that $H^k(\fn_{I},M^{w_2})$ has finite length (hence is in $\cO^{\fb_{I}}_{\fl_{I},\rm{alg}}$) is parallel to the one for the proof of Proposition~\ref{lem: n coh O}, using as (crucial) inputs Lemma \ref{lem: HC n coh} and $\Dim_E H^k(\fn_{I},M^{w_2})_{\mu}<\infty$ for $\mu\in\Lambda$. Note that $M^{w_2}$ itself is not in general in $\cO^{\fb}_{\rm{alg}}$ as it might not be locally $\fb$-finite

We prove \ref{it: H0 conjugate 2}. By the second statement in \ref{it: H0 conjugate 1}, it suffices to treat the case when $w=w_2\in W^{I',I}$. If $w_2=1$, we have $H^0(\fn_{I},L(\mu))\cong L^I(\mu)\neq 0$ from \ref{it: dominance 2} of Lemma \ref{lem: n coh dominance}. We assume from now on $w_2\neq 1$ and $H^0(\fn_{I},L(\mu)^{w_2})\neq 0$ and seek a contradiction.
Let $M_I$ be any non-zero simple $U(\fl_{I})$-module that embeds into $H^0(\fn_{I},L(\mu)^{w_2})$ ($M_I$ exists and belongs to $\cO_{\fl_I,\rm{alg}}^{\fb_I}$ by the first statement in \ref{it: H0 conjugate 1}). Then the injection $M_I\hookrightarrow H^0(\fn_{I},L(\mu)^{w_2})$ induces a non-zero map $U(\fg)\otimes_{U(\fp_I)}M_I\rightarrow L(\mu)^{w_2}$ by (\ref{equ: g spectral seq 0}) which has to be a surjection as $L(\mu)^{w_2}$ is irreducible. As $U(\fg)\otimes_{U(\fp_I)}M_I$ is in $\cO_{\rm{alg}}^{\fb}$ (see below (\ref{belongtoOb})), $L(\mu)^{w_2}$ is also an object of $\cO_{\rm{alg}}^{\fb}$, and in particular $L(\mu)$ is locally $w_2 \fb w_2^{-1}$-finite. As $1\neq w_2$, there exists $\al\in \Delta$ (recall it is a positive simple root for $\fb^+$) such that $w_2^{-1}(\al)\in\fb$, or equivalently $\fu_\al\subseteq w_2 \fb w_2^{-1}$ (where $\fu_\al\subseteq \fu\subset \fb^+$ is the one dimensional root subspace corresponding to $\alpha$), which implies that $L(\mu)$ is locally $\fu_\al$-finite. Note that $\alpha\notin I'$ because $w_2^{-1}(I')\subseteq \Phi^+$. Let $v\in M(\mu)_\mu$ be a highest weight vector of $M(\mu)$ and $0\neq x_\al\in \fu_\al$. As $L(\mu)$ is locally $\fu_\al$-finite, there exists $N\geq 1$ such that $x_\al^N\cdot v\in N(\mu)_{\mu+N\al}\subseteq M(\mu)_{\mu+N\al}$, which implies $N(\mu)\ne 0$ as $M(\mu)\cong U(\fu^+)$. As each Jordan-H\"older factor of $N(\mu)$ has highest weight $w'\cdot\mu$ for some $w'\in W(G)$ such that $w'\cdot\mu-\mu\in \Z_{> 0}\Phi^+$ (\cite[Thm.~5.1]{Hum08}, in particular $w'\ne 1$), we deduce $\mu+N\al-w'\cdot\mu\in \Z_{\geq 0}\Phi^+$ for such a $w'$. Since $w'\ne 1$, there also exists $\beta\in \Delta$ such that $w'\cdot\mu-\mu=(w'\cdot\mu-s_\beta\cdot \mu)+(s_\beta\cdot \mu - \mu)\in \beta + \Z_{\geq 0}\Phi^+$. We can choose $\beta\ne \alpha$, as otherwise this would mean $w'\cdot\mu-\mu\in \Z_{> 0}\alpha$, hence $w'=s_\alpha$, but this is impossible since $\alpha\notin I'$ and $I'$ is maximal such that $L(\mu)\in \cO_{\rm{alg}}^{\fp_{I'}}$, equivalently $I'$ is maximal such that $\mu$ is dominant with respect to $\fb_{I'}$, which implies $s_\alpha\cdot\mu-\mu\in \Z_{\leq 0}\Phi^+$. It follows that $(\mu+N\al-w'\cdot\mu) + (w'\cdot\mu-\mu) = N\al \in \beta + \Z_{\geq 0}\Phi^+$ which is impossible as $\beta\ne \al$.

Finally we prove \ref{it: H0 conjugate 3}. As in \ref{it: H0 conjugate 1} above, $w\in W(L_{I'})$ implies $L(\mu)^w\cong L(\mu)$ which implies $L(\mu)^w\in \cO_{\rm{alg}}^{\fb}$. Conversely, assume $L(\mu)^w\in \cO_{\rm{alg}}^{\fb}$, which implies $H^0(\fu, L(\mu)^w)\neq 0$. By \ref{it: H0 conjugate 2} applied with $I=\emptyset$ we get $w\in W(L_{I'})W(L_{I})=W(L_{I'})$.
\end{proof}

The following consequence of Lemma \ref{lem: H0 Weyl conjugate} will be used later.

\begin{lem}
For $k\geq 0$, \ $M$,\ $M'$\ in \ $\cO^{\fb}_{\rm{alg}}$ \ and \ $w\in W(G)$, \ the \ $E$-vector \ space \ $\mathrm{Ext}_{U(\fg)}^k(M',M^{w})$ is finite dimensional.
\end{lem}
\begin{proof}
It follows from the first statement in \ref{it: H0 conjugate 1} of Lemma~\ref{lem: H0 Weyl conjugate} applied with $I=\emptyset$ that $H^{\ell}(\fu,L(\mu)^{w})$ is a finite dimensional semi-simple $U(\ft)$-module for $\ell \geq 0$, $\mu\in \Lambda$ and $w\in W(G)$, and thus $\mathrm{Ext}_{U(\ft)}^k(\mu', H^{\ell}(\fu,L(\mu)^{w}))$ is finite dimensional for $k,\ell\geq 0$, $\mu,\mu'\in\Lambda$. By (\ref{equ: g spectral seq}) (applied with $I=\emptyset$) we deduce
\begin{equation}\label{fdM}
\Dim_E\mathrm{Ext}_{U(\fg)}^k(M(\mu'),L(\mu)^{w})<+\infty
\end{equation}
for $k\geq 0$, $\mu,\mu'\in\Lambda$ and $w\in W(G)$. Now, let $\mu'\in \Lambda$, if there is no $\mu''\ne \mu'$ such that $\mu'' \uparrow \mu'$ (where $\uparrow$ is the strong linkage relation from \cite[\S 5.1]{Hum08}) then $M(\mu')\buildrel\sim\over\rightarrow L(\mu')$ by \cite[Thm.~5.1(b)]{Hum08} and hence $\mathrm{Ext}_{U(\fg)}^k(L(\mu'),L(\mu)^{w})$ is finite dimensional for $k\geq 0$ by (\ref{fdM}). Assume by induction that $\mathrm{Ext}_{U(\fg)}^k(L(\mu''),L(\mu)^{w})$ is finite dimensional for $k\geq 0$ and any $\mu'' \uparrow \mu'$, $\mu''\ne \mu'$. Then $\mathrm{Ext}_{U(\fg)}^k(N(\mu'),L(\mu)^{w})$ and $\mathrm{Ext}_{U(\fg)}^k(M(\mu'),L(\mu)^{w})$ are both finite dimensional, the first by d\'evissage and induction using \cite[Thm.~5.1(a),(b)]{Hum08}, and the second by (\ref{fdM}). We deduce $\Dim_E\mathrm{Ext}_{U(\fg)}^k(L(\mu'),L(\mu)^{w})<+\infty$ for $k\geq 0$ by an obvious d\'evissage. The statement of the lemma then follows by (another) obvious d\'evissage on the constituents of $M'$ and $M^w$.
\end{proof}

Recall from \S\ref{subsec: category} that, for $x\in W(G)$ and $I\subseteq \Delta$, $L^I(x)$ is the unique simple quotient of $U(\fl_I)\otimes_{U(\fb_I)}x\cdot \mu_0$ (and lies in $\cO_{\fl_I,\rm{alg}}^{\fb_I}$).

\begin{lem}\label{lem: rank one case}
Let $x\in W(G)$ and $j\in\Delta$.
\begin{enumerate}[label=(\roman*)]
\item \label{it: rank one 1} We have $L^{\{j\}}(x)\cong L^{\{j\}}(x)^{s_j}$ if and only if $j\notin D_L(x)$, in which case we have $U(\ft)$-equivariant isomorphisms $H^0(\fu_{\{j\}},L^{\{j\}}(x))\cong x\cdot\mu_0$ and $H^1(\fu_{\{j\}},L^{\{j\}}(x))\cong s_jx\cdot\mu_0$.
\item \label{it: rank one 2} If $j\in D_L(x)$, we have $U(\ft)$-equivariant isomorphisms $H^0(\fu_{\{j\}},L^{\{j\}}(x))\cong x\cdot\mu_0$, $H^1(\fu_{\{j\}},L^{\{j\}}(x))=0$, $H^0(\fu_{\{j\}},L^{\{j\}}(x)^{s_j})=0$ and $H^1(\fu_{\{j\}},L^{\{j\}}(x)^{s_j})\cong s_jx\cdot\mu_0$.
\end{enumerate}
\end{lem}
\begin{proof}
We prove \ref{it: rank one 1}. The isomorphism $H^0(\fu_{\{j\}}, L^{\{j\}}(x))\cong x\cdot\mu_0$ always holds by \ref{it: dominance 2} of Lemma \ref{lem: n coh dominance}. By \ref{it: coh 3} of Lemma~\ref{lem: H0 and H1} (applied with the reductive group $L_{\{j\}}$ instead of $G$ and with $I=\emptyset$) and the fact that $H^1(\fu_{\{j\}},L^{\{j\}}(x))$ is a semi-simple $U(\ft)$-module (cf.Proposition~\ref{lem: n coh O}) we deduce that $H^1(\fu_{\{j\}},L^{\{j\}}(x))\cong s_jx\cdot\mu_0$ if $x<s_jx$ (equivalently if $j\notin D_L(x)$), and $H^1(\fu_{\{j\}},L^{\{j\}}(x))=0$ if $s_jx<x$ (equivalently if $j\in D_L(x)$). By \ref{it: H0 conjugate 3} of Lemma~\ref{lem: H0 Weyl conjugate} (applied with $L_{\{j\}}$ and $I=\emptyset$), we have $L^{\{j\}}(x)\cong L^{\{j\}}(x)^{s_j}$ if and only if $j\notin D_L(x)$ (noting that $I'\ne \emptyset$ in \emph{loc.~cit.} if and only if $I'=\{j\}$ if and only if $j\notin D_L(x)$ by Lemma \ref{lem: dominance and left set}).

We prove \ref{it: rank one 2}. By \ref{it: H0 conjugate 2} of Lemma~\ref{lem: H0 Weyl conjugate} (applied with $L_{\{j\}}$ and $I=\emptyset$) if $j\in D_L(x)$ we have $I'=\emptyset$ and $H^0(\fu_{\{j\}},L^{\{j\}}(x)^{s_j})=0$. In view of what was proven before, it remains to show that
$H^1(\fu_{\{j\}},L^{\{j\}}(x)^{s_j})\cong s_jx\cdot\mu_0$ when $j\in D_L(x)$. Assume $j\in D_L(x)$, then $L^{\{j\}}(x)\cong U(\fl_{\{j\}})\otimes_{U(\fb_{\{j\}})}x\cdot\mu_0$ is a free $U(\fu_{\{j\}}^+)$-module with a generator $0\neq v\in L^{\{j\}}(x)_{x\cdot\mu_0}$. So $L^{\{j\}}(x)^{s_j}$ is a free $U(\fu_{\{j\}})=U((\fu_{\{j\}}^+)^{s_j})$-module with a generator $0\neq v^{s_j}\in (L^{\{j\}}(x)^{s_j})_{s_j(x\cdot\mu_0)}$. Recall that $H^1(\fu_{\{j\}},L^{\{j\}}(x)^{s_j})$ is by definition the cokernel of the map of $U(\ft)$-modules
\[L^{\{j\}}(x)^{s_j}\rightarrow L^{\{j\}}(x)^{s_j}\otimes_E \fu_{\{j\}}^\vee,\ m\mapsto u(m)\otimes u^\vee\]
where $\fu_{\{j\}}^\vee$ is the $E$-vector space dual to $\fu_{\{j\}}$, $u\neq 0$ a fixed element of $\fu_{\{j\}}$ and $u^\vee$ the dual basis (recall $\dim_E\fu_{\{j\}}=1$). By an easy computation (we are with $\GL_2$) we have that $v^{s_j}\otimes u^\vee$ spans this cokernel as $E$-vector space. In particular as $U(\ft)$-module the cokernel is isomorphic to $s_j(x\cdot\mu_0)+\al_j=s_j\cdot(x\cdot\mu_0)=s_jx\cdot \mu_0$ where $\al_j=e_j-e_{j+1}\in\Delta$ is the positive simple root also denoted $j$. This finishes the proof of \ref{it: rank one 2}.
\end{proof}

\begin{lem}\label{lem: H1 conjugate}
Let $I,I'\subseteq \Delta$, $1\neq w_1\in W^{I',I}$ and $w\in W(G)$ such that $I'=\Delta\setminus D_L(w)$. Let $x\in W(G)$, then the following statements are equivalent.
\begin{enumerate}[label=(\roman*)]
\item \label{it: H1 conj 1} $\Hom_{U(\fl_I)}(L^I(x), H^1(\fn_I, L(w)^{w_1}))\neq 0$;
\item \label{it: H1 conj 2} $\mathrm{Ext}_{U(\fg)}^1(L(x),L(w)^{w_1})\neq 0$;
\item \label{it: H1 conj 3} $\mathrm{Ext}_{U(\fg)}^1(L(x),L(w))\neq 0$ and $w_1\in W(L_{I'})W(L_{\Delta\setminus D_L(x)})$.
\end{enumerate}
Moreover, if these statements hold, we have a canonical isomorphism of $E$-vector spaces of dimension $\Dim_E \mathrm{Ext}_{U(\fg)}^1(L(x),L(w))$ induced by applying the functor $H^0(\fu_I,-)$:
\begin{equation}\label{equ: H1 conj socle}
\Hom_{U(\fl_I)}(L^I(x), H^1(\fn_I, L(w)^{w_1}))\buildrel\sim\over\longrightarrow \Hom_{U(\ft)}\big(x\cdot\mu_0, H^0(\fu_I, H^1(\fn_I, L(w)^{w_1}))\big).
\end{equation}
\end{lem}
\begin{proof}
By Lemma \ref{lem: dominance and left set} if $w\in W(G)$ and $I'=\Delta\setminus D_L(w)$ then $L(w)$ is in $\cO_{\rm{alg}}^{\fp_{I'}}$ and $I'$ is maximal for that condition. In particular $L(w)^{w'}\cong L(w)$ for $w'\in W(L_{I'})$ (see just before Lemma \ref{lem: H0 Weyl conjugate}) and, if $w_1=w_2w_3$ for some $w_2\in W(L_{\Delta\setminus D_L(w)})$ and $w_3\in W(L_{\Delta\setminus D_L(x)})$, then $L(w)^{w_2}\cong L(w)$ and $L(x)^{w_3^{-1}}\cong L(x)$. Hence, we deduce from (\ref{equ: twist isom}) (with $M_1=L(x)$, $M_2=L(w)^{w_1}$ and $g=w_3^{-1}$):
\begin{equation}\label{chiant}
\mathrm{Ext}_{U(\fg)}^1(L(x),L(w)^{w_1})\cong \mathrm{Ext}_{U(\fg)}^1(L(x)^{w_3^{-1}},L(w)^{w_2})\cong \mathrm{Ext}_{U(\fg)}^1(L(x),L(w)).
\end{equation}
In particular, \ref{it: H1 conj 3} implies \ref{it: H1 conj 2}.

Since $w_1\ne 1$ by \ref{it: H0 conjugate 2} of Lemma~\ref{lem: H0 Weyl conjugate} we have $H^0(\fn_I, L(w)^{w_1})=0$, which together with (\ref{equ: 5terms}) (applied with $M_I=L^I(x)$ and $M=L(w)^{w_1}$) gives
\begin{equation}\label{equ: Ext 1 Weyl conj}
\mathrm{Ext}_{U(\fg)}^1(M^I(x),L(w)^{w_1})\cong \Hom_{U(\fl_I)}(L^I(x), H^1(\fn_I, L(w)^{w_1})).
\end{equation}
By a parallel argument, we have $H^0(\fu, L(w)^{w_1})=0$ and a canonical isomorphism
\begin{equation}\label{equ: Ext 1 Weyl conj prime}
\mathrm{Ext}_{U(\fg)}^1(M(x),L(w)^{w_1})\cong \Hom_{U(\ft)}(x\cdot\mu_0, H^1(\fu, L(w)^{w_1})).
\end{equation}
The vanishing $H^0(\fn_I, L(w)^{w_1})=0$ together with (\ref{equ: double n coh}) (applied with $I'=\emptyset$) also implies
\begin{equation*}
H^1(\fu, L(w)^{w_1})\cong H^0(\fu_{I}, H^1(\fn_{I}, L(w)^{w_1})),
\end{equation*}
hence we have
\begin{equation}\label{equ: u coh conj}
\Hom_{U(\ft)}\big(x\cdot\mu_0, H^0(\fu_{I}, H^1(\fn_{I}, L(w)^{w_1}))\big)\cong \Hom_{U(\ft)}(x\cdot\mu_0, H^1(\fu, L(w)^{w_1})).
\end{equation}
By \ref{it: H0 conjugate 3} of Lemma~\ref{lem: H0 Weyl conjugate} we have $L(w)^{w_1}\notin \cO_{\rm{alg}}^{\fb}$, and thus $\Hom_{U(\fg)}(L(x'),L(w)^{w_1})=0$ for any $L(x')\in \cO_{\rm{alg}}^{\fb}$. Hence, the surjections $M(x)\twoheadrightarrow M^I(x)\twoheadrightarrow L(x)$ induce injections
\begin{equation}\label{equ: Ext1 conj injection}
\mathrm{Ext}_{U(\fg)}^1(L(x),L(w)^{w_1})\hookrightarrow \mathrm{Ext}_{U(\fg)}^1(M^I(x),L(w)^{w_1})\buildrel q_1 \over \hookrightarrow \mathrm{Ext}_{U(\fg)}^1(M(x),L(w)^{w_1}).
\end{equation}
It follows that \ref{it: H1 conj 2} implies $\mathrm{Ext}_{U(\fg)}^1(M^I(x),L(w)^{w_1})\neq 0$ which is equivalent to \ref{it: H1 conj 1} by (\ref{equ: Ext 1 Weyl conj}). Thus \ref{it: H1 conj 2} implies \ref{it: H1 conj 1}.

The injection $q_1$ corresponds under (\ref{equ: Ext 1 Weyl conj}), (\ref{equ: Ext 1 Weyl conj prime}) and (\ref{equ: u coh conj}) to the injection induced by applying the functor $H^0(\fu_I,-)$
\begin{equation}\label{equ: Ext1 conj coh}
\Hom_{U(\fl_I)}(L^I(x), H^1(\fn_I, L(w)^{w_1})) \buildrel {q_2} \over \hookrightarrow \Hom_{U(\ft)}\big(x\cdot\mu_0, H^0(\fu_{I}, H^1(\fn_{I}, L(w)^{w_1}))\big).
\end{equation}
Note that by (\ref{chiant}), (\ref{equ: Ext 1 Weyl conj}), (\ref{equ: Ext 1 Weyl conj prime}) and (\ref{equ: Ext1 conj injection}), in order to prove that $q_2$ (i.e.~(\ref{equ: H1 conj socle})) is an isomorphism it is enough to prove that (\ref{equ: u coh conj}) has dimension $\Dim_E \mathrm{Ext}_{U(\fg)}^1(L(x),L(w))$. Moreover, if \ref{it: H1 conj 1} holds, then (\ref{equ: u coh conj}) must be non-zero by (\ref{equ: Ext1 conj coh}). Consequently, in order to prove that \ref{it: H1 conj 1} implies \ref{it: H1 conj 3} and the last statement of the lemma (and hence the lemma), it suffices to prove:
\begin{claim}\label{claim: H1 conj}
Assume \ that \ (\ref{equ: u coh conj}) \ is \ non-zero, \ then \ \ref{it: H1 conj 3} \ holds \ and \ (\ref{equ: u coh conj}) \ has \ dimension $\Dim_E \mathrm{Ext}_{U(\fg)}^1(L(x),L(w))$.
\end{claim}

We now prove Claim~\ref{claim: H1 conj} by an increasing induction on $\ell(w_1)\geq 1$. Note first that, if we take $k=1$ and $M=L(w)^{w_1}$ in (\ref{equ: special double n coh}), we have the short exact sequence for any $j\in \Delta$
\begin{multline}\label{equ: conj key sequence}
0\rightarrow H^1(\fu_{\{j\}}, H^0(\fn_{\{j\}}, L(w)^{w_1}))\rightarrow H^1(\fu, L(w)^{w_1})\\
\rightarrow H^0(\fu_{\{j\}}, H^1(\fn_{\{j\}}, L(w)^{w_1}))\rightarrow 0.
\end{multline}

We first prove the case $\ell(w_1)=1$, i.e.\ $w_1=s_j$ for some $j\in \Delta$. As $w_1=s_j\in W^{I',I}$ we have $j\notin I\cup I'$ and thus $j\in D_L(w)$. We also have $s_j\fn_{\{j\}}s_j=\fn_{\{j\}}$ and thus (arguing as in (\ref{twisttwist}))
\begin{equation}\label{forgottentoquote}
H^k(\fn_{\{j\}}, L(w)^{s_j})\cong H^k(s_j\fn_{\{j\}}s_j, L(w)^{s_j}) \cong H^k(\fn_{\{j\}}, L(w))^{s_j}
\end{equation}
for $k\geq 0$. By \ref{it: dominance 2} of Lemma \ref{lem: n coh dominance} this implies $H^0(\fn_{\{j\}},L(w)^{s_j})\cong L^{\{j\}}(w)^{s_j}$. Using (\ref{equ: conj key sequence}), the non-vanishing of $\Hom_{U(\ft)}(x\cdot\mu_0, H^1(\fu, L(w)^{w_1}))$ forces either
\begin{equation}\label{equ: conj 1 0 term}
\Hom_{U(\ft)}(x\cdot\mu_0,H^1(\fu_{\{j\}}, L^{\{j\}}(w)^{s_j}))\ne 0
\end{equation}
or
\begin{equation}\label{equ: conj 0 1 term}
\Hom_{U(\ft)}(x\cdot\mu_0,H^0(\fu_{\{j\}},H^1(\fn_{\{j\}},L(w)^{s_j})))\ne 0.
\end{equation}
We treat these two possibilities separately.
\begin{itemize}
\item Since $j\in D_L(w)$, by \ref{it: rank one 2} of Lemma~\ref{lem: rank one case} (applied with $x=w$!) we see that (\ref{equ: conj 1 0 term}) is non-zero if and only if $x=s_{j}w (<w)$, in which case (\ref{equ: conj 1 0 term}) is one dimensional.
\item Now we consider (\ref{equ: conj 0 1 term}). We deduce from \ref{it: coh 3} of Lemma~\ref{lem: H0 and H1} (applied with $I=\{j\}$) that $\mathrm{soc}_{U(\fl_{\{j\}})}(H^1(\fn_{\{j\}},L(w)))\cong \bigoplus_{x'}L^{\{j\}}(x')^{\oplus \mu(w,x')}$ and therefore
\begin{equation}\label{socforgotten}
\mathrm{soc}_{U(\fl_{\{j\}})}(H^1(\fn_{\{j\}},L(w))^{s_j})\cong \bigoplus_{x'}(L^{\{j\}}(x')^{s_j})^{\oplus \mu(w,x')}
\end{equation}
where $x'$ runs through elements of $W(G)$ such that $w\prec x'$ (note that, as $W(L_{\{j\}})=\{1,s_{j}\}$ and $j\in D_L(w)$, $w\prec x'$ also ensures $x'\notin W(L_{\{j\}})w$). It follows from \ref{it: H0 conjugate 2} of Lemma~\ref{lem: H0 Weyl conjugate} (applied with $L_{\{j\}}$ and $I=\emptyset$) that $H^0(\fu_{\{j\}}, L^{\{j\}}(x')^{s_j})$ is non-zero if and only if $s_j\in W(L_{\Delta\setminus D_L(x')})$ if and only if $j\notin D_L(x')$, in which case it is isomorphic to $x'\cdot \mu_0$ by \ref{it: rank one 1} of Lemma \ref{lem: rank one case}. Using Proposition~\ref{lem: n coh O}, each indecomposable direct summand of the $U(\fl_{\{j\}})$-module $H^1(\fn_{\{j\}},L(w))$ has one of the five forms described in \cite[Prop.~3.12]{Hum08}. But we also have
\[H^0\big(\fu_{\{j\}}, \mathrm{soc}_{U(\fl_{\{j\}})}(H^1(\fn_{\{j\}},L(w)))\big)\cong H^0(\fu_{\{j\}}, H^1(\fn_{\{j\}},L(w)))\]
by \ref{it: coh 5} of Lemma~\ref{lem: H0 and H1} (applied with $I=\{j\}$). It is then easy to check that an indecomposable direct summand of $H^1(\fn_{\{j\}},L(w))$ is either irreducible of the form $L^{\{j\}}(x')$ (for $w\prec x'$) or uniserial of length two with socle $L^{\{j\}}(x')$ and cosocle $L^{\{j\}}(s_jx')$ (for $w\prec x'$ and $j\notin D_L(x')$). Since $H^0(\fu_{\{j\}}, L^{\{j\}}(s_jx')^{s_j})=0$ if $j\notin D_L(x')$ (as seen above), this in turns implies using (\ref{forgottentoquote}) (for $k=1$)
\[H^0\big(\fu_{\{j\}}, \mathrm{soc}_{U(\fl_{\{j\}})}(H^1(\fn_{\{j\}},L(w)^{s_j}))\big)\cong H^0(\fu_{\{j\}}, H^1(\fn_{\{j\}},L(w)^{s_j})).\]
By (\ref{socforgotten}) and the discussion just after it, we deduce that (\ref{equ: conj 0 1 term}) is non-zero if and only if $w\prec x$ and $j\notin D_L(x)$, in which case (\ref{equ: conj 0 1 term}) is $\mu(w,x)$-dimensional.
\end{itemize}
Note that in both cases $w_1=s_j\in W(L_{\Delta\setminus D_L(x)})$ and $x\ne w$. Moreover when either of (\ref{equ: conj 1 0 term}), (\ref{equ: conj 0 1 term}) is non-zero, then the other is zero. Using \ref{it: rabiotext 2} of Lemma \ref{rabiotext} (together with the $U(\ft)$-semi-simplicity of $H^1(\fu, L(w)^{w_1})$) we deduce that (\ref{equ: u coh conj}) always has dimension $\Dim_E\mathrm{Ext}_{U(\fg)}^1(L(x),L(w))$. This finishes the proof of Claim~\ref{claim: H1 conj} when $\ell(w_1)=1$.\bigskip

We now assume $\ell(w_1)>1$ and prove the induction step. We choose an arbitrary $j\in D_R(w_1)$ and set $w_2\defeq w_1s_{j}<w_1$. As $s_j\fn_{\{j\}}s_j=\fn_{\{j\}}$, we have as in (\ref{forgottentoquote})
\begin{equation}\label{equ: composed Weyl conj}
H^k(\fn_{\{j\}}, L(w)^{w_1})\cong H^k(s_j\fn_{\{j\}}s_j, L(w)^{w_2s_{j}}) \cong H^k(\fn_{\{j\}}, L(w)^{w_2})^{s_{j}}.
\end{equation}
As $w_1\in W^{I',I}\subseteq W^{I',\emptyset}$ and $w_1=w_2s_{j}>w_2$, we have $w_2\in W^{I',\emptyset}$, which together with $j\in D_R(w_1)\setminus D_R(w_2)$ implies $w_2\in W^{I',\{j\}}$. Note that $w_2\neq 1$ since $\ell(w_1)>1$, and therefore $H^0(\fn_{\{j\}}, L(w)^{w_1})\cong H^0(\fn_{\{j\}}, L(w)^{w_2})^{s_{j}}=0$ by (\ref{equ: composed Weyl conj}) and \ref{it: H0 conjugate 2} of Lemma~\ref{lem: H0 Weyl conjugate}. By (\ref{equ: conj key sequence}) we deduce an isomorphism (using again (\ref{equ: composed Weyl conj})) $H^1(\fu, L(w)^{w_1})\cong H^0(\fu_{\{j\}}, H^1(\fn_{\{j\}}, L(w)^{w_2})^{s_{j}})$. In particular, $\Hom_{U(\ft)}(x\cdot\mu_0, H^1(\fu, L(w)^{w_1}))$ is non-zero if and only if
\begin{equation}\label{equ: composed conj 0 1 term}
\Hom_{U(\ft)}\big(x\cdot\mu_0, H^0(\fu_{\{j\}}, H^1(\fn_{\{j\}}, L(w)^{w_2})^{s_{j}})\big)\ne 0.
\end{equation}
Moreover, as $\ell(w_2)=\ell(w_1)-1$ and $w_2\in W^{I',\{j\}}$, the induction assumption and the discussion before Claim \ref{claim: H1 conj} imply that
\begin{equation}\label{equ: conj induction term}
\Hom_{U(\fl_{\{j\}})}(L^{\{j\}}(x'), H^1(\fn_{\{j\}}, L(w)^{w_2})) \neq 0
\end{equation}
if and only if $\mathrm{Ext}_{U(\fg)}^1(L(x'),L(w))\neq 0$ and $w_2\in W(L_{I'})W(L_{\Delta\setminus D_L(x')})$, in which case (\ref{equ: conj induction term}) has the same dimension as $\mathrm{Ext}_{U(\fg)}^1(L(x'),L(w))$.

Now let $x\in W(G)$ satisfying (\ref{equ: composed conj 0 1 term}), we have the following two cases.
\begin{itemize}
 \item Assume $j\in D_L(x)$. Then $L^{\{j\}}(x)^{s_j}$ is not locally $\fb_{\{j\}}$-finite (\ref{it: H0 conjugate 3} of Lemma \ref{lem: H0 Weyl conjugate} and \ref{it: rank one 1} of Lemma \ref{lem: rank one case}). Let $M_{\{j\}}\defeq U(\fl_{\{j\}})\otimes_{U(\fb_{\{j\}})}x\cdot \mu_0$ which is isomorphic to $L^{\{j\}}(x)$ since $j\in D_L(x)$. By (\ref{equ: composed conj 0 1 term}) there exists a non-zero map $M_{\{j\}}\rightarrow H^1(\fn_{\{j\}}, L(w)^{w_2})^{s_{j}}$, and thus a non-zero map $M_{\{j\}}^{s_j}\rightarrow H^1(\fn_{\{j\}}, L(w)^{w_2})$. This is impossible as $H^1(\fn_{\{j\}}, L(w)^{w_2})$ is locally $\fb_{\{j\}}$-finite (use $w_2\in W^{I',\{j\}}$ and \ref{it: H0 conjugate 1} of Lemma \ref{lem: H0 Weyl conjugate}) but $M_{\{j\}}^{s_j}\cong L^{\{j\}}(x)^{s_j}$ is not. Hence this case can't happen.
\item Assume $j\notin D_L(x)$. Consider the natural embedding
\begin{multline}\label{embeddingrab}
\Hom_{U(\ft)}\big(x\cdot\mu_0, H^0(\fu_{\{j\}}, \mathrm{soc}_{U(\fl_{\{j\}})}(H^1(\fn_{\{j\}}, L(w)^{w_2})^{s_{j}}))\big)\\
\hookrightarrow \Hom_{U(\ft)}(x\cdot\mu_0, H^0(\fu_{\{j\}}, H^1(\fn_{\{j\}}, L(w)^{w_2})^{s_{j}})).
\end{multline}
Assume that (\ref{embeddingrab}) is not an isomorphism. Then $H^1(\fn_{\{j\}}, L(w)^{w_2})^{s_{j}}$ must contain $M_{\{j\}}=U(\fl_{\{j\}})\otimes_{U(\fb_{\{j\}})}x\cdot \mu_0$, and thus $H^1(\fn_{\{j\}}, L(w)^{w_2})$ must contain $M_{\{j\}}^{s_j}$. This is impossible as $H^1(\fn_{\{j\}}, L(w)^{w_2})$ is locally $\fb_{\{j\}}$-finite but $\mathrm{soc}_{U(\fl_{\{j\}})}(M_{\{j\}}^{s_j})\cong L^{\{j\}}(s_jx)^{s_j}$ is not (same argument as in the previous case). Hence (\ref{embeddingrab}) is an isomorphism, and in particular its left hand side is non-zero by (\ref{equ: composed conj 0 1 term}). But since $j\notin D_L(x)$, we have $s_j\in W(L_{\Delta\setminus D_L(x)})$ and thus $L^{\{j\}}(x)^{s_j}\cong L^{\{j\}}(x)$ (\ref{it: H0 conjugate 3} of Lemma \ref{lem: H0 Weyl conjugate}). Hence the two (non-zero) $U(\fl_{\{j\}})$-modules in (\ref{embeddingrab}), once ``untwisted'' by $s_j$, become isomorphic to $\Hom_{U(\fl_{\{j\}})}(L^{\{j\}}(x), H^1(\fn_{\{j\}}, \!L(w)^{w_2}))$, which is thus also non-zero and has the same dimension as $\mathrm{Ext}_{U(\fg)}^1(L(x),L(w))$ by what follows (\ref{equ: conj induction term}) applied to $x'=x$.
\end{itemize}
We have shown that (\ref{equ: composed conj 0 1 term}) implies $\mathrm{Ext}_{U(\fg)}^1(L(x),L(w))\neq 0$, $w_2\in W(L_{I'})W(L_{\Delta\setminus D_L(x)})$ (using the induction) and $j\notin D_L(x)$, which implies $w_1=w_2s_j\in W(L_{I'})W(L_{\Delta\setminus D_L(x)})$. Moreover, we have also seen that (\ref{equ: composed conj 0 1 term}) has the same dimension as $\mathrm{Ext}_{U(\fg)}^1(L(x),L(w))$. The proof of the induction step is thus finished.
\end{proof}

\begin{rem}\label{mergerem}
Let $I\subseteq \Delta$, $w\in W(G)$, $I'=\Delta\setminus D_L(w)$ and $1\neq w_1\in W^{I',I}$.
\begin{enumerate}[label=(\roman*)]
\item \label{rem: H1 conj u coh} Let $I\subseteq \Delta$, $w\in W(G)$, $I'\defeq\Delta\setminus D_L(w)$ and $1\neq w_1\in W^{I',I}$, it follows from (\ref{equ: H1 conj socle}), the fact that $L(w)$ and $L(w)^{w_1}$ have the same infinitesimal character (which is the one of $L(\mu_0)\cong L(\mu_0)^{w_1}$) and Lemma \ref{lem: HC n coh} (together with Harish-Chandra's theorem) that the inclusion $\mathrm{soc}_{U(\fl_I)}(H^1(\fn_I,L(w)^{w_1}))\subseteq H^1(\fn_I,L(w)^{w_1})$ induces an isomorphism
\[H^0(\fu_I, \mathrm{soc}_{U(\fl_I)}(H^1(\fn_I,L(w)^{w_1})))\buildrel\sim\over\longrightarrow H^0(\fu_I, H^1(\fn_I,L(w)^{w_1})).\]
\item \label{rem: H1 conj Ext1 vanishing}
We deduce from Lemma~\ref{lem: H1 conjugate} and (\ref{equ: Ext 1 Weyl conj}) that, for $x\in W(G)$ such that $I=\Delta\setminus D_L(x)$, we must have
\begin{multline*}
\mathrm{Ext}_{U(\fg)}^1(L(x),L(w)^{w_1})=\mathrm{Ext}_{U(\fg)}^1(M^I(x),L(w)^{w_1})\\
=\Hom_{U(\fl_I)}(L^I(x),H^1(\fn_I, L(w)^{w_1}))=0.
\end{multline*}
\end{enumerate}
\end{rem}

\begin{lem}\label{lem: special H2}
Let $j\in\Delta$ and $x,w\in W(G)$ with $D_L(x)=D_L(w)=\{j\}$.
Then we have
\begin{equation}\label{equ: H2 conj vanishing}
\Hom_{U(\ft)}(x\cdot\mu_0, H^2(\fu, L(w)^{s_j}))=0.
\end{equation}
\end{lem}
\begin{proof}
We deduce from (\ref{equ: special double n coh}) the following exact sequence
\begin{equation}\label{equ: conj u coh 1}
0\rightarrow H^1(\fu_{\{j\}},H^1(\fn_{\{j\}}, L(w)^{s_j}))\rightarrow H^2(\fu, L(w)^{s_j})\rightarrow H^0(\fu_{\{j\}},H^2(\fn_{\{j\}}, L(w)^{s_j}))\rightarrow 0
\end{equation}
and recall that $s_j\fn_{\{j\}}s_j=\fn_{\{j\}}$ implies $H^k(\fn_{\{j\}}, L(w)^{s_j})\cong H^k(\fn_{\{j\}}, L(w))^{s_j}$ for $k\geq 0$ (see (\ref{twisttwist})). By Lemma~\ref{lem: rank one case}, for $x'\in W(G)$, we have $H^0(\fu_{\{j\}}, L^{\{j\}}(x')^{s_j})\neq 0$ if and only if $j\notin D_L(x')$, in which case $H^0(\fu_{\{j\}}, L^{\{j\}}(x')^{s_j})\cong H^0(\fu_{\{j\}}, L^{\{j\}}(x'))\cong x'\cdot\mu_0$. As $j\in D_L(x)$, $H^2(\fn_{\{j\}}, L(w)^{s_j})\cong H^2(\fn_{\{j\}}, L(w))^{s_j}$ and $H^2(\fn_{\{j\}}, L(w) \in \cO^{\fb_{\{j\}}}_{\fl_{\{j\}},\rm{alg}}$ by Proposition~\ref{lem: n coh O},
we deduce by d\'evissage on the constituents of $H^2(\fn_{\{j\}}, L(w)^{s_j})$
\[\Hom_{U(\ft)}(x\cdot\mu_0, H^0(\fu_{\{j\}},H^2(\fn_{\{j\}}, L(w)^{s_j})))=0.\]
Assume as a contradiction that (\ref{equ: H2 conj vanishing}) fails. Then we deduce by \ref{equ: conj u coh 1}
\begin{equation}\label{equ: H2 conj 1 1}
\Hom_{U(\ft)}(x\cdot\mu_0, H^1(\fu_{\{j\}},H^1(\fn_{\{j\}}, L(w)^{s_j})))\neq 0.
\end{equation}
Using \cite[Thm.~1.10]{Hum08} and a standard argument, one sees that any irreducible constituent of $H^1(\fn_{\{j\}}, L(w))$ has the form $L^{\{j\}}(x'')$ with $x''\in\{x',s_jx'\}$ for some $L^{\{j\}}(x')$ showing up in the socle of $H^1(\fn_{\{j\}}, L(w))$. By \ref{it: coh 5} of Lemma~\ref{lem: H0 and H1} we have
\[H^0\big(\fu_{\{j\}}, \mathrm{soc}_{U(\fl_{\{j\}})}(H^1(\fn_{\{j\}}, L(w)))\big)\cong H^0(\fu_{\{j\}}, H^1(\fn_{\{j\}}, L(w)))\]
and by Proposition~\ref{lem: n coh O} $H^1(\fn_{\{j\}}, L(w))$ is in $\cO^{\fb_{\{j\}}}_{\fl_{\{j\}},\rm{alg}}$. Using the explicit list of all indecomposable objects in $\cO^{\fb_{\{j\}}}_{\fl_{\{j\}},\rm{alg}}$ (see \cite[\S 3.12]{Hum08}) we easily deduce that each indecomposable direct summand of $H^1(\fn_{\{j\}}, L(w))$ is either irreducible of the form $L^{\{j\}}(x')$, or uniserial of length two with socle $L^{\{j\}}(x')$, cosocle $L^{\{j\}}(s_jx')$ and $j\notin D_L(x')$. Let $M_{\{j\}}$ be such an indecomposable direct summand. It follows from Lemma~\ref{lem: rank one case} (and a straightforward d\'evissage using $H^2(\fu_{\{j\}},-)=0$) that $H^1(\fu_{\{j\}}, M_{\{j\}}^{s_j})\cong s_jx'\cdot \mu_0$ if $M_{\{j\}}=L^{\{j\}}(x')$, and $H^1(\fu_{\{j\}}, M_{\{j\}}^{s_j})\cong s_jx'\cdot \mu_0\oplus x'\cdot \mu_0$ if $M_{\{j\}}$ has length two with socle $L^{\{j\}}(x')$ (and $j\notin D_L(x')$). Since $j\in D_L(x)$, it follows that (\ref{equ: H2 conj 1 1}) forces $x=s_jx'>x'$ for some $x'$ such that $L^{\{j\}}(x')$ shows up in the socle of $H^1(\fn_{\{j\}}, L(w))$ and $j\notin D_L(x')$ (so $x'\prec x$). By \ref{it: coh 3} of Lemma~\ref{lem: H0 and H1} such an $x'$ satisfies $w\prec x'$.
The existence of such a triple $w\prec x'\prec x$ with $D_L(x)=D_L(w)=\{j\}$ and $j\notin D_L(x')$ contradicts
Lemma~\ref{lem: saturated triple}.
\end{proof}

\begin{lem}\label{lem: special H2 I}
Let $j\in\Delta$, $x,w\in W(G)$ with $D_L(x)=D_L(w)=\{j\}$ and $I\defeq \Delta\setminus\{j\}$.
Then we have
\begin{equation}\label{equ: special H2 I}
\Hom_{U(\fl_I)}(L^I(x), H^2(\fn_I, L(w)^{s_j}))=0.
\end{equation}
\end{lem}
\begin{proof}
As $j\in D_L(w)$, it follows from \ref{it: H0 conjugate 1}, \ref{it: H0 conjugate 2} of Lemma~\ref{lem: H0 Weyl conjugate} (applied with $I'= \Delta\setminus D_L(w)=I$) that $H^k(\fn_I, L(w)^{s_j})\in\cO^{\fb_I}_{\fl_I,\rm{alg}}$ for $k\geq 1$ and that $H^0(\fn_I, L(w)^{s_j})=0$. As $I\cap D_L(x)=\emptyset$ and thus $x$ has minimal length in $W(L_I)x$, we deduce from Lemma~\ref{lem: n coh wt} (applied with $\fl_I$ instead of $\fg$ and with $I$ there being $\emptyset$) that $H^{\ell}(\fu_I,L^I(w'))_{x\cdot\mu_0}=0$ for $\ell\geq 1$ and any $w'\in W(G)$. By a d\'evissage on the constituents of $H^k(\fn_I, L(w)^{s_j})$ we obtain for $\ell\geq 1$ and $k\geq 0$
\begin{equation}\label{equ: higher u coh}
\Hom_{U(\ft)}\big(x\cdot\mu_0, H^{\ell}(\fu_I, H^k(\fn_I, L(w)^{s_j}))\big)=0.
\end{equation}
By (\ref{equ: double n coh}) (applied with $I'=\emptyset$) we have a spectral sequence of semi-simple $U(\ft)$-modules
\[H^{\ell}(\fu_I, H^k(\fn_I, L(w)^{s_j})) \implies H^{\ell+k}(\fu, L(w)^{s_j}).\]
Applying the exact functor $\Hom_{U(\ft)}(x\cdot\mu_0,-)$ to this spectral sequence, we deduce in particular from (\ref{equ: higher u coh}) and Lemma~\ref{lem: special H2}:
\begin{equation}\label{equ: special H2 n u}
\Hom_{U(\ft)}(x\cdot\mu_0, H^0(\fu_I, H^2(\fn_I, L(w)^{s_j})))\cong \Hom_{U(\ft)}(x\cdot\mu_0, H^2(\fu, L(w)^{s_j}))=0.
\end{equation}
Now, from the surjection $U(\fl_I)\otimes_{U(\fb_I)}x\cdot\mu_0\twoheadrightarrow L^I(x)$ and (\ref{equ: g spectral seq 0}) we have
\begin{multline*}
\Hom_{U(\fl_I)}(L^I(x), H^2(\fn_I, L(w)^{s_j}))\hookrightarrow \Hom_{U(\fl_I)}(U(\fl_I)\otimes_{U(\fb_I)}x\cdot\mu_0, H^2(\fn_I, L(w)^{s_j}))\\
\cong \Hom_{U(\ft)}(x\cdot\mu_0, H^0(\fu_I, H^2(\fn_I, L(w)^{s_j}))),
\end{multline*}
which together with (\ref{equ: special H2 n u}) gives (\ref{equ: special H2 I}).
\end{proof}

\begin{prop}\label{prop: typical Ext2}
Let $j\in\Delta$, $w\in W(G)$ such that $D_L(w)=\{j\}$ and $I\defeq \Delta\setminus\{j\}$.
Let $S_0\defeq \{x'\mid x'\in W(L_I)w, \ell(x')=\ell(w)+1, j\notin D_L(x')\}$.
Then we have
\begin{equation*}
\Dim_E \mathrm{Ext}_{U(\fg)}^2(M^{I}(w),L(w)^{s_j})=\#S_0
\end{equation*}
and
\begin{equation*}
\Dim_E \mathrm{Ext}_{U(\fg)}^2(M^{I}(x),L(w)^{s_j})=0
\end{equation*}
for each $x\neq w$ satisfying $D_L(x)=\{j\}$.
\end{prop}
\begin{proof}
Let $x\in W(G)$ such that $D_L(x)=\{j\}$ (allowing $x=w$). As $j\in D_L(w)$, it follows from \ref{it: H0 conjugate 2} of Lemma~\ref{lem: H0 Weyl conjugate} that $H^0(\fn_I, L(w)^{s_j})=0$, which together with Lemma~\ref{lem: special H2 I} and (\ref{equ: g spectral seq}) give an isomorphism
\begin{equation}\label{equ: Ext2 conj H1}
\mathrm{Ext}_{U(\fg)}^2(M^I(x),L(w)^{s_j})\cong \mathrm{Ext}_{U(\fl_I)}^1(L^I(x), H^1(\fn_I, L(w)^{s_j})).
\end{equation}
For any $x'\in W(G)$ we have by Lemma~\ref{lem: H1 conjugate} that
\[\Hom_{U(\fl_I)}(L^I(x'), \mathrm{soc}_{U(\fl_I)}H^1(\fn_I, L(w)^{s_j}))=\Hom_{U(\fl_I)}(L^I(x'), H^1(\fn_I, L(w)^{s_j}))\neq 0\]
if and only if $\mathrm{Ext}_{U(\fg)}^1(L(x'),L(w))\neq 0$ and $j\notin D_L(x')$, in which case both dimensions are $\Dim_E\mathrm{Ext}_{U(\fg)}^1(L(x'),L(w))$. As $j\in D_L(x)$, the vector space in (\ref{equ: H1 conj socle}) is zero, which together with \ref{rem: H1 conj u coh} of Remark \ref{mergerem} and Lemma~\ref{lem: Ext1 with socle} (which can then be applied to $\fl_I$ instead of $\fg$, noting that $L^I(x)$ is a twist of $L^I(\mu_0)=L^I(1)$ since $D_L(x)=\{j\}$) give an isomorphism by (\ref{equ: Ext2 conj H1})
\begin{equation}\label{equ: Ext2 conj H1 socle}
\mathrm{Ext}_{U(\fg)}^2(M^I(x),L(w)^{s_j})\cong \mathrm{Ext}_{U(\fl_I)}^1(L^I(x), \mathrm{soc}_{U(\fl_I)}(H^1(\fn_I, L(w)^{s_j}))).
\end{equation}
By the discussion following (\ref{equ: Ext2 conj H1}), we see that the dimension of (\ref{equ: Ext2 conj H1 socle}) is
\begin{equation}\label{sumextext}
\sum_{x'}\Dim_E \mathrm{Ext}_{U(\fl_I)}^1(L^I(x),L^I(x')) \Dim_E \mathrm{Ext}_{U(\fg)}^1(L(x'),L(w))
\end{equation}
where $x'$ runs through the elements of $W(G)$ such that $j\notin D_L(x')$ and
\[\mathrm{Ext}_{U(\fl_I)}^1(L^I(x),L^I(x'))\neq 0\neq \mathrm{Ext}_{U(\fg)}^1(L(x'),L(w)).\]
Let $x'\in W(G)$ such that $x'\ne x$ and $x'\ne w$. As $x$ is minimal in $W(L_I)x$, by \ref{it: rabiotext 2} of Lemma \ref{rabiotext} (applied to $\fl_I$ instead of $\fg$) and \ref{it: block 3} of Lemma \ref{lem: central component} (with Harish-Chandra's theorem), we have $\mathrm{Ext}_{U(\fl_I)}^1(L^I(x),L^I(x'))\neq 0$ if and only if $x'\in W(L_I)x$ and $\ell(x')=\ell(x)+1$, in which case it has dimension $\mu(x,x')=1$. Likewise by \ref{it: rabiotext 2} of Lemma \ref{rabiotext} $\mathrm{Ext}_{U(\fg)}^1(L(x'),L(w))\neq 0$ if and only if either $x'\prec w$ or $w\prec x'$, in which case it has dimension $\mu(x',w)$ or $\mu(w,x')$. So in (\ref{sumextext}) we sum up over those $x'\in W(G)$ such that $j\notin D_L(x')$, $x'\in W(L_I)x$ with $\ell(x')=\ell(x)+1$ (which implies $x\prec x'$), and either $x'\prec w$ or $w\prec x'$.

Since $x\prec x'$, we can't have $x'\prec w$ because this would contradict Lemma~\ref{lem: saturated triple}. Hence we have $w\prec x'$. Using Lemma \ref{lem: dominance control}, $D_L(w)=\{j\}$ and $j\notin D_L(x')$ then force $x'\in W(L_I)w$ and $\ell(x')=\ell(w)+1$ (hence $\mu(w,x')=1$), which together with $x'\in W(L_I)x$ and $D_L(x)=D_L(w)=\{j\}$ force $x=w$. Consequently, we conclude that (\ref{equ: Ext2 conj H1 socle}) is non-zero only when $x=w$, in which case it has dimension $\#S_0$.
\end{proof}

\begin{lem}\label{lem: H0 of Ext}
Let $x,w\in W(G)$ and $I\subseteq \Delta$ such that $x\in W(L_I)w$, $x>w$ and $\ell(x)=\ell(w)+1$. Let $M\in \cO_{\rm{alg}}^{\fb}$ be the (unique) length two object with socle $L(x)$ and cosocle $L(w)$ (\ref{it: rabiotext 2} of Lemma \ref{rabiotext}). Then $H^0(\fn_I, M)$ is the unique length two object in $\cO_{\fl_I,\rm{alg}}^{\fb_I}$ with socle $L^I(x)$ and cosocle $L^I(w)$.
\end{lem}
\begin{proof}
Using \ref{it: dominance 2} of Lemma \ref{lem: n coh dominance} the short exact sequence $0\rightarrow L(x)\rightarrow M\rightarrow L(w)\rightarrow 0$ induces an exact sequence
\begin{equation*}
0\rightarrow L^I(x)\rightarrow H^0(\fn_I, M)\xrightarrow{q} L^I(w)\rightarrow H^1(\fn_I, L(x)).
\end{equation*}
Let $\xi:Z(\fl_I)\rightarrow E$ be the unique character such that $L^I(w)_\xi\ne 0$, or equivalently $L^I(x)_\xi\ne 0$. Then by \ref{it: dominance 3} of Lemma \ref{lem: n coh dominance} $H^1(\fn_I, L(x))_\xi=0$, which implies that the map $L^I(w)=L^I(w)_\xi\rightarrow H^1(\fn_I, L(x))$ is $0$ and hence that $q$ is surjective.

It remains to prove that the short exact sequence $0\rightarrow L^I(x)\rightarrow H^0(\fn_I, M)\xrightarrow{q} L^I(w)\rightarrow 0$ is non-split. For $M'$ in $\cO_{\rm{alg}}^{\fb}$ write $[M':L(x)]\in \Z_{\geq 0}$ for the multiplicity of $L(x)$ in $M'$. Since $M(w)\cong U(\fg)\otimes_{U(\fp_{I})}(U(\fl_{I})\otimes_{U(\fb_{I})}w\cdot \mu_0)$ we have the following obvious equality (see (\ref{belongtoOb}) for $M^I(x')$)
\begin{equation}\label{equ: mult one in Verma}
[M(w):L(x)]=\sum_{x'\in W(G)}e_I(w,x')[M^I(x'):L(x)]
\end{equation}
where $e_I(w,x')\in \Z_{\geq 0}$ is the multiplicity of $L^I(x')$ in $U(\fl_I)\otimes_{U(\fb_I)}w\cdot\mu_0$.
As $x\in W(L_I)w$, $x>w$ and $\ell(x)=\ell(w)+1$, we have $e_I(w,x)=1$. As $x>w$ and $\ell(x)=\ell(w)+1$ we also have $[M(w):L(x)]=1$. As $[M^I(x):L(x)]=1$ and $e_I(w,w)=1$, we see that (\ref{equ: mult one in Verma}) and $[M(w):L(x)]=1$ force $[M^I(w):L(x)]=0$. Now, assume on the contrary that the above short exact sequence splits. Then $\Hom_{U(\fg)}(M^I(w),M)\cong \Hom_{U(\fl_I)}(L^I(w),H^0(\fn_I,M))\neq 0$, and in particular $L(x)$ must be a constituent of $M^I(w)$ which contradicts $[M^I(w):L(x)]=~0$.
\end{proof}

The following lemma is essentially proved in \cite[(4.71)]{Schr11}, we provide a self-contained proof for the reader's convenience.

\begin{lem}\label{lem: nonsplit rank two}
Let $w\in W(G)$ and $j_1,j_2\in\Delta$ such that $|j_1-j_2|=1$, $s_{j_1}w>w$ and $D_L(w)\cap \{j_1,j_2\}=\emptyset$. Let $M$ be the unique length two $U(\fl_{\{j_1,j_2\}})$-module with socle $L^{\{j_1,j_2\}}(s_{j_1}w)$ and cosocle $L^{\{j_1,j_2\}}(w)$ (\ref{it: rabiotext 2} of Lemma \ref{rabiotext}). Then we have
\begin{enumerate}[label=(\roman*)]
\item \label{it: rank two H1} the $U(\fl_{\{j_2\}})$-module $H^1(\fl_{\{j_1,j_2\}}\cap \fn_{\{j_2\}}, M^{s_{j_1}})$ has length two with socle $L^{\{j_2\}}(s_{j_2}s_{j_1}w)$ and cosocle $L^{\{j_2\}}(s_{j_1}w)$;
\item the $U(\fl_{\{j_2\}})$-module $H^2(\fl_{\{j_1,j_2\}}\cap \fn_{\{j_2\}}, M^{s_{j_1}})$ has length two with socle $L^{\{j_2\}}(s_{j_2}s_{j_1}s_{j_2}w)$ and cosocle $L^{\{j_2\}}(s_{j_1}s_{j_2}w)$.
\end{enumerate}
\end{lem}
\begin{proof}
Since $D_L(w)\cap \{j_1,j_2\}=\emptyset$, $L^{\{j_1,j_2\}}(w)$ is a twist of $L^{\{j_1,j_2\}}(\mu_0)=L^{\{j_1,j_2\}}(1)$, so that ``untwisting everything'', we can reduce to the case $\fg=\fg\fl_3$, $\Delta=\{j_1,j_2\}$, $w=1$ and $x=s_{j_1}$, in which case we have natural $U(\fl_{\{j_2\}})$-equivariant isomorphisms
\[M=M^{\{j_2\}}(1)=U(\fg)\otimes_{U(\fp_{\{j_2\}})}L^{\{j_2\}}(1)\cong U(\fn_{\{j_2\}}^+)\otimes_E L^{\{j_2\}}(1).\]
Write $\al_s$ for the positive simple root corresponding to $j_s$ for $s=1,2$, then $\fn_{\{j_2\}}\cong \fn_{-\al_1}\oplus \fn_{-\al_1-\al_2}$ where $\fn_{\al}\subseteq \fg$ is the one dimensional subspace corresponding to a root $\al\in\Phi$. Then $M^{s_{j_1}}\cong U(s_{j_1}\fn_{\{j_2\}}^+s_{j_1})\otimes_E L^{\{j_2\}}(1)^{s_{j_1}}\cong U(\fn_{-\al_1})\otimes_E U(\fn_{\al_2})\otimes_E L^{\{j_2\}}(1)^{s_{j_1}}$ is a free $U(\fn_{-\al_1})$-module (of infinite rank). We have $H^0(\fn_{-\al_1}, U(\fn_{-\al_1}))=0$ and a $U(\ft)$-equivariant isomorphism $H^1(\fn_{-\al_1}, U(\fn_{-\al_1}))\cong \al_1$ (use $U(\fn_{-\al_1})\cong E[X]$). As $\fn_{\{j_2\}}^+$ and hence $s_{j_1}\fn_{\{j_2\}}^+s_{j_1}$ are commutative Lie algebras, we deduce $H^0(\fn_{-\al_1}, M^{s_{j_1}})=0$. From the (analogue of) (\ref{equ: special double n coh}) for the spectral sequence $H^{\ell_1}(\fn_{-\al_1-\al_2}, H^{\ell_2}(\fn_{-\al_1}, M^{s_{j_1}}))\implies H^{\ell_1+\ell_2}(\fn_{\{j_2\}}, M^{s_{j_1}})$ we deduce then $H^k(\fn_{\{j_2\}}, M^{s_{j_1}})=0$ for $k\notin \{1,2\}$ and $U(\fb_{\{j_2\}}^+)$-equivariant isomorphisms for $k\in \{1,2\}$:
\begin{eqnarray*}
H^k(\fn_{\{j_2\}}, M^{s_{j_1}})&\cong &H^{k-1}(\fn_{-\al_1-\al_2}, H^1(\fn_{-\al_1},M^{s_{j_1}}))\\
&\cong &H^{k-1}\big(\fn_{-\al_1-\al_2}, U(\fn_{\al_2})\otimes_E L^{\{j_2\}}(1)^{s_{j_1}}\otimes_E H^1(\fn_{-\al_1},U(\fn_{-\al_1}))\big)\\
&\buildrel q \over \cong& U(\fn_{\al_2})\otimes_E\left(H^{k-1}(\fn_{-\al_1-\al_2}, L^{\{j_2\}}(1)^{s_{j_1}})\otimes_E \al_1\right)\\
&\cong& U(\fn_{\al_2})\otimes_E \left(H^{k-1}(\fn_{-\al_2}, L^{\{j_2\}}(1))^{s_{j_1}}\otimes_E \al_1\right)\\
&\cong& U(\fn_{\al_2})\otimes_E s_{j_1}\cdot\mu_k = U(\fu_{\{j_2\}}^+)\otimes_E s_{j_1}\cdot\mu_k\\
&\cong& U(\fl_{\{j_2\}})\otimes_{U(\fb_{\{j_2\}})} s_{j_1}\cdot\mu_k
\end{eqnarray*}
where $\mu_k\defeq H^{k-1}(\fn_{-\al_2}, L^{\{j_2\}}(1))$, where by definition $\fu_{\{j_2\}}^+=\fn_{\al_2}$, and where the isomorphism $q$ above is checked to be $U(\fb_{\{j_2\}}^+)$-equivariant using that $\fn_{-\al_1}$ acts trivially on $L^{\{j_2\}}(1)^{s_{j_1}}\otimes_E H^1(\fn_{-\al_1},U(\fn_{-\al_1}))$. By \ref{it: rank one 1} of Lemma~\ref{lem: rank one case} we have $\mu_1=\mu_0$ and $\mu_2=s_{j_2}\cdot\mu_0$, and thus $U(\fl_{\{j_2\}})$-equivariant isomorphisms
\[H^1(\fn_{\{j_2\}}, M^{s_{j_1}})\cong U(\fl_{\{j_2\}})\otimes_{U(\fb_{\{j_2\}})} s_{j_1}\cdot\mu_0,\ \ H^2(\fn_{\{j_2\}}, M^{s_{j_1}})\cong U(\fl_{\{j_2\}})\otimes_{U(\fb_{\{j_2\}})} s_{j_1}s_{j_2}\cdot\mu_0,\]
where we use the standard fact that a $U(\fl_{\{j_2\}})$-module which is isomorphic to $U(\fl_{\{j_2\}})\otimes_{U(\fb_{\{j_2\}})}\mu$ as $U(\fb_{\{j_2\}}^+)$-module must in fact be isomorphic to $U(\fl_{\{j_2\}})\otimes_{U(\fb_{\{j_2\}})}\mu$ as $U(\fl_{\{j_2\}})$-module. This finishes the proof.
\end{proof}

\begin{lem}\label{lem: nonsplit H1 conj}
Let $w\in W(G)$, $j\in D_L(w)$ and $I\defeq \Delta\setminus D_L(w)$. Let $M\in \cO_{\rm{alg}}^{\fb}$ be the unique length two $U(\fg)$-module with socle $L(w)$ and cosocle $L(s_jw)$. The $U(\fl_I)$-module $L^I(w)$ appears with multiplicity one in $H^1(\fn_I, M^{s_j})\in {\cO}_{\fl_I,\rm{alg}}^{\fb_I}$, and $H^1(\fn_I, M^{s_j})$ contains a (unique) $U(\fl_I)$-submodule $M_I^1$ with cosocle $L^I(w)$ which fits into a short exact sequence:
\begin{equation}\label{equ: nonsplit H1 conj}
0\longrightarrow\!\!\bigoplus_{\overset{x\in W(L_I)w}{\overset{\ell(x)=\ell(w)+1}{j\notin D_L(x)}}}\!\!\!\!\!\!L^I(x)\longrightarrow M_I^1\longrightarrow L^I(w)\longrightarrow 0.
\end{equation}
\end{lem}
\begin{proof}
Note first that the existence of $M$ is clear from \ref{it: Ext O 3} of Lemma~\ref{lem: Ext 1 category O} and that $H^1(\fn_I, M^{s_j})$ is in ${\cO}_{\fl_I,\rm{alg}}^{\fb_I}$ by \ref{it: H0 conjugate 1} of Lemma \ref{lem: H0 Weyl conjugate} applied with $I'=\emptyset$.

As $j\notin D_L(s_jw)$, by Lemma \ref{lem: dominance and left set} and \ref{it: H0 conjugate 3} of Lemma \ref{lem: H0 Weyl conjugate} we have $L(s_jw)^{s_j}\cong L(s_jw)\in \cO_{\rm{alg}}^{\fb}$ and $L(w)^{s_j}\notin \cO_{\rm{alg}}^{\fb}$. Since $L(w)^{s_j}$ is the socle of $M^{s_j}$, we deduce
\begin{equation}\label{equ: H0 conj nonsplit}
H^0(\fn_I, M^{s_j})=0,
\end{equation}
otherwise $M^{s_j}$ would contain a non-zero $U(\fg)$-submodule in $\cO_{\rm{alg}}^{\fb}$ using (\ref{equ: g spectral seq 0}). By \ref{it: dominance 2} of Lemma \ref{lem: n coh dominance} we have $H^0(\fn_I, L(s_jw)^{s_j})\cong L^I(s_jw)$. The short exact sequence $0\rightarrow L(w)^{s_j}\rightarrow M^{s_j}\rightarrow L(s_jw)^{s_j}\rightarrow 0$ then induces the long exact sequence of $U(\fl_I)$-modules
\begin{multline}\label{lesI}
0\rightarrow H^1(\fn_I, L(w)^{s_j})/L^I(s_jw)\rightarrow H^1(\fn_I, M^{s_j})\rightarrow H^1(\fn_I, L(s_jw))\\
\rightarrow H^2(\fn_I, L(w)^{s_j})\rightarrow H^2(\fn_I, M^{s_j})\rightarrow H^2(\fn_I, L(s_jw)).
\end{multline}
If $L^I(w)$ appears with multiplicity one in $H^1(\fn_I, M^{s_j})$, we define $M_I^1$ to be the unique $U(\fl_I)$-submodule of $H^1(\fn_I, M^{s_j})$ with cosocle $L^I(w)$.\bigskip

We first prove that $L^I(w)$ appears with multiplicity one in $H^1(\fn_I, L(s_jw))$. It follows from \ref{it: coh 5} of Lemma~\ref{lem: H0 and H1} (applied with $x=s_jw$) that $L^I(w)$ appears with multiplicity one in the socle of $H^1(\fn_I, L(s_jw))$. As $w$ is minimal in $W(L_I)w$ (since $I=\Delta\setminus D_L(w)$) and $H^1(\fn_I, L(s_jw))\in {\cO}_{\fl_I,\rm{alg}}^{\fb_I}$ is semi-simple as $U(\ft)$-module (Proposition~\ref{lem: n coh O}), we must have the inclusion
\[H^1(\fn_I, L(s_jw))_{w\cdot\mu_0}\subseteq H^0(\fu_I, H^1(\fn_I, L(s_jw)))\]
(otherwise, apply $\fu_I$ to any vector in $H^1(\fn_I, L(s_jw))_{w\cdot\mu_0}$). Together with (\ref{equ: coh of socle}), this implies that any copy of $L^I(w)$ that appears in $H^1(\fn_I, L(s_jw))$ must appear in its socle, and thus $L^I(w)$ has multiplicity one in $H^1(\fn_I, L(s_jw))$.\bigskip

We now prove that $L^I(w)$ is not a Jordan-H\"older factor of $H^1(\fn_I, L(w)^{s_j})$. By \ref{rem: H1 conj u coh} of Remark \ref{mergerem} we have an isomorphism
\[H^0(\fu_I, \mathrm{soc}_{U(\fl_I)}(H^1(\fn_I, L(w)^{s_j})))\buildrel\sim\over\longrightarrow H^0(\fu_I, H^1(\fn_I, L(w)^{s_j})).\]
Hence, by the same argument as in the previous paragraph, it is enough to prove that $L^I(w)$ is not in the socle of $H^1(\fn_I, L(w)^{s_j})$. But this follows from the first part of Lemma~\ref{lem: H1 conjugate} (as $w_1=s_j\notin W(L_I)$). Note that we can then apply Lemma~\ref{lem: Ext1 with socle} (with $\fl_I$ instead of $\fg$ and noting that $I=\Delta\setminus D_L(w)$ implies $w\cdot \mu_0 \in \Lambda_I^{\dom}$ hence $L^I(w)\cong L^I(1)$ up to twist) and deduce
\begin{equation}\label{soctonosoc}
\mathrm{Ext}_{U(\fl_I)}^1(L^I(w), \mathrm{soc}_{U(\fl_I)}(H^1(\fn_I, L(w)^{s_j})))\buildrel\sim\over\longrightarrow \mathrm{Ext}_{U(\fl_I)}^1(L^I(w), H^1(\fn_I, L(w)^{s_j})).
\end{equation}

Assuming that $M_I^1$ exists (i.e.\ that $L^I(w)$ ``survives'' in $H^1(\fn_I, M^{s_j})$), we prove that the radical $\mathrm{rad}(M_I^1)$ (see \S\ref{generalnotation}) is contained in the kernel of the short exact sequence (\ref{equ: nonsplit H1 conj}). Note first that, since $s_j\notin W(L_I)$, by Lemma \ref{lem: HC} we have $\mathrm{Ext}_{U(\fl_I)}^k(L^I(w), L^I(s_jw))=0$ for $k\geq 0$ and thus
\[\mathrm{Ext}_{U(\fl_I)}^1(L^I(w), H^1(\fn_I, L(w)^{s_j})) \buildrel\sim\over\longrightarrow \mathrm{Ext}_{U(\fl_I)}^1(L^I(w), H^1(\fn_I, L(w)^{s_j})/L^I(s_jw)).\]
Hence we can forget about the quotient by $L^I(s_jw)$ in (\ref{lesI}). By the above results, in particular (\ref{soctonosoc}) and the fact $L^I(w)$ appears in $\soc_{U(\fl_I)}H^1(\fn_I, L(s_jw))$, we must have
\begin{equation}\label{equ: H1 radical bound}
\mathrm{rad}(M_I^1)\subseteq \mathrm{soc}_{U(\fl_I)}(H^1(\fn_I, L(w)^{s_j})),
\end{equation}
so that $\mathrm{rad}(M_I^1)$ is semi-simple and its constituents $L^I(x)$ satisfy $\mathrm{Ext}_{U(\fl_I)}^1(L^I(w),L^I(x))\neq 0$ with $x\ne w$.
Since $w$ is minimal in $W(L_I)w$, any $x\in W(L_I)w$ is such that $w\leq x$, and if moreover $w\prec x$ it follows from Lemma \ref{lem: dominance control} that $\ell(x)=\ell(w)+1$. Then by Lemma~\ref{lem: HC} and \ref{it: rabiotext 2} of Lemma \ref{rabiotext} we see that, for $x\in W(G)$ with $x\ne w$, $\mathrm{Ext}_{U(\fl_I)}^1(L^I(w),L^I(x))$ is non-zero if and only if $x\in W(L_I)w$ and $\ell(x)=\ell(w)+1$. Combining this with the equivalence between conditions (i) and (iii) in Lemma~\ref{lem: H1 conjugate}, we obtain from (\ref{equ: H1 radical bound})
\begin{equation}\label{equ: H1 radical sharp bound}
\mathrm{rad}(M_I^1)\subseteq \bigoplus_{\overset{x\in W(L_I)w}{\overset{\ell(x)=\ell(w)+1}{j\notin D_L(x)}}}\!\!\!\!\!\!L^I(x).
\end{equation}

Finally we prove that $L^I(w)$ is a constituent of $H^1(\fn_I, M^{s_j})$ and that (\ref{equ: H1 radical sharp bound}) is an isomorphism. It is enough to prove that for any $x\in W(L_I)w$ with $\ell(x)=\ell(w)+1$ and $j\notin D_L(x)$, $H^1(\fn_I, M^{s_j})$ contains a subquotient with socle $L^I(x)$ and cosocle $L^I(w)$. As $x\in W(L_I)w$ and $\ell(x)=\ell(w)+1$, there exists $j_1\in I$ such that $x=s_{j_1}w$. As $j\in D_L(w)\setminus D_L(x)$, we must have $j_1\in \{j-1,j+1\}$. As $H^0(\fn_I, M^{s_j})=0$ by (\ref{equ: H0 conj nonsplit}), we obtain the following isomorphism from the spectral sequence (\ref{equ: double n coh}) (applied with $I'\defeq \{j,j_1\}$ and thus $I\cap I'=\{j_1\}$)
\begin{equation}\label{equ: pass H1 conj to rank one}
H^0(\fl_I\cap \fn_{I'}, H^1(\fn_I, M^{s_j}))\cong H^1(\fn_{I\cap I'}, M^{s_j}).
\end{equation}
Switching the roles of $I$ and $I'$, another application of (\ref{equ: double n coh}) gives an injection
\begin{equation}\label{equ: pass H1 conj to rank two}
H^1(\fl_{I'}\cap \fn_I, H^0(\fn_{I'}, M^{s_j}))\hookrightarrow H^1(\fn_{I\cap I'}, M^{s_j}).
\end{equation}
It follows from Lemma~\ref{lem: H0 of Ext} (applied with $x=w$, $w=s_jw$ and $I=I'$ there) and $s_j\fn_{I'}{s_j}\cong \fn_{I'}$ that $H^0(\fn_{I'}, M^{s_j})\cong H^0(\fn_{I'}, M)^{s_j}$ is a length two $U(\fl_{I'})$-module with socle $L^{I'}(w)^{s_j}$ and cosocle $L^{I'}(s_jw)^{s_j}$ (and thus isomorphic to $M_{\{j_1\}}^{s_j}$ where $M_{\{j_1\}}\defeq U(\fl_{I'})\otimes_{U(\fl_{I'}\cap \fp_{\{j_1\}})}L^{\{j_1\}}(s_jw)$). Now by \ref{it: rank two H1} of Lemma~\ref{lem: nonsplit rank two} applied with $M$ there being the $U(\fl_{I'})$-module $H^0(\fn_{I'}, M)$ we have
\[H^1(\fl_{I'}\cap \fn_I, H^0(\fn_{I'}, M^{s_j}))\cong U(\fl_{\{j_1\}})\otimes_{U(\fb_{\{j_1\}})}w\cdot\mu_0,\] which is a $U(\fl_{\{j_1\}})$-module of length two with socle $L^{\{j_1\}}(s_{j_1}w)$ and cosocle $L^{\{j_1\}}(w)$. By (\ref{equ: pass H1 conj to rank two}) and (\ref{equ: pass H1 conj to rank one}) we deduce that $U(\fl_{\{j_1\}})\otimes_{U(\fb_{\{j_1\}})}w\cdot\mu_0$ \emph{embeds into} $H^0(\fl_I\cap \fn_{I'}, H^1(\fn_I, M^{s_j}))$. By (\ref{equ: g spectral seq 0}) (applied with $\fl_I$ instead of $\fg$ and $\fl_I\cap \fn_{I'}$ instead of $\fn_I$) this first forces the constituent $L^I(w)$ (the cosocle of $U(\fl_I)\otimes_{U(\fb_I)}w\cdot\mu_0\cong U(\fl_I)\otimes_{U(\fl_I\cap \fp_{I'})} (U(\fl_{\{j_1\}})\otimes_{U(\fb_{\{j_1\}})}w\cdot\mu_0)$) to show up in $H^1(\fn_I, M^{s_j})$. Then together with (\ref{equ: H1 radical sharp bound}) this also forces $H^1(\fn_I, M^{s_j})$ to have a length two subquotient with socle $L^I(s_{j_1}w)=L^I(x)$ and cosocle $L^I(w)$.
\end{proof}

\subsection{\texorpdfstring{$\mathrm{Ext}$}{Ext}-squares of \texorpdfstring{$U(\fg)$}{U(g)}-modules}\label{U(g)square}

We use all previous results of \S\ref{sec: n coh} to construct important finite length $U(\fg)$-modules which are uniserial (Lemma \ref{lem: simple length three g}) or ``squares'' (Proposition \ref{lem: explicit g square}, Lemma \ref{lem: Ext square type 3}).

\begin{defn}\label{def: Ext hypercube}
Let $w,x,w_1\in W(G)$, $I\defeq \Delta\setminus D_L(w)$ and $I'\defeq \Delta\setminus D_L(x)$.
\begin{enumerate}[label=(\roman*)]
\item A $U(\fg)$-module $Q_{w_1}(x,w)$ is an \emph{$\mathrm{Ext}$-hypercube} if the following properties hold
\begin{itemize}
\item $w_1\in W^{I,I'}$;
\item $Q_{w_1}(x,w)$ is semi-simple as $U(\ft)$-module;
\item $Q_{w_1}(x,w)$ is finite length multiplicity free and rigid as $U(\fg)$-module;
\item $Q_{w_1}(x,w)$ has socle $L(w)^{w_1}$ and cosocle $L(x)$.
\end{itemize}
\item An $\mathrm{Ext}$-hypercube is an \emph{$\mathrm{Ext}$-square} if it has Loewy length three, and an \emph{$\mathrm{Ext}$-cube} if it has Loewy length four.
\item \label{ext3} An $\mathrm{Ext}$-hypercube $Q_{w_1}(x,w)$ is \emph{minimal} if for any $U(\fg)$-submodules $M_4\subsetneq M_3\subseteq M_2\subsetneq M_1\subseteq Q_{w_1}(x,w)$, we have $\mathrm{Ext}_{U(\fg)}^1(M_1/M_2,M_3/M_4)\neq 0$ if and only if $M_2=M_3$, in which case $\mathrm{Ext}_{U(\fg)}^1(M_1/M_2,M_3/M_4)$ has dimension $1$ and the sequence $0\rightarrow M_3/M_4=M_2/M_4 \rightarrow M_1/M_4 \rightarrow M_1/M_2 \rightarrow 0$ is non-split.
\item When $w_1=1$, an $\mathrm{Ext}$-hypercube $Q_{1}(x,w)$ is \emph{minimal in $\cO_{\rm{alg}}^{\fb}$} if $Q_{1}(x,w)$ lies in $\cO_{\rm{alg}}^{\fb}$ and if \ref{ext3} holds with $\mathrm{Ext}_{\cO_{\rm{alg}}^{\fb}}^1$ instead of $\mathrm{Ext}_{U(\fg)}^1$.
\end{enumerate}
\end{defn}

If an $\mathrm{Ext}$-hypercube $Q_{w_1}(x,w)$ exists, then $Q_{w_1^{-1}}(w,x)$ exists for the triple $(w,x,w_1^{-1})$ by setting $Q_{w_1^{-1}}(w,x)\defeq (Q_{w_1}(x,w)^\tau)^{w_1^{-1}}$.

\begin{rem}
It is not true that $Q_{w_1}(x,w)$ is uniquely determined by the triple $(x,w,w_1)$. For example, if $\fg=\fg\fl_3$, $x=s_1s_2$, $w=s_2s_1$ and $w_1=1$, one can check that there are three different choices of $Q_{w_1}(x,w)$ with minimal possible length, with middle layer being respectively $L(s_1)\oplus L(s_2)$, $L(s_1)\oplus L(s_1s_2s_1)$ and $L(s_2)\oplus L(s_1s_2s_1)$.
\end{rem}

Recall that $\mathrm{Rad}^k$ and $\mathrm{Rad}_k$ for $k\geq 0$ are defined in \S\ref{generalnotation}.

\begin{lem}\label{lem: squareminimal}
An $\mathrm{Ext}$-square $Q_{w_1}(x,w)$ is minimal if and only if it satisfies the following conditions:
\begin{enumerate}[label=(\roman*)]
\item \label{it: minimal sq 1} $\mathrm{Ext}_{U(\fg)}^1(L(x),L(w)^{w_1})=0$;
\item \label{it: minimal sq 2} $\Dim_E \mathrm{Ext}_{U(\fg)}^1(C,L(w)^{w_1})=1=\Dim_E \mathrm{Ext}_{U(\fg)}^1(L(x),C)$ for any irreducible constituent $C$ of $\mathrm{Rad}_1(Q_{w_1}(x,w))$;
\item \label{it: minimal sq 3} there exists an irreducible constituent $C$ of $\mathrm{Rad}_1(Q_{w_1}(x,w))$ such that
\[\mathrm{Ext}_{U(\fg)}^1(L(x),\mathrm{Rad}^1(Q_{w_1}(x,w))^C)=0\]
where $\mathrm{Rad}^1(Q_{w_1}(x,w))^C\subset \mathrm{Rad}^1(Q_{w_1}(x,w))$ is the unique subobject not containing $C$ (in its cosocle).
\end{enumerate}
Moreover, if $w_1=1$, $Q_{1}(x,w)$ is minimal in $\cO_{\rm{alg}}^{\fb}$ if and only if the same conditions as above hold replacing everywhere $\mathrm{Ext}_{U(\fg)}^1$ by $\mathrm{Ext}_{\cO_{\rm{alg}}^{\fb}}^1$.
\end{lem}
\begin{proof}
Note that, contrary to what the terminology ``square'' may suggest, $\mathrm{Rad}_1(Q_{w_1}(x,w))$ can contain more than $2$ constituents. We prove the $U(\fg)$-module case, the proof for $\cO_{\rm{alg}}^{\fb}$ being the same. Note first that \ref{it: minimal sq 1} and \ref{it: minimal sq 2} are contained in \ref{ext3} of Definition \ref{def: Ext hypercube}. Let $C$ be an irreducible constituent of $\mathrm{Rad}_1(Q_{w_1}(x,w))$ and assume that \ref{it: minimal sq 1} and \ref{it: minimal sq 2} hold. Since $\Dim_E \mathrm{Ext}_{U(\fg)}^1(L(x),C)=1$ and $Q_{w_1}(x,w)$ is an $\mathrm{Ext}$-square, we have a short exact sequence
\begin{multline*}
0\longrightarrow \mathrm{Ext}_{U(\fg)}^1(L(x),\mathrm{Rad}^1(Q_{w_1}(x,w))^C)\longrightarrow \mathrm{Ext}_{U(\fg)}^1(L(x),\mathrm{Rad}^1(Q_{w_1}(x,w)))\\
\longrightarrow \mathrm{Ext}_{U(\fg)}^1(L(x),C)\longrightarrow 0.
\end{multline*}
So we see that \ref{it: minimal sq 3} holds if and only if $\Dim_E \mathrm{Ext}_{U(\fg)}^1(L(x),\mathrm{Rad}^1(Q_{w_1}(x,w)))=1$ if and only if \ref{it: minimal sq 3} holds for all constituents $C$ if and only if $\mathrm{Ext}_{U(\fg)}^1(L(x),M)=0$ for any $M\subsetneq \mathrm{Rad}^1(Q_{w_1}(x,w))$ if and only if there is no $\mathrm{Ext}$-square with socle $L(w)^{w_1}$, cosocle $L(x)$ and middle layer strictly contained in $\mathrm{Rad}_1(Q_{w_1}(x,w))$ if and only if \ref{ext3} of Definition \ref{def: Ext hypercube} holds.
\end{proof}

\begin{rem}\label{minimalbis}
The proof of Lemma \ref{lem: squareminimal} shows that, in the presence of \ref{it: minimal sq 1} and \ref{it: minimal sq 2} of Lemma \ref{lem: squareminimal}, condition \ref{it: minimal sq 3} of \emph{loc.\ cit.}\ is equivalent to: there exists an irreducible constituent $C$ of $\mathrm{Rad}_1(Q_{w_1}(x,w))$ such that $\mathrm{Ext}_{U(\fg)}^1((Q_{w_1}(x,w)/L(w)^{w_1})/C, L(w)^{w_1})=0$. Likewise with $\cO_{\rm{alg}}^{\fb}$ instead of $U(\fg)$-modules.
\end{rem}

In the rest of this section, we construct several minimal $\mathrm{Ext}$-squares. Our main tool to do that are \emph{wall-crossing functors}.\bigskip

For $\lambda,\mu\in \Lambda$ we first have an exact translation functor (see for instance \cite[\S 7.1]{Hum08})
\begin{equation}\label{equ: translation functor}
T^\mu_{\lambda}: \cO^{\fb}_{\rm{alg}}\rightarrow \cO^{\fb}_{\rm{alg}}: M\mapsto \mathrm{pr}_{\mu}(L\otimes_E \mathrm{pr}_{\lambda}(M))
\end{equation}
where $L$ is the unique finite dimensional $U(\fg)$-module with highest weight in the $W(G)$-orbit of $\mu-\lambda$ (for the naive action), and $\mathrm{pr}_{\mu}$, $\mathrm{pr}_{\lambda}$ is the projection onto the generalized eigenspace for the infinitesimal character associated with $L(\mu)$, $L(\lambda)$ respectively. Let $j\in \{1,\dots,n-1\}$ and $\mu\in \Lambda$ such that $\langle \mu + \rho, \alpha^\vee\rangle \geq 0$ for all $\alpha\in \Phi^+$ and the stabilizer of $\mu$ in $W(G)$ for the dot action is $\{1, s_j\}$. We define $\Theta_{s_j}\defeq T_\mu^{w_0\cdot \mu_0} \circ T_{w_0\cdot \mu_0}^\mu: \cO^{\fb}_{\rm{alg}}\rightarrow \cO^{\fb}_{\rm{alg}}$ which doesn't depend on the choice of $\mu$ as above (\cite[Example~10.8]{Hum08}) and is called a wall-crossing functor (\cite[\S 7.15]{Hum08}, the $w_0$ comes from the conventions of \emph{loc.\ cit.}\ which uses antidominant weights). For any $M$ in $\cO^{\fb}_{\rm{alg}}$ there are two canonical adjunction maps $\Theta_{s_j}(M)\rightarrow M$ and $M\rightarrow\Theta_{s_j}(M)$ which are non-zero as soon as both $M$ and $\Theta_{s_j}(M)$ are non-zero (\cite[Prop.~7.2(a)]{Hum08}).

\begin{prop}\label{prop: Jantzen middle}
Let $w\in W(G)$ and $j\in\Delta$. Then we have $\Theta_{w_0s_jw_0}(L(w))\neq 0$ if and only if $j\in D_R(w)$, in which case $\Theta_{w_0s_jw_0}(L(w))$ has Loewy length three with both socle and cosocle isomorphic to $L(w)$ and middle layer isomorphic to
\begin{equation}\label{equ: wall crossing middle}
L(ws_j)\oplus\oplus_{x \in S}L(x)^{\mu(x,w)}
\end{equation}
where $S\defeq \{x\mid w\prec x, j\notin D_R(x)\}$.
\end{prop}
\begin{proof}
It follows from \cite[Thm.~7.14(c),(f),(g)]{Hum08} together with \cite[Thm.~7.9]{Hum08} that $\Theta_{w_0s_jw_0}(L(w))\neq 0$ if and only if $j\in D_R(w)$ (i.e.~$ws_j<w$), in which case $L(w)$ is both the socle and cosocle of $\Theta_{w_0s_jw_0}(L(w))$, that each $L(x)\in\mathrm{JH}_{U(\fg)}(\Theta_{w_0s_jw_0}(L(w)))$ with $x\neq w$ satisfies $j\notin D_R(x)$ (i.e.~$x<xs_j$), and that if $j\notin D_R(x)$ we have an isomorphism
\begin{equation}\label{equ: wall crossing middle Ext}
\Hom_{\cO^{\fb}_{\rm{alg}}}(\mathrm{rad}(\Theta_{w_0s_jw_0}(L(w))),L(x))\cong \mathrm{Ext}_{\cO^{\fb}_{\rm{alg}}}^1(L(w),L(x)).
\end{equation}
By Vogan's Conjecture (which follows from the proof of the Kazhdan-Lusztig Conjecture) $\mathrm{rad}(\Theta_{w_0s_jw_0}(L(w)))/\mathrm{soc}(\Theta_{w_0s_jw_0}(L(w)))$ is semi-simple (cf.~\cite[\S 7.15, \S 8.10]{Hum08}), which together with \ref{it: rabiotext 1} of Lemma~\ref{rabiotext} and (\ref{equ: wall crossing middle Ext}) gives (\ref{equ: wall crossing middle}). Here we use the fact that the only $x$ satisfying $x\prec w$ and $j\in D_R(w)\setminus D_R(x)$ is $x=ws_j$ (see Lemma~\ref{lem: dominance control}).
\end{proof}

For $j_0,j_1\in \Delta$ we define
\begin{equation}\label{wj1j2}
w_{j_1,j_0}\defeq s_{j_1}s_{j_1-1}\cdots s_{j_0}\in W(G)\text{ if }j_1\geq j_0,\ \ w_{j_1,j_0}\defeq s_{j_1}s_{j_1+1}\cdots s_{j_0}\in W(G)\text{ if }j_1\leq j_0
\end{equation}
(with $w_{j_0,j_0}=s_{j_0}$). It is clear that $D_L(w_{j_1,j_0})=\{j_1\}$ and $D_R(w_{j_1,j_0})=\{j_0\}$, and one can check that $w_{j_1,j_0}$ is the unique partial-Coxeter element satisfying these two properties.

\begin{rem}\label{rem: explicit middle}
Let $j_0,j_1\in \Delta$, $w=w_{j_1,j_0}$ and $S$ as in Proposition \ref{prop: Jantzen middle}, we deduce from \ref{it: special descent 4} of Lemma~\ref{lem: special descent change} that
\[S=\{w_{j_1,j_0'}\mid j_0'\in\Delta,\ |j_0'-j_0|=1,\ w_{j_1,j_0'}>w_{j_1,j_0}\}\]
and that $\mu(w,x)=1$ for $x\in S$. More precisely, $w_{j_1,j_0+1}\in S$ if and only if $j_1\leq j_0<n-1$, and $w_{j_1,j_0-1}\in S$ if and only if $j_1\geq j_0>1$.
\end{rem}

\begin{lem}\label{lem: g square as wall crossing}
Let $j_0,j_1,j_0',j_1'\in\Delta$ with $|j_0-j_0'|=1$ and $|j_1-j_1'|=1$. Let $M_0$ be the unique length $2$ object in $\cO^{\fb}_{\rm{alg}}$ with socle $L(w_{j_1',j_0})$ and cosocle $L(w_{j_1,j_0})$ (see \ref{it: rabiotext 1} of Lemma~\ref{rabiotext} and Lemma~\ref{lem: special descent change}). Then $L(w_{j_1',j_0'})$ occurs with multiplicity one in $\Theta_{w_0s_{j_0}w_0}(M_0)$, and the unique quotient of $\Theta_{w_0s_{j_0}w_0}(M_0)$ with socle $L(w_{j_1',j_0'})$ is an $\mathrm{Ext}$-square $Q_1(w_{j_1,j_0},w_{j_1',j_0'})$ in $\cO^{\fb}_{\rm{alg}}$ with socle $L(w_{j_1',j_0'})$, cosocle $L(w_{j_1,j_0})$ and middle layer contained in $L(w_{j_1,j_0'})\oplus L(w_{j_1',j_0})\oplus L(1)$ if $j_0=j_1$ and $j_0'=j_1'$, and contained in $L(w_{j_1,j_0'})\oplus L(w_{j_1',j_0})$ otherwise.
\end{lem}
\begin{proof}
We write $S\defeq \{x\mid w_{j_1,j_0}\prec x, \ j_0\notin D_R(x)\}$ and $S'\defeq \{x\mid w_{j_1',j_0}\prec x, \ j_0\notin D_R(x)\}$. By Lemma~\ref{lem: dominance control} we have $\ell(x)=\ell(w_{j_1,j_0})+1$ and thus $\mu(w_{j_1,j_0},x)=1$ (resp.~$\ell(x)=\ell(w_{j_1',j_0})+1$ and thus $\mu(w_{j_1',j_0},x)=1$) for each $x\in S$ (resp.~for each $x\in S'$). We write $L_S\defeq \bigoplus_{x\in S}L(x)$ and similarly for $L_{S'}$.
By Proposition~\ref{prop: Jantzen middle} and Lemma~\ref{lem: coxeter prec} we know that $\Theta_{w_0s_{j_0}w_0}(L(w_{j_1,j_0}))$ (resp.~$\Theta_{w_0s_{j_0}w_0}(L(w_{j_1',j_0}))$) has Loewy length three with both socle and cosocle $L(w_{j_1,j_0})$ (resp.~$L(w_{j_1',j_0})$) and with middle layer $L(w_{j_1,j_0}s_{j_0})\oplus L_S$ (resp.~$L(w_{j_1',j_0}s_{j_0})\oplus L_{S'}$). Recall that for any non-zero $M$ in $\cO^{\fb}_{\rm{alg}}$ we have adjunction maps $\Theta_{w_0s_{j_0}w_0}(M)\rightarrow M$, $M\rightarrow\Theta_{w_0s_{j_0}w_0}(M)$ which are non-zero if $\Theta_{w_0s_{j_0}w_0}(M)\ne 0$. Together with the previous discussion and the exactness of $\Theta_{w_0s_{j_0}w_0}$, we easily deduce that $\Theta_{w_0s_{j_0}w_0}(M_0)\rightarrow M_0$ is surjective and $M_0\rightarrow\Theta_{w_0s_{j_0}w_0}(M_0)$ is injective. As $w_{j_1',j_0'}\in S'\setminus S$ (see Remark~\ref{rem: coxeter pair}) and $w_{j_1',j_0'}\neq w_{j_1,j_0},w_{j_1',j_0}$, we see that $L(w_{j_1',j_0'})$ appears with multiplicity one in $\Theta_{w_0s_{j_0}w_0}(L(w_{j_1',j_0}))$ and $\Theta_{w_0s_{j_0}w_0}(M_0)$. We define $M$ (resp.~$Q$) as the unique quotient of $\Theta_{w_0s_{j_0}w_0}(L(w_{j_1',j_0}))$ (resp.~$\Theta_{w_0s_{j_0}w_0}(M_0)$) with socle $L(w_{j_1',j_0'})$ and cosocle $L(w_{j_1',j_0})$. Then $M$ has length $2$ and the (non-zero) composition
\[\Theta_{w_0s_{j_0}w_0}(L(w_{j_1',j_0}))\hookrightarrow \Theta_{w_0s_{j_0}w_0}(M_0)\twoheadrightarrow Q\]
must have image $M$. In particular $Q/M$ is a quotient of $\Theta_{w_0s_{j_0}w_0}(L(w_{j_1,j_0}))$. Since $L(w_{j_1',j_0'})$ does not occur in $M_0$ the composition
\[M_0\hookrightarrow\Theta_{w_0s_{j_0}w_0}(M_0)\twoheadrightarrow Q\]
must be zero, and thus $Q$ is a quotient of $\Theta_{w_0s_{j_0}w_0}(M_0)/M_0$. As $L(w_{j_1',j_0})=\mathrm{soc}(M_0)$ does not occur in $\Theta_{w_0s_{j_0}w_0}(L(w_{j_1,j_0}))$ and $L(w_{j_1,j_0})=\mathrm{cosoc}(M_0)$ does not occur in $\Theta_{w_0s_{j_0}w_0}(L(w_{j'_1,j_0}))$ (using the above description), the composition
\[M_0\hookrightarrow \Theta_{w_0s_{j_0}w_0}(M_0)\twoheadrightarrow \Theta_{w_0s_{j_0}w_0}(L(w_{j_1,j_0}))\]
factors through $M_0\twoheadrightarrow L(w_{j_1,j_0})\hookrightarrow \Theta_{w_0s_{j_0}w_0}(L(w_{j_1,j_0}))$. It follows that $Q/M$ is a quotient of $\Theta_{w_0s_{j_0}w_0}(L(w_{j_1,j_0}))/L(w_{j_1,j_0})$. Consequently, $Q$ is an $\mathrm{Ext}$-square (see Definition~\ref{def: Ext hypercube}) with socle $L(w_{j_1',j_0'})$, cosocle $L(w_{j_1,j_0})$ and middle layer
\begin{equation}\label{equ: middle bound g crossing}
\mathrm{rad}_1(Q)\subseteq L(w_{j_1',j_0})\oplus L(w_{j_1,j_0}s_{j_0})\oplus\bigoplus_{x\in S''}L(x)
\end{equation}
where $S''=\{x\in S\mid \mathrm{Ext}_{U(\fg)}^1(L(x), L(w_{j_1',j_0'}))\neq 0\}$.

By \ref{it: special descent 4} of Lemma~\ref{lem: special descent change} we know that $x\in S$ if and only if $x=w_{j_1, j_0-1}$ with $j_1\geq j_0$ or $x=w_{j_1, j_0+1}$ with $j_1\leq j_0$. In both cases $D_L(x)=\{j_1\}$ (and $j_1\ne j'_1$), hence by \ref{it: special descent 3} of Lemma~\ref{lem: special descent change} and (\ref{equ: dual Ext1}) we deduce that $x\in S''$ if and only if $x=w_{j_1,j'_0}$ (we already know $j_1=j'_1\pm 1$). Since $w_{j_1,j_0}s_{j_0}$ is obviously not in $S$, we also see that the right hand side of (\ref{equ: middle bound g crossing}) is multiplicity free.

Assume that $L(w_{j_1,j_0}s_{j_0})$ occurs in $\mathrm{rad}_1(Q)$. Then $\mathrm{Ext}_{U(\fg)}^1(L(w_{j_1,j_0}s_{j_0}), L(w_{j_1',j_0'}))\neq 0$. Assume first $j_0\ne j_1$ so that $D_L(w_{j_1,j_0}s_{j_0})=\{j_1\}$, then by \ref{it: special descent 3} of Lemma~\ref{lem: special descent change} again (and (\ref{equ: dual Ext1})) we have $w_{j_1,j_0}s_{j_0}=w_{j_1,j_0'}$. Assume now $j_0= j_1$, then $w_{j_1,j_0}s_{j_0}=1\ne w_{j_1',j_0'}$ and \ref{it: rabiotext 2} of Lemma \ref{rabiotext} implies $1\prec w_{j_1',j_0'}$ which forces $j'_1=j'_0$ by Lemma \ref{lem: dominance control}.

Summing up, we have shown that $\mathrm{rad}_1(Q)$ is contained in $L(w_{j_1,j_0'})\oplus L(w_{j_1',j_0})\oplus L(1)$ if $j_0=j_1$ and $j_0'=j_1'$, and is contained in $L(w_{j_1,j_0'})\oplus L(w_{j_1',j_0})$ otherwise.
\end{proof}

\begin{lem}\label{lem: minimal g square}
Let $j_0,j_1,j_0',j_1'\in\Delta$ with $|j_0-j_0'|=1$ and $|j_1-j_1'|=1$. Let $Q_1(w_{j_1,j_0},w_{j_1',j_0'})$ be an $\mathrm{Ext}$-square with socle $L(w_{j_1',j_0'})$, cosocle $L(w_{j_1,j_0})$ and middle layer contained in $L(w_{j_1,j_0'})\oplus L(w_{j_1',j_0})\oplus L(1)$ if $j_0=j_1$ and $j_0'=j_1'$, and contained in $L(w_{j_1,j_0'})\oplus L(w_{j_1',j_0})$ otherwise. Then $Q_1(w_{j_1,j_0},w_{j_1',j_0'})$ is minimal, and in particular unique up to isomorphism, and is in $\cO^{\fb}_{\rm{alg}}$.
\end{lem}
\begin{proof}
We write $Q\defeq Q_1(w_{j_1,j_0},w_{j_1',j_0'})$ for short and note that the very last statement follows from unicity and Lemma \ref{lem: g square as wall crossing}. Moreover unicity follows from \ref{it: minimal sq 3} and \ref{it: minimal sq 2} of Lemma \ref{lem: squareminimal}, hence we only need to prove minimality.

By \ref{it: special descent 3} of Lemma~\ref{lem: special descent change}, we see that \ref{it: minimal sq 1} and \ref{it: minimal sq 2} of Lemma \ref{lem: squareminimal} hold. Hence it suffices to prove \ref{it: minimal sq 3} of Lemma \ref{lem: squareminimal}. More precisely it suffices to show $\mathrm{Ext}_{U(\fg)}^1(L(w_{j_1,j_0}),M)=0$ where $M\subseteq Q$ is the maximal $U(\fg)$-submodule such that $L(w_{j_1,j_0}),L(w_{j_1',j_0})\notin\mathrm{JH}_{U(\fg)}(M)$. Note that the vanishing $\mathrm{Ext}_{U(\fg)}^1(L(w_{j_1,j_0}),M)=0$ actually forces $L(w_{j_1',j_0})\in\mathrm{JH}_{U(\fg)}(\mathrm{rad}_1(Q))$.

Assume on the contrary $\mathrm{Ext}_{U(\fg)}^1(L(w_{j_1,j_0}),M)\neq 0$, then an arbitrary object $M^+$ that fits into a non-split extension $0\rightarrow M\rightarrow M^+\rightarrow L(w_{j_1,j_0})\rightarrow 0$ contains a unique submodule $Q^-$ with cosocle $L(w_{j_1,j_0})$. Since $L(w_{j_1,j_0}),L(w_{j_1',j_0})\notin\mathrm{JH}_{U(\fg)}(M)$, any constituent $L(x)$ of $Q^-$ distinct from its cosocle is such that $j_0\notin D_R(x)$, hence is killed by $\Theta_{w_0s_{j_0}w_0}$ by the first statement of Proposition~\ref{prop: Jantzen middle}. Thus the surjection $Q^-\twoheadrightarrow L(w_{j_1,j_0})$ induces an isomorphism
\begin{equation}\label{equ: minimal g square isom}
\Theta_{w_0s_{j_0}w_0}(Q^-)\buildrel\sim\over\longrightarrow \Theta_{w_0s_{j_0}w_0}(L(w_{j_1,j_0}))
\end{equation}
and in particular $\Theta_{w_0s_{j_0}w_0}(Q^-)$ has cosocle $L(w_{j_1,j_0})$ by the second statement of Proposition~\ref{prop: Jantzen middle}. Since $Q^-$ is multiplicity free (as $Q$ is) with cosocle $L(w_{j_1,j_0})$, the canonical (non-zero) adjunction map $\Theta_{w_0s_{j_0}w_0}(Q^-)\rightarrow Q^-$ must be a surjection. Note that we have $L(w_{j_1',j_0'})\notin\mathrm{JH}_{U(\fg)}(\Theta_{w_0s_{j_0}w_0}(L(w_{j_1,j_0})))$ by Proposition~\ref{prop: Jantzen middle} and Remark~\ref{rem: explicit middle}, and thus $L(w_{j_1',j_0'})\notin\mathrm{JH}_{U(\fg)}(\Theta_{w_0s_{j_0}w_0}(Q^-))$ by (\ref{equ: minimal g square isom}). However $L(w_{j_1',j_0'})$ is the socle of $M$ and hence occurs in $Q^-$. This contradicts the surjection $\Theta_{w_0s_{j_0}w_0}(Q^-)\twoheadrightarrow Q^-$.
\end{proof}

\begin{prop}\label{lem: explicit g square}
Let $j_0,j_1,j_0',j_1'\in\Delta$ with $|j_0-j_0'|=1$ and $|j_1-j_1'|=1$. Then there exists a unique minimal $\mathrm{Ext}$-square $Q_1(w_{j_1,j_0},w_{j_1',j_0'})$ with socle $L(w_{j_1',j_0'})$, cosocle $L(w_{j_1,j_0})$ and middle layer $L(w_{j_1,j_0'})\oplus L(w_{j_1',j_0})\oplus L(1)$ if $j_0=j_1$ and $j_0'=j_1'$, and $L(w_{j_1,j_0'})\oplus L(w_{j_1',j_0})$ otherwise. Moreover $Q_1(w_{j_1,j_0},w_{j_1',j_0'})$ is in $\cO^{\fb}_{\rm{alg}}$.
\end{prop}
\begin{proof}
By Lemma~\ref{lem: minimal g square} we know that $Q$ is minimal, unique up to isomorphism and is in $\cO^{\fb}_{\rm{alg}}$. Switching $w_{j_1,j_0}$ and $w_{j_1',j_0'}$, we have similar statements for $Q_1(w_{j_1',j_0'},w_{j_1,j_0})$, which implies $Q^{\tau}\cong Q_1(w_{j_1',j_0'},w_{j_1,j_0})$ by unicity.
Moreover, the proof of Lemma~\ref{lem: minimal g square} implies $L(w_{j_1',j_0})\in\mathrm{JH}_{U(\fg)}(\mathrm{rad}_1(Q))$, and similarly $L(w_{j_1,j_0'})\in\mathrm{JH}_{U(\fg)}(\mathrm{rad}_1(Q^{\tau}))$. Since $\mathrm{JH}_{U(\fg)}(\mathrm{rad}_1(Q^{\tau}))=\mathrm{JH}_{U(\fg)}(\mathrm{rad}_1(Q))$, we deduce $L(w_{j_1,j_0'})\oplus L(w_{j_1',j_0})\subseteq \mathrm{rad}_1(Q)$. It remains to prove $L(1)\subseteq \mathrm{rad}_1(Q)$ when $j_0=j_1$ and $j_0'=j_1'$. Assume on the contrary that $L(1)\notin\mathrm{JH}_{U(\fg)}(\mathrm{rad}_1(Q))$ and recall (from our convention) that if a weight in $w\cdot \mu_0-\Z_{\geq 0}\Phi^+$ occurs in $L(w')=L(w'\cdot\mu_0)$ then $w'\leq w$. It follows that none of $L(s_{j_0'})$, $L(w_{j_0,j_0'})$, $L(w_{j_0',j_0})$ contain weights in $s_{j_0}\cdot\mu_0-\Z_{>0}\Phi^+$ and hence that $Q_{s_{j_0}\cdot\mu_0}$ is one dimensional and $Q_{s_{j_0}\cdot\mu_0-\nu}$ is $0$ for $\nu\in \Z_{>0}\Phi^+$. Since the action of a non-zero element of $\fu$ modifies a weight by a character in $-\Z_{>0}\Phi^+$, this forces $0\neq Q_{s_{j_0}\cdot\mu_0}\subseteq H^0(\fu,Q)_{s_{j_0}\cdot\mu_0}$, which together with (\ref{equ: g spectral seq 0}) (for $I=\emptyset$) gives a non-zero map $M(s_{j_0})\rightarrow Q$. As $Q$ is multiplicity free with cosocle $L(s_{j_0})$, it has to be a surjection. But this contradicts the fact that $L(s_{j_0'})\in\mathrm{JH}_{U(\fg)}(Q)\setminus \mathrm{JH}_{U(\fg)}(M(s_{j_0}))$. Hence, we must have $\mathrm{rad}_1(Q)\cong L(w_{j_1,j_0'})\oplus L(w_{j_1',j_0})\oplus L(1)$ when $j_0=j_1$ and $j_0'=j_1'$.
\end{proof}

\begin{rem}\label{rem: other simple square}
By a similar argument as in Lemma~\ref{lem: g square as wall crossing}, Lemma~\ref{lem: minimal g square} and Lemma~\ref{lem: explicit g square}, one can prove that, for $j_0,j_1\in\Delta$ with $|j_0-j_1|=1$, there exists a unique minimal $\mathrm{Ext}$-square $Q_{1}(1,w_{j_1,j_0})$ (resp.~$Q_{1}(w_{j_1,j_0},1)$) with middle layer $L(s_{j_0})\oplus L(s_{j_1})$, and that $Q_{1}(1,w_{j_1,j_0})$, $Q_{1}(w_{j_1,j_0},1)$ are in $\cO^{\fb}_{\rm{alg}}$.
\end{rem}

For $S\subseteq W(G)$, we write $L_S\defeq \bigoplus_{x\in S}L(x)$.

\begin{lem}\label{lem: Ext square type 3}
Let $j\in\Delta$, $I=\Delta\setminus\{j\}$ and $w\in W(G)$ such that $D_L(w)=\{j\}$. Assume $S_0\defeq \{x'\mid x'\in W(L_I)w, \ell(x')=\ell(w)+1, j\notin D_L(x')\}\ne \emptyset$. Then there exists a unique minimal $\mathrm{Ext}$-square $Q_{s_j}(w,w)$ such that
\[\mathrm{Rad}_1(Q_{s_j}(w,w))\cong L(s_jw)\oplus L_{S_0}.\]
\end{lem}
\begin{proof}
There exists a unique $U(\fg)$-module $M_0$ with socle $L(w)$ and cosocle $L(s_jw)\oplus L_{S_0}$. Since $L(x')^{s_j}\cong L(x')$ if $j\notin D_L(x')$ (\ref{it: H0 conjugate 3} of Lemma \ref{lem: H0 Weyl conjugate}), $M_0^{s_j}$ has socle $L(w)^{s_j}$ and cosocle $L(s_jw)\oplus L_{S_0}$. By \ref{it: rabiotext 2} of Lemma \ref{rabiotext} we have
\begin{equation}\label{equ: square conj middle dim}
\Dim_E \mathrm{Ext}_{U(\fg)}^1(L(w), M_0^{s_j}/L(w)^{s_j})=\Dim_E \mathrm{Ext}_{U(\fg)}^1(L(w), L(s_jw)\oplus L_{S_0})=1+\#S_0.
\end{equation}
By the equivalence between conditions (ii) and (iii) in Lemma~\ref{lem: H1 conjugate} (applied with $I=I'=\Delta\setminus\{j\}$ and $w_1=s_j$) we have $\mathrm{Ext}_{U(\fg)}^1(L(x),L(w)^{s_j})=0$ for $x\in W(G)$ such that $D_L(x)=\{j\}$, and in particular for each irreducible constituent $L(x)$ of $M^I(w)$.
This together with Proposition~\ref{prop: typical Ext2} and a d\'evissage using $0\rightarrow N^I(w)\rightarrow M^I(w)\rightarrow L(w)\rightarrow 0$ implies
\begin{equation}\label{equ: square conj Ext2 dim}
\Dim_E \mathrm{Ext}_{U(\fg)}^2(L(w), L(w)^{s_j})\leq\#S_0.
\end{equation}
The short exact sequence $0\rightarrow L(w)^{s_j}\rightarrow M_0^{s_j}\rightarrow M_0^{s_j}/L(w)^{s_j}\rightarrow 0$ yields the exact sequence
\[0\rightarrow \mathrm{Ext}_{U(\fg)}^1(L(w), M_0^{s_j})
\rightarrow \mathrm{Ext}_{U(\fg)}^1(L(w), M_0^{s_j}/L(w)^{s_j})
\rightarrow \mathrm{Ext}_{U(\fg)}^2(L(w), L(w)^{s_j})\]
which together with (\ref{equ: square conj middle dim}) and (\ref{equ: square conj Ext2 dim}) implies
\begin{equation*}
\Dim_E \mathrm{Ext}_{U(\fg)}^1(L(w), M_0^{s_j})\geq 1.
\end{equation*}
In particular, there exists an $\mathrm{Ext}$-square $Q_{s_j}(w,w)$ with socle $L(w)^{s_j}$, cosocle $L(w)$ and middle layer
\[\mathrm{Rad}_1(Q_{s_j}(w,w))\subseteq M_0^{s_j}/L(w)^{s_j}=L(s_jw)\oplus L_{S_0}.\]
We now prove $L(s_jw)\oplus L_{S_0}\subseteq \mathrm{Rad}_1(Q_{s_j}(w,w))$ and the minimality of $Q_{s_j}(w,w)$.\bigskip

\textbf{Step $1$}: Let $M$ be any $U(\fg)$-module with socle $L(w)$ and $M'$ a $U(\fg)$-module of finite length with all irreducible constituents in $\cO_{\rm{alg}}^{\fb}$. As $L(w)^{s_j}$ is not in $\cO_{\rm{alg}}^{\fb}$ (by \ref{it: H0 conjugate 3} of Lemma \ref{lem: H0 Weyl conjugate} and Lemma \ref{lem: dominance and left set}), we have $\Hom_{U(\fg)}(M', M^{s_j})=0$.\bigskip

\textbf{Step $2$}: We show that $Q_{s_j}(w,w)$ is minimal and is unique such that $\mathrm{Rad}_1(Q_{s_j}(w,w))\subseteq L(s_jw)\oplus L_{S_0}$.\\
Using again \ref{it: rabiotext 2} of Lemma \ref{rabiotext}, there exists a unique $U(\fg)$-module $M_1$ with cosocle $L(w)$ and socle $L(s_jw)\oplus L_{S_0}$. By unicity of $M_1$ there is a surjection $M_1\twoheadrightarrow Q_{s_j}(w,w)/L(w)^{s_j}$. Moreover it follows from \ref{it: rabiotext 1} of Lemma \ref{rabiotext} and (\ref{equ: O Hom radical}) (applied with $x=w$ and $w=x'\in S_0$) that $M_1/L(s_jw)$ is a quotient of $M(w)$. As $\Hom_{U(\fg)}(\mathrm{ker}(M(w)\twoheadrightarrow M_1/L(s_jw)),L(w)^{s_j})=0$ since $L(w)^{s_j}$ is not in $\cO_{\rm{alg}}^{\fb}$, the surjection $M(w)\twoheadrightarrow M_1/L(s_jw)$ induces an embedding $\mathrm{Ext}_{U(\fg)}^1(M_1/L(s_jw), L(w)^{s_j})\hookrightarrow \mathrm{Ext}_{U(\fg)}^1(M(w), L(w)^{s_j})$, which using \ref{rem: H1 conj Ext1 vanishing} of Remark \ref{mergerem} (applied with $I=\emptyset$, $x=w$ and $w_1=s_j$) implies
\begin{equation}\label{zerozero}
\mathrm{Ext}_{U(\fg)}^1(M_1/L(s_jw), L(w)^{s_j})~=~0.
\end{equation}
Then from $0\rightarrow L(s_jw)\rightarrow M_1\rightarrow M_1/L(s_jw)\rightarrow 0$ we deduce
\begin{equation}\label{embed1}
\mathrm{Ext}_{U(\fg)}^1(M_1, L(w)^{s_j})\hookrightarrow \mathrm{Ext}_{U(\fg)}^1(L(s_jw), L(w)^{s_j}).
\end{equation}
But $\Dim_E\mathrm{Ext}_{U(\fg)}^1(L(s_jw), L(w)^{s_j})=1$ by condition (iii) in Lemma~\ref{lem: H1 conjugate} with (\ref{chiant}) (applied with $I=I'=\Delta\setminus D_L(w)$ and $w_1=s_j$) and \ref{it: rabiotext 2} of Lemma \ref{rabiotext}. Since we have $\mathrm{Ext}_{U(\fg)}^1(Q_{s_j}(w,w)/L(w)^{s_j}, L(w)^{s_j})\ne 0$ (by definition of $Q_{s_j}(w,w)$) and
\begin{equation}\label{embed2}
\mathrm{Ext}_{U(\fg)}^1(Q_{s_j}(w,w)/L(w)^{s_j}, L(w)^{s_j})\hookrightarrow \mathrm{Ext}_{U(\fg)}^1(M_1, L(w)^{s_j})
\end{equation}
(by Step $1$ applied with $M'=\mathrm{ker}(M_1\twoheadrightarrow Q_{s_j}(w,w)/L(w)^{s_j})$ and $M=L(w)$), we deduce from (\ref{embed1}) and (\ref{embed2}) isomorphisms of $1$-dimensional vector spaces
\begin{equation*}
\begin{aligned}
\mathrm{Ext}_{U(\fg)}^1(M_1, L(w)^{s_j})&\buildrel\sim\over\longrightarrow \mathrm{Ext}_{U(\fg)}^1(L(s_jw), L(w)^{s_j})\\
\mathrm{Ext}_{U(\fg)}^1(Q_{s_j}(w,w)/L(w)^{s_j}, L(w)^{s_j})&\buildrel\sim\over\longrightarrow \mathrm{Ext}_{U(\fg)}^1(M_1, L(w)^{s_j}).
\end{aligned}
\end{equation*}
The composition gives an isomorphism of $1$-dimensional $E$-vector spaces
\[\mathrm{Ext}_{U(\fg)}^1(Q_{s_j}(w,w)/L(w)^{s_j}, L(w)^{s_j})\buildrel\sim\over\longrightarrow\mathrm{Ext}_{U(\fg)}^1(L(s_jw), L(w)^{s_j})\]
which implies that $L(s_jw)$ appears in $\mathrm{Rad}_1(Q_{s_j}(w,w))$. We now prove minimality. Using again Step $1$, we deduce from (\ref{zerozero})
\[\mathrm{Ext}_{U(\fg)}^1((Q_{s_j}(w,w)/L(w)^{s_j})/L(s_jw), L(w)^{s_j})=0.\]
By Lemma \ref{lem: squareminimal} and Remark \ref{minimalbis} this implies that $Q_{s_j}(w,w)$ is minimal, and is actually the unique $\mathrm{Ext}$-square with socle $L(w)^{s_j}$, cosocle $L(w)$ and middle layer $\mathrm{Rad}_1(Q_{s_j}(w,w))\subseteq L(s_jw)\oplus L_{S_0}$.\bigskip

\textbf{Step $3$}: We finally show that $L_{S_0}\subseteq \mathrm{Rad}_1(Q_{s_j}(w,w))$.\\
Let $M_1'$ be the unique $U(\fg)$-module with socle $L(w)$ and cosocle $L(s_jw)\oplus L_{S_0}$, which is in $\cO_{\rm{alg}}^{\fb}$ (use once more \ref{it: rabiotext 2} of Lemma \ref{rabiotext}). Then $Q_{s_j}(w,w)$ contains a unique maximal $U(\fg)$-submodule $Q_1$ with socle $L(w)^{s_j}$ and cosocle contained in $L(s_jw)\oplus L_{S_0}$ and containing $L(s_jw)$. We have $Q_1\subseteq (M_1')^{s_j}$, $Q_1/L(w)^{s_j}\in \cO_{\rm{alg}}^{\fb}$ and $\mathrm{Ext}_{U(\fg)}^1(L(w), Q_1)\ne 0$. In order to prove the statement, it suffices to show $Q_1=(M_1')^{s_j}$. Assume on the contrary $Q_1\subsetneq (M_1')^{s_j}$ and let $S_0'\subsetneq S_0$ such that $Q_1$ has cosocle $L(s_jw)\oplus L_{S_0'}$ (and socle $L(w)^{s_j}$). The surjection $M^I(w)\twoheadrightarrow L(w)$ together with $\Hom_{U(\fg)}(\mathrm{ker}(M^I(w)\rightarrow L(w)),Q_1)=0$ gives an embedding $\mathrm{Ext}_{U(\fg)}^1(L(w), Q_1)\hookrightarrow \mathrm{Ext}_{U(\fg)}^1(M^I(w), Q_1)$, which forces
\begin{equation}\label{equ: sq conj Verma}
\mathrm{Ext}_{U(\fg)}^1(M^I(w), Q_1)\neq 0.
\end{equation}
By \ref{it: H0 conjugate 2} of Lemma \ref{lem: H0 Weyl conjugate} (applied with $I=I'$ and $w=w_2=s_j$) we have $H^0(\fn_I,L(w)^{s_j})=0$, from which we deduce by d\'evissage (and \ref{it: dominance 2} of Lemma \ref{lem: n coh dominance}) that $H^0(\fn_I,Q_1)$ is in $\cO_{\rm{alg}}^{\fp_{I}}$. Since the socle $L(w)^{s_j}$ of $Q_1$ is not in $\cO_{\rm{alg}}^{\fb}$, (\ref{equ: g spectral seq 0}) implies $H^0(\fn_I,Q_1)=0$. Hence (\ref{equ: 5terms}) (applied with $M=Q_1$ and $M_I=L^I(w)$) and (\ref{equ: sq conj Verma}) give
\begin{equation}\label{equ: sq conj 0 1}
\Hom_{U(\fl_I)}(L^I(w), H^1(\fn_I,Q_1)\neq 0.
\end{equation}
Let $M_2$ be the unique $U(\fg)$-module with socle $L(w)$ and cosocle $L(s_jw)$. The short exact sequence $0\rightarrow M_2^{s_j}\rightarrow Q_1\rightarrow Q_1/M_2^{s_j}\rightarrow 0$ and $H^0(\fn_I,Q_1)=0$ give an exact sequence
\begin{equation}\label{longexactI}
0\rightarrow H^0(\fn_I, Q_1/M_2^{s_j})\rightarrow H^1(\fn_I,M_2^{s_j})\rightarrow H^1(\fn_I,Q_1)\rightarrow H^1(\fn_I, Q_1/M_2^{s_j}).
\end{equation}
Since $Q_1/M_2^{s_j}\cong L_{S_0'}^{s_j}\cong L_{S_0'}$, we have $H^0(\fn_I, Q_1/M_2^{s_j})\cong L^I_{S_0'}$ by \ref{it: dominance 2} of Lemma \ref{lem: n coh dominance}. Let $\xi: Z(\fl_I) \rightarrow E$ such that $L^I(w)_\xi\neq 0$. Since $x'\in W(L_I)w$, we have $L^I(x')_\xi\neq 0$, and thus $H^1(\fn_I, L(x'))_\xi=0$ by \ref{it: dominance 3} of Lemma \ref{lem: n coh dominance}. In particular $H^1(\fn_I, Q_1/M_2^{s_j})\cong \oplus_{x'\in S_0'} H^1(\fn_I, L(x'))$ does not have $L^I(w)$ as Jordan-H\"older factor. By Lemma~\ref{lem: nonsplit H1 conj} $L^I(w)$ has multiplicity one in $H^1(\fn_I,M_2^{s_j})$ and $H^1(\fn_I,M_2^{s_j})$ contains a $U(\fl_I)$-submodule $M_I$ with socle $L^I_{S_0}$ and cosocle $L^I(w)$. From (\ref{longexactI}) we deduce that $L^I(w)$ appears with multiplicity one in $H^1(\fn_I,Q_1)$ and that $H^1(\fn_I,Q_1)$ contains a submodule with socle $L^I_{S_0\setminus S_0'}$ and cosocle $L^I(w)$ (the image of $M_I$ in $H^1(\fn_I,Q_1)$). Since by assumption $S_0\setminus S_0'\ne \emptyset$, we must have $\Hom_{U(\fl_I)}(L^I(w), H^1(\fn_I,Q_1))=0$, which contradicts (\ref{equ: sq conj 0 1}). It follows that $S_0'=S_0$ and $Q_1=(M_1')^{s_j}$, which finishes the proof.
\end{proof}

\begin{rem}\label{rem: Ext 1 conjugate vanshing}
Keep the notation of Lemma \ref{lem: Ext square type 3} and let $M\subsetneq M_1'$ where $M_1'$ is the unique $U(\fg)$-module with socle $L(w)$ and cosocle $L(s_jw)\oplus L_{S_0}$, and assume that $L(s_jw)$ appears in (the cosocle of) $M$. Then we proved in Step $3$ of the above proof that $H^0(\fn_I,M^{s_j})=0$ and $\Hom_{U(\fl_I)}(L^I(w), H^1(\fn_I,M^{s_j}))=0$. By (\ref{equ: 5terms}) (applied with $M$ there being $M^{s_j}$) it follows that we have for $\ell\leq 1$:
\[\mathrm{Ext}_{U(\fg)}^{\ell}(M^I(w), M^{s_j})=0.\]
\end{rem}

\begin{ex}
The following special cases of Lemma~\ref{lem: Ext square type 3} will be useful. Let $j,j'\in\Delta$ and $w\defeq w_{j,j'}$ (so $D_L(w)=\{j\}$). If $j>j'$, then we have $S_0=\{w_{j+1,j'}\}$ when $j<n-1$ and $S_0=\emptyset$ when $j=n-1$. If $j<j'$, then we have $S_0=\{w_{j-1,j'}\}$ when $j>1$ and $S_0=\emptyset$ when $j=1$. If $j=j'$, then we have $S_0=\{w_{j'',j'}\mid j''\in\Delta, |j-j''|=1\}$.
\end{ex}

\begin{lem}\label{lem: simple length three g}
Let $j\in\Delta$. There exists a $\fz$-semi-simple uniserial $U(\fg)$-module of length $3$ with both socle and cosocle $L(1)$ and middle layer $L(s_{j})$.
\end{lem}
\begin{proof}
Upon applying $T_{w_0\cdot0}^{w_0\cdot\mu_0}$ (and using \cite[Thm.~7.8]{Hum08}), it is harmless to assume in the rest of the proof that $\mu_0=0$. We write $I\defeq \widehat{j}$, $Z$ (resp.~$Z_{I}$) for the center of $G$ (resp.~$L_{I}$) and $\fz$ (resp.~$\fz_{I}$) their associated Lie algebras. We write $1_{I}$ for the trivial object of $\mathrm{Mod}_{U(\fz_{I}/\fz)}$. Recall we have (for instance using (\ref{equ: CE complex}))
\[\mathrm{Ext}_{U(\fz_{I}/\fz)}^1(1_{I},1_{I})\cong \Hom_{E}(\fz_{I}/\fz, E)\neq 0.\]
Choose an arbitrary non-split extension $0\rightarrow 1_{I} \rightarrow \tld{1}_{I} \rightarrow 1_{I}\rightarrow 0$ in $\mathrm{Mod}_{U(\fz_{I}/\fz)}$. We define $\tld{L}^{I}(1)\defeq L^{I}(1)\otimes_E \tld{1}_{I}$, which we see in $\mathrm{Mod}_{U(\fl_{I}/\fz)}$ (and thus in $\mathrm{Mod}_{U(\fl_{I})}$) writing $\fl_{I}/\fz\cong \fl_{I}/\fz_I\times \fz_{I}/\fz$ and noting that $L^{I}(1)$ is in $\mathrm{Mod}_{U(\fl_{I}/\fz_{I})}$. It is a non-split extension
\begin{equation}\label{equ: LI deform}
0\longrightarrow L^{I}(1) \longrightarrow \tld{L}^{I}(1) \longrightarrow L^{I}(1) \longrightarrow 0,
\end{equation}
and $\tld{M}^{I}(1)\defeq U(\fg)\otimes_{U(\fp_{I})}\tld{L}^{I}(1)$, which is in the category $\cO^{\fp_{I},\infty}_{\rm{alg}}$ (see the beginning of \S \ref{subsec: category}) and is an extension of $U(\fg)$-modules
\begin{equation}\label{equ: Verma deform}
0\longrightarrow M^{I}(1) \longrightarrow \tld{M}^{I}(1) \longrightarrow M^{I}(1) \longrightarrow 0.
\end{equation}
We denote by $M^+$ the pushforward of $\tld{M}^{I}(1)$ along the surjection $M^{I}(1)\twoheadrightarrow L(1)$ (on the left). Then $M^+$ is a quotient of $\tld{M}^{I}(1)$ which fits into an exact sequence $0\rightarrow L(1)\rightarrow M^+\rightarrow M^{I}(1)\rightarrow 0$. Recall that by \cite[Thm.~9.4(c)]{Hum08} and Lemma \ref{lem: dominance and left set} $M^{I}(1)$ is the maximal quotient of $M(1)$ with constituents $L(y)$ such that $D_L(y)\subseteq \Delta\setminus I=\{j\}$. By \ref{it: Ext O 3} of Lemma~\ref{lem: Ext 1 category O} and (\ref{equ: O Hom radical}) $M(1)$ admits a (unique) length $2$ quotient $M^-$ with socle $L(s_{j})$ and cosocle $L(1)$. Hence $M^-$ is also a quotient of $M^{I}(1)$. Since the two conditions $D_L(y)\subseteq \{j\}$ and $\ell(y)\leq 1$ force $y\in\{1,s_{j}\}$, and since $L(1)$, $L(s_{j})$ occur with multiplicity $1$ in $M(1)$, it follows that any constituent $L(y)$ of $\mathrm{ker}(M^{I}(1)\twoheadrightarrow M^-)$ satisfies $\ell(y)\geq 2$. As $\mathrm{Ext}_{U(\fg)}^1(L(y),L(1))=0$ by \ref{it: rabiotext 2} of Lemma~\ref{rabiotext}, a d\'evissage yields $\mathrm{Ext}_{U(\fg)}^1(M^-,L(1))\buildrel\sim\over\rightarrow \mathrm{Ext}_{U(\fg)}^1(M^{I}(1),L(1))$. In particular $M^+$ admits a unique length $3$ quotient $M$ that fits into an exact sequence of $U(\fg)$-modules
\begin{equation}\label{equ: simple length three object}
0\longrightarrow L(1)\longrightarrow M \longrightarrow M^-\longrightarrow 0.
\end{equation}

Let $\xi: Z(\fl_I)\rightarrow E$ be the unique infinitesimal character such that $L^I(1)_{\xi}\neq 0$. By \ref{it: dominance 2} and \ref{it: dominance 3} of Lemma~\ref{lem: n coh dominance} we have $H^0(\fn_{I},L(1))_{\xi}\cong L^I(1)_{\xi}=L^I(1)\neq 0$ and $H^k(\fn_{I},L(1))_{\xi}=0$ for $k\geq 1$. Given $L(y)\in\mathrm{JH}_{U(\fg)}(M^{I}(1))\setminus\{L(1)\}$, we have $y\neq 1$ and $D_L(y)\subseteq \{j\}$ (Lemma \ref{lem: dominance and left set}), hence $D_L(y)= \{j\}$ and thus $y\notin W(L_I)$. Let $k\geq 0$ and assume $H^k(\fn_{I},L(y))_{\xi}\ne 0$. Let $L^I(z)$ be a constituent of $H^k(\fn_{I},L(y))_{\xi}$ for some $z\in W(L_I)$ (using Proposition \ref{lem: n coh O} and \cite[Thm.~1.10]{Hum08} for $L_I$). As $L^I(z)_{z\cdot \mu_0}\ne 0$, we have $H^k(\fn_{I},L(y))_{z\cdot \mu_0}\ne 0$ (weight spaces). By Lemma~\ref{lem: n coh wt} this implies $z\cdot \mu_0-y\cdot \mu_0\in \Z_{\geq 0}\Phi^+$, which implies $z\geq y$ in $W(G)$ (in view of our conventions), i.e.~$y$ is a subword of $z$ (in $W(G)$). But since $y\notin W(L_I)$ this contradicts $z\in W(L_I)$. Hence we have $H^k(\fn_{I},L(y))_{\xi}= 0$ for $k\geq 0$, and by d\'evissage we obtain in particular $H^1(\fn_{I},M^{I}(1))_\xi=0$. It follows from all this that the surjection $M^{1}(1)\twoheadrightarrow L(1)$ induces an isomorphism $H^0(\fn_{I},M^{I}(1))_{\xi}\buildrel\sim\over\rightarrow H^0(\fn_{I},L(1))_{\xi}\cong L^I(1)$, and that (\ref{equ: Verma deform}) induces a short exact sequence
\begin{equation}\label{equ: n coh deform}
0 \rightarrow H^0(\fn_{I},M^{I}(1))_{\xi} \rightarrow H^0(\fn_{I}, \tld{M}^{I}(1))_{\xi} \rightarrow H^0(\fn_{I},M^{I}(1))_{\xi} \rightarrow 0.
\end{equation}
By (\ref{equ: g spectral seq 0}) (and \ref{it: block 2} of Lemma~\ref{lem: central component}) we have canonical isomorphisms
\begin{multline*}
0\neq \Hom_{U(\fg)}(\tld{M}^{I}(1),\tld{M}^{I}(1))\cong \Hom_{U(\fl_{I})}(\tld{L}^{I}(1),H^0(\fn_{I},\tld{M}^{I}(1)))\\
\cong \Hom_{U(\fl_{I})}(\tld{L}^{I}(1),H^0(\fn_{I},\tld{M}^{I}(1))_{\xi}).
\end{multline*}
As the identity map on $\tld{M}^{I}(1)$ does not factor through $\tld{M}^{I}(1)\twoheadrightarrow M^{I}(1)$, the corresponding map $\tld{L}^{I}(1)\rightarrow H^0(\fn_{I},\tld{M}^{I}(1))_{\xi}$ (via (\ref{equ: g spectral seq 0})) does not factor through $\tld{L}^{I}(1)\twoheadrightarrow L^{I}(1)$. Hence it is injective and thus must be an isomorphism using (\ref{equ: n coh deform}) and $H^0(\fn_{I},M^{I}(1))_{\xi}\cong L^I(1)$. As $H^k(\fn_{I}, L(y))_{\xi}=0$ for $k\geq 0$ (in fact $k=0,1$ is enough) and $L(y)\in\mathrm{JH}_{U(\fg)}(M^{I}(1))\!\setminus\{L(1)\}$, the surjection $\tld{M}^{I}(1)\twoheadrightarrow M$ induces an isomorphism $\tld{L}^{I}(1) \cong H^0(\fn_{I},\tld{M}^{I}(1))_{\xi}\buildrel\sim\over\rightarrow H^0(\fn_{I},M)_{\xi}$. In particular, the short exact sequence (\ref{equ: simple length three object}) is non-split, otherwise we would have $\tld{L}^{I}(1)\cong H^0(\fn_{I},M)_{\xi}\cong L^{I}(1)\oplus H^0(\fn_{I},M^-)_{\xi}\cong L^{I}(1)\oplus L^{I}(1)$ which contradicts (\ref{equ: LI deform}). Let us prove that $M$ is actually uniserial (of length $3$). If it were not uniserial, it would contain as subquotient a non-split extension of $L(1)$ by $L(1)$. As $U(\fz)$ acts on $\tld{L}^{I}(1)$ by scalars (in $E$), it also acts on $\tld{M}^{I}(1)$ and $M$ by (the same) scalars. Since $\mathrm{Ext}_{U(\fg/\fz)}^1(L(1),L(1))=0$ (cf.~for instance \cite[(3.27)]{Schr11}), this would yield a contradiction. So $M$ must be uniserial of length $3$, with both socle and cosocle $L(1)$ and middle layer $L(s_{j})$. In particular, the unique length $2$ $U(\fg)$-submodule of $M$ is isomorphic to $(M^-)^\tau$ (see (\ref{tauduality} for the notation) which is the unique length $2$ $U(\fg)$-module with socle $L(1)$ and cosocle $L(s_{j})$ by \ref{it: rabiotext 2} of Lemma~\ref{rabiotext}.
\end{proof}

\begin{lem}\label{lem: special Verma wall crossing}
Let $I\defeq \Delta\setminus\{n-1\}$.
\begin{enumerate}[label=(\roman*)]
\item \label{it: special Verma 1} We have $M^I(w_{n-1,1})\cong L(w_{n-1,1})$, and $M^I(w_{n-1,n-k})$ is the unique length $2$ $U(\fg)$-module with socle $L(w_{n-1,n-k-1})$ and cosocle $L(w_{n-1,n-k})$ for each $0\leq k\leq n-2$ (with the convention $w_{n-1,n}=1$).
\item \label{it: special Verma 2} Let $1\leq k\leq n-1$. Then $\Theta_{s_k}(M^I(w_{n-1,n-k}))$ admits a subquotient (resp.~quotient) of the form $M^I(\mu)$ for some $\mu\in \Lambda_J^{\dom}$ if and only if $\mu\in\{w_{n-1,n-k}\cdot\mu_0,w_{n-1,n-k+1}\cdot\mu_0\}$ (resp.~if and only if $\mu=w_{n-1,n-k}\cdot\mu_0$). Moreover, $\Theta_{s_k}(M^I(w_{n-1,n-k}))$ fits into the following short exact sequence
\begin{equation}\label{equ: special Verma crossing seq}
0\rightarrow M^I(w_{n-1,n-k+1})\rightarrow \Theta_{s_k}(M^I(w_{n-1,n-k}))\rightarrow M^I(w_{n-1,n-k})\rightarrow 0.
\end{equation}
\end{enumerate}
\end{lem}
\begin{proof}
We have $W^{I,\emptyset}=\{x\in W(G)\mid D_L(x)\subseteq \{n-1\}\}=\{1\}\sqcup \{w_{n-1,n-k}\mid 1\leq k\leq n-1\}$ where the second equality follows from Lemma~\ref{lem: maximal dominance} and Remark~\ref{rem: 1 dominance}.
We fix $w\in W^{I,\emptyset}$, hence we have $w=w_{n-1,n-k}$ (with $w_{n-1,n}=1$) for some $k\in \{0,\dots,n-1\}$.

We prove \ref{it: special Verma 1}. By \cite[Thm.~9.4(c)]{Hum08} $M^I(w)$ is the maximal length quotient of $M(w)$ which belongs to $\cO^{\fp_I}_{\rm{alg}}$. We consider an arbitrary constituent $L(x)$ of $M^I(w)$ with $x>w$. By \cite[Thm.~9.4(c)]{Hum08} $L(x)\in \cO^{\fp_I}_{\rm{alg}}$ and thus by Lemma~\ref{lem: dominance and left set} $x\in W^{I,\emptyset}=\{1\}\sqcup \{w_{n-1,n-k}\mid 1\leq k\leq n-1\}$. As $x$ is partial-Coxeter (see above Lemma \ref{lem: saturated expansion}), by \ref{it: coxeter mult 1} of Lemma~\ref{lem: KL coxeter} $L(x)$ has multiplicity one in $M(x')\subseteq M(w)$ for each $x'$ satisfying $w<x'<x$ (using \cite[Thm.~5.1(a)]{Hum08}), so the quotient $M^I(w)$ of $M(w)$ which admits $L(x)$ as a constituent must also admit $L(x')$ as a constituent (use that the composition $M(x')\hookrightarrow M(w)\twoheadrightarrow M^I(w)$ is non-zero as its image contains the constituent $L(x)$). If $\ell(x)>\ell(w)+1$, we may always choose $w<x'<x$ such that $D_L(x')\not\subseteq \{n-1\}$, but $L(x')\in \mathrm{JH}_{U(\fg)}(M^I(w))$ forces $D_L(x')\subseteq \{n-1\}$ by \cite[Thm.~9.4(c)]{Hum08} (and Lemma~\ref{lem: dominance and left set}), a contradiction. Hence, we must have $\ell(x)=\ell(w)+1$, which together with $w=w_{n-1,n-k}$ forces $x=w_{n-1,n-k-1}$ (and thus $k<n-1$). When $0\leq k<n-1$, by (\ref{equ: O Hom radical}) we know that a length $2$ quotient of $M(w)$ with socle $L(x)$ and cosocle $L(w)$ exists, and is unique by \ref{it: rabiotext 2} of Lemma~\ref{rabiotext}. This finishes the proof of \ref{it: special Verma 1}.

We prove \ref{it: special Verma 2}. Assume now $1\leq k\leq n-1$. By the description of $M^I(w)$ in \ref{it: special Verma 1}, Proposition~\ref{prop: Jantzen middle} and Remark~\ref{rem: explicit middle} (and the exactness of $\Theta_{s_k}$) we deduce that $\Theta_{s_k}(M^I(w))\cong \Theta_{s_k}(L(w))$ has Loewy length $3$ with socle and cosocle $L(w)$ and middle layer $L(w_{n-1,n-k+1})\oplus L(w_{n-1,n-k-1})$ if $k<n-1$ (resp.~$L(w_{n-1,n-k+1})$ if $k=n-1$). Any $M^I(\mu)$ that appears as a subquotient of $\Theta_{s_k}(L(w))$ satisfies $\mathrm{JH}_{U(\fg)}(M^I(\mu))\subseteq \mathrm{JH}_{U(\fg)}(\Theta_{s_k}(L(w)))$, which by the previous description of $\Theta_{s_k}(L(w))$ and \ref{it: special Verma 1} forces $\mu\in\{w\cdot\mu_0, w_{n-1,n-k+1}\cdot\mu_0\}$. Moreover, as $\Theta_{s_k}(L(w))$ has cosocle $L(w)$ and $M^I(\mu)$ has cosocle $L(\mu)$, we see that $M^I(\mu)$ can be a quotient of $\Theta_{s_k}(L(w))$ only if $\mu=w\cdot\mu_0$. Finally (\ref{equ: special Verma crossing seq}) follows from the above explicit structure of $\Theta_{s_k}(M^I(w))\cong \Theta_{s_k}(L(w))$ and from \ref{it: special Verma 1}.
\end{proof}

\newpage

\section{Computing Ext groups of locally analytic representations}\label{sec: spectral seq}

We prove an important result (Corollary \ref{cor: Ext P graded}, which follows from Theorem \ref{prop: p coh graded}) which enables us to compute the $\mathrm{Ext}$ groups (\ref{extdef}) of certain locally analytic representations of $G$ by a computation of $\mathrm{Ext}$ groups purely on the smooth side and purely on the Lie algebra side.

\subsection{Fr\'echet spaces with \texorpdfstring{$U(\ft)$}{U(t)}-action}\label{subsec: Frechet t action}

We define and study certain (left) $U(\ft)$-modules over $E$ and canonical Fr\'echet completions of them. This section has an intersection with \cite{Schm13}, but our treatment is self-contained.

\begin{defn}\label{def: small t mod}
A (semi-simple) $U(\ft)$-module $M$ over $E$ is \emph{small} if the following two conditions hold.
\begin{itemize}
\item There is a $U(\ft)$-equivariant isomorphism $M\cong \bigoplus_{\mu\in\Lambda}M_{\mu}$ where $M_{\mu}$ is the eigenspace attached to the weight $\mu\in \Lambda=X(T)$, i.e.~$M$ is semi-simple with \emph{integral} weights.
\item There exist finitely many $\mu_1,\dots,\mu_k\in\Lambda$ such that $M_{\mu}\neq 0$ only if $\mu-\mu_{k'}\in\Z_{\geq 0}\Phi^+$ for some $1\leq k'\leq k$, and $M_{\mu}$ is always a finite dimensional $E$-vector space.
\end{itemize}
\end{defn}

Note that by \cite[Lemma~3.6.1]{Schm13} the second condition in Definition \ref{def: small t mod} implies that the set of weights of $M$ is relatively compact (\cite[\S 2]{Schm13}). For each $\mu\in\Lambda$, we write $\langle\mu\rangle$ for the one dimensional $U(\ft)$-module such that $\langle\mu\rangle_{\mu}\neq 0$. If $M$ is a small $U(\ft)$-module $M$, then so is $\langle\mu\rangle\otimes_E M$ for each $\mu\in\Lambda$. When $\mu\in \Z\Phi^+$, we write $\mu=\sum_{\al\in\Delta}\mu_{\al}\alpha$ for some $\mu_\alpha\in \Z$ and set $|\mu|\defeq \sum_{\al\in\Delta}\mu_\alpha \in\Z$.

\begin{rem}\label{rem: wt reduction}
\hspace{2em}
\begin{enumerate}[label=(\roman*)]
\item \label{rem: splitwt}
For $\mu,\mu'\in\Lambda$ such that $\mu-\mu'\in\Z\Phi^+$, there always exist $\mu''\in\Lambda$ such that $(\mu + \Z_{\geq 0}\Phi^+) \cup (\mu' + \Z_{\geq 0}\Phi^+)\subseteq \mu'' + \Z_{\geq 0}\Phi^+$. Putting ``together'' the weights in $\mu_i+\Z_{\geq 0}\Phi^+$ and $\mu_j+\Z_{\geq 0}\Phi^+$ for all $i,j$ such that $\mu_i - \mu_j \in \Z\Phi^+$, it follows that a small $U(\ft)$-module admits a \emph{canonical} decomposition into \emph{small} direct summands indexed by finitely many cosets of $\Lambda/\Z\Phi^+$. By considering each direct summand (with weights contained in some $\mu''+\Z_{\geq 0}\Phi^+$) separately and twisting it by $\langle-\mu''\rangle$, we will often reduce ourselves to small $U(\ft)$-modules with weights in $\Z_{\geq 0}\Phi^+$ in the sequel.
\item 
Note that the Verma module $M(\mu)=U(\fg)\otimes_{U(\fb)}\mu$, which is isomorphic to the twist $\langle\mu\rangle\otimes_E U(\fu^+)$ as $U(\ft)$-module, is a small $U(\ft)$-module, and thus so is any object in $\cO^{\fb}_{\rm{alg}}$. Conversely, any $U(\fg)$-module which is small as a $U(\ft)$-module is necessarily such that $\fu^-$ acts nilpotently, hence is the union of its finitely generated $U(\fg)$-submodules, which are all in $\cO^{\fb}_{\rm{alg}}$. An instructive example of such an $M$ which is not in $\cO^{\fb}_{\rm{alg}}$ is an infinite direct sum of Verma modules $M(\mu)$ where $\mu\in -\Lambda^+$ and $|\mu|$ tends to $+\infty$ (because for a given weight, only a finite number of such $M(\mu)$ will contribute).
\end{enumerate}
\end{rem}

Now we explicitly define a Fr\'echet topology on a small $U(\ft)$-module $M$ as follows. By \ref{rem: splitwt} of Remark \ref{rem: wt reduction}, we choose $\mu_1,\dots,\mu_k$ in $\Lambda$ which are distinct in $\Lambda/\Z\Phi^+$ and such that each weight $\mu$ of $M$ is in $\mu_{i(\mu)}+\Z_{\geq 0}\Phi^+$ for a unique $i(\mu)\in \{1,\dots,k\}$. For each weight $\mu$ of $M$, we fix a choice of a $p$-adic norm $|\cdot|_{\mu}$ on $M_{\mu}$, which is equivalent to the choice of an $\cO_E$-lattice $\{x\in M_{\mu}\mid |x|_{\mu}\leq 1\}$ in the finite dimensional $E$-vector space $M_{\mu}$. We write $\cL=\{|\cdot|_{\mu}\}_{\mu}$ for this collection of norms.
Then we define a semi-norm $|\cdot|_{\cL,r}$ on $M$ for $r\in \bQ_{>0}$ by
\begin{equation}\label{equ: general semi norm}
|x|_{\cL,r}\defeq \max_{\mu}|x_{\mu}|_{\mu}r^{|\mu-\mu_{i(\mu)}|}
\end{equation}
where $x=\sum_{\mu}x_{\mu}\in M$. It is clear that $|x|_{\cL,r}\leq |x|_{\cL,r'}$ for $r'\geq r>0$. If we modify the $\mu_i$ by elements in $\Z\Phi^+$, we get an equivalent semi-norm using that, for each $r\in \bQ_{>0}$ and $\nu\in \Z\Phi^+$, there exist $C_{\nu,r},C_{\nu,r}'\in \bQ_{>0}$ such that $C_{\nu,r}'r^{|\mu'-\nu|}\leq r^{|\mu'|}\leq C_{\nu,r}r^{|\mu'-\nu|}$ for $\mu'\in \Z_{\geq 0}\Phi^+$. The countable family of semi-norms $\{|\cdot|_{\cL,r}\}_{r\in \bQ_{>0}}$ (or equivalently the countable family $\{|\cdot|_{\cL,r}\}_{r\in p^{\bQ}}$) defines a Fr\'echet topology on $M$ which we denote by $\cT_{\cL}$. We write $\widehat{M}_{\cL}$ for the completion of $M$ under $\cT_{\cL}$. For each $U(\ft)$-submodule $M'\subseteq M$, the induced subspace topology on $M'$ can be defined by the family of semi-norms $\{|\cdot|_{\cL',r}\}_{r\in \bQ_{>0}}$ associated with the collection $\cL'=\{|\cdot|_{\mu}'\}_{\mu}$ where $|\cdot|_{\mu}'$ is the restriction of $|\cdot|_{\mu}$ to $M'_{\mu}$.\bigskip

For later convenience, we introduce the following definition.

\begin{defn}\label{def: standard semi norm}
Given a small $U(\ft)$-module $M$, a semi-norm $|\cdot|$ on $M$ is called \emph{standard} if, up to equivalence of semi-norms, it satisfies $|x|=\max_{\mu}|x_{\mu}|$ for $x=\sum_\mu x_{\mu}\in M$.
\end{defn}

Fort instance the semi-norm (\ref{equ: general semi norm}) is standard for any $r\in \bQ_{>0}$.

\begin{rem}\label{rem: wt direct sum}
Let $M$ be a small $U(\ft)$-module with a Fr\'echet topology $\cT_{\cL}$ as above. We set $M_{\Lambda'}\defeq \bigoplus_{\mu\in\Lambda'}M_{\mu}$ for any subset $\Lambda'\subseteq \Lambda$, i.e.~we only keep in $M_{\Lambda'}$ the weights of $M$ that are in $\Lambda'$. We clearly have $M\cong M_{\Lambda'}\oplus M_{\Lambda\setminus\Lambda'}$. We equip $M_{\Lambda'}$ and $M_{\Lambda\setminus\Lambda'}$ with the subspace topology (again denoted by $\cT_{\cL}$) induced from $M$. It is then clear from (\ref{equ: general semi norm}) that the original Fr\'echet topology $\cT_{\cL}$ on $M$ is equivalent to the direct sum topology on $M_{\Lambda'}\oplus M_{\Lambda\setminus\Lambda'}$, and thus we have a $U(\ft)$-equivariant topological isomorphism $\widehat{M}_{\cL}\cong \widehat{M}_{\Lambda',\cL}\oplus \widehat{M}_{\Lambda\setminus\Lambda',\cL}$.
\end{rem}

\begin{lem}\label{lem: diagonal map}
Let $M$ be a small $U(\ft)$-module and let $\cL=\{|\cdot|_{\mu}\}_{\mu}$ be as above. Let $q\in\mathrm{End}_{U(\ft)}(M)$ with $q|_{M_\mu}=a_{\mu}\in\Z$ for each $\mu\in \Lambda$ (with $a_{\mu}\defeq 0$ if $M_\mu=0$). Then $q$ is continuous for the topology on $M$ given by the semi-norms (\ref{equ: general semi norm}) and extends to a continuous endomorphism $\widehat{q}\in\mathrm{End}_{U(\ft)}\widehat{M}_{\cL}$. Assume moreover that there exists $C\in\Z_{\geq 0}$ such that $|a_{\mu}|_{\infty}\leq |\mu - \mu_{i(\mu)}|^C$ for each weight $\mu$ of $M$ (with the notation of (\ref{equ: general semi norm})), then $\widehat{q}$ is strict, and is invertible if and only if $a_{\mu}\neq 0$ for each $\mu$ such that $M_{\mu}\neq 0$.
\end{lem}
\begin{proof}
The continuity of $q$ is clear from the definition of $|\cdot|_{\cL,r}$ as well as the fact that $|a_{\mu}|_p\leq 1$.

Now we assume that there exists $C\in\Z_{\geq 0}$ such that $|a_{\mu}|_{\infty}\leq |\mu - \mu_{i(\mu)}|^C$ for any $\mu$ such that $M_\mu\ne 0$. We set $\Lambda_0\defeq \{\mu\in\Lambda\mid M_{\mu}\neq 0, a_{\mu}=0\}$, $\Lambda_1\defeq \{\mu\in\Lambda\mid M_{\mu}\neq 0, a_{\mu}\neq 0\}$ and $M_i\defeq \bigoplus_{\mu\in\Lambda_i}M_{\mu}$ for $i=0,1$. We have a $U(\ft)$-equivariant topological isomorphism $\widehat{M}_{\cL}\cong \widehat{M}_{0,\cL}\oplus \widehat{M}_{1,\cL}$ by Remark~\ref{rem: wt direct sum}.
As $\widehat{q}(x)=\sum_{\mu\in\Lambda}a_{\mu}x_{\mu}$ for each convergent sum $x=\sum_{\mu\in\Lambda}x_{\mu}\in\widehat{M}_{\cL}$, it is clear that $\mathrm{ker}(\widehat{q})=\widehat{M}_{0,\cL}$ and thus we have a $U(\ft)$-equivariant topological isomorphism $\widehat{M}_{\cL}/\mathrm{ker}(\widehat{q})\cong \widehat{M}_{1,\cL}$.
It suffices to show that $q': M_1\rightarrow M_1, \sum_{\mu\in\Lambda_1}x_{\mu}\mapsto \sum_{\mu\in\Lambda_1}a_{\mu}^{-1}x_{\mu}$ extends continuously to $\widehat{M}_{1,\cL}$ (this will necessarily give the inverse of $\widehat{q}|_{\widehat{M}_{1,\cL}}$). For each $r'>r>0$ in $\bQ$, our assumption $|a_{\mu}|_{\infty}\leq |\mu-\mu_{i(\mu)}|^C$ (together with $a_\mu \in \Z$ and $a_{\mu}\neq 0$ for $\mu\in\Lambda_1$) implies the existence of $C_{r,r'}\in\bQ_{>0}$ such that $|a_{\mu}|_p^{-1}(r/r')^{|\mu-\mu_{i(\mu)}|}\leq C_{r,r'}$ for each $\mu\in\Lambda_1$, which together with (\ref{equ: general semi norm}) gives the continuity of $q'$ on $M_1$.
\end{proof}

\begin{rem}
It follows from Lemma~\ref{lem: diagonal map} that the $U(\ft)$-action on $M$ extends to a continuous $U(\ft)$-action on $\widehat{M}_{\cL}$ (i.e.~any element of $U(\ft)$ acts continuously): an arbitrary element $x\in \widehat{M}_{\cL}$ can be expressed as a convergent infinite sum $\sum_{\mu\in\Lambda}x_{\mu}$ with $x_{\mu}\in M_{\mu}$ for $\mu\in\Lambda$, and $t\cdot x=\sum_{\mu\in\Lambda}\mu(t)x_{\mu}$. In particular, we have $x\in(\widehat{M}_{\cL})_{\mu}$, i.e.~$t\cdot x=\mu(t)x$ for $t\in\ft$, if and only if $x=x_{\mu}$. In other words, the embedding $M\hookrightarrow \widehat{M}_{\cL}$ induces an equality $M_{\mu}=(\widehat{M}_{\cL})_{\mu}$ for each $\mu\in\Lambda$.
\end{rem}

It is convenient to introduce the following definitions.

\begin{defn}\label{def: Frechet t module}
\hspace{2em}
\begin{enumerate}[label=(\roman*)]
\item Let $V$ be a Fr\'echet space over $E$. We say that $V$ is a \emph{Fr\'echet $U(\ft)$-module} if it is equipped with a continuous $U(\ft)$-action.
\item \label{it: Frechet t 2} We say that a Fr\'echet $U(\ft)$-module $V$ is \emph{small} if there exists a (semi-simple) small $U(\ft)$-module $M$ (see Definition~\ref{def: small t mod}) and a collection of norms $\cL=\{|\cdot|_{\mu}\}_{\mu}$ such that we have a $U(\ft)$-equivariant topological isomorphism $V\cong \widehat{M}_{\cL}$.
\end{enumerate}
\end{defn}

\begin{rem}
Given a small Fr\'echet $U(\ft)$-module $V\cong \widehat{M}_{\cL}$, it is clear that the Fr\'echet $U(\ft)$-module $\langle\mu_1\rangle\otimes_E V$ is small for each $\mu_1\in\Lambda$.
\end{rem}

We refer to \cite[\S 17.B]{S02} for the definition and properties of the semi-norm tensor product of two semi-norms.

\begin{lem}\label{lem: tensor with fin dim}
Let $V\cong \widehat{M}_{\cL}$ be a small Fr\'echet $U(\ft)$-module and $N$ a finite dimensional $U(\ft)$-module (with its canonical Banach topology). We fix an arbitrary norm $|\cdot|_{N,\nu}$ on $N_{\nu}$ for each $\nu\in \Lambda$ such that $N_{\nu}\neq 0$.
\begin{enumerate}[label=(\roman*)]
\item \label{it: tensor with fin dim 1}
For $\mu\in\Lambda$ and $\sum_{\nu}y_{\nu}\otimes_E x_{\mu-\nu}\in (N\otimes_E M)_{\mu}$ with $y_{\nu}\in N_{\nu}\neq 0$ and $x_{\mu-\nu}\in M_{\mu-\nu}$ we define
\begin{equation}\label{equ: tensor wt space}
\Big|\sum_{\nu}y_{\nu}\otimes_E x_{\mu-\nu}\Big|_{N\otimes_E M,\mu}\defeq \max_{\nu}|y_{\nu}|_{N,\nu}|x_{\mu-\nu}|_{\mu-\nu}.
\end{equation}
Then $N\otimes_E V\cong (\widehat{N\otimes_E M})_{\cL'}$ where the collection of norms $\cL'$ on the weight spaces $(N\otimes_E M)_{\mu}$ is uniquely determined by (\ref{equ: tensor wt space}). In particular $N\otimes_E V$ is a small Fr\'echet $U(\ft)$-module.
\item \label{it: tensor with fin dim 2}
For any standard semi-norm $|\cdot|$ on $M$ (Definition \ref{def: standard semi norm}) the semi-norm $(\max_{\nu}|\cdot|_{N,\nu})\otimes_E |\cdot|$ on $N\otimes_E M$ is standard.
\end{enumerate}
\end{lem}
\begin{proof}
We prove \ref{it: tensor with fin dim 1}. Note first that $N\otimes_EM$ is a small $U(\ft)$-module and that the Fr\'echet topology on $N\otimes_E V$ is defined by the family of semi-norms $\{(\max_{\nu}|\cdot|_{N,\nu})\otimes_E |\cdot|_{\cL,r}\}_{r\in \bQ_{>0}}$. Let $z=\bigoplus_{\mu\in\Lambda}z_{\mu}\in N\otimes_E M$ with $z_{\mu}=\sum_{\nu}y_{\nu}\otimes_E x_{\mu-\nu}$ for some $y_{\nu}$ and $x_{\mu-\nu}$ as in (\ref{equ: tensor wt space}). We have by (\ref{equ: general semi norm})
\begin{equation*}
|z|_{\cL',r}=\max_{\mu,\nu}|y_{\nu}|_{N,\nu}|x_{\mu-\nu}|_{\mu-\nu}r^{|\mu-\mu_{i(\mu)}|}.
\end{equation*}
The evaluation of the semi-norm $(\max_{\nu}|\cdot|_{N,\nu})\otimes_E |\cdot|_{\cL,r}$ on $z$ is
\begin{equation}\label{equ: tensor semi norm 2}
\max_{\mu,\nu}|y_{\nu}|_{N,\nu}|x_{\mu-\nu}|_{\mu-\nu}r^{|\mu-\nu-\mu_{i(\mu-\nu)}|}.
\end{equation}
As there exists only finitely many $\nu$ such that $N_{\nu}\neq 0$, we conclude from (\ref{equ: tensor semi norm 2}) that for $r\in \bQ_{>0}$ the semi-norms $|\cdot|_{\cL',r}$ and $(\max_{\nu}|\cdot|_{N,\nu})\otimes_E |\cdot|_{\cL,r}$ are equivalent on $N\otimes_E M$. It follows that $N\otimes_E V\cong (\widehat{N\otimes_E M})_{\cL'}$.

The proof of \ref{it: tensor with fin dim 2} is an easy exercise that is left to the reader.
\end{proof}

The following result is contained in \cite[\S 1.3.1]{Lac99}, we reproduce a proof for the reader's convenience.

\begin{lem}\label{lem: t closed}
Let $V_0$ be a small Fr\'echet $U(\ft)$-module and $V\subseteq V_0$ a closed Fr\'echet $U(\ft)$-submodule. Then $V$ is again small. In particular $\bigoplus_{\mu\in\Lambda}V_{\mu}$ is a small $U(\ft)$-module which is dense in $V$.
\end{lem}
\begin{proof}
We fix throughout this proof a small $U(\ft)$-module $M$ and a collection of norms $\cL=\{|\cdot|_{\mu}\}_{\mu}$ such that we have a $U(\ft)$-equivariant topological isomorphism $V_0\cong \widehat{M}_{\cL}$.
We may write $x=\sum_{\mu\in\Lambda}x_{\mu}$ as a convergent infinite sum with $x_{\mu}\in M_{\mu}$ for each $\mu\in\Lambda$. Note that $V_{\mu}\subseteq (V_0)_\mu=M_{\mu}$ is finite dimensional for each $\mu\in\Lambda$. Hence to prove that $V$ is small, it suffices to show that $\bigoplus_{\mu\in\Lambda}V_{\mu}$ is dense in $V$.
In other words, for each convergent infinite sum $x=\sum_{\mu\in\Lambda}x_{\mu}\in V\subseteq V_0$ with $x_\mu\in M_{\mu}$, we need to show that
\begin{equation}\label{equ: weight vector}
x_{\mu}\in V.
\end{equation}

Let $e\in\ft$ be an element such that $\mu(e)\in \Z$ for all $\mu\in \Lambda$ and $\alpha(e)\in\Z_{>0}$ for all $\alpha\in\Delta$. Replacing $V_0$ by $\langle-\mu''\rangle\otimes_E V_0$ for some $\mu''\in \Lambda$, we can moreover assume that $\mu(e)\in \Z_{\geq 0}$ for all the weights $\mu$ of $V_0$.
For each $N\in\Z_{\geq 0}$, we set
\[\Lambda_N\defeq \{\mu\mathrm{\ weight\ of\ }V_0\mid \mu(e)=N\}\]
and similarly for $\Lambda_{\leq N}$, $\Lambda_{<N}$ and $\Lambda_{>N}$. Writing $\mu=\mu_{i(\mu)}+\sum_{\al\in\Delta}\mu_{\al}\alpha$ with $\mu_{i(\mu)}$ as in (\ref{equ: general semi norm}) and $\mu_\alpha\in \Z_{\geq 0}$, we have $\mu(e)=\mu_{i(\mu)}(e)+\sum_{\al\in\Delta}\mu_{\al}\alpha(e)$. Since $\alpha(e)\in\Z_{>0}$ for all $\alpha$ and since there is a finite number of $\mu_i$, there is only a finite number (possibly $0$) of tuples $(\mu_\alpha)_{\alpha\in \Delta}$ in $\Z_{\geq 0}^{\Delta}$ such that $\sum_{\al\in\Delta}\mu_{\al}\alpha(e)=N-\mu_i(e)$ for some $i$, i.e.~$\Lambda_N$ is always a finite set. Therefore, to prove (\ref{equ: weight vector}), it suffices to prove by an increasing induction on $N\in\Z_{\geq 0}$ that
\begin{equation}\label{equ: weight vector e}
\sum_{\mu\in\Lambda_N}x_{\mu}\in V
\end{equation}
for each $N\in\Z_{\geq 0}$ (use the action of $\U(\ft)$ on $\sum_{\mu\in\Lambda_N}x_{\mu}$ and the fact the $x_\mu$ are weight vectors with distinct weights to deduce that each $x_\mu$ is then also in $V$). We fix $N\in \Z_{\geq 0}$, assume inductively that (\ref{equ: weight vector e}) holds if $N$ is replaced by any $N'<N$, and prove the statement for $N$. Upon replacing $x$ by $x-\sum_{N'<N}\sum_{\mu\in\Lambda_{N'}}x_{\mu}=x-\sum_{\mu\in\Lambda_{<N}}x_{\mu}$, we can assume $x_{\mu}=0$ for $\mu\in\Lambda_{<N}$.

For each $N_1\in \Z_{> 0}$, we set
\begin{equation*}
\Theta_{e}^{N,N_1}\defeq \frac{\prod_{k=1}^{N_1}(N+k-e)}{N_1!}\in U(\ft).
\end{equation*}
Then $x_{\mu}$ for $\mu$ a weight of $V_0$ is an eigenvector for $\Theta_{e}^{N,N_1}$ with eigenvalue
\begin{equation*}
C_{\mu,N,N_1}\defeq \frac{\prod_{k=1}^{N_1}(N-\mu(e)+k)}{N_1!}\in \Z.
\end{equation*}
In particular, we have $C_{\mu,N,N_1}=1$ if $\mu\in\Lambda_N$, and $C_{\mu,N,N_1}=0$ if and only if $\mu\in\Lambda_{N'}$ for some $N<N'\leq N+N_1$. As $\Theta_{e}^{N,N_1}\in U(\ft)$ and $V$ is $U(\ft)$-stable, we have $\Theta_{e}^{N,N_1}(x)\in V$. Since $V$ is closed in $V_0$, to prove (\ref{equ: weight vector e}), it suffices to show that the sequence $\{\Theta_{e}^{N,N_1}(x)\}_{N_1\in \Z_{\geq 1}}$ of vectors of $V$ converges to $\Theta_{e}^{N,N_1}(\sum_{\mu\in\Lambda_N}x_{\mu})=\sum_{\mu\in\Lambda_N}x_{\mu}$ inside $V_0$ when $N_1\rightarrow+\infty$, or equivalently that $\{\Theta_{e}^{N,N_1}(x-\sum_{\mu\in\Lambda_N}x_{\mu})\}_{N_1\in \Z_{\geq 1}}$ converges to zero in $V_0$ when $N_1\rightarrow+\infty$. As $C_{\mu,N,N_1}\in \Z$, we have $|C_{\mu,N,N_1}|_p\leq 1$ and thus for $r\in \bQ_{>0}$:
\begin{equation*}
\left\vert\Theta_{e}^{N,N_1}\Big(x-\sum_{\mu\in\Lambda_N}x_{\mu}\Big)\right\vert_{\cL,r}=\ \ \left\vert\sum_{\mu\in\Lambda_{>N+N_1}}C_{\mu,N,N_1}x_{\mu}\right\vert_{\cL,r}\leq \ \ \left\vert\sum_{\mu\in\Lambda_{>N+N_1}}x_{\mu}\right\vert_{\cL,r}.
\end{equation*}
Since $|\mu-\mu_{i(\mu)}|= \sum_{\al\in\Delta}\mu_\alpha$ tends to $+\infty$ if and only if $\mu(e)=\mu_{i(\mu)}(e)+\sum_{\al\in\Delta}\mu_{\al}\alpha(e)$ tends to $+\infty$, we see that the convergence of the infinite sum $x=\sum_{\mu\in\Lambda}x_{\mu}$ in $V_0$ implies $\varinjlim_{N_1\rightarrow+\infty}|\sum_{\mu\in\Lambda_{>N+N_1}}x_{\mu}|_{\cL,r}=0$ for $r\in \bQ_{>0}$, which gives the desired result.
\end{proof}

\begin{rem}\label{rem: closure bijection}
By Lemma~\ref{lem: t closed}, there exists a natural bijection between $U(\ft)$-submodu\-les of $M$ and closed Fr\'echet $U(\ft)$-submodules of $\widehat{M}_{\cL}$, given by sending $M'\subseteq M$ to its closure in $\widehat{M}_{\cL}$, with inverse given by sending a closed Fr\'echet $U(\ft)$-submodule $V\subseteq \widehat{M}_{\cL}$ to its dense $U(\ft)$-submodule
\[\bigoplus_{\mu\in\Lambda}V_{\mu}=\bigoplus_{\mu\in\Lambda}V\cap M_{\mu}=V\cap\Big(\bigoplus_{\mu\in\Lambda}M_{\mu}\Big)=V\cap M.\]
This is special case of \cite[Satz~1.3.19]{Lac99}, see also \cite[Prop.~2.0.1]{Schm13}.
\end{rem}

\begin{rem}\label{rem: dense in Banach}
Let $M$ be a small $U(\ft)$-module, $|\cdot|$ a standard semi-norm on $M$ (see Definition~\ref{def: standard semi norm}) and $\widehat{M}$ the corresponding completion. The proof of Lemma~\ref{lem: t closed} actually shows that, for any closed $U(\ft)$-stable Banach subspace $V\subseteq \widehat{M}$, the subspace $V\cap M$ is dense in $V$. In particular, there exists a natural bijection between $U(\ft)$-submodules of $M$ and closed $U(\ft)$-stable Banach subspaces of $\widehat{M}$ which sends $M'\subseteq M$ to its closure $\widehat{M}'$ in $\widehat{M}$. Note that the induced semi-norm on the submodule $M'$ (resp.~on the quotient $M/M'$) of $M$ is again standard, with $\widehat{M}'$ (resp.~$\widehat{M}/\widehat{M}'$) being the corresponding completion.
\end{rem}

\begin{lem}\label{lem: closed complement}
Let $V_0$ be a small Fr\'echet $U(\ft)$-module and $V\subseteq V_0$ a closed Fr\'echet $U(\ft)$-submodule. Then there exists another closed Fr\'echet $U(\ft)$-submodule $V'\subseteq V_0$ such that the natural map $V\oplus V'\rightarrow V_0$ is a $U(\ft)$-equivariant topological isomorphism.
\end{lem}
\begin{proof}
We fix a small $U(\ft)$-module $M$ and a collection of norms $\cL=\{|\cdot|_{\mu}\}_{\mu}$ such that $V_0\cong \widehat{M}_{\cL}$. By Remark~\ref{rem: closure bijection}, the choice of $V'$ is equivalent to the choice of a subspace $V'_{\mu}\subseteq M_{\mu}$ for each $\mu\in\Lambda$. From the definition of the semi-norms $|\cdot|_{\cL,r}$ for $r\in \bQ_{>0}$ in (\ref{equ: general semi norm}), we see that it suffices to construct $V'_{\mu}\subseteq M_{\mu}$ such that $V_{\mu}\oplus V'_{\mu}\cong M_{\mu}$ and
\begin{equation}\label{equ: opposite subspace}
|x|_{\mu}=\max\{|y_{\mu}|_{\mu}, |z_{\mu}|_{\mu}\}
\end{equation}
for $x=y+z\in M_{\mu}$ with $y\in V_{\mu}$ and $z\in V'_{\mu}$. If we set $M_{\mu}^{\circ}\defeq \{x\in M_{\mu}\mid |x|_{\mu}\leq 1\}$, which is an $\cO_E$-lattice in the $E$-vector space $M_{\mu}$, then (\ref{equ: opposite subspace}) is equivalent to the equality (inside $M_\mu$)
\begin{equation*}
M_{\mu}^{\circ}=(V_{\mu}\cap M_{\mu}^{\circ})+(V'_{\mu}\cap M_{\mu}^{\circ}).
\end{equation*}
But such $V'_{\mu}\subseteq M_{\mu}$ exists for each given $V_{\mu}$ (and our fixed $M_{\mu}^{\circ}\subseteq M_{\mu}$) because we can extend an arbitrary $\cO_E$-basis of $V_{\mu}\cap M_{\mu}^{\circ}$ into one of $M_{\mu}^{\circ}$, and define $V'_{\mu}$ as the $E$-span of the new basis elements we added. This finishes the proof.
\end{proof}

\begin{rem}
We consider a ($U(\ft)$-equivariant) strict exact sequence of Fr\'echet $U(\ft)$-modules
\begin{equation}\label{equ: split seq}
0\longrightarrow V_1\longrightarrow V_2\longrightarrow V_3\longrightarrow 0
\end{equation}
with $V_2$ small (and thus $V_1$ is also small by Lemma~\ref{lem: t closed}). By Lemma~\ref{lem: closed complement} there exists another (small) Fr\'echet $U(\ft)$-submodule $V_3'\subseteq V_2$ such that the natural map $V_1\oplus V_3'\rightarrow V_2$ is a $U(\ft)$-equivariant topological isomorphism, which together with our assumption $V_2/V_1\buildrel\sim\over\rightarrow V_3$ gives a $U(\ft)$-equivariant topological isomorphism $V_3'\buildrel\sim\over\rightarrow V_3$, and thus a $U(\ft)$-equivariant topological isomorphism $V_2\cong V_1\oplus V_3$. In other words, the short exact sequence (\ref{equ: split seq}) always splits (non-canonically). Moreover, since $V'_3$ is small so is $V_3$, and we have canonical isomorphisms of small $U(\ft)$-modules
\begin{equation}\label{equ: quotient wt}
\Big(\bigoplus_{\mu\in\Lambda}V_{2,\mu}\Big)\Big/\Big(\bigoplus_{\mu\in\Lambda}V_{1,\mu}\Big)\cong \bigoplus_{\mu\in\Lambda}V_{2,\mu}/V_{1,\mu}\cong \bigoplus_{\mu\in\Lambda}V_{3,\mu}.
\end{equation}
\end{rem}

\begin{lem}\label{lem: invertible action}
Let $e\in\ft$ such that $\al(e)\in\Z_{>0}$ for each $\al\in\Delta$ and let $V$ be a small Fr\'echet $U(\ft)$-module. Then the action of $e-N$ on
\begin{equation*}
V\big/\Big(\bigoplus_{\mu\in\Lambda}V_{\mu}\Big)
\end{equation*}
is invertible for any $N\in\Z$.
\end{lem}
\begin{proof}
By \ref{rem: splitwt} of Remark~\ref{rem: wt reduction} we can assume $V_{\mu}\neq 0$ only if $\mu\in\Z_{\geq 0}\Phi^+$. As in the proof of Lemma \ref{lem: t closed}, the set $\Lambda_N=\{\mu\in \Z_{\geq 0}\Phi^+\mid \mu(e)=N\}$ is finite and we write $\Lambda_{\neq N}$ for its complement in $\Z_{\geq 0}\Phi^+$. Let $V_{\neq N}$ be the closure of $\bigoplus_{\mu\in\Lambda_{\neq N}}V_{\mu}$ in $V$, then by Remark~\ref{rem: wt direct sum} we have a canonical isomorphism of Fr\'echet $U(\ft)$-modules
\begin{equation}\label{equ: remove kernel}
V\cong V_{\neq N}\oplus \Big(\bigoplus_{\mu\in\Lambda_N}V_{\mu}\Big),
\end{equation}
with $\bigoplus_{\mu\in\Lambda_N}V_{\mu}$ being (finite dimensional and) exactly the kernel of $e-N$ in $V$.
The isomorphism (\ref{equ: remove kernel}) induces the following isomorphism of $U(\ft)$-modules
\begin{equation}\label{equ: remove finite wt}
V\Big/\Big(\bigoplus_{\mu\in\Lambda}V_{\mu}\Big)\cong V_{\neq N}\Big/\Big(\bigoplus_{\mu\in\Lambda}(V_{\neq N})_{\mu}\Big)=V_{\neq N}\Big/\Big(\bigoplus_{\mu\in\Lambda_{\neq N}}V_{\mu}\Big).
\end{equation}
By Lemma~\ref{lem: diagonal map} applied to $M=\bigoplus_{\mu\in\Lambda}(V_{\neq N})_{\mu}$ and $q=e-N$, we know that $e-N$ has invertible action on $V_{\neq N}$, and thus on \emph{any} $U(\ft)$-equivariant quotient of $V_{\neq N}$. In particular, $e-N$ has invertible action on (\ref{equ: remove finite wt}).
\end{proof}

\begin{rem}\label{rem: non semi simple case}
The results in this section can be generalized to certain $U(\ft)$-modules $M$ which are not necessarily semi-simple as follows. In the first part of Definition~\ref{def: small t mod}, we can replace ``where $M_{\mu}$ is the eigenspace attached to the weight $\mu$'' by ``where $M_{\mu}$ is the generalized eigenspace attached to the weight $\mu$ such that there is a fixed $N\in \Z_{\geq 1}$ satisfying $(t-\mu(t))^N=0$ on $M_{\mu}$ for all $\mu$ and all $t\in \ft$'' (the second part of Definition~\ref{def: small t mod} being unchanged). Then such an $M$ contains a sequence of $U(\ft)$-submodules $M=M_N\supseteq M_{N-1}\supseteq \cdots \supseteq M_1\supseteq M_0=0$ such that $M_{k}/M_{k-1}$ is the maximal semi-simple $U(\ft)$-submodule of $M/M_{k-1}$ for $1\leq k\leq N$. We can also fix a collection of norms $\cL=\{|\cdot|_{\mu}\}_{\mu}$ on each $M_\mu$ and take the completion $V\defeq \widehat{M}_{\cL}$ for the same semi-norms (\ref{equ: general semi norm}). For $0\leq k\leq N$ we define $V_k$ as the closure of $M_k$ in $V$, and we easily check that $V_k/V_{k-1}$ for $1\leq k\leq N$ is the completion of $M_k/M_{k-1}$ for the quotient topology (which is given by (\ref{equ: general semi norm}) for the quotient norms of the norms $|\cdot|_{\mu}$).

Let $V'\subseteq V$ be a closed Fr\'echet $U(\ft)$-submodule. For $0\leq k\leq N$ we have a closed Fr\'echet $U(\ft)$-submodule $V'_k\defeq V'\cap V_k\subseteq V_k$ . As $M_k/M_{k-1}$ is small (semi-simple) with completion $V_k/V_{k-1}$, by Lemma~\ref{lem: t closed} the closed Fr\'echet $U(\ft)$-submodule $V'_k/V'_{k-1}\subseteq V_k/V_{k-1}$ is also small with
\[\bigoplus_{\mu\in\Lambda}(V'_k/V'_{k-1})_{\mu}=(V'_k/V'_{k-1})\cap (M_k/M_{k-1})\subseteq M_k/M_{k-1}.\]
By an increasing induction on $k\in \{1,\dots, N\}$ one easily checks that $V'\cap M_k$ is dense in $V'_k$ with $V'\cap M_k=\bigoplus_{\mu\in\Lambda}(V'_k)_{\mu}$, and hence satisfies the generalized Definition~\ref{def: small t mod} (as defined above).

For such non-semi-simple $M$ as above, Lemma~\ref{lem: closed complement} clearly fails in general, and likewise a strict exact sequence as in (\ref{equ: split seq}) does not split in general. However, one can prove (using the above filtration $V=V_N\supseteq \cdots \supseteq V_0$) that (\ref{equ: quotient wt}) still holds. Lemma~\ref{lem: invertible action} also admits a straightforward generalization to such non-semi-simple $M$.
\end{rem}

\subsection{Preliminaries on locally analytic distributions}\label{subsec: loc an dist}

We prove several important (technical) results on locally analytic and locally constant distribution algebras with prescribed support. All these results are used in the next sections.\bigskip

For a paracompact locally $K$-analytic manifold $M$ of finite dimension (\cite[\S 8]{S11}), we write $C^{\rm{an}}(M)$ for the space of $E$-valued locally $K$-analytic functions on $M$ equipped with its usual locally convex topology (\cite[\S 12]{S11}). We write $D(M)\defeq C^{\rm{an}}(M)_b^\vee$ for its continuous $E$-dual equipped with the strong topology (\cite[\S 9]{S02}, \cite[\S 2]{ST02a}). For a closed subset $C\subseteq M$, we let $D(M)_C\subseteq D(M)$ be the closed subspace of distributions supported on $C$, which is the strong dual of the quotient $C^{\rm{an}}(M)/C^{\rm{an}}(M)_{M\setminus C}$ where $C^{\rm{an}}(M)_{M\setminus C}$ is the closed subspace of $C^{\rm{an}}(M)$ of functions with support contained in the open $M\setminus C$ (\cite[\S 1.2]{Koh07}).\bigskip

Similarly, we write $C^{\infty}(M)$ for the space of $E$-valued locally constant functions on $M$ equipped with its usual locally convex topology and $D^{\infty}(M)\defeq C^{\infty}(M)^\vee_b$ for its strong continuous $E$-dual. We have a closed embedding $C^{\infty}(M)\hookrightarrow C^{\rm{an}}(M)$ (use that $C^{\infty}(M)$ is the kernel of the continuous map $C^{\rm{an}}(M)\rightarrow C^{\rm{an}}(TM)$ where $TM$ is the tangent bundle, see \cite[\S 9]{S11}, \cite[Def.~12.4.i]{S11} and \cite[Rem.~6.2]{S11}) which induces a strict continuous surjection $D(M)\twoheadrightarrow D^{\infty}(M)$ (\cite[\S 2]{ST01}). We write $D^{\infty}(M)_C$ for the image of $D(M)_C$ under $D(M)\twoheadrightarrow D^{\infty}(M)$.

\begin{lem}\label{lem: closedinfini}
The subspace $D^{\infty}(M)_C$ is closed in $D^{\infty}(M)$.
\end{lem}
\begin{proof}
Define $\widetilde D^{\infty}(M)_C$ as the closed subspace of $D^{\infty}(M)$ which is the strong dual of $C^{\infty}(M)/C^{\infty}(M)_{M\setminus C}$ where $C^{\infty}(M)_{M\setminus C}$ is the closed subspace of $C^{\infty}(M)$ of functions with support contained in $M\setminus C$. The composition $D(M)_C\hookrightarrow D(M) \twoheadrightarrow D^{\infty}(M)$ factors through a continuous map $D(M)_C\rightarrow \widetilde D^{\infty}(M)_C$, and it is enough to prove that this map is surjective. Arguing exactly as in the first part of the proof of \cite[Lemme~3.12]{BD23} for both $C^{\infty}(M)$ and $C^{\rm{an}}(M)$, we can assume that the closed subset $C$ is compact. Then using \cite[(3.3)]{BD23} for both $\widetilde D^{\infty}(M)_C$ and $D(M)_C$, it is enough to prove that the natural injection of locally convex spaces of compact type $\varinjlim_U C^{\infty}(U) \hookrightarrow \varinjlim_U C^{\rm{an}}(U)$ is a closed embedding, where $U$ runs among the compact open subsets of $M$ containing $C$ and the transition maps are the restrictions. But this easily follows from the short exact sequence
\[0\longrightarrow \varinjlim_U C^{\infty}(U) \longrightarrow \varinjlim_U C^{\rm{an}}(U) \longrightarrow \varinjlim_U C^{\rm{an}}(TU)\]
where $TU$ is the tangent bundle of $U$ (see \cite[\S 9]{S11}).
\end{proof}

If $M\cong M_1\times M_2$ and the closed subset $C\subseteq M$ has the form $C_1\times C_2$ with $C_i\subseteq M_i$ for $i=1,2$, then we have a canonical topological isomorphism (\cite[Lemme~3.12]{BD23})
\begin{equation}\label{equ: loc an product}
D(M)_C\cong D(M_1)_{C_1}\widehat{\otimes}_{E,\iota} D(M_2)_{C_2}
\end{equation}
and by the same proof as for \emph{loc.~cit.}
\begin{equation}\label{equ: smooth product}
D^{\infty}(M)_C\cong D^{\infty}(M_1)_{C_1}\widehat{\otimes}_{E,\iota} D^{\infty}(M_2)_{C_2}.
\end{equation}

\begin{lem}\label{lem: smooth restriction}
Let $M$ be a paracompact locally $K$-analytic manifold and $C\subseteq M$ a closed paracompact locally $K$-analytic submanifold. The embedding $C\hookrightarrow M$ induces a closed topological embedding $D^{\infty}(C)\hookrightarrow D^{\infty}(M)$, which induces a topological isomorphism
\begin{equation}\label{equ: smooth restriction}
D^{\infty}(C)\xrightarrow{\sim} D^{\infty}(M)_C.
\end{equation}
\end{lem}
\begin{proof}
The first statement follows from \cite[Prop.~1.1.2]{Koh07} (and its proof). It suffices to show that (\ref{equ: smooth restriction}) is surjective. If $(M_i)_{i\in I}$ is a covering of $M$ by compact open disjoint subsets, then $(C_i)_{i\in I}\defeq (M_i\cap C)_{i\in I}$ is a covering of $C$ by compact open disjoint subsets, and using $D^{\infty}(C)\cong \bigoplus_i D^{\infty}(C_i)$, $D(M)_C\cong \bigoplus_i D(M_i)_{C_i}$ (see e.g.~the proof of \cite[Lemme~3.12]{BD23}), we have $D^\infty(M)_C\cong \bigoplus_i D^\infty(M_i)_{C_i}$ which shows that we can assume $M$ and $C$ compact. Let $\delta\in D^{\infty}(M)_C$, we define below an element $\overline{\delta}\in D^{\infty}(C)$ which maps to $\delta$. For each $f\in C^{\infty}(C)$, there exists a finite partition $\cP$ of $C$ into compact open subsets $U$ such that $f|_{U}$ is a constant for each $U\in \cP$. We can moreover choose an arbitrary partition $\tld{\cP}$ of $M$ into compact open subsets $\tld{U}$ such that we have either $\tld{U}\cap C\in\cP$ or $\tld{U}\cap C=\emptyset$ for $\tld{U}\in\tld{\cP}$. Then we define $\tld{f}\in C^{\infty}(M)$ by requiring that $\tld{f}$ is a constant on each $\tld{U}\in\tld{\cP}$ and that $\tld{f}|_{C}=f$. Now we define $\overline{\delta}(f)\defeq \delta(\tld{f})\in E$ and observe that $\delta(\tld{f})$ is independent of the choice of $\tld{f}$ (and the partitions ${\cP}$, $\tld{\cP}$) as above. Indeed, if $\tld{f}'$ is another choice, then by construction $\tld{f}-\tld{f}'$ is zero in an open neighborhood of $C$ and thus $\delta(\tld{f}-\tld{f}')=0$ since $\delta\in D^{\infty}(M)_C$. It is then clear that $\overline{\delta}\in \Hom_E(C^\infty(M), E)\cong D^{\infty}(C)$ has image $\delta$ under $D^{\infty}(C)\hookrightarrow D^{\infty}(M)$.
\end{proof}

If $M=G$ is a locally $K$-analytic group (automatically paracompact by \cite[Cor.~18.8]{S11}), then $D(G)$ is a unital associative algebra (with multiplication being the convolution, see \cite[Prop.~2.3]{ST02a}). For a locally $K$-analytic closed subgroup $H\subseteq G$, the closed subspace $D(G)_H$ is closed under convolution (\cite[Cor.~1.2.6]{Koh07}), making $D(G)_H$ a unital associative subalgebra of $D(G)$. Similarly, $D^{\infty}(G)$ is also a unital associative algebra isomorphic to the quotient $D(G)\otimes_{U(\fg)}E$ of $D(G)$ by the closed two-sided ideal generated by the Lie algebra $\fg$ of $G$ (cf.~\cite[\S 13]{S11} and \cite[Rem.~1.1(iii)]{ST05}). The surjective algebra homomorphism $D(G)\twoheadrightarrow D^{\infty}(G)$ induces a surjective homomorphism $D(G)_H\twoheadrightarrow D^{\infty}(G)_H$ showing that $D^{\infty}(G)_H$ is a closed unital associative subalgebra of $D^{\infty}(G)$ isomorphic to $D^{\infty}(H)$ (using Lemma \ref{lem: closedinfini} and Lemma~\ref{lem: smooth restriction}). When $H$ is the trivial group, we write $D(G)_1$ and $D^{\infty}(G)_1$ instead of $D(G)_{\{1\}}$, $D^{\infty}(G)_{\{1\}}$.\bigskip

For $A$ a Fr\'echet-Stein algebra we let $\cC_A$ be the abelian category of coadmissible left $A$-modules (see \cite[\S 3]{ST03}), which is a full subcategory of $\mathrm{Mod}_{A}$ (\cite[Cor.~3.5]{ST03}). By the discussion before and after \cite[Lemma~3.6]{ST03}, each $M\in \cC_A$ carries a canonical topology as a $E$-Fr\'echet space and any $A$-linear map in $\cC_A$ is continuous and strict.\bigskip

If $G$ is a compact locally $K$-analytic group, by \cite[Thm.~5.1]{ST03} the algebra $D(G)$ is Fr\'echet-Stein, and so is its quotient $D^{\infty}(G)$ by \cite[Prop.~3.7]{ST03}. Moreover the continuous surjection $D(G)\twoheadrightarrow D^{\infty}(G)$ induces a fully faithful embedding $\cC_{D^{\infty}(G)}\hookrightarrow \cC_{D(G)}$. When $G$ is not necessarily compact, one defines the category $\cC_{D(G)}$ of coadmissible $D(G)$-modules over $E$ as the full subcategory of $\mathrm{Mod}_{D(G)}$ of $D(G)$-modules which are coadmissible as $D(H)$-module for any - equivalently one - compact open subgroup $H$ of $G$ (see \cite[\S 6]{ST03} and the references there). Replacing $\mathrm{Mod}_{D(G)}$ by $\mathrm{Mod}_{D^{\infty}(G)}$, one defines in a similar way the full subcategory $\cC_{D^{\infty}(G)}$ of $\cC_{D(G)}$.\bigskip

Let $G$ be an arbitrary locally $K$-analytic group, then an admissible locally $K$-analytic - or just locally analytic - representation of $G$ over $E$ is a locally $K$-analytic $G$-representation on a $E$-vector space of compact type $V$ such that the strong dual $V^\vee_b$ is in $\cC_{D(G)}$. Here the left $D(G)$-action on $V^\vee$ is given by $(\delta_g\cdot \delta)(x)\defeq \delta(g^{-1}\cdot x)$ for $g\in G$, $x\in V$ and $\delta\in V^\vee$. We write $\mathrm{Rep}^{\rm{an}}_{\rm{adm}}(G)$ for the abelian category of admissible locally analytic representations of $G$ over $E$. The strong dual gives an anti-equivalence $\mathrm{Rep}^{\rm{an}}_{\rm{adm}}(G)\xrightarrow{\sim}\cC_{D(G)}$ (\cite[Thm.~6.3]{ST03}). Via this equivalence, the ``inverse image'' of the full subcategory $\cC_{D^{\infty}(G)}$ recovers the abelian category $\mathrm{Rep}^{\infty}_{\rm{adm}}(G)$ (see \S\ref{subsec: BZ}) of admissible smooth representations of $G$ over $E$ (\cite[Thm.~6.6]{ST03}).\bigskip

Given a locally $K$-analytic group $G$, we write $\mathrm{Ext}_{D(G)}^\bullet(-,-)$ for the extension groups computed in the category $\mathrm{Mod}_{D(G)}$. For $V_0,V_1$ in $\mathrm{Rep}^{\rm{an}}_{\rm{adm}}(G)$, we use the notation
\begin{equation}\label{extdef}
\mathrm{Ext}_G^\bullet(V_0,V_1)\defeq \mathrm{Ext}_{D(G)}^\bullet(V_1^\vee,V_0^\vee).
\end{equation}
We will need the following lemmas.

\begin{lem}\label{erratum}
Let $G$ be a locally $K$-analytic group and $V_0,V_1$ in $\mathrm{Rep}^{\infty}(G)$. Assume either that both $V_0,V_1$ are admissible, or that the map $V_1\rightarrow (V_1^\vee)^\vee$ is an isomorphism (i.e.~$V_1$ is finite dimensional). Then we have isomorphisms $\mathrm{Ext}_{D^{\infty}(G)}^\bullet(V_1^\vee,V_0^\vee)\cong \mathrm{Ext}_{G}^\bullet(V_0,V_1)^{\infty}$.
\end{lem}
\begin{proof}
The first case is proven in Corollary 0.2 of the erratum to \cite{Schr11} (\url{http://math.univ-lyon1.fr/homes-www/schraen/Erratum_GL3.pdf}), we prove the second. First recall that for any $M_0, M_1$ in $\mathrm{Mod}_{D^{\infty}(G)}$ one has a functorial isomorphism
\begin{equation}\label{adjunctM}
\Hom_{D^{\infty}(G)}(M_1, M_0^\vee)\cong \Hom_{D^{\infty}(G)}(M_0, M_1^\vee)
\end{equation}
where $M^\vee=\Hom_E(M,E)$ is the algebraic dual of $M$ with $\delta\in D^{\infty}(G)$ acting on $f\in M^\vee$ by $(\delta\cdot f)(m)\defeq f(\tilde\delta\cdot m)$ where $m\in M$ and $\tilde{}$ is the unique anti-involution on $D^{\infty}(G)$ extending $g\mapsto g^{-1}$ on $G$. In particular for $V_0,V_1$ in $\mathrm{Rep}^{\infty}(G)$ we have a functorial isomorphism $\Hom_{D^{\infty}(G)}(V_1^\vee,V_0^\vee) \cong \Hom_{D^{\infty}(G)}(V_0,(V_1^\vee)^\vee)$ (recall that any smooth representation of $G$ over $E$ is also a module over $D^{\infty}(G)$, see \cite[p.~300]{ST05}). Under the assumption on $V_1$, we thus have a functorial isomorphism in $V_0$:
\begin{equation}\label{adjunct}
\Hom_{D^{\infty}(G)}(V_1^\vee,V_0^\vee) \cong \Hom_{D^{\infty}(G)}(V_0,V_1)=\Hom_{G}(V_0,V_1).
\end{equation}
Take any projective resolution $P^\cdot$ of $V_0$ in $\mathrm{Rep}^{\infty}(G)$ (which exists by \cite[\S 2]{Be92}). Using (\ref{adjunctM}) together with the fact that the functor from $\mathrm{Mod}_{D^{\infty}(G)}$ to $\mathrm{Rep}^{\infty}(G)$ sending a module $M$ to the subspace of its smooth vectors (under the action of $G$) is exact (as follows from \cite[Lemma 1.3]{ST05}), one easily deduces that the algebraic dual ${P^\cdot}^\vee$ is an injective resolution of $V_0^\vee$ in $\mathrm{Mod}_{D^{\infty}(G)}$. It then follows from (\ref{adjunct}) that $\mathrm{Ext}_{D^{\infty}(G)}^i(V_1^\vee,V_0^\vee)\cong \mathrm{Ext}_{G}^i(V_0,V_1)^{\infty}$ for $i\geq 0$.
\end{proof}

\begin{lem}\label{lem: top ML}
Let $I$ be a countable directed index set and $(0\rightarrow A_i\rightarrow B_i\rightarrow C_i\rightarrow 0)_{i\in I}$ a inverse system of strict exact sequences of Banach $E$-spaces with continuous transition maps. Assume that the image of $A_{i'}$ is dense in $A_i$ for each pair of indices $i\leq i'$. Then the inverse limit
\[0\rightarrow \varprojlim_{i\in I}A_i \rightarrow \varprojlim_{i\in I}B_i \rightarrow \varprojlim_{i\in I}C_i \rightarrow 0\]
is a strict exact sequence of Fr\'echet $E$-spaces.
\end{lem}
\begin{proof}
Continuous maps between $E$-Banach spaces are uniformly continuous. So this is a special case of \cite[Chap.~II, \S~3.5, Th.~1]{Bo} (see also \cite[Rem.~13.2.4(i)]{EGAIII}).
\end{proof}

We now introduce some notation largely following \cite[\S 5.5]{OS15}. Let $\mathbf{G}_0$ be a split reductive algebraic group scheme over $\cO_K$ and $\mathbf{T}_0\subseteq \mathbf{G}_0$ a maximal split torus. Let $\mathbf{P}_0\subseteq \mathbf{G}_0$ be a parabolic subgroup scheme that contains $\mathbf{T}_0$, and $\mathbf{P}_0^-$ the opposite parabolic with unipotent radical $\mathbf{N}_0^-$. We write $G_0\defeq \mathbf{G}_0(\cO_K)$, $G\defeq \mathbf{G}_0(K)$, $\fg_0\defeq \mathrm{Lie}({G}_0)$, $\fg\defeq \mathrm{Lie}({G})$ and use similar notation for the other subgroup schemes (e.g.~$P_0$, $P$, $\fp_0$, $\fp$, $T_0$, etc.). Note that $\fg_0$ is an $\cO_K$-lattice in the $K$-vector space $\fg$. We fix an integer $m_0\geq 1$ if $p>2$, $m_0\geq 2$ if $p=2$, and set $\kappa\defeq 1$ if $p>2$, $\kappa\defeq 2$ if $p=2$. By \cite[\S 9.4]{DDMS99}, the $\cO_K$-lattices $p^{m_0}\fg_0\subseteq \fg$, $p^{m_0}\fp_0\subseteq \fp$ and $p^{m_0}\fn_0^-$ are powerful $\Z_p$-Lie algebras, and thus the exponential map $\exp_{G}: \fg \dashrightarrow G$ converges on these $\cO_K$-lattices. Hence we may define $G_1\defeq \exp_G(p^{m_0}\fg_0)$, $P_1\defeq \exp_G(p^{m_0}\fp_0)$ and $N_1^-\defeq \exp_G(p^{m_0}\fn_0^-)$ which are uniform pro-$p$ groups by \cite[Thm.~9.10]{DDMS99}. Since the adjoint action of $G_0$ leaves $\fg_0$ invariant, $G_1$ is normal in $G_0$. Moreover using coordinates ``of the second kind'' (\cite[\S 34]{S11}, \cite[Thm.~4.9]{DDMS99}) one checks that $G_1=N_1^-P_1$, $G_1\cap P_0=P_1$ and $G_1\cap N_0^-=N_1^-$. In particular the embedding $N_1^-\hookrightarrow G_1$ gives a (locally $K$-analytic) section $G_1/P_1\cong N_1^-\hookrightarrow G_1$ of the natural surjection $\mathrm{pr}_1: G_1\twoheadrightarrow G_1/P_1$. Using $G_0/P_0=\bigsqcup_{g\in G_0/G_1P_0}g G_1/P_1$ and choosing a system of representatives of $G_0/G_1P_0$, we extend it to a section $s: G_0/P_0\hookrightarrow G_0$ of $\mathrm{pr}_0: G_0\twoheadrightarrow G_0/P_0$. As $G_0/P_0=G/P$ (by the Iwasawa decomposition of Bruhat-Tits, see e.g.~\cite[Lemma 3.4]{He11}), $s$ determines a section (again denoted by) $s:G/P\hookrightarrow G$ of the surjection $\mathrm{pr}: G\twoheadrightarrow G/P$.\bigskip

We let $\cI\subseteq ]0,1[\cap p^{\bQ}$ be the subset of $r$ satisfying \cite[(5.5.3)]{OS15} and $|\cdot|_r$ the norm on $D(G_1)$ associated to $r\in \cI$ and the canonical $p$-valuation on the uniform pro-$p$ group $G_1$, see \cite[\S\S 2.2.3, 2.2.6]{OS10} (with our $G_1$ and $|\cdot|_r$ denoted by $H$ and $\overline{q}_r$ there, see \cite[\S 2.2.6 Step $3$]{OS10}). We write $D(G_1)_r$ for the completion of $D(G_1)$ with respect to $|\cdot|_r$. Then $D(G_1)_r$ is an $E$-Banach algebra and $D(G_1)\cong \varprojlim_{r\in\cI} D(G_1)_r$ gives the Fr\'echet-Stein structure on $D(G_1)$ (the projective limit is for $r\rightarrow 1$ in $\cI$). As in \cite[(5.5.4), (5.5.5)]{OS15}, we extend the norm $|\cdot|_r$ on $D(G_1)$ to a maximum norm (still denoted) $|\cdot|_r$ on $D(G_0)\cong \bigoplus_{g\in G_0/G_1}\delta_g D(G_1)$. Then $D(G_0)_r\defeq D(G_0)\otimes_{D(G_1)}D(G_1)_r$ is isomorphic to the completion of $D(G_0)$ with respect to $|\cdot|_r$ and $D(G_0)\cong \varprojlim_{r\in\cI} D(G_0)_r$ is a Fr\'echet-Stein structure on $D(G_0)$. Similarly, we equip $D(P_1)$ (resp.~$D(N_1^-)$) with the norm $|\cdot|_r$ attached to the canonical $p$-valuation on $P_1$ (resp.~$N_1^-$), extend it to a norm $|\cdot|_r$ on $D(P_0)$ (resp.~$D(N_0^-)$) and write $D(P_0)_r$ (resp.~$D(N_0^-)_r$) for the corresponding completion. Since the canonical $p$-valuation of $G_1$ restricts to the one of $P_1$ and $N_1^-$, the norm $|\cdot|_r$ on $G_1$ restricts to the norm $|\cdot|_r$ on $P_1$ and $N_1^-$, and thus the norm $|\cdot|_r$ on $G_0$ restricts to the norm $|\cdot|_r$ on $P_0$ and $N_0^-$ (cf.~the discussion before \cite[(5.5.8)]{OS15}). Finally, using $D^{\infty}(G_0)\cong D(G_0)\otimes_{U(\fg)}E$ (cf.~the discussion below Lemma \ref{lem: smooth restriction}), we deduce from \cite[Prop.~3.7]{ST03} (and its proof) that $D^{\infty}(G_0)$ admits the Fr\'echet-Stein structure $D^{\infty}(G_0)\cong \varprojlim_{r\in\cI}D^{\infty}(G_0)_r$ where $D^{\infty}(G_0)_r\defeq D(G_0)_r\otimes_{U(\fg)}E\cong E\otimes_{U(\fg)}D(G_0)_r$.\bigskip

Let $(P_m(G_1))_{m\geq 1}$ be the lower $p$-series of $G_1$ (cf.~\cite[Def.~1.15]{DDMS99}). Following the notation of \cite{OS15}, we write $G_1^{m}\defeq P_{m+1}(G_1)$. By the proof of \cite[Thm.~4.2]{DDMS99} we have that $G_1^{m}$ is a uniform pro-$p$ group and using \cite[Thm.~3.6(iii)]{DDMS99} we have that its $\Z_p$-Lie algebra is $p^m\mathrm{Lie}_{\Z_p}(G_1)$ (in fact $G_1^{m}\simeq G_1^{p^m}\defeq\{x^{p^m},x\in G_1\}$). Similarly, we define $P_1^{m}\defeq P_{m+1}(P_1)$ for $m\geq 0$ and as before we have $P_1^{m}=G_1^m\cap P_1=G_1^{m}\cap P_0$. Let $s=r^{p^m}$ with $s>\frac{1}{p}$ and $s^{\kappa}<p^{-1/(p-1)}$. Following \cite[\S 6]{Schm08}, we write $|\cdot|_s^{(m)}$ for the norm on $D(G_1^m)$ attached to $s$ and the canonical $p$-valuation on $G_1^m$, and $D(G_1^m)_s$ for the corresponding completion. As before $|\cdot|_s^{(m)}$ restricts to the norm on $D(P_1^m)$ defined similarly using the canonical $p$-valuation on $P_1^m$. We have the following result from \cite{Schm08} (refining results by \cite{Fro03} and \cite{Koh07}).

\begin{lem}\label{lem: free mod coset basis}
For $r,m,s$ as above the restriction of the norm $|\cdot|_r$ of $D(G_1)$ to $D(G_1^m)$ is equivalent to $|\cdot|_s^{(m)}$, and $D(G_1)_r$ is a finite free right $D(G_1^m)_s$-module with a basis given by any set of coset representatives for $G_1/G_1^m$.
\end{lem}
\begin{proof}
This is \cite[Prop.~6.2, Cor.~6.4]{Schm08}.
\end{proof}

For $r\in \cI$, $i=0,1$ and a closed subset $C\subseteq G_i$, we write $D(G_i)_{C,r}$ (resp.~$D^{\infty}(G_i)_{C,r}$) for the closure of $D(G_i)_C$ in $D(G_i)_r$ (resp.~of $D^{\infty}(G_i)_C$ in $D^{\infty}(G_i)_r$). As $U(\fg)$ is dense in $D(G_1)_1$ (\cite[Prop.~1.2.8]{Koh07}), in the particular case $C=\{1\}$ we see that $D(G_1)_{1,r}$ is also the closure of $U(\fg)$ in $D(G_1)_r$.

\begin{lem}\label{lem: envelop dense}
If $r^{\kappa}<p^{-1/(p-1)}$, then $U(\fg)$ is dense in $D(G_1)_r$. In particular $D(G_1)_{1,r}=D(G_1)_r$ and $D(G_1)_{1,r}\fg=\fg D(G_1)_{1,r}=\fg D(G_1)_r=D(G_1)_r\fg$ is the (two-sided) augmentation ideal of $D(G_1)_r$.
\end{lem}
\begin{proof}
This is \cite[(5.5.6)]{OS15}, which follows from \cite[Prop.~5.6]{Schm08}.
\end{proof}

The following lemma will be used many times in the sequel.

\begin{lem}\label{lem: support basis}
Let $r\in\cI$ and $m\geq 0$ such that $s=r^{p^m}$ satisfies $s>p^{-1}$ and $s^{\kappa}<p^{-1/(p-1)}$.
\begin{enumerate}[label=(\roman*)]
\item \label{it: support basis 1} For each closed subset $C\subseteq G_0$, $D(G_0)_{C,r}$ (resp.~$D^{\infty}(G_0)_{C,r}$) is a finite free right $D(G_0)_{1,r}$-module (resp.~a finite dimensional $E$-vector space) with a basis given by any set of coset representatives of $CG_1^m/G_1^m$.
\item \label{it: support basis 2} Let $H_0\subseteq G_0$ be a closed subgroup, for each closed subset $C\subseteq G_0$ such that $CH_0=C$ in $G_0$, $D(G_0)_{C,r}$ (resp.~$D^{\infty}(G_0)_{C,r}$) is a finite free right $D(G_0)_{H_0,r}$-module (resp.~a finite free right $D^{\infty}(G_0)_{H_0,r}$-module) with a basis given by any set of coset representatives of $CG_1^m/H_0G_1^m$.
\end{enumerate}
\end{lem}
\begin{proof}
We fix $r$, $m$ and $s=r^{p^m}$ as in the statement. As $s^{\kappa}<p^{-1/(p-1)}$, by Lemma~\ref{lem: envelop dense} $D(G_1^m)_{1,s}=D(G_1^m)_s$. As $G_1^m$ is compact open in $G_1$, we have $D(G_1^m)\buildrel\sim\over\rightarrow D(G_1)_{G_1^m}$ and $D(G_1^m)_1\buildrel\sim\over\rightarrow D(G_1)_1$ and likewise with $D(G_0)$ instead of $D(G_1)$. Together with Lemma~\ref{lem: free mod coset basis} this implies $D(G_1)_{G_1^m,r}=D(G_1)_{1,r}=D(G_1^m)_{1,s}=D(G_1^m)_s$ (likewise with $D(G_0)$ instead of $D(G_1)$) and
\[D(G_1)_r\cong \bigoplus_{g\in G_1/G_1^m}\delta_g D(G_1)_{G_1^m,r}\cong \bigoplus_{g\in G_1/G_1^m}D(G_1)_{gG_1^m,r},\]
which together with $D(G_0)_r\cong \bigoplus_{g\in G_0/G_1}\delta_g D(G_1)_r$ implies
\begin{equation}\label{equ: support basis}
D(G_0)_r\cong \bigoplus_{g\in G_0/G_1^m}\delta_g D(G_0)_{G_1^m,r}\cong \bigoplus_{g\in G_0/G_1^m}D(G_0)_{gG_1^m,r}.
\end{equation}
By Lemma~\ref{lem: envelop dense} (and the above discussion) $D(G_1^m)_{1,s} \fg =D(G_1^m)_s \fg=D(G_0)_{G_1^m,r} \fg$ is the augmentation ideal of $D(G_1^m)_{1,s} =D(G_1^m)_s=D(G_0)_{G_1^m,r}$, which together with the first equality in (\ref{equ: support basis}) implies
\begin{equation}\label{dinfinityr}
D^{\infty}(G_0)_r=D(G_0)_r/D(G_0)_r \fg\cong E[G_0/G_1^m].
\end{equation}

Now we fix a closed subset $C\subseteq G_0$ and note that $D(G_0)\cong \bigoplus_{g\in G_0/G_1^m}\delta_g D(G_1^m)\cong \bigoplus_{g\in G_0/G_1^m}D(G_0)_{gG_1^m}$ implies $D(G_0)_{C}\cong \bigoplus_{g\in G_0/G_1^m}D(G_0)_{C\cap gG_1^m}$ (see e.g.~the proof of \cite[Lemma 1.2.5]{Koh07}) which in turn implies
\begin{equation}\label{equ: support coset decomp}
D(G_0)_{C,r}=\bigoplus_{g\in G_0/G_1^m}D(G_0)_{C\cap gG_1^m,r}.
\end{equation}

\textbf{Step $1$}: We prove $D(G_0)_{C\cap gG_1^m,r}\!=\!D(G_0)_{gG_1^m,r}$ if $g\!\in \!CG_1^m$ and $D(G_0)_{C\cap gG_1^m,r}\!=\!0$ otherwise.\\
Note that $g\in CG_1^m$ if and only if $C\cap gG_1^m\neq \emptyset$. If $C\cap gG_1^m=\emptyset$, we obviously have $D(G_0)_{C\cap gG_1^m,r}=0$. If $C\cap gG_1^m\neq \emptyset$, for any $h\in C\cap gG_1^m$ we have closed embeddings of $E$-Banach spaces
\[D(G_0)_{gG_1^m,r}\supseteq D(G_0)_{C\cap gG_1^m,r}\supseteq D(G_0)_{\{h\},r}.\]
Writing $D(G_0)_{\{h\},r}=\delta_h D(G_0)_{1,r}=\delta_h D(G_0)_{G_1^m,r}=D(G_0)_{hG_1^m,r}=D(G_0)_{gG_1^m,r}$ (where the last equality follows from $hG_1^m=gG_1^m$), we deduce $D(G_0)_{\{h\},r}=D(G_0)_{C\cap gG_1^m,r}=D(G_0)_{gG_1^m,r}$.\bigskip

\textbf{Step $2$}: We prove \ref{it: support basis 1} and \ref{it: support basis 2}.\\
We combine Step $1$ with (\ref{equ: support coset decomp}) and first deduce
\begin{equation}\label{equ: support full coset decomp}
D(G_0)_{C,r}=\!\!\bigoplus_{g\in CG_1^m/G_1^m}D(G_0)_{C\cap gG_1^m,r}=\bigoplus_{g\in CG_1^m/G_1^m}\!\!\delta_g D(G_0)_{G_1^m,r}=\bigoplus_{g\in CG_1^m/G_1^m}\!\!\delta_g D(G_0)_{1,r}
\end{equation}
which proves the first statement in \ref{it: support basis 1}. Define $\widetilde D^{\infty}(G_0)_{C,r}\defeq D(G_0)_{C,r}\otimes_{U(\fg)}E$ (noting that $D(G_0)_{C,r}$ is a right $D(G_0)_{1,r}$-module) and note that $D^{\infty}(G_0)_{C,r}$ is the closure of the image of $\widetilde D^{\infty}(G_0)_{C,r}$ in $D^{\infty}(G_0)_r$. Applying $(-)\otimes_{U(\fg)}E$ to (\ref{equ: support full coset decomp}), we deduce
\begin{equation*}
\widetilde D^{\infty}(G_0)_{C,r}=E[CG_1^m/G_1^m],
\end{equation*}
which together with (\ref{dinfinityr}) shows that $\widetilde D^{\infty}(G_0)_{C,r}$ is already closed in $D^{\infty}(G_0)_r$. This implies $\widetilde D^{\infty}(G_0)_{C,r}\buildrel\sim\over\rightarrow D^{\infty}(G_0)_{C,r}$ and also gives the second statement of \ref{it: support basis 1}. Note that $G_1^m$ is normal in $G_0$ and thus $G_1^mH_0=H_0G_1^m$ is a compact open subgroup of $G_0$. If furthermore $CH_0=C$, then we may rewrite (\ref{equ: support full coset decomp}) as
\begin{multline*}
D(G_0)_{C,r}=\bigoplus_{g\in CG_1^m/G_1^m}\delta_g D(G_0)_{G_1^m,r}=\bigoplus_{g\in CG_1^m/H_0G_1^m,h\in H_0G_1^m/G_1^m} \delta_g\delta_h D(G_0)_{G_1^m,r}\\
=\bigoplus_{g\in CG_1^m/H_0G_1^m} \delta_g D(G_0)_{H_0G_1^m,r}=\bigoplus_{g\in CG_1^m/H_0G_1^m} \delta_g D(G_0)_{H_0,r}
\end{multline*}
where the third equality follows from (\ref{equ: support full coset decomp}) applied with $C=H_0G_1^m$ and the last equality follows from \ref{it: support basis 1} applied with $C=H_0$ and $C=H_0G_1^m$ (noting that $H_0G_1^m/G_1^m=(H_0G_1^m)G_1^m/G_1^m$). Applying $(-)\otimes_{U(\fg)}E$ and using $\widetilde D^{\infty}(G_0)_{C,r}\buildrel\sim\over\rightarrow D^{\infty}(G_0)_{C,r}$ (for both $C$ and $H_0$, see Step $1$) we deduce
\[D^{\infty}(G_0)_{C,r}=\bigoplus_{g\in CG_1^m/H_0G_1^m} \delta_g D^{\infty}(G_0)_{H_0,r}\]
which finishes the proof of \ref{it: support basis 2}.
\end{proof}

The following result is well-known, we provide a proof for lack of a precise reference.

\begin{lem}\label{lem: Frechet Stein}
Let $H_0\subseteq G_0$ be a closed subgroup. The closed subalgebra $D(G_0)_{H_0}$ of $D(G_0)$ is Fr\'echet-Stein via $D(G_0)_{H_0}\cong \varprojlim_{r\in\cI}D(G_0)_{H_0,r}$.
\end{lem}
\begin{proof}
We consider $r,r'\in\cI$ with $r\leq r'$. We note that the image of $D(G_0)_{H_0,r'}$ is dense in $D(G_0)_{H_0,r}$ (both contain as a dense subspace the image of $D(G_0)_{H_0}$). It suffices to check that $D(G_0)_{H_0,r}$ is flat over $D(G_0)_{H_0,r'}$. By \ref{it: support basis 2} of Lemma~\ref{lem: support basis} (applied with $C=G_0$), we know that $D(G_0)_r$ (resp.~$D(G_0)_{r'}$) is finite free as a $D(G_0)_{H_0,r}$-module (resp.~$D(G_0)_{H_0,r'}$-module) with $D(G_0)_{H_0,r}$ being a direct summand. Since $D(G_0)_r$ is flat over $D(G_0)_{r'}$ (cf.~\cite[Thm.~5.1]{ST03} which uses \cite[Thm.~4.10, Prop.~3.7]{ST03}) and $D(G_0)_{r'}$ is free as a $D(G_0)_{H_0,r'}$-module, we see that $D(G_0)_r$ is flat over $D(G_0)_{H_0,r'}$. Since $D(G_0)_{H_0,r}$ is a direct summand of $D(G_0)_r$ as $D(G_0)_{H_0,r'}$-modules, this implies that $D(G_0)_{H_0,r}$ is also flat over $D(G_0)_{H_0,r'}$.
\end{proof}

Let $H\subseteq G$ be a closed subgroup, we say that a left $D(G)_H$-module is coadmissible if it is coadmissible as a $D(G_0)_{H\cap G_0}$-module with the Fr\'echet-Stein structure given by Lemma \ref{lem: Frechet Stein}. If $G'_0\subseteq G_0$ is another compact open subgroup, then $D(G_0)_{H\cap G_0}$ is free of finite rank over $D(G'_0)_{H\cap G'_0}$ (see the argument for \cite[(1.7)]{Koh07}) and arguing as in \cite[\S 6]{ST05} we see that being a coadmissible $D(G)_H$-module doesn't depend on the compact open subgroup $G_0$. We let $\cC_{D(G)_H}$ be the abelian category of coadmissible (left) $D(G)_H$-modules. Recall that each coadmissible $D(G)_H$-module carries a canonical Fr\'echet topology, and any $D(G)_H$-linear map between coadmissible $D(G)_H$-modules is continuous and strict. See \cite[\S 7.7]{AS22} for a similar discussion (also) based on \cite[Corollary~1.4.3]{Koh07}.\bigskip

Let $H_0\subseteq G_0$ be a closed subgroup, as $D(G_0)_{H_0}$ is a Fr\'echet-Stein $U(\fg)$-module, it follows from \cite[Prop.~3.7]{ST03} and its proof that $D(G_0)_{H_0}\otimes_{U(\fg)}E$ is also Fr\'echet-Stein with $D(G_0)_{H_0}\otimes_{U(\fg)}E\buildrel\sim\over\rightarrow \varprojlim_{r\in \cI} (D(G_0)_{H_0,r}\otimes_{U(\fg)}E)$. Since $D(G_0)_{H_0,r}\otimes_{U(\fg)}E\buildrel\sim\over \rightarrow D^\infty(G_0)_{H_0,r}\subseteq D^\infty(G_0)_{r}$ (see Step $2$ in the proof of Lemma \ref{lem: support basis}), we have $D(G_0)_{H_0}\otimes_{U(\fg)}E\buildrel\sim\over \rightarrow D^{\infty}(G_0)_{H_0}\subseteq D^\infty(G_0)$ (recall $D^{\infty}(G_0)_{H_0}$ is by definition the image of $D(G_0)_{H_0}$ in $D^\infty(G_0)$). This statement can in fact be generalized to any closed subset $C$ in $G_0$.

\begin{lem}\label{lem: pass to smooth support}
For any closed subset $C\subseteq G_0$ we have a topological isomorphism
\begin{equation}\label{equ: pass to smooth support}
D(G_0)_C\otimes_{U(\fg)}E \buildrel\sim\over\longrightarrow D(G_0)_C\otimes_{D(G_0)_1}E \buildrel\sim\over\longrightarrow D^{\infty}(G_0)_C.
\end{equation}
\end{lem}
\begin{proof}
From the definitions each map is surjective and continuous, hence it is enough to prove injectivity of the composition. Since $D^{\infty}(G_0)_C$ is closed in $D^{\infty}(G_0)$ (Lemma \ref{lem: closedinfini}), we have $D^{\infty}(G_0)_C\cong \varprojlim_{r\in\cI}D^{\infty}(G_0)_{C,r}\subseteq \varprojlim_{r\in\cI}D^{\infty}(G_0)_{r}\cong D^{\infty}(G_0)$. Since $D(G_0)_{C,r}\otimes_{U(\fg)}E\buildrel\sim\over\rightarrow D^{\infty}(G_0)_{C,r}$ (Step $2$ in the proof of Lemma \ref{lem: support basis}) we have exact sequences $0\rightarrow D(G_0)_{C,r}\fg \rightarrow D(G_0)_{C,r} \rightarrow D^{\infty}(G_0)_{C,r}\rightarrow 0$ and since the maps $D(G_0)_{C,r'}\fg\rightarrow D(G_0)_{C,r}\fg$ for $r\leq r'$ have dense image (as the image of $D(G_0)_C\fg$ is dense everywhere by construction), Lemma \ref{lem: top ML} gives a short exact sequence of Fr\'echet $E$-spaces
\begin{equation}\label{newses}
0\longrightarrow \varprojlim_{r\in\cI} (D(G_0)_{C,r}\fg) \longrightarrow D(G_0)_C \longrightarrow D^{\infty}(G_0)_{C}\longrightarrow 0.
\end{equation}
Since $D(G_0)_1\buildrel\sim\over\rightarrow \varprojlim_{r\in\cI}D(G_0)_{1,r}$ is a Fr\'echet-Stein algebra (e.g.~by Lemma \ref{lem: Frechet Stein}), both $D(G_0)_{1}\otimes_E \fg$ and $D(G_0)_{1}\fg$ are coadmissible $D(G_0)_1$-modules (by \cite[Cor.~3.4.iv]{ST03} for the latter). By \cite[Cor.~3.4.ii]{ST03} we have in particular short exact sequences of finitely generated $D(G_0)_{1,r}$-modules $0\rightarrow M_r \rightarrow D(G_0)_{1,r}\otimes_E \fg \rightarrow D(G_0)_{1,r}\fg\rightarrow 0$ where the maps $M_{r'}\rightarrow M_{r}$ have dense image. By \ref{it: support basis 1} of Lemma \ref{lem: support basis} $D(G_0)_{C,r}$ is a free $D(G_0)_{1,r}$-module of finite rank (this rank growing when $r$ tends to $1$), hence tensoring by $D(G_0)_{C,r}$ over $D(G_0)_{1,r}$ again gives short exact sequences of Banach spaces
\[0\longrightarrow D(G_0)_{C,r}\otimes_{D(G_0)_{1,r}} M_r \longrightarrow D(G_0)_{C,r}\otimes_E \fg \longrightarrow D(G_0)_{C,r}\fg\longrightarrow 0\]
where the maps $D(G_0)_{C,r'}\otimes_{D(G_0)_{1,r'}} M_{r'}\rightarrow D(G_0)_{C,r}\otimes_{D(G_0)_{1,r}} M_r$ have dense image. By Lemma \ref{lem: top ML} (and since $\fg$ is finite dimensional over $E$) we deduce in particular a surjection $D(G_0)_{C}\otimes_E \fg \twoheadrightarrow \varprojlim_{r\in\cI} (D(G_0)_{C,r}\fg)$. By (\ref{newses}) this implies $D(G_0)_C\otimes_{U(\fg)}E\cong D(G_0)_C/D(G_0)_C\fg \buildrel\sim\over\rightarrow D^{\infty}(G_0)_C$.
\end{proof}

Recall the surjections $\mathrm{pr}_0: G_0\twoheadrightarrow G_0/P_0$ and $\mathrm{pr}: G\twoheadrightarrow G/P$.
For each subset $C\subseteq G_0$ (resp.~$C\subseteq G$), we have $CP_0=\mathrm{pr}_0^{-1}(\mathrm{pr}_0(C))$ (resp.~$CP=\mathrm{pr}^{-1}(\mathrm{pr}(C))$).
As the inclusion $G_0\subseteq G$ induces an isomorphism $G_0/P_0\buildrel\sim\over\rightarrow G/P$, we have $C\cap G_0=\mathrm{pr}_0^{-1}\mathrm{pr}(C)$ for each subset $C\subseteq G$ such that $CP=C$.
The fixed section $s: G_0/P_0\hookrightarrow G_0$ of $\mathrm{pr}_0$ induces an isomorphism $G_0/P_0\times P_0\buildrel\sim\over\longrightarrow G_0$ and thus an isomorphism $C_0/P_0\times P_0\buildrel\sim\over\longrightarrow C_0$ for each closed subset $C_0\subseteq G_0$ such that $C_0P_0=C_0$. We deduce a topological isomorphism for such closed subsets (by applying (\ref{equ: loc an product}) with $M_1=G_0/P_0$, $C_1=C_0/P_0$ and $M_2=C_2=P_0$)
\begin{equation}\label{equ: loc an coset section}
D(G_0/P_0)_{C_0/P_0}\widehat{\otimes}_E D(P_0) \cong D(G_0)_{C_0}.
\end{equation}
Similarly, using (\ref{equ: smooth product}) $s$ also induces a topological isomorphism
\begin{equation}\label{equ: smooth coset section}
D^{\infty}(G_0/P_0)_{C_0/P_0}\widehat{\otimes}_E D^\infty(P_0) \cong D^{\infty}(G_0)_{C_0}.
\end{equation}

We now fix a locally closed subset $X\subseteq G$ (i.e.~$X$ is the intersection of a closed subset and an open subset of $G$, equivalently $X$ is open in its closure $\overline X\subseteq G$) such that $XP=X$ and we set $X_0\defeq X\cap G_0$ and $Y\defeq \overline X$. Note that $X_0$ is also open in its closure $Y_0\defeq \overline X_0\cong Y\cap G_0$ and that $X_0=X_0P_0$. A compact open subset of $X_0$ is the same thing as a compact open subset of $Y_0$ which is contained in $X_0$, in particular it is a compact subset of $G_0$ of the form $Y_0\cap U_0$ where $U_0$ is a compact open subset of $G_0$ (of which there are only countably many). For each compact open subset $C_0$ of $X_0$, recall that $D(G_0)_{C_0}$ is a Fr\'echet space equipped with a separately continuous right action of $D(G_0)_{P_0}$. If $C_0\subseteq C'_0$ are two compact open subsets of $X_0$, one can easily find compact open subsets $U_0\subseteq U'_0$ of $G_0$ such that $C_0=X_0\cap U_0=Y_0\cap U_0$ and $C'_0=X_0\cap U'_0=Y_0\cap U'_0$, in particular $C_0=C'_0\cap U_0$, $D(G_0)_{C_0}=D(U_0)_{C_0}$ and $D(G_0)_{C'_0}=D(U'_0)_{C'_0}$. Writing $U_0=U_0 \amalg U'_0\!\setminus \!U_0$ and noting that $U'_0\!\setminus \!U_0$ is also a compact open subset of $G_0$, we have $D(U'_0)_{C'_0}\cong D(U_0)_{C_0}\oplus D(U'_0\!\setminus \! U_0)_{C'_0\cap U'_0\setminus U_0}$, and hence deduce a canonical projection of Fr\'echet spaces $D(G_0)_{C'_0}\twoheadrightarrow D(G_0)_{C_0}$. We then consider the projective limit
\begin{equation}\label{projc0}
\widehat{D}(G_0)_{X_0}\defeq \varprojlim_{C_0}D(G_0)_{C_0}
\end{equation}
over the compact open subsets $C_0$ of $X_0$ such that $C_0P_0=C_0$ with transition maps given by the above projections for $C_0\subseteq C'_0$. This is still a Fr\'echet space equipped with a separately continuous right $D(G_0)_{P_0}$-action and we have
\begin{equation}\label{equ: support limit}
\widehat{D}(G_0)_{X_0}\cong \varprojlim_{C_0}\big(D(G_0/P_0)_{C_0/P_0}\widehat{\otimes}_E D(P_0)\big)
\cong \big(\varprojlim_{C_0} D(G_0/P_0)_{C_0/P_0}\big)\widehat{\otimes}_E D(P_0)
\end{equation}
where the first isomorphism comes from (\ref{equ: loc an coset section}) and the second from \cite[Prop.~1.1.29]{Em17}. We also deduce a right $D(G)_P$-action on
\begin{multline}\label{equ: loc an hat}
\widehat{D}(G)_X\defeq \widehat{D}(G_0)_{X_0}\otimes_{D(G_0)_{P_0}}D(G)_P \cong \widehat{D}(G_0)_{X_0}\otimes_{D(P_0)}D(P)\\
 \cong (\varprojlim_{C_0}D(G_0/P_0)_{C_0/P_0})\widehat{\otimes}_{E,\iota}D(P)
\end{multline}
where the first isomorphism follows from $D(G)_P\cong D(G_0)_{P_0}\otimes_{D(P_0)}D(P)$ (see the argument for \cite[(1.7)]{Koh07}) and the second from (\ref{equ: support limit}) and the same argument as at the end of the proof of \cite[Prop.~A.3]{ST05}. Similarly, we consider the right $D^{\infty}(G_0)_{P_0}\cong D^{\infty}(P_0)$-module (see Lemma \ref{lem: smooth restriction} for the latter isomorphism)
\[\widehat{D}^{\infty}(G_0)_{X_0}\defeq \varprojlim_{C_0}D^{\infty}(G_0)_{C_0}\cong (\varprojlim_{C_0}D^{\infty}(G_0/P_0)_{C_0/P_0})\widehat{\otimes}_E D^{\infty}(P_0)\]
where we have used (\ref{equ: smooth coset section}) for the last isomorphism, and the right $D^{\infty}(G)_P\cong D^{\infty}(P)$-module
\begin{equation}\label{equ: smooth hat}
\widehat{D}^{\infty}(G)_X\defeq \widehat{D}^{\infty}(G_0)_{X_0} \otimes_{D^{\infty}(P_0)} D^{\infty}(P)\cong (\varprojlim_{C_0}D^{\infty}(G_0/P_0)_{C_0/P_0})\widehat{\otimes}_{E,\iota}D^{\infty}(P)
\end{equation}
(note that $D^\infty(P)\cong D^{\infty}(P_0)\otimes_{D(P_0)}D(P)$).

\begin{lem}
For any compact open subset $C_0$ of $X_0$ such that $C_0=C_0P_0$ we have topological isomorphisms
\begin{eqnarray}
\label{equ: pass to smooth limit}
\widehat{D}(G_0)_{X_0}\otimes_{U(\fg)}E &\buildrel\sim\over\longrightarrow &\widehat{D}(G_0)_{X_0}\otimes_{D(G_0)_1}E \buildrel\sim\over\longrightarrow \widehat{D}^{\infty}(G_0)_{X_0}\\
\label{equ: hat loc an to smooth}
\widehat{D}(G)_{X}\otimes_{U(\fg)}E &\buildrel\sim\over\longrightarrow &\ \widehat{D}(G)_{X}\otimes_{D(G)_1}E \ \ \buildrel\sim\over\longrightarrow \ \widehat{D}^{\infty}(G)_{X}.
\end{eqnarray}
\end{lem}
\begin{proof}
The proof of (\ref{equ: pass to smooth limit}) is analogous to the proof of Lemma \ref{lem: pass to smooth support} using twice Lemma \ref{lem: top ML} (first with $\varprojlim_{r}$ for a given $C_0$ as in the proof of \emph{loc.~cit.}, then with $\varprojlim_{C_0}$) and noting that all transition maps always have dense image. Writing $D(P)=\bigoplus_{h\in P_0\backslash P}D(P_0)\delta_h$, we have and by (\ref{equ: loc an hat}) and (\ref{equ: smooth hat})
\begin{equation}\label{P0toP-}
\begin{gathered}
\widehat{D}(G)_X\cong \bigoplus_{h\in P_0\backslash P}\!\widehat{D}(G_0)_{X_0}\delta_h\cong \widehat{D}(G_0)_{X_0}\otimes_{D(P_0)}D(P)\\
\widehat{D}^\infty(G)_X\cong \bigoplus_{h\in P_0\backslash P}\!\widehat{D}^\infty(G_0)_{X_0}\delta_h\cong \widehat{D}^\infty(G_0)_{X_0}\otimes_{D(P_0)}D(P)
\end{gathered}
\end{equation}
from which we easily deduce (\ref{equ: hat loc an to smooth}) using (\ref{equ: pass to smooth limit}).
\end{proof}

As the morphisms in (\ref{equ: pass to smooth support}) are right $D(G_0)_{P_0}$-equivariant, so are the morphisms in (\ref{equ: pass to smooth limit}), and thus the morphisms (\ref{equ: hat loc an to smooth}) are also right $D(G)_P$-equivariant by definition. Note that the topological isomorphisms $G_0/P_0 \times P_0 \buildrel\sim\over\rightarrow G_0$, $G/P \times P \rightarrow G$, $G_0/P_0\buildrel\sim\over\rightarrow G/P$ (given by the section $s$ and the Iwasawa decomposition) induce topological isomorphisms $X_0/P_0 \times P_0 \buildrel\sim\over\rightarrow X_0$, $X/P \times P \rightarrow X$, $X_0/P_0\buildrel\sim\over\rightarrow X/P$. If $X$ is moreover such that $QX=X$ for some locally $K$-analytic closed subgroup $Q\subseteq G$, we then have a natural separately continuous left $D(Q)$-action on $\varprojlim_{C}D(G/P)_{C}$ and $\varprojlim_{C}D^{\infty}(G/P)_{C}$, where $C$ runs among the compact open subsets of $X/P=X_0/P_0$ and the transition maps are defined as previously, which is uniquely determined by ($h\in Q$, $(\delta_C)_{C}\in \varprojlim D(G/P)_{C}$ or $\varprojlim D^{\infty}(G/P)_{C}$):
\[\delta_{h}\cdot(\delta_C)_{C}=(\delta_{h}\cdot\delta_C)_{h C}.\]
Via the last isomorphism in (\ref{equ: loc an hat}) (resp.~(\ref{equ: smooth hat})), we deduce a separately continuous left $D(Q)$-action on $\widehat{D}(G)_X$ (resp.~$\widehat{D}^{\infty}(G)_X$), which makes (\ref{equ: pass to smooth limit}) and (\ref{equ: hat loc an to smooth}) left $D(Q)$-equivariant.

\begin{lem}\label{lem: limit vers quotient}
Let $X\subseteq G$ be a locally closed subset such that $XP=X$. We have canonical right topological $D(G)_P$-equivariant isomorphisms
\begin{eqnarray}
\label{equ: limit vers quotient}
D(G)_{\overline{X}}/D(G)_{\overline{X}\setminus X}&\buildrel\sim\over\longrightarrow &\widehat{D}(G)_X\\
\label{equ: limit vers quotient smooth}
D^{\infty}(G)_{\overline{X}}/D^{\infty}(G)_{\overline{X}\setminus X} &\buildrel\sim\over\longrightarrow &\widehat{D}^{\infty}(G)_X.
\end{eqnarray}
If moreover $QX=X$ for some locally $K$-analytic subgroup $Q\subseteq G$, then (\ref{equ: limit vers quotient}) and (\ref{equ: limit vers quotient smooth}) are $D(Q)$-equivariant for the natural left $D(Q)$-action on both sides.
\end{lem}
\begin{proof}
We only give the proofs for (\ref{equ: limit vers quotient}), leaving the case of (\ref{equ: limit vers quotient smooth}) to the reader (arguing as in the proof of (\ref{equ: hat loc an to smooth}) above). As before we note $Y=\overline{X}$, $X_0= X\cap G_0$, $Y_0= Y\cap G_0$, and we also define $Z\defeq \overline{X}\setminus X$ and $Z_0\defeq Z\cap G_0$. Hence $Z_0\subseteq Y_0$ are closed subspaces of the compact group $G_0$ with $X_0=Y_0\setminus Z_0$ and we have $YP=Y$, $ZP=Z$, $Y_0P_0=Y_0$, $Z_0P_0=Z_0$. Each compact open subset $C_0\subseteq X_0$ is compact open in $Y_0$, and writing $Y_0=C_0 \amalg Y_0\!\setminus \!C_0$ we have as in (\ref{projc0}) a surjection of Fr\'echet spaces $D(G_0)_{Y_0}\twoheadrightarrow D(G_0)_{C_0}$ with kernel $D(G_0)_{Y_0\!\setminus \!C_0}$ containing $D(G_0)_{Z_0}$ as closed subspace. Taking the projective limit over those $C_0$ such that $C_0P_0=C_0$ we deduce a canonical morphism of Fr\'echet spaces
\begin{equation}\label{equ: limit as quotient prime}
D(G_0)_{Y_0}/D(G_0)_{Z_0}\longrightarrow \varprojlim_{C_0}D(G_0)_{C_0}=\widehat{D}(G_0)_{X_0}.
\end{equation}

\textbf{Step $1$}: We prove that (\ref{equ: limit as quotient prime}) is a topological isomorphism.\\
By \ref{it: support basis 2} of Lemma~\ref{lem: support basis}, $D(G_0)_{Z_0,r}$, $D(G_0)_{Y_0,r}$ and $D(G_0)_{Y_0,r}/D(G_0)_{Z_0,r}$ are finite free right $D(G_0)_{P_0,r}$-modules with a basis given by $\{\delta_g\}$ with $g$ running through coset representatives of $Z_0G_1^m/P_0G_1^m$, $Y_0G_1^m/P_0G_1^m$ and $(Y_0G_1^m\!\setminus \!Z_0G_1^m)/P_0G_1^m$ respectively (using the notation of \emph{loc.~cit.}). Since $(Y_0G_1^m\!\setminus \!Z_0G_1^m)/P_0G_1^m=(Y_0\!\setminus \!(Z_0G_1^m))G_1^m/P_0G_1^m$ (recall $G_1^mP_0=P_0G_1^m$) it follows from \emph{loc.~cit.} again that we have short exact sequence of $E$-Banach spaces (which are finite free right $D(G_0)_{P_0,r}$-modules)
\begin{equation}\label{equ: Banach seq}
0\rightarrow D(G_0)_{Z_0,r} \rightarrow D(G_0)_{Y_0,r}\rightarrow D(G_0)_{Y_0\setminus (Z_0G_1^m),r}\rightarrow 0.
\end{equation}
Applying Lemma~\ref{lem: top ML}, we obtain a short exact sequence of $E$-Fr\'echet spaces
\begin{equation}\label{equ: limit Frechet seq}
0\rightarrow D(G_0)_{Z_0} \rightarrow D(G_0)_{Y_0}\rightarrow \varprojlim_{r\in \cI}D(G_0)_{Y_0\setminus (Z_0G_1^m),r}\rightarrow 0.
\end{equation}
As $Z_0$ is compact, the $(G_1^m)_m$ form a system of neighborhood of $1$ in $G_0$ when $r\rightarrow 1$ in $\cI$ and as $G_1^mP_0=P_0G_1^m$, it is easy to check that $(Y_0\setminus (Z_0G_1^m))_m$ is cofinal among compact open subsets of $X_0=Y_0\!\setminus \!Z_0$ stable under right multiplication by $P_0$, and thus
\begin{equation}\label{equ: Banach seqlim}
\varprojlim_{r}D(G_0)_{Y_0\setminus (Z_0G_1^m),r}\cong \varprojlim_{C_0}\varprojlim_{r}D(G_0)_{C_0,r}\cong \varprojlim_{C_0}D(G_0)_{C_0}\cong \widehat{D}(G_0)_{X_0}.
\end{equation}

\textbf{Step $2$}: We prove that (\ref{equ: limit as quotient prime}) induces a topological isomorphism (\ref{equ: limit vers quotient}) which is $D(G)_P$-equivariant.\\
Writing $D(P)=\bigoplus_{h\in P_0\backslash P}D(P_0)\delta_h$, we have by the argument at the end of the proof of \cite[Lemma 1.2.5]{Koh07} for $S\in \{Z,Y\}$
\begin{equation}\label{P0toP}
D(G)_S\cong \bigoplus_{h\in P_0\backslash P}\!D(G_0)_{S_0}\delta_h\cong D(G_0)_{S_0}\otimes_{D(P_0)}D(P).
\end{equation}
Together with the first line in (\ref{P0toP-}), it follows that (\ref{equ: limit as quotient prime}) induces a topological isomorphism (\ref{equ: limit vers quotient}). As the projection $D(G_0)_{Y_0}\twoheadrightarrow D(G_0)_{C_0}$ is right $D(G_0)_{P_0}$-equivariant for each compact open $C_0\subseteq X_0$ such that $C_0P_0=C_0$, we deduce that (\ref{equ: limit as quotient prime}) is right $D(G_0)_{P_0}$-equivariant, and therefore (\ref{equ: limit vers quotient}) is right $D(G)_P\cong D(G_0)_{P_0}\otimes_{D(P_0)}D(P)$-equivariant.\bigskip

\textbf{Step $3$}: If furthermore $QX=X$, we prove that (\ref{equ: limit vers quotient}) is left $D(Q)$-equivariant.\\
As $X$ is left $Q$-stable, $Y$ and $Z$ are also left $Q$-stable, so we have natural left $D(Q)$-action on $D(G)_Y$ and $D(G)_Z$. Let $C_0\subseteq X_0$ be a compact open subset such that $C_0P_0=C_0$, the above short exact sequence $0\rightarrow D(G_0)_{Y_0\!\setminus \!C_0} \rightarrow D(G_0)_{Y_0}\rightarrow D(G_0)_{C_0}\rightarrow 0$ induces a short exact sequence
\begin{multline}\label{equ: limit vers quotient projection}
0\longrightarrow D(G_0)_{Y_0\!\setminus \!C_0}\otimes_{D(P_0)}D(P)\longrightarrow D(G)_Y\cong D(G_0)_{Y_0}\otimes_{D(P_0)}D(P)\\
\longrightarrow D(G_0)_{C_0}\otimes_{D(P_0)}D(P)\cong D(G)_{C_0P}\longrightarrow 0
\end{multline}
where the kernel contains the closed subspace $D(G)_Z\cong D(G_0)_{Z_0}\otimes_{D(P_0)}D(P)$ of $D(G)_Y$ and where the last isomorphism follows again from the proof of \cite[Lemma 1.2.5]{Koh07}. We thus have a continuous morphism $D(G)_Y/D(G)_Z\rightarrow D(G)_{C_0P}$. Using $Z_0=\bigcap_{C_0}Y_0\!\setminus \!C_0$, and hence $D(G_0)_{Z_0}=\cap_{C_0}D(G_0)_{Y_0\!\setminus \!C_0}\cong \varprojlim_{C_0}D(G_0)_{Y_0\!\setminus \!C_0}$, we deduce from (\ref{equ: limit vers quotient projection}) an embedding
\begin{equation}\label{equ: naive limit}
D(G)_Y/D(G)_Z \hookrightarrow \varprojlim_{C_0}D(G)_{C_0P}.
\end{equation}
Note that $D(G)_{C_0P}$ does not have a left action of $D(Q)$. However, for each $\delta\in D(Q)$ and $\delta_{C_0P}\in D(G)_{C_0P}$, $\mathrm{Supp}(\delta\delta')\subseteq QC_0P\subseteq X$ is compact (\cite[Rem.~1.2.3]{Koh07}), so there always exists a compact open subset $C_0'\subseteq X_0$ such that $C'_0P_0=C'_0$ and $\mathrm{Supp}(\delta\delta')\subseteq C_0'P$ (use $X/P\cong X_0/P_0$). Using this, we can equip $\varprojlim_{C_0}D(G)_{C_0P}$ with a natural left $D(Q)$-action so that (\ref{equ: naive limit}) is $D(Q)$-equivariant. From the first isomorphism in (\ref{equ: loc an hat}) we have a natural map
\[\widehat{D}(G)_X\rightarrow D(G_0)_{C_0}\otimes_{D(P_0)}D(P)\cong D(G)_{C_0P}\]
for each compact open $C_0\subseteq X_0$ such that $C_0P_0=C_0$, hence a map
\begin{equation}\label{equ: limit dense image}
\widehat{D}(G)_X\rightarrow \varprojlim_{C_0}D(G)_{C_0P}
\end{equation}
which is $D(Q)$-equivariant by the definition of the $D(Q)$-actions on both sides. We see that the $D(Q)$-equivariant embedding (\ref{equ: naive limit}) factors as the composition of the isomorphism (\ref{equ: limit vers quotient}) with the $D(Q)$-equivariant map (\ref{equ: limit dense image}). This forces (\ref{equ: limit vers quotient}) to be also $D(Q)$-equivariant.
\end{proof}

If $A$ is a locally convex $E$-vector space endowed with a structure of a separately continuous algebra and $V$ (resp.~$W$) is a locally convex $E$-vector space endowed with a structure of a separately continuous right (resp.~left) $A$-module, we define $V\widehat\otimes_A W$ as the quotient of $V\widehat\otimes_{E,\iota}W$ by the closure of the $E$-vector subspace generated by elements $v a\otimes w - v\otimes a w$, $(a,v,w) \in A\times V\times W$ with the quotient topology. Let $V\otimes_{A} W$ be the obvious quotient of $V\otimes_{E,\iota} W$ (without completing) endowed with the quotient topology, it is then not difficult to check using the various universal properties that $V\widehat\otimes_A W$ is also the universal Hausdorff completion of $V\otimes_{A} W$ (\cite[\S 7]{S02}). The following lemma will be very useful.

\begin{lem}\label{limprojtensor}
Let $V = \varprojlim_{r} V_r$ and $W= \varprojlim_{r} W_r$ be Fr\'echet spaces written as countable projective limits of $E$-Banach spaces $V_r$, $W_r$ and assume that the transition maps have dense image. Let $A =\varprojlim_{r} A_r$ be a Fr\'echet algebra which is a countable projective limit of Noetherian $E$-Banach algebras $A_r$ with transition maps having dense image. Assume that $V$ and $W$ admit a separately continuous $A$-action induced by a separately continuous action of $A_r$ on $V_r$ and $W_r$. Assume finally that $V_r$ is a finitely generated $A_r$-module for each $r$. Then we have a canonical isomorphism of Fr\'echet spaces $V \widehat{\otimes}_A W \buildrel\sim\over\rightarrow \varprojlim_{r} (V_r \otimes_{A_r} W_r)$.
\end{lem}
\begin{proof}
By \cite[Prop.~1.1.29]{Em17} we have an isomorphism $V \widehat{\otimes}_E W \buildrel\sim\over\rightarrow \varprojlim_{r} (V_r \widehat{\otimes}_E W_r)$ and by \cite[Chap.~II, \S~3.5, Th.~1]{Bo} the image of $V$, $W$, $A$ in respectively $V_r$, $W_r$, $A_r$ has dense image. Let $C\subseteq V \widehat{\otimes}_E W$ (resp.~$C_r\subseteq V_r \widehat{\otimes}_E W_r$) be the closure of the $E$-vector subspace generated by elements $va \otimes w-v \otimes aw$ for $(v,w)\in V\times W$ and $a\in A$ (resp.~$(v,w)\in V_r\times W_r$ and $a\in A_r$), then $C_r$ is also the closure of the image of $C$ in $V_r \widehat{\otimes}_E W_r$ and hence $C\buildrel\sim\over\rightarrow \varprojlim_{r}C_r$. Applying Lemma \ref{lem: top ML} to the short exact sequences $0\rightarrow C_r\rightarrow V_r \widehat{\otimes}_E W_r\rightarrow V_r \widehat{\otimes}_{A_r} W_r\rightarrow 0$ we deduce an isomorphism $V \widehat{\otimes}_A W\cong \varprojlim_{r} (V_r \widehat{\otimes}_{A_r} W_r)$. But since $V_r$ is finitely generated over the Noetherian algebra $A_r$, hence finitely presented, using for instance \cite[Prop.~2.1.iii]{ST03} we see that $V_r \otimes_{A_r} W_r$ is already complete, hence $V_r \widehat{\otimes}_A W_r\cong V_r \otimes_{A_r} W_r$.
\end{proof}

Recall that a locally closed locally $K$-analytic submanifold of $G$ is a closed locally $K$-analytic submanifold of an open subset of $G$ (with its induced structure of locally $K$-analytic manifold). We need to slightly generalize the definition of smooth compact induction given before Lemma \ref{lem: induction coinvariant}. We fix $X\subseteq G$ a locally closed locally $K$-analytic submanifold of $G$ such that $XP=X$. We also fix a representation $\pi^{\infty}$ of $P$ in $\mathrm{Rep}^{\infty}_{\rm{adm}}(P)$ and recall that $(\pi^{\infty})^\vee=\Hom_E(\pi^{\infty},E)$ is in $\cC_{D^{\infty}(P)}$. We define $(\mathrm{ind}_{P}^{X}\pi^{\infty})^{\infty}$ to be the set of locally constant functions $f: X\rightarrow \pi^{\infty}$ such that
\begin{itemize}
\item $f(xh)=h^{-1}\cdot f(x)$ for $x\in X$ and $h\in P$;
\item there exists a compact open subset $C_f$ of $X$ such that $f(x)=0$ for $x\notin C_fP$.
\end{itemize}
If moreover $QX=X$ for some locally $K$-analytic closed subgroup $Q\subseteq G$, then $(\mathrm{ind}_{P}^{X}\pi^{\infty})^{\infty}$ is naturally a (left) smooth $Q$-representation via $(h'(f))(x)\defeq f((h')^{-1}x)$ ($h'\in Q$, $x\in X$, $f\in (\mathrm{ind}_{P}^{X}\pi^{\infty})^{\infty}$).
As $D(G)_P\otimes_{U(\fg)}E\cong (D(G_0)_{P_0}\otimes_{D(P_0)}D(P))\otimes_{U(\fg)}E\cong D^{\infty}(G_0)_{P_0}\otimes_{D^{\infty}(P_0)}D^\infty(P)\cong D^{\infty}(P)$ (Lemma~\ref{lem: smooth restriction}) we see that $(\pi^{\infty})^\vee$ is in $\cC_{D(G)_P}$ by \cite[Prop.~3.7]{ST03} and \cite[Lemma 3.8]{ST03}.

\begin{lem}\label{lem: dual of compact induction}
With the above notation $(\mathrm{ind}_{P}^{X}\pi^{\infty})^{\infty}$ is a locally convex $E$-vector space of compact type and we have a canonical isomorphism of Fr\'echet $E$-spaces
\begin{equation}\label{equ: dual of compact induction}
\widehat{D}(G)_X\widehat{\otimes}_{D(G)_P}(\pi^{\infty})^\vee\cong
\widehat{D}^{\infty}(G)_X\widehat{\otimes}_{D^{\infty}(G)_P}(\pi^{\infty})^\vee
\cong \big((\mathrm{ind}_{P}^{X}\pi^{\infty})^{\infty}\big)^\vee.
\end{equation}
If moreover $QX=X$ for some locally $K$-analytic closed subgroup $Q\subseteq G$, then the isomorphisms in (\ref{equ: dual of compact induction}) are left $D(Q)$-equivariant with the $D(Q)$-actions factoring through $D^{\infty}(Q)$.
\end{lem}
\begin{proof}
It follows from (\ref{equ: smooth hat}) and (\ref{equ: hat loc an to smooth}) that
\begin{equation}\label{equ: pass to smooth}
\widehat{D}(G)_X\otimes_{D(G)_1}E\cong \widehat{D}^{\infty}(G)_X\cong (\varprojlim D^{\infty}(G_0)_{C_0})\otimes_{D^{\infty}(P_0)}D^{\infty}(P)
\end{equation}
where $C_0$ runs through the compact open subsets of $X_0\defeq X\cap G_0$ such that $C_0P_0=C_0$. As the left $D(G)_P$-action on $(\pi^{\infty})^\vee$ factors through $D^{\infty}(G)_P$ and as $D^{\infty}(G)_P\cong D^{\infty}(P)$ (Lemma \ref{lem: smooth restriction}), by (\ref{equ: pass to smooth}) we have topological isomorphisms
\[\widehat{D}(G)_X{\otimes}_{D(G)_P}(\pi^{\infty})^\vee\cong \widehat{D}(G)_X{\otimes}_{D^\infty(G)_P}(\pi^{\infty})^\vee\cong (\varprojlim_{C_0}D^{\infty}(G_0)_{C_0}){\otimes}_{D^{\infty}(P_0)}(\pi^{\infty})^\vee.\]
Taking universal Hausdorff completion, we obtain topological isomorphisms
\[\widehat{D}(G)_X\widehat{\otimes}_{D(G)_P}(\pi^{\infty})^\vee\cong \widehat{D}^{\infty}(G)_X\widehat{\otimes}_{D^{\infty}(G)_P}(\pi^{\infty})^\vee\cong \big(\varprojlim_{C_0}D^{\infty}(G_0)_{C_0}\big)\widehat{\otimes}_{D^{\infty}(P_0)}(\pi^{\infty})^\vee\]
which the last space is a Fr\'echet space since both $\varprojlim_{C_0}D^{\infty}(G_0)_{C_0}$ and $(\pi^{\infty})^\vee$ are. We now prove the last isomorphism in (\ref{equ: dual of compact induction}). We have topological isomorphisms of Fr\'echet spaces
\begin{multline*}
\big(\varprojlim_{C_0}D^{\infty}(G_0)_{C_0}\big)\widehat{\otimes}_{D^{\infty}(P_0)}(\pi^{\infty})^\vee\cong
\Big(\varprojlim_{C_0}\big(D^{\infty}(G_0/P_0)_{C_0/P_0}\widehat{\otimes}_E D^\infty(P_0)\big)\Big)\widehat{\otimes}_{D^{\infty}(P_0)}(\pi^{\infty})^\vee\\
\cong \Big(\big(\varprojlim_{C_0}D^{\infty}(G_0/P_0)_{C_0/P_0}\big)\widehat{\otimes}_E D^\infty(P_0)\Big)\widehat{\otimes}_{D^{\infty}(P_0)}(\pi^{\infty})^\vee\cong \big(\varprojlim_{C_0} D^{\infty}(G_0/P_0)_{C_0/P_0}\big)\widehat{\otimes}_E (\pi^{\infty})^\vee\\
\cong \varprojlim_{C_0}\big(D^{\infty}(G_0/P_0)_{C_0/P_0}\widehat{\otimes}_E(\pi^{\infty})^\vee\big)\cong \varprojlim_{C_0}\Big(\big(D^{\infty}(G_0/P_0)_{C_0/P_0}\widehat{\otimes}_E D^\infty(P_0)\big)\widehat{\otimes}_{D^{\infty}(P_0)}(\pi^{\infty})^\vee\Big)\\
\cong \varprojlim_{C_0} \big(D^{\infty}(G_0)_{C_0}\widehat{\otimes}_{D^{\infty}(P_0)}(\pi^{\infty})^\vee\big)
\end{multline*}
where the first and last isomorphisms follow from (\ref{equ: smooth coset section}), the second and fourth from \cite[Prop.~1.1.29]{Em17} (recall there are countably many $C_0$) and the third and fifth from \cite[Lemme 3.1]{BD23}. For $r\in \cI$ let $(\pi^{\infty})^\vee_r\defeq D^{\infty}(P_0)_r\otimes_{D^{\infty}(P_0)}(\pi^{\infty})^\vee$. By the same proof as for \ref{it: support basis 1} of Lemma \ref{lem: support basis}~replacing $G_0$ by $P_0$ and $G_1^m$ by $P_1^m=G_1^m\cap P_0$ (see the discussion above Lemma \ref{lem: free mod coset basis}) shows that $D^\infty(P_0)_{C,r}$ is finite dimensional. Applying this with $C=P_0$ and since $(\pi^{\infty})^\vee$ is in $\cC_{D(G)_P}$, it follows that $(\pi^{\infty})^\vee_r$ is a finite dimensional $E$-vector space. Then by Lemma \ref{limprojtensor} we have an isomorphism of Fr\'echet spaces
\begin{equation*}
D^{\infty}(G_0)_{C_0}\widehat{\otimes}_{D^{\infty}(P_0)}(\pi^{\infty})^\vee\cong \varprojlim_{r}\big(D^{\infty}(G_0)_{C_0,r}{\otimes}_{D^{\infty}(P_0)_r}(\pi^{\infty})^\vee_r\big)
\end{equation*}
Putting everything together we obtain an isomorphism of Fr\'echet spaces
\begin{equation}\label{equ: loc an smooth tensor}
\widehat{D}(G)_X\widehat{\otimes}_{D(G)_P}(\pi^{\infty})^\vee\cong \varprojlim_{C_0,}\varprojlim_{r}(D^{\infty}(G_0)_{C_0,r}{\otimes}_{D^{\infty}(P_0)_r}(\pi^{\infty})^\vee_r).
\end{equation}
Given $r\in \cI$ let $m\geq 0$, $s=r^{p^m}$ as in Lemma~\ref{lem: support basis}, and define $\pi^{\infty}_r\defeq (\pi^{\infty})^{P_1^m}$ which is a finite dimensional representation of the finite group $P_0/P_1^m$. As $\{P_1^m\}_{r\in \cI}$ is a system of open neighborhoods of $1$ inside $P_0$, we have $\pi^{\infty}=\varinjlim_{r\in \cI}\pi^{\infty}_r$.
Recall from \ref{it: support basis 1} of Lemma~\ref{lem: support basis} (with $G_0$ there replaced with $P_0$) that
\begin{equation}\label{equ: p0r}
D^{\infty}(P_0)_r\cong E[P_0/P_1^m]\cong E\otimes_{D^{\infty}(P_1^m)}D^{\infty}(P_0),
\end{equation}
which implies
\begin{equation}\label{equ: smooth dual of inv}
(\pi^{\infty})^\vee_r=D^{\infty}(P_0)_r\otimes_{D^{\infty}(P_0)}(\pi^{\infty})^\vee\cong E\otimes_{D^{\infty}(P_1^m)}(\pi^{\infty})^\vee \cong (\pi^{\infty}_r)^\vee.
\end{equation}
For any compact open subset $C_0\subseteq X_0$ such that $C_0P_0=C_0$, by \ref{it: support basis 2} of Lemma~\ref{lem: support basis} we have
$D^{\infty}(G_0)_{C_0,r}\cong E[C_0G_1^m/G_1^m]$, which together with (\ref{equ: smooth dual of inv}) and (\ref{equ: p0r}) implies
\begin{equation}\label{equ: smooth induction basis}
D^{\infty}(G_0)_{C_0,r}\otimes_{D^{\infty}(P_0)_r}(\pi^{\infty})^\vee_r\cong E[C_0G_1^m/G_1^m]\otimes_{E[P_0/P_1^m]}(\pi^{\infty}_r)^\vee\cong \big(\mathrm{ind}_{P_0/P_1^m}^{C_0G_1^m/G_1^m}\pi^{\infty}_r\big)^\vee
\end{equation}
where $\mathrm{ind}_{P_0/P_1^m}^{C_0G_1^m/G_1^m}\pi^{\infty}_r$ is the set of functions $\phi: C_0G_1^m/G_1^m\rightarrow \pi^{\infty}_r$ satisfying $\phi(x\overline{h})=\overline{h^{-1}}\cdot \phi(x)$ for $x\in C_0G_1^m/G_1^m$ and $\overline{h}\in P_0/P_1^m$.
If we lift such $\phi: C_0G_1^m/G_1^m\rightarrow \pi^{\infty}_r$ to a locally constant function on $C_0G_1^m$, take its restriction to $C_0$ and then extend it by zero on the open subset $X_0\setminus C_0$ of $X_0$, we obtain a locally constant function $f: X_0\rightarrow \pi^{\infty}_r\subseteq \pi^{\infty}$ which is supported on the compact open subset $C_0$ of $X_0$, constant on the compact open subsets $yG_1^m\cap C_0\subseteq C_0$ for $y\in C_0$ and which satisfies $f(yh)=h^{-1}\cdot f(y)$ for $y\in C_0$ and $h\in P_0$. Using $X_0/P_0\buildrel\sim\over\rightarrow X/P$ (see the paragraph before Lemma \ref{lem: limit vers quotient}) there exists a unique extension of $f: X_0\rightarrow \pi^{\infty}$ to a locally constant function (still denoted) $f:X\rightarrow \pi^{\infty}$ such that $f(yh)=h^{-1}\cdot f(y)$ for $y\in X$ and $h\in P$. This defines an injection
\begin{equation}\label{equ: smooth induction embedding}
\mathrm{ind}_{P_0/P_1^m}^{C_0G_1^m/G_1^m}\pi^{\infty}_r \hookrightarrow (\mathrm{ind}_P^{X}\pi^{\infty})^{\infty}
\end{equation}
with image consisting of those $f:X\rightarrow \pi^{\infty}$ such that
\begin{itemize}
\item $f$ is supported on $C_0P;$
\item the restriction$f|_{C_0}$ has image in $\pi^{\infty}_r$ and is constant on $yG_1^m\cap C_0$ for $y\in C_0$.
\end{itemize}
As each element of $(\mathrm{ind}_P^{X}\pi^{\infty})^{\infty}$ obviously satisfies the above two conditions for some compact open subset $C_0\subseteq X_0$ and some $r\in \cI$, we deduce from (\ref{equ: smooth induction embedding}) a topological isomorphism with both sides endowed with the finest locally convex topology
\begin{equation}\label{equ: smooth induction as colimit}
\varinjlim_{C_0}\varinjlim_{r}\mathrm{ind}_{P_0/P_1^m}^{C_0G_1^m/G_1^m}\pi^{\infty}_r \buildrel\sim\over\longrightarrow (\mathrm{ind}_P^{X}\pi^{\infty})^{\infty},
\end{equation}
and in particular $(\mathrm{ind}_P^{X}\pi^{\infty})^{\infty}$ is of compact type.
Taking the dual of (\ref{equ: smooth induction as colimit}), by \cite[Prop.~16.10]{S02} and (\ref{equ: smooth induction basis}), we obtain topological isomorphisms of Fr\'echet spaces
\begin{multline*}
\big((\mathrm{ind}_P^{X}\pi^{\infty})^{\infty}\big)^\vee \buildrel\sim\over\longrightarrow \big(\varinjlim_{C_0}\varinjlim_{r}\mathrm{ind}_{P_0/P_1^m}^{C_0G_1^m/G_1^m}\!\pi^{\infty}_r\big)^\vee \cong \varprojlim_{C_0}\varprojlim_{r}\big(\mathrm{ind}_{P_0/P_1^m}^{C_0G_1^m/G_1^m}\pi^{\infty}_r\big)^\vee\\
\cong \varprojlim_{C_0}\varprojlim_{r}\big(D^{\infty}(G_0)_{C_0,r}{\otimes}_{D^{\infty}(P_0)_r}(\pi^{\infty})^\vee_r\big).
\end{multline*}
Together with (\ref{equ: loc an smooth tensor}) this finishes the proof of (\ref{equ: dual of compact induction}). The left $D(Q)$-equivariance of the isomorphisms in (\ref{equ: dual of compact induction}) is easy and left to the reader.
\end{proof}

\subsection{Fr\'echet completion of objects of \texorpdfstring{$\cO^{\fb}_{\rm{alg}}$}{Obalg}}\label{subsec: completion O}

We prove several statements on the canonical Fr\'echet completions of $U(\fg)$-modules in $\cO^{\fb}_{\rm{alg}}$ defined in \S\ref{subsec: Frechet t action} and use them to give a useful description of the continuous dual of Orlik-Strauch representations (Proposition \ref{prop: OS via completion}).\bigskip

We fix $I\subseteq \Delta$ and use the notation of \S\ref{subsec: loc an dist} with $\mathbf{G}_0\defeq \GL_n/{\cO_K}$ and $\mathbf{P}_0\defeq P_I/{\cO_K}$, so that we are back with $P_I=\mathbf{P}_0(K)\subseteq G=\mathbf{G}_0(K)=\GL_n(K)$, $\fp_I= \mathrm{Lie}(P_I)$ and we have $G_0=\GL_n(\cO_K)$, $P_{I,0}=P_I(\cO_K)$, etc. We use without further ado that a $U(\fp_I)$-module which is finite dimensional over $E$ is an algebraic representation of $P_I$ (\cite[Lemma 3.2]{OS15}), and hence in particular is a $D(P_I)$-module, and thus also a $D(P_I)_1$-module. We start by recalling Schmidt's result on the canonical Fr\'echet completion of objects in $\cO^{\fb}_{\rm{alg}}$ (\cite{Schm13}). We define $\widehat{\cO}^{\fb}_{\rm{alg}}\subseteq \cC_{D(G)_1}$ as the full subcategory consisting of those coadmissible left $D(G)_1$-modules $D$ such that $D\vert_{U(\ft)}$ is a small Fr\'echet $U(\ft)$-module in the sense of \ref{it: Frechet t 2} of Definition~\ref{def: Frechet t module} and the (left) $U(\fg)$-module $\bigoplus_{\mu\in\Lambda}D_{\mu}$ is an object of $\cO^{\fb}_{\rm{alg}}$ (recall that $D_{\mu}\subseteq D$ is the eigenspace of $D$ for the weight $\mu$).

\begin{prop}\label{prop: category O completion}
\hspace{2em}
\begin{enumerate}[label=(\roman*)]
\item \label{it: O completion 1} The algebra $D(G)_1$ is flat over $U(\fg)$.
\item \label{it: O completion 2} The functor $M\mapsto \cM\defeq D(G)_1\otimes_{U(\fg)}M$ induces an equivalence of (abelian) categories between $\cO^{\fb}_{\rm{alg}}$ and $\widehat{\cO}^{\fb}_{\rm{alg}}$, with a quasi-inverse given by $D\mapsto \bigoplus_{\mu\in\Lambda}D_{\mu}$.
\item \label{it: O completion 3} Let $M$ in $\cO^{\fp_I}_{\rm{alg}}$ and $X$ a finite dimensional $U(\fp_I)$-module such that one has a surjection $q: U(\fg)\otimes_{U(\fp_I)}X\twoheadrightarrow M$. Then $D(G)_1\mathrm{ker}(q)\subseteq D(G)_1\otimes_{U(\fp_I)}X$ is a coadmissible (left) $D(G)_1$-submodule of $D(G)_1\otimes_{U(\fp_I)}X\buildrel\sim\over\rightarrow D(G)_1\otimes_{D(P_I)_1}X$ and we have an isomorphism in $\widehat{\cO}^{\fb}_{\rm{alg}}\subseteq\cC_{D(G)_1}$:
\begin{equation}\label{equ: presentation completion}
\cM=D(G)_1\otimes_{U(\fg)}M\cong (D(G)_1\otimes_{D(P_I)_1}X)/(D(G)_1\mathrm{ker}(q)).
\end{equation}
\end{enumerate}
\end{prop}
\begin{proof}
Part \ref{it: O completion 1} is \cite[Thm.~4.3.3]{Schm13} and part \ref{it: O completion 2} is the case $I=\emptyset$ of \cite[Thm.~4.3.1]{Schm13}, with the harmless difference that we have the extra condition that the weights of $M$ are integral.
We prove \ref{it: O completion 3}. By applying \ref{it: O completion 1} and \ref{it: O completion 2} to $0\rightarrow \mathrm{ker}(q)\rightarrow U(\fg)\otimes_{U(\fp_I)}X\rightarrow M\rightarrow 0$, we obtain a short exact sequence in $\widehat{\cO}^{\fb}_{\rm{alg}}\subseteq\cC_{D(G)_1}$:
\begin{equation}\label{equ: O completion seq}
0\longrightarrow D(G)_1\otimes_{U(\fg)}\mathrm{ker}(q)\longrightarrow D(G)_1\otimes_{U(\fp_I)}X \longrightarrow D(G)_1\otimes_{U(\fg)}M\longrightarrow 0.
\end{equation}
By the density of $U(\fg)$ in $D(G)_1$ and \ref{equ: O completion seq} we know that $D(G)_1\otimes_{U(\fg)}\mathrm{ker}(q)$ is the closure of $\mathrm{ker}(q)$ in $D(G)_1\otimes_{U(\fp_I)}X$ and by \ref{it: O completion 2} we have
\[\bigoplus_{\mu\in\Lambda}(D(G)_1\otimes_{U(\fg)}\mathrm{ker}(q))_{\mu}\cong \mathrm{ker}(q).\]
As $\mathrm{ker}(q)$ is a finitely generated $U(\fg)$-submodule of $U(\fg)\otimes_{U(\fp_I)}X$, $D(G)_1\mathrm{ker}(q)$ is a finitely generated $D(G)_1$-submodule of $D(G)_1\otimes_{U(\fp_I)}X$, hence is coadmissible by \cite[Cor.~3.4.iv]{ST03} and thus closed in $D(G)_1\otimes_{U(\fp_I)}X$ by \cite[Lemma~3.6]{ST03}. As it contains $\mathrm{ker}(q)$ as dense subspace, we deduce $D(G)_1\otimes_{U(\fg)}\mathrm{ker}(q)\buildrel\sim\over\rightarrow D(G)_1\mathrm{ker}(q)$. Comparing (\ref{equ: presentation completion}) with (\ref{equ: O completion seq}), it remains to show that the natural map $D(G)_1\otimes_{U(\fp_I)}X\rightarrow D(G)_1\otimes_{D(P_I)_1}X$ is an isomorphism. It is enough to prove that we have an isomorphism
\begin{equation}\label{equ: O completion tensor}
X\buildrel\sim\over\longrightarrow D(P_I)_1\otimes_{U(\fp_I)}X.
\end{equation}
Since the action of $U(\fp_I)$ on $X$ extends to $D(P_I)_1$, there is a natural surjection $D(P_I)_1\otimes_{U(\fp_I)}X\twoheadrightarrow X$, $\mathfrak d \otimes x\mapsto {\mathfrak d}x$ which is the identity of $X$ when composed with (\ref{equ: O completion tensor}). Hence the map (\ref{equ: O completion tensor}) is injective, and it is enough to prove its surjectivity. But using \cite[Cor.~3.4]{ST03} applied to the Fr\'echet-Stein algebra $D(P_I)_1$, one can check that $D(P_I)_1\otimes_{U(\fp_I)}X$ is a coadmissible $D(P_I)_1$-module, hence is a Fr\'echet space. Since the map $X\rightarrow D(P_I)_1\otimes_{U(\fp_I)}X$ has dense image (as $U(\fp_I)$ is dense in $D(P_I)_1$) with its source finite dimensional, it follows that it is surjective (and that $D(P_I)_1\otimes_{U(\fp_I)}X$ is finite dimensional).
\end{proof}

Let $I\subseteq \Delta$, following \cite[\S 3]{ST02a} we say that a smooth representation $\pi^{\infty}$ of $P_{I,0}$ over $E$ is strongly admissible if there exists a $P_{I,0}$-equivariant embedding $\pi^{\infty}\hookrightarrow C^{\infty}(P_{I,0})^m$ for some $m\geq 0$. Since $P_{I,0}$ is compact, $C^{\infty}(P_{I,0})$ is a direct sum of (finite dimensional) irreducible smooth representations of $P_{I,0}$ over $E$. In particular, $\pi^{\infty}$ is a direct summand of $C^{\infty}(P_{I,0})^m$. The following lemma can be extracted from the beginning of \cite[\S 4.4]{Schr11}, we provide a proof for the reader's convenience.

\begin{lem}\label{lem: loc alg FIN}
Let $I\subseteq \Delta$, $\pi^{\infty}$ a strongly admissible smooth representation of $P_{I,0}$ over $E$ (see \cite[\S 2]{ST01}) and $X$ a finite dimensional left $D(P_{I,0})$-module. Then $X\otimes_E (\pi^{\infty})^\vee$ admits a resolution by finite free $D(P_{I,0})$-modules. In other words, $X\otimes_E (\pi^{\infty})^\vee$ satisfies the condition (FIN) from \cite[p.321]{ST05}.
\end{lem}
\begin{proof}
Since $\pi^{\infty}$ is a strongly admissible smooth representation of $P_{I,0}$ over $E$, taking duals there exists $m\geq 0$ such that $(\pi^{\infty})^{\vee}$ is a direct summand of $D^{\infty}(P_{I,0})^m$, and in particular $(\pi^{\infty})^{\vee}$ admits a resolution
\[[\cdots\rightarrow Y_{i}\rightarrow \cdots\rightarrow Y_{0}]\buildrel\sim\over\longrightarrow (\pi^{\infty})^{\vee}\]
with $Y_{i}=D^{\infty}(P_{I,0})^m$ for each $i\geq 0$ (writing $D^{\infty}(P_{I,0})^m \cong D^\infty \oplus (\pi^{\infty})^{\vee}$, the differential maps alternatively project onto $(\pi^{\infty})^{\vee}$ or onto $D^\infty$). Recall from \cite[p.307]{ST05} that $D^{\infty}(P_{I,0})$ admits the following (standard) resolution by finite free $D(P_{I,0})$-modules
\[D(P_{I,0})\otimes_E\wedge^{\bullet}\fp_{I} \rightarrow D^{\infty}(P_{I,0}).\]
We set $X_{i,j}\defeq Y_{i}\otimes_E\wedge^{j}\fp_{I}=(D(P_{I,0})\otimes_E\wedge^{j}\fp_{I})^m$ for each $i,j\geq 0$, so that $X_{i,\bullet}$ is a resolution of $Y_{i}$ by finite free $D(P_{I,0})$-modules for $i\geq 0$. Thanks to a construction of C.T.C Wall (see \cite[Chapitre V~(3.1),(3.2)]{La65} and \cite[Thm.~6.1]{Koh11}), there exists a map $d_{i,j}^{k}: X_{i,j}\rightarrow X_{i-k,j+k-1}$ for each $i,j,k\geq 0$ such that the data $(X_{i,j},d_{i,j}^{k})_{i,j,k\geq 0}$ altogether form a Wall complex, whose total complex $(X_{r},\Delta_{r})_{r\geq 0}$ (defined as $X_{r}\defeq\bigoplus_{i+j=r}X_{i,j}$ and $\Delta_{r}\defeq \sum_{i=0}^{r}\sum_{k=0}^{i}d_{i,r-i}^{k}: X_{r}\rightarrow X_{r-1}$) is a resolution of $(\pi^{\infty})^{\vee}$. Since $X_{i,j}$ is a finite free $D(P_{I,0})$-module which is non-zero only if $i,j\geq 0$, we conclude that $X_{r}$ is a finite free $D(P_{I,0})$-module for each $r\geq 0$ and thus that $(\pi^{\infty})^\vee$ satisfies (FIN). Since $X$ is a finite dimensional left $D(P_{I,0})$-module, we deduce from \cite[Lemma~3.5]{Schr11} that $X\otimes_E X_{r}$ is a finite free $D(P_{I,0})$-module for each $r\geq 0$. We thus obtain the resolution $X\otimes_EX_{\bullet}$ of $X\otimes_E (\pi^{\infty})^\vee$ by finite free $D(P_{I,0})$-modules, namely $X\otimes_E (\pi^{\infty})^\vee$ also satisfies the condition (FIN).
\end{proof}

\begin{lem}\label{lem: support variant}
Let $M$ in $\cO^{\fp_I}_{\rm{alg}}$ and $X$ a finite dimensional $\ft$-semi-simple $U(\fp_I)$-module such that one has a surjection $q: U(\fg)\otimes_{U(\fp_I)}X\twoheadrightarrow M$ in $\cO^{\fp_I}_{\rm{alg}}$. Then the left $D(G)_1$-action on $\cM\defeq D(G)_1\otimes_{U(\fg)}M\in \cC_{D(G)_1}$ (see \ref{it: O completion 2} of Proposition \ref{prop: category O completion}) naturally extends to a left $D(G)_{P_I}$-action which fits into an isomorphism of coadmissible $D(G)_{P_I}$-modules:
\begin{equation*}
\cM\cong (D(G)_{P_I}\otimes_{D(P_I)} X)/(D(G)_{P_I}\mathrm{ker}(q))\buildrel\sim\over\longrightarrow (D(G)_{P_I}\widehat\otimes_{D(P_I)} X)/(D(G)_{P_I}\mathrm{ker}(q)).
\end{equation*}
\end{lem}
\begin{proof}
First note that, by Lemma~\ref{lem: loc alg FIN}, the finite dimensional representation $X$ fits into an exact sequence of the form
\begin{equation}\label{equ: support variant FIN}
\cdots \rightarrow D^{k}\rightarrow \cdots \rightarrow D^0\rightarrow X \rightarrow 0
\end{equation}
where $D^k$ is a finite free $D(P_{I,0})$-modules for $k\leq 0$. Moreover, it follows from the argument in the proof of \cite[Lemma 6.3i]{ST05} with $D(G)_{P_I}$ instead of $D(G)$ that applying $D(G_0)_{P_{I,0}} \otimes_{D(P_{I,0})}(-)$ to (\ref{equ: support variant FIN}) gives another exact sequence
\begin{equation*}
\cdots \rightarrow D(G_0)_{P_{I,0}} \otimes_{D(P_{I,0})}D^{k}\rightarrow \cdots \rightarrow D(G_0)_{P_{I,0}} \otimes_{D(P_{I,0})}D^0\rightarrow D(G_0)_{P_{I,0}} \otimes_{D(P_{I,0})}X \rightarrow 0
\end{equation*}
(using $D(G)_{P_I}\cong D(G_0)_{P_{I,0}}\otimes_{D(P_{I,0})}D(P_I)$). Note that the argument of \emph{loc.~cit.}~indeed extends since $D(G_0)_{P_{I,0},r}$ is flat over $D(P_{I,0})_r$ for $r\in \cI$, which follows from the flatness of $D(G_0)_r$ over $D(P_{I,0})_r$ (\cite[Prop.~2.6]{Schm09}) and the fact that $D(G_0)_{P_{I,0},r}$ is a direct summand of $D(G_0)_r$ as $D(P_{I,0})_r$-module by \ref{it: support basis 2} of Lemma \ref{lem: support basis} (applied with $C=G_0$ and $H_0=P_{I,0}$). This implies that $D(G)_{P_I} \otimes_{D(P_I)} X$ is a $D(G_0)_{P_{I,0}}$-module of finite presentation, hence is coadmissible by \cite[Cor.~3.4.v]{ST03}. In particular we have $D(G)_{P_I} \otimes_{D(P_I)} X \buildrel\sim\over\rightarrow D(G)_{P_I} \widehat{\otimes}_{D(P_I)} X$ (compare with \cite[Rem.~5.1.3(ii)]{Bre19}). Since we have seen in the proof of \ref{it: O completion 3} of Proposition \ref{prop: category O completion} that $D(G)_{1} \otimes_{D(P_I)_1} X\cong D(G)_{1} \otimes_{U(\fp_I)} X$ is a coadmissible, hence complete, $D(G)_1$-module, we also have $D(G_0)_{1} \otimes_{D(P_{I,0})_1} X\cong D(G)_{1} \otimes_{D(P_I)_1} X \buildrel\sim\over\rightarrow D(G)_{1} \widehat{\otimes}_{D(P_I)_1} X$.

By \ref{it: O completion 3} of Proposition~\ref{prop: category O completion} we have an isomorphism of coadmissible left $D(G)_1$-modules
\begin{equation*}
\cM=D(G)_1\otimes_{U(\fg)}M\cong (D(G)_1\otimes_{D(P_I)_1} X)/(D(G)_1\mathrm{ker}(q)),
\end{equation*}
hence it is enough to prove that the natural map
\[D(G_0)_1\otimes_{D(P_{I,0})_1} X\longrightarrow D(G_0)_{P_{I,0}}\otimes_{D(P_{I,0})} X\cong D(G)_{P_{I}}\otimes_{D(P_{I})} X\]
is a topological isomorphism which sends $D(G)_1\mathrm{ker}(q)$ to $D(G)_{P_I}\mathrm{ker}(q)$. For the first statement, using Lemma \ref{limprojtensor} (noting that $X$ is a finitely generated $D(P_I)_{1,r}$-module or $D(P_{I,0})_r$-module) and the above discussion, it suffices to prove that we have canonical isomorphisms of $E$-Banach spaces for $r\in \cI$
\begin{equation*}
D(G)_{1,r}\otimes_{D(P_I)_{1,r}} X\buildrel\sim\over\longrightarrow D(G_0)_{P_{I,0},r} \otimes_{D(P_{I,0})_r} X.
\end{equation*}
But this follows from
\begin{equation}\label{equ: reduce support r}
D(G)_{1,r}{\otimes}_{D(P_I)_{1,r}}D(P_{I,0})_r\cong D(G_0)_{P_{I,0},r},
\end{equation}
which itself follows from \ref{it: support basis 1} of Lemma \ref{lem: support basis} (applied with $C=P_{I,0}$ and both $G_0$ and $P_{I,0}$, noting that $P_{I,0}/P_{I,1}^m\buildrel\sim\over\rightarrow P_{I,0}G_1^m/G_1^m$). Thus $D(G)_1\otimes_{D(P_I)_1} X$ is a coadmissible left $D(G)_{P_I}$-module. As the left $U(\fp_I)$-action on $\mathrm{ker}(q)$ extends to a left $D(P_I)$-action (this is so for any object of $\cO^{\fp_I}_{\rm{alg}}$, see \cite[\S 3.4]{OS15}), we deduce that $D(G)_1\mathrm{ker}(q)$ (which is the closure of $\mathrm{ker}(q)$ in $D(G)_1\otimes_{D(P_I)_1} X$, see (\ref{equ: presentation completion})) is a closed subspace of $D(G)_{P_I}\otimes_{D(P_I)} X$ which is both $D(P_I)$-stable and $D(G)_1$-stable. As $D(P_I)$ and $D(G)_1$ generate a dense subalgebra of $D(G)_{P_I}$ (see \cite[Prop.~1.2.12]{Koh07}), we deduce that $D(G)_1\mathrm{ker}(q)$ is a closed $D(G)_{P_I}$-submodule of $D(G)_{P_I}\otimes_{D(P_I)} X$. In particular $D(G)_1\mathrm{ker}(q)\cong D(G)_{P_I} \mathrm{ker}(q)$ since both spaces are the closure of $\mathrm{ker}(q)$ in $D(G)_{P_I}\otimes_{D(P_I)} X$. This finishes the proof.
\end{proof}

We will need the following result.

\begin{lem}\label{lem: family of standard semi norms}
Let $M$ in $\cO_{\rm{alg}}^{\fp_I}$ and $\cM$ the canonical Fr\'echet completion of $M$ as a coadmissible $D(G_0)_{P_{I,0}}$-module (Lemma~\ref{lem: support variant}) with $\cM \cong \varprojlim_{r\in\cI} \cM_r$ where $\cM_r\!\defeq \! D(G_0)_{P_{I,0},r}\otimes_{D(G_0)_{P_{I,0}}}\cM$. Then there exists a family of standard semi-norms $\{|\cdot|_r\}_{r\in \cI}$ on $M$ (see Definition~\ref{def: standard semi norm}) such that $\cM_r$ is the completion of $M$ under $|\cdot|_r$ for $r\in \cI$. Moreover the natural map $\cM\rightarrow \cM_r$ is a continuous injection.
\end{lem}
\begin{proof}
We first prove that the first statement implies the second. Assume that the (continuous) map $\cM\rightarrow \cM_r$ is not injective, then its non-zero kernel is a closed Fr\'echet $U(\ft)$-submodule of $\cM$. By Lemma \ref{lem: t closed} its intersection with $M$ is non-zero. By Remark \ref{rem: dense in Banach} applied to the standard semi-norm $|\cdot|_r$, the completion under $|\cdot|_r$ of this intersection is also non-zero in $\cM_r$, which is a contradiction. Hence the map $\cM\rightarrow \cM_r$ is injective.\bigskip

Let $X$ be a finite dimensional $\ft$-semi-simple $U(\fp_I)$-module such that one has a surjection $q: U(\fg)\otimes_{U(\fp_I)}X\twoheadrightarrow M$ in $\cO^{\fp_I}_{\rm{alg}}$. We fix $r\in\cI$ and let $m,s$ be as before Lemma~\ref{lem: free mod coset basis}. The isomorphism $\fg\cong \fn_I^+\oplus \fp_I$ induces an isomorphism of locally $K$-analytic manifolds $G_1^m=(G_1^m\cap N_I^+)\times (G_1^m\cap P_I)$ (see the paragraph before Lemma~\ref{lem: free mod coset basis}), which together with \cite[Prop.~5.9, Prop.~6.2]{Schm08} (see also Lemma~\ref{lem: free mod coset basis}) implies that
\begin{multline*}
D(G)_{1,r}=D(G)_{G_1^m,r}\cong D(G_1^m)_s\cong D(G_1^m\cap N_I^+)_s\widehat{\otimes}_E D(G_1^m\cap P_I)_s\\
\cong D(N_I^+)_{G_1^m\cap N_I^+,r}\widehat{\otimes}_E D(P_I)_{G_1^m\cap P_I,r}=D(N_I^+)_{1,r}\widehat{\otimes}_E D(P_I)_{1,r}.
\end{multline*}
This together with (\ref{equ: reduce support r}) gives the following isomorphisms
\begin{multline*}
D(G_0)_{P_{I,0},r}\otimes_{D(P_{I,0})}X=D(G_0)_{P_{I,0},r}\otimes_{D(P_{I,0})_r}(D(P_{I,0})_r\otimes_{D(P_{I,0})}X)\\
\cong D(G_0)_{P_{I,0},r}\otimes_{D(P_{I,0})_r}X \cong (D(G)_{1,r}\otimes_{D(P_I)_{1,r}}D(P_{I,0})_r)\otimes_{D(P_{I,0})_r}X\\
\cong D(G)_{1,r}\otimes_{D(P_{I})_{1,r}}X \cong D(N_I^+)_{1,r}\otimes_E X,
\end{multline*}
and taking projective limit over $r\in\cI$ recovers $D(G_0)_{P_{I,0}}\otimes_{D(P_{I,0})}X \cong D(N_I^+)_{1}\otimes_E X$. By replacing $\fg$ in \cite[Thm.~1.4.2]{Koh07} with $\fn_I^+$, we know that $U(\fn_I^+)$ admits a standard semi-norm $|\cdot|_r'$ for which its completion is $D(N_I^+)_{1,r}$ (and thus $D(N_I^+)_{1}\cong \varprojlim_{r\in\cI}D(N_I^+)_{1,r}$ is the completion of $U(\fn_I^+)$ under the Fr\'echet topology defined by $\{|\cdot|_r'\}_{r\in\cI}$). By \ref{it: tensor with fin dim 2} of Lemma~\ref{lem: tensor with fin dim} we know that $U(\fg)\otimes_{U(\fp_I)}X\cong U(\fn_I^+)\otimes_E X$ admits a standard semi-norm $|\cdot|_r''$ under which its completion is $D(G_0)_{P_{I,0},r}\otimes_{D(P_{I,0})}X\cong D(N_I^+)_{1,r}\otimes_E X$. In particular, still denoting by $|\cdot|_r''$ the induced semi-norm on $D(G_0)_{P_{I,0},r}\otimes_{D(P_{I,0})}X$, the Fr\'echet topology on the small Fr\'echet $U(\ft)$-module $D(G_0)_{P_{I,0}}\otimes_{D(P_{I,0})}X \cong D(N_I^+)_{1}\otimes_E X$ can be defined by the family of semi-norms $\{|\cdot|_r''\}_{r\in\cI}$. By Lemma~\ref{lem: support variant} $\cM$ is a quotient of the coadmissible $D(G_0)_{P_{I,0}}$-module $D(G_0)_{P_{I,0}}\otimes_{D(P_{I,0})}X$ and for $r\in \cI$ the $D(G_0)_{P_{I,0},r}$-module $\cM_r$ is a quotient of $D(G_0)_{P_{I,0},r}\otimes_{D(P_{I,0})}X$. In particular, $\cM_r$ is a quotient of $D(G_0)_{P_{I,0},r}\otimes_{D(P_{I,0})}X$ by a (closed) Banach subspace stable under the action of $U(\ft)$, so by Remark~\ref{rem: dense in Banach} the standard semi-norm $|\cdot|_r''$ on $U(\fg)\otimes_{U(\fp_I)}X$ induces a standard semi-norm $|\cdot|_r$ on its quotient $M$ so that $\cM_r$ is the completion of $M$ under $|\cdot|_r$.
\end{proof}

For $M$ in $\cO_{\rm{alg}}^{\fp_I}$, by Lemma~\ref{lem: support variant} the canonical completion $\cM$ of $M$ is a coadmissible left $D(G)_{P_I}$-module. For $\pi^{\infty}$ in $\mathrm{Rep}^{\infty}_{\rm{adm}}(L_I)$, recall that $(\pi^{\infty})^\vee$ is a coadmissible left $D(G)_{P_I}$-module via the surjections (where the middle isomorphism follows from Lemma~\ref{lem: smooth restriction})
\begin{equation*}
D(G)_{P_I}\twoheadrightarrow D^{\infty}(G)_{P_I}\cong D^{\infty}(P_I)\twoheadrightarrow D^{\infty}(L_I).
\end{equation*}
We will use the following lemma.

\begin{lem}\label{mistakecorrected}
Let $M$ in $\cO_{\rm{alg}}^{\fp_I}$, $\cM=D(G)_1\otimes_{U(\fg)}M$ its canonical Fr\'echet completion, $X$ a finite dimensional $U(\fp_I)$-module with a surjection $q: U(\fg) \otimes_{U(\fp_I)}X \twoheadrightarrow M$. Let $\pi^\infty$ in $\mathrm{Rep}^{\infty}_{\rm{adm}}(L_I)$ and $D\defeq (\pi^\infty)^\vee$.
\begin{enumerate}[label=(\roman*)]
\item \label{it: mistake 1} The $D(P_I)$-module $X \otimes_E D$ (with the diagonal action) is coadmissible and we have an isomorphism of coadmissible $D(G)_{P_I}$-modules
\begin{equation}\label{mistake1addendum}
D(G)_{P_I}\widehat{\otimes}_{D(P_I)} (X \otimes_E D) \buildrel\sim\over\longrightarrow (D(G)_{P_I} \otimes_{D(P_I)} X) \widehat{\otimes}_E D
\end{equation}
with the diagonal action of $D(G)_{P_I}$ on the right hand side.
\item \label{it: mistake 2} The $D(G)_{P_I}$-module $\cM\widehat{\otimes}_E D$ (with the diagonal action) is coadmissible and we have a short exact sequence of coadmissible $D(G)_{P_I}$-modules (with the diagonal action)
\begin{equation}\label{mistake2addendum}
0 \longrightarrow \big(D(G)_{P_I} \Ker(q)\big) \widehat{\otimes}_E D \longrightarrow (D(G)_{P_I} \otimes_{D(P_I)} X)\widehat\otimes_E D \longrightarrow \cM \widehat{\otimes}_E D \longrightarrow 0.
\end{equation}
\end{enumerate}
\end{lem}
\begin{proof}
Note first that for the diagonal actions, one uses the comultiplication map $D(G)_{P_I} \rightarrow D(G)_{P_I}\widehat\otimes_E D(G)_{P_I}$ deduced from (\ref{equ: loc an product}) as in \cite[\S A]{ST05}. We prove \ref{it: mistake 1}. By the proof of \cite[Lemma 2.4(i)]{OS15}, the $D(P_I)$-module $X\otimes_E D$ is coadmissible, which implies that $(X\otimes_E D)_r\defeq D(P_{I,0})_r\otimes_{D(P_{I,0})}(X\otimes_E D)$ is a finitely generated $D(P_{I,0})_r$-module for $r\in \cI$ (see \cite[\S 3]{OS15}). Moreover, using the universal property of the inductive tensor product (\cite[\S 17.A]{S02}), one has a topological isomorphism $(D(G_0)_{P_{I,0}}\otimes_{D(P_{I,0})}D(P_I))\otimes_{D(P_I)}(X \otimes_E D)\cong D(G_0)_{P_{I,0}}\otimes_{D(P_{I,0})}(X \otimes_E D)$. Then we have topological isomorphisms
\begin{eqnarray*}
D(G)_{P_I} \widehat\otimes_{D(P_I)} (X \otimes_E D)&\cong &D(G_0)_{P_{I,0}} \widehat\otimes_{D(P_{I,0})} (X \otimes_E D)\\
&\buildrel\sim\over\longrightarrow & \varprojlim_{r}\big(D(G_0)_{P_{I,0},r} \otimes_{D(P_{I,0})_r}(X\otimes_E D)_r\big)
\end{eqnarray*}
where the first isomorphism follows from $D(G)_{P_I}\cong D(G_0)_{P_{I,0}}\otimes_{D(P_{I,0})}D(P_I)$ and the fact that $D(G)_{P_I} \widehat\otimes_{D(P_I)} (X \otimes_E D)$ is the completion of $D(G)_{P_I} \otimes_{D(P_I)} (X \otimes_E D)\cong D(G_0)_{P_{I,0}}\otimes_{D(P_{I,0})}(X \otimes_E D)$, and the second from Lemma \ref{limprojtensor} and the beginning of the proof. Now, since we have for $r\leq r'$ in $\cI$
\begin{multline*}
D(G_0)_{P_{I,0},r}\otimes_{D(G_0)_{P_{I,0},r'}}\big(D(G_0)_{P_{I,0},r'} \otimes_{D(P_{I,0})_{r'}}(X\otimes_E D)_{r'}\big)\\
\begin{array}{cl}
\cong &D(G_0)_{P_{I,0},r}\otimes_{D(P_{I,0})_{r'}}(X\otimes_E D)_{r'}\\
\cong &D(G_0)_{P_{I,0},r}\otimes_{D(P_{I,0})_r}\big(D(P_{I,0})_r\otimes_{D(P_{I,0})_{r'}}(X\otimes_E D)_{r'}\big)\\
\cong &D(G_0)_{P_{I,0},r}\otimes_{D(P_{I,0})_r}(X\otimes_E D)_{r}
\end{array}
\end{multline*}
which is a finitely generated $D(G_0)_{P_{I,0},r}$-module, it follows from the definition of coadmissibility (\cite[p.~152]{OS15}) that $D(G)_{P_I} \widehat\otimes_{D(P_I)} (X \otimes_E D)$ is a coadmissible $D(G)_{P_I}$-module. The $D(P_I)$-equivariant embedding $X \hookrightarrow D(G)_{P_I} \otimes_{D(P_I)} X$, $x\mapsto 1 \otimes x$ induces a continuous $D(P_I)$-equivariant embedding $X \otimes_E D \hookrightarrow (D(G)_{P_I} \otimes_{D(P_I)} X) \widehat{\otimes}_E D$ (with diagonal action of $D(P_I)$ on both sides via \cite[\S A]{ST05}), which itself induces a continuous $D(G)_{P_I}$-equivariant~map
\begin{equation}\label{mistakeaddendum}
D(G)_{P_I} \widehat\otimes_{D(P_I)} (X \otimes_E D) \longrightarrow (D(G)_{P_I} \otimes_{D(P_I)} X) \widehat{\otimes}_E D
\end{equation}
(where we have used that the right hand side is complete thanks to the first paragraph of the proof of Lemma \ref{lem: support variant}). Hence it suffices to show that (\ref{mistakeaddendum}) is a topological isomorphism. Note that we also have a topological isomorphism $X \otimes_E D\buildrel\sim\over\rightarrow X \otimes_E (\varprojlim_{r}D_r)\cong \varprojlim_{r}(X \otimes_E D_r)$ where $D_r\defeq D(P_{I,0})_r\otimes_{D(P_{I,0})}D\cong D^{\infty}(P_0)_r\otimes_{D^{\infty}(P_0)}D$ (a finite dimensional $E$-vector space). Using Lemma \ref{limprojtensor} (noting that both $X$ and $X \otimes_E D_r$ are finitely generated $D(P_{I,0})_r$-modules) it suffices to prove that the natural morphism of $E$-Banach spaces
\begin{equation}\label{mistake1addendum2}
D(G_0)_{P_{I,0},r} \otimes_{D(P_{I,0})_r} (X \otimes_E D_r) \longrightarrow (D(G_0)_{P_{I,0},r} \otimes_{D(P_{I,0})_r} X) \otimes_E D_r
\end{equation}
is a topological isomorphism. By \ref{it: support basis 1} of Lemma \ref{lem: support basis} we have an isomorphism
\begin{equation}\label{DPR}
D(G_0)_{1,r} \otimes_{D(P_{I,0})_{1,r}} D(P_{I,0})_r\buildrel\sim\over\longrightarrow D(G_0)_{P_{I,0},r}
\end{equation}
(using \emph{loc.~cit.}~and the discussion above Lemma \ref{lem: free mod coset basis} one checks that both are free over $D(G_0)_{1,r}$ with same basis). So (\ref{mistake1addendum2}) can be rewritten as
\begin{equation}\label{mistake1addendum3}
D(G_0)_{1,r} \otimes_{D(P_{I,0})_{1,r}} (X \otimes_E D_r) \longrightarrow (D(G_0)_{1,r} \otimes_{D(P_{I,0})_{1,r}} X) \otimes_E D_r.
\end{equation}
But since $D(P_{I,0})_{1,r}$ acts on $D_r$ via the surjection $D(P_{I,0})_{1,r} \twoheadrightarrow D^\infty(P_{I,0})_{1,r}\cong E$, one trivially checks that (\ref{mistake1addendum3}) is an isomorphism of $E$-Banach spaces.

We prove \ref{it: mistake 2}. By Lemma \ref{lem: support variant} we have a short exact sequence in $\cC_{D(G)_{P_I}}
$
\[0 \longrightarrow D(G)_{P_I}\Ker(q) \longrightarrow D(G)_{P_I} \otimes_{D(P_I)} X \longrightarrow \cM \longrightarrow 0\]
which gives a short exact sequence of Fr\'echet spaces as in (\ref{mistake2addendum}) by \cite[Lemme 4.13]{Schr11}. So it suffices to show that $\cM \widehat{\otimes}_E D$ and $(D(G)_{P_I} \Ker(q)) \widehat{\otimes}_E D$ (with the diagonal $D(G)_{P_I}$-action) are coadmissible $D(G)_{P_I}$-modules. By \ref{it: mistake 1} we know that $(D(G)_{P_I} \otimes_{D(P_I)} X) \widehat{\otimes}_E D$ is a coadmissible $D(G)_{P_I}$-module. Since $(D(G)_{P_I} \Ker(q)) \widehat{\otimes}_E D$ is closed in $(D(G)_{P_I} \otimes_{D(P_I)} X) \widehat{\otimes}_E D$ (by the above (\ref{mistake2addendum})) and stable under $D(G)_{P_I}$, it is a coadmissible $D(G)_{P_I}$-module by \cite[Lemma 3.6]{ST03}, hence $\cM \widehat{\otimes}_E D$ is a coadmissible $D(G)_{P_I}$-module by \cite[Lemma 3.4.ii]{ST03}.
\end{proof}

Given $M$ in $\cO_{\rm{alg}}^{\fp_I}$ and $\pi^{\infty}$ in $\mathrm{Rep}^{\infty}_{\rm{adm}}(L_I)$, Orlik and Strauch define a representation $\cF_{P_I}^G(M,\pi^{\infty})$ in $\mathrm{Rep}^{\rm{an}}_{\rm{adm}}(G)$, we refer the reader to \cite{OS15} for details (see also Theorem \ref{prop: OS property} below). The following result gives a convenient description of the continuous dual $\cF_{P_I}^G(M,\pi^{\infty})^\vee$ of $\cF_{P_I}^G(M,\pi^{\infty})$ using the completion functor $M\mapsto \cM$ of Proposition~\ref{prop: category O completion}.

\begin{prop}\label{prop: OS via completion}
For $M$ in $\cO_{\rm{alg}}^{\fp_I}$, $\cM= D(G)_1\otimes_{U(\fg)}M$ and $\pi^{\infty}$ a smooth admissible representation of $L_I$ over $E$ we have a canonical isomorphism of coadmissible left $D(G)$-modules
\begin{equation}\label{equ: OS completion}
\cF_{P_I}^G(M,\pi^{\infty})^\vee \cong D(G)\widehat{\otimes}_{D(G)_{P_I}}(\cM\widehat{\otimes}_E (\pi^{\infty})^\vee).
\end{equation}
\end{prop}
\begin{proof}
We let $X$ be a finite dimensional left $U(\fp_I)$-module $X$ which is $\ft$-semi-simple and such that we have a surjection $q: U(\fg)\otimes_{U(\fp_I)}X\twoheadrightarrow M$ in $\cO_{\rm{alg}}^{\fp_I}$, and we define $D_1\defeq D(G)_{P_I}\mathrm{ker}(q)$, $D_2\defeq D(G)_{P_I}\otimes_{D(P_I)} X$ and $D_3\defeq \cM$. By (\ref{mistake2addendum}) (applied with $D=(\pi^{\infty})^\vee$) we have a short exact sequence in $\cC_{D(G)_{P_I}}$
\[0\longrightarrow D_1\widehat{\otimes}_E (\pi^{\infty})^\vee\longrightarrow D_2\widehat{\otimes}_E (\pi^{\infty})^\vee\longrightarrow D_3\widehat{\otimes}_E (\pi^{\infty})^\vee \longrightarrow 0.\]
Let $R_k\defeq D(G)\widehat{\otimes}_{D(G)_{P_I}}(D_k\widehat{\otimes}_E (\pi^{\infty})^\vee)$ for $1\leq k\leq 3$, then by \cite[Lemma 4.27]{Schr11} (the proof of which extends to our setting) we deduce a short exact sequence of (left) $D(G)$-modules
\begin{equation}\label{equ: OS def seq}
0\longrightarrow R_1\longrightarrow R_2 \longrightarrow R_3\longrightarrow 0
\end{equation}
(one could also again argue using Lemma \ref{limprojtensor} and \ref{it: support basis 2} of Lemma \ref{lem: support basis}). Let $(\pi^{\infty})^\vee_r\defeq D(P_0)_r\otimes_{D(P_0)}(\pi^{\infty})^\vee\cong D^{\infty}(P_0)_r\otimes_{D^{\infty}(P_0)}(\pi^{\infty})^\vee$, we have isomorphisms of coadmissible $D(G)_{P_I}$-modules:
\begin{equation*}
D_2\widehat{\otimes}_E (\pi^{\infty})^\vee\cong D(G)_{P_I}\widehat{\otimes}_{D(P_I)}(X\otimes_E (\pi^{\infty})^\vee)\buildrel\sim\over\longrightarrow
\varprojlim_{r\in\cI} \big(D(G_0)_{P_{I,0},r} \otimes_{D(P_{I,0})_r} (X \otimes_E (\pi^{\infty})^\vee_r)\big)
\end{equation*}
where the first isomorphism follows from (\ref{mistake1addendum}) and the second from the proof of \ref{it: mistake 1} of Lemma \ref{mistakecorrected}. Thus (using Lemma \ref{limprojtensor} again with \ref{it: support basis 2} of Lemma \ref{lem: support basis}) we have isomorphisms of topological $D(G)$-modules
\begin{eqnarray}\label{r2r2}
\nonumber R_2 &\cong &\varprojlim_{r\in\cI} \Big(D(G_0)_r\otimes_{D(G_0)_{P_{I,0},r}}\big(D(G_0)_{P_{I,0},r} \!\otimes_{D(P_{I,0})_r} (X \otimes_E (\pi^{\infty})^\vee_r)\big)\Big)\\
\nonumber &\cong & \varprojlim_{r\in\cI} \big(D(G_0)_r\otimes_{D(P_{I,0})_r} (X \otimes_E (\pi^{\infty})^\vee_r)\big)\\
&\cong & D(G)\widehat{\otimes}_{D(P_I)}(X\otimes_E (\pi^{\infty})^\vee)
\end{eqnarray}
where the last $D(G)$-module is coadmissible by \cite[Lemma 2.4(i)]{OS15} (recall we have an isomorphism of Fr\'echet spaces $R_2\cong ((\mathrm{Ind}_{P_I}^{G}(X^\vee\otimes_E \pi^{\infty}))^{\rm{an}})^\vee$ using e.g.~\cite[(52), (56), Rem.~5.4]{Koh11} with \cite[Lemme 3.1]{BD23}). It then follows from \cite[Lemma 3.6]{ST03} that $R_1$ and $R_3$ are also coadmissible $D(G)$-modules.

By \cite[(3.2.2),(4.4.1)]{OS15} we have a natural pairing
\begin{equation}\label{equ: OS pairing}
(U(\fg)\otimes_{U(\fp_I)}X)\otimes_E \big(\mathrm{Ind}_{P_I}^{G}(X^\vee\otimes_E \pi^{\infty})\big)^{\rm{an}} \longrightarrow C^{\rm{an}}(G,\pi^{\infty})
\end{equation}
which sends $(\delta\otimes_E x)\otimes_E f$ (for $\delta\in U(\fg)$, $x\in X$ and $f\in(\mathrm{Ind}_{P_I}^{G}(X^\vee\otimes_E \pi^{\infty}))^{\rm{an}}$) to
\begin{equation*}
\big[g\mapsto \delta\big(h\mapsto f(gh)(x)\big)\big]\in C^{\rm{an}}(G,\pi^{\infty})
\end{equation*}
with $f(gh)(x)\in\pi^{\infty}$, $[h\mapsto f(gh)(x)]\in C^{\rm{an}}(G,\pi^{\infty})$ and $\delta(h\mapsto f(gh)(x))\in \pi^{\infty}$.
Let $g\in G$, composing (\ref{equ: OS pairing}) with the evaluation map $C^{\rm{an}}(G,\pi^{\infty})\rightarrow \pi^{\infty}, \varphi\mapsto \varphi(g)$ induces a pairing
\begin{equation}\label{equ: OS pairing 1}
(U(\fg)\otimes_{U(\fp_I)}X)\otimes_E \big(\mathrm{Ind}_{P_I}^{G}(X^\vee\otimes_E \pi^{\infty})\big)^{\rm{an}} \longrightarrow \pi^{\infty}
\end{equation}
which sends $(\delta\otimes_E x)\otimes_E f$ to $\delta(h\mapsto f(gh)(x))$. By pairing with $(\pi^{\infty})^\vee$ on both sides of (\ref{equ: OS pairing 1}), we obtain the pairing
\begin{equation}\label{equ: OS pairing 2}
\big((U(\fg)\otimes_{U(\fp_I)}X)\otimes_E (\pi^{\infty})^\vee\big)\otimes_E \big(\mathrm{Ind}_{P_I}^{G}(X^\vee\otimes_E \pi^{\infty})\big)^{\rm{an}} \longrightarrow E
\end{equation}
which sends $((\delta\otimes_E x)\otimes_E\delta^{\infty})\otimes_E f$ to $\delta^{\infty}(\delta(h\mapsto f(gh)(x)))$ for $\delta,x,f$ as above and $\delta^{\infty}\in(\pi^{\infty})^\vee$. From the isomorphism $R_2\cong ((\mathrm{Ind}_{P_I}^{G}(X^\vee\otimes_E \pi^{\infty}))^{\rm{an}})^\vee$, we have a perfect pairing
\begin{equation}\label{equ: OS induction dual}
R_2\otimes_E (\mathrm{Ind}_{P_I}^{G}X^\vee\otimes_E \pi^{\infty})^{\rm{an}}\longrightarrow E.
\end{equation}
Then (\ref{equ: OS pairing 2}) can be reinterpreted as the pairing (recall $\delta_g\in D(G)$ is the Dirac distribution)
\begin{equation*}
(\delta_g\cdot((U(\fg)\otimes_{U(\fp_I)}X)\otimes_E (\pi^{\infty})^\vee))\otimes_E (\mathrm{Ind}_{P_I}^{G}X^\vee\otimes_E \pi^{\infty})^{\rm{an}} \longrightarrow E
\end{equation*}
obtained from (\ref{equ: OS induction dual}) via composition with the natural map
\[\delta_g\cdot((U(\fg)\otimes_{U(\fp_I)}X)\otimes_E (\pi^{\infty})^\vee)\longrightarrow R_2\cong D(G)\widehat{\otimes}_{D(P_I)}(X\otimes_E (\pi^{\infty})^\vee).\]
The representation $\cF_{P_I}^G(M,\pi^{\infty})$ is defined in \cite[(4.4.1)]{OS15} as the closed (invariant) subspace of $(\mathrm{Ind}_{P_I}^{G}X^\vee\otimes_E \pi^{\infty})^{\rm{an}}$ which is the orthogonal of $\mathrm{ker}(q)$ under the pairing (\ref{equ: OS pairing}). Since $\varphi\in C^{\rm{an}}(G,\pi^{\infty})$ is zero if and only if $\varphi(g)=0$ for all $g\in G$, we deduce that $\cF_{P_I}^G(M,\pi^{\infty})$ is the closed subspace of $(\mathrm{Ind}_{P_I}^{G}X^\vee\otimes_E \pi^{\infty})^{\rm{an}}$ which is the orthogonal of the image of $\sum_{g\in G}\delta_g\cdot(\mathrm{ker}(q)\otimes_E (\pi^{\infty})^\vee)$ in $R_2$ under the pairing (\ref{equ: OS induction dual}), that is, $\cF_{P_I}^G(M,\pi^{\infty})^\vee$ is the quotient of $R_2$ by the closure of
\begin{equation}\label{equ: sum of translation}
\sum_{g\in G}\delta_g\cdot(\mathrm{ker}(q)\otimes_E (\pi^{\infty})^\vee).
\end{equation}
As $\mathrm{ker}(q)$ is dense in $D_1$ (see the last sentence of the proof of Lemma \ref{lem: support variant}), $\mathrm{ker}(q)\otimes_E (\pi^{\infty})^\vee$ is dense in $D_1\widehat{\otimes}_E (\pi^{\infty})^\vee$, and as $\{\delta_g\mid g\in G\}$ is dense in $D(G)$, one easily checks that (\ref{equ: sum of translation}) is dense inside $R_1=D(G)\widehat{\otimes}_{D(G)_{P_I}}(D_1\widehat{\otimes}_E (\pi^{\infty})^\vee)\subseteq R_2$ (see (\ref{equ: OS def seq})). Hence, we deduce from \emph{loc.~cit.}~that
\[\cF_{P_I}^G(M,\pi^{\infty})^\vee\cong R_2/R_1\cong R_3=D(G)\widehat{\otimes}_{D(G)_{P_I}}(\cM\widehat{\otimes}_E (\pi^{\infty})^\vee).\qedhere\]
\end{proof}

We recall for later convenience the main theorem of \cite{OS15} in the following form.

\begin{thm}\label{prop: OS property}
\hspace{2em}
\begin{enumerate}[label=(\roman*)]
\item \label{it: OS property 1} The functor $\cF_{P_I}^G(-,-)$ is contravariant (resp.~covariant) in the first (resp.~the second) argument and is exact for both arguments.
\item \label{it: OS property 2} For $I_1\subseteq I\subseteq \Delta$, $M$ in $\cO^{\fp_I}_{\rm{alg}}$ and $\pi_1^{\infty}$ in $\mathrm{Rep}^{\infty}_{\rm{adm}}(L_{I_1})$, we have a canonical isomorphism
\[\cF_{P_{I_1}}^G(M,\pi_1^{\infty})\cong \cF_{P_I}^G(M,i_{I_1,I}^{\infty}(\pi_1^{\infty})).\]
\item \label{it: OS property 3} Let $I\subseteq \Delta$, $M$ in $\cO^{\fp_I}_{\rm{alg}}$ with $I$ maximal for $M$ and $\pi^{\infty}$ in $\mathrm{Rep}^{\infty}_{\rm{adm}}(L_I)$. If both $M$ and $\pi^{\infty}$ are simple, then so is $\cF_{P_I}^G(M,\pi^{\infty})$.
\item \label{it: OS property 4} Let $I\subseteq \Delta$, $M$ in $\cO^{\fp_I}_{\rm{alg}}$ and $\pi^{\infty}$ in $\mathrm{Rep}^{\infty}_{\rm{adm}}(L_I)$. Then a Jordan-H\"older factor of $\cF_{P_I}^G(M,\pi^{\infty})$ has the form $\cF_{P_{I_1}}^G(M_1,\pi_1^{\infty})$ for some Jordan-H\"older factor $M_1$ of $M$ such that $I_1\supseteq I$ is maximal for $M_1\in\cO^{\fp_{I_1}}_{\rm{alg}}$, and some Jordan-H\"older factor $\pi_1^{\infty}$ of $i_{I,I_1}^{\infty}(\pi^{\infty})$.
\end{enumerate}
\end{thm}

We finish this section with several remarks.

\begin{rem}\label{rem: non t semisimple}
For $I\subseteq \Delta$, recall the subcategory $\cO^{\fp_I,\infty}_{\rm{alg}}\subseteq \tld{\cO}_{\rm{alg}}^{\fb}\cap \cC_{\rm{alg}}^{\rm{fin}}$ from the discussion at the beginning of \S \ref{subsec: category}. Recall also from \cite{AS22} that, at least for smooth strongly admissible representations $\pi^{\infty}$ of $L_I$ over $E$, the construction of $\cF_{P_I}^G(M,\pi^{\infty})$ can be extended to $M$ in $\cO^{\fp_I,\infty}_{\rm{alg}}$ (modulo a choice of logarithm that we ignore here). For $M$ in $\cO^{\fp_I,\infty}_{\rm{alg}}$, we can still define its canonical Fr\'echet completion as $\cM\defeq D(G)_1\otimes_{U(\fg)}M$. Let $\widehat{\cO}^{\fp_I,\infty}_{\rm{alg}}\subseteq \cC_{D(G)_1}$ be the abelian full subcategory of coadmissible $D(G)_1$-modules $D$ such that $D$ is a (generalized) small Fr\'echet $U(\ft)$-module in the sense of Remark~\ref{rem: non semi simple case} and $\bigoplus_{\mu\in\Lambda} D_{\mu}$ lies in $\cO^{\fp_I,\infty}_{\rm{alg}}$. Then using Remark~\ref{rem: non semi simple case} one can check that the functor $M\mapsto \cM$ induces an equivalence of categories $\cO^{\fp_I,\infty}_{\rm{alg}}\xrightarrow{\sim}\widehat{\cO}^{\fp_I,\infty}_{\rm{alg}}$. Lemma~\ref{lem: support variant}, Proposition~\ref{prop: OS via completion} and Theorem~\ref{prop: OS property} also remain true for $M$ in $\cO^{\fp_I,\infty}_{\rm{alg}}$.
\end{rem}

\begin{rem}\label{addendumstrongly}
If $I\subseteq \Delta$, $\pi^{\infty}$ is a smooth strongly admissible representations of $L_I$ over $E$ and $X$ is a finite dimensional left $U(\fp_I)$-module, then by Lemma~\ref{lem: loc alg FIN} $X\otimes_E (\pi^\infty)^\vee$ admits a resolution by finite free $D(P_{I,0})$-modules (we use here the strong admissibility of $\pi^{\infty}$). Arguing as in the first paragraph of the proof of Lemma \ref{lem: support variant} we then deduce that $D(G)_{P_I} \otimes_{D(P_I)} (X \otimes_E D)$ is a coadmissible $D(G)_{P_{I}}$-module, in particular we have a topological isomorphism $D(G)_{P_I} \otimes_{D(P_I)} (X \otimes_E (\pi^{\infty})^\vee)\buildrel\sim\over\rightarrow D(G)_{P_I} \widehat{\otimes}_{D(P_I)} (X \otimes_E (\pi^{\infty})^\vee)$ (so no need to complete on the right). By the same kind of arguments we also have $R_2\cong D(G) {\otimes}_{D(P_I)} (X \otimes_E (\pi^{\infty})^\vee)$ in (\ref{r2r2}). However, it is not clear to us if we can replace $D(G)\widehat{\otimes}_{D(G)_{P_I}}(-)$ by $D(G){\otimes}_{D(G)_{P_I}}(-)$ in (\ref{equ: OS completion}).
\end{rem}

\begin{rem}
Let $E'$ be a finite extension of $E$. Given two $D(G)=D(G,E)$-modules $D_1$ and $D_2$, we may consider the $D(G,E')=D(G)\otimes_E E'$-modules $D_1\otimes_E E'$ and $D_2\otimes_E E'$. Note that any injective resolution of $D_1\otimes_E E'$ in $\mathrm{Mod}_{D(G,E')}$ restricts to an injective resolution of $D_1\otimes_E E'$ in $\mathrm{Mod}_{D(G)}$, and any injective resolution of $D_1$ in $\mathrm{Mod}_{D(G)}$ remains by scalar extension to $E'$ an injective resolution of $D_1\otimes_E E'$ seen in $\mathrm{Mod}_{D(G)}$. We thus have canonical isomorphisms for $k\geq 0$
\begin{equation*}
\mathrm{Ext}_{D(G)}^k(D_2,D_1)\otimes_E E'\buildrel\sim\over\rightarrow \mathrm{Ext}_{D(G)}^k(D_2,D_1\otimes_E E')\buildrel\sim\over\rightarrow\mathrm{Ext}_{D(G,E')}^k(D_2\otimes_E E', D_1\otimes_E E').
\end{equation*}
Moreover, one has the following results.
\begin{enumerate}[label=(\roman*)]
\item \label{it: loc an rational 1} Let $I\subseteq \Delta$, $M$ in $\cO^{\fp_I}_{\rm{alg}}$, $\pi^{\infty}$ in $\mathrm{Rep}^{\infty}_{\rm{adm}}(L_I)$ and $D\defeq \cF_{P_I}^{G}(M,\pi^{\infty})^\vee$, then we have $D\otimes_E E'\cong \cF_{P_I}^{G}(M\otimes_E E', \pi^{\infty}\otimes_E E')^\vee$ as coadmissible $D(G,E')$-modules. If moreover $\pi^{\infty}\in\cB^I_{\Sigma}$ for some left $W(L_I)$-coset $\Sigma\subseteq \widehat{T}^{\infty}$ (see above (\ref{equ: sm block decomposition})), then it follows from Remark \ref{rem: rational JH} (using the last statement of Remark \ref{rem: support of PS}) that $D$ is a simple $D(G)$-module if and only if $D\otimes_E E'$ is a simple $D(G,E')$-module.
\item Let $D$ be a multiplicity free finite length $D(G)$-module with irreducible constituents all of the form $\cF_{P_I}^{G}(L(\mu),\pi^{\infty})^\vee$ for some $I\subseteq \Delta$, $\mu\in\Lambda_I^{\dom}$, $\pi^{\infty}\in\cB^I_{\Sigma}$ for some left $W(L_I)$-coset $\Sigma\subseteq \widehat{T}^{\infty}$. Then from \ref{it: loc an rational 1} the scalar extension $(-)\otimes_E E'$ induces a bijection of partially-ordered sets
\[\mathrm{JH}_{D(G)}(D)\buildrel\sim\over\longrightarrow \mathrm{JH}_{D(G,E')}(D\otimes_E E').\]
\end{enumerate}
\end{rem}

\subsection{Orlik-Strauch representations and the Bruhat filtration}\label{gradedpiecessubsection}

Using the results of \S\ref{subsec: completion O} we describe explicitly the duals of the graded pieces of the Bruhat filtration on Orlik-Strauch representations.\bigskip

For $i=0,1$ we fix $I_i\subseteq \Delta$ and recall that $W^{I_0,I_1}$ is the set of minimal length representatives of $W(L_{I_0})\backslash W(G)/W(L_{I_1})$. We have the Bruhat decomposition (\cite[Lemma 5.5]{DM91})
\begin{equation}\label{equ: Bruhat stratification}
G=\bigsqcup_{w\in W^{I_0,I_1}} P_{I_1}w^{-1}P_{I_0}.
\end{equation}
For $w\in W^{I_0,I_1}$ we write $S_w\defeq P_{I_1}w^{-1}P_{I_0}$ and $\overline{S_w}$ its closure in $G$. The following lemma is surely well-known, but we couldn't find a proof.

\begin{lem}\label{closurepar}
We have
\begin{equation*}
\overline{S_w}=\bigsqcup_{w,w'\in W^{I_0,I_1},\ w'\leq w}S_{w'}.
\end{equation*}
\end{lem}
\begin{proof}
By \cite[Lemma 2.12(1)]{EK23} $W(L_{I_1})w^{-1}W(L_{I_0})$ is an interval in $W(G)$, hence there is a unique maximal element $w_{\rm max}$ in $W(L_{I_1})w^{-1}W(L_{I_0})$ for the Bruhat order (and also a unique minimal element which is $w^{-1}$). It follows that all cosets $Bw'B$ appearing in $S_w$ are such that $w'\leq w_{\rm max}$ and hence are in the closure $\overline {Bw_{\rm max}B}$ of $Bw_{\rm max}B$. This implies ${S_w}\subseteq \overline {Bw_{\rm max}B}$ and hence $\overline{S_w}\subseteq \overline {Bw_{\rm max}B}$. But since $Bw_{\rm max}B\subseteq S_w$, we have $\overline{S_w}=\overline {Bw_{\rm max}B}=\bigsqcup_{w'\leq w_{\rm max}}Bw'B$. Since $P_{I_1}\overline{S_w}P_{I_0}=\overline{S_w}$, it follows that $\overline{S_w}=\bigcup_{w'\leq w_{\rm max}}P_{I_1}w'P_{I_0}$. Writing $w_{\rm max}=x_{{\rm max}, 1}w^{-1}x_{{\rm max}, 0}$ with $x_{{\rm max}, i}\in W(L_{I_i})$ such that $x_{{\rm max}, 1}w^{-1}$ is minimal in $x_{{\rm max}, 1}w^{-1}W(L_{I_0})$ (\cite[Lemma 5.4(iii)]{DM91}), by the subword property of the Bruhat order and \cite[Prop.~2.4.4]{BB05} any $w'\leq w_{\rm max}$ can be written $w'=x'_1w''x'_0$ with $x'_i\leq x_{{\rm max}, i}$ in $W(L_{I_i})$ and $w''\leq w^{-1}$. Thus $\overline{S_w}=\bigcup_{w''\leq w^{-1}}P_{I_1}w''P_{I_0}$, or equivalently $\overline{S_w}=\bigcup_{w''\leq w}P_{I_1}{w''}^{-1}P_{I_0}$. The result follows replacing $w''$ by its minimal representative in $P_{I_0}{w''}P_{I_1}$.
\end{proof}

In particular it follows from Lemma \ref{closurepar} and its proof that $S_w$ is open in its closure $\overline{S_w}$, i.e.~$S_w$ is locally closed in $G$. The closed subsets $(\overline{S_w})_{w\in W^{I_0,I_1}}$ of $G$ induce an exhaustive (by (\ref{equ: Bruhat stratification})) left $D(P_{I_1})$-equivariant and right $D(P_{I_0})$-equivariant filtration
\begin{equation*}
\mathrm{Fil}_w(D(G))\defeq D(G)_{\overline{S_w}}\subseteq D(G)
\end{equation*}
on $D(G)$ indexed by $W^{I_0,I_1}$ (using Lemma \ref{closurepar}). Note that by \emph{loc.~cit.}~we have $\mathrm{Fil}_1(D(G))=D(G)_{P_{I_1}P_{I_0}}\subseteq \mathrm{Fil}_w(D(G))$ for any $w\in W^{I_0,I_1}$. The terminology ``filtration'' comes from $\mathrm{Fil}_{w'}(D(G))\subseteq \mathrm{Fil}_w(D(G))$ whenever $w'\leq w$. For $w\in W^{I_0,I_1}$ recall that $\widehat{D}(G)_{S_w}$ was defined in (\ref{equ: loc an hat}) (using (\ref{projc0})), and we have
\begin{equation*}
\mathrm{gr}_w(D(G))\defeq \mathrm{Fil}_w(D(G))\big/\!\sum_{w'<w}\mathrm{Fil}_{w'}(D(G))
\cong D(G)_{\overline{S_w}}/D(G)_{\overline{S_w}\setminus S_w}\cong \widehat{D}(G)_{S_w}
\end{equation*}
where the first isomorphism follows from Lemma \ref{closurepar} and the second from (\ref{equ: limit vers quotient}).\bigskip

Let $M_0$ in $\cO_{\rm{alg}}^{\fp_{I_0}}$, $\pi_0^{\infty}$ a smooth admissible representation of $L_{I_0}$ over $E$ and let $V_0\defeq \cF_{P_{I_0}}^G(M_0,\pi_0^{\infty})$. Recall from Proposition~\ref{prop: OS via completion} that
\begin{equation}\label{equ: OS completion 0}
V_0^\vee\cong D(G)\widehat{\otimes}_{D(G)_{P_{I_0}}}(\cM_0\widehat{\otimes}_E(\pi_0^{\infty})^\vee)
\end{equation}
where $\cM_0=D(G)_1\otimes_{U(\fg)}M_0$ is a coadmissible $D(G)_{P_{I_0}}$-module (Lemma~\ref{lem: support variant}).
We define the following filtration indexed by $w\in W^{I_0,I_1}$ on the $D(P_{I_1})$-module $V_0^\vee$:
\begin{equation}\label{equ: Bruhat filtration OS}
\mathrm{Fil}_w(V_0^\vee)\defeq D(G)_{\overline{S_w}}\widehat{\otimes}_{D(G)_{P_{I_0}}}(\cM_0\widehat{\otimes}_E(\pi_0^{\infty})^\vee).
\end{equation}

\begin{lem}\label{closedw}
For $w\in W^{I_0,I_1}$ the inclusion $D(G)_{\overline{S_w}}\subseteq D(G)$ induces a closed embedding of left $D(G)_{P_{I_1}}$-modules $\mathrm{Fil}_w(V_0^\vee)\hookrightarrow V_0^\vee$.
\end{lem}
\begin{proof}
It is clear from (\ref{equ: OS completion 0}) and (\ref{equ: Bruhat filtration OS}) that the natural map $\mathrm{Fil}_w(V_0^\vee)\rightarrow V_0^\vee$ is $D(G)_{P_{I_1}}$-equivariant (as $\overline{S_w}=P_{I_1}\overline{S_w}$). Hence we have to prove that it is a closed embedding. Since $\cM_0\widehat{\otimes}_E(\pi_0^{\infty})^\vee$ is a coadmissible $D(G)_{P_{I_0}}$-module by \ref{it: mistake 2} of Lemma \ref{mistakecorrected}, we have
\begin{equation}\label{limprojmpi0}
\cM_0\widehat{\otimes}_E(\pi_0^{\infty})^\vee\buildrel\sim\over\longrightarrow \varprojlim_{r\in \cI}(\cM_0\widehat{\otimes}_E(\pi_0^{\infty})^\vee)_r
\end{equation}
where $(\cM_0\widehat{\otimes}_E(\pi_0^{\infty})^\vee)_r\defeq D(G_0)_{P_{I_0,0},r}\otimes_{D(G_0)_{P_{I_0,0}}}(\cM_0\widehat{\otimes}_E(\pi_0^{\infty})^\vee)$ and $P_{I_0,0}\defeq P_{I_0}\cap G_0$. Set $Y\defeq \overline{S_w}$ and $Y_0\defeq Y\cap G_0$, we have a topological isomorphisms by (\ref{P0toP})
\begin{multline*}
D(G)_Y\cong D(G_0)_{Y_0}\otimes_{D(P_{I_0,0})}D(P_{I_0})\cong D(G_0)_{Y_0}\otimes_{D(G_0)_{P_{I_0,0}}}\big(D(G_0)_{P_{I_0,0}}\otimes_{D(P_{I_0,0})}D(P_{I_0})\big)\\
 \cong D(G_0)_{Y_0}\otimes_{D(G_0)_{P_{I_0,0}}}D(G)_{P_{I_0}}.
 \end{multline*}
From the universal property of the inductive tensor product, we obtain a topological isomorphism
\begin{equation}\label{forlater}
D(G)_{Y}{\otimes}_{D(G)_{P_{I_0}}}(\cM_0\widehat{\otimes}_E(\pi_0^{\infty})^\vee)\cong D(G_0)_{Y_0}{\otimes}_{D(G_0)_{P_{I_0,0}}}(\cM_0\widehat{\otimes}_E(\pi_0^{\infty})^\vee),
\end{equation}
hence (taking completions) a topological isomorphism
\begin{equation}\label{GtoG0}
\mathrm{Fil}_w(V_0^\vee)\cong D(G_0)_{Y_0}{\widehat\otimes}_{D(G_0)_{P_{I_0,0}}}(\cM_0\widehat{\otimes}_E(\pi_0^{\infty})^\vee).
\end{equation}
By \ref{it: support basis 2} of Lemma \ref{lem: support basis} $D(G_0)_r$ and $D(G_0)_{Y_0,r}$ are finite free over $D(G_0)_{P_{I_0,0},r}$. From (\ref{limprojmpi0}), (\ref{GtoG0}) and Lemma \ref{limprojtensor} we deduce isomorphisms of Fr\'echet spaces
\begin{eqnarray}\label{Filwlimproj}
\mathrm{Fil}_w(V_0^\vee)&\cong &\varprojlim_{r}\big(D(G_0)_{Y_0,r}\otimes_{D(G_0)_{P_{I_0,0},r}}(\cM_0\widehat{\otimes}_E(\pi_0^{\infty})^\vee)_r)\\
\nonumber V_0^\vee&\cong &\varprojlim_{r}\big(D(G_0)_{r}\otimes_{D(G_0)_{P_{I_0,0},r}}(\cM_0\widehat{\otimes}_E(\pi_0^{\infty})^\vee)_r).
\end{eqnarray}
The result follows from the fact that $D(G_0)_{Y_0,r}$ is a direct summand of $D(G_0)_r$ as $D(G_0)_{P_{I_0,0},r}$-module by \ref{it: support basis 2} of Lemma \ref{lem: support basis} again.
\end{proof}

For $w\in W^{I_0,I_1}$, we define $P_w\defeq w^{-1}P_{I_0}w$ and $L_w\defeq w^{-1}L_{I_0}w$. Recall that $\pi_0^{\infty,w}$ is the representation of $L_w$ with the same underlying vector space as $\pi_0^{\infty}$ but where $h\in L_w=w^{-1}L_{I_0}w$ acts by $whw^{-1}$ (by inflation $\pi_0^{\infty,w}$ is also a representation of $P_w$). Likewise we write $\cM_0^w$ for the coadmissible $D(G)_{P_w}$-module with the same underlying (topological) vector space as $\cM_0$ but where $\delta\in D(G)_{P_w}$ acts by $w \delta w^{-1}\in D(G)_{P_{I_0}}$.

\begin{prop}\label{lem: graded OS tensor}
For $w\in W^{I_0,I_1}$ let $\mathrm{gr}_w(V_0^\vee)\defeq \mathrm{Fil}_w(V_0^\vee)/\sum_{w'<w}\mathrm{Fil}_{w'}(V_0^\vee)$ (via Lemma \ref{closedw}), we have a canonical (left) $D(P_{I_1})$-equivariant isomorphism of Fr\'echet spaces
\begin{equation}\label{equ: graded OS tensor}
\mathrm{gr}_w(V_0^\vee)\cong \widehat{D}(G)_{P_{I_1}P_w}\widehat{\otimes}_{D(G)_{P_w}}(\cM_0^w\widehat{\otimes}_E(\pi_0^{\infty,w})^\vee)
\end{equation}
where $\widehat{D}(G)_{P_{I_1}P_w}$ is defined in (\ref{equ: loc an hat}) with its left $D(P_{I_1})$-action as before lemma \ref{lem: limit vers quotient}.
\end{prop}
\begin{proof}
For $w\in W^{I_0,I_1}$ let $P_{w,0}\defeq P_w\cap G_0$, $S\subseteq G$ a closed subset such that $S=SP_w=P_{I_1}S$ and $S_0\defeq S\cap G_0$. Arguing as for (\ref{GtoG0}) and (\ref{Filwlimproj}) we have $D(G)_{P_{I_1}}$-equivariant isomorphisms of Fr\'echet spaces (where $(\cM_0^w\widehat{\otimes}_E(\pi_0^{\infty,w})^\vee)_r\defeq D(G_0)_{P_{w,0},r}\otimes_{D(G_0)_{P_{w,0}}}(\cM_0^w\widehat{\otimes}_E(\pi_0^{\infty,w})^\vee)$)
\begin{eqnarray}\label{equ: grade OS tensor compact}
D(G)_{S\cdot w^{-1}}\widehat{\otimes}_{D(G)_{P_{I_0}}}(\cM_0\widehat{\otimes}_E(\pi_0^{\infty})^\vee)
&\cong &D(G)_{S}\widehat{\otimes}_{D(G)_{P_w}}(\cM_0^w\widehat{\otimes}_E(\pi_0^{\infty,w})^\vee)\\
\nonumber &\cong &D(G_0)_{S_0}\widehat{\otimes}_{D(G_0)_{P_{w,0}}}(\cM_0^w\widehat{\otimes}_E(\pi_0^{\infty,w})^\vee)\\
\nonumber &\cong &\varprojlim_{r}\big(D(G_0)_{S_0,r}\otimes_{D(G_0)_{P_{w,0},r}}(\cM_0^w\widehat{\otimes}_E(\pi_0^{\infty,w})^\vee)_r).
\end{eqnarray}
Let $X\defeq S_w w=P_{I_1}P_w\subseteq G$, $Y\defeq \overline{X}$ and $Z\defeq \overline{X}\setminus X$, applying these isomorphisms with $S=Y$ and $S=Z$ and using (\ref{equ: Banach seq}), (\ref{equ: limit Frechet seq}), (\ref{equ: Banach seqlim}) and Lemma~\ref{lem: top ML} we obtain an exact sequence of Fr\'echet spaces
\begin{multline*}
0\longrightarrow D(G_0)_{Z_0}\widehat{\otimes}_{D(G_0)_{P_{w,0}}}(\cM_0^w\widehat{\otimes}_E(\pi_0^{\infty,w})^\vee)\longrightarrow D(G_0)_{Y_0}\widehat{\otimes}_{D(G_0)_{P_{w,0}}}(\cM_0^w\widehat{\otimes}_E(\pi_0^{\infty,w})^\vee)\\
\longrightarrow \widehat{D}(G_0)_{X_0}\widehat{\otimes}_{D(G_0)_{P_{w,0}}}(\cM_0^w\widehat{\otimes}_E(\pi_0^{\infty,w})^\vee)\longrightarrow 0.
\end{multline*}
Using the definition of $\widehat{D}(G)_{X}$ (see (\ref{equ: loc an hat})) and arguing as before (\ref{GtoG0}), we deduce a $D(P_I)$-equivariant exact sequence of Fr\'echet spaces
\begin{multline*}
0\longrightarrow D(G)_{Z}\widehat{\otimes}_{D(G)_{P_{w}}}(\cM_0^w\widehat{\otimes}_E(\pi_0^{\infty,w})^\vee)\longrightarrow D(G)_{Y}\widehat{\otimes}_{D(G)_{P_{w}}}(\cM_0^w\widehat{\otimes}_E(\pi_0^{\infty,w})^\vee)\\
\longrightarrow \widehat{D}(G)_{X}\widehat{\otimes}_{D(G)_{P_{w}}}(\cM_0^w\widehat{\otimes}_E(\pi_0^{\infty,w})^\vee)\longrightarrow 0.
\end{multline*}
By (\ref{GtoG0}) and using (\ref{equ: grade OS tensor compact}) ``backwards'' together with Lemma \ref{closurepar}, it remains to prove
\[D(G_0)_{Z'_0}\widehat{\otimes}_{D(G_0)_{P_{I_0,0}}}(\cM_0\widehat{\otimes}_E(\pi_0^{\infty})^\vee)\cong \sum_{w'<w}\mathrm{Fil}_{w'}(V_0^\vee)\]
where $Z'_0\defeq \cup_{w'<w}(\overline{S_{w'}}\cap G_0)$. Let $D\defeq \cM_0\widehat{\otimes}_E(\pi_0^{\infty})^\vee$, it is enough to prove that if $C_1, C_2$ are two closed subsets of $G_0$ such that $C_i=C_iP_{I_0}$ then we have inside $D(G_0)\widehat{\otimes}_{D(G_0)_{P_{I_0,0}}}D$:
\[D(G_0)_{C_1}\widehat{\otimes}_{D(G_0)_{P_{I_0,0}}}D + D(G_0)_{C_2}\widehat{\otimes}_{D(G_0)_{P_{I_0,0}}}D=D(G_0)_{C_1\cup C_2}\widehat{\otimes}_{D(G_0)_{P_{I_0,0}}}D.\]
Using \ref{it: support basis 2} of Lemma \ref{lem: support basis}, for $r\in \cI$ we have an isomorphism of finite free $D(G_0)_{P_{I_0,0},r}$-modules $D(G_0)_{C_1,r}+ D(G_0)_{C_2,r}\buildrel\sim\over\rightarrow D(G_0)_{C_1\cup C_2,r}$ which gives an isomorphism of Banach spaces where $D_r\defeq D(G_0)_{P_{I_0,0},r}\otimes_{D(G_0)_{P_{I_0,0}}}D$
\[(D(G_0)_{C_1,r}+ D(G_0)_{C_2,r})\otimes_{D(G_0)_{P_{I_0,0},r}}D_r\buildrel\sim\over \longrightarrow D(G_0)_{C_1\cup C_2,r}\otimes_{D(G_0)_{P_{I_0,0},r}}D_r.\]
Using Lemma \ref{limprojtensor}, it is enough to prove
\[\varprojlim_{r\in \cI}(D(G_0)_{C_1,r}+ D(G_0)_{C_2,r})\cong D(G_0)_{C_1}+D(G_0)_{C_2}\]
(inside $\varprojlim_{r}D(G_0)_{r} \cong D(G_0)$). The image of $D(G_0)_{C_i}$ is dense in $D(G_0)_{C_i,r}$, hence the image of $D(G_0)_{C_1}+D(G_0)_{C_2}$ is dense in $D(G_0)_{C_1,r}+ D(G_0)_{C_2,r}$ (inside $D(G_0)_{r}$). As $D(G_0)_{C_1,r}+ D(G_0)_{C_2,r}$ is closed in $D(G_0)_r$ (see above), it follows that $\varprojlim_{r}(D(G_0)_{C_1,r}+ D(G_0)_{C_2,r})$ is the closure of $D(G_0)_{C_1}+D(G_0)_{C_2}$ inside $D(G_0)$. Hence it is enough to prove $D(G_0)_{C_1}+D(G_0)_{C_2}\cong D(G_0)_{C_1\cap C_2}$ (the latter being closed in $D(G_0)$). Let $U_1, U_2$ be compact open subsets of $G_0$ containing respectively $C_1, C_2$, then we have a short exact sequence of locally convex $E$-vector spaces of compact type
\[0\longrightarrow C^{\rm{an}}(U_1\cup U_2) \longrightarrow C^{\rm{an}}(U_1) \oplus C^{\rm{an}}(U_2) \longrightarrow C^{\rm{an}}(U_1\cap U_2)\longrightarrow 0\]
where the maps are the restrictions. Taking the colimit over such $U_1,U_2$, we deduce a short exact sequence of locally convex $E$-vector spaces of compact type
\begin{equation*}
0\longrightarrow \varinjlim_{U_1,U_2} C^{\rm{an}}(U_1\cup U_2) \longrightarrow \varinjlim_{U_1} C^{\rm{an}}(U_1) \oplus \varinjlim_{U_2} C^{\rm{an}}(U_2) \longrightarrow \varinjlim_{U_1,U_2} C^{\rm{an}}(U_1\cup U_2) \longrightarrow 0.
\end{equation*}
In particular the injection on the left is a closed embedding. Noting that compact open subsets of $G_0$ of the form $U_1 \cup U_2$, where $U_i$ are compact open subsets of $G_0$ containing $C_i$, are cofinal among compact open subsets of $G_0$ containing $C_1\cup C_2$, by \cite[(3.3)]{BD23} and \cite[Cor.~9.4]{S02} we deduce a surjection of Fr\'echet spaces $D(G_0)_{C_1}\oplus D(G_0)_{C_2}\twoheadrightarrow D(G_0)_{C_1\cup C_2}$. In particular we have $D(G_0)_{C_1}+D(G_0)_{C_2}= D(G_0)_{C_1\cap C_2}$ in $D(G_0)$.
\end{proof}

\subsection{Ext groups of Orlik-Strauch representations}\label{subsec: loc an spectral seq}

We prove several results on the $\mathrm{Ext}$ groups of Orlik-Strauch representations, in particular that they are finite dimensional when their smooth entries are finite length representations (Theorem \ref{thm: finitedim}). The most important statements are Corollary \ref{cor: Ext P graded} (which follows from Theorem \ref{prop: p coh graded}) and Corollary \ref{prop: OS spectral seq}.\bigskip

We keep the notation of \S\ref{gradedpiecessubsection}, in particular we have $I_0,I_1\subseteq \Delta$, $M_0\in \cO_{\rm{alg}}^{\fp_{I_0}}$, $\pi_0^{\infty}$ a smooth strongly admissible representation of $L_{I_0}$ over $E$ and $V_0=\cF_{P_{I_0}}^G(M_0,\pi_0^{\infty})$. Note that we assume $\pi_0^{\infty}$ strongly admissible (\cite[\S 3]{ST02a}) instead of just admissible and we recall that if $\pi_0^{\infty}$ is of finite length then it is strongly admissible (\cite[Prop.~2.2]{ST01}). From now until Corollary~\ref{prop: OS spectral seq} (included), we let $M_1\defeq M^{I_1}(\mu)\in\cO_{\rm{alg}}^{\fp_{I_1}}$ for some $\mu\in \Lambda_{I_1}^{\dom}$ (a generalized Verma module, see (\ref{belongtoOb}) and the lines below it), $\pi_1^{\infty}$ a smooth strongly admissible representation of $L_{I_1}$ over $E$ and $V_1\defeq \cF_{P_{I_1}}^G(M_1,\pi_1^{\infty})$. Our main aim in this section is to study the $E$-vector spaces $\mathrm{Ext}_{D(G)}^k(V_1^\vee,V_0^\vee)$ for $k\geq 0$.\bigskip

Recall from Remark \ref{addendumstrongly} and (\ref{r2r2}) (noting that $R_1=0$ in (\ref{equ: OS def seq})) that we have a $D(G)$-equivariant isomorphism $V_1^\vee \cong D(G)\otimes_{D(P_{I_1})}(L^{I_1}(\mu)\otimes_E (\pi_1^{\infty})^\vee)$. By Lemma~\ref{lem: loc alg FIN} and the proof of \cite[Lemma 6.3(ii)]{ST05} (where we use \cite[Prop.~2.6]{Schm09} instead of \cite[Lemma~6.2]{ST05} as we are locally $K$-analytic) we deduce isomorphisms for $k\geq 0$
\begin{equation}\label{isov1v0}
\mathrm{Ext}_{D(G)}^k(V_1^\vee,V_0^\vee)\cong \mathrm{Ext}_{D(P_{I_1})}^k(L^{I_1}(\mu)\otimes_E (\pi_1^{\infty})^\vee, V_0^\vee).
\end{equation}
Thus our main aim is to compute $\mathrm{Ext}_{D(P_{I_1})}^k(L^{I_1}(\mu)\otimes_E (\pi_1^{\infty})^\vee, V_0^\vee)$ for $k\geq 0$.\bigskip

For $M$ in $\mathrm{Mod}_{U(\fp_{I_1})}$, we endow $\Hom_E(L^{I_1}(\mu),M)$ with a structure of $U(\fp_{I_1})$-module by
\[(u\cdot f)(x)= u\cdot f(x)-f(u\cdot x)\]
($u\in \fp_{I_1}$, $f\in \Hom_E(L^{I_1}(\mu),M)$, $x\in L^{I_1}(\mu)$). In particular
\begin{equation}\label{equ: p H0pI}
H^0(\fp_{I_1},\Hom_E(L^{I_1}(\mu),M)) \cong \Hom_{U(\fp_{I_1})}(L^{I_1}(\mu),M).
\end{equation}
It is easy to check for $M'$ in $\mathrm{Mod}_{U(\fp_{I_1})}$:
\begin{equation}\label{equ: p Hom adjunction}
\Hom_{U(\fp_{I_1})}(M', \Hom_E(L^{I_1}(\mu),M))\cong \Hom_{U(\fp_{I_1})}(M'\otimes_E L^{I_1}(\mu), M)
\end{equation}
with $U(\fp_{I_1})$ acting diagonally on $M'\otimes_E L^{I_1}(\mu)$. It follows from (\ref{equ: p Hom adjunction}) and the exactness of $\Hom_E(L^{I_1}(\mu),-)$ that if $M^\bullet$ is an injective resolution of $M$ in $\mathrm{Mod}_{U(\fp_{I_1})}$, then $\Hom_E(L^{I_1}(\mu),M^\bullet)$ is an injective resolution of $\Hom_E(L^{I_1}(\mu),M)$ in $\mathrm{Mod}_{U(\fp_{I_1})}$. Using (\ref{equ: p H0pI}) this implies canonical isomorphisms for $\ell\geq 0$
\begin{equation}\label{equ: Ext g vers coh}
\mathrm{Ext}_{U(\fp_{I_1})}^{\ell}(L^{I_1}(\mu), M)\cong H^{\ell}(\fp_{I_1}, \Hom_E(L^{I_1}(\mu), M)).
\end{equation}
Recall that $H^{\ell}(\fp_{I_1}, \Hom_E(L^{I_1}(\mu), M))$ can also be computed by the Chevalley-Eilenberg complex $C^\bullet$ where for $\ell\geq 0$
\begin{equation}\label{equ: CE conj}
C^{\ell}\defeq\Hom_E(\wedge^{\ell} \fp_{I_1}, \Hom_E(L^{I_1}(\mu),M))\cong \Hom_E(\wedge^{\ell} \fp_{I_1}\otimes_E L^{I_1}(\mu),M)
\end{equation}
(see also (\ref{equ: CE complex})). The left (continuous) action of $P_I$ on $L^{I_1}(\mu)$ induces a left action on $\Hom_E(L^{I_1}(\mu),E)$ defined by $(g\cdot f)(x)=f(g^{-1}x)$ ($g\in P_{I_1}$, $f\in \Hom_E(L^{I_1}(\mu),E)$, $x\in L^{I_1}(\mu)$). As $\Hom_E(L^{I_1}(\mu),E)$ is finite dimensional, this action extends to a left $D(P_{I_1})$-action. Let $D$ be any $D(P_{I_1})$-module, we endow $\Hom_E(L^{I_1}(\mu),D)\cong \Hom_E(L^{I_1}(\mu),E)\otimes_E D$ with the diagonal (left) action of $D(P_{I_1})$ (via \cite[\S A]{ST05}). It is easy to check that for any $D(P_{I_1})$-module $D'$:
\begin{equation}\label{equ: P Hom adjunction}
\Hom_{D(P_{I_1})}(D', \Hom_E(L^{I_1}(\mu),D))\cong \Hom_{D(P_{I_1})}(D'\otimes_E L^{I_1}(\mu), D)
\end{equation}
with $D(P_{I_1})$ acting diagonally on $D'\otimes_E L^{I_1}(\mu)$. Hence, if $D^\bullet$ is an injective resolution of $D$ in $\mathrm{Mod}_{D(P_{I_1})}$, then $\Hom_E(L^{I_1}(\mu),D^\bullet)$ is an injective resolution of $\Hom_E(L^{I_1}(\mu),D)$ satisfying
\[\Hom_{D(P_{I_1})}(L^{I_1}(\mu)\otimes_E (\pi_1^{\infty})^\vee, D^k)\cong \Hom_{D(P_{I_1})}((\pi_1^{\infty})^\vee, \Hom_E(L^{I_1}(\mu),D^k))\]
for $k\geq 0$. We thus obtain canonical isomorphisms for $k\geq 0$
\begin{equation}\label{equ: P move alg}
\mathrm{Ext}_{D(P_{I_1})}^k(L^{I_1}(\mu)\otimes_E (\pi_1^{\infty})^\vee, D)\cong \mathrm{Ext}_{D(P_{I_1})}^k((\pi_1^{\infty})^\vee, \Hom_E(L^{I_1}(\mu),D)).
\end{equation}
It follows from \cite[\S 3]{ST05} (more precisely from \cite[page 307 line 5]{ST05} together with \cite[page 306 line -11]{ST05} applied with $X=X^\bullet=\Hom_E(L^{I_1}(\mu),D)$ and $Y^\bullet=(\pi_1^{\infty})^\vee$, both in degree $0$) that we have a spectral sequence
\[\mathrm{Ext}_{D^{\infty}(P_{I_1})}^k\big((\pi_1^{\infty})^\vee, H^{\ell}(\fp_{I_1}, \Hom_E(L^{I_1}(\mu),D))\big) \implies \mathrm{Ext}_{D(P_{I_1})}^{k+\ell}\big((\pi_1^{\infty})^\vee, \Hom_E(L^{I_1}(\mu),D)\big),\]
(in particular $H^{\ell}(\fp_{I_1}, \Hom_E(L^{I_1}(\mu),D))$ is naturally a $D^{\infty}(P_{I_1})$-module) which together with (\ref{equ: P move alg}) and (\ref{equ: Ext g vers coh}) (applied with $M=D$) gives a spectral sequence
\begin{equation}\label{equ: ST seq}
\mathrm{Ext}_{D^{\infty}(P_{I_1})}^k\big((\pi_1^{\infty})^\vee, \mathrm{Ext}_{U(\fp_{I_1})}^{\ell}(L^{I_1}(\mu),D)\big)\implies \mathrm{Ext}_{D(P_{I_1})}^{k+\ell}\big(L^{I_1}(\mu)\otimes_E (\pi_1^{\infty})^\vee, D\big).
\end{equation}
The spectral sequence (\ref{equ: ST seq}) applied to graded pieces of $V_0^\vee$ will be our primary means of accessing $\mathrm{Ext}_{D(P_{I_1})}^\bullet(L^{I_1}(\mu)\otimes_E (\pi_1^{\infty})^\vee, V_0^\vee)$ (and hence $\mathrm{Ext}_{D(G)}^k(V_1^\vee,V_0^\vee)$ by (\ref{isov1v0})).\bigskip

From now until Theorem~\ref{prop: p coh graded} (included) we fix $w\in W^{I_0,I_1}$. We write $C^\bullet$ (resp.~$\cC^\bullet$) for the Chevalley-Eilenberg complex (\ref{equ: CE conj}) with $M=M_0^w$ (resp.~$M=\cM_0^w$) where $M_0^w$ is defined above (\ref{equ: twist isom}) and $\cM_0^w$ as in Proposition \ref{lem: graded OS tensor}. Note that as $\wedge^{\ell} \fp_{I_1}\otimes_E L^{I_1}(\mu)$ is finite dimensional each $\cC^{\ell}$ is a Fr\'echet space. The algebra $U(\fp_{I_1})$ acts on both $\wedge^{\ell} \fp_{I_1}\otimes_E L^{I_1}(\mu)$ and $M_0^w$ (resp.~$\cM_0^w$), and thus acts diagonally on $C^{\ell}$ (resp.~$\cC^{\ell}$). We write $d^{\ell}: C^{\ell}\rightarrow C^{\ell+1}$ (resp.~${\delta}^{\ell}: \cC^{\ell}\rightarrow \cC^{\ell+1}$) for the differential maps of the complex $C^\bullet$ (resp.~$\cC^\bullet$), see \cite[p.~305]{ST05} for instance. The $U(\fg)$-equivariant embedding $1\otimes \Id:M_0\hookrightarrow \cM_0=D(G)_1\otimes_{U(\fg)}M_0$ induces an $U(\fg)$-equivariant embedding $M_0^w\hookrightarrow \cM_0^w$ which induces a map of complexes $C^\bullet \rightarrow \cC^\bullet$ which induces for $\ell\geq 0$
\begin{equation}\label{equ: p coh M}
\kappa^{\ell}: \mathrm{Ext}_{U(\fp_{I_1})}^{\ell}(L^{I_1}(\mu), M_0^w)\longrightarrow \mathrm{Ext}_{U(\fp_{I_1})}^{\ell}(L^{I_1}(\mu), \cM_0^w).
\end{equation}

\begin{lem}\label{lem: p coh M isom}
For $\ell\!\geq \!0$ $\kappa^{\ell}$ in (\ref{equ: p coh M}) is an isomorphism of finite dimensional $E$-vector~spaces.
\end{lem}
\begin{proof}
By d\'evissage using the short exact sequence $0\rightarrow M_0^w\rightarrow \cM_0^w\rightarrow (\cM_0/M_0)^w\rightarrow 0$, it suffices to show for $\ell\geq 0$:
\begin{equation}\label{equ: p coh M vanishing}
\mathrm{Ext}_{U(\fp_{I_1})}^{\ell}(L^{I_1}(\mu), (\cM_0/M_0)^w)=0.
\end{equation}
Let $e\in\ft$ such that $\alpha(e)\in\Z_{>0}$ for all $\alpha\in\Phi^+$. We divide the proof into three steps.\bigskip

\textbf{Step $1$}: We \ prove \ that \ $\mathrm{ad}(w^{-1})(e)-N$ \ acts \ bijectively \ on \ the \ $U(\ft)$-module \ $H^{\ell_1}(\fu, (\cM_0/M_0)^w)$ for $N\in\Z$ and $\ell_1\geq 0$.\\
By \ref{it: O completion 2} of Proposition \ref{prop: category O completion} and Lemma~\ref{lem: invertible action} $e-N$ for $N\in\Z$ acts bijectively on $\cM_0/M_0$ and thus $\mathrm{ad}(w^{-1})(e)-N$ acts bijectively on $(\cM_0/M_0)^w$. As $\wedge^{\ell_1}\fu$ is $\ft$-semi-simple we have $\wedge^{\ell_1}\fu=\bigoplus_{M\in\Z}\wedge^{\ell_1}\fu|_{\mathrm{ad}(w^{-1})(e)=M}$, and we see that $\mathrm{ad}(w^{-1})(e)-N$ acts invertibly on $\Hom_E(\wedge^{\ell_1}\fu, (\cM_0/M_0)^w)$, and thus on any of its $U(\ft)$-subquotient, in particular on its $U(\ft)$-subquotient $H^{\ell_1}(\fu, (\cM_0/M_0)^w)$.\bigskip

\textbf{Step $2$}: We prove $\mathrm{Ext}_{U(\ft)}^{\ell_2}(\mu_1, H^{\ell_1}(\fu, (\cM_0/M_0)^w))=0$ for $\ell_1,\ell_2\geq 0$ and $\mu_1\in\Lambda$.\\
Let $\ft'\defeq E (\mathrm{ad}(w^{-1})(e))\subseteq \ft$, for $\mu_1\in\Lambda$ and any $U(\ft)$-module we have a Hochschild-Serre spectral sequence
\[H^{\ell_2''}(\ft/\ft', H^{\ell_2'}(\ft', \Hom_E(\mu_1,D))) \implies H^{\ell_2'+\ell_2''}(\ft, \Hom_E(\mu_1, D))\cong \mathrm{Ext}_{U(\ft)}^{\ell_2'+\ell_2''}(\mu_1, D),\]
(where the last isomorphism is proved as (\ref{equ: Ext g vers coh})). Taking $D=H^{\ell_1}(\fu, (\cM_0/M_0)^w))$ it suffices to prove for $\ell_2' \geq 0$:
\begin{equation}\label{equ: p coh M t vanishing prime}
H^{\ell_2'}\big(\ft', \Hom_E(\mu_1, H^{\ell_1}(\fu, (\cM_0/M_0)^w))\big)=0.
\end{equation}
But since $\dim_E\ft' = 1$ the Chevalley-Eilenberg complex that computes (\ref{equ: p coh M t vanishing prime}) is just
\[\Hom_E(\mu_1, H^{\ell_1}(\fu, (\cM_0/M_0)^w)) \buildrel \mathrm{ad}(w^{-1})(e) \over \longrightarrow \Hom_E(\mu_1, H^{\ell_1}(\fu, (\cM_0/M_0)^w))\]
and by Step $1$ the unique differential map is an isomorphism, whence the result.\bigskip

\textbf{Step $3$}: We prove for $\ell\geq 0$ and $\mu_1\in\Lambda$:
\begin{equation}\label{equ: p coh M vanishing general}
\mathrm{Ext}_{U(\fp_{I_1})}^{\ell}(L^{I_1}(\mu_1), (\cM_0/M_0)^w)=0.
\end{equation}
Let $\mu_1\in\Lambda$, by (\ref{equ: g spectral seq}) applied with $I=\emptyset$ we have the spectral sequence
\[\mathrm{Ext}_{U(\ft)}^{\ell_2}(\mu_1, H^{\ell_1}(\fu, (\cM_0/M_0)^w)) \implies \mathrm{Ext}_{U(\fb)}^{\ell_1+\ell_2}(\mu_1, (\cM_0/M_0)^w),\]
which together with Step $2$ implies for $\ell\geq 0$
\begin{equation}\label{equ: p coh M vanishing Verma}
\mathrm{Ext}_{U(\fp_{I_1})}^{\ell}(U(\fp_{I_1})\otimes_{U(\fb)}\mu_1, (\cM_0/M_0)^w)\cong \mathrm{Ext}_{U(\fb)}^{\ell}(\mu_1, (\cM_0/M_0)^w)=0
\end{equation}
where the isomorphism in (\ref{equ: p coh M vanishing Verma}) is Shapiro's lemma. Recall $U(\fp_{I_1})\otimes_{U(\fb)}\mu_1\cong U(\fl_{I_1})\otimes_{U(\fb_{I_1})}\mu_1$. Hence if $U(\fl_{I_1})\otimes_{U(\fb_{I_1})}\mu_1\cong L^{I_1}(\mu_1)$, we are done. In general, we argue by induction using \cite[Thm.~5.1]{Hum08} applied to $\fl_{I_1}$ and a d\'evissage on the constituents of $U(\fl_{I_1})\otimes_{U(\fb_{I_1})}\mu_1$. Finally, applying (\ref{equ: p coh M vanishing general}) with $\mu_1=\mu$ gives (\ref{equ: p coh M vanishing}).
\end{proof}

For $\ell\geq 0$ the map ${\delta}^{\ell}: \cC^{\ell}\rightarrow \cC^{\ell+1}$ is a continuous map between Fr\'echet spaces, in particular $\mathrm{ker}({\delta}^{\ell})\subseteq \cC^{\ell}$ is a closed subspace. We endow $\mathrm{ker}({\delta}^{\ell})$ with the subspace topology of $\cC^{\ell}$ and $H^{\ell}(\cC^\bullet)=\mathrm{ker}({\delta}^{\ell})/\mathrm{im}({\delta}^{\ell-1})$ with the quotient topology (\cite[\S 5.B]{S02}).

\begin{lem}\label{lem: p coh M separated}
For $\ell\geq 0$ the differential map ${\delta}^{\ell}: \cC^{\ell}\rightarrow \cC^{\ell+1}$ has closed image and $H^{\ell}(\cC^\bullet)$ is a finite dimensional separated $E$-vector space (with its natural Banach topology).
\end{lem}
\begin{proof}
Let $\overline{\mathrm{im}({\delta}^{\ell-1})}$ be the closure of $\mathrm{im}({\delta}^{\ell-1})$ in $\cC^{\ell}$, which is still contained in $\mathrm{ker}({\delta}^{\ell})$, then $H^{\ell}(\cC^\bullet)$ is separated if and only if $\mathrm{im}({\delta}^{\ell-1})=\overline{\mathrm{im}({\delta}^{\ell-1})}$. The identification $M_0=\bigoplus_{\mu_1\in\Lambda}(\cM_0)_{\mu_1}$ from \ref{it: O completion 2} of Proposition~\ref{prop: category O completion} implies $M_0^w=\bigoplus_{\mu_1\in\Lambda}(\cM_0^w)_{\mu_1}$ and thus
$C^{\ell}=\bigoplus_{\mu_1\in\Lambda}(\cC^{\ell})_{\mu_1}$ for $\ell\geq 0$ with each $(\cC^{\ell})_{\mu_1}$ being finite dimensional (recall $\wedge^{\ell} \fp_{I_1}\otimes_E L^{I_1}(\mu)$ is finite dimensional and $U(\ft)$-semi-simple). In particular, together with the $U(\ft)$-equivariance of the differential maps $d^\bullet$ and ${\delta}^\bullet$ we deduce
\[\mathrm{im}(d^{\ell-1})=\bigoplus_{\mu_1\in\Lambda}\big(\mathrm{im}({\delta}^{\ell-1})\big)_{\mu_1}=\bigoplus_{\mu_1\in\Lambda}\big(\overline{\mathrm{im}({\delta}^{\ell-1})}\big)_{\mu_1}\ \ {\rm and}\ \ \mathrm{ker}(d^{\ell})=\bigoplus_{\mu_1\in\Lambda}\big(\mathrm{ker}({\delta}^{\ell})\big)_{\mu_1}.\]
and from (\ref{equ: quotient wt})
\[
\bigoplus_{\mu_1\in\Lambda}\big(\mathrm{ker}({\delta}^{\ell})/\overline{\mathrm{im}({\delta}^{\ell})}\big)_{\mu_1}\cong \Big(\bigoplus_{\mu_1\in\Lambda}\big(\mathrm{ker}({\delta}^{\ell})\big)_{\mu_1}\Big)/\Big(\bigoplus_{\mu_1\in\Lambda}\overline{\mathrm{im}\big({\delta}^{\ell})}\big)_{\mu_1}\Big)\\
\cong \mathrm{ker}(d^{\ell})/\mathrm{im}(d^{\ell}).
\]
This forces the composition $\mathrm{ker}(d^{\ell})/\mathrm{im}(d^{\ell})\xrightarrow{\kappa^{\ell}} \mathrm{ker}({\delta}^{\ell})/\mathrm{im}({\delta}^{\ell})\xrightarrow{\theta^{\ell}} \mathrm{ker}({\delta}^{\ell})/\overline{\mathrm{im}({\delta}^{\ell})}$ to be an injection. But $\kappa^{\ell}$ is a bijection by Lemma~\ref{lem: p coh M isom} and $\theta^{\ell}$ is a surjection by definition, so $\theta^{\ell}$ must also be an isomorphism, which means $\overline{\mathrm{im}({\delta}^{\ell})}=\mathrm{im}({\delta}^{\ell})$. The rest of the statement follows from Lemma \ref{lem: p coh M isom} (and (\ref{equ: Ext g vers coh})).
\end{proof}

We write $D^{\infty}\defeq (\pi_0^{\infty,w})^\vee$ and $D\defeq \cM_0^w\widehat{\otimes}_E D^{\infty}$ for short.
For $g\in G$ we write $\cM_0^{wg^{-1}}$, $D^{\infty,g^{-1}}$, $D^{g^{-1}}$ for the $D(G)_{gP_wg^{-1}}$-module with the same underlying space as $\cM_0^{w}$, $D^\infty$, $D$ (respectively) but with $\delta\in D(G)_{gP_wg^{-1}}$ acting by $\delta_g^{-1}\delta\delta_g$.
For $g\in G$ and $h\in gP_wg^{-1}$, the map $v\mapsto \delta_h v$ gives a $D(G)_{gP_wg^{-1}}=D(G)_{hgP_w(hg)^{-1}}$-equivariant topological isomorphism
\begin{equation}\label{equ: isom modulo coset}
\cM_0^{wg^{-1}}\buildrel\sim\over\longrightarrow \cM_0^{w(hg)^{-1}}
\end{equation}
and likewise $D(G)_{gP_wg^{-1}}$-equivariant isomorphisms $D^{\infty,g^{-1}}\buildrel\sim\over\rightarrow D^{\infty,(hg)^{-1}}$, $D^{g^{-1}}\buildrel\sim\over\rightarrow D^{(hg)^{-1}}$.\bigskip

For $g\in G$, we consider the Chevalley-Eilenberg complex attached to $\cM_0^{wg^{-1}}$:
\begin{equation*}
\cC_{g}^{\bullet}\defeq \Hom_E(\wedge^{\bullet} \fp_{I_1}\otimes_E L^{I_1}(\mu), \cM_0^{wg^{-1}})
\end{equation*}
and we denote by $\delta_{g}^{\ell}: \cC_{g}^{\ell}\rightarrow \cC_{g}^{\ell+1}$ the differential maps (not to be confused with the Dirac distribution $\delta_g$!).
For $\ell\geq 0$ we write $\cD_{g}^{\ell}\defeq \delta_{g}^{\ell-1}(\cC_{g}^{\ell-1})$ and
\begin{multline*}
\cH_{g}^{\ell}\defeq H^{\ell}(\cC_{g}^{\bullet})=\mathrm{ker}(\delta_{g}^\ell)/\cD_{g}^{\ell}= H^{\ell}(\fp_{I_1}, \Hom_E(L^{I_1}(\mu), \cM_0^{wg^{-1}}))\\
\buildrel (\ref{equ: Ext g vers coh}) \over \cong \mathrm{Ext}_{U(\fp_{I_1})}^{\ell}(L^{I_1}(\mu), \cM_0^{wg^{-1}}).
\end{multline*}
For $\ell\geq 0$, $g\in G$ and $h\in gP_wg^{-1}$, the $D(G)_{gP_wg^{-1}}$-equivariant topological isomorphism (\ref{equ: isom modulo coset}) induces a $U(\fp_{I_1})$-equivariant topological isomorphism
\begin{equation}\label{equ: twist CE complex modulo coset}
\cC_{g}^\ell\buildrel\sim\over\longrightarrow \cC_{hg}^\ell,
\end{equation}
which further induces topological isomorphisms $\cD_{g}^\ell\buildrel\sim\over\rightarrow \cD_{hg}^\ell$, $\mathrm{ker}(\delta_{g}^\ell)\buildrel\sim\over\rightarrow \mathrm{ker}(\delta_{hg}^\ell)$ and
\begin{equation}\label{equ: twist CE coh modulo coset}
\cH_{g}^\ell\buildrel\sim\over\longrightarrow \cH_{hg}^\ell.
\end{equation}
It follows from (\ref{equ: twist CE complex modulo coset}) and (\ref{equ: twist CE coh modulo coset}) that the complex $\cC_{g}^\bullet$ and the cohomology space $\cH_{g}^\ell$ only depend on the coset $(gP_wg^{-1})g=gP_w$ up to natural $U(\fp_{I_1})$-equivariant topological isomorphisms. We tacitly use this in the sequel.\bigskip

Let $\ell\geq 0$, $g\in G$, $h \in P_{I_1}$.
As $\wedge^{\ell} \fp_{I_1}\otimes_E L^{I_1}(\mu)$ is a finite dimensional $D(P_{I_1})$-module, the map $v\mapsto \delta_{h} v$ induces an isomorphism of $U(\fp_{I_1})$-modules $\wedge^{\bullet} \fp_{I_1}\otimes_E L^{I_1}(\mu)\buildrel\sim\over\longrightarrow (\wedge^{\bullet} \fp_{I_1}\otimes_E L^{I_1}(\mu))^{h^{-1}}$, which induces a topological isomorphism
\begin{multline}\label{equ: twist CE isom}
\theta_{g,h}^\ell: \cC_{g}^{\ell} = \Hom_E(\wedge^{\ell} \fp_{I_1}\otimes_E L^{I_1}(\mu), \cM_0^{wg^{-1}})\\
\cong \Hom_E((\wedge^{\ell} \fp_{I_1}\otimes_E L^{I_1}(\mu))^{h^{-1}}, \cM_0^{wg^{-1}h^{-1}})\\
\buildrel\sim\over\longrightarrow \Hom_E(\wedge^{\ell} \fp_{I_1}\otimes_E L^{I_1}(\mu), \cM_0^{w(hg)^{-1}}) = \cC_{hg}^{\ell}.
\end{multline}
Moreover under (\ref{equ: twist CE isom}) the differential map $\delta_{g}^{\ell}$ corresponds to $\delta_{hg}^{\ell}$, hence we deduce a topological isomorphism for $\ell\geq 0$, $g\in G$ and $h \in P_{I_1}$
\begin{equation}\label{equ: twist CE coh isom}
\omega_{g,h}^\ell: \cH_{g}^{\ell} = \mathrm{Ext}_{U(\fp_{I_1})}^{\ell}(L^{I_1}(\mu), \cM_0^{wg^{-1}})
\buildrel\sim\over\longrightarrow \mathrm{Ext}_{U(\fp_{I_1})}^{\ell}(L^{I_1}(\mu), \cM_0^{w(hg)^{-1}}) = \cH_{hg}^{\ell}.
\end{equation}
By a direct check for $g\in G$ and $h, h'\in P_{I_1}$ we have $\theta_{hg,h'}^\ell\circ\theta_{g,h}^\ell=\theta_{g,h'h}^\ell$ and therefore
\begin{equation}\label{equ: CE coh isom relation}
\omega_{hg,h'}^\ell\circ\omega_{g,h}^\ell=\omega_{g,h'h}^\ell.
\end{equation}
Note that $\cC_{g}^\ell$ (for $g\in G$) contains the ($U(\fp_{I_1})$-equivariant) subcomplex
\begin{equation*}
C_g^\bullet\defeq \Hom_E(\wedge^{\bullet} \fp_{I_1}\otimes_E L^{I_1}(\mu), M_0^{wg^{-1}})
\end{equation*}
with $C_g^\ell$ being a dense subspace of $\cC_{g}^\ell$ for each $\ell\geq 0$. The map between complexes $C_g^\bullet\rightarrow \cC_{g}^\bullet$ induces a natural map for $\ell\geq 0$
\begin{equation}\label{equ: p coh M isom g}
H^\ell(C_g^\bullet)\longrightarrow H^\ell(\cC_g^\bullet)=\cH_{g}^\ell.
\end{equation}

\begin{lem}\label{lem: p coh M separated g}
Assume that $g\in X=P_{I_1}P_w$. Then for $\ell\geq 0$, the morphism (\ref{equ: p coh M isom g}) is an isomorphism, the image $\cD_{g}^\ell$ of the differential map $\delta_{g}^{\ell-1}: \cC_{g}^{\ell-1}\rightarrow \cC_{g}^\ell$ is closed and $\cH_{g}^\ell=H^{\ell}(\cC_{g}^\bullet)=\mathrm{ker}(\delta_{g}^\ell)/\cD_{g}^\ell$ is a finite dimensional separated $E$-vector space (with its natural Banach topology).
\end{lem}
\begin{proof}
As $g\in X=P_{I_1}P_w$, there exists $h\in gP_wg^{-1}$ such that $hg\in P_{I_1}$.
Using the $U(\fp_{I_1})$-equivariant topological isomorphism (\ref{equ: twist CE complex modulo coset}) (which restricts to an isomorphism $C_{g}^\ell\buildrel\sim\over\longrightarrow C_{hg}^\ell$) and upon replacing $g$ with $hg$, we can assume $g\in P_{I_1}$. As $g\in P_{I_1}$, we know that $\delta_{g}^{\ell}$ corresponds to $\delta_{1}^\ell=\delta^{\ell}$ under (\ref{equ: twist CE isom}) for $h=g^{-1}$. Consequently, the desired results for $\delta_{g}^{\ell-1}$, $\cD_{g}^\ell$ and $\cH_{g}^\ell$ follow from those for $\delta_{1}^{\ell-1}=\delta^{\ell-1}$, $\cD_{1}^\ell$ and $\cH_{1}^\ell$, which are proven in Lemma~\ref{lem: p coh M isom} and Lemma~\ref{lem: p coh M separated}. In particular, it follows from (\ref{equ: twist CE isom}) (for $h=g^{-1}\in P_{I_1}$) and (\ref{equ: CE coh isom relation}) that $\omega_{g,g^{-1}}^\ell=(\omega_{1,g}^\ell)^{-1}: \cH_{g}^\ell\buildrel\sim\over\rightarrow \cH_{1}^\ell$ is a topological isomorphism of finite dimensional $E$-vector spaces (with their natural Banach topology).
\end{proof}

\begin{lem}\label{lem: coset dependence}
For $g,h\in P_{I_1}$ such that $h\in g(P_{I_1}\cap P_w)g^{-1}$ we have
\begin{equation}\label{equ: coset dependence}
\omega_{g,h}^\ell=\mathrm{Id}_{\cH_{g}^\ell}.
\end{equation}
In particular, the map $\omega_{g,h}^\ell$ only depends on the cosets $gP_w$ and $hgP_w$.
\end{lem}
\begin{proof}
For $\ell\geq 0$ and $g\in P_{I_1}$ by Lemma~\ref{lem: p coh M separated g} the embedding $M_0^{wg^{-1}}\hookrightarrow \cM_0^{wg^{-1}}$ induces an isomorphism $\mathrm{Ext}_{U(\fp_{I_1})}^{\ell}(L^{I_1}(\mu), M_0^{wg^{-1}})\buildrel\sim\over\rightarrow \mathrm{Ext}_{U(\fp_{I_1})}^{\ell}(L^{I_1}(\mu), \cM_0^{wg^{-1}})\cong \cH_{g}^{\ell}$. If $g^{-1}hg\in P_{I_1}\cap P_w$, the $D(G)_{gP_wg^{-1}}=D(G)_{hgP_w(hg)^{-1}}$-equivariant isomorphism (\ref{equ: isom modulo coset}) of Fr\'echet spaces restricts to an isomorphism of $U(\fg)$-modules $M_0^{wg^{-1}}\buildrel\sim\over\rightarrow M_0^{w(hg)^{-1}}$.
In the following we identify $\cM_0^{w(hg)^{-1}}$ with $\cM_0^{wg^{-1}}$ and $M_0^{w(hg)^{-1}}$ with $M_0^{wg^{-1}}$ via (\ref{equ: isom modulo coset}). Hence (under (\ref{equ: isom modulo coset})) $h\mapsto \theta_{g,h}^\ell$ gives an action of $g(P_{I_1}\cap P_w)g^{-1}$ on $\cC_{g}^{\ell}$ which preserves the subspace
\[C_g^\ell=\Hom_E(\wedge^{\ell} \fp_{I_1}\otimes_E L^{I_1}(\mu), M_0^{wg^{-1}})\subseteq \cC_{g}^{\ell}.\]
By the definition of $\theta_{g,h}^\ell$, we see that the action of $g(P_{I_1}\cap P_w)g^{-1}$ on $C_g^\ell$ is algebraic, and its derivative at $1\in g(P_{I_1}\cap P_w)g^{-1}$ recovers the natural action of $g(\fp_{I_1}\cap \fp_w)g^{-1}$ on $C_g^\ell$.
Since $\cH_{g}^\ell\cong \mathrm{Ext}_{U(\fp_{I_1})}^{\ell}(L^{I_1}(\mu), M_0^{wg^{-1}})$ is a $g(P_{I_1}\cap P_w)g^{-1}$-subquotient of $C_g^\ell$ on which $g(\fp_{I_1}\cap \fp_w)g^{-1}$ (and even $g\fp_{I_1}g^{-1}=\fp_{I_1}$) acts trivially, we deduce that $g(P_{I_1}\cap P_w)g^{-1}$ also acts trivially on $\cH_{g}^\ell$, i.e.~$\omega_{g,h}^\ell=\mathrm{Id}_{\cH_{g}^\ell}$ for $h\in g(P_{I_1}\cap P_w)g^{-1}$. This proves (\ref{equ: coset dependence}).

Together with (\ref{equ: CE coh isom relation}) this implies that the map $\omega_{g,h}^\ell$ for $g,h\in P_{I_1}$ only depends on the cosets $g(P_{I_1}\cap P_w)$ and $hg(P_{I_1}\cap P_w)g^{-1}$. For the final statement, note that, as $g,h\in P_{I_1}$, the cosets $g(P_{I_1}\cap P_w)$ and $hg(P_{I_1}\cap P_w)g^{-1}$ determine uniquely the cosets $gP_w$ and $hgP_w$ and vice versa.
\end{proof}

Let $P_{w,0}\defeq P_w\cap G_0$ and $P_{I_1,0}\defeq P_{I_1}\cap G_0$. For $r\in\cI$ we let $(\cM_0^w)_r\defeq D(G_0)_{P_{w,0},r}\otimes_{D(G_0)_{P_{w,0}}}\cM_0^w$ and consider the Chevalley-Eilenberg complex
\begin{equation*}
\cC_{r}^{\ell}\defeq \Hom_E(\wedge^{\ell} \fp_{I_1}\otimes_E L^{I_1}(\mu), (\cM_0^w)_r)
\end{equation*}
where we denote by $\delta_{r}^{\ell}: \cC_{r}^{\ell}\rightarrow \cC_{r}^{\ell+1}$ the differential maps.
We let $\cD_{r}^{\ell}$ be the closure of $\cD^\ell\defeq \delta^{\ell-1}(\cC^{\ell-1})\subseteq \cC^\ell$ in $\cC_{r}^{\ell}$. Then $\cD_{r}^{\ell}$ is a closed subspace of $\mathrm{ker}(\delta_{r}^{\ell})\subseteq \cC_{r}^{\ell}$ and we define $\cH_{r}^{\ell}\defeq \mathrm{ker}(\delta_{r}^{\ell})/\cD_{r}^{\ell}$ (with the quotient topology). The natural continuous map $\cM_0^w\rightarrow (\cM_0^w)_r$ induces continuous maps $\cC^\ell\rightarrow \cC_{r}^\ell$, $\cD^\ell\rightarrow \cD_{r}^\ell$, $\mathrm{ker}(\delta^\ell)\rightarrow \mathrm{ker}(\delta_{r}^\ell)$ and
\begin{equation}\label{equ: CE coh completion isom}
\cH^\ell\longrightarrow \cH_{r}^\ell.
\end{equation}
As $\cM_0^w\rightarrow (\cM_0^w)_r$ is a continous injection with dense image (see Lemma \ref{lem: family of standard semi norms}), so are $\cC^\ell\rightarrow \cC_{r}^\ell$ and $\cD^\ell\rightarrow \cD_{r}^\ell$.

\begin{lem}\label{lem: isom finite Banach}
For $r\in\cI$ the map $\mathrm{ker}(\delta^\ell)\rightarrow \mathrm{ker}(\delta_{r}^\ell)$ has dense image and the map (\ref{equ: CE coh completion isom}) is a topological isomorphism of finite dimensional $E$-Banach spaces.
\end{lem}
\begin{proof}
It follows from Lemma~\ref{lem: family of standard semi norms} (after conjugation by $w$) that there exists a family of standard semi-norms $\{|\cdot|_r\}_{r\in\cI}$ (see Definition~\ref{def: standard semi norm}) on $M_0^w$ such that $(\cM_0^w)_r$ is the completion of $M_0^w$ under $|\cdot|_r$. By \ref{it: tensor with fin dim 2} of Lemma~\ref{lem: tensor with fin dim}, there exists a norm $|\cdot|$ on the finite dimensional Fr\'echet $U(\ft)$-module $\Hom_E(\wedge^\ell\fp_{I_1}\otimes_E L^{I_1}(\mu),E)$ such that the family of semi-norms $|\cdot|\otimes_E |\cdot|_r$ on $C^\ell=\Hom_E(\wedge^\ell\fp_{I_1}\otimes_E L^{I_1}(\mu),M_0^w)$ is standard with $\cC_{r}^\ell$ being the completion of $C^\ell$ under $|\cdot|\otimes_E |\cdot|_r$. By Remark~\ref{rem: dense in Banach} we know that $\mathrm{ker}(\delta_{r}^\ell)\cap C^\ell$ is dense in $\mathrm{ker}(\delta_{r}^\ell)$. Since $\mathrm{ker}(\delta_{r}^\ell)\cap C^\ell=\mathrm{ker}(\delta_{r}^\ell|_{C^\ell})=\mathrm{ker}(\delta^\ell|_{C^\ell})\subseteq \mathrm{ker}(\delta^\ell)$, we deduce that $\mathrm{ker}(\delta^\ell)$ is dense in $\mathrm{ker}(\delta_{r}^\ell)$ (here we use $\cC^\ell\hookrightarrow \cC_{r}^\ell$), and thus (\ref{equ: CE coh completion isom}) has dense image. By Lemma~\ref{lem: p coh M isom} the injection $M_0^w\hookrightarrow \cM_0^w$ induces an isomorphism $H^\ell(C^\bullet)\buildrel\sim\over\rightarrow H^\ell(\cC^\bullet)=\cH^\ell$ and by Lemma~\ref{lem: p coh M separated} $\cH^\ell$ is a finite dimensional $E$-Banach space. As $\delta^{\ell-1}$ is continuous, $\delta^{\ell-1}(C^{\ell-1})$ is dense in $\delta^\ell(\cC^{\ell-1})=\cD^\ell$ and hence in $\cD_{r}^\ell$, which forces $\delta^{\ell-1}(C^{\ell-1})=\cD_{r}^\ell\cap C^\ell$ by the bijection statement in Remark~\ref{rem: dense in Banach}. It follows that the map
$H^\ell(C^\bullet) \rightarrow \cH^\ell_r$ is an injection, and hence (with the isomorphism $H^\ell(C^\bullet)\cong \cH^\ell$) that (\ref{equ: CE coh completion isom}) is a (continuous) injection. Hence (\ref{equ: CE coh completion isom}) is a continuous injection with dense image from a finite dimensional $E$-Banach space to an $E$-Banach space, it must therefore be a topological isomorphism of finite dimensional $E$-Banach spaces.
\end{proof}

Let us fix $r\in\cI$. For $g\in G_0$ we define $(\cM_0^w)_r^{g^{-1}}$ as $D(G_0)_{gP_{w,0}g^{-1}}$-modules. Using $D(G_0)_{P_{w,0},r}=D(G_0)_{P_{w,0}G_1^m,r}$ (see \ref{it: support basis 1} of Lemma \ref{lem: support basis} with $m$ defined from $r$ as in \emph{loc.~cit.}), for $g\in G_0$ and $h\in gP_{w,0}G_1^mg^{-1}$ the map $v\mapsto \delta_hv$ gives a $D(G_0)_{gP_{w,0}g^{-1},r}=D(G_0)_{hgP_{w,0}(hg)^{-1},r}$-equivariant topological isomorphism
\begin{equation}\label{equ: isom modulo coset r}
(\cM_0^w)_r^{g^{-1}}\buildrel\sim\over\longrightarrow (\cM_0^w)_r^{(hg)^{-1}}.
\end{equation}
We consider the following Chevalley-Eilenberg complex
\begin{equation*}
\cC_{g,r}^{\ell}\defeq \Hom_E(\wedge^{\ell} \fp_{I_1}\otimes_E L^{I_1}(\mu), (\cM_0^w)_r^{g^{-1}})
\end{equation*}
and denote by $\delta_{g,r}^{\ell}: \cC_{g,r}^{\ell}\rightarrow \cC_{g,r}^{\ell+1}$ the differential maps.
We write $\cD_{g,r}^{\ell}$ for the closure of the image of $\cD_{g}^{\ell}=\delta_{g}^{\ell-1}(\cC_{g}^{\ell-1})$ in $\cC_{g,r}^{\ell}$ and define $\cH_{g,r}^{\ell}\defeq \mathrm{ker}(\delta_{g,r}^{\ell})/\cD_{g,r}^{\ell}$.
For $g\in G_0$ and $h\in gP_{w,0}G_1^mg^{-1}$, (\ref{equ: isom modulo coset r}) induces a $U(\fp_{I_1})$-equivariant topological isomorphism $\cC_{g,r}^\ell\buildrel\sim\over\rightarrow \cC_{hg,r}^\ell$, which further induces topological isomorphisms $\cD_{g,r}^\ell\buildrel\sim\over\rightarrow \cD_{hg,r}^\ell$, $\mathrm{ker}(\delta_{g,r}^{\ell})\buildrel\sim\over\rightarrow \mathrm{ker}(\delta_{hg,r}^{\ell})$ and $\cH_{g,r}^\ell\buildrel\sim\over\rightarrow \cH_{hg,r}^\ell$. Hence $\cC_{g,r}^\ell$, $\cD_{g,r}^\ell$, $\mathrm{ker}(\delta_{g,r}^{\ell})$ and $\cH_{g,r}^\ell$ only depend on the coset $gP_{w,0}G_1^m$, or equivalently on the group $gP_{w,0}G_1^mg^{-1}$ (writing $gh=ghg^{-1}g$ for $h\in P_{w,0}G_1^m$), up to canonical $U(\fp_{I_1})$-equivariant topological isomorphisms.\bigskip

For $g\in P_{I_1,0}$, using Lemma~\ref{lem: top ML} and that $\cD_{g}^\ell$ is closed in $\cC_{g}^\ell$ (Lemma~\ref{lem: p coh M separated g}), the projective limit over $r\in \cI$ of the strict exact sequence of $E$-Banach spaces $0\rightarrow \cD_{g,r}^\ell \rightarrow \mathrm{ker}(\delta_{g,r}^{\ell}) \rightarrow \cH_{g,r}^\ell \rightarrow 0$ gives back the strict exact sequence of $E$-Fr\'echet spaces $0\rightarrow \cD_{g}^\ell \rightarrow \mathrm{ker}(\delta_{g}^\ell) \rightarrow \cH_{g}^\ell\rightarrow 0$.\bigskip

Similar to (\ref{equ: twist CE isom}), for $r\in \cI$, $g\in G_0$, $h\in P_{I_1,0}$ and $\ell\geq 0$, we have topological isomorphisms of $E$-Banach spaces
\begin{multline}\label{equ: twist CE isom r}
\theta_{g,h,r}^\ell: \cC_{g,r}^{\ell} = \Hom_E(\wedge^{\ell} \fp_{I_1}\otimes_E L^{I_1}(\mu), (\cM_0^w)_r^{g^{-1}})\\
\cong \Hom_E((\wedge^{\ell} \fp_{I_1}\otimes_E L^{I_1}(\mu))^{h^{-1}}, (\cM_0^w)_r^{g^{-1}h^{-1}})\\
\buildrel\sim\over\longrightarrow \Hom_E(\wedge^{\ell} \fp_{I_1}\otimes_E L^{I_1}(\mu), (\cM_0^w)_r^{(hg)^{-1}}) = \cC_{hg,r}^{\ell}
\end{multline}
under which the differential map $\delta_{g,r}^{\ell}$ corresponds to $\delta_{hg,r}^{\ell}$. Hence we deduce a topological isomorphism of $E$-Banach spaces ($r\in \cI$, $g\in G_0$, $h\in P_{I_1,0}$, $\ell\geq 0$)
\begin{equation}\label{equ: twist CE coh isom r}
\omega_{g,h,r}^\ell: \cH_{g,r}^{\ell} \buildrel\sim\over\longrightarrow \cH_{hg,r}^{\ell}.
\end{equation}
By Lemma~\ref{lem: isom finite Banach} $\mathrm{ker}(\delta_{1}^\ell)$ is dense in $\mathrm{ker}(\delta_{1,r}^\ell)$ and (\ref{equ: CE coh completion isom}) is a topological isomorphism of finite dimensional $E$-Banach spaces. For $g\in P_{I_1,0}$, it then follows from (\ref{equ: twist CE isom r}) (with $1,g$ instead of $g,h$) that $\mathrm{ker}(\delta_{g}^\ell)$ is dense in $\mathrm{ker}(\delta_{g,r}^\ell)$, and from (\ref{equ: twist CE coh isom r}) and (\ref{equ: twist CE coh isom}) (also applied with $1,g$) that the natural map
\begin{equation}\label{equ: CE coh completion isom g}
\cH_{g}^\ell\longrightarrow \cH_{g,r}^\ell
\end{equation}
is a topological isomorphism of finite dimensional $E$-Banach space for $r\in \cI$.\bigskip

We write for $r\in \cI$
\begin{multline}\label{varprojlimr sm}
D^{\infty}_r\defeq D(G_0)_{P_{w,0},r}\otimes_{D(G_0)_{P_{w,0}}}D^{\infty}\cong \big(D(G_0)_{P_{w,0},r}\otimes_{D(G_0)_{P_{w,0}}}D^\infty(G_0)_{P_{w,0}}\big)\otimes_{D^\infty(G_0)_{P_{w,0}}}D^{\infty}\\
\cong D^\infty(G_0)_{P_{w,0},r}\otimes_{D^\infty(P_{w,0})}D^{\infty}\cong D^{\infty}(P_{w,0})_r\otimes_{D^{\infty}(P_{w,0})}D^{\infty}
\end{multline}
where we recall that the second isomorphism follows from $D(G_0)_{P_{w,0},r}\otimes_{U(\fg)}E\cong D^\infty(G_0)_{P_{w,0},r}$ (see Step 2 in the proof of Lemma \ref{lem: support basis}) and $D^\infty(G_0)_{P_{w,0}}\cong D^\infty(P_{w,0})$ (see Lemma \ref{lem: smooth restriction}) and that the last follows using (\ref{DPR}). Recall also that \begin{equation}\label{varprojlimr}
D=\cM_0^w\widehat{\otimes}_E D^{\infty} \cong \varprojlim_{r\in\cI}D_r
\end{equation}
where $D_r\defeq (\cM_0^w)_r\otimes_E D^{\infty}_r$ with $D^{\infty}_r$ finite dimensional (use the coadmissibility of $\cM_0$ (Lemma \ref{lem: support variant}) together with Lemma \ref{limprojtensor}). For $g\in G_0$ and $r\in\cI$ we define $D^{\infty,g^{-1}}_r$ and $D_r^{g^{-1}}$ as $D(G_0)_{gP_{w,0}g^{-1}}$-modules.\bigskip

Recall that $X=P_{I_1}P_w$ and $X_0=X\cap G_0$. We write $\cX_0$ for the set of compact open subsets of $X_0$ stable under right multiplication by $P_{w,0}$, and $\cX$ for the set of open subsets of $X$ which are stable under right multiplication by $P_w$ and which have compact image in $X/P_w$. We denote by $U_0$ and $U$ a general element of respectively $\cX_0$ and $\cX$ (do not confuse $U$ here with the unipotent radical of the Borel subgroup $B\subseteq G$!). Since $X_0/P_{w,0}=X/P_w$, the map $\cX_0\rightarrow \cX: U_0\mapsto U_0P_w$ is a bijection with inverse given by $U\mapsto U\cap G_0$.\bigskip

Recall from (\ref{projc0}) that we have topological isomorphisms
\begin{equation*}
\widehat{D}(G_0)_{X_0}\cong \varprojlim_{U_0\in \cX_0}D(G_0)_{U_0} \cong \varprojlim_{U_0\in \cX_0,r\in\cI}D(G_0)_{U_0,r},
\end{equation*}
and from \ref{it: support basis 2} of Lemma~\ref{lem: support basis} that $D(G_0)_{U_0,r}$ is a finite free $D(G_0)_{P_{w,0},r}$-module. Applying Lemma~\ref{limprojtensor} to $V=\widehat{D}(G_0)_{X_0}$, $W=D$, $A=D(G_0)_{P_{w,0}}$ (with the index $r$ replaced by $U_0,r$) and again to $V=D(G_0)_{U_0}$, $W=D$, $A=D(G_0)_{P_{w,0}}$, we obtain topological isomorphisms (where $\widehat{\otimes}_{D(G_0)_{P_{w,0}}}$ is defined before Lemma \ref{limprojtensor})
\begin{multline}\label{equ: graded short term}
\mathrm{gr}_w(V_0^\vee)\cong \widehat{D}(G_0)_{X_0}\widehat{\otimes}_{D(G_0)_{P_{w,0}}}D
\cong \varprojlim_{U_0\in\cX_0,r\in\cI}\big(D(G_0)_{U_0,r}\otimes_{D(G_0)_{P_{w,0},r}}D_r\big)\\
\cong \varprojlim_{U_0\in\cX_0}\big(D(G_0)_{U_0}\widehat\otimes_{D(G_0)_{P_{w,0}}}D\big)\cong \varprojlim_{U\in\cX}\big(D(G)_{U}\widehat\otimes_{D(G)_{P_w}}D\big).
\end{multline}
where the first isomorphism follows from Proposition \ref{lem: graded OS tensor} and its proof, and the last isomorphism follows from
\begin{equation}\label{equ: G0 to G}
D(G_0)_{U_0}\widehat\otimes_{D(G_0)_{P_{w,0}}}D \cong D(G)_{U}\widehat\otimes_{D(G)_{P_w}}D
\end{equation}
which is analogous to (\ref{forlater}).\bigskip

By (\ref{equ: CE conj}) for $\ell\geq 0$ we have $\mathrm{Ext}_{U(\fp_{I_1})}^{\ell}(L^{I_1}(\mu), \mathrm{gr}_w(V_0^\vee))\cong H^{\ell}(\cC_w^\bullet)$ where
\begin{equation*}
\cC_w^{\ell}\defeq \Hom_E(\wedge^{\ell} \fp_{I_1}, \Hom_E(L^{I_1}(\mu), \mathrm{gr}_w(V_0^\vee)))\cong \Hom_E(\wedge^{\ell} \fp_{I_1}\otimes_E L^{I_1}(\mu), \mathrm{gr}_w(V_0^\vee))
\end{equation*}
and where $\cC_w^{\ell}$ is a left $D(P_{I_1})$-module by the discussion before (\ref{equ: P Hom adjunction}) replacing $L^{I_1}(\mu)$ there by $\wedge^{\ell} \fp_{I_1}\otimes_E L^{I_1}(\mu)$ with the diagonal $D(P_{I_1})$-action (via \cite[\S A]{ST05} and with the adjoint action of $P_{I_1}$ on $\wedge^{\ell} \fp_{I_1}$).
We write $\delta_w^{\ell}: \cC_w^{\ell}\rightarrow \cC_w^{\ell+1}$ for the differential maps.\bigskip

The topological isomorphisms (\ref{equ: graded short term}) induce topological isomorphisms for $\ell\geq 0$
\begin{equation}\label{equ: CE graded limit}
\cC_w^{\ell}\cong \varprojlim_{U_0\in\cX_0,r\in\cI} \cC_{U_0,r}^{\ell} \cong \varprojlim_{U_0\in\cX_0} \cC_{U_0}^{\ell} \cong \varprojlim_{U\in\cX} \cC_{U}^{\ell}
\end{equation}
where
\[\cC_{U_0,r}^{\ell}\defeq \Hom_E\big(\wedge^{\ell} \fp_{I_1}\otimes_E L^{I_1}(\mu), D(G_0)_{U_0,r}\otimes_{D(G_0)_{P_{w,0},r}}D_r\big)\]
and similarly with $\cC_{U_0}^{\ell}$, $\cC_{U}^{\ell}$. Note that (\ref{equ: G0 to G}) gives a natural identification $\cC_{U_0}^\ell=\cC_{U}^\ell$. We write $\delta_{U_0,r}^\ell$, $\delta_{U_0}^\ell$ and $\delta_{U}^\ell$ for the corresponding differential maps.
Under (\ref{equ: CE graded limit}), we have $\delta_w^{\ell}=\varprojlim_{U_0,r}\delta_{U_0,r}^{\ell}=\varprojlim_{U_0}\delta_{U_0}^{\ell}=\varprojlim_{U}\delta_{U}^{\ell}$.
Recall from \ref{it: support basis 2} of Lemma~\ref{lem: support basis} that $D(G_0)_{U_0,r}$ is a finite free $D(G_0)_{P_{w,0},r}$-module with a basis given by $\{\delta_g\}_{g\in U_0G_1^m/P_{w,0}G_1^m}$, so that we have
\begin{multline}\label{equ: tensor as sum}
D(G_0)_{U_0,r}\otimes_{D(G_0)_{P_{w,0},r}}D_r \cong \bigoplus_{g\in U_0G_1^m/P_{w,0}G_1^m} D_r^{g^{-1}}\\
\cong \bigoplus_{g\in U_0G_1^m/P_{w,0}G_1^m} \big(((\cM_0^w)_r)^{g^{-1}}\otimes_E D^{\infty,g^{-1}}_r\big).
\end{multline}
This induces $U(\fp_{I_1})$-equivariant isomorphisms for $\ell\geq 0$
\begin{equation}\label{equ: differential decomposition}
\cC_{U_0,r}^{\ell}\cong \bigoplus_{g\in U_0G_1^m/P_{w,0}G_1^m} \Hom_E(\wedge^{\ell} \fp_{I_1}\otimes_E L^{I_1}(\mu), D_r^{g^{-1}})\cong \bigoplus_{g\in U_0G_1^m/P_{w,0}G_1^m} \big(\cC_{g,r}^{\ell}\otimes_E D^{\infty,g^{-1}}_r\big)
\end{equation}
and the differential maps satisfy
\begin{equation}\label{equ: differential component}
\delta_{U_0,r}^{\ell}=\bigoplus_{g\in U_0G_1^m/P_{w,0}G_1^m} (\delta_{g,r}^{\ell}\otimes_E \mathrm{Id}_{D^{\infty,g^{-1}}_r}).
\end{equation}

We write $\cD_w^{\ell}$ for the closure of $\delta_w^{\ell-1}(\cC_w^{\ell-1})$ in $\cC_w^{\ell}$ and define $\cH_w^{\ell}\defeq \mathrm{ker}(\delta_w^\ell)/\cD_w^\ell$, so that we have a strict short exact sequence
\begin{equation}\label{equ: p coh graded seq}
0\longrightarrow \cD_w^{\ell} \longrightarrow \mathrm{ker}(\delta_w^{\ell}) \longrightarrow \cH_w^{\ell} \longrightarrow 0.
\end{equation}
For $U_0\in \cX_0$ and $r\in\cI$, we similarly define $\cD_{\ast}^\ell$ and $\cH_{\ast}^\ell$ with $\ast$ being $U_0,r$ or $U_0$ or $U=U_0P_w$. In particular, we have a short exact sequence of $E$-Banach spaces
\begin{equation}\label{equ: p coh graded seq r}
0\longrightarrow \cD_{U_0,r}^{\ell} \longrightarrow \mathrm{ker}(\delta_{U_0,r}^{\ell}) \longrightarrow \cH_{U_0,r}^{\ell} \longrightarrow 0.
\end{equation}
The projective limit over $r\in \cI$ of (\ref{equ: p coh graded seq r}) gives a strict exact sequence of $E$-Fr\'echet spaces
\begin{equation}\label{equ: p coh graded seq compact}
0\longrightarrow \cD_{U_0}^{\ell} \longrightarrow \mathrm{ker}(\delta_{U_0}^{\ell}) \longrightarrow \cH_{U_0}^{\ell} \longrightarrow 0,
\end{equation}
and the projective limit of (\ref{equ: p coh graded seq compact}) over $U_0\in \cX_0$ gives back (\ref{equ: p coh graded seq}) (using (\ref{equ: CE graded limit}) and Lemma~\ref{lem: top ML}).
Moreover we have topological isomorphisms by (\ref{equ: differential decomposition}) and (\ref{equ: differential component})
\[\mathrm{ker}(\delta_{U_0,r}^{\ell})\cong \!\!\!\!\bigoplus_{g\in U_0G_1^m/P_{w,0}G_1^m} \big(\mathrm{ker}(\delta_{g,r}^\ell) \otimes_E D^{\infty,g^{-1}}_r\big)\ \ \mathrm{and}\ \ \cD_{U_0,r}^{\ell}\cong \!\!\!\!\bigoplus_{g\in U_0G_1^m/P_{w,0}G_1^m} \big(\cD_{g,r}^\ell\otimes_E D^{\infty,g^{-1}}_r\big),\]
which together with (\ref{equ: p coh graded seq compact}) give a topological isomorphism
\begin{equation}\label{equ: p coh decomposition compact}
\cH_{U_0,r}^\ell \cong \bigoplus_{g\in U_0G_1^m/P_{w,0}G_1^m} \big(\cH_{g,r}^\ell\otimes_E D^{\infty,g^{-1}}_r\big).
\end{equation}

For $g\in X/P_w=X_0/P_{w,0}$ there is a natural continuous map $\cC_g^\ell\rightarrow \cC_w^\ell$ induced by the continuous map $D^{g^{-1}}=\delta_g D\rightarrow \mathrm{gr}_w(V_0^\vee)$ (see (\ref{equ: graded OS tensor})).
Taking the direct sum over $g\in X/P_w$, we obtain a map
\begin{equation}\label{equ: CE graded embedding}
\cC_{w,\flat}^\ell\defeq \bigoplus_{g\in X/P_w} \cC_g^\ell\widehat{\otimes}_E (D^{\infty})^{g^{-1}}\longrightarrow \cC_w^\ell.
\end{equation}
Similarly, for $U\in \cX$, we have a map
\begin{equation}\label{equ: CE graded embedding compact}
\cC_{U,\flat}^\ell\defeq \bigoplus_{g\in U/P_w} \cC_g^\ell\widehat{\otimes}_E (D^{\infty})^{g^{-1}}\longrightarrow \cC_{U}^\ell.
\end{equation}

\begin{lem}\label{lem: pointed embedding}
The maps (\ref{equ: CE graded embedding}) and (\ref{equ: CE graded embedding compact}) are injective with dense image.
\end{lem}
\begin{proof}
Let $U\in \cX$, $U_0\defeq U\cap G_0\in \cX_0$ and note that $U_0/P_{w,0}=U/P_w$. For $g,g'\in U_0$ and $r\in \cI$ one checks from (\ref{equ: differential decomposition}) that $\cC_{g,r}^\ell=\cC_{g',r}^\ell\subseteq \cC_{U_0,r}^\ell$ if and only if $g P_{w,0}G_1^m=g' P_{w,0}G_1^m$. Let $S\subseteq U_0$ be a finite subset such that the cosets $gP_{w,0}$ are distinct for $g\in S$. (Since $X_0/P_{w,0}=X/P_w$, this is equivalent to $gP_w$ being distinct for $g\in S$.) Then for $r\in \cI$ sufficiently close to $1$ the cosets $gP_{w,0}G_1^m$ are distinct for $g\in S$ and thus from (\ref{equ: differential decomposition}) $\cC_{U_0,r}^\ell$ (for such an $r$) contains as closed subspace the direct sum $\bigoplus_{g\in S}\cC_{g,r}^\ell\otimes_E (D^{\infty}_r)^{g^{-1}}$. Taking $\varprojlim_{r\in\cI}$ we obtain by Lemma \ref{limprojtensor} a closed embedding
\begin{equation}\label{equ: finite term embedding compact}
\bigoplus_{g\in S}\cC_{g}^\ell\widehat{\otimes}_E (D^{\infty})^{g^{-1}}\hookrightarrow \cC_{U_0}^\ell\cong \cC_{U}^\ell,
\end{equation}
and taking the projective limit over $U_0\in \cX_0$ a closed embedding (using (\ref{equ: CE graded limit}))
\begin{equation}\label{equ: finite term embedding}
\bigoplus_{g\in S}\cC_{g}^\ell\widehat{\otimes}_E (D^{\infty})^{g^{-1}}\hookrightarrow \cC_w^\ell.
\end{equation}
Taking the colimit over $S$ in (\ref{equ: finite term embedding compact}) and (\ref{equ: finite term embedding}) gives the maps (\ref{equ: CE graded embedding compact}) and (\ref{equ: CE graded embedding}) which are thus injective. As the map (\ref{equ: CE graded embedding compact}) induces a surjection $\cC_{U_0,\flat}^\ell\twoheadrightarrow \cC_{U_0,r}^\ell$ for each $r\in\cI$ (see (\ref{equ: differential decomposition})), it has dense image. The inclusion $U_0\subseteq X_0$ induces a projection $\cC_{w,\flat}^\ell\twoheadrightarrow \cC_{U_0,\flat}^\ell$ which is compatible with the natural projection $\cC_w^\ell\twoheadrightarrow \cC_{U_0}^\ell$. Since $\cC_w^\ell\cong \varprojlim_{U_0\in \cX_0}\cC_{U_0}^\ell$ and since (\ref{equ: CE graded embedding compact}) has dense image for each $U_0$, it follows that (\ref{equ: CE graded embedding}) also has dense image.
\end{proof}

In particular, the natural map $\cC_g^\ell\rightarrow \cC_w^\ell$ for $g\in X$ (resp.~$\cC_g^\ell\rightarrow \cC_{U}^\ell$ for $U\in\cX$ and $g\in U$) is injective.\bigskip

Note that the dense subspace $\cC_{w,\flat}^\ell\subseteq \cC_w^\ell$ is $P_{I_1}$-stable with $h\in P_{I_1}$ acting (on the left) via the collection of maps (for $g\in X/P_w$)
\[\theta_{g,h}^\ell\widehat{\otimes}_E \kappa_{g,h}: \cC_{g}^\ell\widehat{\otimes}_E (D^{\infty})^{g^{-1}} \longrightarrow \cC_{hg}^\ell\widehat{\otimes}_E (D^{\infty})^{(hg)^{-1}}\]
where
\begin{equation}\label{equ: CE conj smooth}
\kappa_{g,h}: (D^{\infty})^{g^{-1}}\longrightarrow (D^{\infty})^{(hg)^{-1}}
\end{equation}
sends $v$ to $\delta_h v$. As $\delta_{hg}^\ell=\delta_h\delta_g^\ell\delta_{h^{-1}}$ for $\ell\geq 0$, $g\in X/P_w$ and $h\in P_{I_1}$, the action of $P_{I_1}$ on $\cC_{w,\flat}^\ell$ stabilizes the two (closed) subspaces
\begin{equation*}
\bigoplus_{g\in X/P_w} \cD_{g}^{\ell}\widehat{\otimes}_E (D^{\infty})^{g^{-1}} \ \subseteq \ \mathrm{ker}(\delta_w^\ell)\cap \cC_{w,\flat}^\ell=\!\!\bigoplus_{g\in X/P_w} \mathrm{ker}(\delta_g^\ell)\widehat{\otimes}_E (D^{\infty})^{g^{-1}}
\end{equation*}
and thus induces a left action of $P_{I_1}$ on
\[\cH_{w,\flat}^{\ell}\defeq \bigoplus_{g\in X/P_w} \cH_{g}^\ell\otimes_E (D^{\infty})^{g^{-1}}\]
with $h\in P_{I_1}$ acting via the collection of maps (for $g\in X/P_w$)
\begin{equation}\label{equ: CE coh dense P action}
\omega_{g,h}^\ell\otimes_E \kappa_{g,h}: \cH_{g}^\ell\otimes_E (D^{\infty})^{g^{-1}} \longrightarrow \cH_{hg}^\ell\otimes_E (D^{\infty})^{(hg)^{-1}}.
\end{equation}
The maps (\ref{equ: CE graded embedding}) and (\ref{equ: CE graded embedding compact}) induce maps
\begin{equation*}
\cH_{w,\flat}^{\ell} \longrightarrow \cH_w^\ell
\end{equation*}
and for $U\in \cX$
\begin{equation}\label{equ: CE graded coh embedding compact}
\cH_{U,\flat}^{\ell}\defeq \bigoplus_{g\in U/P_w} \cH_{g}^\ell\otimes_E (D^{\infty})^{g^{-1}} \longrightarrow \cH_{U}^\ell.
\end{equation}

\begin{lem}\label{lem: CE coh total small}
For $U\in\cX$ such that $U\subseteq P_{I_1,0}P_w$ the map (\ref{equ: CE graded coh embedding compact}) is injective with dense image.
\end{lem}
\begin{proof}
Since $U\subseteq P_{I_1,0}P_w$ and thus $U_0\subseteq P_{I_1,0}P_{w,0}$, each (representative) $g$ in the decomposition (\ref{equ: p coh decomposition compact}) can be chosen in $P_{I_1,0}$, in which case we have the topological isomorphism (\ref{equ: CE coh completion isom g}). This forces the composition $\cH_{U_0,\flat}^{\ell}\rightarrow \cH_{U_0}^\ell\rightarrow \cH_{U_0,r}^\ell$ to be a surjection for each $r\in \cI$, and thus (\ref{equ: CE graded coh embedding compact}) has dense image. By a similar argument as in the proof of Lemma~\ref{lem: pointed embedding} using (\ref{equ: p coh decomposition compact}), (\ref{equ: CE coh completion isom g}), (\ref{equ: p coh graded seq r}) and (\ref{equ: p coh graded seq compact}) the map
\[\bigoplus_{g\in S} \cH_{g}^\ell\otimes_E (D^{\infty})^{g^{-1}} \longrightarrow \cH_{U}^\ell=\cH_{U_0}^\ell\]
is injective for each finite subset $S\subseteq U_0$ such that the cosets $gP_{w,0}$ for $g\in S$ are distinct. Taking the colimit over such $S$, it follows that (\ref{equ: CE graded coh embedding compact}) is injective.
\end{proof}

By Lemma~\ref{lem: coset dependence} applied with $1,g$ ($g\in P_{I_1}$) instead of $g,h$ the map $\omega^\ell_{1,g}$ in (\ref{equ: twist CE coh isom}) for $g\in P_{I_1}$ only depends on the coset $gP_w$. Consequently, we can define the map
\begin{equation}\label{equ: CE coh trivialization}
\zeta^\ell: \cH_{1}^\ell\otimes_E \Big(\bigoplus_{g\in X/P_w} (D^{\infty})^{g^{-1}}\Big) \longrightarrow \cH_{w,\flat}^{\ell}
\end{equation}
by sending $(x,y)\in \cH_{1}^\ell\otimes_E (D^{\infty})^{g^{-1}}$ to $(\omega^\ell_{1,g}(x),y)\in \cH_{g}^\ell\otimes_E (D^{\infty})^{g^{-1}}$ for $g\in P_{I_1}$. Note that since $\omega^\ell_{1,g}$ is bijective, the map $\zeta^\ell$ is also bijective. We defined above a left action of $P_{I_1}$ on $\cH_{w,\flat}^{\ell}$. We define one on the left hand side of (\ref{equ: CE coh trivialization}) by letting $P_{I_1}$ act trivially on $\cH_{1}^\ell$ and by (\ref{equ: CE conj smooth}) on $\bigoplus_{g\in X/P_w} (D^{\infty})^{g^{-1}}$.

\begin{lem}\label{lem: equivariance on dense subspace}
The map $\zeta^\ell$ is $P_{I_1}$-equivariant.
\end{lem}
\begin{proof}
We have for $h\in P_{I_1}$, $g\in X/P_w$ and $(x,y)\in \cH_{1}^\ell\otimes_E (D^{\infty})^{g^{-1}}$:
\[\zeta^\ell(\delta_h\cdot (x,y))=\zeta^\ell((x,\kappa_{g,h}(y)))=(\omega^\ell_{1,hg}(x),\kappa_{g,h}(y))=(\omega_{g,h}^\ell(\omega^\ell_{1,g}(x)),\kappa_{g,h}(y))=\delta_h\cdot\zeta^\ell(x,y)\]
where the third equality follows from (\ref{equ: CE coh isom relation}) and the last from (\ref{equ: CE coh dense P action}).
\end{proof}

For $U\in \cX$ the map $\zeta^\ell$ in (\ref{equ: CE coh trivialization}) restricts to a bijection
\begin{equation}\label{equ: CE coh trivialization compact}
\zeta_{U}^\ell: \cH_{1}^\ell\otimes_E \Big(\bigoplus_{g\in U/P_w} (D^{\infty})^{g^{-1}}\Big) \buildrel\sim\over\longrightarrow \cH_{U,\flat}^{\ell}.
\end{equation}
We recall that $\cH_{1}^\ell$ is finite dimensional by Lemma \ref{lem: isom finite Banach}.

\begin{lem}\label{lem: p coh graded small}
For $U\in\cX$ such that $U\subseteq P_{I_1,0}P_w$ the differential map $\delta_{U}^\ell$ has closed image and the bijection $\zeta_{U}^\ell$ uniquely extends into an isomorphism of $E$-Fr\'echet spaces
\begin{equation}\label{equ: p coh graded small}
\cH_{1}^\ell\otimes_E \big({D}^{\infty}(G)_{U}\widehat{\otimes}_{D^{\infty}(G)_{P_w}}D^{\infty}\big) \buildrel\sim\over\longrightarrow \cH_{U}^{\ell}.
\end{equation}
\end{lem}
\begin{proof}
Note that, although $U$ is not compact, $U/P_w=U_0/P_{w,0}$ is and by (\ref{equ: smooth hat}) we have ${D}^{\infty}(G)_{U}\cong \widehat{D}^{\infty}(G)_{U}$. Recall that given $U,U'\in \cX$ satisfying $U\cap U'=\emptyset$, we have $D(G)_{U\sqcup U'}\cong D(G)_U \oplus D(G)_{U'}$ and $D^\infty(G)_{U\sqcup U'}\cong D^\infty(G)_U \oplus D^\infty(G)_{U'}$, and thus $\cC_{U\sqcup U'}^\ell\cong \cC_{U}^\ell\oplus \cC_{U'}^\ell$, $\delta_{U\sqcup U'}^\ell=\delta_{U}^\ell\oplus \delta_{U'}^\ell$, $\cD_{U\sqcup U'}^\ell\cong \cD_{U}^\ell\oplus \cD_{U'}^\ell$ and $\cH_{U\sqcup U'}^\ell\cong \cH_{U}^\ell\oplus \cH_{U'}^\ell$. Moreover we have $P_{I_1,0}P_w\setminus U\in \cX$ for $U\in \cX$ such that $U\subseteq P_{I_1,0}P_w$ (writing $U=U_0P_w$ where $U_0=U\cap G_0\subseteq P_{I_1,0}P_{w,0}$, one checks that $P_{I_1,0}P_w\setminus U = (P_{I_1,0}P_{w,0}\setminus U_0)P_w$). It follows from all this that the statement of the lemma for $U=P_{I_1,0}P_w$ is the ``direct sum'' of the statement for $U$ and for $P_{I_1,0}P_w\setminus U$. Hence it is enough to prove (\ref{equ: p coh graded small}) for $U=P_{I_1,0}P_w$ or equivalently $U_0=P_{I_1,0}P_{w,0}$.\bigskip

Given a norm $|\cdot|$ on a $E$-vector space $V$ equipped with a left action by a group $H$, we say that $|\cdot|$ is $H$-invariant if $|h\cdot x|=|x|$ for each $h\in H$ and $x\in V$.
For $r\in\cI$, let $m,s$ be as before Lemma~\ref{lem: free mod coset basis}. By \cite[Prop.~5.6]{Schm08} the natural norm on the Banach algebra $D(G_1^m)_s$ is multiplicative, which together with \cite[Prop.~6.2]{Schm08} (cf.~Lemma~\ref{lem: free mod coset basis}) implies that the natural norm on $D(G_0)_{1,r}=D(G_0)_{G_1^m,r}$ is multiplicative and in particular $G_1^m$-invariant (for the natural left action of $G_1^m$ on $D(G_0)_{G_1^m,r}$). This together with \ref{it: support basis 1} of Lemma~\ref{lem: support basis} implies that there exists a $P_{w,0}G_1^m$-invariant norm on $D(G_0)_{P_{w,0},r}=D(G_0)_{P_{w,0}G_1^m,r}$ which defines its Banach topology. As $(\cM_0^w)_r$ is a finitely generated $D(G_0)_{P_{w,0},r}=D(G_0)_{P_{w,0}G_1^m,r}$-module, the $P_{w,0}G_1^m$-invariant norm on $D(G_0)_{P_{w,0}G_1^m,r}$ induces a $P_{w,0}G_1^m$-invariant norm on $(\cM_0^w)_r$ that defines its Banach topology (\cite[Prop.~2.1.i]{ST03}). Similarly, we can choose a $P_{w,0}G_1^m$-invariant norm $|\cdot|_r^{\infty}$ on $D^{\infty}_r$ (which defines its Banach topology). We write $|\cdot|_r$ (resp.~$|\cdot|_r^{\infty}$) for the induced semi-norm on $\cM_0^w$ (resp.~on $D^{\infty}$) and $|\cdot|_{D,r}\defeq |\cdot|_r\otimes_E |\cdot|_r^{\infty}$ for the induced semi-norm on $D$ (under which the completion of is $D_r$). (Note that this semi-norm $|\cdot|_r$ on $\cM_0^w$ might be different than the one in the proof of Lemma \ref{lem: isom finite Banach}.)\bigskip

For $r\in \cI$ we let $J_{r}\defeq P_{I_1,0}P_{w,0}G_1^m/P_{w,0}G_1^m$ (a finite set). We fix a choice of representatives $\tld{J}_{r}\subseteq P_{I_1,0}$ for $J_{r}$, and choose them in a compatible way so that we have a surjection $\tld{J}_{r'}\twoheadrightarrow \tld{J}_{r}$ for $r\leq r'$.\bigskip

Given $g\in P_{I_1,0}$, we write $|\cdot|_r^{g^{-1}}$ (resp.~$|\cdot|_r^{\infty,g^{-1}}$, $|\cdot|_{D^{g^{-1}},r}$) for the corresponding semi-norm on $\cM_0^{wg^{-1}}$ (resp.~$D^{\infty,g^{-1}}$, $D^{g^{-1}}$), and $(\cM_0^{wg^{-1}})_r$ (resp.~$(D^{\infty,g^{-1}})_r$, $(D^{g^{-1}})_r$) for the corresponding completion. We have obvious identifications $(\cM_0^{wg^{-1}})_r=(\cM_0^w)_r^{g^{-1}}$, $(D^{\infty,g^{-1}})_r=D^{\infty,g^{-1}}_r$ and $(D^{g^{-1}})_r=D_r^{g^{-1}}$. Since we have $|\delta_h\cdot x|_r=|x|_r$ for $h\in P_{w,0}G_1^m$ and $x\in (\cM_0^w)_r$, the norm $|\cdot|_r^{g^{-1}}$ on $(\cM_0^{wg^{-1}})_r$ only depends on the coset $gP_{w,0}G_1^m$. Similar facts hold for the norms $|\cdot|_r^{\infty,g^{-1}}$ and $|\cdot|_{D^{g^{-1}},r}$ on respectively $(D^{\infty,g^{-1}})_r$ and $(D^{g^{-1}})_r$.\bigskip

We fix a norm $|\cdot|_{\fp_{I_1}}$ on $\fp_{I_1}$. For $x\in \fp_{I_1}$ the map $g\mapsto |\mathrm{Ad}(g)(x)|_{\fp_{I_1}}$ is a locally constant function on $P_{I_1,0}$. Hence, replacing $|\cdot|_{\fp_{I_1}}$ by the norm $x\mapsto |\int_{P_{I_1,0}}\mathrm{Ad}(g)(x)|_{\fp_{I_1}}dg$ for some Haar measure $dg$ on $P_{I_1,0}$, we can assume $|\mathrm{Ad}(g)(x)|_{\fp_{I_1}}=|x|_{\fp_{I_1}}$ for $g\in P_{I_1,0}$ and $x\in \fp_{I_1}$. For $\ell\geq 0$, the norm $|\cdot|_{\fp_{I_1}}$ induces a norm on $\otimes^\ell_E\fp_{I_1}$ (the tensor product norm) and then a norm $|\cdot|_{\wedge^\ell \fp_{I_1}}$ on $\wedge^\ell \fp_{I_1}$ (the quotient norm) which satisfies ($g\in P_{I_1,0}$, $x\in \wedge^\ell \fp_{I_1}$)
\begin{equation}\label{equ: inv norm}
|\mathrm{Ad}(g)(x)|_{\wedge^\ell \fp_{I_1}}=|x|_{\wedge^\ell \fp_{I_1}}.
\end{equation}
For $g\in P_{I_1,0}$, the norm $|\cdot|_{\wedge^\ell \fp_{I_1}}$ on $\wedge^\ell \fp_{I_1}$ with the semi-norm $|\cdot|_r^{g^{-1}}$ on $\cM_0^{wg^{-1}}$ induce a semi-norm $|\cdot|_{\cC_g^{\ell},r}$ on $\cC_g^{\ell}$ such that the corresponding completion is $\cC_{g,r}^{\ell}$ (in particular $\{|\cdot|_{\cC_g^{\ell},r}\}_{r\in\cI}$ defines the Fr\'echet topology on $\cC_g^{\ell}$). It follows from (\ref{equ: inv norm}) and from (\ref{equ: twist CE isom}) that for $g,h\in P_{I_1,0}$, $r\in \cI$ and $x\in \cC_{g}^\ell$
\begin{equation}\label{equ: CE inv norm}
|\theta_{g,h}^\ell(x)|_{\cC_{hg}^\ell,r}=|x|_{\cC_{g}^\ell,r}.
\end{equation}
The (induced) norm $|\cdot|_{\cC_g^{\ell},r}$ on $\cC_{g,r}^{\ell}$ induces a norm on the closed subspaces $\cD_{g,r}^{\ell}$, $\mathrm{ker}(\delta_{g,r}^{\ell})$ and thus a norm $|\cdot|_{\cH_g^{\ell},r}$ on $\cH_{g,r}^{\ell}$, which satisfies for $g,h\in P_{I_1,0}$, $r\in \cI$ and $x\in \cH_{g}^\ell$ (using (\ref{equ: twist CE coh isom}))
\begin{equation}\label{equ: CE inv norm coh}
|\omega_{g,h}^\ell(x)|_{\cH_{hg}^\ell,r}=|x|_{\cH_{g}^\ell,r}.
\end{equation}
As the norm $|\cdot|_r^{g^{-1}}$ on $(\cM_0^{wg^{-1}})_r$ only depends on $gP_{w,0}G_1^m$, so does the norm $|\cdot|_{\cC_g^{\ell},r}$ on $\cC_{g,r}^{\ell}$ and the norm $|\cdot|_{\cH_g^{\ell},r}$ on $\cH_{g,r}^{\ell}$.\bigskip

Using (\ref{equ: differential decomposition}), we obtain a norm $|\cdot|_{\cC_{U_0}^\ell,r}$ on $\cC_{U_0,r}^\ell$ by taking the maximum of $|\cdot|_{\cC_g^{\ell},r}\otimes_E (|\cdot|_r^{\infty})^{g^{-1}}$ on the direct summands $\cC_{g,r}^{\ell}\otimes_E (D^{\infty}_r)^{g^{-1}}$ for $g\in \tld{J}_{r}$ (recall $U_0=P_{I_1,0}P_{w,0}$). As both $|\cdot|_{\cC_g^{\ell},r}$ and $(|\cdot|_r^{\infty})^{g^{-1}}$ only depend on $gP_{w,0}G_1^m$, so does $|\cdot|_{\cC_g^{\ell},r}\otimes_E (|\cdot|_r^{\infty})^{g^{-1}}$, and thus $|\cdot|_{\cC_{U_0}^\ell,r}$ does not depend on the choice of $\tld{J}_{r}$. Since $\cC_{U_0,r}^\ell$ is a Banach space with norm $|\cdot|_{\cC_{U_0}^\ell,r}$ and $\cC_{U_0}^\ell=\varprojlim_{r\in\cI}\cC_{U_0,r}^\ell$, the family of semi-norms $|\cdot|_{\cC_{U_0}^\ell,r}$ defines the Fr\'echet topology on $\cC_{U_0}^\ell$. Note that for $g\in P_{I_1,0}$ the restriction of $|\cdot|_{\cC_{U_0}^\ell,r}$ to $\cC_g^\ell\widehat{\otimes}_E (D^{\infty})^{g^{-1}}$ via (\ref{equ: CE graded embedding compact}) (and Lemma \ref{lem: pointed embedding}) is $|\cdot|_{\cC_g^\ell,r}\otimes_E (|\cdot|_r^{\infty})^{g^{-1}}$. Finally the norm $|\cdot|_{\cC_{U_0}^\ell,r}$ on $\cC_{U_0,r}^\ell$ induces a norm on the closed subspaces $\cD_{U_0,r}^\ell$, $\mathrm{ker}(\delta_{U_0,r}^\ell)$ and thus a norm on $\cH_{U_0,r}^\ell$ which, under (\ref{equ: p coh decomposition compact}), is explicitly
\begin{equation}\label{equ: CE coh norm}
|\cdot|_{\cH_{U_0}^\ell,r}\defeq \max_{g\in \tld{J}_{r}}\big(|\cdot|_{\cH_g^\ell,r}\otimes_E (|\cdot|_r^{\infty})^{g^{-1}}\big).
\end{equation}

We now prove that the image of the differential map $\delta_{U_0}^{\ell}$ is closed for $\ell\geq 0$. Recall that the Fr\'echet topology on $\cC_{U_0}^\ell$ is defined by the family of semi-norms $\{|\cdot|_{\cC_{U_0}^{\ell},r}\}_{r\in \cI}$. We write $|\cdot|_{\cC_{U_0}^{\ell},r}'$ for the induced semi-norm on the quotient $\cC_{U_0}^{\ell}/\mathrm{ker}(\delta_{U_0}^{\ell})$. Since $\delta_{1}^\ell:\cC_1^{\ell}/\mathrm{ker}(\delta_1^{\ell})\buildrel\sim\over\rightarrow \cD_1^{\ell+1}\subseteq \cC_1^{\ell+1}$ is a topological isomorphism by Lemma \ref{lem: p coh M separated g} (applied to $g=1$), given $r\in \cI$, there exists $r'\geq r$ and $A_{r,r'}\in\bQ_{>0}$ such that
\begin{equation}\label{equ: semi norm lower bound 1}
|x|_{\cC_{1}^{\ell},r}'\leq A_{r,r'}|\delta_{1}^{\ell}(x)|_{\cC_{1}^{\ell+1},r'}
\end{equation}
for $x\in \cC_{1}^{\ell}/\mathrm{ker}(\delta_{1}^{\ell})$. Let $g\in P_{I_1,0}$, by (\ref{equ: twist CE isom}) (applied with $1,g$ instead of $h,g$), we have a topological isomorphism $\theta_{1,g}^\ell: \cC_{1}^\ell\buildrel\sim\over\longrightarrow \cC_{g}^\ell$. This together with $\theta_{1,g}^{\ell+1}\circ\delta_{1}^\ell=\delta_{g}^\ell\circ\theta_{1,g}^\ell$, (\ref{equ: semi norm lower bound 1}) and (\ref{equ: CE inv norm}) implies
\begin{equation}\label{equ: semi norm lower bound g}
|y|_{\cC_{g}^{\ell},r}'\leq A_{r,r'}|\delta_{g}^{\ell}(y)|_{\cC_{g}^{\ell+1},r'}
\end{equation}
for $y\in \cC_{g}^{\ell}/\mathrm{ker}(\delta_{g}^{\ell})$. By continuity (\ref{equ: semi norm lower bound g}) holds for any $g\in P_{I_1,0}$ and any $y\in \cC_{g,r'}^{\ell}/\mathrm{ker}(\delta_{g,r'}^{\ell})$. Using (\ref{equ: differential decomposition}) (for $r$ and $r'$), we deduce for $y\in \cC_{U_0,r'}^{\ell}/\mathrm{ker}(\delta_{U_0,r'}^{\ell})$
\begin{equation}\label{equ: semi norm lower bound total}
|y|_{\cC_{U_0}^{\ell},r}'\leq A_{r,r'}|\delta_{U_0}^{\ell}(y)|_{\cC_{U_0}^{\ell+1},r'}.
\end{equation}
As (\ref{equ: semi norm lower bound total}) \emph{a fortiori} holds for any $y\in \cC_{U_0}^{\ell}/\mathrm{ker}(\delta_{U_0}^{\ell})$, we see that $\delta_{U_0}^\ell:\cC_{U_0}^{\ell}/\mathrm{ker}(\delta_{U_0}^{\ell})\hookrightarrow \cC_{U_0}^{\ell+1}$ is a closed embedding for $\ell\geq 0$ (in particular $\delta_{U_0}^{\ell}(\cC_{U_0}^{\ell})=\cD_{U_0}^{\ell+1}$).\bigskip

We finally prove the isomorphism (\ref{equ: p coh graded small}). Recall first that by Lemma \ref{limprojtensor} and \ref{it: support basis 2} of Lemma \ref{lem: support basis} we have isomorphisms
\begin{multline}\label{equ: CE coh smooth 2}
{D}^{\infty}(G_0)_{U_0}\widehat{\otimes}_{D^{\infty}(G_0)_{P_{w,0}}}\!D^{\infty}
\cong \varprojlim_{r\in\cI} \big(\!D^{\infty}(G_0)_{U_0,r}\!\otimes_{D^{\infty}(G_0)_{P_{w,0},r}}\!D^{\infty}_r\big)
\cong \varprojlim_{r\in\cI} \!\Big(\!\!\bigoplus_{g\in \tld{J}_{r}} (D^{\infty}_r)^{g^{-1}}\!\Big).
\end{multline}
By an argument similar to the proof of the injectivity of (\ref{equ: CE graded embedding compact}) in Lemma \ref{lem: pointed embedding}, we see that (\ref{equ: CE coh smooth 2}) contains as a dense subspace
\begin{equation}\label{equ: CE coh smooth 1}
\bigoplus_{g\in U/P_w} (D^{\infty})^{g^{-1}}=\bigoplus_{g\in U_0/P_{w,0}} (D^{\infty})^{g^{-1}}.
\end{equation}
The Fr\'echet topology on (\ref{equ: CE coh smooth 2}) is defined by $\{|\cdot|_{U_0,r}^{\infty}\}_{r\in \cI}$ where $|\cdot|_{U_0,r}^{\infty}\defeq \max_{g\in \tld{J}_{r}}(|\cdot|_r^{\infty})^{g^{-1}})$, and (\ref{equ: CE coh smooth 2}) is the completion of (\ref{equ: CE coh smooth 1}) under this topology. In particular, the left hand side of (\ref{equ: p coh graded small}) is the completion of the left hand side of (\ref{equ: CE coh trivialization compact}) under the Fr\'echet topology defined by the family of semi-norms
\begin{equation}\label{equ: family of semi norm 2}
\{|\cdot|_{\cH_{1}^\ell,r}\otimes_E|\cdot|_{U_0,r}^{\infty}\}_{r\in \cI}.
\end{equation}

On the other hand the Fr\'echet topology on $\cH_{U_0}^\ell$ is defined by the family of semi-norms
\begin{equation}\label{equ: family of semi norm 1}
\{|\cdot|_{\cH_{U_0}^\ell,r}\}_{r\in \cI}
\end{equation}
from (\ref{equ: CE coh norm}).
Since the map (\ref{equ: CE graded coh embedding compact}) is injective with dense image by Lemma~\ref{lem: CE coh total small}, $\cH_{U_0}^\ell$ can be identified with the completion of $\cH_{U_0,\flat}^\ell$ under $\{|\cdot|_{\cH_{U_0}^\ell,r}\}_{r\in \cI}$.\bigskip

By (\ref{equ: CE inv norm coh}) (applied with $1,g$ instead of $h,g$), for $g\in P_{I_1,0}$ the norm $|\cdot|_{\cH_{1}^\ell,r}$ on $\cH_{1,r}^\ell$ is identical to the norm $|\cdot|_{\cH_{g}^\ell,r}$ on $\cH_{g,r}^\ell$ under (\ref{equ: twist CE coh isom r}). Consequently, the two families of semi-norms (\ref{equ: family of semi norm 1}) and (\ref{equ: family of semi norm 2}) are identical under the isomorphism $\zeta_{U_0}^\ell$. It follows that $\zeta_{U_0}^\ell$ uniquely extends by completion into a topological isomorphism of $E$-Fr\'echet spaces as in (\ref{equ: p coh graded small}).
\end{proof}

We can finally prove the first important result of this section.

\begin{thm}\label{prop: p coh graded}
We have a (left) $D^{\infty}(P_{I_1})$-equivariant isomorphism of Fr\'echet spa\-ces for $\ell\geq 0$
\begin{multline}\label{equ: p coh graded}
\mathrm{Ext}_{U(\fp_{I_1})}^{\ell}(L^{I_1}(\mu), M_0^w)\otimes_E\big(\widehat{D}^{\infty}(G)_{P_{I_1}P_w}\widehat{\otimes}_{D^{\infty}(G)_{P_w}}(\pi_0^{\infty,w})^\vee\big)\\
\buildrel\sim\over\longrightarrow
\mathrm{Ext}_{U(\fp_{I_1})}^{\ell}(L^{I_1}(\mu), \mathrm{gr}_w(V_0^\vee))
\end{multline}
with trivial $D^{\infty}(P_{I_1})$-action on the finite dimensional $E$-vector space $\mathrm{Ext}_{U(\fp_{I_1})}^{\ell}(L^{I_1}(\mu), M_0^w)$ and where $\widehat{D}^{\infty}(G)_{P_{I_1}P_w}$ is defined in (\ref{equ: smooth hat}).
\end{thm}
\begin{proof}
Recall that $\mathrm{Ext}_{U(\fp_{I_1})}^{\ell}(L^{I_1}(\mu), M_0^w)\buildrel\sim\over\rightarrow \mathrm{Ext}_{U(\fp_{I_1})}^{\ell}(L^{I_1}(\mu), \cM_0^w)$ is finite dimensional by Lemma \ref{lem: p coh M isom}. Let $U\in \cX$. If there exists $t\in T$ such that $U\subseteq tP_{I_1,0}P_w=tP_{I_1,0}t^{-1}P_w$, then replacing $G_0$ (resp.~$P_{I_1,0}$, $P_{w,0}$) by $tG_0t^{-1}$ (resp.~$tP_{I_1,0}t^{-1}$, $tP_{w,0}t^{-1}$), a parallel argument as in the proof of Lemma~\ref{lem: p coh graded small} shows that $\delta_{U}^\ell$ is strict and that $\zeta_{U}^\ell$ uniquely extends into a topological isomorphism as in (\ref{equ: p coh graded small}). In general, since $U/P_w$ is compact there exists a finite partition $U=\bigsqcup_{i}U_i$ and elements $t_i\in T$ such that $U_i\in \cX$ and $U_i\subseteq t_iP_{I_1,0}P_w=t_iP_{I_1,0}t_i^{-1}P_w$ for each $i$ (use the open covering $(tP_{I_1,0}t^{-1}/(tP_{I_1,0}t^{-1}\cap P_w))_{t\in T}$ of $P_{I_1}/(P_{I_1}\cap P_w)$). Since $\cC_{U}^\ell\cong \bigoplus_{i}\cC_{U_i}^\ell$, $\delta_{U}^\ell=\bigoplus_{i} \delta_{U_i}^\ell$ and $\zeta_{U}^\ell=\bigoplus_{i} \zeta_{U_i}^\ell$, we deduce that $\delta_{U}^\ell$ is strict and that $\zeta_{U}^\ell$ uniquely extends into a topological isomorphism (\ref{equ: p coh graded small}). Since $\cC_{w}^\ell\cong \varprojlim_{U\in\cX}\cC_{U}^\ell$ with $\delta_{w}^\ell=\varprojlim_{U\in\cX}\delta_{U}^\ell$ (see (\ref{equ: CE graded limit})), we deduce that $\delta_{w}^{\ell-1}$ is strict with closed image $\cD_w^\ell$ and thus $H^\ell(\cC_{w}^\bullet)=\cH_{w}^\ell$. Taking $\varprojlim_{U\in\cX}$ on (\ref{equ: p coh graded small}), we obtain a topological isomorphism
\begin{equation}\label{equ: p coh graded as limit}
\varprojlim_{U\in\cX}\zeta_{U}^\ell: \cH_{1}^\ell\otimes_E \big(\widehat{D}^\infty(G)_{X}\widehat{\otimes}_{D^\infty(G)_{P_w}}D^{\infty}\big) \buildrel\sim\over\longrightarrow \cH_{w}^\ell=H^\ell(\cC_{w}^\bullet)
\end{equation}
as the composition of the topological isomorphisms
\begin{multline*}
\cH_{1}^\ell\otimes_E \big(\widehat{D}^{\infty}(G)_{X}\widehat{\otimes}_{D^{\infty}(G)_{P_w}}D^{\infty}\big)
\cong \cH_{1}^\ell\otimes_E \big((\varprojlim_{U \in \cX}D^{\infty}(G)_{U})\widehat{\otimes}_{D^{\infty}(G)_{P_w}}D^{\infty}\big)\\
\cong \varprojlim_{U \in \cX} \big(\cH_{1}^\ell\otimes_E D^{\infty}(G)_{U}\widehat{\otimes}_{D^{\infty}(G)_{P_w}}D^{\infty}\big)
\cong \varprojlim_{U \in \cX} \cH_{U}^\ell \cong \cH_{w}^\ell.
\end{multline*}
Here we have used the topological isomorphisms
\begin{multline*}
\widehat{D}^{\infty}(G)_{X}\widehat{\otimes}_{D^{\infty}(G)_{P_w}}D^{\infty}\cong \widehat{D}^{\infty}(G_0)_{X_0}\widehat{\otimes}_{D^{\infty}(G_0)_{P_{w,0}}}D^{\infty}\\
\cong \varprojlim_{U_0 \in \cX_0}\big(D^{\infty}(G_0)_{U_0}\widehat{\otimes}_{D^{\infty}(G_0)_{P_{w,0}}}D^{\infty}\big)\cong
\varprojlim_{U \in \cX}\big(D^{\infty}(G)_{U}\widehat{\otimes}_{D^{\infty}(G)_{P_w}}D^{\infty}\big)
\end{multline*}
which follows by similar argument as for (\ref{equ: graded short term}). The source (resp.~target) of (\ref{equ: p coh graded as limit}) contains $\cH_{1}^\ell\otimes_E \big(\bigoplus_{g\in X/P_w} (D^{\infty})^{g^{-1}}\big)$ (resp.~$\cH_{w,\flat}^\ell$) as a $P_{I_1}$-stable dense subspace and the isomorphism (\ref{equ: p coh graded as limit}) restricts to the isomorphism $\zeta^\ell$ of (\ref{equ: CE coh trivialization}) on these dense subspaces. Then Lemma~\ref{lem: equivariance on dense subspace} forces (\ref{equ: p coh graded as limit}) to be $D(P_{I_1})$-equivariant (hence $D^\infty(P_{I_1})$-equivariant). Finally, by (\ref{equ: CE conj}) and (\ref{equ: Ext g vers coh}) the isomorphism (\ref{equ: p coh graded as limit}) is (\ref{equ: p coh graded}).
\end{proof}

We write $\Sigma_i\defeq W(L_{I_i})\cdot\cJ(\pi_i^{\infty})\subseteq \widehat{T}^{\infty}$ for $i=0,1$ where $\cJ(\pi_i^{\infty})$ is in Definition \ref{def: basic rep} and the dot action $\cdot$ in (\ref{equ: smooth dot action}). Recall that the categories $\cB^{I_i}_{\Sigma_i}$ are defined above (\ref{equ: sm block decomposition}) and that any $\pi_i^{\infty}$ in $\cB^{I_i}_{\Sigma_i}$ is of finite length. The following corollary to Theorem \ref{prop: p coh graded} is crucial.

\begin{cor}\label{cor: Ext P graded}
Keep the setting of Theorem \ref{prop: p coh graded} and assume moreover that the smooth representations $\pi_i^{\infty}$ of $L_{I_i}$ are of finite length for $i=0,1$. We have canonical isomorphisms of finite dimensional $E$-vector spaces for $w\in W^{I_0,I_1}$ and $k,\ell\geq 0$
\begin{multline}\label{equ: extp coh graded}
\mathrm{Ext}_{D^{\infty}(P_{I_1})}^k\big((\pi_1^{\infty})^\vee, \mathrm{Ext}_{U(\fp_{I_1})}^{\ell}(L^{I_1}(\mu), \mathrm{gr}_w(V_0^\vee))\big)\\
\cong \mathrm{Ext}_{U(\fp_{I_1})}^{\ell}(L^{I_1}(\mu), M_0^w)\otimes_E \mathrm{Ext}_{L_{I_1}}^k(i_{I_0,I_1,w}^{\infty}(J_{I_0,I_1,w}(\pi_0^{\infty})), \pi_1^{\infty})^{\infty}
\end{multline}
where $i_{I_0,I_1,w}^{\infty}$ is defined in (\ref{equ: Bruhat induction}) and $J_{I_0,I_1,w}$ in (\ref{equ: Bruhat Jacquet}). If moreover $\Sigma_i$ is a single $G$-regular $W(L_{I_i})$-coset, $\pi_i^{\infty}$ is in $\cB^{I_i}_{\Sigma_i}$ for $i=0,1$ and $\Sigma_1\cap W(G)\cdot\Sigma_0\neq\emptyset$, then (\ref{equ: extp coh graded}) is non-zero for at most one $w\in W^{I_0,I_1}$, which is the unique $w\in W^{I_0,I_1}$ such that $\Sigma_1\cap w^{-1}\cdot\Sigma_0\neq \emptyset$.
\end{cor}
\begin{proof}
Recall from Lemma~\ref{lem: dual of compact induction} and its proof (applied with $P=P_w$, $X=P_{I_1}P_w$ and $Q=P_{I_1}$) that $(\mathrm{ind}_{P_w}^{P_{I_1}P_w}\pi_0^{\infty,w})^{\infty}$ (with the finest locally convex topology) is a smooth representation of $P_{I_1}$ on a vector space of compact type and that we have isomorphisms of $E$-Fr\'echet spaces with separately continuous $D^\infty(P_{I_1})$-actions
\begin{equation}\label{equ: dual sm compact induction}
\widehat{D}^\infty(G)_{P_{I_1}P_w}\widehat{\otimes}_{D^\infty(G)_{P_w}}(\pi_0^{\infty,w})^\vee\cong \big((\mathrm{ind}_{P_w}^{P_{I_1}P_w}\pi_0^{\infty,w})^{\infty}\big)^\vee \cong \big((\mathrm{ind}_{P_{I_1}\cap P_w}^{P_{I_1}}\pi_0^{\infty,w})^{\infty}\big)^\vee.
\end{equation}
From (\ref{equ: dual sm compact induction}) and Theorem~\ref{prop: p coh graded} we deduce isomorphisms for $k,\ell\geq 0$ (using that $D^{\infty}(P_{I_1})$ acts trivially on $\mathrm{Ext}_{U(\fp_{I_1})}^{\ell}(L^{I_1}(\mu), M_0^w)$)
\begin{multline}\label{equ: Ext graded sm compact induction}
\mathrm{Ext}_{D^{\infty}(P_{I_1})}^k\big((\pi_1^{\infty})^\vee, \mathrm{Ext}_{U(\fp_{I_1})}^{\ell}(L^{I_1}(\mu), \mathrm{gr}_w(V_0^\vee))\big)\\
\cong \mathrm{Ext}_{U(\fp_{I_1})}^{\ell}(L^{I_1}(\mu), M_0^w)\otimes_E \mathrm{Ext}_{D^{\infty}(P_{I_1})}^k\big((\pi_1^{\infty})^\vee, ((\mathrm{ind}_{P_{I_1}\cap P_w}^{P_{I_1}}\pi_0^{\infty,w})^{\infty})^\vee\big).
\end{multline}
By Lemma \ref{erratum} we have for $k_1\geq 0$ (writing $1$ for the trivial representation of $N_{I_1}$)
\begin{equation}\label{equ: sm coh equal}
\mathrm{Ext}_{D^{\infty}(N_{I_1})}^{k_1}\Big(1, \big((\mathrm{ind}_{P_{I_1}\cap P_w}^{P_{I_1}}\pi_0^{\infty,w})^{\infty}\big)^\vee\Big)\cong \mathrm{Ext}_{N_{I_1}}^{k_1}\Big(\big((\mathrm{ind}_{P_{I_1}\cap P_w}^{P_{I_1}}\pi_0^{\infty,w})^{\infty}\big), 1\Big)^\infty.
\end{equation}
Since Jacquet functors are exact on smooth representations, we deduce from (\ref{equ: sm coh equal}) that we have for $k_1>0$
\begin{equation}\label{equ: sm coh vanishing}
\mathrm{Ext}_{D^{\infty}(N_{I_1})}^{k_1}\Big(1, \big((\mathrm{ind}_{P_{I_1}\cap P_w}^{P_{I_1}}\pi_0^{\infty,w})^{\infty}\big)^\vee\Big)=0
\end{equation}
and a canonical $D^{\infty}(L_{I_1})$-equivariant isomorphism of $E$-Fr\'echet spaces
\begin{equation}\label{equ: sm coinv dual}
\Hom_{D^{\infty}(N_{I_1})}\Big(1, \big((\mathrm{ind}_{P_{I_1}\cap P_w}^{P_{I_1}}\pi_0^{\infty,w})^{\infty}\big)^\vee\Big)\cong \Big(\big((\mathrm{ind}_{P_{I_1}\cap P_w}^{P_{I_1}}\pi_0^{\infty,w})^{\infty}\big)_{N_{I_1}}\Big)^\vee.
\end{equation}
Combining (\ref{equ: sm coh vanishing}) and (\ref{equ: sm coinv dual}) with the (standard ``Hochschild-Serre type'') spectral sequence
\[\mathrm{Ext}_{D^{\infty}(L_{I_1})}^{k_2}\big((\pi_1^{\infty})^\vee, \mathrm{Ext}_{D^{\infty}(N_{I_1})}^{k_1}(1,-)\big)\implies \mathrm{Ext}_{D^{\infty}(P_{I_1})}^{k_1+k_2}((\pi_1^{\infty})^\vee,-)\]
we deduce isomorphisms for $k\geq 0$
\begin{multline}\label{equ: sm pass to Levi}
\mathrm{Ext}_{D^{\infty}(P_{I_1})}^k\Big((\pi_1^{\infty})^\vee, \big((\mathrm{ind}_{P_{I_1}\cap P_w}^{P_{I_1}}\pi_0^{\infty,w})^{\infty}\big)^\vee\Big)\\
\cong \mathrm{Ext}_{D^{\infty}(L_{I_1})}^k\Big((\pi_1^{\infty})^\vee, \big(((\mathrm{ind}_{P_{I_1}\cap P_w}^{P_{I_1}}\pi_0^{\infty,w})^{\infty})_{N_{I_1}}\big)^\vee\Big).
\end{multline}
Recall from (\ref{equ: N coinvariant induction}) the isomorphism of smooth admissible representations of $L_{I_1}$
\begin{equation*}
\big((\mathrm{ind}_{P_{I_1}\cap P_w}^{P_{I_1}}\pi_0^{\infty,w})^{\infty}\big)_{N_{I_1}}\cong i_{I_0,I_1,w}^{\infty}(J_{I_0,I_1,w}(\pi_0^{\infty})),
\end{equation*}
which together with (\ref{equ: sm pass to Levi}), (\ref{equ: Ext graded sm compact induction}) and Lemma \ref{erratum} gives the isomorphism (\ref{equ: extp coh graded}). The finite dimensionality of $\mathrm{Ext}_{L_{I_1}}^k(i_{I_0,I_1,w}^{\infty}(J_{I_0,I_1,w}(\pi_0^{\infty})), \pi_1^{\infty})^{\infty}$ comes from the finite length of both representations $i_{I_0,I_1,w}^{\infty}(J_{I_0,I_1,w}(\pi_0^{\infty}))$ and $\pi_1^{\infty}$ (\cite[\S3 Cor.~3]{SS93} noting that $i_{I_0,I_1,w}^{\infty}(J_{I_0,I_1,w}(\pi_0^{\infty}))$ is of finite length as $\pi_0^{\infty}$ is by \cite[\S VI.6.4]{Re10} and \cite[\S VI.6.2]{Re10}). For the last statement, note first that by the definition of regularity (Definition \ref{def: basic rep}), we have $\Sigma_1\cap W(G)\cdot\Sigma_0\neq\emptyset$ if and only if there exists $w\in W(G)$ such that $\Sigma_1\cap w^{-1}\cdot\Sigma_0\neq \emptyset$, and we can take $w$ in $ W^{I_0,I_1}$ which is then unique. The last statement then follows from the first and \ref{it: sml2} of Lemma \ref{lem: smooth geometric lemma} (applied with $L_I=G$).
\end{proof}

\begin{rem}\label{rem: functorial isom}
For $i=0,1$ and $\Sigma_i$ as before Corollary \ref{cor: Ext P graded}, assume that $\pi_i^{\infty}$ is in $\cB^{I_i}_{\Sigma_i}$. For $w\in W^{I_0,I_1}$, the isomorphism (\ref{equ: p coh graded}) is functorial in $\pi_0^{\infty}$ and the isomorphism (\ref{equ: extp coh graded}) is functorial in both $\pi_0^{\infty}$ and $\pi_1^{\infty}$.
\end{rem}

We now use Corollary \ref{cor: Ext P graded} to derive several results on the groups $\mathrm{Ext}_{G}^\bullet(V_0,V_1)$ of (\ref{extdef}).

\begin{cor}\label{prop: OS spectral seq}
Keep the setting of Theorem \ref{prop: p coh graded}, let $\Sigma_i= W(L_{I_i})\cdot\cJ(\pi_i^{\infty})\subseteq \widehat{T}^{\infty}$ and assume moreover that the smooth representations $\pi_i^{\infty}$ of $L_{I_i}$ are of finite length ($i=0,1$).
\begin{enumerate}[label=(\roman*)]
\item \label{it: OS spectral seq 0} The $E$-vector space $\mathrm{Ext}_{G}^{k}(V_0,V_1)$ is finite dimensional for $k\geq 0$.
\item \label{it: OS spectral seq 1} Assume that $\pi_i^{\infty}$ is in $\cB^{I_i}_{\Sigma_i}$ for $i=0,1$. If $\Sigma_1\cap W(G)\cdot\Sigma_0=\emptyset$, then $\mathrm{Ext}_{G}^{k}(V_0,V_1)=0$ for $k\geq 0$.
\item \label{it: OS spectral seq 2} Assume that $\pi_i^{\infty}$ is in $\cB^{I_i}_{\Sigma_i}$ for $i=0,1$. Let $\xi: Z(\fl_{I_1})\rightarrow E$ be the unique infinitesimal character such that $L^{I_1}(\mu)_\xi\neq 0$. If $\mathrm{Ext}_{G}^{k}(V_0,V_1)\neq 0$ for some $k\geq 0$, then there exists $w\in W^{I_0,I_1}$ and $\ell\leq k$ such that $\Sigma_1\cap w^{-1}\cdot\Sigma_0\neq \emptyset$ and $H^{\ell}(\fn_{I_1},M_0^w)_\xi\neq 0$ (see before Lemma \ref{lem: central component} for the notation).
\item \label{it: OS spectral seq 3} Assume that $\Sigma_i$ is a single $G$-regular $W(L_{I_i})$-coset and that $\pi_i^{\infty}$ is in $\cB^{I_i}_{\Sigma_i}$ for $i=0,1$. If $\Sigma_1\cap W(G)\cdot\Sigma_0\neq\emptyset$, then there exists a unique $w\in W^{I_0,I_1}$ such that $\Sigma_1\cap w^{-1}\cdot\Sigma_0\neq \emptyset$, and we have a spectral sequence
\begin{equation}\label{equ: OS main seq}
\mathrm{Ext}_{U(\fg)}^{\ell}(M_1, M_0^w)\otimes_E \mathrm{Ext}_{L_{I_1}}^k\big(i_{I_0,I_1,w}^{\infty}(J_{I_0,I_1,w}(\pi_0^{\infty})), \pi_1^{\infty}\big)^{\infty}\implies \mathrm{Ext}_{G}^{k+\ell}(V_0,V_1).
\end{equation}
\end{enumerate}
\end{cor}
\begin{proof}
We \ prove \ \ref{it: OS spectral seq 0}. \ By \ the \ first \ statement \ of \ Corollary~\ref{cor: Ext P graded} \ the \ $E$-vector \ space $\mathrm{Ext}_{D^{\infty}(P_{I_1})}^k((\pi_1^{\infty})^\vee, \mathrm{Ext}_{U(\fp_{I_1})}^{\ell}(L^{I_1}(\mu), \mathrm{gr}_w(V_0^\vee)))$ is finite dimensional for $k\geq 0$, $\ell\geq 0$ and $w\in W^{I_0,I_1}$. By the spectral sequence (\ref{equ: ST seq}) this implies that $\mathrm{Ext}_{D(P_{I_1})}^k(L^{I_1}(\mu)\otimes_E (\pi_1^{\infty})^\vee, \mathrm{gr}_w(V_0^\vee))$ is finite dimensional for $k\geq 0$ and $w\in W^{I_0,I_1}$. By an obvious d\'evissage and (\ref{isov1v0}) we deduce that $\mathrm{Ext}_{G}^k(V_0,V_1)$ is finite dimensional for $k\geq 0$.

We prove \ref{it: OS spectral seq 1}. As $\Sigma_1\cap W(G)\cdot\Sigma_0=\emptyset$ and by (\ref{equ: support bound}), we have $i_{I_0,I_1,w}^{\infty}(J_{I_0,I_1,w}(\pi_0^{\infty}))_{\cB^{I_1}_{\Sigma_1}}=0$ using the property of the Bernstein block $\cB^{I_1}_{\Sigma_1}$ (see the paragraph after Remark \ref{rem: twist Jacquet support}). Since $i_{I_0,I_1,w}^{\infty}(J_{I_0,I_1,w}(\pi_0^{\infty}))$ and $\pi_1^{\infty}$ live in different Bernstein blocks, we deduce for $k\geq 0$:
\[\mathrm{Ext}_{L_{I_1}}^k(i_{I_0,I_1,w}^{\infty}(J_{I_0,I_1,w}(\pi_0^{\infty})), \pi_1^{\infty})^{\infty}=0.\]
Together with (\ref{equ: ST seq}) applied with $D=\mathrm{gr}_w(V_0^\vee)$, (\ref{equ: extp coh graded}) and a d\'evissage with respect to the $W^{I_0,I_1}$-filtration $\mathrm{Fil}_{\bullet}(V_0^\vee)$ on $V_0^\vee$, we deduce $\mathrm{Ext}_{D(P_{I_1})}^{k}(L^{I_1}(\mu)\otimes_E (\pi_1^{\infty})^\vee, V_0^\vee)=0$ for $k\geq 0$ and hence by (\ref{isov1v0}) $\mathrm{Ext}_{G}^{k}(V_0,V_1)=0$ for $k\geq 0$.

We prove \ref{it: OS spectral seq 2}. If $\mathrm{Ext}_{G}^{k}(V_0,V_1)\neq 0$ we have $\mathrm{Ext}_{D(P_{I_1})}^k(L^{I_1}(\mu)\otimes_E (\pi_1^{\infty})^\vee, V_0^\vee)\ne 0$ by (\ref{isov1v0}), and thus there exists $w\in W^{I_0,I_1}$ such that $\mathrm{Ext}_{D(P_{I_1})}^k(L^{I_1}(\mu)\otimes_E (\pi_1^{\infty})^\vee, \mathrm{gr}_w(V_0^\vee))\neq 0$. By (\ref{equ: ST seq}) (with $D=\mathrm{gr}_w(V_0^\vee)$), (\ref{equ: extp coh graded}) and Lemma \ref{lem: p coh M isom} (together with Shapiro's lemma for Lie algebra cohomology) there exist $\ell,k'\leq k$ such that $\mathrm{Ext}_{U(\fg)}^{\ell}(M_1, M_0^w)\neq 0$ and $\mathrm{Ext}_{L_{I_1}}^{k'}(i_{I_0,I_1,w}^{\infty}(J_{I_0,I_1,w}(\pi_0^{\infty})), \pi_1^{\infty})^{\infty}\ne 0$. The latter implies $\Sigma_1\cap w^{-1}\cdot\Sigma_0\neq \emptyset$ (otherwise $i_{I_0,I_1,w}^{\infty}(J_{I_0,I_1,w}(\pi_0^{\infty}))$ and $\pi_1^{\infty}$ would live in distinct blocks by (\ref{equ: support bound}) and the paragraph after Remark \ref{rem: twist Jacquet support}). By (\ref{equ: g spectral seq}) the former implies $\mathrm{Ext}_{U(\fl_{I_1})}^{\ell_1}(L^{I_1}(\mu), H^{\ell_2}(\fn_{I_1}, M_0^w))\neq 0$ for some $\ell_1,\ell_2\geq 0$ such that $\ell_1+\ell_2=\ell$. We thus deduce $H^{\ell_2}(\fn_{I_1}, M_0^w)_\xi\neq 0$ from \ref{it: block 3} of Lemma~\ref{lem: central component}.

We prove \ref{it: OS spectral seq 3}. By the last statement in Corollary~\ref{cor: Ext P graded} we have that $w$ is unique, and together with (\ref{equ: ST seq}) applied with $D=\mathrm{gr}_{w'}(V_0^\vee)$ for $w'\neq w\in W^{I_0,I_1}$ we deduce $\mathrm{Ext}_{D(P_{I_1})}^k(L^{I_1}(\mu)\otimes_E (\pi_1^{\infty})^\vee, \mathrm{gr}_{w'}(V_0^\vee))=0$ for $k\geq 0$ and such $w'$. Via a d\'evissage with respect to the $W^{I_0,I_1}$-filtration $\mathrm{Fil}_{\bullet}(V_0^\vee)$ on $V_0^\vee$, we deduce for $k\geq 0$
\begin{equation}\label{equ: unique graded}
\mathrm{Ext}_{D(P_{I_1})}^k(L^{I_1}(\mu)\otimes_E (\pi_1^{\infty})^\vee, V_0^\vee)\cong \mathrm{Ext}_{D(P_{I_1})}^k(L^{I_1}(\mu)\otimes_E (\pi_1^{\infty})^\vee, \mathrm{gr}_w(V_0^\vee)).
\end{equation}
Finally, we apply (\ref{equ: ST seq}) with $D=\mathrm{gr}_w(V_0^\vee)$, which together with (\ref{equ: extp coh graded}) and the isomorphisms (\ref{equ: unique graded}), (\ref{isov1v0}) give the spectral sequence (\ref{equ: OS main seq}).
\end{proof}

By the last statement in \ref{it: PS 1} of Lemma \ref{lem: Jacquet of PS} and the comment below Lemma \ref{lem: Jacquet of PS}, Corollary \ref{cor: Ext P graded} and \ref{it: OS spectral seq 3} of Corollary \ref{prop: OS spectral seq} can in particular be applied to $G$-basic $\pi_i^\infty$.

\begin{rem}
If $I_1=\Delta$ (with $\mu_1\in\Lambda^{\dom}$ and $V_1\cong L(\mu_1)^\vee\otimes_E \pi_1^{\infty}$), then $W^{I_0,\Delta}=\{1\}$ and the same argument as at the end of the proof of \ref{it: OS spectral seq 3} of Corollary~\ref{prop: OS spectral seq} (for $\pi_0^\infty,\pi_1^\infty$ of finite length) gives a spectral sequence of finite dimensional $E$-vector spaces
\begin{equation*}
\mathrm{Ext}_{U(\fg)}^{\ell}(L(\mu_1), M_0)\otimes_E\mathrm{Ext}_{G}^k(i_{I_0}^{\infty}(\pi_0^{\infty}), \pi_1^{\infty})^{\infty}\implies \mathrm{Ext}_{G}^{k+\ell}(V_0,V_1).
\end{equation*}
If moreover $I_0=\Delta$ (with $\mu_0\in\Lambda^{\dom}$ and $V_0\cong L(\mu_0)^\vee\otimes_E \pi_0^{\infty}$), we obtain a spectral sequence
\begin{equation}\label{equ: Ext loc alg 2}
\mathrm{Ext}_{U(\fg)}^{\ell}(L(\mu_1), L(\mu_0))\otimes_E\mathrm{Ext}_{G}^k(\pi_0^{\infty}, \pi_1^{\infty})^{\infty} \implies \mathrm{Ext}_{G}^{k+\ell}(V_0,V_1).
\end{equation}
Note that, as $\mu_0,\mu_1\in\Lambda^{\dom}$, we have $L(\mu_0)_{\xi}\neq 0\neq L(\mu_1)_{\xi}$ for some $\xi: Z(\fg)\rightarrow E$ if and only if $\mu_0=\mu_1$. Hence by \ref{it: block 3} of Lemma~\ref{lem: central component} we have $\mathrm{Ext}_{U(\fg)}^{\ell}(L(\mu_1), L(\mu_0))\neq 0$ for some $\ell\geq 0$ if and only if $\mu_0=\mu_1$, in which case a translation functor argument gives a canonical isomorphism $\mathrm{Ext}_{U(\fg)}^{\ell}(L(\mu_0), L(\mu_0))\cong H^{\ell}(\fg,1_{\fg})$ for $\ell \geq 0$.
Consequently, when $I_0=I_1=\Delta$, $\mathrm{Ext}_{G}^k(V_0,V_1)\neq 0$ for some $k\geq 0$ only if $\mu_0=\mu_1$, in which case $\mathrm{Ext}_{G}^k(V_0,V_1)=0$ for $k<d(\pi_0^{\infty},\pi_1^{\infty})$ (see \ref{it: basic 3} of Definition \ref{def: basic rep} for $d(\pi_0^{\infty},\pi_1^{\infty})$) and we have a canonical isomorphism
\begin{equation}\label{equ: reduce to sm Ext}
\mathrm{Ext}_{G}^{d(\pi_0^{\infty},\pi_1^{\infty})}(\pi_0^{\infty},\pi_1^{\infty})^{\infty}\buildrel\sim\over\longrightarrow \mathrm{Ext}_{G}^{d(\pi_0^{\infty},\pi_1^{\infty})}(V_0,V_1).
\end{equation}
\end{rem}

\begin{lem}\label{lem: distinct block}
For $i=0,1$ let $V_i= \cF_{P_{I_i}}^G(M_i,\pi_i^{\infty})$ with $I_i\subseteq \Delta$, $M_i$ in $\cO^{\fp_{I_i}}_{\rm{alg}}$ and $\pi_i^{\infty}$ smooth finite length representations of $L_{I_i}$ such that $\pi_i^{\infty}$ is in $\cB^{I_i}_{\Sigma_i}$ with $\Sigma_i= W(L_{I_i})\cdot\cJ(\pi_i^{\infty})\subseteq \widehat{T}^{\infty}$. Assume that
\begin{equation}\label{equ: distinct block}
\Sigma_0\cap W(G)\cdot \Sigma_1=\emptyset.
\end{equation}
Then we have for $k\geq 0$
\begin{equation}\label{equ: block vanishing}
\mathrm{Ext}_{G}^k(V_0,V_1)=0.
\end{equation}
\end{lem}
\begin{proof}
By d\'evissage we can assume $M_1$ simple of the form $L(\mu)$ for some $\mu\in\Lambda_{I_1}^{\dom}$. If $M_1=L(\mu)=M^{I_1}(\mu)$, then (\ref{equ: block vanishing}) follows directly from (\ref{equ: distinct block}) and \ref{it: OS spectral seq 1} of Corollary~\ref{prop: OS spectral seq}. In general, we assume inductively that
\begin{equation}\label{equ: distinct block lower wt}
\mathrm{Ext}_{G}^k(V_0,\cF_{P_{I_1}}^{G}(L(\mu'),\pi_1^{\infty}))=0
\end{equation}
for $\mu'\in\Lambda_{I_1}^{\dom}$ such that $\mu'-\mu\in\Z_{\geq 0}\Phi^+$ and $\mu'\neq \mu$. Recall $M^{I_1}(\mu)$ fits into $0\rightarrow N^{I_1}(\mu)\rightarrow M^{I_1}(\mu)\rightarrow L(\mu)\rightarrow 0$ with all Jordan-H\"older factors of $N^{I_1}(\mu)$ of the form $L(\mu')$ for some $\mu'\in\Lambda_{I_1}^{\dom}$ such that $\mu'-\mu\in\Z_{\geq 0}\Phi^+$ and $\mu'\neq \mu$ (use \cite[Thm.~5.1]{Hum08}). This together with (\ref{equ: distinct block lower wt}) (and \ref{it: OS property 1} of Theorem~\ref{prop: OS property}) implies $\mathrm{Ext}_{G}^k(V_0,\cF_{P_{I_1}}^{G}(N^{I_1}(\mu),\pi_1^{\infty}))=0$ for $k\geq 0$. Using $0\rightarrow V_1\rightarrow \cF_{P_{I_1}}^{G}(M^{I_1}(\mu),\pi_1^{\infty})\rightarrow \cF_{P_{I_1}}^{G}(N^{I_1}(\mu),\pi_1^{\infty})\rightarrow 0$ we obtain an isomorphism $\mathrm{Ext}_{G}^k(V_0,V_1)\buildrel\sim\over\rightarrow \mathrm{Ext}_{G}^k(V_0,\cF_{P_{I_1}}^{G}(M^{I_1}(\mu),\pi_1^{\infty}))$ for $k\geq 0$. We then again deduce (\ref{equ: block vanishing}) from (\ref{equ: distinct block}) and \ref{it: OS spectral seq 1} of Corollary~\ref{prop: OS spectral seq}, which finishes the proof by induction.
\end{proof}

We isolate the following result because it has its own interest.

\begin{thm}\label{thm: finitedim}
For $i=0,1$ let $V_i= \cF_{P_{I_i}}^G(M_i,\pi_i^{\infty})$ with $I_i\subseteq \Delta$, $M_i$ in $\cO^{\fp_{I_i}}_{\rm{alg}}$ and $\pi_i^{\infty}$ smooth finite length representations of $L_{I_i}$ over $E$. Then the $E$-vector space $\mathrm{Ext}_{G}^{k}(V_0,V_1)$ is finite dimensional for $k\geq 0$.
\end{thm}
\begin{proof}
By d\'evissage we can assume $M_1$ simple of the form $L(\mu)$ for some $\mu\in\Lambda_{I_1}^{\dom}$. If $M_1=L(\mu)=M^{I_1}(\mu)$, then the result is \ref{it: OS spectral seq 0} of Corollary \ref{prop: OS spectral seq}. In general, we argue by induction as in the proof of Lemma \ref{lem: distinct block} using the short exact sequence $0\rightarrow N^{I_1}(\mu)\rightarrow M^{I_1}(\mu)\rightarrow L(\mu)\rightarrow 0$.
\end{proof}

We can obtain better vanishing results when $\Sigma_i$ is a single $G$-regular left $W(L_{I_i})$-coset for $i=0,1$.

\begin{lem}\label{lem: improved vanishing}
For $i=0,1$ let $V_i= \cF_{P_{I_i}}^G(M_i,\pi_i^{\infty})$ with $I_i\subseteq \Delta$, $M_i$ in $\cO^{\fp_{I_i}}_{\rm{alg}}$ and $\pi_i^{\infty}$ smooth finite length representations of $L_{I_i}$ such that $\pi_i^{\infty}$ is in $\cB^{I_i}_{\Sigma_i}$ with $\Sigma_i= W(L_{I_i})\cdot\cJ(\pi_i^{\infty})\subseteq \widehat{T}^{\infty}$. Assume that $\Sigma_i$ is a single $G$-regular left $W(L_{I_i})$-coset and recall that $d(\pi_0^{\infty},\pi_1^{\infty})$ is defined in \ref{it: basic 3} of Definition \ref{def: basic rep}.
\begin{enumerate}[label=(\roman*)]
\item \label{it: Ext vanishing 1} We have $\mathrm{Ext}_{G}^k(V_0, V_1)=0$ for $k<d(\pi_0^{\infty},\pi_1^{\infty})$.
\item \label{it: Ext vanishing 2} Let $w\in W^{I_0,I_1}$ such that $\Sigma\defeq \Sigma_1\cap w^{-1}\cdot\Sigma_0\ne \emptyset$ and $I\defeq w^{-1}(I_0)\cap I_1$. Then $\Sigma$ is a single $G$-regular $W(L_I)$-coset and, if $J_{I_0,I_1,w}(\pi_0^{\infty})_{\cB^{I}_{\Sigma}}=0$, we have $\mathrm{Ext}_{G}^k(V_0, V_1)=0$ for $k\geq 0$.
\end{enumerate}
\end{lem}
\begin{proof}
We prove \ref{it: Ext vanishing 1}.
By (\ref{equ: first adjunction}) and \ref{it: sml1} of Lemma~\ref{lem: smooth geometric lemma} (applied with $I=\Delta$) we have for $k\geq 0$
\[\mathrm{Ext}_{G}^k(i_{I_0,\Delta}^{\infty}(\pi_0^{\infty}),i_{I_1,\Delta}^{\infty}(\pi_1^{\infty}))^{\infty}\cong \bigoplus_{w\in W^{I_0,I_1}}\mathrm{Ext}_{L_{I_1}}^k(i_{I_0,I_1,w}^{\infty}(J_{I_0,I_1,w}(\pi_0^{\infty})),\pi_1^{\infty})\]
and thus $\mathrm{Ext}_{L_{I_1}}^k(i_{I_0,I_1,w}^{\infty}(J_{I_0,I_1,w}(\pi_0^{\infty})),\pi_1^{\infty})=0$ for $w\in W^{I_0,I_1}$ and $k<d(\pi_0^{\infty},\pi_1^{\infty})$. This together with (\ref{equ: OS main seq}) gives \ref{it: Ext vanishing 1} when $M_1=M^{I_1}(\mu_1)$ for some $\mu_1\in\Lambda_{I_1}^{\dom}$. The result for general $M_1$ follows from the same induction as in the proof of Lemma~\ref{lem: distinct block}.

We prove \ref{it: Ext vanishing 2}. The fact that $\Sigma$ is a single $G$-regular $W(L_I)$-coset follows easily from the fact that $\Sigma_i$, for $i=0,1$, is also a single $G$-regular $W(L_I)$-coset (last statement in \ref{it: PS 1} of Lemma \ref{lem: Jacquet of PS}). We have for $k\geq 0$ (see (\ref{extranot}) for $J_{I_1,I}'(\pi_1^{\infty})$)
\begin{multline*}
\mathrm{Ext}_{L_{I_1}}^k(i_{I_0,I_1,w}^{\infty}(J_{I_0,I_1,w}(\pi_0^{\infty})), \pi_1^{\infty})^{\infty}\cong \mathrm{Ext}_{L_I}^k(J_{I_0,I_1,w}(\pi_0^{\infty}), J_{I_1,I}'(\pi_1^{\infty}))^{\infty}\\
\cong \mathrm{Ext}_{L_I}^k(J_{I_0,I_1,w}(\pi_0^{\infty})_{\cB^{I}_{\Sigma}}, J_{I_1,I}'(\pi_1^{\infty})_{\cB^{I}_{\Sigma}})^{\infty}=0
\end{multline*}
where the first isomorphism is (\ref{equ: Hom induction adjunction}), the second follows from Lemma~\ref{lem: Jacquet basic} and Remark \ref{rem: twist Jacquet basic} (arguing as above (\ref{equ: Hom vers induction block})) and where the last equality follows from the assumption $J_{I_0,I_1,w}(\pi_0^{\infty})_{\cB^{I}_{\Sigma}}=0$. Together with (\ref{equ: OS main seq}) this implies \ref{it: Ext vanishing 2} when $M_1=M^{I_1}(\mu_1)$ for some $\mu_1\in\Lambda_{I_1}^{\dom}$. The result for general $M_1$ again follows from the same induction as in the proof of Lemma~\ref{lem: distinct block}.
\end{proof}

\begin{lem}\label{lem: distinct inf char}
For $i=0,1$ let $V_i= \cF_{P_{I_i}}^G(M_i,\pi_i^{\infty})$ with $I_i\subseteq \Delta$, $M_i$ in $\cO^{\fp_{I_i}}_{\rm{alg}}$ and $\pi_i^{\infty}$ smooth finite length representations of $L_{I_i}$ such that $\pi_i^{\infty}$ is in $\cB^{I_i}_{\Sigma_i}$ with $\Sigma_i= W(L_{I_i})\cdot\cJ(\pi_i^{\infty})\subseteq \widehat{T}^{\infty}$. If $\mathrm{Ext}_{G}^k(V_0,V_1)\neq 0$ for some $k\geq 0$ then there exists $\xi: Z(\fg)\rightarrow E$ such that $M_{0,\xi}\neq 0\neq M_{1,\xi}$.
\end{lem}
\begin{proof}
By d\'evissage we can assume $M_1$ simple of the form $L(\mu)$ for some $\mu\in\Lambda_{I_1}^{\dom}$. Moreover, we may choose $\mu$ so that any $\mu'\in\Lambda_{I_1}^{\dom}$ such that $\mu'\in W(G)\cdot\mu\setminus\{\mu\}$ and $\mu'-\mu\in\Z_{\geq 0}\Phi^+$ must also satisfy for $k\geq 0$
\begin{equation}\label{equ: distinct inf char vanishing}
\mathrm{Ext}_{G}^k(V_0,\cF_{P_{I_1}}^{G}(L(\mu'),\pi_1^{\infty}))=0.
\end{equation}
Define $N^{I_1}(\mu)$ as in the proof of Lemma \ref{lem: distinct block}. As the irreducible constituents $L(\mu')$ of $N^{I_1}(\mu)$ satisfy $\mu'\in\Lambda_{I_1}^{\dom}$, $\mu'\in W(G)\cdot\mu\setminus\{\mu\}$ and $\mu'-\mu\in\Z_{\geq 0}\Phi^+$ (use \cite[Thm.~5.1]{Hum08}), we have $\mathrm{Ext}_{G}^k(V_0,\cF_{P_{I_1}}^{G}(N^{I_1}(\mu),\pi_1^{\infty}))=0$ for $k\geq 0$ by (\ref{equ: distinct inf char vanishing}), and thus the surjection $M^{I_1}(\mu)\twoheadrightarrow L(\mu)$ induces an isomorphism for $k\geq 0$
\[\mathrm{Ext}_{G}^k(V_0,\cF_{P_{I_1}}^{G}(L(\mu),\pi_1^{\infty}))\buildrel\sim\over\longrightarrow \mathrm{Ext}_{G}^k(V_0,\cF_{P_{I_1}}^{G}(M^{I_1}(\mu),\pi_1^{\infty}))\]
(which is non-zero for some $k\geq 0$ by assumption). By the proof of \ref{it: OS spectral seq 2} of Corollary~\ref{prop: OS spectral seq} we have $\mathrm{Ext}_{U(\fg)}^{\ell}(M^{I_1}(\mu),M_0^w)\neq 0$ for some $w\in W^{I_0,I_1}$ and some $\ell\geq 0$, which together with \ref{it: block 3} of Lemma~\ref{lem: central component} implies $(M_0^w)_{\xi}\neq 0$ for the unique infinitesimal character $\xi: Z(\fg)\rightarrow E$ such that $M^{I_1}(\mu)_{\xi}\neq 0\neq L(\mu)_{\xi}$. As the adjoint action of $G$, and in particular of $w\in W(G)$, on $Z(\fg)$ (inside $U(\fg)$) is the identity, this is equivalent to $(M_{0,\xi})^w\neq 0$, i.e.~to $M_{0,\xi}\neq 0$.
\end{proof}

\begin{rem}\label{rem: OS isotypic}
Let $V$ be a finite length object in $\mathrm{Rep}^{\rm{an}}_{\rm{adm}}(G)$ with each constituent of the form $\cF_{P_I}^{G}(M,\pi^{\infty})$ for some $I\subseteq \Delta$, some $M$ in $\cO^{\fp_I}_{\rm{alg}}$ and some smooth (finite length) representation $\pi^{\infty}$ of $L_I$ in $\cB^I_{W(L_I)\cdot\cJ(\pi^{\infty})}$. Given $\xi: Z(\fg)\rightarrow E$ and a $W(G)$-coset $\Sigma\subseteq \widehat{T}^{\infty}$, there exists a maximal closed subrepresentation $V_{\xi,\Sigma}\subseteq V$ such that each constituent of $V_{\xi,\Sigma}$ has the form $\cF_{P_I}^{G}(M,\pi^{\infty})$ for $I$, $M$, $\pi^{\infty}$ with $M=M_{\xi}$ and $\cJ(\pi^{\infty})\subseteq \Sigma$. It follows from Lemma~\ref{lem: distinct block} and Lemma~\ref{lem: distinct inf char} that $V$ is the direct sum of $V_{\xi,\Sigma}$ over all pairs $\xi,\Sigma$. Similarly, we can define $V_{\xi}$ just taking $\xi$ into account, as well as $V_{\Sigma}$. Taking continuous duals we can define $D_{\xi,\Sigma}$, $D_{\xi}$ and $D_{\Sigma}$ for a finite length coadmissible $D(G)$-module with constituents of the form $\cF_{P_I}^{G}(M,\pi^{\infty})^\vee$ (with $M,\pi^{\infty}$ as above).
\end{rem}

\begin{rem}
Let $I_0\subseteq \Delta$, $x_0\in W(G)$ with $I_0=\Delta\setminus D_L(x_0)$, $M_0\defeq L(x_0)\in\cO_{\rm{alg}}^{\fp_{I_0}}$ (Lemma \ref{lem: dominance and left set}), $\pi_0^{\infty}$ a smooth strongly admissible representation of $L_{I_0}$ over $E$ and $V_0\defeq \cF_{P_{I_0}}^G(M_0,\pi_0^{\infty})$. We have $H^0(\fu,M_0^w)=0$ for $1\neq w\in W^{I_0,\emptyset}$ by \ref{it: H0 conjugate 2} of Lemma~\ref{lem: H0 Weyl conjugate} (applied with $I'=I_0$ and $I=\emptyset$) and $H^0(\fu,M_0)$ is finite dimensional. Since $H^0(\fu,\cM_0^w)\subseteq \cM_0^w$ is a closed Fr\'echet $U(\ft)$-submodule, it is small by Lemma~\ref{lem: t closed}, and thus contains $H^0(\fu,M_0^w)=H^0(\fu,\cM_0^w)\cap M_0^w$ as a dense subspace. In particular $H^0(\fu,\cM_0^w)=0$ for $1\neq w\in W^{I_0,\emptyset}$ and the injection $M_0\hookrightarrow \cM_0$ induces a (topological) isomorphism of finite dimensional $E$-vector spaces $H^0(\fu,M_0)\buildrel\sim\over\rightarrow H^0(\fu,\cM_0)$.
We fix a compact open subgroup $G_0\subseteq G$ and for $w\in W^{I_0,\emptyset}$, we write $P_w\defeq w^{-1}P_{I_0}w$, $P_{w,0}\defeq P_w\cap G_0$, and $\cX_w$ for the set of compact open subsets of $(BP_w)\cap G_0$ stable under right multiplication by $P_{w,0}$. We write $U_w$ for a general element of $\cX_w$. For $r\in\cI$, it follows from Lemma~\ref{lem: family of standard semi norms} that $(\cM_0^w)_r=D(G_0)_{P_{w,0},r}\otimes_{D(G_0)_{P_{w,0}}}\cM_0^w$ is the completion of $\cM_0^w$ under a standard semi-norm, which by Remark~\ref{rem: dense in Banach} implies that $H^0(\fu,M_0^w)=H^0(\fu,(\cM_0^w)_r)\cap M_0^w$ is dense inside the closed Banach subspace $H^0(\fu,(\cM_0^w)_r)\subseteq (\cM_0^w)_r$.
Consequently, for $1\neq w\in W^{I_0,\emptyset}$, we have for $g\in B$ and $r\in\cI$ (since $H^0(\fu,M_0^w)=0$)
\[H^0(\fu,((\cM_0^w)_r)^{g^{-1}})\cong H^0(g\fu g^{-1},((\cM_0^w)_r)^{g^{-1}})\cong H^0(\fu,(\cM_0^w)_r)=0,\]
and therefore
\[\begin{array}{lll}
H^0(\fu, \mathrm{gr}_w(V_0^\vee))&\cong&
\big(\varprojlim_{U_w\in\cX_w,r\in\cI}(D(G_0)_{U_w,r}\otimes_{D(G_0)_{P_{w,0},r}}D_r)\big)[\fu]\\
&\cong&
\varprojlim_{U_w\in\cX_w,r\in\cI}\big((D(G_0)_{U_w,r}\otimes_{D(G_0)_{P_{w,0},r}}D_r)[\fu]\big)\\
&\cong& \varprojlim_{U_w\in\cX_w,r\in\cI}\big(\bigoplus_{g\in U_wG_1^m/P_{w,0}G_1^m}\big(((\cM_0^w)_r)^{g^{-1}}[\fu]\otimes_E D_r^{\infty,g^{-1}}\big)\big)\\
&=&0
\end{array}\]
where the first isomorphism uses (\ref{equ: graded short term}) and the third isomorphism uses (\ref{equ: tensor as sum}), and where $D_r$ and $D_r^{\infty}$ are the same as in \emph{loc.~cit.}~and depend on the choice of $w$ (see (\ref{varprojlimr}), (\ref{varprojlimr sm}) and the paragraph before (\ref{equ: isom modulo coset})).
For $w=1$, we have
\[\begin{array}{lllll}
H^0(\fu, \mathrm{gr}_1(V_0^\vee))&\cong &\!\big(\varprojlim_{r\in\cI} ((\cM_0)_r\otimes_E (\pi_0^{\infty})^\vee_r)\big)[\fu]&\cong &\varprojlim_{r\in\cI}\big(((\cM_0)_r\otimes_E (\pi_0^{\infty})^\vee_r)[\fu]\big)\\
&\cong &\varprojlim_{r\in\cI}( (\cM_0)_r[\fu]\otimes_E (\pi_0^{\infty})^\vee_r)
&\cong &\varprojlim_{r\in\cI} ((\cM_0)_r[\fu])\widehat{\otimes}_E \varprojlim_{r\in\cI}(\pi_0^{\infty})^\vee_r)\\
&\cong &\cM_0[\fu] \widehat\otimes_E(\pi_0^{\infty})^\vee&& \\
&\cong &\cM_0[\fu]\otimes_E(\pi_0^{\infty})^\vee&&
\end{array}\]
where the first isomorphism follows from Proposition \ref{lem: graded OS tensor} (applied with $I_1=\emptyset$) and from (\ref{varprojlimr}), the fourth follows from Lemma \ref{limprojtensor}, and the last follows from the finite di\-mensionality of $\cM_0[\fu]$. We deduce from all this $H^0(U, \mathrm{gr}_w(V_0^\vee))=0$ for $1\neq w\in W^{I_0,\emptyset}$ and $H^0(U, \mathrm{gr}_1(V_0^\vee))\cong H^0(\fu, M_0)\otimes_E (J_{I_0,\emptyset}(\pi_0^{\infty}))^\vee$ (cf.~(\ref{equ: sm coinv dual}) with $I_1=\emptyset$ and $w=1$). By d\'evissage on the $W^{I_0,\emptyset}$-filtration $\mathrm{Fil}_{\bullet}(V_0^\vee)$ we finally obtain the following isomorphism of $D(T)$-modules (which will be used in the proof of Theorem \ref{prop: top deg} below)
\begin{equation}\label{equ: N coh H0}
H^0(U, V_0^\vee)\cong H^0(\fu, M_0)\otimes_E (J_{I_0,\emptyset}(\pi_0^{\infty}))^\vee.
\end{equation}
Note that (\ref{equ: N coh H0}) is a finite dimensional (coadmissible) $D(T)$-module which is non-zero when $\pi_0^{\infty}$ is $G$-basic (and $M_0=L(x_0)$).
\end{rem}

\newpage

\section{Complexes of locally analytic representa\-tions}\label{sec: extension}

We describe the global sections of the de Rham complex of the Drinfeld space in dimension $n-1$ as a complex of explicit finite length coadmissible $D(G)$-modules, and we describe an explicit quasi-isomorphism with the direct sum of its (shifted) cohomology groups. The proof works essentially \emph{verbatim} for the complex of holomorphic discrete series of \cite{S92}.

\subsection{Results on locally analytic \texorpdfstring{$\mathrm{Ext}^0$}{Ext0}, \texorpdfstring{$\mathrm{Ext}^1$}{Ext1} and \texorpdfstring{$\mathrm{Ext}^2$}{Ext2} groups}\label{subsec: Ext OS}

We prove various useful results on $\mathrm{Hom}$, $\mathrm{Ext}^1$ and $\mathrm{Ext}^2$ groups between certain Orlik-Strauch representations.\bigskip

We use the notation of \S\ref{sec: smooth rep} and \S\ref{sec: n coh}, in particular $\mu_0\in \Lambda^{\dom}$ is a fixed weight, $L(w)=L(w\cdot \mu_0)$ for $w\in W(G)$, $N^{I}(w)= \Ker(M^{I}(w)\twoheadrightarrow L(w))$, etc. For $x\in W(G)$ we write $I_x\defeq \Delta\setminus D_L(x)$. Throughout this section for $i=0,1$, $I_i$ is a subset of $\Delta$, $M_i$ is a $U(\fg)$-module in $\cO^{\fp_{I_i}}_{\rm{alg}}$ (see below (\ref{equ: simple self dual})), $\Sigma_i$ is a finite subset of $\widehat{T}^{\infty}$ preserved under the left action (\ref{equ: smooth dot action}) of $W(L_{I_i})$ and $\pi_i^{\infty}$ is a smooth (finite length) representation of $L_{I_i}$ in the category $\cB^{I_i}_{\Sigma_i}$ (see above (\ref{equ: sm block decomposition})). We write $V_i\defeq \cF_{P_{I_i}}^{G}(M_i,\pi_i^{\infty})$ for $i=0,1$. Depending on the statements, we will add various assumptions on $I_i$, $M_i$, $\Sigma_i$ or $\pi_i^\infty$ (and thus $V_i$).

\begin{lem}\label{lem: Hom OS}
For $i=0,1$ assume that $I_i=\Delta\setminus D_L(x_i)$ and $M_i=L(x_i)$ for some $x_i\in W(G)$. Then we have $\Hom_{G}(V_0,V_1)\neq 0$ only if $x_0=x_1$ (and $I_0=I_1$), in which case we have a canonical isomorphism
\begin{equation}\label{equ: Hom OS}
\Hom_{G}(V_0,V_1) \cong \Hom_{U(\fg)}(M_1,M_0)\otimes_E \Hom_{L_{I_1}}(\pi_0^{\infty}, \pi_1^{\infty}).
\end{equation}
In particular, if $V_0$ and $V_1$ are irreducible (which forces $\pi_i^{\infty}$ to be irreducible for $i=0,1$), they are isomorphic if and only if $x_0=x_1$ and $\pi_0^{\infty}\cong \pi_1^{\infty}$.
\end{lem}
\begin{proof}
Let $V_2\defeq \cF_{P_{I_1}}^{G}(M^{I_1}(x_1),\pi_1^{\infty})$, the surjection $M^{I_1}(x_1)\twoheadrightarrow L(x_1)$ induces an injection $V_1\hookrightarrow V_2$ and thus an injection
\begin{equation}\label{equ: Hom OS embedding}
\Hom_{G}(V_0,V_1)\hookrightarrow \Hom_{G}(V_0,V_2).
\end{equation}
Recall from \ref{it: H0 conjugate 3} of Lemma~\ref{lem: H0 Weyl conjugate} that $L(x_0)^w$ is in $\cO^{\fb}_{\rm{alg}}$ for some $w\in W^{I_0,I_1}$ if and only if $w=1$. Hence for $1\neq w\in W^{I_0,I_1}$ we have (using Shapiro's lemma for the first isomorphism)
\[\Hom_{U(\fg)}(M^{I_1}(x_1),L(x_0)^w)\cong \Hom_{U(\fp_{I_1})}(L^{I_1}(x_1),L(x_0)^w)=0.\]
By (\ref{isov1v0}) (applied with $\mathrm{gr}_w(V_0^\vee)$ instead of $V_0^\vee$ there, which doesn't change the argument), (\ref{equ: ST seq}) (applied with $D=\mathrm{gr}_w(V_0^\vee)$ and $k=\ell=0$) and Corollary~\ref{cor: Ext P graded} (for $k=\ell=0$), we deduce $\Hom_{D(G)}(V_2^\vee,\mathrm{gr}_w(V_0^\vee))=0$ for $1\neq w\in W^{I_0,I_1}$ and for $w=1$:
\[\Hom_{D(G)}(V_2^\vee,\mathrm{gr}_1(V_0^\vee))\cong\Hom_{U(\fg)}(M^{I_1}(x_1),L(x_0))\otimes_E \Hom_{L_{I_1}}(i_{I_0,I_1,1}^{\infty}(J_{I_0,I_1,1}(\pi_0^{\infty})),\pi_1^{\infty}).\]
By a d\'evissage with respect to $(\mathrm{Fil}_w(V_0^\vee))_{w\in W^{I_0,I_1}}$ we obtain
\begin{multline}\label{equ: Hom OS seq}
\Hom_{G}(V_0,V_2)\cong \Hom_{D(G)}(V_2^\vee,V_0^\vee)\\
\cong\Hom_{U(\fg)}(M^{I_1}(x_1),L(x_0))\otimes_E \Hom_{L_{I_1}}(i_{I_0,I_1,1}^{\infty}(J_{I_0,I_1,1}(\pi_0^{\infty})),\pi_1^{\infty}).
\end{multline}
If $\Hom_{G}(V_0,V_1)\neq 0$ then $\Hom_{G}(V_0,V_2)\neq 0$ by (\ref{equ: Hom OS embedding}), thus $\Hom_{U(\fg)}(M^{I_1}(x_1),L(x_0))\neq 0$ by (\ref{equ: Hom OS seq}) and thus
\begin{equation}\label{isosocle}
\Hom_{U(\fg)}(L(x_1),L(x_0))\buildrel\sim\over\longrightarrow \Hom_{U(\fg)}(M^{I_1}(x_1),L(x_0))\neq 0
\end{equation}
(as $M^{I_1}(x_1)$ has irreducible cosocle $L(x_1)$) which forces $x_0=x_1$ and hence $I_0=I_1$. As any irreducible constituent $L(x')$ of $N^{I_1}(x_1)$ satisfies $x'>x_1=x_0$, the argument above with $M^{I_1}(x')$ and $L(x')$ instead of $M^{I_1}(x_1)$ and $L(x_1)$ (note that $I_1=\Delta\setminus D_L(x')$ by Lemma \ref{lem: dominance and left set}) shows that $\Hom_{G}(V_0,\cF_{P_{I_1}}^{G}(L(x'),\pi_1^{\infty}))=0$. As $V_2/V_1\cong \cF_{P_{I_1}}^{G}(N^{I_1}(x_1),\pi_1^{\infty})$ (see \ref{it: OS property 1} of Theorem \ref{prop: OS property}), we deduce by d\'evissage $\Hom_{G}(V_0,V_2/V_1)=0$ and hence $\Hom_{G}(V_0,V_1)\buildrel\sim\over\rightarrow \Hom_{G}(V_0,V_2)$. Then (\ref{equ: Hom OS}) follows from (\ref{equ: Hom OS seq}) (with $I_0=I_1$) and (\ref{isosocle}).
\end{proof}

\begin{lem}\label{lem: sm to OS}
Let $x\in W(G)$, $I\defeq\Delta\setminus D_L(x)$, $\Sigma$ a finite subset of $\widehat{T}^{\infty}$ preserved under the left action (\ref{equ: smooth dot action}) of $W(L_{I})$, $\pi^{\infty}$ a smooth (finite length) multiplicity free representation of $L_{I}$ in $\cB^I_{\Sigma}$ and $V\defeq \cF_{P_I}^{G}(L(x),\pi^{\infty})$. Then $V$ is of finite length and multiplicity free, and the functor $\cF_{P_I}^{G}(L(x),-)$ induces a bijection $\mathrm{JH}_{L_I}(\pi^{\infty})\buildrel\sim\over\rightarrow \mathrm{JH}_{G}(V)$ between partially ordered sets (see \S\ref{generalnotation} for the definition of this partial order).
\end{lem}
\begin{proof}
It follows from Theorem~\ref{prop: OS property} and Lemma~\ref{lem: Hom OS} that $V$ is finite length and multiplicity free
with $\mathrm{JH}_{G}(V)=\{\cF_{P_I}^{G}(L(x),\tau^{\infty})\mid \tau^{\infty}\in\mathrm{JH}_{L_I}(\pi^{\infty})\}$. Let $\tau_0^{\infty},\tau_1^{\infty}$ be two distinct irreducible constituents of $\pi^{\infty}$ and define $\sigma^{\infty}$ as the unique subrepresentation of $\pi^{\infty}$ with cosocle $\tau_0^{\infty}$. Recall that, by definition of the partial order on $\mathrm{JH}_{L_I}(\pi^{\infty})$, $\tau_1^{\infty}<\tau_0^{\infty}$ if and only if $\tau_1^{\infty}\in \mathrm{JH}_{L_I}(\sigma^{\infty})$. Let $W\defeq \cF_{P_I}^{G}(L(x),\sigma^{\infty})$ and $W_i\defeq \cF_{P_I}^{G}(L(x),\tau_i^{\infty})$ for $i=0,1$. By Theorem~\ref{prop: OS property} and the last claim in Lemma~\ref{lem: Hom OS} $\tau_1^{\infty}\in \mathrm{JH}_{L_I}(\sigma^{\infty})$ if and only if $W_1\in\mathrm{JH}_{G}(W)$.
By (\ref{equ: Hom OS}) applied with $x_0=x_1=x$, $\pi_0^{\infty}=\sigma^{\infty}$ and $\pi_1^{\infty}$ any constituent of $\sigma^{\infty}$, we deduce that $W$ has simple cosocle $W_0$, and thus $W_1<W_0$ if and only if $W_1\in \mathrm{JH}_{G}(W)$. We have shown $\tau_1^{\infty}<\tau_0^{\infty}$ if and only if $W_1<W_0$, which proves the lemma.
\end{proof}

Recall that $G$-regular and $G$-basic smooth representations of $L_I$ over $E$ are defined in \ref{it: basic 1} of Definition \ref{def: basic rep}. Recall also that a $G$-basic representation of $L_I$ is in $\cB^I_{\Sigma}$ for $\Sigma=W(L_I)\cdot\cJ(\pi^{\infty})$ (see below Lemma \ref{lem: Jacquet of PS}) and is multiplicity free (\ref{it: PS 1} of \emph{loc.~cit.}), in particular we can apply Lemma \ref{lem: sm to OS} when $\pi^{\infty}$ is $G$-basic.

\begin{lem}\label{lem: Ext1 from sm}
For $i=0,1$ assume that $I\defeq I_0=I_1=\Delta\setminus D_L(x)$ and $M_0=M_1=L(x)$ for some $x\in W(G)$. Assume moreover that $\pi_0^{\infty}$ and $\pi_1^{\infty}$ have no common Jordan-H\"older factor. Then the functor $\cF_{P_I}^{G}(L(x),-)$ induces a canonical isomorphism
\begin{equation}\label{equ: Ext1 from sm}
\mathrm{Ext}_{L_I}^1(\pi_0^{\infty},\pi_1^{\infty})^\infty\buildrel\sim\over\longrightarrow \mathrm{Ext}_{G}^1(V_0,V_1).
\end{equation}
If moreover $\pi_0^{\infty}$ and $\pi_1^{\infty}$ are $G$-basic, then (\ref{equ: Ext1 from sm}) is one dimensional if non-zero.
\end{lem}
\begin{proof}
Recall first that both spaces in (\ref{equ: Ext1 from sm}) are finite dimensional over $E$, as follows from Theorem \ref{thm: finitedim} and (the references in) its proof.\bigskip

\textbf{Step $1$}: We prove that (\ref{equ: Ext1 from sm}) is injective.\\
If $\mathrm{Ext}_{L_I}^1(\pi_0^{\infty},\pi_1^{\infty})^\infty=0$, there is nothing to prove. Otherwise, let $\pi^{\infty}$ in $\mathrm{Rep}^{\infty}_{\rm{adm}}(L_I)$ that fits into a non-split extension
\begin{equation}\label{equ: Ext1 from sm seq}
0\rightarrow \pi_1^{\infty} \rightarrow \pi^{\infty} \rightarrow \pi_0^{\infty} \rightarrow 0.
\end{equation}
and note that $\pi^{\infty} $ is in $\cB^I_{\Sigma}$ for $\Sigma=\Sigma_0\cup \Sigma_1$ and is multiplicity free. In particular, there exists $\sigma_i^{\infty}\in\mathrm{JH}_{L_I}(\pi_i^{\infty})$ such that $\pi^{\infty}$ admits a length $2$ subquotient $\sigma^{\infty}$ with socle $\sigma_1^{\infty}$ and cosocle $\sigma_0^{\infty}$. Applying $\cF_{P_I}^{G}(L(x),-)$ to (\ref{equ: Ext1 from sm seq}) and using Lemma~\ref{lem: sm to OS}, it follows that $V\defeq \cF_{P_I}^{G}(L(x),\pi^{\infty})$ admits a length $2$ subquotient $\cF_{P_I}^{G}(L(x),\sigma^{\infty})$ with socle $\cF_{P_I}^{G}(L(x),\sigma_1^{\infty})$ and cosocle $\cF_{P_I}^{G}(L(x),\sigma_0^{\infty})$. In particular the short exact sequence $0\rightarrow V_1\rightarrow V\rightarrow V_0\rightarrow 0$ is non-split. This proves the injectivity of (\ref{equ: Ext1 from sm}).\bigskip

\textbf{Step $2$}: We prove $\Dim_E \mathrm{Ext}_{G}^1(V_0,V_1)\leq \Dim_E \mathrm{Ext}_{L_I}^1(\pi_0^{\infty},\pi_1^{\infty})^\infty$.\\
We can assume $\mathrm{Ext}_{G}^1(V_0,V_1)\neq 0$. Let $V_2\defeq \cF_{P_I}^{G}(M^I(x),\pi_1^{\infty})$, we have $V_1\hookrightarrow V_2$ with $V_2/V_1\cong \cF_{P_I}^{G}(N^I(x),\pi_1^{\infty})$. Since $L(x)\notin \mathrm{JH}_{U(\fg)}(N^I(x))$, by Lemma~\ref{lem: Hom OS} (using Lemma \ref{lem: dominance and left set}) we know that $V_2/V_1$ and $V_0$ share no common Jordan-H\"older factor and in particular $\Hom_{G}(V_0,V_2/V_1)=0$. A d\'evissage on $0\rightarrow V_1\rightarrow V_2\rightarrow V_2/V_1\rightarrow 0$ then gives an injection
\begin{equation}\label{equ: Ext1 from sm embedding}
\mathrm{Ext}_{G}^1(V_0,V_1)\hookrightarrow \mathrm{Ext}_{G}^1(V_0,V_2)
\end{equation}
and thus $\mathrm{Ext}_{G}^1(V_0,V_2)\neq 0$ since $\mathrm{Ext}_{G}^1(V_0,V_1)\neq 0$. By a d\'evissage on $(\mathrm{Fil}_w(V_0^\vee))_{w\in W^{I,I}}$ there exists $w\in W^{I,I}$ such that $\mathrm{Ext}_{D(G)}^1(V_2^\vee,\mathrm{gr}_w(V_0^\vee))\neq 0$, which together with (\ref{isov1v0}) (applied with $\mathrm{gr}_w(V_0^\vee)$ instead of $V_0^\vee$), (\ref{equ: ST seq}) (applied with $D=\mathrm{gr}_w(V_0^\vee)$) and Corollary~\ref{cor: Ext P graded} implies
\begin{equation}\label{equ: Ext1 from sm spectral seq}
\mathrm{Ext}_{L_I}^k(i_{I,I,w}^{\infty}(J_{I,I,w}(\pi_0^{\infty})), \pi_1^{\infty})^{\infty}\otimes_E \mathrm{Ext}_{U(\fg)}^{\ell}(M^I(x), L(x)^w)\neq 0
\end{equation}
for some $k,\ell\geq 0$ such that $k+\ell=1$. By \ref{it: H0 conjugate 3} of Lemma~\ref{lem: H0 Weyl conjugate} (for $\ell=0$) and \ref{rem: H1 conj Ext1 vanishing} of Remark~\ref{mergerem} (for $\ell =1$) we have $\mathrm{Ext}_{U(\fg)}^{\ell}(M^I(x), L(x)^{w'})=0$ for $\ell\leq 1$ and $1\neq w'\in W^{I,I}$. Hence, (\ref{equ: Ext1 from sm spectral seq}) can hold only when $w=1$. Since $\pi_0^{\infty}$ and $\pi_1^{\infty}$ have no common Jordan-H\"older factor, we have $\Hom_{L_I}(\pi_0^{\infty},\pi_1^{\infty})=0$, so that the only possibly non-zero term in (\ref{equ: Ext1 from sm spectral seq}) is for $k=1$ and $\ell=0$. The spectral sequence (\ref{equ: ST seq}) (with Corollary~\ref{cor: Ext P graded}) and (\ref{equ: Ext1 from sm embedding}) then yield
\begin{multline}\label{equ: Ext1 from sm upper bound seq}
\Dim_E \mathrm{Ext}_{G}^1(V_0,V_1)\leq \Dim_E \mathrm{Ext}_{G}^1(V_0,V_2)\\
\leq \Dim_E \mathrm{Ext}_{L_I}^1(\pi_0^{\infty}, \pi_1^{\infty})^{\infty} \Dim_E \Hom_{U(\fg)}(M^I(x), L(x)).
\end{multline}
As $\Hom_{U(\fg)}(M^I(x), L(x))$ is one dimensional, (\ref{equ: Ext1 from sm upper bound seq}) implies the statement.\bigskip

{Step $1$} and {Step $2$} imply that (\ref{equ: Ext1 from sm}) is an isomorphism. When $\pi_0^{\infty}$ and $\pi_1^{\infty}$ are moreover $G$-basic and (\ref{equ: Ext1 from sm}) is non-zero, then $d_I(\pi_0^{\infty},\pi_1^{\infty})= 1$ (see \ref{it: basic 3} of Definition \ref{def: basic rep}) and Lemma~\ref{lem: dim one} then implies that (\ref{equ: Ext1 from sm}) is one dimensional.
\end{proof}

\begin{lem}\label{lem: OS basic unique}
Let $I\subseteq \Delta$, $L(x)$ in $\cO^{\fp_I}_{\rm{alg}}$, $\pi^{\infty}$ a smooth $G$-basic representation of $L_I$ and $V\defeq \cF_{P_I}^{G}(L(x),\pi^{\infty})$. Let $V'$ be a smooth multiplicity free finite length representation of $G$ such that $\mathrm{JH}_{G}(V')=\mathrm{JH}_{G}(V)$ as partially ordered sets. Then we have $V\cong V'$.
\end{lem}
\begin{proof}
Using \ref{it: OS property 2} of Theorem \ref{prop: OS property} we can assume $I=\Delta\setminus D_L(x)$ and by \ref{it: PS 2} of Lemma \ref{lem: Jacquet of PS} $\pi^{\infty}$ is (finite length) multiplicity free. By Lemma~\ref{lem: sm to OS} $V$ is multiplicity free and the functor $\cF_{P_I}^{G}(L(x),-)$ induces a bijection of partially ordered sets $\mathrm{JH}_{L_I}(\pi^{\infty})\buildrel\sim\over\longrightarrow \mathrm{JH}_{G}(V)$. We prove the statement by induction on the length of $\pi^{\infty}$. If $\pi^{\infty}$ is irreducible, then $V$ is irreducible by \ref{it: OS property 3} of Theorem~\ref{prop: OS property}, so $\mathrm{JH}_{G}(V)=\mathrm{JH}_{G}(V')$ forces $V\cong V'$ and there is nothing to prove. If $\pi^{\infty}$ is reducible, by Corollary~\ref{cor: basic reducible} there exist smooth $G$-basic representations $\pi_0^{\infty},\pi_1^{\infty}$ of $L_I$ and a non-split short exact sequence $0\rightarrow \pi_1^{\infty}\rightarrow \pi^{\infty}\rightarrow \pi_0^{\infty}\rightarrow 0$. Applying $\cF_{P_I}^{G}(L(x),-)$ we obtain a short exact sequence $0\rightarrow V_1\rightarrow V \rightarrow V_0\rightarrow 0$ with $V_i\defeq \cF_{P_I}^{G}(L(x),\pi_i^{\infty})$, which is non-split by Step 1 of the proof of Lemma~\ref{lem: Ext1 from sm}. Since $\mathrm{JH}_{G}(V)=\mathrm{JH}_{G}(V')$ as partially ordered sets, $V'$ also fits into a non-split short exact sequence $0\rightarrow V_1'\rightarrow V' \rightarrow V_0'\rightarrow 0$ with $\mathrm{JH}_{G}(V_i)=\mathrm{JH}_{G}(V_i')$ as partially ordered sets for $i=0,1$. As $\pi_i^{\infty}$ has strictly smaller length than $\pi^{\infty}$, by induction we have $V_i\cong V_i'$ for $i=0,1$. As both $V$ and $V'$ fit into non-split short exact sequences $0\rightarrow V_1\rightarrow \ast\rightarrow V_0\rightarrow 0$ and $\Dim_E \mathrm{Ext}_{G}^1(V_0,V_1)=1$ by Lemma~\ref{lem: Ext1 from sm}, we deduce $V\cong V'$.
\end{proof}

\begin{lem}\label{lem: Ext1 g vers OS}
For $i=0,1$ assume that $I_i=\Delta\setminus D_L(x_i)$ and $M_i=L(x_i)$ for some $x_i\in W(G)$. Assume moreover that $x_0\neq x_1$, that $\pi_0^{\infty}$ and $\pi_1^{\infty}$ are $G$-basic and that
\begin{equation}\label{equ: Ext1 g vers OS sm}
\Hom_{L_{I_1}}(i_{I,I_1}^{\infty}(J_{I_0,I}(\pi_0^{\infty})),\pi_1^{\infty})\neq 0
\end{equation}
where $I\defeq I_0\cap I_1$. Then there exists a unique $G$-basic representation $\pi^{\infty}$ of $L_I$ which is both a subrepresentation of $J_{I_1,I}'(\pi_1^{\infty})$(see (\ref{extranot})) and a quotient of $J_{I_0,I}(\pi_0^{\infty})$.
Moreover, we have a canonical injection
\begin{equation}\label{equ: Ext1 g vers OS}
\Hom_{L_{I_0}}(\pi_0^{\infty},i_{I,I_0}^{\infty}(\pi^{\infty}))\otimes_E\mathrm{Ext}_{U(\fg)}^1(M_1,M_0)\otimes_E\Hom_{L_{I_1}}(i_{I,I_1}^{\infty}(\pi^{\infty}),\pi_1^{\infty})\hookrightarrow \mathrm{Ext}_{G}^1(V_0,V_1).
\end{equation}
\end{lem}
\begin{proof}
\textbf{Step $1$}: We construct the desired $G$-basic $\pi^{\infty}$.\\
By (\ref{equ: second adjunction}), (\ref{equ: Ext1 g vers OS sm}) is equivalent to
\begin{equation}\label{equ: Ext1 g vers OS sm prime}
\Hom_{L_I}(J_{I_0,I}(\pi_0^{\infty}),J_{I_1,I}'(\pi_1^{\infty}))\neq 0.
\end{equation}
For $i=0,1$ let $\Sigma_i\defeq W(L_{I_i})\cdot\cJ(\pi_i^{\infty})$. By (\ref{equ: twist Jacquet support}) and \ref{it: PS 1} of Lemma~\ref{lem: Jacquet of PS} $W(L_{I_1})\cdot\cJ(J_{I_1,I}'(\pi_1^{\infty}))\!=\Sigma_1$ is a single left $W(L_{I_1})$-coset and $\Sigma_0$ is a single left $W(L_{I_0})$-coset, and thus (using the regularity of the characters, see \ref{it: PS 1} of Lemma \ref{lem: Jacquet of PS})
\begin{equation}\label{equ: Ext1 g coset}
\Sigma\defeq \Sigma_0\cap W(L_{I_1})\cdot\cJ(J_{I_1,I}'(\pi_1^{\infty}))=\Sigma_0\cap \Sigma_1
\end{equation}
(which is non-empty using (\ref{equ: Ext1 g vers OS sm prime})) is a single left $W(L_I)$-coset. In particular, (\ref{equ: Ext1 g vers OS sm prime}) is equivalent to
\begin{equation*}
\Hom_{L_I}\big(J_{I_0,I}(\pi_0^{\infty})_{\cB^I_{\Sigma}},J_{I_1,I}'(\pi_1^{\infty})_{\cB^I_{\Sigma}}\big)\neq 0
\end{equation*}
with both $J_{I_0,I}(\pi_0^{\infty})_{\cB^I_{\Sigma}}$ and $J_{I_1,I}'(\pi_1^{\infty})_{\cB^I_{\Sigma}}$ non-zero and thus $G$-basic by Lemma~\ref{lem: Jacquet basic} and Remark~\ref{rem: twist Jacquet basic}. The desired $\pi^{\infty}$ is necessarily the image of a non-zero map $J_{I_0,I}(\pi_0^{\infty})\rightarrow J_{I_1,I}'(\pi_1^{\infty})$, which itself is necessarily the image of the unique (up to scalar) non-zero map $J_{I_0,I}(\pi_0^{\infty})_{\cB^I_{\Sigma}}\rightarrow J_{I_1,I}'(\pi_1^{\infty})_{\cB^I_{\Sigma}}$ (the unicity follows from the fact $G$-basic representations are multiplicity free with simple socle (and cosocle), see \ref{it: basic as image} of Remark \ref{rem: basic PS intertwine}). In particular $\pi^{\infty}$ is a quotient of $J_{I_0,I}(\pi_0^{\infty})_{\cB^I_{\Sigma}}$ with simple socle and cosocle and thus is also $G$-basic by Corollary~\ref{cor: basic subquotient}. The definition of $\pi^{\infty}$ together with (\ref{equ: first adjunction}) and (\ref{equ: second adjunction}) implies
\begin{equation}\label{equ: Ext1 g vers OS adjunction}
\Hom_{L_{I_0}}(\pi_0^{\infty}, i_{I,I_0}^{\infty}(\pi^{\infty}))\neq 0 \ \ \mathrm{and}\ \ \Hom_{L_{I_1}}(i_{I,I_1}^{\infty}(\pi^{\infty}),\pi_1^{\infty}) \neq 0
\end{equation}
with both spaces being one dimensional (using \ref{it: basic as image} of Remark \ref{rem: basic PS intertwine} as above).\bigskip

\textbf{Step $2$}: We construct the map (\ref{equ: Ext1 g vers OS}).\\
The exact functor $\cF_{P_I}^{G}(-,\pi^{\infty})$ induces a canonical map
\begin{equation}\label{equ: Ext1 g map}
\mathrm{Ext}_{U(\fg)}^1(M_1,M_0)\longrightarrow \mathrm{Ext}_{G}^1(\cF_{P_I}^G(M_0,\pi^{\infty}),\cF_{P_I}^G(M_1,\pi^{\infty})).
\end{equation}
From \ref{it: OS property 2} of Theorem~\ref{prop: OS property} and Lemma~\ref{lem: Hom OS} we have canonical isomorphisms
\[\Hom_{L_{I_0}}(\pi_0^{\infty},i_{I,I_0}^{\infty}(\pi^{\infty}))\cong \Hom_{G}(V_0, \cF_{P_I}^G(M_0,\pi^{\infty}))\]
and
\[\Hom_{L_{I_1}}(i_{I,I_1}^{\infty}(\pi^{\infty}),\pi_1^{\infty})\cong \Hom_{G}(\cF_{P_I}^G(M_1,\pi^{\infty}), V_1),\]
which together with (\ref{equ: Ext1 g map}) and obvious functorialities give a canonical map (\ref{equ: Ext1 g vers OS}). In the rest of the proof we assume that the left hand side of (\ref{equ: Ext1 g vers OS}) is non-zero, equivalently $\mathrm{Ext}_{U(\fg)}^1(M_1,M_0)\ne 0$ by the end of Step $1$ (otherwise (\ref{equ: Ext1 g vers OS}) is obviously an injection).\bigskip

\textbf{Step $3$}: We reduce to the case $\pi_0^{\infty}$ is a subobject of $i_{I,I_0}^{\infty}\!(\pi^{\infty})$ and $\pi_1^{\infty}$ a quotient~of~$i_{I,I_1}^{\infty}\!(\pi^{\infty})$.\\
Let $\pi_{0,-}^{\infty}$ (resp.~$\pi_{1,-}^{\infty}$) be the image of the unique (up to scalar) non-zero map $\pi_0^{\infty}\rightarrow i_{I,I_0}^{\infty}(\pi^{\infty})$ (resp.~$i_{I,I_1}^{\infty}(\pi^{\infty})\rightarrow \pi_1^{\infty}$) which is $G$-basic by Corollary~\ref{cor: basic subquotient}. As $\pi_0^{\infty}$ and $\pi_1^{\infty}$ are multiplicity free and $x_0\neq x_1$, by the last statement in Lemma \ref{lem: Hom OS} $V_0$ and $V_1$ are multiplicity free with no common Jordan-H\"older factor. In particular
\[\Hom_{G}\big(V_0,V_1/\cF_{P_{I_1}}^{G}(M_1,\pi_1^{\infty}/\pi_{1,-}^{\infty})\big)=0=\Hom_{G}\big(\cF_{P_{I_0}}^{G}(M_0,\mathrm{ker}(\pi_0^{\infty}\twoheadrightarrow\pi_{0,-}^{\infty})),\cF_{P_{I_1}}^{G}(M_1,\pi_{1,-}^{\infty})\big)\]
hence we have injections (using $V_0\twoheadrightarrow \cF_{P_{I_0}}^{G}(M_0,\pi_{0,-}^{\infty})$ and $\cF_{P_{I_1}}^{G}(M_1,\pi_{1,-}^{\infty})\hookrightarrow V_1$)
\[\mathrm{Ext}_{G}^1(\cF_{P_{I_0}}^{G}(M_0,\pi_{0,-}^{\infty}),\cF_{P_{I_1}}^{G}(M_1,\pi_{1,-}^{\infty})) \hookrightarrow
\mathrm{Ext}_{G}^1(V_0,\cF_{P_{I_1}}^{G}(M_1,\pi_{1,-}^{\infty})) \hookrightarrow
\mathrm{Ext}_{G}^1(V_0,V_1).\]
It easily follows that, to prove the injectivity of (\ref{equ: Ext1 g vers OS}), we can assume in the rest of the proof $\pi_0^{\infty}=\pi_{0,-}^{\infty}$ and $\pi_1^{\infty}=\pi_{1,-}^{\infty}$. We fix an injection $f_0: \pi_0^{\infty}\hookrightarrow i_{I,I_0}^{\infty}(\pi^{\infty})$ and a surjection $f_1: i_{I,I_1}^{\infty}(\pi^{\infty})\twoheadrightarrow \pi_1^{\infty}$ (unique up to scalar).\bigskip

\textbf{Step $4$}: We reduce to the case $\pi^{\infty},\pi_0^{\infty}$ and $\pi_1^{\infty}$ are all simple.\\
We choose an arbitrary non-split extension $0\rightarrow M_0\rightarrow M\rightarrow M_1\rightarrow 0$, which induces an extension $0\rightarrow R_1\rightarrow R\rightarrow R_0\rightarrow 0$ where $R\defeq \cF_{P_I}^G(M,\pi^{\infty})$ and $R_i\defeq \cF_{P_I}^G(M_i,\pi^{\infty})\cong \cF_{P_{I_i}}^G(M_i,i_{I,I_i}^{\infty}(\pi^{\infty}))$ for $i=0,1$. We write $V$ for the extension
\begin{equation}\label{equ: nonsplit Ext OS 2}
0\rightarrow V_1\rightarrow V\rightarrow V_0\rightarrow 0
\end{equation}
induced from $R$ by pullback and pushforward along $f_0$ and $f_1$. We want to prove that (\ref{equ: nonsplit Ext OS 2}) is non-split. Let $\sigma^{\infty}$ be an arbitrary irreducible ($G$-regular) constituent of $\pi^{\infty}$, then $\cJ(\sigma^{\infty})\subseteq \cJ(\pi^{\infty})\subseteq \cJ(J_{I_0,I}(\pi_0^{\infty})_{\cB^I_{\Sigma}})$ using the first statement in \ref{it: PS 1} of Lemma \ref{lem: Jacquet of PS} applied to $J_{I_0,I}(\pi_0^{\infty})_{\cB^I_{\Sigma}}$. By the last statement in Lemma~\ref{lem: Jacquet basic} the irreducible ($G$-regular) constituents of $J_{I_0,I}(\pi_0^{\infty})_{\cB^I_{\Sigma}}$ are the $J_{I_0,I}(\sigma_0^{\infty})_{\cB^I_{\Sigma}}$ for $\sigma_0^{\infty}$ an irreducible ($G$-regular) constituent of $\pi_0^{\infty}$. Hence it follows from the last statement in \ref{it: PS 3} of Lemma~\ref{lem: Jacquet of PS} that there exists a unique irreducible constituent $\sigma_0^{\infty}$ of $\pi_0^{\infty}$ such that $\sigma^{\infty}\cong J_{I_0,I}(\sigma_0^{\infty})_{\cB^I_{\Sigma}}$. A similar argument using moreover Remark~\ref{rem: twist Jacquet basic} gives a unique irreducible constituent $\sigma_1^{\infty}$ of $\pi_1^{\infty}$ such that $\sigma^{\infty}\cong J_{I_1,I}'(\sigma_1^{\infty})_{\cB^I_{\Sigma}}$. By (\ref{equ: first adjunction}) and $J_{I_0,I}(\sigma_0^{\infty})_{\cB^I_{\Sigma}}\buildrel\sim\over\rightarrow \sigma^{\infty}$ (resp.~(\ref{equ: second adjunction}) and $\sigma^{\infty}\buildrel\sim\over\rightarrow J_{I_1,I}'(\sigma_1^{\infty})_{\cB^I_{\Sigma}}$), we deduce $\sigma_0^{\infty}\buildrel\sim\over\rightarrow \mathrm{soc}_{L_{I_0}}(i_{I,I_0}^{\infty}(\sigma^{\infty}))$ (resp.~$\mathrm{cosoc}_{L_{I_1}}(i_{I,I_1}^{\infty}(\sigma^{\infty}))\buildrel\sim\over\rightarrow \sigma_1^{\infty}$).

Now, let us assume that $\cF_{P_I}^{G}(M,\sigma^{\infty})$ admits a uniserial length $2$ subquotient with socle $\cF_{P_{I_1}}^G(M_1,\sigma_1^{\infty})$ and cosocle $\cF_{P_{I_0}}^G(M_0,\sigma_0^{\infty})$, then so does $R=\cF_{P_I}^G(M,\pi^{\infty})$. But such a uniserial length $2$ subquotient of $R$ must also be a subquotient of $V$ as $R$ is multiplicity free (using Lemma \ref{lem: Hom OS}) and $\cF_{P_{I_i}}^G(M_i,\sigma_i^{\infty})$ for $i=0,1$ is a constituent of $V_i$, hence of $V$. In particular (\ref{equ: nonsplit Ext OS 2}) is then non-split. Hence to prove the injectivity of (\ref{equ: Ext1 g vers OS}), we see that can we replace $\pi^{\infty}$, $\pi_0^{\infty}$ and $\pi_1^{\infty}$ by $\sigma^{\infty}$, $\sigma_0^{\infty}$ and $\sigma_1^{\infty}$ respectively, and it is enough to prove that $\cF_{P_I}^{G}(M,\pi^{\infty})$ admits such a subquotient.\bigskip

\textbf{Step $5$}: We assume $\pi^{\infty}$, $\pi_0^{\infty}$, $\pi_1^{\infty}$ ($G$-regular) irreducible and prove that $\cF_{P_I}^{G}(M,\pi^{\infty})$ has a uniserial length $2$ subquotient with socle $\cF_{P_{I_1}}^G(M_1,\pi_1^{\infty})$ and cosocle $\cF_{P_{I_0}}^G(M_0,\pi_0^{\infty})$.\\
Note \ that \ we \ then \ have \ $\pi_0^{\infty}\buildrel\sim\over\rightarrow \mathrm{soc}_{L_{I_0}}(i_{I,I_0}^{\infty}(\pi^{\infty}))$ \ and \ $\mathrm{cosoc}_{L_{I_1}}(i_{I,I_1}^{\infty}(\pi^{\infty}))\buildrel\sim\over\rightarrow \pi_1^{\infty}$. As Lemma~\ref{lem: sm to OS} (together with \ref{it: basic as image} of Remark \ref{rem: basic PS intertwine}) implies that $\cF_{P_I}^{G}(M_i,\pi^{\infty})$ has simple socle $\cF_{P_{I_i}}^{G}(M_i,\mathrm{soc}_{L_{I_i}}(i_{I,I_i}^{\infty}(\pi^{\infty})))$ and simple cosocle $\cF_{P_{I_i}}^{G}(M_i,\mathrm{cosoc}_{L_{I_i}}(i_{I,I_i}^{\infty}(\pi^{\infty})))$ for $i=0,1$, we see that the existence of such a subquotient forces $\cF_{P_I}^{G}(M,\pi^{\infty})$ to have simple socle $\cF_{P_{I_1}}^{G}(M_1,\mathrm{soc}_{L_{I_1}}(i_{I,I_1}^{\infty}(\pi^{\infty})))$ and simple cosocle $\cF_{P_{I_0}}^{G}(M_0,\mathrm{soc}_{L_{I_0}}(i_{I,I_0}^{\infty}(\pi^{\infty})))$.
Also, as $\mathrm{Ext}_{U(\fg)}^1(M_1,M_0)\neq 0$ we have either $x_0\prec x_1$ or $x_1\prec x_0$ by \ref{it: rabiotext 2} of Lemma~\ref{rabiotext}. We have the following two possibilities.

\textbf{Case $5.1$}: Assume $\Hom_{U(\fg)}(M^I(x_1),M)=0$. Since $I\subseteq I_1$ $M^{I_1}(x_1)$ is a quotient of $M^I(x_1)$ (\cite[Thm.~9.4(c)]{Hum08}) and thus we \emph{a fortiori} have $\Hom_{U(\fg)}(M^{I_1}(x_1),M)=0$. Let $V_2\defeq \cF_{P_{I_1}}^{G}(M^{I_1}(x_1),\pi_1^{\infty})$, by (\ref{equ: OS main seq}) for $k=\ell=0$ and $w=1$ this implies $\Hom_{G}(R, V_2)=0$. The surjection $M^{I_1}(x_1)\twoheadrightarrow L(x_1)$ induces an injection $V_1\hookrightarrow V_2$ and thus an injection $\Hom_{G}(R, V_1)\hookrightarrow \Hom_{G}(R, V_2)$ which thus implies $\Hom_{G}(R, V_1)=0$. As $\mathrm{cosoc}_{G}(R_1)\buildrel\sim\over\rightarrow V_1$ by Lemma~\ref{lem: sm to OS} and $\mathrm{cosoc}_{L_{I_1}}(i_{I,I_1}^{\infty}(\pi^{\infty}))\buildrel\sim\over\rightarrow \pi_1^{\infty}$, there exists an irreducible constituent $V_0'=\cF_{P_{I_0}}^{G}(L(x_0),\tau_0^{\infty})$ of $R_0$ (with $\tau_0^{\infty}$ an irreducible constituent of $i_{I,I_0}^{\infty}(\pi^{\infty})$) such that $R$ admits a subquotient $V'$ with socle $V_1$ and cosocle $V_0'$. In particular $\mathrm{Ext}_{G}^1(V_0',V_1)\neq 0$, which together with \ref{it: Ext vanishing 1} of Lemma~\ref{lem: improved vanishing} implies $d(\tau_0^{\infty},\pi_1^{\infty})\leq 1$. Now by (\ref{equ: first adjunction}) we have for $k\geq 0$
\[\mathrm{Ext}_{G}^k(i_{I_0, \Delta}^\infty(\tau_0^{\infty}), i_{I_1,\Delta}^{\infty}(\pi_1^{\infty}))^{\infty}\cong \mathrm{Ext}_{L_{I_1}}^k(J_{\Delta,I_1}(i_{I_0, \Delta}^\infty(\tau_0^{\infty})), \pi_1^{\infty})^{\infty}\]
and thus from \ref{it: sml1} of Lemma \ref{lem: smooth geometric lemma} applied with $I=\Delta$ we deduce
\[\mathrm{Ext}_{G}^k(i_{I_0, \Delta}^\infty(\tau_0^{\infty}), i_{I_1,\Delta}^{\infty}(\pi_1^{\infty}))^{\infty} \cong \bigoplus_{w\in W^{I_0,I_1}} \mathrm{Ext}_{L_{I_1}}^k\big(i_{I_0,I_1,w}^{\infty}(J_{I_0,I_1,w}(\tau_0^{\infty})), \pi_1^{\infty}\big)^{\infty}.\]
But from (\ref{equ: support bound}) we have $\cJ(i_{I_0,I_1,w}^{\infty}(J_{I_0,I_1,w}(\pi_0^{\infty})))\subseteq W(L_{I_1})w^{-1}\Sigma_0$, and from (\ref{equ: Ext1 g coset}) and the $G$-regularity (\ref{it: PS 1} of Lemma \ref{lem: Jacquet of PS}) we have $W(L_{I_1})w^{-1}\Sigma_0 \cap \Sigma_1=\emptyset$ if $w\ne 1$. It follows that $(i_{I_0,I_1,w}^{\infty}(J_{I_0,I_1,w}(\tau_0^{\infty})))_{\cB^{I_1}_{\Sigma_1}}=0$ if $w\ne 1$ and thus $\mathrm{Ext}_{G}^k(i_{I_0, \Delta}^\infty(\tau_0^{\infty}), i_{I_1,\Delta}^{\infty}(\pi_1^{\infty}))^{\infty} \cong \mathrm{Ext}_{L_{I_1}}^k\big(i_{I,I_1}^{\infty}(J_{I_0,I}(\tau_0^{\infty})), \pi_1^{\infty}\big)^{\infty}$ for $k\geq 0$. Then the last statement of \ref{it: Hom vers induction2} of Lemma~\ref{lem: Hom vers induction} (applied with $w=1$, $\sigma_0^\infty=\tau_0^\infty$ and $\sigma_1^\infty=\pi_1^\infty$) implies $\tau_0^{\infty}\cong \pi_0^{\infty}\cong \mathrm{soc}_{L_{I_0}}(i_{I,I_0}^{\infty}(\pi^{\infty}))$. In other words, we must have $V_0'\cong V_0$ and $V'$ is the desired $V$.

\textbf{Case $5.2$}: Assume $\Hom_{U(\fg)}(M^I(x_1),M)\ne 0$. As $x_0\neq x_1$ we must have a surjection $M^I(x_1)\twoheadrightarrow M$, which induces an injection $R\hookrightarrow R'\defeq \cF_{P_I}^{G}(M^I(x_1),\pi^{\infty})$ and thus an injection
\begin{equation}\label{equ: Hom vanishing cosocle embedding}
\Hom_{G}(V_0, R)\hookrightarrow \Hom_{G}(V_0, R').
\end{equation}
By (\ref{equ: OS main seq}) (for $k=\ell=0$ and $w=1$) we have $\Hom_{G}(V_0, R')=0$ as $\Hom_{U(\fg)}(M^I(x_1),L(x_0))=0$, and by (\ref{equ: Hom vanishing cosocle embedding}) we deduce $\Hom_{G}(V_0,R)=0$. As $V_0\cong \mathrm{soc}_{G}(R_0)$ by Lemma~\ref{lem: sm to OS}, there exists an irreducible constituent $V_1''=\cF_{P_{I_1}}^{G}(L(x_1),\tau_1^{\infty})$ of $R_1$ such that $R$ admits a subquotient $V''$ with socle $V_1''$ and cosocle $V_0$. In particular $\mathrm{Ext}_{G}^1(V_0,V_1'')\neq 0$, which together with \ref{it: Ext vanishing 1} of Lemma~\ref{lem: improved vanishing} implies $d(\pi_0^{\infty},\tau_1^{\infty})\leq 1$. By an argument analogous to the one in Case $5.1$, \ref{it: Hom vers induction2} of Lemma~\ref{lem: Hom vers induction} (applied with $w=1$, $\sigma_0^\infty=\pi_0^\infty$ and $\sigma_1^\infty=\tau_1^\infty$) implies $\tau_1^{\infty}\cong \pi_1^{\infty}\cong \mathrm{cosoc}_{L_{I_1}}(i_{I,I_1}^{\infty}(\pi^{\infty}))$. In other words, we must have $V_1''\cong V_1$ and $V''$ is the desired $V$. This finishes the proof.
\end{proof}

\begin{lem}\label{lem: Lie length two OS}
For $i=0,1$ let $x_i\in W(G)$ such that $x_0\ne x_1$, $I_i\defeq \Delta\setminus D_L(x_i)$, $I\subseteq I_0\cap I_1$ and $M$ a uniserial length $2$ object in $\cO^{\fp_I}_{\rm{alg}}$ with socle $L(x_0)$ and cosocle $L(x_1)$. Let $\Sigma$ be a $G$-regular finite subset of $\widehat{T}^{\infty}$ preserved under the left action (\ref{equ: smooth dot action}) of $W(L_{I})$, $\pi^{\infty}$ a smooth (finite length) multiplicity free representation of $L_I$ in $\cB^I_{\Sigma}$ and $V\defeq \cF_{P_I}^{G}(M,\pi^{\infty})$. Then we have $\mathrm{soc}_{G}(V)=\cF_{P_{I_1}}^{G}(L(x_1),\mathrm{soc}_{L_{I_1}}(i_{I,I_1}^{\infty}(\pi^{\infty})))$ and $\mathrm{cosoc}_{G}(V)=\cF_{P_{I_0}}^{G}(L(x_0),\mathrm{cosoc}_{L_{I_0}}(i_{I,I_0}^{\infty}(\pi^{\infty})))$.
\end{lem}
\begin{proof}
Replacing $I$ by $I_0\cap I_1$ and $\pi^{\infty}$ by $i_{I,I_0\cap I_1}^{\infty}(\pi^{\infty})$, we can assume $I=I_0\cap I_1$. Using (\ref{equ: sm block decomposition}) we can assume that $\Sigma$ is a single ($G$-regular) $W(L_{I})$-coset, i.e.~$\Sigma = W(L_I)\cdot \chi$ for some regular $\chi$. Then using Remark \ref{rem: support of PS} the irreducible constituents of $\pi^{\infty}$ are (some of) the irreducible constituents of $i_{\emptyset,I}^{\infty}(\chi)$, and using \ref{it: PS 2} of Lemma \ref{lem: Jacquet of PS} for $I_i$ ($i=0,1$) we deduce that $i_{I,I_i}^{\infty}(\pi^{\infty})$ is still multiplicity free. Hence by Lemma~\ref{lem: sm to OS} (and \ref{it: OS property 2} of Theorem~\ref{prop: OS property}) $V_i\defeq \cF_{P_I}^{G}(L(x_i),\pi^{\infty})\cong \cF_{P_{I_i}}^{G}(L(x_i),i_{I,I_i}^{\infty}(\pi^{\infty}))$ is multiplicity free with socle $\cF_{P_{I_i}}^{G}(L(x_i),\mathrm{soc}_{L_{I_i}}(i_{I,I_i}^{\infty}(\pi^{\infty})))$ and cosocle $\cF_{P_{I_i}}^{G}(L(x_i),\mathrm{cosoc}_{L_{I_i}}(i_{I,I_i}^{\infty}(\pi^{\infty})))$ for $i=0,1$. Fix $\sigma_0^{\infty}\in\mathrm{JH}_{L_{I_0}}(\mathrm{soc}_{L_{I_0}}(i_{I,I_0}^{\infty}(\pi^{\infty})))$, $\sigma_1^{\infty}\in\mathrm{JH}_{L_{I_1}}(\mathrm{cosoc}_{L_{I_1}}(i_{I,I_1}^{\infty}(\pi^{\infty})))$ (both irreducible $G$-regular) and let $W_i\defeq \cF_{P_{I_i}}^{G}(L(x_i),\sigma_i^{\infty})$ for $i=0,1$. It suffices to show that $W_0$ (resp.~$W_1$) does not show up in the socle (resp.~cosocle) of $V$. By (\ref{equ: first adjunction}) and the assumption on $\pi^{\infty}$ we have isomorphisms
\begin{equation}\label{equ: Lie length two OS 1}
0\neq \Hom_{L_{I_0}}(\sigma_0^{\infty},i_{I,I_0}^{\infty}(\pi^{\infty}))\cong \Hom_{L_I}(J_{I_0,I}(\sigma_0^{\infty}),\pi^{\infty})\cong \Hom_{L_I}(J_{I_0,I}(\sigma_0^{\infty})_{\cB^I_{\Sigma}},\pi^{\infty}),
\end{equation}
which by the last statement in Lemma~\ref{lem: Jacquet basic} (and our assumption on $\Sigma$) implies that $\tau_0^{\infty}\defeq J_{I_0,I}(\sigma_0^{\infty})_{\cB^I_{\Sigma}}$ is $G$-regular irreducible. Another application of (\ref{equ: first adjunction}) (with $\tau_0^{\infty}$ instead of $\pi^{\infty}$) together with the last statement in \ref{it: PS 3} of Lemma \ref{lem: Jacquet of PS} show that $\sigma_0^{\infty}\buildrel\sim\over\rightarrow \mathrm{soc}_{L_{I_0}}(i_{I,I_0}^{\infty}(\tau_0^{\infty}))$. Note that any injection $W_0\hookrightarrow V_0$ determines by (\ref{equ: Hom OS}) an injection $\sigma_0^{\infty}\hookrightarrow i_{I,I_0}^{\infty}(\pi^{\infty})$, then by (\ref{equ: Lie length two OS 1}) an injection $\tau_0^{\infty}\hookrightarrow \pi^{\infty}$, and finally an injection $\tld{W}_0\defeq \cF_{P_I}^{G}(M,\tau_0^{\infty})\hookrightarrow V$. Moreover $\tld{W}_0$ is multiplicity free as $V$ is by Lemma~\ref{lem: Hom OS} (and using that the $V_i$ are multiplicity free). It follows from Step $5$ of the proof of Lemma~\ref{lem: Ext1 g vers OS} (replacing $\pi^{\infty}$ there by $\tau_0^{\infty}$) that $\tld{W}_0$ contains a unique length $2$ subquotient with socle $\cF_{P_{I_1}}^{G}(L(x_1),\mathrm{cosoc}_{L_{I_1}}(i_{I,I_1}^{\infty}(\tau_0^{\infty})))$ and cosocle $W_0$, and thus the pullback of $0\rightarrow V_1\rightarrow V \rightarrow V_0\rightarrow 0$ along $W_0\hookrightarrow V_0$ is non-split. In particular, $W_0$ does not show up in $\mathrm{soc}_{G}(V)$.
Similarly, we can define $\tau_1^{\infty}\defeq J_{I_1,I}'(\sigma_1^{\infty})_{\cB^I_{\Sigma}}$ and deduce as above from Remark~\ref{rem: twist Jacquet basic} and (\ref{equ: second adjunction}) (together with \ref{it: PS 3} of Lemma \ref{lem: Jacquet of PS} and (\ref{equ: twist Jacquet support})) that $\tau_1^{\infty}$ is $G$-regular irreducible with $\sigma_1^{\infty}\cong \mathrm{cosoc}_{L_{I_1}}(i_{I,I_1}^{\infty}(\tau_1^{\infty}))$. Parallel to the above argument for $W_0$, any surjection $V_1\twoheadrightarrow W_1$ determines a surjection $V\twoheadrightarrow \tld{W}_1\defeq \cF_{P_I}^{G}(M,\tau_1^{\infty})$ with $\tld{W}_1$ admitting a unique length $2$ subquotient with socle $W_1$ and cosocle $\cF_{P_{I_0}}^{G}(L(x_0),\mathrm{soc}_{L_{I_0}}(i_{I,I_0}^{\infty}(\tau_1^{\infty})))$, forcing the pushforward of $0\rightarrow V_1\rightarrow V \rightarrow V_0\rightarrow 0$ along $V_1\twoheadrightarrow W_1$ to be non-split. In particular $W_1$ does not show up in $\mathrm{cosoc}_{G}(V)$, which finishes the proof.
\end{proof}

Keeping the notation at the very beginning of this section, we now assume till its end that, for $i=0,1$, $\Sigma_i$ is a single $G$-regular $W(L_{I_i})$-coset and that the $L_{I_i}$-representation $\pi_i^\infty$ (in $\cB^{I_i}_{\Sigma_i}$) is multiplicity free. We will add other assumptions, depending on the statements.

\begin{lem}\label{lem: Hom OS socle cosocle}
Assume that $I_0=\Delta\setminus D_L(x_0)$ and $M_0=L(x_0)$ for some $x_0\in W(G)$.
\begin{enumerate}[label=(\roman*)]
\item \label{it: Hom OS socle} If $\Hom_{G}(V_0,V_1)\neq 0$, then we have $I_1\subseteq I_0$ and a canonical isomorphism
\begin{equation}\label{equ: Hom OS socle}
\Hom_{G}(V_0,V_1)\cong \Hom_{U(\fg)}(M_1,M_0)\otimes_E \Hom_{L_{I_0}}(\pi_0^{\infty},i_{I_1,I_0}^{\infty}(\pi_1^{\infty})).
\end{equation}
\item \label{it: Hom OS cosocle} If $\Hom_{G}(V_1,V_0)\neq 0$, then we have $I_1\subseteq I_0$ and a canonical isomorphism
\begin{equation}\label{equ: Hom OS cosocle}
\Hom_{G}(V_1,V_0)\cong \Hom_{U(\fg)}(M_0,M_1)\otimes_E \Hom_{L_{I_0}}(i_{I_1,I_0}^{\infty}(\pi_1^{\infty}),\pi_0^{\infty}).
\end{equation}
\end{enumerate}
\end{lem}
\begin{proof}
By Lemma \ref{lem: distinct inf char} (and Harish-Chandra's theorem) we can assume that the irreducible constituents of the $U(\fg)$-module $M_1$ are of the form $L(w)$ for some $w\in W(G)$. We only prove \ref{it: Hom OS cosocle} as \ref{it: Hom OS socle} is symmetric. The assumption $\Hom_{G}(V_1,V_0)\ne 0$ implies that there is an irreducible constituent $L(x_1)$ of $M_1$ such that $\Hom_{G}(\cF_{P_{I_1}}^{G}(L(x_1),\pi_1^{\infty}),V_0)\ne 0$. By \ref{it: OS property 2} of Theorem \ref{prop: OS property} and Lemma \ref{lem: dominance and left set} we can replace $I_1$ by the larger $\Delta\setminus D_L(x_1)$. Then by the first statement of Lemma \ref{lem: Hom OS} we deduce $x_0=x_1$, and in particular $I_1\subseteq I_0$ (note that the smooth induction of $\pi_1^{\infty}$ satisfies the condition at the beginning of this section, for instance using \ref{it: PS 2} of Lemma \ref{lem: Jacquet of PS} for $I=\Delta\setminus D_L(x_1)$).

Considering each irreducible constituent of $\mathrm{soc}_{U(\fg)}(M_1)$ and arguing as above, it follows from (\ref{equ: Hom OS}) that we have an isomorphism (where both sides could \emph{a priori} be $0$)
\[\Hom_{U(\fg)}(M_0,\mathrm{soc}_{U(\fg)}(M_1))\!\otimes_E \!\Hom_{L_{I_0}}(i_{I_1,I_0}^{\infty}(\pi_1^{\infty}),\pi_0^{\infty})\buildrel\sim\over\longrightarrow \Hom_{G}(\cF_{P_{I_1}}^{G}(\mathrm{soc}_{U(\fg)}(M_1),\pi_1^{\infty}),V_0).\]
Hence, as $\Hom_{U(\fg)}(M_0,\mathrm{soc}_{U(\fg)}(M_1))\cong \Hom_{U(\fg)}(M_0,M_1)$ (since $M_0=L(x_0)$), it is enough to show that the map
\begin{equation}\label{equ: Hom OS cosocle 1}
\Hom_{G}(V_1,V_0)\longrightarrow \Hom_{G}(\cF_{P_{I_1}}^{G}(M_1/\mathrm{soc}_{U(\fg)}(M_1),\pi_1^{\infty}),V_0)
\end{equation}
induced by the surjection $M_1\twoheadrightarrow M_1/\mathrm{soc}_{U(\fg)}(M_1)$ is zero. Assume on the contrary that (\ref{equ: Hom OS cosocle 1}) is non-zero and let $M_2$ be a quotient of $M_1/\mathrm{soc}_{U(\fg)}(M_1)$ of minimal length such that the composition $M_1\twoheadrightarrow M_1/\mathrm{soc}_{U(\fg)}(M_1)\twoheadrightarrow M_2$ induces a non-zero map
\begin{equation*}
\Hom_{G}(V_1,V_0)\longrightarrow \Hom_{G}(\cF_{P_{I_1}}^{G}(M_2,\pi_1^{\infty}),V_0).
\end{equation*}
The minimality of $M_2$ easily implies that $\mathrm{soc}_{U(\fg)}(M_2)$ is irreducible, isomorphic to $L(x)$ for some $x\in W(G)$, and that the map $\Hom_{G}(V_1,V_0)\rightarrow \Hom_{G}(\cF_{P_{I_1}}^{G}(M_2/L(x),\pi_1^{\infty}),V_0)$ (induced by $M_1\twoheadrightarrow M_2/L(x)$) is zero. Using the exact sequence
\begin{multline*}
0\rightarrow \Hom_{G}\big(\cF_{P_{I_1}}^{G}(\mathrm{ker}(M_1\rightarrow M_2/L(x)),\pi_1^{\infty}),V_0\big) \rightarrow \Hom_{G}(V_1,V_0)\\
\rightarrow \Hom_{G}\big(\cF_{P_{I_1}}^{G}(M_2/L(x),\pi_1^{\infty}),V_0\big)
\end{multline*}
we deduce that $\mathrm{ker}(M_1\rightarrow M_2/L(x))\hookrightarrow M_1$ induces an isomorphism
\[\Hom_{G}\big(\cF_{P_{I_1}}^{G}(\mathrm{ker}(M_1\rightarrow M_2/L(x)),\pi_1^{\infty}),V_0\big)\buildrel\sim\over\longrightarrow \Hom_{G}(V_1,V_0).\]
Hence we can replace $M_1$ by $\mathrm{ker}(M_1\rightarrow M_2/L(x))$ and $M_2$ by $L(x)$. The surjection $M_1\twoheadrightarrow M_2$ (which is $0$ on $\mathrm{soc}_{U(\fg)}(M_1)$) then factors through surjections $M_1\twoheadrightarrow M_1'\twoheadrightarrow M_2$ where $M_1'$ is uniserial of length $2$ with cosocle $L(x)$ and socle $L(w)$ for some $w\in W(G)$. Moreover we have $w\neq x$ by \ref{it: rabiotext 1} of Lemma~\ref{rabiotext}. Hence $M_1'\twoheadrightarrow M_2$ induces a non-zero map
\begin{equation}\label{equ: Hom OS cosocle 3}
\Hom_{G}(\cF_{P_{I_1}}^{G}(M_1',\pi_1^{\infty}),V_0)\longrightarrow \Hom_{G}(\cF_{P_{I_1}}^{G}(M_2,\pi_1^{\infty}),V_0).
\end{equation}
The target of (\ref{equ: Hom OS cosocle 3}) being non-zero forces $x=x_0$ by Lemma~\ref{lem: Hom OS}.
But by Lemma~\ref{lem: Lie length two OS} (applied with $x_0=w$, $x_1=x$, $I=I_1$, $\pi^\infty=\pi_1^\infty$ and $M=M'_1$) $\cF_{P_{I_1}}^{G}(M_1',\pi_1^{\infty})$ has cosocle $\cF_{P_{I_w}}^{G}(L(w),\mathrm{cosoc}_{L_{I_w}}(i_{I_1,I_w}^{\infty}(\pi_1^{\infty})))$ where $I_w\defeq \Delta\setminus D_L(w)\supseteq I_1$ (by Lemma \ref{lem: dominance and left set}). Lemma~\ref{lem: Hom OS} then forces $\Hom_{G}(\cF_{P_{I_1}}^{G}(M_1',\pi_1^{\infty}),V_0)=0$, a contradiction to (\ref{equ: Hom OS cosocle 3}) being non-zero. This finishes the proof.
\end{proof}

\begin{rem}\label{rem: Hom functorial isom}
Using Remark~\ref{rem: functorial isom}, we can check through the proof of Lemma~\ref{lem: Hom OS} that the isomorphism (\ref{equ: Hom OS}) is functorial in $\pi_0^{\infty}$ and $\pi_1^{\infty}$. Similarly, the proof of Lemma~\ref{lem: Hom OS socle cosocle} shows that both (\ref{equ: Hom OS socle}) and (\ref{equ: Hom OS cosocle}) are functorial in $\pi_0^{\infty}$ and $\pi_1^{\infty}$.
\end{rem}

\begin{lem}
Let $I\subseteq \Delta$, $\Sigma$ a $G$-regular $W(L_{I})$-coset and $\pi^{\infty}$ a smooth (finite length) representation of $L_{I}$ in $\cB^I_{\Sigma}$ such that $i_{I,\Delta}^{\infty}(\pi^{\infty})$ is multiplicity free. Let $M,M'$ in $\cO^{\fp_I}_{\rm{alg}}$ with all irreducible constituents of the form $L(w)$ for some $w\in W(G)$ and $V\defeq \cF_{P_I}^{G}(M,\pi^{\infty})$, $V'\defeq \cF_{P_I}^{G}(M',\pi^{\infty})$. If $V'$ is a subquotient of $V$, then there exists a subquotient $M''$ of $M$ such that $V'\cong \cF_{P_I}^{G}(M'',\pi^{\infty})$.
\end{lem}
\begin{proof}
Note that, although it is possible that $M''=M'$, we actually do not need that (and we do not prove it below). Replacing $M$ by $Q$, we can assume that there does not exist a strict subquotient $Q$ of $M$ such that $V'$ is a subquotient of $\cF_{P_I}^{G}(Q,\pi^{\infty})$. Let $V_{-}'$ (resp.~$V_{+}'$) a (closed) subrepresentation (resp.~a quotient) of $V$ such that the composition $V_{-}'\hookrightarrow V \twoheadrightarrow V_{+}'$ has image $V'$. It is an easy exercise to check that we can always choose $V_{-}'$ and $V_{+}'$ such that the injection $V'\hookrightarrow V_{+}'$ (resp.~the surjection $V_{-}'\twoheadrightarrow V'$) induces an isomorphism on socles (resp.~on cosocles).

Under the first assumption, it is enough to show that $V'=V$, or equivalently $V_{-}'=V=V_{+}'$. We consider an arbitrary surjection $q: M\twoheadrightarrow L(x)$, which induces an injection $\cF_{P_I}^{G}(L(x),\pi^{\infty})\hookrightarrow V$. The composition
\begin{equation}\label{equ: socle to quotient}
\cF_{P_I}^{G}(L(x),\pi^{\infty})\hookrightarrow V\twoheadrightarrow V_{+}'
\end{equation}
must be non-zero, otherwise the surjection $V\twoheadrightarrow V_{+}'$ factors through $\cF_{P_I}^{G}(\mathrm{ker}(q),\pi^{\infty})\twoheadrightarrow V_{+}'$ which forces $V'$ to be a subquotient of $\cF_{P_I}^{G}(\mathrm{ker}(q),\pi^{\infty})$ and thus contradicts our first assumption. Note that $I\subseteq I_x$ by Lemma \ref{lem: dominance and left set} (recall $I_x= \Delta\setminus D_L(x)$). As $i_{I_x,\Delta}^{\infty}(i_{I,I_x}^{\infty}(\pi^{\infty}))\cong i_{I,\Delta}^{\infty}(\pi^{\infty})$ is multiplicity free, so is $i_{I,I_x}^{\infty}(\pi^{\infty})$. By Lemma~\ref{lem: sm to OS} (which can be applied since $i_{I,I_x}^{\infty}(\pi^{\infty})$ satisfies the assumption there using the (second statement in) \ref{it: PS 2} of Lemma \ref{lem: Jacquet of PS}) we deduce that $\cF_{P_I}^{G}(L(x),\pi^{\infty})\cong \cF_{P_{I_x}}^{G}(L(x),i_{I,I_x}^{\infty}(\pi^{\infty}))$ is multiplicity free with socle $\cF_{P_{I_x}}^{G}(L(x),\sigma_0^{\infty})$ where $\sigma_0^{\infty}\defeq \mathrm{soc}_{L_{I_x}}(i_{I,I_x}^{\infty}(\pi^{\infty}))$.
If the non-zero composition (\ref{equ: socle to quotient}) is not injective, then $\mathrm{soc}_{G}(V_{+}')$ has an irreducible constituent $W$ of the form $\cF_{P_{I_x}}^{G}(L(x),\sigma^{\infty})$ for some constituent $\sigma^{\infty}$ of $i_{I,I_x}^{\infty}(\pi^{\infty})/\sigma_0^{\infty}$. Since $\mathrm{soc}_{G}(V')\buildrel\sim\over\rightarrow \mathrm{soc}_{G}(V'_+)$, $W$ is also a constituent of $\mathrm{soc}_{G}(V')$, in particular $\Hom_{G}(W,V')\neq 0$. But by \ref{it: Hom OS socle} of Lemma~\ref{lem: Hom OS socle cosocle} (applied with $V_0=W$ and $V_1=V'$) $\Hom_{G}(W,V')\neq 0$ implies $\Hom_{L_{I_x}}(\sigma^{\infty},i_{I,I_x}^{\infty}(\pi^{\infty}))\neq 0$, which contradicts $\sigma^{\infty}\in \mathrm{JH}_{L_{I_x}}(i_{I,I_x}^{\infty}(\pi^{\infty})/\sigma_0^{\infty})$ (as $i_{I,I_x}^{\infty}(\pi^{\infty})$ is multiplicity free). Hence, the composition (\ref{equ: socle to quotient}) must be injective. Since $q: M\twoheadrightarrow L(x)$ was arbitrary, the composition
\[\cF_{P_I}^{G}(\mathrm{cosoc}_{U(\fg)}(M),\pi^{\infty})\hookrightarrow V\twoheadrightarrow V_{+}'\]
(induced by $M\twoheadrightarrow \mathrm{cosoc}_{U(\fg)}(M)$) has to be injective. Applying Lemma \ref{lem: Lie length two OS} to $\cF_{P_I}^{G}(N,\pi^{\infty})$ where $N$ is any uniserial length $2$ subquotient of $M$ (using \ref{it: rabiotext 1} of Lemma~\ref{rabiotext}), we obtain (using \ref{it: OS property 4} of Theorem \ref{prop: OS property}) that the injection $\cF_{P_I}^{G}(\mathrm{cosoc}_{U(\fg)}(M),\pi^{\infty})\hookrightarrow V$ contains $\mathrm{soc}_{G}(V)$. It follows that the composition $\mathrm{soc}_{G}(V)\hookrightarrow V\twoheadrightarrow V_{+}'$ is also injective, and thus the surjection $V\twoheadrightarrow V_{+}'$ is an isomorphism. A symmetric argument (using \ref{it: Hom OS cosocle} of Lemma~\ref{lem: Hom OS socle cosocle}) shows that the injection $V_{-}'\hookrightarrow V$ must also be an isomorphism, hence we have $V'=V$.
\end{proof}

\begin{lem}\label{lem: Ext1 OS Verma}
For $i=0,1$ assume that $I_i= \Delta\setminus D_L(x_i)$, $M_0=L(x_0)$ and $M_1=M^{I_1}(x_1)$ for some $x_i\in W(G)$ such that $x_0\ne x_1$. Then we have a canonical isomorphism
\begin{equation}\label{equ: Ext1 OS Verma}
\mathrm{Ext}_{G}^1(V_0,V_1)\cong \mathrm{Ext}_{U(\fg)}^1(M_1,M_0)\otimes_E \Hom_{L_{I_1}}(i_{I_0\cap I_1,I_1}^{\infty}(J_{I_0,I_0\cap I_1}(\pi_0^{\infty})),\pi_1^{\infty}).
\end{equation}
\end{lem}
\begin{proof}
Assume $\Sigma_1\cap W(G)\cdot \Sigma_0= \emptyset$, then $\mathrm{Ext}_{G}^1(V_0,V_1)=0$ by \ref{it: OS spectral seq 2} of Corollary~\ref{prop: OS spectral seq} and $\Hom_{L_{I_1}}(i_{I_0\cap I_1,I_1}^{\infty}(J_{I_0,I_0\cap I_1}(\pi_0^{\infty})),\pi_1^{\infty})=0$ by (\ref{equ: refine geometric lemma}) for $w=1$ (using $\cJ(\pi_0^\infty)\subseteq \Sigma_0$), so (\ref{equ: Ext1 OS Verma}) holds (both sides being $0$). Assume $\Sigma_1\cap W(G)\cdot \Sigma_0\ne \emptyset$, then by \ref{it: OS spectral seq 3} of Corollary~\ref{prop: OS spectral seq} and as $\Hom_{U(\fg)}(M_1,M_0)=0$ (since $x_0\neq x_1$), it suffices to show $\mathrm{Ext}_{U(\fg)}^k(M_1,M_0^w)=0$ for $k\leq 1$ and $1\neq w\in W^{I_0,I_1}$. As $L(x_0)^w$ is not in $\cO^{\fb}_{\rm{alg}}$ for $1\neq w\in W^{I_0,I_1}$ by \ref{it: H0 conjugate 3} of Lemma~\ref{lem: H0 Weyl conjugate}, we have $\Hom_{U(\fg)}(M_1,M_0^w)=0$. And we have $\mathrm{Ext}_{U(\fg)}^1(M_1,M_0^w)=0$ for $1\neq w\in W^{I_0,I_1}$ by \ref{rem: H1 conj Ext1 vanishing} of Remark~\ref{mergerem}.
\end{proof}

\begin{lem}\label{lem: Ext1 upper bound}
For $i=0,1$ assume that $I_i= \Delta\setminus D_L(x_i)$ and $M_i=L(x_i)$ for some $x_i\in W(G)$ such that $x_0\ne x_1$. Then we have a canonical injection
\begin{equation}\label{equ: Ext1 upper bound}
\mathrm{Ext}_G^1(V_0,V_1)\hookrightarrow \mathrm{Ext}_{U(\fg)}^1(M_1,M_0)\otimes_E \Hom_{L_{I_1}}(i_{I_0\cap I_1,I_1}^{\infty}(J_{I_0,I_0\cap I_1}(\pi_0^{\infty})),\pi_1^{\infty}).
\end{equation}
In particular $\mathrm{Ext}_G^1(V_0,V_1)\neq 0$ implies $\mathrm{Ext}_{U(\fg)}^1(M_1,M_0)\neq 0$ and
\begin{equation}\label{equ: Ext1 upper bound sm}
\Hom_{L_{I_1}}(i_{I_0\cap I_1,I_1}^{\infty}(J_{I_0,I_0\cap I_1}(\pi_0^{\infty})),\pi_1^{\infty})\neq 0.
\end{equation}
\end{lem}
\begin{proof}
It is enough to prove (\ref{equ: Ext1 upper bound}) when $\mathrm{Ext}_G^1(V_0,V_1)\neq 0$, which we now assume. Let $V_2\defeq \cF_{P_{I_1}}^G(M^{I_1}(x_1),\pi_1^{\infty})$ and note that there is an injection $V_1\hookrightarrow V_2$. As $x_0\neq x_1$, we have $\Hom_{U(\fg)}(M^{I_1}(x_1),L(x_0))=0$. If $\Hom_{G}(V_0,V_2)\ne 0$, then $\Hom_{U(\fg)}(M^{I_1}(x_1),L(x_0))\ne 0$ by \ref{it: Hom OS socle} of Lemma~\ref{lem: Hom OS socle cosocle} which is a contradiction, hence $\Hom_{G}(V_0,V_2)=0$. This together with the exact sequence $0\rightarrow V_1\rightarrow V_2\rightarrow V_2/V_1\rightarrow 0$ induce an exact sequence
\begin{equation}\label{equ: Ext1 OS seq}
0\rightarrow\Hom_{G}(V_0,V_2/V_1)\rightarrow \mathrm{Ext}_{G}^1(V_0,V_1) \rightarrow \mathrm{Ext}_{G}^1(V_0,V_2).
\end{equation}

\textbf{Case $1$}: $\mathrm{Ext}_{G}^1(V_0,V_2)=0$.\\
Then $\Hom_{G}(V_0,V_2/V_1)\buildrel\sim\over\rightarrow \mathrm{Ext}_{G}^1(V_0,V_1)\neq 0$ by (\ref{equ: Ext1 OS seq}). Using $V_2/V_1\cong \cF_{P_{I_1}}^G(N^{I_1}(x_1),\pi_1^{\infty})$, by \ref{it: Hom OS socle} of Lemma~\ref{lem: Hom OS socle cosocle} we have $I_1\subseteq I_0$ and an isomorphism
\begin{equation}\label{equ: Ext1 OS case 1}
\Hom_{G}(V_0,V_2/V_1)\cong \Hom_{U(\fg)}(N^{I_1}(x_1),L(x_0))\otimes_E \Hom_{L_0}(\pi_0^{\infty},i_{I_1,I_0}^{\infty}(\pi_1^{\infty})).
\end{equation}
As $N^{I_1}(x_1)$ is the image of $N(x_1)$ under $M(x_1)\twoheadrightarrow M^{I_1}(x_1)$ and $I_1\subseteq I_0$, we deduce from \cite[Thm.~9.4(c)]{Hum08} and Lemma \ref{lem: dominance and left set} that
\[\Hom_{U(\fg)}(N^{I_1}(x_1),L(x_0))\buildrel\sim\over\longrightarrow \Hom_{U(\fg)}(N(x_1),L(x_0))\ne 0,\]
which by (\ref{equ: O Hom radical}) and \ref{it: Ext O 2} of Lemma \ref{lem: Ext 1 category O} has dimension $\mu(x_1,x_0)$ (and $x_1\prec x_0$). As $x_1<x_0$ and $M^{I_1}(x_1)$ has cosocle $L(x_1)$ we have $\Hom_{U(\fg)}(M^{I_1}(x_1),L(x_0))=0$ which implies
\begin{equation}\label{equ: Ext1 OS case 1 embedding}
\Hom_{U(\fg)}(N^{I_1}(x_1),L(x_0))\hookrightarrow \mathrm{Ext}_{U(\fg)}^1(L(x_1),L(x_0)).
\end{equation}
By \ref{it: rabiotext 2} of Lemma \ref{rabiotext} the target of (\ref{equ: Ext1 OS case 1 embedding}) also has dimension $\mu(x_1,x_0)$, hence (\ref{equ: Ext1 OS case 1 embedding}) is an isomorphism and (\ref{equ: Ext1 OS case 1}) together with $\Hom_{G}(V_0,V_2/V_1)\buildrel\sim\over\rightarrow \mathrm{Ext}_{G}^1(V_0,V_1)$ give
\begin{equation}\label{equ: Ext1 OS case 1 isom}
\mathrm{Ext}_G^1(V_0,V_1)\cong \mathrm{Ext}_{U(\fg)}^1(L(x_1),L(x_0))\otimes_E \Hom_{L_0}(\pi_0^{\infty},i_{I_1,I_0}^{\infty}(\pi_1^{\infty})).
\end{equation}
The right hand side of (\ref{equ: Ext1 OS case 1 isom} is exactly the right hand side of (\ref{equ: Ext1 upper bound}) by (\ref{equ: first adjunction}) and $I_0\cap I_1=I_1$.\bigskip

\textbf{Case $2$}: $\mathrm{Ext}_{G}^1(V_0,V_2)\neq 0$.\\
The isomorphism from Lemma~\ref{lem: Ext1 OS Verma}
\begin{equation}\label{equ: Ext1 OS case 2}
\mathrm{Ext}_{G}^1(V_0,V_2)\cong \mathrm{Ext}_{U(\fg)}^1(M^{I_1}(x_1), L(x_0))\otimes_E \Hom_{L_{I_1}}(i_{I_0\cap I_1,I_1}^{\infty}(J_{I_0,I_0\cap I_1}(\pi_0^{\infty})),\pi_1^{\infty})
\end{equation}
implies $\mathrm{Ext}_{U(\fg)}^1(M^{I_1}(x_1), L(x_0))\neq 0$. Hence the end of Remark~\ref{rem: Levi Ext} (which can be applied since $\mathrm{Ext}_G^1(V_0,V_1)\neq 0$) implies $\mathrm{Ext}_{U(\fg)}^1(L(x_1),L(x_0))\buildrel\sim\over\rightarrow \mathrm{Ext}_{U(\fg)}^1(M^{I_1}(x_1), L(x_0))$ and $\Hom_{U(\fg)}(N^{I_1}(x_1),L(x_0))=0$. In particular (\ref{equ: Ext1 OS case 2}) gives
\begin{equation}\label{equ: Ext1 OS case 2 isom}
\mathrm{Ext}_{G}^1(V_0,V_2)\cong \mathrm{Ext}_{U(\fg)}^1(L(x_1),L(x_0))\otimes_E \Hom_{L_{I_1}}(i_{I_0\cap I_1,I_1}^{\infty}(J_{I_0,I_0\cap I_1}(\pi_0^{\infty})),\pi_1^{\infty}).
\end{equation}
Moreover $\Hom_{U(\fg)}(N^{I_1}(x_1),L(x_0))=0$ implies $\Hom_{G}(V_0,V_2/V_1)=0$ (otherwise use (\ref{equ: Ext1 OS case 1})) and thus $\mathrm{Ext}_{G}^1(V_0,V_1) \hookrightarrow \mathrm{Ext}_{G}^1(V_0,V_2)$ by (\ref{equ: Ext1 OS seq}). With (\ref{equ: Ext1 OS case 2 isom}) this finishes the proof.
\end{proof}

\begin{rem}\label{rem: Ext1 OS distance}
As $i_{I_0\cap I_1,I_1}^{\infty}(J_{I_0,I_0\cap I_1}(\pi_0^{\infty}))\cong i_{I_0,I_1,1}^{\infty}(J_{I_0,I_1,1}(\pi_0^{\infty}))$ is a direct summand of $J_{\Delta, I_1}(i_{I_0,\Delta}^{\infty}(\pi_0^{\infty}))$ by (\ref{equ: smooth geometric lemma}), (\ref{equ: Ext1 upper bound}) being non-zero forces (using (\ref{equ: first adjunction}))
\[0\neq \Hom_{L_{I_1}}(J_{\Delta, I_1}(i_{I_0,\Delta}^{\infty}(\pi_0^{\infty})),\pi_1^{\infty})\cong \Hom_{G}(i_{I_0,\Delta}^{\infty}(\pi_0^{\infty}),i_{I_1,\Delta}^{\infty}(\pi_1^{\infty})).\]
In particular $\mathrm{Ext}_G^1(V_0,V_1)\neq 0$ and $x_0\neq x_1$ imply $d(\pi_0^{\infty},\pi_1^{\infty})=0$.
\end{rem}

\begin{rem}\label{rem: Ext1 functorial}
Using Remark~\ref{rem: functorial isom}, an examination of the proof of Lemma~\ref{lem: Ext1 OS Verma} shows that the isomorphism (\ref{equ: Ext1 OS Verma}) is functorial in $\pi_0^{\infty}$ and $\pi_1^{\infty}$. Since the proof of Lemma~\ref{lem: Ext1 upper bound} is based on Lemma~\ref{lem: Ext1 OS Verma} and Lemma~\ref{lem: Hom OS socle cosocle}, we deduce from Remark~\ref{rem: Hom functorial isom} that the injection (\ref{equ: Ext1 upper bound}) is also functorial in $\pi_0^{\infty}$ and $\pi_1^{\infty}$.
\end{rem}

\begin{prop}\label{prop: Ext1 OS}
For $i=0,1$ assume that $I_i=\Delta\setminus D_L(x_i)$ and $M_i=L(x_i)$ for some $x_i\in W(G)$. Assume moreover $x_0\neq x_1$ and that $\pi_0^{\infty}$, $\pi_1^{\infty}$ are $G$-basic. Then the injection (\ref{equ: Ext1 upper bound}) is an isomorphism
\begin{equation}\label{equ: Ext1 OS}
\mathrm{Ext}_G^1(V_0,V_1)\buildrel\sim\over\longrightarrow \mathrm{Ext}_{U(\fg)}^1(M_1,M_0)\otimes_E \Hom_{L_{I_1}}(i_{I_0\cap I_1,I_1}^{\infty}(J_{I_0,I_0\cap I_1}(\pi_0^{\infty})),\pi_1^{\infty}).
\end{equation}
\end{prop}
\begin{proof}
It suffices to show that both sides of (\ref{equ: Ext1 OS}) have the same dimension when the right hand side is non-zero. It follows from Lemma \ref{lem: Jacquet basic} and Lemma~\ref{lem: smooth geometric lemma} that
\[i_{I_0\cap I_1,I_1}^{\infty}(J_{I_0,I_0\cap I_1}(\pi_0^{\infty}))_{\cB^{I_1}_{\Sigma_1}}=i_{I_0,I_1,1}^{\infty}(J_{I_0,I_1,1}(\pi_0^{\infty}))_{\cB^{I_1}_{\Sigma_1}}\]
is either zero or $G$-basic, and by \ref{it: basic as image} of Remark \ref{rem: basic PS intertwine} that
\begin{equation}\label{equ: Ext1 OS Hom}
\Hom_{L_{I_1}}(i_{I_0\cap I_1,I_1}^{\infty}(J_{I_0,I_0\cap I_1}(\pi_0^{\infty})),\pi_1^{\infty})\cong \Hom_{L_{I_1}}(i_{I_0\cap I_1,I_1}^{\infty}(J_{I_0,I_0\cap I_1}(\pi_0^{\infty}))_{\cB^{I_1}_{\Sigma_1}},\pi_1^{\infty})
\end{equation}
is one dimensional if non-zero. Consequently the right hand side of (\ref{equ: Ext1 OS}), if non-zero, has dimension $\Dim_E \mathrm{Ext}_{U(\fg)}^1(M_1,M_0)$. However, under the assumption that (\ref{equ: Ext1 OS Hom}) is non-zero, it follows from (\ref{equ: Ext1 g vers OS}) and (\ref{equ: Ext1 g vers OS adjunction}) that
\[\Dim_E \mathrm{Ext}_G^1(V_0,V_1)\geq \Dim_E \mathrm{Ext}_{U(\fg)}^1(M_1,M_0),\]
which forces (\ref{equ: Ext1 OS}) to be an isomorphism.
\end{proof}

\begin{lem}\label{lem: Ext1 OS socle cosocle}
For $i=0,1$ assume that $I_i=\Delta\setminus D_L(x_i)$ and $M_i=L(x_i)$ for some $x_i\in W(G)$. Assume moreover $x_0\neq x_1$ and that $\pi_0^{\infty}$, $\pi_1^{\infty}$ are $G$-basic. Let $V$ in $\mathrm{Rep}^{\rm{an}}_{\rm{adm}}(G)$ which fits into a non-split short exact sequence $0\rightarrow V_1\rightarrow V\rightarrow V_0\rightarrow 0$.
\begin{enumerate}[label=(\roman*)]
\item \label{it: Ext1 OS cosocle} We have $\mathrm{cosoc}_{G}(V)\cong \mathrm{cosoc}_{G}(V_0)$ if and only if $V$ has simple cosocle if and only if $\pi_1^{\infty}$ is a quotient of $J_{\Delta,I_1}(i_{I_0,\Delta}^{\infty}(\pi_0^{\infty}))$.
\item \label{it: Ext1 OS socle} We have $\mathrm{soc}_{G}(V)\cong \mathrm{soc}_{G}(V_1)$ if and only if $V$ has simple socle if and only if $\pi_0^{\infty}$ is a subrepresentation of $J_{\Delta,I_0}'(i_{I_1,\Delta}^{\infty}(\pi_1^{\infty}))$.
\end{enumerate}
\end{lem}
\begin{proof}
Recall that, since $\Sigma_i$ is a single $W(L_{I_i})$-coset, we have $\Sigma_i=W(L_{I_i})\cdot\cJ(\pi_i^{\infty})$ for $i=0,1$. The assumptions imply that we have (\ref{equ: Ext1 upper bound sm}) and by (\ref{equ: refine geometric lemma}) (for $w=1$) we deduce $\Sigma_0\cap \Sigma_1\neq \emptyset$. By Lemma~\ref{lem: sm to OS} and (the last statement in) \ref{it: basic as image} of Remark \ref{rem: basic PS intertwine} $V_i$ has simple socle $\cF_{P_{I_i}}^{G}(M_i,\mathrm{soc}_{L_{I_i}}(\pi_i^{\infty}))$ and simple cosocle $\cF_{P_{I_i}}^{G}(M_i,\mathrm{cosoc}_{L_{I_i}}(\pi_i^{\infty}))$ for $i=0,1$. Hence $V$ has simple socle (resp.~cosocle) if and only if $\mathrm{soc}_{G}(V)=\mathrm{soc}_{G}(V_1)$ (resp.~$\mathrm{cosoc}_{G}(V)=\mathrm{cosoc}_{G}(V_0)$) if and only if the map $\mathrm{Ext}_{G}^1(V_0,V_1)\rightarrow \mathrm{Ext}_{G}^1(\mathrm{soc}_{G}(V_0),V_1)$ (resp.~$\mathrm{Ext}_{G}^1(V_0,V_1)\rightarrow \mathrm{Ext}_{G}^1(V_0,\mathrm{cosoc}_{G}(V_1))$) is non-zero. Since the isomorphism (\ref{equ: Ext1 OS}) is functorial in $\pi_0^{\infty}$ and $\pi_1^{\infty}$ by Remark~\ref{rem: Ext1 functorial}, $\mathrm{Ext}_{G}^1(V_0,V_1)\rightarrow \mathrm{Ext}_{G}^1(\mathrm{soc}_{G}(V_0),V_1)$ is non-zero if and only if the map
\begin{equation}\label{equ: Ext1 functorial sm 1}
\Hom_{L_{I_1}}(i_{I_0\cap I_1,I_1}^{\infty}(J_{I_0,I_0\cap I_1}(\pi_0^{\infty})),\pi_1^{\infty})\rightarrow \Hom_{L_{I_1}}(i_{I_0\cap I_1,I_1}^{\infty}(J_{I_0,I_0\cap I_1}(\mathrm{soc}_{L_{I_0}}(\pi_0^{\infty}))),\pi_1^{\infty})
\end{equation}
is non-zero, and $\mathrm{Ext}_{G}^1(V_0,V_1)\rightarrow \mathrm{Ext}_{G}^1(V_0,\mathrm{cosoc}_{G}(V_1))$) is non-zero if and only if the map
\begin{equation}\label{equ: Ext1 functorial sm 2}
\Hom_{L_{I_1}}(i_{I_0\cap I_1,I_1}^{\infty}(J_{I_0,I_0\cap I_1}(\pi_0^{\infty})),\pi_1^{\infty})\rightarrow \Hom_{L_{I_1}}(i_{I_0\cap I_1,I_1}^{\infty}(J_{I_0,I_0\cap I_1}(\pi_0^{\infty})),\mathrm{cosoc}_{L_{I_1}}(\pi_1^{\infty}))
\end{equation}
is non-zero. As $\mathrm{cosoc}_{L_{I_1}}(\pi_1^{\infty})$ is simple, by (both parts of) Lemma \ref{lem: smooth geometric lemma} the map (\ref{equ: Ext1 functorial sm 2}) is non-zero if and only if $\pi_1^{\infty}$ is a quotient of $J_{\Delta,I_1}(i_{I_0,\Delta}^{\infty}(\pi_0^{\infty}))_{\cB^{I_1}_{\Sigma_1}}$ if and only if $\pi_1^{\infty}$ is a quotient of $J_{\Delta,I_1}(i_{I_0,\Delta}^{\infty}(\pi_0^{\infty}))$. This proves \ref{it: Ext1 OS cosocle}.

By (\ref{equ: second adjunction}) followed by (\ref{equ: first adjunction}) and by both parts of Lemma \ref{lem: smooth geometric lemma} we have functorial isomorphisms in representations $\ast$ in $\mathrm{Rep}^{\infty}_{\rm{adm}}(L_{I_0})$
\begin{multline*}
\Hom_{L_{I_0}}(\ast,J_{\Delta,I_0}'(i_{I_1,\Delta}^{\infty}(\pi_1^{\infty})))\cong \Hom_{L_{I_1}}(J_{\Delta,I_1}(i_{I_0,\Delta}^{\infty}(\ast)),\pi_1^{\infty})\\
\cong \Hom_{L_{I_1}}(i_{I_0\cap I_1,I_1}^{\infty}(J_{I_0,I_0\cap I_1}(\ast)),\pi_1^{\infty}).
\end{multline*}
Hence (\ref{equ: Ext1 functorial sm 1}) is the map
\begin{equation}\label{equ: Ext1 functorial sm 3}
\Hom_{L_{I_0}}(\pi_0^{\infty},J_{\Delta,I_0}'(i_{I_1,\Delta}^{\infty}(\pi_1^{\infty})))\longrightarrow \Hom_{L_{I_0}}(\mathrm{soc}_{L_{I_0}}(\pi_0^{\infty}),J_{\Delta,I_0}'(i_{I_1,\Delta}^{\infty}(\pi_1^{\infty}))).
\end{equation}
Since $\mathrm{soc}_{L_{I_0}}(\pi_0^{\infty})$ is simple, (\ref{equ: Ext1 functorial sm 3}) is non-zero if and only if $\pi_0^{\infty}$ is a subrepresentation of $J_{\Delta,I_0}'(i_{I_1,\Delta}^{\infty}(\pi_1^{\infty}))_{\cB^{I_1}_{\Sigma_0}}$ if and only if $\pi_0^{\infty}$ is a subrepresentation of $J_{\Delta,I_0}'(i_{I_1,\Delta}^{\infty}(\pi_1^{\infty}))$. This proves \ref{it: Ext1 OS socle}.
\end{proof}

\begin{lem}\label{lem: simple Ext2 OS}
For $i=0,1$ assume that $I_i=\Delta\setminus D_L(x_i)$ and $M_i=L(x_i)$ for some $x_i\in W(G)$. Assume moreover that $M_0$ is not a constituent of $M^{I_1}(x_1)$, that $d(\pi_0^{\infty},\pi_1^{\infty})\geq 1$ and that $\mathrm{Ext}_{G}^2(V_0,V_1)\neq 0$. Then we have $\mathrm{Ext}_{U(\fg)}^1(M_1,M_0)\neq 0$ and
\begin{equation}\label{equ: simple Ext2 OS}
\mathrm{Ext}_{L_{I_1}}^1(i_{I_0\cap I_1,I_1}^{\infty}(J_{I_0,I_0\cap I_1}(\pi_0^{\infty})),\pi_1^{\infty})^{\infty}\neq 0.
\end{equation}
\end{lem}
\begin{proof}
Let $V_2\defeq \cF_{P_{I_1}}^{G}(M^{I_1}(x_1),\pi_1^{\infty})$, the surjection $M^{I_1}(x_1)\twoheadrightarrow M_1$ induces an injection $V_1\hookrightarrow V_2$ with $V_2/V_1\cong \cF_{P_{I_1}}^{G}(N^{I_1}(x_1),\pi_1^{\infty})$. Since $M_0$ is not a constituent of $M^{I_1}(x_1)$, we have in particular $x_0\ne x'_1$ for any constituent $L(x_1')$ of $M^{I_1}(x_1)$. We then deduce from $d(\pi_0^{\infty},\pi_1^{\infty})\geq 1$ and Remark~\ref{rem: Ext1 OS distance} that $\mathrm{Ext}_{G}^1(V_0,\cF_{P_{I_1}}^{G}(L(x_1'),\pi_1^{\infty}))=0$ for any constituent $L(x_1')$ of $N^{I_1}(x_1)$ (note that we also implicitly use Lemma \ref{lem: dominance and left set} applied to $L(x'_1)$ together with \ref{it: OS property 2} of Theorem \ref{prop: OS property}). Therefore we have $\mathrm{Ext}_{G}^1(V_0,V_2/V_1)=0$ and the exact sequence $0\rightarrow V_1\rightarrow V_2\rightarrow V_2/V_1\rightarrow 0$ induces an embedding
\begin{equation*}
0\neq \mathrm{Ext}_{G}^2(V_0,V_1)\rightarrow \mathrm{Ext}_{G}^2(V_0,V_2).
\end{equation*}
By a d\'evissage on $(\mathrm{Fil}_w(V_0^\vee))_{w\in W^{I_0,I_1}}$ there is $w\in W^{I_0,I_1}$ such that $\mathrm{Ext}_{D(G)}^2(V_2^\vee,\mathrm{gr}_w(V_0^\vee))\neq 0$, which together with (\ref{isov1v0}) (applied with $\mathrm{gr}_w(V_0^\vee)$ instead of $V_0^\vee$), (\ref{equ: ST seq}) (applied with $D=\mathrm{gr}_w(V_0^\vee)$) and Corollary~\ref{cor: Ext P graded} implies
\begin{equation}\label{equ: simple Ext2 OS seq}
\mathrm{Ext}_{L_{I_1}}^k(i_{I_0,I_1,w}^{\infty}(J_{I_0,I_1,w}(\pi_0^{\infty})),\pi_1^{\infty})^{\infty}\otimes_E \mathrm{Ext}_{U(\fg)}^{\ell}(M^{I_1}(x_1), M_0^w)\neq 0
\end{equation}
for some $k,\ell\geq 0$ such that $k+\ell=2$. As $\Hom_{U(\fg)}(M^{I_1}(x_1),M_0^w)=0$ using $M_0\notin \mathrm{JH}_{U(\fg)}(M^{I_1}(x_1))$ and \ref{it: H0 conjugate 3} of Lemma~\ref{lem: H0 Weyl conjugate}, and as $\Hom_{L_{I_1}}(i_{I_0,I_1,w}^{\infty}(J_{I_0,I_1,w}(\pi_0^{\infty})),\pi_1^{\infty})=0$ using $d(\pi_0^{\infty},\pi_1^{\infty})\geq 1$ and (\ref{equ: smooth geometric lemma}) for $I=\Delta$ (with (\ref{equ: first adjunction})), we see that (\ref{equ: simple Ext2 OS seq}) can hold only when $k=\ell=1$. As $\mathrm{Ext}_{U(\fg)}^1(M^{I_1}(x_1),M_0^w)=0$ for $1\neq w\in W^{I_0,I_1}$ by \ref{rem: H1 conj Ext1 vanishing} of Remark~\ref{mergerem}, we must have $w=1$. Then (\ref{equ: simple Ext2 OS seq}) implies (\ref{equ: simple Ext2 OS}) and $\mathrm{Ext}_{U(\fg)}^1(M^{I_1}(x_1),M_0)\neq 0$. Finally, it follows from the discussion in the paragraph below (\ref{equ: Levi Ext H1}) (applied with $w,x$ being $x_1,x_0$) that if $\mathrm{Ext}_{U(\fg)}^1(M^{I_1}(x_1),M_0)\ne 0$ then $\mathrm{Ext}_{U(\fg)}^1(M^{I_1}(x_1),M_0)$ and $\mathrm{Ext}_{U(\fg)}^1(M_1,M_0)$ have the same dimension. In particular $\mathrm{Ext}_{U(\fg)}^1(M_1,M_0)\neq 0$.
\end{proof}

\begin{lem}\label{lem: extra Ext1 Ext2}
For $i=0,1$ assume that $I_i=\Delta\setminus D_L(x_i)$ and $M_i=L(x_i)$ for some $x_i\in W(G)$. Assume moreover that $\Sigma_0\cap \Sigma_1=\emptyset$. Then the following results hold.
\begin{enumerate}[label=(\roman*)]
\item \label{it: extra Ext 1} We have $\mathrm{Ext}_{G}^1(V_0,V_1)=0$.
\item \label{it: extra Ext 2} If $d(\pi_0^{\infty},\pi_1^{\infty})\geq 1$, then we have $\mathrm{Ext}_{G}^2(V_0,V_1)=0$.
\end{enumerate}
\end{lem}
\begin{proof}
We let $V_2\defeq \cF_{P_{I_1}}^{G}(M^{I_1}(x_1),\pi_1^{\infty})$, which contains $V_1$.

We prove \ref{it: extra Ext 1}. Assume $\Hom_{G}(V_0,V_2/V_1)\ne 0$, then by \ref{it: Hom OS socle} of Lemma \ref{lem: Hom OS socle cosocle} we have $I_1\subseteq I_0$ and (using (\ref{equ: first adjunction}))
\[\Hom_{L_0}(\pi_0^{\infty},i_{I_1,I_0}^{\infty}(\pi_1^{\infty}))\cong \Hom_{L_{I_1}}(i_{I_0\cap I_1,I_1}^{\infty}(J_{I_0,I_0\cap I_1}(\pi_0^{\infty})),\pi_1^{\infty})\ne 0.\]
But by (\ref{equ: refine geometric lemma}) for $w=1$ this contradicts the assumption $\Sigma_0\cap \Sigma_1=\emptyset$. Hence $\Hom_{G}(V_0,V_2/V_1)=0$, so the injection $V_1\hookrightarrow V_2$ induces an injection $\mathrm{Ext}_{G}^1(V_0,V_1)\hookrightarrow \mathrm{Ext}_{G}^1(V_0,V_2)$. Assume on the contrary $\mathrm{Ext}_{G}^1(V_0,V_1)\neq 0$ and thus $\mathrm{Ext}_{G}^1(V_0,V_2)\neq 0$. By a d\'evissage on $(\mathrm{Fil}_w(V_0^\vee))_{w\in W^{I_0,I_1}}$ there exists $w\in W^{I_0,I_1}$ such that $\mathrm{Ext}_{D(G)}^2(V_2^\vee,\mathrm{gr}_w(V_0^\vee))\neq 0$, which as in the previous proof implies
\begin{equation}\label{equ: extra Ext spectral seq}
\mathrm{Ext}_{L_{I_1}}^k(i_{I_0,I_1,w}^{\infty}(J_{I_0,I_1,w}(\pi_0^{\infty})),\pi_1^{\infty})^\infty\otimes_E \mathrm{Ext}_{U(\fg)}^{\ell}(M^{I_1}(x_1),L(x_0)^w)\neq 0
\end{equation}
for some $k,\ell\geq 0$ such that $k+\ell=1$. Then (\ref{equ: refine geometric lemma}) and $\Sigma_0\cap \Sigma_1=\emptyset$ force $w\neq 1$. But then we have $\mathrm{Ext}_{U(\fg)}^{\ell}(M^{I_1}(x_1),L(x_0)^w)=0$ for $\ell\leq 1$ by \ref{it: H0 conjugate 3} of Lemma~\ref{lem: H0 Weyl conjugate} ($\ell=0$) and \ref{rem: H1 conj Ext1 vanishing} of Remark~\ref{mergerem} ($\ell =1$), which contradicts (\ref{equ: extra Ext spectral seq}).

We prove \ref{it: extra Ext 2}. By \ref{it: extra Ext 1}, $V_2/V_1\cong \cF_{P_{I_1}}^{G}(N^{I_1}(x_1),\pi_1^{\infty})$ and a d\'evissage on the constituents of $N^{I_1}(x_1)$, we obtain $\mathrm{Ext}_{G}^1(V_0,V_2/V_1)=0$. Assume on the contrary $\mathrm{Ext}_{G}^2(V_0,V_1)\neq 0$, then we must have $\mathrm{Ext}_{G}^2(V_0,V_2)\neq 0$, and by the same d\'evissage as in \ref{it: extra Ext 1} there exists $w\in W^{I_0,I_1}$ such that (\ref{equ: extra Ext spectral seq}) holds for some $k,\ell\geq 0$ with $k+\ell=2$. As $\Sigma_0\cap \Sigma_1=\emptyset$, we have again $w\neq 1$ and thus $\mathrm{Ext}_{U(\fg)}^{\ell}(M^{I_1}(x_1),L(x_0)^w)=0$ for $\ell\leq 1$. But we also have
\[\Hom_{L_{I_1}}(i_{I_0,I_1,w}^{\infty}(J_{I_0,I_1,w}(\pi_0^{\infty})),\pi_1^{\infty})=0\]
by (\ref{equ: smooth geometric lemma}) (for $I=\Delta$), (\ref{equ: first adjunction}) and the assumption $d(\pi_0^{\infty},\pi_1^{\infty})\geq 1$. This contradicts (\ref{equ: extra Ext spectral seq}) for all $k,\ell\geq 0$ such that $k+\ell=2$ and finishes the proof.
\end{proof}

\begin{prop}\label{prop: typical Ext2 OS}
For $i=0,1$ assume that $I\defeq I_i=\Delta\setminus \{j\}$ and $M_i=L(x_i)$ for some $j\in \Delta$ and some $x_i\in W(G)$ such that $D_L(x_i)=\{j\}$. Assume moreover that $\pi_0^{\infty}$, $\pi_1^{\infty}$ are $G$-basic, that $\mathrm{Ext}_G^2(V_0,V_1)\neq 0$ and that
\begin{equation}\label{equ: special cell}
i_{I,I,s_j}^{\infty}(J_{I,I,s_j}(\pi_0^{\infty}))_{\cB^{I}_{\Sigma_1}}\neq 0.
\end{equation}
Then we have $x_0=x_1$,
\begin{equation}\label{equ: special cell Hom}
\Hom_{L_I}(i_{I,I,s_j}^{\infty}(J_{I,I,s_j}(\pi_0^{\infty})),\pi_1^{\infty})\neq 0
\end{equation}
and
\begin{equation}\label{equ: special cell upper bound}
\Dim_E \mathrm{Ext}_G^2(V_0,V_1)\leq \#S_0
\end{equation}
where $S_0=\{x'\mid x'\in W(L_I)x_0, \ \ell(x')=\ell(x_0)+1, \ j\notin D_L(x')\}$.
\end{prop}
\begin{proof}
We have canonical isomorphisms
\begin{multline*}
\Hom_{G}(i_{I,\Delta}^{\infty}(\pi_0^{\infty}),i_{I,\Delta}^{\infty}(\pi_1^{\infty}))\cong \Hom_{L_{I}}(J_{\Delta,I}(i_{I,\Delta}^{\infty}(\pi_0^{\infty})),\pi_1^{\infty})\\
\cong \Hom_{L_{I}}(i_{I,I,s_j}^{\infty}(J_{I,I,s_j}(\pi_0^{\infty})),\pi_1^{\infty})\cong \Hom_{L_{I}}(i_{I,I,s_j}^{\infty}(J_{I,I,s_j}(\pi_0^{\infty}))_{\cB^{I}_{\Sigma_1}},\pi_1^{\infty})
\end{multline*}
where the first isomorphism is (\ref{equ: first adjunction}) and the other two follow from both parts of Lemma~\ref{lem: smooth geometric lemma} together with (\ref{equ: special cell}) remembering that $\Sigma_1=W(L_I)\cdot \cJ(\pi_1^\infty)$ is a $G$-regular left $W(L_I)$-coset. In particular, assuming (\ref{equ: special cell}), (\ref{equ: special cell Hom}) is equivalent to $d(\pi_0^{\infty},\pi_1^{\infty})=0$.

Let $V_2\defeq \cF_{P_{I}}^{G}(M^{I}(x_1),\pi_1^{\infty})$ and recall $M^{I}(x_1)\twoheadrightarrow L(x_1)$ induces $V_1\hookrightarrow V_2$. The exact sequence $0\rightarrow V_1\rightarrow V_2\rightarrow V_2/V_1\rightarrow 0$ induces an exact sequence
\begin{equation}\label{equ: special cell seq}
\mathrm{Ext}_{G}^1(V_0,V_2/V_1)\rightarrow \mathrm{Ext}_{G}^2(V_0,V_1) \rightarrow \mathrm{Ext}_{G}^2(V_0,V_2).
\end{equation}

\textbf{Step $1$}: We prove $\mathrm{Ext}_{G}^1(V_0,V_2/V_1)=0$.\\
We have $V_2/V_1\cong \cF_{P_{I}}^{G}(N^{I}(x_1),\pi_1^{\infty})$ which admits a filtration with graded pieces of the form $\cF_{P_{I}}^{G}(L(x'),\pi_1^{\infty})$ for some $x'>x_1$ with $D_L(x')\subseteq D_L(x_1)$ (by Lemma \ref{lem: dominance and left set}) and hence $D_L(x')=\{j\}$. We deduce from (\ref{equ: special cell}) and Lemma~\ref{lem: smooth geometric lemma} that
\begin{equation}\label{equ: internonempty}
\Sigma_1=W(L_{I})\cdot\cJ(\pi_1^{\infty})\subseteq W(L_{I})s_jW(L_{I})\cdot\cJ(\pi_0^{\infty}),
\end{equation}
which forces
\begin{equation}\label{equ: typical Ext2 disjoint support}
\Sigma_1\cap W(L_{I})\cdot\cJ(\pi_0^{\infty})=\Sigma_1\cap \Sigma_0=\emptyset.
\end{equation}
Then (\ref{equ: typical Ext2 disjoint support}) together with \ref{it: extra Ext 1} of Lemma~\ref{lem: extra Ext1 Ext2} applied to each $\cF_{P_{I}}^{G}(L(x'),\pi_1^{\infty})$ and a d\'evissage imply the statement.\bigskip

\textbf{Step $2$}: We prove the proposition.\\
Note first that Step $1$, (\ref{equ: special cell seq}) and the assumption $\mathrm{Ext}_G^2(V_0,V_1)\neq 0$ imply $\mathrm{Ext}_G^2(V_0,V_2)\neq 0$. By (\ref{equ: OS main seq}) and (\ref{equ: internonempty}) we have a spectral sequence
\[\mathrm{Ext}_{U(\fg)}^{\ell}(M^{I}(x_1),L(x_0)^{s_j})\otimes_E\mathrm{Ext}_{L_{I}}^k(i_{I,I,s_j}^{\infty}(J_{I,I,s_j}(\pi_0^{\infty})),\pi_1^{\infty})^\infty \implies \mathrm{Ext}_{G}^{k+\ell}(V_0,V_2).\]
Since $\mathrm{Ext}_{U(\fg)}^{\ell}(M^{I}(x_1),L(x_0)^{s_j})=0$ for $\ell\leq 1$ by \ref{it: H0 conjugate 3} of Lemma~\ref{lem: H0 Weyl conjugate} ($\ell=0$) and \ref{rem: H1 conj Ext1 vanishing} of Remark~\ref{mergerem} ($\ell=1$), we deduce
\begin{equation}\label{equ: special cell Ext2 reduction}
0\ne \mathrm{Ext}_{G}^2(V_0,V_2)\cong \Hom_{L_{I}}(i_{I,I,s_j}^{\infty}(J_{I,I,s_j}(\pi_0^{\infty})),\pi_1^{\infty})\otimes_E \mathrm{Ext}_{U(\fg)}^2(M^{I}(x_1),L(x_0)^{s_j}).
\end{equation}
From (\ref{equ: special cell Ext2 reduction}) we have (\ref{equ: special cell Hom}), which is one dimensional by \ref{it: Hom vers induction1} of Lemma~\ref{lem: Hom vers induction}, and thus $\Dim_E \mathrm{Ext}_{G}^2(V_0,V_2)=\Dim_E\mathrm{Ext}_{U(\fg)}^2(M^{I}(x_1),L(x_0)^{s_j})\neq 0$, which by Proposition~\ref{prop: typical Ext2} gives $x_0=x_1$ and
\begin{equation}\label{equ: special cell Ext2 dim}
\Dim_E \mathrm{Ext}_{G}^2(V_0,V_2)=\#S_0.
\end{equation}
Finally, we deduce (\ref{equ: special cell upper bound}) from (\ref{equ: special cell seq}), Step $1$ and (\ref{equ: special cell Ext2 dim}).
\end{proof}

Let $I\subseteq \Delta$, $M$ multiplicity free in $\cO_{\rm{alg}}^{\fp_I}$ and recall the set $\mathrm{JH}_{U(\fg)}(M)$ is equipped with a partial order (see \S\ref{generalnotation}). For $\pi^{\infty}$ $G$-basic in $\mathrm{Rep}^{\infty}_{\rm{adm}}(L_I)$ it follows that $\cF_{P_I}^{G}(M,\pi^{\infty})$ is also multiplicity free using Lemma~\ref{lem: Hom OS} and \ref{it: basic as image} of Remark \ref{rem: basic PS intertwine} (and Theorem \ref{prop: OS property}). As $M$ is multiplicity free, for each $L(x)\in \mathrm{JH}_{U(\fg)}(M)$ there is a unique subobject of $M$ with cosocle $L(x)$. This defines an increasing filtration on $M$ indexed by $\mathrm{JH}_{U(\fg)}(M)$, where ``increasing'' means that the inclusions respect the partial order on $\mathrm{JH}_{U(\fg)}(M)$. By the exactness of the contravariant functor $\cF_{P_I}^{G}(-,\pi^{\infty})$ this in turn defines a decreasing filtration on $\cF_{P_I}^{G}(M,\pi^{\infty})$ indexed by $\mathrm{JH}_{U(\fg)}(M)$. We call it the $\mathrm{JH}_{U(\fg)}(M)$-filtration on $\cF_{P_I}^{G}(M,\pi^{\infty})$.

\begin{lem}\label{lem: general construction}
With the above assumptions assume that for each $L(x)\in \mathrm{JH}_{U(\fg)}(M)$ there exists a $G$-basic subquotient $\sigma_x^{\infty}$ of $i_{I,I_x}^{\infty}(\pi^{\infty})$ such that $d(\sigma_x^{\infty},\sigma_{x'}^{\infty})=0$ when $L(x')\leq L(x)$ in $\mathrm{JH}_{U(\fg)}(M)$. Then $\cF_{P_I}^{G}(M,\pi^{\infty})$ contains a unique subquotient $V$ such that
\begin{enumerate}[label=(\roman*)]
\item \label{it: general construction 1} the $\mathrm{JH}_{U(\fg)}(M)$-filtration on $\cF_{P_I}^{G}(M,\pi^{\infty})$ induces a decreasing filtration on $V$ indexed by $\mathrm{JH}_{U(\fg)}(M)$ with $L(x)$-graded piece $V_x\defeq \cF_{P_{I_x}}^{G}(L(x),\sigma_x^{\infty})$ for $L(x)\in \mathrm{JH}_{U(\fg)}(M)$;
\item \label{it: general construction 2} for each uniserial length $2$ subquotient of $M$ with socle $L(x')$ and cosocle $L(x)$, $V$ contains a (unique) subquotient $V_{x,x'}$ which fits into a non-split extension $0\rightarrow V_x\rightarrow V_{x,x'}\rightarrow V_{x'}\rightarrow 0$;
\item \label{it: general construction 3} if $\sigma_x^{\infty}$ is furthermore a subrepresentation (resp.~a quotient) of $i_{I,I_x}^{\infty}(\pi^{\infty})$ for each $L(x)\in \mathrm{JH}_{U(\fg)}(M)$, then $V$ is a subrepresentation (resp.~a quotient) of $\cF_{P_I}^{G}(M,\pi^{\infty})$.
\end{enumerate}
\end{lem}
\begin{proof}
Recall that $i_{I,I_x}^{\infty}(\pi^{\infty})$ and $\sigma_x^{\infty}$ are $G$-basic and thus multiplicity free with simple socle and cosocle (\ref{it: basic as image} of Remark \ref{rem: basic PS intertwine}). For $L(x)\in \mathrm{JH}_{U(\fg)}(M)$, let $\sigma_{x,+}^{\infty}$ (resp.~$\sigma_{x,-}^{\infty}$) be the unique quotient (resp.~subrepresentation) of $i_{I,I_x}^{\infty}(\pi^{\infty})$ with socle $\mathrm{soc}_{L_{I_x}}(\sigma_x^{\infty})$ (resp.~cosocle $\mathrm{cosoc}_{L_{I_x}}(\sigma_x^{\infty})$), which is multiplicity free with simple socle and cosocle.
By Corollary~\ref{cor: basic subquotient} both $\sigma_{x,+}^{\infty}$ and $\sigma_{x,-}^{\infty}$ are $G$-basic. In fact, $\sigma_{x,+}^{\infty}$ (resp.~$\sigma_{x,-}^{\infty}$) is the unique quotient (resp.~subrepresentation) of minimal length of $i_{I,I_x}^{\infty}(\pi^{\infty})$ that admits $\sigma_x^{\infty}$ as a subrepresentation (resp.~a quotient).

We first prove that \ref{it: general construction 1} implies \ref{it: general construction 3}. It follows from Lemma~\ref{lem: sm to OS} and \ref{it: Hom OS socle} of Lemma~\ref{lem: Hom OS socle cosocle} that $\cF_{P_I}^{G}(M,\pi^{\infty})$ has socle
\begin{equation}\label{soc5.1.19}
\bigoplus_{L(x)\in\mathrm{JH}_{U(\fg)}(\mathrm{cosoc}_{U(\fg)}(M))}\cF_{P_{I_x}}^{G}(L(x),\mathrm{soc}_{L_{I_x}}(i_{I,I_x}^{\infty}(\pi^{\infty}))),
\end{equation}
and from Lemma~\ref{lem: sm to OS} and \ref{it: Hom OS cosocle} of Lemma~\ref{lem: Hom OS socle cosocle} that $\cF_{P_I}^{G}(M,\pi^{\infty})$ has cosocle
\begin{equation}\label{cosoc5.1.19}
\bigoplus_{L(x)\in\mathrm{JH}_{U(\fg)}(\mathrm{soc}_{U(\fg)}(M))}\cF_{P_{I_x}}^{G}(L(x),\mathrm{cosoc}_{L_{I_x}}(i_{I,I_x}^{\infty}(\pi^{\infty}))).
\end{equation}
If $\sigma_x^{\infty}$ is a subrepresentation of $i_{I,I_x}^{\infty}(\pi^{\infty})$ for each $L(x)\in \mathrm{JH}_{U(\fg)}(M)$, we have $\mathrm{soc}_{L_{I_x}}(\sigma_x^{\infty})=\mathrm{soc}_{L_{I_x}}(i_{I,I_x}^{\infty}(\pi^{\infty}))$ for each $L(x)\in \mathrm{JH}_{U(\fg)}(M)$ as $i_{I,I_x}^{\infty}(\pi^{\infty})$ is $G$-basic with simple socle by \ref{it: basic as image} of Remark~\ref{rem: basic PS intertwine}. Hence by Lemma~\ref{lem: sm to OS} and \ref{it: general construction 1} we have
\[\mathrm{soc}_{G}(V_x)\cong \cF_{P_{I_x}}^{G}(L(x),\mathrm{soc}_{L_{I_x}}(\sigma_x^{\infty}))= \cF_{P_{I_x}}^{G}(L(x),\mathrm{soc}_{L_{I_x}}(i_{I,I_x}^{\infty}(\pi^{\infty}))),\]
which implies by (\ref{soc5.1.19}) and \ref{it: Hom OS socle} of Lemma~\ref{lem: Hom OS socle cosocle}
\[\mathrm{soc}_{G}(V)\supseteq\bigoplus_{L(x)\in\mathrm{JH}_{U(\fg)}(\mathrm{cosoc}_{U(\fg)}(M))}\mathrm{soc}_{G}(V_x)=\mathrm{soc}_{G}(\cF_{P_I}^{G}(M,\pi^{\infty})).\]
Since $\cF_{P_I}^{G}(M,\pi^{\infty})$ is multiplicity free, this forces $V$ to be a subrepresentation of $\cF_{P_I}^{G}(M,\pi^{\infty})$. If $\sigma_x^{\infty}$ is a quotient of $i_{I,I_x}^{\infty}(\pi^{\infty})$ for each $L(x)\in \mathrm{JH}_{U(\fg)}(M)$, a symmetric argument using (\ref{cosoc5.1.19}) shows $\mathrm{cosoc}_{G}(V)\supseteq \mathrm{cosoc}_{G}(\cF_{P_I}^{G}(M,\pi^{\infty}))$, which forces $V$ to be a quotient of $\cF_{P_I}^{G}(M,\pi^{\infty})$.

We prove \ref{it: general construction 1}, \ref{it: general construction 2} and \ref{it: general construction 3} by increasing induction on the length $\ell(M)\geq 1$. If $\ell(M)=1$ with $\mathrm{JH}_{U(\fg)}(M)=\{L(x)\}$, then $V=V_x$, which is a subquotient of $\cF_{P_I}^{G}(L(x),\pi^{\infty})\cong \cF_{P_{I_x}}^{G}(L(x),i_{I,I_x}^{\infty}(\pi^{\infty}))$. Assume from now on $\ell(M)\geq 2$. If $\ell(\mathrm{cosoc}_{U(\fg)}(M))\geq 2$, then there exists $M',M''\subsetneq M$ such that $M=M'+M''$.
In this case, by Theorem~\ref{prop: OS property} we have
\begin{equation*}
\cF_{P_I}^{G}(M,\pi^{\infty})\cong \cF_{P_I}^{G}(M',\pi^{\infty})\times_{\cF_{P_I}^{G}(M'\cap M'',\pi^{\infty})}\cF_{P_I}^{G}(M'',\pi^{\infty}).
\end{equation*}
By our induction hypothesis there is a subquotient $V'$ (resp.~$V''$) of $\cF_{P_I}^{G}(M',\pi^{\infty})$ (resp.~of $\cF_{P_I}^{G}(M'',\pi^{\infty})$) which satisfies \ref{it: general construction 1} and \ref{it: general construction 2} (with $M$ replaced by $M'$, $M''$ respectively). Moreover, still by induction $V'$ and $V''$ admit a common quotient $V'''$ which is a subquotient of $\cF_{P_I}^{G}(M'\cap M'',\pi^{\infty})$ and satisfies \ref{it: general construction 1} and \ref{it: general construction 2} (with $M$ replaced by $M'\cap M''$). Then we can take $V\defeq V'\times_{V'''}V''$.

We assume from now on $\mathrm{cosoc}_{U(\fg)}(M)\cong L(x)$ for some $x\in W(G)$ and we write $M'\defeq \mathrm{Rad}^1(M)\subsetneq M$. By induction hypothesis there is a subquotient $V'$ of $\cF_{P_I}^{G}(M',\pi^{\infty})$ that satisfies \ref{it: general construction 1} and \ref{it: general construction 2} (with $M$ replaced by $M'$).
By Lemma~\ref{lem: sm to OS} (and the definition of $\sigma_{x',-}^{\infty}$) for $L(x')\in\mathrm{JH}_{U(\fg)}(M')$, $\cF_{P_{I_{x'}}}^{G}(L(x'),\sigma_{x',-}^{\infty})$ is the unique minimal length subrepresentation of $\cF_{P_I}^{G}(L(x'),\pi^{\infty})\cong \cF_{P_{I_{x'}}}^{G}(L(x'),i_{I,I_{x'}}^{\infty}(\pi^{\infty}))$ with $V_{x'}=\cF_{P_{I_{x'}}}^{G}(L(x'),\sigma_{x'}^{\infty})$ as quotient. As $\sigma_{x',-}^{\infty}$ is a subrepresentation of $i_{I,I_{x'}}^{\infty}(\pi^{\infty})$ which admits $\sigma_{x'}^{\infty}$ as quotient, $i_{I_{x'},\Delta}^{\infty}(\sigma_{x',-}^{\infty})$ is a subrepresentation of $i_{I_{x'},\Delta}^{\infty}(i_{I,I_{x'}}^{\infty}(\pi^{\infty}))\cong i_{I,\Delta}^{\infty}(\pi^{\infty})$ which admits $i_{I_{x'},\Delta}^{\infty}(\sigma_{x'}^{\infty})$ as quotient. Let $L(x'),L(x'')\in\mathrm{JH}_{U(\fg)}(M')$ such that $L(x')\leq L(x'')$. The assumption $d(\sigma_{x'}^{\infty},\sigma_{x''}^{\infty})=0$, i.e.~$\Hom_{G}(i_{I_{x'},\Delta}^{\infty}(\sigma_{x'}^{\infty}),i_{I_{x''},\Delta}^{\infty}(\sigma_{x''}^{\infty}))\neq 0$, implies (with \ref{it: basic as image} of Remark \ref{rem: basic PS intertwine})
\[\mathrm{cosoc}_{G}(i_{I_{x'},\Delta}^{\infty}(\sigma_{x',-}^{\infty}))=\mathrm{cosoc}_{G}(i_{I_{x'},\Delta}^{\infty}(\sigma_{x'}^{\infty}))\in \mathrm{JH}_{G}(i_{I_{x''},\Delta}^{\infty}(\sigma_{x''}^{\infty}))\subseteq \mathrm{JH}_{G}(i_{I_{x''},\Delta}^{\infty}(\sigma_{x'',-}^{\infty})),\]
which forces $i_{I_{x'},\Delta}^{\infty}(\sigma_{x',-}^{\infty})\subseteq i_{I_{x''},\Delta}^{\infty}(\sigma_{x'',-}^{\infty})\subseteq i_{I,\Delta}^{\infty}(\pi^{\infty})$ (using that all these representations are $G$-basic hence multiplicity free) and thus $d(\sigma_{x',-}^{\infty},\sigma_{x'',-}^{\infty})=0$. As $\ell(M')<\ell(M)$, by induction there is subquotient $V'_{-}$ of $\cF_{P_I}^{G}(M',\pi^{\infty})$ with a decreasing filtration indexed by $\mathrm{JH}_{U(\fg)}(M')$ such that its $L(x')$-graded piece is $\cF_{P_{I_{x'}}}^{G}(L(x'),\sigma_{x',-}^{\infty})$ for each $L(x')\in\mathrm{JH}_{U(\fg)}(M')$.
By Lemma~\ref{lem: sm to OS} again (and the definition of $\sigma_{x',-}^{\infty}$), for $L(x')\in\mathrm{JH}_{U(\fg)}(M')$, $\cF_{P_{I_{x'}}}^{G}(L(x'),\sigma_{x',-}^{\infty})$ is the unique minimal length subrepresentation of $\cF_{P_I}^{G}(L(x'),\pi^{\infty})\cong \cF_{P_{I_{x'}}}^{G}(L(x'),i_{I,I_{x'}}^{\infty}(\pi^{\infty}))$ with $\cF_{P_{I_{x'}}}^{G}(L(x'),\sigma_{x'}^{\infty})$ as quotient. Hence by \ref{it: general construction 3} applied to $V'_{-}$ (which holds by the induction hypothesis and the beginning of the proof) $V'_{-}$ is the minimal length subrepresentation of $\cF_{P_I}^{G}(M',\pi^{\infty})$ with $V'$ as quotient. As $\cF_{P_I}^{G}(M,\pi^{\infty})$ fits into an exact sequence
\[0\rightarrow \cF_{P_I}^{G}(L(x),\pi^{\infty}) \rightarrow \cF_{P_I}^{G}(M,\pi^{\infty})\rightarrow \cF_{P_I}^{G}(M',\pi^{\infty})\rightarrow 0,\]
the subrepresentation $V'_{-}$ of $\cF_{P_I}^{G}(M',\pi^{\infty})$ and the quotient $\cF_{P_{I_{x}}}^{G}(L(x),\sigma_{x,+}^{\infty})$ of $\cF_{P_I}^{G}(L(x),\pi^{\infty})$ uniquely determine a subquotient $W$ of $\cF_{P_I}^{G}(M,\pi^{\infty})$ that fits into
\begin{equation}\label{equ: subquotient step 1}
0\rightarrow \cF_{P_{I_{x}}}^{G}(L(x),\sigma_{x,+}^{\infty})\rightarrow W \rightarrow V'_{-}\rightarrow 0.
\end{equation}
We construct the desired $V$ as a subquotient of $W$ through the following steps.\bigskip

\textbf{Step $1$}: We prove that, for each length $2$ quotient $Q$ of $M$ (with cosocle $L(x)$ and socle some $L(x')$), the exact sequence
\begin{equation}\label{equ: subquotient step 1 prime}
0\rightarrow \cF_{P_{I_{x}}}^{G}(L(x),\sigma_{x,+}^{\infty})\rightarrow R \rightarrow \cF_{P_{I_{x'}}}^{G}(L(x'),\sigma_{x',-}^{\infty}) \rightarrow 0
\end{equation}
induced from (\ref{equ: subquotient step 1}) is non-split.\\
As $\sigma_{x}^{\infty}$ is a subrepresentation of $\sigma_{x,+}^{\infty}$ and $\sigma_{x'}^{\infty}$ is a quotient of $\sigma_{x',-}^{\infty}$, we deduce from the assumption $d(\sigma_{x'}^{\infty},\sigma_x^{\infty})=0$ that
\begin{equation}\label{equ: subquotient step 1 Hom}
\Hom_{G}(i_{I_{x'}}^{\infty}(\sigma_{x',-}^{\infty}),i_{I_{x}}^{\infty}(\sigma_{x,+}^{\infty}))\neq 0.
\end{equation}
By Lemma~\ref{lem: Lie length two OS} $\cF_{P_I}^{G}(Q,\pi^{\infty})$ has simple socle $\cF_{P_{I_x}}^{G}(L(x),\mathrm{soc}_{L_{I_x}}(i_{I,I_x}^{\infty}(\pi^{\infty})))$ and simple cosocle $\cF_{P_{I_{x'}}}^{G}(L(x'),\mathrm{cosoc}_{L_{I_{x'}}}(i_{I,I_{x'}}^{\infty}(\pi^{\infty})))$. In particular no constituent of $\cF_{P_{I_{x'}}}^{G}(L(x'),\sigma_{x',-}^{\infty})$ shows up in the socle of $\cF_{P_I}^{G}(Q,\pi^{\infty})$. Hence, if (\ref{equ: subquotient step 1 prime}) splits, an easy diagram chase shows that there must exist a constituent $\cF_{P_{I_{x}}}^{G}(L(x),\tau^{\infty})$ of $\cF_{P_I}^{G}(L(x),\pi^{\infty})$ not in $\cF_{P_{I_{x}}}^{G}(L(x),\sigma_{x,+}^{\infty})$ such that $\mathrm{Ext}_{G}^1(\cF_{P_{I_{x'}}}^{G}(L(x'),\sigma_{x',-}^{\infty}),\cF_{P_{I_{x}}}^{G}(L(x),\tau^{\infty}))\neq 0$. By Remark~\ref{rem: Ext1 OS distance} this implies
\begin{equation}\label{equ: subquotient step 1 Hom prime}
\Hom_{G}(i_{I_{x'},\Delta}^{\infty}(\sigma_{x',-}^{\infty}),i_{I_{x},\Delta}^{\infty}(\tau^{\infty}))\neq 0.
\end{equation}
Note that $i_{I_{x},\Delta}^{\infty}(\sigma_{x,+}^{\infty})$ and $i_{I_{x},\Delta}^{\infty}(\tau^{\infty})$ have no common constituent since $i_{I_{x},\Delta}^{\infty}(i_{I,I_{x}}^{\infty}(\pi^{\infty}))\cong i_{I,\Delta}^{\infty}(\pi^{\infty})$ is $G$-basic, thus multiplicity free, and $\tau^{\infty}\in\mathrm{JH}_{L_{I_x}}(i_{I,I_{x}}^{\infty}(\pi^{\infty}))\setminus \mathrm{JH}_{L_{I_x}}(\sigma_{x,+}^{\infty})$. But (\ref{equ: subquotient step 1 Hom}) and (\ref{equ: subquotient step 1 Hom prime}) force both $i_{I_{x}}^{\infty}(\sigma_{x,+}^{\infty})$ and $i_{I_{x}}^{\infty}(\tau^{\infty})$ to have the (simple) cosocle of $i_{I_{x'}}^{\infty}(\sigma_{x',-}^{\infty})$ as a Jordan-H\"older factor, which is a contradiction and shows that (\ref{equ: subquotient step 1 prime}) is non-split.\bigskip

\textbf{Step $2$}: We construct $V$ as a subquotient of $W$.\\
Since $\sigma_{x,+}^{\infty}$ is a quotient of $i_{I,I_x}^{\infty}(\pi^{\infty})$ which contains $\sigma_{x}^{\infty}$ as a subrepresentation, $i_{I_x,\Delta}^{\infty}(\sigma_{x,+}^{\infty})$ is a quotient of $i_{I,\Delta}^{\infty}(\pi^{\infty})\cong i_{I_x,\Delta}^{\infty}(i_{I,I_x}^{\infty}(\pi^{\infty}))$ which contains $i_{I_x,\Delta}^{\infty}(\sigma_{x}^{\infty})$ as a subrepresentation. In fact it is the unique quotient of $i_{I,\Delta}^{\infty}(\pi^{\infty})$ with the same (simple) socle as $i_{I_x,\Delta}^{\infty}(\sigma_{x}^{\infty})$. Likewise, for $L(y)\in\mathrm{JH}_{U(\fg)}(M')$, $i_{I_y,\Delta}^{\infty}(\sigma_{y,-}^{\infty})$ is the unique subrepresentation of $i_{I,\Delta}^{\infty}(\pi^{\infty})$ with the same (simple) cosocle as $i_{I_y,\Delta}^{\infty}(\sigma_{y}^{\infty})$ (recall all these objects are $G$-basic and thus multiplicity free with simple socle and cosocle by \ref{it: basic as image} of Remark \ref{rem: basic PS intertwine}). It follows from $d(\sigma_{y}^{\infty},\sigma_{x}^{\infty})=0$ that the injection $\sigma_{x}^{\infty}\hookrightarrow \sigma_{x,+}^{\infty}$ and the surjection $\sigma_{y,-}^{\infty}\twoheadrightarrow\sigma_{y}^{\infty}$ induce an isomorphism of $1$-dimensional $E$-vector spaces
\[0\neq \Hom_{G}(i_{I_y,\Delta}^{\infty}(\sigma_{y}^{\infty}),i_{I_x,\Delta}^{\infty}(\sigma_{x}^{\infty}))\buildrel\sim\over\longrightarrow \Hom_{G}(i_{I_y,\Delta}^{\infty}(\sigma_{y,-}^{\infty}),i_{I_x,\Delta}^{\infty}(\sigma_{x,+}^{\infty})).\]
In particular, the unique (up to scalar) non-zero map $i_{I_y,\Delta}^{\infty}(\sigma_{y,-}^{\infty})\rightarrow i_{I_x,\Delta}^{\infty}(\sigma_{x,+}^{\infty})$ factors through $i_{I_y,\Delta}^{\infty}(\sigma_{y}^{\infty})\rightarrow i_{I_x,\Delta}^{\infty}(\sigma_{x}^{\infty})$ and we have
\begin{equation}\label{equ: subquotient step 2 JH}
\mathrm{JH}_{G}(i_{I_y,\Delta}^{\infty}(\sigma_{y,-}^{\infty}))\cap \mathrm{JH}_{G}(i_{I_x,\Delta}^{\infty}(\sigma_{x,+}^{\infty}))=\mathrm{JH}_{G}(i_{I_y,\Delta}^{\infty}(\sigma_{y}^{\infty}))\cap \mathrm{JH}_{G}(i_{I_x,\Delta}^{\infty}(\sigma_{x}^{\infty}))
\end{equation}
(which is the set of constituents of $i_{I,\Delta}^{\infty}(\pi^{\infty})$ ``between'' the socle of $i_{I_x,\Delta}^{\infty}(\sigma_{x}^{\infty})$ and the cosocle of $i_{I_y,\Delta}^{\infty}(\sigma_{y}^{\infty})$ for the partial order on $\mathrm{JH}_{G}(i_{I,\Delta}^{\infty}(\pi^{\infty}))$). For $\tau_x^{\infty}\in\mathrm{JH}_{L_{I_x}}(\sigma_{x,+}^{\infty})\setminus \mathrm{JH}_{L_{I_x}}(\sigma_{x}^{\infty})$ and $\tau_{y}^{\infty}\in\mathrm{JH}_{L_{I_y,\Delta}}(\sigma_{y,-}^{\infty})\setminus \mathrm{JH}_{L_{I_y,\Delta}}(\sigma_{y}^{\infty})$ we also have using (\ref{equ: subquotient step 2 JH})
\begin{equation}\label{equ: subquotient step 2 Hom 1}
\Hom_{G}(i_{I_y,\Delta}^{\infty}(\tau_{y}^{\infty}),i_{I_x,\Delta}^{\infty}(\tau_{x}^{\infty}))=0
\end{equation}
and
\begin{equation}\label{equ: subquotient step 2 Hom 2}
\Hom_{G}(i_{I_y,\Delta}^{\infty}(\tau_{y}^{\infty}),i_{I_x,\Delta}^{\infty}(\sigma_{x}^{\infty}))=0=\Hom_{G}(i_{I_y,\Delta}^{\infty}(\sigma_{y}^{\infty}),i_{I_x,\Delta}^{\infty}(\tau_{x}^{\infty})).
\end{equation}
It follows from (\ref{equ: subquotient step 2 Hom 1}) and (\ref{equ: subquotient step 2 Hom 2}) together with Lemma \ref{lem: Hom OS} (for $\ell=0$), Remark~\ref{rem: Ext1 OS distance} and Proposition~\ref{prop: Ext1 OS} (for $\ell=1$) that for $\ell\leq 1$ and $\tau_x$, $\tau_y$ as above
\begin{equation}\label{equ: subquotient step 2 Ext 1}
\mathrm{Ext}_{G}^\ell(\cF_{P_{I_y}}^{G}(L(y),\tau_{y}^{\infty}), \cF_{P_{I_x}}^{G}(L(x),\tau_{x}^{\infty}))=0
\end{equation}
and
\begin{equation}\label{equ: subquotient step 2 Ext 2}
\mathrm{Ext}_{G}^\ell(\cF_{P_{I_y}}^{G}(L(y),\sigma_{y}^{\infty}), \cF_{P_{I_x}}^{G}(L(x),\tau_{x}^{\infty}))=0=\mathrm{Ext}_{G}^\ell(\cF_{P_{I_y}}^{G}(L(y),\tau_{y}^{\infty}), \cF_{P_{I_x}}^{G}(L(x),\sigma_{x}^{\infty})).
\end{equation}
Then by d\'evissage from (\ref{equ: subquotient step 2 Ext 1}) and (\ref{equ: subquotient step 2 Ext 2}) we deduce an isomorphism for $L(y)\in \mathrm{JH}_{U(\fg)}(M')$
\begin{equation}\label{equ: subquotient step 2 factor}
\mathrm{Ext}_{G}^1(\cF_{P_{I_y}}^{G}(L(y),\sigma_{y}^{\infty}), \cF_{P_{I_x}}^{G}(L(x),\sigma_{x}^{\infty}))\buildrel\sim\over\longrightarrow \mathrm{Ext}_{G}^1(\cF_{P_{I_y}}^{G}(L(y),\sigma_{y,-}^{\infty}), \cF_{P_{I_x}}^{G}(L(x),\sigma_{x,+}^{\infty})).
\end{equation}
As for (\ref{equ: subquotient step 2 Ext 1}) and (\ref{equ: subquotient step 2 Ext 2}) we deduce from (\ref{equ: subquotient step 2 Hom 1}) and (\ref{equ: subquotient step 2 Hom 2}) for $\ell\leq 1$
\[\mathrm{Ext}_{G}^\ell(V', \cF_{P_{I_x}}^{G}(L(x),\sigma_{x,+}^{\infty}/\sigma_{x}^{\infty}))=\mathrm{Ext}_{G}^\ell(\mathrm{ker}(V'_-\twoheadrightarrow V'), \cF_{P_{I_x}}^{G}(L(x),\sigma_{x,+}^{\infty}))=0\]
which implies an isomorphism
\begin{equation*}
\mathrm{Ext}_{G}^1(V', V_x)=\mathrm{Ext}_{G}^1(V', \cF_{P_{I_x}}^{G}(L(x),\sigma_{x}^{\infty}))\buildrel\sim\over\longrightarrow \mathrm{Ext}_{G}^1(V'_{-}, \cF_{P_{I_x}}^{G}(L(x),\sigma_{x,+}^{\infty})).
\end{equation*}
In other words, the object $W$ from (\ref{equ: subquotient step 1}) admits a unique subquotient $V$ that fits into a short exact sequence $0\rightarrow V_x\rightarrow V\rightarrow V'\rightarrow 0$. It is then clear that $V$ satisfies \ref{it: general construction 1}.\bigskip

\textbf{Step $3$}: We check \ref{it: general construction 2} for $V$ as in {Step $2$}.\\
Let $Q'$ be a subquotient of $M$ which is uniserial of length $2$ with socle $L(y)$ and cosocle $L(y')$, we want to show that $V$ in Step $2$ admits a (unique) subquotient $V_{y',y}$ which fits into a non-split short exact sequence $0\rightarrow V_{y'}\rightarrow V_{y',y}\rightarrow V_y\rightarrow 0$. If $y'\neq x$, then $Q'$ is a subquotient of $M'$ and the existence of $V_{y',y}$ follows from the induction hypothesis on $V'$. If $y'=x$, then it follows from Step $1$ that $W$ contains a subrepresentation $R$ which fits into a non-split extension (\ref{equ: subquotient step 1 prime}). But (\ref{equ: subquotient step 2 factor}) ensures that $R$ contains a subquotient $V_{x,y}$ as desired.
\end{proof}

Let $I\subseteq \Delta$, $M$ multiplicity free in $\cO_{\rm{alg}}^{\fp_I}$ and $\pi^{\infty}$ $G$-basic in $\mathrm{Rep}^{\infty}_{\rm{adm}}(L_I)$. For $L(x)\in\mathrm{JH}_{U(\fg)}(M)$ let $\sigma_x^{\infty}$ be a $G$-basic subquotient of $i_{I,I_x}^{\infty}(\pi^{\infty})$ as in Lemma~\ref{lem: general construction} and let $V$ as in \emph{loc.~cit.} Let $V_{-}$ be the minimal length subrepresentation of $\cF_{P_I}^G(M,\pi^{\infty})$ with $V$ as quotient. We define $\sigma_{x,-}^{\infty}$ as in the first paragraph of the proof of Lemma~\ref{lem: general construction}, and as in \emph{loc.~cit.} (see the paragraph before (\ref{equ: subquotient step 1}) with $M'$ there replaced by $M$), $V_{-}$ has a decreasing filtration (indexed by $\mathrm{JH}_{U(\fg)}(M)$) induced by the one on $\cF_{P_I}^G(M,\pi^{\infty})$ with $L(x)$-graded piece $\cF_{P_{I_x}}^{G}(L(x),\sigma_{x,-}^{\infty})$ for $L(x)\in\mathrm{JH}_{U(\fg)}(M)$. Now let $x_0\in W(G)$, $I_0\defeq \Delta\setminus D_L(x_0)$, $\pi_0^{\infty}$ a smooth $G$-basic representation of $L_{I_0}$, $V_0\defeq \cF_{P_{I_0}}^{G}(L(x_0),\pi_0^{\infty})$, $\Sigma_0\defeq W(L_0)\cdot\cJ(\pi_0^{\infty})$ and $\Sigma_{0,x}\defeq \Sigma_0\cap W(L_{I_x})\cdot\cJ(\sigma_x^{\infty})$ for $L(x)\in\mathrm{JH}_{U(\fg)}(M)$. Note that $\Sigma_{0,x}$, if non-empty, is a single left $W(L_{I_0\cap I_x})$-coset by \ref{it: PS 1} of Lemma \ref{lem: Jacquet of PS} and $G$-regularity.

\begin{lem}\label{lem: minimal OS}
With the above assumptions assume $d(\pi_0^{\infty},\sigma_x^{\infty})=0$ for $L(x)\in\mathrm{JH}_{U(\fg)}(M)$, $J_{I_x,I_0\cap I_x}'(\tau_{x}^{\infty})_{\cB^{I_0\cap I_x}_{\Sigma_{0,x}}}=0$ for $L(x)\in\mathrm{JH}_{U(\fg)}(M)\setminus\{L(x_0)\}$ and $\tau_{x}^{\infty}\in\mathrm{JH}_{L_{I_x}}(\sigma_{x,-}^{\infty})\setminus \mathrm{JH}_{L_{I_x}}(\sigma_{x}^{\infty})$, and $\Sigma_{0,x_0}=\emptyset$ if $L(x_0)\in\mathrm{JH}_{U(\fg)}(M)$. Assume moreover that $L(x_0)$ is not a constituent of $M^{I_x}(x)$ for $L(x)\!\in\! \mathrm{JH}_{U(\fg)}(M)$ with $x \neq x_0$. Then the unique (up to scalar) injection $V_{-}\hookrightarrow \cF_{P_I}^G(M,\pi^{\infty})$ and surjection $V_{-}\twoheadrightarrow V$ induce isomorphisms
\begin{equation}\label{equ: minimal OS}
\mathrm{Ext}_{G}^1(V_0,V) \buildrel{\sim} \over\longleftarrow\mathrm{Ext}_{G}^1(V_0,V_{-}) \buildrel{\sim} \over\longrightarrow \mathrm{Ext}_{G}^1(V_0,\cF_{P_I}^G(M,\pi^{\infty})).
\end{equation}
\end{lem}
\begin{proof}
Note first that, for any $L(x)\in\mathrm{JH}_{U(\fg)}(M)$ and any constituent $\tau_x^{\infty}$ of $i_{I,I_x}^{\infty}(\pi^{\infty})$, it follows from \ref{it: PS 1} and \ref{it: PS 2} of Lemma \ref{lem: Jacquet of PS} that $W(L_{I_x})\cdot\cJ(\tau_{x}^{\infty})=W(L_{I_x})\cdot\cJ(\sigma_{x}^{\infty})=W(L_{I_x})\cdot\cJ(\pi^{\infty})$.\bigskip

\textbf{Step $1$}: We prove
\begin{equation}\label{equ: minimal OS vanishing 1}
\mathrm{Ext}_{G}^1(V_0,\cF_{P_{I_x}}^{G}(L(x),\tau_x^{\infty}))=0
\end{equation}
for $L(x)\in\mathrm{JH}_{U(\fg)}(M)$ and $\tau_x^{\infty}\in\mathrm{JH}_{L_{I_x}}(i_{I,I_x}^{\infty}(\pi^{\infty}))\setminus \mathrm{JH}_{L_{I_x}}(\sigma_{x}^{\infty})$.\\
As $d(\pi_0^{\infty},\sigma_{x}^{\infty})=0$, $i_{I_x,\Delta}^{\infty}(\sigma_{x}^{\infty})$ contains $\mathrm{cosoc}_{G}(i_{I_0,\Delta}^{\infty}(\pi_0^{\infty}))$ as an (irreducible) constituent. As $i_{I_x,\Delta}^{\infty}(i_{I,I_x}^{\infty}(\pi^{\infty}))\cong i_{I,\Delta}^{\infty}(\pi^{\infty})$ is multiplicity free (since $G$-basic), for $\tau_x$ as above $i_{I_x,\Delta}^{\infty}(\tau_x^{\infty})$ cannot have $\mathrm{cosoc}_{G}(i_{I_0,\Delta}^{\infty}(\pi_0^{\infty}))$ as a constituent, which forces
\begin{equation}\label{equ: minimal OS sm}
\Hom_{G}(i_{I_0,\Delta}^{\infty}(\pi_0^{\infty}), i_{I_x,\Delta}^{\infty}(\tau_x^{\infty}))=0.
\end{equation}
If $x\neq x_0$, then (\ref{equ: minimal OS sm}) and Remark~\ref{rem: Ext1 OS distance} imply (\ref{equ: minimal OS vanishing 1}). If $x=x_0$, then $\Sigma_{0,x_0}=\emptyset$ and \ref{it: extra Ext 1} of Lemma~\ref{lem: extra Ext1 Ext2} imply (\ref{equ: minimal OS vanishing 1}).\bigskip

\textbf{Step $2$}: We prove
\begin{equation}\label{equ: minimal OS vanishing 2}
\mathrm{Ext}_{G}^2(V_0,\cF_{P_{I_x}}^{G}(L(x),\tau_{x}^{\infty}))=0
\end{equation}
for $L(x)\in\mathrm{JH}_{U(\fg)}(M)$ and $\tau_{x}^{\infty}\in\mathrm{JH}_{L_{I_x}}(\sigma_{x,-}^{\infty})\setminus \mathrm{JH}_{L_{I_x}}(\sigma_{x}^{\infty})$.\\
It follows from (\ref{equ: minimal OS sm}) that $d(\pi_0^{\infty},\tau_{x}^{\infty})\geq 1$. If $x=x_0$, then $\Sigma_{0,x_0}=\emptyset$ together with $d(\pi_0^{\infty},\tau_{x}^{\infty})\geq 1$ and \ref{it: extra Ext 2} of Lemma~\ref{lem: extra Ext1 Ext2} imply (\ref{equ: minimal OS vanishing 2}). Assume on the contrary that (\ref{equ: minimal OS vanishing 2}) fails for some $\tau_{x}^{\infty}$ with $x\neq x_0$, then we deduce from $d(\pi_0^{\infty},\tau_{x}^{\infty})\geq 1$, Lemma~\ref{lem: simple Ext2 OS} and (\ref{equ: second adjunction})
\begin{equation}\label{equ: minimal OS sm Ext}
\mathrm{Ext}_{L_{I_0\cap I_x}}^1(J_{I_0,I_0\cap I_x}(\pi_0^{\infty}),J_{I_x,I_0\cap I_x}'(\tau_{x}^{\infty}))^{\infty}\cong \mathrm{Ext}_{L_{I_x}}^1(i_{I_0\cap I_x,I_x}^{\infty}(J_{I_0,I_0\cap I_x}(\pi_0^{\infty})),\tau_{x}^{\infty})^{\infty}\neq 0.
\end{equation}
We have $\cJ(J_{I_0,I_0\cap I_x}(\pi_0^{\infty}))=\cJ(\pi_0^{\infty})\subseteq \Sigma_0$. We also have $\cJ(J_{I_x,I_0\cap I_x}'(\tau_{x}^{\infty}))\subseteq W(L_{I_x})\cdot\cJ(\tau_{x}^{\infty})=W(L_{I_x})\cdot\cJ(\pi^{\infty})$ using (\ref{equ: twist Jacquet support}) for the first inclusion. Hence we deduce from (\ref{equ: minimal OS sm Ext})
\[\mathrm{Ext}_{L_{I_0\cap I_x}}^1(J_{I_0,I_0\cap I_x}(\pi_0^{\infty})_{\cB^{I_0\cap I_x}_{\Sigma_{0,x}}},J_{I_x,I_0\cap I_x}'(\tau_{x}^{\infty})_{\cB^{I_0\cap I_x}_{\Sigma_{0,x}}})^{\infty}\neq 0,\]
which contradicts the assumption $J_{I_x,I_0\cap I_x}'(\tau_{x}^{\infty})_{\cB^{I_0\cap I_x}_{\Sigma_{0,x}}}=0$.\bigskip

\textbf{Step $3$}: We prove the isomorphisms in (\ref{equ: minimal OS}).\\
As we have assumed $\Sigma_{0,x_0}=\emptyset$ when $L(x_0)\in\mathrm{JH}_{U(\fg)}(M)$, we deduce from Lemma~\ref{lem: Hom OS} that $V_0$ has no constituent in commun with $\cF_{P_I}^G(M,\pi^{\infty})$, and in particular the injection $V_{-}\hookrightarrow \cF_{P_I}^G(M,\pi^{\infty})$ induces an exact sequence
\begin{equation}\label{equ: minimal OS seq 1}
0\rightarrow \mathrm{Ext}_{G}^1(V_0,V_{-}) \rightarrow \mathrm{Ext}_{G}^1(V_0,\cF_{P_I}^G(M,\pi^{\infty})) \rightarrow \mathrm{Ext}_{G}^1(V_0,\cF_{P_I}^G(M,\pi^{\infty})/V_{-}).
\end{equation}
The surjection $V_{-}\twoheadrightarrow V$ induces an exact sequence
\begin{equation}\label{equ: minimal OS seq 2}
\mathrm{Ext}_{G}^1(V_0,V_{--})\rightarrow \mathrm{Ext}_{G}^1(V_0,V_{-}) \rightarrow \mathrm{Ext}_{G}^1(V_0,V) \rightarrow \mathrm{Ext}_{G}^2(V_0,V_{--})
\end{equation}
where we write $V_{--}\subseteq V_{-}$ for the unique subrepresentation such that $V_{-}/V_{--}\cong V$. By d\'evissage using (\ref{equ: minimal OS vanishing 1}) we have $\mathrm{Ext}_{G}^1(V_0,\cF_{P_I}^G(M,\pi^{\infty})/V_{-})=0=\mathrm{Ext}_{G}^1(V_0,V_{--})$, and using (\ref{equ: minimal OS vanishing 2}) we have $\mathrm{Ext}_{G}^2(V_0,V_{--})=0$. Together with (\ref{equ: minimal OS seq 1}) and (\ref{equ: minimal OS seq 2}) we obtain (\ref{equ: minimal OS}).
\end{proof}

\begin{rem}\label{rem: dual minimal OS}
There exists a ``dual'' version of Lemma~\ref{lem: minimal OS} with a parallel proof. Let $V_{+}$ be the minimal length quotient of $\cF_{P_I}^G(M,\pi^{\infty})$ which has $V$ as a subrepresentation and define $\sigma_{x,+}^{\infty}$ as in the proof of Lemma~\ref{lem: general construction} for $L(x)\in\mathrm{JH}_{U(\fg)}(M)$. Similar to $V_{-}$, $V_{+}$ is equipped with a decreasing filtration indexed by $\mathrm{JH}_{U(\fg)}(M)$ which is induced from $\cF_{P_I}^G(M,\pi^{\infty})$, with $L(x)$-graded piece $\cF_{P_{I_x}}^{G}(L(x),\sigma_{x,+}^{\infty})$. Now let $x_0$, $\pi_0^{\infty}$, $V_0$, $\Sigma_0$ and $\Sigma_{0,x}$ (for $L(x)\in\mathrm{JH}_{U(\fg)}(M)$) as before Lemma \ref{lem: minimal OS}. Assume $d(\sigma_x^{\infty},\pi_0^{\infty})=0$ for $L(x)\in\mathrm{JH}_{U(\fg)}(M)$, $J_{I_x,I_0\cap I_x}(\tau_{x}^{\infty})_{\cB^{I_0\cap I_x}_{\Sigma_{0,x}}}=0$ for $L(x)\in\mathrm{JH}_{U(\fg)}(M)\setminus\{L(x_0)\}$ and $\tau_{x}^{\infty}\in\mathrm{JH}_{L_{I_x}}(\sigma_{x,+}^{\infty})\setminus \mathrm{JH}_{L_{I_x}}(\sigma_{x}^{\infty})$, $\Sigma_{0,x_0}=\emptyset$ if $L(x_0)\in\mathrm{JH}_{U(\fg)}(M)$ and $\mathrm{JH}_{U(\fg)}(M) \cap \mathrm{JH}_{U(\fg)}(M^{I_0}(x_0)) \subseteq \{L(x_0)\}$. Then as in (\ref{equ: minimal OS}) the unique (up to scalar) surjection $\cF_{P_I}^G(M,\pi^{\infty})\twoheadrightarrow V_{+}$ and injection $V\hookrightarrow V_{+}$ induce isomorphisms
\begin{equation*}
\mathrm{Ext}_{G}^1(V,V_0) \xleftarrow{\sim} \mathrm{Ext}_{G}^1(V_{+},V_0) \xrightarrow{\sim} \mathrm{Ext}_{G}^1(\cF_{P_I}^G(M,\pi^{\infty}),V_0).
\end{equation*}
\end{rem}

\subsection{\texorpdfstring{$\mathrm{Ext}$}{Ext}-squares and \texorpdfstring{$\mathrm{Ext}$}{Ext}-cubes of Orlik-Strauch representations}\label{subsec: square}

We use (essentially) all previous results to construct locally analytic representations of $G$ which are either uniserial (Lemma \ref{lem: simple length three}) or ``squares'' (Proposition \ref{prop: easy square}, Proposition \ref{prop: hard square}) and ``cubes'' (Proposition \ref{prop: factor cube}) of irreducible Orlik-Strauch representations.\bigskip

We keep the notation of \S\ref{subsec: Ext OS}. We first define some irreducible locally analytic representations via Orlik-Strauch's construction (Theorem \ref{prop: OS property}).\bigskip

For $j,j'\in\Delta$, since $D_L(w_{j,j'})=\{j\}$ we have $L(w_{j,j'})\in \cO_{\rm{alg}}^{\fp_{\widehat{j}}}$ by Lemma~\ref{lem: dominance and left set}, where $\widehat{j}=\Delta\setminus\{j\}$ (see the beginning of \S \ref{subsec: sm example}) and $w_{j,j'}$ is defined in (\ref{wj1j2}) (in particular we always have $w_{j,j'}\ne 1$). We set
\[\mathbf{J}=\mathbf{J}_n\defeq \{\un{j}=(j_0,j_1,j_2)\mid 1\leq j_0\leq n-1,~1\leq j_1\leq n-1,~1\leq j_2\leq n,~0\leq j_2-j_1\leq n-1\}\]
that we equip with the partial order $\un{j}=(j_0,j_1,j_2)\leq \un{j}'=(j_0',j_1',j_2')$ if and only if $j_0\leq j_0'$, $j_2\leq j_2'$ and $j_2-j_1\leq j_2'-j_1'$. In particular forgetting $j_0$ gives a surjection $\mathbf{J}\twoheadrightarrow \mathbf{J}^{\infty}$ which respects the partial orders (see the beginning of \S \ref{subsec: sm example} for the set $\mathbf{J}^{\infty}$). We will see later in \S \ref{subsec: final} that this partial order on $\mathbf{J}$ is motivated by the layer structure of certain admissible finite length multiplicity free locally analytic representations. For $\un{j},\un{j}'\in\mathbf{J}$ we write $d(\un{j},\un{j}')\defeq |j_0-j_0'|+|j_2-j_2'|+|(j_2-j_1)-(j_2'-j_1')|$ for the \emph{distance} between $\un{j}$ and $\un{j}'$. Finally, for $\un{j}=(j_0,j_1,j_2)\in \mathbf{J}$, we set
\begin{equation}\label{cj}
C_{\un{j}}=C_{(j_0,j_1,j_2)}\defeq \cF_{P_{\widehat{j}_1}}^{G}(L(w_{j_1,j_0}),\pi_{j_1,j_2}^{\infty})
\end{equation}
where $\pi_{j_1,j_2}^{\infty}$ is the irreducible $G$-regular smooth representation of $L_{\widehat{j}_1}$ defined in (\ref{pij_1,_2}).

\begin{lem}\label{lem: Ext1 factor 1}
Let $\un{j},\un{j}'\in\mathbf{J}$ such that $\un{j}'\not \leq \un{j}$. Then $\mathrm{Ext}_G^1(C_{\un{j}'},C_{\un{j}})\neq 0$ if and only if $d(\un{j},\un{j}')=1$, in which case it is one dimensional.
\end{lem}
\begin{proof}
Note first that $\un{j}'\not \leq \un{j}$ with $j_0=j'_0$ and $|j_1-j_1'|\leq 1$ force $\un{j}<\un{j}'$ in $\mathbf{J}$.

Assume first $w_{j_1,j_0}=w_{j_1',j_0'}$, equivalently $(j_0, j_1)=(j'_0, j'_1)$. Hence we have $\un{j}<\un{j}'$ by the previous sentence, hence $j_2<j'_2$ since $j_1=j'_1$, and thus $d(\un{j},\un{j}')\geq 2$. We also have $\Sigma_{j_1,j_2}\cap \Sigma_{j_1',j_2'}=\emptyset$ by \ref{it: coset 1} of Lemma~\ref{lem: coset intersection} (see above (\ref{equ: special Weyl element}) for the notation), which by \ref{it: extra Ext 1} of Lemma~\ref{lem: extra Ext1 Ext2} implies $\mathrm{Ext}_G^1(C_{\un{j}'},C_{\un{j}})=~0$. In particular the lemma trivially holds when $w_{j_1,j_0}=w_{j_1',j_0'}$ since both assertions in the statement never happen.

We assume from now on $w_{j_1,j_0}\neq w_{j_1',j_0'}$ and write $I=\Delta\setminus\{j_1, j'_1\}$. It follows from Proposition~\ref{prop: Ext1 OS} that $\mathrm{Ext}_G^1(C_{\un{j}'},C_{\un{j}})\neq 0$ if and only if $\mathrm{Ext}_{U(\fg)}^1(L(w_{j_1,j_0}),L(w_{j_1',j_0'}))\neq 0$ and
\begin{equation}\label{equ: Ext1 factor 1 Hom sm}
\Hom_{L_{\widehat{j}_1}}(i_{I,\widehat{j}_1}^{\infty}(J_{\widehat{j}_1',I}(\pi_{j_1',j_2'}^{\infty})),\pi_{j_1,j_2}^{\infty})\neq0.
\end{equation}
It follows from the last assertion in \ref{it: rabiotext 2} of Lemma~\ref{rabiotext} that $\mathrm{Ext}_{U(\fg)}^1(L(w_{j_1,j_0}),L(w_{j_1',j_0'}))\neq 0$ if and only if either $w_{j_1,j_0}\prec w_{j_1',j_0'}$ or $w_{j_1',j_0'}\prec w_{j_1,j_0}$. By Remark~\ref{rem: coxeter pair} this is equivalent to $w_{j_1,j_0}< w_{j_1',j_0'}$ with $\ell(w_{j_1',j_0'})=\ell(w_{j_1,j_0})+1$ or $w_{j_1',j_0'}< w_{j_1,j_0}$ with $\ell(w_{j_1,j_0})=\ell(w_{j_1',j_0'})+1$, which is easily checked to be equivalent to $|j_0-j_0'|+|j_1-j_1'|=1$. Hence $\mathrm{Ext}_G^1(C_{\un{j}'},C_{\un{j}})\neq 0$ if and only if $|j_0-j_0'|+|j_1-j_1'|=1$ and (\ref{equ: Ext1 factor 1 Hom sm}) holds. In that case, note that $\mathrm{Ext}_{U(\fg)}^1(L(w_{j_1,j_0}),L(w_{j_1',j_0'}))$ is one dimensional by the last assertion in \ref{it: rabiotext 2} of Lemma~\ref{rabiotext}.

We distinguish the two cases $j_0=j'_0$ and $j_0\ne j'_0$. 

Assume $j_0=j'_0$, then $|j_0-j_0'|+|j_1-j_1'|=1$ is equivalent to $|j_1-j_1'|=1$ which implies $\un{j}<\un{j}'$ by the first sentence and also $(j_1, j_2)< (j'_1, j'_2)$ in $\mathbf{J}^{\infty}$. When $(j_1, j_2)< (j'_1, j'_2)$ in $\mathbf{J}^{\infty}$, (\ref{equ: Ext1 factor 1 Hom sm}) is equivalent to $|j_2-j_2'|+|(j_2-j_1)-(j_2'-j_1')|\leq 1$ by \ref{it: sm cell 1} of Lemma~\ref{lem: sm cell}. It follows that, when $j_0=j'_0$ (and $\un{j}'\not \leq \un{j}$, $w_{j_1,j_0}\neq w_{j_1',j_0'}$), $\mathrm{Ext}_G^1(C_{\un{j}'},C_{\un{j}})\neq 0$ is equivalent to $|j_1-j_1'|=1$ and $|j_2-j_2'|+|(j_2-j_1)-(j_2'-j_1')|\leq 1$, which is easily checked to be equivalent to $|j_2-j_2'|+|(j_2-j_1)-(j_2'-j_1')|=1$ (using $\un{j}'\ne \un{j}$), i.e.~$d(\un{j},\un{j}')=1$.

Assume $j_0\ne j'_0$, then $|j_0-j_0'|+|j_1-j_1'|=1$ implies $j_1=j'_1$ and $I=\widehat{j}_1$, and thus (\ref{equ: Ext1 factor 1 Hom sm}) is equivalent to $j_2=j'_2$. Thus $\mathrm{Ext}_G^1(C_{\un{j}'},C_{\un{j}})\neq 0$ is equivalent to $|j_0-j_0'|=1$, $j_1=j'_1$, $j_2=j'_2$, which is equivalent to $d(\un{j},\un{j}')=1$ (when $j_0\ne j'_0$).
\end{proof}

Recall from the paragraph below Lemma~\ref{lem: explicit smooth induction} that, for $(j_1,j_2)\in\mathbf{J}^{\infty}$, we have defined $I_{j_1,j_2}^+, I_{j_1,j_2}^-\subseteq \Delta$ by $\mathrm{soc}_{G}(i_{\widehat{j}_1,\Delta}^{\infty}(\pi_{j_1,j_2}^{\infty}))\cong V_{I_{j_1,j_2}^+,\Delta}^{\infty}$ and $\mathrm{cosoc}_{G}(i_{\widehat{j}_1,\Delta}^{\infty}(\pi_{j_1,j_2}^{\infty}))\cong V_{I_{j_1,j_2}^-,\Delta}^{\infty}$ (see (\ref{gensteinberg}) for the $G$-regular irreducible representation $V_{I,\Delta}^{\infty}$). Recall also that $L(1)=L(\mu_0)$.

\begin{lem}\label{lem: Ext1 factor 2}
Let $\un{j}\in\mathbf{J}$ and $I\subseteq \Delta$.
\begin{enumerate}[label=(\roman*)]
\item \label{it: Ext1 with alg 1} We have $\mathrm{Ext}_G^1(L(1)^\vee\otimes_E V_{I,\Delta}^{\infty},C_{\un{j}})\neq 0$ if and only if $j_0=j_1$ and $I=I_{j_1,j_2}^+$, in which case it is one dimensional.
\item \label{it: Ext1 with alg 2} We have $\mathrm{Ext}_G^1(C_{\un{j}},L(1)^\vee\otimes_E V_{I,\Delta}^{\infty})\neq 0$ if and only if $j_0=j_1$ and $I=I_{j_1,j_2}^-$, in which case it is one dimensional.
\end{enumerate}
\end{lem}
\begin{proof}
We prove \ref{it: Ext1 with alg 1}. As $1\neq w_{j_1,j_0}$, we deduce from Proposition~\ref{prop: Ext1 OS} that $\mathrm{Ext}_G^1(L(1)^\vee\otimes_E V_{I,\Delta}^{\infty},C_{\un{j}})\neq 0$ if and only if $\mathrm{Ext}_{U(\fg)}^1(L(w_{j_1,j_0}),L(1))\neq 0$ and $\Hom_{G}(V_{I,\Delta}^{\infty},i_{\widehat{j}_1,\Delta}^{\infty}(\pi_{j_1,j_2}^{\infty}))\neq 0$, if and only if $\ell(w_{j_1,j_0})=1$ i.e.~$j_0=j_1$ (using \ref{it: rabiotext 2} of Lemma~\ref{rabiotext}), and $I=I_{j_1,j_2}^+$ (using \ref{it: distance St 1} of Lemma~\ref{lem: distance from St}). The proof of \ref{it: Ext1 with alg 2} is completely analogous.
\end{proof}

\begin{lem}\label{lem: Ext2 factor 1}
Let $\un{j},\un{j}'\in\mathbf{J}$ such that $\un{j}<\un{j}'$ and $j_0'\leq j_0+1$.
\begin{enumerate}[label=(\roman*)]
\item \label{it: Ext2 factor 1 1} If $\mathrm{Ext}_G^2(C_{\un{j'}},C_{\un{j}})\neq 0$ then $(j_1',j_2')\in\{(j_1+1,j_2+1),(j_1-1,j_2),(j_1,j_2+1),(j_1,j_2)\}$.
\item \label{it: Ext2 factor 1 2} Assume \ $(j_1',j_2')=(j_1,j_2+1)$. \ Then \ $\mathrm{Ext}_G^2(C_{\un{j'}},C_{\un{j}})\neq 0$ \ implies \ $j_0'=j_0$, \ and \ $\mathrm{Ext}_G^2(C_{\un{j'}},C_{\un{j}})$ is at most one dimensional except when $2\leq j_0=j_1\leq n-2$ in which case it is at most two dimensional.
\end{enumerate}
\end{lem}
\begin{proof}
We first prove \ref{it: Ext2 factor 1 1}. Assume on the contrary that there exists $\un{j},\un{j}'\in\mathbf{J}$ that satisfies $\un{j}< \un{j}'$, $j_0'\leq j_0+1$, $(j_1',j_2')\notin\{(j_1+1,j_2+1),(j_1-1,j_2),(j_1,j_2+1),(j_1,j_2)\}$ and $\mathrm{Ext}_G^2(C_{\un{j'}},C_{\un{j}})\neq 0$. The condition $(j_1',j_2')\notin\{(j_1+1,j_2+1),(j_1-1,j_2),(j_1,j_2+1),(j_1,j_2)\}$ together with \ref{it: connect 1} of Lemma~\ref{lem: connect Hom} implies $d(\pi_{j_1',j_2'}^{\infty},\pi_{j_1,j_2}^{\infty})\geq 1$. We have two cases.\bigskip

\textbf{Case $1$}: $j_1\neq j_1'$. Then $w_{j_1,j_0}\neq w_{j_1',j_0'}$ and $L(w_{j_1',j_0'})$ is not a constituent of $M^{\widehat{j}_1}(w_{j_1,j_0})$ using Lemma \ref{lem: dominance and left set} and \cite[Thm.~9.4(c)]{Hum08}. We can then apply Lemma~\ref{lem: simple Ext2 OS} which gives $\mathrm{Ext}_{U(\fg)}^1(L(w_{j_1,j_0}),L(w_{j_1',j_0'}))\neq 0$ and $\mathrm{Ext}_{L_{\widehat{j}_1}}^1(i_{\widehat{j'}_1\cap \widehat{j}_1,\widehat{j}_1}^{\infty}(J_{\widehat{j'}_1,\widehat{j'}_1\cap \widehat{j}_1}(\pi_{j_1',j_2'}^{\infty})),\pi_{j_1,j_2}^{\infty})^{\infty}\neq 0$. The first inequality implies $|j_1-j_1'|\leq 1$ by \ref{it: rabiotext 2} of Lemma~\ref{rabiotext} and Remark~\ref{rem: coxeter pair}. By (\ref{equ: first adjunction}) we have $\mathrm{Ext}_{G}^1(i_{\widehat{j'}_1,\Delta}^{\infty}(\pi_{j_1',j_2'}^{\infty}),i_{\widehat{j}_1,\Delta}^{\infty}(\pi_{j_1,j_2}^{\infty}))^{\infty} = \mathrm{Ext}^1_{L_{\widehat{j}_1}}(J_{\Delta, \widehat{j}_1}(i_{\widehat{j'}_1,\Delta}^{\infty}(\pi_{j_1',j_2'}^{\infty})),\pi_{j_1,j_2}^{\infty})^{\infty}$, which is non-zero by (\ref{equ: smooth geometric lemma}) (applied with $I=\Delta$ and taking $w=1$) and the second inequality. Hence $d(\pi_{j_1',j_2'}^{\infty},\pi_{j_1,j_2}^{\infty})=1$, which by \ref{it: connect 2} of Lemma~\ref{lem: connect Hom} implies $(j_1',j_2')\in \{(j_1+2,j_2+2),(j_1-2,j_2)\}$ and thus $|j_1-j_1'|=2$. This contradicts $|j_1-j_1'|\leq 1$.\bigskip

\textbf{Case $2$}: $j_1=j_1'$. Then $j_2< j'_2$ and \ref{it: coset 1} of Lemma~\ref{lem: coset intersection} forces $\Sigma_{j_1,j_2}\cap \Sigma_{j_1',j_2'}=\emptyset$, which by \ref{it: extra Ext 2} of Lemma~\ref{lem: extra Ext1 Ext2} implies $\mathrm{Ext}_G^2(C_{\un{j'}},C_{\un{j}})=0$, a contradiction.

We prove \ref{it: Ext2 factor 1 2}. As $(j_1',j_2')=(j_1,j_2+1)$ by \ref{it: sm cell 2} of Lemma~\ref{lem: sm cell} we have
\[\Hom_{L_{\widehat{j}_1}}(i_{\widehat{j}_1',\widehat{j}_1, s_{j_1}}^{\infty}(J_{\widehat{j}_1',\widehat{j}_1, s_{j_1}}(\pi_{j_1',j_2'}^{\infty})),\pi_{j_1,j_2}^{\infty})\neq 0.\]
By Proposition~\ref{prop: typical Ext2 OS} if $\mathrm{Ext}_G^2(C_{\un{j'}},C_{\un{j}})\ne 0$ then $j_0\neq j_0'$, and when $j_0=j_0'$ then
\begin{equation*}
\Dim_E \mathrm{Ext}_G^2(C_{\un{j'}},C_{\un{j}})\leq \#S_0
\end{equation*}
where $S_0= \{x'\mid x'\in W(L_{\widehat{j}_1})w_{j_1,j_0}, \ \ell(x')=\ell(w_{j_1,j_0})+1, \ j_1\notin D_L(x')\}$. Let $x'\in S_0$ and $j\in D_L(x')$. The condition $j_1\notin D_L(x')$ forces $j=j_1+1$ if $j_1>j_0$, $j=j_1-1$ if $j_1<j_0$, and $j\in\{j_1-1,j_1+1\}$ if $j_1=j_0$. In particular, we always have $S_0\subseteq\{w_{j_1+1,j_0},w_{j_1-1,j_0}\}$, with $w_{j_1+1,j_0}\in S_0$ if and only if $j_0\leq j_1<n-1$, and $w_{j_1-1,j_0}\in S_0$ if and only if $j_0\geq j_1>1$. This finishes the proof.
\end{proof}


\begin{lem}\label{lem: Ext2 factor 2}
Let $\un{j}\in\mathbf{J}$ and $I\subseteq \Delta$.
\begin{enumerate}[label=(\roman*)]
\item \label{it: Ext2 with alg 1} If $I\neq I_{j_1,j_2}^+$ and if $\mathrm{Ext}_G^2(L(1)^\vee\otimes_E V_{I,\Delta}^{\infty},C_{\un{j}})\neq 0$ then $j_0=j_1$, $I\notin[I_{j_1,j_2}^+,I_{j_1,j_2}^-]$ and $d(I,I_{j_1,j_2}^+)=1$ (see (\ref{equ: same interval}) and (\ref{di0i1}) for the notation).
\item \label{it: Ext2 with alg 2} If $I\neq I_{j_1,j_2}^-$ and if $\mathrm{Ext}_G^2(C_{\un{j}},L(1)^\vee\otimes_E V_{I,\Delta}^{\infty})\neq 0$ then $j_0=j_1$, $I\notin[I_{j_1,j_2}^+,I_{j_1,j_2}^-]$ and $d(I_{j_1,j_2}^-,I)=1$.
\end{enumerate}
\end{lem}
\begin{proof}
We first prove \ref{it: Ext2 with alg 1}. As $I\neq I_{j_1,j_2}^+$, we have $d(V_{I,\Delta}^{\infty},\pi_{j_1,j_2}^{\infty})\geq 1$ by \ref{it: distance St 1} of Lemma \ref{lem: distance from St}. By Lemma~\ref{lem: simple Ext2 OS} (note that $L(1)$ is not a constituent of $M^{\widehat{j}_1}(w_{j_1,j_0})$) $\mathrm{Ext}_G^2(L(1)^\vee\otimes_E V_{I,\Delta}^{\infty},C_{\un{j}})\neq 0$ implies $\mathrm{Ext}_{U(\fg)}^1(L(w_{j_1,j_0}),L(1))\neq 0$ and $\mathrm{Ext}_{L_{\widehat{j}_1}}^1(i_{I\cap \widehat{j}_1,\widehat{j}_1}^{\infty}(J_{I,I\cap \widehat{j}_1}(V_{I,\Delta}^{\infty})),\pi_{j_1,j_2}^{\infty})^{\infty}\neq 0$. The first inequality is equivalent to $j_0=j_1$ by \ref{it: rabiotext 2} of Lemma~\ref{rabiotext} and Lemma \ref{lem: dominance control}. Arguing as in Case $1$ of the proof of Lemma \ref{lem: Ext2 factor 1} using (\ref{equ: smooth geometric lemma}), (\ref{equ: refine geometric lemma}) and that $\pi_{j_1,j_2}^{\infty}$ is in a single Bernstein block (as it is irreducible), the second inequality is equivalent to
\[\mathrm{Ext}_{G}^1(V_{I,\Delta}^{\infty}, i_{\widehat{j}_1,\Delta}^{\infty}(\pi_{j_1,j_2}^{\infty}))^{\infty}\buildrel (\ref{equ: first adjunction}) \over\cong \mathrm{Ext}_{L_{\widehat{j}_1}}^1(J_{\Delta,\widehat{j}_1}(V_{I,\Delta}^{\infty}),\pi_{j_1,j_2}^{\infty})^\infty\neq 0,\]
which is equivalent to $d(V_{I,\Delta}^{\infty},\pi_{j_1,j_2}^{\infty})= 1$ (since $d(V_{I,\Delta}^{\infty},\pi_{j_1,j_2}^{\infty})\geq 1$), which is equivalent to $I\notin[I_{j_1,j_2}^+,I_{j_1,j_2}^-]$ and $d(I,I_{j_1,j_2}^+)=1$ by \ref{it: distance St 2} of Lemma~\ref{lem: distance from St}. The proof of \ref{it: Ext2 with alg 2} is similar.
\end{proof}

\begin{lem}\label{lem: Ext2 with St}
Let $\un{j}\in\mathbf{J}$.
\begin{enumerate}[label=(\roman*)]
\item \label{it: Ext2 with St 1} If $(j_1,j_2)\neq (1,n)$ and $\un{j}\ne (2,2,n)$, then we have $\mathrm{Ext}_G^2(L(1)^\vee\otimes_E \mathrm{St}_n^{\infty},C_{\un{j}})=0$.
\item \label{it: Ext2 with St 2} If $(j_1,j_2)\neq (1,1)$ and $\un{j}\ne (2,2,2)$, then we have $\mathrm{Ext}_G^2(C_{\un{j}},L(1)^\vee\otimes_E \mathrm{St}_n^{\infty})=0$.
\end{enumerate}
\end{lem}
\begin{proof}
It follows from Lemma~\ref{lem: explicit smooth induction} that $I_{j_1,j_2}^+\neq \emptyset$ if and only if $(j_1,j_2)\neq (1,n)$, $I_{j_1,j_2}^-\neq \emptyset$ if and only if $(j_1,j_2)\neq (1,1)$, and $\emptyset\notin[I_{j_1,j_2}^+,I_{j_1,j_2}^-]$ (equivalently $I_{j_1,j_2}^+\cap I_{j_1,j_2}^-\ne \emptyset$) if and only if $j_1>1$.

If $(j_1,j_2)\neq (1,n)$ and $\un{j}\ne (2,2,n)$, we have either $j_1=1$ in which case $\emptyset \in [I_{j_1,j_2}^+,I_{j_1,j_2}^-]$, or $j_1>2$ in which case $d(\emptyset,I_{j_1,j_2}^+)=\#I_{j_1,j_2}^+>1$, or $j_1=2$, $j_2<n$ in which case $d(\emptyset,I_{j_1,j_2}^+)>1$ (again), or $(j_1,j_2)=(2,n)$ in which case $j_0\ne j_1$ by assumption. In all these cases \ref{it: Ext2 with alg 1} of Lemma~\ref{lem: Ext2 factor 2} implies $\mathrm{Ext}_G^2(L(1)^\vee\otimes_E \mathrm{St}_n^{\infty},C_{\un{j}})=0$, which gives \ref{it: Ext2 with St 1}.

If $(j_1,j_2)\neq (1,1)$ and $\un{j}\ne (2,2,2)$, we have either $j_1=1$ in which case $\emptyset \in [I_{j_1,j_2}^+,I_{j_1,j_2}^-]$, or $j_1>2$ in which case $d(\emptyset,I_{j_1,j_2}^-)=\#I_{j_1,j_2}^->1$, or $j_1=2$, $j_2>2$ in which case $d(\emptyset,I_{j_1,j_2}^-)>1$ (again), or $(j_1,j_2)=(2,2)$ in which case $j_0\ne j_1$ by assumption. In all these cases \ref{it: Ext2 with alg 2} of Lemma~\ref{lem: Ext2 factor 2} implies $\mathrm{Ext}_G^2(C_{\un{j}},L(1)^\vee\otimes_E \mathrm{St}_n^{\infty})=0$, which gives \ref{it: Ext2 with St 2}.
\end{proof}

Let $\Gamma_{\!\rm{OS}}$ be the set of pairs $(x,\pi^{\infty})$ with $x\in W(G)$ and $\pi^{\infty}$ an isomorphism class of $G$-basic representation in $\mathrm{Rep}^{\infty}_{\rm{adm}}(L_{I_x})$ where $I_x\defeq \Delta\setminus D_L(x)$ (``OS'' for ``Orlik-Strauch''). For $(x,\pi^{\infty})\in \Gamma_{\!\rm{OS}}$ we write $V_{x,\pi^{\infty}}\defeq \cF_{P_{I_x}}^{G}(L(x),\pi^{\infty})$ (which is not assumed to be irreducible). We consider a finite length object $V$ in $\mathrm{Rep}^{\rm{an}}_{\rm{adm}}(G)$ equipped with a decreasing filtration (for some $d\geq 0$)
\begin{equation}\label{equ: OS cube filtration}
0=\mathrm{Fil}^{d+1}(V)\subsetneq \mathrm{Fil}^d(V)\subsetneq \cdots\subsetneq \mathrm{Fil}^1(V) \subsetneq \mathrm{Fil}^0(V)=V,
\end{equation}
and a finite subset $\Gamma(V)\subseteq \Gamma_{\!\rm{OS}}$ equipped with a partition $\Gamma(V)=\bigsqcup_{k=0}^d\Gamma_k(V)$, such that there is an isomorphism for $0\leq k\leq d$
\begin{equation}\label{equ: OS cube graded}
0\neq \mathrm{gr}^k(V)\defeq \mathrm{Fil}^k(V)/\mathrm{Fil}^{k+1}(V)\cong \!\!\bigoplus_{(x,\pi^{\infty})\in \Gamma_k(V)}V_{x,\pi^{\infty}}.
\end{equation}
We will \emph{always} assume that the pairs $(x,W(L_{I_x})\cdot\cJ(\pi^{\infty}))$ are distinct for different choices of $(x,\pi^{\infty})\in \Gamma(V)$. Under this assumption, Lemma~\ref{lem: Hom OS} and Lemma \ref{lem: Jacquet of PS} imply that $V$ is multiplicity free and that the set $\Gamma(V)$ is uniquely determined by $V$ (but not necessarily the partition $\bigsqcup_{k=0}^d\Gamma_k(V)$ of $\Gamma(V)$ nor the filtration (\ref{equ: OS cube filtration}) of $V$). In particular, the constituents of the subquotient $V_{x,\pi^{\infty}}$ of $V$ are the constituents of $V$ of the form $\cF_{P_{I_x}}^{G}(L(x),\sigma^{\infty})$ with $\sigma^{\infty}\in\cB^{I_x}_{W(L_{I_x})\cdot\cJ(\pi^{\infty})}$ and two distinct $(x,\pi^{\infty})$ in $\Gamma(V)$ lead to two multiplicity free $V_{x,\pi^{\infty}}$ which have no constituent in commun.\bigskip

For $V$ fixed as above and $V'$ a subquotient of $V$, we define $\mathrm{Fil}^k(V')\subseteq V'$ for $k\in \{0,\dots, d+1\}$ as the maximal (for inclusion) subrepresentation of $V'$ such that its constituents are constituents of $\mathrm{Fil}^k(V)$ (recall that $V$ is multiplicity free of finite length). We say that a subquotient $V'$ of $V$ is a \emph{basic subquotient} of $V$ if there exists $(x,\pi^{\infty})\in\Gamma(V)$ such that $V'\cong V_{x,\pi^{\infty}}$, and is a \emph{good subquotient} of $V$ if there exists $\Gamma_k(V')\subseteq \Gamma_k(V)$ for every $k\in \{0,\dots, d\}$ such that
\begin{equation*}
\mathrm{gr}^k(V')=\mathrm{Fil}^k(V')/\mathrm{Fil}^{k+1}(V')\cong \bigoplus_{(x,\pi^{\infty})\in \Gamma_k(V')}V_{x,\pi^{\infty}}.
\end{equation*}
In particular every subquotient of $V$ is good if the $V_{x,\pi^{\infty}}$ are irreducible for every $(x,\pi^{\infty})\in \Gamma(V)$. (Note that ``basic'' here is relative to a fixed representation $V$, and is quite different from ``$G$-basic'' as in Definition \ref{def: basic rep}.)\bigskip

The partial order on $\mathrm{JH}_{G}(V)$ induces a partial order on $\Gamma(V)$ as follows. Given two distinct $(x_0,\pi_0^{\infty}), (x_1,\pi_1^{\infty})\in \Gamma(V)$, we write $(x_1,\pi_1^{\infty})<(x_0,\pi_0^{\infty})$ if there exists $W_i\in \mathrm{JH}_{G}(V_{x_i,\pi_i^{\infty}})\subseteq \mathrm{JH}_{G}(V)$ for $i=0,1$ such that $W_1<W_0$ in the sense of \S\ref{generalnotation}. If $d_i\in \{0,\dots,d\}$ is such that $(x_i,\pi_i^{\infty})\in \Gamma_{d_i}(V)$, we note from (\ref{equ: OS cube filtration}) and (\ref{equ: OS cube graded}) that $(x_1,\pi_1^{\infty})<(x_0,\pi_0^{\infty})$ implies $d_0<d_1$. If there does not exist $(x,\pi^{\infty})\in \Gamma(V)$ such that $(x_1,\pi_1^{\infty})<(x,\pi^{\infty})<(x_0,\pi_0^{\infty})$, then there is a good subquotient $V'$ of $V$ that fits into a
non-split short exact sequence
\[0\rightarrow V_{x_1,\pi_1^{\infty}} \rightarrow V' \rightarrow V_{x_0,\pi_0^{\infty}} \rightarrow 0.\]
The following definition is reminiscent of (though different from) Definition \ref{def: Ext hypercube}.

\begin{defn}\label{def: OS cube}
Let $d$ be some integer in $\Z_{\geq 0}$.
\begin{enumerate}[label=(\roman*)]
\item \label{it: OS cube 1}
A finite length multiplicity free representation $V$ as in (\ref{equ: OS cube filtration}) and (\ref{equ: OS cube graded}) is an \emph{$\mathrm{Ext}$-hypercube} of rank $d$ if the following properties hold
\begin{itemize}
\item the partially ordered set $\Gamma(V)$ admits a unique maximal element which is (necessarily) in $\Gamma_0(V)$ and a unique minimal element which is (necessarily) in $\Gamma_d(V)$;
\item for $0\leq d_0,d_1\leq d$ and $(x_i,\pi_i^{\infty})\in\Gamma_{d_i}(V)$, $i=0,1$, we have
\begin{equation}\label{equ: OS cube adjacent}
\mathrm{Ext}_{G}^1(V_{x_0,\pi_0^{\infty}}, V_{x_1,\pi_1^{\infty}})\neq 0
\end{equation}
\emph{if and only if} $d_1=d_0+1$ and $(x_1,\pi_1^{\infty})<(x_0,\pi_0^{\infty})$, in which case (\ref{equ: OS cube adjacent}) is one dimensional.
\end{itemize}
\item 
An $\mathrm{Ext}$-hypercube is an \emph{$\mathrm{Ext}$-square} if $d=2$, and an \emph{$\mathrm{Ext}$-cube} if $d=3$.
\item \label{it: OS cube 3}
An $\mathrm{Ext}$-hypercube $V$ is \emph{strict} if it has simple socle and cosocle, in which case $\mathrm{soc}_{G}(V)=\mathrm{soc}_{G}(\mathrm{gr}^d(V))$ and $\mathrm{cosoc}_{G}(V)=\mathrm{cosoc}_{G}(\mathrm{gr}^0(V))$.
\item \label{it: OS cube 4}
An $\mathrm{Ext}$-hypercube $V$ of rank $d$ is \emph{minimal} if for any good subquotient $V'$ of $V$ which is an $\mathrm{Ext}$-hypercube of rank $d'\leq d$, there does not exist an $\mathrm{Ext}$-hypercube $V''$ of rank $d'$ such that $\mathrm{gr}^0(V'')\cong \mathrm{gr}^0(V')$, $\mathrm{gr}^{d'}(V'')\cong \mathrm{gr}^{d'}(V')$ and $\mathrm{gr}^k(V'')\subseteq \mathrm{gr}^k(V')$ is a good direct summand of $\mathrm{gr}^k(V')$ for $1\leq k\leq d'-1$ with at least one inclusion being strict.
\end{enumerate}
\end{defn}

If $V$ is an $\mathrm{Ext}$-hypercube of rank $d$ the conditions in \ref{it: OS cube 1} imply that both $\Gamma_0(V)$, $\Gamma_d(V)$ are singletons. In fact, it is not difficult to see that the filtration (\ref{equ: OS cube filtration}) is uniquely determined by (the isomorphism class of) $V$, and that the second condition in \ref{it: OS cube 1} above implies that, for $0\leq d_0,d_1\leq d$ and $(x_i,\pi_i^{\infty})\in\Gamma_{d_i}(V)$ ($i=0,1$) such that $d_1>d_0+1$ and $(x_1,\pi_1^{\infty})<(x_0,\pi_0^{\infty})$, there exists $(y_k,\sigma_k^{\infty})\in\Gamma_k(V)$ for $d_0\leq k\leq d_1$ such that $(y_{d_i},\sigma_{d_i}^{\infty})=(x_i,\pi_i^{\infty})$ for $i=0,1$ and $(y_{\ell},\sigma_{\ell}^{\infty})<(y_{\ell+1},\sigma_{\ell+1}^{\infty})$ for $d_0\leq \ell\leq d_1-1$. Moreover for $(x_0,\pi_0^{\infty}),(x_1,\pi_1^{\infty})\in \Gamma(V)$ such that $(x_1,\pi_1^{\infty})<(x_0,\pi_0^{\infty})$, one easily checks that $V$ admits a unique good subquotient $V_{x_1,\pi_1^{\infty}}^{x_0,\pi_0^{\infty}}$ with $\Gamma(V_{x_1,\pi_1^{\infty}}^{x_0,\pi_0^{\infty}})$ consisting exactly of those $(y,\sigma^{\infty})\in\Gamma(V)$ such that $(x_1,\pi_1^{\infty})\leq (y,\sigma^{\infty})\leq (x_0,\pi_0^{\infty})$, and that $V_{x_1,\pi_1^{\infty}}^{x_0,\pi_0^{\infty}}$ is an $\mathrm{Ext}$-hypercube. Conversely any good subquotient of $V$ which is an $\mathrm{Ext}$-hypercube has the form $V_{x_1,\pi_1^{\infty}}^{x_0,\pi_0^{\infty}}$ for some $(x_1,\pi_1^{\infty})\leq (x_0,\pi_0^{\infty})\in \Gamma(V)$.

\begin{lem}\label{lem: minimal vers Ext}
Let $V$ be an $\mathrm{Ext}$-hypercube of rank $d$. Assume that there does not exist another $\mathrm{Ext}$-hypercube $V'$ of rank $d$ such that $\mathrm{gr}^0(V')\cong \mathrm{gr}^0(V)$, $\mathrm{gr}^d(V')\cong \mathrm{gr}^d(V)$ and $\mathrm{gr}^k(V')\subseteq \mathrm{gr}^k(V)$ is a good direct summand for $1\leq k\leq d-1$ with at least one inclusion being strict. Then $\mathrm{Ext}_{G}^1(\mathrm{gr}^0(V),\mathrm{Fil}^1(V))$ is one dimensional.
\end{lem}
\begin{proof}
We can assume $d\geq 2$. The existence of $V$ forces $\mathrm{Ext}_{G}^1(\mathrm{gr}^0(V),\mathrm{Fil}^1(V))\ne 0$. Assume on the contrary that $\mathrm{Ext}_{G}^1(\mathrm{gr}^0(V),\mathrm{Fil}^1(V))$ has dimension $\geq 2$. Choose a good subrepresentation $V_1\subseteq \mathrm{Fil}^1(V)$ such that $\mathrm{Fil}^1(V)/V_1$ is a basic (hence non-zero) subquotient of $V$, then we have $\Dim_E \mathrm{Ext}_{G}^1(\mathrm{gr}^0(V),\mathrm{Fil}^1(V)/V_1)=1$ by the second condition in \ref{it: OS cube 1} of Definition~\ref{def: OS cube}, which by d\'evissa\-ge from $0\rightarrow V_1\rightarrow \mathrm{Fil}^1(V)\rightarrow \mathrm{Fil}^1(V)/V_1\rightarrow 0$ forces $\mathrm{Ext}_{G}^1(\mathrm{gr}^0(V),V_1)\neq 0$. Hence, there exists $V''$ that fits into a non-split short exact sequence $0\rightarrow V_1\rightarrow V''\rightarrow \mathrm{gr}^0(V)\rightarrow 0$. We equip $V''$ with the filtration $\mathrm{Fil}^0(V'')\defeq V'$ and $\mathrm{Fil}^k(V'')\defeq \mathrm{Fil}^k(V)\cap V_1$ for $k\geq 1$. We now define $V'$ as the minimal length good subrepresentation of $V''$ such that $\mathrm{gr}^0(V')=\mathrm{gr}^0(V'')=\mathrm{gr}^0(V)$ (note that $V'$ can be strictly smaller than $V''$ since $\mathrm{gr}^0(V'')$ may have more than one maximal element in $\Gamma(V'')$). Then the conditions in \ref{it: OS cube 1} of Definition~\ref{def: OS cube} for $V$ imply the similar conditions for $V'$, in particular $V'$ is an $\mathrm{Ext}$-hypercube of rank $d$, which contradicts the minimality of $V$ as $\mathrm{gr}^1(V')\subsetneq \mathrm{gr}^1(V)$.
\end{proof}

The following formal lemma gives a rigidity property of $\mathrm{Ext}$-hypercubes.

\begin{lem}\label{lem: minimal is unique}
Let $V$ be a minimal $\mathrm{Ext}$-hypercube and $V'$ a finite length multiplicity free representation in $\mathrm{Rep}^{\rm{an}}_{\rm{adm}}(G)$. If $\mathrm{JH}_{G}(V)= \mathrm{JH}_{G}(V')$ as partially ordered sets, then $V\cong V'$.
\end{lem}
\begin{proof}
Let $d\geq 0$ be the rank of $V$. If $d=0$, Lemma~\ref{lem: OS basic unique} implies $V\cong V'$. We assume from now $d\geq 1$ and prove $V\cong V'$ by an increasing induction on $d\geq 0$. The equality $\mathrm{JH}_{G}(V)=\mathrm{JH}_{G}(V')$ as partially ordered sets implies that $V'$ admits a unique quotient $V_0'$ such that $\mathrm{JH}_{G}(V_0')=\mathrm{JH}_{G}(\mathrm{gr}^0(V))$ as partially ordered sets, which by Lemma~\ref{lem: OS basic unique} implies $V_0'\cong \mathrm{gr}^0(V)$. For $k\geq 1$ and $(x,\pi^{\infty})\in\Gamma_k(V)$, $V$ admits a unique good subrepresentation $\tld{V}_{x,\pi^{\infty}}$ such that $\Gamma(\tld{V}_{x,\pi^{\infty}})$ consists of all pairs $(w,\sigma^{\infty})\in\Gamma(V)$ such that $(w,\sigma^{\infty})\leq (x,\pi^{\infty})$. It is clear from \ref{it: OS cube 1} of Definition~\ref{def: OS cube} that $\tld{V}_{x,\pi^{\infty}}$ is itself an $\mathrm{Ext}$-hypercube of rank $d-k\geq 0$. Let $S\subseteq \Gamma(V)\setminus \Gamma_0(V)$ be a subset such that any $(w,\sigma^{\infty})\in\Gamma(V)$ satisfying $(w,\sigma^{\infty})\leq (x,\pi^{\infty})$ for some $(x,\pi^{\infty})\in S$ also satisfies $(w,\sigma^{\infty})\in S$ and define the good subrepresentation $\mathrm{Fil}_{S}(V)\defeq \bigcup_{(x,\pi^{\infty})\in S}\tld{V}_{x,\pi^{\infty}}$.
The equality $\mathrm{JH}_{G}(V)=\mathrm{JH}_{G}(V')$ as partially ordered sets implies that $V'$ admits unique subrepresentations $\tld{V}_{x,\pi^{\infty}}'$ for $(x,\pi^{\infty})\in \Gamma(V)\setminus \Gamma_0(V)$ and $\mathrm{Fil}_S(V')$ for $S\subseteq \Gamma(V)$ as above such that $\mathrm{JH}_{G}(\tld{V}_{x,\pi^{\infty}})=\mathrm{JH}_{G}(\tld{V}_{x,\pi^{\infty}}')$ and $\mathrm{JH}_{G}(\mathrm{Fil}_{S}(V))=\mathrm{JH}_{G}(\mathrm{Fil}_S(V'))$ as partially ordered sets. By induction we have in $\mathrm{Rep}^{\rm{an}}_{\rm{adm}}(G)$ for $(x,\pi^{\infty})\in \Gamma(V)\setminus \Gamma_0(V)$
\begin{equation}\label{equ: minimal OS induction}
\tld{V}_{x,\pi^{\infty}}\cong \tld{V}_{x,\pi^{\infty}}'.
\end{equation}
We prove $\mathrm{Fil}_{S}(V)\cong\mathrm{Fil}_S(V')$ by induction on $\#S\geq 1$. Let $(x,\pi^{\infty})$ be a maximal element of $S$, $S'\defeq S\setminus \{(x,\pi^{\infty})\}$ and $S''$ the set of $(w,\sigma^{\infty})\in \Gamma(V)$ such that $(w,\sigma^{\infty})\leq (x,\pi^{\infty})$. If $S''=S$ (i.e.~$(x,\pi^{\infty})$ is the unique maximal element of $S$), then $\mathrm{Fil}_{S}(V)\cong \tld{V}_{x,\pi^{\infty}}\cong \tld{V}_{x,\pi^{\infty}}' \cong \mathrm{Fil}_{S}(V')$ by (\ref{equ: minimal OS induction}). Otherwise we have $S',S''\subsetneq S$ which implies $\mathrm{Fil}_{S'}(V)\cong \mathrm{Fil}_{S'}(V')$, $\mathrm{Fil}_{S''}(V)\cong \mathrm{Fil}_{S''}(V')$ and $\mathrm{Fil}_{S'\cap S''}(V)\cong \mathrm{Fil}_{S'\cap S''}(V')$ by induction on $S$. It follows from the first condition in \ref{it: OS cube 1} of Definition~\ref{def: OS cube} (applied to $V$) that $\mathrm{Fil}_{S'}(V)$, $\mathrm{Fil}_{S''}(V)$ and $\mathrm{Fil}_{S'\cap S''}(V)$ are indecomposable representations, and hence that $\mathrm{Fil}_{S}(V)$ is the amalgamate sum of $\mathrm{Fil}_{S'}(V)$ and $\mathrm{Fil}_{S''}(V)$ over $\mathrm{Fil}_{S'\cap S''}(V)$ (recall all representations are multiplicity free). Likewise $\mathrm{Fil}_{S'}(V')$, $\mathrm{Fil}_{S''}(V')$, $\mathrm{Fil}_{S'\cap S''}(V')$ are indecomposable and $\mathrm{Fil}_{S}(V')$ is the amalgamate sum of $\mathrm{Fil}_{S'}(V')$ and $\mathrm{Fil}_{S''}(V')$ over $\mathrm{Fil}_{S'\cap S''}(V')$. Hence we deduce $\mathrm{Fil}_{S}(V)\cong \mathrm{Fil}_{S}(V')$. Now take $S=\Gamma(V)\setminus \Gamma_0(V)$, so that $\mathrm{Fil}_S(V)=\mathrm{Fil}^1(V)$. From $\mathrm{JH}_{G}(V)= \mathrm{JH}_{G}(V')$ and the previous paragraph $V'$ fits into a non-split short exact sequence
\begin{equation}\label{equ: minimal OS ext}
0\longrightarrow \mathrm{Fil}^1(V)\cong \mathrm{Fil}_S(V')\longrightarrow V' \longrightarrow V_0'\cong \mathrm{gr}^0(V)\longrightarrow 0.
\end{equation}
The minimality of $V$ and Lemma~\ref{lem: minimal vers Ext} imply $\Dim_E\mathrm{Ext}_{G}^1(\mathrm{gr}^0(V),\mathrm{Fil}^1(V))=1$, and thus it follows from (\ref{equ: minimal OS ext}) that we must have $V'\cong V$.
\end{proof}

\begin{lem}\label{lem: OS minimal criterion}
Let $V$ be an $\mathrm{Ext}$-square. Assume that there exist good subrepresentations $\mathrm{gr}^2(V)\subseteq V_2\subseteq V_1\subseteq \mathrm{Fil}^1(V)$ such that $V_1/V_2$ is a basic subquotient of $V$ and $\mathrm{Ext}_{G}^1(V/V_1,V_2)\!=0$.
Then $V$ is minimal.
\end{lem}
\begin{proof}
We fix throughout the proof $V_1,V_2$ as in the statement. Let $V_2'$ and $V_1'$ be good subrepresentations of $V$ such that $\mathrm{gr}^2(V)\subseteq V_2' \subseteq V_2\subseteq V_1\subseteq V_1' \subseteq \mathrm{Fil}^1(V)$. As $V$ is multiplicity free, the injection $V_2'\hookrightarrow V_2$ and the surjection $V_1\twoheadrightarrow V_1'$ induce an injection $\mathrm{Ext}_{G}^1(V/V_1',V_2')\hookrightarrow \mathrm{Ext}_{G}^1(V/V_1,V_2)$ which implies
\begin{equation}\label{equ: OS minimal square Ext 0}
\mathrm{Ext}_{G}^1(V/V_1',V_2')=0.
\end{equation}

\textbf{Step $1$}: We prove that $\dim_E\mathrm{Ext}_{G}^1(\mathrm{gr}^0(V), \mathrm{Fil}^1(V))=1$.\\
The existence of the $\mathrm{Ext}$-square $V$ ensures $\mathrm{Ext}_{G}^1(\mathrm{gr}^0(V), \mathrm{Fil}^1(V))\ne 0$. As $V_1/V_2$ is a basic direct summand of $\mathrm{gr}^1(V)$ by (\ref{equ: OS cube graded}), we have a surjection $q: \mathrm{Fil}^1(V)\twoheadrightarrow V_1/V_2$ which induces an exact sequence
\begin{equation}\label{equ: OS minimal square Ext 2}
\mathrm{Ext}_{G}^1(\mathrm{gr}^0(V), \mathrm{ker}(q))\rightarrow \mathrm{Ext}_{G}^1(\mathrm{gr}^0(V), \mathrm{Fil}^1(V))\rightarrow \mathrm{Ext}_{G}^1(\mathrm{gr}^0(V), V_2/V_1).
\end{equation}
Assume $\dim_E\mathrm{Ext}_{G}^1(\mathrm{gr}^0(V), \mathrm{Fil}^1(V))\geq 2$. Then (\ref{equ: OS minimal square Ext 2}) and $\Dim_E \mathrm{Ext}_{G}^1(\mathrm{gr}^0(V), V_2/V_1)=1$ (see the second condition in \ref{it: OS cube 1} of Definition~\ref{def: OS cube}) imply $\mathrm{Ext}_{G}^1(\mathrm{gr}^0(V), \mathrm{ker}(q))\neq 0$. Similar arguments as in the proof of Lemma~\ref{lem: minimal vers Ext} with $V_1$ there replaced by $\mathrm{ker}(q)$ then show that there exists an $\mathrm{Ext}$-square $V'$ with $\mathrm{gr}^0(V')=\mathrm{gr}^0(V)$, $\mathrm{gr}^2(V')=\mathrm{gr}^2(V)$ and $\mathrm{gr}^1(V')$ a good direct summand of $\mathrm{ker}(q)/\mathrm{gr}^2(V)\cong \mathrm{gr}^1(V)/(V_1/V_2)$. But such a $V'$ necessarily fits into a non-split extension $0\rightarrow V_2'\rightarrow V' \rightarrow V/V_1'\rightarrow 0$ for some good subrepresentations $V_1',V_2'$ of $V$ such that $\mathrm{gr}^2(V)\subseteq V_2' \subseteq V_2\subseteq V_1\subseteq V_1' \subseteq \mathrm{Fil}^1(V)$, which contradicts (\ref{equ: OS minimal square Ext 0}).\bigskip

\textbf{Step $2$}: We prove that $V$ is minimal.\\
By \ref{it: OS cube 4} of Definition \ref{def: OS cube}, it suffices to show that there does not exist an $\mathrm{Ext}$-square $V'$ such that $\mathrm{gr}^0(V')=\mathrm{gr}^0(V)$, $\mathrm{gr}^2(V')=\mathrm{gr}^2(V)$ and $\mathrm{gr}^1(V')$ is a proper good direct summand of $\mathrm{gr}^1(V)$.
Assume on the contrary that such $V'$ exists. Let $V_0$ be a basic subquotient of $V$ which is a direct summand of $\mathrm{gr}^1(V)/\mathrm{gr}^1(V')$ and fix a surjection $q': \mathrm{Fil}^1(V)\twoheadrightarrow V_0$. Using that $V$ is multiplicity free, the injection $\mathrm{Fil}^1(V')\subseteq \mathrm{ker}(q')$ induces an injection
\begin{equation}\label{equ: OS minimal square Ext 3}
0\neq \mathrm{Ext}_{G}^1(\mathrm{gr}^0(V), \mathrm{Fil}^1(V')) \rightarrow \mathrm{Ext}_{G}^1(\mathrm{gr}^0(V), \mathrm{ker}(q')),
\end{equation}
and the short exact sequence $0\rightarrow \mathrm{ker}(q')\rightarrow \mathrm{Fil}^1(V)\rightarrow V_0\rightarrow 0$ induces an exact sequence
\begin{equation}\label{equ: OS minimal square Ext 4}
0\rightarrow \mathrm{Ext}_{G}^1(\mathrm{gr}^0(V), \mathrm{ker}(q'))\rightarrow \mathrm{Ext}_{G}^1(\mathrm{gr}^0(V), \mathrm{Fil}^1(V))\rightarrow \mathrm{Ext}_{G}^1(\mathrm{gr}^0(V), V_0).
\end{equation}
By \ref{it: OS cube 1} of Definition~\ref{def: OS cube} $\dim_E\mathrm{Ext}_{G}^1(\mathrm{gr}^0(V), V_0)=1$ and the last map in (\ref{equ: OS minimal square Ext 4}) is non-zero. By Step $1$ it is an isomorphism, which forces $\mathrm{Ext}_{G}^1(\mathrm{gr}^0(V), \mathrm{ker}(q'))\!=0$, contradicting~(\ref{equ: OS minimal square Ext 3}).
\end{proof}

Recall that for $\un{j}\in\mathbf{J}$ the representation $C_{\un{j}}$ is defined in (\ref{cj}).

\begin{prop}\label{prop: easy square}
Let $\un{j},\un{j}'\in\mathbf{J}$ and write $\underline j=(j_0,j_1,j_2)$, $\underline j'=(j'_0,j'_1,j'_2)$.
\begin{enumerate}[label=(\roman*)]
\item \label{it: easy square 1} If $\un{j}'=(j_0+1,j_1-1,j_2)$ there exists a unique $\mathrm{Ext}$-square $V_{\un{j},\un{j}'}$ such that
\[\left\{\begin{array}{rcl}
\mathrm{gr}^0(V_{\un{j},\un{j}'})&\cong &C_{\un{j}'}\\
\mathrm{gr}^1(V_{\un{j},\un{j}'})&\cong &C_{(j_0+1,j_1,j_2)}\oplus C_{(j_0,j_1-1,j_2)}\\
\mathrm{gr}^2(V_{\un{j},\un{j}'})&\cong &C_{\un{j}}.
\end{array}\right.\]
\item \label{it: easy square 2} If $\un{j}'=(j_0+1,j_1+1,j_2+1)$ there exists a unique $\mathrm{Ext}$-square $V_{\un{j},\un{j}'}$ such that
\[\left\{\begin{array}{rcl}
\mathrm{gr}^0(V_{\un{j},\un{j}'})&\cong &C_{\un{j}'}\\
\mathrm{gr}^1(V_{\un{j},\un{j}'})&\cong &C_{(j_0+1,j_1,j_2)}\oplus C_{(j_0,j_1+1,j_2+1)}\oplus L(1)^\vee\otimes_E \pi_{\un{j},\un{j}'}^{\infty}\\
\mathrm{gr}^2(V_{\un{j},\un{j}'})&\cong &C_{\un{j}}
\end{array}\right.\]
where $\pi_{\un{j},\un{j}'}^{\infty}\cong V_{[1,j_1],\Delta}^{\infty}$ if $j_0=j_1=j_2$, $\pi_{\un{j},\un{j}'}^{\infty}$ is the unique $G$-basic length $2$ representation of $G$ with socle $V_{[j_2-j_1+1,j_2],\Delta}^{\infty}$ and cosocle $V_{[j_2-j_1,j_2],\Delta}^{\infty}$ if $j_0=j_1<j_2$ (first statement in Lemma \ref{lem: length two trivial block}), and $\pi_{\un{j},\un{j}'}^{\infty}$ is zero otherwise.
\item \label{it: easy square 3} If $j_1=j_0+1$ there exists a unique $\mathrm{Ext}$-square $V_{\un{j},\infty}$ such that
\[\left\{\begin{array}{rcl}
\mathrm{gr}^0(V_{\un{j},\infty})&\cong &L(1)^\vee\otimes_E \pi_{\un{j},\infty}^{\infty}\\
\mathrm{gr}^1(V_{\un{j},\infty})&\cong &C_{(j_0,j_0,j_2)}\oplus C_{(j_0+1,j_0+1,j_2)}\\
\mathrm{gr}^2(V_{\un{j},\infty})&\cong &C_{\un{j}}
\end{array}\right.\]
where \ $\pi_{\un{j},\infty}^{\infty}$ \ is \ the \ unique \ $G$-basic \ length \ two \ representation \ of \ $G$ \ with \ socle $V_{[j_2-j_0+1,j_2]\cap\Delta,\Delta}^{\infty}$ and cosocle $V_{[j_2-j_0,j_2]\cap\Delta,\Delta}^{\infty}$ (Lemma \ref{lem: length two trivial block}).
\item \label{it: easy square 4} If $j_1=j_0-1$ there exists a unique $\mathrm{Ext}$-square $V_{\infty,\un{j}}$ such that
\[\left\{\begin{array}{rcl}
\mathrm{gr}^0(V_{\infty,\un{j}})&\cong &C_{\un{j}}\\
\mathrm{gr}^1(V_{\infty,\un{j}})&\cong &C_{(j_0,j_0,j_2)}\oplus C_{(j_0-1,j_0-1,j_2)}\\
\mathrm{gr}^2(V_{\infty,\un{j}})&\cong &L(1)^\vee\otimes_E \pi_{\infty,\un{j}}^{\infty}
\end{array}\right.\]
where $\pi_{\infty,\un{j}}^{\infty}\cong V_{[1,j_2-1],\Delta}^{\infty}$ if $j_1+1=j_0=j_2$, and $\pi_{\infty,\un{j}}^{\infty}$ is the unique $G$-basic length $2$ representation of $G$ with socle $V_{[j_2-j_0+1,j_2-1],\Delta}^{\infty}$ and cosocle $V_{[j_2-j_0,j_2-1],\Delta}^{\infty}$ (Lemma \ref{lem: length two trivial block}) if $j_1+1=j_0<j_2$.
\end{enumerate}
Moreover, all $\mathrm{Ext}$-squares above are minimal and strict (see Definition \ref{def: OS cube}).
\end{prop}
\begin{proof}
We make crucial use of the Ext-squares of $U(\fg)$-modules constructed in \S\ref{U(g)square}. We divide the proof into two steps.\bigskip

\textbf{Step $1$}: We construct $I\subseteq \Delta$, $Q$ in $\cO^{\fp_I}_{\rm{alg}}$ (see \S\ref{subsec: category}) and a $G$-regular irreducible $\pi^{\infty}$ in $\mathrm{Rep}^{\infty}_{\rm{adm}}(L_I)$ such that $\cF_{P_I}^G(Q,\pi^{\infty})$ contains a unique subquotient of the form $V_{\un{j},\un{j}'}$ (resp.~$V_{\un{j},\infty}$, $V_{\infty,\un{j}}$) which is a strict $\mathrm{Ext}$-square as described in \ref{it: easy square 1} and \ref{it: easy square 2} (resp.~\ref{it: easy square 3}, \ref{it: easy square 4}).\\
We define $I\defeq \Delta\setminus\{j_1,j_1'\}=\widehat{j}_1\cap \widehat{j}_1'$ in \ref{it: easy square 1} and \ref{it: easy square 2}, and $I\defeq \Delta\setminus\{j_0,j_1\}=\widehat{j}_0\cap \widehat{j}_1$ in \ref{it: easy square 3} and \ref{it: easy square 4}. In case \ref{it: easy square 1} and \ref{it: easy square 2} we set $Q\defeq Q_1(w_{j_1,j_0},w_{j_1',j_0'})$ using Proposition~\ref{lem: explicit g square}. In case \ref{it: easy square 3} we set $Q\defeq Q_1(w_{j_1,j_0},1)$ and in case \ref{it: easy square 4} we set $Q\defeq Q_1(1,w_{j_1,j_0})$ using Remark~\ref{rem: other simple square}. All these $\mathrm{Ext}$-squares of $U(\fg)$-modules are actually in $\cO^{\fp_I}_{\rm{alg}}$ using Lemma \ref{lem: dominance and left set} and \cite[Prop.~9.3(c)]{Hum08}.

We consider the $W(L_I)$-coset $\Sigma$ given by $\Sigma_{j_1,j_2}\cap \Sigma_{j_1',j_2'}$ in \ref{it: easy square 1} and \ref{it: easy square 2}, and by $\Sigma_{j_1,j_2}\cap \Sigma_{j_0,j_2}$ in \ref{it: easy square 3} and \ref{it: easy square 4}. Thanks to \ref{it: sm cell 1} of Lemma~\ref{lem: sm cell} and \ref{it: Hom vers induction1} of Lemma~\ref{lem: Hom vers induction}, we can define a $G$-regular irreducible smooth representation $\pi^{\infty}$ of $L_I$ in \ref{it: easy square 1} and \ref{it: easy square 2} as
\begin{equation}\label{equ: easy square smooth part 1 2}
\pi^{\infty}\defeq J_{\widehat{j}_1',I}(\pi_{j_1',j_2'}^{\infty})_{\cB^{I}_{\Sigma}}\cong J_{\widehat{j}_1,I}'(\pi_{j_1,j_2}^{\infty})_{\cB^{I}_{\Sigma}},
\end{equation}
and note that $\pi_{j_1,j_2}^{\infty}$ (resp.~$\pi_{j_1',j_2'}^{\infty}$) is a subquotient of $i_{I,\widehat{j}_1}^{\infty}(\pi^{\infty})$ (resp.~$i_{I,\widehat{j}_1'}^{\infty}(\pi^{\infty})$) by (\ref{equ: second adjunction}) (resp.~by (\ref{equ: first adjunction})). Similarly, using \ref{it: sm cell 1} of Lemma~\ref{lem: sm cell} (applied with $(j_1,j'_1,j_2,j'_2)$ there being $(j_1,j_0,j_2,j_2)$) and \ref{it: Hom vers induction1} of Lemma~\ref{lem: Hom vers induction}) we have in \ref{it: easy square 3}
\begin{equation*}
\pi^{\infty}\defeq J_{\widehat{j}_0,I}(\pi_{j_0,j_2}^{\infty})_{\cB^{I}_{\Sigma}}\cong J_{\widehat{j}_1,I}'(\pi_{j_1,j_2}^{\infty})_{\cB^{I}_{\Sigma}}
\end{equation*}
and using \ref{it: sm cell 1} of Lemma~\ref{lem: sm cell} (applied with $(j_1,j'_1,j_2,j'_2)$ there being $(j_0,j_1,j_2,j_2)$) and \ref{it: Hom vers induction1} of Lemma~\ref{lem: Hom vers induction}) we have in \ref{it: easy square 4}
\begin{equation}\label{equ: easy square smooth part 4}
\pi^{\infty}\defeq J_{\widehat{j}_1,I}(\pi_{j_1,j_2}^{\infty})_{\cB^{I}_{\Sigma}}\cong J_{\widehat{j}_0,I}'(\pi_{j_0,j_2}^{\infty})_{\cB^{I}_{\Sigma}}.
\end{equation}
Moreover in \ref{it: easy square 3} and \ref{it: easy square 4} $\pi_{j_1,j_2}^{\infty}$ (resp.~$\pi_{j_0,j_2}^{\infty}$) is a subquotient of $i_{I,\widehat{j}_1}^{\infty}(\pi^{\infty})$ (resp.~$i_{I,\widehat{j}_0}^{\infty}(\pi^{\infty})$) using again (\ref{equ: first adjunction}) and (\ref{equ: second adjunction}). Using Lemma~\ref{lem: explicit smooth induction} applied to $\sigma_{j_1,j_2}^{\infty}=\pi^\infty$ one easily checks that:\\
$\bullet$ in \ref{it: easy square 2} $\pi_{\un{j},\un{j}'}^{\infty}$ is a subquotient of $i_{I,\Delta}^{\infty}(\pi^{\infty})$ which is a subrepresentation of $i_{\widehat{j}_1,\Delta}^{\infty}(\pi_{j_1,j_2}^{\infty})$ and a quotient of $i_{\widehat{j}_1',\Delta}^{\infty}(\pi_{j_1',j_2'}^{\infty})$;\\
$\bullet$ in \ref{it: easy square 3} $\pi_{\un{j},\infty}^{\infty}$ is a subquotient of $i_{I,\Delta}^{\infty}(\pi^{\infty})$ which is a subrepresentation of $i_{\widehat{j}_0,\Delta}^{\infty}(\pi_{j_0,j_2}^{\infty})$ such that $\Hom_{G}(\pi_{\un{j},\infty}^{\infty},i_{\widehat{j}_1,\Delta}^{\infty}(\pi_{j_1,j_2}^{\infty}))\neq 0$;\\
$\bullet$ in \ \ref{it: easy square 4} \ $\pi_{\infty,\un{j}}^{\infty}$ \ is \ a \ subquotient \ of \ $i_{I,\Delta}^{\infty}(\pi^{\infty})$ \ which \ is \ a \ quotient \ of \ $i_{\widehat{j}_0,\Delta}^{\infty}(\pi_{j_0,j_2}^{\infty})$ \ such that $\Hom_{G}(i_{\widehat{j}_1,\Delta}^{\infty}(\pi_{j_1,j_2}^{\infty}),\pi_{\infty,\un{j}}^{\infty})~\neq~0$.\\
We see that all assumptions in Lemma~\ref{lem: general construction} are satisfied and thus $\cF_{P_I}^G(Q,\pi^{\infty})$ contains a unique subquotient of the form $V_{\un{j},\un{j}'}$ (resp.~$V_{\un{j},\infty}$, $V_{\infty,\un{j}}$) in \ref{it: easy square 1} and \ref{it: easy square 2} (resp.~in \ref{it: easy square 3}, in \ref{it: easy square 4}).
More precisely the first condition in \ref{it: OS cube 1} of Definition~\ref{def: OS cube} is obvious and the second condition for $V_{\un{j},\un{j}'}$, $V_{\un{j},\infty}$ and $V_{\infty,\un{j}}$ respectively can be checked using Lemma~\ref{lem: Ext1 factor 1}, Lemma~\ref{lem: Ext1 factor 2} and \ref{it: general construction 2} of Lemma~\ref{lem: general construction}.
In \ref{it: easy square 1} and \ref{it: easy square 2}, since $\mathrm{gr}^0(V_{\un{j},\un{j}'})$ and $\mathrm{gr}^2(V_{\un{j},\un{j}'})$ are simple, the $\mathrm{Ext}$-square $V_{\un{j},\un{j}'}$ is strict. In \ref{it: easy square 3}, as $\mathrm{gr}^0(V_{\un{j},\infty})$ is the only reducible basic subquotient of $V_{\un{j},\infty}$ and $\pi_{\un{j},\infty}^{\infty}$ is a subrepresentation of $i_{\widehat{j}_0,\Delta}^{\infty}(\pi_{j_0,j_2}^{\infty})$, we deduce from \ref{it: Ext1 OS socle} of Lemma~\ref{lem: Ext1 OS socle cosocle} (applied with $I_0=\Delta$) that $V_{\un{j},\infty}$ is strict.
Similarly, in \ref{it: easy square 4}, as $\mathrm{gr}^2(V_{\infty,\un{j}})$ is the only reducible basic subquotient of $V_{\infty,\un{j}}$ and $\pi_{\infty,\un{j}}^{\infty}$ is a quotient of $i_{\widehat{j}_0,\Delta}^{\infty}(\pi_{j_0,j_2}^{\infty})$, we deduce from \ref{it: Ext1 OS cosocle} of Lemma~\ref{lem: Ext1 OS socle cosocle} (applied with $I_1=\Delta$) that $V_{\infty,\un{j}}$ is strict.\bigskip

\textbf{Step $2$}: We prove the minimality of the $\mathrm{Ext}$-squares in the statement.\\
In each case, we choose a subrepresentation $M\subseteq Q$ and a quotient $M'$ of $Q/M$ as follows.
\begin{itemize}
\item In \ref{it: easy square 1} and \ref{it: easy square 2}, $M'$ has length $2$ with socle $L(w_{j_1,j_0+1})$ and cosocle $L(w_{j_1,j_0})$ if $j_1\leq j_0$, and $M'\defeq L(w_{j_1,j_0})$ if $j_1\geq j_0+1$.
\item In \ref{it: easy square 1}, $M$ has length $2$ with socle $L(w_{j_1-1,j_0+1})$ and cosocle $L(w_{j_1-1,j_0})$ if $j_1\geq j_0+1$, and $M\defeq L(w_{j_1-1,j_0+1})$ if $j_1\leq j_0$. In \ref{it: easy square 2}, $M$ has length $2$ with socle $L(w_{j_1+1,j_0+1})$ and cosocle $L(w_{j_1+1,j_0})$ if $j_1\geq j_0$, and $M\defeq L(w_{j_1+1,j_0+1})$ if $j_1\leq j_0-1$.
\item In \ref{it: easy square 3}, $M'\defeq L(w_{j_1,j_0})$ and $M$ has length $2$ with socle $L(1)$ and cosocle $L(s_{j_0})$.
\item In \ref{it: easy square 4}, $M'\defeq L(1)$ and $M$ has length $2$ with socle $L(w_{j_1,j_0})$ and cosocle $L(s_{j_1})$.
\end{itemize}
Since $M'$ is a proper quotient of $Q/M$, we have for $k\leq 1$
\begin{equation}\label{equ: easy square g minimal}
\mathrm{Ext}_{U(\fg)}^k(M',M)=0
\end{equation}
by minimality of the $\mathrm{Ext}$-square $Q$, see Proposition~\ref{lem: explicit g square} and Remark~\ref{rem: other simple square}. By \ref{it: general construction 2} of Lemma~\ref{lem: general construction} $M$ uniquely determines a quotient $V$ of $V_{\un{j},\un{j}'}$ in \ref{it: easy square 1} and \ref{it: easy square 2} (resp.~$V_{\un{j},\infty}$ in \ref{it: easy square 3},~$V_{\infty,\un{j}}$ in \ref{it: easy square 4}) and $M'$ uniquely determines a subrepresentation $V'$ of $V_{\un{j},\un{j}'}$ (resp.~$V_{\un{j},\infty}$, $V_{\infty,\un{j}}$). We prove the minimality of the $\mathrm{Ext}$-square $V_{\un{j},\un{j}'}$, $V_{\un{j},\infty}$ and $V_{\infty,\un{j}}$ using Lemma~\ref{lem: OS minimal criterion}. To check the assumption of Lemma~\ref{lem: OS minimal criterion}, it suffices to prove in each case above
\begin{equation}\label{equ: easy square OS minimal}
\mathrm{Ext}_{G}^1(V,V')=0.
\end{equation}
The idea is to apply (\ref{equ: OS main seq}) and use (\ref{equ: easy square g minimal}). But in order to do so, we need to replace $V'$ by a parabolic Verma module as this is part of the assumptions in (\ref{equ: OS main seq}).

\textbf{Case $2.1$}: In \ref{it: easy square 1} and \ref{it: easy square 2}, we have $V\cong \cF_{P_{\widehat{j}_1'}}^{G}(M,\pi_{j_1',j_2'}^{\infty})$ and $V'\cong \cF_{P_{\widehat{j}_1}}^{G}(M',\pi_{j_1,j_2}^{\infty})$. We first note that $M'$ is a quotient of $M(w_{j_1,j_0})$ (using for instance (\ref{equ: O Hom radical})), hence of $M^{\widehat{j}_1}(w_{j_1,j_0})$ by \cite[Thm.~9.4(c)]{Hum08}. By \emph{loc.~cit.}~and Lemma \ref{lem: dominance and left set}, $M^{\widehat{j}_1}(w_{j_1,j_0})$ and $M$ do not share any constituent, and thus $\cF_{P_{\widehat{j}_1}}^{G}(M^{\widehat{j}_1}(w_{j_1,j_0}),\pi_{j_1,j_2}^{\infty})$ and $V$ do not share any constituent by Lemma~\ref{lem: Hom OS}. Hence the surjection $M^{\widehat{j}_1}(w_{j_1,j_0})\twoheadrightarrow M'$ induces an embedding
\begin{equation}\label{equ: easy square case 1 embedding}
\mathrm{Ext}_{G}^1(V,V')\hookrightarrow \mathrm{Ext}_{G}^1\big(V,\cF_{P_{\widehat{j}_1}}^{G}(M^{\widehat{j}_1}(w_{j_1,j_0}),\pi_{j_1,j_2}^{\infty})\big).
\end{equation}
Let $M''\subseteq M^{\widehat{j}_1}(w_{j_1,j_0})$ such that $M^{\widehat{j}_1}(w_{j_1,j_0})/M''\cong M'$ and let $L(x)\in\mathrm{JH}_{U(\fg)}(M'')$. We have $x\geq w_{j_1,j_0}$ (from the structure of $M(w_{j_1,j_0})$), $D_L(x)=\{j_1\}$ (from Lemma \ref{lem: dominance and left set}) and $x\notin\{w_{j_1,j_0},w_{j_1,j_0+1}\}$ (from $L(x)\notin\mathrm{JH}_{U(\fg)}(M')$ and the fact both $L(w_{j_1,j_0})$ and $L(w_{j_1,j_0+1})$ have multiplicity $\leq 1$ in $M(w_{j_1,j_0})$ by \ref{it: coxeter mult 1} of Lemma \ref{lem: KL coxeter} for instance). Since each $L(y)\in\mathrm{JH}_{U(\fg)}(M)$ satisfies $y\in\{w_{j_1',j_0},w_{j_1',j_0+1}\}$ we deduce $\mathrm{Ext}_{U(\fg)}^1(L(x),L(y))=0$ from \ref{it: special descent 3} of Lemma~\ref{lem: special descent change}. Using this and $L(x)\notin \mathrm{JH}_{U(\fg)}(M)$ we deduce by an obvious d\'evissage on $0\rightarrow M''\rightarrow M^{\widehat{j}_1}(w_{j_1,j_0})\rightarrow M'\rightarrow 0$:
\[\mathrm{Ext}_{U(\fg)}^k(M',M)\buildrel\sim\over\longrightarrow \mathrm{Ext}_{U(\fg)}^k(M^{\widehat{j}_1}(w_{j_1,j_0}),M)\ \ \mathrm{for}\ \ k\leq 1,\]
which together with (\ref{equ: easy square g minimal}) implies $\mathrm{Ext}_{U(\fg)}^k(M^{\widehat{j}_1}(w_{j_1,j_0}),M)=0$ for $k\leq 1$. As $\Sigma_{j_1,j_2}\cap \Sigma_{j_1',j_2'}\neq \emptyset$ (\ref{it: coset 1} of Lemma~\ref{lem: coset intersection}), we can apply (\ref{equ: OS main seq}) with $w=1$ and obtain
\[\mathrm{Ext}_{G}^1(V,\cF_{P_{\widehat{j}_1}}^{G}(M^{\widehat{j}_1}(w_{j_1,j_0}),\pi_{j_1,j_2}^{\infty}))=0.\]
By (\ref{equ: easy square case 1 embedding}) this implies (\ref{equ: easy square OS minimal}).

\textbf{Case $2.2$}: In \ref{it: easy square 3}, we have $V'\cong C_{\un{j}}$ and $V$ is a subrepresentation of $\cF_{P_{\widehat{j}_0}}^{G}(M,\pi_{j_0,j_2}^{\infty})$ as $\pi_{\un{j},\infty}^{\infty}$ is a subrepresentation of $i_{\widehat{j}_0,\Delta}^{\infty}(\pi_{j_0,j_2}^{\infty})$. As $(j_0,j_2)=(j_1-1,j_2)$, we have by \ref{it: sm cell 1} of Lemma~\ref{lem: sm cell}
\[\Hom_{L_{\widehat{j}_1}}(i_{\widehat{j}_0\cap\widehat{j}_1}^{\infty}(J_{\widehat{j}_0,\widehat{j}_0\cap\widehat{j}_1}(\pi_{j_0,j_2}^{\infty})),\pi_{j_1,j_2}^{\infty})\neq 0,\]
which together with \ref{it: Hom vers induction2} of Lemma~\ref{lem: Hom vers induction} implies $J_{\widehat{j}_1}(\tau^{\infty})_{\cB^{\widehat{j}_1}_{\Sigma_{j_1,j_2}}}\!\!\!\!=0$ for any constituent $\tau^{\infty}$ of $i_{\widehat{j}_0,\Delta}^{\infty}(\pi_{j_0,j_2}^{\infty})/\pi_{\un{j},\infty}^{\infty}$. Note that $\mathrm{JH}_{U(\fg)}(M)\cap \mathrm{JH}_{U(\fg)}(M^{\widehat{j}_1}(w_{j_1,j_0}))=\emptyset$ (using $j_1>j_0$). Hence, we deduce from Remark~\ref{rem: dual minimal OS} (applied when $V_+=V$) an isomorphism
\begin{equation}\label{equ: easy square case 2 sub}
\mathrm{Ext}_{G}^1(\cF_{P_{\widehat{j}_0}}^G(M,\pi_{j_0,j_2}^{\infty}),V')\buildrel\sim\over\longrightarrow \mathrm{Ext}_{G}^1(V,V').
\end{equation}
Note that $\cF_{P_{\widehat{j}_0}}^G(M,\pi_{j_0,j_2}^{\infty})$ and $\cF_{P_{\widehat{j}_1}}^{G}(M^{\widehat{j}_1}(w_{j_1,j_0}),\pi_{j_1,j_2}^{\infty})$ also have no common constituent by Lemma~\ref{lem: Hom OS}. So the surjection $M^{\widehat{j}_1}(w_{j_1,j_0})\twoheadrightarrow L(w_{j_1,j_0})$ induces an embedding
\begin{equation}\label{equ: easy square case 2 embedding}
\mathrm{Ext}_{G}^1(\cF_{P_{\widehat{j}_0}}^G(M,\pi_{j_0,j_2}^{\infty}),V')\hookrightarrow \mathrm{Ext}_{G}^1(\cF_{P_{\widehat{j}_0}}^G(M,\pi_{j_0,j_2}^{\infty}),\cF_{P_{\widehat{j}_1}}^{G}(M^{\widehat{j}_1}(w_{j_1,j_0}),\pi_{j_1,j_2}^{\infty})).
\end{equation}
By a similar argument as in Case $2.1$ using \ref{it: special descent 3} of Lemma~\ref{lem: special descent change} and (\ref{equ: easy square g minimal}), one shows $\mathrm{Ext}_{U(\fg)}^k(M^{\widehat{j}_1}(w_{j_1,j_0}),M)=0$ for $k\leq 1$. As $\Sigma_{j_1,j_2}\cap \Sigma_{j_0,j_2}\neq \emptyset$ (from \ref{it: coset 1} of Lemma~\ref{lem: coset intersection} and $j_0=j_1-1$ here) we can apply (\ref{equ: OS main seq}) with $w=1$ and deduce
\[\mathrm{Ext}_{G}^1(\cF_{P_{\widehat{j}_0}}^G(M,\pi_{j_0,j_2}^{\infty}),\cF_{P_{\widehat{j}_1}}^{G}(M^{\widehat{j}_1}(w_{j_1,j_0}),\pi_{j_1,j_2}^{\infty}))=0.\]
Then (\ref{equ: easy square case 2 embedding}) and (\ref{equ: easy square case 2 sub}) imply (\ref{equ: easy square OS minimal}).

\textbf{Case $2.3$}: In \ref{it: easy square 4}, we have $V\cong \cF_{P_{\widehat{j}_1}}^{G}(M,\pi_{j_1,j_2}^{\infty})$ and $V'\cong L(1)^\vee\otimes_E \pi_{\infty,\un{j}}^{\infty}$, and the assumption $\Sigma_{j_1,j_2}\cap W(G)\cdot \cJ(\pi_{\infty,\un{j}}^{\infty})\neq\emptyset$ is satisfied since both cosets contain $\cJ(\pi^\infty)$ (see case (iv) below (\ref{equ: easy square smooth part 4})). Thus we can apply directly (\ref{equ: OS main seq}) (with $w=1$) and obtain (\ref{equ: easy square OS minimal}) from (\ref{equ: easy square g minimal}).
\end{proof}

A $D(G)$-module is said to be $Z(\fg)$-finite if every element is killed by an ideal of finite codimension in $Z(\fg)$ (recall that $Z(\fg)$ lies in the center of $D(G)$ by \cite[Prop.~3.7]{ST02a}). For $\lambda, \mu\in \Lambda=\Lambda$, \cite[(2)]{JLS24} defines (by the formula (\ref{translationintro})) an exact endofunctor $\cT^{\mu}_{\lambda}$ on the abelian category of $Z(\fg)$-finite $D(G)$-modules such that for $I\subseteq \Delta$, $M$ in $\cO^{\fp_I}_{\rm{alg}}$ and $\pi^{\infty}$ a strongly admissible smooth representation of $L_I$ over $E$ we have an isomorphism of $D(G)$-modules (see \cite[Thm.~2]{JLS24})
\begin{equation}\label{equ: OS translation}
\cT^{\mu}_{\lambda}(\cF_{P_I}^{G}(M,\pi^{\infty})^\vee)\cong \cF_{P_I}^{G}(T^{\mu}_{\lambda}(M),\pi^{\infty})^\vee
\end{equation}
where $T^{\mu}_{\lambda}$ is defined in (\ref{equ: translation functor}) (note that $T^{\mu}_{\lambda}(M)$ remains in $\cO^{\fp_I}_{\rm{alg}}$ when $M$ is in $\cO^{\fp_I}_{\rm{alg}}$ using the argument in the proof of \cite[Thm.~1.1(d)]{Hum08}). Recall from \S\ref{subsec: loc an dist} that $\cC_{D(G)}$ is the abelian category of coadmissible $D(G)$-modules over $E$. If $V$ is an admissible locally analytic representation of $G$ over $E$ such that $V^\vee$ is $Z(\fg)$-finite, it follows from the discussion below \cite[Def.~2.4.5]{JLS24} that $\cT^\mu_{\lambda}(V^\vee)$ is again in $\cC_{D(G)}$, hence we can define an admissible locally analytic representation by $\cT^\mu_{\lambda}(V)\defeq \cT^\mu_{\lambda}(V^\vee)^\vee$. Note also that if $V$ has a central character, then (using that an irreducible algebraic representation of $G$ always has a central character), we deduce from \cite[(2)]{JLS24} that $\cT^\mu_{\lambda}(V)$ also has a central character.\bigskip

Let $j\in \{1,\dots,n-1\}$ and $\mu\in \Lambda$ such that $\langle \mu + \rho, \alpha^\vee\rangle \geq 0$ for $\alpha\in \Phi^+$ and the stabilizer of $\mu$ in $W(G)$ for the dot action is $\{1, s_j\}$. By the same formula as below (\ref{equ: translation functor}) we define a wall-crossing functor $\Theta_{\mu}\defeq \cT_\mu^{w_0\cdot \mu_0} \circ \cT_{w_0\cdot \mu_0}^\mu$ which is an exact endofunctor on the abelian category of $Z(\fg)$-finite $D(G)$-modules. If $V$ is an admissible locally analytic representation of $G$ over $E$ such that $V^\vee$ is $Z(\fg)$-finite, we define $\Theta_{\mu}(V)\defeq \Theta_{\mu}(V^\vee)^\vee$. By (\ref{equ: OS translation}) for $I\subseteq \Delta$, $M$ in $\cO^{\fp_I}_{\rm{alg}}$ and $\pi^{\infty}$ a strongly admissible smooth representation of $L_I$ over $E$ we have
\begin{equation}\label{equ: OS wall crossing}
\Theta_{\mu}(\cF_{P_I}^{G}(M,\pi^{\infty})^\vee)\cong \cF_{P_I}^{G}(\Theta_{s_j}(M),\pi^{\infty})^\vee.
\end{equation}
As $\Theta_{\mu}(\cF_{P_I}^{G}(M,\pi^{\infty})^\vee)$ only depends on $s_j$ (and on $\mu_0$) by (\ref{equ: OS wall crossing}), we write $\Theta_{s_j}(\cF_{P_I}^{G}(M,\pi^{\infty})^\vee)$ in that case.

\begin{rem}\label{independence}
Though we do not need it, it is possible that the endofunctor $\Theta_{\mu}=\cT_\mu^{w_0\cdot \mu_0} \circ \cT_{w_0\cdot \mu_0}^\mu$ only depends on $s_j$ (and on $\mu_0$) up to isomorphism. Let $\mu'\in \Lambda$ such that $\langle \mu' + \rho, \alpha^\vee\rangle \geq 0$ for $\alpha\in \Phi^+$ and the stabilizer of $\mu'$ in $W(G)$ for the dot action is $\{1, s_j\}$. By \cite[Prop.~5.0.8(c)]{Bez13} we have isomorphisms of functors $T_{\mu}^{\mu'}\circ T_{w_0\cdot \mu_0}^{\mu}\buildrel\sim\over\rightarrow T_{w_0\cdot \mu_0}^{\mu'}$ and $T_{\mu}^{w_0\cdot \mu_0}\circ T_{\mu'}^{\mu}\buildrel\sim\over\rightarrow T_{\mu'}^{w_0\cdot \mu_0}$. It would be enough to prove the analogous statements with $\cT_{\mu}^{\mu'}\circ \cT_{w_0\cdot \mu_0}^{\mu}$ and $\cT_{\mu}^{w_0\cdot \mu_0}\circ \cT_{\mu'}^{\mu}$. Indeed, we would then have isomorphisms of endofunctors
\[\cT_{\mu'}^{w_0\cdot \mu_0} \circ \cT_{w_0\cdot \mu_0}^{\mu'} \cong (\cT_{\mu}^{w_0\cdot \mu_0}\circ \cT_{\mu'}^{\mu})\circ (\cT_{\mu}^{\mu'}\circ \cT_{w_0\cdot \mu_0}^{\mu})\cong \cT_{\mu}^{w_0\cdot \mu_0}\circ (\cT_{\mu'}^{\mu}\circ \cT_{\mu}^{\mu'})\circ \cT_{w_0\cdot \mu_0}^{\mu} \cong \cT_{\mu}^{w_0\cdot \mu_0}\circ \cT_{w_0\cdot \mu_0}^{\mu}\]
where the last isomorphism follows from $\cT_{\mu'}^{\mu}\circ \cT_{\mu}^{\mu'}\cong \Id$ (see \cite[Thm.~3.2.1]{JLS24}). If $M$ is a $Z(\fg)$-finite $D(G)$-module, seeing $M$ as a $U(\fg)$-module we have an isomorphism of $U(\fg)$-modules $(T_{\mu}^{\mu'}\circ T_{w_0\cdot \mu_0}^{\mu})(M)\buildrel\sim\over\rightarrow T_{w_0\cdot \mu_0}^{\mu'}(M)$, and since the functors $\cT^{\mu}_{\lambda}$ in \cite[\S 1]{JLS24} are just $T^{\mu}_{\lambda}$ on the underlying $U(\fg)$-modules, it would be enough to prove that this isomorphism is $D(G)$-equivariant. An examination of the proof of \cite[Thm.~3.2.1]{JLS24} shows that it would even be enough to prove this $D(G)$-equivariance for $M$ the form $D(G)\otimes_{U(\fg)}N$ where $N$ is in $\mathrm{Mod}_{U(\fg)}$. However, we couldn't find a quick argument for this $D(G)$-equivariance (if true).
\end{rem}

\begin{lem}\label{lem: wall crossing of simple}
Let $x\in W(G)$ and $\pi^{\infty}$ an irreducible $G$-regular smooth representation of $L_{I_x}$ such that $V=\cF_{P_{I_x}}^G(L(x),\pi^{\infty})$ is irreducible (see \ref{it: OS property 3} of Theorem~\ref{prop: OS property}).
\begin{enumerate}[label=(\roman*)]
\item \label{it: wall crossing simple 1} If $j\notin D_R(x)$ then $\Theta_{w_0s_jw_0}(V)=0$.
\item \label{it: wall crossing simple 2} If $j\in D_R(x)$ then $\Theta_{w_0s_jw_0}(V)$ has $V$ as both socle and cosocle.
\end{enumerate}
\end{lem}
\begin{proof}
By Proposition~\ref{prop: Jantzen middle} $\Theta_{w_0s_jw_0}(L(x))$ is zero if $j\notin D_R(x)$, and has Loewy length $3$ with both socle and cosocle isomorphic to $L(x)$ (and middle layer not containing $L(x)$) if $j\in D_R(x)$. Then \ref{it: wall crossing simple 1} follows from (\ref{equ: OS wall crossing}). Assume $j\in D_R(x)$. Let $W$ be an irreducible constituent of $\Theta_{w_0s_jw_0}(V)$, which we can write $W=\cF_{I_w}^G(L(w),\sigma^{\infty})$ for some $w\in W(G)$ such that $I_x\subseteq I_w$ and some irreducible $G$-regular smooth representation $\sigma^{\infty}$ which is a subquotient of $i_{I_x,I_w}^\infty(\pi^\infty)$ (using \ref{it: OS property 2} of Theorem~\ref{prop: OS property} and \ref{it: PS 1}, \ref{it: PS 2} of Lemma~\ref{lem: Jacquet of PS}). By the first sentence of this proof together with \ref{it: Hom OS socle} of Lemma~\ref{lem: Hom OS socle cosocle} applied with $I_0=I_w$, $I_1=I_x$, $M_0=L(w)$, $M_1=\Theta_{w_0s_jw_0}(L(x))$, $\pi_0^{\infty}=\sigma^{\infty}$ and $\pi_1^{\infty}=\pi^{\infty}$, we have $\Hom_{G}(W,\Theta_{w_0s_jw_0}(V))\neq 0$ if and only if $w=x$ and $\sigma^{\infty}=\pi^{\infty}$ if and only if $W=V$ (using Lemma \ref{lem: Hom OS} for the last equivalence). Moreover the space of homomorphisms is then $1$-dimensional (still by \ref{it: Hom OS socle} of Lemma~\ref{lem: Hom OS socle cosocle}). It follows that $V$ is the socle of $\Theta_{w_0s_jw_0}(V)$. An analogous argument using \ref{it: Hom OS cosocle} of Lemma~\ref{lem: Hom OS socle cosocle} gives that $V$ is also the cosocle of $\Theta_{w_0s_jw_0}(V)$.
\end{proof}

\begin{lem}\label{lem: square as wall crossing}
Let $\un{j}=(j_0,j_1,j_2)\in\mathbf{J}$ and $\mu\in \Lambda$ such that $\langle \mu + \rho, \alpha^\vee\rangle \geq 0$ for $\alpha\in \Phi^+$ and the stabilizer of $\mu$ in $W(G)$ for the dot action is $\{1, w_0s_{j_0}w_0\}$.
\begin{enumerate}[label=(\roman*)]
\item \label{it: square crossing 1} The representation $\Theta_{w_0s_{j_0}w_0}(C_{\un{j}})$ of $G$ has simple socle and cosocle $C_{\un{j}}$ and middle layer
\begin{equation*}
\mathrm{rad}_{G}(\Theta_{w_0s_{j_0}w_0}(C_{\un{j}}))/\mathrm{soc}_{G}(\Theta_{w_0s_{j_0}w_0}(C_{\un{j}}))
\cong C_{(j_0-1,j_1,j_2)}\oplus C_{(j_0+1,j_1,j_2)}\oplus L(1)^\vee\otimes_E \pi^{\infty}
\end{equation*}
where $\pi^{\infty}$ is non-zero if and only if $j_0=j_1$, in which case $\pi^{\infty}\cong i_{\widehat{j}_0,\Delta}^{\infty}(\pi_{j_0,j_2}^{\infty})$, and where we omit $C_{(j_0-1,j_1,j_2)}$ when $j_0=1$ and $C_{(j_0+1,j_1,j_2)}$ when $j_0=n-1$.
\item \label{it: square crossing 2} Assume $j_0<n-1$ and let $\un{j}'\in\{(j_0,j_1+1,j_2+1), (j_0,j_1-1,j_2)\}$, $\un{j}''\defeq (j_0+1,j_1',j_2')$ and $V_{\un{j},\un{j}''}$ the minimal $\mathrm{Ext}$-square with socle $C_{\un{j}}$ and cosocle $C_{\un{j}''}$ constructed in \ref{it: easy square 1} or \ref{it: easy square 2} of Proposition~\ref{prop: easy square}. Let $V$ be the unique length $2$ representation of $G$ with socle $C_{\un{j}}$ and cosocle $C_{\un{j}'}$ defined from Lemma~\ref{lem: Ext1 factor 1}. Then $\Theta_{\mu}(V)$ admits a unique subquotient isomorphic to $V_{\un{j},\un{j}''}$ which is moreover a subrepresentation.
\end{enumerate}
\end{lem}
\begin{proof}
We prove \ref{it: square crossing 1}.
It follows from Lemma~\ref{lem: wall crossing of simple} that $\Theta_{w_0s_{j_0}w_0}(C_{\un{j}})$ has simple socle and cosocle $C_{\un{j}}$.
By Proposition~\ref{prop: Jantzen middle} and Remark~\ref{rem: explicit middle} $\Theta_{w_0s_{j_0}w_0}(L(w_{j_1,j_0}))$ has Loewy length $3$ with both socle and cosocle $L(w_{j_1,j_0})$ and with middle layer $\mathrm{rad}_1(\Theta_{w_0s_{j_0}w_0}(L(w_{j_1,j_0})))$ semi-simple and multiplicity free. More precisely, from Remark~\ref{rem: explicit middle} $L(x)$ is a constituent of $\mathrm{rad}_1(\Theta_{w_0s_{j_0}w_0}(L(w_{j_1,j_0})))$ if and only if either $j_0=j_1$ and $x=1$, or $j_0>1$ and $x=w_{j_1,j_0-1}$, or $j_0<n-1$ and $x=w_{j_1,j_0+1}$. If $x=w_{j_1,j}$ for $j\in\{j_0-1,j_0+1\}$, we have $\cF_{P_{\widehat{j}_1}}^{G}(L(x),\pi_{j_1,j_2}^{\infty})=C_{(j,j_1,j_2)}$ by definition (see (\ref{cj})). If $x=1$ (with $j_0=j_1$), we have $\cF_{P_{\widehat{j}_1}}^{G}(L(1),\pi_{j_1,j_2}^{\infty})=L(1)^\vee\otimes_E i_{\widehat{j}_0,\Delta}^{\infty}(\pi_{j_0,j_2}^{\infty})$ using \ref{it: OS property 2} of Theorem~\ref{prop: OS property}. By (\ref{equ: OS wall crossing}) this finishes the proof of \ref{it: square crossing 1}.

We prove \ref{it: square crossing 2}. Let $I\defeq \Delta\setminus\{j_1,j'_1\}=\widehat{j}_1\cap \widehat{j}_1'$ and $M_0$ the (unique) length $2$ object in $\cO^{\fp_I}_{\rm{alg}}$ with socle $L(w_{j_1',j_0})$ and cosocle $L(w_{j_1,j_0})$. Let $\Sigma\defeq \Sigma_{j_1,j_2}\cap \Sigma_{j_1',j_2'}$ and $\pi^{\infty}$ as in (\ref{equ: easy square smooth part 1 2}). Then by Step $5$ in the proof of Lemma~\ref{lem: Ext1 g vers OS} we know that $V$ (as in the statement) is a subquotient of $\cF_{P_I}^{G}(M_0,\pi^{\infty})$. By Lemma~\ref{lem: g square as wall crossing} $L(w_{j_1',j_0+1})$ appears with multiplicity one in $\Theta_{w_0s_{j_0}w_0}(M_0)$ and the unique quotient of $\Theta_{w_0s_{j_0}w_0}(M_0)$ with socle $L(w_{j_1',j_0+1})$ is isomorphic to $Q_1(w_{j_1,j_0},w_{j_1',j_0+1})$ (see Proposition \ref{lem: explicit g square}). Moreover $\pi_{j_1'',j_2''}^{\infty}=\pi_{j_1',j_2'}^{\infty}$ appears with multiplicity one in the $G$-basic, hence multiplicity free, representation $i_{I,\widehat{j}_1'}^\infty(\pi^{\infty})$ (see below (\ref{equ: easy square smooth part 1 2})). Using Theorem~\ref{prop: OS property} and Lemma~\ref{lem: Hom OS} this implies that $C_{\un{j}''}$ appears with multiplicity one in $\Theta_{w_0s_{j_0}w_0}(\cF_{P_I}^{G}(M_0,\pi^{\infty}))\cong \cF_{P_I}^{G}(\Theta_{w_0s_{j_0}w_0}(M_0),\pi^{\infty})$. Moreover the unique subrepresentation of $\Theta_{w_0s_{j_0}w_0}(\cF_{P_I}^{G}(M_0,\pi^{\infty}))$ with cosocle $C_{\un{j}''}$ is a subrepresentation of $\cF_{P_I}^{G}(Q_1(w_{j_1,j_0},w_{j_1',j_0+1}),\pi^{\infty})$ and thus is multiplicity free (using Lemma~\ref{lem: Hom OS} and the facts that $Q_1(w_{j_1,j_0},w_{j_1',j_0+1})$ is multiplicity free and $\pi^{\infty}$ irreducible $G$-regular).
By the last paragraph in Step $1$ in the proof of Proposition~\ref{prop: easy square}, we know that $V_{\un{j},\un{j}''}$ is the unique subquotient of $\cF_{P_I}^{G}(Q_1(w_{j_1,j_0},w_{j_1',j_0+1}),\pi^{\infty})$ with socle $C_{\un{j}}$ and cosocle $C_{\un{j}''}$. Since any subquotient of $\Theta_{w_0s_{j_0}w_0}(\cF_{P_I}^{G}(M_0,\pi^{\infty}))$ with cosocle $C_{\un{j}''}$ is a subquotient of $\cF_{P_I}^{G}(Q_1(w_{j_1,j_0},w_{j_1',j_0+1}),\pi^{\infty})$ by the previous discussion, $V_{\un{j},\un{j}''}$ is also the unique subquotient of $\Theta_{w_0s_{j_0}w_0}(\cF_{P_I}^{G}(M_0,\pi^{\infty}))$ with socle $C_{\un{j}}$ and cosocle $C_{\un{j}''}$. As $\Theta_{\mu}(V)$ is a subquotient of $\Theta_{w_0s_{j_0}w_0}(\cF_{P_I}^{G}(M_0,\pi^{\infty}))$ and $C_{\un{j}}, C_{\un{j}''}\in \mathrm{JH}_{G}(\Theta_{\mu}(V))$ (by \ref{it: square crossing 1} applied to $C_{\un{j}}$ and $C_{\un{j}'}$), $V_{\un{j},\un{j}''}$ is the unique subquotient of $\Theta_{\mu}(V)$ with socle $C_{\un{j}}$ and cosocle $C_{\un{j}''}$. By the explicit description of $V_{\un{j},\un{j}''}$ in \ref{it: easy square 1} or \ref{it: easy square 2} of Proposition~\ref{prop: easy square} we know that $V$ injects into $V_{\un{j},\un{j}''}$. By the first statement in Proposition~\ref{prop: Jantzen middle} and (\ref{equ: OS wall crossing}) (and the exactness of $\Theta_\mu$) we have $\Theta_{\mu}(V_{\un{j},\un{j}''}/V)=0$ and thus the injection $V\hookrightarrow V_{\un{j},\un{j}''}$ induces an isomorphism $\Theta_{\mu}(V)\buildrel\sim\over\rightarrow \Theta_{\mu}(V_{\un{j},\un{j}''})$. Moreover, the canonical maps $C_{\un{j}}\rightarrow \Theta_{\mu}(C_{\un{j}})$ and $C_{\un{j}'}\rightarrow \Theta_{\mu}(C_{\un{j}'})$ are injective (since non-zero), hence so is the canonical map $V\rightarrow\Theta_{\mu}(V)\cong \Theta_{\mu}(V_{\un{j},\un{j}''})$. We deduce that the restriction of $V_{\un{j},\un{j}''}\rightarrow \Theta_{\mu}(V_{\un{j},\un{j}''})$ to $V$ is injective. Since $V$ and $V_{\un{j},\un{j}''}$ have same socle $C_{\un{j}}$, it follows that $V_{\un{j},\un{j}''}$ injects into $\Theta_{\mu}(V_{\un{j},\un{j}''})\cong \Theta_{\mu}(V)$.
\end{proof}

\begin{rem}\label{rem: opposite crossing}
Similarly, we can prove that $\Theta_{w_0s_{j_0+1}w_0}(C_{(j_0+1,j_1,j_2)})$ has Loewy length $3$ with both socle and cosocle $C_{(j_0+1,j_1,j_2)}$, and with middle layer not containing $C_{(j_0+1,j_1,j_2)}$ but containing $C_{(j_0,j_1,j_2)}$ with multiplicity one. Keeping the notation and assumption in \ref{it: square crossing 2} of Lemma \ref{lem: square as wall crossing}, let $V'$ be the unique length $2$ representation of $G$ with socle $C_{(j_0+1,j_1,j_2)}$ and cosocle $C_{\un{j}''}$ (using Lemma~\ref{lem: Ext1 factor 1}), then $\Theta_{\mu}(V')$ also uniquely admits $V_{\un{j},\un{j}''}$ as a subquotient (with $V_{\un{j},\un{j}''}$ constructed in \ref{it: easy square 1} or \ref{it: easy square 2} of Proposition~\ref{prop: easy square}), which is moreover a quotient.
\end{rem}

\begin{rem}\label{rem: loc alg crossing}
Using similar arguments as in the proof of \ref{it: square crossing 2} of Lemma~\ref{lem: square as wall crossing}, we can also prove the following results (where we choose $\mu\in \Lambda$ such that $\langle \mu + \rho, \alpha^\vee\rangle \geq 0$ for $\alpha\in \Phi^+$ and the stabilizer of $\mu$ in $W(G)$ for the dot action is $\{1, w_0s_jw_0\}$).
\begin{enumerate}[label=(\roman*)]
\item \label{it: loc alg crossing 1} Let $V$ be the (unique) length $2$ representation of $G$ with socle $C_{(j_0,j_0+1,j_2)}$ and cosocle $C_{(j_0,j_0,j_2)}$ for some $1\leq j_0<j_2\leq n$ (Lemma~\ref{lem: Ext1 factor 1}), then $V_{(j_0,j_0+1,j_2),\infty}$ in \ref{it: easy square 3} of Proposition \ref{prop: easy square} is isomorphic to the unique subrepresentation of $\Theta_{\mu}(V)$ with cosocle $L(1)^\vee\otimes_E V_{[j_2-j_0,j_2]\cap\Delta,\Delta}^{\infty}$.
\item Let $V'$ be the (unique) length $2$ representation of $G$ with socle $C_{(j_0,j_0,j_2)}$ and cosocle $C_{(j_0,j_0-1,j_2)}$ for some $1< j_0\leq j_2\leq n$ (Lemma~\ref{lem: Ext1 factor 1}), then $V_{\infty, (j_0,j_0-1,j_2)}$ in \ref{it: easy square 4} of Proposition \ref{prop: easy square} is isomorphic to the unique quotient of $\Theta_{\mu}(V')$ with socle $L(1)^\vee\otimes_E V_{[j_2-j_0+1,j_2-1],\Delta}^{\infty}$.
\end{enumerate}
\end{rem}

We now introduce some technical but useful notation, which will be used in Lemma \ref{lem: square Ext1 vanishing} and Lemma \ref{lem: conj Ext1 vanishing OS} below. For a fixed $\un{j}\in\mathbf{J}$ with $(j_0,j_1,j_2+1)\in\mathbf{J}$ we define
\begin{equation}\label{equ: Sjj}
S_{j_0,j_1}\defeq \{x\in W(G) \mid w_{j_1,j_0}\prec x, j_1\notin D_L(x)\}.
\end{equation}
It follows from \ref{it: special descent 2} of Lemma \ref{lem: special descent change} that $S_{j_0,j_1}=\{w_{j_1+1,j_0}\}$ if $j_1>j_0$ and $j_1+1\in \Delta$, $S_{j_0,j_1}=\{w_{j_1-1,j_0}\}$ if $j_1<j_0$ and $j_1-1\in \Delta$, $S_{j_0,j_1}=\{w_{j,j_0}\mid j\in \{j_1-1,j_1+1\}\cap\Delta\}$ if $j_1=j_0$, and $S_{j_0,j_1}=\emptyset$ otherwise.
We define
\[\left\{\begin{array}{rclcc}
A_{\un{j}}&\defeq &C_{(j_0,j_1-1,j_2)} &\mathrm{if}& j_1>j_0\\
A_{\un{j}}&\defeq & C_{(j_0,j_1+1,j_2+1)}&\mathrm{if}& j_1<j_0\\
A_{\un{j}}&\defeq &L(1)^\vee\otimes_E V_{[j_2-j_0+1,j_2],\Delta}^{\infty}& \mathrm{if} & j_1=j_0.
\end{array}\right.\]
We also define
\[B_{\un{j}}\defeq \oplus_{\un{j'}}C_{\un{j'}},\]
the direct sum being over those $\un{j'}\in \mathbf{J}$ such that $\un{j}\leq \un{j'}$, $j_0=j_0'$, $d(\un{j}, \un{j'})=|j_2'-j_2|+|(j_2'-j_1')-(j_2-j_1)|=1$ and $w_{j_1',j_0}\in S_{j_0,j_1}$. An easy check shows that there is in fact a bijection between $S_{j_0,j_1}$ and $\mathrm{JH}_{G}(B_{\un{j}})$ given by $w_{j_1+1,j_0}\mapsto C_{(j_0,j_1+1,j_2+1)}$, $w_{j_1-1,j_0}\mapsto C_{(j_0,j_1-1,j_2)}$. For $S\subseteq S_{j_0,j_1}$ we define $B_{\un{j},S}$ as the direct summand of $B_{\un{j}}$ corresponding to $S$ under this bijection (with $B_{\un{j},\emptyset}\defeq 0$). We finally define $W_{\un{j}}$ as the unique (up to isomorphism) representation with socle $C_{\un{j}}$ which fits into an exact sequence $0\rightarrow C_{\un{j}}\rightarrow W_{\un{j}}\rightarrow B_{\un{j}}\rightarrow 0$. Note that the existence and unicity of $W_{\un{j}}$ follows from Lemma~\ref{lem: Ext1 factor 1}.

\begin{lem}\label{lem: square Ext1 vanishing}
Let $\un{j}\in\mathbf{J}$ with $(j_0,j_1,j_2+1)\in\mathbf{J}$. Then we have
\begin{equation}\label{equ: square Ext1 vanishing}
\mathrm{Ext}_{G}^1(C_{(j_0,j_1,j_2+1)}, W_{\un{j}})=0.
\end{equation}
\end{lem}
\begin{proof}
If $S_{j_0,j_1}=\emptyset$, i.e.~$B_{\un{j}}=0$, we have $W_{\un{j}}=C_{\un{j}}$ and (\ref{equ: square Ext1 vanishing}) follows from Lemma~\ref{lem: Ext1 factor 1} (note that $d(\un{j},(j_0,j_1,j_2+1))=2$). We assume from now on $S_{j_0,j_1}\neq \emptyset$ and write $M$ for the unique $U(\fg)$-module with cosocle $L(w_{j_1,j_0})$ which fits into an exact sequence $0\rightarrow L_{S_{j_0,j_1}}\rightarrow M\rightarrow L(w_{j_1,j_0})\rightarrow 0$ (recall that for $S\subseteq W(G)$ we define $L_S=\oplus_{x\in S}L(x)$). We let $I\subseteq \Delta$ be the maximal subset such that $M$ is in $\cO_{\rm{alg}}^{\fp_I}$. In the following, we only prove (\ref{equ: square Ext1 vanishing}) when $j_0=j_1$ and $S_{j_1,j_0}=\{j_1-1,j_1+1\}$, the other cases being simpler and left to the (interested) reader. In particular we have $2\leq j_0=j_1\leq n-2$ and $I=\Delta\setminus\{j_1-1,j_1,j_1+1\}$. We recall the notation (above) Lemma~\ref{lem: coset square}: $\Sigma_{\pm}=\Sigma_{j_1,j_2}\cap s_{j_1}\cdot \Sigma_{j_1,j_2+1}$ (a left $W(L_I)$-coset) and
\[\pi_{\pm}^{\infty}=J_{\widehat{j}_1,\widehat{j}_1,s_{j_1}}(\pi_{j_1,j_2+1}^{\infty})_{\cB^{I}_{\Sigma_{\pm}}}\cong J_{\widehat{j}_1,I}'(\pi_{j_1,j_2}^{\infty})_{\cB^{I}_{\Sigma_{\pm}}}\]
which is an irreducible $G$-regular representation of $L_I$ over $E$ (note that $I_{\pm}=\Delta\setminus\{j_1-1,j_1,j_1+1\}=I$). In particular $\cF_{P_{I}}^{G}(M,\pi_{\pm}^{\infty})$ is multiplicity free by Lemma \ref{lem: Hom OS}.\bigskip

\textbf{Step $1$}: We prove that $W_{\un{j}}$ is a subquotient of $\cF_{P_{I}}^{G}(M,\pi_{\pm}^{\infty})$.\\
By (\ref{equ: second adjunction}) applied to (v) above Lemma \ref{lem: coset square}, $\pi_{j_1,j_2}^{\infty}$ is in the cosocle of $i_{I,\widehat{j}_1}^{\infty}(\pi_{\pm}^{\infty})$. By \ref{it: coset square 2} of Lemma~\ref{lem: coset square} (and the definition of $\pi_{+,1}^{\infty}$, $\pi_{-,1}^{\infty}$ above \emph{loc.~cit.}) $\pi_{\ast,1}^{\infty}\cong \mathrm{cosoc}_{L_{I_{\ast}}}(i_{I,I_{\ast}}^{\infty}(\pi_{\pm}^{\infty}))$ for $\ast\in\{+,-\}$. By (\ref{equ: first adjunction}) applied to (iii) and (iv) above Lemma \ref{lem: coset square}, $\pi_{j_1+1,j_2+1}^{\infty}$, $\pi_{j_1-1,j_2}^{\infty}$ is in the socle of $i_{I_+,\Delta\setminus\{j_1+1\}}^{\infty}(\pi_{+,1}^{\infty})$, $i_{I_-,\Delta\setminus\{j_1-1\}}^{\infty}(\pi_{-,1}^{\infty})$ respectively. Thus $\pi_{j_1+1,j_2+1}^{\infty}$, $\pi_{j_1-1,j_2}^{\infty}$ is a subquotient of $i_{I,\Delta\setminus\{j_1+1\}}^{\infty}(\pi_{\pm}^{\infty})$, $i_{I,\Delta\setminus\{j_1-1\}}^{\infty}(\pi_{\pm}^{\infty})$ respectively. Moreover we have $d(\pi_{j_1+1,j_2+1}^{\infty},\pi_{j_1,j_2}^{\infty})=d(\pi_{j_1-1,j_2}^{\infty},\pi_{j_1,j_2}^{\infty})=0$ by \ref{it: connect 1} of Lemma~\ref{lem: connect Hom}. We can thus apply Lemma~\ref{lem: general construction} to $M$ as above and $\pi_{\pm}^{\infty}$, which gives that $W_{\un{j}}$ is a subquotient of $\cF_{P_{I}}^{G}(M,\pi_{\pm}^{\infty})$. We define $W_{\un{j},-}$ as the (unique) minimal length subrepresentation of $\cF_{P_{I}}^{G}(M,\pi_{\pm}^{\infty})$ which admits $W_{\un{j}}$ as a quotient.\bigskip

\textbf{Step $2$}: We prove that the injection $W_{\un{j},-}\hookrightarrow \cF_{P_{I}}^{G}(M,\pi_{\pm}^{\infty})$ and surjection $W_{\un{j},-}\twoheadrightarrow W_{\un{j}}$ induce isomorphisms
\begin{equation}\label{equ: square Ext1 subquotient}
\mathrm{Ext}_{G}^1(C_{(j_0,j_1,j_2+1)}, W_{\un{j}}) \xleftarrow{\sim} \mathrm{Ext}_{G}^1(C_{(j_0,j_1,j_2+1)}, W_{\un{j},-})\xrightarrow{\sim} \mathrm{Ext}_{G}^1(C_{(j_0,j_1,j_2+1)}, \cF_{P_{I}}^{G}(M,\pi_{\pm}^{\infty})).
\end{equation}
We wish to apply Lemma~\ref{lem: minimal OS} to $V_0=C_{(j_0,j_1,j_2+1)}$, $\pi^{\infty}=\pi_{\pm}^{\infty}$ and $V=W_{\un{j}}$, so we check that all assumptions there are satisfied (using the notation of \emph{loc.~cit.}). We recall that $j_0=j_1$, so $w_{j_1,j_0}=s_{j_1}$. First, for $L(x)\in\mathrm{JH}_{U(\fg)}(M)$ with $x\neq w_{j_1,j_0}=s_{j_1}$ we have $j_1\notin D_L(x)$, and hence $L(s_{j_1})$ is not a constituent of $M^{I_x}(x)$ by \cite[Thm.~9.4(c)]{Hum08} and Lemma \ref{lem: dominance and left set}. Secondly we have $d(\pi_{j_1,j_2+1}^{\infty},\pi_{j_1,j_2}^{\infty})=d(\pi_{j_1,j_2+1}^{\infty},\pi_{j_1+1,j_2+1}^{\infty})=d(\pi_{j_1,j_2+1}^{\infty},\pi_{j_1-1,j_2}^{\infty})=0$ by \ref{it: connect 1} of Lemma~\ref{lem: connect Hom}. Thirdly, writing $\Sigma_{0,x}= \Sigma_{j_1,j_2+1}\cap W(L_{I_x})\cdot\cJ(\sigma_x^{\infty})$ for $L(x)\in \mathrm{JH}_{U(\fg)}(M)$, we have $\Sigma_{0,s_{j_1}}=\Sigma_{j_1,j_2+1}\cap \Sigma_{j_1,j_2} =\emptyset$ by \ref{it: coset 1} of Lemma \ref{lem: coset intersection}. Finally we need to check that for $L(x)\in \mathrm{JH}_{U(\fg)}(M)\setminus\{L(s_{j_1})\}$ and $\tau_x^{\infty}\in\mathrm{JH}_{L_{I_x}}(\sigma_{x,-}^{\infty})\setminus\{\sigma_x^{\infty}\}$ we have
\begin{equation*}
J_{I_x,I_x\cap \widehat{j}_1}'(\tau_x^{\infty})_{\cB^{I_x\cap \widehat{j}_1}_{\Sigma_{0,x}}}=0
\end{equation*}
where $\sigma_x^{\infty}\defeq \pi_{j_1+1,j_2+1}^{\infty}$ if $x=w_{j_1+1,j_0}$ and $\sigma_x^{\infty}\defeq \pi_{j_1-1,j_2}^{\infty}$ if $x=w_{j_1-1,j_0}$. But this follows from \ref{it: coset square 3} of Lemma \ref{lem: coset square}, noting that the constituents in $\mathrm{JH}_{L_{I_x}}(\sigma_{x,-}^{\infty})\setminus\{\sigma_x^{\infty}\}$ are exactly the constituents $\tau_x^{\infty}$ in $i_{I,I_x}^{\infty}(\pi_{\pm}^{\infty})$ such that $\tau_x^{\infty}< \sigma_x^{\infty}$.\bigskip

\textbf{Step $3$}: We prove (for $j_0=j_1$)
\begin{equation}\label{equ: square Ext1 vanishing prime}
\mathrm{Ext}_{G}^1(C_{(j_0,j_1,j_2+1)}, \cF_{P_{I}}^{G}(M,\pi_{\pm}^{\infty}))=0.
\end{equation}
As $M$ has cosocle $L(s_{j_1})$ and socle $L(w_{j_1-1,j_1})\oplus L(w_{j_1+1,j_1})$, it is a quotient of $M(s_{j_1})$ (using (\ref{equ: O Hom radical})), which by \cite[Thm.~9.4(c)]{Hum08} implies that $M$ is a quotient of $M^{I}(s_{j_1})$. Let $W\defeq \cF_{P_{I}}^{G}(M^{I}(s_{j_1}),\pi_{\pm}^{\infty})$, then $W/\cF_{P_{I}}^{G}(M,\pi_{\pm}^{\infty})\cong \cF_{P_{I}}^{G}(Q,\pi_{\pm}^{\infty})$ where $Q\defeq \mathrm{ker}(M^{I}(s_{j_1})\twoheadrightarrow M)$ (\ref{it: OS property 1} of Theorem~\ref{prop: OS property}). It is clear that $L(s_{j_1})\notin\mathrm{JH}_{U(\fg)}(Q)$ and thus we have by Lemma~\ref{lem: Hom OS}
\[\Hom_{G}(C_{(j_0,j_1,j_2+1)}, W/\cF_{P_{I}}^{G}(M,\pi_{\pm}^{\infty}))=0.\]
Hence, the injection $\cF_{P_{I}}^{G}(M,\pi_{\pm}^{\infty})\hookrightarrow W$ induces an embedding
\begin{equation}\label{equ: square Ext1 embedding}
\mathrm{Ext}_{G}^1(C_{(j_0,j_1,j_2+1)}, \cF_{P_{I}}^{G}(M,\pi_{\pm}^{\infty}))\hookrightarrow \mathrm{Ext}_{G}^1(C_{(j_0,j_1,j_2+1)}, W).
\end{equation}
By \ref{it: H0 conjugate 2} of Lemma~\ref{lem: H0 Weyl conjugate} applied with $I'=\widehat{j}_1$ and $j_1\notin I$ we have $H^0(\fn_{I},L(s_{j_1})^{s_{j_1}})=0$. By \ref{it: H1 conj 1} of Lemma~\ref{lem: H1 conjugate} applied with $I'=\widehat{j}_1$ we also have $\Hom_{U(\fl_I)}(L^I(s_{j_1}), H^1(\fn_I, L(s_{j_1})^{s_{j_1}}))= 0$. By (\ref{equ: g spectral seq}) applied with $M_I=L^I(s_{j_1})$ and $M=L(s_{j_1})^{s_{j_1}}$ we deduce for $\ell\leq 1$
\begin{equation}\label{equ: Extell}
\mathrm{Ext}_{U(\fg)}^{\ell}(M^{I}(s_{j_1}), L(s_{j_1})^{s_{j_1}})=0.
\end{equation}
Since $\Sigma_{\pm}\cap s_{j_1}\cdot\Sigma_{j_1,j_2+1}=\Sigma_{\pm}\ne \emptyset$ we can apply (\ref{equ: OS main seq}) with $w=s_{j_1}$ which gives using (\ref{equ: Extell})
\[\mathrm{Ext}_{G}^1(C_{(j_0,j_1,j_2+1)},W)=0.\]
This together with (\ref{equ: square Ext1 embedding}) gives (\ref{equ: square Ext1 vanishing prime}). Then (\ref{equ: square Ext1 vanishing}) follows from (\ref{equ: square Ext1 vanishing prime}) and (\ref{equ: square Ext1 subquotient}).
\end{proof}

For $S\subseteq S_{j_0,j_1}$ we define $\tld{B}_{\un{j},S}$ as the unique (up to isomorphism) representation with cosocle $C_{(j_0,j_1,j_2+1)}$ which fits into an exact sequence $0\rightarrow A_{\un{j}}\oplus B_{\un{j},S}\rightarrow \tld{B}_{\un{j},S}\rightarrow C_{(j_0,j_1,j_2+1)}\rightarrow 0$. Note that the existence and unicity of $\tld{B}_{\un{j},S}$ follows from Lemma~\ref{lem: Ext1 factor 1} and \ref{it: Ext1 with alg 2} of Lemma~\ref{lem: Ext1 factor 2}.

\begin{lem}\label{lem: conj Ext1 vanishing OS}
Let $\un{j}\in\mathbf{J}$ with $(j_0,j_1,j_2+1)\in\mathbf{J}$. Then we have for $S\subsetneq S_{j_0,j_1}$
\begin{equation}\label{equ: conj Ext1 vanishing OS}
\mathrm{Ext}_{G}^1(\tld{B}_{\un{j},S}, C_{\un{j}})=0.
\end{equation}
\end{lem}
\begin{proof}
Let $S\subsetneq S_{j_0,j_1}$. By \ref{it: Ext O 3} of Lemma~\ref{lem: Ext 1 category O} we have $\Dim_E\mathrm{Ext}_{\cO^{\fb}_{\rm{alg}}}^1(L(s_{j_1}w_{j_1,j_0}),L(w_{j_1,j_0}))\!=1$ (recall $j_1\in D_L(w_{j_1,j_0})$). By \ref{it: rabiotext 1} of Lemma \ref{rabiotext} and the explicit description of $S_{j_0,j_1}$ below (\ref{equ: Sjj}) we have $\Dim_E\mathrm{Ext}_{\cO^{\fb}_{\rm{alg}}}^1(L(x),L(w_{j_1,j_0}))=1$ for $x\in S_{j_0,j_1}$. Hence there is a unique $M$ in $\cO^{\fb}_{\rm{alg}}$ with socle $L(w_{j_1,j_0})$ which fits into an exact sequence (where $L_S=\oplus_{x\in S}L(x)$)
\[0\longrightarrow L(w_{j_1,j_0})\longrightarrow M\longrightarrow L(s_{j_1}w_{j_1,j_0})\oplus L_S\longrightarrow 0.\]
Let $I\subseteq \Delta$ be the maximal subset such that $M$ is in $\cO^{\fp_I}_{\rm{alg}}$. An explicit check using Lemma \ref{lem: dominance and left set} and the explicit description of $S_{j_0,j_1}$ below (\ref{equ: Sjj}) shows that we have the following cases (with the notation in (\ref{equ: I+-}) and using that $S$ is \emph{strictly} smaller than $S_{j_0,j_1}$):
\begin{itemize}
\item$I=I_{-}$ if either $j_1>j_0$ with $S=\emptyset$, or $j_1=j_0$ with $S=\{w_{j_1-1,j_0}\}$;
\item$I=I_{+}$ if either $j_1<j_0$ with $S=\emptyset$, or $j_1=j_0$ with $S=\{w_{j_1+1,j_0}\}$;
\item$I=\widehat{j}_1$ if $j_1=j_0$ with $S=\emptyset$.
\end{itemize}
We let $\Sigma\defeq \Sigma_{\ast,0}$ if $I=I_{\ast}$ with $\ast\in\{+,-\}$ and $\Sigma\defeq \Sigma_{j_1,j_2+1}$ if $I=\widehat{j}_1$ (see (\ref{equ: sigmainter})). From (\ref{equ: coset square intersection}) and (\ref{equ: sigmainter}) we have
\begin{equation}\label{equ: cube coset intersection}
\Sigma_{j_1,j_2}\cap s_{j_1}\cdot \Sigma=\Sigma_{j_1,j_2}\cap s_{j_1}\cdot\Sigma_{j_1,j_2+1}=\Sigma_{\pm}\neq \emptyset.
\end{equation}
We define $\pi^{\infty}\defeq J_{\widehat{j}_1,I}(\pi_{j_1,j_2+1}^{\infty})_{\cB^I_{\Sigma}}$, an irreducible $G$-regular representation of $L_I$ (so $\pi^{\infty}=\pi_{j_1,j_2+1}^{\infty}$ if $I=\widehat{j}_1$ and $\pi^{\infty}=\pi_{\ast,0}^{\infty}$ if $I=I_{\ast}$ for $\ast\in\{+,-\}$ with the notation above Lemma \ref{lem: coset square}). By (\ref{equ: first adjunction}) and using \ref{it: basic as image} of Remark \ref{rem: basic PS intertwine} $\pi_{j_1,j_2+1}^{\infty}$ is the socle of $i_{I,\widehat{j}_1}^{\infty}(\pi^{\infty})$, and by (\ref{equ: second adjunction}) with the isomorphism in (i) (resp.~(ii)) above Lemma \ref{lem: coset square} $\pi_{j_1+1,j_2+1}^{\infty}$ (resp.~$\pi_{j_1-1,j_2}^{\infty}$) is the cosocle of $i_{I,\Delta\setminus\{j_1+1\}}^{\infty}(\pi^{\infty})$ if $I=I_{+}$ (resp.~of $i_{I,\Delta\setminus\{j_1-1\}}^{\infty}(\pi^{\infty})$ if $I=I_{-}$). By Lemma~\ref{lem: explicit smooth induction} (and $j_1\leq j_2<j_2+1$) we have
\begin{equation}\label{equ: cube Ext1 vanishing socle}
\mathrm{cosoc}_{G}(i_{\widehat{j}_1,\Delta}^{\infty}(\pi_{j_1,j_2+1}^{\infty}))\cong V_{[j_2-j_1+1,j_2],\Delta}^{\infty}\cong \mathrm{soc}_{G}(i_{\widehat{j}_1,\Delta}^{\infty}(\pi_{j_1,j_2}^{\infty})),
\end{equation}
hence $V_{[j_2-j_1+1,j_2],\Delta}^{\infty}$ is a subquotient of $i_{I,\Delta}^{\infty}(\pi^{\infty})$. By \ref{it: connect 1} of Lemma~\ref{lem: connect Hom} and \ref{it: distance St 1} Lemma~\ref{lem: distance from St} we have
\[d(\pi_{j_1,j_2+1}^{\infty},\pi_{j_1-1,j_2}^{\infty})=d(\pi_{j_1,j_2+1}^{\infty},\pi_{j_1+1,j_2+1}^{\infty})=d(\pi_{j_1,j_2+1}^{\infty},V_{[j_2-j_1+1,j_2],\Delta}^{\infty})=0.\]
Hence we can apply Lemma~\ref{lem: general construction} with $\sigma_x^{\infty}\defeq \pi_{j_1,j_2+1}^{\infty}$ if $x=w_{j_1,j_0}$, $\sigma_x^{\infty}\defeq\pi_{j_1-1,j_2}^{\infty}$ if $x=w_{j_1-1,j_0}$, $\sigma_x^{\infty}\defeq\pi_{j_1+1,j_2+1}^{\infty}$ if $x=w_{j_1+1,j_0}$ and $\sigma_x^{\infty}\cong V_{[j_2-j_1+1,j_2],\Delta}^{\infty}$ if $j_0=j_1$ and $x=1$. This gives that $\tld{B}_{\un{j},S}$ is a subquotient of $\cF_{P_I}^G(M,\pi^{\infty})$. Since $M$ is multiplicity free and $i_{I,I_x}^{\infty}(\pi^{\infty})$ is multiplicity free for any constituent $L(x)$ of $M$ (see \ref{it: basic as image} of Remark \ref{rem: basic PS intertwine}), $\cF_{P_I}^G(M,\pi^{\infty})$ is also multiplicity free by Lemma \ref{lem: sm to OS} and Lemma \ref{lem: Hom OS}. We can thus define $\tld{B}_{\un{j},S}^+$ as the minimal length quotient of $\cF_{P_I}^G(M,\pi^{\infty})$ which admits $\tld{B}_{\un{j},S}$ as a subrepresentation. By a similar argument as in the paragraph before {Step $1$} of the proof of Lemma~\ref{lem: general construction}, $\tld{B}_{\un{j},S}^+$ admits a decreasing filtration indexed by $\mathrm{JH}_{U(\fg)}(M)$ with $L(x)$-graded piece given by $\cF_{P_{I_x}}^G(L(x),\sigma_{x,+}^{\infty})$, where $\sigma_{x,+}^{\infty}$ is the minimal length quotient of $i_{I,I_x}^{\infty}(\pi^{\infty})$ which admits $\sigma_x^{\infty}$ as a subrepresentation. For $L(x)\in\mathrm{JH}_{U(\fg)}(M)$, we have from the previous discussion
\begin{itemize}
\item if $x=w_{j_1-1,j_0}$, $I=I_{-}$ and $\sigma_x^{\infty}\cong \pi_{j_1-1,j_2}^{\infty}\cong \mathrm{cosoc}(i_{I,\Delta\setminus\{j_1-1\}}^{\infty}(\pi^{\infty}))$, thus $\sigma_{x,+}^{\infty}=\sigma_x^{\infty}$;
\item if $x=w_{j_1+1,j_0}$, $I=I_{+}$ and $\sigma_x^{\infty}\cong \pi_{j_1+1,j_2+1}^{\infty}\cong \mathrm{cosoc}(i_{I,\Delta\setminus\{j_1+1\}}^{\infty}(\pi^{\infty}))$, thus $\sigma_{x,+}^{\infty}=\sigma_x^{\infty}$;
\item if $x=w_{j_1,j_0}$, $\sigma_x^{\infty}\cong \pi_{j_1,j_2+1}^{\infty}\cong \mathrm{soc}_{L_{\widehat{j}_1}}(i_{I,\widehat{j}_1}^{\infty}(\pi^{\infty}))$, thus $\sigma_{x,+}^{\infty}=i_{I,\widehat{j}_1}^{\infty}(\pi^{\infty})$;
\item If $x=1$, $\sigma_{x,+}^{\infty}$ is the unique quotient of $i_{I,\Delta}^{\infty}(\pi^{\infty})$ with socle $\sigma_x^{\infty}\cong V_{[j_2-j_1+1,j_2],\Delta}^{\infty}$.
\end{itemize}
We need to make $\sigma_{x,+}^{\infty}$ a bit more explicit when $x=1$. We have isomorphisms
\begin{multline*}
\Hom_{G}(i_{I,\Delta}^{\infty}(\pi^{\infty}),i_{\widehat{j}_1,\Delta}^{\infty}(\pi_{j_1,j_2}^{\infty}))\cong \Hom_{L_{\widehat{j}_1}}(J_{\Delta,\widehat{j}_1}(i_{I,\Delta}^{\infty}(\pi^{\infty})), \pi_{j_1,j_2}^{\infty})\\
\cong \Hom_{L_{\widehat{j}_1}}(i_{I,\widehat{j}_1,s_{j_1}}^{\infty}(J_{I,\widehat{j}_1,s_{j_1}}(\pi^{\infty})), \pi_{j_1,j_2}^{\infty})\cong \Hom_{L_{I}}(J_{I,\widehat{j}_1,s_{j_1}}(\pi^{\infty}),J_{\widehat{j}_1,I_{\pm}}'(\pi_{j_1,j_2}^{\infty}))
\end{multline*}
where the first isomorphism follows from (\ref{equ: first adjunction}), the second from \ref{it: sml1} and \ref{it: sml2} of Lemma~\ref{lem: smooth geometric lemma} and from $\Sigma_{j_1,j_2}\subseteq W(L_{\widehat{j}_1})s_{{j}_1}W(L_I)\cdot \cJ(\pi^\infty)$ (which follows from (\ref{equ: cube coset intersection}) and $\Sigma= W(L_I)\cdot \cJ(\pi^\infty)$), and the last from $s_{j_1}(I)\cap \widehat{j}_1=s_{j_1}(\widehat{j}_1)\cap \widehat{j}_1=I_{\pm}$ (see (\ref{equ: Bruhat induction}) and (\ref{equ: I+-})) followed by (\ref{equ: second adjunction}). Moreover these spaces are all non-zero since we have
\[J_{\widehat{j}_1,I_{\pm}}'(\pi_{j_1,j_2}^{\infty})_{\cB^{I}_{\Sigma_{\pm}}}\cong J_{\widehat{j}_1,\widehat{j}_1,s_{j_1}}(\pi_{j_1,j_2+1}^{\infty})_{\cB^{I}_{\Sigma_{\pm}}}\cong J_{I_{\ast},\widehat{j}_1,s_{j_1}}(\pi_{\ast,0}^{\infty})_{\cB^{I}_{\Sigma_{\pm}}}\ne 0\]
by (v) above Lemma \ref{lem: coset square} and \ref{it: coset square 2} of Lemma~\ref{lem: coset square}. Since $\sigma_x^{\infty}\cong V_{[j_2-j_1+1,j_2],\Delta}^{\infty}$ (when $x=1$) is the socle of $i_{\widehat{j}_1,\Delta}^{\infty}(\pi_{j_1,j_2}^{\infty})$ by (\ref{equ: cube Ext1 vanishing socle}) (and since all representations are multiplicity free by \ref{it: basic as image} of Remark \ref{rem: basic PS intertwine}), we finally deduce that $\sigma_{x,+}^{\infty}$ injects into $i_{\widehat{j}_1,\Delta}^{\infty}(\pi_{j_1,j_2}^{\infty})$.\bigskip

We divide the rest of the proof into two steps.\bigskip

\textbf{Step $1$}: We prove that
\begin{equation}\label{equ: conj Ext1 vanishing OS OS}
\mathrm{Ext}_{G}^1(\cF_{P_I}^G(M,\pi^{\infty}), C_{\un{j}})=0.
\end{equation}
We first claim that $\cF_{P_I}^G(M,\pi^{\infty})$ has no constituent isomorphic to $\cF_{P_{\widehat{j}_1}}^{G}(L(x),\pi_{j_1,j_2}^{\infty})$ for any $x$ such that $D_L(x)=\{j_1\}$. Assume on the contrary that $W$ is such a constituent, then by Lemma~\ref{lem: Hom OS} (and \ref{it: OS property 2} of Theorem \ref{prop: OS property}) $W$ must be a constituent of $\cF_{P_{\widehat{j}_1}}^G(L(x),i_{I,\widehat{j}_1}^\infty(\pi^{\infty}))$, and thus $\pi_{j_1,j_2}^{\infty}$ is a constituent of $i_{I,\widehat{j}_1}^\infty(\pi^{\infty})$ (\ref{it: OS property 4} of Theorem \ref{prop: OS property}). This implies $\Sigma_{j_1,j_2}=W(L_{\widehat{j}_1})\cdot \cJ(i_{I,\widehat{j}_1}^\infty(\pi^{\infty}))$ since both are single $W(L_{\widehat{j}_1})$-cosets by the last statement in \ref{it: PS 1} of Lemma \ref{lem: Jacquet of PS}. By (\ref{equ: support of induction}) this implies $\Sigma=W(L_I)\cdot \cJ(\pi^\infty)\subseteq \Sigma_{j_1,j_2}$. But by (\ref{equ: sigmainter}) we have $\Sigma\subseteq \Sigma_{j_1,j_2+1}$, and since $\Sigma_{j_1,j_2+1}\cap \Sigma_{j_1,j_2}=\emptyset$ by \ref{it: coset 1} of Lemma~\ref{lem: coset intersection}, we derive a contradiction. Now, since any constituent $L(x)$ of $N^{\widehat{j}_1}(w_{j_1,j_0})$ is such that $D_L(x)=\{j_1\}$ by Lemma \ref{lem: dominance and left set}, we deduce in particular
\[\Hom_{G}\big(\cF_{P_I}^G(M,\pi^{\infty}),\cF_{P_{\widehat{j}_1}}^{G}(N^{\widehat{j}_1}(w_{j_1,j_0}),\pi_{j_1,j_2}^{\infty})\big)=0,\]
which together with $0\rightarrow N^{\widehat{j}_1}(w_{j_1,j_0})\rightarrow M^{\widehat{j}_1}(w_{j_1,j_0}) \rightarrow L(w_{j_1,j_0})\rightarrow 0$ gives an embedding
\begin{equation}\label{equ: conj Ext1 vanishing OS embedding}
\mathrm{Ext}_{G}^1(\cF_{P_I}^G(M,\pi^{\infty}), C_{\un{j}})\hookrightarrow \mathrm{Ext}_{G}^1\big(\cF_{P_I}^G(M,\pi^{\infty}), \cF_{P_{\widehat{j}_1}}^{G}(M^{\widehat{j}_1}(w_{j_1,j_0}),\pi_{j_1,j_2}^{\infty})\big).
\end{equation}
By Remark~\ref{rem: Ext 1 conjugate vanshing} (with $j$, $w$, $S_0$ there being $j_1$, $w_{j_1,j_0}$, $S_{j_0, j_1}$) we have for $\ell\leq 1$
\begin{equation}\label{equ: cube Ext1 vanishing g}
\mathrm{Ext}_{U(\fg)}^{\ell}(M^{\widehat{j}_1}(w_{j_1,j_0}), M^{s_{j_1}})=0.
\end{equation}
By (\ref{equ: cube coset intersection}) we can apply (\ref{equ: OS main seq}), which gives using (\ref{equ: cube Ext1 vanishing g})
\[\mathrm{Ext}_{G}^1\big(\cF_{P_I}^G(M,\pi^{\infty}), \cF_{P_{\widehat{j}_1}}^{G}(M^{\widehat{j}_1}(w_{j_1,j_0}),\pi_{j_1,j_2}^{\infty})\big)=0.\]
By (\ref{equ: conj Ext1 vanishing OS embedding}) this gives (\ref{equ: conj Ext1 vanishing OS OS}).\bigskip

\textbf{Step $2$}: We prove that the injection $\tld{B}_{\un{j},S}\hookrightarrow \tld{B}_{\un{j},S}^+$ and the surjection $\cF_{P_I}^G(M,\pi^{\infty})\twoheadrightarrow \tld{B}_{\un{j},S}^+$ induce isomorphisms
\begin{equation}\label{equ: conj Ext1 vanishing OS subquotient}
\mathrm{Ext}_{G}^1(\tld{B}_{\un{j},S}, C_{\un{j}})\xrightarrow{\sim} \mathrm{Ext}_{G}^1(\tld{B}_{\un{j},S}^+, C_{\un{j}}) \xleftarrow{\sim} \mathrm{Ext}_{G}^1(\cF_{P_I}^G(M,\pi^{\infty}), C_{\un{j}}).
\end{equation}
We wish to apply Remark~\ref{rem: dual minimal OS} to $V_0=C_{\un{j}}$ and $V=\tld{B}_{\un{j},S}$, so we check that all assumptions there are satisfied (using the notation of \emph{loc.~cit.}). From the definition of $M$, Lemma \ref{lem: dominance and left set} and \cite[Thm.~9.4(c)]{Hum08}, we have $\mathrm{JH}_{U(\fg)}(M)\cap \mathrm{JH}_{U(\fg)}(M^{\widehat{j}_1}(w_{j_1,j_0}))=\{L(w_{j_1,j_0})\}$. By \ref{it: connect 1} of Lemma~\ref{lem: connect Hom} and \ref{it: distance St 1} Lemma~\ref{lem: distance from St} we have $d(\sigma_x^\infty,\pi_{j_1,j_0}^\infty)=0$ for all constituents $L(x)$ of $M$. By \ref{it: coset 1} of Lemma~\ref{lem: coset intersection} we have $\Sigma_{j_1,j_2+1}\cap \Sigma_{j_1,j_2}=\emptyset$. Hence it remains to check $J_{\widehat{j}_1}(\tau_x^{\infty})_{\cB^{\widehat{j}_1}_{\Sigma_{j_1,j_2}}}\!\!=0$ for each $L(x)\in\mathrm{JH}_{U(\fg)}(M)$ with $x\neq w_{j_1,j_0}$ and each $\tau_x^{\infty}\in\mathrm{JH}_{L_{I_x}}(\sigma_{x,+}^{\infty})\setminus\{\sigma_x^{\infty}\}$. Since $\sigma_{x,+}^{\infty}=\sigma_x^{\infty}$ when $x\in \{w_{j_1-1,j_0},w_{j_1+1,j_0}\}$ by the explicit description of $\sigma_{x,+}^{\infty}$ before Step $1$, such a $\tau_x^{\infty}$ only exists when $x=1$. In this case $\sigma_x^{\infty}\cong V_{[j_2-j_1+1,j_2],\Delta}^{\infty}$ and $\tau_x^{\infty}\in \mathrm{JH}_{G}(i_{\widehat{j}_1,\Delta}^{\infty}(\pi_{j_1,j_2}^{\infty}))\setminus \{V_{[j_2-j_1+1,j_2],\Delta}^{\infty}\}$ by the discussion just before Step $1$. We claim that we have isomorphisms
\begin{equation}\label{equ: conj Ext1 vanishing OS Jacquet}
J_{\Delta,\widehat{j}_1}(V_{[j_2-j_1+1,j_2],\Delta}^{\infty})_{\cB^{\widehat{j}_1}_{\Sigma_{j_1,j_2}}}\cong \pi_{j_1,j_2}^{\infty}\cong J_{\Delta,\widehat{j}_1}(i_{\widehat{j}_1,\Delta}^{\infty}(\pi_{j_1,j_2}^{\infty}))_{\cB^{\widehat{j}_1}_{\Sigma_{j_1,j_2}}}.
\end{equation}
Applying (\ref{equ: first adjunction}) to $V_{[j_2-j_1+1,j_2],\Delta}^{\infty}\hookrightarrow i_{\widehat{j}_1,\Delta}^{\infty}(\pi_{j_1,j_2}^{\infty})$ and using the last assertion of Lemma \ref{lem: Jacquet basic} gives the first isomorphism. The second isomorphism is the very last isomorphism in \ref{it: sml2} of Lemma~\ref{lem: smooth geometric lemma} applied with $I$, $I_1$, $I_0$, $w$, $\pi_0^\infty$ there being $\Delta$, $\widehat{j}_1$, $\widehat{j}_1$, $1$, $\pi_{j_1,j_2}^{\infty}$. By exactness of $J_{\Delta,\widehat{j}_1}$ we deduce from (\ref{equ: conj Ext1 vanishing OS Jacquet}) $J_{\Delta,\widehat{j}_1}(\tau_x^{\infty})_{\cB^{\widehat{j}_1}_{\Sigma_{j_1,j_2}}}\!\!=0$ for $\tau_x$ as above (and $x=1$). By Remark~\ref{rem: dual minimal OS} we thus have (\ref{equ: conj Ext1 vanishing OS subquotient}). Finally, (\ref{equ: conj Ext1 vanishing OS}) follows from (\ref{equ: conj Ext1 vanishing OS subquotient}) and (\ref{equ: conj Ext1 vanishing OS OS}).
\end{proof}

\begin{prop}\label{prop: hard square}
Let $\un{j}\in\mathbf{J}$ such that $\un{j}'=(j'_0, j'_1, j'_2)\defeq (j_0, j_1, j_2+1)$ is still in $\mathbf{J}$.
\begin{enumerate}[label=(\roman*)]
\item \label{it: hard square 1} If $j_0\ne j_1$ there exists a unique minimal $\mathrm{Ext}$-square $V_{\un{j},\un{j}'}$ such that
\[\left\{\begin{array}{rcl}
\mathrm{gr}^0(V_{\un{j},\infty})&\cong &C_{(j_0,j_1,j_2+1)}\\
\mathrm{gr}^1(V_{\un{j},\infty})&\cong &C_{(j_0,j_1-1,j_2)}\oplus C_{(j_0,j_1+1,j_2+1)}\\
\mathrm{gr}^2(V_{\un{j},\infty})&\cong &C_{(j_0,j_1,j_2)}.
\end{array}\right.\]
\item \label{it: hard square 2} If $j_0= j_1$ there exists a unique minimal $\mathrm{Ext}$-square $V_{\un{j},\un{j}'}$ such that
\[\left\{\begin{array}{rcl}
\mathrm{gr}^0(V_{\un{j},\infty})&\cong &C_{(j_0,j_1,j_2+1)}\\
\mathrm{gr}^1(V_{\un{j},\infty})&\cong &C_{(j_0,j_1-1,j_2)}\oplus C_{(j_0,j_1+1,j_2+1)}\oplus L(1)^\vee\otimes_EV_{[j_2-j_0+1,j_2],\Delta}^{\infty}\\
\mathrm{gr}^2(V_{\un{j},\infty})&\cong &C_{(j_0,j_1,j_2)}.
\end{array}\right.\]
\end{enumerate}
In both \ref{it: hard square 1} and \ref{it: hard square 2} we omit $C_{(j_0,j_1+1,j_2+1)}$ if $j_1=j_2=n-1$ and $C_{(j_0,j_1-1,j_2)}$ if $j_1=1$.
\end{prop}
\begin{proof}
Note that, with the notation before Lemma~\ref{lem: square Ext1 vanishing}, we have $\mathrm{gr}^1(V_{\un{j},\un{j}'})\cong A_{\un{j}}\oplus B_{\un{j}}$, and we prove \ref{it: hard square 1} and \ref{it: hard square 2} simultaneously. We define $W_{\un{j}}^+$ as the unique (up to isomorphism) representation with socle $C_{\un{j}}$ which fits into an exact sequence $0\rightarrow C_{\un{j}}\rightarrow W_{\un{j}}^+\rightarrow A_{\un{j}}\oplus B_{\un{j}}\rightarrow 0$. Note that the existence and unicity of $W_{\un{j}}^+$ follows from Lemma~\ref{lem: Ext1 factor 1} and \ref{it: Ext1 with alg 1} of Lemma~\ref{lem: Ext1 factor 2}, and that we have an exact sequence $0\rightarrow W_{\un{j}}\rightarrow W_{\un{j}}^+\rightarrow A_{\un{j}}\rightarrow 0$ (see before Lemma \ref{lem: square Ext1 vanishing} for $W_{\un{j}}$). An obvious d\'evissage gives
\begin{equation}\label{equ: hard square upper bound}
\Dim_E \mathrm{Ext}_{G}^1(C_{\un{j}'}, W_{\un{j}}^+)\leq \Dim_E \mathrm{Ext}_{G}^1(C_{\un{j}'}, W_{\un{j}})+ \Dim_E \mathrm{Ext}_{G}^1(C_{\un{j}'}, A_{\un{j}})=0+1=1
\end{equation}
where the last equality follows from (\ref{equ: square Ext1 vanishing}), Lemma~\ref{lem: Ext1 factor 1} and \ref{it: Ext1 with alg 2} of Lemma~\ref{lem: Ext1 factor 2}. The short exact sequence $0\rightarrow C_{\un{j}}\rightarrow W_{\un{j}}^+\rightarrow A_{\un{j}}\oplus B_{\un{j}}\rightarrow 0$ yields an exact sequence
\begin{equation}\label{equ: hard square seq}
\mathrm{Ext}_{G}^1(C_{\un{j}'},C_{\un{j}})\rightarrow \mathrm{Ext}_{G}^1(C_{\un{j}'},W_{\un{j}}^+)\rightarrow \mathrm{Ext}_{G}^1(C_{\un{j}'},A_{\un{j}}\oplus B_{\un{j}})\rightarrow \mathrm{Ext}_{G}^2(C_{\un{j}'},C_{\un{j}}).
\end{equation}
We have by Lemma~\ref{lem: Ext1 factor 1}
\begin{equation}\label{equ: hard square vanishing simple}
\mathrm{Ext}_{G}^1(C_{\un{j}'},C_{\un{j}})=0
\end{equation}
and by Lemma~\ref{lem: Ext1 factor 1} and \ref{it: Ext1 with alg 2} of Lemma~\ref{lem: Ext1 factor 2}
\begin{equation}\label{equ: hard square Ext1 middle}
\Dim_E\mathrm{Ext}_{G}^1(C_{\un{j}'}, A_{\un{j}}\oplus B_{\un{j}})=1+\#S_{j_0,j_1}.
\end{equation}
By \ref{it: Ext2 factor 1 2} of Lemma~\ref{lem: Ext2 factor 1} and the explicit description of $S_{j_0,j_1}$ below (\ref{equ: Sjj}) we have
\[\Dim_E \mathrm{Ext}_{G}^2(C_{\un{j}'},C_{\un{j}})\leq \#S_{j_0,j_1},\]
which together with (\ref{equ: hard square vanishing simple}), (\ref{equ: hard square seq}) and (\ref{equ: hard square Ext1 middle}) implies
\begin{equation*}
\Dim_E \mathrm{Ext}_{G}^1(C_{\un{j}'}, W_{\un{j}}^+)\geq 1+\#S_{j_0,j_1}-\#S_{j_0,j_1}=1.
\end{equation*}
With (\ref{equ: hard square upper bound}), we deduce $\Dim_E \mathrm{Ext}_{G}^1(C_{\un{j}'}, W_{\un{j}}^+)=1$. The unique (up to scalar) non-zero class in $\mathrm{Ext}_{G}^1(C_{\un{j}'}, W_{\un{j}}^+)$ determines a unique (up to isomorphism) representation of $G$ over $E$, which is moreover multiplicity free using Lemma \ref{lem: Hom OS}. We define $V_{\un{j},\un{j}'}$ as its unique subrepresentation with cosocle $C_{\un{j}'}$. It is clear from its definition and from Lemma~\ref{lem: Ext1 factor 1} and Lemma~\ref{lem: Ext1 factor 2} that $V_{\un{j},\un{j}'}$ is an $\mathrm{Ext}$-square. Moreover we have by (\ref{equ: square Ext1 vanishing}) (and since $A_{\un{j}}$ is irreducible)
\begin{equation}\label{equ: hard square middle layer bound}
A_{\un{j}}\subseteq \mathrm{gr}^1(V_{\un{j},\un{j}'})\subseteq A_{\un{j}}\oplus B_{\un{j}}.
\end{equation}
If the second inclusion in (\ref{equ: hard square middle layer bound}) is strict, there exists $S\subsetneq S_{j_0,j_1}$ such that $\mathrm{gr}^1(V_{\un{j},\un{j}'})\cong A_{\un{j}}\oplus B_{\un{j},S}$, and the existence of $V_{\un{j},\un{j}'}$ forces $\mathrm{Ext}_{G}^1(\tld{B}_{\un{j},S}, C_{\un{j}})\neq 0$, which contradicts (\ref{equ: conj Ext1 vanishing OS}). Hence we have $W_{\un j}\hookrightarrow V_{\un{j},\un{j}'}$. Finally, the minimality of $V_{\un{j},\un{j}'}$ follows from Lemma~\ref{lem: OS minimal criterion} (applied with $V_1= \mathrm{Fil}^1(V_{\un{j},\un{j}'})$, $V_2=W_{\un j}$) and Lemma \ref{lem: square Ext1 vanishing}.
\end{proof}

Until the end of the proof of Proposition~\ref{prop: factor cube}, we fix $\un{j}\in\mathbf{J}$ such that $\un{j}'=(j'_0, j'_1, j'_2)\!\defeq (j_0+1, j_1, j_2+1)$ is still in $\mathbf{J}$. Note that we have $1\leq j_0\leq n-2$, $1\leq j_1\leq j_2\leq n-1$ and $j_2-j_1\leq n-2-j_0$. The following notation is convenient
\[\begin{array}{cclcclccl}
C^{0,0,0}&\defeq &C_{(j_0,j_1,j_2)}&&&&&&\\
C^{1,0,0}&\defeq &C_{(j_0+1,j_1,j_2)}&C^{0,1,0}&\defeq &C_{(j_0,j_1+1,j_2+1)}& C^{0,0,1}&\defeq &C_{(j_0,j_1-1,j_2)}\\
C^{1,1,0}&\defeq &C_{(j_0+1,j_1+1,j_2+1)}& C^{0,1,1}&\defeq &C_{(j_0,j_1,j_2+1)}& C^{1,0,1}&\defeq &C_{(j_0+1,j_1-1,j_2)}\\
C^{1,1,1}&\defeq &C_{(j_0+1,j_1,j_2+1)}&&&&&&
\end{array}\]
(so $C^{0,0,0}=C_{\un j}$ and $C^{1,1,1}=C_{\un{j}'}$). By Proposition~\ref{prop: hard square} there exists a minimal $\mathrm{Ext}$-square $V_0$ with socle $C^{0,0,0}$, cosocle $C^{0,1,1}$ and middle layer containing $C^{0,1,0}\oplus C^{0,0,1}$. By \emph{loc.~cit.}~applied with $\un{j}$, $\un{j}'$ there being $(j_0+1,j_1,j_2)$, $(j_0+1,j_1,j_2+1)$, there exists a minimal $\mathrm{Ext}$-square $V_1$ with socle $C^{1,0,0}$, cosocle $C^{1,1,1}$ and middle layer containing $C^{1,1,0}\oplus C^{1,0,1}$. By \emph{loc.cit.}~$V_0$ (resp.~$V_1$) is multiplicity free and admits a locally algebraic constituent if and only if $j_1=j_0$ (resp.~$j_1=j_0+1$).\bigskip

Until Proposition~\ref{prop: factor cube} we also fix $\mu\in \Lambda$ such that $\langle\mu+\rho,\alpha^\vee\rangle\geq 0$ for all $\al\in \Phi^+$ and the stabilizer of $\mu$ in $W(G)$ for the dot action is $\{1,w_0s_{j_0}w_0\}$. The following lemma gives technical results on $\Theta_{w_0s_{j_0}w_0}(C^{0,b,c})$ for $b,c\in\{0,1\}$ which will be used in various proofs, in particular the proof of Proposition~\ref{prop: factor cube}.

\begin{lem}\label{lem: cube step 1}
For each $(b,c)\in\{0,1\}^2$ we have:
\begin{enumerate}[label=(\roman*)]
\item \label{it: cube step 1 1} $C^{a',b',c'}$ is a constituent of $\Theta_{w_0s_{j_0}w_0}(C^{0,b,c})$ for some $a',b',c'\in\{0,1\}$ if and only if $b'=b$ and $c'=c$;
\item \label{it: cube step 1 2} $C^{1,b,c}$ (resp.~$C^{0,b,c}$) appears with multiplicity $1$ (resp.~$2$) in $\Theta_{w_0s_{j_0}w_0}(C^{0,b,c})$ and in $\Theta_{\mu}(V_0)$;
\item \label{it: cube step 1 3} $\Theta_{w_0s_{j_0}w_0}(C^{0,b,c})$ has simple socle and cosocle $C^{0,b,c}$, and admits a unique length $2$ subrepresentation (resp.~quotient) with socle $C^{0,b,c}$ and cosocle $C^{1,b,c}$ (resp.~with socle $C^{1,b,c}$ and cosocle $C^{0,b,c}$).
\end{enumerate}
\end{lem}
\begin{proof}
We fix $b,c\in\{0,1\}$ and write $C^{0,b,c}=\cF_{P_{I_x}}^G(L(x),\sigma^{\infty})$ for some $x\in W(G)$ and some irreducible $G$-regular smooth $\sigma^{\infty}$, which can be made explicit using (\ref{cj}). In particular by Proposition \ref{prop: Jantzen middle} we have $\Theta_{w_0s_{j_0}w_0}(L(x))\ne 0$.

We prove \ref{it: cube step 1 1}. For $a',b',c'\in\{0,1\}$ we write $C^{a',b',c'}=\cF_{P_{I_w}}^G(L(w),\pi^{\infty})$ for some $w\in W(G)$ and some irreducible $G$-regular smooth $\pi^{\infty}$. Note that $\#D_L(x)=\#D_L(w)=1$, and thus $I_w\supseteq I_x$ if and only if $I_x=I_w$. Hence by (\ref{equ: OS wall crossing}) and Lemma~\ref{lem: Hom OS} $C^{a',b',c'}$ is a constituent of $\Theta_{w_0s_{j_0}w_0}(C^{0,b,c})$ if and only if $L(w)$ is a constituent of $\Theta_{w_0s_{j_0}w_0}(L(x))$, $I_x=I_w$ and $\sigma^{\infty}\cong \pi^{\infty}$. On one hand, $I_x=I_w$ and $\sigma^{\infty}\cong \pi^{\infty}$ force $b=b'$ and $c=c'$ (by definition of $C^{0,b,c}$ and $C^{a',b',c'}$). On the other hand, if $b=b'$ and $c=c'$, then $\sigma^{\infty}\cong \pi^{\infty}$ and a case by case check from the definition of $C^{0,b,c}$, $C^{a',b',c'}$ gives that $x=w_{j,j_0}$ for some $j\in\{j_1-1,j_1,j_1+1\}\cap\Delta$ and either $w=x$ or $w=w_{j,j_0+1}$. In all these cases $C^{a',b',c'}$ is a constituent of $\Theta_{w_0s_{j_0}w_0}(C^{0,b,c})$ by Proposition \ref{prop: Jantzen middle} (with Remark~\ref{rem: explicit middle}).

We \ prove \ \ref{it: cube step 1 2}. \ By \ the \ first \ statement \ in \ Proposition~\ref{prop: Jantzen middle} \ and \ (\ref{equ: OS wall crossing}) \ we \ have $\Theta_{w_0s_{j_0}w_0}(W)=0$ if $W$ is a locally algebraic constituent of $V_0$. As $V_0$ is multiplicity free with non locally algebraic constituents exactly given by $C^{0,b,c}$ for $b,c\in\{0,1\}$ (see Proposition~\ref{prop: hard square}), we deduce from the exactness of $\Theta_{\mu}$ that $\Theta_{\mu}(V_0)$ admits a filtration with graded pieces given by $\Theta_{w_0s_{j_0}w_0}(C^{0,b,c})$ for $b,c\in\{0,1\}$. This together with \ref{it: cube step 1 1} implies that $C^{1,b,c}$ (resp.~$C^{0,b,c}$) appears in $\Theta_{w_0s_{j_0}w_0}(C^{0,b,c})$ and $\Theta_{\mu}(V_0)$ with the same multiplicity. It then follows from the (second half of) the first sentence in \ref{it: square crossing 1} of Lemma~\ref{lem: square as wall crossing} that $C^{1,b,c}$ (resp.~$C^{0,b,c}$) appears with multiplicity one (resp.~two) in $\Theta_{w_0s_{j_0}w_0}(C^{0,b,c})$ (we apply \emph{loc.~cit.}~with $\un{j}$ there being $(j_0,j_1,j_2)$, $(j_0,j_1+1,j_2+1)$, $(j_0,j_1-1,j_2)$ and $(j_0,j_1,j_2+1)$ which corresponds to $(b,c)$ being $(0,0)$, $(1,0)$, $(0,1)$ and $(1,1)$).

Finally, \ref{it: cube step 1 3} follows from the rest of the statement in \ref{it: square crossing 1} of Lemma~\ref{lem: square as wall crossing}.
\end{proof}

For $\lambda, \lambda'\in \Lambda$ and $C, D$ any $Z(\fg)$-finite $D(G)$-modules, by \cite[Thm.~2.4.7]{JLS24} we have canonical isomorphisms
\begin{equation}\label{adjunctionlocadm}
\Hom_{D(G)}(\cT^{\lambda'}_{\lambda}(C), D)\cong \Hom_{D(G)}(C, \cT^{\lambda}_{\lambda'}(D)),
\end{equation}
from which we (formally) obtain canonical functorial adjunction maps $\Theta_{\lambda}(D)\rightarrow D$, $D\rightarrow\Theta_{\lambda}(D)$ for $D$ any $Z(\fg)$-finite $D(G)$-module and $\lambda\in \Lambda$ such that $\langle \lambda + \rho, \alpha^\vee\rangle \geq 0$ for $\alpha\in \Phi^+$ and the stabilizer of $\lambda$ in $W(G)$ for the dot action is $\{1, s_{j}\}$ for some $j\in \{0,\dots,n-1\}$. It is clear that these two adjunction maps are non-zero when both $D$ and $\Theta_{\lambda}(D)$ are non-zero. If $V$ is an admissible locally analytic representation of $G$ over $E$ such that $V^\vee$ is $Z(\fg)$-finite, we therefore have canonical functorial adjunctions maps $V\rightarrow\Theta_{\lambda}(V)$ and $\Theta_{\lambda}(V)\rightarrow V$. Note that, when $D=\cF_{P_I}^{G}(M,\pi^{\infty})^\vee$ for $I\subseteq \Delta$, $M$ in $\cO^{\fp_I}_{\rm{alg}}$ and $\pi^{\infty}$ a strongly admissible smooth representation of $L_I$, an examination of the proof of (\ref{equ: OS wall crossing}) in \cite[Thm.~4.1.12]{JLS24} shows that the adjunction maps $\Theta_{\lambda}(D)=\Theta_{s_{j}}(D)\rightarrow D$ coming from (i) the above argument and (ii) the adjunction map $\Theta_{s_{j}}(M)\rightarrow M$ (see below (\ref{equ: translation functor})) by functoriality of the Orlik-Strauch functor are the same. Likewise with the adjunction maps $D\rightarrow\Theta_{s_{j}}(D)$.\bigskip

Going back to our previous running notation, since $C^{1,1,1}$ occurs with multiplicity $1$ in $\Theta_{\mu}(V_0)$ by \ref{it: cube step 1 2} of Lemma~\ref{lem: cube step 1}, we can (and do) define $V_{\un{j},\un{j}'}$ as the unique subrepresentation of $\Theta_{\mu}(V_0)$ with cosocle $C^{1,1,1}$.

\begin{lem}\label{lem: cube step 2}
The representations $\Theta_{\mu}(V_0)$ and $V_{\un{j},\un{j}'}$ have socle $C^{0,0,0}$. Moreover, $V_{\un{j},\un{j}'}$ satisfies the following properties:
\begin{enumerate}[label=(\roman*)]
\item \label{it: cube step 2 1} $V_0$ injects into $V_{\un{j},\un{j}'}$;
\item \label{it: cube step 2 2} for each $(b,c), (b',c')\in \{0,1\}^2$ such that $b\leq b'$, $c\leq c'$ and $b'+c'=b+c+1$, $\Theta_{\mu}(V_0)$ admits a subquotient which is the minimal $\mathrm{Ext}$-square (with socle $C^{0,b,c}$, cosocle $C^{1,b',c'}$ and middle layer containing $C^{0,b',c'}\oplus C^{1,b,c}$) constructed in \ref{it: easy square 1} or \ref{it: easy square 2} of Proposition~\ref{prop: easy square};
\item \label{it: cube step 2 3} for each $(a,b,c)\in \{0,1\}^3$, $C^{a,b,c}$ appears in $V_{\un{j},\un{j}'}$ with multiplicity $1$;
\item \label{it: cube step 2 4} for each $(a,b,c), (a',b',c')\in \{0,1\}^3$ such that $a\leq a'$, $b\leq b'$, $c\leq c'$ and $a'+b'+c'=a+b+c+1$, $V_{\un{j},\un{j}'}$ admits a unique length $2$ subquotient with socle $C^{a,b,c}$ and cosocle $C^{a',b',c'}$.
\end{enumerate}
\end{lem}
\begin{proof}
By the first statement in \ref{it: square crossing 1} of Lemma~\ref{lem: square as wall crossing}, for $b,c\in\{0,1\}$ we have (up to non-zero scalars) a canonical injection $C^{0,b,c}\hookrightarrow \Theta_{w_0s_{j_0}w_0}(C^{0,b,c})$ and a canonical surjection $\Theta_{w_0s_{j_0}w_0}(C^{0,b,c})\twoheadrightarrow C^{0,b,c}$. Since $V_0$ has socle $C^{0,0,0}$ and cosocle $C^{0,1,1}$, it follows from the exactness and functoriality of $\Theta_{\mu}$ that the adjunction map $V_0\rightarrow \Theta_{\mu}(V_0)$ is injective and the adjunction map $\Theta_{\mu}(V_0)\rightarrow V_0$ is surjective. As the composition $C^{0,b,c}\hookrightarrow \Theta_{w_0s_{j_0}w_0}(C^{0,b,c})\twoheadrightarrow C^{0,b,c}$ is (obviously) zero for $b,c\in\{0,1\}$ and as $\Theta_{\mu}$ kills any locally algebraic constituent of $V_0$ (by the first statement in Proposition~\ref{prop: Jantzen middle} and (\ref{equ: OS wall crossing})), by functoriality and exactness again of $\Theta_{\mu}$ the composition $V_0\rightarrow \Theta_{\mu}(V_0)\rightarrow V_0$ is also zero. As $C^{1,1,1}$ is not a constituent of $V_0$ and as $C^{0,b,c}$ appears in $\Theta_{\mu}(V_0)$ (resp.~$V_0$) with multiplicity $2$ (resp.~$1$) for each $(b,c)\in\{0,1\}^2$ by \ref{it: cube step 1 2} of Lemma \ref{lem: cube step 1}, we deduce that $V_{\un{j},\un{j}'}$ injects into the kernel of the surjection $\Theta_{\mu}(V_0)\twoheadrightarrow V_0$ and that each $C^{0,b,c}$ appears in $V_{\un{j},\un{j}'}$ with multiplicity at most one.

We prove that $\Theta_{\mu}(V_0)$ has socle $C^{0,0,0}$ (and thus its subrepresentation $V_{\un{j},\un{j}'}$ also has socle $C^{0,0,0}$).
We choose an arbitrary irreducible $W\hookrightarrow \mathrm{soc}_{G}(\Theta_{\mu}(V_0))$. By the discussion before this lemma, we have canonical isomorphisms
\begin{equation}\label{equ: wall crossing adjunction}
0\neq \Hom_{G}(W,\Theta_{\mu}(V_0))\cong \Hom_{G}(\cT^{\mu}_{w_0\cdot\mu_0}(W), \cT^{\mu}_{w_0\cdot\mu_0}(V_0))\cong \Hom_{G}(\Theta_{w_0s_{j_0}w_0}(W),V_0).
\end{equation}
In particular, we have $\Theta_{w_0s_{j_0}w_0}(W)\neq 0$ and $W$ is not locally algebraic by Theorem~\ref{prop: OS property} and (\ref{equ: OS wall crossing}). By Lemma~\ref{lem: wall crossing of simple} $\Theta_{w_0s_{j_0}w_0}(W)$ has simple cosocle $W$, which together with (\ref{equ: wall crossing adjunction}) forces $W$ to be a non locally algebraic factor of $V_0$, and thus $W=C^{0,b,c}$ for some $b,c\in\{0,1\}$. But by \ref{it: cube step 1 1} of Lemma~\ref{lem: cube step 1} we know that the socle $C^{0,0,0}$ of $V_0$ is a constituant of $\Theta_{w_0s_{j_0}w_0}(C^{0,b,c})$ if and only if $b=c=0$. Hence, a non-zero map $\Theta_{w_0s_{j_0}w_0}(W)\rightarrow V_0$ exists only if $W=C^{0,0,0}$, in which case this map must factor through $\Theta_{w_0s_{j_0}w_0}(W)\twoheadrightarrow C^{0,0,0}=\mathrm{soc}_{G}(V_0)$ and (\ref{equ: wall crossing adjunction}) is one dimensional by \ref{it: cube step 1 3} of Lemma~\ref{lem: cube step 1}. It follows from this discussion that $\Theta_{\mu}(V_0)$ has simple socle $C^{0,0,0}$.

We prove \ref{it: cube step 2 1}. The surjection $V_0\twoheadrightarrow C^{0,1,1}$ induces a surjection $\Theta_{\mu}(V_0)\twoheadrightarrow \Theta_{w_0s_{j_0}w_0}(C^{0,1,1})$ such \ that \ (by \ functoriality \ of \ the \ adjunction \ map) \ the \ composition \ $V_0\hookrightarrow \Theta_{\mu}(V_0)\twoheadrightarrow \Theta_{w_0s_{j_0}w_0}(C^{0,1,1})$ coincides with the composition
\[V_0\twoheadrightarrow C^{0,1,1}\cong \soc_G(\Theta_{w_0s_{j_0}w_0}(C^{0,1,1})\hookrightarrow \Theta_{w_0s_{j_0}w_0}(C^{0,1,1})\]
(see \ref{it: cube step 1 3} of Lemma~\ref{lem: cube step 1} for the above isomorphism). By \emph{loc.~cit.}~$\Theta_{w_0s_{j_0}w_0}(C^{0,1,1})$ admits a unique length $2$ subrepresentation with socle $C^{0,1,1}$ and cosocle $C^{1,1,1}$, and we denote by $\widetilde V$ its inverse image in $\Theta_{\mu}(V_0)$ via $\Theta_{\mu}(V_0)\twoheadrightarrow \Theta_{w_0s_{j_0}w_0}(C^{0,1,1})$. Then $C^{0,1,1}$ has multiplicity $1$ in $\widetilde V$ since it has multiplicity $2$ in $\Theta_{\mu}(V_0)$ and $\Theta_{w_0s_{j_0}w_0}(C^{0,1,1})$ (by \ref{it: cube step 1 2} of Lemma~\ref{lem: cube step 1}). As $V_0$ has cosocle $C^{0,1,1}$ and $V_0\subseteq \widetilde V$ by the above discussion, it follows that $V_0$ is the unique subrepresentation of $\widetilde V$ with cosocle $C^{0,1,1}$. Moreover any subrepresentation of $\widetilde V$ which contains $C^{0,1,1}$ as a constituent also contains the subrepresentation $V_0$. Since $V_{\un{j},\un{j}'}$ has cosocle $C^{1,1,1}$ (which appears with multiplicity $1$ in $\Theta_{\mu}(V_0)$), its image in $\Theta_{w_0s_{j_0}w_0}(C^{0,1,1})$ is the unique subrepresentation of cosocle $C^{1,1,1}$, in particular $V_{\un{j},\un{j}'}\subseteq \widetilde V$, and thus $V_0\subseteq V_{\un{j},\un{j}'}$. Hence, by the end of the last sentence of the first paragraph of the proof, we also deduce that each $C^{0,b,c}$ for $b,c\in \{0,1\}$ appears in $V_{\un{j},\un{j}'}$ with multiplicity exactly one.

We prove \ref{it: cube step 2 2}. Let $(b,c), (b',c')\in \{0,1\}^2$ such that $b\leq b'$, $c\leq c'$ and $b'+c'=b+c+1$, then $V_0$ admits a unique length $2$ subquotient $V$ with socle $C^{0,b,c}$ and cosocle $C^{0,b',c'}$. It follows from \ref{it: square crossing 2} of Lemma~\ref{lem: square as wall crossing} that $\Theta_{\mu}(V)$ admits a unique subquotient which is the $\mathrm{Ext}$-square constructed in \ref{it: easy square 1} or \ref{it: easy square 2} of Proposition~\ref{prop: easy square} with socle $C^{0,b,c}$, cosocle $C^{1,b',c'}$ and middle layer containing $C^{1,b,c}\oplus C^{0,b',c'}$. Since $\Theta_{\mu}(V)$ is a subquotient of $\Theta_{\mu}(V_0)$, this proves \ref{it: cube step 2 2}. In particular $\Theta_{\mu}(V_0)$ admits a length $2$ subquotient with socle $C^{1,b,c}$ and cosocle $C^{1,b',c'}$.

We prove \ref{it: cube step 2 3}. Let $(b,c), (b',c')\in \{0,1\}^2$ such that $b\leq b'$, $c\leq c'$, and recall from \ref{it: cube step 1 2} of Lemma~\ref{lem: cube step 1} that $C^{1,b,c}$ and $C^{1,b',c'}$ have multiplicity one in $\Theta_{\mu}(V_0)$. Hence, the end of the previous paragraph implies that $\Theta_{\mu}(V_0)$ admits a unique length $2$ subquotient with socle $C^{1,b,c}$ and cosocle $C^{1,b',c'}$ when $b'+c'=b+c+1$. Since $C^{1,1,1}$ appears in $V_{\un{j},\un{j}'}$, we deduce that $C^{1,0,1}$ and $C^{1,1,0}$, and then $C^{1,0,0}$, must all appear in $V_{\un{j},\un{j}'}$. So $V_{\un{j},\un{j}'}$ admits a unique length $2$ subquotient with socle $C^{1,b,c}$ and cosocle $C^{1,b',c'}$ when $b'+c'=b+c+1$. Together with the last sentence of the proof of \ref{it: cube step 2 1}, we deduce \ref{it: cube step 2 3}.

We prove \ref{it: cube step 2 4}. As $V_0$ injects into $V_{\un{j},\un{j}'}$ by \ref{it: cube step 2 1} and $V_0$ admits a unique length $2$ subquotient with socle $C^{0,b,c}$ and cosocle $C^{0,b',c'}$ when $b'+c'=b+c+1$, such a subquotient uniquely appears in $V_{\un{j},\un{j}'}$ (unicity uses multiplicity $1$ in \ref{it: cube step 2 3}). As $C^{1,b,c}$ has multiplicity one in $\Theta_{\mu}(V_0)$ and $V_{\un{j},\un{j}'}$ by the previous paragraph, $V_{\un{j},\un{j}'}$ contains the unique subrepresentation of $\Theta_{\mu}(V_0)$ with cosocle $C^{1,b,c}$. As $\Theta_{w_0s_{j_0}w_0}(C^{0,b,c})$ admits a length $2$ subrepresentation with socle $C^{0,b,c}$ and cosocle $C^{1,b,c}$ by \ref{it: cube step 1 3} of Lemma~\ref{lem: cube step 1}, this subrepresentation of $\Theta_{\mu}(V_0)$, and hence $V_{\un{j},\un{j}'}$, also admit a length $2$ subquotient with socle $C^{0,b,c}$ and cosocle $C^{1,b,c}$. Moreover, for $V_{\un{j},\un{j}'}$ this subquotient is unique again by multiplicity $1$ in \ref{it: cube step 2 3}. Together with the one but last sentence of the proof of \ref{it: cube step 2 3}, this gives all cases of \ref{it: cube step 2 4}.
\end{proof}

\begin{rem}\label{rem: one dim Hom}
Let $q: V_0\rightarrow \Theta_{\mu}(V_0)$ be any non-zero map, which is necessarily injective as $V_0$ is multiplicity free with socle $C^{0,0,0}$ and $\Theta_{\mu}(V_0)$ has socle $C^{0,0,0}$ by (the first statement of) Lemma~\ref{lem: cube step 2}, and denote by $p: \Theta_{\mu}(V_0)\twoheadrightarrow V_0$ the canonical surjection (cf.~the first paragraph of the proof of Lemma~\ref{lem: cube step 2}). If the composition $p\circ q: V_0\rightarrow V_0$ is non-zero, then it has to be an isomorphism as $V_0$ is multiplicity free with socle $C^{0,0,0}$. By exactness of $\Theta_\mu$ and \ref{it: cube step 1 3} of Lemma \ref{lem: cube step 1}, the restriction $(p\circ q)|_{\mathrm{soc}_{G}(V_0)}$ is the composition $C^{0,0,0}\hookrightarrow \Theta_{\mu}(C^{0,0,0})\twoheadrightarrow C^{0,0,0}$ which is zero, a contradiction to $p\circ q$ being an isomorphism. Hence, we have $p\circ q=0$ and thus $\mathrm{im}(q)\subseteq \mathrm{ker}(p)$. But from \ref{it: cube step 1 2} of Lemma~\ref{lem: cube step 1} $C^{0,1,1}$ has multiplicity one in $\mathrm{ker}(p)$, so $\mathrm{im}(q)$ has to be the unique subrepresentation of $\mathrm{ker}(p)$ with cosocle $C^{0,1,1}$. It follows that $q$ is unique up to a scalar, or equivalently $\Dim_E\Hom_{G}(V_0, \Theta_{\mu}(V_0))=1$.
\end{rem}

For each $(b,c)\in\{0,1\}^2$, we write $V_{b,c}$ for the unique subrepresentation of $V_0$ with cosocle $C^{0,b,c}$. By exactness of $\Theta_{\mu}$ and since $\Theta_{\mu}$ kills any locally algebraic constituent of $V_0$, $\Theta_{\mu}(V_0)$ admits a natural increasing filtration $\mathrm{Fil}_{b,c}(\Theta_{\mu}(V_0))\defeq \Theta_{\mu}(V_{b,c})$ with graded piece
\[\mathrm{gr}_{b,c}(\Theta_{\mu}(V_0))\defeq \mathrm{Fil}_{b,c}(\Theta_{\mu}(V_0))\Big/\!\!\!\sum_{b'+c'<b+c}\mathrm{Fil}_{b',c'}(\Theta_{\mu}(V_0))\cong \Theta_{w_0s_{j_0}w_0}(C^{0,b,c}).\]
The filtration $\{\mathrm{Fil}_{b,c}(\Theta_{\mu}(V_0))\}_{0\leq b,c\leq 1}$ induces a filtration $\{\mathrm{Fil}_{b,c}(V_{\un{j},\un{j}'})\}_{0\leq b,c\leq 1}$ on $V_{\un{j},\un{j}'}\subseteq \Theta_{\mu}(V_0)$ with graded piece $\mathrm{gr}_{b,c}(V_{\un{j},\un{j}'})\subseteq \Theta_{w_0s_{j_0}w_0}(C^{0,b,c})$. We deduce from \ref{it: cube step 1 2} and \ref{it: cube step 1 3} of Lemma~\ref{lem: cube step 1} combined with \ref{it: cube step 2 3} of Lemma \ref{lem: cube step 2} that we have inclusions for $(b,c)\in\{0,1\}^2$
\begin{equation}\label{equ: grade in rad}
\mathrm{gr}_{b,c}(V_{\un{j},\un{j}'})\subseteq \mathrm{rad}_{G}(\Theta_{w_0s_{j_0}w_0}(C^{0,b,c}))=\mathrm{ker}(\Theta_{w_0s_{j_0}w_0}(C^{0,b,c})\twoheadrightarrow C^{0,b,c}).
\end{equation}

\begin{lem}\label{lem: cube step loc alg}
Assume $j_1\neq j_0+1$. Then $V_{\un{j},\un{j}'}$ admits a locally algebraic constituent if and only if $j_0=j_1$. Moreover, if $j_0=j_1$ then $\mathrm{gr}_{b,c}(V_{\un{j},\un{j}'})$ (for $(b,c)\in\{0,1\}^2$) admits a locally algebraic constituent if and only if $b=c=0$.
\end{lem}
\begin{proof}
Assume that $V_{\un{j},\un{j}'}$ admits a locally algebraic constituent $W$, then there exists $(b,c)\in\{0,1\}^2$ such that $W$ occurs in $\mathrm{gr}_{b,c}(V_{\un{j},\un{j}'})\subseteq \Theta_{w_0s_{j_0}w_0}(C^{0,b,c})$.
Writing $C^{0,b,c}=\cF_{P_{I_x}}^G(L(x),\pi^{\infty})$ with $x\in W(G)$ and $\pi^{\infty}$ irreducible $G$-regular, it follows from (\ref{equ: OS wall crossing}) and Lemma~\ref{lem: Hom OS} (and \ref{it: OS property 2} of Theorem~\ref{prop: OS property}) that $L(1)$ is a constituent of $\Theta_{w_0s_{j_0}w_0}(L(x))$, which by Proposition~\ref{prop: Jantzen middle} forces $x=s_{j_0}$. By checking the definition of $C^{0,b,c}$ for $b,c\in\{0,1\}$ (see above Lemma \ref{lem: cube step 1} with (\ref{cj})) and using $j_1\neq j_0+1$, we have either $b=c$ with $j_1=j_0$, or $b=1$, $c=0$ with $j_0=j_1+1$.\bigskip

\textbf{Case $1$}: Assume $W$ occurs in $\mathrm{gr}_{1,1}(V_{\un{j},\un{j}'})$. By \ref{it: cube step 1 3} of Lemma~\ref{lem: cube step 1} $\Theta_{w_0s_{j_0}w_0}(C^{0,1,1})$ admits a unique length $2$ subrepresentation $V$ with socle $C^{0,1,1}$ and cosocle $C^{1,1,1}$, and by the proof of \ref{it: cube step 2 1} of Lemma \ref{lem: cube step 2}, $V$ is the image of $V_{\un{j},\un{j}'}$ via $\Theta_{\mu}(V_0)\twoheadrightarrow \Theta_{w_0s_{j_0}w_0}(C^{0,1,1})$, i.e.~$\mathrm{gr}_{1,1}(V_{\un{j},\un{j}'})\cong V$. In particular, $\mathrm{gr}_{1,1}(V_{\un{j},\un{j}'})$ has no locally algebraic constituent, a contradiction.\bigskip

\textbf{Case $2$}: Assume $W$ occurs in $\mathrm{gr}_{1,0}(V_{\un{j},\un{j}'})$ (and $j_0=j_1+1$). Then by Proposition~\ref{prop: Jantzen middle} and (\ref{equ: grade in rad}), $W$ occurs in the cosocle of $\mathrm{gr}_{1,0}(V_{\un{j},\un{j}'})$. As $V_{\un{j},\un{j}'}$ has cosocle $C^{1,1,1}$, there must exist $W'\in\mathrm{JH}_{G}(\mathrm{gr}_{1,1}(V_{\un{j},\un{j}'}))=\{C^{0,1,1},C^{1,1,1}\}$ (see \textbf{Case $1$} for the equality) such that $V_{\un{j},\un{j}'}$ admits a length $2$ subquotient $V'$ with socle $W$ and cosocle $W'$. Recall that $V_0$ is the unique subrepresentation of $V_{\un{j},\un{j}'}$ with cosocle $C^{0,1,1}$ (use \ref{it: cube step 2 1} and \ref{it: cube step 2 3} of Lemma~\ref{lem: cube step 2}). Thus, if $W'=C^{0,1,1}$ then $W$ must occur in $V_0$. But $V_0$ admits a locally algebraic factor if and only if $j_0=j_1$ (see Proposition~\ref{prop: hard square}), which contradicts $j_0=j_1+1$. Hence we must have $W'=C^{1,1,1}$. Since $\mathrm{Ext}_{U(\fg)}^1(L(1),L(w_{j_1,j_0+1}))=0$ when $j_0=j_1+1$ by \ref{it: rabiotext 2} of Lemma~\ref{rabiotext}, Proposition~\ref{prop: Ext1 OS} forces $\mathrm{Ext}_{G}^1(C^{1,1,1},W)=0$, which contradicts the existence of $V'$ above.

So the only remaining case is $b=c=0$ with $j_1=j_0$. This proves the ``only if'' of the two statements of the lemma. For the ``if'', we can use that $V_0$ has a locally algebraic constituent when $j_1=j_0$ by \ref{it: hard square 2} of Proposition \ref{prop: hard square}, hence so does $V_{\un{j},\un{j}'}$ by \ref{it: cube step 2 1} of Lemma~\ref{lem: cube step 2}.
\end{proof}

\begin{lem}\label{lem: cube step 3}
Assume $j_1\neq j_0+1$. The non-locally algebraic constituents of $V_{\un{j},\un{j}'}$ are the $C^{a,b,c}$ for $(a,b,c)\in\{0,1\}^3$.
\end{lem}
\begin{proof}
Assume on the contrary that there exists a non-locally algebraic constituent $W_0$ of $V_{\un{j},\un{j}'}$ such that $W_0\neq C^{a,b,c}$ for any $(a,b,c)\in\{0,1\}^3$. Then there exists $(b,c)\in\{0,1\}^2$ such that $W_0$ occurs in $\mathrm{gr}_{b,c}(V_{\un{j},\un{j}'})\subseteq \Theta_{w_0s_{j_0}w_0}(C^{0,b,c})$. We take $(b,c)$ such that $b+c$ is maximal among those $(b,c)\in\{0,1\}^2$ such that $\mathrm{gr}_{b,c}(V_{\un{j},\un{j}'})$ has a non-locally algebraic constituent distinct from the $C^{a,b,c}$. By definition we can write $C^{0,b,c}=\cF_{P_{\widehat{j}}}^G(L(w_{j,j_0}),\pi^{\infty})$ for some $j\in\{j_1-1,j_1,j_1+1\}\cap\Delta$ and some irreducible $G$-regular $\pi^{\infty}$ (see above Lemma \ref{lem: cube step 1}), and note that, by definition again, we have $C^{1,b,c}=\cF_{P_{\widehat{j}}}^G(L(w_{j,j_0+1}),\pi^{\infty})$. It then follows from Proposition~\ref{prop: Jantzen middle}, Remark~\ref{rem: explicit middle}, (\ref{equ: OS wall crossing}) and Lemma~\ref{lem: Hom OS} (and the above assumptions on $W_0$) that we must have $W_0\cong \cF_{P_{\widehat{j}}}^G(L(w_{j,j_0-1}),\pi^{\infty})$, and using moreover (\ref{equ: grade in rad}) and the first statement of \ref{it: cube step 1 3} of Lemma \ref{lem: cube step 1} that $W_0$ must be in the cosocle of $\mathrm{gr}_{b,c}(V_{\un{j},\un{j}'})$. Since $V_{\un{j},\un{j}'}$ has cosocle $C^{1,1,1}\ne W_0$ there must exist a constituent $W_1$ of $\mathrm{gr}_{b',c'}(V_{\un{j},\un{j}'})$ for some $b',c'\in\{0,1\}$ with $b'+c'>b+c$ such that $V_{\un{j},\un{j}'}$ admits a length $2$ subquotient $V$ with socle $W_0$ and cosocle $W_1$. As $b'+c'>b+c\geq 0$, it follows from Lemma~\ref{lem: cube step loc alg} that $W_1$ is not locally algebraic. By maximality of $b+c$, we deduce $W_1=C^{a',b',c'}$ for some $(a',b',c')\in\{0,1\}^3$. If $a'=0$, then $W_1$ occurs in $V_0$ by \ref{it: cube step 2 1} and \ref{it: cube step 2 3} of Lemma~\ref{lem: cube step 2}, but $W_0$ does not, a contradiction. Assume $a'=1$, and thus $W_1=C^{1,b',c'}=\cF_{P_{\widehat{j}'}}^G(L(w_{j',j_0+1}),\sigma^{\infty})$ for some $j'\in\{j_1-1,j_1,j_1+1\}\cap\Delta$ and some irreducible $G$-regular $\sigma^{\infty}$ (see above Lemma \ref{lem: cube step 1}). By \ref{it: rabiotext 2} of Lemma~\ref{rabiotext} and (the very last statement in) Remark~\ref{rem: coxeter pair}, we have $\mathrm{Ext}_{U(\fg)}^1(L(w_{j,j_0-1}),L(w_{j',j_0+1}))=0$. By Proposition~\ref{prop: Ext1 OS} this forces $\mathrm{Ext}_{G}^1(W_1,W_0)=0$, a contradiction to the existence of $V$.
\end{proof}

\begin{lem}\label{lem: cube step 4}
Assume $j_1\neq j_0+1$. Then $V_{\un{j},\un{j}'}$ admits a unique quotient isomorphic to $V_1$ (see above Lemma \ref{lem: cube step 1} for $V_1$).
\end{lem}
\begin{proof}
Using \ref{it: cube step 2 3} of Lemma~\ref{lem: cube step 2} we define $V_1'$ as the unique quotient of $V_{\un{j},\un{j}'}$ with socle $C^{1,0,0}$. We let $\mathrm{Fil}_{b,c}(V_1')$ be the image of $\mathrm{Fil}_{b,c}(V_{\un{j},\un{j}'})$ via $V_{\un{j},\un{j}'}\twoheadrightarrow V_1'$. It follows from \ref{it: cube step 1 2} of Lemma \ref{lem: cube step 1} and from (\ref{equ: grade in rad}) that $C^{1,0,0}$ occurs in the cosocle of $\mathrm{gr}_{0,0}(V_{\un{j},\un{j}'})$, which forces $\mathrm{gr}_{0,0}(V_1')\cong C^{1,0,0}$. Since $\mathrm{gr}_{b,c}(V_{\un{j},\un{j}'})$, and hence $\mathrm{gr}_{b,c}(V_1')$, have no locally algebraic constituent when $b+c>0$ by Lemma~\ref{lem: cube step loc alg}, it follows that $V_1'$ has no locally algebraic constituent. By Lemma~\ref{lem: cube step 3}, all constituents of $V'_1$ are therefore of the form $C^{a,b,c}$ for some $(a,b,c)\in\{0,1\}^3$. In particular $V_1'$ is multiplicity free by \ref{it: cube step 2 3} of Lemma~\ref{lem: cube step 2} with socle $C^{1,0,0}$ and cosocle $C^{1,1,1}$. Since $V_0$ maps to $0$ via $V_0\hookrightarrow V_{\un{j},\un{j}'}\twoheadrightarrow V_1'$ (using \ref{it: cube step 2 1} of Lemma~\ref{lem: cube step 2} and the fact $V'_1$ has cosocle $C^{1,0,0}$ which doesn't occur in $V_0$), it follows from \ref{it: cube step 2 3} of Lemma~\ref{lem: cube step 2} that the constituents of $V'_1$ are of the form $C^{1, b', c'}$ for some $(b',c')\in \{0,1\}^2$. Since $C^{1,b,c}$ is a constituent of $\mathrm{gr}_{b',c'}(V_{\un{j},\un{j}'})$ if and only if $b=b'$ and $c=c'$ by \ref{it: cube step 1 1} of Lemma~\ref{lem: cube step 1}, we deduce with \ref{it: cube step 2 4} of Lemma \ref{lem: cube step 2} that $V_1'$ has Loewy length $3$ with socle $C^{1,0,0}$, cosocle $C^{1,1,1}$ and middle layer $C^{1,0,1}\oplus C^{1,1,0}$. By Lemma~\ref{lem: minimal is unique} and the definition of $V_1$, this implies $V_1'\cong V_1$.
\end{proof}

For $I\subseteq \Delta$ we write (see also (\ref{gensteinberg}))
\begin{equation}\label{vIalg}
V_{I,\Delta}^{\rm{alg}}\defeq L(1)^\vee\otimes_E V_{I,\Delta}^{\infty} = L(\mu_0)^\vee\otimes_E V_{I,\Delta}^{\infty}.
\end{equation}

\begin{lem}\label{lem: cube step 5}
Assume $j_1\neq j_0+1$. Then $V_{\un{j},\un{j}'}$ is multiplicity free. Moreover, when $j_1=j_0$, $V_{\un{j},\un{j}'}$ admits a unique quotient isomorphic to $V_{\infty,\un{j}'}$ (see \ref{it: easy square 4} of Proposition \ref{prop: easy square}) which contains all its locally algebraic constituents.
\end{lem}
\begin{proof}
If $j_1\neq j_0$ (and $j_1\neq j_0+1$), then by Lemma~\ref{lem: cube step loc alg} $V_{\un{j},\un{j}'}$ has no locally algebraic constituent, which by Lemma~\ref{lem: cube step 3} and \ref{it: cube step 2 3} of Lemma~\ref{lem: cube step 2} implies that $V_{\un{j},\un{j}'}$ is multiplicity free. We assume $j_1=j_0$ in the rest of the proof. Note that, by \emph{loc.~cit.}~each non locally algebraic constituent of $V_{\un{j},\un{j}'}$ occurs with multiplicity $1$. Recall from Proposition~\ref{prop: Jantzen middle} that $L(1)$ has multiplicity $1$ in $\Theta_{w_0s_{j_0}w_0}(L(s_{j_0}))$. It then follows from (\ref{equ: OS wall crossing}), \ref{it: OS property 2} of Theorem~\ref{prop: OS property} and Lemma~\ref{lem: Hom OS} that $L(1)^\vee\otimes_E i_{\widehat{j}_0,\Delta}^{\infty}(\pi_{j_0,j_2}^{\infty})$ contains all the locally algebraic constituents of $\Theta_{w_0s_{j_0}w_0}(C^{0,0,0})$, and therefore of $V_{\un{j},\un{j}'}$ by Lemma~\ref{lem: cube step loc alg}. As $i_{\widehat{j}_0,\Delta}^{\infty}(\pi_{j_0,j_2}^{\infty})$ is multiplicity free by \ref{it: basic as image} of Remark~\ref{rem: basic PS intertwine}, we deduce that $V_{\un{j},\un{j}'}$ is multiplicity free.

It remains to prove that $V_{\un{j},\un{j}'}$ admits a (unique) quotient isomorphic to $V_{\infty,\un{j}'}$ which contains all its locally algebraic constituents. As $\mathrm{gr}_{0,0}(V_{\un{j},\un{j}'})\subseteq \Theta_{w_0s_{j_0}w_0}(C^{0,0,0})$, and the subquotient $L(1)^\vee\otimes_E i_{\widehat{j}_0,\Delta}^{\infty}(\pi_{j_0,j_2}^{\infty})$ of $\Theta_{w_0s_{j_0}w_0}(C^{0,0,0})$ contains all its locally algebraic constituents, there exists a (unique) subrepresentation $\pi^{\infty}\subseteq i_{\widehat{j}_0,\Delta}^{\infty}(\pi_{j_0,j_2}^{\infty})$ such that $L(1)^\vee\otimes_E \pi^{\infty}$ is a subquotient of $V_{\un{j},\un{j}'}$ which contains its locally algebraic constituents. Moreover by (\ref{equ: grade in rad}), Proposition \ref{prop: Jantzen middle} and Lemma~\ref{lem: Hom OS} we know that $L(1)^\vee\otimes_E \pi^{\infty}$ is a quotient of $\mathrm{gr}_{0,0}(V_{\un{j},\un{j}'})$.

Now we study all possible length $2$ subquotients $V$ of $V_{\un{j},\un{j}'}$ with socle $W_0$ and cosocle $W_1$ such that $W_0$ is locally algebraic and $W_1$ is not. By Lemma~\ref{lem: cube step 3} $W_1=C^{a,b,c}$ for some $(a,b,c)\in\{0,1\}^3$ and by \ref{it: cube step 1 1} of Lemma~\ref{lem: cube step 1} $W_1$ is a constituent of $\mathrm{gr}_{b,c}(V_{\un{j},\un{j}'})$. As $W_0$ is necessarily a constituent of $L(1)^\vee\otimes_E \pi^{\infty}$ which is a quotient of $\mathrm{gr}_{0,0}(V_{\un{j},\un{j}'})$, we must have $b+c>0$. If $a=0$, then $W_1$ occurs in $V_0$ (using \ref{it: cube step 2 1} and \ref{it: cube step 2 3} of Lemma~\ref{lem: cube step 2}), which together with the layer structure of $V_0$ (see \ref{it: hard square 2} of Proposition~\ref{prop: hard square}) forces $W_1=C^{0,1,1}$ and $W_0=L(1)^\vee\otimes_E V_{[j_2-j_0+1,j_2],\Delta}^{\infty}$.
Assume $a=1$ and thus $W_1=\cF_{P_{\widehat{j}}}^G(L(w_{j,j_0+1}),\tau^{\infty})$ for some $j\in\{j_1-1,j_1,j_1+1\}\cap\Delta$ and some irreducible $G$-regular $\tau^{\infty}$ (see above Lemma \ref{lem: cube step 1}). The existence of $V$ implies $\mathrm{Ext}_{G}^1(W_1,W_0)\neq 0$, which together with Proposition~\ref{prop: Ext1 OS} forces $\mathrm{Ext}_{U(\fg)}^1(L(1),L(w_{j,j_0+1}))\neq 0$. By \ref{it: rabiotext 2} of Lemma~\ref{rabiotext} and our assumption $j_0=j_1$, we deduce $j=j_0+1=j_1+1$ and thus $W_1=C^{1,1,0}$ (with $w_{j,j_0+1}=s_{j_0+1}$ and $\tau^{\infty}=\pi_{j_0+1,j_2+1}^{\infty}$). We can write $W_0=L(1)^\vee\otimes_E \sigma^{\infty}$ for some irreducible $G$-regular constituent $\sigma^{\infty}$ of $\pi^{\infty}$. Then $\mathrm{Ext}_{G}^1(W_1,W_0)\neq 0$ and the last statement in Remark~\ref{rem: Ext1 OS distance} imply $\Hom_{G}(\sigma^{\infty}, i_{\Delta\setminus\{j_0+1\},\Delta}^{\infty}(\pi_{j_0+1,j_2+1}^{\infty}))\neq 0$. By Lemma~\ref{lem: explicit smooth induction} we deduce $\sigma^{\infty}=V_{[j_2-j_0,j_2],\Delta}^{\infty}$ if $j_0=j_1<j_2$ and $\sigma^{\infty}=V_{[j_2-j_0+1,j_2],\Delta}^{\infty}$ if $j_0=j_1=j_2$. We have thus shown that all possible pairs $(W_0,W_1)$ are $(V_{[j_2-j_0+1,j_2],\Delta}^{\rm{alg}},C^{0,1,1})$, $(V_{[j_2-j_0,j_2],\Delta}^{\rm{alg}}, C^{1,1,0})$ when $j_0=j_1<j_2$, and $(V_{[j_2-j_0+1,j_2],\Delta}^{\rm{alg}}, C^{1,1,0})$ when $j_0=j_1=j_2$.

Let $\sigma^{\infty}$ be a constituent of $\mathrm{cosoc}_{G}(\pi^{\infty})$, then since $L(1)^\vee\otimes_E \pi^{\infty}$ is a quotient of $\mathrm{gr}_{0,0}(V_{\un{j},\un{j}'})$ there must exist $(W_0,W_1)$ in the above list with $W_0=L(1)^\vee\otimes_E \sigma^{\infty}$. An examination of Lemma~\ref{lem: explicit smooth induction} shows that $\pi^{\infty}\subseteq i_{\widehat{j}_0,\Delta}^{\infty}(\pi_{j_0,j_2}^{\infty})$ must be the length $2$ subrepresentation with socle $V_{[j_2-j_0+1,j_2],\Delta}^{\infty}$ and cosocle $V_{[j_2-j_0,j_2],\Delta}^{\infty}$ when $j_0=j_1<j_2$, and the irreducible representation $V_{[j_2-j_0+1,j_2],\Delta}^{\infty}$ when $j_0=j_1=j_2$. From the description of the possible pairs above, from Lemma~\ref{lem: cube step 3} and \ref{it: cube step 2 3}, \ref{it: cube step 2 4} of Lemma~\ref{lem: cube step 2}, and from the fact $V_{\un{j},\un{j}'}$ has (irreducible) cosocle $C^{1,1,1}$, we can deduce that $V_{\un{j},\un{j}'}$ admits a unique (multiplicity free) quotient with constituents
\[\{V_{[j_2-j_0+1,j_2],\Delta}^{\rm{alg}},\ V_{[j_2-j_0,j_2],\Delta}^{\rm{alg}},\ C^{0,1,1},\ C^{1,1,0},\ C^{1,1,1}\}\]
with partial order (in the sense of \S\ref{generalnotation}) which can only be: $V_{[j_2-j_0+1,j_2],\Delta}^{\rm{alg}}\leq V_{[j_2-j_0,j_2],\Delta}^{\rm{alg}}$, $V_{[j_2-j_0+1,j_2],\Delta}^{\rm{alg}}\leq C^{0,1,1}$, $V_{[j_2-j_0,j_2],\Delta}^{\rm{alg}}\leq C^{1,1,0}$, $C^{0,1,1}\leq C^{1,1,1}$ and $C^{1,1,0}\leq C^{1,1,1}$ (with $V_{[j_2-j_0,j_2],\Delta}^{\rm{alg}}$ omitted when $j_0=j_1=j_2$). It then follows from \ref{it: easy square 4} of Proposition~\ref{prop: easy square} (with $\un{j}$ there replaced by $\un{j}'$ here) together with the minimality of the $\mathrm{Ext}$-square $V_{\infty,\un{j}'}$ in \emph{loc.~cit.}~(last statement of Proposition~\ref{prop: easy square}) and with Lemma~\ref{lem: minimal is unique} that this quotient must be $V_{\infty,\un{j}'}$.
\end{proof}

\begin{lem}\label{lem: cube step 6}
Assume $j_1\neq j_0+1$. Then $V_{\un{j},\un{j}'}$ is an $\mathrm{Ext}$-cube such that
\begin{enumerate}[label=(\roman*)]
\item \label{it: cube step 6 1} $\mathrm{gr}^0(V_{\un{j},\un{j}'})\cong C_{\un{j}'}$ and $\mathrm{gr}^3(V_{\un{j},\un{j}'})\cong C_{\un{j}}$;
\item \label{it: cube step 6 2} $V_{\un{j},\un{j}'}$ contains a unique subquotient of the form $V_{\un{j}'',\un{j}'''}$ for each pair $(\un{j}'',\un{j}''')\in \mathbf{J}^2$ satisfying $\un{j}\leq \un{j}''\leq \un{j}'''\leq \un{j}'$ and $d(\un{j}'',\un{j}''')=2$ (these $V_{\un{j}'',\un{j}'''}$ are defined in \ref{it: easy square 1}, \ref{it: easy square 2} of Proposition \ref{prop: easy square} and \ref{it: hard square 1}, \ref{it: hard square 2} of Proposition \ref{prop: hard square}).
\end{enumerate}
\end{lem}
\begin{proof}
We first prove \ref{it: cube step 6 2}. Note that $\un{j}''\in \mathbf{J}$ satisfying $\un{j}\leq \un{j}'' \leq \un{j}'$ is equivalent to the choice of $(a,b,c)\in \{0,1\}^3$ (writing $\un{j}''=(j_0+a,j_1+c-b,j_2+c)$). In particular, the choice of $\un{j}'',\un{j}'''\in \mathbf{J}$ such that $\un{j}\leq \un{j}''\leq \un{j}''' \leq \un{j}'$ and $d(\un{j''},\un{j}''')=2$ is equivalent to the choice of $(a,b,c), (a',b',c')\in \{0,1\}^3$ such that $a\leq a'$, $b\leq b'$, $c\leq c'$ and $a'+b'+c'=a+b+c+2$. It follows from \ref{it: cube step 2 3} of Lemma~\ref{lem: cube step 2} that if $V_{\un{j}'',\un{j}'''}$ is a subquotient of $V_{\un{j},\un{j}'}$, then it is necessarily the unique subquotient of $V_{\un{j},\un{j}'}$ with socle $C_{\un{j}''}$ and cosocle $C_{\un{j}'''}$.
If $b=b'$ or $c=c'$, then the existence of a subquotient $V_{\un{j}'',\un{j}'''}$ of $V_{\un{j},\un{j}'}$ follows from \ref{it: cube step 2 2}, \ref{it: cube step 2 3} of Lemma~\ref{lem: cube step 2} and \ref{it: cube step 1 2} of Lemma~\ref{lem: cube step 1}. If $a=a'=0$, then $V_{\un{j}'',\un{j}'''}=V_0$ and the existence of a subquotient $V_0$ of $V_{\un{j},\un{j}'}$ follows from \ref{it: cube step 2 1} of Lemma~\ref{lem: cube step 2}. If $a=a'=1$, then $V_{\un{j}'',\un{j}'''}=V_1$ and the existence of a subquotient $V_1$ of $V_{\un{j},\un{j}'}$ follows from Lemma~\ref{lem: cube step 4}. This proves \ref{it: cube step 6 2}.

Now we construct a filtration on $V_{\un{j},\un{j}'}$ that makes it an $\mathrm{Ext}$-cube (see Definition \ref{def: OS cube}). Recall that $V_{\un{j},\un{j}'}$ has socle $C^{0,0,0}$ and cosocle $C^{1,1,1}$. It follows from \ref{it: cube step 6 2} above, from \ref{it: cube step 2 3} and \ref{it: cube step 2 4} of Lemma~\ref{lem: cube step 2}, from Lemma~\ref{lem: cube step 3} and from Lemma~\ref{lem: cube step 5} (with \ref{it: easy square 4} and \ref{it: easy square 2} of Proposition \ref{prop: easy square}) that $V_{\un{j},\un{j}'}$ admits a unique decreasing filtration $(\mathrm{Fil}^k(V_{\un{j},\un{j}'}))_{0\leq k\leq 3}$ such that for $0\leq k\leq 3$
\begin{equation}\label{grjj'}
\mathrm{gr}^k(V_{\un{j},\un{j}'})\cong C_{\infty,k}\bigoplus\bigg(\bigoplus_{a+b+c=3-k}C^{a,b,c}\bigg),
\end{equation}
where $C_{\infty,k}\defeq \mathrm{gr}^2(V_{\infty,\un{j}'})$ when $k=2$ and $j_1=j_0$, and $C_{\infty,k}\defeq 0$ otherwise. In particular, \ref{it: cube step 6 1} holds for the filtration $(\mathrm{Fil}^k(V_{\un{j},\un{j}'}))_k$.
It remains to check that $V_{\un{j},\un{j}'}$ is an $\mathrm{Ext}$-cube for this filtration. But the first condition in \ref{it: OS cube 1} of Definition~\ref{def: OS cube} holds by definition of $V_{\un{j},\un{j}'}$ and of the filtration $(\mathrm{Fil}^k(V_{\un{j},\un{j}'}))_k$, and the second condition in \ref{it: OS cube 1} of Definition~\ref{def: OS cube} holds by Lemma~\ref{lem: Ext1 factor 1} and Lemma~\ref{lem: Ext1 factor 2}.
\end{proof}

\begin{lem}\label{lem: cube step 7}
Assume $j_1\neq j_0+1$. Then the $\mathrm{Ext}$-cube $V_{\un{j},\un{j}'}$ is strict and minimal (see \ref{it: OS cube 3} and \ref{it: OS cube 4} of Definition \ref{def: OS cube}).
\end{lem}
\begin{proof}
The $\mathrm{Ext}$-cube $V_{\un{j},\un{j}'}$ is strict as it has simple socle ($C^{0,0,0}=C_{\un{j}}$) and simple cosocle ($C^{1,1,1}=C_{\un{j}'}$).
By \ref{it: easy square 1}, \ref{it: easy square 2} of Proposition~\ref{prop: easy square} and Proposition \ref{prop: hard square}, $V_{\un{j}'',\un{j}'''}$ is minimal for each $(\un{j}'',\un{j}''')\in \mathbf{J}^2$ such that $\un{j}\leq \un{j}''\leq \un{j}'''\leq \un{j}'$ and $d(\un{j}'',\un{j}''')=2$, and by \ref{it: easy square 4} of Proposition~\ref{prop: easy square} $V_{\infty,\un{j}'}$ is minimal when $j_1=j_0$. Hence to check the minimality of $V_{\un{j},\un{j}'}$, it suffices to show that any $\mathrm{Ext}$-cube $V'$ such that $\mathrm{gr}^0(V')=C_{\un{j}'}$, $\mathrm{gr}^3(V')=C_{\un{j}}$ and $\mathrm{gr}^k(V')$ is a good direct summand of $\mathrm{gr}^k(V_{\un{j},\un{j}'})$ for $k=1,2$ must satisfy $\mathrm{gr}^k(V')=\mathrm{gr}^k(V_{\un{j},\un{j}'})$ for $k=1,2$. Let $V'$ be such an $\mathrm{Ext}$-cube, then $\mathrm{gr}^k(V')\neq 0$ for $k=1,2$. Thus by (\ref{grjj'}) $\mathrm{gr}^1(V')$ contains a constituent $C_{\un{j}''}$ for some $\un{j}< \un{j}''< \un{j}'$ such that $d(\un{j},\un{j}'')=2$. It follows from Lemma \ref{lem: cube step 6} and its proof that the unique subrepresentation $V''$ of $V'$ with cosocle $C_{\un{j}''}$ is an $\mathrm{Ext}$-square such that $\mathrm{gr}^0(V'')=C_{\un{j}''}$, $\mathrm{gr}^2(V'')=C_{\un{j}}$ and $\mathrm{gr}^1(V'')$ is a good direct summand of $\mathrm{gr}^1(V_{\un{j},\un{j}''})$. The minimality of $V_{\un{j},\un{j}''}$ forces $\mathrm{gr}^1(V'')\cong\mathrm{gr}^1(V_{\un{j},\un{j}''})$, and thus for each $\un{j}'''$ such that $\un{j}< \un{j}'''< \un{j}''$ the constituent $C_{\un{j}'''}$ shows up in $\mathrm{gr}^1(V'')$, and therefore in $\mathrm{gr}^2(V')$. A similar argument using the minimality of $V_{\un{j}''',\un{j}'}$ for $\un{j}'''$ such that $\un{j}< \un{j}'''< \un{j}''$ gives $\mathrm{gr}^1(V')\cong \mathrm{gr}^1(V_{\un{j},\un{j}'})$. In particular $C_{\un{j}''}$ occurs in $\mathrm{gr}^1(V')$ for \emph{each} $\un{j}< \un{j}''< \un{j}'$ with $d(\un{j},\un{j}'')=2$, and we can repeat the previous argument with the subrepresentation $V''$ of $V'$ with cosocle $C_{\un{j}''}$. We now deduce that $C_{\un{j}'''}$ occurs in $\mathrm{gr}^1(V'')$, and therefore in $\mathrm{gr}^2(V')$, for \emph{any} $\un{j}'''$ such that $\un{j}< \un{j}'''< \un{j}'$ and $d(\un{j},\un{j}''')=1$. If $j_0\ne j_1$, there are no constituents left and we have proven $\mathrm{gr}^k(V')=\mathrm{gr}^k(V_{\un{j},\un{j}'})$ for $k=1,2$. We now assume $j_1=j_0$. We already have $\mathrm{gr}^1(V')\cong \mathrm{gr}^1(V_{\un{j},\un{j}'})$ but it remains to prove $\mathrm{gr}^2(V')\cong \mathrm{gr}^2(V_{\un{j},\un{j}'})$, and for that we have to prove that the locally algebraic constituents of $V_{\un{j},\un{j}'}$, equivalently of $\mathrm{gr}^2(V_{\un{j},\un{j}'})$ by Lemma \ref{lem: cube step 5}, are also in $\mathrm{gr}^2(V')$. Let $\un{j}''=(j_0+1,j_1+1,j_2+1)$, from the minimality of $V_{\un{j},\un{j}''}$ in \ref{it: easy square 2} of Proposition \ref{prop: easy square} and the fact $C_{\un{j}''}$ occurs in $\mathrm{gr}^1(V')$, we deduce $\mathrm{gr}^2(V_{\infty,\un{j}'})\subseteq \mathrm{gr}^2(V')$. This finally implies $\mathrm{gr}^2(V')\cong \mathrm{gr}^2(V_{\un{j},\un{j}'})$ and finishes the proof.
\end{proof}

We now note that Lemma \ref{lem: cube step 1} has a symmetric statement replacing $V_0$ by $V_1$, $\Theta_{w_0s_{j_0}w_0}$ (resp.~$\Theta_\mu$) by $\Theta_{w_0s_{j_0+1}w_0}$ (resp.~by $\Theta_\mu$ where $\mu\in \Lambda$ is such that $\langle\mu+\rho,\alpha^\vee\rangle\geq 0$ for $\al\in \Phi^+$ and the stabilizer of $\mu$ in $W(G)$ for the dot action is $\{1,w_0s_{j_0+1}w_0\}$) and switching $C^{0,b,c}$ and $C^{1,b,c}$ everywhere. Then we can define $\widetilde V_{\un{j},\un{j}'}$ as the unique quotient of $\Theta_{\mu}(V_1)$ with socle $C_{\un{j}}$. Then Lemma \ref{lem: cube step 2} also has a symmetric version ($\Theta_{\mu}(V_1)$ and $\widetilde V_{\un{j},\un{j}'}$ have cosocle $C^{1,1,1}$, $\widetilde V_{\un{j},\un{j}'}$ surjects onto $V_1$, etc.) and likewise all statements from Lemma \ref{lem: cube step loc alg} to Lemma \ref{lem: cube step 7} have symmetric versions replacing the assumption $j_1\ne j_0+1$ by the assumption $j_1\ne j_0$, $V_0$ by $V_1$ and the case $j_1=j_0$ by the case $j_1=j_0+1$. We let the reader work out by himself the symmetric statements. For instance the symmetric Lemma \ref{lem: cube step 5} is: assume $j_1\neq j_0$, then $V_{\un{j},\un{j}'}$ is multiplicity free, and moreover when $j_1=j_0+1$, $\widetilde V_{\un{j},\un{j}'}$ admits a unique subrepresentation isomorphic to $V_{\un{j},\infty}$ which contains all its locally algebraic constituents. Note that all minimal $\mathrm{Ext}$-squares that are used in the (symmetric) proofs are still provided by Proposition \ref{prop: easy square} and Proposition \ref{prop: hard square}. It also follows from Lemma \ref{lem: cube step 6} and its proof, from the minimality in Lemma \ref{lem: cube step 7}, from their symmetric versions, and from Lemma \ref{lem: minimal is unique} that $V_{\un{j},\un{j}'}\cong \widetilde V_{\un{j},\un{j}'}$ when $j_1\notin \{j_0,j_0+1\}$.\bigskip

In order to avoid too much notation, we will now denote by $V_{\un{j},\un{j}'}$ the Ext-cube previously (also) denoted $V_{\un{j},\un{j}'}$ when $j_1\ne j_0+1$, and by $V_{\un{j},\un{j}'}$ the Ext-cube previously denoted $\widetilde V_{\un{j},\un{j}'}$ when $j_1\ne j_0$. Note that $V_{\un{j},\un{j}'}$ is well defined by the above isomorphism. The following proposition sums up some of the previous results on $V_{\un{j},\un{j}'}$ which are proven in Lemma~\ref{lem: cube step 7}, Lemma~\ref{lem: cube step 3}, Lemma~\ref{lem: cube step 6}, Lemma~\ref{lem: cube step 5} and in their symmetric versions.

\begin{prop}\label{prop: factor cube}
Let $\un{j}\in\mathbf{J}$ such that $\un{j}'=(j'_0, j'_1, j'_2)\defeq (j_0+1, j_1, j_2+1)$ is still in $\mathbf{J}$. There exists a minimal $\mathrm{Ext}$-cube $V_{\un{j},\un{j}'}$ such that
\begin{enumerate}[label=(\roman*)]
\item the non-locally algebraic constituents of $V_{\un{j},\un{j}'}$ are the $C_{\un{j}''}$ for $\un{j}\leq \un{j}''\leq \un{j}'$;
\item $\mathrm{gr}^0(V_{\un{j},\un{j}'})\cong C_{(j_0+1,j_1,j_2+1)}$ and $\mathrm{gr}^3(V_{\un{j},\un{j}'})\cong C_{(j_0,j_1,j_2)}$;
\item $V_{\un{j},\un{j}'}$ contains a unique subquotient of the form $V_{\un{j}'',\un{j}'''}$ for each pair $(\un{j}'',\un{j}''')\in \mathbf{J}^2$ satisfying $\un{j}\leq \un{j}''\leq \un{j}'''\leq \un{j}'$ and $d(\un{j}'',\un{j}''')=2$ (\ref{it: easy square 1}, \ref{it: easy square 2} of Proposition \ref{prop: easy square} and \ref{it: hard square 1}, \ref{it: hard square 2} of Proposition \ref{prop: hard square});
\item $V_{\un{j},\un{j}'}$ admits a locally algebraic constituent if and only if one of the following holds: either $j_1=j_0$ and $V_{\un{j},\un{j}'}$ admits a unique quotient isomorphic to $V_{\infty,\un{j}'}$ (\ref{it: easy square 4} of Proposition \ref{prop: easy square}) which contains all its locally algebraic constituents, or $j_1=j_0+1$ and $V_{\un{j},\un{j}'}$ admits a unique subrepresentation isomorphic to $V_{\un{j},\infty}$ (\ref{it: easy square 3} of Proposition \ref{prop: easy square}) which contains all its locally algebraic constituents.
\end{enumerate}
\end{prop}

\begin{rem}
In fact one can prove that the representations $V_{\un{j},\un{j}'}$ and $\widetilde V_{\un{j},\un{j}'}$ just below Lemma \ref{lem: cube step 7} are \emph{always} isomorphic (even if $j_1\in \{j_0, j_0+1\}$). But we won't need that result.
\end{rem}

We end up this section with three lemmas, two of which construct other finite length representations $V_{\un{j},\infty}^+$ and $V_{\un{j}}$ of $G$ in $\mathrm{Rep}^{\rm{an}}_{\rm{adm}}(G)$ which will be important to define key representations in \S\ref{subsec: final}.

\begin{lem}\label{lem: sm top}
Let $1\leq j_0=j_1\leq j_2<n-1$ and let $\un{j}\defeq (j_0,j_0,j_2)$, $\un{j}'\defeq (j_0,j_0,j_2+1)$, $\un{j}''\defeq (j_0+1,j_0+1,j_2+1)$. Then there exists a unique multiplicity free finite length representation $V_{\un{j},\infty}^+$ in $\mathrm{Rep}^{\rm{an}}_{\rm{adm}}(G)$ such that (see (\ref{vIalg}) for the notation)
\begin{enumerate}[label=(\roman*)]
\item \label{it: sm top 0} $V_{\un{j},\infty}^+$ has socle $C_{\un{j}}$ and cosocle $V_{[j_2-j_0+1,j_2+1],\Delta}^{\rm alg}$;
\item \label{it: sm top 1} both $V_{\un{j},\un{j}'}$ and $V_{\un{j},\un{j}''}$ inject into $V_{\un{j},\infty}^+$ (see \ref{it: hard square 2} of Proposition \ref{prop: hard square} and \ref{it: easy square 2} of Proposition \ref{prop: easy square} respectively);
\item \label{it: sm top 2} the quotient $V_{\un{j},\infty}^+/(V_{\un{j},\un{j}'}+V_{\un{j},\un{j}''})$ is uniserial of length $2$ with socle $V_{[j_2-j_0+2,j_2+1],\Delta}^{\rm alg}$ and cosocle $V_{[j_2-j_0+1,j_2+1],\Delta}^{\rm alg}$.
\end{enumerate}
\end{lem}
\begin{proof}
We first prove unicity of $V_{\un{j},\infty}^+$. We fix $\mu\in \Lambda$ such that $\langle\mu+\rho,\alpha^\vee\rangle\geq 0$ for $\al\in \Phi^+$ and the stabilizer of $\mu$ in $W(G)$ for the dot action is $\{1,w_0s_{j_0}w_0\}$. Condition \ref{it: sm top 2} together with Proposition~\ref{prop: Jantzen middle} and (\ref{equ: OS wall crossing}) imply $\Theta_{\mu}(V_{\un{j},\infty}^+/V_{\un{j},\un{j}'})=0$. The exactness of $\Theta_{\mu}$ then imply $\Theta_{\mu}(V_{\un{j},\un{j}'})\buildrel\sim\over\rightarrow \Theta_{\mu}(V_{\un{j},\infty}^+)$. Since $V_{\un{j},\infty}^+$ has socle $C_{\un{j}}$ by condition \ref{it: sm top 1}, it follows from the exactness and functoriality of $\Theta_{\mu}$ that the adjunction map $V_{\un{j},\infty}^+\rightarrow \Theta_{\mu}(V_{\un{j},\infty}^+)$ is injective. We thus deduce an injection $V_{\un{j},\infty}^+\hookrightarrow \Theta_{\mu}(V_{\un{j},\un{j}'})$. It follows from Proposition~\ref{prop: Jantzen middle} that $L(1)$ is a constituent of $\Theta_{w_0s_{j_0}w_0}(L(x))$ for some $x\in W(G)$ if and only if $x=s_{j_0}$, in which case $L(1)$ occurs with multiplicity $1$. Together with \ref{it: hard square 2} of Proposition \ref{prop: hard square}, (\ref{equ: OS wall crossing}) and the fact $\Theta_\mu$ kills locally algebraic constituents of $V_{\un{j},\un{j}'}$ (Proposition~\ref{prop: Jantzen middle}), we see that if $W$ is a constituent of $V_{\un{j},\un{j}'}$, then $\Theta_{\mu}(W)$ admits locally algebraic constituents if and only if $W\in \{C_{\un{j}},C_{\un{j}'}\}$, and together with \ref{it: OS property 2} of Theorem \ref{prop: OS property} that $\Theta_{\mu}(C_{\un{j}})$ (resp.~$\Theta_{\mu}(C_{\un{j}'})$) admits a subquotient $L(1)^\vee\otimes_E i_{\widehat{j}_0,\Delta}^{\infty}(\pi_{j_0,j_2}^{\infty})$ (resp.~$L(1)^\vee\otimes_E i_{\widehat{j}_0,\Delta}^{\infty}(\pi_{j_0,j_2+1}^{\infty})$) which contains all its locally algebraic constituents. By Lemma~\ref{lem: explicit smooth induction} $i_{\widehat{j}_0,\Delta}^{\infty}(\pi_{j_0,j_2}^{\infty})$ and $i_{\widehat{j}_0,\Delta}^{\infty}(\pi_{j_0,j_2+1}^{\infty})$ are multiplicity free,
and $V_{[j_2-j_0+1,j_2+1],\Delta}^{\infty},V_{[j_2-j_0+2,j_2+1],\Delta}^{\infty}\in \mathrm{JH}_{G}(i_{\widehat{j}_0,\Delta}^{\infty}(\pi_{j_0,j_2+1}^{\infty}))\setminus \mathrm{JH}_{G}(i_{\widehat{j}_0,\Delta}^{\infty}(\pi_{j_0,j_2}^{\infty}))$. This implies that $V_{[j_2-j_0+1,j_2+1],\Delta}^{\rm{alg}},V_{[j_2-j_0+2,j_2+1],\Delta}^{\rm{alg}}\in \mathrm{JH}_{G}(\Theta_{\mu}(C_{\un{j}'}))\setminus \mathrm{JH}_{G}(\Theta_{\mu}(C_{\un{j}}))$, that $V_{[j_2-j_0+1,j_2+1],\Delta}^{\rm{alg}}$, $V_{[j_2-j_0+2,j_2+1],\Delta}^{\rm{alg}}$ both appear with multiplicity $1$ in $\Theta_{\mu}(V_{\un{j},\un{j}'})$, and that $\Theta_{\mu}(V_{\un{j},\un{j}'})$ contains as a subquotient the unique non-split extension of $V_{[j_2-j_0+1,j_2+1],\Delta}^{\rm{alg}}$ by $V_{[j_2-j_0+2,j_2+1],\Delta}^{\rm{alg}}$ (see Lemma~\ref{lem: explicit smooth induction} and Lemma \ref{lem: Ext1 from sm}). In particular the second statement in condition \ref{it: sm top 0} and the beginning of the proof force $V_{\un{j},\infty}^+$, if it exists, to be the unique subrepresentation of $\Theta_{\mu}(V_{\un{j},\un{j}'})$ with cosocle $V_{[j_2-j_0+1,j_2+1],\Delta}^{\rm{alg}}$. This proves unicity of $V_{\un{j},\infty}^+$.

We now prove that $V_{\un{j},\infty}^+$, defined as the unique subrepresentation of $\Theta_{\mu}(V_{\un{j},\un{j}'})$ with cosocle $V_{[j_2-j_0+1,j_2+1],\Delta}^{\rm{alg}}$, is multiplicity free and satisfies \ref{it: sm top 0}, \ref{it: sm top 1}, \ref{it: sm top 2}. We already note that \ref{it: sm top 0} follows from the fact that $\Theta_{\mu}(V_{\un{j},\un{j}'})$ has socle $C_{\un{j}}$ (see the first statement of Lemma~\ref{lem: cube step 2} and note that $V_{\un{j},\un{j}'}$ is the representation denoted $V_0$ there in the case $j_1=j_0$).\bigskip

\textbf{Step $1$}: We prove that $V_{\un{j},\un{j}'}$ injects into $V_{\un{j},\infty}^+$.\\
Let $\pi^{\infty}$ be the unique smooth length $2$ representation of $G$ with socle $V_{[j_2-j_0+2,j_2+1],\Delta}^{\infty}$ and cosocle $V_{[j_2-j_0+1,j_2+1],\Delta}^{\infty}$, which is $G$-basic by Lemma~\ref{lem: length two trivial block}. By \ref{it: explicit induction 1} of Lemma~\ref{lem: explicit smooth induction} (applied with $j_1,j_2$ there being $j_0,j_2+1$) $\pi^{\infty}$ injects into $i_{\widehat{j}_0,\Delta}^{\infty}(\pi_{j_0,j_2+1}^{\infty})$, and thus $d(\pi^{\infty}, \pi_{j_0,j_2+1}^{\infty})=0$. Let $M$ be the unique length $2$ $U(\fg)$-module with socle $L(1)$ and cosocle $L(s_{j_0})$ (\ref{it: rabiotext 2} of Lemma~\ref{rabiotext}), then \ref{it: general construction 2} of Lemma~\ref{lem: general construction} implies that $\cF_{P_{\widehat{j}_0}}^G(M,\pi_{j_0,j_2+1}^{\infty})$ contains a unique subrepresentation $V$ which fits into a non-split extension $0\rightarrow C_{\un{j}'}\rightarrow V\rightarrow L(1)^\vee\otimes_E \pi^{\infty}\rightarrow 0$. Moreover, as $\pi^{\infty}$ injects into $i_{\widehat{j}_0,\Delta}^{\infty}(\pi_{j_0,j_2+1}^{\infty})$, we deduce from \ref{it: Ext1 OS socle} of Lemma~\ref{lem: Ext1 OS socle cosocle} (applied with $V_1=C_{\un{j}'}$ and $V_0=L(1)^\vee\otimes_E \pi^{\infty}$) that $V$ is uniserial of length $3$, with socle $C_{\un{j}'}$, cosocle $V_{[j_2-j_0+1,j_2+1],\Delta}^{\rm{alg}}$ and middle layer $V_{[j_2-j_0+2,j_2+1],\Delta}^{\rm{alg}}$. As $M$ is a quotient of $\Theta_{w_0s_{j_0}w_0}(L(s_{j_0}))$ by Proposition~\ref{prop: Jantzen middle}, we deduce from (\ref{equ: OS wall crossing}) that $\cF_{P_{\widehat{j}_0}}^G(M,\pi_{j_0,j_2+1}^{\infty})$ is a subrepresentation of $\Theta_{\mu}(C_{\un{j}'})\cong \cF_{P_{\widehat{j}_0}}^G(\Theta_{w_0s_{j_0}w_0}(L(s_{j_0})),\pi_{j_0,j_2+1}^{\infty})$. Hence $V$ is also a subrepresentation of $\Theta_{\mu}(C_{\un{j}'})$, and is necessarily the unique subrepresentation with cosocle $V_{[j_2-j_0+1,j_2+1],\Delta}^{\rm{alg}}$ (note that $V_{[j_2-j_0+1,j_2+1],\Delta}^{\rm{alg}}$ occurs with multiplicity $1$ in $\Theta_{\mu}(C_{\un{j}'})$ arguing as in the first paragraph of the proof of Lemma \ref{lem: cube step 5}). Let us consider the surjection
\begin{equation}\label{equ: cosocle wall crossing}
\Theta_{\mu}(V_{\un{j},\un{j}'})\twoheadrightarrow \Theta_{\mu}(C_{\un{j}'})
\end{equation}
induced by $V_{\un{j},\un{j}'}\twoheadrightarrow C_{\un{j}'}$. As $V_{[j_2-j_0+1,j_2+1],\Delta}^{\rm{alg}}$ occurs with multiplicity $1$ in both $\Theta_{\mu}(V_{\un{j},\un{j}'})$ and $\Theta_{\mu}(C_{\un{j}'})$, it follows from the definition of $V_{\un{j},\infty}^+$ that $V$ is the image of $V_{\un{j},\infty}^+$ under (\ref{equ: cosocle wall crossing}). In particular $V_{\un{j},\infty}^+$ contains the constituent $C_{\un{j}'}$. As $C_{\un{j}'}$ does not occur in the kernel of (\ref{equ: cosocle wall crossing}) by \ref{it: cube step 1 1} of Lemma~\ref{lem: cube step 1}, it follows from \ref{it: cube step 1 2} of Lemma~\ref{lem: cube step 1} that $C_{\un{j}'}$ has multiplicity $1$ in the inverse image $\widetilde V$ of $V$ via (\ref{equ: cosocle wall crossing}) (which contains $V_{\un{j},\infty}^+$). Note that, by the same proof as for $V_{\un{j},\infty}^+$, the adjunction map $V_{\un{j},\un{j}'}\rightarrow \Theta_{\mu}(V_{\un{j},\un{j}'})$ is injective. By functoriality of $\Theta_\mu$ the composition of (\ref{equ: cosocle wall crossing}) with the injection $V_{\un{j},\un{j}'}\hookrightarrow \Theta_{\mu}(V_{\un{j},\un{j}'})$ factors through $V_{\un{j},\un{j}'}\twoheadrightarrow C_{\un{j}'}\hookrightarrow \Theta_{\mu}(C_{\un{j}'})$. Hence the image of $V_{\un{j},\un{j}'}$ under (\ref{equ: cosocle wall crossing}) is just its cosocle $C_{\un{j}'}$, and we deduce that $V_{\un{j},\un{j}'}$ is the unique subrepresentation of $\widetilde V$ with cosocle $C_{\un{j}'}$. Since $V_{\un{j},\infty}^+$ is a subrepresentation of $\widetilde V$ which contains $C_{\un{j}'}$, it follows that $V_{\un{j},\un{j}'}$ is also the unique subrepresentation of $V_{\un{j},\infty}^+$ with cosocle $C_{\un{j}'}$. In particular $V_{\un{j},\un{j}'}$ injects into $V_{\un{j},\infty}^+$.\bigskip

\textbf{Step $2$}: We prove that $V_{\un{j},\un{j}''}$ injects into $V_{\un{j},\infty}^+$.\\
Let $V'$ be the unique length $2$ representation of $G$ with socle $C_{(j_0,j_0+1,j_2+1)}$ and cosocle $C_{(j_0,j_0,j_2+1)}$ (Lemma~\ref{lem: Ext1 factor 1}). Note that we have a canonical surjection $V_{\un{j},\un{j}'}\twoheadrightarrow V'$ by \ref{it: hard square 2} of Proposition \ref{prop: hard square}. By \ref{it: loc alg crossing 1} of Remark~\ref{rem: loc alg crossing} (applied with $j_0,j_2$ there being $j_0,j_2+1$), the representation $V_{(j_0,j_0+1,j_2+1),\infty}$ in \ref{it: easy square 3} of Proposition \ref{prop: easy square} is isomorphic to the unique subrepresentation of $\Theta_{\mu}(V')$ with cosocle $V_{[j_2-j_0+1,j_2+1],\Delta}^{\rm{alg}}$. Let us consider the surjection
\begin{equation}\label{equ: length two quotient crossing}
\Theta_{\mu}(V_{\un{j},\un{j}'})\twoheadrightarrow \Theta_{\mu}(V')
\end{equation}
induced from $V_{\un{j},\un{j}'}\twoheadrightarrow V'$. As $V_{[j_2-j_0+1,j_2+1],\Delta}^{\rm{alg}}$ occurs with multiplicity $1$ in both $\Theta_{\mu}(V_{\un{j},\un{j}'})$ and $\Theta_{\mu}(V')$, it follows from the definition of $V_{\un{j},\infty}^+$ that $V_{(j_0,j_0+1,j_2+1),\infty}$ is the image of $V_{\un{j},\infty}^+$ under (\ref{equ: length two quotient crossing}). As $C_{\un{j}''}$ occurs in $V_{(j_0,j_0+1,j_2+1),\infty}$ and $\Theta_{\mu}(V_{\un{j},\un{j}'})$ with multiplicity $1$ (by \ref{it: easy square 3} of Proposition~\ref{prop: easy square} and \ref{it: cube step 1 2} of Lemma~\ref{lem: cube step 1} respectively), we deduce that $C_{\un{j}''}$ occurs in $V_{\un{j},\infty}^+$ with multiplicity $1$. Now let $V''$ be the unique length $2$ representation of $G$ with socle $C_{(j_0,j_0,j_2)}$ and cosocle $C_{(j_0,j_0+1,j_2+1)}$ (Lemma~\ref{lem: Ext1 factor 1}). Note that we have a canonical injection $V''\hookrightarrow V_{\un{j},\un{j}'}$ by \ref{it: hard square 2} of Proposition \ref{prop: hard square}. By \ref{it: square crossing 2} of Lemma~\ref{lem: square as wall crossing} $V_{\un{j},\un{j}''}$ is isomorphic to the unique subrepresentation of $\Theta_{\mu}(V'')$ with cosocle $C_{\un{j}''}$. The injection $V''\hookrightarrow V_{\un{j},\un{j}'}$ induces an injection $\Theta_{\mu}(V'')\hookrightarrow \Theta_{\mu}(V_{\un{j},\un{j}'})$, which therefore allows us to identify $V_{\un{j},\un{j}''}$ with the unique subrepresentation of $\Theta_{\mu}(V_{\un{j},\un{j}'})$ with cosocle $C_{\un{j}''}$. Since $V_{\un{j},\infty}^+$ is a subrepresentation of $\Theta_{\mu}(V_{\un{j},\un{j}'})$ which contains $C_{\un{j}''}$, it follows that $V_{\un{j},\un{j}''}$ is also the unique subrepresentation of $V_{\un{j},\infty}^+$ with cosocle $C_{\un{j}''}$. In particular $V_{\un{j},\un{j}''}$ injects into $V_{\un{j},\infty}^+$.\bigskip

\textbf{Step $3$}: We prove that, when $j_0>1$, there does not exist $W\in\mathrm{JH}_{G}(\Theta_{\mu}(C_{(j_0,j_0-1,j_2)}))$ such that $\mathrm{Ext}_{G}^1(V_{I}^{\rm{alg}},W)\neq 0$ for some $I\in\{[j_2-j_0+2,j_2+1],[j_2-j_0+1,j_2+1]\}$.\\
By \ (\ref{equ: OS wall crossing}) \ and \ Lemma~\ref{lem: Hom OS} \ a \ constituent \ of \ $\Theta_{\mu}(C_{(j_0,j_0-1,j_2)})$ \ has \ the \ form $W=\cF_{P_{\Delta\setminus\{j_0-1\}}}^G(L(x),\pi_{j_0-1,j_2}^{\infty})$ where $L(x)$ is a constituent of $\Theta_{w_0s_{j_0}w_0}(L(w_{j_0-1,j_0}))$, and by Proposition~\ref{prop: Jantzen middle} and Remark~\ref{rem: explicit middle} we have $D_L(x)=\{j_0-1\}$. Assume on the contrary that $\mathrm{Ext}_{G}^1(V_{I}^{\rm{alg}},W)\neq 0$ for some $I\in\{[j_2-j_0+2,j_2+1],[j_2-j_0+1,j_2+1]\}$, then we have $d(V_{I,\Delta}^{\infty},\pi_{j_0-1,j_2}^{\infty})=0$ by the last statement in Remark \ref{rem: Ext1 OS distance}, which implies $I=[j_2-j_0+2,j_2]$ by \ref{it: explicit induction 1} of Lemma~\ref{lem: explicit smooth induction}, a contradiction.\bigskip

\textbf{Step $4$}:
We prove that the constituents in the following list
\begin{equation}\label{list}
\big(\mathrm{JH}_{G}(V_{\un{j},\un{j}'})\cup \mathrm{JH}_{G}(V_{\un{j},\un{j}''})\big)\amalg\left\{V_{[j_2-j_0+2,j_2+1],\Delta}^{\rm{alg}},V_{[j_2-j_0+1,j_2+1],\Delta}^{\rm{alg}}\right\}
\end{equation}
occur with multiplicity $1$ in $V_{\un{j},\infty}^+$.\\
Note first that, from Step $1$, Step $2$ and the definition of $V_{\un{j},\infty}^+$, all these constituents occur in $V_{\un{j},\infty}^+$. By the end of the first paragraph of the proof of Lemma \ref{lem: cube step 2}, we have a surjective adjunction map
\begin{equation}\label{equ: length two quotient crossingbis}
\Theta_{\mu}(V_{\un{j},\un{j}'})\twoheadrightarrow V_{\un{j},\un{j}'}.
\end{equation}
Since the cosocle $V_{[j_2-j_0+1,j_2+1],\Delta}^{\rm{alg}}$ of $V_{\un{j},\infty}^+$ doesn't occur in $V_{\un{j},\un{j}'}$ by \ref{it: hard square 2} of Proposition \ref{prop: hard square}, $V_{\un{j},\infty}^+$ is contained in the kernel of (\ref{equ: length two quotient crossingbis}). Hence it is enough to prove that all constituents in (\ref{list}) occur with multiplicity $\leq 1$ in this kernel. For $V_{[j_2-j_0+2,j_2+1],\Delta}^{\rm{alg}}$ and $V_{[j_2-j_0+1,j_2+1],\Delta}^{\rm{alg}}$, this holds because they occur with multiplicity $1$ in $\Theta_{\mu}(V_{\un{j},\un{j}'})$ (see above Step $1$). For the non-locally algebraic constituents in $\mathrm{JH}_{G}(V_{\un{j},\un{j}'})\cup \mathrm{JH}_{G}(V_{\un{j},\un{j}''})$, this also holds using \ref{it: cube step 1 2} of Lemma \ref{lem: cube step 1}. The discussion in the first paragraph of the proof shows that only $\Theta_{\mu}(C_{\un{j}})$ and $\Theta_{\mu}(C_{\un{j}'})$ have locally algebraic constituents given by the constituents of $L(1)^\vee\otimes_E i_{\widehat{j}_0,\Delta}^{\infty}(\pi_{j_0,j_2}^{\infty})$ and $L(1)^\vee\otimes_E i_{\widehat{j}_0,\Delta}^{\infty}(\pi_{j_0,j_2+1}^{\infty})$ respectively. By Lemma \ref{lem: explicit smooth induction} the only common constituent of these two locally algebraic representations is $V_{[j_2-j_0+1,j_2],\Delta}^{\rm{alg}}$, which occurs in $V_{\un{j},\un{j}'}$ by \ref{it: hard square 2} of Proposition \ref{prop: hard square}. Hence all these locally algebraic constituents finally occur with multiplicity $1$ in the kernel of (\ref{equ: length two quotient crossingbis}).\bigskip

\textbf{Step $5$}: We prove that $V_{\un{j},\infty}^+$ is multiplicity free and satisfies \ref{it: sm top 2}.\\
For $V_{\un{j},\infty}^+$ multiplicity free, by Step $4$ it is enough to prove
\begin{equation}\label{equ: sm top JH}
\mathrm{JH}_{G}(V_{\un{j},\infty}^+)=\big(\mathrm{JH}_{G}(V_{\un{j},\un{j}'})\cup \mathrm{JH}_{G}(V_{\un{j},\un{j}''})\big)\amalg\left\{V_{[j_2-j_0+2,j_2+1],\Delta}^{\rm{alg}},V_{[j_2-j_0+1,j_2+1],\Delta}^{\rm{alg}}\right\}.
\end{equation}
By Step $4$ the right hand side of (\ref{equ: sm top JH}) is contained in $\mathrm{JH}_{G}(V_{\un{j},\infty}^+)$. Assume on the contrary that (\ref{equ: sm top JH}) is a strict inclusion. Then $V_{\un{j},\infty}^+/(V_{\un{j},\un{j}'}+V_{\un{j},\un{j}''})$ admits a quotient containing a constituent $W$ which is not in the list (\ref{equ: sm top JH}). Taking such a quotient of \emph{minimal} length and using Step $4$, we can assume that it contains a length $2$ subrepresentation with socle $W$ and cosocle $V_{I}^{\rm{alg}}$ for some $I\in\{[j_2-j_0+2,j_2+1],[j_2-j_0+1,j_2+1]\}$ (since any quotient of $V_{\un{j},\infty}^+/(V_{\un{j},\un{j}'}+V_{\un{j},\un{j}''})$ has cosocle $V_{[j_2-j_0+1,j_2+1],\Delta}^{\rm{alg}}$ and since $V_{I}^{\rm{alg}}$ for $I\in\{[j_2-j_0+2,j_2+1],[j_2-j_0+1,j_2+1]\}$ are the only ``remaining'' constituents in (\ref{equ: sm top JH})). This forces $\mathrm{Ext}_{G}^1(V_{I}^{\rm{alg}},W)\neq 0$. Define $V'''$ as $C_{\un{j}}$ if $j_0=1$, and as the unique length $2$ representation of $G$ with socle $C_{\un{j}}$ and cosocle $C_{(j_0,j_0-1,j_2)}$ if $j_0>1$ (using Lemma~\ref{lem: Ext1 factor 1}). By \ref{it: hard square 2} of Proposition~\ref{prop: hard square} $\mathrm{ker}(V_{\un{j},\un{j}'}\rightarrow V')/V'''$ is locally algebraic and thus $\Theta_{\mu}(\mathrm{ker}(V_{\un{j},\un{j}'}\rightarrow V')/V''')=0$ by Proposition~\ref{prop: Jantzen middle} and (\ref{equ: OS wall crossing}). This together with the exactness of $\Theta_{\mu}$ allows us to identify $\Theta_{\mu}(V''')$ with the kernel of the surjection (\ref{equ: length two quotient crossing}).
We have seen in {Step $2$} that the image of $V_{\un{j},\infty}^+$ under (\ref{equ: length two quotient crossing}) is isomorphic to $V_{(j_0,j_0+1,j_2+1),\infty}$. By \ref{it: easy square 3} of Proposition~\ref{prop: easy square} (with $j_2$ there being $j_2+1$) $\mathrm{JH}_{G}(V_{(j_0,j_0+1,j_2+1),\infty})$ is contained in the right and side of (\ref{equ: sm top JH}). Since $\Theta_{\mu}(V''')$ is the kernel of (\ref{equ: length two quotient crossing}), we deduce that $W\in\mathrm{JH}_{G}(\Theta_{\mu}(V'''))\setminus (\mathrm{JH}_{G}(V_{\un{j},\un{j}'})\cup \mathrm{JH}_{G}(V_{\un{j},\un{j}''}))$. Moreover by Step $3$ $W\notin \mathrm{JH}_{G}(\Theta_{\mu}(C_{(j_0,j_0-1,j_2)}))$ when $j_0>1$. Therefore, by definition of $V'''$ (and the exactness of $\Theta_{\mu}$), we must have
\begin{equation*}
W\in\mathrm{JH}_{G}\big(\Theta_{\mu}(C_{\un{j}})\big)\setminus \big(\mathrm{JH}_{G}(V_{\un{j},\un{j}'})\cup \mathrm{JH}_{G}(V_{\un{j},\un{j}''})\big).
\end{equation*}
Since $W$ is a constituent of $\Theta_{\mu}(C_{\un{j}})$, by (\ref{equ: OS wall crossing}), \ref{it: OS property 2} of Theorem~\ref{prop: OS property} and Lemma~\ref{lem: Hom OS} it has the form $W=\cF_{P_x}^G(L(x),\sigma^{\infty})$ where $L(x)$ is a constituent of $\Theta_{w_0s_{j_0}w_0}(L(s_{j_0}))$ and $\sigma^{\infty}$ a constituent of $i_{\widehat{j}_0,I_x}^{\infty}(\pi_{j_0,j_2}^{\infty})$. If $x\neq 1$, then $\mathrm{Ext}_{G}^1(V_{I}^{\rm{alg}},W)\neq 0$ together Proposition~\ref{prop: Ext1 OS} and Remark~\ref{rem: Ext1 OS distance} force $\mathrm{Ext}_{U(\fg)}^1(L(x),L(1))\neq 0$ and $d(V_{I,\Delta}^{\infty},\sigma^{\infty})=0$. By Proposition~\ref{prop: Jantzen middle} and \ref{it: rabiotext 2} of Lemma~\ref{rabiotext} we deduce $x=s_{j_0}$, which forces $W=C_{\un{j}}$ and $\sigma^{\infty}=\pi_{j_0,j_2}^{\infty}$, a contradiction as $d(V_{I,\Delta}^{\infty},\pi_{j_0,j_2}^{\infty})=0$ if and only if $I=[j_2-j_0+1,j_2]$ (Lemma~\ref{lem: explicit smooth induction}) but we have $I\in\{[j_2-j_0+2,j_2+1],[j_2-j_0+1,j_2+1]\}$. Hence we have $x=1$, i.e.~$W=L(1)^\vee\otimes_E\sigma^{\infty}$. Since $W$ is not a constituent of $V_{\un{j},\un{j}''}$ and since $L(1)^\vee\otimes_E \pi_{\un{j},\un{j}''}^{\infty}$ is a subquotient of $V_{\un{j},\un{j}''}$ by \ref{it: easy square 2} of Proposition~\ref{prop: easy square}, it follows that
\[\sigma^{\infty}\in\mathrm{JH}_{G}(i_{\widehat{j}_0,\Delta}^{\infty}(\pi_{j_0,j_2}^{\infty}))\setminus\mathrm{JH}_{G}(\pi_{\un{j},\un{j}''}^{\infty})=\{V_{[j_2-j_0+1,j_2-1],\Delta}^\infty, V_{[j_2-j_0,j_2-1],\Delta}^\infty\}\]
with $V_{[j_2-j_0,j_2-1],\Delta}^\infty$ omitted when $j_0=j_2$ (the equality follows from Lemma \ref{lem: explicit smooth induction}). Write $\sigma^{\infty}=V_{I',\Delta}^\infty$, then one can check for each $I\in\{[j_2-j_0+2,j_2+1],[j_2-j_0+1,j_2+1]\}$ that $d(V_{I,\Delta}^{\infty},V_{I',\Delta}^\infty)>1$ (for instance using \cite{Or05}) and thus $\mathrm{Ext}_{G}^1(V_{I}^{\rm{alg}},W)=0$ by Lemma~\ref{lem: Ext1 from sm}, another contradiction. We conclude that $W$ cannot exist, and thus (\ref{equ: sm top JH}) holds. Hence $V_{\un{j},\infty}^+$ is multiplicity free. The statement \ref{it: sm top 2} comes from the fact $\Theta_{\mu}(V_{\un{j},\un{j}'})$ contains as a subquotient the unique non-split extension of $V_{[j_2-j_0+1,j_2+1],\Delta}^{\rm{alg}}$ by $V_{[j_2-j_0+2,j_2+1],\Delta}^{\rm{alg}}$ (see above Step $1$).
\end{proof}

\begin{rem}\label{rem: explicit sm top}
An inspection of Proposition~\ref{prop: hard square} and \ref{it: easy square 2} of Proposition~\ref{prop: easy square} shows that the (multiplicity free) subrepresentation $V_{\un{j},\un{j}'}+V_{\un{j},\un{j}''}$ of $V_{\un{j},\infty}^+$ must be the amalgamate sum of $V_{\un{j},\un{j}'}$ and $V_{\un{j},\un{j}''}$ over the length $3$ subrepresentation of $V_{\un{j},\infty}^+$ with socle $C_{\un{j}}$ and cosocle $C_{(j_0,j_0+1,j_2+1)}\oplus V_{[j_2-j_0+1,j_2],\Delta}^{\rm alg}$. This determines the partially ordered set $\mathrm{JH}_{G}(V_{\un{j},\un{j}'}+V_{\un{j},\un{j}''})$ completely. The partial order on $\mathrm{JH}_{G}(V_{\un{j},\infty}^+)$ is then determined by the one on $\mathrm{JH}_{G}(V_{\un{j},\un{j}'}+V_{\un{j},\un{j}''})$ and the relations $C_{\un{j}'}\leq V_{[j_2-j_0+2,j_2+1],\Delta}^{\rm alg}\leq V_{[j_2-j_0+1,j_2+1],\Delta}^{\rm alg}$ and $C_{\un{j}''}\leq V_{[j_2-j_0+1,j_2+1],\Delta}^{\rm alg}$, using that the representation $V_{(j_0,j_0+1,j_2+1),\infty}$ of \ref{it: easy square 3} of Proposition \ref{prop: easy square} is a quotient of $V_{\un{j},\infty}^+$ (see {Step $2$} of the proof of Lemma~\ref{lem: sm top}) and using $\mathrm{Ext}_{G}^1(V_{[j_2-j_0+2,j_2+1],\Delta}^{\rm alg},C_{\un{j}''})=0$ which follows from \ref{it: Ext1 with alg 1} of Lemma \ref{lem: Ext1 factor 2}.
\end{rem}

Recall from \ref{it: sm top 0} of Lemma~\ref{lem: sm top} that $V_{\un{j},\infty}^+$ has cosocle $V_{[j_2-j_0+1,j_2+1],\Delta}^{\rm alg}$ and socle $C_{\un{j}}$, and thus fits into a non-split extension $0\rightarrow \mathrm{rad}_{G}(V_{\un{j},\infty}^+)\rightarrow V_{\un{j},\infty}^+\rightarrow V_{[j_2-j_0+1,j_2+1],\Delta}^{\rm{alg}}\rightarrow 0$. In particular
\begin{equation}\label{equ: sm top Ext}
\mathrm{Ext}_{G}^1(V_{[j_2-j_0+1,j_2+1],\Delta}^{\rm alg},\mathrm{rad}_{G}(V_{\un{j},\infty}^+))\ne 0.
\end{equation}

\begin{lem}\label{lem: sm top without top}
The vector space in (\ref{equ: sm top Ext}) is one dimensional, and any multiplicity free finite length $V^{\flat}$ in $\mathrm{Rep}^{\rm{an}}_{\rm{adm}}(G)$ which satisfies $\mathrm{JH}_{G}(V^{\flat})=\mathrm{JH}_{G}(\mathrm{rad}_{G}(V_{\un{j},\infty}^+))$ as partially ordered sets must satisfy $V^{\flat}\cong \mathrm{rad}_{G}(V_{\un{j},\infty}^+)$.
\end{lem}
\begin{proof}
We borrow all notation from the proof of Lemma~\ref{lem: sm top}.

We prove that the vector space in (\ref{equ: sm top Ext}) is one dimensional, and it suffices to show that it has dimension at most one. Using Remark~\ref{rem: explicit sm top}, we check that there is a unique increasing filtration on $\mathrm{rad}_{G}(V_{\un{j},\infty}^+)$ whose only reducible graded pieces are the length $2$ quotient $U_0$ of $\mathrm{rad}_{G}(V_{\un{j},\infty}^+)$ with socle $C_{(j_0,j_0+1,j_2+1)}$ and cosocle $C_{\un{j}''}$, and the length $2$ subrepresentation $U_1$ of $\mathrm{rad}_{G}(V_{\un{j},\infty}^+)$ with socle $C_{\un{j}}$ and cosocle $V_{[j_2-j_0+1,j_2],\Delta}^{\rm{alg}}$.
By Lemma~\ref{lem: Ext1 from sm} and Lemma~\ref{lem: Ext sm St} we have $\Dim_E\mathrm{Ext}_{G}^1(V_{[j_2-j_0+1,j_2+1],\Delta}^{\rm{alg}},V_{[j_2-j_0+2,j_2+1],\Delta}^{\rm{alg}})=1$, and by \ref{it: Ext1 with alg 1} of Lemma~\ref{lem: Ext1 factor 2} (and Remark~\ref{rem: explicit sm top}) we have $\mathrm{Ext}_{G}^1(V_{[j_2-j_0+1,j_2+1],\Delta}^{\rm{alg}},W)=0$ for $W\in\mathrm{JH}_{G}(\mathrm{rad}_{G}(V_{\un{j},\infty}^+))\setminus (\mathrm{JH}_{G}(U_0)\sqcup \mathrm{JH}_{G}(U_1)\sqcup\{V_{[j_2-j_0+2,j_2+1],\Delta}^{\rm{alg}}\})$. Hence by an easy d\'evissage it is enough to prove
\begin{equation}\label{equ: sm top Ext vanishing}
\mathrm{Ext}_{G}^1(V_{[j_2-j_0+1,j_2+1],\Delta}^{\rm{alg}},U_i)=0\ \ \mathrm{for}\ i=0,1.
\end{equation}
The surjection $\mathrm{gr}^0(V_{(j_0,j_0+1,j_2+1),\infty})\twoheadrightarrow V_{[j_2-j_0+1,j_2+1],\Delta}^{\rm{alg}}$ induces an embedding
\[\mathrm{Ext}_{G}^1(V_{[j_2-j_0+1,j_2+1],\Delta}^{\rm{alg}},U_0)\hookrightarrow \mathrm{Ext}_{G}^1(\mathrm{gr}^0(V_{(j_0,j_0+1,j_2+1),\infty}),U_0).\]
But the minimality of $V_{(j_0,j_0+1,j_2+1),\infty}$ in \ref{it: easy square 3} of Proposition~\ref{prop: easy square} (applied with $j_2$ there being $j_2+1$) implies $\mathrm{Ext}_{G}^1(\mathrm{gr}^0(V_{(j_0,j_0+1,j_2+1),\infty}),U_0)=0$. Hence we deduce (\ref{equ: sm top Ext vanishing}) for $i=0$. Assume on the contrary that (\ref{equ: sm top Ext vanishing}) fails for $i=1$. As $\mathrm{Ext}_{G}^1(V_{[j_2-j_0+1,j_2+1],\Delta}^{\rm{alg}},C_{\un{j}})=0$ by \ref{it: Ext1 with alg 1} of Lemma~\ref{lem: Ext1 factor 2}, there must exist a uniserial length $3$ representation $\tld{U}_1$ containing $U_1$ with socle $C_{\un{j}}$, cosocle $V_{[j_2-j_0+1,j_2+1],\Delta}^{\rm{alg}}$ and middle layer $V_{[j_2-j_0+1,j_2],\Delta}^{\rm{alg}}$. By Lemma~\ref{lem: Ext1 from sm} and Lemma~\ref{lem: Ext sm St} there exists a unique length $2$ representation $L(1)\otimes_E \tau^{\infty}$ with $\tau^{\infty}$ of length $2$ with socle $V_{[j_2-j_0+1,j_2],\Delta}^{\infty}$ and cosocle $V_{[j_2-j_0+1,j_2+1],\Delta}^{\infty}$ (and thus $G$-basic by Lemma~\ref{lem: length two trivial block}). By Lemma~\ref{lem: explicit smooth induction} $\Hom_{G}(\tau^{\infty},i_{\widehat{j}_0,\Delta}^{\infty}(\pi_{j_0,j_2}^{\infty}))=0$, which together with the last statement in Remark~\ref{rem: Ext1 OS distance} forces $\mathrm{Ext}_{G}^1(L(1)^\vee\otimes_E \tau^{\infty},C_{\un{j}})=0$, contradicting the existence of $\tld{U}_1$.
Hence (\ref{equ: sm top Ext vanishing}) holds and thus the vector space in (\ref{equ: sm top Ext}) has dimension $1$.

Now we let $V^{\flat}$ in $\mathrm{Rep}^{\rm{an}}_{\rm{adm}}(G)$ be a multiplicity free finite length representation such that $\mathrm{JH}_{G}(V^{\flat})=\mathrm{JH}_{G}(\mathrm{rad}_{G}(V_{\un{j},\infty}^+))$ as partially ordered sets. By Remark~\ref{rem: one dim Hom} (recall $V_{\un{j},\un{j}'}$ is the representation $V_0$ there for $j_0=j_1$) $\Theta_{\mu}(V_{\un{j},\un{j}'})$ contains a unique subrepresentation isomorphic to $V_{\un{j},\un{j}'}$. By \ref{it: cube step 1 2} of Lemma~\ref{lem: cube step 1} and \ref{it: cube step 2 2} of Lemma~\ref{lem: cube step 2} (and the fact $\Theta_{\mu}(V_{\un{j},\un{j}'})$ has socle $C_{\un{j}}$) we see that $\Theta_{\mu}(V_{\un{j},\un{j}'})$ also contains a unique subrepresentation isomorphic to $V_{\un{j},\un{j}''}$. Finally, by the first paragraph of the proof of Lemma~\ref{lem: sm top}, recall that $V_{[j_2-j_0+1,j_2+1],\Delta}^{\rm{alg}}$ has multiplicity one in $\Theta_{\mu}(V_{\un{j},\un{j}'})$.

Given $V^{\flat}$ as above, the equality of partially ordered sets $\mathrm{JH}_{G}(V^{\flat})=\mathrm{JH}_{G}(\mathrm{rad}_{G}(V_{\un{j},\infty}^+))$ forces $V^{\flat}$ to have socle $C_{\un{j}}$. Lemma~\ref{lem: minimal is unique} then forces $V^{\flat}$ to contain both $V_{\un{j},\un{j}'}$ and $V_{\un{j},\un{j}''}$. Using the first statement of Proposition~\ref{prop: Jantzen middle}, (\ref{equ: OS wall crossing}) and the explicit description of $\mathrm{JH}_{G}(V^{\flat}/V_{\un{j},\un{j}'})=\mathrm{JH}_{G}(\mathrm{rad}_{G}(V_{\un{j},\infty}^+)/V_{\un{j},\un{j}'})$ (cf.~Remark~\ref{rem: explicit sm top}), we have $\Theta_{\mu}(V^{\flat}/V_{\un{j},\un{j}'})=0$ and thus $\Theta_{\mu}(V_{\un{j},\un{j}'})\buildrel\sim\over\rightarrow \Theta_{\mu}(V^{\flat})$. As $V^{\flat}$ and $\Theta_{\mu}(V_{\un{j},\un{j}'})$ have socle $C_{\un{j}}$, the canonical (non-zero) map $V^{\flat}\rightarrow \Theta_{\mu}(V^{\flat})\cong \Theta_{\mu}(V_{\un{j},\un{j}'})$ is injective, hence $V^{\flat}$ is a subrepresentation of $\Theta_{\mu}(V_{\un{j},\un{j}'})$. Finally the equality $\mathrm{JH}_{G}(V^{\flat})=\mathrm{JH}_{G}(\mathrm{rad}_{G}(V_{\un{j},\infty}^+))$ and the above unicity statements in $\Theta_{\mu}(V_{\un{j},\un{j}'})$ force $V^{\flat}\cong \mathrm{rad}_{G}(V_{\un{j},\infty}^+)$.
\end{proof}

\begin{lem}\label{lem: simple length three}
Let $\un{j}\in\mathbf{J}$ with $1\leq j_0=j_1<j_2<n$. Then there exists a unique uniserial length $3$ representation $V_{\un{j}}$ in $\mathrm{Rep}^{\rm{an}}_{\rm{adm}}(G)$ with socle $L(1)^\vee\otimes_E V_{[j_2-j_1,j_2-1],\Delta}^{\infty}$, cosocle $L(1)^\vee\otimes_E V_{[j_2-j_1+1,j_2],\Delta}^{\infty}$ and middle layer $C_{\un{j}}$.
\end{lem}
\begin{proof}
We write $I\defeq \Delta\setminus\{j_1\}$. By Lemma~\ref{lem: simple length three g} there exists a $\fz$-semi-simple uniserial $U(\fg)$-module $M$ of length $3$ with both socle and cosocle $L(1)$ and middle layer $L(s_{j_1})$. Let $M^-$ be the unique length $2$ $U(\fg)$-module with socle $L(s_{j_1})$ and cosocle $L(1)$ (using \ref{it: rabiotext 2} of Lemma~\ref{rabiotext}), then $M^-$ is isomorphic to the unique length $2$ quotient of $M$, and the unique length $2$ $U(\fg)$-submodule of $M$ is isomorphic to $(M^-)^\tau$ (see (\ref{tauduality} for the notation) which is the unique length $2$ $U(\fg)$-module with socle $L(1)$ and cosocle $L(s_{j_1})$ (again by \ref{it: rabiotext 2} of Lemma~\ref{rabiotext}).

We define $W\defeq \cF_{P_{I}}^{G}(M,\pi_{j_1,j_2}^{\infty})$, which is well defined by Remark~\ref{rem: non t semisimple} and comes with an injection $q_1:\cF_{P_{I}}^{G}(M^-,\pi_{j_1,j_2}^{\infty})\hookrightarrow W$ and a surjection $q_2: W\twoheadrightarrow \cF_{P_{I}}^{G}((M^-)^\tau,\pi_{j_1,j_2}^{\infty})$ such that $q_2\circ q_1$ has (simple) image $C_{\un{j}}\cong\cF_{P_{I}}^{G}(L(s_{j_1}),\pi_{j_1,j_2}^{\infty})$. As $j_1<j_2<n$, recall that $V_{[j_2-j_1,j_2-1],\Delta}^{\infty}\cong \mathrm{cosoc}_{G}(i_{I}^{\infty}(\pi_{j_1,j_2}^{\infty}))$ and $V_{[j_2-j_1+1,j_2],\Delta}^{\infty}\cong \mathrm{soc}_{G}(i_{I}^{\infty}(\pi_{j_1,j_2}^{\infty}))$ by \ref{it: explicit induction 1} of Lemma~\ref{lem: explicit smooth induction}. Now by Lemma~\ref{lem: general construction} $\cF_{P_{I}}^{G}(M^-,\pi_{j_1,j_2}^{\infty})$ (resp.~$\cF_{P_{I}}^{G}((M^-)^\tau,\pi_{j_1,j_2}^{\infty})$) admits a quotient $W_1$ (resp.~a subrepresentation $W_2$) which is uniserial of length $2$ with socle $L(1)^\vee\otimes_E \mathrm{cosoc}_{G}(i_{I}^{\infty}(\pi_{j_1,j_2}^{\infty}))\cong V_{[j_2-j_1,j_2-1],\Delta}^{\rm alg}$ and cosocle $C_{\un{j}}$ (resp.~with socle $C_{\un{j}}$ and cosocle $L(1)^\vee\otimes_E \mathrm{soc}_{G}(i_{I}^{\infty}(\pi_{j_1,j_2}^{\infty}))\cong V_{[j_2-j_1+1,j_2],\Delta}^{\rm alg}$). Taking the pushforward of $W$ along $\cF_{P_{I}}^{G}(M^-,\pi_{j_1,j_2}^{\infty})\twoheadrightarrow W_1$, and then the pullback along $W_2\hookrightarrow \cF_{P_{I}}^{G}((M^-)^\tau,\pi_{j_1,j_2}^{\infty})$ gives a length $3$ subquotient $V_{\un{j}}$ of $W$ which admits $W_1$ as a subrepresentation and $W_2$ as a quotient. As we have $d(V_{[j_2-j_1+1,j_2],\Delta}^{\infty},V_{[j_2-j_1,j_2-1],\Delta}^{\infty})=2>1$ by \cite{Or05}, we have by (the sentence before) (\ref{equ: reduce to sm Ext})
\[\mathrm{Ext}_{G}^1(V_{[j_2-j_1+1,j_2],\Delta}^{\rm alg}, V_{[j_2-j_1,j_2-1],\Delta}^{\rm alg})=0.\]
We also have by Lemma~\ref{lem: Ext1 factor 2}
\[\Dim_E \mathrm{Ext}_{G}^1(V_{[j_2-j_1+1,j_2],\Delta}^{\rm alg}, C_{\un{j}})=1=\Dim_E \mathrm{Ext}_{G}^1(C_{\un{j}},V_{[j_2-j_1,j_2-1],\Delta}^{\rm alg}).\]
All this implies (by an easy d\'evissage) that the representation $V_{\un{j}}$ is uniserial and unique.
\end{proof}

\subsection{Hooking all constituents together}\label{subsec: final}

We hook together the various $\mathrm{Ext}$-squares, $\mathrm{Ext}$-cubes and uniserial representations constructed in \S\ref{subsec: square} to define two complexes of explicit finite length coadmissible $D(G)$-modules $\mathbf{D}^\bullet$, $\tld{\mathbf{D}}^\bullet$ and explicit quasi-isomorphisms $\bigoplus_{\ell=0}^{n-1}H^{\ell}(\mathbf{D}^\bullet)[-\ell]\leftarrow\tld{\mathbf{D}}^\bullet\rightarrow \mathbf{D}^\bullet$ (Theorem \ref{thm: main split}).\bigskip

From now on we tacitly identify the set $\mathbf{J}$ of \S\ref{subsec: square} with the set $\{C_{\un{j}}\mid \un{j}\in\mathbf{J}\}$ of irreducible admissible representations of $G$ over $E$ defined in (\ref{cj}), and we recall that $\mathbf{J}$ is equipped with the partial order $\un{j}\leq \un{j}'$ if and only if $j_0\leq j_0'$, $j_2\leq j_2'$ and $j_2-j_1\leq j_2'-j_1'$. In order to take into account locally algebraic constituents, it is convenient to enlarge $\mathbf{J}$ and define
\[\tld{\mathbf{J}}\defeq \mathbf{J}\sqcup \{V_{[j,j'],\Delta}^{\rm{alg}}\mid 1\leq j\leq j'\leq n-1\}\]
(see (\ref{vIalg})) equipped with the (unique) weakest partial order such that
\begin{itemize}
\item the partial order on $\tld{\mathbf{J}}$ restricts to the one on $\mathbf{J}$;
\item $\left\{\begin{array}{lclccl}
C_{(j_0,j_0,j_2)}&< &V_{[j_2-j_0+1,j_2],\Delta}^{\rm{alg}}&< &C_{(j_0,j_0,j_2+1)}&\mathrm{for\ }1\leq j_0\leq j_2\leq n-1\\
V_{[1,j_0],\Delta}^{\rm{alg}}& < &C_{j_0+1,j_0+1,j_0+1}&&&\mathrm{for\ }1\leq j_0\leq n-2\\
V_{[j_2-j_0+1,j_2],\Delta}^{\rm{alg}}&< &V_{[j_2-j_0,j_2],\Delta}^{\rm{alg}}&&&\mathrm{for\ }1\leq j_0< j_2\leq n-1.
\end{array}\right.$
\end{itemize}
In particular we have $V_{[j_2-j_0+1,j_2],\Delta}^{\rm alg} \leq C_{(j_0,j_0,j_2+1)} \leq V_{[j_2-j_0+2,j_2+1],\Delta}^{\rm alg}$, and an easy induction shows that if $C_{(j_0,j_1,j_2)} \leq V_{[j'_2-j'_0+1,j'_2],\Delta}^{\rm{alg}}$ then $j_0\leq j'_0$. For $V,W\in\tld{\mathbf{J}}$ such that $V\leq W$, we then define
\begin{equation}\label{intervalVW}
[V,W]\defeq \{V'\in\tld{\mathbf{J}}\mid V\leq V'\leq W\}\subseteq \tld{\mathbf{J}},
\end{equation}
and for $1\leq j_0\leq n-1$ we define
\begin{equation}\label{tldJj0def}
\tld{\mathbf{J}}_{j_0}\defeq [C_{(j_0,1,1)},C_{(j_0,1,n)}]\subseteq \tld{\mathbf{J}},
\end{equation}
which can be more explicitly described as
\begin{equation}\label{tldJj0}
\tld{\mathbf{J}}_{j_0}=\{C_{(j_0,j_1,j_2)}\mid (j_0,j_1,j_2)\in \mathbf{J}\} \sqcup \{V_{[j_2-j_0+1,j_2],\Delta}^{\rm{alg}}\mid j_0\leq j_2\leq n-1\}.
\end{equation}
Let $\un{j},\un{j}'\in\mathbf{J}$ such that $\un{j}\leq \un{j}'$ and
\begin{equation}\label{equ: cube condition}
\max\{j_0'-j_0, j_2'-j_2, (j_2'-j_1')-(j_2-j_1)\}\leq 1,
\end{equation}
and note that the corresponding $C_{\un{j}'}$ are exactly the $8$ constituents above Lemma \ref{lem: cube step 1}. When $d(\un{j},\un{j}')= |j_0-j_0'|+|j_2-j_2'|+|(j_2-j_1)-(j_2'-j_1')|\geq 2$ we have defined in \ref{it: easy square 1}, \ref{it: easy square 2} of Proposition~\ref{prop: easy square}, Proposition~\ref{prop: hard square} and Proposition~\ref{prop: factor cube} a minimal $\mathrm{Ext}$-hypercube $V_{\un{j},\un{j}'}$ of socle $C_{\un{j}}$ and cosocle $C_{\un{j}'}$. When $d(\un{j},\un{j}')=1$ we let $V_{\un{j},\un{j}'}$ be the unique non-split extension of $C_{\un{j}'}$ by $C_{\un{j}}$ in Lemma \ref{lem: Ext1 factor 1}, and when $d(\un{j},\un{j}')=0$ we set $V_{\un{j},\un{j}'}\defeq C_{\un{j}}$. From the definition of the partial order on $\tld{\mathbf{J}}$, we immediately check:

\begin{lem}
Let $\un{j},\un{j}'\in\mathbf{J}$ satisfying $\un{j}\leq \un{j}'$ and the bound (\ref{equ: cube condition}), we have $\mathrm{JH}_{G}(V_{\un{j},\un{j}'})=[C_{\un{j}},C_{\un{j}'}]$ as partially ordered sets.
\end{lem}

For $\un{j}'\in\mathbf{J}$, we set $\tld{\mathbf{J}}(C_{\un{j}'})\defeq [C_{\un{j}},C_{\un{j}'}]\subseteq \tld{\mathbf{J}}$ where $\un{j}\in\mathbf{J}$ is the (unique) minimal element that satisfies $\un{j}\leq \un{j}'$ and (\ref{equ: cube condition}) (for instance if $j'_0, j'_2\geq 2$ we have $\un{j}=(j'_0-1, j'_1, j'_2-1)$, the remaining cases are easily worked out). For $2\leq j_0\leq j_2\leq n-1$, we set
\[\tld{\mathbf{J}}(V_{[j_2-j_0+1,j_2],\Delta}^{\rm{alg}})\defeq [C_{(j_0-1,j_0-1,j_2-1)},V_{[j_2-j_0+1,j_2],\Delta}^{\rm{alg}}]=\mathrm{JH}_{G}(V_{(j_0-1,j_0-1,j_2-1),\infty}^+)\]
where the last identification (as partially ordered sets) follows from Lemma~\ref{lem: sm top} and Remark \ref{rem: explicit sm top}. For $2\leq j_2\leq n-1$ (and $j_0=1$) we set $\tld{\mathbf{J}}(V_{\{j_2\},\Delta}^{\rm{alg}})\defeq [V_{\{j_2-1\},\Delta}^{\rm{alg}},V_{\{j_2\},\Delta}^{\rm{alg}}]$, and (when $j_2=j_0=1$) $\tld{\mathbf{J}}(V_{\{1\},\Delta}^{\rm{alg}})\defeq [C_{(1,1,1)},V_{\{1\},\Delta}^{\rm{alg}}]$ (in fact, in each case we have $\tld{\mathbf{J}}(V_{[j_2-j_0+1,j_2],\Delta}^{\rm{alg}})= \mathrm{JH}_{G}(V_{(j_0-1,j_0-1,j_2-1),\infty}^+)$ where we keep the constituents of $V_{(j_0-1,j_0-1,j_2-1),\infty}^+$ which ``remain'').

\begin{lem}\label{lem: lower rep}
Let $\un{j}\in\mathbf{J}$, $j'_0\in \{1,\dots,n-1\}$ and $W\in\tld{\mathbf{J}}_{j_0'}$ such that $C_{\un{j}}\leq W$.
\begin{enumerate}[label=(\roman*)]
\item \label{it: lower rep 1} There exists a unique multiplicity free finite length representation $V_{\un{j},W}$ in $\mathrm{Rep}^{\rm{an}}_{\rm{adm}}(G)$ such that $\mathrm{JH}_{G}(V_{\un{j},W})=[C_{\un{j}},W]\cap\tld{\mathbf{J}}(W)$ as partially ordered sets.
\item \label{it: lower rep 2} Let $V_{\un{j},<W}$ be the unique subrepresentation of $V_{\un{j},W}$ such that $V_{\un{j},W}/V_{\un{j},<W}\cong W$ and assume $V_{\un{j},<W}\neq 0$, then $V_{\un{j},<W}$ is the unique multiplicity free finite length representation in $\mathrm{Rep}^{\rm{an}}_{\rm{adm}}(G)$ such that $\mathrm{JH}_G(V_{\un{j},<W})=\mathrm{JH}_G(V_{\un{j},W})\setminus\{W\}$ as partially ordered sets. Moreover we have
\begin{equation}\label{equ: lower rep dim one}
\Dim_E \mathrm{Ext}_{G}^1(W,V_{\un{j},<W})=1.
\end{equation}
\end{enumerate}
\end{lem}
\begin{proof}
Assume first $W=C_{\un{j}'}$ for some $\un{j}'\in\tld{\mathbf{J}}_{j_0'}$ (so we have $\un{j}\leq \un{j}'$). Then $[C_{\un{j}},W]\cap\tld{\mathbf{J}}(W)=[C_{\un{j}''},W]$ for the (unique) minimal $\un{j}''=(j_0'',j_1'',j_2'')\in\mathbf{J}$ such that $\un{j}\leq \un{j}''\leq \un{j}'$ and $\max\{j_0'-j_0'',j_2'-j_2'',(j_2'-j_1')-(j_2''-j_1'')\}\leq 1$. It follows from the definition of the partial order on $\mathbf{J}$ that $\un{j}''$ is the unique element in $\mathbf{J}$ such that
\begin{multline*}
j_0'-j_0''=\min\{1,j_0'-j_0\},\ j_2'-j_2''=\min\{1,j_2'-j_2\},\\
(j_2'-j_1')-(j_2''-j_1'')=\min\{1,(j_2'-j_1')-(j_2-j_1)\}.
\end{multline*}
In particular $d(\un{j}'',\un{j}')\in \{0,1,2,3\}$. We have $d(\un{j}'',\un{j}')=0$ if and only if $\un{j}''=\un{j}'$ if and only if $\un{j}=\un{j}'$, in which case $V_{\un{j},W}\cong C_{\un{j}}$. If $d(\un{j}'',\un{j}')=1$ then $V_{\un{j},W}$ is the unique length $2$ representation with socle $C_{\un{j}''}$ and cosocle $C_{\un{j}'}$ (Lemma~\ref{lem: Ext1 factor 1}). If $d(\un{j}'',\un{j}')=2$ then $V_{\un{j},W}$ is the unique minimal $\mathrm{Ext}$-square $V_{\un{j}'',\un{j}'}$ constructed in \ref{it: easy square 1}, \ref{it: easy square 2} of Proposition~\ref{prop: easy square} and Proposition~\ref{prop: hard square}. Finally if $d(\un{j}'',\un{j}')=3$ then $V_{\un{j},W}$ is the unique minimal $\mathrm{Ext}$-cube $V_{\un{j}'',\un{j}'}$ constructed in Proposition~\ref{prop: factor cube}). This proves \ref{it: lower rep 1} in that case. Assume now $V_{\un{j},<W}\neq 0$, i.e.~$\un{j}''<\un{j}'$. The unicity of $V_{\un{j},<W}$ is obvious when $d(\un{j}'',\un{j}')=1$ (as $V_{\un{j},<W}\cong C_{\un{j}''}$), and follows from the minimality of $V_{\un{j}'',\un{j}'}$ and Lemma~\ref{lem: minimal is unique} when $d(\un{j}'',\un{j}')\geq 2$. Finally (\ref{equ: lower rep dim one}) follows from Lemma~\ref{lem: Ext1 factor 1} when $d(\un{j}'',\un{j}')=1$, and from Lemma~\ref{lem: minimal vers Ext} (together with the minimality of $V_{\un{j}'',\un{j}'}$) when $d(\un{j}'',\un{j}')\geq 2$. This proves \ref{it: lower rep 2} in the case $W=C_{\un{j}'}$.\bigskip

Assume now $W=V_{[j_2'-j_0'+1,j_2'],\Delta}^{\rm{alg}}$ for some $j_0'\leq j_2'\leq n-1$. If $j_0=j_0'=1$, then $[C_{\un{j}},W]\cap\tld{\mathbf{J}}(W)$ is either $\tld{\mathbf{J}}(W)$ or $[C_{(1,1,j_2')},W]$, in which case the existence and unicity of $V_{\un{j},W}$ follows from Lemma~\ref{lem: simple length three} and \ref{it: Ext1 with alg 1} of Lemma~\ref{lem: Ext1 factor 2} respectively. Assume from now $j_0'>1$ so that $\tld{\mathbf{J}}(W)=[C_{(j_0'-1,j_0'-1,j_2'-1)},W]=\mathrm{JH}_{G}(V_{(j_0'-1,j_0'-1,j_2'-1),\infty}^+)$. From the definition of the partial order on $\tld{\mathbf{J}}$ one checks that $C_{\un{j}}\leq V_{[j_2'-j_0'+1,j_2'],\Delta}^{\rm{alg}}$ if and only if either $\un{j}\leq (j_0',j_0',j_2'-1)$ or $\un{j}\leq (j_0'-1,j_0'-1,j_2'-1)$. Hence the intersection $[C_{\un{j}},W]\cap\tld{\mathbf{J}}(W)=[C_{\un{j}},W]\cap [C_{(j_0'-1,j_0'-1,j_2'-1)},W]$ is equal to $[C_{\un{j}''},W]$ for some $\un{j}''\in \mathbf{J}$ such that $\un{j}\leq \un{j}''$, $(j_0'-1,j_0'-1,j_2'-1)\leq \un{j}''$ and either $\un{j}''\leq (j_0',j_0',j_2'-1)$ or $\un{j}''\leq (j_0'-1,j_0'-1,j_2'-1)$. We see that we necessarily have:
\begin{multline*}
\un{j}''\in \{(j_0'-1,j_0'-1,j_2'-1),(j_0'-1,j_0',j_2'),(j_0'-1,j_0'-2,j_2'-1),(j_0'-1,j_0'-1,j_2'),\\
(j_0',j_0'-1,j_2'-1),(j_0',j_0',j_2')\}.
\end{multline*}
But from Remark \ref{rem: explicit sm top} any constituent $C_{\un{j}''}$ for $\un{j}''$ in the above set is a constituent of $V_{(j_0'-1,j_0'-1,j_2'-1),\infty}^+$, hence it follows that $[C_{\un{j}},W]\cap\tld{\mathbf{J}}(W)=[C_{\un{j}''},W]=\mathrm{JH}_{G}(Q_{\un{j}''})$ (as partially ordered sets) where $Q_{\un{j}''}$ is the unique quotient of $V_{(j_0'-1,j_0'-1,j_2'-1),\infty}^+$ with socle $C_{\un{j}''}$ (and cosocle $W=V_{[j_2'-j_0'+1,j_2'],\Delta}^{\rm{alg}}$). We thus take $V_{\un{j},W}= Q_{\un{j}''}$ and $V_{\un{j},<W}= \mathrm{rad}_{G}(Q_{\un{j}''})$. It remains to prove \ref{it: lower rep 2} (which implies the unicity in \ref{it: lower rep 1} once we know the existence in \ref{it: lower rep 1}). If $\un{j}''=(j_0'-1,j_0'-1,j_2'-1)$ then $\mathrm{JH}_{G}(Q_{\un{j}''})=\mathrm{JH}_{G}(V_{(j_0'-1,j_0'-1,j_2'-1),\infty}^+)$, and \ref{it: lower rep 2} follows from Lemma~\ref{lem: sm top without top}. If $\un{j}''=(j_0'-1,j_0',j_2')$ then $\mathrm{JH}_{G}(Q_{\un{j}''})=\mathrm{JH}_{G}(V_{(j_0'-1,j_0',j_2'),\infty})$ (see Step $2$ of the proof of Lemma~\ref{lem: sm top} and recall that $V_{(j_0'-1,j_0',j_2'),\infty}$ is defined in \ref{it: easy square 3} of Proposition~\ref{prop: easy square}), and \ref{it: lower rep 2} follows from the existence and minimality of the $\mathrm{Ext}$-square $V_{(j_0'-1,j_0',j_2'),\infty}$ (\ref{it: easy square 3} of Proposition~\ref{prop: easy square}) and from Lemma~\ref{lem: minimal is unique}. In the remaining $4$ cases for $\un{j}''$, $Q_{\un{j}''}$ and $\mathrm{rad}_{G}(Q_{\un{j}''})$ are uniserial, and \ref{it: lower rep 2} follows by d\'evissage from Lemma \ref{lem: Ext1 factor 1}, \ref{it: Ext1 with alg 1} of Lemma~\ref{lem: Ext1 factor 2}, Lemma~\ref{lem: Ext1 from sm} and Lemma~\ref{lem: Ext sm St}.
\end{proof}

It follows from (\ref{tldJj0}) that $\tld{\mathbf{J}}_{1}=\mathbf{J} \sqcup \{V_{\{j\},\Delta}^{\rm{alg}}\mid 1\leq j\leq n-1\}$.

\begin{lem}\label{lem: vanishing outside cube}
Let $j_0, j_0'\in \{1,\dots,n-1\}$ such that $j_0\leq j_0'\leq j_0+1$ and $V\in \tld{\mathbf{J}}_{j_0}$, $V'\in\tld{\mathbf{J}}_{j_0'}$ such that $V<V'$.
\begin{enumerate}[label=(\roman*)]
\item \label{it: vanishing outside cube 1}
If $V\notin \tld{\mathbf{J}}(V')$, then we have $\mathrm{Ext}_{G}^1(V',V)=0$.
\item \label{it: vanishing outside cube 2}
If $V\notin \tld{\mathbf{J}}(V')$ and $\mathrm{Ext}_{G}^2(V',V)\neq 0$, then either $j_0=j'_0=1$, or $V'\cong C_{(j_0+1,j_0+1,j_2+2)}$ and $V\cong V_{[j_2-j_0+1,j_2],\Delta}^{\rm{alg}}$ for some $1\leq j_0\leq j_2\leq n-2$.
\end{enumerate}
\end{lem}
\begin{proof}
Assume first that $V$ and $V'$ are both not locally algebraic. Then $V=C_{\un{j}}$ and $V'=C_{\un{j}'}$ with ${\un{j}}<{\un{j}'}$ in ${\mathbf{J}}$. The assumption $V\notin \tld{\mathbf{J}}(V')$ then implies that (\ref{equ: cube condition}) fails, and it follows from Lemma~\ref{lem: Ext1 factor 1} and \ref{it: Ext2 factor 1 1} of Lemma \ref{lem: Ext2 factor 1} that 
\begin{equation}\label{equ: vanishing outside cube}
\mathrm{Ext}_{G}^1(V',V)=\mathrm{Ext}_{G}^2(V',V)=0.
\end{equation}

Assume $V=C_{\un{j}}=C_{(j_0,j_1,j_2)}$ and $V'$ is locally algebraic. Then there exists $j_2'\in \{j'_0, \dots,n-1\}$ such that $V'\cong V_{[j_2'-j_0'+1,j_2'],\Delta}^{\rm{alg}}$. An examination of the partial order on $\tld{\mathbf{J}}$ shows that $V<V'$ implies $j_2\leq j'_2\leq n-1$, and thus $I_{j_1,j_2}^+=[j_2-j_1+1,j_2]$ by Lemma \ref{lem: explicit smooth induction} (see (\ref{i+-chiant}) for $I_{j_1,j_2}^+$). If $I_{j_1,j_2}^+=[j_2'-j_0'+1,j_2']$ then $j_2=j'_2$ and $j_1=j'_0$, hence (using $j'_0\in\{j_0, j_0+1\}$) $\un{j}\in \{(j'_0,j'_0,j'_2), (j'_0-1,j'_0,j'_2)\}$ which implies $V=C_{\un{j}}\in \tld{\mathbf{J}}(V')$ and contradicts $V\notin \tld{\mathbf{J}}(V')$. Hence $I_{j_1,j_2}^+\ne [j_2'-j_0'+1,j_2']$ and thus $\mathrm{Ext}_{G}^1(V',V)=0$ by \ref{it: Ext1 with alg 1} of Lemma~\ref{lem: Ext1 factor 2}. If $\mathrm{Ext}_{G}^2(V',V)\ne 0$, then by \ref{it: Ext2 with alg 1} of Lemma \ref{lem: Ext2 factor 2} we have $j_0=j_1$ and $d([j_2'-j_0'+1,j_2'],[j_2-j_1+1,j_2])=d([j_2'-j_0'+1,j_2'],[j_2-j_0+1,j_2])=1$. Since $j_2\leq j'_2$ and $j'_0\in\{j_0, j_0+1\}$, the latter implies either $j'_2=j_2$ and $j'_0=j_0+1$, or $j'_2=j_2+1$ and $j'_0=j_0+1$. In the first case we have $V=C_{(j'_0-1,j'_0-1,j'_2)}\in \tld{\mathbf{J}}(V')$ and in the second $V=C_{(j'_0-1,j'_0-1,j'_2-1)}\in \tld{\mathbf{J}}(V')$. Hence both cases contradict $V\notin \tld{\mathbf{J}}(V')$, and we must have $\mathrm{Ext}_{G}^2(V',V)= 0$.\bigskip

Assume $V$ is locally algebraic, i.e.~$V\cong V_{[j_2-j_0+1,j_2],\Delta}^{\rm{alg}}$ for some $j_2\in \{j_0, \dots,n-1\}$, and $V'=C_{\un{j}'}=C_{(j'_0,j'_1,j'_2)}$. Then one checks that $V<V'$ implies $j'_2>j'_1$ and $j'_2\geq j_2+1$, and thus $I_{j'_1,j'_2}^-=[j'_2-j'_1,j'_2-1]$ by Lemma \ref{lem: explicit smooth induction} (see (\ref{i+-chiant}) for $I_{j'_1,j'_2}^-$). If $I_{j'_1,j'_2}^-=[j_2-j_0+1,j_2]$ then $j_2=j'_2-1$ and $j'_1=j_0$, hence (using $j'_0\in\{j_0, j_0+1\}$) $\un{j}'\in \{(j_0,j_0,j_2+1), (j_0+1,j_0,j_2+1)\}$ which implies $V=C_{\un{j}}\in \tld{\mathbf{J}}(V')$ and contradicts $V\notin \tld{\mathbf{J}}(V')$. Hence $I_{j'_1,j'_2}^-\ne [j_2-j_0+1,j_2]$ and thus $\mathrm{Ext}_{G}^1(V',V)=0$ by \ref{it: Ext1 with alg 2} of Lemma~\ref{lem: Ext1 factor 2}. If $\mathrm{Ext}_{G}^2(V',V)\ne 0$, then by \ref{it: Ext2 with alg 2} of Lemma \ref{lem: Ext2 factor 2} we have $j'_0=j'_1$ and $d([j'_2-j'_1,j'_2-1],[j_2-j_0+1,j_2])=d([j'_2-j'_0,j'_2-1],[j_2-j_0+1,j_2])=1$. Since $j_2+1\leq j'_2$ and $j'_0\in\{j_0, j_0+1\}$, the latter implies either $j'_2=j_2+1$ and $j'_0=j_0+1$, or $j'_2=j_2+2$ and $j'_0=j_0+1$. In the first case we have $V'\cong C_{(j_0+1,j_0+1,j_2+1)}$ and thus $V\cong V_{[j_2-j_0+1,j_2],\Delta}^{\rm{alg}}\in \tld{\mathbf{J}}(V')$ which contradicts $V\notin \tld{\mathbf{J}}(V')$. So the only possible case is $V\cong V_{[j_2-j_0+1,j_2],\Delta}^{\rm{alg}}$ and $V'\cong C_{(j_0+1,j_0+1,j_2+2)}$.
\bigskip

We now finally assume that $V$ and $V'$ are both locally algebraic, i.e.~$V\cong V_{[j_2-j_0+1,j_2],\Delta}^{\rm{alg}}$ and $V'\cong V_{[j_2'-j_0'+1,j_2'],\Delta}^{\rm{alg}}$ for some $j_2\in \{j_0, \dots,n-1\}$ and some $j_2'\in \{j'_0, \dots,n-1\}$. The assumption $V< V'$ forces $j_0 \leq j_0'$ and $j_2 \leq j_2'$ by an examination of the partial order on $\tld{\mathbf{J}}$. If $j_0'=1$ (and hence $j_0=1$), then we have $V\cong V_{\{j_2\},\Delta}^{\rm{alg}}$ and $V'\cong V_{\{j_2'\},\Delta}^{\rm{alg}}$ with $j_2<j_2'$, in which case we have
\[\mathrm{Ext}_{G}^1(V',V)=0\neq \mathrm{Ext}_{G}^2(V',V)\]
by Lemma~\ref{lem: Ext sm St} and the sentence before (\ref{equ: reduce to sm Ext}). In particular this finishes the proof of \ref{it: vanishing outside cube 1}. Assume from now $j_0'>1$, which forces $\tld{\mathbf{J}}(V')=\mathrm{JH}_{G}(V_{(j_0'-1,j_0'-1,j_2'-1),\infty}^+)$. If $j_2\in \{j_2',j'_2-1\}$ then the assumption $j_0\leq j_0'\leq j_0+1$ implies
\[V\in \{V_{[j_2'-j_0'+1,j_2'],\Delta}^{\rm{alg}},V_{[j_2'-j_0'+2,j_2'],\Delta}^{\rm{alg}}, V_{[j_2'-j_0',j_2'-1],\Delta}^{\rm{alg}},V_{[j_2'-j_0'+1,j_2'-1],\Delta}^{\rm{alg}}\}\]
which are precisely the $4$ locally algebraic constituents of $\tld{\mathbf{J}}(V')=\mathrm{JH}_{G}(V_{(j_0'-1,j_0'-1,j_2'-1),\infty}^+)$ by Remark \ref{rem: explicit sm top}, a contradiction with $V\notin \tld{\mathbf{J}}(V')$. Hence we have $j_2<j_2'-1$. The bounds $j_0'\leq j_0+1$ and $j_2<j_2'-1$ imply $j_2'-j_0'+1 \not\leq j_2-j_0+1$, hence in particular $j_2-j_0+1\in [j_2-j_0+1,j_2]\setminus [j_2'-j_0'+1,j_2']$, and the bounds $j_0'>1$, $j_2<j_2'-1$ imply $j_2',j_2'-1\in [j_2'-j_0'+1,j_2']\setminus [j_2-j_0+1,j_2]$. Thus $d([j_2'-j_0'+1,j_2'],[j_2-j_0+1,j_2])\geq 3$ by (\ref{di0i1}). Hence by Lemma~\ref{lem: Ext sm St} we have $d(V_{[j_2'-j_0'+1,j_2'],\Delta}^{\infty},V_{[j_2-j_0+1,j_2],\Delta}^{\infty})\geq 3$, which by the sentence before (\ref{equ: reduce to sm Ext}) implies (\ref{equ: vanishing outside cube}) in this case. By summarizing all the above cases, we have \ref{it: vanishing outside cube 2}.
\end{proof}

\begin{lem}\label{lem: unique rep reduction}
Let $\un{j}=(j_0,j_1,j_2)\in\mathbf{J}$, $j'_0\in \{1,\dots,n-1\}$ such that $j_0\leq j_0'\leq j_0+1$ and $W\in\tld{\mathbf{J}}_{j_0'}$ such that $C_{\un{j}}\leq W$. There exists at most one multiplicity free finite length representation $V$ in $\mathrm{Rep}^{\rm{an}}_{\rm{adm}}(G)$ such that $\mathrm{JH}_{G}(V)=[C_{\un{j}},W]$ as partially ordered sets. Moreover, if $j_0=j_0'>1$ then such a $V$ exists.
\end{lem}
\begin{proof}
We fix $\un{j}\in\mathbf{J}$ and let ${\mathcal C}_{\un{j}}\subseteq \tld{\mathbf{J}}$ be the subset of constituents $W$ such that $C_{\un{j}}\leq W$ and there exists $j'_0\in \{j_0, j_0+1\}$ with $j'_0\leq n-1$ such that $W\in\tld{\mathbf{J}}_{j_0'}$. As $W\in\tld{\mathbf{J}}_{j_0'}$, using (\ref{tldJj0def}) we have $C_{(j_0,1,1)}\leq C_{\un{j}}\leq W\leq C_{(j_0',1,n)}$. By (\ref{tldJj0def}) we have $[C_{(j_0,1,1)}, C_{(j'_0,1,n)}] = [C_{(j_0,1,1)},C_{(j_0,1,n)}]=\tld{\mathbf{J}}_{j_0}$ when $j_0=j_0'$, and an examination of the partial order on $\tld{\mathbf{J}}$ gives $[C_{(j_0,1,1)}, C_{(j'_0,1,n)}] = [C_{(j_0,1,1)},C_{(j_0,1,n)}] \sqcup [C_{(j'_0,1,1)}, C_{(j'_0,1,n)}]=\tld{\mathbf{J}}_{j_0}\sqcup \tld{\mathbf{J}}_{j'_0}$ when $j_0'=j_0+1$. Hence, we always have $[C_{\un{j}},W]\subseteq [C_{j_0,1,1},C_{j_0',1,n}]=\tld{\mathbf{J}}_{j_0}\cup \tld{\mathbf{J}}_{j'_0}$ for each $W\in\tld{\mathbf{J}}_{j'_0}$ such that $C_{\un{j}}\leq W$. Note that we have ${\mathcal C}_{\un{j}}=[C_{\un{j}},C_{j_0+1,1,n}]$ when $j_0<n-1$, and ${\mathcal C}_{\un{j}}=[C_{\un{j}},C_{n-1,1,n}]$ when $j_0=n-1$.

We now prove the statement by an increasing induction on $W\in {\mathcal C}_{\un{j}}$ (for the partial order on ${\mathcal C}_{\un{j}}$ induced by $\tld{\mathbf{J}}$). If $W=C_{\un{j}}$, we have $V=C_{\un{j}}$ and there is nothing to prove. We assume from now on $C_{\un{j}}<W$. By induction, for each $W'\in [C_{\un{j}},W]\setminus \{W\}$ there exists at most one multiplicity free finite length representation $V'$ in $\mathrm{Rep}^{\rm{an}}_{\rm{adm}}(G)$ such that $\mathrm{JH}_{G}(V')=[C_{\un{j}},W']$ as partially ordered sets. Taking the amalgamate sum of all such $V'$ for $W'\in [C_{\un{j}},W]\setminus \{W\}$ (noting that we need here the unicity in the induction hypothesis) we obtain a representation $\tld{V}_{\un{j},<W}$ in $\mathrm{Rep}^{\rm{an}}_{\rm{adm}}(G)$ which is the unique (if it exists) multiplicity free finite length representation such that $\mathrm{JH}_{G}(\tld{V}_{\un{j},<W})=[C_{\un{j}},W]\setminus \{W\}$ (as partially ordered sets). Replacing $[C_{\un{j}},W]\setminus \{W\}$ by $[C_{\un{j}},W]\setminus \tld{\mathbf{J}}(W)$ and noting that any $W''\in [C_j,W]$ such that $W'' \leq W'$ for some $W' \in [C_j,W]\setminus \tld{\mathbf{J}}(W)$ is still in $[C_j,W]\setminus \tld{\mathbf{J}}(W)$, we also obtain a unique (if it exists) multiplicity free finite length representation $\tld{V}_{\un{j},<W}'$ in $\mathrm{Rep}^{\rm{an}}_{\rm{adm}}(G)$ such that $\mathrm{JH}_{G}(\tld{V}_{\un{j},<W}')=[C_{\un{j}},W]\setminus \tld{\mathbf{J}}(W)$ and such that $\tld{V}_{\un{j},<W}'$ is a subrepresentation of $\tld{V}_{\un{j},<W}$. It then follows from \ref{it: lower rep 2} of Lemma \ref{lem: lower rep} that we have a short exact sequence (if $\tld{V}_{\un{j},<W}$ exists)
\[0\longrightarrow \tld{V}_{\un{j},<W}'\longrightarrow \tld{V}_{\un{j},<W}\longrightarrow V_{\un{j},<W}\rightarrow 0,\]
which in turn induces a long exact sequence
\begin{equation}\label{equ: unique rep seq}
\mathrm{Ext}_{G}^1(W, \tld{V}_{\un{j},<W}') \rightarrow \mathrm{Ext}_{G}^1(W, \tld{V}_{\un{j},<W}) \xrightarrow{q} \mathrm{Ext}_{G}^1(W, V_{\un{j},<W}) \rightarrow \mathrm{Ext}_{G}^2(W, \tld{V}_{\un{j},<W}').
\end{equation}
It follows from \ref{it: vanishing outside cube 1} of Lemma~\ref{lem: vanishing outside cube} that $\mathrm{Ext}_{G}^1(W, W')=0$ for each constituent $W'$ of $\tld{V}_{\un{j},<W}'$, which by d\'evissage implies $\mathrm{Ext}_{G}^1(W, \tld{V}_{\un{j},<W}')=0$. Then (\ref{equ: unique rep seq}) and (\ref{equ: lower rep dim one}) imply $\Dim_E\mathrm{Ext}_{G}^1(W, \tld{V}_{\un{j},<W})\leq 1$, which shows that $V$ as in the statement is unique if it exists. Note that if $j_0=j_0'>1$, then \ref{it: vanishing outside cube 2} of Lemma~\ref{lem: vanishing outside cube} implies $\mathrm{Ext}_{G}^2(W, \tld{V}_{\un{j},<W}')=0$ by an analogous d\'evissage and thus $q$ is an isomorphism. Going back through the induction above, we see that $V$ exists if $j_0=j_0'>1$.
\end{proof}

\begin{prop}\label{prop: general wall crossing}
Let $1\leq j_0\leq n-1$ and $\mu\in \Lambda$ such that $\langle\mu+\rho,\alpha^\vee\rangle\geq 0$ for $\al\in \Phi^+$ and the stabilizer of $\mu$ in $W(G)$ for the dot action is $\{1,w_0s_{j_0}w_0\}$.
\begin{enumerate}[label=(\roman*)]
\item \label{it: general wall crossing 1} For each representation $V\in\tld{\mathbf{J}}\setminus (\tld{\mathbf{J}}_{j_0}\cap\mathbf{J})$, we have $\Theta_{w_0s_{j_0}w_0}(V)=0$.
\item \label{it: general wall crossing 3} Assume $j_0>1$ and set $j_0^+\defeq \min\{j_0+1,n-1\}$. Let $(j_1,j_2)\in\{(n-1,n-1),(1,n)\}$ and $V_0$ be the unique multiplicity free finite length representation in $\mathrm{Rep}^{\rm{an}}_{\rm{adm}}(G)$ such that $\mathrm{JH}_{G}(V_0)=[C_{(j_0,1,1)},C_{(j_0,j_1,j_2)}]$ as partially ordered sets (see Lemma~\ref{lem: unique rep reduction}). Then we have a short exact sequence 
\begin{equation*}
0\longrightarrow V_0^+\longrightarrow \Theta_{\mu}(V_0)\longrightarrow V_0^-\longrightarrow 0, 
\end{equation*}
where $V_0^-$ (resp.~$V_0^+$) is the unique multiplicity free finite length representation in $\mathrm{Rep}^{\rm{an}}_{\rm{adm}}(G)$ such that $\mathrm{JH}_{G}(V_0^-)=[C_{(j_0-1,1,1)},C_{(j_0,j_1,j_2)}]$ (resp.~such that $\mathrm{JH}_{G}(V_0^+)=[C_{(n-1,1,1)},L(1)^\vee]$ when $j_0=j_1=j_2=n-1$, and $\mathrm{JH}_{G}(V_0^+)=[C_{(j_0,1,1)},C_{(j_0^+,j_1,j_2)}]$ otherwise) as partially ordered sets. 
\end{enumerate}
\end{prop}
\begin{proof}
\ref{it: general wall crossing 1} follows from \ref{it: wall crossing simple 1} of Lemma~\ref{lem: wall crossing of simple} and the fact that each $V\in\tld{\mathbf{J}}\setminus (\tld{\mathbf{J}}_{j_0}\cap\mathbf{J})$ is either locally algebraic or of the form $\cF_{P_{\widehat{j}_1}}^{G}(L(w_{j_1,j_0'}),\pi_{j_1,j_2}^{\infty})$ with $D_R(w_{j_1,j_0'})=\{j_0'\}$ and $j_0'\neq j_0$.\bigskip

We prove \ref{it: general wall crossing 3}. 
It follows from Lemma~\ref{lem: unique rep reduction} that $V_0$ exists (using $j_0>1$) and is unique, that $V_0^-$ and $V_0^+$ are unique if they exist, and that $V_0^+$ exists if $j_0=n-1=j_0^+$. By a decreasing induction on $1<j_0\leq n-1$, it suffices to assume the existence of $V_0^+$, construct an embedding $V_0^+\hookrightarrow \Theta_{\mu}(V_0)$, and then prove that $\Theta_{\mu}(V_0)/V_0^+$ is a multiplicity free finite length representation that satisfies $\mathrm{JH}_G(\Theta_{\mu}(V_0)/V_0^+)=[C_{(j_0-1,1,1)},C_{(j_0,j_1,j_2)}]$ as partially ordered sets (for instance, when $j_0=n-1$, then $\Theta_{\mu}(V_0)/V_0^+$ is actually the $V_0^+$ for $j_0=n-2$ and the induction goes on). Assume from now that $V_0^+$ exists. As $C_{(j_0,1,1)}\leq C_{(j_0,j_1,j_2)}\leq C_{(j_0^+,j_1,j_2)}$, $V_0$ is the unique subrepresentation of $V_0^+$ with cosocle $C_{(j_0,j_1,j_2)}$. We divide the rest of the proof of \ref{it: general wall crossing 3} into the following steps.\bigskip

\textbf{Step $1$}: We prove that the injection $V_0\hookrightarrow V_0^+$ induces an isomorphism $\Theta_{\mu}(V_0)\buildrel\sim\over\longrightarrow \Theta_{\mu}(V_0^+)$, and that the adjunction map $V_0^+\rightarrow \Theta_{\mu}(V_0^+)$ is injective.\\
Note that we have
\[\begin{array}{rclllll}
V_0^+/V_0&\cong &L(1)^\vee&\text{ if }&j_0=j_0^+=n-1&\text{ and }&(j_1,j_2)=(n-1,n-1)\\
V_0^+& = &V_0&\text{ if }&j_0=j_0^+=n-1&\text{ and }&(j_1,j_2)=(1,n),
\end{array}\]
and when $j_0^+=j_0+1\leq n-1$:
\[\mathrm{JH}_{G}(V_0^+/V_0)=[C_{(j_0,1,1)},C_{(j_0^+,j_1,j_2)}]\setminus [C_{(j_0,1,1)},C_{(j_0,j_1,j_2)}]=[C_{(j_0^+,1,1)},C_{(j_0^+,j_1,j_2)}]\subseteq \tld{\mathbf{J}}_{j_0'}.\]
Since in the third case $\tld{\mathbf{J}}_{j_0^+}\cap \tld{\mathbf{J}}_{j_0}=\emptyset$, we have $\mathrm{JH}_{G}(V_0^+/V_0)\subseteq \tld{\mathbf{J}}\setminus \tld{\mathbf{J}}_{j_0}$, and \emph{a fortiori} $\mathrm{JH}_{G}(V_0^+/V_0)\subseteq\tld{\mathbf{J}}\setminus (\tld{\mathbf{J}}_{j_0}\cap\mathbf{J})$. The first two cases also clearly satisfy $\mathrm{JH}_{G}(V_0^+/V_0)\subseteq\tld{\mathbf{J}}\setminus (\tld{\mathbf{J}}_{j_0}\cap\mathbf{J})$. Hence we deduce from \ref{it: general wall crossing 1} that $\Theta_{\mu}(V_0^+/V_0)=0$ and thus that the injection $V_0\hookrightarrow V_0^+$ induces an isomorphism $\Theta_{\mu}(V_0)\buildrel\sim\over\longrightarrow \Theta_{\mu}(V_0^+)$. We define $q^+$ as the composition $V_0^+\rightarrow \Theta_{\mu}(V_0^+)\cong \Theta_{\mu}(V_0)$ (the first map being the adjunction map). As $q^+$ restricts to an injection $C_{(j_0,1,1)}\hookrightarrow \Theta_{w_0s_{j_0}w_0}(C_{(j_0,1,1)})$ by (the first statement in) \ref{it: square crossing 1} of Lemma \ref{lem: square as wall crossing} and $V_0^+$ has socle $C_{(j_0,1,1)}$, the map $q^+$ is injective.\bigskip

\textbf{Step $2$}: We prove that $\Theta_{\mu}(V_0)/V_0^+$ is a multiplicity free finite length representation such that there is an equality of \emph{sets}
\begin{equation}\label{equ: wall crossing cube quotient JH}
\mathrm{JH}_{G}(\Theta_{\mu}(V_0)/V_0^+)=[C_{(j_0-1,1,1)},C_{(j_0,j_1,j_2)}].
\end{equation}
It suffices to show that we have an equality in the \emph{Grothendieck group} of finite length admissible representations of $G$
\begin{equation}\label{equ: wall crossing mult}
[\Theta_{\mu}(V_0)]=[C_{(j_0,1,1)},C_{(j_0^+,j_1,j_2)}]+[C_{(j_0-1,1,1)},C_{(j_0,j_1,j_2)}]
\end{equation}
(where now $[C_{(j_0,1,1)},C_{(j_0^+,j_1,j_2)}]$, $[C_{(j_0-1,1,1)},C_{(j_0,j_1,j_2)}]$ mean the direct sums of the corresponding constituents in this Grothendieck group, but this does not lead to confusions). By \ref{it: wall crossing simple 1} of Lemma~\ref{lem: wall crossing of simple} we have
\begin{equation*}
[\Theta_{\mu}(V_0)]=\sum_{(j_1',j_2')\leq (j_1,j_2)}[\Theta_{w_0s_{j_0}w_0}(C_{(j_0,j_1',j_2')})],
\end{equation*}
and we can compute $[\Theta_{w_0s_{j_0}w_0}(C_{(j_0,j_1',j_2')})]$ using \ref{it: square crossing 1} of Lemma~\ref{lem: square as wall crossing}. For each $(j_1',j_2')\leq (j_1,j_2)$, we check that $C_{(j_0,j_1',j_2')}$ (resp.~$C_{(j_0-1,j_1',j_2')}$) appears with multiplicity $2$ (resp.~$1$) on both sides of (\ref{equ: wall crossing mult}), and that $C_{(j_0+1,j_1',j_2')}$ appears with multiplicity $1$ on both sides of (\ref{equ: wall crossing mult}) when $j_0'=j_0+1\leq n-1$. We now treat the locally algebraic constituents. Assume first $(j_1,j_2)=(1,n)$. Using \ref{it: square crossing 1} of Lemma~\ref{lem: square as wall crossing} and Lemma~\ref{lem: explicit smooth induction}, we check that the contribution of locally algebraic constituents on the left hand side of (\ref{equ: wall crossing mult}) is
\begin{multline}\label{equ: LHS loc alg mult}
\sum_{j_2'=j_0}^{n}[L(1)^\vee\otimes_Ei_{\widehat{j}_0,\Delta}^{\infty}(\pi_{j_0,j_2'}^{\infty})]
=2\sum_{j_2'=j_0}^{n-1}[V_{[j_2'-j_0+1,j_2'],\Delta}^{\rm{alg}}]+\sum_{j_2'=j_0+1}^{n-1}[V_{[j_2'-j_0,j_2'],\Delta}^{\rm{alg}}]\\
+\sum_{j_2'=j_0}^{n}[V_{[j_2'-j_0+1,j_2'-1],\Delta}^{\rm{alg}}].
\end{multline}
Using the relations
\[\left\{\begin{array}{rclclll}
C_{(j_0,j_0,j_2')}&\leq &V_{[j_2'-j_0+1,j_2'],\Delta}^{\rm{alg}}&\leq &C_{(j_0,j_0,j_2'+1)}&\ \mathrm{for}\ &j_0\leq j_2'\leq n-1;\\C_{(j_0+1,j_0+1,j_2')}&\leq &V_{[j_2'-j_0,j_2'],\Delta}^{\rm{alg}}&\leq &C_{(j_0+1,j_0+1,j_2'+1)}&\ \mathrm{for}\ &j_0+1\leq j_2'\leq n-1;\\
C_{j_0-1,j_0-1,j_2'-1}&\leq &V_{[j_2'-j_0+1,j_2'-1],\Delta}^{\rm{alg}}&\leq &C_{j_0-1,j_0-1,j_2'}&\ \mathrm{for}\ &j_0\leq j_2'\leq n,
\end{array}\right.\]
we check that the contribution of the locally algebraic constituents in $[V_0^+]+[V_0^-]$ is given by exactly the same formula (\ref{equ: LHS loc alg mult}). Finally, when $(j_1,j_2)=(n-1,n-1)$, using again \ref{it: square crossing 1} of Lemma~\ref{lem: square as wall crossing}, Lemma~\ref{lem: explicit smooth induction} and the above relations, one checks that the contribution of locally algebraic constituents on both sides of (\ref{equ: wall crossing mult}) is $[L(1)\otimes_E i_{\hat{j}_0,\Delta}^{\infty}(\pi_{j_0,j_0}^{\infty})] =[V_{[1,j_0],\Delta}^{\rm alg}]+[V_{[1,j_0-1],\Delta}^{\rm alg}]$.\bigskip

\textbf{Step $3$}: We assume $(j_1,j_2)=(1,n)$ and prove that the partial order on $\mathrm{JH}_{G}(\Theta_{\mu}(V_0)/V_0^+)$ is at least as strong as the one on $[C_{(j_0-1,1,1)},C_{(j_0,1,n)}]$ induced by $\tld{\mathbf{J}}$.\\
Note that the underlying set of (\ref{equ: wall crossing cube quotient JH}) is explicitly given by
\begin{multline*}
\{C_{\un{j}'}\mid (j_0-1,1,1)\leq \un{j}'\leq (j_0,1,n)\}\sqcup \{V_{[j_2'-j_0+2,j_2'],\Delta}^{\rm{alg}}\mid 1\leq j_2'\leq n-1\}\\
\sqcup \{V_{[j_2'-j_0+1,j_2'],\Delta}^{\rm{alg}}\mid 2\leq j_2'\leq n-1\}.
\end{multline*}
Recall from \ref{it: square crossing 1} of Lemma~\ref{lem: square as wall crossing} that the adjunction map $\Theta_{\mu}(C_{(j_0,1,n)})\rightarrow C_{(j_0,1,n)}$ is surjective, and the composition $C_{(j_0,1,n)}\rightarrow \Theta_{\mu}(C_{(j_0,1,n)})\rightarrow C_{(j_0,1,n)}$ with the other adjunction map is zero. Since $V_0$ has cosocle $C_{(j_0,1,n)}$ and the adjunction maps $V_0\rightarrow \Theta_{\mu}(V_0)\rightarrow V_0$ are compatible with the adjunction maps $C_{(j_0,1,n)}\rightarrow \Theta_{\mu}(C_{(j_0,1,n)})\rightarrow C_{(j_0,1,n)}$ under the surjection $V_0\twoheadrightarrow C_{(j_0,1,n)}$, we deduce that the adjunction map $\Theta_{\mu}(V_0)\rightarrow V_0$ is surjective and the composition of $V_0\rightarrow \Theta_{\mu}(V_0)\rightarrow V_0$ is zero. If $j_0=j_0^+=n-1$, then $V_0=V_0^+$ and we obtain a surjection $\Theta_{\mu}(V_0)/V_0^+\twoheadrightarrow V_0$. If $j_0<n-1$, then as the cosocle $C_{(j_0+1,1,n)}$ of $V_0^+$ does not occur in $V_0$, this forces the composition $V_0^+\hookrightarrow \Theta_{\mu}(V_0) \twoheadrightarrow V_0$ to be zero and we also obtain a surjection $\Theta_{\mu}(V_0)/V_0^+\twoheadrightarrow V_0$. It follows that the partial order on $\mathrm{JH}_{G}(\Theta_{\mu}(V_0)/V_0^+)$ restricts to the one on $\mathrm{JH}_{G}(V_0)=[C_{(j_0,1,1)},C_{(j_0,1,n)}]$. In particular, inside $\mathrm{JH}_{G}(\Theta_{\mu}(V_0)/V_0^+)$ we have $C_{(j_0,j_1',j_2')}<C_{(j_0,j_1'',j_2'')}$ for $(j_1',j_2')<(j_1'',j_2'')\in\mathbf{J}^{\infty}$, and $C_{j_0,j_0,j_2'}<V_{[j_2'-j_0+1,j_2'],\Delta}^{\rm{alg}}<C_{j_0,j_0,j_2'+1}$ for $2\leq j_2'\leq n-1$ (as these relations occur in $[C_{(j_0,1,1)},C_{(j_0,1,n)}]$). 
Now, let $(j_1',j_2'),(j_1'',j_2'')\in\mathbf{J}^{\infty}$ with $(j_1'',j_2'')\in\{(j_1'+1,j_2'+1),(j_1'-1,j_2')\}$, then $V_0$ admits a unique length $2$ subquotient $V_1$ with socle $C_{(j_0,j_1',j_2')}$ and cosocle $C_{(j_0,j_1'',j_2'')}$ (Lemma \ref{lem: Ext1 factor 1}). 
By the last statement in Remark \ref{rem: opposite crossing} applied to $j_0$ (there) being $j_0-1$ (here), $\Theta_{\mu}(V_1)$ admits a unique subquotient isomorphic to $V_{(j_0-1,j_1',j_2'),(j_0,j_1'',j_2'')}$ as in \ref{it: easy square 1} or \ref{it: easy square 2} of Proposition~\ref{prop: easy square}. Since $\mathrm{soc}_{G}(V_{(j_0-1,j_1',j_2'),(j_0,j_1'',j_2'')})\cong C_{(j_0-1,j_1',j_2')}\notin\mathrm{JH}_{G}(V_0^+)$, it follows that $V_{(j_0-1,j_1',j_2'),(j_0,j_1'',j_2'')}$ is still a subquotient of $\Theta_{\mu}(V_0)/V_0^+$. In particular, from the structure of $V_{(j_0-1,j_1',j_2'),(j_0,j_1'',j_2'')}$ in \emph{loc.~cit.}, inside $\mathrm{JH}_{G}(\Theta_{\mu}(V_0)/V_0^+)$ we have $C_{(j_0-1,j_1',j_2')}< C_{(j_0-1,j_1'',j_2'')}< C_{(j_0,j_1'',j_2'')}$, $C_{(j_0-1,j_1',j_2')}< C_{(j_0,j_1',j_2')}$, and
\begin{equation}\label{equ: relation with loc alg 1}
\left\{\begin{array}{ccccc}
\scriptstyle{C_{(j_0-1,j_0-1,j_2')}<V_{[j_2'-j_0+2,j_2'],\Delta}^{\rm{alg}}<V_{[j_2'-j_0+1,j_2'],\Delta}^{\rm{alg}}<C_{(j_0,j_0,j_2'+1)}}&{\mathrm{ if }}&\scriptstyle{j_1=1}&{\mathrm{ and }}&\scriptstyle{j_2''=j_2'+1>j_0}\\
\scriptstyle{C_{(j_0-1,j_0-1,j_0-1)}<V_{[1,j_0-1],\Delta}^{\rm{alg}}<C_{(j_0,j_0,j_0)}}&{\mathrm{ if }}&j_1=1\scriptstyle{j_1=1}&{\mathrm{ and }}&\scriptstyle{j_2''=j_2'+1=j_0}.
\end{array}\right.
\end{equation}
In particular the partial order on $\mathrm{JH}_{G}(\Theta_{\mu}(V_0)/V_0^+)$ and on $[C_{(j_0-1,1,1)},C_{(j_0,1,n)}]$ have the same restriction to $[C_{(j_0-1,1,1)},C_{(j_0,1,n)}]\cap\mathbf{J}$ (i.e.~non locally algebraic constituents). We now deal with locally algebraic constituents. We prove by an increasing induction on $1\leq j_2'\leq n-1$ that, inside $\mathrm{JH}_{G}(\Theta_{\mu}(V_0)/V_0^+)$, $V_{[j_2'-j_0+2,j_2'],\Delta}^{\rm{alg}}$ is the only locally algebraic constituent $V$ such that:
\begin{equation}\label{equ: relation with loc alg 2}
C_{(j_0-1,j_0-1,j_2')}<V<C_{(j_0-1,j_0-1,j_2'+1)}.
\end{equation}
Note first that all constituents of the form $V_{[j_2'-j_0+1,j_2'],\Delta}^{\rm alg}$ occur in $V_0$ and hence can't lie below $C_{(j_0-1,j_0-1,j_2'+1)}$ in $\mathrm{JH}_{G}(\Theta_{\mu}(V_0)/V_0^+)$ (since $C_{(j_0-1,j_0-1,j_2'+1)}\in \mathrm{JH}_{G}(\Theta_{\mu}(V_0)/V_0^+)$ maps to $0$ in the quotient $V_0$ of $\Theta_{\mu}(V_0)/V_0^+$, the same holds for any constituent below $C_{(j_0-1,j_0-1,j_2'+1)}$ in $\mathrm{JH}_{G}(\Theta_{\mu}(V_0)/V_0^+)$). It follows from (\ref{equ: relation with loc alg 1}) that $C_{(j_0-1,j_0-1,j_2'+1)}\leq C_{(j_0-1,j_0-1,j)}< V_{[j-j_0+2,j],\Delta}^{\rm{alg}}$ for $j>j_2'$. Our induction hypothesis implies $V_{[j-j_0+2,j],\Delta}^{\rm{alg}}<C_{(j_0-1,j_0-1,j+1)}\leq C_{(j_0-1,j_0-1,j_2')}$ for $j<j_2'$ (and is empty when $j_2'=j_0-1$). Hence, the only locally algebraic constituent $V$ that could satisfy (\ref{equ: relation with loc alg 2}) is $V=V_{[j_2'-j_0+2,j_2'],\Delta}^{\rm{alg}}$. To prove that $V_{[j_2'-j_0+2,j_2'],\Delta}^{\rm{alg}}$ is indeed there (in $\mathrm{JH}_{G}(\Theta_{\mu}(V_0)/V_0^+)$), one has to check that the unique subquotient $V_2$ of $\Theta_{\mu}(V_0)/V_0^+$ with socle $C_{(j_0-1,j_0-1,j_2')}$ and cosocle $C_{(j_0-1,j_0-1,j_2'+1)}$ has Loewy length $3$ with middle layer either $C_{(j_0-1,j_0,j_2'+1)}$ or $C_{(j_0-1,j_0,j_2'+1)}\oplus V_{[j_2'-j_0+2,j_2'],\Delta}^{\rm{alg}}$. But by minimality of the $\mathrm{Ext}$-square $V_{(j_0-1,j_0-1,j_2'),(j_0-1,j_0-1,j_2'+1)}$ (\ref{it: hard square 2} of Proposition~\ref{prop: hard square} applied with $j_0, j_1$ there both being $j_0-1$) we see that $V_{[j_2'-j_0+2,j_2'],\Delta}^{\rm{alg}}$ must appear in the middle layer of $V_2$ and thus we indeed have (\ref{equ: relation with loc alg 2}) for $V=V_{[j_2'-j_0+2,j_2'],\Delta}^{\rm{alg}}$. Since the partial order on $[C_{(j_0-1,1,1)},C_{(j_0,1,n)}]$ is generated by the relations (\ref{equ: relation with loc alg 1}), (\ref{equ: relation with loc alg 2}) and $C_{\un{j}}<C_{\un{j}'}$ for $\un{j}<\un{j}'$, we have shown that the partial order on $\mathrm{JH}_{G}(\Theta_{\mu}(V_0)/V_0^+)$ is at least as strong as the one on $[C_{(j_0-1,1,1)},C_{(j_0,1,n)}]$ induced by $\tld{\mathbf{J}}$.\bigskip

\textbf{Step $4$}: We assume $(j_1,j_2)=(1,n)$ and prove that the equality of sets (\ref{equ: wall crossing cube quotient JH}) holds as \emph{partially ordered sets}.\\
We consider an arbitrary length $2$ subquotient $V_3$ of $\Theta_{\mu}(V_0)/V_0^+$ with socle $V$ and cosocle $V'$. We have $\mathrm{Ext}_{G}^1(V',V)\neq 0$ and $V$ is right below $V'$ for the partial order on $\mathrm{JH}_{G}(\Theta_{\mu}(V_0)/V_0^+)$. Let us prove $V<V'$ for the partial order on $\tld{\mathbf{J}}$. If $V$ and $V'$ are not both locally algebraic, then this follows from Lemma \ref{lem: Ext1 factor 1} and Lemma \ref{lem: Ext1 factor 2}. Assume that $V$ and $V'$ are both locally algebraic. Since $\mathrm{JH}_{G}(V_0)=[C_{(j_0,1,1)},C_{(j_0,1,n)}]$ (as partially ordered sets), one easily checks that $V_0$ does not contain any locally algebraic length $2$ subquotient. Hence $V$ is not a constituent of $V_0$, and we therefore have $V=V_{[j_2'-j_0+2,j_2'],\Delta}^{\rm alg}$ for some $1 \leq j_2' \leq n-1$. Then by Lemma \ref{lem: Ext1 from sm} and Lemma \ref{lem: Ext sm St} we have $V'=V_{[j_2''-j_0+1,j_2''],\Delta}^{\rm alg}$ with $d([j_2'-j_0+2,j_2'],[j_2''-j_0+1,j_2''])=1$. If $V'=V_{[j_2'-j_0+1,j_2'],\Delta}^{\rm alg}$ with $j_2'>1$, then $V$ lies right below $V'$ in $\tld{\mathbf{J}}$. Otherwise $V'=V_{[j_2'-j_0+2,j_2'+1],\Delta}^{\rm alg}$ and we have for the partial order on $\mathrm{JH}_{G}(\Theta_{\mu}(V_0)/V_0^+)$
\[V=V_{[j_2'-j_0+2,j_2'],\Delta}^{\rm alg}<C_{(j_0-1,j_0-1,j_2'+1)}<V_{[j_2'-j_0+3,j_2'+1],\Delta}^{\rm alg}<V'=V_{[j_2'-j_0+2,j_2'+1],\Delta}^{\rm alg}\]
(as these relations hold in $\tld{\mathbf{J}}$ and the partial order on $\mathrm{JH}_{G}(\Theta_{\mu}(V_0)/V_0^+)$ is at least as strong). But this contradicts the existence of $V_3$. We conclude that (\ref{equ: wall crossing cube quotient JH}) is an equality of partially ordered sets.\bigskip

\textbf{Step $5$}: We finish the proof of \ref{it: general wall crossing 3}.\\
If $(j_1,j_2)=(1,n)$, we conclude \ref{it: general wall crossing 3} from {Step $4$} together with a decreasing induction on $1<j_0\leq n-1$ (as explained above Step $1$). In particular, for each $1<j_0\leq n-1$ there exists a unique multiplicity free finite length representation $V$ in $\mathrm{Rep}^{\rm{an}}_{\rm{adm}}(G)$ such that $\mathrm{JH}_{G}(V)=[C_{(j_0-1,1,1)},C_{(j_0,1,n)}]$ as partially ordered sets. Then its unique subrepresentation $V'$ with cosocle $C_{(j_0,n-1,n-1)}$ satisfies $\mathrm{JH}_{G}(V')=[C_{(j_0-1,1,1)},C_{(j_0,n-1,n-1)}]$ as partially ordered sets. Assume $(j_1,j_2)=(n-1,n-1)$ and let $V_0$, $V_0^-$ as in \ref{it: general wall crossing 3} (both of which exist now). We have a surjection $V_0\twoheadrightarrow V_0^-$. By a symmetric argument as for $V_0^+$ in Step $1$, this induces an isomorphism $\Theta_{\mu}(V_0^-)\cong \Theta_{\mu}(V_0)$ and the adjunction map $\Theta_{\mu}(V_0^-)\rightarrow V_0^-$ is surjective. We define $q^-$ as the composition $\Theta_{\mu}(V_0)\cong \Theta_{\mu}(V_0^-)\twoheadrightarrow V_0^-$. As $\mathrm{soc}_{G}(V_0^-)=C_{(j_0-1,1,1)}$ is not a constituent of $V_0^+$, the composition $q^-\circ q^+$ (see Step $1$ for $q^+$) is necessarily zero, hence we have a complex $[V_0^+\buildrel q^+\over\longrightarrow \Theta_{\mu}(V_0)\buildrel q^-\over\longrightarrow V_0^-]$ which is exact on the left by Step $1$ and on the right by above. Then {Step $2$} (applied with $(j_1,j_2)=(n-1,n-1)$) gives exactness in the middle, and also that (\ref{equ: wall crossing cube quotient JH}) is again an equality of partially ordered sets. This finishes the proof of \ref{it: general wall crossing 3}.
\end{proof}

\begin{thm}\label{thm: unique rep}
For each $\un{j}=(j_0,j_1,j_2)\in\mathbf{J}$, $j'_0\in \{1,\dots,n-1\}$ such that $j_0\leq j_0'\leq j_0+1$ and $W\in\tld{\mathbf{J}}_{j_0'}$ such that $C_{\un{j}}\leq W$, there exists a unique multiplicity free finite length representation $V$ in $\mathrm{Rep}^{\rm{an}}_{\rm{adm}}(G)$ such that $\mathrm{JH}_{G}(V)=[C_{\un{j}},W]$ as partially ordered sets.
\end{thm}
\begin{proof}
The unicity of $V$ follows from Lemma~\ref{lem: unique rep reduction}. By \ref{it: general wall crossing 3} of Proposition~\ref{prop: general wall crossing} there exists a unique multiplicity free finite length representation $V_0$ in $\mathrm{Rep}^{\rm{an}}_{\rm{adm}}(G)$ such that $\mathrm{JH}_{G}(V_0)=[C_{(j_0,1,1)},C_{(j_0',1,n)}]$ as partially ordered sets. Since we have $C_{\un{j}}, W\in [C_{(j_0,1,1)},C_{(j_0',1,n)}]$ with $C_{\un{j}}\leq W$, $V_0$ admits a unique subquotient with socle $C_{\un{j}}$ and cosocle $W$, which gives the existence of $V$.
\end{proof}

\begin{rem}
It follows from Theorem~\ref{thm: unique rep} that there exists a unique multiplicity free finite length representation $V_0$ in $\mathrm{Rep}^{\rm{an}}_{\rm{adm}}(G)$ such that $\mathrm{JH}_{G}(V_0)=[C_{(1,1,1)},C_{(1,1,n)}]$ as partially ordered sets.
When $n=2$ the representation $V_0$ is well-known, it is (up to unramified twist) the representation denoted $\Pi^{1,\sigma}(\underline k_\sigma, \underline D)$ at the bottom of \cite[p.638]{Bre19}. 
When $n=3$ and $K=\Qp$, $V_0$ is the representation denoted $\Pi^1(\underline k, D)$ in \cite[Thm.~1.2]{Bre19} (the value of $\underline k=(k_1,k_2,k_3)$ being determined by $\mu_0$).
\end{rem}

Recall that $\mathrm{St}_n^{\rm{alg}}=V_{\emptyset,\Delta}^{\rm alg}=L(1)^\vee\otimes_E V_{\emptyset,\Delta}^{\infty}=L(\mu_0)^\vee\otimes_E\mathrm{St}_n^{\infty}$ (with $\mathrm{St}_n^{\infty}$ being the smooth Steinberg representation of $G$). Let $Z(G)$ be the center of $G$ and $\chi:Z(G)\rightarrow E^\times$ the central character of $\mathrm{St}_n^{\rm{alg}}$ (which depends on $\mu_0$). We consider the extension groups $\mathrm{Ext}_{D(G),\chi^\vee}^\bullet$ computed in the full subcategory of $\mathrm{Mod}_{D(G)}$ of $D(G)$-modules where $D(Z(G))$ acts by $\chi^\vee:D(Z(G))\rightarrow E^\times$ (see e.g.~\cite[Rem.~5.1.3(i)]{Bre19}). It follows from the analogue of the spectral sequence (\ref{equ: Ext loc alg 2}) (with (\ref{extdef})) where we fix central characters everywhere and from $\mathrm{Ext}_{U(\fg),\chi^\vee}^{\ell}(L(1), L(1))=0$ for $\ell=1,2$ (Whitehead's lemma) that we have isomorphisms for $\ell\leq 2$ and $I,I'\subseteq \Delta$
\begin{equation}\label{equ: reduce to sm Extcentral}
\mathrm{Ext}_{G, 1}^{\ell}(V_{I,\Delta}^{\infty}, V_{I',\Delta}^{\infty})^\infty\buildrel\sim\over\longrightarrow \mathrm{Ext}_{D(G),\chi^\vee}^{\ell}((V_{I',\Delta}^{\rm{alg}})^\vee, (V_{I,\Delta}^{\rm{alg}})^\vee)
\end{equation}
(where $\mathrm{Ext}_{G, 1}^{\ell}(-, -)^\infty$ means smooth extensions with trivial central character). By \cite[Thm.~2]{Or05} and (\ref{equ: reduce to sm Extcentral}) we deduce
\begin{equation}\label{equ: vanishing St central}
\mathrm{Ext}_{D(G),\chi^\vee}^1((\mathrm{St}_n^{\rm{alg}})^\vee, (\mathrm{St}_n^{\rm{alg}})^\vee)=0.
\end{equation}

\begin{lem}\label{lem: wall crossing first branch}
Let $j_0^+=\min\{2,n-1\}$, $(j_1,j_2)\in\{(n-1,n-1),(1,n)\}$ and $V_0$ (resp.~$V_0^+$) be the unique multiplicity free finite length representation in $\mathrm{Rep}^{\rm{an}}_{\rm{adm}}(G)$ such that $\mathrm{JH}_{G}(V_0)=[C_{(1,1,1)},C_{(1,j_1,j_2)}]$ (resp.~$\mathrm{JH}_{G}(V_0^+)=[C_{(1,1,1)},L(1)^\vee]$ if $n=2$ with $(j_1,j_2)=(1,1)$, and $\mathrm{JH}_{G}(V_0^+)=[C_{(1,1,1)},C_{(j_0^+,j_1,j_2)}]$ otherwise) as partially ordered sets (see Theorem~\ref{thm: unique rep}). 
Let $\mu\in \Lambda$ such that $\langle\mu+\rho,\alpha^\vee\rangle\geq 0$ for $\al\in \Phi^+$ and the stabilizer of $\mu$ in $W(G)$ for the dot action is $\{1,w_0s_{1}w_0\}$. We set $m=1$ if $(j_1,j_2)=(n-1,n-1)$ and $m=n$ if $(j_1,j_2)=(1,n)$.
\begin{enumerate}[label=(\roman*)]
\item \label{it: wall crossing first branch 1} The injection $V_0\hookrightarrow V_0^+$ induces an isomorphism $\Theta_{\mu}(V_0)\buildrel\sim\over\longrightarrow\Theta_{\mu}(V_0^+)$, and the canonical adjunction map $V_0^+\rightarrow \Theta_{\mu}(V_0^+)$ is an injection.
\item \label{it: wall crossing first branch 2} The representation $\Theta_{\mu}(V_0)/V_0^+$ admits a central character, has socle $(\mathrm{St}_n^{\rm{alg}})^{\oplus n}$ and fits into a short exact sequence
\begin{equation}\label{equ: wall crossing first branch}
0\longrightarrow (\mathrm{St}_n^{\rm{alg}})^{\oplus m} \longrightarrow \Theta_{\mu}(V_0)/V_0^+ \longrightarrow V_0\longrightarrow 0.
\end{equation}
\end{enumerate}
\end{lem}
\begin{proof}
The proof of \ref{it: wall crossing first branch 1} is similar to {Step $1$} of the proof of \ref{it: general wall crossing 3} of Proposition~\ref{prop: general wall crossing}. We prove \ref{it: wall crossing first branch 2}. As $V_0$ is multiplicity free, it admits a central character, and so do $\Theta_{\mu}(V_0)$ and $\Theta_{\mu}(V_0)/V_0^+$. By {Step $3$} and {Step $5$} of the proof of \ref{it: general wall crossing 3} of Proposition~\ref{prop: general wall crossing}, the adjunction $\Theta_{\mu}(V_0)\rightarrow V_0$ induces a surjection $q:\Theta_{\mu}(V_0)/V_0^+\twoheadrightarrow V_0$. Similar to {Step $2$} of the proof of \ref{it: general wall crossing 3} of Proposition~\ref{prop: general wall crossing}, we can check using \ref{it: wall crossing simple 1} of Lemma~\ref{lem: wall crossing of simple} and \ref{it: square crossing 1} of Lemma~\ref{lem: square as wall crossing} that we have in the Grothendieck group of finite length admissible representations of $G$
\begin{equation}\label{equ: first branch mult}
[\Theta_{\mu}(V_0)]-[V_0]-[V_0^+]=m[\mathrm{St}_n^{\rm{alg}}],
\end{equation}
which implies $[\mathrm{ker}(q)]=m[\mathrm{St}_n^{\rm{alg}}]$ by the previous sentence. Since $\Theta_{\mu}(V_0)/V_0^+$ has a central character, so does its subrepresentation $\mathrm{ker}(q)$, which together with (\ref{equ: vanishing St central}) forces $\mathrm{ker}(q)\cong (\mathrm{St}_n^{\rm{alg}})^{\oplus m}$ and gives (\ref{equ: wall crossing first branch}). It remains to show that $\Theta_{\mu}(V_0)/V_0^+$ has socle $(\mathrm{St}_n^{\rm{alg}})^{\oplus m}$. Since $V_0$ has socle $C_{(1,1,1)}$, by (\ref{equ: wall crossing first branch}) it suffices to show that $C_{(1,1,1)}$ does not appear in the socle of $\Theta_{\mu}(V_0)/V_0^+$. It follows from \ref{it: square crossing 1} of Lemma~\ref{lem: square as wall crossing} that $\Theta_{\mu}(C_{(1,1,1)})$ admits a quotient $V$ with socle $\mathrm{St}_n^{\rm{alg}}$ and cosocle $C_{(1,1,1)}$. Since $\mathrm{St}_n^{\rm{alg}}$ is not constituent of $V_0^+$, the subquotient $V$ of $\Theta_{\mu}(C_{(1,1,1)})$ and thus of $\Theta_{\mu}(V_0)$ must be a subquotient of $\Theta_{\mu}(V_0)/V_0^+$. Since $C_{(1,1,1)}$ has multiplicity $1$ in $\Theta_{\mu}(V_0)/V_0^+$ by (\ref{equ: first branch mult}), we deduce that $C_{(1,1,1)}$ cannot show up in the socle of $\Theta_{\mu}(V_0)/V_0^+$. 
\end{proof}

We now define several important finite length multiplicity free coadmissible $D(G)$-modu\-les. By Theorem~\ref{thm: unique rep} we only need to specify the corresponding partially ordered set of irreducible constituents of $\tld{\mathbf{J}}$.\bigskip

We set $X_0=Y_0=Z_0=0$, and for $1\leq k\leq n-1$, we define (with the notation (\ref{intervalVW}))
\begin{equation}\label{definitionofxyz}
\left\{\begin{array}{cclll}
X_k&\ {\rm such\ that}\ &\mathrm{JH}_{G}(X_k^\vee)&=&[C_{(n-k,1,1)},C_{(n-k,n-1,n-1)}]\\
Y_k&\ {\rm such\ that}\ &\mathrm{JH}_{G}(Y_k^\vee)&=&[C_{(n-k,1,2)},C_{(n-k,1,n)}]\\
Z_k&\ {\rm such\ that}\ &\mathrm{JH}_{G}(Z_k^\vee)&=&[C_{(n-k,1,1)},C_{(n-k,1,n)}].
\end{array}\right.
\end{equation}
We can check that $Y_1^\vee$ has no locally algebraic constituent but that for $k\geq 2$ the locally algebraic constituents of $Y_k^\vee$ are
\[V_{[2,n-k+1],\Delta}^{\rm{alg}},V_{[3,n-k+2],\Delta}^{\rm{alg}},\dots, V_{[k,n-1],\Delta}^{\rm{alg}},\]
and that for $k\geq 1$ the locally algebraic constituents of $Z_k^\vee$ are
\[V_{[1,n-k],\Delta}^{\rm{alg}},V_{[2,n-k+1],\Delta}^{\rm{alg}},V_{[3,n-k+2],\Delta}^{\rm{alg}},\dots, V_{[k,n-1],\Delta}^{\rm{alg}}.\]
We define
\begin{equation}\label{definitionofd0}
\left\{\begin{array}{cclll}
D_0&\ {\rm such\ that}\ &\mathrm{JH}_{G}(D_0^\vee)&=&[C_{(n-1,1,1)},V_{\Delta,\Delta}^{\rm{alg}}]\ \ =\ \ [C_{(n-1,1,1)},L(1)^\vee]\\
\tld{D}_0&\ {\rm such\ that}\ &\mathrm{JH}_{G}(\tld{D}_0^\vee)&=&[C_{(n-1,1,1)},C_{(n-1,1,n)}],
\end{array}\right.
\end{equation}
(note that $\tld{D}_0=Z_1$) and for $1\leq k\leq n-2$
\begin{equation}\label{definitionofdk}
\left\{\begin{array}{cclll}
D_k&\ {\rm such\ that}\ &\mathrm{JH}_{G}(D_k^\vee)&=&[C_{(n-k-1,1,1)},C_{(n-k,n-1,n-1)}]\\
\tld{D}_k&\ {\rm such\ that}\ &\mathrm{JH}_{G}(\tld{D}_k^\vee)&=&[C_{(n-k-1,1,1)},C_{(n-k,1,n)}].
\end{array}\right.
\end{equation}
From their definition (and the definition of $\tld{\mathbf{J}}$), we see that all coadmissible $D(G)$-modules $X_k$, $Y_k$, $Z_k$, $D_k$ and $\tld{D}_k$ are indecomposable multiplicity free with irreducible socle and cosocle.\bigskip

The following remark will be useful.

\begin{rem}\label{rem: unique sub quotient}
Let $V_0$ be a multiplicity free finite length representation in $\mathrm{Rep}^{\rm{an}}_{\rm{adm}}(G)$, and let $S_1,S_2\subseteq \mathrm{JH}_{G}(V_0)$ be subsets such that $S_1\cap S_2=\emptyset$. The partial order on $\mathrm{JH}_{G}(V_0)$ restricts to a partial order on $S_i$ for $i=1,2$. Assume that the following conditions hold for each $S_i$:
\begin{itemize}
\item the partially ordered set $S_i$ admits a unique minimal element $V_i'$ and a unique maximal element $V_i''$;
\item each $V\in \mathrm{JH}_{G}(V_0)$ such that $V_i'\leq V\leq V_i''$ is in $S_i$.
\end{itemize}
Then the following results are easily checked.
\begin{enumerate}[label=(\roman*)]
\item For $i=1,2$ $V_0$ admits a unique subquotient $V_i$ such that $\mathrm{JH}_{G}(V_i)=S_i$.
\item \label{it: unique subquotient 2} If $V_2'\not\leq V_1''$ in $\mathrm{JH}_{G}(V_0)$ and $\mathrm{JH}_{G}(V_0)=S_1\sqcup S_2$ as sets, then $V_0$ fits into a (possibly split) short exact sequence $0\rightarrow V_1\rightarrow V_0\rightarrow V_2\rightarrow 0$.
\item \label{it: unique subquotient 3} If $V_2'\not\leq V_1''$ and $V_1'\not\leq V_2''$ in $\mathrm{JH}_{G}(V_0)$, then $V_0$ admits a unique subquotient isomorphic to $V_1\oplus V_2$.
\end{enumerate}
\end{rem}

We now sum up the main properties of the above coadmissible $D(G)$-modules.

\begin{thm}\label{thm: main rep}
\hspace{2em}
\begin{enumerate}[label=(\roman*)]
\item \label{it: main rep 1} The finite length coadmissible $D(G)$-modules $X_k$, $Y_k$, $Z_k$, $D_k$ and $\tld{D}_k$ are multiplicity free with simple socle and cosocle and are uniquely determined (up to isomorphism) by their set of constituents endowed with the partial order of \S\ref{generalnotation}.
\item \label{it: main rep 3} For $1\leq k\leq n-1$, the coadmissible $D(G)$-module $Z_k$ admits a unique increasing $3$-stage filtration by (closed) $D(G)$-submodules with subrepresentation $Y_k$, middle graded piece $(V_{[1,n-k],\Delta}^{\rm{alg}})^\vee$ and quotient $X_k$.
\item \label{it: main rep 2} For $0\leq k\leq n-1$, the coadmissible $D(G)$-module $D_k$ admits a unique increasing $3$-stage filtration by (closed) $D(G)$-submodules with subrepresentation $X_k$, middle graded piece $(V_{[1,n-k-1],\Delta}^{\rm{alg}})^\vee$ and quotient $X_{k+1}$ (with $X_n\defeq 0$ if $k=n-1$).
\item \label{it: main rep 4} For $0\leq k\leq n-2$ we have a short exact sequence $0\rightarrow Z_k\rightarrow \tld{D}_k\rightarrow Z_{k+1}\rightarrow 0$ and a surjection $\tld{D}_k\twoheadrightarrow D_k$. More precisely, for $1\leq k\leq n-2$, the coadmissible $D(G)$-module $\tld{D}_k$ admits a unique increasing $5$-stage filtration by (closed) $D(G)$-submodules with subrepresentation $Y_k$, second graded piece $(V_{[1,n-k],\Delta}^{\rm{alg}})^\vee$, third graded piece $X_k\oplus Y_{k+1}$, fourth graded piece $(V_{[1,n-k-1],\Delta}^{\rm{alg}})^\vee$ and quotient $X_{k+1}$.
\end{enumerate}
\end{thm}
\begin{proof}
\ref{it: main rep 1} follows directly from Theorem~\ref{thm: unique rep}. \ref{it: main rep 3}, \ref{it: main rep 2} and \ref{it: main rep 4} follow from Remark~\ref{rem: unique sub quotient} and corresponding decompositions inside $\tld{\mathbf{J}}$ of the respective partially ordered sets of constituents. Let us prove \ref{it: main rep 4} and leave the other (easier) cases to the reader. For $1\leq k\leq n-2$ we have
\begin{multline*}
\mathrm{JH}_{G}(\tld{D}_k^\vee)=[C_{(n-k-1,1,1)}, C_{(n-k,1,n)}]=[C_{(n-k-1,1,1)}, C_{(n-k-1,1,n)}]\sqcup [C_{(n-k,1,1)}, C_{(n-k,1,n)}]\\
=\mathrm{JH}_{G}(Z_{k+1}^\vee)\sqcup \mathrm{JH}_{G}(Z_k^\vee),
\end{multline*}
which by \ref{it: unique subquotient 2} of Remark~\ref{rem: unique sub quotient} and the fact $C_{(n-k,1,1)}\not\leq C_{(n-k-1,1,n)}$ in $[C_{(n-k-1,1,1)}, C_{(n-k,1,n)}]$ gives the first exact sequence in \ref{it: main rep 4}. We also have
\begin{multline*}
\mathrm{JH}_{G}(\tld{D}_k^\vee)=[C_{(n-k-1,1,1)},C_{(n-k,n-1,n-1)}]\sqcup [C_{(n-k-1,1,2)}, C_{(n-k,1,n)}]\\
=\mathrm{JH}_{G}(D_k^\vee) \sqcup [C_{(n-k-1,1,2)}, C_{(n-k,1,n)}],
\end{multline*}
which \ by \ \ref{it: unique subquotient 2} \ of \ Remark~\ref{rem: unique sub quotient} \ and \ the \ fact \ that\ $C_{(n-k-1,1,2)}\not\leq C_{(n-k,n-1,n-1)} $ \ in \ $[C_{(n-k-1,1,1)}, C_{(n-k,1,n)}]$ gives an injection $D_k^\vee\hookrightarrow \tld{D}_k^\vee$, and thus a surjection $\tld{D}_k\twoheadrightarrow D_k$. Finally, by \ref{it: unique subquotient 3} of Remark~\ref{rem: unique sub quotient} and the fact $C_{(n-k-1,1,2)}\not\leq C_{(n-k,n-1,n-1)}$ and $C_{(n-k,1,1)}\not\leq C_{(n-k-1,1,n)}$ (with $\mathrm{JH}_{G}(X_k^\vee)=[C_{(n-k,1,1)},C_{(n-k,n-1,n-1)}]$ and $\mathrm{JH}_{G}(Y_{k+1}^\vee)=[C_{(n-k-1,1,2)},C_{(n-k-1,1,n)}]$), we deduce that $X_k\oplus Y_{k+1}$ is a subquotient of $\tld{D}_k$.
\end{proof}

\begin{thm}\label{thm: main Ext}
\hspace{2em}
\begin{enumerate}[label=(\roman*)]
\item \label{it: main Ext 0} The $E$-vector space $\mathrm{Ext}_{D(G)}^1((\mathrm{St}_n^{\rm{alg}})^\vee, X_{n-1})$ has dimension $1$.
\item \label{it: main Ext 1} The $E$-vector spaces $\mathrm{Ext}_{D(G)}^1((\mathrm{St}_n^{\rm{alg}})^\vee,Z_{n-1})$ and $\mathrm{Ext}_{D(G)}^1(Z_{n-1}, (\mathrm{St}_n^{\rm{alg}})^\vee)$ have dimension $n$.
\item \label{it: main Ext 3} For $1\leq k\leq n-2$ we have $\mathrm{Ext}_{D(G)}^1(D_k, (V_{[1,n-k],\Delta}^{\rm{alg}})^\vee)=0$.
\item \label{it: main Ext 2} For $1\leq k\leq n-1$, the $E$-vector space $\mathrm{Ext}_{D(G)}^1(X_k, (V_{[1,n-k],\Delta}^{\rm{alg}})^\vee)$ has dimension $1$ and the $E$-vector space $\mathrm{Ext}_{D(G)}^1((V_{[1,n-k-1],\Delta}^{\rm{alg}})^\vee, \tld{X}_k)$ has dimension $2$ where $\tld{X}_k$ is the unique non-split extension of $X_k$ by $(V_{[1,n-k],\Delta}^{\rm{alg}})^\vee$.
\end{enumerate}
\end{thm}
\begin{proof}
We prove \ref{it: main Ext 0}. By \ref{it: Ext1 with alg 2} of Lemma~\ref{lem: Ext1 factor 2} applied with $I=\emptyset$ and by \ref{it: Ext2 with St 2} of Lemma~\ref{lem: Ext2 with St} we have $\Dim_E \mathrm{Ext}_{G}^1(C_{(1,1,1)},\mathrm{St}_n^{\rm{alg}})=1$ and $\mathrm{Ext}_{G}^i(V,\mathrm{St}_n^{\rm{alg}})=0$ for $i=\{1,2\}$ and $V=C_{(1,j_1,j_1)}\in \mathrm{JH}_{G}(X_{n-1}^\vee)=[C_{(1,1,1)},C_{(1,n-1,n-1)}]$ with $2\leq j_1\leq n-1$ (note that $I^-_{j_1,j_1}=\emptyset$ if and only if $j_1=1$ by Lemma~\ref{lem: explicit smooth induction} and (\ref{i+-chiant})). The statement then follows by d\'evissage.

We prove \ref{it: main Ext 1}. Recall the functors $\mathrm{Ext}_{D(G),\chi^\vee}^\bullet$ from the discussion above Lemma~\ref{lem: wall crossing first branch}. Note \ first \ that, \ since \ $\mathrm{St}_n^{\rm{alg}}$ \ is \ not \ a \ constituent \ of \ $Z_{n-1}^\vee$, \ we \ have \ isomorphisms $\mathrm{Ext}_{D(G),\chi}^1((\mathrm{St}_n^{\rm{alg}})^\vee,Z_{n-1})\buildrel\sim\over\rightarrow \mathrm{Ext}_{D(G)}^1((\mathrm{St}_n^{\rm{alg}})^\vee,Z_{n-1})$ and $\mathrm{Ext}_{D(G),\chi^\vee}^1(Z_{n-1}, (\mathrm{St}_n^{\rm{alg}})^\vee)\buildrel\sim\over\rightarrow \mathrm{Ext}_{D(G)}^1(Z_{n-1}, (\mathrm{St}_n^{\rm{alg}})^\vee)$. By \cite[Thm.~2]{Or05} and (\ref{equ: reduce to sm Extcentral}) we deduce for $j\in \{1,\dots,n-1\}$
\begin{equation}\label{equ: exticentral}
\Dim_E \mathrm{Ext}_{D(G),\chi^\vee}^1((\mathrm{St}_n^{\rm{alg}})^\vee,(V_{\{j\}}^{\rm{alg}})^\vee)=1,\ \mathrm{Ext}_{D(G),\chi^\vee}^2((\mathrm{St}_n^{\rm{alg}})^\vee, (V_{\{j\}}^{\rm{alg}})^\vee)=0.
\end{equation}
Moreover, by \ref{it: Ext1 with alg 2} of Lemma~\ref{lem: Ext1 factor 2} applied with $I=\emptyset$ (where we recall that $\mathrm{Ext}_{D(G),\chi^\vee}^1=\mathrm{Ext}_{D(G)}^1$ there) and by the analogue of \ref{it: Ext2 with St 2} of Lemma~\ref{lem: Ext2 with St} with $\mathrm{Ext}_{D(G),\chi^\vee}^2$ instead of $\mathrm{Ext}_{D(G)}^2$ (the proof of which is the same), we deduce for $\ell\in \{1,2\}$ and $V=C_{(1,j_1,j_2)}\in \mathrm{JH}_{G}(Z_{n-1}^\vee)\setminus\{C_{(1,1,1)}\}=[C_{(1,1,1)},C_{(1,1,n)}]\setminus\{C_{(1,1,1)}\}$
\begin{equation}\label{equ: exticentralbis}
\mathrm{Ext}_{D(G),\chi^\vee}^\ell((\mathrm{St}_n^{\rm{alg}})^\vee, V^\vee)=0.
\end{equation}
Then combining $\Dim_E\mathrm{Ext}_{D(G),\chi^\vee}^1((\mathrm{St}_n^{\rm{alg}})^\vee, C_{(1,1,1)}^\vee)=1$ in \ref{it: main Ext 0} with (\ref{equ: exticentral}), (\ref{equ: exticentralbis}) and a d\'evissa\-ge, we obtain the first statement in \ref{it: main Ext 1}.
Similarly, by \cite[Thm.~2]{Or05} and (\ref{equ: reduce to sm Extcentral}) we deduce for $j\in \{1,\dots,n-1\}$
\begin{equation}\label{equ: exticentra2}
\Dim_E \mathrm{Ext}_{D(G),\chi^\vee}^1((V_{\{j\}}^{\rm{alg}})^\vee, (\mathrm{St}_n^{\rm{alg}})^\vee)=1,\ \mathrm{Ext}_{D(G),\chi^\vee}^2((V_{\{j\}}^{\rm{alg}})^\vee, (\mathrm{St}_n^{\rm{alg}})^\vee)=0.
\end{equation}
Recall from Lemma~\ref{lem: explicit smooth induction} and (\ref{i+-chiant}) (see also the proof of Lemma~\ref{lem: Ext2 with St}) that $I_{j_1,j_2}^+=\emptyset$ if and only if $(j_1,j_2)=(1,n)$. Then by \ref{it: Ext1 with alg 1} of Lemma~\ref{lem: Ext1 factor 2} applied with $I=\emptyset$ (noting that $\mathrm{Ext}_{D(G),\chi^\vee}^1=\mathrm{Ext}_{D(G)}^1$ there) and by the analogue of \ref{it: Ext2 with St 1} of Lemma~\ref{lem: Ext2 with St} with $\mathrm{Ext}_{D(G),\chi^\vee}^2$ instead of $\mathrm{Ext}_{D(G)}^2$ (the proof of which is the same), we deduce for $\ell\in \{1,2\}$ and $V=C_{(1,j_1,j_2)}\in \mathrm{JH}_{G}(Z_{n-1}^\vee)\setminus\{C_{(1,1,n)}\}=[C_{(1,1,1)},C_{(1,1,n)}]\setminus\{C_{(1,1,n)}\}$
\begin{equation}\label{equ: exticentralbisbis}
\mathrm{Ext}_{D(G),\chi^\vee}^\ell(V^\vee, (\mathrm{St}_n^{\rm{alg}})^\vee)=0.
\end{equation}
Then combining $\Dim_E\mathrm{Ext}_{D(G),\chi^\vee}^1(C_{(1,1,n)}^\vee, (\mathrm{St}_n^{\rm{alg}})^\vee)=1$ from \ref{it: Ext1 with alg 1} of Lemma~\ref{lem: Ext1 factor 2} with (\ref{equ: exticentra2}), (\ref{equ: exticentralbisbis}) and a d\'evissa\-ge, we obtain the second statement in \ref{it: main Ext 1}.

We prove \ref{it: main Ext 3}. By construction $D_k$ admits a unique subquotient $D$ with $\mathrm{JH}_{G}((D)^\vee)=[C_{(n-k-1,n-k-1,n-k-1)},C_{(n-k,n-k,n-k)}]$ as partially ordered sets. More precisely $D^\vee$ is the $\mathrm{Ext}$-square constructed in \ref{it: easy square 2} of Proposition~\ref{prop: easy square} with socle $C_{(n-k-1,n-k-1,n-k-1)}$, cosocle $C_{(n-k,n-k,n-k)}$, and middle layer $C_{(n-k-1,n-k,n-k)}\oplus C_{(n-k,n-k-1,n-k-1)}\oplus V_{[1,n-k-1],\Delta}^{\rm{alg}}$. By \ref{it: Ext1 with alg 1} of Lemma~\ref{lem: Ext1 factor 2} and Lemma~\ref{lem: Ext1 from sm} together with Lemma~\ref{lem: Ext sm St} we know that, for $V\in\mathrm{JH}_{G}(D_k^\vee)$, we have $\mathrm{Ext}_{D(G)}^1(V^\vee,(V_{[1,n-k],\Delta}^{\rm{alg}})^\vee)\neq 0$ if and only if $V\in\{C_{n-k,n-k,n-k},V_{[1,n-k-1],\Delta}^{\rm{alg}}\}$. So by d\'evissage it suffices to show
\begin{equation}\label{equ: main Ext 3}
\mathrm{Ext}_{G}^1(V_{[1,n-k],\Delta}^{\rm{alg}},D^\vee)=0.
\end{equation}
Assume (\ref{equ: main Ext 3}) does not hold. Then there exists a non-split extension $0\rightarrow D^\vee\rightarrow W\rightarrow V_{[1,n-k],\Delta}^{\rm{alg}}\rightarrow 0$ in $\mathrm{Rep}^{\rm{an}}_{\rm{adm}}(G)$. As $D^\vee$ has socle $V_1\defeq C_{(n-k-1,n-k-1,n-k-1)}$, so does $W$. Let $W^-$ be the unique length $2$ subrepresentation of $D^\vee\subseteq W$ with socle $V_1$ and cosocle $V_0\defeq C_{(n-k-1,n-k,n-k)}$. By \ref{it: wall crossing simple 1} of Lemma~\ref{lem: wall crossing of simple} we have $\Theta_{s_k}(V)=0$ for each constituent $V$ of $W/W^-$, and thus the injection $W^-\hookrightarrow W$ induces an isomorphism $\Theta_{s_k}(W^-)\buildrel\sim\over\longrightarrow\Theta_{s_k}(W)$. As the canonical adjunction map $W\rightarrow \Theta_{s_k}(W)$ restricts to an injection $V_1\hookrightarrow \Theta_{s_k}(V_1)$ (as $V_1$ is irreducible and the adjunction map is non-zero), and as $W$ has socle $V_1$, we obtain an injection $W\hookrightarrow \Theta_{s_k}(W)\cong \Theta_{s_k}(W^-)$. But from \ref{it: square crossing 1} of Lemma~\ref{lem: square as wall crossing} we deduce $V_{[1,n-k],\Delta}^{\rm{alg}}\notin\mathrm{JH}_{G}(\Theta_{s_k}(W^-))=\mathrm{JH}_{G}(\Theta_{s_k}(V_1))\cup \mathrm{JH}_{G}(\Theta_{s_k}(V_0))$, which is a contradiction.

We prove \ref{it: main Ext 2}. Let $1\leq k\leq n-1$ and recall that $\mathrm{JH}_{G}(X_k^\vee)=[C_{(n-k,1,1)},C_{(n-k,n-1,n-1)}]=\{C_{(n-k,j,j)}\mid 1\leq j\leq n-1\}$. We have $d([1,n-k-1],[1,n-k])=1$ (with $[1,n-k-1]=\emptyset$ when $k=n-1$), which by (\ref{equ: reduce to sm Extcentral}) and \cite[Thm.~1]{Or05} gives 
\begin{multline}\label{cancellation0}
\mathrm{Ext}_{D(G),\chi^\vee}^2((V_{[1,n-k-1],\Delta}^{\rm{alg}})^\vee,(V_{[1,n-k],\Delta}^{\rm{alg}})^\vee)=0 \ \ \mathrm{ and } \ \\\Dim_E\mathrm{Ext}_{D(G),\chi^\vee}^1((V_{[1,n-k-1],\Delta}^{\rm{alg}})^\vee,(V_{[1,n-k],\Delta}^{\rm{alg}})^\vee)=1.
\end{multline} 
By \ref{it: Ext1 with alg 1} of Lemma~\ref{lem: Ext1 factor 2} (with (\ref{i+-chiant})), we have
\begin{equation}\label{cancellation1}
\mathrm{Ext}_{D(G),\chi^\vee}^1(C_{(n-k,j,j)}^\vee, (V_{[1,n-k],\Delta}^{\rm{alg}})^\vee)\neq 0\ {\rm if\ and\ only\ if}\ j=n-k
\end{equation}
and $\Dim_E \mathrm{Ext}_{D(G),\chi^\vee}^1(C_{(n-k,n-k,n-k)}
^\vee, (V_{[1,n-k],\Delta}^{\rm{alg}})^\vee)=1$. 
By \ref{it: Ext1 with alg 2} of Lemma~\ref{lem: Ext1 factor 2}, we have
\begin{equation}\label{cancellation2}
\mathrm{Ext}_{D(G),\chi^\vee}^1((V_{[1,n-k-1],\Delta}^{\rm{alg}})^\vee,C_{(n-k,j,j)}^\vee)\neq 0\ {\rm if\ and\ only\ if}\ j=n-k
\end{equation}
and \ $\Dim_E \mathrm{Ext}_{D(G),\chi^\vee}^1((V_{[1,n-k-1],\Delta}^{\rm{alg}})^\vee,C_{(n-k,n-k,n-k)}
^\vee)=1$.\ By \ the \ analogue \ of \ \ref{it: Ext2 with alg 1} of Lemma~\ref{lem: Ext2 factor 2} (see also the second paragraph in the proof of Lemma~\ref{lem: vanishing outside cube}) with $\mathrm{Ext}_{D(G),\chi^\vee}^2$ instead of $\mathrm{Ext}_{D(G)}^2$ (the proof of which is the same), we have for $j<n-k$
\begin{equation}\label{cancellation3}
\mathrm{Ext}_{D(G),\chi^\vee}^2(C_{(n-k,j,j)}^\vee, (V_{[1,n-k],\Delta}^{\rm{alg}})^\vee)=0.
\end{equation}
By the analogue of \ref{it: Ext2 with alg 2} of Lemma~\ref{lem: Ext2 factor 2} (see also the third paragraph in the proof of Lemma~\ref{lem: vanishing outside cube}) with $\mathrm{Ext}_{D(G),\chi^\vee}^2$ instead of $\mathrm{Ext}_{D(G)}^2$, we have for $j>n-k$
\begin{equation}\label{cancellation4}
\mathrm{Ext}_{D(G),\chi^\vee}^2((V_{[1,n-k-1],\Delta}^{\rm{alg}})^\vee,C_{(n-k,j,j)}^\vee)=0.
\end{equation}
Let \ $V_0$ \ (resp.~$V_1$) \ be \ the \ unique \ quotient \ of \ $X_k^\vee$ \ with \ socle \ $C_{(n-k,n-k,n-k)}$ \ (resp. \ $C_{(n-k,n-k+1,n-k+1)}$) with $V_1=0$ when $k=1$. By a d\'evissage using (\ref{cancellation1}) and (\ref{cancellation3}), we deduce $\mathrm{Ext}_{D(G),\chi^\vee}^1(V_1^\vee, (V_{[1,n-k],\Delta}^{\rm{alg}})^\vee)=0$ and $\mathrm{Ext}_{D(G),\chi^\vee}^{\ell}(X_k/V_0^\vee, (V_{[1,n-k],\Delta}^{\rm{alg}})^\vee)=0$ for $\ell=1,2$. Hence, the surjection $X_k\twoheadrightarrow X_k/V_1^\vee$ and the injection $C_{(n-k,n-k,n-k)}^\vee\hookrightarrow X_k/V_1^\vee$ (with cokernel $X_k/V_0^\vee$) induce isomorphisms between $1$-dimensional $E$-vector spaces
\begin{multline*}
\mathrm{Ext}_{D(G),\chi^\vee}^1(X_k, (V_{[1,n-k],\Delta}^{\rm{alg}})^\vee) \buildrel\sim\over\longleftarrow\mathrm{Ext}_{D(G),\chi^\vee}^1(X_k/V_1^\vee, (V_{[1,n-k],\Delta}^{\rm{alg}})^\vee)\\
\buildrel\sim\over\longrightarrow \mathrm{Ext}_{D(G),\chi^\vee}^1(C_{(n-k,n-k,n-k)}^\vee, (V_{[1,n-k],\Delta}^{\rm{alg}})^\vee).
\end{multline*}
In particular there exists a unique (up to isomorphism) $D(G)$-module that fits into a non-split extension \ $0\rightarrow (V_{[1,n-k],\Delta}^{\rm{alg}})^\vee\rightarrow \tld{X}_k\rightarrow X_k\rightarrow 0$. \ A \ symmetric \ argument \ gives $\Dim_E \mathrm{Ext}_{D(G),\chi^\vee}^1((V_{[1,n-k-1],\Delta}^{\rm{alg}})^\vee, X_k/V_1^\vee)=1$. By (\ref{cancellation0}) we have a short exact sequence
\begin{multline*}
0\rightarrow \mathrm{Ext}_{D(G),\chi^\vee}^1((V_{[1,n-k-1],\Delta}^{\rm{alg}})^\vee, (V_{[1,n-k],\Delta}^{\rm{alg}})^\vee)\rightarrow \mathrm{Ext}_{D(G),\chi^\vee}^1((V_{[1,n-k-1],\Delta}^{\rm{alg}})^\vee, \tld{X}_k/V_1^\vee)\\
\rightarrow \mathrm{Ext}_{D(G),\chi^\vee}^1((V_{[1,n-k-1],\Delta}^{\rm{alg}})^\vee, X_k/V_1^\vee)\rightarrow 0
\end{multline*}
where the left hand side and the right hand side both have dimension $1$. By another d\'evissage using (\ref{cancellation2}) and (\ref{cancellation4}) we deduce $\mathrm{Ext}_{D(G),\chi^\vee}^1((V_{[1,n-k-1],\Delta}^{\rm{alg}})^\vee, X_k/V_0^\vee)=0$ and $\mathrm{Ext}_{D(G),\chi^\vee}^{\ell}((V_{[1,n-k-1],\Delta}^{\rm{alg}})^\vee,V_1^\vee)=0$ for $\ell=1,2$, which gives an isomorphism
\[\mathrm{Ext}_{D(G),\chi^\vee}^1((V_{[1,n-k-1],\Delta}^{\rm{alg}})^\vee, \tld{X}_k) \buildrel\sim\over\longrightarrow \mathrm{Ext}_{D(G),\chi^\vee}^1((V_{[1,n-k-1],\Delta}^{\rm{alg}})^\vee, \tld{X}_k/V_1^\vee).\]
In particular, we conclude that $\Dim_E \mathrm{Ext}_{D(G),\chi^\vee}^1((V_{[1,n-k-1],\Delta}^{\rm{alg}})^\vee, \tld{X}_k)=2$. Finally, since all irreducible constituents in the various $\mathrm{Ext}_{D(G),\chi^\vee}^1$ are actually distinct, we can replace $\mathrm{Ext}_{D(G),\chi^\vee}^1$ by $\mathrm{Ext}_{D(G)}^1$.
\end{proof}

\begin{rem}
Note that, when $n=2$, the first statement in \ref{it: main Ext 1} of Theorem \ref{thm: main Ext} proves \cite[Conj.~3.2.1]{Bre19}.
\end{rem}

For $1\leq j\leq n-1$, define $Z_{n-1,\geq j}$ (resp.~$Z_{n-1,\leq j}$) as the $D(G)$-module which is the unique (closed) subspace (resp.~(topological) quotient) of $Z_{n-1}$ such that $\mathrm{JH}_{G}(Z_{n-1,\geq j}^\vee)=[V_{\{j\},\Delta}^{\rm{alg}},C_{(1,1,n)}]$ (resp.~$\mathrm{JH}_{G}(Z_{n-1,\leq j}^\vee)=[C_{(1,1,1)},V_{\{j\},\Delta}^{\rm{alg}}]$). The same argument as in the proof of \ref{it: main Ext 1} of Theorem~\ref{thm: main Ext} shows that the injections $Z_{n-1,\geq n-1}\hookrightarrow \cdots\hookrightarrow Z_{n-1,\geq 1}\hookrightarrow Z_{n-1}$ induce a decreasing filtration of subspaces of $\mathrm{Ext}_{D(G)}^1((\mathrm{St}_n^{\rm{alg}})^\vee,Z_{n-1})$
\begin{equation}\label{equ: branch filtration 1}
\left\{\mathrm{Ext}_{D(G)}^1((\mathrm{St}_n^{\rm{alg}})^\vee,Z_{n-1,\geq j})\right\}_{1\leq j\leq n-1}
\end{equation}
where $\Dim_E \mathrm{Ext}_{D(G)}^1((\mathrm{St}_n^{\rm{alg}})^\vee,Z_{n-1,\geq j})=n-j$. Similarly, the surjections $Z_{n-1}\!\twoheadrightarrow \!Z_{n-1,\leq n-1}\!\twoheadrightarrow \cdots \twoheadrightarrow Z_{n-1,\leq 1}$ induce an increasing filtration of subspaces of $\mathrm{Ext}_{D(G)}^1(Z_{n-1}, (\mathrm{St}_n^{\rm{alg}})^\vee)$
\begin{equation*}
\left\{\mathrm{Ext}_{D(G)}^1(Z_{n-1,\leq j}, (\mathrm{St}_n^{\rm{alg}})^\vee)\right\}_{1\leq j\leq n-1}
\end{equation*}
where $\Dim_E \mathrm{Ext}_{D(G)}^1(Z_{n-1,\leq j}, (\mathrm{St}_n^{\rm{alg}})^\vee)=j$.\bigskip

By \ref{it: main rep 4} of Theorem~\ref{thm: main rep}, for $0\leq k\leq n-2$, $\tld{D}_k$ has the form
\begin{equation}\label{equ: rough shape of rep}
\tld{D}_k\ \cong \ \begin{xy}
(-3,0)*+{Y_{k}}="a"; (20,0)*+{(V_{[1,n-k],\Delta}^{\rm{alg}})^\vee}="b"; (40,-10)*+{X_{k}}="c"; (60,0)*+{(V_{[1,n-k-1],\Delta}^{\rm{alg}})^\vee}="d"; (40,10)*+{Y_{k+1}}="e"; (87,0)*+{X_{k+1}}="f";
{\ar@{-}"a";"b"}; {\ar@{-}"b";"c"}; {\ar@{-}"c";"d"}; {\ar@{-}"b";"e"}; {\ar@{-}"e";"d"}; {\ar@{-}"d";"f"};
\end{xy}
\end{equation}
where we write subrepresentations on the left, quotients on the right, where lines represent non-split extensions, and where $Y_k=(V_{[1,n-k],\Delta}^{\rm{alg}})^\vee=X_k=0$ when $k=0$. Let $D_{n-1}$ be the unique (up to isomorphism) non-split extension in \ref{it: main Ext 0} of Theorem~\ref{thm: main Ext}:
\[0\longrightarrow X_{n-1}\longrightarrow D_{n-1} \longrightarrow (\mathrm{St}_n^{\rm{alg}})^\vee \longrightarrow 0.\]
Define $\tld{D}_{n-1}$ as \emph{any} non-split extension \emph{with cosocle} $(\mathrm{St}_n^{\rm{alg}})^\vee$:
\[0\longrightarrow Z_{n-1}\longrightarrow \tld{D}_{n-1} \longrightarrow (\mathrm{St}_n^{\rm{alg}})^\vee \longrightarrow 0.\]
Note that, by \ref{it: main Ext 1} of Theorem~\ref{thm: main Ext} and the discussion around (\ref{equ: branch filtration 1}), the isomorphism class of $\tld{D}_{n-1}$ depends on $n-1$ ``parameters''. More precisely, the set of isomorphism classes of coadmissible $D(G)$-modules $\tld{D}_{n-1}$ is in natural bijection with the set
\[\big(\mathrm{Ext}_{D(G)}^1((\mathrm{St}_n^{\rm{alg}})^\vee,Z_{n-1})\setminus \mathrm{Ext}_{D(G)}^1((\mathrm{St}_n^{\rm{alg}})^\vee,Z_{n-1,\geq 1})\big)/E^\times\]
which is in non-canonical bijection with ${\mathbb A}^{n-1}(E)$. By \ref{it: main rep 3} of Theorem~\ref{thm: main rep} (for $k=n-1$) $\tld{D}_{n-1}$ has the form
\begin{equation}\label{equ: dn-1tilde}
\tld{D}_{n-1}\ \cong\ \begin{xy}
(0,0)*+{Y_{n-1}}="a"; (20,0)*+{(V_{\{1\},\Delta}^{\rm{alg}})^\vee}="b"; (41,0)*+{X_{n-1}}="c"; (61,0)*+{(\mathrm{St}_n^{\rm{alg}})^\vee.}="d";
{\ar@{-}"a";"b"}; {\ar@{-}"b";"c"}; {\ar@{-}"c";"d"};
\end{xy}
\end{equation}

It follows from \ref{it: main rep 1}, \ref{it: main rep 2} of Theorem~\ref{thm: main rep} and the above definition of $D_{n-1}$ that, for $0\leq k\leq n-2$, there exists a unique (up to scalar) non-zero map $d_{\mathbf{D}}^k: D_k \rightarrow D_{k+1}$ whose image is $X_{k+1}$. In particular we can consider the complex of finite length coadmissible $D(G)$-modules (with $D_k$ in degree $k$)
\begin{equation}\label{equ: explicit dR}
\mathbf{D}^\bullet\ \defeq \ [D_0\buildrel {d_{\mathbf{D}}^0} \over \longrightarrow \cdots\buildrel {d_{\mathbf{D}}^{k-1}} \over\longrightarrow D_k \buildrel {d_{\mathbf{D}}^k} \over\longrightarrow D_{k+1}\buildrel {d_{\mathbf{D}}^{k+1}} \over\longrightarrow \cdots \buildrel {d_{\mathbf{D}}^{n-2}} \over\longrightarrow D_{n-1}].
\end{equation}
For $0\leq k\leq n-1$ we define $H^k(\mathbf{D}^\bullet)\defeq \mathrm{ker}(d_{\mathbf{D}}^k)/\mathrm{im}(d_{\mathbf{D}}^{k-1})$ (with the convention $d_{\mathbf{D}}^{n-1}=0$). By \ref{it: main rep 2} of Theorem~\ref{thm: main rep} and the above definition of $D_{n-1}$ we have for $0\leq k\leq n-1$
\begin{equation*}
H^k(\mathbf{D}^\bullet)\cong (V_{[1,n-k-1],\Delta}^{\rm{alg}})^\vee.
\end{equation*}

Similarly it follows from \ref{it: main rep 1}, \ref{it: main rep 4} of Theorem~\ref{thm: main rep} and the above definition of $\tld{D}_{n-1}$ that, for $0\leq k\leq n-2$, there exists a unique (up to scalar) non-zero map $d_{\tld{\mathbf{D}}}^k:\tld{D}_k \rightarrow \tld{D}_{k+1}$ whose image is $Z_{k+1}$, and we obtain a complex of finite length coadmissible $D(G)$-modules (with $\tld{D}_k$ in degree $k$)
\begin{equation}\label{equ: explicit dR split top}
\tld{\mathbf{D}}^\bullet\ \defeq \ [\tld{D}_0\buildrel {d_{\tld{\mathbf{D}}}^0} \over \longrightarrow \cdots\buildrel {d_{\tld{\mathbf{D}}}^{k-1}} \over \longrightarrow \tld{D}_k \buildrel {d_{\tld{\mathbf{D}}}^{k}} \over \longrightarrow \tld{D}_{k+1}\buildrel {d_{\tld{\mathbf{D}}}^{k+1}} \over \longrightarrow \cdots \buildrel {d_{\tld{\mathbf{D}}}^{n-2}} \over \longrightarrow \tld{D}_{n-1}].
\end{equation}
Recall that $\tld{\mathbf{D}}^\bullet$ is not unique because $\tld{{D}}_{n-1}$ depends non-canonically on some element in ${\mathbb A}^{n-1}(E)$. For $0\leq k\leq n-1$ we define $H^k(\tld{\mathbf{D}}^\bullet)\defeq \mathrm{ker}(d_{\tld{\mathbf{D}}}^k)/\mathrm{im}(d_{\tld{\mathbf{D}}}^{k-1})$ (with the convention $d_{\tld{\mathbf{D}}}^{n-1}=0$). By \ref{it: main rep 4} of Theorem~\ref{thm: main rep} and the above definition of $\tld{D}_{n-1}$ we have $H^k(\tld{\mathbf{D}}^\bullet)=0$ for $0\leq k\leq n-2$ and $H^{n-1}(\tld{\mathbf{D}}^\bullet)\cong (\mathrm{St}_n^{\rm{alg}})^\vee$. In particular the canonical morphism of complexes $\tld{\mathbf{D}}^\bullet\twoheadrightarrow H^{n-1}(\tld{\mathbf{D}}^\bullet)[-(n-1)]\cong (\mathrm{St}_n^{\rm{alg}})^\vee[-(n-1)]$ is a quasi-isomorphism.\bigskip

For $0\leq k\leq n-2$, by \ref{it: main rep 2}, \ref{it: main rep 4} of Theorem~\ref{thm: main rep} we have a surjection $\tld{D}_k\twoheadrightarrow D_k$ which is unique up to scalar (as $\tld{D}_k$ is multiplicity free with simple cosocle by \ref{it: main rep 1} of Theorem~\ref{thm: main rep}). It follows from (\ref{equ: dn-1tilde}) and the definition of $D_{n-1}$ that there is a unique (up to scalar) surjection $\tld{D}_{n-1}\twoheadrightarrow D_{n-1}$. Consequently, we see from the definition of the complexes (\ref{equ: explicit dR}) and (\ref{equ: explicit dR split top}) that there is a natural morphism of complexes of $D(G)$-modules
\begin{equation*}
\tld{\mathbf{D}}^\bullet\longrightarrow \mathbf{D}^\bullet
\end{equation*}
which is an isomorphism on $H^{n-1}$. We thus have proven the following theorem.

\begin{thm}\label{thm: main split}
The canonical morphism of complexes
\[\mathbf{D}^\bullet\twoheadrightarrow H^{n-1}(\mathbf{D}^\bullet)[-(n-1)]\cong (\mathrm{St}_n^{\rm{alg}})^\vee[-(n-1)]\]
admits an explicit section in the derived category of finite length coadmissible $D(G)$-modules with Orlik-Strauch constituents (Theorem \ref{prop: OS property}) given by
\begin{equation*}
(\mathrm{St}_n^{\rm{alg}})^\vee[-(n-1)]\longleftarrow\tld{\mathbf{D}}^\bullet\longrightarrow \mathbf{D}^\bullet.
\end{equation*}
\end{thm}

Recall that $\tld{\mathbf{D}}^\bullet$ (and hence the corresponding section) depends on a parameter in ${\mathbb A}^{n-1}(E)$.\bigskip

At this point, one can ask the following question. For $0\leq \ell\leq n-1$, define the usual truncated subcomplex $\tau_{\leq\ell}\mathbf{D}^\bullet$ of $\mathbf{D}^\bullet$ (with $\tau_{\leq n-1}\mathbf{D}^\bullet=\mathbf{D}^\bullet$)
\begin{equation*}
\tau_{\leq \ell}\mathbf{D}^\bullet\ \defeq \ [D_0\longrightarrow \cdots\longrightarrow \cdots \longrightarrow D_{\ell-1}\longrightarrow \mathrm{ker}(d_{\mathbf{D}}^{\ell})].
\end{equation*}
Then we again have a canonical morphism of complexes
\begin{equation}\label{forsection}
\tau_{\leq \ell}\mathbf{D}^\bullet \twoheadrightarrow H^{\ell}(\mathbf{D}^\bullet)[-\ell]\cong (V_{[1,n-\ell-1],\Delta}^{\rm{alg}})^\vee[-\ell].
\end{equation}
In view of Theorem \ref{thm: main split}, it is natural to ask if, for $0\leq \ell\leq n-2$ there also exists a section to this morphism in the derived category of finite length coadmissible $D(G)$-modules with Orlik-Strauch constituents. This is obvious when $\ell=0$ (since $\tau_{\leq 0}\mathbf{D}^\bullet\cong H^0(\mathbf{D}^\bullet)[0]$), and not too complicated when $\ell=1$:

\begin{prop}\label{prop: little split}
The canonical morphism of complexes
\[\tau_{\leq 1}\mathbf{D}^\bullet \twoheadrightarrow H^{1}(\mathbf{D}^\bullet)[-1]\cong (V_{[1,n-2],\Delta}^{\rm{alg}})^\vee[-1]\]
admits an explicit section in the derived category of finite length coadmissible $D(G)$-modules with Orlik-Strauch constituents given by
\begin{equation*}
(V_{[1,n-2],\Delta}^{\rm{alg}})^\vee[-1]\longleftarrow[{D}_0\rightarrow \!\begin{xy} (0,0)*+{{D}_0}="a"; (18,0)*+{(V_{[1,n-2],\Delta}^{\rm{alg}})^\vee}="b"; {\ar@{-}"a";"b"}\end{xy}\!]\longrightarrow \tau_{\leq 1}\mathbf{D}^\bullet
\end{equation*}
where $\begin{xy} (0,0)*+{{D}_0}="a"; (19,0)*+{(V_{[1,n-2],\Delta}^{\rm{alg}})^\vee}="b"; {\ar@{-}"a";"b"}\end{xy}\cong \begin{xy} (0,0)*+{(V_{[1,n-1],\Delta}^{\rm{alg}})^\vee}="a"; (19,0)*+{X_1}="b"; (38,0)*+{(V_{[1,n-2],\Delta}^{\rm{alg}})^\vee}="c"; {\ar@{-}"a";"b"}; {\ar@{-}"b";"c"}\end{xy}\!$ is one of the representations in \ref{it: main Ext 2} of Theorem \ref{thm: main Ext} (applied with $k=1$) depending on a parameter in ${\mathbb A}^{1}(E)$, and where the morphisms of complexes are the obvious ones.
\end{prop}

One can prove that the parameter in ${\mathbb A}^{1}(E)$ corresponds to a choice of a $p$-adic logarithm $\log : K^\times \rightarrow E$. Note that, when $n=2$, $\begin{xy} (0,0)*+{{D}_0}="a"; (19,0)*+{(V_{[1,n-2],\Delta}^{\rm{alg}})^\vee}="b"; {\ar@{-}"a";"b"}\end{xy}\!\cong \widetilde D_1/Y_1$. Theorem \ref{thm: main split} and Proposition \ref{prop: little split} have the following consequence.

\begin{cor}\label{cor: split3}
For $n=3$ there exists an explicit splitting in the derived category of finite length coadmissible $D(G)$-modules with Orlik-Strauch constituents $\mathbf{D}^\bullet \cong \oplus_{\ell=0}^2 H^{\ell}(\mathbf{D}^\bullet)[-\ell]$.
\end{cor}

When $2\leq \ell\leq n-2$, finding a ``nice'' explicit section to (\ref{forsection}) becomes more complicated. For instance when $n=4$ and $\ell=2$, we can use a variant of the complex $\tld{\mathbf{D}}^\bullet$ of (\ref{equ: explicit dR split top}) for ${\rm GL}_3(K)$ combined with parabolic induction to ${\rm GL}_4(K)$ to build a complex of finite length coadmissible $D(G)$-modules with Orlik-Strauch constituents which is exact in degrees $0,1$ and maps to $\tau_{\leq 2}\mathbf{D}^\bullet$ with an isomorphism on $H^2$, hence which gives an explicit section. But this complex is not nice (contrary to $\tld{\mathbf{D}}^\bullet$). It just gives us enough confidence to state the following conjecture.

\begin{conj}\label{conj: split}
For $2\leq \ell\leq n-2$ the morphism of complexes (\ref{forsection}) admits a section in the derived category of finite length coadmissible $D(G)$-modules with Orlik-Strauch constituents.
\end{conj}

Recall we know that a section exists in the derived category of \emph{all} (abstract) $D(G)$-modules by \cite[Thm.~6.1]{Schr11} and Dat (\cite[Cor.~A.1.3]{Dat06}.

\subsection{Application to de Rham complex of the Drinfeld space}

We show that, for $E=K$ and $\mu_0=(0,\cdots,0)$, the complex of coadmissible $D(G)$-modules $\mathbf{D}^\bullet$ in (\ref{equ: explicit dR}) is isomorphic to the global sections of the de Rham complex of the rigid analytic Drinfeld space over $K$ of dimension $n-1$, and for arbitrary $\mu_0$ to the complex of holomorphic discrete series of \cite{S92}.\bigskip

Throughout this section, we use the notation $I\defeq \widehat{1}$ and $J\defeq \widehat{n-1}$. Recall from \S\ref{generalnotation} that $w_0$ (resp.~$w_I$) is the longest element of $W(G)$ (resp.~of $W(L_I)$). We check that $w_0=w_{n-1,1}w_I=w_Iw_{1,n-1}$. We also keep the notation of \S\ref{subsec: final} (in particular we have fixed a weight $\mu_0$ in $\Lambda^{\dom}$).\bigskip

We start with two more results on coadmissible $D(G)$-modules which will be used in Theorem \ref{thm: main dR} below. The first statement shows that the $D(G)$-modules $D_k$ and $\tld{D}_{k}$ of (\ref{definitionofd0}), (\ref{definitionofdk}) have a nice behaviour with respect to wall-crossing functors.

\begin{thm}\label{thm: D wall crossing}
Let $1\leq k\leq n-1$ and $\mu\in \Lambda$ such that $\langle\mu+\rho,\alpha^\vee\rangle\geq 0$ for $\al\in \Phi^+$ and the stabilizer of $\mu$ in $W(G)$ for the dot action is $\{1,s_k\}$.
\begin{enumerate}[label=(\roman*)]
\item \label{it: D wall crossing 0} The morphisms ${D}_{k-1}\twoheadrightarrow X_k \hookrightarrow {D}_{k}$ from \ref{it: main rep 2} of Theorem \ref{thm: main rep} induce isomorphisms $\Theta_{\mu}({D}_{k-1})\buildrel\sim\over\rightarrow \Theta_{\mu}(X_k) \buildrel\sim\over\rightarrow \Theta_{\mu}({D}_k)$.
\item \label{it: D wall crossing 1} The morphisms $\tld{D}_{k-1}\twoheadrightarrow Z_k \hookrightarrow \tld{D}_{k}$ from \ref{it: main rep 4} of Theorem \ref{thm: main rep} induce isomorphisms $\Theta_{\mu}(\tld{D}_{k-1})\buildrel\sim\over\rightarrow \Theta_{\mu}(Z_k) \buildrel\sim\over\rightarrow \Theta_{\mu}(\tld{D}_k)$.
\item \label{it: D wall crossing 2} We have non-split short exact sequences $0\rightarrow D_k\rightarrow \Theta_{\mu}(D_k)\rightarrow D_{k-1}\rightarrow 0$ where $D_k\rightarrow \Theta_{\mu}(D_k)$ is the canonical adjunction map.
\item \label{it: D wall crossing 3} If $1\leq k\leq n-2$ we have non-split short exact sequences $0\rightarrow \tld{D}_k\rightarrow \Theta_{\mu}(\tld{D}_k)\rightarrow \tld{D}_{k-1}\rightarrow 0$ where $\tld{D}_k\rightarrow \Theta_{\mu}(\tld{D}_k)$ is the canonical adjunction map.
\item \label{it: D wall crossing 4} If $k=n-1$ the composition $\Theta_{\mu}(\tld{D}_{n-1})\cong \Theta_{\mu}(\tld{D}_{n-2})\rightarrow \tld{D}_{n-2}$ (where the first isomorphism follows from \ref{it: D wall crossing 1} and the second map is the adjunction map) is surjective and its kernel is the ``universal'' extension
\begin{equation}\label{equ: universal first branch}
0\longrightarrow Z_{n-1}\longrightarrow \ast \longrightarrow \big((\mathrm{St}_n^{\rm{alg}})^\vee\big)^{\oplus n}\longrightarrow 0
\end{equation}
deduced from \ref{it: main Ext 1} of Theorem~\ref{thm: main Ext}.
\end{enumerate}
\end{thm}
\begin{proof}
We prove \ref{it: D wall crossing 0} and \ref{it: D wall crossing 1}. One checks that any $V\in \mathrm{JH}_{G}(\mathrm{ker}(\tld{D}_{k}^\vee\twoheadrightarrow Z_k^\vee))$ (resp.~any $V\in \mathrm{JH}_{G}(\tld{D}_{k-1}^\vee/Z_k^\vee)$) either is locally algebraic or is in $\tld{\mathbf{J}}_{n-k-1}$ if $k<n-1$ (resp.~in $\tld{\mathbf{J}}_{n-k+1}$ if $1<k$) (see (\ref{tldJj0def})). Hence we deduce from \ref{it: wall crossing simple 1} of Lemma~\ref{lem: wall crossing of simple} that $\Theta_{\mu}(V)=0$. The exactness of $\Theta_{\mu}(-)$ then gives the isomorphisms in \ref{it: D wall crossing 1}. The proof of \ref{it: D wall crossing 0} is completely analogous.

We prove \ref{it: D wall crossing 2} and \ref{it: D wall crossing 3}. When $1\leq k\leq n-2$, \ref{it: D wall crossing 2} (resp.~\ref{it: D wall crossing 3}) follows from the definition of $D_k$ (resp.~of $\tld{D}_k$) in (\ref{definitionofdk}) and \ref{it: general wall crossing 3} of Proposition~\ref{prop: general wall crossing} applied with $V_0^-=D_k^\vee$ (i.e.~$j_0=n-k$) and $(j_1,j_2)=(n-1,n-1)$ (resp.~and $(j_1,j_2)=(1,n)$), noting that in the proof of \emph{loc.~cit.} we have $\Theta_\mu(V_0)\cong \Theta_\mu(V_0^-)$. When $k=n-1$, \ref{it: D wall crossing 2} follows from the $(j_1,j_2)=(n-1,n-1)$ case of \ref{it: wall crossing first branch 2} of Lemma~\ref{lem: wall crossing first branch} (where $V=X_{n-1}^\vee$ and $V^+=D_{n-2}^\vee$), the fact that $D_{n-1}$ is the unique $D(G)$-module that fits into a non-split extension $0\rightarrow X_{n-1}\rightarrow D_{n-1}\rightarrow (\mathrm{St}_n^{\rm{alg}})^\vee\rightarrow 0$ (\ref{it: main Ext 0} of Theorem \ref{thm: main Ext}), and $\Theta_\mu(X_{n-1})\buildrel\sim\over\rightarrow \Theta_\mu(D_{n-1})$ since $\Theta_\mu((\mathrm{St}_n^{\rm{alg}})^\vee)=0$ by \ref{it: wall crossing simple 1} of Lemma~\ref{lem: wall crossing of simple}. The non-splitness easily follows from \ref{it: wall crossing simple 2} of Lemma~\ref{lem: wall crossing of simple}.

Finally \ref{it: D wall crossing 4} follows from the $(j_1,j_2)=(1,n)$ case of \ref{it: wall crossing first branch 2} of Lemma~\ref{lem: wall crossing first branch} with $V_0=Z_{n-1}^\vee$ and $V_0^+=\tld{D}_{n-2}^\vee$ in \emph{loc.~cit.}~(and noting that $\Theta_{\mu}(V_0)=\Theta_\mu(V_0^+)$).
\end{proof}

\begin{rem}\label{secondadjunction}
With the notation of Theorem \ref{thm: D wall crossing}, for $1\leq k\leq n-1$ the composition $\Theta_{\mu}({D}_k)\twoheadrightarrow D_{k-1} \rightarrow D_k$ where the last map is the differential map $d_{\mathbf{D}}^{k-1}$ (see (\ref{equ: explicit dR})) is nothing else than the (non-zero) canonical adjunction map $\Theta_{\mu}({D}_k)\rightarrow D_k$. Indeed, by functoriality of the adjunction maps and since $\Theta_{\mu}(X_k)\buildrel\sim\over\rightarrow \Theta_{\mu}({D}_k)$ (see \ref{it: D wall crossing 0} of Theorem \ref{thm: D wall crossing}), this adjunction map factors as $\Theta_{\mu}({D}_k)\rightarrow X_k\hookrightarrow D_k$. One easily checks from (\ref{definitionofxyz}) and (\ref{definitionofdk}) that the (irreducible) cosocle of ${D}_k$ does not appear in $X_k$. It then follows from \ref{it: D wall crossing 2} of Theorem \ref{thm: D wall crossing} that the adjunction map $\Theta_{\mu}({D}_k)\rightarrow D_k$ factors through a non-zero map $D_{k-1} \rightarrow D_k$, which must be $d_{\mathbf{D}}^{k-1}$ (up to a non-zero scalar) by unicity of $d_{\mathbf{D}}^{k-1}$ (see the references above (\ref{equ: explicit dR})). A similar proof replacing $X_k$, $X_{k-1}$ by $Z_k$, $Z_{k-1}$ and using \ref{it: D wall crossing 1}, \ref{it: D wall crossing 3} and \ref{it: D wall crossing 4} of Theorem \ref{thm: D wall crossing} gives that, for $1\leq k\leq n-1$, the composition $\Theta_{\mu}(\tld{D}_k)\twoheadrightarrow \tld{D}_{k-1} \rightarrow \tld{D}_k$ (where the last map is $d_{\tld{\mathbf{D}}}^{k-1}$, see (\ref{equ: explicit dR split top})) is the (non-zero) adjunction map $\Theta_{\mu}(\tld{D}_k)\rightarrow \tld{D}_k$ (for $k=n-1$ one has to use (\ref{equ: universal first branch})).
\end{rem}

Using Theorem \ref{thm: D wall crossing} we can prove the following unicity statements which strengthen Theorem \ref{thm: main rep}.

\begin{cor}\label{unicitystuff}
We have the following unicity results.
\begin{enumerate}[label=(\roman*)]
\item \label{it: unicity rep 0} For $1\leq k\leq n-1$, $Z_k$ is the unique coadmissible $D(G)$-module of the form
\[\begin{xy}(0,0)*+{Y_{k}}="a"; (22.5,0)*+{(V_{[1,n-k],\Delta}^{\rm{alg}})^\vee}="b"; (45.5,0)*+{X_{k}}="c";
{\ar@{-}"a";"b"}; {\ar@{-}"b";"c"}\end{xy}\!.\]
\item \label{it: unicity rep 1} For $0\leq k\leq n-1$, $D_k$ is the unique coadmissible $D(G)$-module of the form
\[\begin{xy}(0,0)*+{X_{k}}="a"; (25.2,0)*+{(V_{[1,n-k-1],\Delta}^{\rm{alg}})^\vee}="b"; (52,0)*+{X_{k+1}}="c";
{\ar@{-}"a";"b"}; {\ar@{-}"b";"c"}\end{xy}\!.\]
\item \label{it: unicity rep 2} For $0\leq k\leq n-2$, $\tld{D}_k$ is the unique coadmissible $D(G)$-module of the form
\[\begin{xy}(0,0)*+{Z_{k}}="a"; (17,0)*+{Z_{k+1}}="b"; {\ar@{-}"a";"b"}\end{xy}\!.\]
\end{enumerate}
\end{cor}
\begin{proof}
We prove \ref{it: unicity rep 0}. It follows from \ref{it: Ext1 with alg 1} of Lemma \ref{lem: Ext1 factor 2} that for $j\in \{1,\dots,n-1\}$ we have $\mathrm{Ext}_{D(G)}^1(C_{n-k,j,j}^\vee, (V_{[1,n-k],\Delta}^{\rm{alg}})^\vee) \neq 0$ if and only if $j=n-k$, in which case it is one dimensional. It follows from \ref{it: Ext2 with alg 1} of Lemma \ref{lem: Ext2 factor 2} that $\mathrm{Ext}_{D(G)}^2(C_{n-k,j,j}^\vee, (V_{[1,n-k],\Delta}^{\rm{alg}})^\vee)=0$ for $1\leq j<n-k$. By d\'evissage and (\ref{definitionofxyz}) this implies
\begin{equation}\label{eq: dimext=1pourintro}
\Dim_E \mathrm{Ext}_{D(G)}^1(X_k, (V_{[1,n-k],\Delta}^{\rm{alg}})^\vee)=1,
\end{equation}
in particular it is enough to prove $\mathrm{Ext}_{D(G)}^1(X_k,Y_k)=0$. 

Assume the contrary and let $M$ be a non-split extension of $X_k$ by $Y_k$. Using again \ref{it: Ext1 with alg 1} of Lemma \ref{lem: Ext1 factor 2}, for $j\in \{1,\dots,n-1\}$ and $i\in \{1,\dots, k-1\}$ we have
\[\mathrm{Ext}_{D(G)}^1(C_{n-k,j,j}^\vee, (V_{[1+i,n-k+i],\Delta}^{\rm{alg}})^\vee) = 0,\]
and hence $\mathrm{Ext}_{D(G)}^1(X_k, (V_{[1+i,n-k+i],\Delta}^{\rm{alg}})^\vee) = 0$ for $i\in \{1,\dots, k-1\}$ (see (\ref{definitionofxyz})). It follows that there exist $j\in \{1,\dots,n-1\}$ and $(n-k,j_1,j_2)\in {\mathbf{J}}$ with $j_1<j_2$ such that $M$ has a length $2$ subquotient with socle $C_{n-k,j_1,j_2}^\vee$ in $Y_k$ and cosocle $C_{n-k,j,j}^\vee$ in $X_k$ (see (\ref{definitionofxyz}) and the lines that follow). Moreover by Lemma \ref{lem: Ext1 factor 1} we must have $|j_2-j|+|j_2-j_1|=1$, i.e.~$j_2=j$ and $j_1=j_2-1=j-1$ (which implies $j\geq 2$). In other terms $C_{n-k,j-1,j}^\vee$ is the only constituent of $Y_k$ that can lie ``just below'' the constituent $C_{n-k,j,j}^\vee$ of $X_k$. 

Using this result we now prove that, for \emph{any} $j'\in \{1,\dots,n-1\}$, $M$ must contain as a subquotient the unique non-split extension of $C_{n-k,j',j'}^\vee$ by $C_{n-k,j'-1,j'}^\vee$ (using again Lemma \ref{lem: Ext1 factor 1}). We have just seen this holds for $j'=j$. Assume this fails for some $j'>j$ and take the minimal such $j'$. Then, by definition of $j'$, $M$ contains as a subquotient the unique non-split extension of $C\defeq C_{n-k,j'-1,j'-1}^\vee$ by $B\defeq C_{n-k,j'-2,j'-1}^\vee$. Note that, from (\ref{definitionofxyz}) and the definition of the partial order on ${\mathbf{J}}$, $Y_k$ and hence $M$ also contain as a subquotient the unique non-split extension of $B$ by $A\defeq C_{n-k,j'-1,j'}^\vee$. Moreover there cannot exist a fourth constituent $B'$ of $M$ distinct from $B$ such that $A<B'<C$ for the partial order on $\mathrm{JH}_{D(G)}(M)$. Indeed, if such a $B'$ exists, one can take it such that $\mathrm{Ext}_{D(G)}^1(B',A)\ne 0$. If $B'\in \mathrm{JH}_{D(G)}(Y_k)$, then the partial order on $Y_k$ and (\ref{definitionofxyz}) force $B'=B$, a contradiction. Hence $B'\in \mathrm{JH}_{D(G)}(X_k)$, but then arguing as at the end of the previous paragraph we must have $B'=C_{n-k,j',j'}^\vee$, contradicting the hypothesis on $j'$. It follows that $M$ contains as a subquotient a uniserial $D(G)$-module of the form $\!\begin{xy}(0,0)*+{A}="a"; (10,0)*+{B}="b"; (20,0)*+{C}="c";{\ar@{-}"a";"b"}; {\ar@{-}"b";"c"}\end{xy}\!$, contradicting the minimality of the $\mathrm{Ext}$-square $V_{(n-k,j'-1,j'-1),(n-k,j'-1,j')}$ in Proposition \ref{prop: hard square}. A symmetric argument when $j'<j$ also yields a contradiction. 

Finally, applying the previous result with $j'=n-k$ and $j'=n-k+1$ (and using the structures of $X_k$ and $Y_k$ from (\ref{definitionofxyz})), we deduce that $M$ contains as a subquotient the dual of the minimal $\mathrm{Ext}$-square $V_{(n-k,n-k,n-k),(n-k,n-k,n-k+1)}$ in \ref{it: hard square 2} of Proposition \ref{prop: hard square}, in particular $M$ contains $(V_{[1,n-k],\Delta}^{\rm{alg}})^\vee$ which contradicts the lines below (\ref{definitionofxyz}). Hence $\mathrm{Ext}_{D(G)}^1(X_k,Y_k)=0$, which proves \ref{it: unicity rep 0}.

We prove \ref{it: unicity rep 1}. The case $k=0$ follows from (\ref{eq: dimext=1pourintro}) for $k=1$, while the case $k=n-1$ is \ref{it: main Ext 0} of Theorem \ref{thm: main Ext}. We assume $k\in \{1,\dots, n-2\}$ and let $M$ be any coadmissible $D(G)$-module as in \ref{it: unicity rep 1}. By the same proof as for \ref{it: D wall crossing 0} of Theorem \ref{thm: D wall crossing} we have $\Theta_{\mu}(X_k)\buildrel\sim\over\rightarrow \Theta_{\mu}(M)$ (with $\mu$ as in \emph{loc.~cit.}), and thus by adjunction (for $M$) a non-zero map $M\rightarrow \Theta_{\mu}(X_k)$. Since $C_{(n-k-1,1,1)}^\vee$ is the (irreducible) cosocle of $X_{k+1}$, it also appears in the cosocle of $M$. Assume first that $C_{(n-k-1,1,1)}^\vee$ maps to $0$ in $\Theta_{\mu}(X_k)$. Then from the form of $M$ it follows that a non-zero strict quotient of $X_k$ embeds into $\Theta_{\mu}(X_k)$, and thus (from (\ref{definitionofxyz})) there exists $1\leq j\leq n-2$ such that $C_{(n-k,j,j)}^\vee$ appears in the socle of $\Theta_{\mu}(X_k)$. By \ref{it: D wall crossing 0} and \ref{it: D wall crossing 2} of Theorem \ref{thm: D wall crossing}, this implies that $C_{(n-k,j,j)}^\vee$ embeds into $D_k$ or into $D_{k-1}$, which is impossible by (\ref{definitionofdk}). Hence $C_{(n-k-1,1,1)}^\vee$ still occurs in the (cosocle of the) image of $M$ in $\Theta_{\mu}(X_k)$. By \ref{it: D wall crossing 0}, \ref{it: D wall crossing 2} of Theorem \ref{thm: D wall crossing} and (\ref{definitionofdk}), $C_{(n-k-1,1,1)}^\vee$ has multiplicity $1$ in $\Theta_{\mu}(X_k)$. Hence the image of $M$ contains the unique (closed) $D(G)$-submodule of $\Theta_{\mu}(X_k)$ of cosocle $C_{(n-k-1,1,1)}^\vee$. But this submodule is $D_k$ using \emph{loc.~cit.} again. Since $M$ and $D_k$ are multiplicity free with the same irreducible constituents, if follows that $M\buildrel\sim\over\rightarrow D_k$. 

The proof of \ref{it: unicity rep 2} is analogous to \ref{it: unicity rep 1} replacing $X_k$, $X_{k+1}$ by $Z_k$, $Z_{k+1}$ and \ref{it: D wall crossing 0}, \ref{it: D wall crossing 2} of Theorem \ref{thm: D wall crossing} by \ref{it: D wall crossing 1}, \ref{it: D wall crossing 3} of Theorem \ref{thm: D wall crossing}.
\end{proof}

We now prove another unicity theorem which will play a key role in the comparison of the complex $\mathbf{D}^\bullet$ when $\mu_0=(0,\cdots,0)$ with the (global sections of) the de Rham complex of the Drinfeld space (Theorem \ref{thm: main dR}).

\begin{thm}\label{prop: top deg}
Let $C$ be a finite length coadmissible $D(G)$-module equipped with a decreasing filtration
\[C=\mathrm{Fil}^0(C)\supseteq \mathrm{Fil}^1(C)\supseteq \cdots \supseteq\mathrm{Fil}^{n-1}(C)=0\]
which satisfies the following conditions:
\begin{enumerate}[label=(\roman*)]
\item \label{it: top deg 1} we have $H^0(N_{J},C)\cong L^{J}(w_{n-1,1})$;
\item \label{it: top deg 2} for $0\leq \ell\leq n-2$, there exists a $U(\fg)$-module $M_{\ell}$ in $\cO^{\fp_{\widehat{\ell+1}}}_{\rm{alg}}$ (see the beginning of \S\ref{subsec: category}) such that $\mathrm{gr}^{0}(C)=C/\mathrm{Fil}^1(C)$ fits into a short exact sequence of coadmissible $D(G)$-modules
\[0\rightarrow \cF_{P_{I}}^{G}(M_0,\pi_{1,1}^{\infty})^\vee \rightarrow \mathrm{gr}^{0}(C) \rightarrow (\mathrm{St}_n^{\rm{alg}})^\vee \rightarrow 0,\]
and $\mathrm{gr}^{\ell}(C)\cong \cF_{P_{\widehat{\ell+1}}}^{G}(M_{\ell},\pi_{\ell+1,\ell+1}^{\infty})^\vee$ if $\ell\neq 0$ (see (\ref{pij_1,_2}) for $\pi_{\ell+1,\ell+1}^{\infty}$).
\end{enumerate}
Then we have $C\cong D_{n-1}$ (and $M_\ell=L(w_{\ell+1,1})$ for $0\leq \ell\leq n-2$).
\end{thm}
\begin{proof}
Note that in \ref{it: top deg 2} we do not specify $M_{\ell}$ neither whether the short exact sequence is split or not. But Condition \ref{it: top deg 1} (which was inspired by \cite[Prop.~6.3]{Schr11} when $G=\GL_3(\Qp)$) ``rigidifies'' everything.

As $\cF_{P_{\widehat{\ell+1}}}^{G}(M_{\ell},\pi_{\ell+1,\ell+1}^{\infty})^\vee$ injects into $\mathrm{gr}^{\ell}(C)$ by condition \ref{it: top deg 2}, we define $\mathrm{Fil}^{\ell}(C)'$ for $0\leq \ell\leq n-2$ as the unique (closed) $D(G)$-submodule of $\mathrm{Fil}^{\ell}(C)$ which fits into a short exact sequence
\begin{equation}\label{equ: step 2 filtration seq}
0\rightarrow\mathrm{Fil}^{\ell+1}(C)\rightarrow\mathrm{Fil}^{\ell}(C)'\rightarrow \cF_{P_{\widehat{\ell+1}}}^{G}(M_{\ell},\pi_{\ell+1,\ell+1}^{\infty})^\vee\rightarrow 0
\end{equation}
(in particular $\mathrm{Fil}^{\ell}(C)'=\mathrm{Fil}^{\ell}(C)$ for $1\leq \ell\leq n-2$). For $\ell\geq 1$, we define $\mathrm{Fil}^{\ell}(D_{n-1})'=\mathrm{Fil}^{\ell}(D_{n-1})$ as the unique (closed) $D(G)$-submodule of $D_{n-1}$ such that $\mathrm{JH}_{G}(\mathrm{Fil}^{\ell}(D_{n-1})^\vee)=[C_{(1,\ell+1,\ell+1)},C_{(1,n-1,n-1)}]$. For $\ell=0$ we define $\mathrm{Fil}^{0}(D_{n-1})\defeq D_{n-1}$ and $\mathrm{Fil}^{0}(D_{n-1})'\defeq X_{n-1}\subseteq D_{n-1}$. We check from the definitions of $D_{n-1}$ and $X_{n-1}$ that $\mathrm{gr}^{0}(D_{n-1})$ fits into a non-split extension
\[0\rightarrow \cF_{P_{I}}^{G}(L(s_{1}),\pi_{1,1}^{\infty})^\vee \rightarrow \mathrm{gr}^{0}(D_{n-1}) \rightarrow (\mathrm{St}_n^{\rm{alg}})^\vee \rightarrow 0,\]
and that $\mathrm{gr}^{\ell}(D_{n-1})\cong C_{(1,\ell+1,\ell+1)}^\vee=\cF_{P_{\widehat{\ell+1}}}^{G}(L(w_{\ell+1,1}),\pi_{\ell+1,\ell+1}^{\infty})^\vee$ for $1\leq \ell\leq n-2$ (see (\ref{cj})).\bigskip

\textbf{Step $1$}: We prove $M_{n-2}\cong L(w_{n-1,1})$.\\
From \ref{it: top deg 1} we deduce
\begin{equation}\label{equ: dR U coh}
H^0(U,C)\cong H^0(U_{J}, H^0(N_{J},C))\cong w_{n-1,1}\cdot\mu_0.
\end{equation}
Note that $\cF_{P_{J}}^{G}(M_{n-2},\pi_{n-1,n-1}^{\infty})^\vee\cong \mathrm{gr}^{n-2}(C)$ injects into $C$ by \ref{it: top deg 2}. Since $H^0(U,-)$ is left exact and since we have by (\ref{equ: N coh H0})
\[H^0(U,\cF_{P_{J}}^{G}(M_{n-2},\pi_{n-1,n-1}^{\infty})^\vee)\cong H^0(\fu,M_{n-2})\otimes_E (J_{J,\emptyset}(\pi_{n-1,n-1}^{\infty}))^\vee\neq 0,\]
we deduce the following isomorphism of $D(T)$-modules
\[H^0(\fu,M_{n-2})\otimes_E (J_{J,\emptyset}(\pi_{n-1,n-1}^{\infty}))^\vee \cong H^0(U,C)\cong w_{n-1,1}\cdot\mu_0.\]
In particular we have $J_{J,\emptyset}(\pi_{n-1,n-1}^{\infty})=J_{J,\emptyset}(1_{L_{J}})=1_{T}$ and an isomorphism of $U(\ft)$-modules
\begin{equation}\label{equ: step 1 n coh}
H^0(\fu,M_{n-2})\cong w_{n-1,1}\cdot\mu_0.
\end{equation}
Left exactness of $H^0(\fu,-)$ with \ref{it: dominance 2} of Lemma~\ref{lem: n coh dominance} force $\mathrm{soc}_{U(\fg)}(M_{n-2})\cong L(w_{n-1,1})$. If $M_{n-2}\not\cong L(w_{n-1,1})$, then $M_{n-2}$ contains a length $2$ $U(\fg)$-submodule $M_{n-2}'$ with cosocle $L(x)$ such that $D_L(x)\subseteq \{n-1\}$ (using $M_{n-2}\in\cO^{\fp_{J}}_{\rm{alg}}$ and Lemma~\ref{lem: dominance and left set}), or equivalently $x\in W^{J,\emptyset}$ (see \S\ref{generalnotation} for $W^{J,\emptyset}$). Moreover $x\neq w_{n-1,1}$ by \ref{it: rabiotext 1} of Lemma~\ref{rabiotext}, and hence $x<w_{n-1,1}$ since $x\in W^{J,\emptyset}$ and $w_{n-1,1}$ is the maximal element in $W^{J,\emptyset}$. But $x<w_{n-1,1}$ together with (\ref{equ: O Hom radical}) imply that $M_{n-2}'$ is a quotient of $M(x)$, which by (\ref{equ: g spectral seq 0}) (applied with $I=\emptyset$) implies
\[0\ne \Hom_{U(\fg)}(M(x),M'_{n-2})\hookrightarrow \Hom_{U(\fg)}(M(x),M_{n-2})\cong \Hom_{U(\ft)}(x\cdot\mu_0,H^0(\fu,M_{n-2})),\]
contradicting (\ref{equ: step 1 n coh}). Hence, we have $M_{n-2}\cong L(w_{n-1,1})$. Note that by (\ref{equ: N coh H0}) we have $H^0(U,V^\vee)\neq 0$ for $V=\cF_{P_I}^{G}(M,\pi^{\infty})$ with $M\in\cO^{\fp_I}_{\rm{alg}}$ and $\pi^{\infty}$ $G$-basic (\ref{it: basic 1} of Definition \ref{def: basic rep}). By left exactness of $H^0(U,-)$ we deduce using (\ref{equ: dR U coh}) and $M_{n-2}\cong L(w_{n-1,1})$
\begin{equation}\label{equ: step 1 socle}
\mathrm{soc}_{D(G)}(C)\cong \mathrm{gr}^{n-2}(C)\cong C_{(1,n-1,n-1)}^\vee.
\end{equation}
Note that, by Remark~\ref{rem: OS isotypic}, (\ref{equ: step 1 socle}) and condition \ref{it: top deg 2} imply $(M_{\ell})_{\xi}=M_{\ell}$ for $0\leq \ell\leq n-2$ where $\xi: Z(\fg)\rightarrow E$ is the unique infinitesimal character such that $L(1)_{\xi}\neq 0$. In particular any irreducible constituent of $M_{\ell}$ (for $0\leq \ell\leq n-2$) is of the form $L(x)$ for $x\in W(G)$.\bigskip

\textbf{Step $2$}: We prove that $L(1)\notin\mathrm{JH}_{U(\fg)}(M_{\ell})$ for $1\leq \ell\leq n-2$.\\
We fix $\mu\in \Lambda$ such that $\langle\mu+\rho,\alpha^\vee\rangle\geq 0$ for $\al\in \Phi^+$ and the stabilizer of $\mu$ in $W(G)$ for the dot action is $\{1,s_{n-1}\}$. Note that $\mathrm{Fil}^{1}(C)\subseteq C$ has socle $C_{(1,n-1,n-1)}^\vee$ by (\ref{equ: step 1 socle}) (for $n\geq 3$), and that the canonical adjunction map $\mathrm{Fil}^{1}(C)\rightarrow \Theta_{\mu}(\mathrm{Fil}^{1}(C))$ is injective as it restricts to the canonical injection $C_{(1,n-1,n-1)}^\vee\hookrightarrow \Theta_{s_{n-1}}(C_{(1,n-1,n-1)}^\vee)$ (see \ref{it: wall crossing simple 2} of Lemma~\ref{lem: wall crossing of simple}). The decreasing filtration $\mathrm{Fil}^{\ell}(C)$ on $\mathrm{Fil}^{1}(C)$ for $1\leq \ell\leq n-2$ induces a decreasing filtration $\Theta_{\mu}(\mathrm{Fil}^{\ell}(C))$ on $\Theta_{\mu}(\mathrm{Fil}^{1}(C))$ with graded pieces $\Theta_{s_{n-1}}(\mathrm{gr}^{\ell}(C))\cong \cF_{P_{\widehat{\ell+1}}}^{G}(\Theta_{s_{n-1}}(M_{\ell}),\pi_{\ell+1,\ell+1}^{\infty})^\vee$ by condition \ref{it: top deg 2} and (\ref{equ: OS wall crossing}). It follows from Proposition~\ref{prop: Jantzen middle} that $L(1)\in\mathrm{JH}_{U(\fg)}(\Theta_{s_{n-1}}(L(x)))$ for some $x\in W(G)$ if and only if $x=s_1$. Since any $L(x)\in\mathrm{JH}_{U(\fg)}(M_{\ell})$ satisfies $D_L(x)\subseteq \{\ell+1\}$ by Lemma~\ref{lem: dominance and left set}, we deduce $x\neq s_1$, and thus $L(1)\notin\mathrm{JH}_{U(\fg)}(\Theta_{s_{n-1}}(L(x)))$, for any $L(x)\in\mathrm{JH}_{U(\fg)}(M_{\ell})$ and any $1\leq \ell\leq n-2$. By the above discussion we deduce that $\Theta_{\mu}(\mathrm{Fil}^{1}(C))$, and hence its $D(G)$-submodule $\mathrm{Fil}^{1}(C)$, do not admit constituents which are duals of locally algebraic representations, or equivalently $L(1)\notin\mathrm{JH}_{U(\fg)}(M_{\ell})$ for $1\leq \ell\leq n-2$ by condition \ref{it: top deg 2} (and Lemma~\ref{lem: Hom OS}).\bigskip

\textbf{Step $3$}: We prove that $M_{\ell}\neq 0$ for $0\leq \ell\leq n-2$.\\
From Step $1$ we have $M_{\ell}\neq 0$ for $\ell = n-2$. Assume first $M_{0}=0$. As $(\mathrm{St}_n^{\rm{alg}})^\vee$ shows up in $C$ by \ref{it: top deg 2} but not in the socle of $C$ by (\ref{equ: step 1 socle}), from condition \ref{it: top deg 2} there exists $1\leq \ell\leq n-2$ and $L(x)\in\mathrm{JH}_{U(\fg)}(M_{\ell})$ such that
\begin{equation*}
\mathrm{Ext}_{D(G)}^1((\mathrm{St}_n^{\rm{alg}})^\vee, \cF_{P_{\widehat{\ell+1}}}^{G}(L(x),\pi_{\ell+1,\ell+1}^{\infty})^\vee)\neq 0.
\end{equation*}
As $x\neq 1$ by {Step $2$}, the last statement in Remark~\ref{rem: Ext1 OS distance} implies $d(\pi_{\ell+1,\ell+1}^{\infty},\mathrm{St}_n^{\infty})=0$. But by \ref{it: explicit induction 2} of Lemma \ref{lem: explicit smooth induction}, for $1\leq \ell\leq n-2$ the representation $i_{\widehat{\ell+1},\Delta}^{\infty}(\pi_{\ell+1,\ell+1}^{\infty})$ has length $2$ with socle $V_{[1,\ell+1],\Delta}^{\infty}$ and cosocle $V_{[1,\ell],\Delta}^{\infty}$, and thus we can't have $d(\pi_{\ell+1,\ell+1}^{\infty},\mathrm{St}_n^{\infty})=0$. Hence we have $M_0\ne 0$. Assume now $M_{\ell}=0$ for some $0<\ell<n-2$. Since $M_{0}\neq 0$, we may choose $\ell$ such that $M_{\ell-1}\neq 0=M_{\ell}$. Let $L(y)\subseteq M_{\ell-1}$ be an arbitrary irreducible $U(\fg)$-submodule and recall that $D_L(y)\subseteq \{\ell\}$ by Lemma~\ref{lem: dominance and left set}. As $M_{\ell}=0$ and $\cF_{P_{\widehat{\ell}}}^{G}(L(y),\pi_{\ell,\ell}^{\infty})^\vee$ has no common constituent with $\mathrm{soc}_{D(G)}(C)$ by (\ref{equ: step 1 socle}) and $D_L(y)\subseteq \{\ell\}$ (using Lemma \ref{lem: Hom OS}), there exist $\ell<\ell'\leq n-2$ and $L(z)\in \mathrm{JH}_{U(\fg)}(M_{\ell'})$ such that
\begin{equation*}
\mathrm{Ext}_{D(G)}^1\big(\cF_{P_{\widehat{\ell}}}^{G}(L(y),\pi_{\ell,\ell}^{\infty})^\vee, \cF_{P_{\widehat{\ell'+1}}}^{G}(L(z),\pi_{\ell'+1,\ell'+1}^{\infty})^\vee\big)\neq 0.
\end{equation*}
Note that $z\ne 1$ by Step $2$ (as $\ell'\geq 1$) and hence $D_L(z)=\{\ell'+1\}$, in particular $z\neq y$ (as $D_L(y)\subseteq \{\ell\}$). But then, the last statement of Remark~\ref{rem: Ext1 OS distance} again implies $d(\pi_{\ell'+1,\ell'+1}^{\infty},\pi_{\ell,\ell}^{\infty})\!=0$, which contradicts \ref{it: connect 1} of Lemma~\ref{lem: connect Hom} since $\ell<\ell'$. This finishes the proof of Step $3$.\bigskip

\textbf{Step $4$}: Let $0\leq \ell<\ell'\leq n-2$ and $L(x)\in\mathrm{JH}_{U(\fg)}(M_{\ell})$, we prove that if
\begin{equation}\label{equ: Ext1 step 4}
\mathrm{Ext}_{D(G)}^1(\cF_{P_{\widehat{\ell+1}}}^{G}(L(x),\pi_{\ell+1,\ell+1}^{\infty})^\vee, C_{(1,\ell'+1,\ell'+1)}^\vee)\neq 0
\end{equation}
then $\ell'=\ell+1$ and $x=w_{\ell+1,1}$, in which case (\ref{equ: Ext1 step 4}) has dimension $1$.\\
As $D_L(w_{\ell'+1,1})=\{\ell'+1\}$ and $D_L(x)\subseteq \{\ell+1\}$ (by Lemma~\ref{lem: dominance and left set}), we have $x\neq w_{\ell'+1,1}$. By Proposition \ \ \ref{prop: Ext1 OS} \ \ we \ \ have \ \ $\mathrm{Ext}_{U(\fg)}^1(L(x),L(w_{\ell'+1,1}))\neq 0$, \ \ or \ \ equivalently \ \ $\mathrm{Ext}_{U(\fg)}^1(L(w_{\ell'+1,1}), L(x))\neq 0$ \ by (\ref{equ: dual Ext1}) \ and \ (\ref{equ: simple self dual}). \ As \ $\ell(w_{\ell'+1,1})\geq 2$, \ we \ have \ $\mathrm{Ext}_{U(\fg)}^1(L(w_{\ell'+1,1}), L(1))=0$ by \ref{it: rabiotext 2} of Lemma~\ref{rabiotext} and Lemma~\ref{lem: Ext with dominant}, thus $x\neq 1$ and $D_L(x)=\{\ell+1\}$. By \ref{it: special descent 3} of Lemma~\ref{lem: special descent change} and $\ell<\ell'$, we obtain $\ell'=\ell+1$ and $x=w_{\ell+1,1}$. Hence $\cF_{P_{\widehat{\ell+1}}}^{G}(L(x),\pi_{\ell+1,\ell+1}^{\infty})^\vee=C_{(1,\ell+1,\ell+1)}^\vee$ and (\ref{equ: Ext1 step 4}) is one dimensional by Lemma~\ref{lem: Ext1 factor 1}.\bigskip

\textbf{Step $5$}: We prove $\mathrm{Fil}^{\ell}(C)'\cong \mathrm{Fil}^{\ell}(D_{n-1})'$ for $0\leq \ell\leq n-2$.\\
We proceed by a decreasing induction on $0\leq \ell\leq n-2$. The case $\ell=n-2$ holds by {Step $1$}. We now assume $0\leq \ell\leq n-3$ and $\mathrm{Fil}^{\ell'}(C)'\cong \mathrm{Fil}^{\ell'}(D_{n-1})'$ for all $\ell'>\ell$, or equivalently $\mathrm{Fil}^{\ell'}(C)\cong \mathrm{Fil}^{\ell'}(D_{n-1})$ (since $\ell'>0$). Let $L(x)\subseteq M_{\ell}$ be an arbitrary (irreducible) $U(\fg)$-submodule, which induces a $D(G)$-submodule $\cF_{P_{\widehat{\ell+1}}}^{G}(L(x),\pi_{\ell+1,\ell+1}^{\infty})^\vee\subseteq \cF_{P_{\widehat{\ell+1}}}^{G}(M_{\ell},\pi_{\ell+1,\ell+1}^{\infty})^\vee$. By pullback and the induction hypothesis, (\ref{equ: step 2 filtration seq}) gives a short exact sequence
\begin{equation}\label{equ: step 5 seq 1}
0\rightarrow \mathrm{Fil}^{\ell+1}(D_{n-1})\rightarrow \ast \rightarrow \cF_{P_{\widehat{\ell+1}}}^{G}(L(x),\pi_{\ell+1,\ell+1}^{\infty})^\vee\rightarrow 0.
\end{equation}
It follows from (\ref{equ: step 1 socle}) (and $\ell\leq n-3$) that (\ref{equ: step 5 seq 1}) is non-split and thus
\begin{equation}\label{equ: Ext1 step 5 1}
\mathrm{Ext}_{D(G)}^1(\cF_{P_{\widehat{\ell+1}}}^{G}(L(x),\pi_{\ell+1,\ell+1}^{\infty})^\vee, \mathrm{Fil}^{\ell+1}(D_{n-1}))\neq 0.
\end{equation}
By \ {Step $4$} \ and \ a \ d\'evissage \ on \ $\mathrm{Fil}^{\ell+1}(D_{n-1})$, \ we \ deduce \ $x=w_{\ell+1,1}$ \ (hence $\cF_{P_{\widehat{\ell+1}}}^{G}(L(x),\pi_{\ell+1,\ell+1}^{\infty})^\vee=C_{(1,\ell+1,\ell+1)}^\vee$) and that the surjection $\mathrm{Fil}^{\ell+1}(D_{n-1})\twoheadrightarrow C_{(1,\ell+2,\ell+2)}^\vee$ induces \ \ an \ \ isomorphism \ \ between \ \ (\ref{equ: Ext1 step 5 1}) \ \ and \ \ the \ \ $1$-dimensional \ \ vector \ \ space \ $\mathrm{Ext}_{D(G)}^1(C_{(1,\ell+1,\ell+1)}^\vee, C_{(1,\ell+2,\ell+2)}^\vee)$.
In particular, $\mathrm{soc}_{U(\fg)}(M_{\ell})\cong L(w_{\ell+1,1})$ and the pushout of (\ref{equ: step 5 seq 1}) along $\mathrm{Fil}^{\ell+1}(D_{n-1})\twoheadrightarrow C_{(1,\ell+2,\ell+2)}^\vee$ is a non-split extension
\begin{equation}\label{equ: step 5 seq 2}
0\rightarrow C_{(1,\ell+2,\ell+2)}^\vee\rightarrow \ast \rightarrow C_{(1,\ell+1,\ell+1)}^\vee\rightarrow 0.
\end{equation}
If $M_{\ell}\cong L(w_{\ell+1,1})$, then both $\mathrm{Fil}^{\ell}(C)'$ and $\mathrm{Fil}^{\ell}(D_{n-1})'$ fit into a non-split extension
\[0\rightarrow\mathrm{Fil}^{\ell+1}(D_{n-1})\rightarrow \ast \rightarrow C_{(1,\ell+1,\ell+1)}^\vee\rightarrow 0,\]
which forces $\mathrm{Fil}^{\ell}(C)'\cong\mathrm{Fil}^{\ell}(D_{n-1})'$ as (\ref{equ: Ext1 step 5 1}) is $1$-dimensional for $x=w_{\ell+1,1}$. We now prove that we must have $M_{\ell}\cong L(w_{\ell+1,1})$. Assume on the contrary that $M_{\ell}$ is not $L(w_{\ell+1,1})$, then it contains a length $2$ $U(\fg)$-submodule $M'_{\ell}$ with socle $L(w_{\ell+1,1})$ and cosocle some $L(x')$ with $D_L(x')\subseteq \{\ell+1\}$ (so $x'\neq w_{\ell+2,1}$) and $\mathrm{Ext}_{\cO^{\fb}_{\rm{alg}}}^1(L(x'), L(w_{\ell+1,1}))\neq 0$. By \ref{it: rabiotext 1} of Lemma~\ref{rabiotext} this implies $x'\neq w_{\ell+1,1}$ and $\vert\ell(x')-\ell(w_{\ell+1,1})\vert$ odd. The pullback of (\ref{equ: step 2 filtration seq}) along $\cF_{P_{\widehat{\ell+1}}}^{G}(M_{\ell}',\pi_{\ell+1,\ell+1}^{\infty})^\vee\hookrightarrow \cF_{P_{\widehat{\ell+1}}}^{G}(M_{\ell},\pi_{\ell+1,\ell+1}^{\infty})^\vee$ followed by the pushout along $\mathrm{Fil}^{\ell+1}(D_{n-1})\twoheadrightarrow C_{(1,\ell+2,\ell+2)}^\vee$ gives a short exact sequence
\begin{equation}\label{equ: step 5 seq 3}
0\rightarrow C_{(1,\ell+2,\ell+2)}^\vee \rightarrow \ast \rightarrow \cF_{P_{\widehat{\ell+1}}}^{G}(M_{\ell}',\pi_{\ell+1,\ell+1}^{\infty})^\vee\rightarrow 0.
\end{equation}
Note that (\ref{equ: step 5 seq 3}) is non-split as it restricts to (\ref{equ: step 5 seq 2}) which is non-split.
In particular, we have
\begin{equation}\label{equ: existence of length three}
\mathrm{Ext}_{D(G)}^1(\cF_{P_{\widehat{\ell+1}}}^{G}(M_{\ell}',\pi_{\ell+1,\ell+1}^{\infty})^\vee, C_{(1,\ell+2,\ell+2)}^\vee)\neq 0.
\end{equation}

\textbf{Case $2.1$}:
If $x'=1$, then {Step $2$} forces $\ell=0$. Note that $M_{0}'$ is a quotient of $M^{I}(1)$ (using (\ref{equ: O Hom radical}) and \cite[Thm.~9.4(c)]{Hum08}) and that $L(w_{2,1})\notin\mathrm{JH}_{U(\fg)}(M^{I}(1))$ (by Lemma \ref{lem: dominance and left set}). Thus $C_{(1,2,2)}^\vee$ is not a constituent of $\cF_{P_{I}}^{G}(M^{I}(1),\pi_{1,1}^{\infty})^\vee$ (by Lemma~\ref{lem: Hom OS}). Hence, the surjection $M^{I}(1)\twoheadrightarrow M_{0}'$ induces an injection
\begin{equation}\label{equ: existence of length three prime}
\mathrm{Ext}_{D(G)}^1(\cF_{P_I}^{G}(M_{0}',\pi_{1,1}^{\infty})^\vee, C_{(1,2,2)}^\vee)\hookrightarrow \mathrm{Ext}_{D(G)}^1(\cF_{P_{I}}^{G}(M^{I}(1),\pi_{1,1}^{\infty})^\vee, C_{(1,2,2)}^\vee).
\end{equation}
Let $\xi_I: Z(\fl_{I})\rightarrow E$ be the unique infinitesimal character such that $L^{I}(1)_{\xi_I}\neq 0$. Since $w_{2,1}\notin W(L_I)$, by exactly the same argument as in the paragraph following (\ref{equ: simple length three object}) we have \ $H^k(\fn_{I}, L(w_{2,1}))_{\xi_I}=0$ \ for \ $k\geq 0$. \ By \ (\ref{equ: g spectral seq}) \ and \ Lemma~\ref{lem: central component} \ this \ implies \ $\mathrm{Ext}_{U(\fg)}^k(M^{I}(1), L(w_{2,1}))=0$ for $k\geq 0$. By (\ref{equ: OS main seq}) (applied with $w=1$, note that the assumption there is satisfed using \ref{it: coset 1} of Lemma \ref{lem: coset intersection}) this forces the right hand side of (\ref{equ: existence of length three prime}) to be zero, which contradicts (\ref{equ: existence of length three}) (when $\ell=0$).

\textbf{Case $2.2$}:
If $x'\neq 1$, then $D_L(x')=\{\ell+1\}$. It follows from Lemma~\ref{lem: maximal dominance} and Remark~\ref{rem: 1 dominance} that $D_R(x')=\{1\}$ and $D_L(x')=\{\ell+1\}$ imply $x'=w_{\ell+1,1}$. Since we have $x'\neq w_{\ell+1}$ (and $x'\ne 1$), we deduce $D_R(x')\not\subseteq \{1\}$. Choose $j\in D_R(x')\setminus\{1\}\subseteq \{2,\dots,n-1\}$. As $\vert\ell(x')-\ell(w_{\ell+1,1})\vert$ is odd, it follows that $\vert\ell(x')-\ell(w_{\ell+2,1})\vert$ is even. Since $x'\neq w_{\ell+2,1}$, by \ref{it: rabiotext 1}, \ref{it: rabiotext 2} of Lemma~\ref{rabiotext} we obtain $\mathrm{Ext}_{U(\fg)}^1(L(x'),L(w_{\ell+2,1}))=0$ and thus by Proposition~\ref{prop: Ext1 OS}
\[\mathrm{Ext}_{D(G)}^1(\cF_{P_{\widehat{\ell+1}}}^{G}(L(x'),\pi_{\ell+1,\ell+1}^{\infty})^\vee, C_{(1,\ell+2,\ell+2)}^\vee)=0.\]
Let $W\defeq \cF_{P_{\widehat{\ell+1}}}^{G}(L(x'),\pi_{\ell+1,\ell+1}^{\infty})$, then (\ref{equ: existence of length three}) forces the existence of a uniserial coadmissible $D(G)$-module $D$ of length $3$ with socle $C_{(1,\ell+2,\ell+2)}^\vee$, cosocle $W^\vee$, and middle layer $C_{(1,\ell+1,\ell+1)}^\vee$. Fix $\mu'\in \Lambda$ such that $\langle\mu'+\rho,\alpha^\vee\rangle\geq 0$ for $\al\in \Phi^+$ and the stabilizer of $\mu'$ in $W(G)$ for the dot action is $\{1,w_0s_jw_0\}$. As $j\ne 1$, by \ref{it: wall crossing simple 1} of Lemma~\ref{lem: wall crossing of simple} we have
\begin{equation*}
\Theta_{w_0s_jw_0}(C_{(1,\ell+2,\ell+2)}^\vee)=0=\Theta_{w_0s_jw_0}(C_{(1,\ell+1,\ell+1)}^\vee)
\end{equation*}
and by \ref{it: wall crossing simple 2} of Lemma~\ref{lem: wall crossing of simple} the canonical map $\Theta_{w_0s_jw_0}(W^\vee)\rightarrow W^\vee$ is surjective. Then the exactness of $\Theta_{\mu'}$ forces the canonical map $\Theta_{\mu'}(D)\rightarrow D$ to be a surjection, and we have $\Theta_{\mu'}(D)\cong \Theta_{\mu'}(W^\vee)$. As $x'\neq w_{\ell+2,1}$ and $\vert\ell(x')-\ell(w_{\ell+2,1})\vert$ is even, it follows from Proposition~\ref{prop: Jantzen middle} with \ref{it: rabiotext 1} of Lemma~\ref{rabiotext} and (\ref{equ: OS wall crossing}) that $C_{(1,\ell+2,\ell+2)}^\vee\notin \mathrm{JH}_{D(G)}(\Theta_{\mu'}(W^\vee))=\mathrm{JH}_{D(G)}(\Theta_{\mu'}(D))$, a contradiction to the surjection $\Theta_{\mu'}(D)\twoheadrightarrow D$.\bigskip

From Step $5$ and condition \ref{it: top deg 2}, we deduce that both $C$ and $D_{n-1}$ fit into a short exact sequence $0\rightarrow X_{n-1}=\mathrm{Fil}^0(D_{n-1})' \rightarrow \ast \rightarrow (\mathrm{St}_n^{\rm{alg}})^\vee\rightarrow 0$. Moreover $D_{n-1}$ is the unique such non-split extension by \ref{it: main Ext 0} of Theorem~\ref{thm: main Ext} (and the definition of $D_{n-1}$). But (\ref{equ: step 1 socle}) implies that $C$ is also a non-split extension, hence we finally deduce $C\cong D_{n-1}$.
\end{proof}

From now on we assume $E=K$. We consider the Drinfeld space $\bH_{/K}$ of dimension $n-1$ defined as in (\ref{defDrinfeld}), i.e.~$\bH\defeq \bP^{n-1}_{\rm{rig}}\setminus\bigcup_{\cH}\cH$ where $\cH$ runs through the $K$-rational hyperplanes inside $\bP^{n-1}_{\rm{rig}}$. Recall from \cite[Prop.~1.4]{SS91} that $\bH$ is quasi-Stein and that for any coherent sheaf $\cF$ on $\bH$ we have $H^k(\cF)=0$ for $k>0$. We follow the convention of \cite{ST02b} and \cite{Or08} and endow $\bH$ with the left action of $G=\GL_n(K)$ coming from the left action of $G$ on $\bP^{n-1}_{\rm{rig}}$ given by $(z_0,\dots,z_{n-1})\mapsto (z_0,\dots,z_{n-1})g^{-1}$ (matrix product) for $g\in G$. This action in turn induces an action of $G$ (which we will always take on the left) on the global sections $\cF(\bH)$ for any $G$-equivariant vector bundle $\cF$ on $\bH$. Moreover if $\cF$ is the restriction to $\bH$ of a $G$-equivariant vector bundle on $\bP^{n-1}_{\rm{rig}}$, then by the argument in \cite[p.593]{Or08} $\cF(\bH)$ is naturally a coadmissible $D(G)$-module. For instance this applies to the sheaf of differential $k$-forms on $\bH$.\bigskip

We consider the \emph{global sections} of the de Rham complex of $\bH$, which we denote (see (\ref{history}))
\begin{equation}\label{equ: Drinfeld dR}
\Omega^\bullet=[\Omega^0\longrightarrow \Omega^1\longrightarrow \cdots \longrightarrow \Omega^{n-1}].
\end{equation}
By the above discussion this is a complex of coadmissible $D(G)$-modules. The action of $G$ induces an action on each cohomology group $H^k(\Omega^\bullet)$ for $0\leq k\leq n-1$. We recall the following seminal result of Schneider-Stuhler \cite[Thm.~3.1, Lem.~4.1]{SS91} (see also \S\ref{history}).

\begin{thm}\label{thm: dR coh}
We have a $G$-equivariant isomorphism $H^k(\Omega^\bullet)\cong (V_{[1,n-k-1],\Delta}^{\infty})^\vee$ for $0\leq k\leq n-1$.
\end{thm}

Given a finite dimensional $U(\fp_I)$-module $X$ over $K$, we can lift it to a finite dimensional algebraic representation of $P_I$ (\cite[Lemma~3.2]{OS15}), and then consider the algebraic $G$-equivariant vector bundle over $G/P_I\cong \bP^{n-1}_{/K}$
\[\cF_X\defeq (G\times X)/P_I\]
where $p_I\in P_I$ acts on $(g,x)\in G\times X$ by $p_I(g,x) = (gp_I, p_I^{-1}x)$ (here we view $G$ and $P_I$ as algebraic groups over $K$ rather than their $K$-points). We also denote by $\cF_X$ the analytification of $\cF_X$, which is a $G$-equivariant vector bundle over $\bP^{n-1}_{\rm{rig}}$, hence by restriction over $\bH$, and we denote its global sections (a coadmissible $D(G)$-module) by
\[D_X\defeq \cF_X(\bH).\]
Note that, given a short exact sequence $0\rightarrow X_1\rightarrow X_2\rightarrow X_3\rightarrow 0$ of finite dimensional $U(\fp_I)$-modules, we have a short exact sequence $0\rightarrow \cF_{X_1}\rightarrow \cF_{X_2}\rightarrow \cF_{X_3}\rightarrow 0$ of $G$-equivariant vector bundles over $\bH$. As $\bH$ is quasi-Stein, by taking global sections we obtain a short exact sequence of coadmissible $D(G)$-modules
\begin{equation}\label{exactbundle}
0\rightarrow D_{X_1}\rightarrow D_{X_2}\rightarrow D_{X_3}\rightarrow 0.
\end{equation}

Now we recall Schneider's holomorphic discrete series (\cite{S92}). For $1\leq k\leq n-1$ we define $\mu_k\defeq w_{n-1,n-k}\cdot\mu_0\in \Lambda_{J}^{\dom}$ where we recall that $w_{n-1,n-k}=s_{n-1}s_{n-2}\cdots s_{n-k}$ (see (\ref{wj1j2})). For $\mu\in\Lambda_J^{\dom}$, we have $w_{1,n-1}(\mu)\in\Lambda_I^{\dom}$ (where $w_{1,n-1}=s_1s_2\cdots s_{n-1}$ and we see $L^I(w_{1,n-1}(\mu))$ as a finite dimensional algebraic representation of $P_I$ over $K$ via the inflation $P_I\twoheadrightarrow L_I$. For $\mu\in\Lambda_J^{\dom}$ we define the coadmissible $D(G)$-module
\begin{equation}\label{notationw0}
D_{\mu}\defeq D_{L^I(w_{1,n-1}(\mu))}.
\end{equation}
In \cite[\S 3]{S92}, Schneider defines a complex denoted there (see \cite[Lemma~9]{S92} and its proof):
\begin{equation}\label{equ: discrete series S}
[D_{\lambda(0)}\rightarrow D_{\lambda(1)}\rightarrow \cdots \rightarrow D_{\lambda(n-1)}].
\end{equation}
Unraveling the various conventions between \cite[\S 3]{S92} and this work, we can check that for $0\leq k\leq n-1$ we have $D_{\lambda(k)}=D_{L^I(w_{1,n-1}(\mu_k))}$. Hence by (\ref{notationw0}), (\ref{equ: discrete series S}) is a complex of coadmissible $D(G)$-modules
\begin{equation}\label{equ: discrete series}
\mathbf{D}_{\mu_0}^\bullet\defeq [D_{\mu_0}\longrightarrow D_{\mu_1}\longrightarrow \cdots \longrightarrow D_{\mu_{n-1}}].
\end{equation}

We recall the following results from \cite[\S 3]{S92}.

\begin{prop}\label{prop: discrete series}
\hspace{2em}
\begin{enumerate}[label=(\roman*)]
\item \label{it: discrete 1} For $0\leq k\leq n-1$ we have an isomorphism of coadmissible $D(G)$-modules $\Omega^k\cong D_{w_{n-1,n-k}\cdot 0}$ (with the convention $w_{n-1,n}=1$).
\item \label{it: discrete 2} For any finite dimensional $U(\fp_I)$-module $X$ over $K$ and any $\nu\in\Lambda^{\dom}$, we see $X\otimes_K L(\nu)$ as a (finite dimensional) $U(\fp_I)$-module via the diagonal action of $U(\fp_I)$, then we have an isomorphism of coadmissible $D(G)$-modules
\begin{equation*}
D_{X\otimes_K L(\nu)}\cong D_X\otimes_KL(\nu).
\end{equation*}
\end{enumerate}
\end{prop}
\begin{proof}
We check that \ref{it: discrete 1} (resp.~\ref{it: discrete 2}) is the translation of \cite[Prop.~1]{S92} (resp.~\cite[Lemma~5]{S92}) under our convention. Recall that the ``diagonal'' $D(G)$-action on $D_X\otimes_KL(\nu)$ is not completely straightforward, see \cite[\S 2.3.1]{JLS24}.
\end{proof}

\begin{lem}\label{lem: discrete series inf char}
For $\mu\in\Lambda_J^{\dom}$, $Z(\fg)$ acts on $D_{\mu}$ by $\xi_{\mu}$ where $\xi_{\mu}: Z(\fg)\rightarrow E$ is the unique infinitesimal character such that $L(\mu)_{\xi_\mu}\neq 0$.
\end{lem}
\begin{proof}
We use results from \cite{Schr11}, where $\rho$ there is denoted $-\delta$ and $\xi_{\mu}$ is denoted $\chi_{\mu-\delta}=\chi_{\mu+\rho}$. By the argument of \cite[Lemma~6.4]{Schr11} (which essentially follows from \cite[Prop.~8.22]{Kna01}) this implies that $Z(\fg)$ acts on $D_{\mu}$ (denoted $D_{w_0(\mu)}$ in \emph{loc.~cit.}) by $\chi_{w_0(\mu)+\delta}=\xi_{w_0(\mu)-2\rho}=\xi_{w_0\cdot\mu}=\xi_{\mu}$.
\end{proof}

We will need the two following (somewhat technical) lemmas.

\begin{lem}\label{lem: opposite filtration}
Let $L$ be a finite dimensional $U(\fg)$-module. Then there exist a \emph{decreasing} exhaustive separated filtration $(\mathrm{Fil}^{\ell}(L))_{\ell\in \Z}$ of $U(\fp_I)$-submodules on $L$, and an \emph{increasing} exhaustive separated filtration $(\mathrm{Fil}_{\ell}(L))_{\ell\in \Z}$ of $U(\fp_I^+)$-submodules on $L$, such that for each $\ell\in\Z$ the $U(\fl_I)$-modules $\mathrm{gr}^{\ell}(L)$ and $\mathrm{gr}_{\ell}(L)$ are either both zero or both simple and isomorphic.
\end{lem}
\begin{proof}
We consider the decomposition $L\vert_{U(\fl_I)}\cong \bigoplus_{\nu_1\in\Lambda_I^{\dom}}L_{\nu_1}$ where $L_{\nu_1}$ is the $L^I(\nu_1)$-isotypic component of $L\vert_{U(\fl_I)}$. We fix an arbitrary total order $\leq$ on $\Lambda_I^{\dom}$ such that $\nu_1',\nu_1''\in\Lambda_I^{\dom}$ satisfy $\nu_1'\leq \nu_1''$ only if $\nu_1'-\nu_1''\in\Z_{\geq 0}\Phi^+$. For $\nu_1\in \Lambda_I^{\dom}$ we define
\[\mathrm{Fil}^{\nu_1}(L)\defeq \bigoplus_{\nu_1'\geq \nu_1}L_{\nu_1'}\text{\ \ \ and\ \ \ }\mathrm{Fil}_{\nu_1}(L)\defeq \bigoplus_{\nu_1'\leq \nu_1}L_{\nu_1'}\]
which is a decreasing, respectively increasing, exhaustive separated filtration of $L\vert_{U(\fl_I)}$ (by $U(\fl_I)$-submodules) indexed by the totally ordered set $\Lambda_I^{\dom}$. The key observation (which is easy to check) is that $\mathrm{Fil}^{\nu_1}(L)$ (resp.~$\mathrm{Fil}_{\nu_1}(L)$) is $U(\fp_I)$-stable (resp.~$U(\fp_I^+)$)-stable in $L$ for each $\nu_1\in\Lambda_I^{\dom}$. Then we fix an arbitrary filtration on $L_{\nu_1}$ for each $\nu_1\in\Lambda_I^{\dom}$ with graded pieces being either simple or zero, and we further refine and reindex $(\mathrm{Fil}^{\nu_1}(L))_{\nu_1\in\Lambda_I^{\dom}}$ (resp.~$(\mathrm{Fil}_{\nu_1}(L))_{\nu_1\in\Lambda_I^{\dom}}$) to get a decreasing filtration $(\mathrm{Fil}^{\ell}(L))_{\ell\in \Z}$ by $U(\fp_I)$-submodules (resp. an increasing filtration $(\mathrm{Fil}_{\ell}(L))_{\ell\in \Z}$ by $U(\fp_I^+)$-submodules) as in the statement.
\end{proof}

Recall that, for $\lambda,\mu\in \Lambda$, $\cT^{\mu}_{\lambda}$ is the translation functor on the category of $Z(\fg)$-finite $D(G)$-modules, see above (\ref{equ: OS translation}).

\begin{lem}\label{lem: opposite Verma filtration}
Let $\nu\in\Lambda_J^{\dom}$ and $\lambda,\mu\in \Lambda$. Then there exist:
\begin{itemize}
\item a decreasing exhaustive separated filtration $(\mathrm{Fil}^{i}(\cT^{\mu}_{\lambda}(D_{\nu})))_{i\in \Z}$ on $\cT^{\mu}_{\lambda}(D_{\nu})$ by (clo\-sed) $D(G)$-submodules
\item an increasing exhaustive separated filtration $(\mathrm{Fil}_{i}(T^{\mu}_{\lambda}(M^J(\nu))))_{i\in \Z}$ on $T^{\mu}_{\lambda}(M^J(\nu))$ by $U(\fg)$-submodules
\end{itemize}
such that for $i\in\Z$ we have either $\mathrm{gr}^{i}(\cT^{\mu}_{\lambda}(D_{\nu}))=\mathrm{gr}_{i}(T^{\mu}_{\lambda}(M^J(\nu)))=0$, or $\mathrm{gr}^{i}(\cT^{\mu}_{\lambda}(D_{\nu}))\cong D_{\kappa_i}$ and \ \ $\mathrm{gr}_{i}(T^{\mu}_{\lambda}(M^J(\nu)))\cong M^J(\kappa_i)$ \ \ for \ \ some \ \ $\kappa_i\in\Lambda_J^{\dom}$. \ \ Moreover \ \ $\mathrm{gr}^{i}(\cT^{\mu}_{\lambda}(D_{\nu}))$ \ \ and $\mathrm{gr}_{i}(T^{\mu}_{\lambda}(M^J(\nu)))$ are non-zero only for finitely many $i\in \Z$.
\end{lem}
\begin{proof}
By Lemma~\ref{lem: discrete series inf char} we can assume $\xi_\nu=\xi_{\lambda}$ otherwise we have $\cT^{\mu}_{\lambda}(D_{\nu})=T^{\mu}_{\lambda}(M^J(\nu))=0$ by definition of $\cT^{\mu}_{\lambda}$ and there is nothing to prove. Let $L$ be the unique finite dimensional simple $U(\fg)$-module with highest weight in the $W(G)$-orbit of $\mu-\lambda$ (for the naive action of $W(G)$). By definition of $\cT^{\mu}_{\lambda}$ in \cite[(2)]{JLS24} and using \ref{it: discrete 2} of Proposition~\ref{prop: discrete series} we have
\begin{equation*}
\cT^{\mu}_{\lambda}(D_{\nu})=(D_{L^I(w_{1,n-1}(\nu))}\otimes_K L)_{\xi_{\mu}}\cong (D_{L^I(w_{1,n-1}(\nu))\otimes_K L})_{\xi_{\mu}}.
\end{equation*}
By Lemma \ref{lem: opposite filtration} we fix on $L$ a decreasing exhaustive separated filtration $(\mathrm{Fil}^{\ell}(L))_{\ell\in \Z}$ by $U(\fp_I)$-submodules, and an increasing exhaustive separated filtration $(\mathrm{Fil}_{\ell}(L))_{\ell\in \Z}$ by $U(\fp_I^+)$-submodules.\bigskip

\textbf{Step $1$}: We construct the decreasing filtration $(\mathrm{Fil}^{i}(\cT^{\mu}_{\lambda}(D_{\nu})))_{i\in \Z}$ on $\cT^{\mu}_{\lambda}(D_{\nu})$.\\
We \ \ first \ \ choose \ \ an \ \ arbitrary \ \ decreasing \ \ exhaustive \ \ separated \ \ filtration \ \ $(\mathrm{Fil}^{i}(L^I(w_{1,n-1}(\nu))\otimes_K L))_{i\in \Z}$ of $U(\fp_I)$-submodules on $L^I(w_{1,n-1}(\nu))\otimes_K L$ which refines the filtration $(L^I(w_{1,n-1}(\nu))\otimes_K \mathrm{Fil}^{\ell}(L))_{\ell\in \Z}$ and with each graded piece $\mathrm{gr}^{i}(L^I(w_{1,n-1}(\nu))\otimes_K L)$ either simple (of the form $L^I(w_{1,n-1}(\kappa_i))$ for some $\kappa_i\in\Lambda_J^{\dom}$) or zero. By (\ref{exactbundle}) the filtration $\mathrm{Fil}^{i}(L^I(w_{1,n-1}(\nu))\otimes_K L)$ induces a filtration $\mathrm{Fil}^{i}(D_{L^I(w_{1,n-1}(\nu))\otimes_K L})\defeq D_{\mathrm{Fil}^{i}(L^I(w_{1,n-1}(\nu))\otimes_K L)}$ by closed \ $D(G)$-submodules \ on \ $D_{L^I(w_{1,n-1}(\nu))\otimes_K L}$, \ which \ further \ induces \ a \ filtration \ $\mathrm{Fil}^{i}((D_{L^I(w_{1,n-1}(\nu))\otimes_K L})_{\xi_{\mu}})\defeq \mathrm{Fil}^{i}(D_{L^I(w_{1,n-1}(\nu))\otimes_K L})_{\xi_{\mu}}$ on $(D_{L^I(w_{1,n-1}(\nu))\otimes_K L})_{\xi_{\mu}}$. In particular, for $i\in\Z$, we have \[\mathrm{gr}^{i}((D_{L^I(w_{1,n-1}(\nu))\otimes_K L})_{\xi_{\mu}})=(D_{\mathrm{gr}^{i}(L^I(w_{1,n-1}(\nu))\otimes_K L)})_{\xi_{\mu}},\]
which is non-zero if and only if $\mathrm{gr}^{i}(L^I(w_{1,n-1}(\nu))\otimes_K L)\cong L^I(w_{1,n-1}(\kappa_i))$ for some $\kappa_i\in\Lambda_J^{\dom}$ such that $L(\kappa_i)_{\xi_{\mu}}\neq 0$ (using (\ref{notationw0}) and Lemma~\ref{lem: discrete series inf char}).\bigskip

\textbf{Step $2$}: We \ \ construct \ \ the \ \ increasing \ \ filtration \ \ $(\mathrm{Fil}_{i}(T^{\mu}_{\lambda}(M^J(\nu))))_{i\in \Z}$ \ \ on $T^{\mu}_{\lambda}(M^J(\nu))$.\\
We first note that $L^I(w_{1,n-1}(\nu))\otimes_K\mathrm{gr}^{\ell}(L)$ (resp.~$L^I(w_{1,n-1}(\nu))\otimes_K\mathrm{gr}_{\ell}(L)$) is a semi-simple $U(\fp_I)$-module (resp.~a semi-simple $U(\fp_I^+)$-module) and they are isomorphic as semi-simple $U(\fl_I)$-modules. Hence, it is possible to choose an increasing filtration $(\mathrm{Fil}_{i}(L^I(w_{1,n-1}(\nu))\otimes_{\! K} L))_{i\in \Z}$ of $U(\fp_I^+)$-modules on $L^I(w_{1,n-1}(\nu))\otimes_K L$ such that for each $i\in\Z$ we have an isomorphism of $U(\fl_I)$-modules $\mathrm{gr}_{i}(L^I(w_{1,n-1}(\nu))\otimes_K L)\cong \mathrm{gr}^{i}(L^I(w_{1,n-1}(\nu))\otimes_K L)$ (which by Step $1$ is either $0$ or $L^I(w_{1,n-1}(\kappa_i))$). For $\mu'\in \Lambda$ and $w\in W(G)$, we write $L^I(\mu')^w$ for the finite dimensional simple $U(w^{-1}\fl_Iw)$-module $L^I(\mu')$ where $x\in U(w^{-1}\fl_Iw)$ acts by $w x w^{-1}$.
Using $\fu_J=w_{1,n-1}^{-1}\fu_I w_{1,n-1}$ and $\ft=w_{1,n-1}^{-1}\ft w_{1,n-1}$, we have an isomorphism of $U(\ft)$-modules $((L^I(w_{1,n-1}(\nu))^{w_{1,n-1}})[\fu_J])_{\nu}=((L^I(w_{1,n-1}(\nu))[\fu_I])_{w_{1,n-1}(\nu)})^{w_{1,n-1}}\neq 0$. As $L^I(w_{1,n-1}(\nu))^{w_{1,n-1}}$ is a finite dimensional simple $U(w_{1,n-1}^{-1}\fl_I w_{1,n-1})=U(\fl_J)$-modules, this implies
\begin{equation}\label{equ: Levi twist}
L^I(w_{1,n-1}(\nu))^{w_{1,n-1}}\cong L^J(\nu)
\end{equation}
as $U(\fl_J)$-modules. Writing similarly $(L^I(w_{1,n-1}(\nu))\otimes_K L)^{w_{1,n-1}}$ for the $U(w_{1,n-1}^{-1}\fp_I^+ w_{1,n-1})=U(\fp_J)$-module on which $x\in U(\fp_J)$ acts by $w_{1,n-1} x w_{1,n-1}^{-1}$, we obtain the following isomorphisms of $U(w_{1,n-1}^{-1}\fp_I^+ w_{1,n-1})=U(\fp_J)$-modules
\begin{equation*}
(L^I(w_{1,n-1}(\nu))\otimes_K L)^{w_{1,n-1}}\cong L^I(w_{1,n-1}(\nu))^{w_{1,n-1}}\otimes_K L^{w_{1,n-1}}\cong L^J(\nu)\otimes_K L.
\end{equation*}
In \ particular \ we \ deduce \ from $(\mathrm{Fil}_{i}(L^I(w_{1,n-1}(\nu))\otimes_{\! K} L))_{i\in \Z}$ \ an \ increasing \ filtration $(\mathrm{Fil}_{i}(L^J(\nu)\otimes_K L))_{i\in \Z}$ of $U(\fp_J)$-modules on $L^J(\nu)\otimes_K L$ such that $\mathrm{gr}_{i}(L^J(\nu)\otimes_K L)$ is either $0$ or $L^J(\kappa_i)$. Using the tensor identity (cf.~\cite[\S 3.6]{Hum08})
\[U(\fg)\otimes_{U(\fp_J)}(L^J(\nu)\otimes_K L)\cong (U(\fg)\otimes_{U(\fp_J)}L^J(\nu))\otimes_K L=M^J(\nu)\otimes_K L\]
we obtain an increasing filtration of $U(\fg)$-submodules on $M^J(\nu)\otimes_K L$
\[(\mathrm{Fil}_{i}(M^J(\nu)\otimes_K L))_{i\in \Z}\defeq (U(\fg)\otimes_{U(\fp_J)}\mathrm{Fil}_{i}(L^J(\nu)\otimes_K L))_{i\in \Z}\]
such that $\mathrm{gr}_{i}(M^J(\nu)\otimes_K L)$ is either $0$ of $M^J(\kappa_i)$. By exactness of the translation functor $T^{\mu}_{\lambda}$ (see (\ref{equ: translation functor})), we deduce an increasing filtration $(\mathrm{Fil}_{i}(T^{\mu}_{\lambda}(M^J(\nu))))_{i\in \Z}\defeq (\mathrm{Fil}_{i}(M^J(\nu)\otimes_K L))_{\xi_{\mu}})_{i\in \Z}$ on $T^{\mu}_{\lambda}(M^J(\nu))=(M^J(\nu)\otimes_K L)_{\xi_{\mu}}$ which is non-zero if and only if $\mathrm{gr}_{i}(T^{\mu}_{\lambda}(M^J(\nu)))\cong M^J(\kappa_i)$ for some $\kappa_i\in\Lambda_J^{\dom}$ such that $(M^J(\kappa_i))_{\xi_\mu}\ne 0$, or equivalently $L(\kappa_i)_{\xi_{\mu}}\neq 0$. Comparing with the end of Step $1$, we obtain the statement.\bigskip

The last statement of the lemma follows from the fact $L^I(w_{1,n-1}(\nu))\otimes_K L$ is finite dimensional.
\end{proof}

\begin{rem}
Lemma \ref{lem: opposite Verma filtration} is applied in Lemma \ref{lem: dR to discrete} and Lemma \ref{lem: discrete wall crossing} below to study certain $\cT^{\mu}_{\lambda}(D_{\nu})$. Although more direct proofs probably exist without using Verma modules (since ultimately it is a matter of decomposing $L^I(w_{1,n-1}(\nu))\otimes_K L$, with the notation of the proof of Lemma \ref{lem: opposite Verma filtration}), the present proofs are more convenient for us as we can use standard results on the translation or wall-crossing of Verma modules.
\end{rem}

\begin{lem}\label{lem: dR to discrete}
For $0\leq k\leq n-1$ we have a canonical isomorphism of coadmissible $D(G)$-modules
\begin{equation}\label{equ: dR to discrete}
\cT_{w_0\cdot 0}^{w_0\cdot\mu_0}(\Omega^k) \cong D_{\mu_k}.
\end{equation}
\end{lem}
\begin{proof}
Let \ $\xi_0$ \ (resp. $\xi_{\mu_0}$) \ be \ the \ unique \ infinitesimal \ character \ such \ that \ $L(0)_{\xi_0}\neq 0$ (resp.~$L(\mu_0)_{\xi_{\mu_0}}\neq 0$). By \cite[\S 7.8]{Hum08} the pair of functors $T_{w_0\cdot 0}^{w_0\cdot\mu_0}$ and $T_{w_0\cdot\mu_0}^{w_0\cdot 0}$ (see (\ref{equ: translation functor})) are both left and right adjoint of each other, and gives an equivalence of categories $(\cO^{\fb}_{\rm{alg}})_{\xi_0}\cong (\cO^{\fb}_{\rm{alg}})_{\xi_{\mu_0}}$ (where $(\cO^{\fb}_{\rm{alg}})_{\xi}$ for an infinitesimal character $\xi$ is the subcategory of $\cO^{\fb}_{\rm{alg}}$ of $M$ such that $M=M_\xi$ (see above Lemma \ref{lem: central component} and recall that \cite[\S 7.8]{Hum08} uses antidominant weights, whence this $w_0$ everywhere). By \cite[Thm.~7.6]{Hum08} for $w\in W(G)$ we have $T_{w_0\cdot 0}^{w_0\cdot\mu_0}(M(w\cdot 0))\cong M(w\cdot\mu_0)$ and $T_{w_0\cdot\mu_0}^{w_0\cdot 0}(M(w\cdot\mu_0))\cong M(w\cdot 0)$. Since $\cO^{\fp_J}_{\rm{alg}}$ is stable under subobjects and tensoring by finite dimensional $U(\fg)$-modules, $T_{w_0\cdot 0}^{w_0\cdot\mu_0}$ and $T_{w_0\cdot\mu_0}^{w_0\cdot 0}$ preserve $\cO^{\fp_J}_{\rm{alg}}$, and thus induce an equivalence of categories $(\cO^{\fp_J}_{\rm{alg}})_{\xi_0}\cong (\cO^{\fp_J}_{\rm{alg}})_{\xi_{\mu_0}}$ (with obvious notation). Let $k\in \{0\dots, n-1\}$, since $M^J(w_{n-1,n-k}\cdot 0)$ (resp.~$M^J(w_{n-1,n-k}\cdot \mu_0)$) is the maximal quotient of $M(w_{n-1,n-k}\cdot 0)$ (resp.~$M(w_{n-1,n-k}\cdot \mu_0)$) in $\cO^{\fp_J}_{\rm{alg}}$ by \cite[Thm.~9.4(c)]{Hum08}, the isomorphism $T_{w_0\cdot 0}^{w_0\cdot\mu_0}(M(w_{n-1,n-k}\cdot 0))\cong M(w_{n-1,n-k}\cdot\mu_0)$ necessarily induces an isomorphism $T_{w_0\cdot 0}^{w_0\cdot\mu_0}(M^J(w_{n-1,n-k}\cdot 0))\cong M^J(w_{n-1,n-k}\cdot\mu_0)$. Since $M^J(\mu)$ for $\mu\in\Lambda_J^{\dom}$ has cosocle $L(\mu)$, any quotient of $M^J(w_{n-1,n-k}\cdot\mu_0)$ of the form $M^J(\mu)$ for some $\mu\in\Lambda_J^{\dom}$ is necessarily $M^J(w_{n-1,n-k}\cdot\mu_0)$ itself. In particular any increasing filtration on $T_{w_0\cdot 0}^{w_0\cdot\mu_0}(M^J(w_{n-1,n-k}\cdot 0))$ as in Lemma~\ref{lem: opposite Verma filtration} has only one non-zero graded piece, which is $M^J(w_{n-1,n-k}\cdot\mu_0)$. Then Lemma~\ref{lem: opposite Verma filtration} implies $\cT_{w_0\cdot 0}^{w_0\cdot\mu_0}(D_{w_{n-1,n-k}\cdot 0})\cong D_{\mu_k}$. Thus (\ref{equ: dR to discrete}) follows from \ref{it: discrete 1} of Proposition~\ref{prop: discrete series}.
\end{proof}

Applying $\cT_{w_0\cdot 0}^{w_0\cdot\mu_0}$ to the complex (\ref{equ: Drinfeld dR}) we deduce from Lemma \ref{lem: dR to discrete} an isomorphism of complexes of coadmissible $D(G)$-modules:
\[\cT_{w_0\cdot 0}^{w_0\cdot\mu_0}(\Omega^\bullet)\cong \mathbf{D}_{\mu_0}^\bullet.\]
This is essentially equivalent to what Schneider does in \cite[Thm.~3]{S92}: though there is no mention of infinitesimal characters there, in \cite[p.643]{S92} he first defines a decreasing filtration on the complex $\Omega^\bullet\otimes_K L(1)$ (note that $L(1)=L(\mu_0)$ is the unique finite dimensional simple $U(\fg)$-module with highest weight in the $W(G)$-orbit of $w_0\cdot\mu_0-w_0\cdot 0=\mu_0$) and then in the proof of \cite[Thm.~3]{S92} projects onto the unique graded piece with infinitesimal character $\xi_{\mu_0}$. As $\cT_{w_0\cdot 0}^{w_0\cdot\mu_0}$ is exact, we deduce from Theorem \ref{thm: dR coh} for $0\leq k\leq n-1$ (with the notation of (\ref{vIalg})
\begin{equation}\label{equ: coh translation}
H^k(\mathbf{D}_{\mu_0}^\bullet)\cong T_{w_0\cdot 0}^{w_0\cdot\mu_0}(H^k(\Omega^\bullet))\cong L(1)\otimes_K (V_{[1,n-k-1],\Delta}^{\infty})^\vee=(V_{[1,n-k-1],\Delta}^{\rm{alg}})^\vee.
\end{equation}
We need the following result on the wall-crossing of holomorphic discrete series.

\begin{lem}\label{lem: discrete wall crossing}
Let $1\leq k\leq n-1$ and $\mu\in \Lambda$ such that $\langle\mu+\rho,\alpha^\vee\rangle\geq 0$ for $\al\in \Phi^+$ and the stabilizer of $\mu$ in $W(G)$ for the dot action is $\{1,s_k\}$. We have a short exact sequence of coadmissible $D(G)$-modules
\begin{equation}\label{equ: discrete wall crossing}
0\longrightarrow D_{\mu_k}\longrightarrow \Theta_{\mu}(D_{\mu_k}) \longrightarrow D_{\mu_{k-1}}\longrightarrow 0.
\end{equation}
\end{lem}
\begin{proof}
Recall from above (\ref{equ: OS wall crossing}) that $\Theta_{\mu}(D_{\mu_k})=\cT_{\mu}^{w_0\cdot\mu_0}(\cT_{w_0\cdot\mu_0}^{\mu}(D_{\mu_k}))$. By Lemma~\ref{lem: opposite Verma filtration} applied with $\nu=\mu_k$ and $\lambda=w_0\cdot\mu_0$, we have a decreasing exhaustive separated filtration $(\mathrm{Fil}^{i}(D))_{i\in \Z}$ on $D\defeq \cT_{w_0\cdot\mu_0}^{\mu}(D_{\mu_k})$, and an increasing exhaustive separated filtration $(\mathrm{Fil}_{i}(M))_{i\in \Z}$ on $M\defeq T^{\mu}_{w_0\cdot\mu_0}(M^J(\mu_k))$ such that for $i\in\Z$ we have either $\mathrm{gr}^{i}(D)=\mathrm{gr}_{i}(M)=0$, or $\mathrm{gr}^{i}(D)\cong D_{\kappa_i}$ and $\mathrm{gr}_{i}(M)\cong M^J(\kappa_i)$ for some $\kappa_i\in\Lambda_J^{\dom}$. Then we apply Lemma~\ref{lem: opposite Verma filtration} again to $\cT_{\mu}^{w_0\cdot\mu_0}(\mathrm{gr}^{i}(D))$ (i.e.~with $\nu=\kappa_i$, $\lambda=\mu$, $\mu=w_0\cdot\mu_0$) for $i\in\Z$ such that $0\neq \mathrm{gr}^{i}(D)\cong D_{\kappa_i}$, and obtain a decreasing exhaustive separated filtration $(\mathrm{Fil}^{j}(\Theta_{\mu}(D_{\mu_k})))_{j\in \Z}$ on $\Theta_{\mu}(D_{\mu_k})=\cT_{\mu}^{w_0\cdot\mu_0}(D)$ which refines the decreasing filtration $(\cT_{\mu}^{w_0\cdot\mu_0}(\mathrm{Fil}^{i}(D)))_{i\in \Z}$, and an increasing exhaustive separated filtration $(\mathrm{Fil}_{j}(\Theta_{s_k}(M^J(\mu_k))))_{j\in \Z}$ on $\Theta_{s_k}(M^J(\mu_k))= T^{w_0\cdot\mu_0}_{\mu}(M)$ which refines the increasing filtration $(T^{w_0\cdot\mu_0}_{\mu}(\mathrm{Fil}_{i}(M)))_{i\in \Z}$) such that for each $j\in\Z$ we have either $\mathrm{gr}^{j}(\Theta_{\mu}(D_{\mu_k}))=\mathrm{gr}_{j}(\Theta_{s_k}(M^J(\mu_k)))=0$, or $\mathrm{gr}^{j}(\Theta_{\mu}(D_{\mu_k}))\cong D_{\theta_j}$ and $\mathrm{gr}_{j}(\Theta_{s_k}(M^J(\mu_k)))\cong M^J(\theta_j)$ for some $\theta_j\in\Lambda_J^{\dom}$. However, by \ref{it: special Verma 2} of Lemma~\ref{lem: special Verma wall crossing} and up to some reindexation, the only such filtration on $\Theta_{s_k}(M^J(\mu_k))$ is
\begin{multline*}
0=\mathrm{Fil}_{0}(\Theta_{s_k}(M^J(\mu_k)))\subsetneq \mathrm{Fil}_{1}(\Theta_{s_k}(M^J(\mu_k)))=M^J(\mu_{k-1})\subsetneq \mathrm{Fil}_{2}(\Theta_{s_k}(M^J(\mu_k)))\\
=\Theta_{s_k}(M^J(\mu_k))
\end{multline*}
with \ $\mathrm{gr}_1(\Theta_{s_k}(M^J(\mu_k)))\cong M^J(\mu_{k-1})$ \ and \ $\mathrm{gr}_2(\Theta_{s_k}(M^J(\mu_k)))\cong M^J(\mu_k)$. \ This \ forces $\Theta_{\mu}(D_{\mu_k})$ to admit a filtration
\[\Theta_{\mu}(D_{\mu_k})=\mathrm{Fil}^{1}(\Theta_{\mu}(D_{\mu_k}))\supsetneq \mathrm{Fil}^{2}(\Theta_{\mu}(D_{\mu_k}))=D_{\mu_k}\supsetneq \mathrm{Fil}^{3}(\Theta_{\mu}(D_{\mu_k}))=0\]
with $\mathrm{gr}^1(\Theta_{\mu}(D_{\mu_k}))\cong D_{\mu_{k-1}}$ and $\mathrm{gr}^2(\Theta_{\mu}(D_{\mu_k}))\cong D_{\mu_k}$. Hence we have (\ref{equ: discrete wall crossing}).
\end{proof}

\begin{rem}
We will see in the proof of Theorem \ref{thm: main dR} below that the injection $D_{\mu_k}\rightarrow \Theta_{\mu}(D_{\mu_k})$ in Lemma \ref{lem: discrete wall crossing} is the canonical adjunction map.
\end{rem}

We now slightly reformulate (a weak form of) \cite[Thm.~2.2]{Or15}.

\begin{thm}\label{thm: dR complex filtration}
Let $\mu\in\Lambda_J^{\dom}$. Then $D_{\mu}$ in (\ref{notationw0}) is a finite length coadmissible $D(G)$-module which admits a decreasing filtration
\[D_{\mu}=\mathrm{Fil}^0(D_{\mu})\supseteq \mathrm{Fil}^1(D_{\mu})\supseteq \cdots \supseteq \mathrm{Fil}^{n-1}(D_{\mu})\]
satisfying the following conditions:
\begin{enumerate}[label=(\roman*)]
\item \label{it: Orlik 1} $\mathrm{gr}^{n-1}(D_{\mu})\neq 0$ if and only if $\mu\in\Lambda^{\dom}$, in which case $\mathrm{gr}^{n-1}(D_{\mu})\cong L(\mu)$;
\item \label{it: Orlik 2} for $0\leq \ell\leq n-2$, there exists $M_{\ell}$ in $\cO^{\fp_{\widehat{\ell+1}}}_{\rm{alg}}$ such that $\mathrm{gr}^{\ell}(D_{\mu})$ fits into a short exact sequence of coadmissible $D(G)$-modules if $w_{\ell+1,n-1}\cdot\mu\in\Lambda^{\dom}$
\[0\rightarrow \cF_{P_{\widehat{\ell+1}}}^{G}(M_{\ell},\pi_{\ell+1,\ell+1}^{\infty})^\vee \rightarrow \mathrm{gr}^{\ell}(D_{\mu}) \rightarrow L(w_{\ell+1,n-1}\cdot\mu)\otimes_K(V_{[1,\ell],\Delta}^{\infty})^\vee \rightarrow 0,\]
and $\mathrm{gr}^{\ell}(D_{\mu})\cong \cF_{P_{\widehat{\ell+1}}}^{G}(M_{\ell},\pi_{\ell+1,\ell+1}^{\infty})^\vee$ otherwise.
\end{enumerate}
\end{thm}
\begin{proof}
Recall from (\ref{notationw0}) that $D_{\mu}=\cF_{\mu}(\bH)$ where $\cF_{\mu}$ is the analytification on $\bP^{n-1}_{\rm{rig}}$ of the $G$-equivariant vector bundle $\cF_{L^I(w_{1,n-1}(\mu))}$ on $\bP^{n-1}$. It follows from \cite[Thm.~6.1]{GS69} (see also \cite[Thm.~6.4]{Kos61}) and from the definition of $\Lambda^{\dom}$ (using $w_0=w_I w_{1,n-1}$) 
that $H^{n-1-\ell}(\bP^{n-1},\cF_{L^I(w_{1,n-1}(\mu))})\neq 0$ \ if \ and \ only \ if \ $w_{\ell+1,n-1}\cdot\mu\in\Lambda^{\dom}$, \ in \ which \ case $H^{n-1-\ell}(\bP^{n-1},\cF_{L^I(w_{1,n-1}(\mu))})\cong L(w_{\ell+1,n-1}\cdot\mu)$. Also recall that we have for $0\leq \ell\leq n-2$ (see the line below (\ref{pij_1,_2}))
\[\pi_{\ell+1,\ell+1}^{\infty}\cong 1_{\ell+1}\boxtimes_E \mathrm{St}_{n-1-\ell}^{\infty}=V_{[1,\ell],\widehat {\ell+1}}^{\infty}.\]
Then \ref{it: Orlik 1} and \ref{it: Orlik 2} follow from \cite[Lemma~2.1]{Or15}, \cite[Thm.~2.2]{Or15} (based on \cite[Thm.~1]{Or08}) and the discussion that follows \emph{loc.~cit.} where for $0\leq \ell\leq n-2$ (see \cite[(7)]{Or15} and the notation there)
\[M_{\ell}\defeq \mathrm{ker}\left(H^{n-1-\ell}_{\bP^{\ell}}(\bP^{n-1},\cF_{L^I(w_{1,n-1}(\mu))})\\
\longrightarrow H^{n-1-\ell}(\bP^{n-1},\cF_{L^I(w_{1,n-1}(\mu))})\right).\qedhere\]
\end{proof}

Recall that $J=\widehat{n-1}$ and $\mu_k=w_{n-1,n-k}\cdot\mu_0$ for $0\leq k\leq n-1$. We will need the following technical result.

\begin{lem}\label{lem: dR complex N coh}
For $0\leq k\leq n-1$ (and $\mu_0\in\Lambda^{\dom}$) we have a $D(L_J)$-equivariant isomorphism
\begin{equation}\label{equ: H0 dR complex}
H^0(N_J, D_{\mu_k})\cong L^J(\mu_k) =L^J(w_{n-1,n-k}).
\end{equation}
\end{lem}
\begin{proof}
Note first that we have an equality of closed vector subspaces of $D_{\mu_k}$
\begin{equation}\label{equ: dR N inv}
H^0(N_J, D_{\mu_k})=H^0(N_J, H^0(\fn_J,D_{\mu_k})).
\end{equation}
Recall that $\bH$ is by definition contained in the rigid analytic space associated with the Zariski open \emph{affine} subscheme $P_Jw_0P_I/P_I$ of $G/P_I=\bP^{n-1}_{/K}$. Using $P_J=N_JL_J$, $w_0L_Jw_0=L_I$ and $w_0=w_{n-1,1}w_I$, we have (as affine schemes)
\begin{multline}\label{equ: simple Bruhat}
P_Jw_0P_I=N_JL_Jw_0P_I=N_Jw_0(w_0L_Jw_0)P_I=N_Jw_0P_I\\
=N_Jw_{n-1,1}P_I\cong N_J\times(w_{n-1,1}P_I).
\end{multline} 
Let $L^I\defeq L^I(w_{1,n-1}(\mu_k))$, the restriction $\cF_{L^I(w_{1,n-1}(\mu_k))}|_{P_Jw_0P_I/P_I}=\cF_{L^I}|_{P_Jw_0P_I/P_I}$ is by definition the quotient $(P_Jw_0P_I\times L^I)/P_I$ where $P_I$ acts on $P_Jw_0P_I\times L^I$ by
\[(P_Jw_0P_I\times L^I)\times P_I\rightarrow P_Jw_0P_I\times L^I, \ ((x,v),h)\mapsto (xh,h^{-1}\cdot v).\]
Using (\ref{equ: simple Bruhat}), we write each $x\in P_Jw_0P_I$ as $x=aw_{n-1,1}b$ for some $a\in N_J$ and $b\in P_I$ (uniquely determined by $x$). We consider the following isomorphism
\[P_Jw_0P_I\times L^I\buildrel\sim\over\longrightarrow P_Jw_0P_I\times L^I, \ (aw_{n-1,1}b, v)\mapsto (aw_{n-1,1}b, b\cdot v)\]
which descends to an isomorphism
\begin{equation}\label{equ: trivialization quotient}
(P_Jw_0P_I\times L^I)/P_I\buildrel\sim\over\longrightarrow (P_Jw_0P_I)/P_I \times L^I, \ (aw_{n-1,1}b, v)P_I \mapsto (aw_{n-1,1}bP_I, b\cdot v).
\end{equation}
We have a (left) $P_J$-action on $(P_Jw_0P_I\times L^I)/P_I$ given by
\begin{equation}\label{equ: first action}
P_J\times ((P_Jw_0P_I\times L^I)/P_I)\rightarrow (P_Jw_0P_I\times L^I)/P_I, \ (g, (x,v)P_I)\mapsto (gx,v)P_I.
\end{equation}
For $g=g'g''\in N_JL_J=P_J$ and $x=aw_{n-1,1}b\in N_Jw_{n-1,1}P_I=P_Jw_0P_I$, we have \[gx=g'(g''a(g'')^{-1})w_{n-1,1}(w_{1,n-1}g''w_{1,n-1}^{-1})b\]
with $g'(g''a(g'')^{-1})\in N_J$ and $(w_{1,n-1}g''w_{1,n-1}^{-1})\in L_I$. Hence, the $P_J$-action (\ref{equ: first action}) translates via (\ref{equ: trivialization quotient}) into the $P_J$-action (with the above notation)
\begin{multline}\label{equ: second action}
P_J\times ((P_Jw_0P_I)/P_I \times L^I)\rightarrow (P_Jw_0P_I)/P_I \times L^I, \ (g, (xP_I, v))\\
\mapsto (gxP_I, (w_{1,n-1}g''w_{1,n-1}^{-1})\cdot v).
\end{multline} 
Recall from the notation in \S\ref{generalnotation} that $(L^I)^{w_{1,n-1}}$ is the $L_J=w_{1,n-1}^{-1}L_Iw_{1,n-1}$-representation with same underlying space as $L^I$ and $g''\in L_J$ acting as $w_{1,n-1}g''w_{1,n-1}^{-1}$ on $L^I$. It follows from the above discussion that (\ref{equ: trivialization quotient}) gives a $P_J$-equivariant isomorphism of $P_J$-equivariant algebraic vector bundles on $P_Jw_0P_I/P_I$ (with the $P_J$-action (\ref{equ: first action}) on the left and the $P_J$-action (\ref{equ: second action}) on the right):
\begin{equation*}
\cF_{L^I}|_{P_Jw_0P_I/P_I}\cong \cO_{P_Jw_0P_I/P_I}\times (L^I)^{w_{1,n-1}}.
\end{equation*}
As in (\ref{equ: Levi twist}), we have $(L^I)^{w_{1,n-1}}\cong L^J(\mu_k)$ as $L_J$-representations. Taking rigid analytic sections over $\bH$, we obtain a $P_J$-equivariant topological isomorphism of Fr\'echet spaces
\begin{equation}\label{equ: trivialization open}
D_{\mu_k}\buildrel (\ref{notationw0}) \over = \cF_{L^I}(\bH)\cong \cO(\bH)\otimes_K L^J(\mu_k)=\Omega^0\otimes_K L^J(\mu_k).
\end{equation}
The $N_J$-equivariant isomorphism of schemes $P_Jw_0P_I/P_I=N_Jw_0P_I/P_I\buildrel\sim\over\longrightarrow N_J$ induces a $N_J$-equivariant isomorphism between their algebraic de Rham complexes $\Omega^\bullet_{P_Jw_0P_I/P_I}\cong \Omega^\bullet_{N_J}$. As $N_J$ is the affine space $\mathbb{A}^{\dim N_J}_{/K}$, $\Omega^\bullet_{N_J}$ is naturally identified with the Chevalley-Eilenberg complex $\Hom_K(\wedge^\bullet\fn_J,\cO_{N_J})\cong \Hom_K(\wedge^\bullet\fn_J,\cO_{P_Jw_0P_I/P_I})$. Taking rigid analytic sections over $\bH$, we can identify the de Rham complex (\ref{equ: Drinfeld dR}) with $\Hom_K(\wedge^\bullet\fn_J,\Omega^0)$, which is the Chevalley-Eilenberg complex of $\Omega^0$ (see also \cite[p.635]{S92} for the isomorphism $\Omega^j\cong \Hom_K(\wedge^j\fn_J,\Omega^0)$ for $0\leq j\leq n-1$). In particular, we obtain an equality of closed vector subspaces of $\Omega^0$
\[H^0(\fn_J,\Omega^0)=H^0(\Omega^\bullet)\cong 1_{G}^\vee\]
(where the second isomorphism follows from Theorem \ref{thm: dR coh}). With (\ref{equ: trivialization open}) (and as $\fn_J$ acts trivially on $(L^I)^{w_{1,n-1}}\cong L^J(\mu_k)$) this gives a $P_J$-equivariant isomorphism of Fr\'echet spaces
\begin{equation*}
H^0(\fn_J,D_{\mu_k})\cong H^0(\fn_J,\Omega^0)\otimes_K L^J(\mu_k)=L^J(\mu_k)
\end{equation*}
with $N_J$ acting trivially on $L^J(\mu_k)$. This together with (\ref{equ: dR N inv}) gives (\ref{equ: H0 dR complex}).
\end{proof}

We actually only use Theorem \ref{thm: dR complex filtration} for $\mu=\mu_{n-1}$, and likewise Lemma \ref{lem: dR complex N coh} for $k=n-1$, in order to prove the key result that follows.

\begin{thm}\label{thm: main dR}
For $0\leq k\leq n-1$ we have an isomorphism of finite length coadmissible $D(G)$-modules
\begin{equation}\label{equ: main dR}
D_k\cong D_{\mu_k}.
\end{equation}
Moreover, for $0\leq k\leq n-2$, the differential map $D_{\mu_k}\rightarrow D_{\mu_{k+1}}$ is the unique (up to scalar) non-zero map $D_k\rightarrow D_{k+1}$ with image $X_{k+1}$.
\end{thm}
\begin{proof}
We prove the first statement. Note first that $D_{\mu_{n-1}}=D_{w_{n-1,1}\cdot\mu_0}$ is a finite length coadmissible $D(G)$-module by the first statement of Theorem~\ref{thm: dR complex filtration}. Condition \ref{it: top deg 1} of Theorem~\ref{prop: top deg} holds for $C=D_{\mu_{n-1}}$ by Lemma~\ref{lem: dR complex N coh}. Condition \ref{it: top deg 2} of Theorem~\ref{prop: top deg} holds for $C=D_{\mu_{n-1}}$ by Theorem~\ref{thm: dR complex filtration}, noting that $w_{\ell+1,n-1}\cdot\mu_{n-1}\in\Lambda^{\dom}$ if and only if $\ell=0$ (in which case it is $\mu_0$). Hence $D_{n-1}\cong D_{\mu_{n-1}}$ by Theorem~\ref{prop: top deg}. One easily checks from (\ref{definitionofxyz}) and (\ref{definitionofdk}) that the (irreducible) cosocle of ${D}_{n-1}$ does not appear in $D_{n-2}$, hence $\Hom_G(D_{n-1},D_{n-2})=0$. Then it follows from \ref{it: D wall crossing 2} of Theorem~\ref{thm: D wall crossing} (for $k=n-1$) that there is (up to a non-zero scalar) a unique injection $D_{n-1}\hookrightarrow \Theta_{\mu}(D_{n-1})$ (with the notation of \emph{loc.~cit.}). In particular this must be the injection in (\ref{equ: discrete wall crossing}) for $k=n-1$. Then it follows from \ref{it: D wall crossing 2} of Theorem~\ref{thm: D wall crossing} and (\ref{equ: discrete wall crossing}) again (both applied with $k=n-1$) that we have $D_{n-2}\cong D_{\mu_{n-2}}$. Since we again have $\Hom_G(D_{n-2},D_{n-3})=0$, the same argument for $k=n-2$ instead of $k=n-1$ gives $D_{n-3}\cong D_{\mu_{n-3}}$. By descending induction we deduce (\ref{equ: main dR}) for all $k$.

We prove the second statement. For $0\leq k\leq n-2$, as $H^{k+1}(\mathbf{D}_{\mu_0}^\bullet)\cong(V_{[1,n-k-2],\Delta}^{\rm{alg}})^\vee$ (see (\ref{equ: coh translation})) and $(V_{[1,n-k-2],\Delta}^{\rm{alg}})^\vee$ does not show up in $\mathrm{soc}_{D(G)}(D_{\mu_{k+1}})\cong\mathrm{soc}_{D(G)}(D_{k+1})$ (using (\ref{definitionofdk})), we deduce that the differential map $D_{\mu_k}\rightarrow D_{\mu_{k+1}}$ is non-zero. Therefore it must be the unique (up to scalar) non-zero map $D_k\rightarrow D_{k+1}$ with image $X_{k+1}$ (see the discussion above (\ref{equ: explicit dR})).
\end{proof}

\begin{rem}
By the definition of $D_{n-1}$ below (\ref{equ: rough shape of rep}), we have $\mathrm{cosoc}_{D(G)}(D_{n-1})=(\mathrm{St}_n^{\rm{alg}})^\vee$. When $\mu_0=(0,\cdots,0)$, we have $D_{n-1}\cong \Omega^{n-1}$ by Theorem~\ref{thm: main dR} and \ref{it: discrete 1} of Proposition~\ref{prop: discrete series}. Thus we obtain $\mathrm{cosoc}_{D(G)}(\Omega^{n-1})=(\mathrm{St}_n^{\infty})^\vee$. Then, by an argument parallel to \cite[Cor.~6.11]{Schr11} using \cite{IS01}, one can expect to deduce from this that $\Omega^{n-1}$ is a quotient of the dual of the locally $K$-analytic Steinberg $(\mathrm{St}_n^{\rm{an}})^\vee$, and thus a subquotient of $\cF_{B}^{G}(U(\fg)\otimes_{U(\fb)}0,1_T)^\vee$. Using Theorem~\ref{thm: main dR}, Lemma~\ref{lem: dR to discrete} and (\ref{equ: OS translation}), it would follow that $D_{n-1}\cong D_{\mu_{n-1}}\cong \cT_{w\cdot 0}^{w\cdot\mu_0}(\Omega^{n-1})$ is a subquotient of $\cT_{w\cdot 0}^{w\cdot\mu_0}(\cF_{B}^{G}(U(\fg)\otimes_{U(\fb)}0,1_T)^\vee)\cong \cF_{B}^{G}(U(\fg)\otimes_{U(\fb)}\mu_0,1_T)^\vee$, which by Lemma~\ref{lem: discrete wall crossing} (or \ref{it: D wall crossing 2} of Theorem~\ref{thm: D wall crossing}) would imply that $D_k\cong D_{\mu_k}$ is a subquotient of
\begin{equation*}
\Theta_{s_{k+1}}\cdots\Theta_{s_{n-1}}(\cF_{B}^{G}(U(\fg)\otimes_{U(\fb)}\mu_0,1_T)^\vee)\cong \cF_{B}^{G}(Q_k,1_T)
\end{equation*}
where $Q_k\defeq \Theta_{s_{k+1}}\cdots\Theta_{s_{n-1}}(U(\fg)\otimes_{U(\fb)}\mu_0)\in \cO^{\fb}_{\rm{alg}}$.
\end{rem}

Recall that we have defined an explicit complex of finite length coadmissible $D(G)$-modules with Orlik-Strauch constituents $\tld{\mathbf{D}}^\bullet$ in (\ref{equ: explicit dR split top}).

\begin{cor}
For $\mu_0\in\Lambda^{\dom}$ we have two morphisms of complexes of finite length coadmissible $D(G)$-modules with Orlik-Strauch constituents
\begin{equation*}
H^{n-1}(\mathbf{D}_{\mu_0}^\bullet)[-(n-1)]\longleftarrow\tld{\mathbf{D}}^\bullet\longrightarrow \mathbf{D}_{\mu_0}^\bullet
\end{equation*}
which give an explicit section to the morphism of complexes $\mathbf{D}_{\mu_0}^\bullet\twoheadrightarrow H^{n-1}(\mathbf{D}_{\mu_0}^\bullet)[-(n-1)]$ in the derived category of finite length coadmissible $D(G)$-modules with Orlik-Strauch constituents.
\end{cor}
\begin{proof}
This follows directly from (\ref{equ: explicit dR}), Theorem~\ref{thm: main split}, (\ref{equ: discrete series}) and Theorem~\ref{thm: main dR}.
\end{proof}

One also has statements analogous to Proposition \ref{prop: little split} and Corollary \ref{cor: split3} with $\mathbf{D}_{\mu_0}^\bullet$ instead of $\mathbf{D}^\bullet$. In particular, if we define $\tld{\Omega}^\bullet\defeq \tld{\mathbf{D}}^\bullet$ when $\mu_0=(0,\cdots,0)$, we at last obtain one of our main results.

\begin{cor}\label{cor: split dR}
We have two morphisms of complexes of finite length coadmissible $D(G)$-modules with Orlik-Strauch constituents.
\begin{equation*}
H^{n-1}(\Omega^\bullet)[-(n-1)]\longleftarrow \tld{\Omega}^\bullet\longrightarrow \Omega^\bullet
\end{equation*}
which give an explicit section to the morphism of complexes $\Omega^\bullet\twoheadrightarrow H^{n-1}(\Omega^\bullet)[-(n-1)]$ in the derived category of finite length coadmissible $D(G)$-modules with Orlik-Strauch constituents.
\end{cor}

\begin{rem}
Assume $\mu_0=(0,\cdots,0)$ and, changing notation, denote by $\tld{\Omega}^{n-1}$ the unique coadmissible $D(G)$-module with cosocle $((\mathrm{St}_n^{\infty})^{\oplus n})^\vee$ which sits in a short exact sequence
\[0\longrightarrow Z_{n-1}\longrightarrow \tld{\Omega}^{n-1} \longrightarrow \big((\mathrm{St}_n^{\infty})^{\oplus n}\big)^\vee \longrightarrow 0\]
(see (\ref{equ: universal first branch})). Then the complex $[\widetilde \Omega^0\rightarrow \widetilde \Omega^1\rightarrow \cdots \rightarrow \widetilde \Omega^{n-2} \rightarrow \widetilde \Omega^{n-1}]$ (the same as in (\ref{equ: explicit dR split top}) except that we have modified the term in degree $n-1$) is canonical, exact in degrees $<n-1$ and its $H^{n-1}$ is $((\mathrm{St}_n^{\infty})^{\oplus n})^\vee$. That is, we have ``made exact'' the de Rham complex $\Omega^\bullet$ in degrees $<n-1$ at the expense of replacing $H^{n-1}(\Omega^\bullet)=(\mathrm{St}_n^{\infty})^\vee$ by $((\mathrm{St}_n^{\infty})^{\oplus n})^\vee$. These properties look similar to the properties of the complexes obtained as $\chi$-isotypic direct factors of the global sections of the de Rham complex of the first covering of $\bH$, where $\chi$ is a smooth character of $\cO_D^\times$ and $D$ is the division algebra over $K$ of invariant $1/n$, see for instance Junger's work \cite{Jun24}. Finally, note that, by \ref{it: D wall crossing 4} of Theorem \ref{thm: D wall crossing}, with this modified complex \ref{it: D wall crossing 3} of Theorem \ref{thm: D wall crossing} now also holds when $k=n-1$.
\end{rem}

\begin{rem}
When $n=2$ and $K=\Qp$, the study of the global sections of the de Rham complex of Drinfeld's coverings $\Sigma_n$ ($n\geq 1$) of $\bH$ is very important for the $p$-adic Langlands program for $\GL_2(\Qp)$, see for instance the deep results of \cite{DLB17} or \cite{Pa25}. Moreover, when $n=2$ and $K$ is arbitrary, at least for the first covering $\Sigma_1$ we now know from \cite[Cor.~B]{AW23}, \cite[Thm.~7.2.1, Rem.~7.2.4]{PSS19} that $H^0(\Sigma_1,\mathcal{O}_{\Sigma_1})$ and $H^0(\Sigma_1,\Omega^1_{\Sigma_1})$ are finite length coadmissible $D(G)$-modules (or rather $D(G^0)$-modules here, where $G^0\subset G$ is the subgroup of matrices with determinant in $\cO_K^\times$). It is therefore a natural question to ask for the internal structure of these $D(G^0)$-modules, and more generally of all the $H^0(\Sigma_n,\Omega^j_{\Sigma_n})$, $n\geq 1$, $j\geq 0$. However, it is unlikely that one can apply (at least directly) the techniques of this work which are adapted to Orlik-Strauch representations, i.e.~subquotients of locally anaytic parabolic inductions of locally algebraic representations.
\end{rem}

\newpage

\appendix

\section{Combinatorics in \texorpdfstring{$W(G)$}{W(G)}}\label{sec: appendix}

We prove several technical combinatorial lemmas involving specific elements of $W(G)$ that are used throughout this work.\bigskip

We equip the set $\Z\times\Z$ with the partial order $(a,b)\leq (a',b')$ if and only if $a\leq a'$ and $b\leq b'$. Let $\Sigma\subseteq \Z\times\Z$ be a finite subset such that $1\leq a-b\leq n-1$ for each $(a,b)\in\Sigma$. For $k\in\Z$, we set
\begin{equation}\label{equ: expansion layer}
x_{\Sigma,k}\defeq \prod_{(a,b)\in\Sigma, a+b=k}s_{a-b}.
\end{equation}
One easily checks that $s_{a-b}s_{a'-b'}=s_{a'-b'}s_{a-b}$ for $(a,b),(a',b')\in\Sigma$ satisfying $a+b=a'+b'=k$, so the definition of $x_{\Sigma,k}$ is independent of the order of the $s_{a-b}$ in (\ref{equ: expansion layer}). Note that
\[\ell(x_{\Sigma,k})=\#\{(a,b)\in\Sigma\mid a+b=k\}.\]
As $\Sigma$ is a finite set, we have $x_{\Sigma,k}\neq 1$ only for finitely many $k\in\Z$, so we can define
\begin{equation*}
x_{\Sigma}\defeq \cdots x_{\Sigma,k}\cdot x_{\Sigma,k-1}\cdots \in W(G).
\end{equation*}

\begin{defn0}
Let $\Sigma\subseteq \Z\times\Z$ be a finite subset such that $1\leq a-b\leq n-1$ for $(a,b)\in\Sigma$.
\begin{enumerate}[label=(\roman*)]
\item We say that $\Sigma$ is an \emph{expansion} of $x_{\Sigma}$ if $\ell(x_{\Sigma})=\sum_{k\in\Z}\ell(x_{\Sigma,k})\ (=\#\Sigma)$.
\item We say that $\Sigma$ is \emph{saturated} if the following extra condition holds: for $(a,b),(a',b')\in\Sigma$ and $(a'',b'')\in\Z\times\Z$ such that $(a,b)\leq (a'',b'')\leq (a',b')$, we have $(a'',b'')\in\Sigma$.
\item We say that a saturated $\Sigma$ is \emph{connected} if for each pair of saturated subsets $\Sigma_1$ and $\Sigma_2$ satisfying $\Sigma_1\cap \Sigma_2=\emptyset$ and $\Sigma=\Sigma_1\sqcup\Sigma_2$, there exists $(a_1,b_1)\in\Sigma_1$ and $(a_2,b_2)\in\Sigma_2$ such that $|(a_1-b_1)-(a_2-b_2)|\leq 1$.
\item We say that $x\in W(G)$ is \emph{saturated} if $x=x_{\Sigma}$ for a saturated $\Sigma$.
\item We say a saturated $x_{\Sigma}$ is \emph{connected} if $\Sigma$ can be moreover taken connected.
\end{enumerate}
\end{defn0}

We define the \emph{saturated closure} of any subset $\Sigma\subseteq \Z\times\Z$ to be the minimal saturated subset of $\Z\times\Z$ which contains $\Sigma$, if it exists. Note that the saturated closure can contain elements $(a'',b'')$ such that $a''-b''\notin \{1,\dots,n-1\}$ (see Example \ref{ex: easy elements} below). Note also that there is an obvious way to decompose an arbitrary saturated subset $\Sigma\subseteq \Z\times \Z$ into its connected components.\bigskip

For $x\in W(G)$ we let $\mathrm{Supp}(x)\subseteq \Delta$ be the set of $j\in\{1,\dots,n-1\}$ such that $s_j$ appears in one (and thus all) reduced decomposition of $x$. Recall that $x\in W(G)$ is called Coxeter if it is a product of all distinct simple reflections (each with multiplicity $1$). We say that $x\in W(G)$ is \emph{partial-Coxeter} if there exists a reduced decomposition of $x$ which is multiplicity free, i.e.\ $s_j$ shows up at most once for each $j\in\Delta$. One easily checks from the braid relations that \emph{any} reduced decomposition of a partial-Coxeter element is multiplicity free.

\begin{lem0}\label{lem: saturated expansion}
Let $\Sigma\subseteq \Z\times \Z$ be a saturated subset.
\begin{enumerate}[label=(\roman*)]
\item \label{it: expansion 1} We have $\ell(x_{\Sigma})=\#\Sigma$, i.e.\ $\Sigma$ is an expansion of $x_{\Sigma}$. Moreover, there exists a natural bijection between the set of reduced decompositions of $x_{\Sigma}$ and the set of total orders on $\Sigma$ refining the partial order on $\Sigma$ induced by the one on $\Z\times\Z$.
\item \label{it: expansion 2} We have $j\in D_L(x_{\Sigma})$ (resp.~$j\in D_R(x_{\Sigma})$) if and only if there exists a maximal (resp.~minimal) element $(a,b)\in\Sigma$ such that $j=a-b$.
\item \label{it: expansion 3} If $\Sigma$ is connected and $\Sigma'\subseteq \Z\times\Z$ is another saturated subset such that $x_{\Sigma}=x_{\Sigma'}$, there exists $c\in\Z$ such that $\Sigma'=\{(a-c,b-c)\mid (a,b)\in\Sigma\}$.
\item \label{it: expansion 4} The element $x_{\Sigma}$ is partial-Coxeter if and only if the map
\begin{equation}\label{equ: expansion support}
\Sigma\longrightarrow \mathrm{Supp}(x_{\Sigma}), \ (a,b)\longmapsto a-b
\end{equation}
is bijective. Moreover, all partial-Coxeter elements have the form $x_{\Sigma}$ for a saturated~$\Sigma$.
\end{enumerate}
\end{lem0}
\begin{proof}
We start with several preliminaries. We first fix an arbitrary finite $\Sigma\subseteq \Z\times \Z$ such that $1\leq a-b\leq n-1$ for $(a,b)\in\Sigma$. We write $\leq^\ast$ for a total order on $\Sigma$ which refines the fixed partial order $\leq$.

For each total order $\leq^\ast$ on $\Sigma$ that refines $\leq$, we consider the word
\begin{equation}\label{equ: order product}
x_{\Sigma}(\leq^\ast)\defeq \prod_{(a,b)\in\Sigma,\leq^\ast}s_{a,b}
\end{equation}
where $s_{a,b}$ shows up at the \emph{right} of $s_{a',b'}$ whenever $(a,b)<^\ast (a',b')$. Here we use the symbol $s_{a,b}$ for the copy of $s_{a-b}$ corresponding to $(a,b)$ inside the word (\ref{equ: order product}).

There exists a total order $\leq^\ast$ on $\Sigma$ that refines $\leq$ such that $x_{\Sigma}$ is the image of the word $x_{\Sigma}(\leq^\ast)$ in $W(G)$. This is obvious as we can fix an arbitrary total order on $\{(a,b)\in\Sigma\mid a+b=k\}$ for each $k\in\Z$, and then glue them to a total order $\leq^\ast$ on $\Sigma$ by requiring that $(a,b)\leq^\ast (a',b')$ whenever $a+b\leq a'+b'$.

For two different total orders $\leq_1^\ast$ and $\leq_2^\ast$ on $\Sigma$ which refine $\leq$, we say that $\leq_1^\ast$ and $\leq_2^\ast$ are \emph{adjacent} if there exists exactly one pair $(a,b),(a',b')\in\Sigma$ such that $(a,b)<_1^\ast(a',b')$ and $(a',b')<_2^\ast (a,b)$. If both words (\ref{equ: order product}) have the same image in $W(G)$, we see that $x_{\Sigma}(\leq_1^\ast)$ is reduced if and only if $x_{\Sigma}(\leq_2^\ast)$ is reduced (since then both have length $\#\Sigma$).

Given two total orders $\leq_1^\ast$ and $\leq_2^\ast$ on $\Sigma$ which refine $\leq$, we say that $\leq_1^\ast$ and $\leq_2^\ast$ are \emph{equivalent} if there exist an integer $t\geq 0$ and total orders $\leq_{1,t'}^\ast$ on $\Sigma$ for $0\leq t'\leq t$ which refine $\leq$ such that $\leq_{1,0}^\ast=\leq_1^\ast$, $\leq_{1,t}^\ast=\leq_2^\ast$ and $\leq_{1,t'}^\ast$ is adjacent to $\leq_{1,t'-1}^\ast$ for $1\leq t'\leq t$.

We now prove that all total orders on $\Sigma$ which refine $\leq$ are equivalent. It suffices to prove that any total order $\leq^\ast$ is equivalent to the total order $\leq_0^\ast$ defined by $(a,b)\leq_0^\ast(a',b')$ if and only if $a+b\leq a'+b'$ and $a\leq a'$ (which refines $\leq$). As a total order on $\Sigma$ corresponds to a bijection between $\Sigma$ and the set $\{1,\dots,\#\Sigma\}$, $\leq^\ast$ differs from $\leq_0^\ast$ by a permutation, and we write $\ell(\leq^\ast)$ for the length of this permutation (with $\ell(\leq_0^\ast)=0$). We prove that $\leq^\ast$ is equivalent to $\leq_0^\ast$ by induction on $\ell(\leq^\ast)$. If $\ell(\leq^\ast)=0$, there is nothing to prove. If $\ell(\leq^\ast)>0$, then the subset
\[\{(a,b)\in \Sigma\mathrm{\ such\ that\ }\exists\ (a',b')\in\Sigma\mathrm{\ with\ }(a',b')<_0^\ast(a,b)\mathrm{\ and\ }(a,b)<^\ast(a',b')\}\]
has a maximal element under $\leq_0^\ast$. We choose such a maximal element $(a,b)$, and we choose $(a',b')$ in the above subset to be minimal under $\leq^\ast$. Assume that there exists $(a'',b'')\in\Sigma$ such that $(a,b)<^\ast (a'',b'')<^\ast (a',b')$. We can assume $(a'',b'')$ to be the minimal such element under $\leq^\ast$. The minimality of our choice of $(a',b')$ forces $(a,b)<_0^\ast (a'',b'')$, which in turn contradicts the maximality of $(a,b)$ (as $(a'',b'')<^\ast (a',b')$ and $(a',b')<_0^\ast(a,b)<_0^\ast (a'',b'')$). Hence, there is no element between $(a,b)$ and $(a',b')$ under $\leq^\ast$, and we can define $\leq_1^\ast$ by interchanging $(a,b)$ and $(a',b')$ in the total order $\leq^\ast$. Then $\ell(\leq_1^\ast)=\ell(\leq^\ast)-1$ by \cite[Prop.\ 1.5.8]{BB05}, so $\leq_1^\ast$ is equivalent to $\leq_0^\ast$ by induction. As $\leq_1^\ast$ is adjacent to $\leq^\ast$, it follows that $\leq^\ast$ is also equivalent to $\leq_0^\ast$.

We now assume that $\Sigma$ is saturated. We prove that if we modify the word $x_{\Sigma}(\leq^\ast)$ using a braid relation, the new word is $x_{\Sigma}(\leq'^\ast)$ for another total order $\leq'^\ast$ on $\Sigma$ (refining $\leq$) which is \emph{adjacent} to $\leq^\ast$. As the braid relations are $s_{j}s_{j'}=s_{j'}s_{j}$ for some $|j-j'|\geq 2$ or $s_{j}s_{j-1}s_{j}=s_{j-1}s_{j}s_{j-1}$ for some $2\leq j\leq n-1$, from the definition of adjacent it is enough to rule out the second possibility. Assume on the contrary that there exist $(a,b)<^\ast(a',b')<^\ast(a'',b'')$ in $\Sigma$ such that the only element in $\Sigma$ between $(a,b)$ and $(a'',b'')$ (under $\leq^\ast$) is $(a',b')$ and such that $a-b=j=a''-b''=j$ and $a'-b'=j-1$ for some $2\leq j\leq n-1$. Since $a-b=j=a''-b''=j$ and $a'-b'=j-1$, and since $\leq^\ast$ refines $\leq$, it is easy to check that we must have $(a,b)< (a',b')< (a'',b'')$. Then it is also easy to see that the saturated condition on $\Sigma$ forces the existence of at least another $(a''',b''')\in\Sigma$ such that $(a,b)< (a''',b''')< (a'',b'')$ and $a'''-b'''=j+1$. But this contradicts the fact that $(a',b')$ is the only element between $(a,b)$ and $(a'',b'')$ under $\leq^\ast$. An analogous argument shows that we cannot have elements $(a,b)<^\ast(a',b')<^\ast(a'',b'')$ in $\Sigma$ such that the only element between $(a,b)$ and $(a'',b'')$ is $(a',b')$ and such that $a-b=j-1=a''-b''=j-1$ and $a'-b'=j$ for some $2\leq j\leq n-1$. Hence, the only braid relations that can be used to modify $x_{\Sigma}(\leq^\ast)$ are $s_{j}s_{j'}=s_{j'}s_{j}$ for some $|j-j'|\geq 2$. This implies that the new word has the form $x_{\Sigma}(\leq'^\ast)$ for some $\leq'^\ast$ which is adjacent to $\leq^\ast$.

We now prove \ref{it: expansion 1}. From the previous statements, we only have to prove that $x_{\Sigma}(\leq^\ast)$ is a reduced word for one (equivalently any) total order $\leq^\ast$ on $\Sigma$ refining $\leq$. Note first that if two copies of $s_j$ occur in $x_{\Sigma}(\leq^\ast)$, they come from two distinct elements $(a_1,b_1), (a_2,b_2) \in \Sigma$ such that $a_1-b_1=a_2-b_2=j$. We can assume $a_1<a_2$, and we see that $(a_1+1,b_1)$ also belongs to $\Sigma$ as $\Sigma$ is saturated. We deduce that $s_{j+1}$ occurs in $x_{\Sigma}(\leq^\ast)$ between these two $s_j$. In particular, the two $s_j$ cannot cancel after applying the braid relations. It follows that $x_{\Sigma}(\leq^\ast)$ is reduced.

We prove \ref{it: expansion 2}. An element $(a,b)\in\Sigma$ is the unique maximal (resp.~minimal) element under some total order $\leq^\ast$ if and only if $(a,b)$ is a maximal (resp.~minimal) element under the partial order $\leq$. Then the statement follows from $\ell(x_{\Sigma})=\#\Sigma$ and \cite[Cor.~1.4.6]{BB05}.

We prove \ref{it: expansion 3}. Let $\Sigma'\subseteq \Z\times\Z$ be another saturated subset such that $x_{\Sigma}=x_{\Sigma'}$, we argue by induction on $\#\Sigma$. If $\#\Sigma=1$, the claim is obvious. If $\#\Sigma>1$, there exists $(a_0,b_0)\in\Sigma$ which is either minimal or maximal (for $\leq$), such that $\Sigma\setminus\{(a_0,b_0)\}$ is still saturated and connected. Let $j\defeq a_0-b_0$. If $(a_0,b_0)$ is minimal (resp.~maximal) in $\Sigma$, we have $j\in D_R(x_{\Sigma})=D_R(x_{\Sigma'})$ (resp.~$j\in D_L(x_{\Sigma})=D_L(x_{\Sigma'})$) (using \ref{it: expansion 2}), which gives a unique minimal (resp.~maximal) element $(a_0',b_0')\in\Sigma'$ satisfying $a_0'-b_0'=j$. By our induction assumption, there exists $c\in\Z$ such that
\[\Sigma'\setminus \{(a_0',b_0')\}=\{(a-c,b-c)\mid (a,b)\in\Sigma\setminus\{(a_0,b_0)\}\}.\]
Thus we only need to prove $(a_0',b_0')=(a_0-c,b_0-c)$. Assume that $(a_0,b_0)$ is maximal in $\Sigma$. Then as $\Sigma$ is connected, either $(a_0-1,b_0)$ or $(a_0,b_0-1)$ is in $\Sigma\setminus\{(a_0,b_0)\}$ but none of $(a_0+1,b_0)$, $(a_0,b_0+1)$. Thus $\Sigma'\setminus\{(a_0',b_0')\}$ contains either $(a_0-c-1,b_0-c)$ or $(a_0-c,b_0-c-1)$ but not $(a_0-c+1,b_0-c), (a_0-c,b_0-c+1)$. As $(a_0',b_0')$ is maximal in $\Sigma'$ and $\Sigma'$ is saturated, this forces $(a_0',b_0')=(a_0-c,b_0-c)$. The proof for $(a_0,b_0)$ minimal is analogous.

We prove \ref{it: expansion 4}. As $x_{\Sigma}(\leq_0^\ast)$ is a reduced word, the map (\ref{equ: expansion support}) takes values in $\mathrm{Supp}(x_{\Sigma})\subseteq \{1,\dots, n-1\}$ and is clearly surjective. Hence, (\ref{equ: expansion support}) is a bijection if and only if it is injective if and only if any reduced decomposition of $x_{\Sigma}$ is multiplicity free, which is equivalent to $x_{\Sigma}$ being partial-Coxeter by definition. We now prove that any partial-Coxeter element $x\in W(G)$ has the form $x_{\Sigma}$ for some saturated $\Sigma$ by induction on $\ell(x)$. The case $\ell(x)=0$ is trivial. Let $1\neq x\in W(G)$ be partial-Coxeter and choose $j\in\mathrm{Supp}(x)$ such that $j'\notin\mathrm{Supp}(x)$ for $j'> j$. We set $x'\defeq s_jx$ (resp.~$x'\defeq xs_j$) if $j\in D_L(x)$ (resp.~if $j\in D_R(x)$). We have $j\in D_L(x)\cap D_R(x)$ if $j-1\notin\mathrm{Supp}(x)$, $j\in D_R(x)\setminus D_L(x)$ if $s_{j-1}s_j\leq x$, $j\in D_L(x)\setminus D_R(x)$ if $s_js_{j-1}\leq x$, and $\ell(x')=\ell(x)-1$ in all cases. By our induction assumption, there is a saturated subset $\Sigma'$ such that $x'=x_{\Sigma'}$. When $j-1\in \mathrm{Supp}(x')$ let $(a',b')\in\Sigma'$ be the unique element such that $a'-b'=j-1$. Then we claim that $x=x_{\Sigma}$ for $\Sigma\defeq \Sigma'\sqcup\{(a,b)\}$ with $(a,b)$ described as follows:
\begin{itemize}
\item if $j-1\notin \mathrm{Supp}(x)$, then $x=s_jx'=x's_j$ and we take $(a,b)$ to be any element of $\Z\times\Z$ satisfying $a-b=j$;
\item if $s_{j-1}s_j\leq x$, then $x=x's_j$ and we take $(a,b)\defeq (a',b'-1)\leq (a',b')$;
\item if $s_js_{j-1}\leq x$, then $x=s_jx'$ and we take $(a,b)\defeq (a'+1,b')\geq (a',b')$. \qedhere
\end{itemize}
\end{proof}

\begin{ex0}\label{ex: easy elements}
We take $\fg=\fg\fl_n$. For $n=4$, the partial-Coxeter elements of length $3$ are $s_3s_2s_1$, $s_1s_2s_3$, $s_1s_3s_2$, $s_2s_1s_3$. The element $s_2s_1s_3s_2$ is saturated and connected (take for instance $\Sigma=\{(2,0), (2,1), (3,0), (3,1) \}$) but not partial-Coxeter (and is in fact the only such element for $\fg\fl_4$). The elements $s_1s_2s_1$ and $s_1s_2s_3s_2s_1$ are not saturated (for the first one, one could think of $\Sigma=\{(2,0), (2,1), (3,2)\}$ or $\Sigma=\{(2,1),(3,1),(3,2)\}$ or $\Sigma=\{(2,0),(2,1),(3,1)\}$ but their saturated closure contains the extra element $(2,2)$ for the first two, the extra element $(3,0)$ for the third). More generally, for $n$ sufficiently large, $s_js_{j-1}s_{j+1}s_js_{j+2}s_{j+1}s_{j+3}s_{j+2}$ and its inverse are examples of saturated connected elements.
\end{ex0}

\begin{lem0}\label{lem: weaker dominance}
Let $\Sigma\subseteq \Z\times\Z$ be a saturated subset with a unique minimal element $(a_1,b_1)$ and a unique maximal element $(a_2,b_2)$. Assume there is $j\in\mathrm{Supp}(x_{\Sigma})$ such that $|j-(a_2-b_2)|=1$. Then we have
\begin{equation}\label{equ: weaker dominance}
\{a_1-b_1\}\subsetneq D_R(s_jx_{\Sigma}).
\end{equation}
\end{lem0}
\begin{proof}
It follows from \ref{it: expansion 2} of Lemma~\ref{lem: saturated expansion} that $D_R(x_{\Sigma})=\{a_1-b_1\}$ and $D_L(x_{\Sigma})=\{a_2-b_2\}$. This forces $s_jx_{\Sigma}>x_{\Sigma}>x_{\Sigma}s_{a_1-b_1}$ (recall from the proof of Proposition \ref{lem: saturated expansion} that the braid relations in $x_{\Sigma}$ are only of the type $s_{j}s_{j'}=s_{j'}s_{j}$), which implies $\ell(s_jx_{\Sigma}s_{a_1-b_1})\leq \ell(x_{\Sigma}s_{a_1-b_1})+1=\ell(s_jx_{\Sigma})-1$ and thus $s_jx_{\Sigma}>s_jx_{\Sigma}s_{a_1-b_1}$. Hence $a_1-b_1\in D_R(s_jx_{\Sigma})$. Now we prove that (\ref{equ: weaker dominance}) is a strict inclusion by increasing induction on $\#\Sigma$. If $\#\Sigma=1$, then $\mathrm{Supp}(x_{\Sigma})=\{a_2-b_2\}$ and there is nothing to prove. We can assume $j=a_2-b_2+1$ as the proof for $j=a_2-b_2-1$ is symmetric. We have the following two possibilities.\bigskip

\textbf{Case $1$}: If $a_2-b_2-1\notin \mathrm{Supp}(x_{\Sigma})$, then we must have $a_1=a_2$ and $\Sigma=\{(a_1,b)\mid b_1\leq b\leq b_2\}$. Hence, we have
\[s_jx_{\Sigma}=s_{a_2-b_2+1}s_{a_2-b_2}s_{a_2-b_2+1}x_{\Sigma'}=s_{a_2-b_2}s_{a_2-b_2+1}s_{a_2-b_2}x_{\Sigma'}=s_{a_2-b_2}s_{a_2-b_2+1}x_{\Sigma'}s_{a_2-b_2}\]
with $\Sigma'\defeq \Sigma\setminus\{(a_2,b_2),(a_2,b_2-1)\}$. This implies $\{a_1-b_1\}\subsetneq \{a_2-b_2,a_1-b_1\}\subseteq D_R(s_jx_{\Sigma})$.\bigskip

\textbf{Case $2$}: If $a_2-b_2-1\in \mathrm{Supp}(x_{\Sigma})$, then we must have $(a_1,b_1)\leq (a_2-1,b_2)\in\Sigma$. Let $\Sigma'$ be the saturated closure of $\{(a_1,b_1),(a_2-1,b_2)\}$, which is a (saturated) subset of $\Sigma$. Moreover, $\Sigma\setminus\Sigma'=\{(a_2,b)\mid b_1\leq b\leq b_2\}$ is also a non-empty saturated subset of $\Sigma$. It is not difficult to check that $x_{\Sigma}=x_{\Sigma\setminus \Sigma'}x_{\Sigma'}$. As $j=a_2-b_2+1\in \mathrm{Supp}(x_{\Sigma})$, we also have $(a_1,b_1)\leq (a_2,b_2-1)\in\Sigma$, and we set $\Sigma''\defeq (\Sigma\setminus \Sigma')\setminus\{(a_2,b_2),(a_2,b_2-1)\}$, which is yet another saturated subset of $\Sigma$. We observe that
\begin{multline}\label{equ: dominance induction}
s_jx_{\Sigma}=s_{a_2-b_2+1}s_{a_2-b_2}s_{a_2-b_2+1}x_{\Sigma''}x_{\Sigma'}=s_{a_2-b_2}s_{a_2-b_2+1}s_{a_2-b_2}x_{\Sigma''}x_{\Sigma'}\\
=s_{a_2-b_2}s_{a_2-b_2+1}x_{\Sigma''}s_{a_2-b_2}x_{\Sigma'}
\end{multline}
where each decomposition in \ref{equ: dominance induction} is reduced. As $(a_1,b_1)\leq (a_2-1,b_2)$ and $(a_1,b_1)\leq (a_2,b_2-1)$, we have $(a_1,b_1)\leq (a_2-1,b_2-1)\leq (a_2-1,b_2)$, which implies $(a_2-1,b_2-1)\in\Sigma'$ and $a_2-b_2=(a_2-1)-(b_2-1)\in\mathrm{Supp}(x_{\Sigma'})$. By our induction assumption applied to $\Sigma'$ and $a_2-b_2$, we have $\{a_1-b_1\}\subsetneq D_R(s_{a_2-b_2}x_{\Sigma'})$, which together with (\ref{equ: dominance induction}) finishes the proof.
\end{proof}

\begin{lem0}\label{lem: maximal dominance}
Let $x\in W(G)$ and $j\in\Delta$ with $D_R(x)=\{j\}$. Then there exists a saturated connected subset $\Sigma\subseteq \Z\times\Z$ such that $x=x_{\Sigma}$. Moreover, $\Sigma$ contains a unique minimal element $(a,b)$ that satisfies $a-b=j$.
\end{lem0}
\begin{proof}
Connectedness is clear otherwise we would have $\#D_R(x)>1$. We prove the statement by induction on $\ell(x)$. If $\ell(x)=1$, then we deduce from $D_R(x)=\{j\}$ that $x=s_j$. So we can take $\Sigma=\{(a,b)\}$ for an arbitrary $(a,b)\in\Z\times\Z$ satisfying $a-b=j$. We assume from now on $\ell(x)\geq 2$. We choose an arbitrary $j'\in D_L(x)$ and set $x'\defeq s_{j'}x$ which satisfies $\ell(x')=\ell(x)-1\geq 1$. For each $j''\in D_R(x')$, we have $x=s_{j'}x'>x'>x's_{j''}$ which implies $\ell(s_{j'}x's_{j''})\leq \ell(x's_{j''})+1=\ell(s_{j'}x')-1$ and thus $x=s_{j'}x'>s_{j'}x's_{j''}=xs_{j''}$, equivalently $j''\in D_R(x)$. As $D_R(x)=\{j\}$, we conclude that $D_R(x')=\{j\}$, which together with our induction assumption gives a saturated subset $\Sigma'\subseteq \Z\times\Z$ such that $x'=x_{\Sigma'}$. Our induction assumption also says that $\Sigma'$ contains a unique minimal element $(a_0,b_0)$ that satisfies $a_0-b_0=j$. We now construct $\Sigma$ from $\Sigma'$, according to the following two possibilities.

Assume $\mathrm{Supp}(x')\cap \{j'+1,j'-1\}=\emptyset$. Then $x=s_{j'}x'=x's_{j'}$ and thus $j'\in D_R(x)$, which forces $j=j'$. Thus $x=x's_j$ which contradicts $D_R(x)=D_R(x')=\{j\}$.

We thus have $\mathrm{Supp}(x')\cap \{j'+1,j'-1\}\neq\emptyset$. Since $\Sigma'$ contains a unique minimal element $(a_0,b_0)$, it is easy to see that, for each saturated subset $\Sigma''\subseteq \Sigma'$ containing $(a_0,b_0)$, $\Sigma'\setminus\Sigma''$ is another saturated subset of $\Sigma'$ and that $x_{\Sigma'}=x_{\Sigma'\setminus\Sigma''}x_{\Sigma''}$.
We choose such $\Sigma''$ to be the saturated closure of the following subset
\begin{equation*}
\{(a_0,b_0)\}\sqcup\{(a,b)\in\Sigma'\mid |(a-b)-j'|\leq 1\}\subseteq \Sigma'.
\end{equation*}
As $s_{j'}s_{a'-b'}=s_{a'-b'}s_{j'}$ for each $(a',b')\in\Sigma'\setminus\Sigma''$, we have
\begin{equation}\label{equ: move reflection}
x=s_{j'}x_{\Sigma'}=x_{\Sigma'\setminus\Sigma''}s_{j'}x_{\Sigma''}
\end{equation}
with both sides being reduced (using \ref{it: expansion 1} of Lemma \ref{lem: saturated expansion}). As $s_{j'}x_{\Sigma'}>x_{\Sigma'}$, we have $s_{j'}x_{\Sigma''}>x_{\Sigma''}$ and thus a maximal element $(a,b)\in\Sigma''$ cannot satisfy $a-b=j'$. As $\mathrm{Supp}(x')\cap \{j'+1,j'-1\}\neq\emptyset$, any maximal element $(a,b)\in\Sigma''$ must satisfy $|a-b-j'|=1$. We now have the following $3$ cases.
\begin{itemize}
\item Assume that $\Sigma''$ has two different maximal elements $(a_1,b_1)$, $(a_2,b_2)$, and we can assume $a_1\geq a_2$. If $a_1-b_1\leq a_2-b_2$, then $a_1\geq a_2$ forces $(a_1,b_1)\geq (a_2,b_2)$ which contradicts the fact that both $(a_1,b_1)$ and $(a_2,b_2)$ are maximal in $\Sigma''$. So we must have $a_1-b_1=j'+1$ and $a_2-b_2=j'-1$. If $a_1\geq a_2+2$, this implies $b_1\geq b_2$ and thus $(a_1,b_1)> (a_2,b_2)$, another contradiction. If $a_1=a_2$, then $b_1=b_2-2$ and thus $(a_1,b_1)< (a_2,b_2)$, also a contradiction. So the only possibility is $a_1=a_2+1$ and $b_1=b_2-1$, in which case we simply take $\Sigma\defeq \Sigma'\sqcup\{(a_1,b_2)\}$ which is easily seen to be saturated and satisfies $x=s_{j'}x'=x_{\Sigma}$.
\item Assume that $\Sigma''$ has a unique maximal element $(a_1,b_1)$ and $j'\notin\mathrm{Supp}(x')$. We can assume $a_1-b_1=j'+1$ as the case $a_1-b_1=j'-1$ is similar. Then $\Sigma''$ is the saturated closure of $\{(a_0,b_0),(a_1,b_1)\}$. Since $j'\notin\mathrm{Supp}(x')$, we must have $a_0=a_1$ and $\Sigma''=\{(a_1,b)\mid b_0\leq b\leq b_1\}$. Hence, it is clear that $\Sigma''\sqcup\{(a_1,b_1+1)\}$ is saturated. Using that $\Sigma''$ contains the unique minimal element of $\Sigma'$, we also deduce that $\Sigma\defeq \Sigma'\sqcup\{(a_1,b_1+1)\}$ is saturated. As $(a_1,b_1+1)$ must also be maximal in $\Sigma$ and $a_1-(b_1+1)=j'$, we conclude that $x_{\Sigma}=s_{j'}x_{\Sigma'}=x$.
 \item Assume that $\Sigma''$ has a unique maximal element $(a_1,b_1)$ and $j'\in\mathrm{Supp}(x')$. Here again we can assume $a_1-b_1=j'+1$. As $\Sigma''$ is the saturated closure of $\{(a_0,b_0),(a_1,b_1)\}$ and $j'\in\mathrm{Supp}(x')$, we must have $j'\in\mathrm{Supp}(x_{\Sigma''})$ and $a_0<a_1$ which forces $(a_0,b_0)\leq (a_1-1,b_1)\leq (a_1,b_1)$ and $(a_1-1,b_1)\in \Sigma''$. But then it follows from Lemma~\ref{lem: weaker dominance} that $\{a_0-b_0\}\subsetneq D_R(s_{j'}x_{\Sigma''})$, which together with (\ref{equ: move reflection}) implies that $\{a_0-b_0\}=\{j\}\subsetneq D_R(s_{j'}x')=D_R(x)=\{j\}$, a contradiction.
\end{itemize}
We have constructed the desired $\Sigma$ from $\Sigma'$ in all possible cases, which finishes the proof.
\end{proof}

\begin{rem0}\label{rem: 1 dominance}
It is an easy check that if a (non-empty) saturated subset $\Sigma\subseteq \Z\times\Z$ admits a unique minimal element $(a_0,b_0)$ and moreover $a_0-b_0=1$, then $\Sigma$ must have the form $\{(a_0+c,b_0)\mid 0\leq c\leq c_0\}$ for some $0\leq c_0\leq n-2$.
Similarly, if a (non-empty) saturated subset $\Sigma\subseteq \Z\times\Z$ admits a unique maximal element $(a_1,b_1)$ and moreover $a_1-b_1=n-1$, then $\Sigma$ must have the form $\{(a_1-c,b_1)\mid 0\leq c\leq c_1\}$ for some $0\leq c_1\leq n-2$.
\end{rem0}

In the lemmas below we use without comment notation related to the Kazhdan-Lusztig polynomials, see the beginning of \S\ref{subsec: Ext O}.

\begin{lem0}\label{lem: saturated triple}
Let $j\in\Delta$ and $x,w\in W(G)$ such that $D_L(x)=D_L(w)=\{j\}$. Then there does not exist $x'$ such that $x\prec x'\prec w$ and $j\notin D_L(x')$.
\end{lem0}
\begin{proof}
Assume on the contrary that such a triple $x,w,x'$ exist. By applying Lemma~\ref{lem: maximal dominance} to $x^{-1}$ and $w^{-1}$ (with $D_R(x^{-1})=D_R(w^{-1})=\{j\}$) we know that both $x$ and $w$ are saturated, i.e.~there exist saturated subsets $\Sigma_1,\Sigma_2\subseteq \Z\times\Z$ such that $x=x_{\Sigma_1}$ and $w=x_{\Sigma_2}$. As $j\in D_L(w)\setminus D_L(x')$ and $w\prec x'$, by Lemma~\ref{lem: dominance control} we have $w=s_jx'>x'$ and in particular $x'$ is saturated with $x'=x_{\Sigma_2\setminus\{(a,b)\}}$ where $(a,b)$ is the maximal element of $\Sigma_2$ satisfying $a-b=j$. Note that $D_L(w)=\{j\}$ forces $D_L(x')\subseteq \{j-1,j+1\}\cap \Delta$.
As $D_L(x)=\{j\}$ and $j\notin D_L(x')\neq \emptyset$, by Lemma~\ref{lem: dominance control} we have $x'=s_{j'}x>x$ for some $j'\in D_L(x')\subseteq \{j-1,j+1\}\cap \Delta$. If we write $(a',b')$ for the maximal element of $\Sigma_2\setminus\{(a,b)\}$ that satisfies $a'-b'=j'$, then we can choose $\Sigma_1=\Sigma_2\setminus\{(a,b),(a',b')\}$. Let $(a'',b'')$ be the unique maximal element of $\Sigma_1$ such that $a''-b''=j$ (using $D_L(x)=\{j\}$ and \ref{it: expansion 2} of Lemma~\ref{lem: saturated expansion}). As $(a'',b'')\in \Sigma_1$, we have $(a'',b'')=(a-c,b-c)$ for some $c\geq 1$, and in particular $(a-1,b),(a,b-1)\in\Sigma_2\setminus\{(a,b)\}$ as $\Sigma_2$ is saturated. Hence, either $(a-1,b)$ or $(a,b-1)$ remains in $\Sigma_1=(\Sigma_2\setminus\{(a,b),(a',b')\}$, contradicting the maximality of $(a'',b'')$ in $\Sigma_1$.
\end{proof}

For $j_1,j_2\in\Delta$ with $|j_1-j_2|=1$ and $x\in W(G)$ with $\#(D_L(x)\cap\{j_1,j_2\})=1$, we define $\theta_{j_1,j_2}(x)$ as the unique element in $\{s_{j_1}x,s_{j_2}x\}$ such that $\#(D_L(\theta_{j_1,j_2}(x))\cap\{j_1,j_2\})=1$. It is clear that $j_i\in D_L(x)$ for $i\in \{1,2\}$ if and only if $j_{3-i}\in D_L(\theta_{j_1,j_2}(x))$.

\begin{lem0}\label{lem: KL involution}
Let $x,w\in W(G)$ and $j_1,j_2\in\Delta$ with $|j_1-j_2|=1$. Assume $w\notin W(L_{\{j_1,j_2\}})x$ and $\#(D_L(\theta_{j_1,j_2}(x))\cap\{j_1,j_2\})=1=\#(D_L(\theta_{j_1,j_2}(w))\cap\{j_1,j_2\})$. Then we have $x\prec w$ if and only if $\theta_{j_1,j_2}(x)\prec \theta_{j_1,j_2}(w)$, in which case $\mu(x,w)=\mu(\theta_{j_1,j_2}(x),\theta_{j_1,j_2}(w))$.
\end{lem0}
\begin{proof}
This is \cite[Thm.~4.2(i)]{KL79}.
\end{proof}

For $j,j'\in \Delta$, recall that $w_{j,j'}\in W(G)$ is defined in (\ref{wj1j2}).

\begin{lem0}\label{lem: coxeter prec}
Let $x,w\in W(G)$ with $x\prec w$. If $x=1$, then $w=s_j$ for some $j\in\Delta$. If $x=w_{j,j'}$, then $\ell(w)>\ell(x)+1$ if and only if $x=s_j$ and $w=s_js_{j+1}s_{j-1}s_j$ for some $2\leq j\leq n-2$, in which case we have $\mu(x,w)=1$.
\end{lem0}
\begin{proof}
The claim when $x=1$ is clear. Assume from now that $x=w_{j,j'}$ and $w$ satisfies $x\prec w$ and $\ell(w)>\ell(x)+1$. We can assume $j\geq j'$ (the case $j\leq j'$ is similar). By Lemma~\ref{lem: dominance control} we must have $D_L(w)=D_L(x)=\{j\}$ and $D_R(w)=D_R(x)=\{j'\}$. By Lemma~\ref{lem: maximal dominance} we know that $w$ is saturated, i.e.~$w=x_{\Sigma}$ for some saturated subset $\Sigma\subseteq \Z\times\Z$. By \ref{it: expansion 2} of Lemma~\ref{lem: saturated expansion} the set $\Sigma$ admits a unique maximal element $(a,b)$ with $a-b=j$ and a unique minimal element $(a',b')$ with $a'-b'=j'$. If either $a=a'$ or $b=b'$, then it is easy to check that $w=w_{j,j'}$ which contradicts $\ell(w)>\ell(x)+1$. Hence, we must have $a'<a$ and $b'<b$ and in particular $(a-1,b-1),(a-1,b),(a,b-1)\in\Sigma$ (and $j-1,j+1\in\Delta$).
One checks that the only $x'\in W(L_{\{j,j+1\}})x$ satisfying $D_L(x')=\{j\}$ and $D_R(x')=\{j'\}$ is $x'=x$, and in particular $w\notin W(L_{\{j,j+1\}})x$. This together with $x\prec w$, $D_L(x)=D_L(w)=\{j\}$ and Lemma~\ref{lem: KL involution} implies $\theta_{j,j+1}(x)\prec \theta_{j,j+1}(w)$. We check that $\theta_{j,j+1}(x)=s_{j+1}x=w_{j+1,j'}$ with $D_L(\theta_{j,j+1}(x))=\{j+1\}$ and $\theta_{j,j+1}(w)=s_jw<w$ with $D_L(\theta_{j,j+1}(w))=\{j-1,j+1\}$ (as $\Sigma\setminus\{(a,b)\}$ has maximal elements $(a-1,b)$ and $(a,b-1)$), which together with Lemma~\ref{lem: dominance control} forces $\theta_{j,j+1}(w)=s_{j-1}\theta_{j,j+1}(x)=s_{j-1}s_{j+1}w_{j,j'}$. Equivalently, we have $x=w_{j,j'}=x_{\Sigma'}$ with $\Sigma'\defeq \Sigma\setminus\{(a,b),(a-1,b),(a,b-1)\}$ (which admits $(a-1,b-1)$ as the unique maximal element). This together with $j\geq j'$ forces $b'=b-1$ and $\Sigma'=\{(a'',b-1)\mid a'\leq a''\leq a-1\}$. As $\Sigma=\Sigma'\sqcup\{(a,b),(a-1,b),(a,b-1)\}$ is saturated, we must also have $a'=a-1$, i.e.~$x=s_j$ and $w=s_js_{j-1}s_{j+1}s_j$ (for some $2\leq j\leq n-2$). Finally, note that $\mu(x,w)=\mu(\theta_{j,j+1}(x),\theta_{j,j+1}(w))=1$ by Lemma~\ref{lem: KL involution} (and $\ell(\theta_{j,j+1}(w))=\ell(\theta_{j,j+1}(x))+1$).
\end{proof}

\begin{rem0}\label{rem: coxeter pair}
Let $j_0,j_1,j_0',j_1'\in\Delta$ and consider $w_{j_1,j_0},w_{j_1',j_0'}\in W(G)$. It is easy to check (using \cite[Thm.~2.2.2]{BB05}) that $w_{j_1,j_0}\leq w_{j_1',j_0'}$ if and only if either $j_1'\geq j_1\geq j_0\geq j_0'$ or $j_1'\leq j_1\leq j_0\leq j_0'$. Combined with Lemma~\ref{lem: coxeter prec}, we observe that $w_{j_1,j_0}\prec w_{j_1',j_0'}$ if and only if $w_{j_1,j_0}<w_{j_1',j_0'}$ and $\ell(w_{j_1',j_0'})=\ell(w_{j_1,j_0})+1$, in which case $|j_0-j_0'|+|j_1-j_1'|=1$.
\end{rem0}

\begin{lem0}\label{lem: special descent change}
Let $j,j'\in\Delta$ and $x\in W(G)$ such that $D_L(x)=\{j_1\}$ for some $j_1\neq j$.
\begin{enumerate}[label=(\roman*)]
\item \label{it: special descent 1} We have $x\prec w_{j,j'}$ if and only if $x=w_{j+1,j'}$ with $j'>j$, or $x=w_{j-1,j'}$ with $j'<j$.
\item \label{it: special descent 2} We have $w_{j,j'}\prec x$ if and only if $x=w_{j+1,j'}$ with $j'\leq j$, or $x=w_{j-1,j'}$ with $j'\geq j$.
\item \label{it: special descent 3} We have $\mathrm{Ext}_{U(\fg)}^1(L(w_{j,j'}),L(x))\neq 0$ if and only if $x=w_{j_1,j'}$ with $j_1\in\{j-1,j+1\}\cap\Delta$, in which case $\Dim_E \mathrm{Ext}_{U(\fg)}^1(L(w_{j,j'}),L(x))=1$.
\item \label{it: special descent 4} We have $w_{j,j'}\prec w$ with $j'\notin D_R(w)$ if and only if $w=w_{j,j'+1}$ with $j\leq j'$, or $w=w_{j,j'-1}$ with $j\geq j'$.
\end{enumerate}
\end{lem0}
\begin{proof}
We prove \ref{it: special descent 1}. As $x\prec w_{j,j'}$ and $j\in D_L(w_{j,j'})\setminus D_L(x)$, by Lemma~\ref{lem: dominance control} we have $w_{j,j'}=s_jx$. As $D_L(x)\neq \emptyset$, we have $x\neq 1$ and thus $j\neq j'$. So either we have $j'>j$ and $x=w_{j+1,j'}$, or we have $j'<j$ and $x=w_{j-1,j'}$. The other direction is obvious.

We prove \ref{it: special descent 2}. As $w_{j,j'}\prec x$ and $j_1\in D_L(x)\setminus D_L(w_{j,j'})$, by Lemma~\ref{lem: dominance control} we have $x=s_{j_1}w_{j,j'}$. As $D_L(x)=\{j_1\}$, by Lemma~\ref{lem: maximal dominance} (applied to $x^{-1}$) we know that $x=x_{\Sigma}$ for some saturated subset $\Sigma\subseteq \Z\times\Z$. By \ref{it: expansion 2} of Lemma~\ref{lem: saturated expansion} we know that $\Sigma$ admits a unique maximal element $(a_1,b_1)$ with $a_1-b_1=j_1$, and that $x_{\Sigma\setminus\{(a_1,b_1)\}}=w_{j,j'}$. In particular, $\Sigma\setminus\{(a_1,b_1)\}$ admits a unique maximal element $(a,b)$ with $a-b=j$ and a unique minimal element $(a',b')$ with $a'-b'=j'$, with $a=a'$ and $b\geq b'$ (resp.~$a\geq a'$ and $b=b'$) if $j\geq j'$ (resp.~if $j\leq j'$). As $\Sigma$ is saturated and $(a,b)$ is the unique maximal element of $\Sigma\setminus\{(a_1,b_1)\}$, either we have $(a,b)=(a_1-1,b_1)$ with $(a_1,b_1-1)\notin\Sigma$, or we have $(a,b)=(a_1,b_1-1)$ with $(a_1-1,b_1)\notin\Sigma$. In other words, either we have $(a,b)=(a_1-1,b_1)$, $a\geq a'$ and $b=b'$ (with $j_1=j+1$ and $j\geq j'$), or we have $(a,b)=(a_1,b_1-1)$, $a=a'$ and $b\geq b'$ (with $j_1=j-1$ and $j\leq j'$). This finishes the proof of \ref{it: special descent 2} as the other direction is trivial.

\ref{it: special descent 3} follows directly from \ref{it: special descent 1}, \ref{it: special descent 2} and from \ref{it: rabiotext 2} of Lemma~\ref{rabiotext}.

We prove \ref{it: special descent 4}. First $D_R(w_{j,j'})=\{j'\}$ and $j'\notin D_R(w)$ together with Lemma~\ref{lem: dominance control} force $w=w_{j,j'}s_{j''}$ for some $j''\neq j'$. If $|j''-j'|\neq 1$, then $s_{j'}s_{j''}=s_{j''}s_{j'}$ and thus $j'\in D_R(w)$, a contradiction. Hence $j''\in \{j'-1,j'+1\}$. If $j\leq j'$ (resp.~$j\geq j'$), then we have $D_R(w_{j,j'}s_{j'+1})=\{j'+1\}$ (resp.~$D_R(w_{j,j'}s_{j'-1})=\{j'-1\}$), which gives the two cases of the statement. If $j<j'$ and $j''=j'-1$ (resp.~$j>j'$ and $j''=j'+1$), then we have $w_{j,j'}s_{j'-1}=s_{j'}w_{j,j'-1}s_{j'}$ (resp.~$w_{j,j'}s_{j'+1}=s_{j'}w_{j,j'+1}s_{j'}$) and thus $j'\in D_R(w_{j,j'}s_{j'-1})$ (resp.~$j'\in D_R(w_{j,j'}s_{j'+1})$), a contradiction. Here again the other direction is obvious.
\end{proof}

\begin{lem0}\label{lem: KL coxeter}
Let $x,w\in W(G)$ with $x<w$.
\begin{enumerate}[label=(\roman*)]
\item \label{it: coxeter mult 1} If $w$ is partial-Coxeter, then $L(w)$ has multiplicity $1$ in $M(x)$.
\item \label{it: coxeter mult 2} If $w=s_jw_{n-1,1}$ and $x=w_{j,1}$ for some $j<n-1$, then $L(w)$ has multiplicity $1$ in $M(x)$.
\end{enumerate}
\end{lem0}
\begin{proof}
It is enough to prove $P_{w_0xw_0,w_0ww_0}=1$ for both \ref{it: coxeter mult 1} and \ref{it: coxeter mult 2} (cf.~\cite[\S 8.4]{Hum08}). We write $w'\defeq w_0ww_0$ and $x'\defeq w_0xw_0$ and note that $x'<w'$.

We prove $P_{x',w'}=1$ for \ref{it: coxeter mult 1}. As $w$ is partial-Coxeter, so are $w'$ and $x'$.
Let $j\in D_L(w')$, since $w'$ is partial-Coxeter we have $s_jz>z$ for $z\prec s_jw'<w'$. We deduce from \cite[(2.2.c)]{KL79} the following two possibilities. If $j\in \mathrm{Supp}(x')$, then $j\in D_L(x')$ since $x'<w'$ (and $w'$ is partial-Coxeter), and we have $P_{x',w'}=P_{s_jx',s_jw'}$ in this case. If $j\notin \mathrm{Supp}(x')$, then $x'\leq s_jw'$, $s_jx'>x'$ and we have $P_{x',w'}=P_{x',s_jw'}$ in this case. A simple induction on $\ell(w')$ thus shows $P_{x',w'}=1$.

We prove $P_{x',w'}=1$ for \ref{it: coxeter mult 2}. As $w=s_jw_{n-1,1}$ and $x=w_{j,1}$, we have $w'=s_{j'}w_{1,n-1}$ and $x'=w_{j',n-1}$ with $j'\defeq n-j>1$. We consider any $z$ satisfying $x'\leq z\prec s_{j'}w'=w_{1,n-1}$ and $s_{j'}z<z$. As $w_{1,n-1}$ is partial-Coxeter, the proof of \ref{it: coxeter mult 1} shows that $P_{z,s_{j'}w'}=1$, which together with $z\prec s_{j'}w'=w_{1,n-1}$ forces $\ell(z)=\ell(w_{n-1,1})-1$. Using \cite[Thm.~2.2.2]{Brenti}, it is easy to check that the only $z$ satisfying $z<w_{n-1,1}$, $\ell(z)=\ell(w_{n-1,1})-1$ and $s_{j'}z<z$ is $z=w_{n-1,j'+1}w_{j',1}$. By the proof of \ref{it: coxeter mult 1} we see that $P_{x',s_{j'}w'}=P_{s_{j'}x',s_{j'}w'}=P_{x',z}=1$, which together with \cite[(2.2.c)]{KL79} gives $P_{x',w'}=P_{s_{j'}x',s_{j'}w'}+q P_{x',s_{j'}w'}- q P_{x',z}=1$.
\end{proof}

\begin{lem0}
Let $j_0,j_1\in\Delta$ with $j_1\geq j_0$.
\begin{enumerate}[label=(\roman*)]
\item \label{it: double prec 1} If $x\prec w_{j_1,j_0}$ and $x\prec w_{j_1+1,j_0+1}$, then we have $x=w_{j_1,j_0}s_{j_0}=s_{j_1+1}w_{j_1+1,j_0+1}$.
\item \label{it: double prec 2} If $w_{j_1,j_0}\prec x$ and $w_{j_1+1,j_0+1}\prec x$, then we have $x=w_{j_1+1,j_0}$ when $j_1>j_0$, and $x\in\{w_{j_0+1,j_0},w_{j_0,j_0+1}\}$ when $j_1=j_0$.
\end{enumerate}
\end{lem0}
\begin{proof}
Note \ that \ $w_{j_1,j_0}$ \ and \ $w_{j_1+1,j_0+1}$ \ are \ both \ partial-Coxeter \ and \ satisfy \ $\ell(w_{j_1,j_0})=\ell(w_{j_1+1,j_0+1})\geq 1$.

We prove \ref{it: double prec 1}. By the proof of \ref{it: coxeter mult 1} of Lemma~\ref{lem: KL coxeter} we have $P_{x,w_{j_1,j_0}}=P_{x,w_{j_1+1,j_0+1}}=1$, which together with $x\prec w_{j_1,j_0}$ and $x\prec w_{j_1+1,j_0+1}$ forces $\ell(x)=\ell(w_{j_1,j_0})-1=\ell(w_{j_1+1,j_0+1})-1$. As $\mathrm{Supp}(x)\subseteq \mathrm{Supp}(w_{j_1,j_0})\cap \mathrm{Supp}(w_{j_1+1,j_0+1})$, we deduce from \cite[Thm.~2.2.2]{Brenti} that $x=w_{j_1,j_0}s_{j_0}=s_{j_1+1}w_{j_1+1,j_0+1}$.

We prove \ref{it: double prec 2}. If $\ell(x)>\ell(w_{j_1,j_0})+1=\ell(w_{j_1+1,j_0+1})+1$, then by Lemma~\ref{lem: dominance control} we have $D_R(x)\subseteq \{j_0\}\cap \{j_0+1\}=\emptyset$, a contradiction. So we have $\ell(x)=\ell(w_{j_1,j_0})+1=\ell(w_{j_1+1,j_0+1})+1$ and $\mathrm{Supp}(x)\supseteq \mathrm{Supp}(w_{j_1,j_0})\cup \mathrm{Supp}(w_{j_1+1,j_0+1})=[j_1+1,j_0]$. If $j_1>j_0$, then $x>w_{j_1,j_0}$, $\ell(x)=\ell(w_{j_1,j_0})+1$ and $\mathrm{Supp}(x)\supseteq [j_1+1,j_0]$ together with \cite[Thm.~2.2.2]{Brenti} imply $x\in\{w_{j_1+1,j_0},w_{j_1,j_0}s_{j_1+1}\}$. Together with $x>w_{j_1+1,j_0+1}$ (and \cite[Thm.~2.2.2]{Brenti}) this implies $x=w_{j_1+1,j_0}$. If $j_1=j_0$, a similar argument shows that $x\in\{w_{j_0+1,j_0},w_{j_0,j_0+1}\}$.
\end{proof}

Recall that the partially-ordered set $\mathbf{J}^{\infty}$ is defined at the beginning of \S\ref{subsec: sm example}, and that, for each $(j_1,j_2)\in \mathbf{J}^{\infty}$, $x_{j_1,j_2}$ is the element of maximal length in the set (\ref{equ: special Weyl element}).

\begin{lem0}\label{lem: rectangle element}
Let $(j_1,j_2)\in \mathbf{J}^{\infty}$. The element $x_{j_1,j_2}$ is saturated. More precisely, for each $(a_1,b_1)\in\Z\times\Z$ satisfying $a_1-b_1=j_1$, we have $x_{j_1,j_2}=x_{\Sigma}$ with
\begin{equation}\label{equ: rectangle element}
\Sigma=\{(a,b)\mid a_1-j_1+1\leq a\leq a_1, b_1-j_2+j_1+1\leq b\leq b_1\}\subseteq \Z\times\Z.
\end{equation}
\end{lem0}
\begin{proof}
As $D_L(x_{j_1,j_2})=\{j_1\}$ by definition, we have $D_R(x_{j_1,j_2}^{-1})=\{j_1\}$, so $x_{j_1,j_2}^{-1}$ is saturated by Lemma~\ref{lem: maximal dominance}, and thus $x_{j_1,j_2}$ is saturated. We write $x_{j_1,j_2}=x_{\Sigma}$ and let $(a_1,b_1)$ be the unique maximal element in $\Sigma$ (see \ref{it: expansion 2} of Lemma~\ref{lem: saturated expansion}). It remains to check the equality (\ref{equ: rectangle element}). We observe that $\Sigma$ as in (\ref{equ: rectangle element}) is the maximal possible saturated subset of $\Z\times\Z$ which admits $(a_1,b_1)$ as its unique maximal element and satisfies $a-b\leq j_2-1$ for each $(a,b)\in\Sigma$. In particular, if $(a,b)\leq (a_1,b_1)$ satisfies $1\leq a'-b'\leq j_2-1$ for each $(a,b)\leq (a',b')\leq (a_1,b_1)$, then we have $(a,b)\geq (a_1-j_1+1,b_1-j_2+j_1+1)$.
\end{proof}

\begin{rem0}\label{rem: general rectangle}
Let $j'\leq j\leq j''\in\Delta$. Then by a straightforward generalization of Lemma~\ref{lem: rectangle element}, the element of maximal length $x$ satisfying $D_L(x)=\{j\}$ and $\mathrm{Supp}(x)\subseteq [j',j'']$ is saturated, and we have $x=x_{\Sigma}$ with $\Sigma=\{(a',b')\mid a-j+j'\leq a'\leq a, b-j''+j\leq b'\leq b\}$, for any $(a,b)$ satisfying $a-b=j$.
\end{rem0}

Recall that $w_I$ for $I\subseteq \Delta$ is the element of $W(L_I)$ of maximal length (see \S\ref{sec: intro}).

\begin{lem0}\label{lem: saturated cell}
Let $(j_1,j_2)<(j_1',j_2')$ in $\mathbf{J}^{\infty}$.
\begin{enumerate}[label=(\roman*)]
\item \label{it: saturated cell 1} If $j_1'>j_1$ and $j_2'-j_1'>j_2-j_1$, then $x_{j_1',j_2'}\not\leq xx_{j_1,j_2}$ where $x$ is the element of maximal length in $W^{\widehat{j}_1\cap\widehat{j}_1',\emptyset}(L_{\widehat{j}_1})$ (see above (\ref{caca}) for the notation).
\item \label{it: saturated cell 2} If $j_1=j_1'$, then $w_{[1,j_2'-1]}w_{[1,j_2-1]}^{-1}\in W(L_{\widehat{j}_1})s_{j_1}W(L_{\widehat{j}_1})$ if and only if either $j_1=1$ or $j_2'=j_2+1$.
\end{enumerate}
\end{lem0}
\begin{proof}
We prove \ref{it: saturated cell 1}. Let $(a_1,b_1)\in\Z\times\Z$ with $a_1-b_1=j_1$. We define
\begin{eqnarray*}
\Sigma_1&\defeq &\{(a,b)\!\mid a_1-j_1+1\leq a\leq a_1, b_1-j_2+j_1+1\leq b\leq b_1\}\\
\Sigma_2&\defeq &\{(a,b)\mid a_1-j_1+1\leq a\leq a_1+j_1'-j_1, b_1-j_2'+j_1'+1\leq b\leq b_1\}\\
\Sigma_3&\defeq &\{(a,b)\mid a_1+1\leq a\leq a_1+j_1'-j_1, b_1-n+j_1'+1\leq b\leq b_1\}.
\end{eqnarray*}
By Lemma~\ref{lem: rectangle element} (and Remark~\ref{rem: general rectangle} applied to $x$) we see that $x_{j_1,j_2}=x_{\Sigma_1}$, $x_{j_1',j_2'}=x_{\Sigma_2}$ and $x=x_{\Sigma_3}$. It is easy to check that $\Sigma_1\sqcup\Sigma_3$ is still saturated with a unique maximal element $(a_1+j_1'-j_1,b_1)$. Since $\Sigma_2$ also admits $(a_1+j_1'-j_1,b_1)$ as its unique maximal element, we see that $x_{\Sigma_2}=x_{j_1',j_2'}\leq xx_{j_1,j_2}=x_{\Sigma_3}x_{\Sigma_1}=x_{\Sigma_1\sqcup\Sigma_3}$ if and only if $\Sigma_2\subseteq \Sigma_1\sqcup\Sigma_3$. But $j_2'-j_1'>j_2-j_1$ implies that the minimal element $(a_1-j_1+1,b_1-j_2'+j_1'+1)$ is not inside $\Sigma_1\sqcup\Sigma_3$. So we must have $x_{j_1',j_2'}\not\leq xx_{j_1,j_2}$.

We prove \ref{it: saturated cell 2}. Let $y$ be the element of maximal length such that $D_R(y)=\{j_2\}$ and $\mathrm{Supp}(y)\subseteq [1,j_2'-1]$, which satisfies $w_{[1,j_2'-1]}=yw_{[1,j_2-1]}w_{[j_2+1,j_2'-1]}$ with $\ell(w_{[1,j_2'-1]})=\ell(y)+\ell(w_{[1,j_2-1]})+\ell(w_{[j_2+1,j_2'-1]})$ and in particular $w_{[1,j_2'-1]}w_{[1,j_2-1]}^{-1}\geq y$.
If $j_1=1$, then $W^{\widehat{j}_1,\widehat{j}_1}=\{1,s_{j_1}\}$. As $1=j_1\in\mathrm{Supp}(y)$, we have $y\in W(L_{\widehat{j}_1})s_{j_1}W(L_{\widehat{j}_1})$ which together with $w_{[1,j_2'-1]}w_{[1,j_2-1]}^{-1}\geq y$ (and \cite[Prop.~2.5.1]{BB05}) implies $w_{[1,j_2'-1]}w_{[1,j_2-1]}^{-1}\in W(L_{\widehat{j}_1})s_{j_1}W(L_{\widehat{j}_1})$.
If $j_2'=j_2+1$, then we have $w_{[1,j_2'-1]}w_{[1,j_2-1]}^{-1}=y=w_{1,j_2}\in W(L_{\widehat{j}_1})s_{j_1}W(L_{\widehat{j}_1})$ as $j_1\leq j_2$.
We assume from now $j_1>1$ and $j_2'>j_2+1$. Let $(a,b)\in\Z\times\Z$ with $a-b=j_2$. Applying Remark~\ref{rem: general rectangle} to $y^{-1}$, we see that $y^{-1}$ is saturated and $y^{-1}=x_{\Sigma}$ for $\Sigma=\{(a',b')\mid a-j_2+1\leq a'\leq a, b-j_2'+j_2+1\leq b'\leq b\}$. Since we have $j_2-1>j_2-j_1$ and $j_2'>j_2+1$, we observe that $\Sigma$ contains the saturated subset $\Sigma'=\{(a',b')\mid a-j_2+j_1-1\leq a'\leq a-j_2+j_1, b-1\leq b'\leq b\}$, which implies $y^{-1}\geq x_{\Sigma'}$ and moreover $w_{[1,j_2'-1]}w_{[1,j_2-1]}^{-1}\geq y\geq x_{\Sigma'}^{-1}=x_{\Sigma'}\in W^{\widehat{j}_1,\widehat{j}_1}$ (using \ref{it: expansion 2} of Lemma~\ref{lem: saturated expansion}). By \cite[Prop.~2.5.1]{BB05}, we know that $w_{[1,j_2'-1]}w_{[1,j_2-1]}^{-1}\in W(L_{\widehat{j}_1})wW(L_{\widehat{j}_1})$ for some $w\in W^{\widehat{j}_1,\widehat{j}_1}$ such that $w\geq x_{\Sigma'}>s_{j_1}$, and in particular $w_{[1,j_2'-1]}w_{[1,j_2-1]}^{-1}\notin W(L_{\widehat{j}_1})s_{j_1}W(L_{\widehat{j}_1})$.
\end{proof}

\newpage

\section{Some figures for \texorpdfstring{$\GL_4$}{GL4} and \texorpdfstring{$\GL_5$}{GL5}}\label{sec: appendix2}

For $G={\rm GL}_4(K)$ or $G={\rm GL}_5(K)$ we draw many of the previous finite length $D(G)$-modules.\bigskip

In the drawings, as usual a bullet $\bullet$ means an irreducible constituent (the socle being at the bottom of the page and the cosocle being at the top) and a line between two irreducible constituents, that is between two $\bullet$, means a \emph{non-split} extension as subquotient (whether it is dotted or not). For clarity in the three-dimensional drawings, we draw the non-split extensions which are ``behind'' as dotted lines.\bigskip

We first draw the complex of $D(G)$-modules $\mathbf{D}^\bullet$ of (\ref{equ: explicit dR}) in Figure $1$ (${\rm GL}_4(K)$) and Figure~$7$ (${\rm GL}_5(K)$).\bigskip

Then we draw the $D(G)$-modules $\!\begin{xy}(-3,0)*+{Y_{k}}="a"; (16,0)*+{(V_{[1,4-k],\Delta}^{\rm{alg}})^\vee}="b"; {\ar@{-}"a";"b"}\end{xy}\!$ for $1\leq k \leq 3$ in Figure $2$ ($G={\rm GL}_4(K)$) and the $D(G)$-modules $\!\begin{xy}(-3,0)*+{Y_{k}}="a"; (16,0)*+{(V_{[1,5-k],\Delta}^{\rm{alg}})^\vee}="b"; {\ar@{-}"a";"b"}\end{xy}\!$ for $1\leq k \leq 4$ in Figure $8$ ($G={\rm GL}_5(K)$). Recall that $Y_k$ is defined in (\ref{definitionofxyz}) and $V_{[1,4-k],\Delta}^{\rm{alg}},V_{[1,5-k],\Delta}^{\rm{alg}}$ in (\ref{vIalg}). Moreover the red bullets are the duals of locally algebraic constituents.\bigskip

We draw the $D(G)$-module $\widetilde{D}_3$ of (\ref{equ: dn-1tilde}) for $G=\GL_4(K)$ in Figure $3$ and the $D(G)$-module $\widetilde{D}_4$ of (\ref{equ: dn-1tilde}) for $G=\GL_5(K)$ in Figure $9$. Finally for $G=\GL_4(K)$ we draw the $D(G)$-modules $\widetilde{D}_2$ (Figure $4$), $\widetilde{D}_1$ (Figure $5$) and $\widetilde{D}_0$ (Figure $6$) of (\ref{equ: rough shape of rep}), and for $G=\GL_5(K)$ the $D(G)$-modules $\widetilde{D}_3$ (Figure $10$), $\widetilde{D}_2$ (Figure $11$), $\widetilde{D}_1$ (Figure $12$) and $\widetilde{D}_0$ (Figure $13$) of (\ref{equ: rough shape of rep}). Note that $\widetilde{D}_0$ is just the ``top'' of $\widetilde{D}_1$ (keeping in mind that the drawings are $3$-dimensional). Note also that, to save place, we sometimes write $V_{[1,4-k-1]}^{\rm{alg}}$, $V_{[1,5-k-1]}^{\rm{alg}}$ instead of (respectively) $V_{[1,4-k-1],\Delta}^{\rm{alg}}$, $V_{[1,5-k-1],\Delta}^{\rm{alg}}$.\bigskip

The reader can have (a bit of) fun guessing on the drawings what are the morphisms $d_{\mathbf{D}}^k: {D}_k\rightarrow {D}_{k+1}$, the surjections $\widetilde{D}_k\twoheadrightarrow D_k$ and the morphisms $d_{\tld{\mathbf{D}}}^k: \widetilde{D}_k\rightarrow \widetilde{D}_{k+1}$.

\newpage

\pagestyle{empty}

\newgeometry{top=0mm, bottom=0mm, margin=0mm}

\begin{figure}
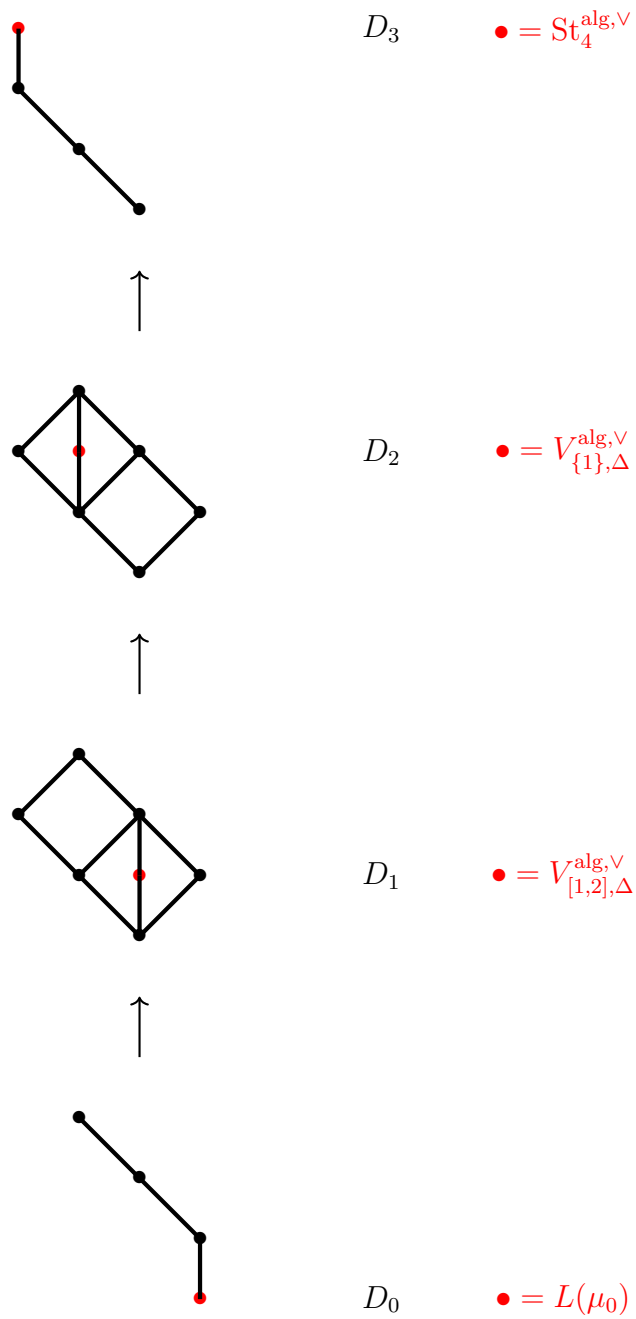

\centering

\vspace{.5cm}
\caption{$\mathbf{D}^\bullet$ for ${\rm GL}_4$}
\end{figure}

\newpage

\begin{figure}
\centering
%
\vspace{1cm}
\caption{$Y_{k}\!-\!V_{[1,4-k],\Delta}^{\rm{alg},\vee}$ for ${\rm GL}_4$ and $1\leq k\leq 3$}
\end{figure}

\newpage

\begin{figure}
\centering
%
\caption{$\widetilde{D}_3$ for ${\rm GL}_4$}
\end{figure}

\newpage

\begin{figure}
\centering
%
\caption{$\widetilde{D}_2$ for ${\rm GL}_4$}
\end{figure}

\newpage

\begin{figure}
\centering
%
\caption{$\widetilde{D}_1$ for ${\rm GL}_4$}
\end{figure}

\newpage

\begin{figure}
\centering
%
\caption{$\mathbf{D}^\bullet$ for ${\rm GL}_5$}
\end{figure}

\newpage

\begin{figure}
\centering
%
\caption{$Y_{k}\!-\!V_{[1,5-k],\Delta}^{\rm{alg},\vee}$ for ${\rm GL}_5$ and $1\leq k\leq 4$}
\end{figure}

\newpage

\begin{figure}
\centering
%
\caption{$\widetilde{D}_0$ for ${\rm GL}_5$}
\end{figure}

\restoregeometry

\newpage

\pagestyle{plain}

\bibliographystyle{alpha}

\end{document}